\documentclass[10pt, leqno]{book}
\usepackage{amsmath,amsfonts,amsthm,amscd,amssymb}
\numberwithin{equation}{section}
\theoremstyle{plain}
\newtheorem{theo}{Theorem}[section]
\newtheorem{prop}[theo]{Proposition}
\newtheorem{coro}[theo]{Corollary} 
\newtheorem{lemm}[theo]{Lemma}
\theoremstyle{definition}
\newtheorem{defi}[theo]{Definition}
\newtheorem{rema}[theo]{Remark}
\newtheorem{theo-defi}[theo]{Theorem-Definition}
\newtheorem{prop-defi}[theo]{Proposition-Definition}
\newtheorem{coro-defi}[theo]{Corollary-Definition}
\newtheorem{rema-defi}[theo]{Remark-Definition}
\newtheorem{exem-defi} [theo]{Example-Definiton}

\newtheorem{exem}[theo]{Example}
\newtheorem{conj}[theo]{Conjecture}
\newtheorem{prob}[theo]{Problem}

\def \al{\alpha}
\def \bet{\beta}
\def \bul{\bullet}
\def \col{\colon}
\def \Del{\Delta}

\def \del{\delta}
\def \eps{\epsilon}
\def \Gam{\Gamma}
\def \gam{\gamma}
\def \inf{\infty}

\def \kap{\kappa}
\def \Lam{\Lambda}
\def \lam{\lambda}

\def \Lo{\Longrightarrow}
\def \lo{\longrightarrow}
\def \lom{\longmapsto}
\def \mab{\mathbb}

\def \nat{\natural} 

\def \Om{\Omega}
\def \om{\omega}
\def \ol{\overline}
\def \os{\overset}
\def \parno{\par\noindent}

\def \Sig{\Sigma}
\def \sig{\sigma}
\def \sq{\square}
\def \sus{\subset}

\def \ul{\underline}
\def \us{\underset}

\def \vp{\varpi}

\def \vil{\varinjlim}
\def \vpl{\varprojlim}

\def \wh{\widehat}
\def \wt{\widetilde}
\bigskip

\newcommand{\getsfrom}
{\ensuremath{\longleftarrow\kern-.
52em\lower-.1ex\hbox%
{$\shortmid\,$}}}

\begin{document}
\title{Limits of weight filtrations and limits of slope filtrations 
on infinitesimal cohomologies in mixed characteristics I
\parno
}
\author{Yukiyoshi Nakkajima
\date{}\thanks{\endgraf}}
\maketitle

\bigskip
\parno
{\bf Abstract.---} 
We construct a fundamental filtered complex 
$(A_{\rm zar},P)$ for the base change of a 
proper simple normal crossing log scheme over a family of log points 
in characteristic $p>0$ with an admissible immersion. 
This complex is a generalization of 
Mokrane's Hyodo-Steenbrink complex modulo torsion. 
As an application of this construction, we construct a theory   
of the limit of the weight filtration and a theory of 
the limit of the slope filtration on the infinitesimal cohomology 
of a proper scheme over the fraction field of 
a complete discrete valuation ring of mixed characteristics 
with perfect residue field by using the log crystalline cohomology 
of the log special fiber of an appropriate model of 
a split proper hypercovering of the scheme. 
By using canonical Hyodo-Kato-Tsuji's isomorphism, 
we also prove that the infinitesimal cohomology 
of a separated scheme of finite type over the field 
has a contravariant canonical action of the crystalline Weil-Deligne group.

\bigskip
\parno



$${\bf Contents}$$
\parno
\S\ref{sec:int}. Introduction
\bigskip
\parno 
Chapter I. Weight filtrations on log crystalline cohomology sheaves
\medskip
\parno 
\S\ref{sec:snclv}. SNCL schemes
\parno 
\S\ref{sec:tdiai}. Review for results in \cite{nh3} 
and bisimplicial (exact) immersions I
\parno 
\S\ref{sec:ldc}. Log crystalline de Rham complexes
\parno 
\S\ref{sec:psc}. 
Zariskian $p$-adic filtered Steenbrink complexes  
\parno 
\S\ref{sec:fcuc}.   
Contravariant functoriality of 
zariskian $p$-adic filtered Steenbrink complexes
\parno
\S\ref{sec:bckf}. Filtered base change theorem I
\parno 
\S\ref{sec:crcks}. 
$p$-adic monodromy operators I and modified $P$-filtered log crystalline complexes 
\parno 
\S\ref{sec:vpmn}. 
Variational $p$-adic monodromy-weight conjecture I and 
variational filtered log $p$-adic hard Lefschetz conjecture I
\bigskip
\parno 
Chapter II. Weight filtrations and slope filtrations on log crystalline cohomologies 
via log de Rham-Witt complexes
\medskip
\parno
\S\ref{sec:ldrwc}. 
Log de Rham-Witt complexes of crystals I
\parno
\S\ref{sec:ldrwii}. 
Log de Rham-Witt complexes of crystals II 
\parno  
\S\ref{sec:flgdw}. Zariskian $p$-adic filtered Steenbrink complexes  
via log de Rham-Witt complexes
\parno 
\S\ref{sec:pmr}. 
$p$-adic monodromy operators 
via log de Rham-Witt complexes  
\parno 
\S\ref{sec:pdwmn}.  
$p$-adic monodromy-weight conjecture 
and filtered log $p$-adic hard Lefschetz conjecture via log de Rham-Witt complexes
\bigskip
\parno
Chapter III. Weight filtrations on log isocrystalline cohomology sheaves 
\medskip
\parno   
\S\ref{sec:bcsncl}.  
Truncated simplicial base changes of SNCL schemes
\parno 
\S\ref{sec:beii}. 
Families of successive split truncated simplicial SNCL schemes and 
bisimplicial (exact) immersions II
\parno 
\S\ref{sec:pgensc}. Iso-zariskian $p$-adic filtered Steenbrink complexes 
\parno
\S\ref{sec:pfnsc}.  
Contravariant functoriality of iso-zariskian $p$-adic filtered Steenbrink complexes 
\parno 
\S\ref{sec:bcakf}. 
Filtered base change theorem II
\parno 
\S\ref{sec:mndomp}. $p$-adic monodromy operators II and 
modified $P$-filtered log isocrystalline complexes 
\bigskip
\parno 
Chapter IV. Weight filtrations and slope filtrations on log isocrystalline cohomologies 
\medskip
\parno  
\S\ref{sec:flgsiidw}.  
Iso-zariskian $p$-adic filtered Steenbrink complexes  
via log de Rham-Witt complexes
\bigskip
\parno
Chapter V. Results for weight filtrations on log isocrystalline cohomology sheaves 
\medskip 
\parno 
\S\ref{sec:cfi}. Log convergent $F^{\infty}$-isospans
\parno 
\S\ref{sec:lcw}.   
Log convergences of the weight filtrations on the log isocrystalline cohomology sheaves 
\parno 
\S\ref{sec:infdi}. 
Infinitesimal deformation invariance of iso-zarisikian 
$p$-adic filtered Steenbrink complexes
\parno 
\S\ref{sec:filbo}. 
The $E_2$-degeneration of 
the $p$-adic weight spectral sequence and strict compatibility 
\parno 
\S\ref{sec:vpmwc}. 
Variational $p$-adic monodromy-weight conjecture II and 
variational log $p$-adic hard Lefschetz conjecture II 
\bigskip
\parno
Chapter VI. Limits of weight filtrations and limits of slope filtrations 
on infinitesimal cohomologies in mixed characteristics
\medskip 
\parno 
\S\ref{sec:ph}. 
Good and split proper hypercoverings 
\parno 
\S\ref{lic}. Log infinitesimal cohomologies 
\parno 
\S\ref{sec:pmm}. 
Monodromy operators on the infinitesimal 
cohomologies of proper schemes in mixed characteristics
\parno 
\S\ref{faiw}. 
Actions of crystalline Weil-Deligne groups on the infinitesimal 
cohomologies of proper schemes in mixed characteristics 
\parno 
\S\ref{sec:pff}. Limits of the weight filtrations 
on the infinitesimal cohomologies of 
proper schemes in mixed characteristics
\parno 
\S\ref{sec:lcwgac}. 
Limits of the slope filtrations on the infinitesimal cohomologies of 
proper schemes in mixed characteristics
\parno 
\S\ref{sec:padicar}. 
$p$-adic arithmetic weight filtrations

\bigskip
\parno
Chapter VII.  Generalized Ogus' conjecture and Fontaine's conjecture
\medskip 
\parno 
\S\ref{sec:ofc}. 
Generalized Ogus's conjecture and Fontaine's conjecture on the action of crystalline 
Weil-Deligne group for varieties

\smallskip
\parno
References

\section{Introduction}\label{sec:int}  
Based on Grothendieck's idea of motives (\cite{grrs}), 
in \cite{dh1} Deligne has made a translation
between $\infty$-adic objects and $l$-adic objects, 
especially the weight filtration 
on the singular cohomology of the analytification of 
a scheme of finite type over ${\mab C}$ 
and the weight filtration on 
the $l$-adic \'{e}tale cohomology of 
a scheme of finite type over a finite field.  
In \cite{dh2} and \cite{dh3} 
Deligne has constructed the theory 
of Hodge-Deligne; especially he has proved 
the existence of the weight filtration 
on the singular cohomology above. 
\par 
In \cite{brc}, \cite{bfi} and \cite{bpre} 
Berthelot has defined the rigid cohomology of a separated scheme 
of finite type over a perfect field of characteristic $p>0$. 
This is an analogue in characteristic $p$ of the singular cohomology above. 
Using theory of filtered derived category of Berthelot (\cite{blec}) and 
our results in \cite{nh2}, in \cite{nh3}   
we have proved the existence of 
the weight filtration on 
the rigid cohomology of a separated scheme 
of finite type over a perfect field of characteristic $p>0$.  
This is an analogue in characteristic $p$ of 
Deligne's result above.  
\par 
On the other hand, in \cite{db} du Bois has constructed the limit of 
the weight filtration on the singular cohomology 
of the nearby cycle complex of the analytification of 
a proper scheme over a smooth curve over ${\mab C}$ 
at the localization of a closed point of the smooth curve by 
using a proper semistable hypercovering of the scheme 
which exists by the semistable reduction theorem in characteristic $0$ 
(\cite{kkms}).  
As a result he has defined the limit of the weight filtration on 
the singular cohomology with coefficient ${\mab Q}$ of 
the analytification of the generic fiber of the proper scheme 
at the localization of the point.  
\par 
In \cite{grcr} Grothendieck has defined the infinitesimal cohomology 
of a noetherian scheme over a noetherian scheme in characteristic $0$ 
to generalize the de Rham cohomology of a smooth scheme 
over a noetherian scheme in characteristic $0$. 
By comparison theorems in \cite{grcr} and \cite{hadr} (see also \cite{cm}),  
we know that the infinitesimal cohomology 
calculates the singular cohomology with coefficient ${\mab C}$. 
\par 
In this book we are interested 
in an analogue of du Bois' result in mixed characteristics. 
Our main technical ingredient is 
the log crystalline cohomology of 
a proper log smooth scheme in characteristic $p$, 
which has been defined by Kato in \cite{klog1} 
by using theory of log structures of Fontaine-Illusie-Kato.  
(See also \cite{fao}.) 
His log crystalline cohomology is a generalization of 
the crystalline cohomology of a proper smooth scheme 
in characteristic $p$ defined by Grothendieck-Berthelot 
(\cite{grcr}, \cite{bb}, \cite{bob}).
\par 
Let ${\cal V}$ be a complete discrete valuation ring 
of mixed characteristics $(0,p)$ with perfect residue field $\kap$. 
Set $K:={\rm Frac}({\cal V})$. 
Let ${\mathfrak X}$ be a proper scheme over $K$. 
(It is better to assume only that ${\mathfrak X}$ 
is separated and finite type over $K$; 
however we make the assumption of the properness of ${\mathfrak X}$ in this book 
except Chapter VII in order to make the length of this book reasonable.) 
One of our main results in this book is to show the existence of 
the ``limit of the weight filtration'' $P$ on the infinitesimal cohomology 
$H^q_{\rm inf}({\mathfrak X}/K)$ $(q\in {\mab N})$ 
by using the log isocrystalline cohomology 
of the log reduction of a certain log smooth model of 
a split proper hypercovering of the base change of 
${\mathfrak X}$ with respect to a finite extension of $K$.  
(We call $P$ the {\it limit weight filtration} of $H^q_{\rm inf}({\mathfrak X}/K)$.)
This is an analogue in mixed characteristics of du Bois' result, 
though our method is not the imitation of du Bois' one 
because the semistable reduction theorem 
has not yet been proved in the mixed characteristics case. 
However each log scheme of the members of 
the split proper hypercovering is the disjoint union of 
the base changes of proper strict semistable families. 
(The existence of the split proper hypercovering is assured 
by de Jong's semistable reduction theorem (\cite{dj}).) 
In this way, proper strict semistable families over ${\cal V}$, 
or more generally, proper strict semistable schemes 
over the log point $s$ of $\kap$ 
(by taking the log special fibers of them), 
are fundamental log schemes in this book. 
(We call the log schemes above 
proper SNCL(=simple normal crossing log) schemes over $s$.)
\par 
Set $K_0:={\cal W}(\kap)$. 
In fact, we define a canonical contravariantly functorial $K_0$-structure 
$H^q_{\rm inf}({\mathfrak X}/K/K_0)$ of  
$H^q_{\rm inf}({\mathfrak X}/K)$ and a $K_0$-form 
$P(0)$ of $P$ on $H^q_{\rm inf}({\mathfrak X}/K/K_0)$. 
Here we note that  $H^q_{\rm inf}({\mathfrak X}/K/K_0)$ is 
only a finite dimensional $K_0$-vector space 
and we have not proved that it is a cohomology of an abelian sheaf in a topos. 
The scheme ${\mathfrak X}$ can be considered as a proper scheme over $K_0$ 
by using the finite morphism ${\rm Spec}(K)\lo {\rm Spec}(K_0)$. 
By the base change theorem for infinitesimal cohomologies, 
the canonical morphism 
\begin{align*} 
H^q_{\rm inf}({\mathfrak X}/K_0)\otimes_{K_0}K
\os{\sim}{\lo} 
H^q_{\rm inf}({\mathfrak X}\otimes_{K_0}K/K)
\end{align*} 
is an isomorphism. However $H^q_{\rm inf}({\mathfrak X}\otimes_{K_0}K/K)\not=
H^q_{\rm inf}({\mathfrak X}/K)$ in general if ${\mathfrak X}\not= \emptyset$.
Hence $H^q_{\rm inf}({\mathfrak X}/K/K_0)\not=H^q_{\rm inf}({\mathfrak X}/K_0)$ in general. 
We also define an action of the monodromy operator on 
$H^q_{\rm inf}({\mathfrak X}/K/K_0)$ 
and the action of the Frobenius endomorphism on it without a choice of 
a uniformizer of $K$. 
This is a generalization of a stronger version of U.~Jannsen's conjecture 
to the case of any proper scheme over ${\rm Frac}({\cal V})$ 
(original his conjecture is a conjecture 
for a proper semistable family over ${\cal V}$ 
and this has been proved by Hyodo and Kato in \cite{hk} 
via famous Hyodo-Kato's isomorphism; see also the introduction of \cite{tsst};  
however see also \cite{ndw} because the proof of Hyodo-Kato isomorphism is not perfect). 
To define the monodromy operator on $H^q_{\rm inf}({\mathfrak X}/K/K_0)$ and 
the action of the Frobenius endomorphism 
on $H^q_{\rm inf}({\mathfrak X}/K/K_0)$, we use 
Tusji's truncated simplicial version of the Hyodo-Kato isomorphism 
(\cite{tsgep}). 
Precisely speaking, 
we define and use the {\it canonical} Hyodo-Kato-Tsuji isomorphism which 
does {\it not} depend on the choice of a uniformizer of ${\cal V}$ by modifying 
the Hyodo-Kato-Tsuji isomorphism. 
(It seems to me that no one has imagined  
that there exists such a canonical Hyodo-Kato-Tsuji isomorphism.) 
Let ${\rm WD}_{\rm crys}(\ol{K})$ be the crystalline Weil-Deligne group 
defined by Ogus in \cite{ollc}; it fits into the following exact sequence 
\begin{align*} 
0\lo \ol{K}(1)\lo {\rm WD}_{\rm crys}(\ol{K})\lo W_{\rm crys}(\ol{K})\lo 0, 
\end{align*} 
where $W_{\rm crys}(\ol{K})$ is 
the crystalline Weil group defined in \cite{boi}: 
$W_{\rm crys}(\ol{K}):=
\{g\in {\rm Aut}(\ol{K})~\vert~$the restriction of $g$ on 
Frac$({\cal W}(\ol{\kap}))$ is the integral power of 
the Frobenius automorphism of Frac$({\cal W}(\ol{\kap}))\}$. 
By using the canonical Hyodo-Kato-Tsuji isomorphism, 
we define the well-defined semi-linear action 
$\rho^q_{{\mathfrak X}_{\ol{K}}/\ol{K}}$ of ${\rm WD}_{\rm crys}(\ol{K})$
on $H^q_{\rm inf}({\mathfrak X}_{\ol{K}}/\ol{K})$. 
That is, we prove that there exists the following canonical semi-linear 
group homomorphism 
\begin{align*} 
\rho^q_{{\mathfrak X}_{\ol{K}}/\ol{K}} \col {\rm WD}_{\rm crys}(\ol{K})\lo {\rm GL}_{{\mab Q}_p}(H^q_{\rm inf}({\mathfrak X}_{\ol{K}}/\ol{K})). 
\end{align*}
In [loc.~cit.] Ogus has conjectured that $\rho^q_{{\mathfrak X}_{\ol{K}}/\ol{K}}$ is contravariantly functorial
with respect to a morphism of the generic fibers of 
proper semistable families over ${\cal V}$.  
(Strictly speaking, he has not defined our $\rho^q_{{\mathfrak X}_{\ol{K}}/\ol{K}}$ in [loc.~cit.] 
even in the proper semistable case.)
In this book we solve this conjecture affirmatively 
in a more generalized form.  
That is, for a separated scheme ${\mathfrak U}$ of finite type over $\ol{K}$, 
we define the following canonical semi-linear 
group homomorphism 
\begin{align*} 
\rho^q_{{\mathfrak U}/\ol{K}} \col {\rm WD}_{\rm crys}(\ol{K})\lo 
{\rm GL}_{{\mab Q}_p}(H^q_{\rm inf}({\mathfrak U}/\ol{K})). 
\end{align*}
and we prove that, 
for a morphism of ${\mathfrak f}\col {\mathfrak U}\lo {\mathfrak V}$ of 
separated schemes of finite type over $\ol{K}$ and for an element 
$\gam$ of ${\rm WD}_{\rm crys}(\ol{K})$, 
the following diagram is commutative: 
\begin{equation*} 
\begin{CD}
H^q_{\rm inf}({\mathfrak U}/\ol{K})@>{\rho^q_{{\mathfrak U}/\ol{K}}(\gam)}>> H^q_{\rm inf}({\mathfrak U}/\ol{K})\\
@A{{\mathfrak f}^*}AA @AA{\mathfrak f}^*A \\
H^q_{\rm inf}({\mathfrak V}/\ol{K})@>{\rho^q_{{\mathfrak V}/\ol{K}}(\gam)}>> 
H^q_{\rm inf}({\mathfrak V}/\ol{K}). 
\end{CD}
\end{equation*} 
Because it is difficult to solve original Ogus' conjecture itself, 
we use Grothendieck's general principle: 
in order to solve a problem in mathematics, 
one has to generalize the problem as possible in order that 
the arguments in the proof of the problem work 
and then one solve it in an evident but often long way. 
In this case we use simplicial method and de Jong's theorem mentioned before. 
Furthermore, for a separated scheme ${\mathfrak U}$ of finite type over $K$, 
by using the natural action $\rho_{\rm inf}$ of ${\rm Gal}(\ol{K}/K)$ 
on $H^q_{\rm inf}({\mathfrak U}/K)\otimes_K\ol{K}$,  
we can define a $\ol{K}$-linear action $\rho_{\rm inf}\circ \rho^{-1}_{{\mathfrak U}_{\ol{K}/\ol{K}}}$ 
of the (usual) Weil-Deligne group ${\rm WD}(\ol{K}/K)$ 
on $H^q_{\rm inf}({\mathfrak U}_{\ol{K}}/\ol{K})$ and 
we can conjecture the compatibility of the standard action of 
${\rm WD}(\ol{K}/K)$ on the $l$-adic \'{e}tale cohomology. 
Here the group ${\rm WD}(\ol{K}/K)$ is, by definition, 
the inverse image of 
$W_{\rm crys}(\ol{K})\cap {\rm Gal}(\ol{K}/K)$ 
in ${\rm W}_{\rm crys}(\ol{K})$ by the morphism 
${\rm WD}_{\rm crys}(\ol{K})\lo {\rm W}_{\rm crys}(\ol{K})$.  
This conjecture is a straightforward generalization of 
the compatibility of the action of ${\rm WD}(\ol{K}/K)$  on 
the log isocrystalline cohomology 
with the action of ${\rm WD}(\ol{K}/K)$ 
obtained by the $l$-adic cohomology 
$H^q_{\rm et}({\cal X}_{\ol{K}},{\mab Q}_l)$ of 
a proper semistable family ${\cal X}$ over ${\cal V}$ in \cite{fos} 
and a special case of 
a more general Fontaine's conjecture about motives in [loc.~cit.]. 
(See also a recent work of Chiarellotto--Lazda \cite{cla} 
for a proper semistable family over a complete discrete valuation ring of 
equal characteristic.) 
In fact, because we construct
the limit weight filtration $P$ on $H^q_{\rm inf}({\mathfrak X}/K)$ 
(and hence on $H^q_{\rm inf}({\mathfrak X}_{\ol{K}}/\ol{K})$)
and because we can construct the limit weight filtration $P$ 
on $H^q_{\rm et}({\mathfrak X}_{\ol{K}},{\mab Q}_l)$ and 
because the actions of ${\rm WD}(\ol{K}/K)$ on the cohomologies 
preserves $P$'s, we can conjecture that 
$P_kH^q_{\rm inf}({\mathfrak X}_{\ol{K}}/\ol{K})$ and 
$P_kH^q_{\rm et}({\mathfrak X}_{\ol{K}},{\mab Q}_l)$ are compatible. 
\par 
We also calculate the slope filtration on 
$H^q_{\rm inf}({\mathfrak X}/K/K_0)$ 
by enlarging $K_0$ by a certain unramified finite extension of $K_0$  
and using the cohomology of the log de Rham-Witt complex 
(\cite{hk}) of the log reduction already stated. 
\par 
Let $U$ be a separated scheme of finite type over $\kap$.
Let $j \col U \os{\sus}{\lo} \ol{U}$ be 
an open immersion into a proper scheme over $\kap$.
Let $(U_{\bul},X_{\bul})$ be a split proper hypercovering
of $(U,\ol{U})$ in the sense of 
\cite{tzp}, that is, $(U_{\bul},X_{\bul})$ is split, 
$U_{\bul}$ is a proper hypercovering of $U$, 
$X_{\bul}$ is a proper simplicial scheme over $\ol{U}$ 
and $U_{\bul}=X_{\bul}\times_{\ol{U}}U$. 
Assume that 
$X_{\bul}$ is a proper smooth simplicial scheme over $\kap$ 
and that $U_{\bul}$ is the complement of 
a simplicial SNCD $D_{\bul}$ 
on  $X_{\bul}$ over $\kap$. 
In \cite{nh3} we have proved 
\begin{equation*}
R\Gam_{\rm rig}(U/K)=
R\Gam((X_{\bul},D_{\bul})/{\cal W}(\kap))_K.  
\end{equation*}  
In particular, we see that $R\Gam((X_{\bul},D_{\bul})/{\cal W}(\kap))_{K_0}$ 
is independent of the choice of $(X_{\bul},D_{\bul})$. 
This proves that de Jong's conjecture in \cite{dj} is true: 
$H^q_{\rm crys}((X_{\bul},D_{\bul})/{\cal W}(\kap))_K$ $(q\in {\mab N})$ 
is independent 
of the choice of the proper hypercovering $U_{\bul}=X_{\bul}\setminus D_{\bul}$ 
of $U$. 
By using $H^q_{\rm inf}({\mathfrak X}/K/K_0)$, we can prove the mixed characteristic 
analogue of this conjecture.  
This will be explained later in this introduction 
(see (\ref{eqn:esc000ffsp}) below for the precise statement). 
\par 
As mentioned above, proper strict semistable families
play important roles in this paper as fundamental ingredients. 
Because there are a lot of important works for 
proper strict semistable families in all characteristics, 
we recall several important results for proper strict semistable families 
over the unit disk over ${\mab C}$ 
and several results for proper strict semistable families in 
mixed characteristics and equal characteristic $p>0$ 
by the use of $p$-adic methods 
before we state our results for proper strict semistable families 
over a complete discrete valuation ring, 
more generally a proper SNCL(=simple normal crossing log) scheme 
over a family of log points.  
(See \cite{rz}, \cite{nd}, \cite{stwsl} and \cite{nlpi} for the $l$-adic case.) 
\par 
In \cite{st1} Steenbrink has 
constructed a bifiltered complex,  
what is called, the bifiltered Steenbrink complex 
of a proper strict semistable analytic family 
over the unit disk over ${\mab C}$.  
When the irreducible components of the special fiber 
are K\"{a}hler or the analytifications of 
proper smooth schemes over ${\mab C}$,   
the bifiltered complex gives the limit of 
the weight filtration on the singular cohomology of 
the generic fiber of the semistable family. 
Though we are interested in the case of characteristic $p>0$ 
or mixed characteristics, 
we are deeply influenced by his ingenious method. 
In \cite{kwn} Kawamata-Namikawa and in \cite{fn} 
Fujisawa-Nakayama have generalized Steenbrink's work as follows. 
\par 
By the theory of log schemes in \cite{klog1}, 
we have the notion of an SNCL(=simple normal crossing log) scheme. 
(See \S\ref{sec:snclv} below for the precise definition of 
an SNCL scheme.) 
In a standard way, a strict semistable family gives an SNCL scheme: 
the special fiber with canonical log structure. 
In [loc.~cit.] Fujisawa and Nakayama 
(see also Kawamata-Namikawa's article for a slightly different formulation)
have constructed a bifiltered complex which gives the weight filtration on 
the singular cohomology of the base change 
of the real blow up (\cite{kn}) of the analytification of 
a proper SNCL scheme over the log point  
$({\rm Spec}({\mab C}),({\mab N}\oplus {\mab C}^*\lo {\mab C}))$ 
by the morphism 
${\mab R}\owns x \lom e^{2\pi \sqrt{-1}x}\in {\mab S}^1$, 
where the morphism ${\mab N}\oplus {\mab C}^*\lo {\mab C}$ is defined by 
$(a,u)\lom 0$ $((a,u)\in {\mab N}\oplus {\mab C}^*, a\not=0)$ and $(0,u)\lom u$.  
(Though Steenbrink has claimed 
that he has constructed a similar bifiltered complex in \cite{st2}, 
he has not proved that his bifiltered complex 
is independent of the choice of local charts in [loc.~cit.]; 
we have already pointed out this in \cite{nlpi}.) 
\par  
Let $\kap$ be a perfect field of characteristic $p>0$.  
Let $s=({\rm Spec}(\kap),({\mab N}\oplus \kap^*\lo \kap))$ 
be the log point of $\kap$, 
where the morphism ${\mab N}\oplus \kap^*\lo \kap$ 
is defined by $(a,u)\lom 0$ $((a,u)\in {\mab N}\oplus {\kap}^*, a\not=0)$ 
and $(0,u)\lom u$. 
Let ${\cal W}(s)$ be a log scheme 
whose underlying scheme is ${\rm Spf}({\cal W})$ 
and whose log structure is 
${\mab N}\oplus {\cal W}^* \lo {\cal W}$ 
defined by $(a,u)\lom 0$ $((a,u)\in {\mab N}\oplus {\cal W}^*, a\not=0)$ 
and $(0,u)\lom u$. 
Let $X$ be an SNCL scheme over $s$. 
As explained in the last paragraph, 
the special fiber with the canonical log structure 
of a strict semistable family 
over a discrete valuation ring with residue field $\kap$  
is an example of an SNCL scheme 
over $s$
(we call this special fiber with the canonical log structure 
the {\it log special fiber of the strict semistable family}).  
Before the works of Kawamata-Namikawa and 
Fujisawa-Nakayama, 
in \cite{msemi} Mokrane has constructed a filtered complex 
$({\cal W}A^{\bul}_X,\{P_k{\cal W}A^{\bul}_X\}_{k\in {\mab Z}})$ on 
the Zariski site on $X$  
by using the log de Rham-Witt complex 
${\cal W}\wt{\Om}^{\bul}_X$ of $X/\kap$ 
essentially defined by Hyodo (\cite{hdw}). 
Let ${\cal W}\Om^{\bul}_X$ be the log de Rham-Witt complex of $X/s$
defined by Hyodo ([loc.~cit.]) and Hyodo-Kato (\cite{hk}). 
(They have denoted ${\cal W}\wt{\Om}^{\bul}_X$ and 
${\cal W}\Om^{\bul}_X$ by ${\cal W}\wt{\om}^{\bul}_X$
and ${\cal W}\om^{\bul}_X$, respectively.) 
In \cite{msemi} and \cite{ndw} Mokrane and I 
have proved that there exists a quasi-isomorphism 
$\theta \wedge \col {\cal W}\Om^{\bul}_X \os{\sim}{\lo} {\cal W}A^{\bul}_X$. 
Furthermore, in \cite{ndw} 
I have proved that $\theta \wedge$ is contravariantly functorial 
with respect to a morphism $f\col X\lo Y$ to 
an SNCL scheme over $s$. 
Let $H^q_{{\rm crys}}(X/{\cal W}(s))$ $(q\in {\mab Z})$ 
be Kato's $q$-th log crystalline cohomology of $X/{\cal W}(s)$. 
(In this book we omit to write ``log'' in the notation 
$H^q_{{\rm log}{\textrm -}{\rm crys}}(X/{\cal W}(s))$ in references.) 
In \cite{hk} Hyodo and Kato and in \cite{ndw} I  
have proved that the complex ${\cal W}\Om^{\bul}_X$ 
calculates $H^q_{{\rm crys}}(X/{\cal W}(s))$: 
$H^q_{{\rm crys}}(X/{\cal W}(s))=H^q(X,{\cal W}\Om^{\bul}_X)$. 
This is the log version of Bloch-Illusie's result (\cite{bl}, \cite{idw}).  
Let $\os{\circ}{X}{}^{(k)}$ $(k\in {\mab Z}_{\geq 0})$ 
be the disjoint union of the $(k+1)$-fold 
intersections of the different irreducible components of 
$\os{\circ}{X}$ as in \cite{stwsl}. 
By calculating ${\rm gr}^P_k{\cal W}A^{\bul}_X$, Mokrane has constructed 
the following $p$-adic weight spectral sequence of 
$H^q_{{\rm crys}}(X/{\cal W}(s))$:  
\begin{align*} 
E_1^{-k,q+k} &=
\bigoplus_{j\geq \max \{-k,0\}} 
H^{q-2j-k}_{\rm crys}(\os{\circ}{X}{}^{(2j+k)}/{\cal W})(-j-k) 
=H^{q+k}(X,{\rm gr}_k^P{\cal W}A^{\bul}_X)
\tag{0.0.1.1}\label{eqn:eswsp} \\ 
& \Lo 
H^q(X,{\cal W}A^{\bul}_X)=H^q(X,{\cal W}\Om^{\bul}_X)=
H^q_{{\rm crys}}(X/{\cal W}(s)) 
\end{align*}  
when the underlying scheme $\os{\circ}{X}$ of $X$ is proper over $\kap$.  
This spectral sequence induces the weight filtration on 
$H^q_{{\rm crys}}(X/{\cal W}(s))$. 
(However there are incomplete and mistaken points 
in \cite{hdw}, \cite{hk} and \cite{msemi}: 
we have completed and corrected all 
in \cite{ndw} except (\ref{theo:ccrw}) (see (\ref{rema:altn}) (2), (3), (4)), 
(\ref{prop:plz}) (see (\ref{rema:qok}) (1), (2), (3)) and (\ref{prop:prm}) 
(see (\ref{rema:lcrdc})) below. See also (\ref{rema:sot}) below.)   
The spectral sequence and the filtered complex 
$({\cal W}A^{\bul}_X,P)$ have many applications and they have been used in, e.~g., 
\cite{msemi}, \cite{stmfph}, \cite{cric}, \cite{ocltm}, \cite{stk}, \cite{kiha}, 
\cite{ndw}, \cite{itup}, \cite{nlpi}, \cite{miepl}, \cite{rx}, \cite{magk3}, 
\cite{mwf}, \cite{zflv}, though the indispensable work \cite{ndw} for 
the constructions of (\ref{eqn:eswsp}) and $({\cal W}A^{\bul}_X,P)$ 
has been ignored in some of them. 
Though we can construct the weight filtration on 
$H^q_{{\rm crys}}(X/{\cal W}(s))$ by using $({\cal W}A^{\bul}_X,P)$, 
it seems impossible to answer 
the following fundamental problem only by using $({\cal W}A^{\bul}_X,P)$.  
Let $f\col X\lo Y$ be a morphism of 
proper SNCL schemes over a morphism of log points 
$s\lo t$ of perfect fields of characteristic $p>0$. 
Then, is the pull-back 
$$f^*\col H^q_{{\rm crys}}(Y/{\cal W}(t))\otimes_{\mab Z}{\mab Q}
\lo H^q_{{\rm crys}}(X/{\cal W}(s))\otimes_{\mab Z}{\mab Q}$$ 
strictly compatible 
with the weight filtration? 
If $s$ and $t$ are log points of finite fields,  
then the answer is ``Yes'' by the purity of the Frobenius weight of
the crystalline cohomology of a proper smooth scheme 
over a finite field (\cite{dw2}, \cite{kme}, \cite{clpu}, \cite{ndw}; 
see also \cite{ny} because the proof of the weak Lefschetz theorem in 
\cite{bwl} has a gap) 
and the spectral sequence (\ref{eqn:eswsp}). 
However there is a serious difficulty to answer this problem 
when $s$ and $t$ are not necessarily log points of finite fields  
because we do not have theory of (log) de Rham-Witt complex 
over a general (log) $p$-adic base formal scheme 
(however see \cite{maoro} for a special (log) $p$-adic base formal scheme, 
which is a log overconvergent version of \cite{lazi}). 
In particular, we do not know whether 
a specialization argument for cohomologies of (log) de Rham-Witt complexes 
over a general perfect field to them over a finite field is possible.  
In this book we give the answer ``Yes'' for this problem by 
constructing another more general filtered complex 
``$(A_{\rm zar},P)$'' using the log crystalline method 
in which the specialization argument is possible 
and by comparing our filtered complex with Mokrane's one 
(see the explanations for (9) and (16) below). 
\par 
In \cite{gkwf} Gro{\ss}e-Kl\"onne has constructed 
the $p$-adic (weight) spectral sequence of 
a (proper) SNCL scheme over $s$ 
by using the sheaves of differential forms on 
local weak formal log schemes. 
(See also \cite{maoro} and \cite{greglan}.)
(However he has not proved that his spectral sequence 
is independent of the choice of an open covering of 
the SNCL scheme nor that it is compatible with 
the Frobenius action.)    
In \cite{tsd} Tsuji has constructed 
the $p$-adic weight spectral sequence 
of the log crystalline cohomology of 
a proper SNCL scheme over $s$  
by using the theory of arithmetic ${\cal D}$-modules of Berthelot 
(\cite{bdmi}, \cite{bdmii}, \cite{bdmitr}). 
Tsuji has also proved important properties of 
his weight spectral sequence, e.~g., 
the compatibility of his weight spectral sequence with  
the contravariant and covariant functorialities. 
His results are $p$-adic analogues of T.~Saito's results 
for $l$-adic \'{e}tale cohomologies (\cite{stwsl}) 
and are useful for the compatibility of local and global 
Langlands correspondences (\cite{tay}, \cite{ytgll}, \cite{lcara}, \cite{pcara}). 
On the other hand, Caraiani has constructed 
the $p$-adic weight spectral sequence of the log crystalline cohomology 
of the product of two proper SNCL schemes over $s$ by using 
the log de Rham-Witt complexes of it (\cite{pcara}). 
As a corollary, she has obtained the limit weight spectral sequence 
of the de Rham cohomology of the generic fiber of the product of 
two proper strict semistable families 
over a complete discrete valuation ring of mixed characteristics 
with residue field $\kap$ by using Hyodo-Kato isomorphism in \cite{hk}. 
\par 
Now we state our results very roughly in each chapter. 
In \cite{nh2} we have constructed the weight spectral sequence of 
the log crystalline cohomology sheaf 
of an open but proper log scheme in characteristic $p>0$ and 
the weight filtration on it, while, in the Chapter I and the Chapter II,
we construct the weight spectral sequence of 
the log crystalline cohomology sheaf 
of a proper SNCL scheme  
in characteristic $p>0$ and 
the weight filtration on it.    
We can consider the Chapter I with the Chapter II of this book 
as a generalization of Mokrane's work to the case where the base scheme is 
any $p$-adic formal family of log points. 
Indeed, in the case where the base log scheme is a  
log point of a perfect field of characteristic $p>0$,  
we prove a comparison theorem between 
the $p$-adic filtered Steenbrink complex in the Chapter I and 
the filtered complex $({\cal W}A^{\bul}_X,P)$ constructed in 
\cite{msemi} and \cite{ndw}. 
The method in this book are different from that in \cite{msemi}: 
we use log crystalline methods mainly as in \cite{nh2}. 
We think that it is the most fundamental method for the construction 
of the $p$-adic weight spectral sequence of a proper SNCL scheme 
over a $p$-adic formal family of log points. 
To develop theory of the $p$-adic weight spectral sequence 
by using log crystalline method is one of our main aims to write this book. 
The point for the construction of the $p$-adic weight spectral sequence 
is the construction of the {\it zariskian $p$-adic filtered Steenbrink complex} 
$(A_{\rm zar},P)$ in a filtered derived category, 
which will be explained later 
(see the explanation for (1) below for more details). 
\par 
More generally, in the Chapter I we construct 
the weight spectral sequence of 
the log crystalline cohomology sheaf 
of a proper split truncated simplicial SNCL scheme 
in characteristic $p>0$ and 
the weight filtration on it. 
More generally again, we construct 
the weight spectral sequence of 
the log crystalline cohomology sheaf 
of a proper truncated simplicial SNCL scheme 
which has the disjoint union of the member of 
an affine truncated simplicial open covering (as in \cite{nh3}) 
and we construct the weight filtration on the log crystalline cohomology sheaf of it 
(the split truncated simplicial SNCL scheme has this disjoint union).  
In this book we call this proper truncated simplicial SNCL scheme 
{\it a proper truncated simplicial SNCL scheme having an  
affine truncated simplicial open covering}. 
This truncated simplicial SNCL scheme  will also be be necessary 
in a future paper when we shall discuss the product structure of 
the log crystalline cohomology sheaf of it with the weight filtration 
because the self-product of 
a split truncated simplicial SNCL scheme is not necessarily split 
but it is a truncated simplicial SNCL scheme which has the disjoint union.  
\par 
From the Chapter III to the Chapter VI, we consider log isocrystalline cohomologies 
and do not consider log crystalline cohomologies; we ignore 
the torsions of log crystalline cohomologies. 
In the Chapter III we define a new proper split truncated simplicial log smooth scheme 
in characteristic $p>0$ which is more general than 
the split truncated simplicial SNCL scheme in the Chapter I.  
This new notion is useful for the construction of the limit of the weight filtration 
on the infinitesimal cohomology of a proper scheme over $K$ 
because we cannot use du Bois' argument in \cite{db} 
(in this book) in which 
he has used the semistable reduction in characteristic $0$; 
if we can use du Bois' argument in any characteristic, 
the split truncated simplicial SNCL scheme in the Chapter I 
is sufficient for the construction. 
We shall give the definition of 
this general split truncated simplicial scheme in the Introduction.  
We call the split truncated simplicial scheme in the Chapter I 
(resp.~that in the Chapter III) 
the ``naive" split truncated  simplicial SNCL scheme  
(resp.~the ``successive" split truncated  simplicial SNCL scheme). 
We construct a key $p$-adic filtered Steenbrink complex modulo torsion 
``$(A_{{\rm zar},{\mab Q}},P)$'' which produces the weight spectral sequence of 
the log isocrystalline cohomology sheaf of 
the successive split truncated  simplicial SNCL scheme 
in characteristic $p>0$ and 
the weight filtration on it.   
\par 
In the Chapter III we also 
give a framework to 
treat the $p$-adic filtered Steenbrink complex modulo torsion
of a proper truncated simplicial SNCL scheme having an  
affine truncated simplicial open covering and that of 
a successive split truncated  simplicial SNCL scheme 
at the same time by introducing new notions 
a {\it truncated simplicial base change of SNCL schemes} 
and an {\it admissible immersion}. 
Though a truncated simplicial SNCL scheme is not stable 
under the base change of families of log points, 
a truncated simplicial base change of SNCL schemes is  
stable under the base change by definition. 
In this chapter we construct 
the $p$-adic filtered Steenbrink complex modulo torsion 
of a truncated simplicial base change of SNCL schemes
with an admissible immersion.
In the Chapter IV we define the $p$-adic filtered Steenbrink complex modulo torsion 
of it over a sequence of log points of 
perfect fields of characteristic $p>0$ by using filtered complexes of 
log de Rham-Witt complexes. 
By using these $p$-adic filtered Steenbrink complexes in the Chapters III and IV, 
we construct the weight filtrations on the log isocrystalline cohomology (sheaf) of 
a (truncated) simplicial base change of SNCL schemes. 
We prove a comparison theorem between 
the $p$-adic filtered Steenbrink complex in the Chapter III and 
that in the Chapter IV in the case where a base sequence of log schemes is 
a sequence of log points of perfect fields of characteristic $p>0$. 
This is a generalization of the comparison theorem already mentioned 
between our $p$-adic Steenbrink complex in the Chapter I modulo torsion 
and the filtered complex $({\cal W}A^{\bul}_X,P)$ modulo torsion 
\par 
In the Chapter V we prove various important properties of 
the weight filtration on the log isocrystalline cohomology sheaf of 
a truncated simplicial base change of SNCL schemes 
which has an admissible immersion. 
We explain these in more details soon later. 
By virtue of fundamental results before 
this chapter, we can obtain fruitful results in this chapter.  
\par  
Using results in the Chapter III and V, in the Chapter VI   
we prove the existence of the limit of the weight filtration 
and the existence of the action of ${\rm WD}_{\rm crys}(\ol{K})$, 
especially  the existences of the $p$-adic monodromy operator and 
the Frobenius endomorphism on the infinitesimal cohomology 
of a proper scheme over $K$.
Using results in the Chapter IV,  
we calculate the limit of the slope filtration on 
the infinitesimal cohomology in a geometric way.  
We also give fundamental properties 
of  the limit of the weight filtration, the limit of the slope filtration 
and the $p$-adic monodromy operator.  
We can consider the Chapter VI 
as a generalization of Caraiani's work because 
we construct the limit weight spectral sequence of the infinitesimal cohomology
(=a generalization of the de Rham cohomology) of 
any proper scheme over $K$ with any singularity. 
\par 
We give the following picture to visualize the relations 
between this book and \cite{nh2} and \cite{nh3}:  
\begin{equation*} 
\begin{CD} 
{\rm \cite{nh2}} @>{{\rm (truncated~simplicial)~SNCL~versions}}>>  
{\rm Chapter~I, II~and~Chapter~III, IV, V}\\ 
@VVV  @VVV \\
{\rm \cite{nh3}}  @>{{\rm a~mixed~characteristics~version}}>> {\rm Chapter~VI}.  
\end{CD}
\end{equation*}  
Here, by the vertical arrows,  
we mean that the targets of them are important applications of 
the sources of them. 
{\rm \cite{nh3}} and Chapter I are related by 
$p$-adic local invariant cycle conjecture by Chiarellotto (\cite{cric}), 
which will be discussed in another paper. 
(See \cite{st1} and \cite{sam} for the local invariant cycle theorem in 
the case of the characteristic $0$.)
\par 
In the Chapter VII we define the action of ${\rm WD}_{\rm crys}(\ol{K})$ 
on the infinitesimal cohomology of a separated scheme of finite type over $\ol{K}$
and we prove the contravariant functoriality of this action and the compatibility with 
the cup product of the infinitesimal cohomology. 
\medskip 
\par
Because we have already stated our main results very roughly, 
we next give explanations for results in this book in more details. 
The reader can use the following explanations as a resume of this book. 
The main results in this book are the following 27 pieces of theorems: 
\smallskip 
\parno
(1): the construction of the $p$-adic weight spectral sequence of  
the log crystalline cohomology sheaf 
(resp.~the log isocrystalline cohomology sheaf)
of a proper naive split truncated  simplicial SNCL scheme 
(resp.~a proper truncated simplicial base change of SNCL schemes 
with admissible immersions) 
over a family of log points,  
\parno
(2): Independence of the (pre)weight filtration with respect to 
the log structure of the base change of a base log scheme, 
\parno
(3): the contravariant functoriality of the weight filtration 
on the log (iso)crystalline cohomology sheaf,  
\parno
(4): the base change theorem of the weight filtration on 
the log (iso)crystalline cohomology sheaf, 
\parno 
(5): the log infinitesimal deformation invariance of 
the pull-back on log isocrystalline cohomology sheaves 
with weight filtrations with respect to 
a morphism of proper truncated simplicial 
base changes of SNCL schemes,  
\parno 
(6): the $E_2$-degeneration of the $p$-adic weight spectral sequence modulo torsion, 
\parno
(7):  the development of a general theory of $F^{\infty}$-isospans,   
\parno
(8): the (log) convergence of the weight filtration,   
\parno
(9): the strict compatibility of the weight filtration modulo torsion 
with respect to the pull-back of a morphism of 
proper truncated simplicial base changes of SNCL schemes,  
\parno 
(10): to give the formulations of 
a variational $p$-adic monodromy-weight conjecture   
and a variational filtered log $p$-adic hard Lefschetz conjecture 
for a projective SNCL scheme over a family of log points,
\parno 
(11): to prove that the variational $p$-adic monodromy-weight conjecture 
is true for the log special fiber of  
a proper strict semistable family over a complete discrete valuation ring $\Delta$ of 
equal positive characteristic $p>0$ 
and a generalization of this result,   
\parno 
(12): to prove that the variational log $p$-adic hard Lefschetz conjecture 
is ture for the log special fiber of  
a projective strict semistable family over $\Delta$ 
and a generalization of this result, 
\parno 
(13): to prove that the variational $p$-adic monodromy-weight conjecture 
and the variational log $p$-adic hard Lefschetz conjecture 
are true for the log special fiber of 
a projective SNCL family over a complete discrete valuation ring of 
mixed characteristics and a generalization of this result, 
\parno 
(14): a log version of 
Etesse's comparison theorem (\cite{et}) 
between the crystalline complex of 
a flat coherent crystal and 
the de Rham-Witt complex of it,  
\parno
(15): a new definition of ``modified'' 
$P$-filtered(=Poincar\'{e} filtered) log de Rham-Witt complexes 
without using the admissible lift in \cite{hdw} and \cite{msemi}, 
\parno
(16): a comparison theorem between 
the $p$-adic filtered Steenbrink complex constructed in this book  
and the filtered 
complex $({\cal W}A_X,P)$ constructed in \cite{msemi} 
and \cite{ndw};   
in fact, we generalize this comparison theorem for flat coherent crystals whose 
corresponding integrable connections have no log poles,   
\parno 
(17): to give the definitions of 
the truncated simplicial base change of SNCL schemes and 
the admissible immersion, 
\parno 
(18): to construct the ``canonical'' Hyodo-Kato isomorphism for 
a proper log smooth $N$-truncated simplicial log scheme ${\cal Y}_{\bul \leq N}$ 
over a log scheme $S$ whose underlying scheme is ${\rm Spec}({\cal V})$ and 
whose log structure is canonical such that the log special fiber of 
${\cal Y}_{\bul \leq N}$ is of Cartier type over the log special fiber of $S$, 
\parno 
(19): 
the construction of the limit weight spectral sequence of 
$H^q_{\rm inf}({\mathfrak X}/K)$ (by extending $K$) 
in a geometric way, 
\parno 
(20): the $E_2$-degeneration of the weight spectral sequence in (19), 
\parno 
(21): the well-definedness of the induced filtration 
on $H^q_{\rm inf}({\mathfrak X}/K)$ by the spectral sequence in (19),  
\parno 
(22): the strict compatibility of the filtration in (21) with respect 
to the pull-back of a morphism of proper schemes ${\mathfrak X}$'s,  
\parno 
(23): to find a canonical $K_0$-structure $H^q_{\rm inf}({\mathfrak X}/K/K_0)$ 
of $H^q_{\rm inf}({\mathfrak X}/K)$,  
\parno 
(24): the calculation of the limit of the slope filtration on 
$H^q_{\rm inf}({\mathfrak X}/K/K_0)$ (by extending $K_0$) 
in a geometric way, 
\parno 
(25): to find a canonical increasing filtration on $H^q_{\rm inf}({\mathfrak X}/K/K_0)$ 
which is a $K_0$-form of the filtration in (21), 
\parno 
(26): to define the action of ${\rm WD}_{\rm crys}(\ol{K})$ on 
$H^q_{\rm inf}({\mathfrak X}_{\ol{K}}/\ol{K})$   
(This leads to the conjecture of the compatibility 
of the linear action of ${\rm WD}(\ol{K}/K)$
on $H^q_{\rm inf}({\mathfrak X}_{\ol{K}}/\ol{K})$ 
with the action of ${\rm WD}(\ol{K}/K)$
on $H^q_{\rm et}({\mathfrak X}_{\ol{K}},{\mab Q}_l)$.),  
\parno
(27): to prove that the action of ${\rm WD}_{\rm crys}(\ol{K})$ on 
$H^q_{\rm inf}({\mathfrak X}_{\ol{K}}/\ol{K})$ is contravariant 
and that the action of ${\rm WD}_{\rm crys}(\ol{K})$ 
preserves the limit of the weight filtration 
$P$ on $H^q_{\rm inf}({\mathfrak X}_{\ol{K}}/\ol{K})$. 
In fact, we prove the existence of the contravariant action of 
${\rm WD}_{\rm crys}(\ol{K})$ on the infinitesimal cohomology of 
a separated scheme of finite type over $K$.  
\par 
\medskip 
\par
In the following we explain each of the main results above. 
In this introduction we first explain (1) 
for the constant simplicial case(=the case of the usual SNCL scheme) 
because this case is more accessible than general truncated simplicial cases. 
We consider a certain nontrivial coefficient(=a flat coherent crystal
whose corresponding integrable connection has no log poles) 
following a suggestion by T.~Tsuji. 
\par 
Let $T_0 \os{\sus}{\lo} T$ be an exact closed immersion defined by 
a quasi-coherent PD-ideal sheaf ${\cal J}$ with PD-structure $\del$. 
For a log scheme $Y$ 
(resp.~a morphism $g\col Y\lo Z$ of log schemes), 
we denote by $\os{\circ}{Y}$ 
(resp.~$\os{\circ}{g}\col \os{\circ}{Y}\lo 
\os{\circ}{Z}$)
the underlying scheme of $Y$ 
(resp.~the underlying morphism of schemes of $g$). 
For a log scheme $Y$ over $T_0$, 
let  $\eps_{Y/T_0}\col Y\lo \os{\circ}{Y}$ 
be the natural morphism over $T_0\lo \os{\circ}{T}_0$ 
forgetting the log structure of $Y$. 
Let $g\col Y\lo T$ be the structural morphism. 
Let $(Y/T)_{\rm crys}$ be 
the log crystalline topos of $Y/(T,{\cal J},\del)$ 
defined by K.~Kato in \cite{klog1}; 
we omit to write the notation 
``log'' and  $\tilde{\quad}$ 
for the log crystalline topos 
$(\wt{Y/T})^{\log}_{\rm crys}$ in \cite{klog1} 
and so on; 
we follow the notation of \cite{bb} and \cite{bob} 
by omitting $\tilde{\quad}$.
Let $((Y/T)_{\rm crys},{\cal O}_{Y/T})$ be 
the log crystalline ringed topos of $Y/(T,{\cal J},\del)$. 
Let 
$$\eps_{Y/T} \col  ((Y/T)_{\rm crys},{\cal O}_{Y/T}) \lo 
((\os{\circ}{Y}/\os{\circ}{T})_{\rm crys},
{\cal O}_{\os{\circ}{Y}/\os{\circ}{T}})$$ 
be the induced morphism by $\eps_{Y/T_0}$.  
Let 
$$u_{Y/T} \col ((Y/T)_{\rm crys},{\cal O}_{Y/T})
\lo (\os{\circ}{Y}_{\rm zar},g^{-1}({\cal O}_T))$$ 
be the canonical projection. 
\par 
Let $S$ be a log scheme on which 
a prime number $p$ is locally nilpotent 
on $\os{\circ}{S}$ and whose log structure is 
Zariski locally isomorphic to ${\mab N}\oplus {\cal O}_S^*$ 
with structural morphism 
${\mab N}\oplus {\cal O}_S^*\lo {\cal O}_S$ 
defined by $(a,u)\lom 0$  
$((a,u)\in {\mab N}\oplus {\cal O}_S^*, a\not=0)$
and $(0,u)\lom u$.  
We call $S$ a {\it family of log points}. 
Let $S(0) \os{\sus}{\lo} S$ be an exact closed immersion defined by 
a quasi-coherent PD-ideal sheaf ${\cal I}$ with PD-structure $\gam$. 
We call $(S,{\cal I},\gam)$ a PD-family of log points. 
Let $f\col X \lo S(0)$ be an SNCL scheme  
(in \S\ref{sec:snclv} in the text 
we will give the definition of the SNCL scheme). 
Intuitively $X/S(0)$ is a family of SNCL schemes over log points. 
By abuse of notation, we denote by $f$ 
the composite structural morphism 
$f\col X\lo S(0)\os{\sus}{\lo}S$. 
Let $E$ be a flat quasi-coherent crystal of  
${\cal O}_{\os{\circ}{X}/\os{\circ}{S}}$-modules. 
Let ${\rm D}^+{\rm F}(f^{-1}({\cal O}_S))$ 
be the derived category of bounded below filtered complexes 
of $f^{-1}({\cal O}_S)$-modules (see, e.~g.,~\cite{dh3}, \cite{blec} 
for the derived category ${\rm D}^+{\rm F}(f^{-1}({\cal O}_S)))$. 
(In \cite{nh2} we explain a part of Berthelot's theory of 
filtered derived category in \cite{blec}.) 
Then we construct a new filtered complex 
\begin{equation*} 
(A_{\rm zar}(X/S,\eps_{X/\os{\circ}{S}}^*(E)),P)
\in {\rm D}^+{\rm F}(f^{-1}({\cal O}_S)) 
\tag{0.0.1.2}\label{eqn:apea}
\end{equation*} 
depending only on $X/(S,{\cal I},\gam)$ and $E$. 
Here we have considered the morphism 
$X\lo \os{\circ}{X}$ over $\os{\circ}{S}$ 
by considering the composite morphism $X\lo S\lo \os{\circ}{S}$. 
Because the notation $\eps_{X/\os{\circ}{S}}^*$ will become too complicated 
soon later (because we consider a base change morphism of $S$), 
we denote $(A_{\rm zar}(X/S,\eps_{X/\os{\circ}{S}}^*(E)),P)$ 
by $(A_{\rm zar}(X/S,E),P)$ by abuse of notation. 
When $E={\cal O}_{\os{\circ}{X}/\os{\circ}{S}}$, 
we denote $(A_{\rm zar}(X/S,E),P)$ by 
$(A_{\rm zar}(X/S),P)$,  
which is a $p$-adic variational analogue of  
Steenbrink's filtered complex in \cite{st1}, \cite{kwn} and \cite{fn}. 
We call $(A_{\rm zar}(X/S,E),P)$ 
the {\it zariskian $p$-adic filtered Steenbrink complex} of 
$\eps_{X/\os{\circ}{S}}^*(E)$ and 
call $(A_{\rm zar}(X/S),P)$ the {\it zariskian $p$-adic filtered Steenbrink complex} 
of $X/S$. 
Though $(A_{\rm zar}(X/S),P)$ is a $p$-adic analogue of 
the Steenbrink filtered complexes in \cite{st1} and \cite{fn}, 
the convention on signs of certain filtered double complexes 
are different from them. 
As in the open log case in \cite{nh2}, the filtered derived category 
${\rm D}^+{\rm F}(f^{-1}({\cal O}_S))$ in which 
$(A_{\rm zar}(X/S,E),P)$ exists 
is an indispensable language in this book (and in \cite{gkwf}, 
though the filtered derived category has not appeared in [loc.~cit.]). 
The filtered complex $(A_{\rm zar}(X/S,E),P)$, 
several generalizations of $(A_{\rm zar}(X/S,E),P)$
and analogous filtered complexes to $(A_{\rm zar}(X/S,E),P)$ 
which will be explained later are core objets in our theory. 
\par 
If we forget the filtration $P$, 
then we prove that there exists a canonical isomorphism 
\begin{equation*} 
Ru_{X/S*}(\eps^*_{X/S}(E)) 
\os{\sim}{\lo} A_{\rm zar}(X/S,E)
\tag{0.0.1.3}\label{eqn:canea}
\end{equation*}
in ${\rm D}^+(f^{-1}({\cal O}_S))$. 
Here we have considered the morphism $X\lo \os{\circ}{X}$ 
over $S\lo \os{\circ}{S}$. 
For a nonnegative integer $k$, 
let $\os{\circ}{X}{}^{(k)}$ 
be the scheme over $\os{\circ}{S}(0)$ 
associated to $\os{\circ}{X}$ 
whose definition will be given in \S\ref{sec:ldc} below. 
(This $\os{\circ}{X}{}^{(k)}$ is a generalization of 
$\os{\circ}{X}{}^{(k)}$ in (\ref{eqn:eswsp}).)  
Let $a^{(k)}\col 
\os{\circ}{X}{}^{(k)} \lo \os{\circ}{X}$ 
be the natural morphism of schemes defined in 
\S\ref{sec:ldc} below. 
Let $a^{(k)}_{\rm crys} \col 
((\os{\circ}{X}{}^{(k)}/\os{\circ}{S})_{\rm crys},
{\cal O}_{\os{\circ}{X}{}^{(k)}/\os{\circ}{S}})
\lo 
((\os{\circ}{X}/\os{\circ}{S})_{\rm crys},
{\cal O}_{\os{\circ}{X}/\os{\circ}{S}})$ 
be the induced morphism of ringed topoi by $a^{(k)}$. 
Set 
$E_{\os{\circ}{X}{}^{(k)}/\os{\circ}{S}}:=
a^{(k)*}_{\rm crys}(E)$. 
Let 
$\vp^{(k)}_{\rm crys}(\os{\circ}{X}/\os{\circ}{S})$ 
be the crystalline orientation sheaf defined in 
\S\ref{sec:ldc} below: 
$\vp^{(k)}_{\rm crys}(\os{\circ}{X}/\os{\circ}{S})$ 
is isomorphic to ${\mab Z}$ in 
$(\os{\circ}{X}{}^{(k)}/\os{\circ}{S})_{\rm crys}$ if 
$\os{\circ}{X}$ is small. 
Let $[-m]$ $(m\in {\mab Z})$ be 
the operation shifting a complex to the right by the $m$ degree with changing 
the signs of the boundary morphisms by $(-1)^m$ as usual.  
Then we prove that, for $k\in {\mab Z}$,  
there exists a canonical isomorphism 
\begin{align*} 
{\rm gr}^P_kA_{\rm zar}(X/S,E)\os{\sim}{\lo}
\bigoplus_{j\geq \max \{-k,0\}} 
a^{(2j+k)}_{*} & 
Ru_{\os{\circ}{X}{}^{(2j+k)}/\os{\circ}{S}*}
(E_{\os{\circ}{X}{}^{(2j+k)}
/\os{\circ}{S}}
\otimes_{\mab Z} 
\tag{0.0.1.4}\label{eqn:caxd} \\
& 
\vp^{(2j+k)}_{\rm crys}(\os{\circ}{X}/\os{\circ}{S}))
(-j-k)[-2j-k]
\end{align*}
in ${\rm D}^+(f^{-1}({\cal O}_S))$. 
Here $(-j-k)$ means the nonstandard Tate twist, 
which we will define in the explanation for (3). 
Set $f_{X/S}:=f\circ u_{X/S}$. 
As a corollary of (\ref{eqn:canea}) and (\ref{eqn:caxd}), 
we obtain the following spectral sequence 
\begin{align*} 
E_1^{-k,q+k}=
\bigoplus_{j\geq \max \{-k,0\}} &
R^{q-2j-k}f_{\os{\circ}{X}{}^{(2j+k)}/\os{\circ}{S}*}
(E_{\os{\circ}{X}{}^{(2j+k)}/\os{\circ}{S}}
\otimes_{\mab Z}
\vp^{(2j+k)}_{\rm crys}(\os{\circ}{X}/\os{\circ}{S}))(-j-k) 
\tag{0.0.1.5}\label{eqn:espsp} \\
&  \Lo 
R^qf_{X/S*}(\eps^*_{X/S}(E)) 
\quad (q\in {\mab Z}). 
\end{align*}  
When $E={\cal O}_{\os{\circ}{X}/\os{\circ}{S}}$, 
(\ref{eqn:espsp}) gives us 
the preweight filtration $P$ on $R^qf_{X/S*}({\cal O}_{X/S})$.
(Here we use the terminology ``preweight filtration'' instead of 
``weight filtration''  because ${\cal O}_S$ is $p$-primary torsion.)  
\par  
More generally, 
we need to construct a generalization of the filtered complex 
$(A_{\rm zar}(X/S),P)$ for any base change 
$(T,{\cal J},\del)\lo (S,{\cal I},\gam)$ of fine log PD-schemes 
such that ${\cal J}$ is quasi-coherent
(we consider $(T,{\cal J},\del)$ as 
an ``extension'' and a ``ramification'' of $(S,{\cal I},\gam)$) 
for applications in the Chapters V and VI.  
(In \cite{nh2} it is not necessary to consider 
this formalism in the open log case 
since the base change of a smooth scheme 
with a relative SNCD by any morphism of base schemes 
is a smooth scheme with a relative SNCD again.) 
\par 
We change notations by introducing 
a ``$p$-primary torsion version'' of a log PD-enlargement
defined by Ogus in \cite{ollc} as follows. 
This variant gives us an essential language. 
For example, we need the notion $S^{[p]}(T)^{\nat}$ which will appear  
in the explanation for (5) later.  
\par 
Let $S$ be a family of log points and 
let $(T,{\cal J},\del)$ be a fine log PD-scheme 
such that ${\cal J}$ is quasi-coherent and 
such that $p$ is locally nilpotent on $T$. 
Let $T_0$ be the exact closed log subscheme of $T$ 
defined by ${\cal J}$. 
Assume that we are given a morphism $T_0\lo S$ of fine log schemes. 
We call $(T,{\cal J},\del)$ with morphism 
$T_0\lo S$ a log PD-enlargement of $S$. 
As in the trivial log case in \cite{bb}, we can define the big log crystalline 
site ${\rm Crys}(S/({\rm Spf}({\mab Z}_p),p{\mab Z}_p,[~]))$.   
Then $(T,{\cal J},\del)$ with morphism 
$T_0\lo S$ is an object of this site. 
Set $S_{\os{\circ}{T}_0}:=S\times_{\os{\circ}{S}}\os{\circ}{T}_0$.  
The log structure of $S_{\os{\circ}{T}_0}$ is nothing but 
the inverse image of the log structure of $S$ by the morphism 
$\os{\circ}{T}_0\lo \os{\circ}{S}$; consequently 
$S_{\os{\circ}{T}_0}$ is a family of log points.   
Then we have a natural morphism $T_0\lo S_{\os{\circ}{T}_0}$ 
of fine log schemes. 
For a fine log scheme $Y$, denote by $M_Y$ the log structure of $Y$. 
It is clear that the natural morphism 
$M_{S_{\os{\circ}{T}_0}}\lo M_{T_0}$ is injective. 
Let $S(T)$ be a fine log scheme whose underlying scheme is $\os{\circ}{T}$ 
and whose log structure $M_{S(T)}$ is characterized 
by the following commutative diagram 
\begin{equation*} 
\begin{CD}
M_T/{\cal O}_T^*@>{\sim}>>M_{T_0}/{\cal O}_{T_0}^*\\
@A{\bigcup}AA @AA{\bigcup}A \\
M_{S(T)}/{\cal O}_T^*@>{\sim}>>M_{S_{\os{\circ}{T}_0}}/{\cal O}_{T_0}^*.
\end{CD}
\end{equation*} 
Obviously we have a natural morphism 
$(T,{\cal J},\del)\lo (S(T),{\cal J},\del)$ of fine log PD-schemes. 
Though 
$S(T)$ is not necessarily hollow in the sense of Ogus' article \cite{ollc}, 
we can consider the hollowing out $S(T)^{\nat}$ of $S(T)$ in [loc.~cit.]  
because $M_{S(T)^{\nat}}/{\cal O}_T^*$ is a constant sheaf of monoids 
on $\os{\circ}{T}$.  
(Recall that, for a log scheme $Z$, 
$M_Z:=(M_Z,\al_Z)$ is said to be hollow in [loc.~cit.] if 
any section of $M_{Z}\setminus {\cal O}^*_Z$ is mapped to $0$ in 
${\cal O}_Z$ by the structural morphism $\al_Z\col M_Z\lo {\cal O}_Z$.) 
Then $S(T)^{\nat}$ is a family of log points and 
$S_{\os{\circ}{T}_0}$ is an exact closed log subscheme 
of $S(T)^{\nat}$ defined by ${\cal J}$. 
Let $U_0$ (resp.~$U$) be $S_{\os{\circ}{T}_0}$ or $T_0$ 
(resp.~$S(T)^{\nat}$ or $T$).  
Set $X_{U_0}:=X\times_{S}U_0$ and 
$\os{\circ}{X}_{U_0}:=\os{\circ}{(X_{U_0})}$ 
(the underlying scheme of $X_{U_0}$). 
Then $X_{S_{\os{\circ}{T}_0}}=X\times_{S}S_{\os{\circ}{T}_0}=
X\times_{\os{\circ}{S}}\os{\circ}{T}_0=:X_{\os{\circ}{T}_0}$ 
(we have considered $X$ as a fine log scheme over $\os{\circ}{S}$) 
and $X_{T_0}= X_{S_{\os{\circ}{T}_0}}\times_{S_{\os{\circ}{T}_0}}T_0$. 
Note that $X_{\os{\circ}{T}_0}/S_{\os{\circ}{T}_0}$ is an SNCL scheme, 
but $X_{T_0}/T_0$ is not necessarily an SNCL scheme. 
Let $f_U\col X_{U_0}\lo U$ be the structural morphism. 
Because it is easy to check an equality $\os{\circ}{X}_{U_0}
=\os{\circ}{X}\times_{\os{\circ}{S}}\os{\circ}{T}_0$, 
we see that $\os{\circ}{X}_{S_{\os{\circ}{T}_0}}=\os{\circ}{X}_{T_0}$ 
and $\os{\circ}{f}_{S(T)^{\nat}}=\os{\circ}{f}_T$.
Hence we have the following morphism 
\begin{align*} 
\eps_{X_{\os{\circ}{T}_0}/(S(T)^{\nat})^{\circ}} \col 
((X_{\os{\circ}{T}_0}/(S(T)^{\nat})^{\circ})_{\rm crys},
{\cal O}_{X_{\os{\circ}{T}_0}/(S(T)^{\nat})^{\circ}}) 
\lo &(({(X_{\os{\circ}{T}_0})}^{\circ}/(S(T)^{\nat})^{\circ})_{\rm crys},
{\cal O}_{{(X_{\os{\circ}{T}_0})^{\circ}}/(S(T)^{\nat})^{\circ}})\\
& =((\os{\circ}{X}_{T_0}/\os{\circ}{T})_{\rm crys},
{\cal O}_{{\os{\circ}{X}_{T_0}}/\os{\circ}{T}}). 
\end{align*}  
Let $E$ be a flat quasi-coherent crystal of  
${\cal O}_{\os{\circ}{X}_{T_0}/\os{\circ}{T}}$-modules. 
Then we construct a filtered complex  
$$(A_{\rm zar}(X_{\os{\circ}{T}_0}/S(T)^{\nat},
\eps^*_{X_{\os{\circ}{T}_0}/(S(T)^{\nat})^{\circ}}(E)),P)=
(A_{\rm zar}(X_{\os{\circ}{T}_0}/S(T)^{\nat},
\eps^*_{X_{\os{\circ}{T}_0}/\os{\circ}{T}}(E)),P),$$ 
which is a generalization of $(A_{\rm zar}(X/S,E),P)$.   
Because this notation is too long for us, we denote it by 
$(A_{\rm zar}(X_{\os{\circ}{T}_0}/S(T)^{\nat},E),P)$ for simplicity of notation. 
When $E={\cal O}_{\os{\circ}{X}_{T_0}/\os{\circ}{T}}$, 
we denote $(A_{\rm zar}(X_{\os{\circ}{T}_0}/S(T)^{\nat},E),P)$ 
by $(A_{\rm zar}(X_{\os{\circ}{T}_0}/S(T)^{\nat}),P)$. 
As a generalization of (\ref{eqn:canea}), 
we prove that there exists a canonical isomorphism 
\begin{equation*} 
Ru_{X_{\os{\circ}{T}_0}/S(T)^{\nat}*}
(\eps^*_{X_{\os{\circ}{T}_0}/S(T)^{\nat}}(E))\os{\sim}{\lo} 
A_{\rm zar}(X_{\os{\circ}{T}_0}/S(T)^{\nat},E)
\tag{0.0.1.6}\label{eqn:catea}
\end{equation*}
in ${\rm D}^+(f^{-1}_T({\cal O}_T))$.   
If $S(T)$ is hollow, then $S(T)=S(T)^{\nat}$ and 
we have the following canonical isomorphism 
\begin{equation*} 
Ru_{X_{T_0}/T*}(\eps^*_{X_{T_0}/T}(E))
\os{\sim}{\longleftarrow}
Ru_{X_{\os{\circ}{T}_0}/S(T)^{\nat}*}(\eps^*_{X_{\os{\circ}{T}_0}/S(T)^{\nat}}(E)) 
\tag{0.0.1.7}\label{eqn:ectsp}
\end{equation*} 
by using the natural morphism $(T,{\cal J},\del)\lo (S(T),{\cal J},\del)$ 
and noting that the morphism $\os{\circ}{T}\lo (S(T))^{\circ}$ 
is the identity of $\os{\circ}{T}$. 
Consequently, if $S(T)$ is hollow, then we have the following 
canonical isomorphism 
\begin{equation*} 
Ru_{X_{T_0}/T*}(\eps^*_{X_{T_0}/T}(E)) \os{\sim}{\lo}
A_{\rm zar}(X_{\os{\circ}{T}_0}/S(T),E). 
\tag{0.0.1.8}\label{eqn:ecetsp}
\end{equation*} 
Set 
$$f_{X_{T_0}/T}:=f_T\circ u_{X_{T_0}/T}\col 
((X_{T_0}/T)_{\rm crys},{\cal O}_{X_{T_0}/T}) 
\lo (\os{\circ}{T}_{\rm zar},{\cal O}_T).$$  
Set also 
$\os{\circ}{X}{}^{(k)}_{T_0}
:=\os{\circ}{X}{}^{(k)}\times_{\os{\circ}{S}}\os{\circ}{T}_0$ $(k\in {\mab N})$. 
Let  
$$a^{(k)}_{T{\rm crys}} \col 
((\os{\circ}{X}{}^{(k)}_{T_0}/\os{\circ}{T})_{\rm crys},
{\cal O}_{\os{\circ}{X}{}^{(k)}_{T_0}/\os{\circ}{T}})
\lo 
((\os{\circ}{X}_{T_0}/\os{\circ}{T})_{\rm crys},
{\cal O}_{\os{\circ}{X}_{T_0}/\os{\circ}{T}})$$ 
be the natural morphism. 
Set  
$E_{\os{\circ}{X}{}^{(k)}_{T_0}/\os{\circ}{T}}
:=
a^{(k)*}_{T{\rm crys}}(E)$ 
and let 
$\vp^{(k)}_{\rm crys}(\os{\circ}{X}_{T_0}/\os{\circ}{T})$ 
be the crystalline orientation sheaf of $\os{\circ}{X}_{T_0}/\os{\circ}{T}$. 
Then the graded complex of the filtered complex 
$(A_{\rm zar}(X_{\os{\circ}{T}_0}/S(T)^{\nat},E),P)$
can be calculated as follows: 
\begin{align*} 
{\rm gr}^P_kA_{\rm zar}(X_{\os{\circ}{T}_0}/S(T)^{\nat},E)
\os{\sim}{\lo}
\bigoplus_{j\geq \max \{-k,0\}} 
a^{(2j+k)}_{*} & 
Ru_{\os{\circ}{X}{}^{(2j+k)}_{T_0}/\os{\circ}{T}*}
(E_{\os{\circ}{X}{}^{(2j+k)}
/\os{\circ}{T}}
\otimes_{\mab Z} 
\tag{0.0.1.9}\label{eqn:catxd} \\
& 
\vp^{(2j+k)}_{\rm crys}(\os{\circ}{X}_{T_0}/\os{\circ}{T}))
(-j-k)[-2j-k] 
\end{align*}
in ${\rm D}^+(f^{-1}_T({\cal O}_T))$. 
Hence we have the following spectral sequence of 
the log crystalline cohomology sheaf of 
$\eps^*_{X_{\os{\circ}{T}_0}/S(T)^{\nat}}(E)$ 
on the SNCL scheme $X_{\os{\circ}{T}_0}/S(T)^{\nat}$: 
\begin{align*} 
E_1^{-k,q+k}=
& \bigoplus_{j\geq \max \{-k,0\}} 
R^{q-2j-k}f_{\os{\circ}{X}{}^{(2j+k)}_{T_0}
/\os{\circ}{T}*}
(E_{\os{\circ}{X}{}^{(2j+k)}_{T_0}/\os{\circ}{T}}
\otimes_{\mab Z}
\vp^{(2j+k)}_{\rm crys}(\os{\circ}{X}_{T_0}/\os{\circ}{T}))(-j-k)  
\tag{0.0.1.10}\label{eqn:esbcpsp} \\
& 
\Lo 
R^qf_{X_{T_0}/T*}(\eps^*_{X_{\os{\circ}{T}_0}/S(T)^{\nat}}(E)) \quad (q\in {\mab Z}).   
\end{align*}  
Consequently, if $S(T)$ is hollow, then we have the following spectral sequence 
of the log crystalline cohomology sheaf of 
$\eps^*_{X_{T_0}/T}(E)$ on 
the not necessarily SNCL scheme $X_{T_0}/T$:
\begin{align*} 
E_1^{-k,q+k}=
& \bigoplus_{j\geq \max \{-k,0\}} 
R^{q-2j-k}f_{\os{\circ}{X}{}^{(2j+k)}_{T_0}
/\os{\circ}{T}*}
(E_{\os{\circ}{X}{}^{(2j+k)}_{T_0}/\os{\circ}{T}}
\otimes_{\mab Z}
\vp^{(2j+k)}_{\rm crys}(\os{\circ}{X}_{T_0}/\os{\circ}{T}))(-j-k)  
\tag{0.0.1.11}\label{eqn:esbwtsp} \\
& 
\Lo 
R^qf_{X_{T_0}/T*}(\eps^*_{X_{T_0}/T}(E)) \quad (q\in {\mab Z}).   
\end{align*}  
As a corollary of the existence of 
$(A_{\rm zar}(X_{\os{\circ}{T}_0}/S(T)^{\nat},E),P)$
and the isomorphism (\ref{eqn:ectsp}), 
if $S$ is an exact closed log subscheme of a fine log PD-scheme 
$(\wt{S},{\cal I},\gam)$ defined by ${\cal I}$, if $\wt{S}$ is a family of log points 
and if there exists a morphism $(T,{\cal J},\del)\lo (\wt{S},{\cal I},\gam)$ 
of fine log PD-schemes, then we see that the induced filtration on 
$R^qf_{X_{T_0}/T*}(\eps^*_{X_{T_0}/T}(E))$ 
by the filtered complex 
$(A_{\rm zar}(X_{\os{\circ}{T}_0}/S(T)^{\nat},E),P)$
depends only on the extension morphism 
$(\os{\circ}{T},{\cal J},\del)\lo (\os{\circ}{\wt{S}},{\cal I},\gam)$ 
of $\os{\circ}{T}_0\lo \os{\circ}{S}$;   
it is independent of the morphism $M_{\wt{S}}\lo M_T$ of log structures. 
This is consistent with Kato's log base change theorem in \cite{klog1}:  
$Rf_{X_{T_0}/T*}(\eps^*_{X_{T_0}/T}(E))=
{\cal O}_T\otimes^L_{{\cal O}_{\wt{S}}}
Rf_{X/\wt{S}*}(\eps^*_{X/\wt{S}}(E))$. 
We can say that the preweight filtration on the log crystalline cohomology sheaf 
(in the case $E={\cal O}_{\os{\circ}{X}/\os{\circ}{\wt{S}}}$)
is independent of the log structure of the base change of a base log scheme ((2)). 
\par 
Next we explain (3). 
\par 
An important property of 
$(A_{\rm zar}(X_{\os{\circ}{T}_0}/S(T)^{\nat},E),P)$ 
is the contravariant functoriality as in the open log case in \cite{nh2}, 
which has not been investigated for analogues of  
$(A_{\rm zar}(X/S),P)$ in any characteristic 
except the Frobenius endomorphism in \cite{ndw} as far as I know.  
(See \cite[\S9]{ndw} for what has been proved before [loc.~cit.] and 
what has not before it. 
The Frobenius action on 
Mokrane's filtered Steenbrink-Hyodo complex 
$({\cal W}A_X,P)$ defined in \cite{cric} 
is mistaken and the Frobenius action 
on the analogous (filtered) complex in \cite{gkwf} has not been defined.) 
The contravariant 
functorialities of 
$$(A_{\rm zar}(X_{\os{\circ}{T}_0}/S(T)^{\nat},E),P)$$ 
and 
$$(A_{\rm zar}(X_{\os{\circ}{T}_0}/S(T)^{\nat},E),P)
\otimes^L_{\mab Z}{\mab Q}$$ 
(in the case where $\os{\circ}{T}$ is a flat $p$-adic formal scheme)
are more complicated than 
the contravariant functorialty of the open log case in \cite{nh2}. 
We need cares for new phenomenons 
of the contravariant functorialty of 
$(A_{\rm zar}(X_{\os{\circ}{T}_0}/S(T)^{\nat},E),P)$, 
which does not appear in the open log case in [loc.~cit.]. 
The precise meaning of the functoriality  
of $(A_{\rm zar}(X_{\os{\circ}{T}_0}/S(T)^{\nat},E),P)$ is as follows 
(we need a nontrivial work to formulate it). 
\par 
Let $S'$ be another PD-family of log points. 
Let $u\col S \lo S'$ be a morphism of log schemes. 
Let $(T',{\cal J}',\del')$ be a log PD-enlargement of $S'$. 
Let $T'_0\os{\sus}{\lo} T'$ be an exact closed immersion defined by ${\cal J}'$. 
Let $v\col (T,{\cal J},\del) \lo (T',{\cal J}',\del')$ be 
a morphism of log PD-enlargements over $u$, that is, a morphism  
$(T,{\cal J},\del) \lo (T',{\cal J}',\del')$ of log PD-schemes 
fitting into the following commutative diagram 
\begin{equation*} 
\begin{CD} 
\os{\circ}{T}_0 @>{v_0}>> \os{\circ}{T}{}'_0 \\ 
@VVV @VVV \\ 
S @>{u}>> S',  
\end{CD} 
\end{equation*} 
where $v_0\col T_0\lo T'_0$ is the induced morphism by $v$. 
Then it is easy to see that exists a natural morphism
$u\col (S(T)^{\nat},{\cal J},\del) \lo (S'(T')^{\nat},{\cal J}',\del')$ of fine log PD-schemes 
fitting into 
the following commutative diagram 
\begin{equation*} 
\begin{CD} 
S_{\os{\circ}{T}_0}@>>> S'_{\os{\circ}{T}{}'_0}\\
@V{\bigcap}VV @VV{\bigcap}V \\
S(T)^{\nat}@>>> S'(T')^{\nat}. 
\end{CD} 
\end{equation*} 
Let $X'/S'$ be another SNCL scheme. 
Let 
\begin{equation*} 
\begin{CD} 
X_{\os{\circ}{T}_0} @>{g}>> X'_{\os{\circ}{T}{}'_0} \\ 
@VVV @VVV \\ 
S_{\os{\circ}{T}_0} @>>> S'_{\os{\circ}{T}{}'_0} 
\end{CD}
\tag{0.0.1.12}\label{eqn:xxpv} 
\end{equation*} 
be a commutative diagram of SNCL schemes. 
Let 
$$\deg(u) \col \os{\circ}{T} \lo {\mab N}$$ 
be the mapping degree function of $u$ 
which will be defined in \S\ref{sec:snclv} in the text. 
(When $S'(T')^{\nat}=S(T)^{\nat}$, 
this is an analogue of the mapping degree function of 
a continuous endomorphism of ${\mab S}^1$ in topology.)
It is easy to see that $\deg(u)$ never take the value $0$. 
Assume that $\deg(u)$ is not divisible by $p$ for any point $x\in \os{\circ}{T}$. 
Let $g_{T'_0T_0} \col X_{T_0} \lo X'_{T'_0}$ 
be the induced morphism by $g$ and $v_0$.  
Let 
$$g_{T'T{\rm crys}}\col 
((X_{T_0}/T)_{\rm crys},{\cal O}_{X_{T_0}/T}) 
\lo 
((X'_{T'_0}/T')_{\rm crys},{\cal O}_{X'_{T'_0}/T'})$$ 
and 
$$\os{\circ}{g}_{T'T}=\os{\circ}{g}\col 
((\os{\circ}{X}_{T_0})_{\rm zar},f^{-1}_T({\cal O}_T)) 
\lo 
((\os{\circ}{X}{}'_{T'_0})_{\rm zar},f'{}^{-1}_{T'}({\cal O}_{T'}))$$  
be the induced morphisms of topoi by $g_{T'_0T_0}$.    
(For simplicity of notations, we denote the latter morphism 
by $g$ or $g_{T'T}$ appropriately depending on situations.)
Let $E$ (resp.~$E'$) be a flat quasi-coherent crystal of 
${\cal O}_{\os{\circ}{X}_{T_0}/\os{\circ}{T}}$-modules  
(resp.~a flat quasi-coherent crystal of 
${\cal O}_{\os{\circ}{X}{}'_{T'_0}/\os{\circ}{T}{}'}$-modules). 
Assume that we are given a morphism 
$E'\lo \os{\circ}{g}_{{\rm crys}*}(E)$ of 
${\cal O}_{\os{\circ}{X}{}'_{T'_0}/\os{\circ}{T}{}'}$-modules. 
By the contravariant functoriality of 
$(A_{\rm zar}(X_{\os{\circ}{T}_0}/S(T)^{\nat},E),P)$, 
we mean that there exists a morphism 
\begin{align*} 
g^*
\col 
(A_{\rm zar}(X'_{\os{\circ}{T}{}'_0}/S'(T')^{\nat},E'),P)
\lo
Rg_{*}
((A_{\rm zar}(X_{\os{\circ}{T}_0}/S(T)^{\nat},E),P))
\tag{0.0.1.13}\label{eqn:agcxds}
\end{align*}
of filtered complexes fitting into 
the following commutative diagram 
\begin{equation*} 
\begin{CD}
A_{\rm zar}(X'_{\os{\circ}{T}{}'_0}/S'(T')^{\nat},E')
@>{g^*}>> 
Rg_{*}(A_{\rm zar}(X_{\os{\circ}{T}_0}/S(T)^{\nat},E))
\\ 
@A{\simeq}AA 
@AA{\simeq}A \\ 
Ru_{X'_{\os{\circ}{T}{}'_0}/S(T')^{\nat}*}
(\eps^*_{X'_{\os{\circ}{T}{}'_0}/S'(T')^{\nat}}(E'))
@>{g^*_{{\rm crys}}}>> Rg_{*}
Ru_{X_{\os{\circ}{T}_0}/S(T)^{\nat}*}(\eps^*_{X_{\os{\circ}{T}_0}/S(T)^{\nat}}(E))
\end{CD}
\tag{0.0.1.14}\label{eqn:auc} 
\end{equation*}  
satisfying a transitive relation 
``$(h\circ g)^*=Rh_*(g^*)\circ h^*$'' 
and 
${\rm id}^*_{X_{\os{\circ}{T}_0}}=
{\rm id}_{(A_{\rm zar}(X_{\os{\circ}{T}_0}/S(T)^{\nat},E),P)}$. 
If $S(T)$ is hollow, 
we obtain the following commutative diagram 
\begin{equation*} 
\begin{CD}
A_{\rm zar}(X'_{\os{\circ}{T}{}'_0}/S'(T')^{\nat},E')
@>{g^*}>> 
Rg_{*}(A_{\rm zar}(X_{\os{\circ}{T}_0}/S(T)^{\nat},E))
\\ 
@A{\simeq}AA 
@AA{\simeq}A \\ 
Ru_{X'_{T'_0}/T'*}
(\eps^*_{X'_{T'_0}/T'}(E'))
@>{g^*_{T'T{\rm crys}}}>> Rg_{*}
Ru_{X_{T_0}/T*}(\eps^*_{X_{T_0}/T}(E))
\end{CD}
\tag{0.0.1.15}\label{eqn:asuc} 
\end{equation*}  
by (\ref{eqn:ectsp}) and (\ref{eqn:auc}). 
Even in the case where 
$\deg(u)$ is divisible by $p$ for points $x$ of $\os{\circ}{T}$, 
we obtain the contravariant functoriality 
in certain cases which include the case 
where the morphism $g$ is 
the ``nonstandard'' relative Frobenius morphism of 
$X_{\os{\circ}{T}_0}/S_{\os{\circ}{T}_0}$ 
which will be explained soon later 
when $X$ is of characteristic $p>0$ and $\os{\circ}{T}$ 
is a flat $p$-adic formal scheme. 
(It seems hopeless to me that there exists the commutative diagram 
(\ref{eqn:auc}) for an unconditional $g$; 
it seems an important problem to define other filtered complexes 
satisfying a similar diagram to (\ref{eqn:auc}) for any $g$.) 
\par 
Now we can explain the Tate twist $(-j-k)$ in (\ref{eqn:caxd}) and (\ref{eqn:espsp})
and more generally in (\ref{eqn:catxd}) and (\ref{eqn:esbcpsp}) as follows. 
\par 
Assume that $\os{\circ}{S}$ is of characteristic $p>0$. 
Let $F_S\col S\lo S$ be the Frobenius endomorphism. 
Set $S^{[p]}:=S\times_{\os{\circ}{S},\os{\circ}{F}_S}\os{\circ}{S}$.  
In \cite{ollc} Ogus has already defined $S^{[p]}$ and 
he has denoted it by $S^{(1)}$. 
In [loc.~cit.] he has used $S^{[p]}$ 
to give a simpler proof of Hyodo-Kato isomorphism than the proof of it in \cite{hk}. 
(See \cite{ollc} when $\os{\circ}{S}$ is not necessarily of characteristic $p>0$.) 
The log scheme $S^{[p]}$ is different from $S^{\{p\}}:=S\times_{S,F_S}S=S$ 
and it has the universality for the log schemes 
$U$ with morphisms $S\lo U$ and $U\lo S$ such that 
the composite morphism $S\lo U\lo S$ is equal to $F_S$ 
and the morphism $U\lo S$ is solid, i.e., 
$U=S\times_{\os{\circ}{S}}\os{\circ}{U}$.  
(We follow [loc.~cit.] for the terminology ``solid'': 
it is usually said to be ``strict''.) 
(This observation is important 
for the proof of the log convergence of the weight filtration in (8).).
Let $F_{S/\os{\circ}{S}}\col S\lo S^{[p]}$ be the relative Frobenius morphism 
of $S$ over $\os{\circ}{S}$. 
Let $(T,{\cal J},\del) \lo (T',{\cal J}',\del')$ be a morphism of 
log PD-enlargements over $F_{S/\os{\circ}{S}}$. 
Then we have the ``base change morphism''  
$(F_{S/\os{\circ}{S}})_{\os{\circ}{T}{}'_0 \os{\circ}{T}_0}\col 
S_{\os{\circ}{T}_0}\lo S^{[p]}_{\os{\circ}{T}{}'_0}$ of $F_{S/\os{\circ}{S}}$. 
The underlying morphism of 
$(F_{S/\os{\circ}{S}})_{\os{\circ}{T}{}'_0 \os{\circ}{T}_0}$ is equal to 
the morphism $\os{\circ}{T}_0\lo \os{\circ}{T}{}'_0$. 
Set also $X^{[p]}:=X\times_{S}S^{[p]}
(=X\times_{\os{\circ}{S},\os{\circ}{F}_S}\os{\circ}{S})$ 
and $X^{\{p\}}:=X\times_{S,F_{S}}S$.  
Note that $X^{[p]}/S^{[p]}$ is an SNCL scheme, but $X^{\{p\}}/S$ is not.   
Then we have the following natural commutative diagram 
\begin{equation*} 
\begin{CD} 
X_{\os{\circ}{T}_0}@>{F}>> X^{[p]}_{\os{\circ}{T}{}'_0}\\
@VVV @VVV \\
S_{\os{\circ}{T}_0}@>>> S^{[p]}_{\os{\circ}{T}{}'_0}\\
@V{\bigcap}VV @VV{\bigcap}V \\
S(T)^{\nat}@>>> S^{[p]}(T')^{\nat}. 
\end{CD} 
\tag{0.0.1.16}\label{eqn:asssc} 
\end{equation*} 
\par
Consider (\ref{eqn:asssc}) in the special case $(T',{\cal J}',\del')=(T,{\cal J},\del)$ 
for the time being for realizing the benefit of the theory of log geometry in the sense of Fontaine-Illusie-Kato. 
(We shall explain another benefit in the explanation for (11) below.) 
We should note that we can take a $p$-th root of a ``prime'' section  
of $M_{S^{[p]}(T)^{\nat}}$ 
without changing the structural sheaf of rings of log schemes 
in the theory of Fontaine-Illusie-Kato. 
Here we mean by the prime section a local section 
of $M_{S^{[p]}(T)^{\nat}}$ whose image 
in $M_{S^{[p]}(T)^{\nat}}/{\cal O}_T^*$ is a local generator. 
(More generally, for a ``free'' hollow log scheme with a global chart 
``${\mab N}^m \lo {\cal O}$'', 
we can take an $n$-th root of 
a ``prime'' section of the log structure by considering 
the global chart $\dfrac{1}{n}({\mab N}^m)\lo {\cal O}$ 
without changing the structural sheaf of rings.)  
This is impossible in the classical definition of the semistableness over a one parameter family. 
We should also note that the morphism 
$F\col X_{\os{\circ}{T}_0} \lo X^{[p]}_{\os{\circ}{T}_0}$ 
is a mixture of a relative Frobenius morphism and an absolute Frobenius morphism: 
it is a relative Frobenius morphism as a morphism of schemes but 
an ``essentially'' absolute Frobenius morphism with respect to the log structure of $X_{\os{\circ}{T}_0}$. We call the morphism 
$X_{\os{\circ}{T}_0} \lo X^{[p]}_{\os{\circ}{T}_0}$ 
the {\it abrelative Frobenius morphism} of $X_{\os{\circ}{T}_0}$. 
(The adjective ``abrelative'' is a coined word; 
it means ``absolute and relative'' or ``far from being relative''.) 
This morphism has an important role in this book,  
especially in the integral theory of the weight filtration on 
the log crystalline cohomology of a proper SNCL scheme; 
it is also useful for the theory of the weight filtration modulo torsion.  
\par 
Now consider again the case where $(T',{\cal J}',\del')$ 
is not necessarily equal to $(T,{\cal J},\del)$ for the explanation 
for the Tate twist $(-j-k)$.  
Assume that we are given a morphism 
$\Phi \col \os{\circ}{F}{}^*_{\rm crys}(E')\lo E$ 
of ${\cal O}_{\os{\circ}{X}_{T_0}/\os{\circ}{T}}$-modules 
and that $\os{\circ}{T}$ 
is a flat $p$-adic formal scheme.   
Then we prove that 
there exists a morphism 
\begin{align*} 
F^*\col (A_{\rm zar}(X^{[p]}_{\os{\circ}{T}{}'_0}/S^{[p]}(T')^{\nat},E'),P) \lo 
RF_*((A_{\rm zar}(X_{\os{\circ}{T}_0}/S(T)^{\nat},E),P))
\tag{0.0.1.17}\label{cd:rfrekt}
\end{align*} 
fitting into the following commutative diagram 
\begin{equation*} 
\begin{CD} 
(A_{\rm zar}(X^{[p]}_{\os{\circ}{T}{}'_0}/S^{[p]}(T')^{\nat},E'),P) 
@>{F^*}>>
RF_*((A_{\rm zar}(X_{\os{\circ}{T}_0}/S(T)^{\nat},E),P)) \\
@A{\simeq}AA @AA{\simeq}A \\
Ru_{X^{[p]}_{\os{\circ}{T}_0}/S(T)^{\nat}*}
(\eps_{X^{[p]}_{\os{\circ}{T}{}'_0}/S^{[p]}(T')^{\nat}}^*(E'))
@>{F^*}>>Ru_{X_{\os{\circ}{T}_0}/S(T)^{\nat}*}
(\eps_{X_{\os{\circ}{T}_0}/S(T)^{\nat}}^*(E)).\\
\end{CD}
\tag{0.0.1.18}\label{cd:refrkt}
\end{equation*} 
The morphism $F\col X\lo X^{[p]}$  induces a natural morphism 
$$F_{2j+k}\col \os{\circ}{X}{}^{(2j+k)}_{T_0}\lo \os{\circ}{X}{}^{(2j+k)}_{T'_0}
\times_{\os{\circ}{T}{}'_0,F_{\os{\circ}{T}{}'_0}}\os{\circ}{T}{}'_0
=:(\os{\circ}{X}{}^{(2j+k)}_{T'_0})^{\{p\}}.$$ 
Let $a'{}^{(j)}\col \os{\circ}{X}{}^{(j)}_{T'_0}\lo \os{\circ}{X}_{T'_0}$ 
be the analogous morphism to $a^{(j)}$. 
We mean the Tate twist $(-j-k)$ for (\ref{eqn:catxd}) by 
the following commutative diagram 
\begin{equation*} 
\begin{CD} 
{\rm gr}^P_kA_{\rm zar}(X^{[p]}_{\os{\circ}{T}{}'_0}/S^{[p]}(T')^{\nat},E')
@=\bigoplus_{j\geq \max \{-k,0\}} \\
@V{F^*}VV \\
F_*({\rm gr}^P_kA_{\rm zar}(X_{\os{\circ}{T}_0}/S(T)^{\nat},E))
@=\bigoplus_{j\geq \max \{-k,0\}} 
\end{CD} 
\end{equation*} 
\begin{equation*} 
\begin{CD} 
a^{(2j+k)}_{*} 
Ru_{(\os{\circ}{X}{}^{(2j+k)}_{T'_0})^{\{p\}}/\os{\circ}{T}{}'*}
(E'_{(\os{\circ}{X}{}^{(2j+k)}_{T'_0})^{\{p\}}/\os{\circ}{T}{}'}\otimes_{\mab Z} 
\vp^{(2j+k)}_{\rm crys}(\os{\circ}{X}{}^{[p]}_{T'_0}/\os{\circ}{T}{}'))[-2j-k] \\
@VV\bigoplus_{j\geq \max \{-k,0\}}p^{j+k}F^*_{2j+k}V \\
a^{(2j+k)}_{*} 
Ru_{\os{\circ}{X}{}^{(2j+k)}_{T_0}/\os{\circ}{T}*}
(E_{\os{\circ}{X}{}^{(2j+k)}_{T_0}/\os{\circ}{T}}\otimes_{\mab Z} 
\vp^{(2j+k)}_{\rm crys}(\os{\circ}{X}_{T_0}/\os{\circ}{T}))[-2j-k]. 
\end{CD} 
\tag{0.0.1.19}\label{eqn:xnt}
\end{equation*} 
We use the Frobenius morphism (\ref{cd:rfrekt}) 
in order to prove fundamental properties of 
(\ref{eqn:esbcpsp}), e.~g., (5), (6), (8), (9), (20), (21) and (22). 
\par 
Let us go back to the situation 
before the explanation for the Tate twist $(-j-k)$. 
In the case where $\os{\circ}{T}$ is a $p$-adic formal scheme, 
we can prove that the commutative diagram 
${\rm (\ref{eqn:auc})}\otimes^L_{\mab Z}{\mab Q}$ exists, 
that is, there exists a morphism 
\begin{equation*} 
g^*
\col 
(A_{\rm zar}(X'_{\os{\circ}{T}{}'_0}/S'(T')^{\nat},E'),P)\otimes^L_{\mab Z}{\mab Q}
\lo
Rg_{*}
((A_{\rm zar}(X_{\os{\circ}{T}_0}/S(T)^{\nat},E),P)\otimes^L_{\mab Z}{\mab Q})
\tag{0.0.1.20}\label{eqn:aucds}
\end{equation*}
of filtered complexes for an unconditional $g$ fitting into 
the following commutative diagram 
\begin{equation*} 
\begin{CD}
A_{\rm zar}(X'_{\os{\circ}{T}{}'_0}/S'(T')^{\nat},E')\otimes^L_{\mab Z}{\mab Q}
@>{g^*}>> 
Rg_{*}
(A_{\rm zar}(X_{\os{\circ}{T}_0}/S(T)^{\nat},E))\otimes^L_{\mab Z}{\mab Q}
\\ 
@A{\simeq}AA 
@AA{\simeq}A \\ 
Ru_{X'_{\os{\circ}{T}{}'_0}/S'(T')^{\nat}*}
(\eps^*_{X'_{\os{\circ}{T}{}'_0}/S'(T')^{\nat}}(E'))
\otimes^L_{\mab Z}{\mab Q}
@>{g^*_{{\rm crys}}}>> Rg_{*}
Ru_{X_{\os{\circ}{T}_0}/S(T)^{\nat}*}
(\eps^*_{X_{\os{\circ}{T}_0}/S(T)^{\nat}}(E))
\otimes^L_{\mab Z}{\mab Q}  
\end{CD}
\tag{0.0.1.21}\label{eqn:aucruc} 
\end{equation*}  
and satisfying 
a transitive relation 
\begin{align*} 
``(h\circ g)^*=Rh_*(g^*)\circ h^*{}{\textrm '}{\textrm '}
\tag{0.0.1.22}\label{ali:aughc} 
\end{align*} 
and ${\rm id}^*_{X_{\os{\circ}{T}_0}}={\rm id}_{(A_{\rm zar}(X_{\os{\circ}{T}_0}/S(T)^{\nat},E),P)
\otimes^L_{\mab Z}{\mab Q}}$. 
We also prove the contravariant functorialty of 
the weight spectral sequence (\ref{eqn:esbcpsp}),  
which is useful for the compatibility of 
local and global Langlands correspondences (\cite{ytgll}). 
\par 
To construct 
$(A_{\rm zar}(X_{\os{\circ}{T}_0}/S(T)^{\nat},E),P)$ 
and to prove the contravariant functoriality of it, 
we need to construct a nice embedding system of 
$X_{\os{\circ}{T}_0}/\ol{S(T)^{\nat}}$, where 
the log scheme $\ol{S(T)^{\nat}}$ is explained in 
the next paragraph. 
We can use two methods for the construction of the embedding system. 
We explain this roughly as follows. 
\par 
Let $\star$ be nothing or $\prime$. 
Let $\ol{S^{\star}(T^{\star})^{\nat}}$ be a fine log scheme defined in 
\S\ref{sec:snclv}; if $S^{\star}(T^{\star})^{\nat}=(\os{\circ}{T}{}^{\star},
{\mab N}\oplus {\cal O}_{T^{\star}}^*\owns 
(1,u)\lom 0\in {\cal O}_{T^{\star}})$ $(u\in {\cal O}_{T^{\star}}^*)$,
then 
$\ol{S^{\star}(T^{\star})^{\nat}}
=(\ul{\rm Spec}_{\os{\circ}{T}{}^{\star}}({\cal O}_{T^{\star}}[t]),
{\mab N}\oplus {\cal O}_{T^{\star}}[t]^*\owns 
(n,u)\lom ut^n\in {\cal O}_{T^{\star}}[t])$) 
$(u\in {\cal O}_{T^{\star}}[t]^*)$.  
For simplicity of the explanation, assume that 
there exists an immersion 
$X^{\star}_{\os{\circ}{T}{}^{\star}_0}\os{\sus}{\lo} 
\ol{\cal P}{}^{\star}$ over $\ol{S^{\star}(T^{\star})^{\nat}}$ 
into a log smooth scheme over $\ol{S^{\star}(T^{\star})^{\nat}}$   
and that $E^{\star}$ is the trivial coefficient.    
Because the immersion $S'(T')^{\nat} \os{\sus}{\lo} \ol{S'(T')^{\nat}}$ 
is not (topologically) nil, 
we do not necessarily have a morphism 
$\ol{\cal P}\lo \ol{\cal P}{}'$ extending a morphism 
$X_{\os{\circ}{T}_0}\lo X'_{\os{\circ}{T}{}'_0}$ even locally on $X_{\os{\circ}{T}_0}$. 
\par 
The first method is a more or less standard method 
in the quite nonstandard situation above. 
First we consider the log PD-envelope $\ol{\mathfrak D}{}^{\star}$ of 
the immersion 
$X^{\star}_{\os{\circ}{T}{}^{\star}_0}\os{\sus}{\lo} \ol{\cal P}{}^{\star}$ 
over $(\os{\circ}{T}{}^{\star},{\cal J}^{\star},\del^{\star})$ 
($\star=$ nothing or $\prime$). 
Let ${\mathfrak D}(\ol{S^{\star}(T^{\star})^{\nat}})$ be the log PD-envelope of 
the immersion $S^{\star}(T^{\star})^{\nat}\os{\sus}{\lo} 
\ol{S^{\star}(T^{\star})^{\nat}}$ 
over $(\os{\circ}{T}{}^{\star},{\cal J}^{\star},\del^{\star})$.   
Set ${\cal P}{}^{\star}:=\ol{\cal P}{}^{\star}
\times_{\ol{S^{\star}(T^{\star})^{\nat}}}S^{\star}(T^{\star})^{\nat}$ 
and 
${\mathfrak D}^{\star}
:=\ol{\mathfrak D}{}^{\star}\times_{{\mathfrak D}(\ol{S^{\star}(T^{\star})^{\nat}})}S^{\star}(T^{\star})^{\nat}$. 
Let $\ol{\cal P}{}^{\star}{}^{\rm ex}$ be the exactification of 
the immersion $X^{\star}_{\os{\circ}{T}{}^{\star}_0}\os{\sus}{\lo} \ol{\cal P}{}^{\star}$. 
(In \cite{s3} Shiho has defined the exactification of the immersion of 
fine log $p$-adic formal schemes; we recall this in \S\ref{sec:tdiai} in the text.) 
Though we do not give the precise definition of 
$(A_{\rm zar}(X^{\star}_{\os{\circ}{T}{}^{\star}_0}/S^{\star}(T^{\star})^{\nat}),P)$ here, 
we only note here that 
$(A_{\rm zar}(X^{\star}_{\os{\circ}{T}{}^{\star}_0}/S^{\star}(T^{\star})^{\nat}),P)$ 
is constructed by the use of the de Rham complex 
${\cal O}^{\star}_{{\mathfrak D}^{\star}}\otimes_{{\cal O}_{{\cal P}^{\star}{}^{\rm ex}}}
\Om^{\bul}_{{\cal P}^{\star}{}^{\rm ex}/\os{\circ}{T}{}^{\star}}$ 
with a filtration $P$ on 
${\cal O}^{\star}_{{\mathfrak D}^{\star}}\otimes_{{\cal O}_{{\cal P}^{\star}{}^{\rm ex}}}
\Om^{\bul}_{{\cal P}^{\star}{}^{\rm ex}/\os{\circ}{T}{}^{\star}}$. 
Because the immersion $X_{\os{\circ}{T}_0}\os{\sus}{\lo} \ol{\mathfrak D}$ is (topologically) nil 
and because $\ol{\cal P}{}'{}^{\rm ex}$ is formally log smooth over $\ol{S(T)^{\nat}}$, 
there exists a morphism 
$\ol{\mathfrak D}\lo \ol{\cal P}{}'{}^{\rm ex}$ extending the composite morphism 
$X\lo X'\os{\sus}{\lo} \ol{\cal P}{}'{}^{\rm ex}$. 
This morphism induces the morphism 
$\ol{\mathfrak D}\lo \ol{\mathfrak D}{}'$ by the universality of the log PD-envelope 
and hence we obtain the morphism ${\mathfrak D}\lo {\mathfrak D}{}'$. 
Using the comparison of a quotient of 
$\Om^{\bul}_{{\mathfrak D}^{\star}/\os{\circ}{T}{}^{\star}}$ 
with 
${\cal O}_{{\mathfrak D}^{\star}}\otimes_{{\cal O}_{{\cal P}^{\star}{}^{\rm ex}}}
\Om^{\bul}_{{\cal P}^{\star}{}^{\rm ex}/\os{\circ}{T}{}^{\star}}$, 
we can prove that the morphism ${\mathfrak D}\lo {\mathfrak D}{}'$ induces 
the following filtered morphism 
\begin{align*} 
g^*\col 
({\cal O}_{{\mathfrak D}'}\otimes_{{\cal O}_{{\cal P}'{}^{\rm ex}}}
\Om^{\bul}_{{\cal P}'{}^{\rm ex}/\os{\circ}{T}{}'},P)
\lo 
g_*(({\cal O}_{{\mathfrak D}}\otimes_{{\cal O}_{{\cal P}^{\rm ex}}}
\Om^{\bul}_{{\cal P}^{\rm ex}/\os{\circ}{T}},P))
\tag{0.0.1.23}\label{ali:dope}
\end{align*}  
and  the modification of this filtered morphism 
induces the $g^*$ in (\ref{eqn:agcxds}) 
for the case of the trivial coefficients. 
(In the case where $E$ and $E'$ are nontrivial, we indeed need the morphism 
$\ol{\mathfrak D}\lo \ol{\mathfrak D}{}'$ instead of the morphism 
${\mathfrak D}\lo {\mathfrak D}{}'$.)
\par 
The second method is extremely simple and necessary in the Chapter III 
but has very useful applications 
even for already known facts; it simplifies the proofs of known facts.   
(I wonder why no one has noticed this simple argument except an argument in 
\cite{bfi}.) 
Set $X'_{\os{\circ}{T}_0}:=X'\times_{\os{\circ}{T}{}'_0}\os{\circ}{T}_0$ and 
$\ol{\cal P}{}'_{\ol{S(T)^{\nat}}}
:=\ol{\cal P}{}'\times_{\ol{S'(T')^{\nat}}}\ol{S(T)^{\nat}}$ and 
consider the fiber product 
\begin{align*} 
\ol{\cal P}{}'':=\ol{\cal P}\times_{\ol{S(T)^{\nat}}}\ol{\cal P}{}'_{\ol{S(T)^{\nat}}}. 
\end{align*} 
Then, by using a natural immersion 
$$X_{\os{\circ}{T}_0}=
X_{\os{\circ}{T}_0}\times_{X'_{\os{\circ}{T}_0}}X'_{\os{\circ}{T}_0}
\os{\sus}{\lo} 
\ol{\cal P}\times_{\ol{S(T)^{\nat}}}\ol{\cal P}{}'_{\ol{S(T)^{\nat}}}=\ol{\cal P}{}'',$$ 
we have the following commutative diagram 
\begin{equation*} 
\begin{CD} 
X_{\os{\circ}{T}_0}@>{\subset}>>\ol{\cal P}{}''\\
@V{g}VV @VVV \\
X'_{\os{\circ}{T}{}'_0}@>{\subset}>>\ol{\cal P}{}'. 
\end{CD} 
\end{equation*} 
Hence we directly have the the following filtered morphism 
\begin{align*} 
g^*\col 
(\Om^{\bul}_{{\cal P}'{}^{\rm ex}/\os{\circ}{T}{}'},P)
\lo 
g_*((\Om^{\bul}_{{\cal P}''{}^{\rm ex}/\os{\circ}{T}},P))
\tag{0.0.1.24}\label{ali:doipe}
\end{align*}  
and the modification of this filtered morphism induces 
the $g^*$ in (\ref{eqn:agcxds}) 
for the case of the trivial coefficients.
\par  
This method also enables us to prove 
the contravariant functorialty of 
Bloch-Illusie's comparison isomorphism 
between de Rham-Witt complexes 
and crystalline complexes in \cite{idw} very simply. 
Indeed, let $\kap/\kap'$ be an extension of perfect fields of characteristic $p>0$. 
Let $h\col Y\lo Y'$ be a morphism of smooth schemes  over 
${\rm Spec}(\kap)\lo {\rm Spec}(\kap')$. 
Let $Y\os{\sus}{\lo} {\cal Q}$ and $Y'\os{\sus}{\lo} {\cal Q}'$ be 
immersions into smooth schemes over $\kap$ and $\kap'$, respectively.  
Assume that $Y$ and $Y'$ are affine. 
Then we have morphisms ${\cal W}_n(Y)\lo {\cal Q}$ and 
${\cal W}_n(Y')\lo {\cal Q}'$ $(n\in {\mab Z}_{\geq 1})$ 
extending $Y\os{\sus}{\lo} {\cal Q}$ 
and $Y'\os{\sus}{\lo} {\cal Q}'$, respectively.  
Set ${\cal Q}'':={\cal Q}\times_{{\rm Spec}({\cal W}_n(\kap))}
({\cal Q}'\otimes_{{\cal W}_n(\kap')}{\cal W}_n(\kap))$. 
Then we have the following commutative diagram 
\begin{equation*} 
\begin{CD} 
{\cal W}_n(Y)@>{\subset}>>{\cal Q}{}''\\
@V{{\cal W}_n(h)}VV @VVV \\
{\cal W}_n(Y')@>{\subset}>>{\cal Q}{}'. 
\end{CD} 
\end{equation*} 
Using this commutative diagram, 
we can dispense with the lemma of Dwork-Dieudnonn\'{e}-Cartier in \cite{idw} 
for the proof of the contravariant functorialty of 
Bloch-Illusie's comparison isomorphism 
between de Rham-Witt complexes if we use the ``reverse'' 
de Rham-Witt complexes defined in \cite{ir} in a sketched way. 
(See also \cite{hk} and \cite{ndw} for the detailed argument 
for the definition of the ``reverse'' de Rham-Witt complexes.)
We can also prove the contravariant functoriality of 
Hyodo-Kato's comparison isomorphism 
between log de Rham-Witt complexes 
and log crystalline complexes in \cite{hk} very simply, 
which has not been even proved 
in literatures (except the Frobenius 
endomorphism in \cite{ndw}) as far as I know. 
Using the technique above and 
the exactification of an embedding system of 
an SNCL scheme (cf.~\cite{ctcs}) 
and a local structure theorem 
of an exact immersion proved in \cite{nh2}, 
we can give a much simpler proof of the existence of 
a nice embedding system 
than the proofs in \cite{gkwf} and \cite{ctcs}. 
Furthermore, using Tsuzuki's simplicial object 
in \cite{ctze}, we can construct the existence of 
a nice embedding system for 
a naive $N$-truncated simplicial SNCL scheme
(cf.~the explanation before the explanation of (19)$\sim$(27) below). 
(It seems hopeless to construct this nice embedding system 
in the $N$-truncated simplicial case 
by using the method in \cite{gkwf} and \cite{ctcs}.) 
\par 
Though $(A_{\rm zar}(X_{\os{\circ}{T}_0}/S(T)^{\nat},E),P)$ has 
the contravariant functoriality as explained, 
it seems to me that $(A_{\rm zar}(X_{\os{\circ}{T}_0}/S(T)^{\nat},E),P)$ 
has a undesirable property: 
the filtered complex seems to have no 
contravariant functorialty 
for a morphism $X_{T_0}\lo X'_{T'_0}$ over $T\lo T'$. 
This is the reason why we do not use the notation 
$(A_{\rm zar}(X_{T_0}/T,E),P)$ for 
$(A_{\rm zar}(X_{\os{\circ}{T}_0}/S(T)^{\nat},E),P)$ 
in spite of the isomorphism (\ref{eqn:ecetsp}).  
Indeed, let us consider the following example. 
Let 
\begin{equation*} 
\begin{CD}   
S_1@>>>S'_1\\
@VVV @VVV \\
S@>>>S' 
\end{CD} 
\tag{0.0.1.25}\label{cd:ssp}
\end{equation*}
be a commutative diagram of  
families of log points.  
Let $(T^1,{\cal J}^1,\del^1)\lo (T^{1}{}',{\cal J}^1{}',\del^1{}')$ 
be a morphism of log PD-enlargements over the morphism $S_1\lo S'_1$.  
Let $\star$ be nothing or $\prime$. 
Let $T_0^{1\star}$ be the exact closed subscheme of $T^{1\star}$ 
defined by ${\cal J}^{1\star}$. 
Set $S^{\star}_{1,{\os{\circ}{T}{}^{\star}_0}}:=
(S^{\star}_1)_{\os{\circ}{T}{}^{1\star}_0}$ and 
$X^{\star}_{S^{\star}_{1,\os{\circ}{T}{}^{1\star}_0}}
:=X^{\star}\times_{S^{\star}}S^{\star}_{1,\os{\circ}{T}{}^{\star}_0}$  
($X_{S^{\star}_{1,\os{\circ}{T}{}^{\star}_0}}/S^{\star}_{1,\os{\circ}{T}{}^{\star}_0}$ 
is a log smooth log scheme but not an SNCL scheme 
if $S^{\star}_1\not=S^{\star}\times_{\os{\circ}{S}{}^{\star}}\os{\circ}{S}{}^{\star}_1$). 
Consider the case $T^{1\star}_0=S^{\star}_1$. 
Then it seems to me that 
$(A_{\rm zar}(X^{\star}_{\os{\circ}{T}{}^{\star}_0}/S^{\star}(T^{1\star})^{\nat},E^{\star}),P)$ 
is not contravariantly functorial for a morphism 
$g\col X_{S_{1,\os{\circ}{T}{}^1_0}}\lo 
X'_{S'_{1,{\os{\circ}{T}{}^{1{}'}_0}}}$ over the morphism 
$(S_1(T^1)^{\nat},{\cal J}^1,\del^1)\lo (S'_1(T^1{}')^{\nat},{\cal J}^1{}',\del^1{}')$. 
To overcome this difficulty, 
we assume that $\os{\circ}{T}{}^1$ is a flat $p$-adic formal scheme 
over ${\rm Spf}({\mab Z}_p)$ 
and we consider a case where 
the morphism $(T^1,{\cal J}^1,\del^1)\lo (T^{1'},{\cal J}^{1'},\del^{1'})$ 
of log PD-enlargements over the morphism $S_1\lo S'_1$ 
fits into the following commutative diagram 
\begin{equation*} 
\begin{CD}   
(T^1,{\cal J}^1,\del^1)@>>>(T^{1'},{\cal J}^{1'},\del^{1'})\\
@VVV @VVV \\
(T,{\cal J},\del)@>>>(T',{\cal J}^{'},\del^{'})
\end{CD} 
\tag{0.0.1.26}\label{ali:qftpdb}
\end{equation*}
of log PD-enlargements over 
(\ref{cd:ssp}). 
Then 
we ignore torsions of 
$(A_{\rm zar}(X^{\star}_{\os{\circ}{T}{}^{\star}}/S^{\star}(T^{\star})^{\nat},E^{\star}),P)$ 
and we define another new filtered complex 
\begin{align*} 
(A_{{\rm zar},{\mab Q}}(X^{\star}_{S^{\star}_{1,\os{\circ}{T}{}^{1\star}_0}}
/S^{\star}_1(T^{1\star})^{\nat},E^{\star}),P)
\in {\rm D}^+{\rm F}(f^{-1}_T({\cal O}_T\otimes_{\mab Z}{\mab Q}))
\tag{0.0.1.27}\label{ali:qftb}
\end{align*} 
(see \S\ref{sec:pgensc} for the precise definition of this filtered complex) 
and in particular, we obtain a filtered complex 
\begin{align*} 
(A_{{\rm zar},{\mab Q}}(X^{\star}_{\os{\circ}{T}{}^{\star}_0}/S^{\star}
(T^{\star})^{\nat},E^{\star}),P):=
(A_{{\rm zar},{\mab Q}}(X^{\star}_{S^{\star}_{\os{\circ}{T}{}^{\star}_0}}/
S^{\star}(T^{\star})^{\nat},E^{\star}),P)
\in {\rm D}^+{\rm F}(f^{-1}_T({\cal O}_T\otimes_{\mab Z}{\mab Q})). 
\end{align*}
As far as I know, an explicit filtered complex such as (\ref{ali:qftb})
has not been constructed for a non SNCL scheme. 
(Strictly speaking, the construction of (\ref{ali:qftb}) depends on the 
mapping degree function obtained by the morphism $S^{\star}_1\lo S^{\star}$; 
however the isomorphism class of (\ref{ali:qftb}) does not depend on it; anyway, 
it is quite surprising to me that the filtered complex 
$(A_{{\rm zar},{\mab Q}}(X^{\star}_{S^{\star}_{1,\os{\circ}{T}{}^{1\star}_0}}
/S^{\star}_1(T^{1\star})^{\nat},E^{\star}),P)$
can be defined.)  
We have the following ``tautological'' isomorphism 
\begin{align*}
(A_{\rm zar}(X^{\star}_{\os{\circ}{T}{}^{\star}_0}/S^{\star}(T^{\star})^{\nat},E^{\star}),P)
\otimes^L_{\mab Z}{\mab Q}
\os{\sim}{\lo} 
(A_{{\rm zar},{\mab Q}}(X^{\star}_{\os{\circ}{T}{}^{\star}_0}/S^{\star}
(T^{\star})^{\nat},E^{\star}),P).  
\tag{0.0.1.28}\label{ali:qqtsb}
\end{align*}
We prove that 
$$(A_{{\rm zar},{\mab Q}}
(X^{\star}_{S^{\star}_{1,\os{\circ}{T}{}^{\star}_0}}/
S^{\star}_1(T^{1\star})^{\nat},E^{\star}),P)$$
is contravariantly functorial for morphisms 
$g\col X_{S_{1,{\os{\circ}{T}{}^1_0}}}\lo X'_{S'_{1,{\os{\circ}{T}{}^{1{}'}_0}}}$ 
and $(S_1(T^1)^{\nat},{\cal J}^1,\del^1)\lo (S'_1(T^1{}')^{\nat},{\cal J}^1{}',\del^1{}')$.  
That is, we prove that there exists a morphism 
\begin{equation*} 
g^*\col 
(A_{{\rm zar},{\mab Q}}(X'_{S'_{1,{\os{\circ}{T}{}^{1{}'}_0}}}/
S'_1(T^1{}')^{\nat},E'),P)
\lo
Rg_{*}
((A_{{\rm zar},{\mab Q}}(X_{S_{1,{\os{\circ}{T}{}^1_0}}}/S_1(T^1)^{\nat},E),P))
\tag{0.0.1.29}\label{eqn:auqqcds}
\end{equation*}
of filtered complexes for an unconditional $g$ fitting into 
the following commutative diagram 
\begin{equation*} 
\begin{CD}
(A_{{\rm zar},{\mab Q}}(X'_{S'_{1,{\os{\circ}{T}{}^{1{}'}_0}}}/
S'_1(T^1{}')^{\nat},E'),P)
@>{g^*}>> 
Rg_{*}((A_{{\rm zar},{\mab Q}}(X_{S_{1,{\os{\circ}{T}{}^1_0}}}/S_1(T^1)^{\nat},E),P))
\\ 
@A{\simeq}AA 
@AA{\simeq}A \\ 
Ru_{X'_{\os{\circ}{T}{}^{1{}'}_0}/S_1{}'(T^1{}')^{\nat}*}
(\eps^*_{X'_{\os{\circ}{T}^{1{}'}_0}/S_1{}'(T^1{}')^{\nat}}(E'))\otimes^L_{\mab Z}{\mab Q}
@>{g^*_{\rm crys}}>> Rg_*
Ru_{X_{\os{\circ}{T}{}^1_0}/S_1(T^1)^{\nat}*}
(\eps^*_{X_{\os{\circ}{T}{}^1_0}/S_1(T^1)^{\nat}}(E))
\otimes^L_{\mab Z}{\mab Q}  
\end{CD}
\tag{0.0.1.30}\label{eqn:aucuc} 
\end{equation*}  
and satisfying a natural transitive relation 
\begin{equation*} 
``(h\circ g)^*=Rh_*(g^*)\circ h^*{\textrm '}{\textrm '}
\tag{0.0.1.31}\label{eqn:hgr} 
\end{equation*}  
and ${\rm id}^*_{X_{S_{1,{\os{\circ}{T}{}_0}}}}
={\rm id}_{(A_{{\rm zar},{\mab Q}}(X_{\os{\circ}{T}_0}/S_1(T^1)^{\nat},E),P)}$. 
We also prove that there exists a comparison isomorphism 
\begin{align*}
(A_{{\rm zar},{\mab Q}}(X^{\star}_{{\os{\circ}{T}{}^{\star}_0}}
/S^{\star}(T^{1\star})^{\nat},E^{\star}),P)
& \os{\sim}{\lo} 
(A_{{\rm zar},{\mab Q}}(X^{\star}_{S^{\star}_{1,\os{\circ}{T}{}^{\star}_0}}/
S^{\star}_1(T^{1\star})^{\nat},E^{\star}),P). 
\tag{0.0.1.32}\label{ali:qqsb}\\
\end{align*}
Using the isomorphisms (\ref{ali:qqtsb}) and (\ref{ali:qqsb}), we see that 
$(A_{\rm zar}(X^{\star}_{\os{\circ}{T}{}^{\star}_0}/S^{\star}(T^{\star})^{\nat},E^{\star}),P)
\otimes^L_{\mab Z}{\mab Q}$  
is contravariantly functorial for morphisms  
$X_{S_{1,{\os{\circ}{T}{}^1_0}}}\lo X'_{S'_{1,{\os{\circ}{T}{}^{1{}'}_0}}}$ 
and $(S_1(T^1)^{\nat},{\cal J},\del)\lo (S'_1(T^{1}{}')^{\nat},{\cal J}',\del')$.  
We also see that the filtered complex 
$(A_{{\rm zar},{\mab Q}}(X_{S_{1,\os{\circ}{T}{}{}^1_0}}/S_1(T^1)^{\nat},E^{\star}),P)$ 
depends only on $X/S$ and $(S(T^1)^{\nat},{\cal J}^1,\del^1)$ 
and the morphism 
$\os{\circ}{T}{}^1_0\lo \os{\circ}{S}$.
\par 
Using the filtered complex (\ref{ali:qftb}), 
we can define a more standard Frobenius action on 
$(A_{\rm zar}(X^{\star}_{S_{\os{\circ}{T}_0}}/S(T)^{\nat},E),P) 
\otimes^L_{\mab Z}{\mab Q}$ than the induced action by (\ref{cd:rfrekt}) 
as follows in a restricted case. (It seems to me that the standard 
Frobenius action on 
$(A_{\rm zar}(X^{\star}_{S_{\os{\circ}{T}_0}}/S(T)^{\nat},E),P) 
\otimes^L_{\mab Z}{\mab Q}$ cannot be defined in general 
even in the case where $S$ is of characteristic $p>0$ 
because $X\times_{S,F_{S}}S$ is not an SNCL scheme.) 

\par 
Assume that $S$ is of characteristic $p>0$ until the explanation for (4). 
Let $(T,{\cal J},\del)\lo (T',{\cal J}',\del')$ be a morphism of 
log PD-enlargements over the morphism 
$F_{S/\os{\circ}{S}}\col S\lo S^{[p]}$.   
Consider the following composite morphism: 
$T_0\lo S_{\os{\circ}{T}_0} \lo S^{[p]}_{\os{\circ}{T}_0}\lo S_{\os{\circ}{T}_0}$. 
Set $X^{\{p\}}_{\os{\circ}{T}_0}:=
(X\times_{S,F_S}S)\times_{\os{\circ}{S}}\os{\circ}{T}_0$. 
Then we have the following two relative Frobenius morphisms 
$F\col X_{\os{\circ}{T}_0}\lo X^{\{p\}}_{\os{\circ}{T}_0}$ over 
$S_{\os{\circ}{T}_0}$ 
and 
$F\col X_{\os{\circ}{T}_0}\lo X^{[p]}_{\os{\circ}{T}_0}$ over 
$S_{\os{\circ}{T}_0}\lo S^{[p]}_{\os{\circ}{T}_0}$. 
We denote the first $F$ by $F^{\{p\}}$ 
to distinguish the first $F$ from the second $F$. 
As before, assume that we are given a morphism 
$\Phi \col \os{\circ}{F}{}^*_{\rm crys}(E')\lo E$ 
of ${\cal O}_{\os{\circ}{X}_{T_0}/\os{\circ}{T}}$-modules 
and that $\os{\circ}{T}$ 
is a flat $p$-adic formal scheme.
Assume also that we are given the following commutative diagram: 
\begin{equation*} 
\begin{CD} 
T@>>> T'@>>> T\\
@A{\bigcup}AA @A{\bigcup}AA @A{\bigcup}AA \\
T_1@>>> T'_1@>{F_{T_1}}>> T_1 \\
@VVV @VVV @VVV\\
S@=S @>{F_S}>> S. 
\end{CD} 
\end{equation*}

Because we can consider 
the morphism $F^{\{p\}}\col X_{\os{\circ}{T}_0}\lo X^{\{p\}}_{\os{\circ}{T}_0}$ 
as a morphism over $S(T)^{\nat}\lo S(T')^{\nat}$, 
we obtain the following pull-back morphism 
\begin{align*} 
F^{\{p\}*}\col &
(A_{{\rm zar}}(X^{[p]}_{\os{\circ}{T}{}'_0}/S^{[p]}(T')^{\nat},E'),P)
\otimes^L_{\mab Z}{\mab Q}
=(A_{{\rm zar},{\mab Q}}(X^{[p]}_{\os{\circ}{T}{}'_0}/S^{[p]}(T')^{\nat},E'),P)
\tag{0.0.1.33}\label{ali:frst}\\
&=(A_{{\rm zar},{\mab Q}}(X^{\{p\}}_{\os{\circ}{T}{}'_0}/S(T')^{\nat},E'),P)\lo 
RF^{\{p\}}_{*}((A_{{\rm zar},{\mab Q}}(X_{\os{\circ}{T}_0}/S(T)^{\nat},E),P))\\
&=
RF^{\{p\}}_{*}((A_{{\rm zar}}(X_{\os{\circ}{T}_0}/S(T)^{\nat},E),P))
\otimes^L_{\mab Z}{\mab Q}. 
\end{align*} 
Here we have used (\ref{ali:qqsb}) for obtaining the second equality in (\ref{ali:frst}). 
We also have the following Frobenius morphism 
\begin{align*} 
F^*\col &
(A_{{\rm zar}}(X^{[p]}_{\os{\circ}{T}{}'_0}/S^{[p]}(T')^{\nat},E'),P)
\otimes^L_{\mab Z}{\mab Q}\lo 
RF_{*}((A_{{\rm zar}}(X_{\os{\circ}{T}_0}/S(T)^{\nat},E),P))
\otimes^L_{\mab Z}{\mab Q}
\tag{0.0.1.34}\label{ali:sptp}
\end{align*} 
induced by the morphism (\ref{cd:rfrekt}). 
The fact that the the pull-backs of the two Frobenius morphisms 
(\ref{ali:frst}) and (\ref{ali:sptp}) are equal is 
a consequence of the relation (\ref{eqn:hgr}). 
\par  
(4) is an SNCL version of 
the filtered base change theorem in the open log case in \cite{nh2}. 
\par 
Let the notations be as in the explanation for (2). 
Assume that $\os{\circ}{T}$ is quasi-compact and 
that $\os{\circ}{f}_{T_0} 
\col \os{\circ}{X}_{T_0}\lo \os{\circ}{T}_0$ 
is quasi-compact and quasi-separated.  
Let $u\col (T',{\cal J}',\del') \lo (T,{\cal J},\del)$ 
be a morphism of fine log PD-schemes.  
Assume that  ${\cal J}'$ is quasi-coherent.  
Let $T'_0$ be the exact closed log subscheme of $T'$ 
defined by ${\cal J}'$. 
Let 
$f' \col X_{\os{\circ}{T}{}'_0}:=
X_{\os{\circ}{T}_0}\times_{\os{\circ}{T}_0}\os{\circ}{T}{}'_0 \lo S(T')^{\nat}$ 
be the base change morphism of $f\col X_{\os{\circ}{T}_0}\lo S(T)^{\nat}$.  
Let $p\col X_{\os{\circ}{T}{}'_0}\lo X_{\os{\circ}{T}_0}$ be the natural morphism and 
let 
$$\os{\circ}{p}_{{\rm crys}} \col 
((\os{\circ}{X}_{T'_0}/\os{\circ}{T}{}')_{\rm crys},{\cal O}_{\os{\circ}{X}_{T'_0}/\os{\circ}{T}{}'})
\lo 
((\os{\circ}{X}_{T_0}/\os{\circ}{T})_{\rm crys},{\cal O}_{\os{\circ}{X}_{T_0}/\os{\circ}{T}})$$ 
be the induced morphism by $\os{\circ}{p}$. 
Then we prove that there exists 
the following canonical filtered isomorphism
\begin{equation*}
Lu^*Rf_*((A_{\rm zar}(X_{\os{\circ}{T}_0}/S(T)^{\nat},E),P)) 
\os{\sim}{\lo} Rf'_*((A_{\rm zar}(X_{\os{\circ}{T}{}'_0}/S(T'),
\os{\circ}{p}{}^*_{{\rm crys}}(E)),P))
\tag{0.0.1.35}\label{eqn:fab}
\end{equation*}
in ${\rm DF}(f'{}^{-1}({\cal O}_{T'}))$ ((4)). 
This is the filtered base change theorem of 
$(A_{\rm zar}(X_{\os{\circ}{T}_0}/S(T)^{\nat},E),P)$, 
which is the SNCL version of the filtered base change 
theorem proved in the open log case in \cite{nh2}. 
The existence of the morphism (\ref{eqn:fab}) follows from 
the contravariant functoriality in (3);  
the fact that the morphism 
(\ref{eqn:fab}) is a filtered isomorphism follows from 
the classical base change formula in \cite{bb} and \cite{bob}. 
\par 
(5) is a log filtered version of the infinitesimal deformation invariance of 
isocrystalline cohomologies in \cite{boi} and an SNCL version of 
the log infinitesimal deformation invariance of log crystalline cohomologies 
of open log schemes in \cite{nh2}. 
The proof of the invariance in this book is much more involved 
than that of the invariance in \cite{nh2}. 
To prove the invariance in this book, 
it is indispensable to consider the framework of a log PD-enlargement 
$(T,{\cal J},\del)$ of $S$ of fine log PD (formal) schemes. 
Even if $S$ is an exact closed subscheme of a scheme defined by a 
PD-ideal sheaf, $S^{[p]}$ is {\it not} necessarily so. 
However we can consider $S^{[p]}(T)$ 
by using the composite morphism $T_0\lo S\lo S^{[p]}$. 
We explain (5) as follows. 
\par 
Assume that $S$ is of characteristic $p>0$.  
Let $S\lo S'$ be a morphism of family of log points. 
Assume that $S'$ is of characteristic $p>0$. 
Let $(T,{\cal J},\del)\lo (T',{\cal J}',\del')$ be a morphism of log PD-enlargements 
over $S\lo S'$.  
Let $\star$ be nothing or $\prime$. 
Let $T^{\star}_0$ be as in the explanation for (3). 
Assume that 
$(\os{\circ}{T}{}^{\star},{\cal J}^{\star},\del^{\star})$ 
is a $p$-adic formal PD-schemes in the sense of {\rm \cite{bob}}. 
Let $X'$ be another SNCL scheme over $S'$. 
Assume that $\os{\circ}{X}{}^{\star}$,  
$\os{\circ}{S}{}^{\star}$ and $\os{\circ}{T}{}^{\star}$ are quasi-compact. 
Let 
$\iota\col T^{\star}_0(0) \os{\subset}{\lo} T^{\star}_0$  
be an exact closed nilpotent immersion.  
Set  
$X^{\star}_{\os{\circ}{T}{}^{\star}_0}(0)
:=X^{\star}\times_{S^{\star}}S^{\star}_{\os{\circ}{T}{}^{\star}_0(0)}$. 
Let 
\begin{align*} 
g_0 \col X_{\os{\circ}{T}_0}(0) \lo X'_{\os{\circ}{T}{}'_0}(0)
\end{align*}   
be a morphism of log schemes over 
$S_{\os{\circ}{T}_0}(0)\lo S'_{\os{\circ}{T}{}'_0}(0)$.  
Then we prove that there exists a canonical filtered morphism
\begin{align*}
g^*_0&: 
(A_{{\rm zar},{\mab Q}}(X'_{\os{\circ}{T}{}'_0}/S'(T')^{\nat}),P) 
\lo 
Rg_{0*}((A_{{\rm zar},{\mab Q}}(X_{\os{\circ}{T}_0}/S(T)^{\nat}),P)). 
\tag{0.0.1.36}\label{eqn:ldeqnvn}
\end{align*}  
The $g^*_0$'s are compatible 
with compositions of $g_0$'s as in (\ref{eqn:hgr}) and 
${\rm id}^*_{X_{\os{\circ}{T}{}_0}(0)}=
{\rm id}_{(A_{{\rm zar},{\mab Q}}(X_{\os{\circ}{T}_0}/S(T)^{\nat}),P)}$. 
If $g_0$ has a lift $g \col X_{\os{\circ}{T}_0} \lo X'_{\os{\circ}{T}{}'_0}$ 
over $S_{\os{\circ}{T}}\lo S'_{\os{\circ}{T}{}'}$, 
then $g^*_0=g^*$ in (\ref{eqn:aucds}) 
via the comparison isomorphism (\ref{ali:qqtsb}) in the case of the trivial coefficients. 
As a corollary of the existence of the morphism (\ref{eqn:ldeqnvn}), if $S(T)^{\nat}=S'(T')^{\nat}$ 
and if $X_{\os{\circ}{T}_0}(0)=X'_{\os{\circ}{T}{}'_0}(0)$, 
then we obtain the following equality: 
\begin{align*} 
(A_{{\rm zar},{\mab Q}}(X_{\os{\circ}{T}_0}/S(T)^{\nat}),P) 
= (A_{{\rm zar},{\mab Q}}(X'_{\os{\circ}{T}{}'_0}/S'(T')^{\nat}),P).  
\tag{0.0.1.37}\label{ali:ldeqquvn}
\end{align*} 
This is an SNCL version of the log deformation invariance for 
open log schemes in \cite{nh2}. 
Set 
$S^{[p^n]}:=(S^{[p^{n-1}]})^{[p]}
=(S\times_{\os{\circ}{S},\os{\circ}{F}{}^n_S}\os{\circ}{S})$, 
$X^{\star[p^n]}:=X^{\star}\times_SS^{[p^n]}$ and  
$X^{\star[p^n]}_{\os{\circ}{T}_0}:=X^{\star}\times_SS^{[p^n]}_{\os{\circ}{T}_0} 
=X^{\star[p^n]}\times_{\os{\circ}{S}{}^{\star}}\os{\circ}{T}{}^{\star}_0$. 
To prove the existence of the morphism (\ref{eqn:ldeqnvn}), 
we use the filtered complex
$(A_{{\rm zar},{\mab Q}}
(X^{\star}{}^{[p^n]}_{\!\! \!\os{\circ}{T}{}^{\star}_0}/S^{\star[p^n]}(T^{\star})^{\nat}),P)$ 
$(n\in {\mab Z}_{\geq 1})$,  
the contravariant functoriality of 
$(A_{{\rm zar},{\mab Q}}
(X^{\star}{}^{[p^n]}_{\!\! \!\os{\circ}{T}{}^{\star}_0}/S^{\star[p^n]}(T^{\star})^{\nat}),P)$ 
for the morphism $X^{[p^n]}_{\os{\circ}{T}_0}\lo X'{}^{[p^n]}_{\!\! \! \os{\circ}{T}{}'_0}$  
and a fact that the iterated abrelative Frobenius morphism  
$X^{\star}{}\lo X^{\star}{}^{[p^n]}$ $(\star=$ nothing or $\prime$) 
over $S^{\star}\lo S^{\star[p^n]}$ 
induces an isomorphism
\begin{align*} 
Rf^{[p^n]}_{T^{\star}*}
(P_kA_{{\rm zar},{\mab Q}}(X^{\star}{}^{[p^n]}_{\!\! \!\os{\circ}{T}_0}/
S^{\star[p^n]}(T^{\star})^{\nat})) 
\os{\sim}{\lo} 
Rf_{T^{\star}*}
(P_kA_{{\rm zar},{\mab Q}}(X^{\star}{}_{\!\! \!\os{\circ}{T}{}^{\star}_0}/
S^{\star}(T^{\star})^{\nat})) 
\quad (k\in {\mab Z}), 
\end{align*} 
which is a log filtered version of an isomorphism in \cite{boi}. 
Here $f^{[p^n]}_{T^{\star}}\col X^{\star}{}^{[p^n]}_{\!\! \!\os{\circ}{T}_0}\lo S^{\star}(T^{\star})$ 
is the structural morphism and 
$S^{\star [p^n]}(T^{\star})$ is obtained by the composite morphism 
$T_0\lo S^{\star}\lo S^{\star [p^n]}$. 
To use the abrelative Frobnenius morphism in order to prove (\ref{ali:ldeqquvn}) 
as above is essentially Dwork's idea (see the introduction of \cite{boi}). 
\par 
We explain (6). 
\par  
Let $S$ be a formal family of log points over ${\cal V}$.  
Let $T$ be a fine log formal scheme over ${\cal V}$. 
Assume that the underlying formal scheme 
$\os{\circ}{T}$ of $T$ is a $p$-adic formal ${\cal V}$-scheme 
in the sense of \cite{of}, that is, 
a noetherian $p$-adic formal scheme which is topologically of finite type 
over ${\rm Spf}({\cal V})$.  
We call the $T$ a log $p$-adic formal ${\cal V}$-scheme. 
Let $T_1$ 
(resp.~$T_0$) 
be the exact closed subscheme of $T$ defined by $p{\cal O}_T$ 
(resp.~the ideal sheaf of local sections $t$ of ${\cal O}_T$ such that 
$t^e\in p{\cal O}_T$ for some $e\geq 1$). 
Let $T_1\lo S$ be a morphism of fine log formal schemes over ${\cal V}$.   
Using the log infinitesimal deformation invariance in (5),  
the specialization argument of Deligne-Illusie 
(\cite{isp})
and the purity of the weight in 
the case where $\kap$ is a finite field 
(\cite{dw2}, \cite{kme}, \cite{clpu}, \cite{ndw}), 
we prove that the spectral sequence 
(\ref{eqn:esbcpsp}) 
in the case $E={\cal O}_{\os{\circ}{X}_{T_1}/\os{\circ}{T}}$ 
degenerates at $E_2$ modulo ${\cal V}$-torsion 
when $\os{\circ}{X}$ is proper over $\os{\circ}{S}$. 
This is an SNCL version of the $E_2$-degeneration of 
the weight spectral sequence of 
an open log scheme in characteristic $p>0$ modulo ${\cal V}$-torsion, 
which has been proved in \cite{nh2}. 
\par  
We would like to explain (8). To explain (8), we need to explain (7).
\par 
Let $S$ be a log $p$-adic formal scheme over ${\cal V}$. 
In \cite{ollc} Ogus has defined the category of 
log $p$-adic enlargements of $S/{\cal V}$ 
and that of log enlargements of $S/{\cal V}$.  
He has assumed the solidness for the definitions 
of log $p$-adic enlargements of $S/{\cal V}$ 
and log enlargements of $S/{\cal V}$. 
A solid log $p$-adic enlargement of $S/{\cal V}$ 
(resp.~solid log enlargement of $S/{\cal V}$) is, by definition, 
a pair $(T/{\cal V},z)$, where $T/{\cal V}$ is 
a fine log $p$-adic formal scheme over ${\cal V}$ 
and $z$ is a solid morphism $T_1\lo S$ (resp.~$T_0\lo S$) 
of log schemes fitting into the following commutative diagram 
\begin{equation*}
\begin{CD} 
T_i@>{\subset}>>T \\
@V{z}VV @VVV \\
S@>>>{\rm Spf}({\cal V})
\end{CD} 
\end{equation*}
$(i=0,1)$. 
A morphism of solid log $p$-adic enlargements of $S/{\cal V}$ 
(resp.~solid log enlargements of $S/{\cal V}$) is defined in an obvious way. 
We denote these categories by 
${\rm Enl}^{\rm sld}_p(S/{\cal V})$ and ${\rm Enl}^{\rm sld}(S/{\cal V})$, respectively. 
We also define the category ${\rm Enl}_p(S/{\cal V})$ 
of log $p$-adic enlargements of $S/{\cal V}$ 
(resp.~the category ${\rm Enl}(S/{\cal V})$ 
of log enlargements of $S/{\cal V}$) without assuming the solidness.  
Let $\star$ be $p$ or nothing and let $\sq$ be ``solid'' or nothing.   
As an obvious generalization of the work in \cite{of}, 
we define a filtered (solid) log $p$-adically convergent isocrystal on 
${\rm Enl}^{\sq}_p(S/{\cal V})$ and 
a filtered (solid) log convergent isocrystal on ${\rm Enl}^{\sq}(S/{\cal V})$. 
Let ${\rm IsocF}^{\sq}_p(S/{\cal V})$ (resp.~${\rm IsocF}^{\sq}(S/{\cal V})$) 
denote the category of filtered (solid) log $p$-adically convergent isocrystals on 
${\rm Enl}^{\sq}_{\star}(S/{\cal V})$ 
(resp.~the category of filtered (solid) log convergent isocrystals on 
${\rm Enl}^{\sq}(S/{\cal V})$). 
Let $F{\textrm -}{\rm IsocF}^{\sq}_p(S/{\cal V})$ 
(resp.~$F{\textrm -}{\rm IsocF}^{\sq}(S/{\cal V})$) 
denote the category of filtered (solid) log $p$-adically convergent $F$-isocrystals on 
${\rm Enl}^{\sq}_{\star}(S/{\cal V})$ 
(resp.~the category of filtered (solid) log convergent $F$-isocrystals on 
${\rm Enl}^{\sq}(S/{\cal V})$). 
Moreover we can define 
the category $F^{\infty}{\textrm -}{\rm IsosF}^{\sq}_{\star}(S/{\cal V})$
of filtered $F^{\infty}$-isospans on ${\rm Enl}^{\sq}_{\star}(S/{\cal V})$ 
($F^{\infty}{\textrm -}{\rm IsocF}^{\sq}_{\star}(S/{\cal V})$ 
is a full subcategory of $F^{\infty}{\textrm -}{\rm IsosF}^{\sq}_{\star}(S/{\cal V})$.)   
This is a generalized filtered $F^{\infty}$-isospan version 
of an $F^{\infty}$-span of crystals in \cite{ollc}. 
We prove the following equivalences of categories: 
\begin{equation*} 
\begin{CD} 
F^{\infty}{\textrm -}{\rm IsosF}(S/{\cal V})@>{\sim}>> 
F^{\infty}{\textrm -}{\rm IsosF}_p(S/{\cal V})\\
@V{\simeq}VV @VV{\simeq}V \\
F^{\infty}{\textrm -}{\rm IsosF}^{\rm sld}(S/{\cal V})@>{\sim}>> 
F^{\infty}{\textrm -}{\rm IsosF}^{\rm sld}_p(S/{\cal V}). 
\end{CD} 
\tag{0.0.1.38}\label{cd:ssv}
\end{equation*} 
Here all the arrows in this commutative diagrams 
are restriction functors. 
(In fact, we obtain the commutative diagram (\ref{cd:ssv}) 
with four equivalences 
for a more general fine log $p$-adic formal scheme over ${\cal V}$.) 
We also prove the following equivalences of categories: 
\begin{equation*} 
\begin{CD} 
F{\textrm -}{\rm IsocF}(S/{\cal V})@>{\sim}>> 
F{\textrm -}{\rm IsocF}_p(S/{\cal V})\\
@V{\simeq}VV @VV{\simeq}V \\
F{\textrm -}{\rm IsocF}^{\rm sld}(S/{\cal V})@>{\sim}>> 
F{\textrm -}{\rm IsocF}^{\rm sld}_p(S/{\cal V}). 
\end{CD} 
\tag{0.0.1.39}\label{cd:ssccv}
\end{equation*} 
Assume that $\os{\circ}{S}$ is of characteristic $p>0$. 
Let $F^*\col {\rm IsocF}^{\rm sld}_{\star}(S/{\cal V})
\lo {\rm IsocF}^{\rm sld}_{\star}(S/{\cal V})$ 
be the pull-back morphism. 
Then we can prove that 
\begin{align*} 
F^*((E,P))_{T}=(E,P)_{S^{[p]}(T)}.
\tag{0.0.1.40}\label{ali:rstnszsfb} 
\end{align*} 
This formula (\ref{ali:rstnszsfb}) is very useful 
for the explicit calculation of the pull-back of a filtered isocrystal  
by $F$ in which we are interested below. (See the explanation for (8) below.) 
\par 
Now we can explain (8). 
\par 
Let $S/{\cal V}$ be as in the explanation for (7). 
Assume that $S$ is a $p$-adic formal family of log points. 
(8) is an SNCL version of 
the log convergence of the weight filtration on 
the log isocrystalline cohomology sheaf of an open log scheme 
in characteristic $p>0$ in \cite{nh2}. 
However the formulation and the proof for (8) are different 
from those in [loc.~cit.] because 
the base change of $X^{\{p\}}=X\times_{S,F_S}S$ 
is not an SNCL scheme when $S$ is of characteristic $p>0$. 
To overcome this difficulty, we consider the category 
$F^{\infty}{\textrm -}{\rm IsocF}^{\rm sld}_{\star}(S/{\cal V})$ 
in the explanation of (7). 
(8) is one of the most important result in this book. 
To prove the log convergence of the weight filtration, 
we proceed as follows as in \cite{of}. 
\par 
First we prove the $p$-adic convergence of the weight filtration. 
That is, we prove that there exists a unique object
\begin{align*} 
(R^qf_{*}(P_kA_{{\rm zar},{\mab Q}}(X/K))^{\sq},P)
\tag{0.0.1.41}\label{ali:pciee}
\end{align*} 
of ${\rm IsocF}^{\sq}_p(S/{\cal V})$ 
for integers $k$ and $q$  such that  
\begin{align*} 
(R^qf_*(P_kA_{{\rm zar},{\mab Q}}(X/K))^{\sq},P)_{T}
=(R^qf_{\os{\circ}{T}*}
(P_kA_{{\rm zar},{\mab Q}}(X_{\os{\circ}{T}_1}/S(T)^{\nat})),P)
\tag{0.0.1.42}\label{ali:pkaiskf} 
\end{align*} 
for any object $T$ of ${\rm Enl}^{\sq}_p(S/{\cal V})$, 
where $f_{\os{\circ}{T}}\col X_{\os{\circ}{T}}\lo S(T)^{\nat}$ 
is the structural morphism.  
In particular, there exists a unique object
\begin{align*} 
(R^qf_{*}({\cal O}_{X/K})^{\nat,\sq},P)
\tag{0.0.1.43}\label{ali:pcniee}
\end{align*}  
of ${\rm IsocF}^{\sq}_p(S/{\cal V})$ 
such that  
\begin{align*} 
(R^qf_{*}({\cal O}_{X/K})^{\nat,\sq},P)_T=
(R^qf_{X_{\os{\circ}{T}_1}/S(T)^{\nat}*}
({\cal O}_{X_{\os{\circ}{T}_1}/S(T)^{\nat}})_{\mab Q},P)
\tag{0.0.1.44}\label{ali:ainpce}
\end{align*} 
for any object $T$ of ${\rm Enl}^{\sq}_p(S/{\cal V})$. 
\par 
Next we prove that there exists an object
\begin{align*} 
((R^q f_{*}(P_kA_{{\rm zar},{\mab Q}}(X/K))^{\rm sld},P),\Phi)
\tag{0.0.1.45}\label{ali:pcinee} 
\end{align*}  
of $F{\textrm -}{\rm IsocF}^{\rm sld}(S/{\cal V})$ 
which is an extension of (\ref{ali:pciee}) in the case $\sq={\rm sld}$. 
In particular, there exists an object
\begin{align*} 
((R^qf_{*}({\cal O}_{X/K})^{\nat,{\rm sld}},P),\Phi)
\tag{0.0.1.46}\label{ali:pcxnkee} 
\end{align*}  
of $F{\textrm -}{\rm IsocF}^{\rm sld}(S/{\cal V})$ 
which is an extension of (\ref{ali:ainpce}) in the case $\sq={\rm sld}$. 
By using the equivalence 
$F{\textrm -}{\rm IsocF}(S/{\cal V})\simeq 
F{\textrm -}{\rm IsocF}^{\rm sld}(S/{\cal V})$ in (\ref{cd:ssccv}), 
we obtain a corresponding object 
$((R^qf_{*}({\cal O}_{X/K})^{\nat},P),\Phi)$ to 
$((R^qf_{*}({\cal O}_{X/K})^{\nat,{\rm sld}},P),\Phi)$    
in $F{\textrm -}{\rm IsocF}(S/{\cal V})$. 
\par 
Moreover, locally on $S$, $P_kR^qf_*({\cal O}_{X/K})$ 
$(k\in {\mab Z})$ also defines 
a convergent isocrystal on ${\rm Enl}(\os{\circ}{S}/{\cal V})$.  
This convergence and Ogus' result  
on convergent isocrystals (\cite{of}) enable us to reduce 
to showing properties of $P$ on 
log isocrystalline cohomology sheaves   
to showing those in the case 
where $\os{\circ}{T}$ is the formal spectrum of 
a complete discrete valuation ring of mixed characteristics 
with perfect residue field.  
We hope that our relative theory of weight filtrations has 
applications for the study of the universal families of 
various modular schemes along boundaries. 
\par 
By using the (log) convergence of the weight filtration 
on log isocrystalline cohomology sheaves and by the argument 
in the previous paragraph, 
we prove the strict compatibility of the weight filtration 
on log isocrystalline cohomology sheaves with respect to the pull-back of 
a morphism of proper SNCL schemes ((9)). 
More generally, we prove the strict compatibility of the weight filtration 
on log isocrystalline cohomology sheaves of 
proper truncated simplicial base changes of SNCL schemes 
with respect to the pull-back of the morphisms of them. 
This strict compatibility has an application for 
the well-definedness of the limit of the weight filtration 
on the infinitesimal cohomology in this book ((21)). 
\par 
\par 
Next we explain (10). 
\par 
Let the notations be as in the explanation for (1). 
\par  
Let $Y$ be a log smooth scheme over $T_0$. 
Let 
\begin{equation*} 
\eps_{Y/\os{\circ}{T}}\col ((Y/\os{\circ}{T})_{\rm crys},{\cal O}_{Y/\os{\circ}{T}})
\lo 
((\os{\circ}{Y}/\os{\circ}{T})_{\rm crys},{\cal O}_{\os{\circ}{Y}/\os{\circ}{T}})
\end{equation*} 
and 
\begin{equation*} 
\eps_{Y/S(T)^{\nat}/\os{\circ}{T}}\col ((Y/S(T)^{\nat})_{\rm crys},{\cal O}_{Y/S(T)^{\nat}})
\lo 
((Y/\os{\circ}{T})_{\rm crys},{\cal O}_{Y/\os{\circ}{T}})
\end{equation*} 
be the morphisms of ringed topoi 
forgetting only the log structure of $Y$ 
and also that of $S(T)^{\nat}$, 
respectively. 
Let $F$ be a flat quasi-coherent crystal of 
${\cal O}_{Y/\os{\circ}{T}}$-modules.  
Then we can define a monodromy operator 
\begin{equation*}
N_{\rm zar} \col 
Ru_{Y/S(T)^{\nat}*}(\eps^*_{Y/S(T)^{\nat}/\os{\circ}{T}}(F)) \lo 
Ru_{Y/S(T)^{\nat}*}(\eps^*_{Y/S(T)^{\nat}/\os{\circ}{T}}(F))(-1) 
\tag{0.0.1.47}\label{eqn:lgm}
\end{equation*}
as in \cite{hk} (see \S\ref{sec:crcks} below). 
We can also define a morphism 
\begin{equation*}
\nu_{\rm zar} \col 
(A_{\rm zar}(X_{\os{\circ}{T}_0}/S(T)^{\nat},E),P) 
\lo 
(A_{\rm zar}(X_{\os{\circ}{T}_0}/S(T)^{\nat},E),P\langle -2 \rangle)(-1) 
\tag{0.0.1.48}\label{eqn:lnu}
\end{equation*} 
(cf.~\cite{st1}, \cite{msemi}).  
Here $\langle \quad \rangle$ 
means the shift of the filtration: 
$P\langle l \rangle_k:=P_{k+l}$ $(k,l\in {\mab Z})$. 
We call $\nu_{\rm zar}$ 
the quasi-monodromy operator of $E$.  
When $Y=X_{\os{\circ}{T}_0}$ and 
$F=\eps^*_{Y/\os{\circ}{T}}(E)$, 
then we prove that 
$N_{\rm zar}=\nu_{\rm zar}$ via 
the canonical isomorphism (\ref{eqn:catea}) 
as in \cite{ndw} (cf.~\cite{msemi}). 
If $X$ is quasi-compact, then the morphism 
\begin{equation*}
\nu_{\rm zar} \col 
A_{\rm zar}(X_{\os{\circ}{T}_0}/S(T)^{\nat},E) 
\lo 
A_{\rm zar}(X_{\os{\circ}{T}_0}/S(T)^{\nat},E)(-1) 
\tag{0.0.1.49}\label{eqn:axe}
\end{equation*} 
is nilpotent. Consequently the morphism 
(\ref{eqn:lgm}) is nilpotent when 
$Y=X_{T_0}$ and 
$F=\eps_{Y/\os{\circ}{T}}^*(E)$. 
This argument is due to Illusie as mentioned in 
the introduction of \cite{hk}, which has been carried out 
in \cite{msemi} (with the correction in \cite{ndw}) 
in the framework of the theory of 
log de Rham-Witt complexes for the trivial coefficient. 
Once one constructs 
$(A_{\rm zar}(X_{\os{\circ}{T}_0}/S(T)^{\nat},E),P)$ 
and obtains an equality 
$N_{\rm zar}=\nu_{\rm zar}$, 
it is natural to give the following conjecture as in 
\cite{dh1}, \cite{st1}, \cite{rz}, \cite{msemi} and  
\cite{nlpi} for the case 
$E={\cal O}_{\os{\circ}{X}_{T_0}/\os{\circ}{T}}$: 

\begin{conj}[{\bf Variational $p$-adic monodromy-weight  conjecture}]
\label{conj:remc}  
Let the notations be as in the explanation for $(6)\sim (9)$. 
Assume that $\os{\circ}{X}_{\os{\circ}{T}_1} \lo \os{\circ}{T}_1$ is projective. 
Let $k$ and $q$ be nonnegative integers. 
Then the induced morphism 
\begin{equation*} 
(N_{\rm zar})^k \col 
{\rm gr}^P_{q+k}R^qf_{X_{\os{\circ}{T}_1}/S(T)^{\nat}*}
({\cal O}_{X_{\os{\circ}{T}_1}/S(T)^{\nat}}) 
\lo {\rm gr}^P_{q-k}
R^qf_{X_{\os{\circ}{T}_1}/S(T)^{\nat}*}({\cal O}_{X_{\os{\circ}{T}_1}/S(T)^{\nat}})(-k)  
\tag{0.1.1.1}\label{eqn:pgme} 
\end{equation*}
by the monodromy operator is an isomorphism modulo torsion.  
\end{conj}

\parno
Mokrane is the first who has given the $p$-adic monodromy-weight conjecture 
for the log special fiber of a projective semistable family 
over ${\cal V}$ by using log de Rham-Witt complexes 
in the case where 
$\os{\circ}{T}={\rm Spf}({\cal W}(\kap))$ (\cite{msemi}). 
The conjecture (\ref{conj:remc}) is equivalent to the coincidence 
of the monodromy filtration and the weight filtration on 
$R^qf_{X_{\os{\circ}{T}_1}/S(T)^{\nat}*}
({\cal O}_{X_{\os{\circ}{T}_1}/S(T)^{\nat}})\otimes_{\mab Z}{\mab Q}$. 
The conjecture (\ref{conj:remc}) 
is a $p$-adic generalized version of 
Kato's conjecture in the $l$-adic case (\cite{kln}, \cite{nlpi}) 
(Kato's conjecture does not assume that 
the projective SNCL scheme is the special fiber of 
a strict semistable family over ${\cal V}$.)  
To prove (\ref{conj:remc}), 
we show that it suffices to prove (\ref{conj:remc}) 
in the case where $\os{\circ}{T}$ is the formal spectrum of 
the Witt ring of a finite field. 
This is a variational $p$-adic version of 
Nakayama's result (\cite{nd}) which says that  
the coincidence of the $l$-adic monodromy filtration 
and the $l$-adic weight filtration 
on the $l$-adic log \'{e}tale cohomology of 
a projective SNCL variety over a general field  
is reduced to the coincidence in  
the case over a finite field (\cite{nd}).  
The result (8), Ogus' theory of convergent $F$-isocrystals  
in \cite{of} and 
the specialization argument of Deligne and Illusie 
(\cite{isp}) play key roles in our reduction. 
In \cite{kaji} Kajiwara has proved the $l$-adic analogue of 
(\ref{conj:remc}) for $q=1$ and $q=2d-1$ for the case of the log point. 
In this paper we prove (\ref{conj:remc}) for $q=1$ and $q=2d-1$ by 
using his result; we also prove (\ref{conj:remc}) for the case $d=2$ 
by using an argument in \cite{rz} and \cite{msemi}. 
\par 
Let the assumption be as in (\ref{conj:remc}). 
Assume that the relative dimension of 
$\os{\circ}{X}_{\os{\circ}{T}_1} \lo \os{\circ}{T}_1$ is of pure dimension $d$. 
Let $L$ be a relatively ample line bundle on 
$\os{\circ}{X}_{T_1}/\os{\circ}{T}_1$. 
As in \cite[\S3]{boi},   
we obtain the log cohomology class 
$\eta=c_{1,{\rm crys},T}(L)$ of $L$
in $R^2f_{X_{\os{\circ}{T}_1}/S(T)^{\nat}*}
({\cal O}_{X_{\os{\circ}{T}_1}/S(T)^{\nat}})(1)$. 
The variational filtered log $p$-adic hard Lefschetz conjecture
is the following: 

\begin{conj}[{\bf Variational filtered 
log $p$-adic hard Lefschetz conjecture}]\label{conj:lhilc}
$(1)$ The following cup product 
\begin{equation*} 
\eta^i \col 
R^{d-i}f_{X_{\os{\circ}{T}_1}/S(T)^{\nat}*}
({\cal O}_{X_{\os{\circ}{T}_1}/S(T)^{\nat}})\otimes_{\mab Z}{\mab Q} 
\lo (R^{d+i}f_{X_{\os{\circ}{T}_1}/
S(T)^{\nat}*}({\cal O}_{X_{\os{\circ}{T}_1}/S(T)^{\nat}})
\otimes_{\mab Z}{\mab Q})(i)
\tag{0.1.2.1}\label{eqn:fcvilpl} 
\end{equation*}
is an isomorphism. 
\par 
$(2)$ 
In fact, $\eta^i_T$ is the following isomorphism of filtered sheaves: 
\begin{equation*} 
\eta^i_T \col 
(R^{d-i}f_{X_{\os{\circ}{T}_1}/S(T)^{\nat}*}
({\cal O}_{X_{\os{\circ}{T}_1}/S(T)^{\nat}})\otimes_{\mab Z}{\mab Q},P) 
\os{\sim}{\lo} ((R^{d+i}f_{X_{\os{\circ}{T}_1}/S(T)^{\nat}*}
({\cal O}_{X_{\os{\circ}{T}_1}/S(T)^{\nat}})
\otimes_{\mab Z}{\mab Q})(i),P). 
\tag{0.1.2.2}\label{eqn:filfilpl} 
\end{equation*}
Here $(i)$ means that 
$$P_k(R^{d+i}f_{X_{\os{\circ}{T}_1}/S(T)^{\nat}*}
({\cal O}_{X_{\os{\circ}{T}_1}/S(T)^{\nat}})\otimes_{\mab Z}{\mab Q})(i))
=P_{k+2i}R^{d+i}f_{_{\os{\circ}{T}_1}/S(T)^{\nat}*}
({\cal O}_{X_{\os{\circ}{T}_1}/S(T)^{\nat}})\otimes_{\mab Z}{\mab Q}.$$ 
\end{conj}

\parno  
The conjecture (\ref{conj:lhilc}) (1) 
is a variational version of the log $p$-adic hard Lefschetz conjecture 
(=the $p$-adic version of the $l$-adic log hard Lefschetz conjecture 
due to K.~Kato) (cf.~\cite{nlpi}). 
We add the conjecture (\ref{conj:lhilc}) (2) to this conjecture. 
To prove (\ref{conj:lhilc}) (1), we show 
that it suffices to prove (\ref{conj:lhilc}) (1) in the case where $\os{\circ}{S}$ is 
the formal spectrum of the Witt ring of a finite field 
as in the case of (\ref{conj:remc}). 
In the near future we would like to discuss (\ref{conj:lhilc}) (2) in another paper. 
\par 
We can also give a formulation of 
the variational filtered log $p$-adic weak Lefschetz conjecture 
as in the trivial log case (cf.~\cite{ny}). 
In the near future we would also like 
to prove this conjecture in another paper. 
\par 
We explain (11). 
\par 
Let $D$ be the spectrum of a complete discrete valuation ring of 
equal characteristic $p>0$ with canonical log structure. 
Let $s$ be the log point of the residue field of the closed point of $D$. 
Assume that $\Gam(s,{\cal O}_s)$ is perfect and 
let ${\cal W}$ be the Witt ring of $\Gam(s,{\cal O}_s)$. 
Let $l\not=p$ be a prime number. 
Let ${\cal X}$ be a proper strict semistable family over $D$. 
Let ${\cal X}_{\ol{\eta}}$ be the geometric generic fiber of ${\cal X}/D$. 
Let ${\cal X}_s$ be the log special fiber of ${\cal X}/D$.  
In \cite{ter} and \cite{itp} Terasoma and Ito 
have proved that the monodromy filtration and 
the weight filtration on the \'{e}tale cohomology 
$H^q_{\rm \acute{e}t}({\cal X}_{\ol{\eta}},{\mab Q}_l)$ 
$(q\in {\mab N})$ 
coincide by using Deligne's fundamental result in \cite{dw2}. 
\par 
On the other hand, in \cite{fao} Faltings has proved that 
the monodromy filtration and the weight filtration 
on the log crystalline cohomology 
$H^q_{\rm crys}({\cal X}_s/{\cal W}(s))\otimes_{\mab Z}{\mab Q}$ 
coincide in the case where ${\rm Spec}(D)$ is 
the localization at a rational point of a smooth curve over $\os{\circ}{s}$ 
when $\os{\circ}{s}$ is the spectrum of a finite field. 
In \cite{ctcs} Chiarellotto and Tsuzuki have also proved this. 
They have used Crew's fundamental result (\cite{cr}) 
which is a $p$-adic analogue of 
Deligne's fundamental result in the previous paragraph. 
Recently C.~Lazda and A.~P\'{a}l have proved the coincidence of 
the monodromy filtration and the weight filtration 
on the log crystalline cohomology 
$H^q_{\rm crys}({\cal X}_s/{\cal W}(s))\otimes_{\mab Z}{\mab Q}$  
in \cite{lp} by using their theory of rigid cohomologies over 
fields of Laurent series, their result:  
``Hyodo-Kato isomorphism in equi-characteristic $p>0$'',  
the $p$-adic local monodromy theorem, 
Marmora's functor (\cite{marmo}) 
and Crew's result above (and more) 
when $\os{\circ}{s}$ is the spectrum of a finite field. 
\par 
Let $S$ be as in the explanation for $(7)\sim (11)$. 
Assume that $\os{\circ}{S}$ is connected. 
For a closed point $\os{\circ}{s}$ of $\os{\circ}{S}_1$, 
let $s$ be the exact closed log subscheme of $S_1$ 
whose underlying scheme is $\os{\circ}{s}$.  
We call $s$ an exact closed point of $S_1$. 
Set $X_s:=X\times_{S_1}s$. 
Assume that there exists an exact closed point $s$ of $S_1$
such that the fiber $X_s$ is the log special fiber of 
a proper strict semistable family ${\cal X}$ over $D$. 
Then we prove that the morphism (\ref{eqn:pgme}) 
is an isomorphism modulo torsion for $X/S_1/S$ 
and for $X_s/s/{\cal W}(s)$ 
for any exact closed point $s$ of $S_1$.  
In particular, the monodromy filtration and the weight filtration 
on the log crystalline cohomology 
$H^q_{\rm crys}({\cal X}_s/{\cal W}(s))$ coincide in the case $S=D$. 
(We do not need to assume that $\Gam(s,{\cal O}_s)$ is a finite field.)
This is a $p$-adic version of Terasoma-Ito's result. 
Thus our general result is a generalization of 
Lazda and P\'{a}l's result. 
\par 
Using our result above, we know the following. 
If one wants to prove that 
the morphism (\ref{eqn:pgme}) is an isomorphism modulo torsion 
for the log special fiber of a proper (not necessarily projective) 
strict semistable family over ${\cal V}$,  
one has only to deform this log special fiber along $S$ 
and to find a fiber of an exact closed point of $S$ 
which lifts to a proper strict semistable family 
over a complete discrete valuation ring of equal characteristic $p>0$ 
with perfect residue field. 
(However I do not know whether this deformation is always possible.)
In particular, 
we can reduce the problem that 
the morphism (\ref{eqn:pgme}) is an isomorphism modulo torsion for 
a projective strict semistable family ${\cal V}$ 
to the problem for 
a proper (not necessarily projective) strict semistable family over $D$ 
when these two families exist at the same time. 
This is a $p$-adic version of K.~Fujiwara's try to 
reduce the $l$-adic monodromy-weight conjecture in mixed characteristics 
to the $l$-adic monodromy-weight conjecture in equal characteristic 
(\cite{fuhs}). 
As pointed out in the $l$-adic case in [loc.~cit.], 
this is a benefit of log geometry in the sense of Fontaine-Illusie-Kato. 
\par 
We explain (12).  
\par 
In \cite{nlpi} we have proved that (\ref{conj:lhilc}) (1) holds in the case 
where $X$ is the log special fiber of a projective strict semistable family 
over a complete discrete valuation ring ${\cal V}$
of mixed characteristics with residue field $\kap$ by 
using Hyodo-Kato's isomorphism and 
the hard Lefschetz theorem in characteristic $0$.  
In this book we prove that (\ref{conj:lhilc}) (1) holds in the case 
where $X$ is the log special fiber of a projective strict semistable family 
over a complete discrete valuation ring of equal characteristic with residue field $\kap$ 
by using Faltings' result for the flatness of 
the log convergent cohomology sheaf of a proper strict semistable family 
over a smooth curve over $\kap$ in \cite{fao}. 
(Chiarellotto and Tsuzuki have also proved this flatness in \cite{ctcs}.) 
\par 
Let $X/S$ be as in the explanation for $(11)$. 
Assume 
that there exists an exact closed point $s$ of $S_1$
such that the fiber $X_s$ is the log special fiber of 
a proper strict semistable family over 
a complete discrete valuation ring of mixed characteristics 
or equal characteristic with perfect residue field. 
Then we prove that (\ref{conj:lhilc}) (1) holds by using our result 
in the previous paragraph.  
\par 
Though we prove nothing about (\ref{conj:lhilc}) (2) in this book, 
we would like to prove (\ref{conj:lhilc}) (2) for 
a proper strict semistable family over a complete discrete valuation ring  
after developing theory of the compatibility of the weight filtrations  
on the log crystalline cohomology sheaves of proper SNCL schemes
with multiplicative structures of the log crystalline cohomology sheaves of them 
in another paper.
\par 
We explain (13).
\par 
We first explain (13) for the simplest case.   
Let $s$ be as in the beginning of this introduction. 
Let $X$ be a projective SNCL scheme over $s$. 
Assume that $X$ has a projective SNCL lift ${\cal X}$ over ${\cal W}(s)$.  
Then we prove that the monodromy-weight conjecture 
for $H^q_{\rm crys}(X/{\cal W}(s))\otimes_{\mab Z}{\mab Q}$ holds. 
To prove this, we use M.~Saito's result which implies 
that the monodromy filtration and 
the weight filtration on the singular cohomology of the analytification of 
a projective SNCL scheme over the log point 
$({\rm Spec}({\mab C}),({\mab N}\oplus {\mab C}^*\lo {\mab C}))$ coincide 
(\cite{sam}). 
\par 
Let $X/S_1$ be as in the explanation for (11). 
Assume that there exists an exact closed point $s$ of $S_1$
such that the fiber $X_s$ is the special fiber of a projective SNCL family over 
a complete discrete valuation ring of mixed characteristics. 
Then we prove that (\ref{conj:remc}) and (\ref{conj:lhilc}) (1) 
hold by using our result in the previous paragraph.  
We also prove the $l$-adic analogues of these results. 
\par 
We explain (14). 
\par 
Let $t$ be a fine log scheme whose underlying scheme 
is the spectrum of a perfect field $\kap'$ 
of characteristic $p>0$. 
Let $n$ be a positive integer. 
Let ${\cal W}_n(t)$ be the canonical lift of $t$ over ${\cal W}_n(\kap')$. 
Let $Y$ be a log smooth log scheme of Cartier type over $t$. 
Let $E$ be a flat coherent log crystal of 
${\cal O}_{Y/{\cal W}_n(t)}$-modules.  
Let $f \col Y \lo t \os{\sus}{\lo} {\cal W}_n(t)$ be 
the structural morphism.  
Set $E_n:=E_{{\cal W}_n(Y)}$. 
Let ${\cal W}_n\Om^{\bul}_Y$ be the log de Rham-Witt complex 
defined and denoted by ${\cal W}_n\om^{\bul}_Y$ in \cite{hk} 
(this is a log version of the de Rham-Witt complex in \cite{ir}) 
and let $({\cal W}_n\Om^{\bul}_Y)'$ be the log de Rham-Witt complex 
which is a nontrivial correction of the log de Rham-Witt complex defined in \cite{hk} 
(the correction has been done in \cite{ndw};  
${\cal W}_n\Om^{\bul}_Y$ and $({\cal W}_n\Om^{\bul}_Y)'$ 
are the reverse log de Rham-Witt complex and 
the obverse log de Rham-Witt complex of $Y/{\cal W}_n(t)$ 
in the sense of \cite{ndw}, which were denoted by 
${\cal W}_n\Lam^{\bul}_Y$ and $({\cal W}_n\Lam^{\bul}_Y)''$ in [loc.~cit.], respectively. 
Then we obtain two isomorphic log de Rham-Witt complexes 
$E_n\otimes_{{\cal W}_n({\cal O}_Y)}{\cal W}_n\Om^{\bul}_Y$ 
and $E_n\otimes_{{\cal W}_n({\cal O}_Y)'}({\cal W}_n\Om^{\bul}_Y)'$. 
We prove that 
there exists a canonical isomorphism 
\begin{equation*} 
Ru_{Y/{\cal W}_n(t)*}(E) \os{\sim}{\lo} 
E_n\otimes_{{\cal W}_n({\cal O}_Y)^{\star}}
({\cal W}_n\Om^{\bul}_Y)^{\star} 
\tag{0.1.2.3}\label{eqn:ywintny}
\end{equation*} 
in ${\rm D}^+(f^{-1}({\cal W}_n(\kap'))$ for 
$\star=$nothing or $\prime$. 
This is the log version of Etesse's comparison theorem (\cite{et}). 
Though he has not proved the contravariant functoriality 
of the isomorphism (\ref{eqn:ywintny}), 
we prove the contravariant functoriality of the isomorphism. 
To prove the contravariant functoriality of the isomorphism 
in the case of the trivial log structure and 
the trivial structure sheaf, Illusie has used 
lifts of Frobenii and the lemma 
of Dwork-Dieudnonn\'{e}-Cartier in \cite{idw}. 
However, by using the very simple argument explained in (3), 
we dispense with lifts of Frobenii and (the log version of) 
the lemma of Dwork-Dieudnonn\'{e}-Cartier in \cite{ndw} 
to prove the contravariant functoriality (at least) 
in the case $\star=$nothing. 
In the case $\star=\prime$, we obtain the functoriality 
by the canonical isomorphism 
$({\cal W}_n\Om^{\bul}_Y)'\os{\sim}{\lo} {\cal W}_n\Om^{\bul}_Y$ 
proved in \cite{hk} and \cite{ndw}. 
As a result, we do not need (the log version of) the lemma of 
Dwork-Dieudnonn\'{e}-Cartier in \cite{idw} and \cite{ndw}.   
When $E$ is a unit root log $F$-crystal on $Y/{\cal W}(t)$, 
then we construct the morphism 
$$V\col E_n\otimes_{{\cal W}_n({\cal O}_Y)^{\star}}
({\cal W}_n{\Om}^{\bul}_Y)^{\star} 
\lo 
E_{n+1}
\otimes_{{\cal W}_{n+1}({\cal O}_Y)^{\star}}
({\cal W}_{n+1}\Om^{\bul}_Y)^{\star}$$ 
following the method of Etesse. 
We also have operators $F$, ${\bf p}$ and $R$ on  
$\{E_n\otimes_{{\cal W}_n({\cal O}_Y)^{\star}}
({\cal W}_n{\Om}^{\bul}_Y)^{\star}\}_{n=1}^{\inf}$ 
and we have the usual relations between 
$F$, $V$, ${\bf p}$ and $R$ on 
$\{E_n\otimes_{{\cal W}_n({\cal O}_Y)^{\star}}
({\cal W}_n
\Om^{\bul}_Y)^{\star}\}_{n\in {\mab Z}_{\geq 1}}$. 
\par 
We explain (15). 
\par 
Let the notations be as in the explanation for (1). 
In \S\ref{sec:crcks} below we construct a certain filtered complex 
in ${\rm D}^+{\rm F}(f_T^{-1}({\cal O}_T))$ by using a method in \cite{hk}.  
We denote it by  
\begin{align*} 
(\wt{R}u_{X_{\os{\circ}{T}_0}/\os{\circ}{T}*}(\eps^*_{X_{\os{\circ}{T}_0}/\os{\circ}{T}}(E)),P). 
\end{align*}
This filtered complex has the following property: 
\begin{align*} 
{\rm gr}^P_k\wt{R}u_{X_{\os{\circ}{T}_0}/\os{\circ}{T}*}
(\eps^*_{X_{\os{\circ}{T}_0}/\os{\circ}{T}}(E))
= 
a^{(k-1)}_{T_0*}
Ru_{\os{\circ}{X}{}^{(k-1)}_{T_0}/\os{\circ}{T}*}
(E_{\os{\circ}{X}{}^{(k-1)}_{T_0}/\os{\circ}{T}}
\otimes_{\mab Z}
\vp^{(k-1)}_{\rm crys}(\os{\circ}{X}_{T_0}/\os{\circ}{T}))  
\tag{0.1.2.4}\label{ali:trgrte}
\end{align*}
for $k\in {\mab Z}_{\geq 1}$ and 
$P_{-1}\wt{R}u_{X_{\os{\circ}{T}_0}/\os{\circ}{T}*}
(\eps^*_{X_{\os{\circ}{T}_0}/\os{\circ}{T}}(E))=0$. 
See (\ref{ali:p0dept}) below 
for the description of 
$P_0\wt{R}u_{X_{\os{\circ}{T}_0}/\os{\circ}{T}*}(\eps^*_{X_{\os{\circ}{T}_0}/\os{\circ}{T}}(E))$.  
The complex 
$\wt{R}u_{X_{\os{\circ}{T}_0}/\os{\circ}{T}*}(\eps^*_{X_{\os{\circ}{T}_0}/\os{\circ}{T}}(E))$ 
fits into the following triangle:  
\begin{align*} 
& Ru_{X_{T_0}/T*}(\eps^*_{X_{T_0}/T}(E))(-1)[-1]\lo 
\wt{R}u_{X_{\os{\circ}{T}_0}/\os{\circ}{T}*}
(\eps^*_{X_{\os{\circ}{T}_0}/\os{\circ}{T}}(E))
\lo \tag{0.1.2.5}\label{ali:tripte} \\
& Ru_{X_{T_0}/T*}
(\eps^*_{X_{T_0}/T}(E)) \os{+1}{\lo}. 
\end{align*} 
Because 
$\wt{R}u_{X_{\os{\circ}{T}_0}/\os{\circ}{T}*}(\eps^*_{X_{\os{\circ}{T}_0}/\os{\circ}{T}}(E))$ 
is not equal to the crystalline complex $Ru_{X_{\os{\circ}{T}_0}/\os{\circ}{T}*}
(\eps^*_{X_{\os{\circ}{T}_0}/\os{\circ}{T}}(E))$ in general, 
we call $(\wt{R}u_{X_{\os{\circ}{T}_0}/\os{\circ}{T}*}
(\eps^*_{X_{\os{\circ}{T}_0}/\os{\circ}{T}}(E)),P)$ the 
{\it modified $P$-filtered log crystalline complex of $E$} 
($P$ is the abbreviation of Poincar\'{e}).
We give the definition of 
$$(\wt{R}u_{X_{\os{\circ}{T}_0}/\os{\circ}{T}*}
(\eps^*_{X_{\os{\circ}{T}_0}/\os{\circ}{T}}(E)),P)$$ 
by using an embedding system of $X$ over $\ol{S(T)^{\nat}}$ 
not over $S(T)^{\nat}$ 
and by using the exactification defined in \cite{s3}. 
We prove that it depends only on $X/S$ and 
$(\os{\circ}{T},{\cal J},\del)/(\os{\circ}{S},{\cal I},\gam)$. 
In the case $E={\cal O}_{\os{\circ}{X}_{T_0}/\os{\circ}{T}}$, 
we denote 
$(\wt{R}u_{X_{\os{\circ}{T}_0}/\os{\circ}{T}*}
(\eps^*_{X_{\os{\circ}{T}_0}/\os{\circ}{T}}(E)),P)$ 
by 
$(\wt{R}u_{X_{\os{\circ}{T}_0}/\os{\circ}{T}*}({\cal O}_{X_{T_0}/\os{\circ}{T}}),P)$. 
In the case $E={\cal O}_{\os{\circ}{X}_{T_0}/\os{\circ}{T}}$ and 
$T=S={\cal W}_n(s)$ ($s$ is the log point of a perfect field $\kap$ of characteristic $p>0$), 
Hyodo and Kato 
have used the triangle (\ref{ali:tripte}) essentially 
to define the $p$-adic monodromy operator of 
$Ru_{X_{T_0}/T*}({\cal O}_{X_{T_0}/T})$ in \cite{hk}. 
\par 
Now assume that $(T,{\cal J},\del)=(S,p{\cal O}_S,[~])$ and 
that $\os{\circ}{S}$ is the spectrum of the Witt ring ${\cal W}_n$ 
of length $n >0$ of a perfect field $\kap$. 
Let $({\cal W}_n\wt{\Om}^{\bul}_X,P)$ be the filtered complex defined 
in \cite{msemi} (the complex ${\cal W}_n\wt{\Om}^{\bul}_X$ 
has been defined essentially in \cite{hdw}).  
In [loc.~cit.] Hyodo and Mokrane have used 
the blow up of the product of two local admissible lifts of $X$ 
and have calculated the log de Rham-complexes of the blow up concretely 
for the proof of the well-definedness of $({\cal W}_n\wt{\Om}^i_X,P)$,  
though I am not sure that this concrete calculation 
in [loc.~cit.] implies 
the well-definedness of $({\cal W}_n\wt{\Om}^i_X,P)$ 
(see (\ref{rema:qok}) below). 
In this book 
we prove that there exists the following canonical filtered isomorphism 
\begin{align*} 
{\cal H}^i
((\wt{R}u_{X/\os{\circ}{s}*}({\cal O}_{X/\os{\circ}{s}}),P))
&:=
({\cal H}^i
(\wt{R}u_{X/\os{\circ}{s}*}({\cal O}_{X/\os{\circ}{s}}),
\{{\cal H}^i
(P_k\wt{R}u_{X/\os{\circ}{s}*}({\cal O}_{X/\os{\circ}{s}}))\}_{k\in {\mab Z}})
\tag{0.1.2.6}\label{ali:unad} \\
& \os{\sim}{\lo} 
({\cal W}_n\wt{\Om}^i_X,P) \quad (i\in {\mab N}). 
\end{align*}
(We can prove that the natural morphism 
\begin{align*} 
{\cal H}^i
(P_k\wt{R}u_{X/\os{\circ}{s}*}({\cal O}_{X/\os{\circ}{s}}))\lo 
{\cal H}^i(\wt{R}u_{X/\os{\circ}{s}*}({\cal O}_{X/\os{\circ}{s}}))
\quad (k\in {\mab Z})
\end{align*} 
is injective.)
Because we do not use the local admissible lift  
for the definition of $(\wt{R}u_{X/\os{\circ}{s}*}({\cal O}_{X/\os{\circ}{s}}),P)$,  
we can consider the left hand side of (\ref{ali:unad}) as a new definition of 
$({\cal W}_n\wt{\Om}^i_X,P)$. 
In this book, the exactification in \cite{s3} replaces the blow up
and it is not necessary to calculate 
the log de Rham-complexes of the blow up concretely 
by using an explicit expression of 
$P_0{\cal W}_n\wt{\Om}^i_X$ and (\ref{ali:trgrte}) 
for the proof of the well-definedness of 
$({\cal W}_n\wt{\Om}^i_X,P)$. 
Our definition is more flexible than the definition in [loc.~cit.] 
because the product of two embedding systems gives us an embedding system 
(the product of two lifts is not a lift). 
\par 
We explain (16). 
\par 
Let $({\cal W}_nA^{\bul}_X,P)$ be 
Mokrane's filtered Steenbrink-Hyodo complex 
constructed in \cite{msemi} (and denoted by 
$(W_nA^{\bul}_X,P)$ in [loc.~cit.]) 
and completed in \cite{ndw} except (\ref{prop:plz}) and (\ref{prop:prm}) below 
(see also (\ref{rema:qok}) (2), (3) below). 
The relation between our zariskian $p$-adic filtered Steenbrink complex 
$(A_{\rm zar}(X/S),P)$ by the use of filtered log crystalline complexes 
and $({\cal W}_nA^{\bul}_X,P)$ 
by the use of filtered log de Rham-Witt complexes is as follows.
\par   
If $\os{\circ}{S}$ is the spectrum of 
the Witt ring ${\cal W}_n$ of $\kap$ of length $n >0$, 
then we prove that 
there exists a canonical filtered isomorphism 
\begin{equation*} 
(A_{\rm zar}(X/S),P) \os{\sim}{\lo} ({\cal W}_nA^{\bul}_X,P).  
\tag{0.1.2.7}\label{eqn:axswp} 
\end{equation*} 
Via this canonical isomorphism,  
the quasi-monodromy operator 
$$\nu_{\rm zar}\col 
(A_{\rm zar}(X/S),P) \lo 
(A_{\rm zar}(X/S),P\langle -2\rangle)(-1)$$ 
is equal to an analogous morphism 
$$\nu_{\rm dRW}\col ({\cal W}_nA^{\bul}_X,P) \lo 
({\cal W}_nA^{\bul}_X,P\langle -2\rangle)(-1),$$ 
which has been constructed in \cite{msemi} 
for the case $\os{\circ}{S}={\rm Spec}({\cal W}_n)$ 
and $E={\cal O}_{\os{\circ}{X}/\os{\circ}{S}}$  
(but $\nu$ in \cite{msemi} is 
mistaken because $\nu$ is not a morphism of complexes 
and corrected in \cite{ndw}).  
In fact, we can generalize (\ref{eqn:axswp}) as follows. 
\par 
Let $s$ be the log point of $\kap$. 
Let ${\cal W}_n$ (resp.~${\cal W}$) 
be the Witt ring of $\kap$ of length $n>0$ (resp.~the Witt ring of $\kap$). 
Let $E$ be a flat coherent crystal of 
${\cal O}_{\os{\circ}{X}/{\cal W}_n(\os{\circ}{s})}$-modules. 
Then we construct a filtered complex 
$$({\cal W}_nA^{\bul}_X(E),P)$$ 
\parno  
and we prove that there exists 
the following canonical filtered isomorphism  
\begin{equation*} 
(A_{\rm zar}(X/{\cal W}_n(s),E),P) 
\os{\sim}{\lo} ({\cal W}_nA^{\bul}_{X}(E),P).      
\tag{0.1.2.8}\label{eqn:axetwp} 
\end{equation*} 
(In the case where $E={\cal O}_{\os{\circ}{X}/{\cal W}_n(\os{\circ}{s})}$, 
(\ref{eqn:axetwp}) is equal to (\ref{eqn:axswp}).) 
This filtered isomorphism is compatible with the projections of 
both hand sides on (\ref{eqn:axetwp}).  
\par
Because the constant simplicial case above is 
a restricted case for the applications $(19)\sim(27)$,  
we generalize the main results above except $(10)\sim (13)$ 
in the constant simplicial case to 
the split truncated simplicial case, 
which is a bridge to the explanation for the Chapter III. 
We explain this briefly as follows. 
\par 
Let $N$ be a nonnegative integer. 
Let $(T,{\cal J},\del)$ be a fine log PD-scheme such that 
${\cal J}$ is quasi-coherent and such that 
$p$ is locally nilpotent on $T$. 
Let $T_0$ be the exact closed log subscheme of 
$T$ defined by ${\cal J}$. 
Let $Y_{\bul \leq N}$ be 
a smooth $N$-truncated simplicial scheme over $T_0$ 
and let $D_{\bul \leq N}$ be 
an $N$-truncated simplicial relative SNCD on 
$Y_{\bul \leq N}/T_0$. 
Then we have the log structure 
$M(D_{\bul \leq N})$ on the Zariski site of 
$Y_{\bul \leq N}$ associated to $D_{\bul \leq N}$ 
(\cite{klog1}, \cite{fao}, \cite{nh2}).  
Let 
$\eps_{(Y_{\bul \leq N},M(D_{\bul \leq N}))/T_0} \col 
(Y_{\bul \leq N},M(D_{\bul \leq N}))\lo Y_{\bul \leq N}$ 
be the natural morphism over $T_0$ 
forgetting the log structure $M(D_{\bul \leq N})$. 
Let 
\begin{align*} 
\eps_{(Y_{\bul \leq N},M(D_{\bul \leq N}))/T}\col &
(((Y_{\bul \leq N},M(D_{\bul \leq N}))/T)_{\rm crys},
{\cal O}_{(Y_{\bul \leq N},M(D_{\bul \leq N}))/T})\\
& \lo ((Y_{\bul \leq N}/T)_{\rm crys},{\cal O}_{Y_{\bul \leq N}/T})
\end{align*} 
be the induced morphism of ringed topoi by  
$\eps_{(Y_{\bul \leq N},M(D_{\bul \leq N}))/T_0}$. 
In \cite{nh3} we have considered the filtered complex 
$(R\eps_{(Y_{\bul \leq N},M(D_{\bul \leq N}))/T*}
({\cal O}_{(Y_{\bul \leq N},M(D_{\bul \leq N}))/T}),\tau)$,  
where $\tau$ is the canonical filtration. 
\par 
Let $(T,{\cal J},\del)$ and $S$ be as in the explanation for (1). 
Let $X_{\bul \leq N}/S$ be a split $N$-truncated simplicial SNCL scheme.
For $U_0:=T_0$ or $S_{\os{\circ}{T}_0}$, 
set $X_{\bul \leq N,U_0}:=X_{\bul \leq N}\times_SU_0$. 
Then  $X_{\bul \leq N,S_{\os{\circ}{T}_0}}=
X_{\bul \leq N}\times_{\os{\circ}{S}}\os{\circ}{T}_0
=:X_{\bul \leq N,\os{\circ}{T}_0}$. 
Let 
$
f\col X_{\bul \leq N,\os{\circ}{T}_0}\lo S_{\os{\circ}{T}_0}$ 
be the structural morphism. 
By abuse of notation,  let us also denote 
the structural morphism 
$X_{\bul \leq N,\os{\circ}{T}_0}\lo  S(T)^{\nat}$ 
by  $f$. 
(We have already called this $N$-truncated simplicial case 
the naive case.) 
Let $E^{\bul \leq N}$ be a flat coherent crystal of      
${\cal O}_{\os{\circ}{X}_{\bul \leq N,T_0}/\os{\circ}{T}}$-modules. 
In this truncated simplicial SNCL case, 
we have not been able to find 
a quickly obtained cosimplicial filtered complex 
as in the previous paragraph.  
To obtain an $N$-truncated cosimplicial filtered complex 
\begin{align*} 
(A_{\rm zar}(X_{\bul \leq N,\os{\circ}{T}_0}/S(T)^{\nat}, E^{\bul \leq N}),P)\in 
{\rm D}^+{\rm F}(f^{-1}({\cal O}_T))
\tag{0.1.2.9}\label{ali:ntcf}
\end{align*}  
which is an $N$-truncated cosimplicial version of 
$(A_{\rm zar}(X_{\os{\circ}{T}_0}/S(T)^{\nat},E),P)$, 
we construct an $(N,\infty)$-truncated bisimplicial  
embedding system and 
we give an explicit filtered complex representing 
$(A_{\rm zar}(X_{\bul \leq N,\os{\circ}{T}_0}/S(T)^{\nat},E^{\bul \leq N}),P)$, 
which is shown to be independent of the system. 
As a generalization of (\ref{eqn:espsp}),
we prove that there exists the following spectral sequence 
\begin{equation*} 
\begin{split} 
{} & 
E_1^{-k,q+k}=
\bigoplus_{m\geq 0}\bigoplus_{j\geq \max \{-(k+m),0\}} 
R^{q-2j-k-2m}
f_{\os{\circ}{X}{}^{(2j+k+m)}_{m,T_0}/\os{\circ}{T}*}
(E^m_{\os{\circ}{X}{}^{(2j+k+m)}_{m,T_0}/\os{\circ}{T}} \\ 
{} & \phantom{R^{q-2j-k-m}f_{(\os{\circ}{X}^{(k)}, 
Z\vert_{\os{\circ}{X}^{(2j+k)}})/S*} 
({\cal O}}
\otimes_{\mab Z}
\varpi^{(2j+k+m)}_{\rm crys}
(\os{\circ}{X}_{m,T_0}/\os{\circ}{T}))(-j-k-m)\\
&\Lo R^qf_{X_{\bul \leq N,\os{\circ}{T}_0}/S(T)^{\nat}*}(
\eps^*_{X_{\bul \leq N,\os{\circ}{T}_0}/S(T)^{\nat}}
(E^{\bul \leq N}))
\quad (q\in {\mab Z}). 
\end{split} 
\tag{0.1.2.10}\label{eqn:despsp}
\end{equation*}   
\par
In the Chapter III, we have to ramify the log structure $M_S$ of $S$ 
for the applications $(19)\sim(27)$. 
Now we change a notation: denote $S$ by $S_0$.  
Consider a sequence 
\begin{equation*} 
S_N \lo S_{N-1} \lo \cdots  \lo S_0
\tag{0.1.2.11}\label{eqn:nrn10}
\end{equation*} 
of morphisms of $p$-adic formal families of log points.   
For simplicity of notation, denote the sequence 
(\ref{eqn:nrn10}) by $\{S_m\}_{m=0}^N$ 
and set $S:=S_N$ (we hope that the reader is not confused by this notation). 
We give the definition of the successive case as follows.
(This notion is necessary for the construction of the limit of 
the weight filtration on the infinitesimal cohomology 
of a proper scheme over $K$ 
because the semistable reduction theorem 
has not been proved in the mixed characteristics case.)   
\par 
For a fine $m$-truncated simplicial log scheme 
$Y_{\bul \leq m}$ over $S_m$,  
let $s^{m-1}_i \col Y_{m-1} \lo Y_m$ 
$(0\leq m\leq N, 0\leq i \leq m-1)$ be 
the degeneracy morphism 
corresponding to the standard degeneracy map 
$\partial^i_m \col [m] \lo [m-1]$:  
$\partial^i_m(j)=j$ $(0\leq j \leq i)$, 
$\partial^i_m(j)=j-1$ $(i< j \leq m)$. 
Let $\{X(m)_{\bul \leq m}/S_m\}_{m=0}^N$ 
be a set of split truncated simplicial log smooth schemes 
such that the intersection of the complements of 
$s^{m-1}_i(X(m)_{m-1})$ $(0 \leq i \leq m-1)$ 
in $X(m)_m$ is an SNCL scheme $N_m$ over $S_m$. 
We say that $\{X(m)_{\bul \leq m}/S_m\}_{m=0}^N$ 
is a {\it family of successive 
split $N$-truncated simplicial SNCL schemes 
with respect to the sequence} 
(\ref{eqn:nrn10}) 
if 
\begin{equation*} 
X(m)_{\bul \leq m-1}=
X(m-1)_{\bul \leq m-1}\times_{S_{m-1}}S_m 
\end{equation*} 
and 
\begin{equation*} 
X(m)_m=
\coprod_{0\leq l \leq m}\coprod_{[m]\twoheadrightarrow [l]}(N_l\times_{S_l}S_m). 
\end{equation*} 
Denote $X(N)_{\bul \leq N}$ simply by $X_{\bul \leq N}$.  
Let $(T,{\cal J},\del)$ be a log PD-enlargement of $S$. 
Assume that $(\os{\circ}{T},{\cal J},\del)$ is a $p$-adic formal PD-schemes 
in the sense of \cite{bob}. 
Then, by (\ref{eqn:nrn10}) we have the following sequence 
\begin{equation*} 
(S(T)^{\nat},{\cal J},\del)=(S_N(T)^{\nat},{\cal J},\del) \lo (S_{N-1}(T)^{\nat},{\cal J},\del) 
\lo \cdots  \lo (S_0(T)^{\nat},{\cal J},\del)
\tag{0.1.2.12}\label{eqn:tsn}
\end{equation*} 
of morphisms of $p$-adic formal PD-families of log points.  
Let $T_0$ be the exact closed log subscheme of $T$ 
defined by ${\cal J}$.  
For $U_0:=T_0$ or $S_{\os{\circ}{T}_0}$, set 
$X(m)_{\bul \leq m,U_0}:=X(m)_{\bul \leq m}\times_{S_m}U_0$, 
$X_{\bul \leq N,U_0}:=X_{\bul \leq N}\times_{S}U_0$, 
$\os{\circ}{N}_{l,U_0}:=\os{\circ}{N}{}_l\times_{\os{\circ}{S}_l}\os{\circ}{U}_0$ 
and  
$\os{\circ}{N}{}^{(k)}_{l,U_0}
:=\os{\circ}{N}{}^{(k)}_l\times_{\os{\circ}{S}_l}\os{\circ}{U}_0$ $(k\in {\mab N})$.  
Let $E^{\bul \leq N}$ be a flat coherent crystal of   
${\cal O}_{\os{\circ}{X}_{\bul \leq N,T_0}/\os{\circ}{T}}$-modules.  
Let $f\col X_{\bul \leq N,\os{\circ}{T}_0}\lo S(T)^{\nat}$  
be the structural morphism. 
Then we construct a filtered complex 
\begin{equation*} 
(A_{{\rm zar},{\mab Q}}(X_{\bul \leq N,\os{\circ}{T}_0}/S(T)^{\nat},E^{\bul \leq N}),P) 
\tag{0.1.2.13}\label{eqn:axtn}
\end{equation*}  
in ${\rm D}^+{\rm F}(f^{-1}({\cal O}_ T)
\otimes_{\mab Z}{\mab Q})$, 
which we call the $p$-adic Steenbrink complex of $E^{\bul \leq N}$. 
(Because we allow the ramification of $M_S$, 
it seems to me that we cannot 
obtain the integral $p$-adic filtered Steenbrink complex 
of $E^{\bul \leq N}$ in ${\rm D}^+{\rm F}(f^{-1}({\cal O}_ T))$ 
for $X_{\bul \leq N}/S$ and $(T,{\cal J},\del)$ with morphism $T_0\lo S$.) 
The filtered complex (\ref{eqn:axtn}) 
is a generalization of the filtered complex 
(\ref{ali:ntcf})$\otimes^L_{\mab Z}{\mab Q}$ 
in the $p$-adic case. 
Set 
$\os{\circ}{X}{}^{(k)}_m:=
\coprod_{0\leq l \leq m}\coprod_{[m] 
\twoheadrightarrow [l]}
\os{\circ}{N}{}^{(k)}_{l}$ and 
$\os{\circ}{X}{}^{(k)}_{m,T_0}:=
\coprod_{0\leq l \leq m}\coprod_{[m] 
\twoheadrightarrow [l]}
\os{\circ}{N}{}^{(k)}_{l,T_0}$ 
for $0\leq m\leq N$. 
Let 
$a^{(k)}_{m,T_0}\col \os{\circ}{X}{}^{(k)}_{m,T_0}\lo \os{\circ}{X}_{m,T_0}$ 
be the natural morphism of schemes. 
(Though we should denote 
$a^{(k)}_{m,T_0}$ by 
$a^{(k)}_{m,\os{\circ}{T}_0}$, 
we do not do so for simplicity of notation.) 
Let 
$$a^{(k)}_{m,T{\rm crys}}
\col 
((\os{\circ}{X}{}^{(k)}_{m,T_0}/\os{\circ}{T})_{\rm crys},
{\cal O}_{\os{\circ}{X}{}^{(k)}_{m,T_0}/\os{\circ}{T}})
\lo 
((\os{\circ}{X}_{m,T_0}/\os{\circ}{T})_{\rm crys},{\cal O}_{\os{\circ}{X}_{m,T_0}/\os{\circ}{T}})$$ 
be the induced morphism of topoi by $a^{(k)}_{m,T_0}$. 
Set 
$E^m\vert_{\os{\circ}{X}{}^{(k)}_{m,T_0}/\os{\circ}{T}}:=a^{(k)*}_{m,T{\rm crys}}(E^m)$. 
Let 
$\os{\circ}{p}{}^{(k)}_{m,T_0}\col 
\os{\circ}{X}{}^{(k)}_{m,T_0}\lo 
\os{\circ}{X}{}^{(k)}_m$ be the first projection 
and let 
$\os{\circ}{p}{}^{(k)}_{m,T{\rm crys}}
\col (\os{\circ}{X}{}^{(k)}_{m,T_0}/\os{\circ}{T})_{\rm crys}
\lo (\os{\circ}{X}{}^{(k)}_{m}/\os{\circ}{S})_{\rm crys}$ 
be the induced morphism of topoi by 
$\os{\circ}{p}{}^{(k)}_{m,T_0}$. 
Set 
$\varpi^{(k)}_{\rm crys}
(\os{\circ}{X}_{m,T_0}/\os{\circ}{T}):=
\os{\circ}{p}{}^{(k),-1}_{m,T{\rm crys}}
\varpi^{(k)}_{\rm crys}
(\os{\circ}{X}_{m}/\os{\circ}{S})$. 
As a generalization of 
(\ref{eqn:despsp})$\otimes_{\mab Z}{\mab Q}$ 
in the $p$-adic case,
we obtain the following spectral sequence 
\begin{align*} 
& E_1^{-k,q+k} = \bigoplus_{m=0}^N
\bigoplus_{j\geq \max \{-(k+m),0\}} 
R^{q-2j-k-2m}
f_{\os{\circ}{X}{}^{(2j+k+m)}_{m,T_0}/\os{\circ}{T}*}
(E^m_{\os{\circ}{X}{}^{(2j+k+m)}_{m,T_0}/\os{\circ}{T}}
\tag{0.1.2.14}\label{eqn:getpsp} \\
&\otimes_{\mab Z}\vp^{(2j+k+m)}_{\rm crys}
(\os{\circ}{X}_{m,\os{\circ}{T}_0}/\os{\circ}{T})) \\
& (-j-k-m)\otimes_{\mab Z}{\mab Q}
\Lo 
R^qf_{X_{\bul \leq N,\os{\circ}{T}_0}/S(T)^{\nat}*}
(\eps^*_{X_{\bul \leq N,\os{\circ}{T}_0}/S(T)^{\nat}}
(E^{\bul \leq N}))\otimes_{\mab Z}{\mab Q}
\quad (q\in {\mab Z}).  
\end{align*}  
As in \cite{nh2}, we prove that 
(\ref{eqn:getpsp}) degenerates at $E_2$ 
for the case $E^{\bul \leq N}=
{\cal O}_{\os{\circ}{X}_{\bul \leq N,T_0}/\os{\circ}{T}}$ ((6)). 
\par 
In the case $T=S={\cal W}(s)$, 
we also construct a filtered complex 
$({\cal W}A_{X_{\bul \leq N},{\mab Q}}(E^{\bul \leq N}),P)$ 
for the successive split truncated simplicial SNCL case, 
which is a generalization of 
$({\cal W}A^{\bul}_{X_{\bul \leq N}}(E^{\bul \leq N})
\otimes_{\mab Z}{\mab Q},P)$ for the naive split truncated simplicial SNCL case. 
We prove that there exists a canonical isomorphism 
\begin{equation*} 
(A_{{\rm zar},{\mab Q}}(X_{\bul \leq N}/{\cal W}(s),E^{\bul \leq N}),P) 
\os{\sim}{\lo} 
({\cal W}A_{X_{\bul \leq N},{\mab Q}}(E^{\bul \leq N}),P)
\tag{0.1.2.15}\label{eqn:wzxytqep} 
\end{equation*} 
in ${\rm D}^+{\rm F}(f^{-1}({\cal W}(\kap))
\otimes_{\mab Z}{\mab Q})$. 
Because log de Rham-Witt complexes are useful for 
the calculation of the slope filtrations on the cohomologies of 
log de Rham-Witt complexes, 
the comparison (\ref{eqn:wzxytqep}) 
is useful for the calculation of the slope filtration on 
the infinitesimal cohomology of a proper scheme over $K$. 
\par 
In \S\ref{sec:bcsncl} we define the notions of 
a truncated simplicial base change 
of SNCL schemes and an admissible immersion ((17)). 
Because the definitions of these notions are involved,  
we do not give them in this introduction; 
see \S\ref{sec:bcsncl} for the precise definitions. 
By virtue of these notions, 
we can treat the log isocrystalline cohomologies in 
the naive truncated simplicial case 
and the successive truncated simplicial case at the same time 
and we can state various results in the Chapters III, IV and V simply. 
\par 
We explain (18). 
\par 
Let $N$ be a nonnegative integer.  
Endow ${\rm Spec}({\cal V})$ with the canonical log structure 
and let $S$ be the resulting log scheme. 
Let $\pi$ be a uniformizer of ${\cal V}$. 
Set $s:=S\otimes_{\cal V}\kap$. 
Let ${\cal Y}_{\bul \leq N}$ 
be a log smooth $N$-truncated simplicial log scheme over $S$.  
Let $Y_{\bul \leq N}/s$ be the log special fiber of ${\cal Y}_{\bul \leq N}/S$. 
Assume that $\os{\circ}{\cal Y}_{\bul \leq N}/\os{\circ}{S}$ 
is proper and that the structural morphism 
$Y_{\bul \leq N} \lo s$ is of Cartier type. 
Let ${\mathfrak Y}_{\bul \leq N}$ be 
the log generic fiber of ${\cal Y}_{\bul \leq N}$. 
In \cite{tsgep} Tsuji has proved that 
there exists the following isomorphism 
\begin{equation*} 
H^q(\Psi_{\pi}) \col H^q_{{\rm crys}}(Y_{\bul \leq N}/{\cal W}(s))\otimes_{{\cal W}}K 
\os{\sim}{\lo} H^q_{\rm dR}({\mathfrak Y}_{\bul \leq N}/K)
\tag{0.1.2.16}\label{eqn:raghk}
\end{equation*} 
depending on $\pi$. 
This is a generalization of the famous Hyodo-Kato isomorphism 
which has been constructed in \cite{hk}. 
(However see (\ref{rema:istp}) (2) below.)
In fact, we can construct an isomorphism 
\begin{equation*} 
\Psi_{\pi} \col R\Gam_{{\rm crys}}(Y_{\bul \leq N}/{\cal W}(s))
\otimes_{{\cal W}}^LK 
\os{\sim}{\lo} R\Gam_{\rm dR}({\mathfrak Y}_{\bul \leq N}/K)
\tag{0.1.2.17}\label{eqn:radghk}
\end{equation*} 
depending on $\pi$. 
By considering the virtual existence of the Hyodo-Kato isomorphism by using 
the element $p$ instead of $\pi$, 
we can get rid of 
the dependence of the Hyodo-Kato isomorphism on $\pi$, 
that is, we construct a canonical isomorphism  
\begin{equation*} 
H^q(\Psi) \col H^q_{{\rm crys}}(Y_{\bul \leq N}/{\cal W}(s))\otimes_{\cal W}K 
\os{\sim}{\lo} H^q_{\rm dR}({\mathfrak Y}_{\bul \leq N}/K)
\tag{0.1.2.18}\label{eqn:raghbk}
\end{equation*} 
which is independent of the choice of $\pi$. 
To prove this, we use Hyodo-Kato's formula 
which tells us a relation between 
$H^q(\Psi_{\pi})$ and $H^q(\Psi_{a\pi})$ for $a\in {\cal V}$. 
Our modified isomorphism (\ref{eqn:raghbk}) is necessary for (26). 
\par 
We conclude this introduction by 
stating (19), (20), (21), (22), (23), (24), (25), (26) and (27) in some details. 
\par 
Let ${\cal X}_K$ be a proper scheme over $K$. 
Let ${\cal X}$ be a proper flat model of 
${\cal X}_K$ over ${\cal V}$. 
Let $H^q_{\rm inf}({\cal X}_K/K)$ $(q\in {\mab N})$ 
be the infinitesimal cohomology of 
${\cal X}_K/K$ defined by the infinitesimal topos 
$({\cal X}_K/K)_{\rm inf}$ and 
the structure sheaf ${\cal O}_{{\cal X}_K/K}$ in 
$({\cal X}_K/K)_{\rm inf}$: 
$H^q_{\rm inf}({\cal X}_K/K)
:=H^q(({\cal X}_K/K)_{\rm inf},{\cal O}_{{\cal X}_K/K})$ 
(\cite{grcr}). 
Let 
\begin{equation*} 
{\cal V}=:{\cal V}_{-1} \subset 
{\cal V}_0 \subset \cdots 
\subset {\cal V}_N \subset \cdots 
\end{equation*} 
be a sequence of finite extensions of  
complete discrete valuation rings of mixed characteristics.   
Set $K_m:={\rm Frac}({\cal V}_m)$. 
Let $\kap_m$ $(-1\leq m\leq N)$ 
be the residue field of ${\cal V}_m$. 
Endow ${\rm Spec}({\cal V}_m)$ 
with the canonical log structure 
$({\cal V}_m\setminus \{0\}) \os{\sus}{\lo} {\cal V}_m$
and let ${\rm Spec}^{\log}({\cal V}_m)$ 
be the resulting log scheme. 
Then we have the following sequence: 
\begin{equation*} 
{\rm Spec}^{\log}({\cal V}_0)\longleftarrow 
{\rm Spec}^{\log}({\cal V}_1)\longleftarrow 
\cdots \longleftarrow {\rm Spec}^{\log}({\cal V}_N)
\longleftarrow \cdots. 
\tag{0.1.2.20}\label{eqn:ssmeq}
\end{equation*} 
Let $s_m$ be the log special fiber of ${\rm Spec}^{\log}({\cal V}_m)$.  
Let 
\begin{equation*} 
s_0\longleftarrow s_1\longleftarrow 
\cdots \longleftarrow s_N 
\longleftarrow \cdots
\tag{0.1.2.21}\label{eqn:sseq}
\end{equation*} 
be the sequence of log points obtained by the sequence (\ref{eqn:ssmeq}).  
Fix a nonnegative integer $q$. 
Set $\kap_m=\Gam(s_m,{\cal O}_{s_m})$. 
Let $N$ be a large positive integer relative to $q$ 
(e.~g., $N\geq 2^{-1}(q+1)(q+2))$. 
Set ${\cal X}_{{\cal V}_m}:={\cal X}\otimes_{\cal V}{\cal V}_m$ 
and ${\cal X}_{K_m}:={\cal X}_K\otimes_KK_m$ 
$(-1\leq m\leq N)$.  
Then, in \S\ref{sec:ph}, 
we construct a proper log smooth 
split $N$-truncated simplicial log scheme 
${\cal X}_{\bul \leq N}$ over 
${\rm Spec}^{\log}({\cal V}_N)$ 
satisfying the following conditions: 
\par 
(a) $\os{\circ}{\cal X}_{\bul \leq N}\otimes_{{\cal V}_N}K_N
\lo {\cal X}_{K_N}$ is an $N$-truncated proper hypercovering.   
\par 
(b) Each $\os{\circ}{\cal X}_m$ $(0\leq m\leq N)$ 
is the disjoint union of the base changes of 
proper strict semistable families over ${\cal V}_l$'s 
$(0\leq l\leq m)$ with respect to the morphisms  
${\rm Spec}({\cal V}_N)\lo {\rm Spec}({\cal V}_l)$'s. 
\par 
(c) The log special fiber $X_{\bul \leq N}$ 
of ${\cal X}_{\bul \leq N}$ is 
the $N$-truncated simplicial log scheme 
associated to the family of 
successive split $N$-truncated simplicial SNCL schemes  
with respect to the sequence 
$s_0 \longleftarrow s_1 \longleftarrow 
\cdots  \longleftarrow s_N$.  
\medskip 
\parno
De Jong's theorem in \cite{dj} about 
the semistable reduction theorem by an alteration 
assures the existence of 
${\cal X}_{\bul \leq N}$.  
Then (19) is the following. 
As a corollary of the existence of the spectral sequence 
(\ref{eqn:getpsp}) and the truncated simplicial version of 
the canonical Hyodo-Kato's isomorphism (\ref{eqn:raghbk}) 
and the proper cohomological descent in characteristic $0$,  
we obtain the following spectral sequence:  
\begin{align*} 
&E_1^{-k,q+k} = \bigoplus_{m=0}^N
\bigoplus_{j\geq \max \{-(k+m),0\}} 
H^{q-2j-k-2m}
((\os{\circ}{X}{}^{(2j+k+m)}_m/
{\cal W}(\os{\circ}{s}_m))_{\rm crys}, 
\tag{0.1.2.22}\label{eqn:esptrfsp}\\
& ({\cal O}_{\os{\circ}{X}{}^{(2j+k+m)}_m
/{\cal W}(\os{\circ}{s}_m)}
\otimes_{\mab Z} 
\vp^{(2j+k+m)}_{\rm crys}(\os{\circ}{X}_m
/{\cal W}(\os{\circ}{s}_m))))
(-j-k-m)\otimes_{{\cal W}(\kap_m)}K_N \\
& \Lo 
H^q_{\rm inf}({\cal X}_{K_N}/K_N)
\end{align*} 
by proving that 
\begin{align*} 
H^q_{\rm crys}(X_{\bul \leq N}/{\cal W}(s_N))\otimes_{{\cal W}(\kap_N)}
K_N=
H^q_{\rm inf}({\cal X}_{K_N}/K_N). 
\tag{0.1.2.23}\label{eqn:escrinffsp}
\end{align*} 
By virtue of (7) this spectral sequence degenerates at $E_2$ ((20)). 
The spectral sequence (\ref{eqn:esptrfsp}) induces an increasing filtration 
$P(N)$ on $H^q_{\rm inf}({\cal X}_{K_N}/K_N)$. 
We prove that $P(N)$ depends only on ${\cal X}_K$ and $K_N$.   
Furthermore we prove that there exists 
a well-defined filtration $P$ 
on $H^q_{\rm inf}({\cal X}_K/K)$ 
which induces $P(N)$ on 
$H^q_{\rm inf}({\cal X}_{K_N}/K_N)
=H^q_{\rm inf}({\cal X}_K/K)\otimes_KK_N$ ((21)). 
Roughly speaking, 
(21) is a $p$-adic analogue of du Bois' result 
for the complex analytic case (\cite{db}), though 
our method is not the same as his method as mentioned before.   
We prove that this filtration is strictly compatible with respect to 
the pull-back of a morphism of proper schemes over $K$ ((22)).  
We prove that there exists a well-defined $K_0$-structure 
$H^q_{\rm inf}({\cal X}_K/K/K_0)$ on 
$H^q_{\rm inf}({\cal X}_K/K)$ with a well-defined Frobenius action ((23)). 
Let $K_{N,0}$ be the fraction field of 
the Witt ring of $\Gam(s_N,{\cal O}_{s_N})$. 
The vector space 
$H^q_{\rm inf}({\cal X}_K/K/K_0)$ over $K_0$ satisfies the following equality: 
\begin{align*} 
H^q_{\rm inf}({\cal X}_K/K/K_0)\otimes_{K_0}K_{N,0}
=
H^q_{\rm crys}(X_{\bul \leq N}/{\cal W}(s_N))\otimes_{{\cal W}(\kap_N)}
K_{N,0}. 
\tag{0.1.2.24}\label{eqn:esc000ffsp}
\end{align*}
We calculate the slope filtration on 
$H^q_{\rm inf}({\cal X}_K/K/K_0)\otimes_{K_0}K_{N,0}$
by using $X_{\bul \leq N}$ ((24)). 
Using the argument of Galois descent, 
we prove that the filtration $P$ 
on $H^q_{\rm inf}({\cal X}_K/K)$ descends to 
a filtration $P$ on $H^q_{\rm inf}({\cal X}_K/K/K_0)$ ((25)). 
The equality (\ref{eqn:escrinffsp}) tells us that 
the log crystalline cohomology 
$H^q_{\rm crys}(X_{\bul \leq N}/{\cal W}(s_N))\otimes_{{\cal W}(\kap_N)}K_N$ 
depends only on ${\cal X}_{K_N}/K_N$; however 
the complex $R\Gam_{\rm crys}(X_{\bul \leq N}/{\cal W}(s_N))\otimes^L_{{\cal W}(\kap_N)}K_N$  
may depend on a uniformizer of ${\cal V}_N$; we conjecture the following 

\begin{conj}\label{conj:comiddd}
\begin{align*} 
R\Gam_{\rm inf}({\cal X}_{K_N}/K_N)=
R\Gam_{\rm crys}(X_{\bul \leq N}/{\cal W}(s_N))\otimes^L_{{\cal W}(\kap_N)}K_N
\tag{0.1.3.1}\label{ali:xknn}
\end{align*} 
In particular, the complex 
$R\Gam_{\rm crys}(X_{\bul \leq N}/{\cal W}(s_N))\otimes^L_{{\cal W}(\kap_N)}K_N$  
depends only on ${\cal X}_{K_N}/K_N$. 
\end{conj} 
Though we can construct an isomorphism 
\begin{align*} 
\Phi_{\pi_N} \col R\Gam_{\rm inf}({\cal X}_{K_N}/K_N)\os{\sim}{\lo} 
R\Gam_{\rm crys}(X_{\bul \leq N}/{\cal W}(s_N))\otimes^L_{{\cal W}(\kap_N)}K_N
\end{align*} 
for a uniformizer $\pi_N$ of ${\cal V}_N$, 
the isomorphism $\Phi_{\pi_N}$ may depend on the choice of the uniformizer $\pi_N$. 
To prove this conjecture, we first prove that 
the monodromy operator 
\begin{align*} 
N_{\rm zar} \col 
R\Gam_{\rm crys}(X_{\bul \leq N}/{\cal W}(s_N))\otimes^L_{{\cal W}(\kap_N)}K_N
\lo R\Gam_{\rm crys}(X_{\bul \leq N}/{\cal W}(s_N))\otimes^L_{{\cal W}(\kap_N)}K_N
\end{align*} 
is nilpotent. We prove this in (\ref{prop:dcis}) below. However 
we cannot prove the ``Hyodo-Kato's formula'' 
\begin{align*}
\Phi_{\pi'_N}=\Phi_{\pi_N}\exp((\log a)N_{\rm zar})
\tag{0.1.3.2}\label{ali:xkpinn}
\end{align*} 
as an equality of endomorphisms of 
$R\Gam_{\rm crys}(X_{\bul \leq N}/{\cal W}(s_N))\otimes^L_{{\cal W}(\kap_N)}K_N$ 
for another uniformizer $\pi'_N$ of ${\cal V}_N$, where $a:=\pi'_N/\pi_N$. 
\par 
By using the canonical Hyodo-Kato isomorphism (\ref{eqn:raghbk}), 
we prove the following theorem ((26), (27)):

\begin{theo}\label{theo:cwsasga}
Let ${\mathfrak U}$ be a separated scheme of finite type over $\ol{K}$. 
There exists a $($canonical$)$ semi-linear action 
\begin{align*} 
\rho^q_{\rm crys}:=
\rho^q_{{\mathfrak U},{\rm crys}}\col 
{\rm WD}_{\rm crys}(\ol{K})\lo 
{\rm GL}_{{\mab Q}_p}(H^q_{\rm inf}({\mathfrak U}/\ol{K}))
\tag{0.1.4.1}\label{ali:xqogr}
\end{align*}  
satisfying the following properties$:$ 
\par 
$(1)$ The action {\rm (\ref{ali:xqogr})} 
is contravariantly functorial with respect to a morphism
${\mathfrak U}\lo {\mathfrak V}$ of separated schemes of finite type over $\ol{K}$. 
That is, for a morphism 
${\mathfrak f}\col {\mathfrak U}\lo {\mathfrak V}$ of 
separated schemes of finite type over $\ol{K}$, 
the following diagram is commutative for an element 
$\gam\in {\rm WD}_{\rm crys}(\ol{K})\!:$
\begin{equation*} 
\begin{CD}
H^q_{\rm inf}({\mathfrak U}/\ol{K})
@>\rho^q_{{\mathfrak U},{\rm crys}}(\gam)>> 
H^q_{\rm inf}({\mathfrak U}/\ol{K})\\
@A{{\mathfrak f}^*}AA @AA{{\mathfrak f}^*}A \\
H^q_{\rm inf}({\mathfrak V}/\ol{K})
@>\rho^q_{{\mathfrak V},{\rm crys}}(\gam)>> 
H^q_{\rm inf}({\mathfrak V}/\ol{K}). 
\end{CD}
\tag{0.1.4.2}\label{cd:xqcgr}
\end{equation*} 
\par 
$(2)$ The action {\rm (\ref{ali:xqogr})} 
is compatible with the cup product of 
$H^{\bul}_{\rm inf}({\mathfrak U}/\ol{K})$. 
That is, for $x\in H^q_{\rm inf}({\mathfrak U}/\ol{K})$, 
$y\in H^{q'}_{\rm inf}({\mathfrak U}/\ol{K})$ and 
$\gam \in {\rm WD}_{\rm crys}(\ol{K})$, 
\begin{align*} 
\rho^{q+q'}_{{\mathfrak U},{\rm crys}}(\gam)(x\cup y)
=
\rho^q_{{\mathfrak U},{\rm crys}}(\gam)(x)\cup 
\rho^{q'}_{{\mathfrak U},{\rm crys}}(\gam)(y).
\tag{0.1.4.3}\label{ali:xxauxy}
\end{align*}  
\par 
$(3)$ When ${\mathfrak U}$ is proper over $\ol{K}$, 
then the action {\rm (\ref{ali:xqogr})} preserves the filtration $P$ on 
$H^q_{\rm inf}({\mathfrak U}/\ol{K})$. 
\end{theo} 
The theorem (1) in (\ref{theo:cwsasga}) is a generalization of Berthelot-Ogus' theorem 
which says that ${\rm W}_{\rm crys}(\ol{K})$ acts 
on $H^q_{\rm dR}({\cal Y}_{\ol{K}}/\ol{K})$ for a 
proper potential smooth family ${\cal Y}$ over ${\cal V}$ 
and that this action is contravariantly functorial (\cite{boi}). 
It also proves Ogus' conjecture for a proper semistable family over ${\cal V}$ (\cite{ollc})
in a more general form. 
Using the action of ${\rm WD}_{\rm crys}(\ol{K})$, 
we can define a linear action of 
${\rm WD}(\ol{K}/K)$ on $P_kH^q_{\rm inf}({\cal X}_{\ol{K}}/\ol{K})$ 
$(k\in {\mab Z})$ and we can conjecture 
the compatibility of this action with the linear action of  ${\rm WD}(\ol{K}/K)$ by 
the analogous linear subspace $P_kH^q_{\rm et}({\cal X}_{\ol{K}},{\mab Q}_l)$ 
of the $l$-adic \'{e}tale cohomology $H^q_{\rm et}({\cal X}_{\ol{K}},{\mab Q}_l)$. 
This is a straightforward generalization of 
the compatibility of the action of the Weil-Deligne group on 
the log isocrystalline cohomology with the action of the Weil-Deligne group 
obtained by the $l$-adic cohomology of 
a proper semistable family over ${\cal V}$ in \cite{fos} ((26)). 
\par 
In a future article we would like to discuss 
the convergent version of this book as 
carried out in \cite{nhw} for an open log case.  
In another future article we would like to discuss 
the themes of this book for an open log semistable case.  
\par 
Lastly we should mention Matsuue's different work 
(\cite{maoro}) from our work in a relative situation. 
He has constructed 
theory of relative log de Rham-Witt complexes, which is a log version 
of theory of relative de Rham-Witt complexes of Langer-Zink (\cite{lazi}). 
Especially he has proved the comparison theorem 
between the log crystalline cohomology of a relative SNCL scheme 
and the cohomology of the relative log de Rham-Witt complex 
of it and he has constructed a weight spectral sequence of it 
and he has proved the degeneration at $E_2$ 
of the weight spectral sequence modulo torsion in the case 
where the base log scheme is the log point of a (not necessarily perfect) field of 
characteristic $p>0$ by using ideas in \cite{ndw}. 
\par
\bigskip
\parno
{\bf Acknowledgment.} 
In this book, Berthelot's theory of filtered derived category (\cite{blec}) 
gives us an essential language as in \cite{nh2}. 
I would like to express my sincere gratitude to P.~Berthelot 
for his permission for using his theory in [loc.~cit.]. 
I would also like to express my sincere gratitude to A.~Shiho for 
teaching me a key idea for the construction of 
an important filtered complex in [loc.~cit.] 
and informing me of T.~Tsuji's result on the Hyodo-Kato isomorphism 
in \cite{tst}.  The filtered complex $(A_{\rm zar},P)$ in the Chapter I 
in the case of the trivial coefficient is 
an $N$-truncated SNCL version of the filtered complex in \cite{nh2}. 
\par 
I would like to express my thanks 
to T.~Tsuji for suggesting to me that we can consider
the nontrivial coefficient in the crystalline topos 
for the construction of filtered complexes $(A_{\rm zar},P)$ and 
to thank T.~Ito and Y.~Mieda 
for teaching me the importance of the contravariant 
(and covariant) functorialities of $p$-adic weight spectral sequences.
I would also like to thank T.~Kajiwara and H.~Matsuue for their kindness for 
informing me of their results \cite{kaji} and \cite{maoro}, respectively.
\par 
Lastly but not least, I would like to thank all mathematicians who appear in 
the references in this book. By virtue of their invaluable efforts for their results, 
I have been able to find a lot of theorems. 
\bigskip
\parno

\bigskip
\par\noindent
{\bf Notations.} 
(1) For a log (formal) scheme $X$ in the sense of 
Fontaine-Illusie-Kato (\cite{klog1}, \cite{klog2}), 
we denote by $\os{\circ}{X}$ (resp.~$M_X:=(M_X,\al_X)$)
the underlying (formal) scheme (resp.~the log structure) of $X$. 
In this book we consider the log structure 
on the Zariski site on $\os{\circ}{X}$. 
For a point $\os{\circ}{x}$ of $\os{\circ}{X}$, 
we endow $\os{\circ}{x}$ with the inverse image of $M_X$ 
and we denote it by $x$ and call it an exact closed point of $X$. 
We denote the stalk 
${\cal O}_{\os{\circ}{X},\os{\circ}{x}}$ by ${\cal O}_{X,x}$ for simplicity of notation. 
For a morphism $f\col X\lo Y$ of log (formal) schemes, 
$\os{\circ}{f}$ denotes the underlying morphism 
$\os{\circ}{X} \lo \os{\circ}{Y}$ of $f$. 
\par 
(2) For morphisms $X\lo S$ and $T\lo S$ of fine log (formal) schemes, 
we denote by $\os{\circ}{X}_{T}$
the underlying (formal) scheme of the fine log (formal) scheme 
$X_T:=X\times_ST$ for simplicity of notation. 
\par
(3) For a (formal) scheme $T$ 
and a commutative monoid $P$ with 
unit element $e$, we denote by 
$P\oplus {\cal O}_T^*$ a log structure on $T$ which is  
obtained by a morphism  
with natural morphism 
$P\oplus {\cal O}_T^* \owns (x,u)\lom 0 \in 
{\cal O}_T$ $(x\not= e)$ and $(e,u)\lom u$ 
of sheaves of monoids in $T_{\rm zar}$. 
This is an example of a constant hollow log structure defined in \cite[Definition 4]{ollc}. 
\par
(4) For a commutative monoid $P$ with unit element and 
for a commutative ring $A$ with unit element, 
${\rm Spec}^{\log}(A[P])$ is, by definition, the log scheme 
whose underlying scheme is ${\rm Spec}(A[P])$ 
and whose log structure is the association of the natural 
inclusion $P \os{\sus}{\lo} A[P]$. 
If $A$ has an $I$-adic topology ($I$ is an ideal of $A$), then 
$A\{P\}$ denotes $\vpl_n(A/I^n[P])$ and 
${\rm Spf}^{\log}(A\{P\})$ denotes the log formal scheme 
whose underlying (formal) scheme is ${\rm Spf}(\vpl_{n}(A/I^n[P]))$ 
and whose log structure is the association of the natural 
inclusion $P \os{\sus}{\lo} A\{P\}$.
\par
(5) 
For a morphism $Y \lo T$ of log (formal) schemes, we denote by 
${\Om}^i_{Y/T}$ ($=\om^i_{Y/T}$ in 
\cite{klog1}) $(i\in {\mab N})$ 
the sheaf of  relative logarithmic differential 
forms on $Y/T$ of degree $i$. 
For a log point $s$ whose underlying scheme is 
the spectrum of a perfect field of characteristic $p>0$ and 
for a log smooth morphism $Y\lo s$ of Cartier type, 
${\cal W}_n{\Om}^i_Y(=W_n\om^i_{Y}$ in 
\cite{hk} and \cite{msemi}) $(n\in {\mab Z}_{\geq 1})$ 
denotes the de Rham-Witt sheaf of logarithmic differential 
forms on $Y/s$ of degree $i$.
\par 
(6) SNCL=simple normal crossing log, 
SNCD=simple normal crossing divisor.
\par 
(7) For a module $M$ over a commutative ring $A$ 
with unit element and for a commutative $A$-algebra 
$B$ with unit element, $M_B$ denotes the tensor product 
$M\us{A}{\otimes}B$. 
\par 
(8) For a formal scheme $T$, 
${\cal K}_T$ denotes ${\cal O}_T\otimes_{\mab Z}{\mab Q}$. 
\par 
(9) For a log (formal) scheme $T$ and 
a quasi-coherent ideal sheaf ${\cal J}$ of 
${\cal O}_T$, we denote by 
$\ul{\rm Spec}^{\log}_T({\cal O}_T/{\cal J})$ 
(resp.~
$\ul{\rm Spec}^{\log}_T({\cal O}_T/{\cal J})_{\rm red}$) 
the log scheme whose underlying scheme is 
$\ul{\rm Spec}_{\os{\circ}{T}}({\cal O}_T/{\cal J})$ 
(resp.~
$\ul{\rm Spec}_{\os{\circ}{T}}({\cal O}_T/{\cal J})_{\rm red})$
and whose log structure is the inverse image of 
the log structure of $T$. 
\par 
(10) By following \cite{of}, 
for a fine log formal scheme $T$ with $p$-adic topology, 
we denote by $T_1$  and $T_0$ the log schemes 
$\ul{\rm Spec}^{\log}_T({\cal O}_T/p{\cal O}_T)$ 
and $\ul{\rm Spec}^{\log}_T({\cal O}_T/p{\cal O}_T)_{\rm red}$, respectively.
\par 
(11) Let $(T,{\cal J},\del)$ be a fine log PD-scheme such that 
${\cal J}$ is a quasi-coherent ideal sheaf of ${\cal O}_T$. 
Set $T_0:=\ul{\rm Spec}^{\log}_T({\cal O}_T/{\cal J})$.   
Let $g\col Y\lo T_0$ be 
a morphism of fine log schemes. 
In the previous books \cite{nh2} and \cite{nh3}, 
following \cite{klog1},  
we denoted the log crystalline site 
(resp.~the log crystalline topos) of $Y/(T,{\cal J},\del)$
by $(Y/T)^{\log}_{\rm crys}$ 
(resp.~$(\wt{Y/T})^{\log}_{\rm crys}$). 
In this book, following \cite{bb} and \cite{bob}, 
we denote the log crystalline site 
(resp.~the log crystalline topos) 
of $Y/(T,{\cal J},\del)$ by ${\rm Crys}(Y/T)$ 
(resp.~$(Y/T)_{\rm crys}$). 
For a (formal) scheme $Z$, we denote 
by ${\rm Zar}(Z)$ the zariski site on 
$Z$ and by $Z_{\rm zar}$ the associated topos to ${\rm Zar}(Z)$. 
For a log (formal) scheme $X$, we denote $\os{\circ}{X}_{\rm zar}$  
by $X_{\rm zar}$ for simplicity of notation. 
\par 
(12) For a ringed topos $({\cal T},{\cal A})$, 
let ${\rm C}^+({\cal A})$ 
(resp.~${\rm D}^+({\cal A})$) be 
the category of bounded below 
complexes of ${\cal A}$-modules 
(resp.~the derived category of bounded below 
complexes of ${\cal A}$-modules); 
let 
${\rm C}^+{\rm F}({\cal A})$ 
(resp.~${\rm D}^+{\rm F}({\cal A})$) 
be the category of 
bounded below filtered complexes of ${\cal A}$-modules 
(resp.~the derived category of 
bounded below filtered complexes of 
${\cal A}$-modules).  
For a bounded above filtered complex $(E^{\bul},P)\in {\rm D}^+{\rm F}({\cal A})$, 
we sometimes denote $(E^{\bul},P)\otimes^L_{\mab Z}{\mab Q}$ 
by $(E^{\bul},P)_{\mab Q}$. 
\par 
(13) For a short exact sequence 
$$0\lo (E^{\bul},d^{\bul}_E) \os{f}{\lo} 
(F^{\bul},d^{\bul}_F) \os{g}{\lo} 
(G^{\bul},d^{\bul}_G) \lo 0$$
of bounded below complexes of objects of ${\rm C}^+({\cal A})$, 
let ${\rm MC}(f):=(E^{\bul}[1],d^{\bul}_E[1])\oplus 
(F^{\bul},d^{\bul}_F)$ be the mapping cone of $f$.
We fix an isomorphism 
``$(E^{\bul}[1],d^{\bul}_E[1])\oplus 
(F^{\bul},d^{\bul}_F)\owns (x,y) \lom g(y)\in 
(G^{\bul},d^{\bul}_G)$'' 
in the derived category ${\rm D}^+({\cal A})$.
\par
Let ${\rm MF}(g):=(F^{\bul},d^{\bul}_F)\oplus 
(G^{\bul}[-1],d^{\bul}_G[-1])$ be 
the mapping fiber of $g$. 
We fix an isomorphism 
``$(E^{\bul},d^{\bul}_E)\owns x \lom (f(x),0)\in 
(F^{\bul},d^{\bul}_F)\oplus 
(G^{\bul}[-1],d^{\bul}_G[-1])$'' 
in the derived category ${\rm D}^+({\cal A})$.
\par
(14) For a morphism 
$f\col (E^{\bul},d^{\bul}_E)\lo (F^{\bul},d^{\bul}_F)$ 
of complexes, we identify 
${\rm MF}(f):=(E^{\bul},d^{\bul}_E)\oplus 
(F^{\bul}[-1],d^{\bul}_F[-1])$ and 
${\rm MC}(f)[-1]:=((E^{\bul},d^{\bul}_E)[1]\oplus 
(F^{\bul},d^{\bul}_F))[-1]$ 
with the following isomorphism 
``${\rm MF}(f) \owns (x,y)\lom (-x,y)\in   
{\rm MC}(f)[-1]$''.   
\bigskip
\par\noindent
{\bf Conventions.}
We omit the second ``log'' in the terminology  
``log smooth log (formal) scheme''. 
Following \cite{ollc}, we say that a morphism $X \lo Y$ 
of log (formal) schemes is {\it solid} if the log 
structure of $X$ is the inverse image of that of $Y$. 
It is obvious that this is equivalent to $X=Y\times_{\os{\circ}{Y}}\os{\circ}{X}$ 
in the category of log schemes. 

\chapter{Weight filtrations on log crystalline cohomology sheaves} 

\parno 
In this chapter we give the definition of an SNCL(=simple normal crossing log) scheme 
over a family of log points. We construct the $p$-adic filtered Steenbrink complex of 
a truncated simplicial SNCL scheme over a family of log points 
by using the modified $P$-filtered(=Poincar\'{e} filtered) log crystalline complex of it. 
We assume that the truncated simplicial SNCL scheme 
has the disjoint union of the member of 
an affine truncated simplicial open covering. 
In order to calculate the graded complex of the $p$-adic filtered Steenbrink complex, 
we calculate the graded complex of the modified $P$-filtered crystalline complex. 
We also define the $p$-adic monodromy operator 
(resp.~the $p$-adic quasi-monodromy operator)
of the log crystalline complex (resp.~the $p$-adic Steenbrink complex) of 
the truncated simplicial SNCL scheme 
and formulate the variational $p$-adic monodromy-weight conjecture 
and the variational $p$-adic filtered log hard Lefschetz conjecture.

\section{SNCL schemes}\label{sec:snclv} 
In this section we generalize the definition of 
a (formal) SNCL scheme in \cite[\S2]{nlk3}. 
\par 
Let $r$ be a fixed positive integer. 
Let $T$ be a log (formal) scheme whose log structure 
is ${\mab N}^r\oplus {\cal O}_T^*$ with a structural morphism  
${\mab N}^r\oplus {\cal O}_T^* \lo {\cal O}_T$ 
defined by 
$(x,u)\lom 0$ for $x\not=(0,\ldots,0)$
$((x,u)\in {\mab N}^r\oplus {\cal O}_T^*)$
and $((0,\ldots,0),u)\lom u$.   
We call this log structure the {\it free hollow log structure of rank} $r$ 
on $\os{\circ}{T}$. This is a special case of 
the constant hollow log structure defined in \cite[Definition 4]{ollc}. 
(Recall that the log structure $M_X=(M_X,\al_X)$ of 
a log (formal) scheme $X$ is hollow  
if $M_X\setminus {\cal O}_X^*$ 
is mapped to $0$ in ${\cal O}_X$ 
by the structural morphism $\al_X\col M_X\lo {\cal O}_X$ 
and that $X$ is constant if $M_X/{\cal O}_X^*$ is a locally constant sheaf.) 
Obviously the free hollow log structure of rank $r$  is split, i.~e.,  
the morphism ${\mab N}^r\lo M_T/{\cal O}_T^*$ is an isomorphism. 
\par   
More generally, we consider a fine log (formal) scheme $T$ 
such that there exists an open covering 
$T=\bigcup_{i\in I}T_i$ such that 
$M_{T_i}$ is a free hollow log structure of rank $r$.
We call $M_T$ a {\it locally free hollow log structure of rank} $r$ on $\os{\circ}{T}$ 
and we call $T$ a {\it $($formal$)$ family of log points of virtual dimension} $r$. 
If $r=1$, then we call $T$ a {\it $($formal$)$ family of log points},  
which we often denote by $S$ in this book. 
When $(\os{\circ}{T},{\cal J},\del)$ is 
a PD-scheme with quasi-coherent PD-ideal sheaf 
and PD-structure, 
we call $(T,{\cal J},\del)$ 
a {\it PD-family of log points of virtual dimension} $r$. 
When $r=1$, we call $(T,{\cal J},\del)$ a {\it PD-family of log points}.  

\begin{exem}\label{exem:flp}
Let $T$ be any scheme and 
let us us consider ${\mab P}^1_T$. 
Set $U_0:=\ul{\rm Spec}_T({\cal O}_T[t])$ 
and $U_1:=\ul{\rm Spec}_T({\cal O}_T[t^{-1}])$ and 
let $\{U_0,U_1\}$ be a standard open covering of 
${\mab P}^1_T$. 
Consider the following two charts 
${\mab N}\owns 1\lom 0\in {\cal O}_{U_0}$ 
and 
${\mab N}\owns 1\lom 0\in {\cal O}_{U_1}$.  
Let $e_0$ and $e_1$ be the standard bases of 
the first ${\mab N}$ and the the second ${\mab N}$, respectively. 
Patch $e_0^{\mab N}\vert_{U_0\cap U_1}$ and 
$e_1^{\mab N}\vert_{U_0\cap U_1}$ on $U_0\cap U_1$ by 
the following equality $e_1=te_0$. 
Then we have the family of log points whose underlying scheme is ${\mab P}^1_T$. 
\end{exem}

\par 
For a log (formal) scheme $X=(X,(M_X,\al_X))$, we call a subsheaf 
$N$ of $M_X$ such that $\al_X\vert_N$ induces an isomorphism 
$(\al_X\vert_N)^{-1}({\cal O}_X^*)\os{\sim}{\lo} {\cal O}_X^*$ 
a {\it sub log structure} of $X$.

\begin{prop}\label{prop:nfi}
Let $X$ and $Y$ be $($not necessarily fine$)$ log $($formal$)$ schemes. 
There exists the following bijection of sets$:$
\begin{align*} 
F\col & \{f\in{\rm Hom}(X,Y)~\vert~f^*_x \col (M_Y/{\cal O}_Y^*)_{f(x)}\lo
(M_X/{\cal O}_X^*)_x~{\rm is}~{\rm injective}\}\lo \tag{1.1.2.1}\label{ali:hxy} \\
&\{(g,\gam)\in {\rm Hom}(\os{\circ}{X},\os{\circ}{Y})\times 
{\rm Hom}(g^*((M_Y,\al_Y)),(M_X,\al_X))
~\vert~
\gam~{\rm is}~{\rm an}~{\rm isomorphism}~\\
&\quad {\rm from}~g^*((M_Y,\al_Y))~
{\rm to}~{\rm a}~{\rm sub}~{\rm log}~{\rm structure}~{\rm of}~(M_X,\al_X)\}. 
\end{align*}
\end{prop}
\begin{proof} 
Let $f$ be an element of the source in (\ref{ali:hxy}). 
Then $I:={\rm Im}(f^*(M_Y)\lo M_X)$ 
with structural morphism $\al_X\vert_I$ 
defines a sub log structure of $(M_X,\al_X)$, 
which is isomorphic to $f^*(M_Y)$.  
Indeed, because $f^*(M_Y)/{\cal O}_X^*=f^{-1}(M_Y/{\cal O}_Y^*)$ (\cite[(1.4.1)]{klog1}), 
the morphism $f^*(M_Y)\lo M_X$ over ${\cal O}_X$ is injective and we have an isomorphism 
$\gam \col f^*(M_Y)\os{\sim}{\lo} I$ over ${\cal O}_X$. 
Set $F(f):=(\os{\circ}{f},\gam)$. 
\par 
Conversely, let $(g,\gam)$ be an element of the target in (\ref{ali:hxy}). 
Then we have the following commutative diagram
\begin{equation*} 
\begin{CD} 
g^*(M_Y)@>{\gam,\sim}>> N\\
@V{g^*(\al_Y)}VV @VV{\al_X\vert_N}V\\
{\cal O}_X@={\cal O}_X, 
\end{CD}
\end{equation*} 
where $N$ is a sub log structure of $(M_X,\al_X)$. 
Hence the composite morphism 
$$g^*(M_Y)/{\cal O}_X^*\os{\sim}{\lo} N/{\cal O}_X^*
\os{\sus}{\lo} M_X/{\cal O}_X^*$$ 
is injective. Therefore we have an element of the source in  (\ref{ali:hxy}). 
This map is the inverse of $F$. 
\end{proof}

\begin{prop}\label{prop:lbl}
Let $T$ and $T'$ be hollow log $($formal$)$ schemes 
such that $M_T/{\cal O}^*_T=:P$ and $M_{T'}/{\cal O}^*_{T'}=:P'$ are 
constant sheaves of monoids defined 
by commutative monoids $P$ and $P'$ with unit elements, respectively.  
Let $\star$ be nothing or $\prime$. 
Assume that 
the sequence
\begin{align*} 
1\lo {\cal O}_{T^{\star}}^*\lo M_{T^{\star}}\lo M_{T^{\star}}/{\cal O}^*_{T^{\star}}\lo 1
\tag{1.1.3.1}\label{ali:otmt} 
\end{align*} 
is split. 
Let $e_P$ and $e_{P'}$ be the unit elements of $P$ and $P'$, respectively. 
Let $u\col \os{\circ}{T}\lo \os{\circ}{T}{}'$ be a morphism of $($formal$)$ schemes.  
Let $g\col P'\lo P$ be a morphism of monoids 
such that $g(m)\not=e_P$ for $\forall m\not=e_{P'}$. 
Then there exists a morphism $v\col T\lo T'$ 
such that $\os{\circ}{v}=u$ and the following pull-back morphism 
$$v^*_x\col P'=(M_{T'}/{\cal O}_{T'}^*)_{v(x)}\lo 
(M_T/{\cal O}_T^*)_{x}=P$$ 
is equal to $g$ for any point $x\in \os{\circ}{T}$. 
\end{prop}
\begin{proof} 
The structural morphism $\al_{T^{\star}}\col M_{T^{\star}}\lo {\cal O}_{T^{\star}}$ 
is identified with $P^{\star}\oplus {\cal O}_{T^{\star}}^*\lo {\cal O}_{T^{\star}}$ 
by using a splitting 
$P^{\star}=M_{T^{\star}}/{\cal O}_{T^{\star}}^*\lo M_{T^{\star}}$ 
of (\ref{ali:otmt}).  
Then we have the following commutative diagram: 
\begin{equation*} 
\begin{CD}
P\oplus u_*({\cal O}_T)^*@>>> u_*({\cal O}_T)\\
@A{g\oplus u^*}AA @AA{u^*}A \\
P'\oplus {\cal O}_{T'}^*@>>> {\cal O}_{T'}. 
\end{CD}
\end{equation*} 
This diagram gives the morphism $v$ which we want. 
\end{proof}


\begin{prop}\label{prop:oot} 
$(1)$ Let $f\col X\lo Y$ be a morphism of log $($formal$)$ schemes.  
Assume that $Y$ is a hollow log $($formal$)$ scheme. 
Then, for any point $x\in \os{\circ}{X}$, 
the image of any nonzero element of $(M_Y/{\cal O}^*_Y)_{f(x)}$ by 
the morphism $(M_Y/{\cal O}^*_Y)_{f(x)}\lo (M_X/{\cal O}^*_X)_x$ 
is a nonzero element of $(M_X/{\cal O}^*_X)_x$. 
\par 
$(2)$ Furthermore, assume that $X$ is integral and that 
$Y$ is a $($formal$)$ family of log poitns. 
Then the morphism $f^*(M_Y)\lo M_X$ is injective. 
\end{prop}
\begin{proof} 
(1): Obvious. 
\par 
(2): This is a local question. We may assume 
that $Y$ has a global chart 
${\mab N}\owns 1\lom 0\in {\cal O}_Y$. 
Let $a$ be the image of $1\in {\mab N}=f^*(M_Y/{\cal O}_Y^*)$ 
in $M_X/{\cal O}^*_X$. 
Assume that $a^n=a^{n'}$ with $n'> n$. Since $M_X$ is integral, 
$a^{n'-n}=1$.  
By the isomorphism $\al_X^{-1}({\cal O}_X^*)\os{\sim}{\lo} {\cal O}_X^*$, 
we see that $a\in {\cal O}_X^*$. 
Since the composite morphism 
$f^{-1}(M_Y)\lo M_X\lo {\cal O}_X$ is equal to 
the composite morphism 
$f^{-1}(M_Y)\lo f^{-1}({\cal O}_Y)\lo {\cal O}_X$, 
the image of $1\in f^{-1}(M_Y)$ in ${\cal O}_X$ is $0$. 
This is a contradiction. 
\end{proof} 

\par 
Let $T=\bigcup_{i\in I}T_i$ be an open covering 
of $T$ such that 
$M_{T_i}/{\cal O}^*_{T_i}\simeq {\mab N}^{r}$.  
Take local sections 
$\tau_{i,1},\ldots,\tau_{i,r}\in \Gam(T_i,M_T)$ such that 
the images of $\tau_{i,1},\ldots,\tau_{i,r}$ in  
$\Gam(T_i,M_T/{\cal O}^*_T)$
is a system of generators. 
In \cite[p.~47]{nh3} we have remarked that 
${\rm Aut}({\mab N}^r)\simeq {\mathfrak S}_r$, 
which is easy to prove. Set $T_{ij}:= T_i\cap T_j$.  
Then, by this remark, 
there exists a unique permutation matrix 
$P_{ji}\in {\rm GL}_r({\mab Z})$ 
and unique sections 
$u_{ji,1},\ldots, u_{ji,r} \in 
\Gam(T_{ij},{\cal O}^*_T)$ such that 
\begin{equation*} 
(\tau_{j,1},\ldots,\tau_{j,r})\vert_{T_{ij}}=
{\rm diag}(u_{ji,1},\ldots, u_{ji,r})
(\tau_{i,1},\ldots,\tau_{i,r})\vert_{T_{ij}}P_{ji}
\tag{1.1.4.1}\label{eqn:tkjd} 
\end{equation*} 
in $\Gam(T_i\cap T_j,M_T)$. 
First consider the case 
where $\os{\circ}{T}$ is a scheme (not a formal scheme). 
Consider the scheme 
${\mab A}^r_{\os{\circ}{T}_i}=
\ul{\rm Spec}_{\os{\circ}{T}_i}
({\cal O}_{T_i}[t_{i,1}\ldots,t_{i,r}])$  
and the log scheme 
$({\mab A}^r_{\os{\circ}{T}_i}, 
({\mab N}^r\owns (0,\ldots,0,\os{m}{1},0,\ldots,0) 
\lom t_{i,m} \in {\cal O}_{T_i}[t_{i,1},\ldots,t_{i,r}])^a)$ 
$(1\leq m \leq r)$. 
Denote this log scheme by $\ol{T}_i$. 
Then, by patching $\ol{T}_i$ and $\ol{T}_j$ along 
$\ol{T}_{ij}:=\ol{T}_i\cap \ol{T}_j$ by the following equation 
$$(t_{j,1},\ldots,t_{j,r})\vert_{\ol{T}_{ij}}=
{\rm diag}(u_{ji,1},\ldots, u_{ji,r})
(t_{i,1},\ldots,t_{i,r})\vert_{\ol{T}_{ij}}P_{ji},$$ 
we have the log scheme $\ol{T}=\bigcup_{i\in I}\ol{T}_i$. 
The ideal sheaves
$\{(t_{i,1},\ldots,t_{i,r}){\cal O}_{\ol{T}_i}\}_{i\in I}$ 
patch together  
and we denote by 
${\cal I}_{\ol{T}}$ 
the resulting ideal sheaf of 
${\cal O}_{\ol{T}}$. 
The isomorphism class of the log scheme $\ol{T}$ 
is independent of the choice of the system of 
generators $\tau_{i,1},\ldots,\tau_{i,r}$'s. 
By using this easy fact and by considering 
the refinement of two open coverings of $T$, 
we see that the isomorphism class
of the log scheme $\ol{T}$ and the ideal sheaf 
${\cal I}_{\ol{T}}$ 
are also independent of 
the choice of the open covering $T=\bigcup_{i\in I}T_i$. 
Because 
$\os{\circ}{\ol{T}}_i
=\os{\circ}{T}_i\otimes_{{\mab Z}}{\mab Z}[{\mab N}^r]$, 
the natural morphism $\ol{T} \lo \os{\circ}{T}$ is log smooth 
by the criterion of the log smoothness (\cite[(3.5)]{klog1}). 
For a log scheme $Y$ over $\ol{T}$, 
we denote 
${\cal I}_{\ol{T}}\otimes_{{\cal O}_{\ol{T}}}{\cal O}_Y$
by ${\cal I}_Y$ by abuse of notation. 
If $\os{\circ}{Y}$ is flat over $\os{\circ}{\ol{T}}$, 
then ${\cal I}_Y$ can be considered as 
an ideal sheaf of ${\cal O}_Y$. 
When $\os{\circ}{T}$ is a formal scheme with a system of ideals of definition 
$\{{\cal J}_{\lam}\}$,  
we obtain the analogue of the object $\ol{T}$ above 
by considering the formal scheme 
$\ul{\rm Spf}_{\os{\circ}{T}_i}
(\vpl_{\lam}({\cal O}_{T_i}/{\cal J}_{\lam}\vert_{T_i}[t_{i,1},\ldots,t_{i,r}]))$ 
instead of ${\mab A}^r_{\os{\circ}{T}_i}$ above.  
We denote this analogue by $\ol{T}$ again. 
For a notational reason in \S\ref{sec:ldrwii} below, we denote 
$t_{i,m}$ by $\tau_{i,m}$ from now on by abuse of notation. 
\par 
By killing ``$\tau_i$'''s, 
we have the following natural exact closed immersion 
\begin{equation*} 
T\os{\sus}{\lo} \ol{T}  
\tag{1.1.4.2}\label{eqn:stas}
\end{equation*} 
over $\os{\circ}{T}$. 

\begin{prop}\label{prop:bar}
Let $v\col T\lo T'$ be a morphism of families of 
log points of virtual dimensions $r$ and $r'$, respectively. 
Then $v$ induces a unique morphism 
$\ol{v}\col \ol{T}\lo \ol{T}{}'$ 
fitting into the following commutative diagram 
\begin{equation*} 
\begin{CD}
T@>{v}>> T'\\ 
@V{\bigcap}VV @VV{\bigcap}V \\
\ol{T} @>{\ol{v}}>> \ol{T}{}'. 
\end{CD}
\end{equation*} 
\end{prop} 
\begin{proof}
This is clear from the definition of $\ol{T}$. 
\end{proof}

\begin{lemm}\label{lemm:etl1}
Let $U$ be a fine log scheme. 
Let $U_0\os{\sus}{\lo} U$ be an exact immersion. 
Let $Z_0\lo U_0$ be a morphism of fine log schemes  
such that the composite morphism 
$Z_0\lo U_0\os{\sus}{\lo} U$ 
has a global chart $Q\lo P$ such that 
the morphism 
$Z_0 \lo 
U_0\times_{{\rm Spec}^{\log}({\mab Z}[Q])}{\rm Spec}^{\log}({\mab Z}[P])$ 
is solid and \'{e}tale. 
Then, Zariski locally on $\os{\circ}{Z}_0$, 
there exists a log scheme $Z$ over $U$ fitting 
into the following cartesian diagram 
\begin{equation*}
\begin{CD} 
Z_0 @>{\sus}>> Z \\
@VVV @VVV \\
U_0\times_{{\rm Spec}^{\log}({\mab Z}[Q])}{\rm Spec}^{\log}({\mab Z}[P])
@>{\sus}>> 
U\times_{{\rm Spec}^{\log}({\mab Z}[Q])}{\rm Spec}^{\log}({\mab Z}[P]) \\ 
@VVV @VVV \\
U_0@>{\sus}>> U,   
\end{CD} 
\tag{1.1.6.1}\label{cd:xwtx} 
\end{equation*}
where the vertical morphism 
$Z \lo U\times_{{\rm Spec}^{\log}({\mab Z}[Q])}{\rm Spec}^{\log}({\mab Z}[P])$ 
is solid and \'{e}tale.  
\end{lemm}
\begin{proof} 
(1): By \cite[I (8.1)]{sga1}  
we may assume that there exists 
an \'{e}tale scheme $\os{\circ}{Z}$ over 
$(U\times_{{\rm Spec}^{\log}({\mab Z}[Q])}{\rm Spec}^{\log}({\mab Z}[P]))^{\circ}$ 
such that 
$\os{\circ}{Z}_0=\os{\circ}{Z}\times_{\os{\circ}{U}}\os{\circ}{U}_0$. 
Endow $\os{\circ}{Z}$ with the inverse image of the log structure of 
$U\times_{{\rm Spec}^{\log}({\mab Z}[Q])}{\rm Spec}^{\log}({\mab Z}[P])$. 
The resulting log scheme $Z$ is a desired log scheme. 
\end{proof} 

\begin{prop}\label{prop:ebtl}
Let $U$ be a fine log scheme. 
Let $U_0 \os{\sus}{\lo} U$ be an exact immersion.  
Let $Z_0$ be a log smooth scheme over $U_0$. 
Let $z$ be a point of $\os{\circ}{Z}_0$. 
Assume that there exists a chart 
$Q\os{\sus}{\lo} P$ of $Z_0/U$ around $z$ such that 
$\sharp  (P^{\rm gp}/Q^{\rm gp})_{\rm tor}$ is invertible on $Z_0$ 
and such that ${\cal O}_{Z_0,z}
\otimes_{\mab Z}(P^{\rm gp}/Q^{\rm gp}) 
\simeq {\Om}^1_{Z_0/U_0,z}$.  
Then, on a neighborhood of $z$, there exists a log smooth scheme $Z/U$ 
fitting into the following cartesian diagram$:$
\begin{equation*}
\begin{CD} 
Z_0 @>{\sus}>> Z \\
@VVV @VVV \\
U_0@>{\sus}>> U. 
\end{CD} 
\tag{1.1.7.1}\label{cd:ytx} 
\end{equation*}
\end{prop} 
\begin{proof} 
We may assume that $\os{\circ}{Z}_0$ is affine. 
By the argument in \cite[(3.5)]{klog1} and 
by the proof of \cite[(2.1.4)]{nh2}, we may assume 
that the morphism $Z_0\lo U_0$ factors through 
a strict \'{e}tale morphism 
$Z_0 \lo U_0\times_{{\rm Spec}^{\log}({\mab Z}[Q])}{\rm Spec}^{\log}({\mab Z}[P])$.  
Now (\ref{prop:ebtl}) follows from (\ref{lemm:etl1}) and 
the criterion of the log smoothness (\cite[(3.5)]{klog1}).  
\end{proof}

\begin{rema}\label{rema:colsm}
Consider 
log structures of schemes in the \'{e}tale topos in this remark. 
Then, by \cite[(3.5)]{klog1} and \cite[I (8.1)]{sga1}, the following holds. 
\par 
Let $U_0 \os{\sus}{\lo} U$ be an exact immersion.  
Let $V_0$ be a log smooth scheme over $U_0$. 
Let $v$ be a point of $\os{\circ}{V}_0$. 
Then there exist an \'{e}tale neighborhood $W_0$ of $v$ and 
a log smooth scheme $W/U$ fitting into the following cartesian diagram$:$  
\begin{equation*}
\begin{CD} 
W_0 @>{\sus}>> W \\
@VVV @VVV \\
U_0@>{\sus}>> U. 
\end{CD} 
\tag{1.1.8.1}\label{cd:eytx} 
\end{equation*} 
\end{rema}

\par 
Let $B$ be a scheme. 
For two nonnegative integers $a$ and $d$ such that 
$a\leq d$,  
consider the following scheme 
\begin{equation*} 
\os{\circ}{\mab A}_B(a,d):=
\ul{{\rm Spec}}_B
({\cal O}_B[x_0, \ldots, x_d]/(\prod_{i=0}^ax_i)). 
\end{equation*}

\begin{defi}\label{defi:zbsnc}
Let $Z$ be a scheme over $B$ 
with structural morphism $g \col Z \lo B$. 
We call $Z$ an 
{\it SNC$($=simple normal crossing$)$ scheme} over 
$B$ if $Z$ is a union of smooth schemes 
$\{Z_{\lam}\}_{\lam \in \Lam}$ over $B$ ($\Lam$ is a set) 
and if, for any point of $z \in Z$, 
there exist an open neighborhood $V$ of 
$z$ and an open neighborhood $W$ of 
$g(z)$ such that  
there exists an \'{e}tale morphism  
$\pi \col V \lo \os{\circ}{\mab A}_W(a,d)$ 
such that 
\begin{equation*} 
\{Z_{\lam}\vert_{V}\}_{\lam \in \Lam} 
=\{\pi^*(x_i)=0\}_{i=0}^a, 
\tag{1.1.9.1}\label{eqn:defilsnc}
\end{equation*}   
where $a$ and $d$ are nonnegative integers 
such that $a\leq d$,  
which depend on 
zariskian local neighborhoods in $Z$, and 
$Z_{\lam}\vert _V:=Z_{\lam}\cap V$.  
The strict meaning of the equality (\ref{eqn:defilsnc})
is as follows.  
There exists a subset $\Lam(V)$ of $\Lam$ such that
there exists a bijection $i\col \Lam(V) \lo \{0,1,\ldots,a\}$ 
such that, if 
$\lam \not\in \Lam(V)$, then 
$Z_{\lam}\vert_{V}=\emptyset$ 
and if $\lam \in \Lam(V)$, then 
$Z_{\lam}$ is a closed subscheme defined by 
the ideal sheaf $(\pi^*(x_{i(\lam)}))$.  
We call the set $\{Z_{\lam}\}_{\lam \in \Lam}$ 
a {\it decomposition of $Z$ by smooth components of} 
$Z$ over $B$. In this case, we call $Z_{\lam}$ a 
{\it smooth component} of $Z$ over $B$. 
\end{defi} 

\parno  
The base change of an SNC scheme is an SNC scheme. 
In the open log case in \cite{nh2}, 
we have been able to find a condition in [loc.~cit., (2.1.7)]
eliminating the corresponding condition to 
(\ref{eqn:defilsnc}). 
However we do not know whether it is possible to do so 
in (\ref{defi:zbsnc}). 
\par 
Set $\Del:=\{Z_{\lam}\}_{\lam \in \Lam}$.  
For an open subscheme $V$ of $Z$, 
set $\Del_{V}:=\{Z_{\lam}\vert_{V}\}_{\lam \in \Lam}$.

\par 
In \cite{nh2} we have proved  the following 
``local uniqueness of smooth components''  
(see \cite[(2.1.7)]{nh2} for the definition of a relative SNCD):

\begin{prop}[{\rm {\bf \cite[Proposition A.0.1]{nh2}}}]\label{prop:mainsncd}
Let $B$ be a scheme. 
Let $g\col (X,D) \lo B$ be a smooth scheme 
with a relative SNCD over $B$. 
Let $\Gam := \{D_{\nu}\}_{\nu \in N}$ be 
a  decomposition of $D$ by smooth components of $D$, 
that is, each $D_{\nu}$ is smooth over $B$ and 
$D = \sum_{\nu \in N}D_{\nu}$ in the monoid of  
effective Cartier divisors on $X/B$ {\rm (\cite[(2.1.7)]{nh2})}.  
Let $z$ be a point of $D$ and assume that 
we are given a cartesian diagram 
\begin{equation*}
\begin{CD}
D @>{\subset}>> X\\ 
@V{}VV  @VV{g}V \\
\ul{\rm Spec}_B({\cal O}_B[y_1, \ldots, y_e]/(y_1\cdots y_s)) 
@>{\sus}>> 
\ul{\rm Spec}_B({\cal O}_B[y_1, \ldots, y_e]) 
\end{CD}
\tag{1.1.10.1}\label{cd:dvs}
\end{equation*}
for some $s\leq e$ such that 
$z \in \bigcap_{i=1}^s\{g^*(y_i) = 0\}$. 
Then, by shrinking $X$, for any $1 \leq i \leq s$, 
there exists a unique element 
$\nu_i \in N$ satisfying $D_{\nu_i} = \{g^*(y_i)=0\}$. 
\end{prop}

\begin{prop}\label{prop:cc}
Let the notations be as in {\rm (\ref{defi:zbsnc})}. 
Let $\Del$ and  $\Del'$ be decompositions of $Z$  
by smooth components of $Z$. 
Then, for any point $z \in Z$, 
there exists an open neighborhood $V$ of 
$z$ in $Z$ such that $\Del_V = \Del'_V$. 
\end{prop}
\begin{proof} 
If $V$ is small enough, 
then we have an \'{e}tale morphism 
$\pi \col V \lo \os{\circ}{\mab A}_W(a,d)$ such that 
$\Del_V=\{\pi^*(x_i)=0\}_{i=0}^a$ in $V$.  
Then we have a natural closed immersion 
\begin{equation*} 
\os{\circ}{\mab A}_W(a,d) \os{\sus}{\lo} 
\ul{{\rm Spec}}_W({\cal O}_W[x_0, \ldots, x_d])=
\os{\circ}{\mab A}{}^{d+1}_W. 
\end{equation*} 
By \cite[I (8.1)]{sga1} 
we may assume that there exists 
an \'{e}tale scheme $U$ over 
$\os{\circ}{\mab A}{}^{d+1}_W$ such that 
$V=U\times_{\os{\circ}{\mab A}{}^{d+1}_W}
\os{\circ}{\mab A}_W(a,d)$. 
Then $V$ is a relative SNCD on $U$ over $W$; 
$\Del_V$ and $\Del'_V$ are decompositions of smooth 
components of $V$ in the sense of \cite[(2.18)]{nh2}. 
By (\ref{prop:mainsncd}) 
we can take a small $U$ such that $\Del_V=\Del'_V$. 
\end{proof}
For a nonnegative integer $k$, 
set 
\begin{equation}
Z_{\{\lam_0, \lam_1,\ldots \lam_k\}} 
:=Z_{\lam_0}\cap \cdots \cap Z_{\lam_k} \quad 
(\lam_i \not= \lam_j~{\rm if}~i\not= j) 
\tag{1.1.11.1}\label{eqn:parlm}
\end{equation}
and set
\begin{equation}
Z^{(k)} =  
\us{\{\lam_0, \ldots,  \lam_{k}~\vert~\lam_i 
\not= \lam_j~(i\not=j)\}}{\coprod}
Z_{\{\lam_0, \lam_1, \ldots, \lam_k\}}.   
\tag{1.1.11.2}\label{eqn:kfdintd}
\end{equation} 
For a negative integer $k$, we set 
$Z^{(k)}=\emptyset$.  

\begin{prop}\label{prop:zki}
The scheme $Z^{(k)}$ 
is independent of the choice of $\Del$.  
\end{prop} 
\begin{proof} 
As in the proof of \cite[(2.2.15)]{nh2}, by using (\ref{prop:cc}), 
we can show that $Z^{(k)}$ is independent of the choice of 
$\Del$.  
\end{proof}
We have the natural morphism 
$b^{(k)}\col Z^{(k)}\lo Z$.

As in \cite[(3.1.4)]{dh2} and \cite[(2.2.18)]{nh2}, 
we have an orientation sheaf $\vp^{(k)}_{\rm zar}(Z/B)$ 
$(k\in {\mab N})$ in $Z^{(k)}_{\rm zar}$ 
associated to the set $\Del$. 
If $B$ is a closed subscheme of $B'$ defined by 
a quasi-coherent nil-ideal sheaf ${\cal J}$ 
which has a PD-structure $\del$, 
then $\vp^{(k)}_{\rm zar}(Z/B)$ 
extends to an abelian sheaf 
$\vp^{(k)}_{\rm crys}(Z/B')$ in 
$(Z^{(k)}/(B',{\cal J},\del))_{\rm crys}$.

\begin{defi}
\label{defi:reft}
Let 
$\Del := \{Z_{\lam}\}_{\lam \in \Lam}$ and 
$\Del' := \{Z'_{\sig}\}_{\sig \in \Sig}$ be decompositions of $Z$ 
by smooth components of $Z/S$. 
Then we say that $\Del'$ is a {\it refinement} of $\Del$ if, 
for any $\lam \in \Lam$, there exists a subset 
$\Sig_{\lam}$ of $\Sig$ with 
$\Sig = \coprod_{\lam \in \Lam}\Sig_{\lam}$ 
such that $Z_{\lam} = \bigcup_{\sig \in \Sig_{\lam}} Z'_{\sig}$. 
\end{defi}

\begin{prop-defi}\label{prop-defi:rf} 
{\it If $Z$ is locally noetherian as a topological space,  
then there exists a unique decomposition of $Z$ 
by smooth components of $Z/S$ 
which is a refinement of any decomposition of $Z$ 
by its smooth components.  
We call this decomposition the {\rm finest decomposition} of $Z$}. 
\end{prop-defi}
\begin{proof}
The proof is the same as that of \cite[A.~0.3]{nh2}. 
The finest decomposition is constructed as follows. 
\par
Let the notations be as in (\ref{defi:zbsnc}). 
Let $Z_{\lam}=\coprod_{\sig \in \Sig_{\lam}}Z_{\sig}$ 
be the disjoint sum of $Z_{\lam}$ by the connected components 
of it. Then  
$\{Z_{\sig}\}_{\lam \in \Lam,\sig\in \Sig_{\lam}}$ 
is the finest decomposition of $Z$ by its smooth components. 
\end{proof}

\par 
In the case where $Z$ is not necessarily locally noetherian, 
we still have the following as in \cite[Proposition A.0.7]{nh2}: 

\begin{prop}
\label{prop:cofinal} 
Let $Z \lo B$ be as in {\rm (\ref{defi:zbsnc})}. 
Let $\Del= \{Z_{\lam}\}_{\lam \in \Lam}$ 
and $\Del'=\{Z'_{\lam'}\}_{\lam' \in \Lam'}$ 
be two decompositions of $Z$ by smooth components. 
Then there exists a 
decomposition $\Del''$ of $Z$ by smooth components 
which is a refinement of $\Del$ and $\Del'$. 
\end{prop}
\begin{proof}
There exists an open covering 
$Z = \bigcup_{j \in J}Z_j$ of $Z$ such that 
$\Del \vert_{Z_j}=\Del' \vert_{Z_j}$ ((\ref{prop:cc})). 
Set 
\begin{equation*} 
(Z_{\lam}\wedge Z'_{\lam'})\vert_{Z_j} :=
\begin{cases} 
Z_{\lam}\vert_{Z_j} & {\rm if}~Z_{\lam}\vert_{Z_j}=Z'_{\lam'}\vert_{Z_j}, \\
\emptyset & {\rm otherwise}.
\end{cases}
\end{equation*}  
Then it is clear that 
$\{(Z_{\lam}\wedge Z'_{\lam'})\vert_{Z_j}\}_j$ patch together and 
we indeed have the closed subscheme 
$Z_{\lam}\wedge Z'_{\lam'}$ whose restriction to $Z_j$ is 
$(Z_{\lam}\wedge Z'_{\lam'})\vert_{Z_j}$. 
(One can easily check that $Z_{\lam}\wedge Z'_{\lam'}$ is independent of 
the choice of the open covering above, though we do not use this fact in this book.) 
It is also clear that 
$Z_{\lam} \wedge Z'_{\lam'}$ is smooth over $B$ (possibly empty). 
We have the following equalities 
$$ Z_{\lam} = \bigcup_{ \lam' \in \Lam',  
Z_{\lam} \wedge Z'_{\lam'} \not= \emptyset} 
Z_{\lam} \wedge Z'_{\lam'}, 
\quad \quad Z'_{\lam'} = \bigcup_{
\lam \in \Lam, Z_{\lam} \wedge Z'_{\lam'} \not= \emptyset} 
Z_{\lam} \wedge Z'_{\lam'}. $$
Set $\Lam'' := \{(\lam,\lam') \in 
\Lam \times \Lam'~\vert~Z_{\lam} \wedge Z'_{\lam'} 
\not= \emptyset\}$.  
By using the equalities above, 
it is easy to see  that $\Del'' := \{Z_{\lam} \wedge 
Z'_{\lam'}\}_{(\lam,\lam') \in \Lam''}$ has a desired property. 
\end{proof}

In the rest of this section, consider the case $r=1$ 
unless stated otherwise. 
First assume that $M_S$ is the free hollow log structure. 
We fix an isomorphism 
\begin{align*} 
(M_S,\al_S)\simeq ({\mab N}\oplus {\cal O}_S^*\lo {\cal O}_S) 
\tag{1.1.15.1}\label{ali:nsos} 
\end{align*} 
globally on $S$. 
Let $n$ be a positive integer.  
Let $n:=n\times$ be the multiplication ${\mab N}\owns j \lom jn\in {\mab N}$ by $n$. 
Denote simply by 
${\mab N}^{{\oplus}(a+1)}\oplus_{{\mab N},n}{\mab N}$ 
the push-out of the following diagram of monoids: 
\begin{equation*}
\begin{CD} 
{\mab N}@>{n}>> {\mab N} \\
@V{{\rm diag.}}VV \\
{\mab N}^{{\oplus}(a+1)}@. @..
\end{CD} 
\end{equation*} 
Let $M_S(a,d,n)$ be the log structure on 
${\mab A}_{\os{\circ}{S}}(a,d)$ associated to the following morphism 
\begin{equation*} 
{\mab N}^{{\oplus}(a+1)}\oplus_{{\mab N},n}{\mab N}\owns 
(0, \ldots,0,\os{i}{1},0,\ldots, 0,j)\lom x_{i-1}\cdot 0^j \in 
{\cal O}_S[x_0, \ldots, x_d]/(\prod_{i=0}^ax_i).
\tag{1.1.15.2}  
\end{equation*} 
Here we set $0^0:=1$. 
Let ${\mab A}_S(a,d,n)$ be the resulting log scheme over $S$.
The morphism 
${\mab N}\owns j \lom (0, \ldots,0,\ldots, 0,j)\in 
{\mab N}^{{\oplus}(a+1)}\oplus_{{\mab N},n}{\mab N}$ 
induces a morphism ${\mab A}_S(a,d,n) \lo S$ of log schemes. 
Set $M_S(a,d):=M_S(a,d,1)$ and ${\mab A}_S(a,d):={\mab A}_S(a,d,1)$.   
We call $M_S(a,d,n)$ the $n$-{\it standard log structure} on ${\mab A}_{\os{\circ}{S}}(a,d)$. 
We call ${\mab A}_S(a,d)$ the {\it standard SNCL scheme} 
and the $1$-standard log structure the {\it standard log structure} simply. 
The following morphism 
\begin{align*}
(n,{\rm id})\col {\mab N}\oplus{\cal O}_S^*\lo {\mab N}\oplus{\cal O}_S^*
\end{align*} 
and the isomorphism (\ref{ali:nsos}) induces the following morphism 
\begin{align*} 
(M_S,\al_S)\lo (M_S,\al_S). 
\tag{1.1.15.3}\label{ali:nsois} 
\end{align*}  
Hence we have the following morphism 
\begin{align*} 
u_n \col S\lo S
\tag{1.1.15.4}\label{ali:nsofs} 
\end{align*}  
of fine log schemes. 
Then we have the following equality 
\begin{align*} 
{\mab A}_S(a,d,n)={\mab A}_S(a,d)\times_{S,u_n}S. 
\tag{1.1.15.5}\label{ali:nstes} 
\end{align*}  
Note that the morphism $u_n$ depends on 
the choice of the isomorphism (\ref{ali:nsos}).  
\par 
Let $S$ be a family of log points 
(such that $M_S$ is not necessarily free). 

\begin{defi}\label{defi:lfac}  
(1) 
Let $f \col X(=(\os{\circ}{X},M_X)) \lo S$ be a morphism of 
log schemes such that $\os{\circ}{X}$ is 
an SNC scheme over $\os{\circ}{S}$ 
with a decomposition 
$\Del:=\{\os{\circ}{X}_{\lam}\}_{\lam \in \Lam}$ 
of $\os{\circ}{X}/\os{\circ}{S}$ by its smooth components. 
We call $f$ (or $X/S$) an 
{\it RSNCL$($=ramified simple normal crossing log$)$ scheme} if, 
for any point of $x \in \os{\circ}{X}$, 
there exist an open neighborhood $\os{\circ}{V}$ of 
$x$ and an open neighborhood $\os{\circ}{W}$ of 
$\os{\circ}{f}(x)$ such that $M_W$ is 
the free hollow log structure of rank $1$ 
and such that  $f\vert_V$ factors through 
a strict \'{e}tale morphism 
$\pi \col V {\lo} {\mab A}_W(a,d,n)$ over $W$
for some $a,d\in {\mab Z}_{\geq 0}$, $n\in {\mab Z}_{\geq 1}$ 
depending on $\os{\circ}{V}$
such that $\Del_{\os{\circ}{V}}
=\{\pi^*(x_i)=0\}_{i=0}^a$ in $\os{\circ}{V}$.  
(Similarly we can give the definition of 
a {\it formal RSNCL scheme} over $S$ 
in the case where $\os{\circ}{S}$ is a formal scheme.) 
When $n\equiv 1$, 
we call a (formal) RSNCL scheme 
a ({\it formal}) {\it SNCL scheme} simply. 
\par 
(2) Let the notations be as in (1). 
Let $S'\lo S$ be a morphism of (formal) families of log points. 
If $S'=S\times_{\os{\circ}{S}}\os{\circ}{S}{}'$ 
(resp.~$S'\not=S\times_{\os{\circ}{S}}\os{\circ}{S}{}'$), 
we call the base change $X\times_SS'/S'$ of $X/S$ 
the {\it nonramified base change} (resp.~{\it ramified base change}) 
of $X/S$ with respect to the morphism $S'\lo S$. 
\end{defi}

\begin{rema}
(1) Let $X/S$ be a (formal) RSNCL scheme. 
Let $S'\lo S$ be a morphism of (formal) families of log points. 
Then $X\times_SS'/S'$ is a (formal) RSNCL scheme. 
\par 
Let $X/S$ be a (formal) SNCL scheme. 
For a nonramified (resp.~ramified) base change morphism $S'\lo S$, 
$X\times_SS'$ is (resp.~is not) an SNCL scheme over $S'$. 
\par 
(2) Let $\os{\circ}{X}$ be an SNC scheme over $\os{\circ}{S}$ 
with union of smooth schemes 
$\{\os{\circ}{X}_{\lam}\}_{\lam \in \Lam}$ over $\os{\circ}{S}$. 
Let $D$ be the non-smooth locus of $\os{\circ}{X}$ over $\os{\circ}{S}$, 
that is, $D$ is the image of $\os{\circ}{X}{}^{(1)}$ in $\os{\circ}{X}$. 
Let ${\cal I}_D$ be the ideal sheaf of $D$ in $\os{\circ}{X}$. 
Let ${\cal I}_{\lam}$ be the defining ideal sheaf of 
$\os{\circ}{X}_{\lam}$ in $\os{\circ}{X}$. 
Set ${\cal O}_D(-X):=\otimes_{\lam}({\cal I}_{\lam}/{\cal I}_{\lam}{\cal I}_D)$ 
(the ideal sheaf ${\cal O}_D(-X)$ is well-defined). 
Set ${\cal O}_D(X):={\cal H}{\it om}_{{\cal O}_X}({\cal O}_D(-X),{\cal O}_X)$. 
Then the existence of the log structure on $\os{\circ}{X}$ 
in (\ref{defi:lfac}) is equivalent to 
the $d$-semistability of Friedmann (\cite[(1.9), (2.3)]{fr1}, \cite[(11.7)]{fkato}): 
\begin{align*} 
{\cal E}{\it xt}^1(\Om^1_{\os{\circ}{X}/\os{\circ}{S}},{\cal O}_{\os{\circ}{X}})=
{\cal O}_D(X)\simeq {\cal O}_D.
\tag{1.1.17.1}\label{ali:eoxb} 
\end{align*}  
\par
(3) Though arithmetic geometers call the SNCL scheme 
the {\it strict semistable scheme}, 
we do not do so following algebraic geometers because 
it seems impossible to have the notion of semistableness 
for a general $S$ and because, even in the case where 
$S$ is the special fiber with the canonical log structure 
of the spectrum of a discrete valuation ring,  
an SNCL scheme is not necessarily the special fiber 
with the canonical log structure of an algebraic strict semistable family 
(see \cite[\S6]{nlpi} for an example). 
\end{rema}

\begin{lemm}\label{lemm:usbc}
Let $X/S$ be as in {\rm (\ref{defi:lfac})}.  
Let $T\lo S$ be a morphism of fine log schemes. 
Then $X_T:=X\times_ST$ is integral. 
$($Note that $\os{\circ}{X}_{T}:=(X_T)^{\circ}=
\os{\circ}{X}\times_{\os{\circ}{S}}\os{\circ}{T}.)$ 
\end{lemm}
\begin{proof} 
Because $M_S/{\cal O}^*_S$ is locally generated 
by a section, the morphism $X\lo S$ is integral 
(\cite[(4.4)]{klog1}). 
(\ref{lemm:usbc}) follows from this. 
\end{proof}

\par 
Assume that $M_S$ is the free hollow log structure of rank 1. 
Then $\ol{S}=({\mab A}^1_{\os{\circ}{S}}, 
({\mab N}\owns 1 \lom \tau \in {\cal O}_S[\tau])^a)$.  
Set 
\begin{equation*} 
\os{\circ}{\mab A}_{\ol{S}}(a,d):=
\ul{\rm Spec}_{\os{\circ}{S}}
({\cal O}_S[x_0, \ldots, x_{d},\tau]/(x_0\cdots x_a-\tau))\simeq {\mab A}^d_{\os{\circ}{S}}. 
\end{equation*} 
Then we have a natural structural morphism 
$\os{\circ}{{\mab A}}_{\ol{S}}(a,d) \lo \os{\circ}{\ol{S}}$. 
Let $\ol{M}_{\ol{S}}(a,d)$ 
be the log structure associated to 
a morphism 
${\mab N}^{a+1} \owns e_i=
(0, \ldots,0,\os{i}{1},0,\ldots, 0) \lom x_{i-1} \in 
{\cal O}_S[x_0, \ldots, x_{d},\tau]/(x_0\cdots x_a-\tau)$. 
Set 
$${\mab A}_{\ol{S}}(a,d)
:=(\ul{\rm Spec}_{\os{\circ}{S}}
({\cal O}_S[x_0, \ldots, x_{d},\tau]
/(x_0\cdots x_a-\tau)),\ol{M}_{\ol{S}}(a,d)).$$  
Then we have the following natural morphism 
\begin{equation*} 
{\mab A}_{\ol{S}}(a,d) \lo \ol{S}.  
\tag{1.1.18.1}
\end{equation*}

\begin{defi} 
Assume that $M_S$ is the free hollow log structure of rank $1$. 
\par 
(1) We call ${\mab A}_{\ol{S}}(a,d)$ 
$(0\leq a \leq d)$ the 
{\it standard semistable log scheme} over $\ol{S}$. 
\par 
(2) We call  $\ol{M}_{\ol{S}}(a,d)$ the 
{\it standard log structure} 
on $\os{\circ}{\mab A}_{\ol{S}}(a,d)$. 
\end{defi} 

\par 
By killing ``$\tau$'', 
we have the following natural exact closed immersion 
\begin{equation*} 
{\mab A}_S(a,d) \os{\sus}{\lo} 
{\mab A}_{\ol{S}}(a,d) 
\tag{1.1.19.1}
\end{equation*}
of fs(=fine and saturated) log schemes 
over $S\lo \ol{S}$ 
if the log structure of $S$ is the free hollow log structure of rank 1:  
we have the following commutative diagram 
of log schemes 
\begin{equation*} 
\begin{CD} 
{\mab A}_S(a,d) 
@>{\sus}>> {\mab A}_{\ol{S}}(a,d) \\ 
@VVV @VVV \\
S @>{\sus}>> \ol{S} \\ 
@VVV @VVV\\
\os{\circ}{S} @= \os{\circ}{S}.  
\end{CD} 
\tag{1.1.19.2}\label{cd:aabss} 
\end{equation*} 


\begin{rema}\label{rem:immsng}
Assume that $\ol{S}=({\mab A}^1_{\os{\circ}{S}}, 
({\mab N}\owns 1 \lom \tau \in {\cal O}_S[\tau])^a)$.  
Let $n$ be a positive integer. 
Set 
\begin{equation*} 
\os{\circ}{\mab A}_{\ol{S}}(a,d,n):=
\ul{\rm Spec}_{\os{\circ}{S}}
({\cal O}_S[x_0, \ldots, x_{d},\tau]/(x_0\cdots x_a-\tau^n)). 
\end{equation*} 
Then we have a natural structural morphism 
$\os{\circ}{{\mab A}}_{\ol{S}}(a,d,n) \lo \os{\circ}{\ol{S}}$. 
Assume that $a\geq 1$ and $n>1$. 
Then the composite morphism 
$\os{\circ}{{\mab A}}_{\ol{S}}(a,d,n) \lo \os{\circ}{\ol{S}}\lo \os{\circ}{S}$ 
is not smooth. Though we consider ${\mab A}_S(a,d,n)$ in this book, 
we do not consider $\os{\circ}{{\mab A}}_{\ol{S}}(a,d,n)$ for the case $n>1$ 
for this reason.  
\end{rema}

\begin{lemm}\label{lemm:lfm} 
The sheaf ${\cal O}_{{\mab A}_{\ol{S}}(a,d)}$ is 
a free ${\cal O}_{\ol{S}}$-module. 
\end{lemm} 
\begin{proof} 
The set $\{\prod_{i=0}^dx_i^{e_i} \vert e_i=0$ for some $0\leq i\leq a\}$ 
gives a basis of ${\cal O}_{{\mab A}_{\ol{S}}(a,d)}$ over 
${\cal O}_{\ol{S}}$. 
\end{proof}  
 
\begin{lemm}\label{lemm:etl}
Let $S$ be a family of log points 
and let $X$ be an SNCL scheme over $S$. 
Zariski locally on $X$, there exists a log scheme 
$\ol{X}$ over $\ol{S}$ fitting into 
the following cartesian diagram 
\begin{equation*}
\begin{CD} 
X @>{\sus}>> \ol{X} \\
@VVV @VVV \\
{\mab A}_S(a,d) 
@>{\sus}>> {\mab A}_{\ol{S}}(a,d) \\ 
@VVV @VVV \\
S@>{\sus}>> \ol{S},  
\end{CD} 
\tag{1.1.22.1}\label{cd:xwatx} 
\end{equation*}
where the vertical morphism 
$\ol{X} \lo {\mab A}_{\ol{S}}(a,d)$ is solid and \'{e}tale.  
\end{lemm} 
\begin{proof} 
This is a special case of (\ref{lemm:etl1}). 
\end{proof}

\begin{defi}\label{defi:lfeac}  
Let $\ol{f} \col \ol{X}\lo \ol{S}$ 
be a morphism of log schemes (on the Zariski sites). 
Set $X:=\ol{X}\times_{\ol{S}}S$.  
We call $\ol{f}$ (or $\ol{X}/\ol{S}$) 
a {\it strict semistable log scheme} over $\ol{S}$ 
if $\os{\circ}{\ol{X}}$ is a smooth scheme over $\os{\circ}{S}$, 
if $\os{\circ}{X}$ is a relative SNCD on 
$\os{\circ}{\ol{X}}/\os{\circ}{S}$ 
(with some decomposition 
$\Del:=\{\os{\circ}{X}_{\lam}\}_{\lam \in \Lam}$ of 
$\os{\circ}{X}$ by smooth components of 
$\os{\circ}{X}$ over $\os{\circ}{S}$) 
and if, for any point of $x \in \os{\circ}{\ol{X}}$, 
there exist an open neighborhood $\os{\circ}{\ol{V}}$ of 
$x$ and an open neighborhood 
$\os{\circ}{\ol{W}}\simeq \ul{\rm Spec}_W({\cal O}_W[\tau])$ 
(where $W:=\ol{W}\times_{\ol{S}}S$) of 
$\os{\circ}{f}(x)$ such that 
$M_{\ol{W}}\simeq  
({\mab N}\owns 1 \lom \tau\in {\cal O}_W[\tau])^a$ and such that  
$\ol{f}\vert_{\ol{V}}$ factors through 
a strict \'{e}tale morphism 
$\pi \col \ol{V} {\lo} {\mab A}_{\ol{W}}(a,d)$ 
such that  
$\Del_{\os{\circ}{\ol{V}}}
=\{\pi^*(x_i)=0\}_{i=0}^a$ in $\os{\circ}{\ol{V}}$.  
(Similarly we can give the definition of 
a formal strict semistable log scheme over $\ol{S}$.) 
\end{defi}

\begin{rema}\label{rema:stsmf}  
For a solid morphism $T\lo S$, the fiber product 
$\ol{X}\times_{\ol{S}}\ol{T}$ is a strict semistable log scheme over 
$\ol{T}$. 
\end{rema}

\par  
We recall the following(=the exactification) 
which has been proved in \cite{s3}.

\begin{prop}\label{prop:zpex}
{\rm {\bf (\cite[Proposition-Definition 2.10]{s3})}}\label{prop:exad}
Let $T$ be a fine log formal ${\mab Z}_p$-scheme.   
Let ${\cal C}^{\rm ex}_{{\rm hom},p}$ $($resp.~${\cal C}_p)$ 
be the category of homeomorphic exact immersions 
$($resp.~the category of immersions$)$ of 
fine log formal ${\mab Z}_p$-schemes with adic topologies 
defined by quasi-coherent ideal sheaves of structure sheaves over $T$.  
Let $\iota_p$ be the natural functor 
${\cal C}^{\rm ex}_{{\rm hom},p} \lo {\cal C}_p$. 
Then $\iota_p$ has a right adjoint functor 
$(~)^{\rm ex}\col {\cal C}_p\lo {\cal C}^{\rm ex}_{{\rm hom},p}$. 
\end{prop}
For a notational reason which will be needed later, 
we recall the construction of $(~)^{\rm ex}$ in \cite{s3}. 
\par 
By the universality, the problem is local. 
We may assume that immersions are closed immersions. 
We may assume that 
an object 
$Y\os{\sus}{\lo} {\cal Q}$ in ${\cal C}_p$ has 
a global chart $(P\lo Q)$ such that the morphism 
$P^{\rm gp}\lo Q^{\rm gp}$ is surjective.  
(As in \cite[Proposition-Definition 2.10]{s3}, we do not assume 
that the morphism $P\lo Q$ is surjective.)
Let $P^{\rm ex}$ be the inverse image of 
$Q$ by the morphism $P^{\rm gp} \lo Q^{\rm gp}$. 
Then the natural morphism $P^{\rm ex}\lo Q$ is surjective. 
We have a fine log formal ${\mab Z}_p$-scheme 
${\cal Q}^{{\rm prex}}:=
{\cal Q}\times_{{\rm Spf}^{\log}({\mab Z}_p\{P\})}
{\rm Spf}^{\log}({\mab Z}_p\{P^{\rm ex}\})$ 
over ${\cal Q}$ with a morphism $Y\lo {\cal Q}^{{\rm prex}}$. 
(The notation ``prex'' is the abbreviation of pre-exactification.) 
By the proof of \cite[(4.10) (1)]{klog1},  
the morphism $Y \lo {\cal Q}^{{\rm prex}}$ is 
an exact immersion. 
Indeed, the element $e$ in the proof of \cite[(4.10) (1)]{klog1} 
comes from an element of $P^{\rm ex}$ 
since the morphism $P^{\rm ex}\lo Q$ is surjective. 
Let ${\cal Q}^{\rm ex}$ be the formal completion of 
${\cal Q}^{{\rm prex}}$ along $Y$. 
Let $g\col Y\lo {\cal Q}^{\rm ex}$ be the natural morphism. 
Since the morphism $g$ is exact, 
the morphism $g^*\col  g^*(M_{{\cal Q}^{\rm ex}})\lo M_Y$ is injective. 
Hence the morphism 
$g^{-1}(M_{{\cal Q}^{\rm ex}}/{\cal O}_{{\cal Q}^{\rm ex}}^*)
=g^*(M_{{\cal Q}^{\rm ex}})/{\cal O}_Y^*\lo M_Y/{\cal O}_Y^*$ 
is injective. 
Since $\os{\circ}{\cal Q}{}^{\rm ex}=\os{\circ}{Y}$ as topological spaces,   
$M_{{\cal Q}^{\rm ex},x}/{\cal O}^*_{{\cal Q}^{\rm ex},x}
=M_{Y,x}/{\cal O}^*_{Y,x}$ for a point $x$ of $\os{\circ}{Y}$. 
Hence we see that the morphism $Y\os{\sus}{\lo} {\cal Q}^{\rm ex}$ 
is the desired homeomorphic exact immersion over $T$. 
We call the immersion 
$Y\os{\sus}{\lo} {\cal Q}^{\rm ex}$ or simply ${\cal Q}^{\rm ex}$ 
the {\it exactification} of the immersion 
$Y\os{\sus}{\lo} {\cal Q}$.  
\par 
More generally, let $T$ be a fine log formal scheme. 
Let ${\cal C}_T$ be the category of immersions of
fine log formal schemes over $T$ with adic topologies 
by quasi-coherent ideal sheaves of structure sheaves 
and let ${\cal C}^{\rm ex}_{{\rm hom},T}$ 
be the category of homeomorphic exact immersions of  
fine log formal schemes over $T$ with adic topologies by 
quasi-coherent ideal sheaves of structure sheaves. 
For an object $Y\os{\sus}{\lo} {\cal Q}$ in ${\cal C}$ 
which has a global chart $(P\lo Q)$ such that the morphism 
$P^{\rm gp}\lo Q^{\rm gp}$ is surjective, 
set ${\cal O}_{\cal Q}\{P\}:=\vpl_{\lam}{\cal O}_{\cal Q}/{\cal I}_{\lam}[P]$ and 
${\cal O}_{\cal Q}\{P^{\rm ex}\}
:=\vpl_{\lam}{\cal O}_{\cal Q}/{\cal I}_{\lam}[P^{\rm ex}]$,  
where $\{{\cal I}_{\lam}\}$ is the system of ideal sheaves of definitions of ${\cal Q}$. 
If we replace ${\mab Z}_p\{P\}$ and ${\mab Z}_p\{P^{\rm ex}\}$ 
with ${\cal O}_{\cal Q}\{P\}$ and ${\cal O}_{\cal Q}\{P^{\rm ex}\}$, respectively, 
then we have the following exactification functor 
\begin{align*} 
(~)^{\rm ex} \col {\cal C}_T \lo {\cal C}^{\rm ex}_{{\rm hom},T}, 
\end{align*} 
which is the right adjoint functor of 
the natural functor 
$\iota \col {\cal C}^{\rm ex}_{{\rm hom},T}\lo {\cal C}_T$.

\begin{rema}\label{rema:pds} 
In \cite[Proposition-Definition 2.10]{s3} Shiho has considered 
immersions and homeomorphic exact immersions in the category 
$({\rm LFS}/{\mathfrak B})$ of 
fine log (not necessarily $p$-adic) formal schemes 
which are separated and topologically of finite type over
a $p$-adic fine log formal scheme ${\mathfrak B}$ which is 
separated and topologically of finite type over ${\rm Spf}(W)$, where 
$W$ is a fixed Cohen ring of a field of characteristic $p > 0$. 
We do not assume that the category in (\ref{prop:zpex}) is 
this restricted one. Instead, 
we assume that immersions and homeomorphic exact immersions of 
fine log formal ${\mab Z}_p$-schemes with adic topologies 
by quasi-coherent ideal sheaves of structure sheaves. 
\end{rema}

\begin{prop}[{\bf Base changes of exactifications}]\label{prop:xpls}  
Let $T'\lo T$ be a morphism of 
fine log formal schemes. 
Let ${\cal Q}$ be a fine log formal scheme over $T$. 
Let $Y\os{\sus}{\lo} {\cal Q}$ be an immersion  
of fine log formal schemes over $T$. 
Let $Y \os{\sus}{\lo} {\cal Q}^{\rm ex}$ be 
the exactification of $Y \os{\sus}{\lo} {\cal Q}$. 
Then the immersion $Y\times_TT' \os{\sus}{\lo} {\cal Q}^{\rm ex}\times_TT'$ 
is the exactification of the immersion 
$Y\times_TT' \os{\sus}{\lo} {\cal Q}\times_TT'$ over $T'$.  
\end{prop}
\begin{proof}  
Because exact morphisms are stable under base changes 
(\cite[p.~209]{klog1}), the immersion 
$Y\times_TT'\os{\sus}{\lo} {\cal Q}^{\rm ex}\times_TT'$ is exact. 
Because $\os{\circ}{{\cal Q}}{}^{\rm ex}$ is the formal completion of 
$\os{\circ}{\cal Q}$ along $\os{\circ}{Y}$, 
it is easy to check that $({\cal Q}^{\rm ex}\times_TT')^{\circ}$ 
is also the formal completion of ${\cal Q}^{\rm ex}\times_TT'$ 
along $(Y\times_TT')^{\circ}$. 
Hence $({\cal Q}^{\rm ex}\times_TT')^{\circ}$ is isomorphic to 
$(Y\times_TT')^{\circ}$ as a topological space. 
We claim that ${\cal Q}^{\rm ex}\times_TT'$ 
has the universal property that $({\cal Q}\times_TT')^{\rm ex}$ has. 
Indeed, let $U \os{\sus}{\lo} V$ be an exact immersion of 
fine log formal schemes over $T'$  
fitting into the following commutative diagram 
\begin{equation*} 
\begin{CD} 
U@>{\sus}>> V\\
@VVV @VVV \\
Y\times_TT'@>{\sus}>>{\cal Q}\times_TT'. 
\end{CD}
\end{equation*}
Then we have the following commutative diagram 
\begin{equation*} 
\begin{CD} 
U@>{\sus}>> V\\
@VVV @VVV \\
Y@>{\sus}>>{\cal Q}. 
\end{CD}
\end{equation*}
Hence we have a natural morphism $V\lo {\cal Q}^{\rm ex}$ 
by the universality of ${\cal Q}^{\rm ex}$. 
Consequently we have the following desired commutative diagram 
\begin{equation*} 
\begin{CD} 
U@>{\sus}>> V\\
@VVV @VVV \\
Y\times_TT'@>{\sus}>>{\cal Q}^{\rm ex}\times_TT', 
\end{CD}
\end{equation*}
which shows that ${\cal Q}^{\rm ex}\times_TT'=({\cal Q}\times_TT')^{\rm ex}$. 
\end{proof} 

\begin{prop}\label{prop:excs}
Let $Y\os{\sus}{\lo}Z\os{\sus}{\lo}{\cal Q}$ be a composite morphism in 
${\cal C}_T$. Let ${\cal Q}^{\rm ex}_Z$ and ${\cal Q}^{\rm ex}_Y$ be 
the exactifications of $Z\os{\sus}{\lo}{\cal Q}$ and $Y\os{\sus}{\lo}{\cal Q}$, respectively. 
Assume that the morphism  $Y\os{\sus}{\lo}Z$ is exact. 
Then ${\cal Q}^{\rm ex}_Y$ is the formal completion of ${\cal Q}^{\rm ex}_Z$ 
along $Y$.
\end{prop}
\begin{proof} 
By the universality of the exactification, we have a natural morphism 
${\cal Q}^{\rm ex}_Y\lo {\cal Q}^{\rm ex}_Z$. The problem is local. 
Take a local chart $P\lo Q$ of $Z\os{\sus}{\lo}{\cal Q}$ such that 
the morphism $P^{\rm gp}\lo Q^{\rm gp}$ is surjective. 
Then the morphism $P\lo Q$ also 
gives a local chart of $Y\os{\sus}{\lo}{\cal Q}$. 
The desired equality follows from 
the local description of the exactification. 
\end{proof} 

\begin{prop}\label{prop:exex}
Let $T$ be a fine log formal scheme. 
Let $Y\os{\sus}{\lo} {\cal Q}$ be an immersion of 
fine log formal schemes over $T$. 
Let $({\cal Q}\times_T{\cal Q})^{\rm ex}_Y$ be the exactification of 
the diagonal immersion $Y\os{\sus}{\lo} {\cal Q}\times_T{\cal Q}$. 
Let $({\cal Q}^{\rm ex}\times_T{\cal Q}^{\rm ex})^{\rm ex}_Y$ 
be also the exactification of the diagonal immersion 
$Y\os{\sus}{\lo} {\cal Q}^{\rm ex}\times_T{\cal Q}^{\rm ex}$. 
Then the following hold$:$ 
\par 
$(1)$ $({\cal Q}^{\rm ex}\times_T{\cal Q}^{\rm ex})^{\rm ex}_Y
=({\cal Q}\times_T{\cal Q})^{\rm ex}_Y$. 
\par 
$(2)$ The log formal scheme 
$({\cal Q}^{\rm ex}\times_T{\cal Q}^{\rm ex})^{\rm ex}_Y$ is 
the formal completion of the exactification of 
the diagonal immersion ${\cal Q}^{\rm ex} \os{\sus}{\lo} 
{\cal Q}^{\rm ex}\times_T{\cal Q}^{\rm ex}$ along $Y$. 
\end{prop}
\begin{proof} 
(1): (1) follows from the universality of the exactificiation. 
\par 
(2): (2) follows from (\ref{prop:excs}).  
\end{proof}

We recall the following (\cite[(2.1.5)]{nh2}) 
which will be used in this section and later sections: 

\begin{prop}[{\bf Local structures of exact closed immersions}]\label{prop:adla}
Let $T_0 \os{\sus}{\lo} T$ be a closed immersion of fine log schemes. 
Let $Y$ $($resp.~${\cal Q})$ be a log smooth scheme over $T_0$ 
$($resp.~$T)$, which can be considered as a log scheme over $T$.  
Let $\iota \col Y \os{\sus}{\lo} {\cal Q}$ 
be an exact closed immersion over $T$. 
Let $y$ be a point of $\os{\circ}{Y}$ and 
assume that there exists a chart $(Q \lo M_T, P \lo M_Y, Q 
\os{\rho}{\lo} P)$ of $Y \lo T_0 \os{\subset}{\lo} T$ 
on a neighborhood of $y$ such that 
$\rho$ is injective, such that ${\rm Coker}(\rho^{\rm gp})$ is torsion free 
and that the natural homomorphism ${\cal O}_{Y,y} \otimes_{{\mab Z}} 
(P^{{\rm gp}}/Q^{{\rm gp}}) \lo \Om^1_{Y/T_0,y}$ is an isomorphism. 
Let ${\mab A}^n_T$ $(n\in {\mab N})$ 
be the log affine space over $T$ 
whose log structure is the inverse image of $T$ by the natural projection 
${\mab A}^n_{\os{\circ}{T}}\lo \os{\circ}{T}$. 
Then, on a neighborhood of $y$, 
there exist a nonnegative integer $c$
and the following cartesian diagram 
\begin{equation}
{\small{\begin{CD}
Y @>>> {\cal Q}' @>>> {\cal Q}\\ 
@VVV  @VVV @VVV\\ 
(T_0\otimes_{{\mab Z}[Q]}{\mab Z}[P],P^a) 
@>{\sus}>> (T\otimes_{{\mab Z}[Q]}{\mab Z}[P],P^a)
@>{\sus}>>  
(T\otimes_{{\mab Z}[Q]}{\mab Z}[P],P^a)
\times_T{\mab A}^c_T, 
\end{CD}}}
\tag{1.1.30.1}\label{eqn:0txda}
\end{equation}
where the vertical morphisms are solid and \'{e}tale and 
the lower second horizontal morphism is the base change of 
the zero section $T \os{\sus}{\lo} {\mab A}^c_T$ 
and ${\cal Q}':={\cal Q}\times_{{\mab A}^c_T}T$. 
\end{prop}

The following is a variant of an enlargement defined in \cite[\S3]{ollc}: 

\begin{defi}\label{defi:lfsz} 
Let $Z$ be a fine log (formal) scheme. 
A {\it log nil enlargement} of $Z$ is, by definition,  
an exact closed immersion $T_0 \os{\sus}{\lo} T$ 
of fine log (formal) schemes
defined by a quasi-coherent nil-ideal sheaf ${\cal J}$ of ${\cal O}_T$ 
with morphism $z\col T_0\lo Z$ of fine log (formal) schemes. 
We denote the log nil enlargement by $(T,T_0,z)$ or simply by $(T,z)$. 
We define a {\it morphism of  log nil enlargements} in an obvious way. 
\end{defi}

Let $T_0\os{\sus}{\lo} T$ be 
an exact nil-immersion of fine log (formal) schemes. 
Let $N=(N,\alpha)$ be a fine sub log structure of 
$T_0=(\os{\circ}{T}_0,M_{T_0}\os{\alpha_{T_0}}{\lo} {\cal O}_{T_0})$. 
Let $N^{\rm inv}$ be the inverse image of 
$N/{\cal O}_{T_0}^*$ by the following morphism$:$
$M_T \os{{\rm proj}.}{\lo} M_T/{\cal O}_T^* 
\os{\sim}{\lo} M_{T_0}/{\cal O}_{T_0}^*$. 
Consider the following composite morphism 
$N^{\rm inv}\os{\subset}{\lo} M_T\os{\alpha_T}{\lo} {\cal O}_T$, 
where $\al_T$ is the structural morphism. 
Then we have essentially proved the following in \cite[(2.3.1)]{nh2}. 

\begin{lemm}\label{lemm:cext}  
Let $n$ be a nonnegative integer and 
let ${\cal J}$ be the ideal sheaf of the immersion $T_0\os{\sus}{\lo} T$. 
Assume that ${\cal J}$ is $n$-torsion and 
that locally on $T$, there exists a chart $P \lo M_{T_0}$ 
such that $P^{\rm gp}$ has no $n$-torsion. 
Then $N^{\rm inv}$ is fine. 
\end{lemm}

\begin{defi}\label{defi:ceivt}  
Let the notations and the assumptions be as in (\ref{lemm:cext}).   
We call $N^{\rm inv}$ the canonical extension of 
$N$ to $\os{\circ}{T}$. 
\end{defi} 

\par 

\begin{defi} 
Let the notations be as in (\ref{defi:lfsz}). 
\par 
(1) Set $Z_{\os{\circ}{T}_0}:=Z\times_{\os{\circ}{Z}}\os{\circ}{T}_0$. 
Let $(M_{Z_{\os{\circ}{T}_0}},
(\al\col M_{Z_{\os{\circ}{T}_0}}\lo {\cal O}_{T_0}))$ be the log structure 
associated to the following composite morphism 
$$M_{Z_{\os{\circ}{T}_0}} \lo M_{T_0}\lo {\cal O}_{T_0}.$$
Then $Z_{\os{\circ}{T}_0}=(\os{\circ}{T}_0,(M_{Z_{\os{\circ}{T}_0}},
(\al\col M_{Z_{\os{\circ}{T}_0}}\lo {\cal O}_{T_0})))$. 
\par 
(2) Set $M:=({\rm Im}(M_{Z_{\os{\circ}{T}_0}} \lo M_{T_0}))^{\rm inv}$. 
Let $Z(T)$ be a log (formal) scheme obtained by $\os{\circ}{T}$ and $M$.  
We call $Z(T)$ the {\it canonical extension} of $Z_{\os{\circ}{T}_0}$.   
(If $M_{Z(T)}/{\cal O}_T^*$ is constant, we can consider the hollowing out 
$Z(T)^{\nat}$ of $Z(T)$ defined in \cite[Remark 7]{ollc}.) 
\end{defi} 

\begin{prop}\label{prop:mzt}
$(1)$ Assume that the morphism 
$M_{Z_{\os{\circ}{T}_0}} \lo M_{T_0}$ is injective. 
Then there exists a natural exact closed immersion 
$Z_{\os{\circ}{T}_0}\os{\sus}{\lo} Z(T)$.  
\par 
$(2)$ Assume that ${\cal J}$ and ${\rm Im}(M_{Z_{\os{\circ}{T}_0}} \lo M_{T_0})$ 
satisfy the assumption in {\rm (\ref{lemm:cext})}. Then $Z(T)$ is fine. 
\end{prop}
\begin{proof} 
(1): Set $N:={\rm Im}(M_{Z_{\os{\circ}{T}_0}} \lo M_{T_0})$. Obviously 
we have a closed immersion $(\os{\circ}{T}_0,N)\os{\sus}{\lo} Z(T)$. 
By the assumption, we have an isomorphism 
$Z_{\os{\circ}{T}_0}\os{\sim}{\lo} (\os{\circ}{T}_0,N)$. 
Consequently we have the closed immersion 
$Z_{\os{\circ}{T}_0}\os{\sus}{\lo} Z(T)$. 
\par 
(2): (2) follows from (\ref{lemm:cext}).     
\end{proof} 

\begin{prop}\label{prop:ostj}
Let $Z\lo Z'$ be a morphism of fine log $($formal$)$ schemes. 
Let $(T,T_0,z)\lo (T',T'_0,z')$ be a morphism of log nil enlargements over $Z\lo Z'$. 
Assume that the morphisms $M_{Z_{\os{\circ}{T}_0}} \lo M_{T_0}$ 
and $M_{Z'_{\os{\circ}{T}{}'_0}} \lo M_{T'_0}$ are injective. 
Then there exist canonical morphisms $Z_{\os{\circ}{T}_0}\lo Z'_{\os{\circ}{T}{}'_0}$ 
and $Z(T)\lo Z'(T')$ of log $($formal$)$ schemes 
fitting into the following commutative diagram 
\begin{equation}
\begin{CD}
Z(T) @>>> Z'(T')\\ 
@A{\bigcup}AA  @AA{\bigcup}A \\ 
Z_{\os{\circ}{T}_0}@>>>  Z'_{\os{\circ}{T}{}'_0}. 
\end{CD}
\tag{1.1.36.1}\label{eqn:vda}
\end{equation} 
Here the vertical morphisms in {\rm (\ref{eqn:vda})} 
are natural exact closed immersions. 
\end{prop}
\begin{proof} 
The proof is straightforward. 
\end{proof} 

The following is a nontrivial statement and 
it is important in this book, 
though the proof of it is not difficult. 

\begin{coro}\label{coro:ns}
Let $Z\lo Z'$ be as in {\rm (\ref{prop:ostj})} and 
let $(T,T_0,z)$ be a log nil enlargement of $Z$. 
Assume that the underlying morphism 
$\os{\circ}{Z} \lo \os{\circ}{Z}{}'$ is equal to ${\rm id}_{\os{\circ}{Z}}$ and 
the morphisms $M_{Z_{\os{\circ}{T}_0}} \lo M_{T_0}$ and 
$M_{Z'_{\os{\circ}{T}_0}}\lo  M_{Z_{\os{\circ}{T}_0}}$ 
are injective. 
Then there exists the following commutative diagram 
\begin{equation}
\begin{CD}
Z(T) @>>> Z'(T)\\ 
@A{\bigcup}AA  @AA{\bigcup}A \\ 
Z_{\os{\circ}{T}_0}@>>>  Z'_{\os{\circ}{T}_0}. 
\end{CD}
\tag{1.1.37.1}\label{eqn:vzda}
\end{equation} 
\end{coro} 
\begin{proof} 
This immediately follows from (\ref{prop:ostj}). 
\end{proof} 

\par 
Now let $S$ be a family of log points.  
Let $T$ be a fine log scheme. 
Let $T_0 \os{\sus}{\lo} T$ be a log nil enlargement of $S$ 
defined by a quasi-coherent nil-ideal sheaf ${\cal J}$ of ${\cal O}_T$. 
By (\ref{prop:oot}) (2) the morphism $M_{{S}_{\os{\circ}{T}_0}} \lo M_{T_0}$ is injective. 
Set $M:=({\rm Im}(M_{S_{\os{\circ}{T}_0}}\lo M_{T_0}))^{\rm inv}$. 
By (\ref{lemm:cext}) this is a fine log structure on $\os{\circ}{T}$. 
(In this case, because the following sequence 
$$0\lo {\cal O}_T^*\lo M\lo M/{\cal O}_T^*\lo 0$$ 
is locally split since there exists an isomorphism 
$M/{\cal O}_T^*\simeq {\mab N}$ locally on $T$,  
we also easily see that 
$M$ is indeed a fine log structure on $\os{\circ}{T}$.)  
Because $M_{S(T)}/{\cal O}^*_{S(T)}$ is constant, we can consider 
the hollowing out $S(T)^{\nat}$ (\cite[Remark 7]{ollc}) of $S(T)$; 
$S(T)^{\nat}$ is a family of log points. 
We obtain equalities 
$S_{\os{\circ}{T}_0}=\ul{\rm Spec}_{S(T)}^{\log}({\cal O}_{S(T)}/{\cal J})
=\ul{\rm Spec}_{S(T)^{\nat}}^{\log}({\cal O}_{S(T)^{\nat}}/{\cal J})$.

The following notion will be necessary in this book. 

\begin{defi}\label{defi:nib} 
We say that a log nil enlargement $T_0\os{\sus}{\lo}T$ with morphism 
$T_0\lo S$ is {\it restrictively hollow} or we say that $T$ is 
{\it restrictively hollow with respective to the morphism} $T_0\lo S$ 
if $S(T)$ is hollow. 
\end{defi} 


\par 
Let $X/S$ be an SNCL scheme.  
Let $T_0\os{\sus}{\lo} T$ and $T_0\lo S$ be before (\ref{defi:nib}). 
Set $X_{\os{\circ}{T}_0}:=X\times_SS_{\os{\circ}{T}_0}=
X\times_{\os{\circ}{S}}\os{\circ}{T}_0$ 
(we can consider $X$ as a fine log scheme over $\os{\circ}{S}$). 
Then $\os{\circ}{(X_{\os{\circ}{T}_0})}=
\os{\circ}{X}\times_{\os{\circ}{S}}\os{\circ}{T}_0$. 
Let 
$X_{\os{\circ}{T}_0} \os{\sus}{\lo} {\cal P}$ be 
an immersion into 
a log smooth scheme over $S(T)^{\nat}$. 
Let 
$X_{\os{\circ}{T}_0} \os{\sus}{\lo} \ol{\cal P}$ be also 
an immersion into 
a log smooth scheme over $\ol{S(T)^{\nat}}$. 
(We do not assume that there exists an immersion 
${\cal P}\os{\sus}{\lo} \ol{\cal P}$). 

\begin{rema}\label{rema:stid} 
Let $S$ and $(T,{\cal J})$ with morphism $T_0\lo S$ be as above. 
Let $T'_0\lo S$ be a morphism of fine log schemes. 
Assume that the morphism $T_0\lo S$ factors through 
a morphism $T_0\lo T'_0$. 
Then it is easy to see that $S_{\os{\circ}{T}{}'_0}(T)=S(T)$. 
\end{rema}

\begin{prop}\label{prop:fsi} 
Assume that  
$X_{\os{\circ}{T}_0}\os{\sus}{\lo} {\cal P}$ 
$($resp.~$X_{\os{\circ}{T}_0}\os{\sus}{\lo} \ol{\cal P})$
has a global chart 
$P\lo Q$ $($resp.~$\ol{P}\lo Q)$ such that 
the morphism $P^{\rm gp}\lo Q^{\rm gp}$ 
$($resp.~$\ol{P}{}^{\rm gp}\lo Q^{\rm gp})$ is surjective. 
Let $P^{\rm ex}$ $($resp.~$\ol{P}{}^{\rm ex})$ 
be the inverse image of 
$Q$ by the morphism $P^{\rm gp}\lo Q^{\rm gp}$ 
$($resp.~$\ol{P}{}^{\rm gp}\lo Q^{\rm gp})$. 
Set ${\cal P}^{\rm prex}
:={\cal P}\times_{{\rm Spec}^{\log}({\mab Z}[P])}
{\rm Spec}^{\log}({\mab Z}[P^{\rm ex}])$ 
$($resp.~$\ol{\cal P}{}^{\rm prex}
:=\ol{\cal P}\times_{{\rm Spec}^{\log}({\mab Z}[\ol{P}])}
{\rm Spec}^{\log}({\mab Z}[\ol{P}{}^{\rm ex}]))$. 
Then, locally on $X_{\os{\circ}{T}_0}$, 
there exists an open neighborhood ${\cal P}^{\rm prex}{}'$ 
$($resp.~$\ol{\cal P}{}^{\rm prex}{}')$ 
of ${\cal P}^{\rm prex}$ $($resp.~$\ol{\cal P}{}^{\rm prex})$ 
fitting into the following cartesian diagram for some $0 \leq a \leq d \leq d':$  
\begin{equation*}
\begin{CD}
X_{\os{\circ}{T}_0} @>{\subset}>> {\cal P}^{\rm prex}{}'\\ 
@VVV  @VVV \\
{\mab A}_{S_{\os{\circ}{T}_0}}(a,d)@>{\sus}>> {\mab A}_{S(T)^{\nat}}(a,d')
\end{CD}
\tag{1.1.40.1}\label{eqn:xdplxda}
\end{equation*}
$($resp.~
\begin{equation*}
\begin{CD}
X_{\os{\circ}{T}_0} @>{\subset}>> \ol{\cal P}{}^{\rm prex}{}'\\ 
@VVV  @VVV \\
{\mab A}_{S_{\os{\circ}{T}_0}}(a,d)
@>{\sus}>> {\mab A}_{\ol{S(T)^{\nat}}}(a,d')~),
\end{CD}
\tag{1.1.40.2}\label{eqn:xdpelda}
\end{equation*}
where the vertical morphisms are solid and \'{e}tale.  
\end{prop}
\begin{proof} 
Because the proof for the case 
$X_{\os{\circ}{T}_0}\os{\sus}{\lo} \ol{\cal P}{}^{{\rm prex}}$  
is the same as that for the case 
$X_{\os{\circ}{T}_0}\os{\sus}{\lo} {\cal P}^{\rm prex}$, 
we give only the proof for $X_{\os{\circ}{T}_0}\os{\sus}{\lo} {\cal P}^{{\rm prex}}$.   
\par 
Let $x$ be a point of $\os{\circ}{X}_{T_0}$. 
Since the immersion 
$X_{\os{\circ}{T}_0}\os{\sus}{\lo} {\cal P}^{{\rm prex}}$ 
is exact, 
we have an isomorphism 
\begin{equation*} 
M_{{\cal P}^{{\rm prex}},x}
/{\cal O}^*_{{\cal P}^{{\rm prex}},x}
\os{\sim}{\lo} 
M_{X_{\os{\circ}{T}_0},x}/{\cal O}^*_{X_{\os{\circ}{T}_0},x}\simeq {\mab N}^{a+1} 
\tag{1.1.40.3}\label{eqn:exemo}
\end{equation*} 
for some $a\in {\mab Z}_{\geq -1}$. 
Now (\ref{prop:fsi}) follows from (\ref{prop:adla}).
\end{proof}

\begin{prop}\label{prop:nexeo}
$(1)$ Let ${\cal P}^{\rm ex}$  
be the exactification of 
the immersion $X_{\os{\circ}{T}_0}\os{\sus}{\lo}{\cal P}$.  
Then ${\cal P}^{\rm ex}$ is a formal SNCL scheme over $S(T)^{\nat}$.  
\par 
$(2)$ 
Let $\ol{\cal P}{}^{\rm ex}$ 
be the exactification of 
the immersion $X\os{\sus}{\lo}\ol{\cal P}$. 
Then $\ol{\cal P}{}^{\rm ex}$ 
is a strict semistable family over $\ol{S(T)^{\nat}};$ 
$\os{\circ}{\ol{\cal P}}{}^{\rm ex}$
is a formally smooth scheme over $\os{\circ}{T}$  
such that 
$\os{\circ}{\cal D}:=\os{\circ}{\ol{\cal P}}{}^{\rm ex}
\times_{\os{\circ}{\ol{S(T)^{\nat}}}}\os{\circ}{T}$ is 
a formal SNCD on 
$\os{\circ}{\ol{\cal P}}{}^{\rm ex}$ and such that 
$M_{{\ol{\cal P}}{}^{\rm ex}}=M(\os{\circ}{\cal D})$, 
where $M(\os{\circ}{\cal D})$ is the associated log structure to 
$\os{\circ}{\cal D}$ on $\os{\circ}{\ol{\cal P}}{}^{\rm ex}/\os{\circ}{T}$ 
$(${\rm \cite[p.~61]{nh2}~(cf.~\cite{fao}, \cite{klog1})}$)$. 
\end{prop}
\begin{proof} 
(1): Since 
$\os{\circ}{\cal P}{}^{\rm ex}$ 
is topologically isomorphic to $\os{\circ}{X}$, 
(\ref{prop:nexeo}) follows from (\ref{prop:fsi}). 
\par 
(2): (2) also immediately follows from (\ref{prop:fsi}) 
as in (1). 
\end{proof}

\par 
By the proof of (\ref{prop:nexeo}), we see the following:  
\par 
For a smooth component $\os{\circ}{X}_{\lam}$ of $\os{\circ}{X}_{T_0}$, 
we can define a closed subscheme 
$\os{\circ}{\ol{\cal P}}{}^{\rm ex}_{\lam}$ of 
$\os{\circ}{\ol{\cal P}}{}^{\rm ex}$ 
(resp.~$\os{\circ}{\cal P}{}^{\rm ex}_{\lam}$ of 
$\os{\circ}{\cal P}{}^{\rm ex}$) 
fitting into the following commutative diagram: 
\begin{equation*} 
\begin{CD} 
\os{\circ}{X}_{\lam} @>{\subset}>>
\os{\circ}{\ol{\cal P}}{}^{\rm ex}_{\lam}\\
@V{\bigcap}VV @VV{\bigcap}V \\
\os{\circ}{X} @>{\subset}>>
\os{\circ}{\ol{\cal P}}{}^{\rm ex}
\end{CD}
\end{equation*} 
(resp. 
\begin{equation*} 
\begin{CD} 
\os{\circ}{X}_{\lam} @>{\subset}>>
\os{\circ}{{\cal P}}{}^{\rm ex}_{\lam}\\
@V{\bigcap}VV @VV{\bigcap}V \\
\os{\circ}{X} @>{\subset}>>
~\os{\circ}{{\cal P}}{}^{\rm ex}~)
\end{CD}
\end{equation*} 
where $\os{\circ}{X}_{\lam} \os{\subset}{\lo} 
\os{\circ}{\ol{\cal P}}{}^{\rm ex}_{\lam}$ 
(resp.~ $\os{\circ}{X}_{\lam} \os{\sus}{\lo} 
\os{\circ}{\cal P}{}^{\rm ex}_{\lam}$) 
is an isomorphism as topological spaces. 
(In \S\ref{sec:ldc} below we give different definitions of 
$\os{\circ}{\ol{\cal P}}{}^{\rm ex}_{\lam}$ and 
$\os{\circ}{\cal P}{}^{\rm ex}_{\lam}$.)

\par 
Let the notations be as in the beginning of this section. 
Let $v \col T' \lo T$ be a morphism of log schemes 
with locally free hollow log structures of ranks $r'$ and $r$. 
Let $x$ be a point of $\os{\circ}{T}{}'$. 
Let 
\begin{align*} 
v^*\col {\mab N}^r= M_{T,v(x)}/{\cal O}^*_{T,v(x)}
\lo M_{T',x}/{\cal O}^*_{T',x}={\mab N}^{r'}
\tag{1.1.41.1}\label{ali:ntvo} 
\end{align*} 
be the induced morphism. 
Set 
$e_i=(0, \ldots,0,\os{i}{1},0,\ldots, 0)\in {\mab N}^r$ 
$(0\leq i\leq r)$. 
Let $A\in M_{r'r}({\mab Z})$ be the representing matrix of 
$v^*{}^{\rm gp}$: 
$v^*{}^{\rm gp}(e_1,\ldots,e_r)=(e_1,\ldots,e_{r'})A$. 
Let 
$d'_1\vert d'_2\vert \cdots \vert d'_{t-1}\vert d'_t$ 
$(0\leq t\leq \min \{r,r'\},d'_1,\ldots,d'_t\in {\mab Z})$ 
be associated divisors to $A$ obtained by 
the theory of elementary divisors.  
Set $d_i:=\vert d'_i\vert$. 

\begin{defi}\label{defi:ddef} 
We denote $(d_1,\ldots, d_t)$ by $\deg(v)_x$;
$\deg(v)_x$ is independent of the choice of $x$. 
We denote it simply by $\deg(v)$ and 
we call $\deg(v)$ the {\it $($mapping$)$ degree} of $v$. 
We set $\deg_{\pi}(v):=d_1\cdot \cdots d_t$. 
(The subscript $\pi$ means the product.) 
\end{defi} 

\begin{rema}\label{rema:rr1}
(1) If $r=r'=1$, then 
it is easy to see that $\deg(v)=1$ if and only if $v$ is solid. 
\par 
(2) 
Let 
\begin{align*} 
v^*_{\mab Z}\col {\mab Z}^r= M_{T,v(x)}/{\cal O}^*_{T,v(x)}
\lo M_{T',x}/{\cal O}^*_{T',x}={\mab Z}^{r'}
\tag{1.1.43.1}\label{ali:ntzro} 
\end{align*} 
be the induced morphism by (\ref{ali:ntvo}). 
It is easy to see that 
\begin{align*} 
\deg_{\pi}(v)={\rm vol}
(({\rm Im}(v^*_{\mab Z})\otimes_{\mab Z}{\mab R})
/{\rm Im}(v^*_{\mab Z})),
\tag{1.1.43.2}\label{ali:volf}
\end{align*} 
where ``${\rm vol}$'' means the volume with respect to the standard Haar measure on 
${\rm Im}(v^*_{\mab Z})\otimes_{\mab Z}{\mab R}$.
\end{rema}


\begin{prop}\label{prop:dvs}
Let $v' \col T'' \lo T'$ be a similar morphism of log schemes 
with locally free hollow log structures of ranks $r''$ and $r'$. 
Then 
\begin{align*} 
\deg_{\pi} (v'\circ v)=\deg_{\pi} (v') \deg_{\pi}(v).
\tag{1.1.44.1}\label{ali:dpvv}
\end{align*} 
\end{prop} 
\begin{proof} 
This immediately follows from the formula (\ref{ali:volf}). 
\end{proof}  

\par 
We conclude this section by proving 
the following useful proposition.

\begin{prop}[{\bf A log version of \cite[N.~B. in 5.27]{bob}}]\label{prop:niip}
Let $Y\lo T$ be a log smooth morphism of fine log schemes. 
Let $\iota \col Z_0\os{\sus}{\lo} Z$ be 
an exact nil closed immersion of fine log affine schemes over $T$. 
Then the natural morphism 
\begin{align*} 
Y(Z):={\rm Hom}_T(Z,Y)\lo Y(Z_0):={\rm Hom}_T(Z_0,Y)
\tag{1.1.45.1}\label{ali:slzy}
\end{align*} 
is surjective. 
\end{prop} 
\begin{proof}  
As in \cite[N.B. in 5.27]{bob}, express 
$\os{\circ}{Z}={\rm Spec}(B)$ and $\os{\circ}{Z}_0={\rm Spec}(B/J)$, 
$J=\bigcup_{\lam} J_{\lam}$, where $J_{\lam}$ is a finitely generated ideal of $B$.  
Set $B_{\lam}:=B/J_{\lam}$. Then we have the natural nilpotent closed immersion 
${\rm Spec}(B_{\lam})\os{\sus}{\lo} \os{\circ}{Z}$. 
Set $Z_{\lam}:=Z\times_{\os{\circ}{Z}}{\rm Spec}(B_{\lam})$.  
Then we have the natural exact closed nilpotent immersion 
$Z_{\lam}\os{\sus}{\lo} Z$ and the natural exact closed nil immersion 
$Z_0\os{\sus}{\lo} Z_{\lam}$. 
Since $\iota$ is an isomorphism as topological spaces, 
we identify a sheaf on $Z$ with a sheaf on $Z_0$ and with a sheaf on $Z_{\lam}$. 
Because $M_{Z_0}$ is the inverse image of $M_Z$, 
$M_{Z_0}=M_Z\oplus_{{\cal O}_Z^*}{\cal O}_{Z_0}^*$. 
Since ${\cal O}_{Z_0}=\vil_{\lam}{\cal O}_{Z_{\lam}}$, 
${\cal O}_{Z_0}^*=\vil_{\lam}{\cal O}_{Z_{\lam}}^*$. 
Hence $M_{Z_0}=M_Z\oplus_{{\cal O}_Z^*}\vil_{\lam}{\cal O}_{Z_{\lam}}^*
=\vil_{\lam}(M_Z\oplus_{{\cal O}_Z^*}{\cal O}_{Z_{\lam}}^*)
=\vil_{\lam}M_{Z_{\lam}}$. 
Let $h\col Z_0\lo Y$ be a morphism over $T$. 
Let $g_0\col Z_0\lo T$ be the structural morphism. 
Then we have the following commutative diagram 
\begin{equation*} 
\begin{CD} 
M_{Z_0}=\vil_{\lam}M_{Z_{\lam}}@<<< h^{-1}(M_Y)\\
@VVV @VVV \\
{\cal O}_{Z_0}=\vil_{\lam}{\cal O}_{Z_{\lam}}@<<< h^{-1}({\cal O}_Y)
\end{CD}
\tag{1.1.45.2}\label{ali:mzmzy}
\end{equation*} 
over $g_0^{-1}(M_T)\lo g_0^{-1}({\cal O}_T)$. 
By \cite[(8.8.2) (i)]{ega43}  
the lower horizontal morphism in (\ref{ali:mzmzy}) factors through a morphism 
$\os{\circ}{h}{}^{-1}_{\lam}({\cal O}_Y)\lo {\cal O}_{Z_{\lam}}$ for some $\lam$, 
where $\os{\circ}{h}_{\lam}\col \os{\circ}{Z}_{\lam}\lo  \os{\circ}{Y}$ is a 
morphism of schemes. 
Hence the morphism $h^{-1}({\cal O}_Y^*)\lo {\cal O}^*_{Z_0}$ 
factors through a morphism 
$\os{\circ}{h}{}^{-1}_{\lam}({\cal O}^*_Y)\lo {\cal O}^*_{Z_{\lam}}$ for some $\lam$.  
Take a local chart $P\lo {\cal O}_Y$ of $M_Y$. 
Since $P$ is finitely generated and since 
$h^{-1}({\cal O}_Y^*)\lo {\cal O}^*_{Z_0}$ factors through 
$\os{\circ}{h}{}^{-1}_{\lam}({\cal O}^*_Y)\lo {\cal O}^*_{Z_{\lam}}$ for some $\lam$, 
the morphism $h^{-1}(M_Y)\lo M_{Z_0}$ 
factors through $\os{\circ}{h}{}^{-1}_{\lam}(M_Y)\lo M_{Z_{\lam}}$ 
for some $\lam$ locally on $\os{\circ}{Z}_0$. 
Since $\os{\circ}{Z}_0$ is quasi-compact,   
the morphism $h^{-1}(M_Y)\lo M_{Z_0}$ 
factors through $\os{\circ}{h}{}^{-1}_{\lam}(M_Y)\lo M_{Z_{\lam}}$ for some $\lam$. 
Hence we have the morphism $h_{\lam}\col Z_{\lam}\lo Y$ over $T$ such that 
the composite morphism $Z_0\os{\sus}{\lo} Z_{\lam}\os{h_{\lam}}{\lo} Y$ is 
equal to $h$.  
By the definition of the log smoothness,  
the morphism $Y(Z)\lo Y(Z_{\lam})$ is surjective. 
Consequently the morphism (\ref{ali:slzy}) is surjective. 
\end{proof}

\section{Review for results in [Nakk4] and 
bisimplicial (exact) immersions I}\label{sec:tdiai}
In this section we review several results 
and additional remarks to \cite{nh3} 
for later sections. 
We give only the statements of the results; 
see [loc.~cit.] for the proofs of them. 
\par
Let $\Del$ be the standard simplicial 
category: an object of $\Del$ is denoted by 
$[n]:=\{0, \ldots, n\}$ $(n\in {\mab N})$; 
a morphism in $\Del$ is 
a non-decreasing function $[n] \lo [m]$ 
$(n,m\in {\mab N})$. 
\par  
Let $Y$ be a fine log $($formal$)$ scheme over   
a fine log $($formal$)$ scheme $T$.
Let $Y_{\bul}$ be a fine simplicial log $($formal$)$ scheme over $Y$. 
Let $s^{m-1}_i \col Y_{m-1} \lo Y_m$ 
$(m\in {\mab Z}_{>0}, 0\leq i \leq m-1)$ be the degeneracy
morphism corresponding to the standard degeneracy map 
$\partial^i_m \col [m] \lo [m-1]$:  
$\partial^i_m(j)=j$ $(0\leq j \leq i)$, 
$\partial^i_m(j)=j-1$ $(i< j \leq m)$. 
Following \cite[${\rm V}^{\rm bis}$ \S5]{sga4-2} 
and \cite[(6.2.1.1)]{dh3},
set $N(Y_0):=Y_0$ 
and let $N(Y_m)$ $(m\in {\mab Z}_{>0})$ be the 
intersection of the complements of 
$s^{m-1}_i(Y_{m-1})$ $(0 \leq i \leq m-1)$. 
The simplicial log $($formal$)$ scheme $Y_{\bul}$ 
is said to be split ([loc.~cit.]) 
if $Y_m=\coprod_{0\leq l \leq m}\coprod_{[m] 
\twoheadrightarrow [l]}N(Y_l)$, 
where the subscripts $[m]\twoheadrightarrow [l]$'s run through the surjective non-decreasing morphisms 
$[m]\lo [l]$'s.

\begin{lemm}[{\rm {\bf \cite[(6.1)]{nh3}}}]\label{lemm:lisj} 
Assume that $Y_{\bul}$ is split.  
Then there exists a fine split simplicial log $($formal$)$ scheme 
$Y'_{\bul}$ with a morphism $Y'_{\bul} \lo Y_{\bul}$ 
of simplicial log $($formal$)$ schemes over $Y$ 
satisfying the following conditions$:$
\medskip 
\parno
$(1.2.1.1)$ $Y'_m$ $(m\in {\mab N})$ is the disjoint union 
of log affine open $($formal$)$ subschemes 
which cover $Y_m$ and whose images in 
$Y$ are contained in 
log affine open $($formal$)$ subschemes of $Y$.
\medskip 
\parno  
$(1.2.1.2)$ If $\os{\circ}{Y}_m$ $(m\in {\mab N})$ 
is quasi-compact, then the cardinarity of 
the log affine open $($formal$)$ subschemes 
in $(1.2.1.1)$ can be assumed to be finite.
\medskip
\parno
Set 
$Y_{m n}:={\rm cosk}_0^{Y_m}(Y_m')_n$ 
$(m,n\in {\mab N})$.
Then there exists a natural morphism
$Y_{\bul \bul} \lo Y_{\bul}$  over $Y$.
For each $n\in {\mab N}$, $Y_{\bul n}$ is split. 
\end{lemm}

\par 
The construction of $Y'_{\bul}$ is as follows. 
\par 
Let $Y'_0$ be the disjoint union of 
log affine open (formal) subschemes which cover
$Y_0$ and whose images in 
$Y$ are contained in 
log affine open (formal) subschemes of $Y$. 
Let $m$ be a positive integer.
Assume that we are given $Y'_{\bul \leq m-1}$. 
The log scheme 
${\rm cosk}^Y_{m-1}
(Y'_{\bul \leq m-1})_m$ 
is the disjoint union of 
the members of a log affine open covering  of 
${\rm cosk}^Y_{m-1}(Y_{\bul \leq m-1})_m$. 
Consider the natural composite morphism 
$N(Y_m) \os{\sus}{\lo} Y_m 
\lo {\rm cosk}^Y_{m-1}
(Y_{\bul \leq m-1})_m$ and 
a log affine open covering of $N(Y_m)$
which refines the inverse image 
of the open covering of 
${\rm cosk}^Y_{m-1}
(Y_{\bul \leq m-1})_m$.  
Let $N(Y_m)'$ be the 
disjoint union of the members of 
this open covering.
Then we have the following 
commutative diagram:
\begin{equation*}
\begin{CD}
N(Y_m)' @>>> {\rm cosk}^Y_{m-1}
(Y'_{\bul \leq m-1})_m \\
@VVV @VVV  \\ 
N(Y_m)@>>> {\rm cosk}^Y_{m-1}
(Y_{\bul \leq m-1})_m.
\end{CD}
\tag{1.2.1.3}\label{cd:spnscyny}
\end{equation*}
Then we set $Y'_m=\coprod_{0\leq l \leq m}\coprod_{[m] 
\twoheadrightarrow [l]}N(Y_l)'$, 
which is the inductive construction of $Y'_{\bul}$.  
\par 
In \cite{nh3} we have called the simplicial log (formal) scheme 
$Y'_{\bul}$ satisfying $(1.2.1.1)$ and $(1.2.1.2)$ 
the {\it disjoint union of the members of 
an affine simplicial open covering} of $Y_{\bul}/Y$. 
In this book, if there exists a morphism $Z'_{\bul}\lo Z_{\bul}$ of 
fine simplicial log (formal) schemes over $Y$ satisfying 
$(1.2.1.1)$ and $(1.2.1.2)$,  then 
we also say that 
$Z_{\bul}$ has an affine simplicial open covering $Z'_{\bul}$ 
over $Y$. 
In [loc.~cit.] we have also called the bisimplicial scheme 
$Z_{\bul \bul}$ in  (\ref{lemm:lisj}) 
the {\it \v{C}ech diagram} of $Z'_{\bul}$
over $Z_{\bul}/Y$. 

\begin{prop}[{\rm {\bf \cite[(6.3)]{nh3}}}]
\label{prop:lissmp}
$(1)$ Let $Z_{\bul}\lo Y_{\bul}$ 
be a morphism of fine split simplicial log $($formal$)$ schemes 
over a morphism $Z\lo Y$ of fine log $($formal$)$ schemes 
over $T$.  
Then there exist the disjoint unions of 
the members of affine simplicial open coverings
$Y'_{\bul}$ and $Z'_{\bul}$
of $Y_{\bul}/Y$ and $Z_{\bul}/Z$, 
respectively, fitting into 
the following commutative 
diagram$:$
\begin{equation*}
\begin{CD}
Z'_{\bul}@>>> Y'_{\bul}  \\
@VVV @VVV  \\
Z_{\bul}@>>>Y_{\bul}.
\end{CD}
\tag{1.2.2.1}\label{cd:clchbs} 
\end{equation*}  
\par
$(2)$ Let $Y_{\bul}/Y/T$ be as in {\rm (\ref{lemm:lisj})}.  
Let $Y'_{\bul}$ and $Y''_{\bul}$ 
be two disjoint unions of the members of 
affine simplicial open coverings of $Y_{\bul}/Y$. 
Then there exists a disjoint union of 
the members of 
an affine simplicial open covering $Y'''_{\bul}$
of $Y_{\bul}/Y$ 
fitting into the following 
commutative diagram$:$
\begin{equation*}
\begin{CD}
Y'''_{\bul} @>>> Y''_{\bul}\\
@VVV @VVV  \\
Y'_{\bul}@>>> Y_{\bul}.
\end{CD}
\tag{1.2.2.2}\label{cd:celccov}
\end{equation*} 
\end{prop}

\par
Next we need the obvious log version 
$\Gam^T_N(Y)^{\leq l}$ $(N, l \in {\mab N})$ of 
\cite[\S11]{ctze} for log schemes. 
\par
Let ${\cal C}$ be a category which has 
finite inverse limits. Let $l$ be a nonnegative integer.
Following \cite[\S11]{ctze},  define a set 
${\rm Hom}_{\Del}^{\leq l}([n], [m])$ 
$(n,m\in {\mab N})$:
the set  ${\rm Hom}_{\Del}^{\leq l}([n], [m])$, 
by definition, consists of 
the morphisms $\gam \col [n] \lo [m]$ in $\Del$ 
such that the cardinality 
of the set $\gam([n])$ is less than or equal to $l$. 
\par 
Let $N$ be a nonnegative integer. 
Let us recall the definition of the simplicial object 
$\Gam:=
\Gam^{\cal C}_{N}(X)^{\leq {}{l}}$ in 
${\cal C}$
for an object $X\in {\cal C}$ 
(\cite[(7.3.1)]{tzcp}):
for an object $[m]\in \Del$, set 
$$\Gam_{{}{m}}:=\prod_{\gam \in 
{\rm Hom}_{\Del}^{\leq {}{l}}
([{}{N}], [{}{m}])}X_{\gam}$$
with $X_{\gam}=X$; 
for a morphism $\al \col [{}{m}'] \lo [{}{m}]$ 
in $\Del$, $\al_{\Gam} \col \Gam_{{}{m}} \lo 
\Gam_{{}{m}'}$ is defined to be the following:
``$(c_\gam) \lom (d_{\bet})$'' with 
``$d_{\bet}=c_{\al \bet}$''. 
In fact, we have a functor
\begin{equation*}
\Gam^{\cal C}_N(?)^{\leq l} \col {\cal C} \lo 
{\cal C}^{\Del}:=
\{\text{simplicial objects of ${\cal C}$}\}.
\tag{1.2.2.3}
\end{equation*}
Set $\Gam^{\cal C}_{{}{N}}(X)=
\Gam^{\cal C}_{{}{N}}(X)^{\leq N}$. 
The functor 
$\Gam^{\cal C}_N(?)^{\leq l}$ commutes 
with finite inverse limits.
For $l\geq N$,   
$\Gam^{\cal C}_{{}{N}}(X)^{\leq {}{l}}=
\Gam^{\cal C}_{{}{N}}(X)$.

\begin{defi}\label{defi:tkifun}
We call the functor $\Gam^{\cal C}_N(?)^{\leq l}$ 
{\it Tsuzuki's functor}.  
For an object $X$ of ${\cal C}$, 
we call $\Gam^{\cal C}_N(X)^{\leq l}$ 
{\it Tsuzuki's simplicial object} of $X$.
\end{defi}

\par
In \cite[(11.2.5)]{ctze} Chiarellotto and Tsuzuki 
have essentially proved the following: 

\begin{lemm}[{\rm {\bf \cite[(7.3.2)]{tzcp}}}]\label{gamfun} 
${\rm cosk}^{\cal C}_l
(\Gam^{\cal C}_N(X)_{\bul \leq l})=
\Gam^{\cal C}_N(X)^{\leq l}$.
\end{lemm}

\par
Let $N$ and $m$ be nonnegative integers.
Let $X_{\bul \leq N}$ be 
an $N$-truncated simplicial object of ${\cal C}$ and 
let $f \col X_N \lo Y$ be a morphism in ${\cal C}$.
Then we have a morphism
\begin{equation*}
X_{\bul \leq N} \lo \Gam^{\cal C}_{N}(Y)_{\bul \leq N}
\tag{1.2.4.1}\label{eqn:cgam}
\end{equation*}
making the following diagram commutative:
\begin{equation*}
\begin{CD}
X_{{}{m}} @>>> \prod_{\gam \in 
{\rm Hom}_{\Del}({}{[{N}]}, {}{[m]})}Y_{\gam}\\ 
@V{X(\del)}VV @VV{\rm proj}V \\
X_N @>{f}>> Y_{\del}=Y
\end{CD}
\tag{1.2.4.2}\label{cd:embdf}
\end{equation*} 
for any morphism $\del \col [N]\lo [m]$ in $\Del$.

\par 
Let $T$ be a fine log (formal) scheme. 
Let ${\cal C}_T$ be the category of 
fine log $($formal$)$ schemes over $T$. 
Set $\Gam^T_N(?):=\Gam^{{\cal C}_T}_N(?)$ 
and ${\rm cosk}_N^T(?)={\rm cosk}_N^{{\cal C}_T}(?)$. 
Now we limit ourselves to the case ${\cal C}={\cal C}_T$. 
By (\cite[(1.6), (2.8)]{klog1}) 
the finite inverse limit exists in ${\cal C}_T$. 
Hence we can apply 
the constructions above to ${\cal C}_T$.
\par
The following are immediate generalizations of 
\cite[(11.2.4), (11.2.6)]{ctze}:

\begin{lemm}[{\rm {\bf \cite[(6.6)]{nh3}}}]\label{lemm:pgamfun} 
Let $N$ be a nonnegative integer. 
Then the following hold$:$
\par
$(1)$ Let $X_{\bul \leq N}$ be a fine 
$N$-truncated simplicial log $($formal$)$ scheme over $T$. 
If $f\col X_N \lo Y$ is 
an immersion of fine log $($formal$)$ schemes over $T$, 
then the morphism 
$X_{\bul \leq N} \lo \Gam^T_N(Y)_{\bul \leq N}$ 
in $(\ref{eqn:cgam})$ is an 
immersion of fine $N$-truncated simplicial 
log $($formal$)$ schemes over $T$. 
\par
$(2)$ Let $X \lo T$ be a  
morphism of fine log $($formal$)$ schemes. 
Assume that the morphism 
$X \lo T$ satisfies a condition {\rm (P)} 
which is stable under the morphism 
$\underset{n~{\rm times}}
{\underbrace{X\times_T
\cdots \times_TX}} 
\os{{\rm proj}.}{\lo}
\underset{m~{\rm times}}
{\underbrace{X\times_T
\cdots \times_TX}}$   
for any $n>m$.
Then the natural morphism
$${\rm cosk}^T_l(\Gam^T_N(X)_{\bul \leq l})_m 
\lo {\rm cosk}^T_{l'}
(\Gam^T_N(X)_{\bul \leq l'})_m$$ 
satisfies {\rm (P)} for $l' < l$ and for any $m\in {\mab N}$.
\end{lemm}

\par 
The following is an obvious generalization of \cite[(6.10)]{nh3} 
in which we have assumed the immersion 
$T \os{\sus}{\lo} {\cal T}$ in the following is a nil-immersion of 
fine log schemes. 
  
\begin{prop}[{\bf The existence of an embedding system}]\label{prop:ytlft} 
Let the notations and the assumptions 
be as in {\rm (\ref{lemm:lisj})}. 
Let $T \os{\sus}{\lo} {\cal T}$ be an immersion of fine log schemes. 
Let $N$ be a nonnegative integer. 
Assume that $Y_N$ is log smooth over $T$. 
For any point $y$ of $Y_N$ over $T$, 
assume that the assumption in {\rm (\ref{prop:ebtl})} holds. 
Then there exists the disjoint union $Y'_{\bul \leq N}$ of 
the members of an affine $N$-truncated simplicial open covering of 
$Y_{\bul \leq N}$ and an immersion 
$Y'_{\bul \leq N} \os{\sus}{\lo} {\cal Q}'_{\bul \leq N}$ 
into 
a log smooth $N$-truncated simplicial log scheme 
over ${\cal T}$. 
\end{prop} 
\begin{proof} 
For the completeness of this book, we give the proof. 
\par 
Because $Y_N$ is log smooth over $T$, 
$Y'_N$ is also log smooth over $T$. 
By (\ref{prop:ebtl}) we may assume that 
there exists a closed immersion 
$Y'_N \os{\sus}{\lo} {\cal Y}'_N$ into a log smooth scheme over ${\cal T}$.  
Set ${\cal Q}'_{\bul}:=\Gam^{\cal T}_N({\cal Y}'_N)$. 
The immersion (\ref{eqn:cgam}) is obtained by  
the following composite morphism 
$$Y'_m  \os{\sus}{\lo} 
\Pi_{\gam \in {\rm Hom}_{\Del}([N],[m])}Y'_m 
\lo \Pi_{\gam \in {\rm Hom}_{\Del}([N],[m])}Y'_N 
\os{\sus}{\lo} 
\Pi_{\gam \in {\rm Hom}_{\Del}([N],[m])}{\cal Y}'_N
={\cal Q}'_m,$$ 
where the morphism 
$Y'_m \os{\sus}{\lo} 
\Pi_{\gam \in {\rm Hom}_{\Del}([N],[m])}Y'_m$ 
is the diagonal immersion.  
In fact, by the following commutative diagram 
\begin{equation*} 
\begin{CD} 
Y'_l @>>> Y'_N @>{\subset}>> {\cal Y}'_N \\ 
@VVV @| @| \\ 
Y'_m @>>> Y'_N @>{\subset}>> {\cal Y}'_N 
\end{CD} 
\end{equation*} 
for morphisms $[m] \lo [l]$ and $[N] \lo [m]$ 
in $\Del$, 
we have the desired $N$-truncated simplicial immersion 
$Y'_{\bul \leq N} \os{\sus}{\lo} {\cal Q}'_{\bul \leq N}$. 
\end{proof}

\par 
Let $Y_{\bul \leq N,\bul}$ be the affine \v{C}ech diagram of 
$Y'_{\bul \leq N}$ over $Y_{\bul \leq N}/Y$: 
$Y_{mn}:={\rm cosk}_0^{Y_m}(Y'_m)_n$ 
$(0\leq m\leq N, n\in {\mab N})$. 
Set ${\cal Q}_{mn}:={\rm cosk}_0^{\cal T}({\cal Q}'_m)_n$.  
Then there exists an immersion 
$Y_{\bul \leq N,\bul} \os{\sus}{\lo} {\cal Q}_{\bul \leq N,\bul}$ 
into a log smooth $(N,\infty)$-truncated bisimplicial scheme 
over ${\cal T}$.


\par  
Let $S$ be a family of log points.  
Let $T_0 \os{\subset}{\lo}T$ be a log nil enlargement of $S$ 
defined by a quasi-coherent nil-ideal sheaf 
${\cal J}$ of ${\cal O}_T$.  
Now we apply the above result to the case of 
a truncated simplicial SNCL scheme over $S_{\os{\circ}{T}_0}$ 
which has an affine truncated simplicial open covering 
and the immersions $S_{\os{\circ}{T}_0}\os{\sus}{\lo} S(T)^{\nat}$ 
and $S_{\os{\circ}{T}_0}\os{\sus}{\lo}\ol{S(T)^{\nat}}$. 
(Recall $S_{\os{\circ}{T}_0}$, $S(T)$ and $S(T)^{\nat}$ before (\ref{defi:nib}) and 
note that the immersion $S_{\os{\circ}{T}_0}\os{\sus}{\lo}\ol{S(T)^{\nat}}$ 
is not nil.)  

\begin{defi}\label{defi:ntrxs} 
Let $N$ be a nonnegative integer. 
Let $f \col X_{\bul \leq N} \lo S$ 
be an $N$-truncated simplicial fine log scheme over $S$. 
We call $f$ or $X_{\bul \leq N}/S$ 
an {\it $N$-truncated simplicial SNCL scheme} over $S$ if 
$X_m$ $(0\leq \forall m\leq N)$ is an SNCL scheme over $S$. 
\end{defi} 

Assume that $X_{\bul \leq N,\os{\circ}{T}_0}:=
X_{\bul \leq N}\times_SS_{\os{\circ}{T}_0}=
X_{\bul \leq N}\times_{\os{\circ}{S}}\os{\circ}{T}_0$ 
has an affine $N$-truncated simplicial open covering 
(we can consider $X_{\bul \leq N}$ as an 
$N$-truncated simplicial fine log scheme over $\os{\circ}{S}$).    
By abuse of notation, we denote by the same symbol $f$ 
the structural morphism 
$X_{\bul \leq N,\os{\circ}{T}_0} \lo S_{\os{\circ}{T}_0}$. 
\par 
Let $X'_{\bul \leq N,\os{\circ}{T}_0}$ be the disjoint union 
of an affine $N$-truncated simplicial open covering of 
$X_{\bul \leq N,\os{\circ}{T}_0}$ over $S_{\os{\circ}{T}_0}$. 
Let $X'_{N,\os{\circ}{T}_0}=\coprod_{i\in I}X_{Ni}$ 
be an expression of $X'_{N,\os{\circ}{T}_0}$, where 
$X_{Ni}$ is a log open subscheme of $X_{N,\os{\circ}{T}_0}$. 
Assume that $f(\os{\circ}{X}_{Ni})$ is contained in 
an affine open subscheme of 
$\os{\circ}{T}=(S_{\os{\circ}{T}_0})^{\circ}$ 
such that the restriction of $M_{S_{\os{\circ}{T}_0}}$ 
to this open subscheme is free of rank $1$. 
Assume also that there exists a solid \'{e}tale morphism 
$X_{Ni}\lo {\mab A}_{S_{\os{\circ}{T}_0}}(a,d)$. 
Then, replacing $X_{Ni}$ by a small log open subscheme of 
$X_{N,\os{\circ}{T}_0}$, we can assume that 
there exists a log smooth scheme 
$\ol{\cal P}{}'_{Ni}/\ol{S(T)^{\nat}}$ 
fitting into the following commutative diagram 
\begin{equation*}
\begin{CD} 
X_{Ni}@>{\subset}>> \ol{\cal P}{}'_{Ni}\\
@VVV @VVV \\ 
{\mab A}_{S_{\os{\circ}{T}_0}}(a,d) @>{\subset}>>
{\mab A}_{\ol{S(T)^{\nat}}}(a,d) \\
@VVV @VVV \\ 
S_{\os{\circ}{T}_0}@>{\subset}>> \ol{S(T)^{\nat}},   
\end{CD}
\tag{1.2.7.1}\label{cd:xnip} 
\end{equation*}
where the morphism 
$\ol{\cal P}{}'_{Ni}\lo {\mab A}_{\ol{S(T)^{\nat}}}(a,d)$ 
is solid and \'{e}tale ((\ref{lemm:etl})). 
Set $\ol{\cal P}{}'_N:=\coprod_{i\in I}\ol{\cal P}{}'_{Ni}$ 
and $\ol{\Gam}{}'_{\bul}
:=\Gam^{\ol{S(T)^{\nat}}}_N(\ol{\cal P}{}'_N)$. 
Set 
\begin{equation*} 
X_{mn,\os{\circ}{T}_0}:={\rm cosk}_0^{X_{m,\os{\circ}{T}_0}}(X'_{m,\os{\circ}{T}_0})_n 
\quad (0\leq m \leq N,n\in {\mab N})  
\tag{1.2.7.2}\label{eqn:iincl}
\end{equation*}   
and 
\begin{equation*} 
\ol{\cal P}_{mn}:={\rm cosk}_0^{\ol{S(T)^{\nat}}}(\ol{\Gam}{}'_m)_n
\quad (0\leq m \leq N,n\in {\mab N}). 
\tag{1.2.7.3}\label{eqn:iigncl}
\end{equation*}  
Then we have an 
$(N,\infty)$-truncated bisimplicial SNCL scheme 
$X_{\bul \leq N,\bul,\os{\circ}{T}_0}$ and an immersion 
\begin{equation*}  
X_{\bul \leq N,\bul,\os{\circ}{T}_0} \os{\sus}{\lo} \ol{\cal P}_{\bul \leq N,\bul} 
\tag{1.2.7.4}\label{eqn:eixd} 
\end{equation*} 
into a log smooth 
$(N,\infty)$-truncated bisimplicial log scheme over $\ol{S(T)^{\nat}}$.  
Thus we have proved the following: 


\par 

\begin{prop}\label{prop:xbn}  
Assume that $X_{\bul \leq N,\os{\circ}{T}_0}$ has 
an affine $N$-truncated simplicial open covering 
of $X'_{\bul \leq N,\os{\circ}{T}_0}$  over $S_{\os{\circ}{T}_0}$ 
for which there exists the commutative diagram {\rm (\ref{cd:xnip})}. 
Then the following hold$:$ 
\par 
$(1)$ There exists an 
immersion 
\begin{equation*}  
\begin{CD} 
X_{\bul \leq N,\bul,\os{\circ}{T}_0} 
@>{\sus}>> \ol{\cal P}_{\bul \leq N,\bul} \\
@VVV @VVV \\
S_{\os{\circ}{T}_0} @>{\subset}>> \ol{S(T)^{\nat}}
\end{CD} 
\tag{1.2.8.1}\label{eqn:eipxd} 
\end{equation*} 
into a log smooth $(N,\infty)$-truncated bisimplicial log scheme 
over $\ol{S(T)^{\nat}}$. 
\par 
$(2)$ There exists an 
immersion 
\begin{equation*}  
\begin{CD} 
X_{\bul \leq N,\bul,\os{\circ}{T}_0} 
@>{\sus}>> {\cal P}_{\bul \leq N,\bul} \\
@VVV @VVV \\
S_{\os{\circ}{T}_0} @>{\subset}>> S(T)^{\nat}
\end{CD} 
\tag{1.2.8.2}\label{eqn:eolnpxd} 
\end{equation*} 
into a log smooth $(N,\infty)$-truncated bisimplicial log scheme over $S(T)^{\nat}$. 
\end{prop} 

\begin{coro}\label{coro:exem} 
Let the notations be as above. 
Then there exist  an 
immersion 
\begin{equation*}  
\begin{CD} 
X_{\bul \leq N,\bul,\os{\circ}{T}_0} @>{\sus}>> 
\ol{\cal P}{}^{\rm ex}_{\bul \leq N,\bul} \\
@VVV @VVV \\
S_{\os{\circ}{T}_0} @>{\subset}>> \ol{S(T)^{\nat}}
\end{CD} 
\tag{1.2.9.1}\label{eqn:eipexxd} 
\end{equation*} 
into an $(N,\infty)$-truncated bisimplicial strict 
semistable log formal scheme over $\ol{S(T)^{\nat}}$ and 
an immersion 
\begin{equation*}  
\begin{CD} 
X_{\bul \leq N,\bul,\os{\circ}{T}_0} @>{\sus}>> 
{\cal P}^{\rm ex}_{\bul \leq N,\bul} \\
@VVV @VVV \\
S_{\os{\circ}{T}_0} @>{\subset}>> S(T)^{\nat}
\end{CD} 
\tag{1.2.9.2}\label{eqn:eoexpxd} 
\end{equation*} 
into 
an $(N,\infty)$-truncated bisimplicial formal SNCL scheme over $S(T)^{\nat}$. 
\end{coro}
\begin{proof} 
(\ref{coro:exem}) follows from 
(\ref{prop:nexeo}). 
\end{proof}

\begin{rema}
(1) Let $(Y_{\bul \leq N},D_{\bul \leq N})$ be 
an $N$-truncated simplicial smooth scheme 
with an $N$-truncated simplicial relative 
SNCD over a scheme $T_0$.  
Let $T_0 \os{\sus}{\lo} T$ be a nil closed immersion.  
Assume that $(Y_{\bul \leq N},D_{\bul \leq N})$ 
has an affine $N$-truncated simplicial open covering.
Then, if we use the same argument as those in 
(\ref{prop:xbn}) and  (\ref{coro:exem}) for 
$(Y_{\bul \leq N},D_{\bul \leq N})$ over $T_0\os{\sus}{\lo} T$,  
we can solve a formal version of \cite[Problem 6.13]{nh3}. That is,  
there exists an admissible immersion 
$(Y_{\bul \leq N,\bul},D_{\bul \leq N,\bul})\os{\sus}{\lo} 
({\cal Y}_{\bul \leq N,\bul},{\cal D}_{\bul \leq N,\bul})$ into 
an $(N,\infty)$-truncated bisimplicial smooth formal scheme 
with an $(N,\infty)$-truncated bisimplicial relative formal SNCD over $T$. 
\par 
(2) In the case $N=0$, 
one can construct an exact immersion 
$X_{0\bul,\os{\circ}{T}_0} \os{\sus}{\lo} \ol{\cal P}{}^{\rm ex}_{0\bul}$ 
by using concrete (log) blow ups 
(cf.~\cite{hdw}, \cite{msemi}, \cite{gkwf}, \cite{nh2}). 
However, as in \cite{ctcs}, 
the concrete construction is not necessary in this book. 
Moreover, one can dispense with 
the complicated local calculations in \cite[\S2]{gkwf} and \cite[\S4]{ctcs} 
if one uses (\ref{prop:adla}) as in (\ref{prop:fsi}), (\ref{prop:nexeo}) 
and (\ref{coro:exem}): (\ref{coro:exem}) simplifies 
the complicated local calculations in [loc.~cit.] remarkably.    
In fact, if we use only concrete log blow ups, 
it seems hopeless to construct the $(N,\infty)$-truncated bisimplicial immersions 
in (\ref{eqn:eipexxd}) and (\ref{eqn:eoexpxd}).  
\end{rema}

\section{Log crystalline de Rham complexes}\label{sec:ldc} 
In this section we give several fundamental properties 
of sheaves of log differential forms 
tensored with the structure sheaves of 
log PD-envelopes. 
First we recall the preweight filtration on 
sheaves of differential forms (\cite{nh3}).
\par 
Let $Y$ be a fine log 
(formal) scheme over a fine log (formal) scheme $T$. 
As in \cite[(4.0.2)]{nh3}, we define the {\it preweight filtration}
$P$ on the sheaf ${\Om}^i_{Y/\os{\circ}{T}}$ $(i\in {\mab N})$ 
of log differential forms on $Y_{\rm zar}$ as follows: 
\begin{equation*} 
P_k{\Om}^i_{Y/\os{\circ}{T}} =
\begin{cases} 
0 & (k<0), \\
{\rm Im}({\Om}^k_{Y/\os{\circ}{T}}{\otimes}_{{\cal O}_Y}
\Om^{i-k}_{\os{\circ}{Y}/\os{\circ}{T}}
\lo {\Om}^i_{Y/\os{\circ}{T}}) & (0\leq k\leq i), \\
{\Om}^i_{Y/\os{\circ}{T}} & (k > i).
\end{cases}
\tag{1.3.0.1}\label{eqn:pkdefpw}
\end{equation*}  
A morphism $g\col Y\lo Z$ of log schemes over 
$\os{\circ}{T}$ 
induces the following morphism of filtered complexes: 
\begin{equation*} 
g^*\col (\Om^{\bul}_{Z/\os{\circ}{T}},P)
\lo 
g_*((\Om^{\bul}_{Y/\os{\circ}{T}},P)).  
\tag{1.3.0.2}\label{eqn:lyzpp}
\end{equation*} 
More generally, for a flat ${\cal O}_Y$-module ${\cal E}$ 
and a flat ${\cal O}_Z$-module ${\cal F}$ with a morphism 
$h \col {\cal F}\lo g_*({\cal E})$ of ${\cal O}_Z$-modules, 
we have the following morphism of filtered complexes: 
\begin{equation*} 
h \col ({\cal F}\otimes_{{\cal O}_Z}\Om^{\bul}_{Z/\os{\circ}{T}},P)
\lo 
g_*(({\cal E}\otimes_{{\cal O}_Y}{\Om}^{\bul}_{Y/\os{\circ}{T}},P)).  
\tag{1.3.0.3}\label{eqn:lyytp}
\end{equation*}

\par 
For a morphism $T'\lo T$ of fine log (formal) schemes, 
set $Y_{T'}:=Y\times_TT'$ and 
let $q\col Y_{T'}\lo Y$ be the first projection. 
For a flat ${\cal O}_{Y_{T'}}$-module ${\cal E}$, 
let $P$ be the induced filtration on 
${\cal E}\otimes_{{\cal O}_{Y_{T'}}}\Om^i_{Y_{T'}/\os{\circ}{T}{}'}$ 
by the filtration $P$ on $\Om^i_{Y_{T'}/\os{\circ}{T}{}'}$: 
\begin{align*} 
P_k({\cal E}\otimes_{{\cal O}_{Y_{T'}}}\Om^i_{Y_{T'}/\os{\circ}{T}{}'})
:={\cal E}\otimes_{{\cal O}_{Y_{T'}}}P_k\Om^i_{Y_{T'}/\os{\circ}{T}{}'}. 
\tag{1.3.0.4}\label{eqn:lyep}
\end{align*}

\begin{prop}\label{prop:bclcd}
Assume that the morphism $T'\lo T$ is solid. 
Then 
\begin{align*} 
P_k\Om^i_{Y_{T'}/\os{\circ}{T}{}'}
=q^*(P_k\Om^i_{Y/\os{\circ}{T}}).
\tag{1.3.1.1}\label{ali:oytt} 
\end{align*}  
\end{prop} 
\begin{proof}
Because the morphism is $T'\lo T$ is solid, 
we have an equality 
$Y_{T'}=Y\times_{\os{\circ}{T}}\os{\circ}{T}{}'$. 
We also have an equality 
$\os{\circ}{Y}_{T'}=\os{\circ}{Y}\times_{\os{\circ}{T}}\os{\circ}{T}{}'$ 
since the morphism $T'\lo T$ is integral (\cite[(4.4)]{klog1}).
Hence (\ref{ali:oytt}) follows from \cite[(1.7) (ii)]{klog1}: 
$\Om^1_{Y_{T'}/\os{\circ}{T}{}'}=q^*(\Om^1_{Y/\os{\circ}{T}})$.   
\end{proof}

\par 
For a finitely generated monoid $Q$ 
and a set $\{q_1,\ldots,q_n\}$ 
$(n\in {\mab Z}_{\geq 1})$ of generators of $Q$, 
we say that $\{q_1,\ldots,q_n\}$ is {\it minimal} 
if $n$ is minimal. 
Assume that $Y=(\os{\circ}{Y},M_Y)$ is 
an fs(=fine and saturated) log (formal) scheme. 
Let $y$ be a point of $\os{\circ}{Y}$. 
Let $m_{1,y},\ldots, m_{r,y}$ be local sections of 
$M_Y$ around $y$ 
whose images in $M_{Y,y}/{\cal O}_{Y,y}^*$ 
form a minimal set of generators of 
$M_{Y,y}/{\cal O}_{Y,y}^*$. 
Let $\os{\circ}{D}(M_Y)_i$ $(1\leq i \leq r)$ 
be the local closed subscheme of 
$\os{\circ}{Y}$ defined by the ideal sheaf generated by 
the image of $m_{i,y}$ in ${\cal O}_{Y,y}$. 
For a nonnegative integer $k$, 
let $\os{\circ}{D}{}^{(k)}(M_Y)$ be the disjoint union of 
the $(k+1)$-fold intersections of 
different $\os{\circ}{D}(M_Y)_i$'s. 
Assume that 
\begin{equation*} 
M_{Y,y}/{\cal O}^*_{Y,y}\simeq {\mab N}^r
\tag{1.3.1.2}\label{eqn:epynr}
\end{equation*}
for any point $y$ of $\os{\circ}{Y}$ 
and for some $r\in {\mab N}$ depending on $y$. 
Since ${\rm Aut}({\mab N}^r)\simeq {\mathfrak S}_r$ 
(as mentioned in \S\ref{sec:snclv}),  
${\rm Aut}(M_{Y,y}/{\cal O}^*_{Y,y})\simeq {\mathfrak S}_r$.  
As a result, the scheme $\os{\circ}{D}{}^{(k)}(M_Y)$ is 
independent of the choice of $m_{1,y},\ldots, m_{r,y}$ 
and it is globalized. 
We denote this globalized scheme 
by the same symbol $\os{\circ}{D}{}^{(k)}(M_Y)$. 
Set $\os{\circ}{D}{}^{(-1)}(M_Y):=\os{\circ}{Y}$ 
and 
$\os{\circ}{D}{}^{(k)}(M_Y):=\emptyset$ for 
$k\in {\mab Z}_{\leq -2}$. 
Let $b^{(k)} \col \os{\circ}{D}{}^{(k)}(M_Y) \lo \os{\circ}{Y}$ 
$(k\in {\mab Z})$ 
be the natural morphism.

\begin{prop}[{\rm {\bf \cite[(4.3)]{nh3}}}]\label{prop:mmoo} 
Let $g\col Y \lo Y'$ be a morphism of 
fs log $($formal$)$ schemes satisfying the condition {\rm (\ref{eqn:epynr})}. 
Let $b'{}^{(k)} \col \os{\circ}{D}{}^{(k)}(M_{Y'}) \lo \os{\circ}{Y}{}'$ 
$(k\in {\mab Z})$ be the analogous morphism to $b^{(k)}$. 
Assume that, for each point $y\in \os{\circ}{Y}$ 
and for each member $m$ of the minimal generators 
of $M_{Y,y}/{\cal O}^*_{Y,y}$, there exists 
a unique member $m'$ of the minimal generators of 
$M_{Y',\os{\circ}{g}(y)}/{\cal O}^*_{Y',\os{\circ}{g}(y)}$ 
such that $g^*(m')\in m^{{\mab Z}_{>0}}$. 
Then there exists a canonical morphism 
$\os{\circ}{g}{}^{(k)}\col \os{\circ}{D}{}^{(k)}(M_Y) \lo 
\os{\circ}{D}{}^{(k)}(M_{Y'})$ fitting 
into the following commutative diagram of schemes$:$
\begin{equation*} 
\begin{CD} 
\os{\circ}{D}{}^{(k)}(M_Y) @>{\os{\circ}{g}{}^{(k)}}>> 
\os{\circ}{D}{}^{(k)}(M_{Y'}) \\ 
@V{b^{(k)}}VV @VV{b'{}^{(k)}}V \\ 
\os{\circ}{Y} @>{\os{\circ}{g}}>> \os{\circ}{Y}{}'. 
\end{CD} 
\tag{1.3.2.1}\label{cd:dmgdm}
\end{equation*} 
\end{prop} 

\begin{rema} 
Because 
the construction of $\os{\circ}{g}{}^{(k)}$ 
is necessary for (\ref{prop:rescos}) below, 
we review it as follows 
(cf.~\cite[p.~614]{stwsl}). 
\par 
First we consider the local case.  
Let $\{m_{i,y}\}_{i=1}^r$ and  
$\{D(M_Y)_i\}_{i=1}^r$ be as above. 
Set $y':=\os{\circ}{g}(y)$. 
Let 
$\{m_{i',y'}\}_{i'=1}^{r'}$ and  
$\{D(M_{Y'})_{i'}\}_{i'=1}^{r'}$  
be the similar objects 
for $Y'$ around $y'$. 
We may assume that 
$g^*(m_{i',y'}) \in m_{i,y}^{{\mab Z}_{>0}}$ 
$(1\leq i\leq r)$.  
This means that $g$ induces a natural local morphism 
$D(M_Y)_{i} \lo D(M_{Y'})_{i}$ $(1\leq i\leq r)$. 
Let $D^{(k)}(M_Y;g)$ be the disjoint union of 
all $(k+1)$-fold intersections of 
$D(M_{Y'})_1,\ldots,D(M_{Y'})_r$. 
Then we have a natural morphism 
$D^{(k)}(M_Y)\lo D^{(k)}(M_Y;g)$. 
By composing this with the natural inclusion morphism 
$D^{(k)}(M_Y;g)\os{\sus}{\lo} D^{(k)}(M_{Y'})$, 
we have a morphism 
$D^{(k)}(M_Y)\lo D^{(k)}(M_{Y'})$. 
This is a desired local morphism. 
This local morphism is compatible with localization. 
Hence it is globalized.  
By this construction, for another analogous morphism 
$g'\col Y'\lo Y''$ to $g$, we have the following equality 
\begin{align*} 
(\os{\circ}{g}{}'\circ \os{\circ}{g})^{(k)}=
(\os{\circ}{g}{}')^{(k)}\circ (\os{\circ}{g})^{(k)}. 
\tag{1.3.3.1}\label{ali:kkrel}
\end{align*} 
We also have an obvious equality 
$\os{\circ}{({\rm id}_Y)}{}^{(k)}={\rm id}_{\os{\circ}{D}{}^{(k)}(M_Y)}$. 
\end{rema}

\par 
As in \cite[(3.1.4)]{dh2} and \cite[(2.2.18)]{nh2}, 
we have an orientation sheaf 
$\vp^{(k)}_{\rm zar}(\os{\circ}{D}(M_Y))$ 
$(k\in {\mab N})$ in $\os{\circ}{D}{}^{(k)}(M_Y)_{\rm zar}$ 
associated to the set $\os{\circ}{D}(M_Y)_i$'s. 
We review it. 
Let $E$ be the set of subsets of $\{\os{\circ}{D}(M_Y)_i\}$ 
consisting of $(k+1)$-elements. 
Set $\vp^{(k)}_{\rm zar}(\os{\circ}{D}(M_Y))
:=\bigwedge^{k+1}({\mab Z}^E)$ if $k\geq 0$, 
$\vp^{(-1)}_{\rm zar}(\os{\circ}{D}(M_Y))={\mab Z}$ 
and $\vp^{(k)}_{\rm zar}(\os{\circ}{D}(M_Y))=0$ 
if $k\leq -2$.  
We call $\vp^{(k)}_{\rm zar}(\os{\circ}{D}(M_Y))$ the 
$(k+1)$-{\it fold} {\it zariskian orientation sheaf} of $Y$. 
We have the following equality:  
\begin{equation*} 
b^{(k)}_*\vp^{(k)}_{\rm zar}(\os{\circ}{D}(M_Y)) 
=\bigwedge^{k+1}(M_Y^{\rm gp}/{\cal O}_Y^*). 
\tag{1.3.3.2}\label{eqn:bkezps}
\end{equation*}

\par 
In the case where $Y=X$ is a (formal) SNCL scheme over 
a (formal) family of log points $S$,  
we denote $\os{\circ}{D}(M_X)_i$, 
$b^{(k)}$ and $\vp^{(k)}_{\rm zar}(\os{\circ}{D}(M_X))$ 
by $\os{\circ}{X}_i$, $a^{(k)}$,  
$\vp^{(k)}_{\rm zar}(\os{\circ}{X}/\os{\circ}{S})$ 
$(k\in {\mab Z}_{\geq 1})$, 
respectively. Note that 
$\vp^{(k)}_{\rm zar}(\os{\circ}{X}/\os{\circ}{S})$ 
depends only on $\os{\circ}{X}$ (not on $X$). 
Note also that $\os{\circ}{D}{}^{(k)}(M_X)$ is nothing but 
$\os{\circ}{X}{}^{(k)}$.  

\par
Let $(T,{\cal J},\del)$ be a fine log PD-scheme  
with quasi-coherent PD-ideal sheaf ${\cal J}$ 
such that a prime number $p$ is locally nilpotent on $T$. 
Set $T_0:=\ul{\rm Spec}^{\log}_T({\cal O}_T/{\cal J})$. 
Let $Y$ (resp.~${\cal Q}$) be 
a log smooth scheme over $T_0$ 
$($resp.~$T)$.  
Let $\iota \col Y \os{\sus}{\lo} {\cal Q}$ be 
an immersion over $T$. 
Let ${\mathfrak E}$ be the log PD-envelope of $\iota$ over $(T,{\cal J},\del)$. 
Let ${\cal Q}^{\rm ex}$ be the exactification of $\iota$. 
Then ${\mathfrak E}$ is the classical PD-envelope of the immersion 
$\os{\circ}{Y}\os{\sus}{\lo}\os{\circ}{\cal Q}{}^{\rm ex}$ over $(\os{\circ}{T},{\cal J},\del)$ 
with the inverse image of the log structure of ${\cal Q}^{\rm ex}$ 
(\cite[(6.6)]{klog1}, \cite[3.20 Remarks 7)]{bob}). 
For any point $y$ of $\os{\circ}{Y}$, 
assume that there exists a chart $(Q \lo M_T, P \lo M_Y, Q 
\os{\rho}{\lo} P)$ of $Y \lo T_0 \os{\subset}{\lo} T$ 
on a neighborhood of $y$ such that 
$\rho$ is injective, 
such that ${\rm Coker}(\rho^{\rm gp})$ is torsion free 
and that the natural homomorphism ${\cal O}_{Y,y} \otimes_{{\mab Z}} 
(P^{{\rm gp}}/Q^{{\rm gp}}) \lo \Om^1_{Y/T_0,y}$ is an isomorphism. 
The following is a generalization of 
\cite[(2.2.17) (1)]{nh2}$:$

\begin{prop}\label{prop:injf}
Set $U:=T$ or $\os{\circ}{T}$. Then the following hold$:$
\par 
$(1)$ The natural morphism 
\begin{equation*} 
{\cal O}_{{\cal Q}^{\rm ex}}
\otimes_{{\cal O}_{\cal Q}} \Om^i_{{\cal Q}/U}\lo 
\Om^i_{{\cal Q}^{\rm ex}/U} 
\quad (i\in {\mab N})
\tag{1.3.4.1}\label{eqn:yxpsnpd}
\end{equation*}
is an isomorphism. 
Consequently the natural morphism 
\begin{equation*} 
{\cal O}_{\mathfrak E}
\otimes_{{\cal O}_{\cal Q}}
\Om^i_{{\cal Q}/U}
\lo {\cal O}_{\mathfrak E}
\otimes_{{\cal O}_{{\cal Q}^{\rm ex}}}
\Om^i_{{\cal Q}^{\rm ex}/U} 
\tag{1.3.4.2}\label{eqn:yxpdsnpd}
\end{equation*}
is an isomorphism. 
\par 
$(2)$ The natural morphism 
\begin{equation*} 
{\cal O}_{\mathfrak E}
\otimes_{{\cal O}_{{\cal Q}^{\rm ex}}}
P_k{\Om}^i_{{\cal Q}^{\rm ex}/\os{\circ}{T}}
\lo 
{\cal O}_{\mathfrak E}
\otimes_{{\cal O}_{{\cal Q}^{\rm ex}}}
{\Om}^i_{{\cal Q}^{\rm ex}/\os{\circ}{T}} 
\quad (i,k\in {\mab Z})
\tag{1.3.4.3}\label{eqn:yxepd}
\end{equation*}
is injective. 
\end{prop}
\begin{proof}
The question is local on $Y$; 
we may assume the existence of 
the commutative diagram (\ref{eqn:0txda}); 
we may also assume the existence of 
the log scheme ${\cal Q}^{\rm prex}$ between 
(\ref{prop:zpex}) and (\ref{rema:pds}). 
\par 
(1):  Because the natural morphism 
${\cal Q}^{\rm prex}\lo {\cal Q}$ is log \'{e}tale, 
${\Om}^i_{{\cal Q}^{\rm prex}/U}
={\cal O}_{{\cal Q}^{\rm prex}}
\otimes_{{\cal O}_{\cal Q}}
{\Om}^i_{{\cal Q}/U}$ $(i\in {\mab N})$. 
Hence 
${\Om}^i_{{\cal Q}^{\rm ex}/U}
={\cal O}_{{\cal Q}^{\rm ex}}
\otimes_{{\cal O}_{{\cal Q}^{\rm prex}}}
{\Om}^i_{{\cal Q}^{\rm prex}/U}
={\cal O}_{{\cal Q}^{\rm ex}}
\otimes_{{\cal O}_{\cal Q}}
{\Om}^i_{{\cal Q}/U}$ and  
${\cal O}_{\mathfrak E}
\otimes_{{\cal O}_{{\cal Q}^{\rm ex}}}
{\Om}^i_{{\cal Q}^{\rm ex}/U} 
=
{\cal O}_{\mathfrak E}\otimes_{{\cal O}_{\cal Q}}
{\Om}^i_{{\cal Q}/U}$.  
\par 
(2): (The proof of (2) is the same 
as that of \cite[(2.2.17) (1)]{nh2}.)   
Let $x_1,\ldots, x_c$ be the coordinates of 
${\mab A}^c_T$ in (\ref{eqn:0txda}). 
Set ${\cal K}:=(x_1,\ldots, x_c){\cal O}_{{\cal Q}^{\rm ex}}$ 
and ${\cal Q}':=\ul{\rm Spec}^{\log}_{{\cal Q}^{\rm ex}}
({\cal O}_{{\cal Q}^{\rm ex}}/{\cal K})$. 
Then ${\cal Q}'$ is a log smooth lift of $Y$ over $T$.  
Since ${\cal Q}'$ is log smooth over $T$, 
there exists a section of the surjection 
${\cal O}_{{\cal Q}^{\rm ex}}/{\cal K}^N \lo 
{\cal O}_{{\cal Q}^{\rm ex}}/{\cal K}={\cal O}_{{\cal Q}'}$. 
Set ${\cal K}_0:=(x_1,\ldots, x_c)$ in 
${\cal O}_{{\cal Q}'}[x_1,\ldots, x_c]$.  
Then, as in \cite[3.32 Proposition]{bob}, we have a morphism 
$${\cal O}_{{\cal Q}'}[x_1,\ldots, x_c] \lo  
{\cal O}_{{\cal Q}^{\rm ex}}/{\cal K}^N$$ 
such that the induced morphism 
${\cal O}_{{\cal Q}'}[x_1,\ldots, x_c]/{\cal K}_0^N 
\lo {\cal O}_{{\cal Q}^{\rm ex}}/{\cal K}^N$ is an isomorphism. 
Set ${\cal Q}'':={\mab A}^c_T$.  
Because $p$ is locally nilpotent on $T$, we may assume that 
there exists a positive integer $N$ such that 
${\cal K}^N{\cal O}_{\mathfrak E}=0$. 
By \cite[3.32 Proposition]{bob}, 
${\cal O}_{\mathfrak E}$ is isomorphic to 
the PD-polynomial algebra 
${\cal O}_{{\cal Q}'}\langle x_1,\ldots, x_c\rangle$. 
Hence we have the following isomorphisms  
\begin{equation*} 
{\cal O}_{\mathfrak E}
\otimes_{{\cal O}_{{\cal Q}^{\rm ex}}}
\Om^i_{{\cal Q}^{\rm ex}/\os{\circ}{T}}\os{\sim}{\lo} 
\bigoplus_{i'+i''=i}
\Om^{i'}_{{\cal Q}'/\os{\circ}{T}}
\otimes_{{\cal O}_T}{\cal O}_T\langle x_1,\ldots, x_c\rangle 
\otimes_{{\cal O}_{{\cal Q}''}}
\Om^{i''}_{\os{\circ}{\cal Q}{}''/\os{\circ}{T}}  
\tag{1.3.4.4}\label{eqn:dopx} 
\end{equation*}
and 
\begin{equation*} 
{\cal O}_{\mathfrak E}
\otimes_{{\cal O}_{{\cal Q}^{\rm ex}}}
P_k{\Om}^i_{{\cal Q}^{\rm ex}/\os{\circ}{T}}
\os{\sim}{\lo} 
\bigoplus_{i'+i''=i}
P_k{\Om}^{i'}_{{\cal Q}'/\os{\circ}{T}}
\otimes_{{\cal O}_T}
{\cal O}_T\langle x_1,\ldots, x_c\rangle 
\otimes_{{\cal O}_{{\cal Q}''}}
\Om^{i''}_{\os{\circ}{\cal Q}{}''/\os{\circ}{T}}.  
\tag{1.3.4.5}\label{eqn:dpopx} 
\end{equation*}
Since the complex 
${\cal O}_T\langle x_1,\ldots, x_c\rangle 
\otimes_{{\cal O}_{{\cal Q}''}}
\Om^{\bul}_{\os{\circ}{\cal Q}{}''/\os{\circ}{T}}$    
consists of free ${\cal O}_T$-modules, 
we obtain the desired injectivity. 
\end{proof}

\begin{coro}[{\bf \cite[Lemma 2.22]{hk}}]\label{coro:gmn}
Assume that the morphism $Y\lo T_0$ is integral. 
Then 
${\cal O}_{\mathfrak E}
\otimes_{{\cal O}_{{\cal Q}^{\rm ex}}}\Om^i_{{\cal Q}^{\rm ex}/T}$ 
is a sheaf of flat ${\cal O}_T$-modules.  
\end{coro}  
\begin{proof}
The question is local on $Y$. 
Hence we can use the argument in the proof of (\ref{prop:injf}) (2).  
Because the morphism $Y\lo T$ is integral, 
so is the morphism ${\cal Q}'\lo T$. 
Because this morphism is log smooth and integral, 
it is flat by \cite[(4.5)]{klog1}. 
Since the complex 
${\cal O}_T\langle x_1,\ldots, x_c\rangle 
\otimes_{{\cal O}_{{\cal Q}''}}\Om^{\bul}_{{\cal Q}{}''/T}$    
consists of free ${\cal O}_T$-modules,  
we obtain (\ref{coro:gmn}) by the analogue of (\ref{eqn:dopx}): 
\begin{equation*} 
{\cal O}_{\mathfrak E}
\otimes_{{\cal O}_{{\cal Q}^{\rm ex}}}
\Om^i_{{\cal Q}^{\rm ex}/T}\os{\sim}{\lo} 
\bigoplus_{i'+i''=i}
\Om^{i'}_{{\cal Q}'/T}
\otimes_{{\cal O}_T}{\cal O}_T\langle x_1,\ldots, x_c\rangle 
\otimes_{{\cal O}_{{\cal Q}''}}\Om^{i''}_{{\cal Q}''/T}.   
\tag{1.3.5.1}\label{eqn:dopcwx} 
\end{equation*}
\end{proof}

\begin{defi} 
(1) We denote by $P$ the filtration 
on ${\cal O}_{\mathfrak E}
\otimes_{{\cal O}_{{\cal Q}^{\rm ex}}}
\Om^{\bul}_{{\cal Q}^{\rm ex}/\os{\circ}{T}}$
defined by the family 
$\{{\cal O}_{\mathfrak E}
\otimes_{{\cal O}_{{\cal Q}^{\rm ex}}}
P_k{\Om}^{\bul}_{{\cal Q}^{\rm ex}/\os{\circ}{T}}\}_{k\in {\mab Z}}$ 
of subcomplexes of 
${\cal O}_{\mathfrak E}\otimes_{{\cal O}_{{\cal Q}^{\rm ex}}}
{\Om}^{\bul}_{{\cal Q}^{\rm ex}/\os{\circ}{T}}:$ 
$P_k({\cal O}_{\mathfrak E}
\otimes_{{\cal O}_{{\cal Q}^{\rm ex}}}
{\Om}^{\bul}_{{\cal Q}^{\rm ex}/\os{\circ}{T}})
={\cal O}_{\mathfrak E}
\otimes_{{\cal O}_{{\cal Q}^{\rm ex}}}
P_k{\Om}^{\bul}_{{\cal Q}^{\rm ex}/\os{\circ}{T}}$. 
We call $P$ the {\it preweight filtration} on 
${\cal O}_{\mathfrak E}
\otimes_{{\cal O}_{{\cal Q}^{\rm ex}}}
{\Om}^{\bul}_{{\cal Q}^{\rm ex}/\os{\circ}{T}}$. 
\par 
(2) Let ${\cal E}$ be a flat ${\cal O}_{\mathfrak E}$-module. 
We denote by $P$ the filtration on 
${\cal E}\otimes_{{\cal O}_{{\cal Q}^{\rm ex}}}
\Om^{\bul}_{{\cal Q}^{\rm ex}/\os{\circ}{T}}$ induced by the filtration $P$ in (1). 
We call $P$ the {\it Poincar\'{e} filtration} on ${\cal E}
\otimes_{{\cal O}_{{\cal Q}^{\rm ex}}}{\Om}^{\bul}_{{\cal Q}^{\rm ex}/\os{\circ}{T}}$. 
\end{defi}

\begin{rema}\label{rema:exp}
(1) Let $Y$ be a fine log scheme over 
a closed log subscheme of $T$. 
Let $Y\os{\sus}{\lo} {\cal Q}$ be an immersion over $T$. 
As in \cite[(6.2)]{nhw}, 
consider the case where the immersion 
$Y\os{\sus}{\lo} {\cal Q}$ is not exact. 
Let $P^{\cal Q}$ (resp.~$P$) be the filtration  
on ${\Om}^i_{{\cal Q}/\os{\circ}{T}}$ 
(resp.~${\Om}^i_{{\cal Q}^{\rm ex}/\os{\circ}{T}}$) 
defined by the formula (\ref{eqn:pkdefpw}). 
It is natural to ask whether 
the following equality  
\begin{equation*} 
P=P^{\cal Q}\otimes_{{\cal O}_{\cal Q}}
{\cal O}_{{\cal Q}^{\rm ex}}
\tag{1.3.7.1}\label{eqn:tensfil}
\end{equation*} 
holds. 
However this does not hold in general. 
Indeed, let $\kap$ be a field and 
let $s:=({\rm Spec}(\kap),{\mab N}\oplus \kap^*)$ 
be the log point. 
Consider the case $T=s$, 
$Y= {\mab A}_s(2,2)$, 
$${\cal Q} = Y\times_sY
=({\rm Spec}(\kap[x,y,x',y']/(xy,x'y')),
({\cal O}^*_{\cal Q}x^{\mab N}y^{\mab N}x'{}^{\mab N}
y'{}^{\mab N})/(xy=x'y'))$$  
and consider the diagonal immersion 
$Y\os{\sus}{\lo} {\cal Q}$. 
Then 
${\cal Q}^{\rm ex}=
({\rm Spec}(\kap[x,y][[u-1,v-1]]/(xy)),
{\cal O}_{{\cal Q}^{\rm ex}}^*x^{\mab N}y^{\mab N})$ 
($x'=xu$, $y'=yv$). 
By the definition of $P$, 
$$P^{\cal Q}_1{\Om}^2_{{\cal Q}/\os{\circ}{T}} = 
{\rm Im}({\Om}^1_{{\cal Q}/\os{\circ}{T}}
{\otimes}_{{\cal O}_{\cal Q}}
{\Om}^1_{\os{\circ}{\cal Q}/\os{\circ}{T}} 
\lo {\Om}^2_{{\cal Q}/\os{\circ}{T}})$$ 
is generated by 
$x d\log x\wedge d\log x'$, 
$x' d\log x\wedge d\log x'$,
$xd\log x \wedge d\log y$, 
$yd\log x \wedge d\log y$, 
$x d\log x\wedge d\log y'$, 
$y' d\log x\wedge d\log y'$,
$x' d\log x'\wedge d\log y$,
$y d\log x'\wedge d\log y$, 
$x' d\log x'\wedge d\log y'$, 
$y' d\log x'\wedge d\log y'$,
$yd\log y \wedge d\log y'$ 
and  
$y'd\log y \wedge d\log y'$. 
By using the relation $x'=xu$, 
we see that 
\begin{align*}
d\log x\wedge d\log u=d\log x\wedge d\log x'\not\in &~{\rm Im}(P^{\cal Q}_1
{\Om}^2_{{\cal Q}/\os{\circ}{T}}
\otimes_{{\cal O}_{{\cal Q}}} 
{\cal O}_{{\cal Q}^{\rm ex}}\lo 
{\Om}^2_{{\cal Q}^{\rm ex}/\os{\circ}{T}})\\
&=P^{\cal Q}_1{\Om}^2_{{\cal Q}/\os{\circ}{T}}
\otimes_{{\cal O}_{{\cal Q}}} 
{\cal O}_{{\cal Q}^{\rm ex}};   
\end{align*} 
but 
\begin{equation*}
d\log x\wedge d\log u\in P_1
{\Om}^2_{{\cal Q}^{\rm ex}/\os{\circ}{T}}. 
\end{equation*}
\parno 
In this way, we think that the filtration 
$P^{\cal Q}$ is not good when 
$Y\os{\sus}{\lo} {\cal Q}$ is not exact. 
Furthermore, in the case where $Y/T$ is an SNCL scheme, 
we do not know whether ${\rm gr}_k^{P^{\cal Q}}(\Om^{\bul}_{{\cal Q}/T})$ 
$(k\in {\mab Z})$ is expressed in terms of 
de Rham complexes of smooth schemes 
over $\os{\circ}{T}$  
when $Y\os{\sus}{\lo} {\cal Q}$ is not exact. 
In these two reasons, we do not consider 
$P^{\cal Q}$ any more in this book 
when this immersion is not exact.  
\par 
(2) If the immersion $Y\os{\sus}{\lo}{\cal Q}$ 
is exact, then it is clear that (\ref{eqn:tensfil}) holds. 
\end{rema}

\begin{defi}\label{defi:fldpd} 
Let $(T,{\cal J},\del)$ be a fine log PD-scheme 
such that $p$ is also locally nilpotent on $T$.  
Assume that ${\cal J}$ is quasi-coherent. 
Set $T_0:=\ul{\rm Spec}^{\log}_T({\cal O}_T/{\cal J})$.  
Let $Z$ be a fine log scheme. 
We call $(T,{\cal J},\del)$ with morphism $T_0\lo Z$ (or simply $(T,{\cal J},\del)$)  
a {\it log PD-enlargement} of $Z$ (cf.~\cite{ollc}). 
Let $Z'\lo Z$ be a morphism from a fine log scheme.  
We define a {\it morphism of log PD-enlargements} over $Z'\lo Z$ in an obvious way. 
\end{defi}

\par 
Let $S$ be a family of log points.  
Let $(T,{\cal J},\del)$ and $T_0$ be as in 
(\ref{defi:fldpd}) for the case $Z=S$.   
Assume that we are given a morphism 
$T_0\lo S$ of fine log schemes.  
Set ${S}_{\os{\circ}{T}_0}:={S}\times_{\os{\circ}{S}}\os{\circ}{T}_0$. 
By abuse of notation, let $M_{S_{\os{\circ}{T}_0}}$ be the inverse image of 
the log structure $M_{{S}_{\os{\circ}{T}_0}}$ 
by the morphism $T_0\lo  {S}_{\os{\circ}{T}_0}$. 
Then the natural morphism $M_{S_{\os{\circ}{T}_0}}\lo M_{T_0}$ is injective 
by (\ref{prop:oot}) (2). 
Let $M$ be the canonical extension of $M_{S_{\os{\circ}{T}_0}}$ to $\os{\circ}{T}$ ((\ref{defi:ceivt})), 
that is, $M$ is the inverse image of 
${\rm Im}(M_{S_{\os{\circ}{T}_0}}\lo M_{T_0})/{\cal O}^*_{T_0}$ 
by the following composite morphism 
\begin{align*} 
M_T\lo M_T/{\cal O}^*_T\os{\sim}{\lo}  M_{T_0}/{\cal O}^*_{T_0}. 
\tag{1.3.8.1}\label{ali:mtt0} 
\end{align*} 
Then $M$ is a fine log structure on $\os{\circ}{T}$ ((\ref{lemm:cext})). 
Let $S(T)$ be the resulting log scheme.  


\parno 
If $T$ is restrictively hollow with respective to the morphism 
$T_0\lo S$ ((\ref{defi:nib})), then $S(T)$ is a family of log points 
and $S_{\os{\circ}{T}}=\ul{\rm Spec}_{S(T)}^{\log}({\cal O}_{S(T)}/{\cal J})$;   
we also have the following natural morphism of log PD-schemes:  
\begin{align*} 
(T,{\cal J},\del)\lo(S(T),{\cal J},\del).
\tag{1.3.8.2}\label{eqn:jsss}
\end{align*}   




\par 
Because we do not assume that $S(T)$ is not hollow,  
we consider the hollowing out $S(T)^{\nat}$ and 
the log scheme $\ol{S(T)^{\nat}}$.  
Let ${\mathfrak D}(\ol{S(T)^{\nat}})$ be the log PD-envelope 
of the immersion $S(T)^{\nat}\os{\sus}{\lo}\ol{S(T)^{\nat}}$ 
over $(\os{\circ}{T},{\cal J},\del)$;   
${\mathfrak D}(\ol{S(T)^{\nat}})$ is the classical PD-envelope of 
the immersion $\os{\circ}{T}\os{\sus}{\lo}\os{\circ}{\ol{S(T)^{\nat}}}$ 
over $(\os{\circ}{T},{\cal J},\del)$ 
endowed with the inverse image of the log structure of $\ol{S(T)^{\nat}}$. 
Locally on $\ol{S(T)^{\nat}}$, ${\cal O}_{{\mathfrak D}(\ol{S(T)^{\nat}})}
\simeq {\cal O}_T\langle \tau \rangle$. 
Hence ${\cal O}_{{\mathfrak D}(\ol{S(T)^{\nat}})}$ is a locally-free 
${\cal O}_T$-module. 
By killing the PD-ideal sheaf of ${\cal O}_{{\mathfrak D}(\ol{S(T)^{\nat}})}$, 
we obtain a well-defined surjective morphism 
${\cal O}_{{\mathfrak D}(\ol{S(T)^{\nat}})}\lo {\cal O}_T$. 
This morphism gives us an exact immersion 
$S(T)^{\nat}\os{\sus}{\lo} {\mathfrak D}(\ol{S(T)^{\nat}})$. 
\par 
Let $Y$ be a log smooth scheme over $S$. 
Let $U_0$ (resp.~$U$) be $S_{\os{\circ}{T}_0}$ or $T_0$ 
(resp.~$S(T)^{\nat}$ or $T$).  
When $U=T$ (and $U\not=S(T)^{\nat}$), 
we always assume that $S(T)$ is hollow. 
Set $Y_{U_0}:=Y\times_{S}U_0$ and 
$\os{\circ}{Y}_{U_0}:=\os{\circ}{(Y_{U_0})}$. 
Then $Y_{S_{\os{\circ}{T}_0}}=Y\times_{S}S_{\os{\circ}{T}_0}=
Y\times_{\os{\circ}{S}}\os{\circ}{T}_0
=:Y_{\os{\circ}{T}_0}$ 
(we can consider $Y$ as a fine log scheme over $\os{\circ}{S}$) 
and $Y_{T_0}= Y_{S_{\os{\circ}{T}_0}}\times_{S_{\os{\circ}{T}_0}}T_0$. 
\par 
Let $Y_{\os{\circ}{T}_0}\os{\sus}{\lo} \ol{\cal Q}$ be 
an immersion into a log smooth scheme over $\ol{S(T)^{\nat}}$ 
with structural morphism $\ol{g}\col \ol{\cal Q}\lo \ol{S(T)^{\nat}}$. 
Let $\ol{\cal Q}{}^{\rm ex}$ be the exactification of the immersion 
$Y_{\os{\circ}{T}_0}\os{\sus}{\lo} \ol{\cal Q}$. 
By abuse of notation, we also denote the structural morphism 
$\ol{\cal Q}{}^{\rm ex}\lo \ol{S(T)^{\nat}}$ by $\ol{g}$.  
Set ${\cal Q}:=\ol{\cal Q}\times_{\ol{S(T)^{\nat}}}S(T)^{\nat}$. 
Let ${\cal Q}^{\rm ex}$ be the exactification of the immersion 
$Y_{\os{\circ}{T}_0}\os{\sus}{\lo} {\cal Q}$.  
Then ${\cal Q}^{\rm ex}=\ol{\cal Q}{}^{\rm ex}\times_{\ol{S(T)^{\nat}}}S(T)^{\nat}$ 
by (\ref{prop:xpls}). 
Let $\ol{\mathfrak E}$ and ${\mathfrak E}$ 
be the log PD-envelopes of 
the immersions $Y_{\os{\circ}{T}_0}\os{\sus}{\lo}\ol{\cal Q}$ and 
$Y_{\os{\circ}{T}_0}\os{\sus}{\lo} {\cal Q}$ over 
$(\os{\circ}{T},{\cal J},\del)$ and $(S(T)^{\nat},{\cal J},\del)$, 
respectively.  
Then there exists a natural morphism 
$\ol{\mathfrak E}\lo {\mathfrak D}(\ol{S(T)^{\nat}})$ by the universality 
of the log PD-envelopes and the following equality holds:
\begin{equation*}
{\mathfrak E}=\ol{\mathfrak E}\times_{{\mathfrak D}(\ol{S(T)^{\nat}})}S(T)^{\nat}.
\tag{1.3.8.3}\label{eqn:ddss}
\end{equation*}  

We restate (\ref{prop:injf}) (2) for the case above because 
the restatement will be important for the definition of $(A_{\rm zar},P)$ 
in the next section. 

\begin{prop}\label{prop:om}
For any point $y$ of $\os{\circ}{Y}_{T_0}$, 
assume that there exists a chart 
$(Q \lo M_{S(T)^{\nat}},P \lo M_{Y_{\os{\circ}{T}_0}},Q \os{\rho}{\lo} P)$ 
of $Y_{\os{\circ}{T}_0} \lo S_{\os{\circ}{T}_0} \os{\subset}{\lo} S(T)^{\nat}$ 
on a neighborhood of $y$ such that 
$\rho$ is injective, 
such that ${\rm Coker}(\rho^{\rm gp})$ is torsion free 
and that the natural homomorphism 
${\cal O}_{Y_{\os{\circ}{T}_0},y} \otimes_{{\mab Z}} 
(P^{{\rm gp}}/Q^{{\rm gp}}) \lo \Om^1_{Y_{\os{\circ}{T}_0}/S_{\os{\circ}{T}_0},y}$ 
is an isomorphism. 
Then the following natural morphism 
\begin{align*} 
{\cal O}_{{\mathfrak E}}
\otimes_{{\cal O}_{{\cal Q}^{\rm ex}}}
P_k\Om^{\bul}_{{\cal Q}{}^{\rm ex}/\os{\circ}{T}}
\lo 
{\cal O}_{{\mathfrak E}}\otimes_{{\cal O}_{{\cal Q}^{\rm ex}}}
\Om^{\bul}_{{\cal Q}{}^{\rm ex}/\os{\circ}{T}}   
\tag{1.3.9.1}\label{eqn:ydpd}
\end{align*} 
is injective. 
\end{prop} 

\begin{rema}  
(1) Assume that $T$ is restrictively hollow with respective to the morphism $T_0\lo S$. 
Then set ${\cal Q}{}^{\rm ex}_T:={\cal Q}{}^{\rm ex}\times_{S(T)}T$. 
The reader should note that, though we shall consider 
$\Om^{\bul}_{{\cal Q}{}^{\rm ex}_T/T}$ 
in this book,  we do not consider 
$\Om^{\bul}_{{\cal Q}{}^{\rm ex}_T/\os{\circ}{T}}$ at all because 
$\Om^i_{{\cal Q}{}^{\rm ex}_T/\os{\circ}{T}}$ $(i\in {\mab Z})$ is a too big 
sheaf in this book by the reason that the log structure of $T$ is 
bigger than that of $S(T)$ in general. 
\par 
(2) In (\ref{rema:noninj}) below we shall show that the following natural morphism 
\begin{align*} 
{\cal O}_{{\mathfrak E}}\otimes_{{\cal O}_{\ol{\cal Q}{}^{\rm ex}}}
P_k\Om^{\bul}_{\ol{\cal Q}{}^{\rm ex}/\os{\circ}{T}} 
\lo 
{\cal O}_{{\mathfrak E}}
\otimes_{{\cal O}_{\ol{\cal Q}{}^{\rm ex}}}
\Om^{\bul}_{\ol{\cal Q}{}^{\rm ex}/\os{\circ}{T}}   
\tag{1.3.10.1}\label{eqn:ydqpd}
\end{align*} 
is not necessarily injective. 
\end{rema}

\begin{prop}\label{prop:ed0}
Let $\star$ be ex or nothing. Then the following hold$:$
\par  
$(1)$ 
The morphism $\ol{\cal Q}{}^{\star}\lo \os{\circ}{T}$ is $($formally$)$ log smooth 
and the morphism $\os{\circ}{\ol{\cal Q}}{}^{\star}\lo \os{\circ}{\ol{T}}$ is flat. 
\par 
$(2)$ The natural morphism 
\begin{equation*} 
{\cal O}_{\ol{\mathfrak E}}
\otimes_{{\cal O}_{\ol{\cal Q}}}
{\Om}^i_{\ol{\cal Q}/\os{\circ}{T}}\lo 
{\cal O}_{\ol{\mathfrak E}}
\otimes_{{\cal O}_{\ol{\cal Q}^{\rm ex}}}
{\Om}^i_{\ol{\cal Q}^{\rm ex}/\os{\circ}{T}} 
\quad (i\in {\mab N})
\tag{1.3.11.1}\label{eqn:yxppnpd}
\end{equation*}
is an isomorphism. 
\par 
$(3)$  
The sheaf ${\Om}^i_{\ol{\cal Q}{}^{\star}/\os{\circ}{T}}$ 
$(i\in {\mab N})$ is a locally free 
${\cal O}_{\ol{\cal Q}{}^{\star}}$-module of finite rank. 
\end{prop} 
\begin{proof} 
(1): 
For any point $t\in \os{\circ}{\ol{S(T)^{\nat}}}$, 
$(M_{\ol{S(T)^{\nat}}}/{\cal O}^*_{\ol{S(T)^{\nat}}})_t$ is generated by an element. 
By \cite[(4.4)]{klog1} the morphism 
$\ol{\cal Q}{}^{\star}\lo \ol{S(T)^{\nat}}$ is integral. 
Since the morphism $\ol{S(T)^{\nat}}\lo \os{\circ}{T}$ is integral,  
the composite morphism 
$\ol{\cal Q}{}^{\star}\lo \os{\circ}{T}$ is also integral (\cite[(4.3.1)]{klog1}). 
Because $\ol{\cal Q}{}^{\star}/\ol{S(T)^{\nat}}$ 
and $\ol{S(T)^{\nat}}\lo \os{\circ}{T}$ 
are  (formally) log smooth,  so is the composite morphism $\ol{\cal Q}{}^{\star}/\os{\circ}{T}$. 
Now (1) follows from \cite[(4.5)]{klog1}. 
\par 
(2):  (2) is a special case of (\ref{prop:injf}) (1).  
\par 
(3): Because $\ol{\cal Q}{}^{\star}/\os{\circ}{T}$ is (formally) log smooth, 
(3) follows from \cite[Proposition 3.10]{klog1}. 
\end{proof}

\par
Next,  in (\ref{prop:incol}) and (\ref{rema:noninj}),  
we discuss delicate problems about the injectivities of 
morphisms of (log) de Rham complexes:

\begin{prop}\label{prop:incol}
Let the notations and the assumptions be as in {\rm (\ref{prop:ed0})}.   
Then the following hold$:$
\par 
$(1)$ Set 
$\Om^{\bul}_{{\ol{\cal Q}{}^{\star}}/\os{\circ}{T}}
(-\os{\circ}{\cal Q}{}^{\star}):={\cal I}_{\ol{\cal Q}{}^{\star}}
\otimes_{{\cal O}_{\ol{\cal Q}{}^{\star}}}
\Om^{\bul}_{\ol{\cal Q}{}^{\star}/\os{\circ}{T}}$ 
$($recall the ideal sheaf ${\cal I}_{\ol{\cal Q}{}^{\star}}$ defined before 
{\rm (\ref{prop:bar})}$)$.  
Then 
$\Om^{\bul}_{{\ol{\cal Q}{}^{\star}}/\os{\circ}{T}}(-\os{\circ}{\cal Q}{}^{\star})$ is 
a subcomplex of $\Om^{\bul}_{{\ol{\cal Q}{}^{\star}}/\os{\circ}{T}}$. 
\par 
$(2)$ 
Assume that, for any point 
$x\in \os{\circ}{\ol{\cal Q}}{}^{\star}$, 
there exists a basis $\{\ol{m}_i\}_{i=1}^k$ 
of $(M_{\ol{\cal Q}{}^{\star}}
/{\cal O}_{\ol{\cal Q}{}^{\star}}^*)_x$ 
and a positive integer $e_i>0$ $(1\leq i\leq k)$
such that $\prod_{i=1}^k\ol{m}{}^{e_i}_i$ is the image of 
the generator 
of $(\ol{g}{}^*(M_{\ol{S(T)^{\nat}}}/{\cal O}^*_{\ol{S(T)^{\nat}}}))_{x}$ 
in $(M_{\ol{\cal Q}{}^{\star}}
/{\cal O}_{\ol{\cal Q}{}^{\star}}^*)_x$. 
Then $\Om^{\bul}_{{\ol{\cal Q}{}^{\star}/\os{\circ}{T}}}
(-\os{\circ}{\cal Q}{}^{\star})$ 
is a subcomplex of 
$P_0\Om^{\bul}_{\ol{\cal Q}{}^{\star}/\os{\circ}{T}}$.  
\par 
$(3)$ Let the assumption be as in $(2)$. 
Then the following formula holds$:$  
\begin{equation*}
P_k\Om^{\bul}_{{\cal Q}^{\rm ex}/\os{\circ}{T}}=
P_k\Om^{\bul}_{\ol{\cal Q}{}^{\rm ex}/\os{\circ}{T}}
/\Om^{\bul}_{\ol{\cal Q}{}^{\rm ex}/\os{\circ}{T}}
(-\os{\circ}{\cal Q}{}^{\rm ex}) \quad (k\in {\mab N}).
\tag{1.3.12.1}\label{eqn:pops}
\end{equation*}   
\par 
$(4)$ Let the assumptions be as in $(2)$ without assuming that 
$\prod_{i=1}^k\ol{m}{}^{e_i}_i$ is the image of 
the generator 
of $(\ol{g}{}^*(M_{\ol{S(T)^{\nat}}}/{\cal O}^*_{\ol{S(T)^{\nat}}}))_{x}$ 
in $(M_{\ol{\cal Q}{}^{\star}}
/{\cal O}_{\ol{\cal Q}{}^{\star}}^*)_x$. 
Let $m_i$ be a lift of $\ol{m}_i$ to $M_{\ol{\cal Q}{}^{\star},x}$ 
$(1\leq i\leq k)$. 
Assume furthermore that 
$\os{\circ}{\ol{\cal Q}}{}^{\star}/\os{\circ}{T}$ 
is smooth and that $\{d\al(m_i)\}_{i=1}^k$ is a part of a basis of 
$\Om^1_{\os{\circ}{\ol{\cal Q}}{}^{\star}/\os{\circ}{T}}$. 
Then the following natural morphism 
\begin{equation*} 
\Om^{\bul}_{\os{\circ}{\ol{\cal Q}}{}^{\star}/\os{\circ}{T}} 
\lo \Om^{\bul}_{\ol{\cal Q}{}^{\star}/\os{\circ}{T}}    
\tag{1.3.12.2}\label{eqn:yxppdppd}
\end{equation*} 
is injective. 
\par 
$(5)$ Let the assumptions be as in $(2)$ and $(4)$.  
Then the injective morphism 
$\Om^{\bul}_{{\ol{\cal Q}{}^{\star}}/\os{\circ}{T}}
(-\os{\circ}{\cal Q}{}^{\star}) \lo 
\Om^{\bul}_{{\ol{\cal Q}{}^{\star}}/\os{\circ}{T}}$ 
factors through the following injective morphism$:$
\begin{equation*}  
{\Om}^{\bul}_{\ol{\cal Q}{}^{\star}/\os{\circ}{T}}
(-\os{\circ}{\cal Q}{}^{\star}) \lo 
\Om^{\bul}_{\os{\circ}{\ol{\cal Q}}{}^{\star}/\os{\circ}{T}}.    
\tag{1.3.12.3}\label{eqn:yxpplnpd}
\end{equation*} 
\par 
\end{prop}
\begin{proof} 
(1): Let $x$ be a point of $\os{\circ}{\ol{\cal Q}}{}^{\star}$.  
Let $\tau$ be a local section of $M_{\ol{S(T)^{\nat}},g(x)}$ which 
gives a generator of $(M_{\ol{S(T)^{\nat}}}/{\cal O}_{\ol{S(T)^{\nat}}}^*)_{g(x)}$. 
By abuse of notation, we denote the image of $\tau$ 
in ${\cal O}_{\ol{S(T)^{\nat}},g(x)}$ by $\tau$.  
The morphism 
$\os{\circ}{\ol{\cal Q}}{}^{\star}\lo  \os{\circ}{\ol{S(T)^{\nat}}}$ is flat. 
Consequently the natural morphism 
${\cal I}_{\ol{\cal Q}{}^{\star}}
\lo {\cal O}_{\ol{\cal Q}{}^{\star}}$ is injective. 
By this injectivity and (\ref{prop:ed0}) (3), the natural morphism 
$\Om^i_{{\ol{\cal Q}}{}^{\star}/\os{\circ}{T}}
(-\os{\circ}{\cal Q}{}^{\star}) 
\lo \Om^i_{\ol{\cal Q}{}^{\star}/\os{\circ}{T}}$ 
$(i\in {\mab N})$ is injective. 
For a local section 
$\om \in {\Om}^{\bul}_{\ol{\cal Q}{}^{\star}/\os{\circ}{T}}$ 
around $x$, 
$d(\tau \om)=\tau d\log \tau\wedge \om+\tau d\om$. 
Hence 
$\Om^{\bul}_{{\ol{\cal Q}}{}^{\star}/\os{\circ}{T}}
(-\os{\circ}{\cal Q}{}^{\star})$ 
is a subcomplex of 
$\Om^{\bul}_{\ol{\cal Q}{}^{\star}/\os{\circ}{T}}$.  
\par 
(2): 
The question is local.
Let $\al \col M_{\ol{\cal Q}{}^{\star},x}\lo 
{\cal O}_{\ol{\cal Q}{}^{\star},x}$ 
be the structural morphism. 
Let $m_i$ be a lift of $\ol{m}_i$ to 
$M_{\ol{\cal Q}{}^{\star},x}$ as stated in (4). 
Consider a section 
$d\log m_{i_1}\wedge \cdots 
\wedge d\log m_{i_l}\wedge \om$ 
$(1\leq i_1< \cdots <i_{l}\leq k)$ 
with $\om\in 
P_0\Om^j_{\ol{\cal Q}{}^{\star}/\os{\circ}{T}}$ 
$(j\in {\mab N})$.  
By the assumption, we may assume that 
$\prod_{i=1}^k\al(m_i^{e_i})=\tau$. 
Then 
$$\tau d\log m_{i_1}\wedge \cdots 
\wedge d\log m_{i_l}\wedge \om
=d\al(m_{i_1})\wedge \cdots 
\wedge d\al(m_{i_l})\wedge \om'$$  
with $\om'\in 
P_0\Om^{j-l}_{\ol{\cal Q}{}^{\star}/\os{\circ}{T}}$. 
Hence 
$\Om^{\bul}_{{\ol{\cal Q}{}^{\star}}/\os{\circ}{T}}
(-\os{\circ}{\cal Q}{}^{\star})$
is a subcomplex of $P_0\Om^{\bul}_{\ol{\cal Q}{}^{\star}/\os{\circ}{T}}$.  
\par 
(3): By (2) the complex $\Om^{\bul}_{{\ol{\cal Q}{}^{\star}}/\os{\circ}{T}}
(-\os{\circ}{\cal Q}{}^{\star})$
is a subcomplex of $P_k\Om^{\bul}_{\ol{\cal Q}{}^{\star}/\os{\circ}{T}}$ $(k\in {\mab N})$.  
By the definition of 
$P_k{\Om}^{\bul}_{{\cal Q}^{\rm ex}/\os{\circ}{T}}$, 
we have the following natural surjective morphism 
\begin{equation*} 
P_k{\Om}^{\bul}_{\ol{\cal Q}{}^{\rm ex}/\os{\circ}{T}}
/{\Om}^{\bul}_{\ol{\cal Q}{}^{\rm ex}/\os{\circ}{T}}
(-\os{\circ}{\cal Q}{}^{\rm ex}) 
 \lo 
P_k{\Om}^{\bul}_{{\cal Q}^{\rm ex}/\os{\circ}{T}}.
\end{equation*}  
The following diagram shows that this morphism 
is injective: 
\begin{equation*} 
\begin{CD} 
P_k{\Om}^{\bul}_{\ol{\cal Q}{}^{\rm ex}/\os{\circ}{T}}/
{\Om}^{\bul}_{\ol{\cal Q}{}^{\rm ex}/\os{\circ}{T}}
(-\os{\circ}{\cal Q}{}^{\rm ex})  
@>>> 
P_k{\Om}^{\bul}_{{\cal Q}^{\rm ex}/\os{\circ}{T}}\\ 
@V{\bigcap}VV  @VV{\bigcap}V \\ 
{\Om}^{\bul}_{\ol{\cal Q}{}^{\rm ex}/\os{\circ}{T}}/
{\Om}^{\bul}_{\ol{\cal Q}{}^{\rm ex}/\os{\circ}{T}}
(-\os{\circ}{\cal Q}{}^{\rm ex})  
@>{\sim}>> 
{\Om}^{\bul}_{{\cal Q}^{\rm ex}/\os{\circ}{T}}. 
\end{CD} 
\tag{1.3.12.4}\label{eqn:icor}
\end{equation*} 
Hence we obtain (\ref{eqn:pops}). 
\par 
(4):  
Take a local chart 
$(\{1\} \lo {\cal O}^*_T,P \lo M_{\ol{\cal Q}{}^{\star}},
\{1\}\os{\subset}{\lo} P)$ of the morphism 
$\ol{\cal Q}{}^{\star} \lo  \os{\circ}{T}$ 
on a neighborhood of $x$.  
Because $\ol{\cal Q}{}^{\star}$ 
is (formally) log smooth over $\os{\circ}{T}$, 
we can take the $P$ such that  
${\cal O}_{\ol{\cal Q}{}^{\star}}$ 
is \'{e}tale over ${\cal O}_{T}[P]$, 
in particular, flat over ${\cal O}_{T}[P]$. 
Since $P$ is integral, 
any element $a\in P$ defines an injective multiplication 
$$a\cdot \col  {\cal O}_{T}[P] \os{\sus}{\lo} {\cal O}_{T}[P].$$
Hence the morphism 
\begin{equation*} 
a\cdot \col {\cal O}_{\ol{\cal Q}{}^{\star}} 
\lo {\cal O}_{\ol{\cal Q}{}^{\star}}
\tag{1.3.12.5}\label{eqn:injm}
\end{equation*} 
is injective.  
We may assume that $m_i$ is the image of an element of $P$ 
and that 
$$\Om^1_{\os{\circ}{\ol{\cal Q}}{}^{\star}/\os{\circ}{T}}
=\bigoplus_{j=1}^k{\cal O}_{\ol{\cal Q}{}^{\star}}
d\al(m_j)
\oplus \bigoplus_{j=1}^l {\cal O}_{\ol{\cal Q}{}^{\star}}d x_j 
\quad (l\in {\mab N}),$$ 
$$\Om^1_{\ol{\cal Q}{}^{\star}/\os{\circ}{T}}
=\bigoplus_{j=1}^k{\cal O}_{\ol{\cal Q}{}^{\star}}
d\log m_j
\oplus \bigoplus_{j=1}^l {\cal O}_{\ol{\cal Q}{}^{\star}}d x_j 
\quad (l\in {\mab N})$$ 
with $x_j \in {\cal O}_{\os{\circ}{\ol{\cal Q}}{}^{\star}}$.  
Let $i$ be a nonnegative integer and consider a local section 
$\om \in \Om^i_{\os{\circ}{\ol{\cal Q}}{}^{\star}/\os{\circ}{T}}$. 
Set $\om_j:=d\log m_j$ $(0\leq j\leq k)$ and 
$\om_j:=dx_{j-k}$ $(k+1\leq j\leq k+l)$.  
Express the image of $\om$ in 
$\Om^i_{\ol{\cal Q}{}^{\star}/\os{\circ}{T}}$ 
by the following form 
$$\om =\sum_{j_1< \cdots <j_i}a_{j_1\cdots j_i}\om_{j_1}
\wedge \cdots \wedge \om_{j_i} \quad (a_{j_1\cdots j_i}\in 
{\cal O}_{\ol{\cal Q}{}^{\star}}).$$ 
For $j_1< \cdots <j_i$, let $n:=n(j_1,\ldots,j_i)$ be an integer 
such that $j_n \leq k$ and $j_{n+1} >k$. 
By the definition of $\om$, we have 
$a_{j_1\cdots j_i}=b_{j_1\cdots j_i}\al(m_{j_1}) \cdots \al(m_{j_n})$ 
for some $b_{j_1\cdots j_i}\in {\cal O}_{\ol{\cal Q}{}^{\star}}$. 
Hence  
$$\om =\sum_{j_1< \cdots <j_i}b_{j_1\cdots j_i}
d\al(m_{j_1}) \wedge \cdots \wedge d\al(m_{j_{n(j_1\cdots j_i)}}) 
\wedge \om_{j_{n(j_1\cdots j_i)+1}}
\wedge \cdots \wedge \om_{j_i}$$
in $\Om^i_{\os{\circ}{\ol{\cal Q}}{}^{\star}/\os{\circ}{T}}$. 
Assume that the image of 
$\om$ in $\Om^i_{\ol{\cal Q}{}^{\star}/\os{\circ}{T}}$ 
is zero. 
Then 
$0=a_{j_1\cdots j_i}=
b_{j_1\cdots j_i}\al(m_{j_1}) \cdots \al(m_{j_n})$ 
by the assumption about the locally freeness of  
$\Om^1_{\os{\circ}{\ol{\cal Q}}{}^{\star}/\os{\circ}{T}}$. 
By the injectivity of the morphism (\ref{eqn:injm}),  
we see that $b_{j_1\cdots j_i}=0$. Hence $\om=0$ 
and 
we have shown that the natural morphism 
$\Om^{\bul}_{\os{\circ}{\ol{\cal Q}}{}^{\star}/\os{\circ}{T}}
\lo \Om^{\bul}_{\ol{\cal Q}{}^{\star}/\os{\circ}{T}}$ 
is injective. 
\par 
(5): By (1) the morphism 
${\Om}^{\bul}_{\ol{\cal Q}{}^{\star}/\os{\circ}{T}}
(-\os{\circ}{\cal Q}{}^{\star}) \lo
\Om^{\bul}_{\ol{\cal Q}{}^{\star}/\os{\circ}{T}}$
is injective. 
By (2) 
the morphism 
${\Om}^{\bul}_{\ol{\cal Q}{}^{\star}/\os{\circ}{T}}
(-\os{\circ}{\cal Q}{}^{\star}_{T}) \lo
\Om^{\bul}_{\ol{\cal Q}{}^{\star}/\os{\circ}{T}}$ 
factors through the morphism 
${\Om}^{\bul}_{\ol{\cal Q}{}^{\star}/\os{\circ}{T}}
(-\os{\circ}{\cal Q}{}^{\star}) \lo
{\rm Im}
(\Om^{\bul}_{\os{\circ}{\ol{\cal Q}}{}^{\star}/\os{\circ}{T}}   
\lo 
\Om^{\bul}_{\ol{\cal Q}{}^{\star}/\os{\circ}{T}})$.  
The target of the last morphism is isomorphic to 
$\Om^{\bul}_{\os{\circ}{\ol{\cal Q}}{}^{\star}/\os{\circ}{T}}$ by (4). 
\end{proof}

\begin{rema}\label{rema:noninj}
Let the notations be as in (\ref{prop:incol}). 
\par 
(1) 
The following natural morphisms are   
not injective in general$:$
\begin{equation*} 
{\cal O}_{\ol{\mathfrak E}}
\otimes_{{\cal O}_{\ol{\cal Q}{}^{\rm ex}}}
{\Om}^{\bul}_{\ol{\cal Q}{}^{\rm ex}/\os{\circ}{T}}
(-\os{\circ}{\cal Q}{}^{\rm ex}) \lo
{\cal O}_{\ol{\mathfrak E}}
\otimes_{{\cal O}_{\ol{\cal Q}{}^{\rm ex}}}
\Om^{\bul}_{\os{\circ}{\ol{\cal Q}}{}^{\rm ex}/\os{\circ}{T}} 
\tag{1.3.13.1}\label{eqn:ypkd}
\end{equation*}
and 
\begin{equation*} 
{\cal O}_{\ol{\mathfrak E}}
\otimes_{{\cal O}_{\ol{\cal Q}{}^{\rm ex}}}
{\Om}^{\bul}_{\ol{\cal Q}{}^{\rm ex}/\os{\circ}{T}}
(-\os{\circ}{\cal Q}{}^{\rm ex}) \lo
{\cal O}_{\ol{\mathfrak E}}
\otimes_{{\cal O}_{\ol{\cal Q}{}^{\rm ex}}}
\Om^{\bul}_{\ol{\cal Q}{}^{\rm ex}/\os{\circ}{T}}.   
\tag{1.3.13.2}\label{eqn:yxnpd}
\end{equation*} 
Indeed, consider the case where $\bul =0$,  
$\ol{\cal Q}=\ol{S(T)^{\nat}}$ and 
${\cal O}_{\ol{S(T)^{\nat}}}={\cal O}_T[\tau]$ for the morphisms 
(\ref{eqn:ypkd}) and (\ref{eqn:yxnpd}). 
Then the morphisms (\ref{eqn:ypkd}) and (\ref{eqn:yxnpd}) 
are equal to 
\begin{equation*} 
{\cal O}_{T}\langle \tau\rangle 
\otimes_{{\cal O}_T[[\tau]]}\tau{\cal O}_{T}[[\tau]]
\lo
{\cal O}_{T}\langle \tau\rangle. 
\tag{1.3.13.3}\label{eqn:yxninpd}
\end{equation*}
Clearly this morphism is not injective. 
Indeed, let $n$ be a positive integer such that 
$p^n{\cal O}_T=0$. 
Consider a section 
$\tau^{[p^n-1]}\otimes \tau$ in 
${\cal O}_{T}\langle \tau\rangle \otimes_{{\cal O}_T[[\tau]]}\tau{\cal O}_{T}[[\tau]]$. 
Then, by the natural isomorphism 
$1\otimes (\tau \times) 
\simeq {\cal O}_{T}\langle \tau\rangle \otimes_{{\cal O}_T[[\tau]]}{\cal O}_{T}[[\tau]]
\os{\sim}{\lo}
{\cal O}_{T}\langle \tau\rangle \otimes_{{\cal O}_T[[\tau]]}\tau{\cal O}_{T}[[\tau]]$, 
this section comes from $\tau^{[p^n-1]}\otimes 1$, which is not zero. 
However the image of $\tau^{[p^n-1]}\otimes \tau$ in ${\cal O}_{T}\langle \tau\rangle$ 
is equal to $\tau \tau^{[p^n-1]}=p^n!/(p^n-1)!\tau^{[p^n]}=p^n\tau^{[p^n]}=0$. 
\par 
(2) The morphism  
\begin{equation*} 
{\cal O}_{\ol{\mathfrak E}}
\otimes_{{\cal O}_{\ol{\cal Q}{}^{\rm ex}}}
\Om^{\bul}_{\os{\circ}{\ol{\cal Q}}{}^{\rm ex}/\os{\circ}{T}} \lo
{\cal O}_{\ol{\mathfrak E}}
\otimes_{{\cal O}_{\ol{\cal Q}{}^{\rm ex}}}
\Om^{\bul}_{{\ol{\cal Q}}{}^{\rm ex}/\os{\circ}{T}} 
\tag{1.3.13.4}\label{eqn:ynkd}
\end{equation*}
is not necessarily injective.  
Consider the case where $\bul =1$,  
$\ol{\cal Q}=\ol{S(T)^{\nat}}$ and 
${\cal O}_{\ol{S(T)^{\nat}}}={\cal O}_T[\tau]$ in the notation (1). 
Then the morphism (\ref{eqn:ynkd}) is equal to 
\begin{equation*} 
{\cal O}_T\langle \tau\rangle 
\otimes_{{\cal O}_T[[\tau]]}{\cal O}_T[[\tau]]d\tau \lo
{\cal O}_T\langle \tau\rangle 
\otimes_{{\cal O}_T[[\tau]]}{\cal O}_T[[\tau]]d\log \tau. 
\end{equation*}
The image of the nonzero section $\tau^{[p^n-1]}\otimes d\tau$  
of ${\cal O}_T\langle \tau\rangle \otimes_{{\cal O}_T[[\tau]]}{\cal O}_T[[\tau]]d\tau$ 
by this morphism is $\tau \tau^{[p^n-1]}\otimes d\log \tau$, which is   
$0$ by the same reason as that of 
the non-injectivity in (1). 
\par 
(3) 
The natural morphisms  
\begin{equation*} 
{\cal O}_{\mathfrak E}
\otimes_{{\cal O}_{\ol{\cal Q}{}^{\rm ex}}}
{\Om}^{\bul}_{\ol{\cal Q}{}^{\rm ex}/\os{\circ}{T}}
(-\os{\circ}{\cal Q}{}^{\rm ex}) \lo
{\cal O}_{\mathfrak E}
\otimes_{{\cal O}_{\ol{\cal Q}{}^{\rm ex}}} 
\Om^{\bul}_{\os{\circ}{\ol{\cal Q}}{}^{\rm ex}/\os{\circ}{T}}
\tag{1.3.13.5}\label{eqn:ydpnpd}
\end{equation*} 
and 
\begin{equation*} 
{\cal O}_{\mathfrak E}
\otimes_{{\cal O}_{\ol{\cal Q}{}^{\rm ex}}}
\Om^{\bul}_{\os{\circ}{\ol{\cal Q}}{}^{\rm ex}/\os{\circ}{T}} 
\lo
{\cal O}_{\mathfrak E}
\otimes_{{\cal O}_{\ol{\cal Q}{}^{\rm ex}}}
\Om^{\bul}_{{\ol{\cal Q}{}^{\rm ex}}/\os{\circ}{T}}   
\tag{1.3.13.6}\label{eqn:ydplnpd}
\end{equation*} 
are not necessarily injective. 
Indeed, consider the case $\bul =0$ and $\ol{\cal Q}=\ol{S(T)^{\nat}}$ 
for the morphism (\ref{eqn:ydpnpd}). 
Then the morphism (\ref{eqn:ydpnpd}) is the zero morphism. 
However the source of the morphism (\ref{eqn:ydpnpd}) is 
${\cal O}_T\otimes_{{\cal O}_{\ol{S(T)^{\nat}}}}\tau{\cal O}_{\ol{S(T)^{\nat}}}
\simeq {\cal O}_T$, 
which is not zero. 
Next consider the case 
$\ol{\cal Q}=\ol{S(T)^{\nat}}$ for the morphism (\ref{eqn:ydplnpd}). 
Then the morphism (\ref{eqn:ydplnpd}) for the case $\bul =1$ is 
the following zero morphism
$$0\col {\cal O}_Td\tau \lo {\cal O}_Td\log \tau.$$
Hence the morphism (\ref{eqn:ydplnpd}) is not injective. 
\end{rema}

\par 
Let $X$ be an SNCL scheme over $S$. 
Set 
\begin{align*} 
X_{\os{\circ}{T}_0}:=X\times_{S}S_{\os{\circ}{T}_0}
=X\times_{\os{\circ}{S}}\os{\circ}{T}_0,\quad   
X_{T_0}:=X\times_{S}T_0,\quad  
\os{\circ}{X}{}^{(k)}_{T_0}
:=\os{\circ}{X}{}^{(k)}\times_{\os{\circ}{S}}\os{\circ}{T}_0\quad 
(k\in {\mab Z}_{\geq 0}).
\end{align*}  
The 
orientation sheaf 
$\vp^{(k)}_{\rm zar}(\os{\circ}{X}_{T_0}/\os{\circ}{T}_0)$ extends to 
an abelian sheaf 
$\vp^{(k)}_{\rm crys}(\os{\circ}{X}_{T_0}/\os{\circ}{T})$ 
on $(\os{\circ}{X}{}^{(k)}_{T_0}/\os{\circ}{T})_{\rm crys}$ 
as in \cite{nh2}. 
We call $\vp^{(k)}_{\rm crys}(\os{\circ}{X}_{T_0}/\os{\circ}{T})$ 
the $(k+1)$-{\it fold} {\it crystalline orientation sheaf} of 
$\os{\circ}{X}_{T_0}/(\os{\circ}{T},{\cal J},\del)$. 
Let $f\col X_{\os{\circ}{T}_0}\lo S_{\os{\circ}{T}_0}$ and 
$\os{\circ}{f}{}^{(k)} \col \os{\circ}{X}{}^{(k)}_{T_0}\lo \os{\circ}{T}_0$ 
be the structural morphisms. 
By abuse of notation, denote by $f$ and $\os{\circ}{f}{}^{(k)}$ 
the structural morphisms 
$X_{\os{\circ}{T}_0}\lo S_{\os{\circ}{T}_0}\lo S(T)^{\nat}$ and 
$\os{\circ}{f}{}^{(k)} \col \os{\circ}{X}{}^{(k)}_{T_0}\lo \os{\circ}{T}$ 
$(k\in {\mab N})$, respectively. 
Let $f_T\col X_{T_0} \lo T_0\os{\sus}{\lo}T$ 
be also the structural morphism. 
Then $\os{\circ}{f}_T=\os{\circ}{f}$. 
\par 
In this section we assume that 
there exists an immersion 
$X_{\os{\circ}{T}_0} \os{\sus}{\lo} {\cal P}$ into 
a log smooth scheme over $S(T)^{\nat}$. 
The log formal schemes 
${\cal P}^{\rm prex}{}$ between (\ref{prop:zpex}) and (\ref{rema:pds}) 
and ${\cal P}^{\rm ex}$ satisfy the condition  
{\rm (\ref{eqn:epynr})} by (\ref{prop:fsi}) and (\ref{prop:nexeo}). 
\par 
Let ${\mathfrak D}$ be the log PD-envelope of 
the immersion $X_{\os{\circ}{T}_0}\os{\sus}{\lo} {\cal P}$ over 
$(S(T)^{\nat},{\cal J},{\del})$. 
By abuse of notation,  
for $1\leq i\leq r$ and 
for any different $i_0,\ldots, i_k$ $(1\leq i_0,\ldots,i_k \leq r)$, 
denote $D(M_{{\cal P}^{\rm ex}})_i$ by 
${\cal P}^{\rm ex}_i$ and 
denote by $\os{\circ}{\cal P}{}^{\rm ex}_{i_0 \cdots i_k}$ 
the local closed subscheme 
$\os{\circ}{\cal P}{}^{\rm ex}_{i_0}
\cap \cdots \cap \os{\circ}{\cal P}{}^{\rm ex}_{i_k}$ 
of $\os{\circ}{\cal P}{}^{\rm ex}$. 
It is easy to check that 
${\cal P}^{\rm ex}_i$ is equal to 
${\cal P}^{\rm ex}_i$ defined in 
\S\ref{sec:snclv}. 
We denote $\os{\circ}{D}{}^{(k)}(M_{{\cal P}^{{\rm ex}}})$  
and 
$\vp^{(k)}_{\rm zar}(\os{\circ}{D}(M_{{\cal P}^{\rm ex}}))$ 
by 
$\os{\circ}{\cal P}{}^{{\rm ex},(k)}$ and 
$\vp^{(k)}_{\rm zar}(\os{\circ}{\cal P}{}^{\rm ex}/\os{\circ}{T})$, respectively. 
We have a natural immersion 
$\os{\circ}{X}{}^{(k)}_{T_0} \os{\sus}{\lo} 
\os{\circ}{\cal P}{}^{{\rm ex},(k)}$ 
of (formal) schemes over $\os{\circ}{T}$. 
Because 
$\os{\circ}{X}{}^{(k)}_{T_0,{\rm zar}}= \os{\circ}{\cal P}{}^{{\rm ex},(k)}_{\rm zar}$ 
as topological spaces, 
we can identify $\vp^{(k)}_{\rm zar}(\os{\circ}{\cal P}{}^{\rm ex}/\os{\circ}{T})$ 
with $\vp^{(k)}_{\rm zar}(\os{\circ}{X}_{T_0}/\os{\circ}{T}_0)$. 
Denote the natural local exact closed immersion 
$\os{\circ}{\cal P}{}^{\rm ex}_{i_0 \cdots i_k}\os{\sus}{\lo} 
\os{\circ}{\cal P}{}^{\rm ex}$ by $b_{i_0\cdots i_k}$. 
Let 
$b^{(k)} \col \os{\circ}{\cal P}{}^{{\rm ex},(k)}\lo 
\os{\circ}{\cal P}{}^{\rm ex}$ $(k\in {\mab N})$ 
be the natural morphism. 
Set 
${\cal P}_T:={\cal P}\times_{S(T)}T$,  
${\cal P}{}^{{\rm ex}}_T:=
{\cal P}{}^{{\rm ex}}\times_{S(T)}T$ 
and ${\mathfrak D}_T={\mathfrak D}\times_{S(T)}T$ 
when $T$ is restrictively hollow 
with respective to the morphism $T_0\lo S$.  
\par 
The following (2) is a generalization of the Poincar\'{e} residue isomorphism, 
which is an SNCL version of \cite[(4.3)]{nh3}: 

\begin{prop}\label{prop:perkl}  
$(1)$ The scheme $\os{\circ}{\cal P}{}^{{\rm ex},(k)}$ 
$(k\in {\mab N})$ is formally smooth over $\os{\circ}{T}$.  
\par 
$(2)$ Identify the points of $\os{\circ}{\mathfrak D}$ 
with those of $\os{\circ}{X}_{T_0}$. 
Identify also the images of the points of $\os{\circ}{X}_{T_0}$ 
in $\os{\circ}{\cal P}{}^{\rm ex}$ with the points of $\os{\circ}{X}_{T_0}$. 
Let $x$ be a point of $\os{\circ}{\mathfrak D}$. 
$($Let $r$ be a nonnegative integer such that 
$M_{X_{\os{\circ}{T},x}}/{\cal O}_{X_{\os{\circ}{T},x}}^*\simeq {\mab N}^r$.$)$ 
Then, for a positive integer $k$, 
the following morphism  
\begin{equation*} 
{\rm Res} \col 
{\cal O}_{\mathfrak D}\otimes_{{\cal O}_{{\cal P}^{\rm ex}}}
P_k{\Om}^{\bul}_{{\cal P}^{\rm ex}/\os{\circ}{T}} \lo 
{\cal O}_{\mathfrak D}\otimes_{{\cal O}_{{\cal P}^{\rm ex}}}
b^{(k-1)}_*(\Om^{\bul}_{\os{\circ}{\cal P}{}^{{\rm ex},(k-1)}
/\os{\circ}{T}}
\otimes_{\mab Z}\vp^{(k-1)}_{\rm zar}
(\os{\circ}{\cal P}{}^{\rm ex}/\os{\circ}{T}))[-k]
\tag{1.3.14.1}\label{eqn:mprrn}
\end{equation*} 
\begin{equation*} 
\sig \otimes d\log m_{i_0}\cdots d\log m_{i_{k-1}} \omega 
\lom 
\sig \otimes b^*_{i_0\cdots i_{k-1}}(\omega)
\otimes({\rm orientation}~(i_0\cdots i_{k-1})) 
\end{equation*} 
$$(\sig \in {\cal O}_{\mathfrak D}, 
\om \in P_0{\Om}^{\bul}_{{\cal P}^{\rm ex}/\os{\circ}{T}})$$
{\rm (cf.~\cite[(1.6)]{gn} (\cite[(3.1.5)]{dh2}))} 
induces the following ``Poincar\'{e} residue isomorphism''  
\begin{align*} 
{\rm gr}^P_k
({\cal O}_{\mathfrak D}\otimes_{{\cal O}_{{\cal P}^{\rm ex}}}
{\Om}^{\bul}_{{\cal P}^{\rm ex}/\os{\circ}{T}}) 
& \os{\sim}{\lo} 
{\cal O}_{\mathfrak D}\otimes_{{\cal O}_{{\cal P}^{\rm ex}}}
b^{(k-1)}_*(\Om^{\bul}_{\os{\circ}{\cal P}{}^{{\rm ex},(k-1)}
/\os{\circ}{T}}
\otimes_{\mab Z}\vp^{(k-1)}_{\rm zar}
(\os{\circ}{\cal P}{}^{\rm ex}/\os{\circ}{T}))[-k]
\tag{1.3.14.2}\label{eqn:prvin} \\ 
& (\os{\sim}{\lo} 
{\cal O}_{\mathfrak D}\otimes_{{\cal O}_{{\cal P}^{\rm ex}}}
b^{(k-1)}_*(\Om^{\bul}_{\os{\circ}{\cal P}{}^{{\rm ex},(k-1)}
/\os{\circ}{T}}
\otimes_{\mab Z}\vp^{(k-1)}_{\rm zar}
(\os{\circ}{X}/\os{\circ}{T}))[-k]).  
\end{align*} 
\end{prop} 
\begin{proof} 
(1): (1) immediately follows from (\ref{prop:fsi}). 
\par 
(2): As in the usual case, we can easily check that 
the morphism (\ref{eqn:mprrn}) is well-defined and surjective. 
The question is local. Replacing ${\cal P}$ by 
${\cal P}^{\rm prex}$, 
we may assume that the immersion 
$X_{\os{\circ}{T}_0}\os{\sus}{\lo} {\cal P}$ is exact.   
We may assume that $M_{S(T)^{\nat}}$ is free of rank $1$.   
Because the question is local on $X_{\os{\circ}{T}_0}$, 
we may assume that 
there exists the following cartesian diagrams   
\begin{equation*}
\begin{CD} 
{\cal P} @>{\sus}>> \ol{\cal P} \\
@VVV @VVV \\
{\mab A}_{S(T)^{\nat}}(a,d') 
@>{\sus}>> {\mab A}_{\ol{S(T)^{\nat}}}(a,d') \\ 
@VVV @VVV \\
S(T)^{\nat}@>{\sus}>> \ol{S(T)^{\nat}},  
\end{CD} 
\tag{1.3.14.3}\label{cd:pwtp} 
\end{equation*}
where the vertical morphism 
$\ol{\cal P} \lo {\mab A}_{\ol{S(T)^{\nat}}}(a,d')$ is solid and \'{e}tale 
((\ref{prop:fsi}),  (\ref{lemm:etl})).  
Let 
$\ol{b}{}^{(k)}\col \os{\circ}{\cal P}{}^{(k)}\lo 
\os{\circ}{\ol{\cal P}}$ 
be the natural morphism. 
Then, by 
the usual Poincar\'{e} residue isomorphism 
\cite[(2.2.21.3)]{nh2} (with logarithmic forms put on the left as in 
(\ref{eqn:mprrn})), 
we have the following isomorphism 
\begin{equation*}  
{\rm Res} \col 
{\rm gr}_k^P{\Om}^{\bul}_{\ol{\cal P}/\os{\circ}{T}}
\os{\sim}{\lo} 
\ol{b}{}^{(k-1)}_*
(\Om^{\bul}_{\os{\circ}{\cal P}{}^{(k-1)}/\os{\circ}{T}}
\otimes_{\mab Z}\vp^{(k-1)}_{\rm zar}(\os{\circ}{\cal P}/\os{\circ}{T}))[-k]
\quad (k\geq 1). 
\tag{1.3.14.4}\label{eqn:gurps}
\end{equation*}  
Because 
${\rm gr}_k^P{\Om}^{\bul}_{\ol{\cal P}/\os{\circ}{T}}
=(P_k{\Om}^{\bul}_{\ol{\cal P}/\os{\circ}{T}}/
{\Om}^{\bul}_{\ol{\cal P}/\os{\circ}{T}}(-\os{\circ}{\cal P}))
/(P_{k-1}{\Om}^{\bul}_{\ol{\cal P}/\os{\circ}{T}}/
{\Om}^{\bul}_{\ol{\cal P}/\os{\circ}{T}}(-\os{\circ}{\cal P}))
={\rm gr}_k^P{\Om}^{\bul}_{{\cal P}/\os{\circ}{T}}$ 
by (\ref{eqn:pops}), 
we have the following isomorphism 
\begin{equation*}  
{\rm Res} \col {\rm gr}_k^P{\Om}^{\bul}_{{\cal P}/\os{\circ}{T}}
\os{\sim}{\lo} 
b^{(k-1)}_*
(\Om^{\bul}_{\os{\circ}{\cal P}{}^{(k-1)}/\os{\circ}{T}}
\otimes_{\mab Z}\vp^{(k-1)}_{\rm zar}(\os{\circ}{\cal P}/\os{\circ}{T}))[-k] 
\quad (k\geq 1).  
\tag{1.3.14.5}\label{eqn:grps}
\end{equation*} 
By (\ref{prop:om}) and (\ref{eqn:grps}), 
we have the following exact sequence 
\begin{align*} 
0 \lo P_{k-1}
({\cal O}_{\mathfrak D}\otimes_{{\cal O}_{\cal P}}
{\Om}^{\bul}_{{\cal P}/\os{\circ}{T}}) 
& \lo 
P_k
({\cal O}_{\mathfrak D}\otimes_{{\cal O}_{\cal P}}
{\Om}^{\bul}_{{\cal P}/\os{\circ}{T}})  
\tag{1.3.14.6}\label{eqn:ypbgd} \\
{} & \lo 
{\cal O}_{\mathfrak D}\otimes_{{\cal O}_{\cal P}}
b^{(k-1)}_*(\Om^{\bul}_{\os{\circ}{\cal P}{}^{(k-1)}/\os{\circ}{T}}
\otimes_{\mab Z}\vp^{(k-1)}_{\rm zar}(\os{\circ}{\cal P}/\os{\circ}{T}))[-k] \lo 0. 
\end{align*}
Let ${\cal J}$ be the defining ideal sheaf of the immersion 
$X\os{\sus}{\lo}{\cal P}$. 
By \cite[3.20 Remarks 7)]{bob}, 
${\cal O}_{\mathfrak D}\otimes_{{\cal O}_{\cal P}}
={\cal O}_{\mathfrak D}\otimes_{{\cal O}_{\cal P}/{\cal J}^n}$ 
for some $n\in {\mab N}$. 
Hence the fourth term of (\ref{eqn:ypbgd}) 
is equal to the right hand side of (\ref{eqn:prvin}). 
\end{proof}

\begin{coro}\label{coro:flt}
The complexes 
${\cal O}_{\mathfrak D}\otimes_{{\cal O}_{{\cal P}^{\rm ex}}}
{\Om}^{\bul}_{{\cal P}^{\rm ex}/\os{\circ}{T}}$,  
${\rm gr}^P_k
({\cal O}_{\mathfrak D}\otimes_{{\cal O}_{{\cal P}^{\rm ex}}}
{\Om}^{\bul}_{{\cal P}^{\rm ex}/\os{\circ}{T}})$ 
$(k\in {\mab N})$ 
and 
$P_k
({\cal O}_{\mathfrak D}\otimes_{{\cal O}_{{\cal P}^{\rm ex}}}
{\Om}^{\bul}_{{\cal P}^{\rm ex}/\os{\circ}{T}})$ 
consist of locally free ${\cal O}_T$-modules.    
\end{coro}
\begin{proof}  
Since the question is local on $X$, 
we may assume
that $X={\mab A}_{S}(a,d)$ and 
${\cal P}^{\rm prex}={\cal P}
= {\mab A}_{S(T)^{\nat}}(a,d')$ 
$(a\in {\mab N}, d\leq d')$ by (\ref{prop:adla}).
Then it is clear that 
${\cal O}_{\mathfrak D}\otimes_{{\cal O}_{{\cal P}^{\rm ex}}}
{\Om}^{\bul}_{{\cal P}^{\rm ex}/\os{\circ}{T}}$ 
consists of locally free ${\cal O}_T$-modules. 
The fourth term of (\ref{eqn:ypbgd}) 
consists of locally free ${\cal O}_T$-modules.
Hence the descending induction on $k$ shows 
the locally freeness of 
$P_k({\cal O}_{\mathfrak D}\otimes_{{\cal O}_{{\cal P}^{\rm ex}}}
{\Om}^{\bul}_{{\cal P}^{\rm ex}/\os{\circ}{T}})$.  
\end{proof} 

\par 
Let $\os{\circ}{\mathfrak D}{}^{(k)}$ 
be the PD-envelope of the immersion 
$\os{\circ}{X}{}^{(k)}_{T_0} \os{\sus}{\lo} \os{\circ}{\cal P}{}^{{\rm ex},(k)}$ 
over 
$(\os{\circ}{T},{\cal J},{\del})$.  
Let $c^{(k)}\col \os{\circ}{\mathfrak D}{}^{(k)}\lo 
\os{\circ}{\mathfrak D}$ be 
the natural morphism of schemes over $\os{\circ}{T}$. 

\begin{lemm}\label{lemm:pde}     
{{\rm {\bf (cf.~\cite[(2.2.16)]{nh2}, 
\cite[(4.5)]{nh3})}}} 
$(1)$ 
$\os{\circ}{\mathfrak D}{}^{(k)}=
\os{\circ}{\mathfrak D}
\times_{\os{\circ}{\cal P}{}^{\rm ex}}
\os{\circ}{\cal P}{}^{{\rm ex},(k)}$. 
\par 
$(2)$ Assume that $T$ is restrictively hollow 
with respective to the morphism $T_0\lo S$. 
Then the log scheme ${\mathfrak D}_T$ 
is the log PD-envelope of the immersion 
$X_{T_0}\os{\sus}{\lo} {\cal P}_T$ over $(T,{\cal J},\del)$. 
\par 
$(2)'$ Let the assumption be as in $(2)$. 
Then $({\mathfrak D}_T)^{\circ}=\os{\circ}{\mathfrak D}$. 
\end{lemm}
\begin{proof}
(1): By the universality of 
the (log) PD-envelope, the question is local. 
We may assume that the immersion 
$X_{\os{\circ}{T}_0}\os{\sus}{\lo} {\cal P}$ is exact. 
Let $x$ be a point of $\os{\circ}{X}_{T_0}$. 
We may assume that 
there exists the following cartesian diagram 
\begin{equation*} 
\begin{CD} 
X_{\os{\circ}{T}_0} @>{\subset}>> {\cal P} \\ 
@VVV @VVV \\ 
{\mab A}_{S_{\os{\circ}{T}_0}}(a,d) @>{\subset}>> 
{\mab A}_{S(T)^{\nat}}(a,d') 
\end{CD}
\tag{1.3.16.1}\label{cd:xpas} 
\end{equation*} 
locally around $x$, where 
the lower horizontal morphism is 
an exact closed immersion defined by the ideal sheaf 
$(x_{d+1},\ldots,x_{d'})$. 
Then we have the equality 
\begin{equation*} 
X_{\os{\circ}{T}_0}=  {\mab A}_{S_{\os{\circ}{T}_0}}(a,d)
\times_{{\mab A}_{S(T)^{\nat}}(a,d')}{\cal P}.
\tag{1.3.16.2}\label{eqn:atad}
\end{equation*}   
The rest of the proof is the same as that of \cite[(2.2.16) (2)]{nh2}. 
\par 
(2), $(2)'$: 
Since $T$ is restrictively hollow with respective to the morphism $T_0\lo S$, 
$S(T)=S(T)^{\nat}$. 
Let ${\mathfrak D}'(T)$ be the log PD-envelope of 
the immersion $X_{T_0}\os{\sus}{\lo} {\cal P}_T$ over $(T,{\cal J},\del)$. 
Then, by \cite[3.21 Proposition]{bob}, 
we have a natural morphism 
${\mathfrak D}'(T) \lo {\mathfrak D}_T$. 
Hence the question is local and 
we may assume the existence of (\ref{cd:xpas}).  
In this case, the natural morphism is an isomorphism 
by the equality (\ref{eqn:atad}). We have also proved $(2)'$.   
\end{proof}

\par 
Now we define various morphisms of log schemes and ringed topoi. 
\par 
Let $p_{T_0}\col X_{T_0}\lo X_{\os{\circ}{T}_0}$ 
(resp.~$p_{\os{\circ}{T}_0}\col X_{T_0}\lo X_{\os{\circ}{T}_0}$) 
be the first projection over 
$T_0\lo S_{\os{\circ}{T}_0}$ (resp.~$\os{\circ}{T}_0$).   
Recall that $\os{\circ}{X}_{T_0}=
\os{\circ}{X}\times_{\os{\circ}{S}}\os{\circ}{T}_0$ ((\ref{lemm:usbc})). 
Hence $\os{\circ}{p}_{T_0}=\os{\circ}{p}_{\os{\circ}{T}_0}
={\rm id}_{\os{\circ}{X}_{T_0}}$. 
Let $p^{\rm PD}_T\col {\mathfrak D}_T\lo{\mathfrak D}$ 
be the first projection over $T\lo S(T)$ when $T$ is restrictively hollow 
with respective to the morphism $T_0\lo S$. 
Then $\os{\circ}{p}{}^{\rm PD}_T={\rm id}_{\os{\circ}{\mathfrak D}}$ 
by (\ref{lemm:pde}) $(2)'$. 
Let $$p_{T{\rm crys}}\col 
((X_{T_0}/T)_{\rm crys},{\cal O}_{X_{T_0}/T})
\lo 
((X_{\os{\circ}{T}_0}/S(T))_{\rm crys},
{\cal O}_{X_{\os{\circ}{T}_0}/S(T)})$$  
and 
$$p_{\os{\circ}{T}{\rm crys}}
\col ((X_{T_0}/\os{\circ}{T})_{\rm crys},
{\cal O}_{X_{T_0}/\os{\circ}{T}})
\lo 
((X_{\os{\circ}{T}_0}/\os{\circ}{T})_{\rm crys},
{\cal O}_{X_{\os{\circ}{T}_0}/\os{\circ}{T}})$$  
be the induced morphism by  $p_{T_0}$ 
and $p_{\os{\circ}{T}_0}$, respectively. 
Set $U_0:=S_{\os{\circ}{T}_0}$ or $T_0$ and $U:=S(T)^{\nat}$ or $T$, respectively. 
Let 
\begin{equation*} 
u_{X_{U_0}/U} \col 
((X_{U_0}/U)_{\rm crys},{\cal O}_{X_{U_0}/U}) 
\lo ((X_{U_0})_{\rm zar},f^{-1}({\cal O}_T))
=((X_{\os{\circ}{T}_0})_{\rm zar},f^{-1}({\cal O}_T)) 
\tag{1.3.16.3}\label{eqn:prjdef}
\end{equation*}   
and 
\begin{align*} 
\os{\circ}{u}_{X_{\os{\circ}{T}_0}/S(T)^{\nat}}=\os{\circ}{u}_{X_{T_0}/T} 
\col ((\os{\circ}{X}_{T_0}/\os{\circ}{T})_{\rm crys},
{\cal O}_{\os{\circ}{X}_{T_0}/\os{\circ}{T}}) 
&\lo ((X_{T_0})_{\rm zar},f^{-1}_T({\cal O}_T))\tag{1.3.16.4}\label{ali:prcuxef}\\
&=((X_{\os{\circ}{T}_0})_{\rm zar},f^{-1}({\cal O}_T))
\end{align*}   
be the natural projections. 
Let 
\begin{align*} 
\eps_{X_{U_0}/U} \col 
((X_{U_0}/U)_{\rm crys},{\cal O}_{X_{U_0}/U}) 
&\lo 
((\os{\circ}{X}_{S_{\os{\circ}{T}_0}}/\os{\circ}{S(T)^{\nat}})_{\rm crys},
{\cal O}_{\os{\circ}{X}_{S_{\os{\circ}{T}_0}}/\os{\circ}{S(T)^{\nat}}})\tag{1.3.16.5}\label{eqn:preoodef}\\
&=
((\os{\circ}{X}_{T_0}/\os{\circ}{T})_{\rm crys},
{\cal O}_{\os{\circ}{X}_{T_0}/\os{\circ}{T}})
\end{align*}   
be the morphism forgetting the log structures
of $X_{U_0}$ and $U$. 
Let 
\begin{equation*} 
\eps_{X_{U_0}/\os{\circ}{T}} \col 
((X_{U_0}/\os{\circ}{T})_{\rm crys},{\cal O}_{X_{U_0}/\os{\circ}{T}}) 
\lo 
((\os{\circ}{X}_{U_0}/\os{\circ}{T})_{\rm crys},
{\cal O}_{\os{\circ}{X}_{U_0}/\os{\circ}{T}})=
((\os{\circ}{X}_{T_0}/\os{\circ}{T})_{\rm crys},
{\cal O}_{\os{\circ}{X}_{T_0}/\os{\circ}{T}})
\tag{1.3.16.6}\label{eqn:prtdef}
\end{equation*}   
be the morphism forgetting the log structure of $X_{U_0}$. 
Let 
\begin{equation*}  
\eps_{X_{U_0}/U/\os{\circ}{T}}\col 
((X_{U_0}/U)_{\rm crys},{\cal O}_{X_{U_0}/U})\lo 
((X_{U_0}/\os{\circ}{T})_{\rm crys},{\cal O}_{X_{U_0}/\os{\circ}{T}})
\tag{1.3.16.7}\label{eqn:tntpi} 
\end{equation*} 
be the morphism forgetting the log structure of $U$. 
Then 
\begin{equation*}  
\eps_{X_{U_0}/U}=\eps_{X_{U_0}/\os{\circ}{T}}\circ \eps_{X_{U_0}/U/\os{\circ}{T}}. 
\tag{1.3.16.8}\label{eqn:tnepi} 
\end{equation*} 
\par 
Let $F$ be a quasi-coherent log crystal of 
${\cal O}_{X_{\os{\circ}{T}_0}/S(T)^{\nat}}$-modules and 
let $({\cal F},\nabla)$ be 
the quasi-coherent ${\cal O}_{\mathfrak D}$-module 
with the integrable connection 
associated to $F$ with respect to the immersion 
$X_{\os{\circ}{T}_0}\os{\sus}{\lo} {\cal P}$ over $S(T)^{\nat}$: 
$$\nabla \col 
{\cal F}\lo {\cal F}\otimes_{{\cal O}_{\cal P}}
{\Om}^1_{{\cal P}/S(T)^{\nat}}.$$ 
Then, by \cite[(6.4)]{klog1} 
(or the log Poincar\'{e} lemma (\cite[(2.2.8)]{nh2})), 
we have the following formula 
\begin{align*}  
Ru_{X_{\os{\circ}{T}_0}/S(T)^{\nat}*}(F)
&= {\cal F}\otimes_{{\cal O}_{\cal P}}
{\Om}^{\bul}_{{\cal P}/S(T)^{\nat}}
={\cal F}\otimes_{{\cal O}_{{\cal P}^{\rm ex}}}
{\Om}^{\bul}_{{\cal P}^{\rm ex}/S(T)^{\nat}}
\tag{1.3.16.9}\label{eqn:eopzdq}\\
\end{align*}
in ${\rm D}^+(f^{-1}({\cal O}_T))$. 
If $T$ is restrictively hollow with respective to the morphism $T_0\lo S$, then 
$S(T)^{\nat}=S(T)$ and 
\begin{align*}  
Ru_{X_{\os{\circ}{T}_0}/S(T)^{\nat}*}(F)
&= {\cal F}\otimes_{{\cal O}_{\cal P}}
{\Om}^{\bul}_{{\cal P}/S(T)^{\nat}}
={\cal F}\otimes_{{\cal O}_{{\cal P}_T}}
{\cal O}_{{\cal P}_T}\otimes_{{\cal O}_{{\cal P}}}
{\Om}^{\bul}_{{\cal P}/S(T)^{\nat}}
\tag{1.3.16.10}\label{eqn:eoprzdq}\\
& =
{\cal F}\otimes_{{\cal O}_{{\cal P}_T}}
{\Om}^{\bul}_{{\cal P}_T/T} =
{\cal F}\otimes_{{\cal O}_{{\cal P}^{\rm ex}_T}}
{\Om}^{\bul}_{{\cal P}^{\rm ex}_T/T}
\end{align*}
in ${\rm D}^+(f^{-1}({\cal O}_T))$. 
We have used (\ref{lemm:pde}) (2), (2)' and 
the formula in the end of \cite[(1.7)]{klog1} 
for the second equality and the third equality, respectively.
\par  
Assume that the following cartesian diagram 
\begin{equation*} 
\begin{CD}
{\cal P} @>{\subset}>> \ol{\cal P}\\
@VVV @VVV \\
S(T)^{\nat} @>{\subset}>> \ol{S(T)^{\nat}}
\end{CD}
\tag{1.3.16.11}\label{cd:wgtx}
\end{equation*} 
exists, where the right vertical morphism is log smooth.  
Let $\ol{\mathfrak D}$ be the log PD-envelope of 
the immersion $X_{\os{\circ}{T}_0}\os{\sus}{\lo} \ol{\cal P}$ 
over $(\os{\circ}{T},{\cal J},\del)$; 
$\os{\circ}{\ol{\mathfrak D}}$ is also the PD-envelope 
the immersion 
$\os{\circ}{X}_{\os{\circ}{T}_0}\os{\sus}{\lo}\os{\circ}{\ol{\cal P}}{}^{\rm ex}$ 
over $(\os{\circ}{T},{\cal J},\del)$.   
Let $\ol{F}$ be a quasi-coherent log crystal of 
${\cal O}_{X_{\os{\circ}{T}_0}/\os{\circ}{T}}$-modules. 
Let $(\ol{\cal F},\ol{\nabla})$ be 
the quasi-coherent ${\cal O}_{\ol{\mathfrak D}}$-module 
with the integrable connection associated to $\ol{F}$: 
$$\ol{\nabla} \col \ol{\cal F}\lo 
\ol{\cal F}\otimes_{{\cal O}_{\ol{\cal P}}}
{\Om}^1_{\ol{\cal P}/\os{\circ}{T}}
=
\ol{\cal F}
\otimes_{{\cal O}_{\ol{\cal P}{}^{\rm ex}}}
{\Om}^1_{\ol{\cal P}{}^{\rm ex}/\os{\circ}{T}}.$$ 
Set ${\cal F}:={\cal O}_{\mathfrak D}
\otimes_{{\cal O}_{\ol{\mathfrak D}}}\ol{\cal F}$ 
and let 
\begin{equation*} 
{\nabla} \col {\cal F} \lo {\cal F}\otimes_{{\cal O}_{\cal P}}
{\Om}^1_{{\cal P}/\os{\circ}{T}}=
{\cal F}\otimes_{{\cal O}_{{\cal P}^{\rm ex}}}
{\Om}^1_{{\cal P}^{\rm ex}/\os{\circ}{T}} 
\tag{1.3.16.12}\label{eqn:nidospltd}
\end{equation*}  
be the induced connection by $\ol{\nabla}$. 
Indeed, we obtain the connection $\nabla$ because 
$df^{[i]}=f^{[i-1]}df$ $(i\in {\mab Z}_{\geq 1})$ for a local section 
$f$ of the PD-ideal sheaf of ${\cal O}_{{\mathfrak D}(\ol{S(T)^{\nat}})}$; 
when $i\geq 2$, $f^{[i-1]}df=0$ in ${\Om}^1_{{\cal P}^{\rm ex}/\os{\circ}{T}}$; 
when $i=1$, $d\tau=\tau d\log \tau=0$ in 
${\Om}^1_{{\cal P}^{\rm ex}/\os{\circ}{T}}$ 
for a local section $\tau$ 
of ${\cal O}_{{\mathfrak D}(\ol{S(T)^{\nat}})}$ such that 
${\cal O}_{{\mathfrak D}(\ol{S(T)^{\nat}})}\simeq 
{\cal O}_T\langle \tau \rangle$. 
These facts tell us the existence of $\nabla$. 
\par 
Let $E$ be a quasi-coherent crystal of 
${\cal O}_{\os{\circ}{X}_{T_0}/\os{\circ}{T}}$-modules. 
Now we assume that 
$\ol{F}=\eps^*_{X_{\os{\circ}{T}_0}/\os{\circ}{T}}(E)$. 
Let $(\ol{\cal E},\ol{\nabla})$ be 
the quasi-coherent ${\cal O}_{\ol{\mathfrak D}}$-module 
with the integrable connection 
associated to $\eps^*_{X_{\os{\circ}{T}_0}/\os{\circ}{T}}(E)$: 
$$\ol{\nabla}\col \ol{\cal E}\lo \ol{\cal E}
\otimes_{{\cal O}_{\ol{\cal P}}}{\Om}^1_{\ol{\cal P}/\os{\circ}{T}}=
\ol{\cal E}
\otimes_{{\cal O}_{\ol{\cal P}{}^{\rm ex}}}
{\Om}^1_{\ol{\cal P}{}^{\rm ex}/\os{\circ}{T}}.$$ 
Set ${\cal E}:={\cal O}_{\mathfrak D}
\otimes_{{\cal O}_{\ol{\mathfrak D}}}\ol{\cal E}$ 
and let 
\begin{equation*} 
\nabla \col {\cal E} \lo {\cal E}\otimes_{{\cal O}_{\cal P}}
{\Om}^1_{{\cal P}/\os{\circ}{T}} 
={\cal E}\otimes_{{\cal O}_{{\cal P}^{\rm ex}}}{\Om}^1_{{\cal P}^{\rm ex}/\os{\circ}{T}} 
\tag{1.3.16.13}\label{eqn:nidpltd}
\end{equation*}  
be the induced connection by $\ol{\nabla}$. 
\par 
Let $k$ be a nonnegative integer. 
Let 
$a^{(k)}\col 
\os{\circ}{X}{}^{(k)}_{T_0}\lo \os{\circ}{X}_{T_0}$ 
be the morphism defined after (\ref{eqn:bkezps}) 
for the SNCL scheme $X_{\os{\circ}{T}_0}/S_{\os{\circ}{T}_0}$ 
and let $a^{(k)}_{{\rm crys}}\col 
(\os{\circ}{X}{}^{(k)}_{T_0}/\os{\circ}{T})_{\rm crys}
\lo (\os{\circ}{X}_{T_0}/\os{\circ}{T})_{\rm crys}$ 
be the morphism of topoi induced by $a^{(k)}$. 
Set 
$E_{\os{\circ}{X}{}^{(k)}_{T_0}/\os{\circ}{T}}
:=a^{(k)*}_{\rm crys}(E)$. 
By the usual Poincar\'{e} lemma and (\ref{lemm:pde}) (1),  
we have the following formula  
\begin{align*}  
Ru_{\os{\circ}{X}{}^{(k)}_{T_0}/\os{\circ}{T}*}
(E_{\os{\circ}{X}{}^{(k)}_{T_0}/\os{\circ}{T}}
\otimes_{\mab Z}
\vp^{(k)}_{\rm crys}(\os{\circ}{X}_{T_0}/\os{\circ}{T}))
& \os{\sim}{\lo} 
c^{(k)*}({\cal E})
\otimes_{{\cal O}_{\os{\circ}{\cal P}{}^{{\rm ex},(k)}}}
\Om^{\bul}_{\os{\circ}{\cal P}{}^{{\rm ex},(k)}/\os{\circ}{T}}
\otimes_{\mab Z}
\vp^{(k)}_{\rm zar}(\os{\circ}{X}_{T_0}/\os{\circ}{T}_0)\tag{1.3.16.14}\label{ali:eopudq}\\
&={\cal E}
\otimes_{{\cal O}_{\os{\circ}{\cal P}{}^{{\rm ex}}}}
\Om^{\bul}_{\os{\circ}{\cal P}{}^{{\rm ex},(k)}/\os{\circ}{T}}
\otimes_{\mab Z}
\vp^{(k)}_{\rm zar}(\os{\circ}{X}_{T_0}/\os{\circ}{T}_0)
\end{align*}
in ${\rm D}^+(\os{\circ}{f}{}^{(k)-1}({\cal O}_T))$.  
\par 


\begin{prop}\label{prop:npwf}
Let the notations be as above. 
Assume that $E$ is a flat quasi-coherent crystal of 
${\cal O}_{\os{\circ}{X}_T/\os{\circ}{T}}$-modules.  
Then the complex 
${\cal E}\otimes_{{\cal O}_{{\cal P}^{\rm ex}}}
{\Om}^{\bul}_{{\cal P}^{\rm ex}/\os{\circ}{T}}$ 
gives the subcomplex 
${\cal E}\otimes_{{\cal O}_{{\cal P}^{\rm ex}}}
P_k{\Om}^{\bul}_{{\cal P}^{\rm ex}/\os{\circ}{T}}$ 
for any $k\in {\mab Z}$. 
\end{prop}
\begin{proof} 
By (\ref{prop:om})  the natural morphism 
\begin{align*}  
{\cal E}\otimes_{{\cal O}_{{\cal P}^{\rm ex}}}
P_k{\Om}^i_{{\cal P}^{\rm ex}/\os{\circ}{T}}
\lo {\cal E}\otimes_{{\cal O}_{{\cal P}^{\rm ex}}}
\Om^i_{{\cal P}^{\rm ex}/\os{\circ}{T}}
\quad (i\in {\mab N}) \tag{1.3.17.1}\label{ali:otsd} 
\end{align*} 
is injective.  
We have only to prove that the connection 
\begin{equation*}
\nabla \col 
{\cal E}\lo  {\cal E}
\otimes_{{\cal O}_{{\cal P}^{\rm ex}}}
{\Om}^1_{{\cal P}^{\rm ex}/\os{\circ}{T}} 
\tag{1.3.17.2}\label{eqn:gtelp}
\end{equation*} 
induces a connection 
\begin{equation*}
{\cal E}
\lo {\cal E}
\otimes_{{\cal O}_{{\cal P}^{\rm ex}}}
P_0{\Om}^1_{{\cal P}^{\rm ex}/\os{\circ}{T}}.  
\tag{1.3.17.3}\label{eqn:gtoelp}
\end{equation*}
The problem is local. 
Replacing ${\cal P}$ with ${\cal P}^{\rm prex}$, 
we may assume that 
there exists the cartesian diagram (\ref{cd:pwtp}). 
Let $\eps \col \ol{\mathfrak D}\lo \os{\circ}{\ol{\mathfrak D}}$
be the natural morphism. Then 
\begin{equation*} 
\ol{\cal E}=\eps^*_{X_{\os{\circ}{T}_0}/\os{\circ}{T}}(E)_{\ol{\mathfrak D}}
=\eps^*(E_{\os{\circ}{\ol{\mathfrak D}}})=E_{\os{\circ}{\ol{\mathfrak D}}} 
\tag{1.3.17.4}\label{eqn:exse}
\end{equation*} 
since $\os{\circ}{\eps}$ is equal to 
${\rm id}_{\os{\circ}{\ol{\mathfrak D}}}$. 
Because $\os{\circ}{\ol{\cal P}}{}^{\rm ex}$ 
is formally smooth over $\os{\circ}{T}$, 
$E$ induces an integrable connection 
\begin{equation*}
\ol{\cal E}\lo \ol{\cal E}
\otimes_{{\cal O}_{\ol{\cal P}^{\rm ex}}}
\Om^1_{\os{\circ}{\ol{\cal P}}{}^{\rm ex}/\os{\circ}{T}}.  
\tag{1.3.17.5}\label{eqn:expe}
\end{equation*}
By (\ref{eqn:pops}) we obtain the connection 
(\ref{eqn:gtoelp}). 
\end{proof}

\par    
By abuse of notation, let us denote by 
$P$ the filtration 
$\{{\cal E}\otimes_{{\cal O}_{{\cal P}^{\rm ex}}}
P_k{\Om}^{\bul}_{{\cal P}^{\rm ex}/\os{\circ}{T}}\}_{k\in {\mab Z}}$ 
on 
${\cal E}\otimes_{{\cal O}_{{\cal P}^{\rm ex}}}
{\Om}^{\bul}_{{\cal P}^{\rm ex}/\os{\circ}{T}}$.  

\begin{defi}\label{defi:kai}
We call $P$ the {\it Poincar\'{e} filtration} on 
${\cal E}\otimes_{{\cal O}_{{\cal P}^{\rm ex}}}
{\Om}^{\bul}_{{\cal P}^{\rm ex}/\os{\circ}{T}}$. 
When $E={\cal O}_{\os{\circ}{X}_{T_0}/\os{\circ}{T}}$, we call 
the Poincar\'{e} filtration the {\it preweight filtration} 
on ${\cal O}_{{\mathfrak D}}
\otimes_{{\cal O}_{{\cal P}^{\rm ex}}}
{\Om}^{\bul}_{{\cal P}^{\rm ex}/\os{\circ}{T}}$. 
\end{defi}

\begin{coro}\label{coro:eqpk} 
Let $E$ be as in {\rm (\ref{prop:npwf})}.  
Let $k$ be a positive integer. 
Then 
\begin{equation*}  
{\rm gr}^P_k
({\cal E}\otimes_{{\cal O}_{{\cal P}^{\rm ex}}}
{\Om}^{\bul}_{{\cal P}^{\rm ex}/\os{\circ}{T}})
\os{\sim}{\lo}  
a^{(k-1)}_*
Ru_{\os{\circ}{X}{}^{(k-1)}_{T_0}/\os{\circ}{T}*}
(E_{\os{\circ}{X}{}^{(k-1)}_{T_0}/\os{\circ}{T}}
\otimes_{\mab Z}
\vp^{(k-1)}_{\rm crys}(\os{\circ}{X}_{T_0}/\os{\circ}{T}))[-k]  
\tag{1.3.19.1}\label{eqn:eoppd}
\end{equation*}
in ${\rm D}^+(f^{-1}({\cal O}_T))$.  
\end{coro}
\begin{proof}
By (\ref{eqn:prvin}) we have the following equality: 
\begin{equation*}  
{\rm gr}^P_k
({\cal E}\otimes_{{\cal O}_{{\cal P}^{\rm ex}}}
{\Om}^{\bul}_{{\cal P}^{\rm ex}/\os{\circ}{T}})
=
c^{(k-1)}_*c^{(k-1)*}({\cal E})
\otimes_{{\cal O}_{\os{\circ}{\cal P}{}^{{\rm ex},(k-1)}}}
\Om^{\bul}_{\os{\circ}{\cal P}{}^{{\rm ex},(k-1)}/\os{\circ}{T}}
\otimes_{\mab Z}
\vp^{(k-1)}_{\rm zar}(\os{\circ}{X}_{T_0}/\os{\circ}{T}_0)[-k]. 
\tag{1.3.19.2}\label{eqn:prvien} 
\end{equation*}  
(We can easily check that (\ref{eqn:prvien}) is indeed the equality 
of complexes.)
(\ref{eqn:eoppd}) immediately follows from (\ref{ali:eopudq}).  
\end{proof}

The following is an SNCL version of \cite[(4.8)]{nh3}: 

\begin{prop}[{\bf The contravariant functoriality of the Poincar\'{e} residue isomorphism}]\label{prop:rescos} 
Let $E$ be as in {\rm (\ref{prop:npwf})}. 
Let $S$ be as above and let 
$u \col S\lo S'$ 
be a morphism of family of log points.  
Let $(T,{\cal J},\del) \lo (T',{\cal J}',\del')$ be a morphism of 
log PD-enlargements over $u$. 
Assume that $p$ is locally nilpotent on $\os{\circ}{T}{}'$.    
Set $T'_0:=\ul{\rm Spec}^{\log}_{T'}({\cal O}_{T'}/{\cal J}')$. 
Let $X'/S'$ be an SNCL scheme.  
Let  $X'_{\os{\circ}{T}{}'_0} \os{\subset}{\lo} \ol{\cal P}{}'$ be  
an immersion into a log smooth scheme over 
$\ol{S'(T')^{\nat}}$ fitting into the following commutative diagram 
over the morphism $\ol{S(T)^{\nat}}\lo \ol{S'(T')^{\nat}}:$ 
\begin{equation*} 
\begin{CD} 
X_{\os{\circ}{T}_0} @>{\subset}>> \ol{\cal P} \\
@V{g\vert_{X_{\os{\circ}{T}_0}}}VV @VV{\ol{g}}V \\ 
X'_{\os{\circ}{T}{}'_0} @>{\subset}>> \ol{\cal P}{}'. \\  
\end{CD} 
\tag{1.3.20.1}\label{eqn:pvpq} 
\end{equation*} 
Set ${\cal P}':=\ol{\cal P}{}'\times_{\ol{S'(T')^{\nat}}}S'(T')^{\nat}$. 
Let $g\col {\cal P}\lo {\cal P}'$ 
be the induced morphism by $\ol{g}$.
Let $g^{\rm ex}\col {\cal P}^{\rm ex}\lo 
{\cal P}'{}^{\rm ex}$ be the induced morphism by $g$. 
Let $a'{}^{(k)}\col \os{\circ}{X}{}'^{(k)}_{T'_0}
\lo \os{\circ}{X}{}'_{T'_0}$ $(k\in {\mab N})$
be the natural morphism of schemes over $\os{\circ}{T}{}'_0$.  
Assume that, for each point 
$x\in \os{\circ}{\cal P}{}^{\rm ex}$ 
and for each member $m$ of the minimal generators 
of $M_{{\cal P}^{\rm ex},x}/{\cal O}^*_{{\cal P}^{\rm ex},x}$, 
there exists a unique member $m'$ of 
the minimal generators of 
$M_{{\cal P}{}'^{\rm ex},\os{\circ}{g}(x)}
/{\cal O}^*_{{\cal P}{}'^{\rm ex},\os{\circ}{g}(x)}$ 
such that $g^*(m')= m$ and such that the image of 
the other minimal generators of 
$M_{{\cal P}{}'^{\rm ex},{\os{\circ}{g}(x)}}
/{\cal O}^*_{{\cal P}{}'^{\rm ex},\os{\circ}{g}(x)}$ by $g^*$ 
are the trivial element of 
$M_{{\cal P}^{\rm ex},{x}}/{\cal O}^*_{{\cal P}^{\rm ex},x}$. 
Let ${\mathfrak D}'$ be the log PD-envelope of 
the immersion $X'_{\os{\circ}{T}{}'_0}\os{\sus}{\lo} {\cal P}'$ over 
$(S'(T')^{\nat},{\cal J}',\del')$. 
Let $\os{\circ}{\mathfrak D}{}'^{(k)}$ 
be the analogous schemes to 
$\os{\circ}{\mathfrak D}{}^{(k)}$. 
Let $E'$ $($resp.~$({\cal E}',\nabla'))$ be 
an analogous sheaf $($resp.~sheaf with integral connection$)$ 
to $E$ $($resp.~$({\cal E},\nabla))$ 
for $X'_{\os{\circ}{T}{}'_0}/(S(T)^{\nat},{\cal J}',\del')$. 
Let 
$a'{}^{(k)}\col \os{\circ}{X}{}'^{(k)}_{T'_0}\lo \os{\circ}{X}{}'_{\os{\circ}{T}{}'_0}$ 
$($resp.~$c'{}^{(k)}\col \os{\circ}{\mathfrak D}{}'^{(k)}\lo \os{\circ}{\mathfrak D}{}')$ 
be an analogous morphism to $a^{(k)}$ $($resp.~$c^{(k)})$ 
for $X'_{\os{\circ}{T}{}'_0}$ 
$($resp.~the PD-immersion 
$X'_{\os{\circ}{T}{}'_0}\os{\sus}{\lo}{\mathfrak D}')$.  
Assume that we are given a morphism 
$E'\lo (\os{\circ}{g}\vert_{X_{\os{\circ}{T}_0}})_{{\rm crys}*}(E)$. 
Then the following diagram is commutative$:$ 
\begin{equation*} 
\begin{CD}  
g^{\rm ex}_*({\rm Res}): 
g^{\rm ex}_*({\rm gr}_k^P
({\cal E}\otimes_{{\cal O}_{{\cal P}^{\rm ex}}}
{\Om}^{\bul}_{{\cal P}^{\rm ex}/\os{\circ}{T}}))
@>{\sim}>> \\ 
@AAA \\
{\rm Res}: 
{\rm gr}_k^P({\cal E}'\otimes_{{\cal O}_{{\cal P}'{}^{\rm ex}}}
{\Om}^{\bul}_{{\cal P}'{}^{\rm ex}/\os{\circ}{T}{}'})
@>{\sim}>> 
\end{CD}
\tag{1.3.20.2}\label{eqn:regdmrp} 
\end{equation*} 
\begin{equation*} 
\begin{CD}  
g^{\rm ex}_*a^{(k-1)}_*
(c^{(k-1)*}({\cal E})
\otimes_{{\cal O}_{{\cal P}^{\rm ex}}}
{\Om}^{\bul}_{\os{\circ}{\cal P}{}^{{\rm ex},(k-1)}/\os{\circ}{T}}
\otimes_{\mab Z}
\vp^{(k-1)}_{\rm zar}(\os{\circ}{X}_{T_0}/\os{\circ}{T}_0)[-k])\\
@AAA \\
a'{}^{(k-1)}_*(c'{}^{(k-1)*}({\cal E}')
\otimes_{{\cal O}_{{\cal P}{}'^{\rm ex}}}
{\Om}^{\bul}_{\os{\circ}{\cal P}{}^{'{\rm ex},(k-1)}/\os{\circ}{T}{}'}
\otimes_{\mab Z}
\vp^{(k-1)}_{\rm zar}(\os{\circ}{X}{}'_{T_0}/\os{\circ}{T}{}'_0)[-k]). 
\end{CD} 
\end{equation*} 
\end{prop}
\begin{proof} 
We immediately obtain (\ref{prop:rescos}) 
from the condition 
$g^*(m')=m$ and by the definition of the morphism 
$g^{{\rm ex},(k)}\col \os{\circ}{\cal P}{}^{{\rm ex},(k)}
\lo \os{\circ}{\cal P}{}'^{{\rm ex},(k)}$ 
in (\ref{prop:mmoo}).  
\end{proof}

\par 
Next we give an intrinsic description of 
$P_0({\cal E}\otimes_{{\cal O}_{{\cal P}^{\rm ex}}}
{\Om}^{\bul}_{{\cal P}^{\rm ex}/\os{\circ}{S}})$ 
in ${\rm D}^+(f^{-1}({\cal O}_T))$.   

\par 
By (\ref{eqn:pops}) 
the natural morphism 
$\Om^{\bul}_{\os{\circ}{\ol{\cal P}}{}^{\rm ex}/\os{\circ}{T}} 
\lo \Om^{\bul}_{\os{\circ}{\cal P}{}^{\rm ex}/\os{\circ}{T}}$ induces a morphism 
$P_0\Om^{\bul}_{{\cal P}^{\rm ex}/\os{\circ}{T}}
\lo c^{(0)}_*(\Om^{\bul}_{\os{\circ}{\cal P}{}^{\rm ex}/\os{\circ}{T}})$. 
Consequently we have the following morphism
\begin{align*} 
P_0
({\cal E}
\otimes_{{\cal O}_{{\cal P}^{\rm ex}}}{\Om}^{\bul}_{{\cal P}^{\rm ex}/\os{\circ}{T}}) 
\lo 
c^{(0)}_*(c^{(0)*}({\cal E})
\otimes_{{\cal O}_{{\cal P}^{\rm ex}}}
\Om^{\bul}_{\os{\circ}{\cal P}{}^{\rm ex}/\os{\circ}{T}}). 
\tag{1.3.20.3}\label{eqn:oppkps}
\end{align*} 
\par 
Let $\{\os{\circ}{X}_{\lam}\}_{\lam \in \Lam}$ be 
the set of the smooth components 
of $\os{\circ}{X}_{T_0}$: 
$\os{\circ}{X}_{T_0}=\bigcup_{\lam \in \Lam}\os{\circ}{X}_{\lam}$; 
$\os{\circ}{X}_{\lam}$ is smooth over $\os{\circ}{S}_{\os{\circ}{T}_0}=\os{\circ}{T}_0$. 
Let $\os{\circ}{\cal P}{}^{\rm ex}_{\lam}$ 
be the smooth closed subscheme of 
$\os{\circ}{\cal P}{}^{\rm ex}$ over $\os{\circ}{T}$ 
which is topologically isomorphic to $\os{\circ}{X}_{\lam}$. 
To define a morphism (\ref{eqn:odpaps}) below, 
we fix a total order on $\Lam$ once and for all. 
For different $\lam_0<\lam_1< \ldots< \lam_k$, 
set $\os{\circ}{X}_{\lam_0\cdots \lam_k}
:=\os{\circ}{X}_{\lam_0}\cap \cdots \cap \os{\circ}{X}_{\lam_k}$ 
and $\os{\circ}{\cal P}{}^{\rm ex}_{\lam_0\cdots \lam_k}
:=\os{\circ}{\cal P}{}^{\rm ex}_{\lam_0}\cap 
\cdots \cap \os{\circ}{\cal P}{}^{\rm ex}_{\lam_k}$. 
Let $\os{\circ}{\mathfrak D}_{\lam_0\cdots \lam_k}$ be 
the PD-envelope of the immersion 
$\os{\circ}{X}_{\lam_0\cdots \lam_k} \os{\sus}{\lo} 
\os{\circ}{\cal P}{}^{\rm ex}_{\lam_0\cdots \lam_k}$ 
over $(\os{\circ}{T},{\cal J},\del)$ and let 
$c_{\lam_0\cdots \lam_k} \col 
\os{\circ}{\mathfrak D}_{\lam_0\cdots \lam_k} 
\lo \os{\circ}{\mathfrak D}$ 
be the natural morphism of schemes. 
Let  
\begin{equation*} 
\iota^{\lam_0\cdots \hat{\lam}_j \cdots \lam_{k+1}}_{\lam_0 \cdots \lam_{k+1}} 
\col 
\os{\circ}{X}_{\lam_0}\cap 
\cdots 
\cap \os{\circ}{X}_{\lam_{k+1}} 
\os{\sus}{\lo} 
\os{\circ}{X}_{\lam_0}\cap \cdots 
\cap \hat{\os{\circ}{X}}_{\lam_j} \cap \cdots 
\cap 
\os{\circ}{X}_{\lam_{k+1}}
\tag{1.3.20.4}\label{eqn:odikps}
\end{equation*} 
be the natural closed immersion.  
Let 
\begin{align*} 
&\iota^{(k)*} \col   
c^{(k)}_*
(c^{(k)*}({\cal E})
\otimes_{{\cal O}_{\os{\circ}{\cal P}{}^{{\rm ex},(k)}}}
\Om^{\bul}_{\os{\circ}{\cal P}{}^{{\rm ex},(k)}/\os{\circ}{T}}
\otimes_{\mab Z}
\vp^{(k)}_{\rm zar}
(\os{\circ}{\cal P}{}^{\rm ex}/\os{\circ}{T})) 
\tag{1.3.20.5}\label{eqn:odpaps}\\
& \lo 
c^{(k+1)}_*
(c^{(k+1)*}({\cal E})
\otimes_{{\cal O}_{\os{\circ}{\cal P}{}^{{\rm ex},(k+1)}}}
\Om^{\bul}_{\os{\circ}{\cal P}{}^{{\rm ex},(k+1)}
/\os{\circ}{T}}
\otimes_{\mab Z}\vp^{(k+1)}_{\rm zar}
(\os{\circ}{\cal P}{}^{\rm ex}/\os{\circ}{T}))  
\end{align*} 
be the natural \v{C}ech morphism, which is the summation of 
the following morphism with respect to 
$\lam_j$ and $0\leq j \leq k+1$ 
(cf.~\cite[(5.0.5)]{ndw}): 
\begin{align*} 
& 
c_{\lam_0\cdots \hat{\lam}_j \cdots \lam_{k+1}*}
(c^*_{\lam_0\cdots \hat{\lam}_j \cdots \lam_{k+1}}({\cal E}) 
\otimes_{
{\cal O}_{{\cal P}{}^{\rm ex}_{\lam_0\cdots \hat{\lam}_j \cdots \lam_{k+1}}}}
\Om^{\bul}_{{\cal P}{}^{\rm ex}_{\lam_0\cdots \hat{\lam}_j \cdots \lam_{k+1}}/\os{\circ}{T}}
\otimes_{\mab Z}
\vp_{{\rm zar}{\lam_0\cdots \hat{\lam}_j \cdots \lam_{k+1}}}
(\os{\circ}{\cal P}{}^{\rm ex}/\os{\circ}{T})) 
\tag{1.3.20.6}\label{eqn:odkps}\\
& 
\lo c_{\lam_0\cdots \lam_{k+1}*}
(c^*_{\lam_0\cdots \lam_{k+1}}({\cal E})\otimes_{
{\cal O}_{{\cal P}{}^{\rm ex}_{\lam_0\cdots \lam_{k+1}}}}
\Om^{\bul}_{{\cal P}^{\rm ex}_{\lam_0\cdots \lam_{k+1}}/\os{\circ}{T}}
\otimes_{\mab Z}
\vp_{{\rm zar}{\lam_0\cdots \lam_{k+1}}}
(\os{\circ}{\cal P}{}^{\rm ex}/\os{\circ}{T})))
\end{align*}  
\begin{equation*}
{\om_{\lam_0\cdots \hat{\lam}_j \cdots \lam_{k+1}}}
\otimes(\lam_0\cdots \hat{\lam}_j \cdots \lam_{k+1})
\lom 
{(-1)^j\iota^{\lam_0\cdots \hat{\lam}_j \cdots \lam_{k+1}*}_{\lam_0\cdots \lam_{k+1}}
(\om_{\lam_0\cdots \hat{\lam}_j \cdots \lam_{k+1}})}
\otimes(\lam_0\cdots \lam_{k+1}).  
\end{equation*}

\begin{prop}
[{\rm {\bf cf.~\cite[Lemma 3.15.1]{msemi}, \cite[(6.29)]{ndw}}}]\label{prop:csrl}  
The following sequence obtained by the morphism {\rm (\ref{eqn:oppkps})} is exact$:$
\begin{equation*} 
0 \lo P_0
({\cal E}
\otimes_{{\cal O}_{{\cal P}^{\rm ex}}}{\Om}^{\bul}_{{\cal P}^{\rm ex}/\os{\circ}{T}}) 
\lo c^{(0)}_*(c^{(0)*}({\cal E})
\otimes_{{\cal O}_{{\os{\circ}{\cal P}{}^{{\rm ex},(0)}}}}
\Om^{\bul}_{\os{\circ}{\cal P}{}^{{\rm ex},(0)}/\os{\circ}{T}}\otimes_{\mab Z}
\vp^{(0)}_{\rm zar}(\os{\circ}{\cal P}{}^{\rm ex}
/\os{\circ}{T})) 
\tag{1.3.21.1}\label{eqn:cptle} 
\end{equation*} 
$$\os{\iota^{(0)*}}{\lo} {c}{}^{(1)}_*(c^{(1)*}({\cal E})
\otimes_{{\cal O}_{\os{\circ}{\cal P}{}^{{\rm ex},(1)}}}
\Om^{\bul}_{\os{\circ}{\cal P}{}^{{\rm ex},(1)}
/\os{\circ}{T}}
\otimes_{\mab Z}
\vp^{(1)}_{\rm zar}(\os{\circ}{\cal P}{}^{\rm ex}
/\os{\circ}{T})) 
\os{\iota^{(1)*}}{\lo} \cdots . $$ 
\end{prop} 
\begin{proof} 
(The idea of the following proof is the same as 
that of (\ref{prop:injf}).)  
The question is local.  
We may assume that $E={\cal O}_{\os{\circ}{X}_{T_0}/\os{\circ}{T}}$.  
We may assume that there exists the cartesian diagram 
(\ref{cd:pwtp}). 
Let the notations be as in the proof of (\ref{prop:injf}) 
for the case ${\cal Q}=\ol{\cal P}$.  
\par   
Set ${\cal P}':=\ol{\cal P}{}'\times_{\ol{S(T)^{\nat}}}S(T)^{\nat}$. 
Let  $b'{}^{(k)} \col \os{\circ}{\cal P}{}'{}^{(k)} \lo 
\os{\circ}{\cal P}{}'$ $(k\in {\mab N})$ 
be the analogous morphism to 
$b^{(k)}\col \os{\circ}{\cal P}{}^{(k)} \lo \os{\circ}{\cal P}$. 
By \cite[(4.2.2) (a), (c)]{di} 
we have the following exact sequence 
\begin{equation*}
0 \lo 
\Om^{\bul}_{\os{\circ}{\ol{\cal P}}{}'/\os{\circ}{T}}/
{\Om}^{\bul}_{{\ol{\cal P}}{}'/\os{\circ}{T}}
(-\os{\circ}{\cal P}{}')  
\lo b'{}^{(0)}_{*}
(\Om^{\bul}_{\os{\circ}{\cal P}{}'{}^{(0)}
/\os{\circ}{T}}
\otimes_{\mab Z}\vp^{(0)}_{\rm zar}
(\os{\circ}{\cal P}{}'/\os{\circ}{T})) 
\lo \cdots.  
\end{equation*} 
Because  
\begin{align*} 
{\cal O}_{\mathfrak D}\otimes_{{\cal O}_{{\cal P}}}
(\Om^{\bul}_{\os{\circ}{\ol{\cal P}}/\os{\circ}{T}}/
{\Om}^{\bul}_{{\ol{\cal P}}/\os{\circ}{T}}
(-\os{\circ}{\cal P}))
\os{\sim}{\lo} & 
(\Om^{\bul}_{\os{\circ}{\ol{\cal P}}{}'/\os{\circ}{T}}/
\Om^{\bul}_{{\ol{\cal P}}{}'/\os{\circ}{T}}
(-\os{\circ}{\cal P}{}'))
\otimes_{{\cal O}_T} \\
{} & 
{\cal O}_T\langle x_{d+1},\ldots, x_{d'} \rangle 
\otimes_{{\cal O}_{{\cal P}''}}
\Om^{\bul}_{{\cal P}''/\os{\circ}{T}},   
\end{align*}
because 
\begin{align*}  
{\cal O}_{\mathfrak D}\otimes_{{\cal O}_{\cal P}}
b{}^{(k)}_{*}(\Om^{\bul}_{\os{\circ}{\cal P}{}^{(k)}
/\os{\circ}{T}}
\otimes_{\mab Z}\vp^{(k)}_{\rm zar}(\os{\circ}{\cal P}/\os{\circ}{T})) 
\os{\sim}{\lo} & 
b'{}^{(k)}_{*}
(\Om^{\bul}_{\os{\circ}{\cal P}{}'{}^{(k)}/\os{\circ}{T}}
\otimes_{\mab Z}\vp^{(k)}_{\rm zar}(\os{\circ}{\cal P}{}'/\os{\circ}{T}))  
\otimes_{{\cal O}_T} \\
{} & 
{\cal O}_T\langle x_{d+1},\ldots, x_{d'} \rangle 
\otimes_{{\cal O}_{{\cal P}''}}
\Om^{\bul}_{{\cal P}''/\os{\circ}{T}}   
\end{align*} 
and because the complex ${\cal O}_T
\langle x_{d+1},\ldots, x_{d} \rangle 
\otimes_{{\cal O}_{{\cal P}''}}
\Om^{\bul}_{{\cal P}''/\os{\circ}{T}}$ 
consists of free ${\cal O}_T$-modules, 
we see that the following sequence is exact: 
\begin{equation*} 
0 \lo  
{\cal O}_{\mathfrak D}
\otimes_{{\cal O}_{{\cal P}}}
(\Om^{\bul}_{\os{\circ}{\ol{\cal P}}/\os{\circ}{T}}/
{\Om}^{\bul}_{{\ol{\cal P}}/\os{\circ}{T}}
(-\os{\circ}{\cal P}))
\lo {\cal O}_{\mathfrak D}\otimes_{{\cal O}_{\cal P}}
b^{(0)}_{*}
(\Om^{\bul}_{\os{\circ}{\cal P}{}^{(0)}/\os{\circ}{T}}
\otimes_{\mab Z}\vp^{(0)}_{\rm zar}
(\os{\circ}{\cal P}/\os{\circ}{T})) 
\lo \cdots.
\tag{1.3.21.2}\label{ali:opot}    
\end{equation*} 
By (\ref{prop:incol}) (4) and (\ref{eqn:pops}) we have 
the following isomorphism: 
\begin{align*}
\Om^{\bul}_{\os{\circ}{\ol{\cal P}}/\os{\circ}{T}}/
\Om^{\bul}_{{\ol{\cal P}}/\os{\circ}{T}}(-\os{\circ}{\cal P})
\simeq P_0\Om^{\bul}_{\ol{\cal P}/\os{\circ}{T}}
/
\Om^{\bul}_{{\ol{\cal P}}/\os{\circ}{T}}(-\os{\circ}{\cal P})
=P_0\Om^{\bul}_{{\cal P}/\os{\circ}{T}}. 
\tag{1.3.21.3}\label{ali:opt} 
\end{align*} 
By (\ref{lemm:pde}) (1), (\ref{ali:opot}) and 
(\ref{ali:opt}) we see that the sequence (\ref{eqn:cptle}) is exact. 
\end{proof}  

The following is an intrinsic description of  
$P_0({\cal E}\otimes_{{\cal O}_{{\cal P}^{\rm ex}}}
{\Om}^{\bul}_{{\cal P}^{\rm ex}/\os{\circ}{T}})$:  

\begin{prop}\label{prop:repmf}
There exists the following isomorphism 
\begin{align*}
P_0({\cal E}
\otimes_{{\cal O}_{{\cal P}^{\rm ex}}}
{\Om}^{\bul}_{{\cal P}^{\rm ex}/\os{\circ}{T}})
\os{\sim}{\lo} 
&{\rm MF}(a^{(0)}_{*}
Ru_{\os{\circ}{X}{}^{(0)}_{T_0}/\os{\circ}{T}*}
(E_{\os{\circ}{X}{}^{(0)}_{T_0}/\os{\circ}{T}}
\otimes_{\mab Z}\vp^{(0)}_{\rm crys}
(\os{\circ}{X}_{T_0}/\os{\circ}{T}))
\tag{1.3.22.1}\label{ali:pdte}\\
&\lo 
{\rm MF}(a^{(1)}_{*}
Ru_{\os{\circ}{X}{}^{(1)}_{T_0}/\os{\circ}{T}*}
(E_{\os{\circ}{X}{}^{(1)}_{T_0}/\os{\circ}{T}}
\otimes_{\mab Z}\vp^{(1)}_{\rm crys}
(\os{\circ}{X}_{T_0}/\os{\circ}{T}))  \\
&\lo \cdots \\
&\lo 
{\rm MF}(a^{(m)}_{*}
Ru_{\os{\circ}{X}{}^{(m)}_{T_0}/\os{\circ}{T}*}
(E_{\os{\circ}{X}{}^{(m)}_{T_0}/\os{\circ}{T}}
\otimes_{\mab Z}\vp^{(m)}_{\rm crys}
(\os{\circ}{X}_{T_0}/\os{\circ}{T}))\\
& \lo \cdots
\cdots)\cdots ). 
\end{align*}
Here ${\rm MF}$ means the  mapping fiber of 
a morphism of complexes of derived categories 
$($cf.~Notations $(14))$. 
The isomorphism {\rm (\ref{ali:pdte})} 
is independent of the choice of the immersion 
$X_{\os{\circ}{T}_0}\os{\sus}{\lo} \ol{\cal P}$ over $\ol{S(T)^{\nat}}$.  
In particular, 
$P_0({\cal E}\otimes_{{\cal O}_{{\cal P}^{\rm ex}}}
{\Om}^{\bul}_{{\cal P}^{\rm ex}/\os{\circ}{T}})$ 
is independent of the choice of the immersion 
$X_{\os{\circ}{T}_0}\os{\sus}{\lo} \ol{\cal P}$.  
\end{prop} 
\begin{proof}
By  (\ref{eqn:cptle}), (\ref{ali:eopudq}) and using mapping fibers repeatedly, 
we obtain the description (\ref{ali:pdte}). 
\par 
Assume that we are given another immersion 
$X\os{\sus}{\lo} \ol{\cal P}{}'$ over $\ol{S(T)^{\nat}}$. 
Then, by considering the product 
$\ol{\cal P}\times_{\ol{S(T)^{\nat}}}\ol{\cal P}{}'$, we may assume 
that there exists a morphism 
$\ol{\cal P}\lo \ol{\cal P}{}'$ extending ${\rm id}_{X_{\os{\circ}{T}_0}}$. 
Let $C^{\bul}$ be the complex on 
the right hand side on (\ref{ali:pdte}). 
Let ${\cal E}'$ be the analogue of ${\cal E}$ for $\ol{\cal P}{}'$. 
Then we have the following commutative diagram: 
\begin{equation*} 
\begin{CD} 
P_0({\cal E}
\otimes_{{\cal O}_{{\cal P}^{\rm ex}}}
{\Om}^{\bul}_{{\cal P}^{\rm ex}/\os{\circ}{T}})
@>{\sim}>>C^{\bul} \\
@A{\simeq}AA @| \\
P_0({\cal E}'\otimes_{{\cal O}_{{\cal P}{}'^{\rm ex}}}
{\Om}^{\bul}_{{\cal P}{}'^{\rm ex}/\os{\circ}{T}})
@>{\sim}>>C^{\bul}
\end{CD}
\end{equation*} 
Hence we obtain the desired independence.  
\end{proof} 

\begin{coro}\label{coro:indfo} 
The filtered complex 
$({\cal E}\otimes_{{\cal O}_{{\cal P}^{\rm ex}}}
{\Om}^{\bul}_{{\cal P}^{\rm ex}/\os{\circ}{T}},P)$ 
is independent of the choice of the immersion 
$X_{\os{\circ}{T}_0}\os{\sus}{\lo} \ol{\cal P}$. 
That is, in the notations in {\rm (\ref{prop:rescos})}, assume that 
$g\vert_{X_{\os{\circ}{T}_0}}={\rm id}_{X_{\os{\circ}{T}_0}}$ 
and that the morphism $E'\lo (\os{\circ}{g}\vert_X)_{{\rm crys}*}(E)$ 
is equal to ${\rm id}_E$.  
Then $({\cal E}\otimes_{{\cal O}_{{\cal P}^{\rm ex}}}
{\Om}^{\bul}_{{\cal P}^{\rm ex}/\os{\circ}{T}},P)
=
({\cal E}'\otimes_{{\cal O}_{{\cal P}'{}^{\rm ex}}}
{\Om}^{\bul}_{{\cal P}'{}^{\rm ex}/\os{\circ}{T}},P)$.  
\end{coro} 
\begin{proof} 
This immediately follow from (\ref{ali:pdte}), (\ref{eqn:eoppd}) 
and the induction on $k$ for $P_k$. 
\end{proof} 

\begin{rema}\label{rema:dily} 
(\ref{eqn:eoppd}), (\ref{prop:hkt}) below and the descending induction 
also show 
(\ref{coro:indfo}). 
\end{rema} 

\par 
Next we give another description of 
$({\cal E}\otimes_{{\cal O}_{{\cal P}^{\rm ex}}}
{\Om}^{\bul}_{{\cal P}^{\rm ex}/\os{\circ}{T}},P)$, 
which is a log and filtered version of \cite[0 (3.1.6)]{idw}. 

\begin{defi}\label{defi:tafs} 
Let $T$ be a fine log $($formal$)$ scheme. 
Let $Y/T$ be a fine log $($formal$)$ scheme over $T$. 
Let $(Y\times_TY)^{\rm ex}$ be the exactification of 
the diagonal immersion 
$Y\os{\sus}{\lo} Y\times_TY$.  
Let $I$ be the defining ideal sheaf of 
the diagonal immersion 
$Y\os{\sus}{\lo}(Y\times_TY)^{\rm ex}$. 
Set 
$P(Y^{\rm ex}/T):=
\ul{\rm Spec}^{\log}_{(Y\times_TY)^{\rm ex}}
({\cal O}_{(Y\times_TY)^{\rm ex}}/I^2)$. 
We call $P(Y^{\rm ex}/T)$ the {\it log infinitesimal neighborhood} of 
$Y^{\rm ex}$ in $(Y\times_TY)^{\rm ex}$. 
\end{defi} 

We recall the following proposition: 

\begin{prop}[{\rm {\bf \cite[(5.8)]{klog1}, \cite[Proposition 3.2.5]{s1}}}]\label{prop:omlinf}  
Let the notations be as in {\rm (\ref{defi:tafs})}.  
Then there exists a canonical isomorphism 
\begin{align*} 
\Om^1_{Y/T}=I/I^2. 
\tag{1.3.26.1}\label{ali:omyti} 
\end{align*} 
\end{prop}
 
\begin{prop}\label{prop:tty}  
Let the notations be as in {\rm (\ref{prop:exex})}. 
Let $J$ be the defining ideal sheaf of the diagonal immersion 
${\cal Q}^{\rm ex}\os{\sus}{\lo}
({\cal Q}^{\rm ex}\times_T{\cal Q}^{\rm ex})^{\rm ex}$.   
Then $P({\cal Q}^{\rm ex}/T)=
\ul{\rm Spec}^{\log}_{({\cal Q}^{\rm ex}\times_T{\cal Q}^{\rm ex})^{\rm ex}}
({\cal O}_{({\cal Q}^{\rm ex}\times_T{\cal Q}^{\rm ex})^{\rm ex}}/J^2)$. 
\end{prop}  
\begin{proof} 
This immediately follows from (\ref{prop:exex}) (1). 
\end{proof}

\par 
Let $(T,{\cal J},\del)$ be a fine log PD-scheme. 
Let $g\col T\lo U$ be a morphism of fine log (formal) schemes. 
Consider the sheaf $\bigoplus_{i\in {\mab N}}\Om^i_{T/U}$ of 
dga's over $g^{-1}({\cal O}_U)$. 
Denote by $\Om^{\bul}_{T/U,\del}$ the quotient of 
$\bigoplus_{i\in {\mab N}}\Om^i_{T/U}$ by the ideal sheaf generated 
by local sections of the form $d(a^{[n]})-\del_{n-1}(a)da$ $(a\in {\cal J}, n\in {\mab Z}_{\geq 1})$.

The following (1) and (2) are log versions of \cite[0 (3.1.4), (3.1.6)]{idw}. 
We also recall Illusie's result itself for our memory. 

\begin{lemm}\label{lemm:esd0}
$(1)$ Let $(U,{\cal K},\eps)$ be a fine log PD-scheme. 
Set $U_0:=\ul{\rm Spec}^{\log}_U({\cal O}_U/{\cal K})$. 
Let $Z$ be a fine log scheme over $U_0$. 
Let $Z\os{\sus}{\lo} {\cal Z}$ be an immersion into 
a $($not necessarily log smooth$)$ fine log $($formal$)$ scheme 
over $U$. 
Let ${\mathfrak E}$ be the log PD-envelope of the immersion 
$Z\os{\sus}{\lo} {\cal Z}$ over $(U,{\cal K},\eps)$. 
Let $g\col Z\lo U$ be the structural morphism. 
Then the following hold$:$ 
\par 
$(1)$ The derivation 
$d\col {\cal O}_{{\cal Z}^{\rm ex}}\lo \Om^1_{{\cal Z}^{\rm ex}/U}$ extends uniquely to 
a derivation $d\col {\cal O}_{\mathfrak E}\lo 
{\cal O}_{\mathfrak E}\otimes_{{\cal O}_{{\cal Z}^{\rm ex}}}
\Om^1_{{\cal Z}^{\rm ex}/U}$ such that 
$d(a^{[n]})=a^{[n-1]}da$ for any local section $a$ of the PD-ideal sheaf of 
${\cal O}_{\mathfrak E}$ and $n\in {\mab Z}_{\geq 1}$. 
\par 
$(2)$ 
There exists the following isomorphism of dga's over $g^{-1}({\cal O}_{U}):$ 
\begin{equation*}
\Om^{\bul}_{{\mathfrak E}/U, [~]}
\os{\sim}{\lo} 
{\cal O}_{\mathfrak E}
\otimes_{{\cal O}_{{\cal Z}^{\rm ex}}}
\Om^{\bul}_{{\cal Z}^{\rm ex}/U}.   
\tag{1.3.28.1}\label{eqn:fwniwu}
\end{equation*}   
\par 
$(3)$ {\rm (\cite[0 (3.1.6)]{idw})} There exists the following isomorphism of 
dga's over $g^{-1}({\cal O}_{U}):$ 
\begin{equation*}
\Om^{\bul}_{\os{\circ}{\mathfrak E}/\os{\circ}{U},[~]}
\os{\sim}{\lo} 
{\cal O}_{\mathfrak E}
\otimes_{{\cal O}_{{\cal Z}^{\rm ex}}}
\Om^{\bul}_{\os{\circ}{\cal Z}{}^{\rm ex}/\os{\circ}{U}}.   
\tag{1.3.28.2}\label{eqn:fwu}
\end{equation*}  
\end{lemm} 
\begin{proof} 
(1): We may assume that the immersion 
$Z\os{\sus}{\lo}{\cal Z}$ is closed. 
Let $i\col Z\os{\sus}{\lo}{\cal Z}^{\rm ex}$ be the exactification of this closed immersion. 
By the log version of \cite[6.2.1 Lemma]{bob}, 
$i_{{\rm crys}*}({\cal O}_{Z/U})$ is a crystal of 
${\cal O}_{{\cal Z}^{\rm ex}/U}$-modules and 
$i_{{\rm crys}*}({\cal O}_{Z/U})_{{\cal Z}^{\rm ex}}={\cal O}_{\mathfrak E}$. 
Hence we have the following log HPD-stratification 
\begin{align*} 
{\cal O}_{\mathfrak E}\otimes_{{\cal O}_{{\cal Z}^{\rm ex}}}
{\cal O}_{{\mathfrak D}_{{\cal Z}^{\rm ex}/U}(1)}
\simeq   
{\cal O}_{{\mathfrak D}_{Z,\eps}({\cal Z}^{\rm ex}\times_U{\cal Z}^{\rm ex})}
\simeq   
{\cal O}_{{\mathfrak D}_{{\cal Z}^{\rm ex}/U}(1)}
\otimes_{{\cal O}_{{\cal Z}^{\rm ex}}}
{\cal O}_{\mathfrak E}.
\end{align*} 
Here ${\mathfrak D}_{{\cal Z}^{\rm ex}/U}(1)$ and 
${\mathfrak D}_{Z,\eps}({\cal Z}^{\rm ex}\times_U{\cal Z}^{\rm ex})$ 
are the log PD-envelopes of the immersions 
${\cal Z}^{\rm ex}\os{\sus}{\lo} {\cal Z}^{\rm ex}\times_U{\cal Z}^{\rm ex}$ 
and $Z\os{\sus}{\lo} {\cal Z}^{\rm ex}\times_U{\cal Z}^{\rm ex}$ over 
$(U,{\cal K},\eps)$, respectively.  
Consequently we have an integrable connection on 
${\cal O}_{\mathfrak E}$. 
This is nothing but the desired connection, which is easily proved 
(see \cite[6.4 Exercise]{bob}). The uniqueness of 
$d\col {\cal O}_{\mathfrak E}\lo 
{\cal O}_{\mathfrak E}\otimes_{{\cal O}_{{\cal Z}^{\rm ex}}}
\Om^1_{{\cal Z}^{\rm ex}/U}$ is obvious. 
\par 
(By the proof of (2) below, 
we can also give another proof of (1) directly as in \cite[0 (3.1.4)]{idw}.)
\par 
(2): (We need a nontrivial additional argument 
to the proof of \cite[0 (3.1.4)]{idw}.)  
\par 
Set ${\cal J}:={\rm Ker}({\cal O}_{{\cal Z}^{\rm ex}}\lo {\cal O}_{Z})$
and $\ol{\cal J}:={\rm Ker}({\cal O}_{\mathfrak E}\lo {\cal O}_{Z})$. 
Consider a sheaf 
${\cal A}:={\cal O}_{\mathfrak E}\oplus
({\cal O}_{\mathfrak E}
\otimes_{{\cal O}_{{\cal Z}^{\rm ex}}}\Om^1_{{\cal Z}^{\rm ex}/U})$ of 
${\cal O}_{\mathfrak E}$-modules. 
Define a multiplicative structure of ${\cal A}$ 
by the following equalities: 
$\om \cdot \om'=0$ and $\om^{[n]}=0$ for 
local sections $\om, \om' \in {\cal O}_{\mathfrak E}
\otimes_{{\cal O}_{{\cal Z}^{\rm ex}}}\Om^1_{{\cal Z}^{\rm ex}/U}$ 
and $n\geq 2$. 
Define also a PD-structure $\del$ on 
$\ol{\cal J}\oplus ({\cal O}_{\mathfrak E}
\otimes_{{\cal O}_{{\cal Z}^{\rm ex}}}\Om^1_{{\cal Z}^{\rm ex}/U})$ 
by following equality: 
$\del_n((a,\om))=(a^{[n]},a^{[n-1]}\om)$ 
($a\in \ol{\cal J}$, $\om \in  {\cal O}_{\mathfrak E}
\otimes_{{\cal O}_{{\cal Z}^{\rm ex}}}\Om^1_{{\cal Z}^{\rm ex}/U}$, $n\in {\mab N})$ 
([loc.~cit.]). 
Because $\os{\circ}{\mathfrak E}$ is the PD-envelope of the immersion 
$Z \os{\sus}{\lo} {\cal Z}^{\rm ex}$ over $(\os{\circ}{U},{\cal K},\eps)$, 
there exists a morphism 
$ \ul{\rm Spec}_{\os{\circ}{\mathfrak E}}({\cal A})
\lo \os{\circ}{\mathfrak E}$ by the universality of 
$\os{\circ}{\mathfrak E}$. A local section $a$ of ${\cal O}_{\mathfrak E}$ 
is mapped to $(a,da)\in {\cal A}$, 
where 
\begin{align*} 
d\col {\cal O}_{\mathfrak E} \lo 
{\cal O}_{\mathfrak E}
\otimes_{{\cal O}_{{\cal Z}^{\rm ex}}}\Om^1_{{\cal Z}^{\rm ex}/U} 
\tag{1.3.28.3}\label{ali:zaedmm}
\end{align*} 
is a derivation of ${\cal O}_{\mathfrak E}$. 
(Since the induced morphism 
${\cal O}_{\mathfrak E}\lo {\cal A}$ is a PD-morphism, we see that 
$d(a^{[n]})=a^{[n-1]}da$ for a local section $a$ of ${\cal J}$ as in [loc.~cit.].  
Hence we see that the derivation (\ref{ali:zaedmm}) is equal to 
the extended derivation in (1).) 
Consequently we have the following morphism of sheaves of abelian groups:
\begin{align*} 
d\log \col {\cal O}_{\mathfrak E}^*\owns u\lom d\log u=u^{-1} du\in 
{\cal O}_{\mathfrak E}
\otimes_{{\cal O}_{{\cal Z}^{\rm ex}}}\Om^1_{{\cal Z}^{\rm ex}/U}.
\tag{1.3.28.4}\label{ali:ezum}
\end{align*} 
\par 
Let $\eta \col {\mathfrak E}\lo {\cal Z}^{\rm ex}$ be the natural morphism.
Let $\al \col M_{{\cal Z}^{\rm ex}}\lo {\cal O}_{{\cal Z}^{\rm ex}}$ and 
$\bet \col M_{\mathfrak E}\lo {\cal O}_{\mathfrak E}$ 
be the structural morphisms. 
We have the following morphism of sheaves of abelian groups: 
\begin{align*} 
\eta^{-1}(M_{{\cal Z}^{\rm ex}})\owns m\lom d\log m\in 
{\cal O}_{\mathfrak E}
\otimes_{{\cal O}_{{\cal Z}^{\rm ex}}}\Om^1_{{\cal Z}^{\rm ex}/U}. 
\tag{1.3.28.5}\label{ali:ezam}
\end{align*} 
Because the restrictions of the morphisms (\ref{ali:ezum}) and (\ref{ali:ezam}) 
to $\eta^{-1}({\cal O}^*_{{\cal Z}^{\rm ex}})$ are the same, 
we have the following morphism 
\begin{align*} 
d\log \col {\cal O}^*_{\mathfrak E}\oplus_{\eta^{-1}({\cal O}^*_{{\cal Z}^{\rm ex}})}
\eta^{-1}(M_{{\cal Z}^{\rm ex}}) \lo {\cal O}_{\mathfrak E}
\otimes_{{\cal O}_{{\cal Z}^{\rm ex}}}\Om^1_{{\cal Z}^{\rm ex}/U}. 
\tag{1.3.28.6}\label{ali:ezamm}
\end{align*} 
Obviously hold the relations   
$d(\al(m))=\al(m)\otimes d\log m$ in 
${\cal O}_{\mathfrak E}
\otimes_{{\cal O}_{{\cal Z}^{\rm ex}}}\Om^1_{{\cal Z}^{\rm ex}/U}$ 
for a local section 
$m$ of $\eta^{-1}(M_{{\cal Z}^{\rm ex}})$ and 
$ud\log u=du$ for a local section 
$u$ of ${\cal O}_{\mathfrak E}^*$. 
Hence holds the relation 
\begin{align*} 
d(\bet(m))=\bet(m)\otimes d\log m
\tag{1.3.28.7}\label{ali:zamm}
\end{align*}  
in 
${\cal O}_{\mathfrak E}
\otimes_{{\cal O}_{{\cal Z}^{\rm ex}}}\Om^1_{{\cal Z}^{\rm ex}/U}$ 
for a local section 
$m$ of 
$M_{\mathfrak E}={\cal O}^*_{\mathfrak E}\oplus_{\eta^{-1}({\cal O}^*_{{\cal Z}^{\rm ex}})}
\eta^{-1}(M_{{\cal Z}^{\rm ex}})$. 

\par 
Because the sheaf 
$\Om^1_{{\mathfrak E}/U}$ is universal for the morphism 
of sheaves of abelian groups from 
${\cal O}_{\mathfrak E}\otimes_{\mab Z}M^{\rm gp}_{\mathfrak E}$  
and the derivation of ${\cal O}_{\mathfrak E}$ 
((\ref{ali:zaedmm})) with the relation (\ref{ali:zamm}),  
we have a natural morphism 
$\Om^1_{{\mathfrak E}/U}{\lo} 
{\cal O}_{\mathfrak E}
\otimes_{{\cal O}_{{\cal Z}^{\rm ex}}}\Om^1_{{\cal Z}^{\rm ex}/U}$.   
By virtue of (1), this morphism induces a morphism 
$\Om^{\bul}_{{\mathfrak E}/U,[~]}{\lo} 
{\cal O}_{\mathfrak E}
\otimes_{{\cal O}_{{\cal Z}^{\rm ex}}}\Om^{\bul}_{{\cal Z}^{\rm ex}/U}$ of dga's 
over $g^{-1}({\cal O}_U)$.  
\par 
On the other hand, the natural composite morphism 
${\cal O}_{\mathfrak E}
\otimes_{{\cal O}_{{\cal Z}^{\rm ex}}}\Om^1_{{\cal Z}^{\rm ex}/U}
\lo \Om^1_{{\mathfrak E}/U} \lo \Om^1_{{\mathfrak E}/U,[~]}$ 
induced by the natural morphism 
${\mathfrak E}\lo {\cal Z}^{\rm ex}$ induce a morphism 
${\cal O}_{\mathfrak E}
\otimes_{{\cal O}_{{\cal Z}^{\rm ex}}}\Om^{\bul}_{{\cal Z}^{\rm ex}/U}
\lo \Om^{\bul}_{{\mathfrak E}/U,[~]}$ of dga's over $g^{-1}({\cal O}_U)$. 
This morphism is the inverse of the morphism 
${\cal O}_{\mathfrak E}
\otimes_{{\cal O}_{{\cal Z}^{\rm ex}}}\Om^{\bul}_{{\cal Z}^{\rm ex}/U}
\lo \Om^{\bul}_{{\mathfrak E}/U,[~]}$. 
\par 
(3): This has been proved in \cite[0 (3.1.4), (3.1.6)]{idw}.
\end{proof}

\parno 
The filtration $P$ on $\Om^{\bul}_{{\mathfrak E}/\os{\circ}{U}}$ 
induces a filtration $P$ on $\Om^{\bul}_{{\mathfrak E}/\os{\circ}{U},[~]}$.

\begin{coro}\label{coro:ecdf}
There exists the following isomorphism of 
filtered dga's over $g^{-1}({\cal O}_{U}):$ 
\begin{equation*}
(\Om^{\bul}_{{\mathfrak E}/\os{\circ}{U},[~]},P)
\os{\sim}{\lo} 
({\cal O}_{\mathfrak E}
\otimes_{{\cal O}_{{\cal Z}^{\rm ex}}}
\Om^{\bul}_{{\cal Z}{}^{\rm ex}/\os{\circ}{U}},P).   
\tag{1.3.29.1}\label{eqn:fwuou}
\end{equation*}  
\end{coro}
\begin{proof} 
(\ref{eqn:fwuou}) immediately follows from (\ref{eqn:fwniwu}) and (\ref{eqn:fwu}). 
\end{proof}

\par 
We can use the following lemma (\ref{lemm:es0}) 
and the following corollary (\ref{coro:edf}) in the proofs of (\ref{theo:funas}) 
and (\ref{theo:funpas}) below (see the explanation for (3) in the Introduction).    

\begin{lemm}\label{lemm:es0}
Let $(U,{\cal K},\eps)\lo (U',{\cal K}',\eps')$ 
be a morphism of fine log PD-schemes. 
Set $U_0:=\ul{\rm Spec}^{\log}_U({\cal O}_U/{\cal K})$ 
and $U'_0:=
\ul{\rm Spec}^{\log}_{U'}({\cal O}_{U'}/{\cal K}')$. 
Let $Z$ $($resp.~$W)$ be a fine log scheme over $U_0$ 
$($resp.~a fine log scheme over $U'_0)$. 
Let $Z\os{\sus}{\lo} {\cal Z}$ $($resp.~$W\os{\sus}{\lo}{\cal W})$  
be an immersion into a fine log $($formal$)$ scheme over $U$ 
$($resp.~a fine log $($formal$)$ scheme over $U')$. 
Let ${\mathfrak E}$ $($resp.~${\mathfrak F})$ 
be the log PD-envelope of the immersion 
$Z\os{\sus}{\lo} {\cal Z}$ over $(U,{\cal K},\eps)$ 
$($resp.~$W\os{\sus}{\lo} {\cal W}$ over $(U',{\cal K}',\eps'))$. 
Then the following hold$:$ 
\par 
$(1)$ Let $\ol{g}\col {\mathfrak E}\lo {\mathfrak F}$ be a morphism of 
fine log PD-schemes over $U\lo U'$ extending a morphism 
$g\col Z\lo W$ over the morphism $U_0\lo U'_0$. 
Then $\ol{g}$ induces 
the following morphism of filtered complexes in 
${\rm C}^+{\rm F}(g^{-1}({\cal O}_{U'}))\!:$
\begin{equation*}
\ol{g}{}^*\col ({\cal O}_{\mathfrak F}
\otimes_{{\cal O}_{{\cal W}^{\rm ex}}}
\Om^{\bul}_{{\cal W}^{\rm ex}/\os{\circ}{U}{}'},P) \lo 
g_*(({\cal O}_{\mathfrak E}
\otimes_{{\cal O}_{{\cal Z}^{\rm ex}}}
\Om^{\bul}_{{\cal Z}^{\rm ex}/\os{\circ}{U}},P)).   
\tag{1.3.30.1}\label{eqn:fwpwu}
\end{equation*}   
\par 
$(2)$ Let the notations be as in $(1)$. 
The composite morphism in 
${\rm D}^+(g^{-1}({\cal O}_{U'}))$ 
\begin{equation*}
{\cal O}_{\mathfrak F}
\otimes_{{\cal O}_{{\cal W}^{\rm ex}}}
\Om^{\bul}_{{\cal W}^{\rm ex}/U'}
\os{\ol{g}{}^*}{\lo}  
g_*({\cal O}_{\mathfrak E}
\otimes_{{\cal O}_{{\cal Z}^{\rm ex}}}
\Om^{\bul}_{{\cal Z}^{\rm ex}/U})\lo 
Rg_*({\cal O}_{\mathfrak E}
\otimes_{{\cal O}_{{\cal Z}^{\rm ex}}}
\Om^{\bul}_{{\cal Z}^{\rm ex}/U})   
\tag{1.3.30.2}\label{eqn:frwu}
\end{equation*}  
is equal to the induced morphism 
${\cal O}_{W/U'}\lo Rg_{{\rm crys}*}({\cal O}_{Z/U})$ 
in ${\rm D}^+({\cal O}_{W/U'})$ 
by the morphism 
${\cal O}_{W/U'}\lo g_{{\rm crys}*}({\cal O}_{Z/U})$ 
in ${\rm C}^+({\cal O}_{W/U'})$ via the log Poincar\'{e} lemma. 
\end{lemm} 
\begin{proof}
\par 
(1): The morphism ${\mathfrak E}\lo {\mathfrak F}$ induces the pull-back morphism 
\begin{align*}
\ol{g}{}^*\col \Om^{\bul}_{{\mathfrak F}/\os{\circ}{U}{}',[~]}
\lo 
g_*(\Om^{\bul}_{{\mathfrak E}/\os{\circ}{U},[~]}).   
\tag{1.3.30.3}\label{eqn:fwewu}
\end{align*} 
Since the following diagram
\begin{equation*} 
\begin{CD} 
{\mathfrak E}@>>> {\mathfrak F}\\
@VVV @VVV \\
\os{\circ}{\mathfrak E}@>>> \os{\circ}{\mathfrak F}
\end{CD} 
\end{equation*} 
is commutative, the morphism (\ref{eqn:fwewu}) gives the following morphism: 
\begin{equation*}
\ol{g}{}^*\col (\Om^{\bul}_{{\mathfrak F}/\os{\circ}{U}{}',[~]},P) \lo 
g_*(\Om^{\bul}_{{\mathfrak E}/\os{\circ}{U},[~]},P).    
\end{equation*}  
By (\ref{eqn:fwuou}) we have the morphism (\ref{eqn:fwpwu}). 
\par 
(2): We have the following commutative diagram 
\begin{equation*} 
\begin{CD} 
{\cal O}_{\mathfrak F}@>>> g_*({\cal O}_{\mathfrak E})\\
@VVV @VVV \\
{\cal O}_{\mathfrak F}
\otimes_{{\cal O}_{{\cal W}^{\rm ex}}}
\Om^{\bul}_{{\cal W}^{\rm ex}/U'} @>>>  
Rg_*({\cal O}_{\mathfrak E}
\otimes_{{\cal O}_{{\cal Z}^{\rm ex}}}\Om^{\bul}_{{\cal Z}^{\rm ex}/U}).     
\end{CD} 
\tag{1.3.30.4}\label{cd:fwdu}
\end{equation*} 
This commutative diagram is nothing but the following commutative diagram 
in ${\rm D}^+(g^{-1}({\cal O}_{U'}))$
by the log Poincar\'{e} lemma: 
\begin{equation*} 
\begin{CD} 
({\cal O}_{W/U'})_{\mathfrak F}
@>>> g_*(({\cal O}_{Z/U})_{\mathfrak E})\\
@VVV @VVV \\
Ru_{W/U'*}({\cal O}_{W/U'})@>>> Rg_*Ru_{Z/U*}({\cal O}_{Z/U})
=Ru_{W/U'*}Rg_{{\rm crys}*}({\cal O}_{Z/U}).     
\end{CD} 
\tag{1.3.30.5}\label{cd:fwrdu}
\end{equation*} 
\end{proof}

\begin{lemm}\label{lemm:obvc} 
Let the notations  be as in {\rm (\ref{lemm:es0})}. 
Assume that ${\cal W}$ is log smooth over $U'$ and 
that $\os{\circ}{Z}$ and $\os{\circ}{\cal Z}$ are affine. 
Let $g\col Z\lo W$ be a morphism over $U_0\lo U'_0$. 
Let $\star$ be nothing or ${\rm ex}$. 
Then there exists a morphism 
$\ol{g}\col {\mathfrak E}\lo {\cal W}^{\star}$ over $U\lo U'$ fitting into 
the following commutative diagram 
\begin{equation*}
\begin{CD} 
Z @>{\sus}>> {\mathfrak E} \\
@V{g}VV @VV{\ol{g}}V \\
W@>{\sus}>> {\cal W}^{\star}.  
\end{CD} 
\tag{1.3.31.1}\label{cd:xdqx} 
\end{equation*}   
$($The morphism $\ol{g}$ induces a morphism 
${\mathfrak E}\lo {\mathfrak F}$, 
which we denote by $\ol{g}$ again by abuse of notation.$)$
\end{lemm} 
\begin{proof} 
Since the immersion $Z\os{\sus}{\lo} {\mathfrak E}$ is nil and 
the morphism ${\cal W}\lo U'$ is formally log smooth,  
we have the desired morphism $\ol{g}\col {\mathfrak E}\lo {\cal W}$ 
by (\ref{prop:niip}). Since the morphism $Z\os{\sus}{\lo} {\mathfrak E}$ 
is exact, we have the morphism $\ol{g}\col {\mathfrak E}\lo {\cal W}^{\rm ex}$ 
by the universality  of the exactification.   
\end{proof} 

\begin{rema}\label{rema:uuc}
Note that we do not necessarily have the following (usual) commutative diagram
\begin{equation*} 
\begin{CD} 
Z@>{\subset}>>{\cal Z}^{\star} \\ 
@VVV @VVV \\ 
W@>{\subset}>>{\cal W}^{\star} 
\end{CD} 
\tag{1.3.32.1}\label{cd:nnc} 
\end{equation*} 
($\star$=nothing or ex) 
because we do not assume that the immersion 
$Z\os{\sus}{\lo} {\cal Z}$ is not nil in (\ref{lemm:es0}) and (\ref{lemm:obvc}). 
As a result, we do not have the morphism $\ol{g}$ in (\ref{cd:xdqx}) a priori. 
However, if we replace ${\cal Z}$ by ${\cal Z}\times_U({\cal W}\times_{U'}U)$, 
we have the following commutative diagram 
\begin{equation*} 
\begin{CD} 
Z@>{\subset}>>{\cal Z}\times_U({\cal W}\times_{U'}U) \\ 
@VVV @VVV \\ 
W@>{\subset}>>{\cal W},  
\end{CD} 
\tag{1.3.32.2}\label{cd:icmm} 
\end{equation*} 
where the right vertical morphism is the following composite morphism 
${\cal Z}\times_U({\cal W}\times_{U'}U)\os{{\rm 2nd. proj.}}{\lo} {\cal W}\times_{U'}U
\os{{\rm 1st. proj.}}{\lo} {\cal W}$. 
This simple observation is very useful for the proof of 
the functoriality of canonical (iso)morphisms 
(e.~g., (\ref{theo:funas}),  (\ref{theo:ccrw}) (2) below).  
I cannot find a reference in which 
the product ${\cal Z}\times_U({\cal W}\times_{U'}U)$ 
and the commutative diagram (\ref{cd:icmm})  
are considered in the proof of the functoriality of a canonical isomorphism 
except the proof of \cite[(1.6)]{bfi}. 
\end{rema} 

\begin{coro}\label{coro:edf}
Assume furthermore that $(U,{\cal K},\eps)=(\os{\circ}{T},{\cal J},\del)$ 
and $(U',{\cal K}',\eps')=(\os{\circ}{T}{}',{\cal J}',\del')$. 
Let $Z$ $($resp.~$W)$ be a fine log affine scheme over $S_{\os{\circ}{T}_0}$ 
$($resp.~a fine log scheme over $S'_{\os{\circ}{T}{}'_0})$. 
Let $Z\os{\sus}{\lo} \ol{\cal Z}$ $($resp.~$W\os{\sus}{\lo}\ol{\cal W})$  
be an immersion into a log $($formal$)$ affine scheme over $\ol{S(T)^{\nat}}$ 
$($resp.~a log smooth $($formal$)$ scheme over $\ol{S'(T')^{\nat}})$. 
Let $\ol{\mathfrak E}$ $($resp.~$\ol{\mathfrak F})$ 
be the log PD-envelope of the immersion 
$Z\os{\sus}{\lo} \ol{\cal Z}$ 
$($resp.~$W\os{\sus}{\lo} \ol{\cal W})$ over 
$(\os{\circ}{T},{\cal J},\del)$ 
$($resp.~$(\os{\circ}{T}{}',{\cal J}',\del'))$. 
Set ${\cal Z}:=\ol{\cal Z}\times_{\ol{S(T)^{\nat}}}S(T)^{\nat}$, 
${\cal W}:=\ol{\cal W}\times_{\ol{S'(T')^{\nat}}}S'(T')^{\nat}$, 
${\mathfrak E}:=\ol{\mathfrak E}\times_{{\mathfrak D}(\ol{S(T)^{\nat}})}S(T)^{\nat}$ 
and 
${\mathfrak F}:=\ol{\mathfrak F}\times_{{\mathfrak D}(\ol{S'(T')^{\nat}})}S'(T')^{\nat}$.  

Then the following hold$:$ 
\par 
$(1)$ 
There exists a morphism 
$\ol{g}\col \ol{\mathfrak E}\lo \ol{\mathfrak F}$ of 
fine log schemes over 
$\ol{S(T)^{\nat}}\lo \ol{S'(T')^{\nat}}$ extending a morphism 
$g\col Z\lo W$ over the morphism $S_{T_0}\lo S'_{T'_0}$. 
\par 
$(2)$ 
The morphism 
\begin{equation*} 
g^*\col (\Om^{\bul}_{W/\os{\circ}{T}{}'_0},P)\lo 
g_*((\Om^{\bul}_{Z/\os{\circ}{T}_0},P))
\end{equation*}
lifts to a morphism 
\begin{equation*} 
\ol{g}{}^*\col 
({\cal O}_{\ol{\mathfrak F}}
\otimes_{{\cal O}_{\ol{\cal W}{}^{\rm ex}}}
\Om^{\bul}_{\ol{\cal W}{}^{\rm ex}/\os{\circ}{T}{}'},P)
\lo 
\ol{g}_*(({\cal O}_{\ol{\mathfrak E}}
\otimes_{{\cal O}_{\ol{\cal Z}{}^{\rm ex}}}
\Om^{\bul}_{\ol{\cal Z}{}^{\rm ex}/\os{\circ}{T}},P)).  
\tag{1.3.33.1}\label{eqn:ofw} 
\end{equation*} 
\par 
$(3)$ 
The morphism 
\begin{equation*} 
g^*\col (\Om^{\bul}_{W/\os{\circ}{T}{}'_0},P)
\lo 
g_*((\Om^{\bul}_{Z/\os{\circ}{T}_0},P)) 
\end{equation*}
lifts to a morphism 
\begin{equation*} 
g^{{\rm PD}*}\col 
({\cal O}_{\mathfrak F}
\otimes_{{\cal O}_{{\cal W}^{\rm ex}}}
\Om^{\bul}_{{\cal W}^{\rm ex}/\os{\circ}{T}{}'_0},P)
\lo 
g_*(({\cal O}_{\mathfrak E}
\otimes_{{\cal O}_{{\cal Z}^{\rm ex}}}
\Om^{\bul}_{{\cal Z}^{\rm ex}/\os{\circ}{T}_0},P)).
\tag{1.3.33.2}\label{eqn:ofpw} 
\end{equation*}
\end{coro}
\begin{proof} 
(1): (1) follows from (\ref{lemm:obvc}).   
\par 
(2): (2) follows from (\ref{lemm:es0}) (1). 
\par 
(3): (3) immediately follows from (2). 
\end{proof}


\begin{rema}\label{rema:tcwew}  
Note that the log formal scheme 
${\cal W}^{\rm ex}$ is not necessarily log smooth 
over $\os{\circ}{T}{}'$. 
\end{rema}

\par 
We finish this section by considering a section of 
$P_1({\cal O}_{\mathfrak D}
\otimes_{{\cal O}_{{\cal P}^{\rm ex}}}
{\Om}^1_{{\cal P}^{\rm ex}/\os{\circ}{T}})$ 
for the next section. 
\par 
Let $\tau$ be a local section of $M_{S(T)^{\nat}}$ such 
that the image of $\tau$ in $M_{S(T)^{\nat}}/{\cal O}_T^*$ is the local generator. 
Then we have a section 
$d\log \tau \in {\Om}^1_{S(T)^{\nat}/\os{\circ}{T}}$. 
Because $d\log u=0$ 
in ${\Om}^1_{S(T)^{\nat}/\os{\circ}{T}}$ 
$(u\in {\cal O}_T^*)$,
$d\log \tau$ is independent of the choice of 
a local section of $M_{S(T)^{\nat}}$
whose image in $M_{S(T)^{\nat}}/{\cal O}_T^*$ 
is the local generator. 
By this observation, 
$\theta_{S(T)^{\nat}}:=d\log \tau \in {\Om}^1_{S(T)^{\nat}/\os{\circ}{T}}$ is 
globalized and it is a well-defined global section of 
${\Om}^1_{S(T)^{\nat}/\os{\circ}{T}}$
for a (general) family $S$ of log points.  
Let 
\begin{equation*} 
\begin{CD} 
Z @>{\sus}>> {\cal Q} \\ 
@VVV @VVV \\ 
S_{\os{\circ}{T}_0} @>{\sus}>> S(T)^{\nat}  
\end{CD} 
\end{equation*} 
be a commutative diagram of fine log schemes, 
where the upper horizontal morphism is an immersion.  
Let ${\mathfrak E}$ be the log PD-envelope of 
the immersion $Z\os{\sus}{\lo}{\cal Q}$ over $(\os{\circ}{T},{\cal J},\del)$. 
Let $\theta_{\cal Q}$ be the image of $\theta_{S(T)^{\nat}}$ in 
$P_1({\cal O}_{\mathfrak E}\otimes_{{\cal O}_{\cal Q}}
{\Om}^1_{{\cal Q}/\os{\circ}{T}})$. 
The section $\theta_{\cal Q}$ 
induces the following morphism 
\begin{equation*} 
\theta_{\cal Q} \wedge \col 
({\cal O}_{\mathfrak E}\otimes_{{\cal O}_{\cal Q}}
{\Om}^i_{{\cal Q}/\os{\circ}{S}})/
P_k({\cal O}_{\mathfrak E}\otimes_{{\cal O}_{\cal Q}}
{\Om}^i_{{\cal Q}/\os{\circ}{T}})
\lo 
({\cal O}_{\mathfrak E}
\otimes_{{\cal O}_{\cal Q}}
{\Om}^{i+1}_{{\cal Q}/\os{\circ}{S}})/
P_{k+1}({\cal O}_{\mathfrak E}\otimes_{{\cal O}_{\cal Q}}
{\Om}^{i+1}_{{\cal Q}/\os{\circ}{T}})  
\quad (i,k\in {\mab Z}). 
\tag{1.3.34.1}\label{eqn:tdpf}
\end{equation*}  
More generally, 
let $E$ be as after (\ref{eqn:nidospltd}).  
Then we obtain the following morphisms 
\begin{align*} 
\theta_{{\cal P}^{\rm ex}} \wedge \col 
({\cal E}\otimes_{{\cal O}_{{\cal P}^{\rm ex}}}
{\Om}^i_{{\cal P}^{\rm ex}/\os{\circ}{S}})/
P_k({\cal E}\otimes_{{\cal O}_{{\cal P}^{\rm ex}}}
{\Om}^i_{{\cal P}^{\rm ex}/\os{\circ}{T}})\lo &
({\cal E}
\otimes_{{\cal O}_{{\cal P}^{\rm ex}}}
{\Om}^{i+1}_{{\cal P}^{\rm ex}/\os{\circ}{S}})/
P_{k+1}({\cal E}\otimes_{{\cal O}_{{\cal P}^{\rm ex}}}
{\Om}^{i+1}_{{\cal P}^{\rm ex}/\os{\circ}{T}})  \tag{1.3.34.2}\label{eqn:tdepf}\\
&\quad (i,k\in {\mab Z}). 
\end{align*}

\section{Zariskian $p$-adic filtered Steenbrink complexes}\label{sec:psc}
Let $S$, $(T,{\cal J},\del)$, $T_0\lo S$, $S_{\os{\circ}{T}_0}$, $S(T)$ and 
$S(T)^{\nat}$  
be as in \S\ref{sec:ldc}. 
Let $X_{\bul \leq N}/S$ be an $N$-truncated simplicial SNCL scheme 
((\ref{defi:ntrxs})). 
Set $X_{\bul \leq N,T_0}:=X_{\bul \leq N}\times_{S}T_0$ 
and $X_{\bul \leq N,\os{\circ}{T}_0}:=X_{\bul \leq N}\times_{S}S_{\os{\circ}{T}_0}
=X_{\bul \leq N}\times_{\os{\circ}{S}}{\os{\circ}{T}_0}$.   
Let 
$f\col X_{\bul \leq N,\os{\circ}{T}_0} \lo S_{\os{\circ}{T}_0}$ 
be the structural morphism. 
By abuse of notation, let us also denote the structural morphism 
$X_{\bul \leq N,\os{\circ}{T}_0}\lo S(T)^{\nat}$ 
by $f$. 
We denote the structural morphism 
$X_{m,\os{\circ}{T}_0} \lo S(T)^{\nat}$ $(0\leq m \leq N)$ by $f_{m}$. 
We would not like to use the notation  
$f_{\bul \leq N}$ for the structural morphism 
$X_{\bul \leq N,\os{\circ}{T}_0} \lo S(T)^{\nat}$ 
because we would like to use this notation for the structural morphism 
$X_{\bul \leq N,\os{\circ}{T}_0} \lo S(T)^{\nat}_{\bul \leq N}$, 
where $S(T)^{\nat}_{\bul \leq N}$ 
is the $N$-truncated simplicial constant log schemes 
obtained by $S(T)^{\nat}$. 
Let $f_T\col X_{\bul \leq N,T_0} \lo T$  
be the structural morphism. 
Then $\os{\circ}{f}=\os{\circ}{f}_{T}$. 
We denote the structural morphism 
$X_{m,T_0} \lo T$ $(0\leq m \leq N)$ by $f_{m,T}$. 
Assume that 
$X_{\bul \leq N,\os{\circ}{T}_0}$ has 
an affine $N$-truncated simplicial open covering 
$X'_{\bul \leq N,\os{\circ}{T}_0}$ of 
$X_{\bul \leq N,\os{\circ}{T}_0}$. 
(We do not assume that $X_{\bul \leq N}$ has 
an affine $N$-truncated simplicial open covering 
of $X_{\bul \leq N}$.) 
Let $E^{\bul \leq N}$ be a flat quasi-coherent crystal of 
${\cal O}_{\os{\circ}{X}_{\bul \leq N,T_0}/\os{\circ}{T}}$-modules.  
\par 
The aim in this section is to construct a filtered complex
\begin{equation*} 
(A_{\rm zar}(X_{\bul \leq N,\os{\circ}{T}_0}/S(T)^{\nat},E^{\bul \leq N}),P)
\in {\rm D}^+{\rm F}(
f^{-1}({\cal O}_T)), 
\end{equation*}  
which we call 
the {\it zariskian $p$-adic filtered  Steenbrink complex} 
of $E^{\bul \leq N}$ for $X_{\bul \leq N,\os{\circ}{T}_0}/S(T)^{\nat}$. 
\par 
For the moment, consider the case $N=0$. 
Set $X:=X_0$ and $E:=E^0$.  
Assume that there exist an immersion 
$X_{\os{\circ}{T}_0}\os{\sus}{\lo} \ol{\cal P}$ over $S(T)^{\nat}$ in 
\S\ref{sec:ldc} and the cartesian diagram (\ref{cd:wgtx}). 
Let $\ol{\mathfrak D}$, ${\mathfrak D}$, 
$(\ol{\cal E},\ol{\nabla})$, 
${\cal P}^{\rm ex}$ and 
$({\cal E},{\nabla})$ 
be as in \S\ref{sec:ldc}.

\par 
Set $\theta:=\theta_{{\cal P}^{\rm ex}}$. 
Set  
\begin{align*} 
A_{\rm zar}({\cal P}^{\rm ex}/S(T)^{\nat},{\cal E})^{ij}
& :={\cal E}\otimes_{{\cal O}_{{\cal P}^{\rm ex}}}
{\Om}^{i+j+1}_{{\cal P}^{\rm ex}/\os{\circ}{T}}/P_j 
\tag{1.4.0.1}\label{cd:accef} \\
& :={\cal E}\otimes_{{\cal O}_{{\cal P}^{\rm ex}}}
{\Om}^{i+j+1}_{{\cal P}^{\rm ex}/\os{\circ}{T}}/
P_j({\cal E}\otimes_{{\cal O}_{{\cal P}^{\rm ex}}}
{\Om}^{i+j+1}_{{\cal P}^{\rm ex}/\os{\circ}{T}})  
\quad (i,j \in {\mab N}). 
\end{align*}   
The sheaf 
$A_{\rm zar}({\cal P}^{\rm ex}/S(T)^{\nat},{\cal E})^{ij}$ 
has a quotient filtration $P$ obtained by the filtration $P$ on 
${\cal E}\otimes_{{\cal O}_{{\cal P}^{\rm ex}}}
{\Om}^{i+j+1}_{{\cal P}^{\rm ex}/\os{\circ}{T}}$. 
We consider the following boundary morphisms of 
double complexes: 
\begin{equation*}
\begin{CD}
A_{\rm zar}({\cal P}^{\rm ex}/S(T)^{\nat}
,{\cal E})^{i,j+1}  @.  \\ 
@A{\theta \wedge}AA  @. \\
A_{\rm zar}({\cal P}^{\rm ex}/S(T)^{\nat},{\cal E})^{ij}
@>{-\nabla}>> 
A_{\rm zar}({\cal P}^{\rm ex}/S(T)^{\nat}
,{\cal E})^{i+1,j}.\\
\end{CD}
\tag{1.4.0.2}\label{cd:lccbd} 
\end{equation*}  
(We think that these are the best boundary morphisms with respect to 
the signs.)
Then we have the double complex 
$A_{\rm zar}({\cal P}^{\rm ex}/S(T)^{\nat},{\cal E})^{\bul \bul}$. 
Indeed, for a section $e\otimes \om \in {\cal E}\otimes_{{\cal O}_{{\cal P}^{\rm ex}}}
{\Om}^{i+j+1}_{{\cal P}^{\rm ex}/\os{\circ}{T}}$, we have the following equalities 
\begin{align*} 
(\nabla \theta +\theta \nabla)(e\otimes \om)&=
\{\nabla(e)\wedge \theta \wedge \om+e\otimes d(\theta \wedge \om)\}
+\{\theta \wedge \nabla(e) \wedge \om+e\otimes (\theta \wedge d\om)\}
\tag{1.4.0.3}\label{ali:eom} \\
&=\nabla(e)\wedge \theta \wedge \om+\theta \wedge \nabla(e) \wedge \om =0
\end{align*} 
since $\theta$ is a closed 1-form. 
The complex 
$A_{\rm zar}({\cal P}^{\rm ex}/S(T)^{\nat}
,{\cal E})^{\bul \bul}$ 
has a filtration $P=\{P_k\}_{k \in {\mab Z}}$ 
defined by the following formula: 
\begin{equation*} 
P_kA_{\rm zar}({\cal P}^{\rm ex}/S(T)^{\nat},{\cal E})^{\bul \bul}
:=(\cdots 
(P_{2j+k+1}+P_j)A_{\rm zar}({\cal P}^{\rm ex}/S(T)^{\nat}
,{\cal E})^{ij}
\cdots)\in {\rm C}^+{\rm F}(f^{-1}({\cal O}_T)).    
\tag{1.4.0.4}\label{eqn:lpcad}
\end{equation*} 
Let $(A_{\rm zar}({\cal P}^{\rm ex}/S(T)^{\nat},{\cal E}),P)$ 
be the filtered single complex of the filtered double complex 
$(A_{\rm zar}({\cal P}^{\rm ex}/S(T)^{\nat},{\cal E})^{\bul \bul},P)$.  

\par 
Let $Y$ be an SNCL scheme over $S$ and 
assume that there exists an immersion 
$Y_{\os{\circ}{T}_0} \os{\subset}{\lo} \ol{\cal Q}$  
into a log smooth scheme over $\ol{S(T)^{\nat}}$. 
Assume that there exist morphisms 
$g\col X_{\os{\circ}{T}_0} \lo Y_{\os{\circ}{T}_0}$ and 
$\ol{g} \col \ol{\cal P} \lo \ol{\cal Q}$ $(i=1,2)$ 
making the following diagram commutative: 
\begin{equation*} 
\begin{CD} 
X_{\os{\circ}{T}_0}@>{\sus}>> \ol{\cal P} \\ 
@V{g}VV @VV{\ol{g}}V \\ 
Y_{\os{\circ}{T}_0}@>{\sus}>> \ol{\cal Q}.  
\end{CD}
\end{equation*} 
Set ${\cal Q}:=\ol{\cal Q}\times_{\ol{S(T)^{\nat}}}S(T)^{\nat}$ and 
let ${\mathfrak E}$ be the log PD-envelope of 
the immersion $Y_{\os{\circ}{T}_0} \os{\subset}{\lo} {\cal Q}$ 
over $(S(T)^{\nat},{\cal J},\del)$. 
Let $\wt{g}\col {\cal P}\lo{\cal Q}$ be 
the base change morphism of $\ol{g}$ 
by the morphism 
$S(T)^{\nat}\os{\sus}{\lo} \ol{S(T)^{\nat}}$. 
Let 
$g^{\rm PD}\col 
{\mathfrak D}\lo {\mathfrak E}$ 
be the natural morphism induced by $\wt{g}$. 
Let 
$\os{\circ}{g}_{\rm crys} \col 
((\os{\circ}{X}_{T_0}/\os{\circ}{T})_{\rm crys},
{\cal O}_{\os{\circ}{X}_{T_0}/\os{\circ}{T}})\lo 
((\os{\circ}{Y}_{T_0}/\os{\circ}{T})_{\rm crys},
{\cal O}_{\os{\circ}{Y}_{T_0}/\os{\circ}{T}})$ 
be the induced morphism of ringed topoi by 
$\os{\circ}{g} \col \os{\circ}{X}_{T_0}\lo \os{\circ}{Y}_{T_0}$. 
Let $E$ (resp.~$F$) be a flat quasi-coherent crystal of  
${\cal O}_{\os{\circ}{X}_{T_0}/\os{\circ}{T}}$-modules  
(resp.~a flat quasi-coherent crystal of 
${\cal O}_{\os{\circ}{Y}_{T_0}/\os{\circ}{T}}$-modules). 
Assume that we are given 
a morphism $F\lo \os{\circ}{g}_{{\rm crys}*}(E)$. 
Let $({\cal F},\nabla)$ be 
the ${\cal O}_{\mathfrak E}$-module 
with integrable connection obtained in (\ref{eqn:nidpltd}) for $F$. 
Then, by (\ref{prop:npwf}), 
we have the following morphism of filtered complexes: 
\begin{equation*}
({\cal F}\otimes_{{\cal O}_{{\cal Q}^{\rm ex}}}
{\Om}^{\bul}_{{\cal Q}^{\rm ex}/\os{\circ}{T}},P)
\lo  
g^{\rm PD}_*(({\cal E}\otimes_{{\cal O}_{{\cal P}^{\rm ex}}}
{\Om}^{\bul}_{{\cal P}^{\rm ex}/\os{\circ}{T}},P)). 
\tag{1.4.0.5}\label{eqn:keyxy}
\end{equation*}
\par 
Now we come back to the situation in the beginning of this section. 
\par 
Let $X'_{\bul \leq N,\os{\circ}{T}_0} 
\os{\sus}{\lo} \ol{\cal P}{}'_{\bul \leq N}$ 
be an immersion into a log smooth scheme over $\ol{S(T)^{\nat}}$ 
(this immersion exists by (\ref{cd:pwtp}) and (\ref{prop:ytlft}) if 
each log affine subscheme of 
$X'_m$ $(0\leq m \leq N)$ is sufficiently small). 
Set $X_{mn}:={\rm cosk}_0^{X_m}(X'_m)_n$ 
$(0\leq m\leq N, n\in{\mab N})$ 
and 
$\ol{\cal P}_{mn}
:={\rm cosk}_0^{\ol{S(T)^{\nat}}}(\ol{\cal P}{}'_m)_n$.
Then we have a natural immersion 
$X_{\bul \leq N,\bul,\os{\circ}{T}_0} \os{\sus}{\lo} \ol{\cal P}_{\bul \leq N,\bul}$ 
over $S_{\os{\circ}{T}_0}\os{\sus}{\lo} \ol{S(T)^{\nat}}$. 
Let $\ol{\mathfrak D}_{\bul \leq N,\bul}$ 
be the log PD-envelope of the immersion 
$X_{\bul \leq N,\bul,\os{\circ}{T}_0} \os{\sus}{\lo} 
\ol{\cal P}_{\bul \leq N,\bul}$ 
over $(\os{\circ}{T},{\cal J},\gam)$. 
Set ${\mathfrak D}_{\bul \leq N,\bul}
:=\ol{\mathfrak D}_{\bul \leq N,\bul}
\times_{{\mathfrak D}(\ol{S(T)^{\nat}})}S(T)^{\nat}$. 
Let $f_{\bul}\col 
X_{\bul \leq N,\bul,\os{\circ}{T}_0} \lo S(T)^{\nat}$ 
and $f_{\bul,T}\col X_{\bul \leq N,\bul,T_0} \lo T$ 
be the structural morphisms. Then 
$\os{\circ}{f}_{\bul} =\os{\circ}{f}_{\bul,T}$. 
We denote the structural morphism 
$X_{m,\bul, T_0} \lo T$ $(0\leq m \leq N)$ by $f_{m,\bul,T}$. 
(We would not like to use the notation  
$f_{\bul \leq N,\bul,T}$ for the structural morphism 
$X_{\bul \leq N,\bul,T_0} \lo T$ as before.)
Set ${\cal P}_{\bul \leq N,\bul}:=
\ol{\cal P}_{\bul \leq N,\bul}\times_{\ol{S(T)^{\nat}}}S(T)^{\nat}$.  
We have a natural immersion
$X_{\bul \leq N,\bul,\os{\circ}{T}_0} 
\os{\sus}{\lo} {\cal P}_{\bul \leq N,\bul}$ 
over $S_{\os{\circ}{T}_0}\os{\sus}{\lo} S(T)^{\nat}$. 
When $T$ is restrictively hollow with respective to the morphism $T_0\lo S$, set 
${\cal P}_{\bul \leq N,\bul,T}:=
{\cal P}_{\bul \leq N,\bul}\times_{S(T)}T$  
and we have a natural immersion
$X_{\bul \leq N,\bul,T_0} 
\os{\sus}{\lo} {\cal P}_{\bul \leq N,\bul}\times_ST$ 
over $T_0\os{\sus}{\lo} T$. 
Let $E^{\bul \leq N,\bul}$ be the flat quasi-coherent crystal 
of ${\cal O}_{\os{\circ}{X}_{\bul \leq N,\bul,T_0}
/\os{\circ}{T}}$-modules obtained by $E^{\bul \leq N}$. 
Let $(\ol{\cal E}{}^{mn},\ol{\nabla}{}^{mn})$ 
$(0\leq m\leq N)$ 
be the quasi-coherent ${\cal O}_{\ol{\mathfrak D}_{mn}}$-module 
with the integrable connection 
associated to $\eps^*_{X_{mn,T_0}/\os{\circ}{T}}(E^{mn})$: 
$\ol{\nabla}{}^{mn}\col \ol{\cal E}{}^{mn}\lo 
\ol{\cal E}{}^{mn}\otimes_{{\cal O}_{\ol{\cal P}_{mn}}}
\Om^1_{\ol{\cal P}_{mn}/\os{\circ}{T}}$. 
Set ${\cal E}^{mn}:={\cal O}_{{\mathfrak D}_{mn}}
\otimes_{{\cal O}_{\ol{\mathfrak D}_{mn}}}\ol{\cal E}{}^{mn}$ 
and let 
\begin{equation*} 
\nabla^{mn} \col {\cal E}^{mn} 
\lo {\cal E}^{mn}\otimes_{{\cal O}_{{\cal P}_{mn}}}
{\Om}^1_{{\cal P}_{mn}/\os{\circ}{T}} 
=
{\cal E}^{mn}\otimes_{{\cal O}_{{\cal P}^{\rm ex}_{mn}}}
{\Om}^1_{{\cal P}^{\rm ex}_{mn}/\os{\circ}{T}} 
\tag{1.4.0.6}\label{eqn:nidopltd}
\end{equation*}  
be the induced connection by $\ol{\nabla}{}^{mn}$. 
By (\ref{eqn:keyxy}) 
we have the following  
$(N,\infty)$-truncated bicosimplicial filtered complex 
\begin{equation*} 
({\cal E}^{\bul \leq N,\bul}
\otimes_{{\cal O}_{{\cal P}^{\rm ex}_{\bul \leq N,\bul}}}
{\Om}^{\bul}_{{\cal P}^{\rm ex}_{\bul \leq N,\bul}
/\os{\circ}{T}},P).   
\tag{1.4.0.7}\label{cd:afil} 
\end{equation*} 
Set  
\begin{equation*} 
A_{\rm zar}({\cal P}^{\rm ex}_{\bul \leq N,\bul}/S(T)^{\nat},
{\cal E}^{\bul \leq N,\bul})^{ij} 
:=({\cal E}^{\bul \leq N,\bul}
\otimes_{{\cal O}_{{\cal P}^{\rm ex}_{\bul \leq N,\bul}}}
{\Om}^{i+j+1}_{{\cal P}^{\rm ex}_{\bul \leq N,\bul}/\os{\circ}{T}})/P_j 
\quad (i,j \in {\mab N}). 
\tag{1.4.0.8}\label{cd:adef} 
\end{equation*}   
We consider the following boundary morphisms of 
the following double complex: 
\begin{equation*}
\begin{CD}
A_{\rm zar}({\cal P}^{\rm ex}_{\bul \leq N,\bul}/S(T)^{\nat},
{\cal E}^{\bul \leq N,\bul})^{i,j+1} @.  \\ 
@A{\theta_{{\cal P}^{\rm ex}_{\bul \leq N,\bul,T}}\wedge}AA 
@. \\
A_{\rm zar}({\cal P}^{\rm ex}_{\bul \leq N,\bul}/S(T)^{\nat},
{\cal E}^{\bul \leq N,\bul})^{ij} 
@>{-\nabla}>> A_{\rm zar}({\cal P}^{\rm ex}_{\bul \leq N,\bul}/S(T)^{\nat},
{\cal E}^{\bul \leq N,\bul})^{i+1,j}.\\
\end{CD}
\tag{1.4.0.9}\label{cd:locstbd} 
\end{equation*}  
Then we have the $(N,\infty)$-truncated 
bicosimplicial double complex 
$A_{\rm zar}({\cal P}^{\rm ex}_{\bul \leq N,\bul}/S(T)^{\nat},
{\cal E}^{\bul \leq N,\bul})^{\bul \bul}$. 
The double complex 
$A_{\rm zar}({\cal P}^{\rm ex}_{\bul \leq N,\bul}/S(T)^{\nat},
{\cal E}^{\bul \leq N,\bul})^{\bul \bul}$ 
has a filtration $P=\{P_k\}_{k \in {\mab Z}}$ 
defined by the following formula: 
\begin{equation*} 
P_kA_{\rm zar}({\cal P}^{\rm ex}_{\bul \leq N,\bul}/S(T)^{\nat},
{\cal E}^{\bul \leq N,\bul})
:=(\cdots 
P_{2j+k+1}
A_{\rm zar}({\cal P}^{\rm ex}_{\bul \leq N,\bul}/S(T)^{\nat},
{\cal E}^{\bul \leq N,\bul})^{ij} 
\cdots).    
\tag{1.4.0.10}\label{eqn:dblad}
\end{equation*} 
Let $(A_{\rm zar}({\cal P}^{\rm ex}_{\bul \leq N,\bul}/S(T)^{\nat},
{\cal E}^{\bul \leq N,\bul}),P)$ 
be the single filtered complex of the filtered double complex 
$(A_{\rm zar}({\cal P}^{\rm ex}_{\bul \leq N,\bul}/S(T)^{\nat},
{\cal E}^{\bul \leq N,\bul})^{\bul \bul},P)$.  
\par 
Let 
\begin{equation*} 
a^{(k)}_{m} \col 
\os{\circ}{X}{}^{(k)}_{m,T_0} \lo 
\os{\circ}{X}_{m,T_0}  \quad (0\leq m\leq N, k\in {\mab N})
\tag{1.4.0.11}\label{eqn:atn}
\end{equation*} 
and 
\begin{equation*} 
a^{(k)}_{m\bul} \col 
\os{\circ}{X}{}^{(k)}_{m\bul,T_0} \lo 
\os{\circ}{X}_{m\bul,T_0}  \quad 
(0\leq m\leq N, k\in {\mab N})
\tag{1.4.0.12}\label{eqn:antb}
\end{equation*} 
be the natural morphisms of 
schemes and the natural morphism of simplicial schemes, 
respectively.   
Let 
\begin{equation*} 
a^{(k)}_{m,{\rm crys}} \col 
((\os{\circ}{X}{}^{(k)}_{m,T_0}/\os{\circ}{T})_{\rm crys},
{\cal O}_{\os{\circ}{X}{}^{(k)}_{m,T_0}/\os{\circ}{T}}) 
\lo ((\os{\circ}{X}_{m,T_0}/\os{\circ}{T})_{\rm crys},
{\cal O}_{\os{\circ}{X}_{m,T_0}/\os{\circ}{T}})  
\quad (0\leq m\leq N, k\in {\mab N})
\tag{1.4.0.13}\label{eqn:atnt}
\end{equation*} 
and 
\begin{equation*} 
a^{(k)}_{m\bul,{\rm crys}} \col 
((\os{\circ}{X}{}^{(k)}_{m\bul,T_0}/\os{\circ}{T})_{\rm crys},
{\cal O}_{\os{\circ}{X}{}^{(k)}_{m\bul,T_0}
/\os{\circ}{T}})  \lo 
((\os{\circ}{X}_{m\bul,T_0}/\os{\circ}{T})_{\rm crys},
{\cal O}_{\os{\circ}{X}_{m\bul,T_0}/\os{\circ}{T}}) 
\quad (0\leq m\leq N, k\in {\mab N})
\tag{1.4.0.14}\label{eqn:antbt}
\end{equation*} 
be the morphisms of ringed topoi 
obtained by (\ref{eqn:atn}) and (\ref{eqn:antb}), 
respectively. 
\par 
The following (1) is an SNCL version of \cite[(4.14)]{nh3}. 
The following is a key lemma for (\ref{theo:indcr}) below:  

\begin{lemm}\label{lemm:knit}
Let $m\leq N$ and $k$ be nonnegative integers. 
Then the following hold$:$ 
\par 
$(1)$ For the morphism 
$g \col {\cal P}^{\rm ex}_{mn}\lo {\cal P}^{\rm ex}_{mn'}$ 
corresponding to a morphism $[n']\lo [n]$ in $\Del$, 
the assumption in {\rm (\ref{prop:mmoo})} is satisfied  
for $Y={\cal P}^{\rm ex}_{mn}$ 
and $Y'={\cal P}^{\rm ex}_{mn'}$. 
Consequently there exists a morphism 
$\os{\circ}{g}{}^{(k)}
\col \os{\circ}{\cal P}{}^{{\rm ex},(k)}_{mn}\lo 
\os{\circ}{\cal P}{}^{{\rm ex},(k)}_{mn'}$ 
$(k\in {\mab N})$ fitting into 
the commutative diagram {\rm (\ref{cd:dmgdm})} 
for $Y={\cal P}^{\rm ex}_{mn}$ and 
$Y'={\cal P}^{\rm ex}_{mn'};$ 
the family 
$\{\os{\circ}{\cal P}{}^{{\rm ex},(k)}_{mn}\}_{n\in {\mab N}}$ 
gives a simplicial formal scheme 
$\os{\circ}{\cal P}{}^{{\rm ex},(k)}_{m\bul}$ 
with natural morphism 
$b^{(k)}_{m\bul} \col 
\os{\circ}{\cal P}{}^{{\rm ex},(k)}_{m\bul}
\lo \os{\circ}{\cal P}{}^{\rm ex}_{m\bul}$ 
of simplicial formal schemes. 
\par 
$(2)$ 
Let $\os{\circ}{\mathfrak D}{}^{(k)}_{m\bul}$ $(k\in {\mab N})$ 
be the PD-envelope of the simplicial immersion 
$\os{\circ}{X}{}^{(k)}_{m\bul,T_0}\os{\sus}{\lo} 
\os{\circ}{\cal P}{}^{{\rm ex},(k)}_{m\bul}$. 
Let $c^{(k)}_{m\bul} \col \os{\circ}{\mathfrak D}{}^{(k)}_{m\bul}
\lo \os{\circ}{\mathfrak D}{}_{m\bul}$ 
be the natural morphism. 
Then the following holds$:$
\begin{align*} 
{\rm gr}^P_k
({\cal E}^{m\bul}
\otimes_{{\cal O}_{{\cal P}^{\rm ex}_{m\bul}}}
{\Om}^{\bul}_{{\cal P}^{\rm ex}_{m\bul}/\os{\circ}{T}})= &
c^{(k-1)}_{m\bul*} 
(c^{(k-1)*}_{m\bul}({\cal E}^{m\bul})  \tag{1.4.1.1}\label{ali:grpmd}   \\
& \otimes_{{\cal O}_{{\cal P}{}^{{\rm ex},(k-1)}_{m\bul}}}
\Om^{\bul}_{\os{\circ}{\cal P}{}^{{\rm ex},(k-1)}_{m\bul}/\os{\circ}{T}}
\otimes_{\mab Z}\vp_{\rm zar}^{(k-1)}
(\os{\circ}{\cal P}{}^{\rm ex}_{m\bul}/\os{\circ}{T}))[-k]
\end{align*}
in 
${\rm C}^+(f^{-1}_{m\bul,T}({\cal O}_T))$. 
\par 
$(3)$ The following holds$:$
\begin{align*} 
& {\rm gr}^P_kA_{\rm zar}
({\cal P}^{\rm ex}_{m\bul}/S(T)^{\nat},{\cal E}^{m\bul})
\tag{1.4.1.2}\label{ali:grm2d} \\
& =\bigoplus_{j\geq \max \{-k,0\}}
c^{(2j+k)}_{m\bul*}
(c^{(2j+k)*}_{m\bul}({\cal E}^{m\bul})
\otimes_{{\cal O}_{\os{\circ}{\cal P}{}^{{\rm ex},(2j+k)}_{m\bul}}}
{\Om}^{\bul}_{\os{\circ}{\cal P}{}^{{\rm ex},(2j+k)}_{m\bul}/\os{\circ}{T}} \\
&\otimes_{\mab Z} \vp_{\rm zar}^{(2j+k)}
(\os{\circ}{\cal P}{}^{\rm ex}_{m\bul}/\os{\circ}{T}))[-2j-k]
\end{align*}
in ${\rm C}^+(f^{-1}_{m\bul}({\cal O}_T))$. 
Consequently  
\begin{align*} 
{\rm gr}^P_kA_{\rm zar}({\cal P}^{\rm ex}_{m\bul}/S(T)^{\nat},{\cal E}^{m\bul})  
=& \bigoplus_{j\geq \max \{-k,0\}} 
a^{(2j+k)}_{m\bul*} 
(Ru_{\os{\circ}{X}{}^{(2j+k)}_{m\bul,T_0}
/\os{\circ}{T}*}
(E^{m\bul}_{\os{\circ}{X}{}^{(2j+k)}_{m\bul,T_0}/\os{\circ}{T}}
\otimes_{\mab Z} \tag{1.4.1.3}\label{ali:grmpd} \\
& \vp_{\rm crys}^{(2j+k)}
(\os{\circ}{X}_{m\bul,T_0}/\os{\circ}{T})))[-2j-k]
\end{align*} 
in ${\rm D}^+(f^{-1}_{m\bul}({\cal O}_T))$.  
\end{lemm}
\begin{proof}
(1): 
Fix $0\leq m\leq N$.  
Set ${\cal M}_{mn}:=M_{{\cal P}^{\rm ex}_{mn}}$ 
$(n\in {\mab N})$ 
and $M_{mn}:=M_{X_{mn}}$. 
Let $e_n$ be a local section of a member of 
the local minimal generators of 
${\cal M}_{mn}/{\cal O}_{{\cal P}^{\rm ex}_{mn}}^*$.  
Then the image $\ol{e}_n$ of 
$e_n$ in $M_{mn}/{\cal O}^*_{X_{mn}}$ 
$(n\in {\mab N})$ 
is also a local section of a member of 
the local minimal generators since the immersion 
$X_{mn}\os{\sus}{\lo}{\cal P}^{\rm ex}_{mn}$ is exact. 
Because a standard degeneracy morphism  
$s_i\col X_{ml}\lo X_{m,l+1}$ 
$(l\in {\mab N},0\leq i\leq l)$ and 
a standard face morphism 
$p_i\col X_{ml}\lo X_{m,l-1}$ 
$(l>0,0\leq i\leq l)$ are obtained by local open immersions, 
we can easily check that there exists a unique member 
$\ol{e}_{l\pm 1}$ of the local minimal generators 
of $M_{m,l\pm 1}/{\cal O}_{X_{m,l\pm 1}}^*$ 
such that $s_i^*(\ol{e}_{l+1})=\ol{e}_l$ and 
$p_i^*(\ol{e}_{l-1})=\ol{e}_l$. 
Let $\ol{g}\col X_{mn}\lo X_{mn'}$ 
be the morphism 
corresponding to a morphism $[n']\lo [n]$. 
Since the morphism 
$[n']\lo [n]$ is a composite morphism of 
standard degeneracy morphisms and 
standard face morphisms, 
we see that there exists a unique member 
$\ol{e}_{n'}$ of the local minimal generators 
of $M_{mn'}/{\cal O}_{X_{mn'}}^*$ 
such that $\ol{g}{}^*(\ol{e}_{n'})=\ol{e}_n$. 
Let $\iota_{m\bul} \col X_{m\bul}
\os{\subset}{\lo} {\cal P}^{\rm ex}_{m\bul}$ 
be the immersion of simplicial log (formal) schemes.  
Then we have the following commutative diagram:  
\begin{equation*} 
\begin{CD} 
g_*({\cal M}_{mn}/{\cal O}_{{\cal P}^{\rm ex}_{mn}}^*) 
@>{\sim}>> 
\iota_{mn'*}\ol{g}_*(M_{mn}/{\cal O}_{X_{mn}}^*) \\
@A{g^*}AA @AA{\iota_{mn'*}(\ol{g}{}^*})A \\
{\cal M}_{mn'}/{\cal O}_{{\cal P}^{\rm ex}_{mn'}}^* 
@>{\sim}>> \iota_{mn'*}(M_{mn'}/{\cal O}_{X_{mn'}}^*). 
\end{CD} 
\tag{1.4.1.4}\label{cd:nnp}
\end{equation*} 
Since the immersion 
$X_{ml}\os{\sus}{\lo}{\cal P}^{\rm ex}_{ml}$ 
$(l=n, n')$ is exact and since 
$\os{\circ}{\cal P}{}^{\rm ex}_{ml}$ is topologically 
isomorphic to $\os{\circ}{X}_{ml}$, 
the horizontal morphisms in (\ref{cd:nnp}) 
are isomorphisms. 
Hence there exists a unique member $e_{n'}$ of 
the local minimal generators of 
${\cal M}_{mn'}/{\cal O}_{{\cal P}_{mn'}}^*$ 
such that $g^*(e_{n'})= e_n$. 
Because we have checked the condition in 
(\ref{prop:mmoo}), we have the simplicial formal scheme 
$\os{\circ}{\cal P}{}^{{\rm ex},(k)}_{m\bul}$ with the natural morphism 
$\os{\circ}{\cal P}{}^{{\rm ex},(k)}_{m\bul}\lo \os{\circ}{\cal P}{}^{\rm ex}_{m\bul}$. 
\par 
(2): Let the notations be as in (1). 
Since $s_i$ and $p_i$ are obtained by open immersions, 
any member of the local minimal generators 
of $M_{mn'}/{\cal O}_{X_{mn'}}^*$ which is different from $\ol{e}_{n'}$'s  
is mapped to 
$\ol{g}_*({\cal O}_{X_{mn}}^*)$ by the pull-back $\ol{g}{}^*$. 
The similar relation holds for $g^*$. 
Hence, by (1), (\ref{eqn:prvien}) and (\ref{prop:rescos}), 
we have the formula (\ref{ali:grpmd}). 
\par 
(3): We have the following formula by (2): 
\begin{align*} 
{\rm gr}^P_kA_{\rm zar}({\cal P}^{\rm ex}_{m\bul}/S(T)^{\nat},
{\cal E}^{m\bul}) & =
\bigoplus_{j\geq \max \{-k,0\}} 
({\rm gr}_{2j+k+1}^P
A_{\rm zar}({\cal P}^{\rm ex}_{m\bul}/S(T)^{\nat},
{\cal E}^{m\bul})^{\bul j}\{-j\},-\nabla)
\tag{1.4.1.5}\label{ali:grmd}\\
& =\bigoplus_{j\geq \max \{-k,0\}} 
({\rm gr}_{2j+k+1}^P({\cal E}^{m\bul}
\otimes_{{\cal O}_{{\cal P}^{\rm ex}_{m\bul}}}
{\Om}^{\bul +j+1}_{{\cal P}^{\rm ex}_{m\bul}/\os{\circ}{T}})
\{-j\},-\nabla) \\ 
& =\bigoplus_{j\geq \max \{-k,0\}} 
(c^{(2j+k)}_{m\bul*}(c^{(2j+k)*}_{m\bul}
{\cal E}^{m\bul}
\otimes_{
{\cal O}_{\os{\circ}{\cal P}{}^{{\rm ex},(2j+k)}_{m\bul}}}
{\Om}^{\bul}_{\os{\circ}{\cal P}{}^{{\rm ex},(2j+k)}_{m\bul}
/\os{\circ}{T}} \\
&\otimes_{\mab Z}
\vp^{(2j+k)}_{\rm zar}
(\os{\circ}{\cal P}{}^{\rm ex}_{m\bul}/\os{\circ}{T}))
[-2j-k],\nabla).  
\end{align*} 
The last statement in (3) follows from these equalities and (\ref{ali:eopudq}). 
\end{proof}

\begin{lemm}\label{lemm:ti} 
Let $k$ be a positive integer. 
Let 
\begin{align*} 
\iota^{(k-1)*}_{mn} \col  & 
c^{(k-1)}_{mn*}
(c^{(k-1)*}_{mn}({\cal E}^{mn})
\otimes_{{\cal O}_{{\cal P}^{\rm ex}_{mn}}}
\Om^{\bul}_{\os{\circ}{\cal P}{}^{{\rm ex},(k-1)}_{mn}/\os{\circ}{T}}
\otimes_{\mab Z}\vp^{(k-1)}_{\rm zar}(\os{\circ}{\cal P}{}^{\rm ex}_{mn}/\os{\circ}{T})) 
\lo 
\tag{1.4.2.1}\label{eqn:odips} \\
& 
c^{(k)}_{mn*}
(c^{(k)*}_{mn}({\cal E}^{mn})
\otimes_{{\cal O}_{{\cal P}^{\rm ex}_{mn}}}
\Om^{\bul}_{\os{\circ}{\cal P}{}^{{\rm ex},(k)}_{mn}
/\os{\circ}{T}}
\otimes_{\mab Z}
\vp^{(k)}_{\rm zar}(\os{\circ}{\cal P}{}^{\rm ex}_{mn}/\os{\circ}{T}))
\end{align*}
be the morphism {\rm (\ref{eqn:odpaps})}.  
Then the following diagram is commutative$:$
\begin{equation*} 
\begin{CD} 
{\rm gr}^P_{k+1}
({\cal E}^{mn}
\otimes_{{\cal O}_{{\cal P}^{\rm ex}_{mn}}}
{\Om}^{i+1}_{{\cal P}^{\rm ex}_{mn}/\os{\circ}{T}})
@>{\simeq}>>  
\\
@A{\theta_{{\cal P}^{\rm ex}_{mn}}\wedge}AA  \\
{\rm gr}^P_k
({\cal E}^{mn}
\otimes_{{\cal O}_{{\cal P}^{\rm ex}_{mn}}}
{\Om}^{i}_{{\cal P}^{\rm ex}_{mn}/\os{\circ}{T}})
@>{\simeq}>> 
\end{CD}
\tag{1.4.2.2}\label{eqn:xdxgras}
\end{equation*} 
\begin{equation*} 
\begin{CD} 
(c^{(k)}_{mn*}
(c^{(k)*}_{mn}
({\cal E}^{mn})\otimes_{
{\cal O}_{\os{\circ}{\cal P}{}^{{\rm ex},(k)}_{mn}}}
\Om^{i-k}_{\os{\circ}{\cal P}{}^{{\rm ex},(k)}_{mn}
/\os{\circ}{T}} \otimes_{\mab Z}
\vp^{(k)}_{\rm zar}
(\os{\circ}{\cal P}{}^{\rm ex}_{mn}/\os{\circ}{T}))[-k-1]) 
\\
@AA{\iota^{(k-1)*}_{mn}}A \\
(c^{(k-1)}_{mn*}
(c^{(k-1)*}_{mn}({\cal E}^{mn})\otimes_{
{\cal O}_{\os{\circ}{\cal P}{}^{{\rm ex},(k-1)}_{mn}}}
\Om^{i-k}_{\os{\circ}{\cal P}{}^{{\rm ex},(k-1)}_{mn}
/\os{\circ}{T}} \otimes_{\mab Z}
\vp^{(k-1)}_{\rm zar}
(\os{\circ}{\cal P}{}^{\rm ex}_{mn}/\os{\circ}{T}))[-k]). 
\end{CD}
\end{equation*} 
\end{lemm}
\begin{proof} 
The proof is the same as that of \cite[4.12]{msemi} (cf.~\cite[(10.1.16)]{ndw}). 
\end{proof}

\par 
Set $U_0:=S_{\os{\circ}{T}_0}$ or $T_0$ 
and $U:=S(T)^{\nat}$ or $T$,  respectively. 
\par 
Let 
\begin{align*} 
\pi_{{\rm zar}} \col 
&((X_{\bul \leq N,\bul,T_0})_{\rm zar},
f^{-1}_{\bul,T}({\cal O}_T))
=((X_{\bul \leq N,\bul,\os{\circ}{T}_0})_{\rm zar},
f^{-1}_{\bul}({\cal O}_T))
\tag{1.4.2.3}\label{ali:pzd} \\
& \lo 
((X_{\bul \leq N,\os{\circ}{T}_0})_{\rm zar},
f^{-1}({\cal O}_T)) 
=
((X_{\bul \leq N,T_0})_{\rm zar},
f^{-1}_{T}({\cal O}_T)) 
\end{align*} 
and 
\begin{align*} 
\pi_{m,{\rm zar}} \col 
&((X_{m\bul,T_0})_{\rm zar},f^{-1}_{m\bul,T}({\cal O}_T))
=((X_{m\bul,\os{\circ}{T}_0})_{\rm zar},f^{-1}_{m\bul}({\cal O}_T))
\tag{1.4.2.4}\label{ali:pmzd} \\ 
&\lo((X_{m,\os{\circ}{T}_0})_{\rm zar},f^{-1}_{m}({\cal O}_T))
=((X_{m,T_0})_{\rm zar},f^{-1}_{m,T}({\cal O}_T)) 
\end{align*} 
be the natural morphisms of ringed topoi.  
Let 
\begin{align*} 
u_{X_{\bul \leq N,U_0}/U} \col 
((X_{\bul \leq N,U_0}/U)_{\rm crys},{\cal O}_{X_{\bul \leq N,U_0}/U}) 
&\lo ((X_{\bul \leq N,U_0})_{\rm zar},f^{-1}({\cal O}_T))
\tag{1.4.2.5}\label{eqn:prxudef}\\
&=((X_{\bul \leq N,\os{\circ}{T}_0})_{\rm zar},f^{-1}({\cal O}_T))
\end{align*}    
and 
\begin{align*} 
\os{\circ}{u}_{X_{\bul \leq N,\os{\circ}{T}_0}/S(T)^{\nat}}=\os{\circ}{u}_{X_{\bul \leq N,T_0}/T} 
\col 
((\os{\circ}{X}_{\bul \leq N,T_0}/\os{\circ}{T})_{\rm crys},
{\cal O}_{\os{\circ}{X}_{\bul \leq N,T_0}/\os{\circ}{T}}) 
&\lo ((X_{\bul \leq N,T_0})_{\rm zar},f^{-1}_T({\cal O}_T))
\tag{1.4.2.6}\label{eqn:pxuxxef}\\
&=((X_{\bul \leq N,T_0})_{\rm zar},f^{-1}({\cal O}_T))
\end{align*}   
be the natural projections. 
Let 
\begin{align*} 
\eps_{X_{\bul \leq N,U_0}/U} \col 
((X_{\bul \leq N,U_0}/U)_{\rm crys},{\cal O}_{X_{\bul \leq N,U_0}/U}) 
&\lo 
((\os{\circ}{X}_{\bul \leq N,U_0}/\os{\circ}{U})_{\rm crys},
{\cal O}_{\os{\circ}{X}_{\bul \leq N,U_0}/\os{\circ}{U}})
\tag{1.4.2.7}\label{eqn:preoxdef}\\
&=((\os{\circ}{X}_{\bul \leq N,T_0}/\os{\circ}{T})_{\rm crys},
{\cal O}_{\os{\circ}{X}_{\bul \leq N,T_0}/\os{\circ}{T}})
\end{align*}   
be the morphism forgetting the log structures of $X_{\bul \leq N,U_0}$ and $U$. 
Let 
\begin{align*} 
\eps_{X_{\bul \leq N,U_0}/\os{\circ}{T}} \col 
((X_{\bul \leq N,U_0}/\os{\circ}{T})_{\rm crys},{\cal O}_{X_{\bul \leq N,U_0}/\os{\circ}{T}}) 
&\lo 
((\os{\circ}{X}_{\bul \leq N,U_0}/\os{\circ}{T})_{\rm crys},
{\cal O}_{\os{\circ}{X}_{\bul \leq N,U_0}/\os{\circ}{T}})\tag{1.4.2.8}\label{ali:prtsdef}\\
&=
((\os{\circ}{X}_{\bul \leq N,T_0}/\os{\circ}{T})_{\rm crys},
{\cal O}_{\os{\circ}{X}_{\bul \leq N,T_0}/\os{\circ}{T}})
\end{align*}   
be the morphism forgetting the log structure of $X_{\bul \leq N,U_0}$. 
Let 
\begin{equation*}  
\eps_{X_{\bul \leq N,U_0}/U/\os{\circ}{T}}\col 
((X_{\bul \leq N,U_0}/U)_{\rm crys},{\cal O}_{X_{\bul \leq N,U_0}/U})\lo 
((X_{\bul \leq N,U_0}/\os{\circ}{T})_{\rm crys},
{\cal O}_{X_{\bul \leq N,U_0}/\os{\circ}{T}})
\tag{1.4.2.9}\label{eqn:txntpi} 
\end{equation*} 
be the morphism forgetting the log structure of $U$. 
Then 
\begin{equation*}  
\eps_{X_{\bul \leq N,U_0}/U}=\eps_{X_{\bul \leq N,U_0}/\os{\circ}{T}}
\circ \eps_{X_{\bul \leq N,U_0}/U/\os{\circ}{T}}. 
\tag{1.4.2.10}\label{eqn:tneupi} 
\end{equation*} 

\begin{prop}\label{prop:tefc}
There exists the following isomorphism 
\begin{align*} 
\theta:=\theta_{X_{\bul \leq N,\os{\circ}{T}_0}/S(T)^{\nat}}
\wedge \col 
&Ru_{X_{\bul \leq N,\os{\circ}{T}_0}/S(T)^{\nat}*}
(\eps^*_{X_{\bul \leq N,\os{\circ}{T}_0/S(T)^{\nat}}}(E^{\bul \leq N}))
\tag{1.4.3.1}\label{eqn:uz} \\
& \os{\sim}{\lo}R\pi_{{\rm zar}*}
(A_{\rm zar}({\cal P}^{\rm ex}_{\bul \leq N,\bul}/S(T)^{\nat},
{\cal E}^{\bul \leq N,\bul}))
\end{align*} 
in ${\rm D}^+(
f^{-1}({\cal O}_T))$. 
This isomorphism is independent of the choice of 
an affine $N$-truncated simplicial open covering of 
$X_{\bul \leq N,\os{\circ}{T}_0}$ 
and the choice of an $(N,\infty)$-truncated bisimplicial immersion 
$X_{\bul \leq N,\bul,\os{\circ}{T}_0} \os{\sus}{\lo} 
\ol{\cal P}_{\bul \leq N,\bul}$ over $\ol{S(T)^{\nat}}$.
In particular, the complex 
$R\pi_{{\rm zar}*}
(A_{\rm zar}({\cal P}^{\rm ex}_{\bul \leq N,\bul}/S(T)^{\nat},
{\cal E}^{\bul \leq N,\bul}))$ 
is independent of the choices above. 
\end{prop}
\begin{proof} 
First we claim that there exists an isomorphism from the source 
to the target of (\ref{eqn:uz}). 
Let 
\begin{align*} 
\pi_{S(T)^{\nat}{\rm crys}} 
\col &((X_{\bul \leq N,\bul,\os{\circ}{T}_0}/S(T)^{\nat})_{\rm crys},
{\cal O}_{X_{\bul \leq N,\bul,\os{\circ}{T}_0}/S(T)^{\nat}})   
\tag{1.4.3.2}\label{eqn:pio} \\
&\lo ((X_{\bul \leq N,\os{\circ}{T}_0}/S(T)^{\nat})_{\rm crys},
{\cal O}_{X_{\bul \leq N,\os{\circ}{T}_0}/S(T)^{\nat}})
\end{align*} 
be the natural morphism of ringed topoi. 
Then, by the cohomological descent, 
we have 
$\eps^*_{X_{\bul \leq N,\os{\circ}{T}_0}/S(T)^{\nat}}(E^{\bul \leq N})=
R\pi_{S(T)^{\nat}{\rm crys}*}
(\eps^*_{X_{\bul \leq N,\bul,\os{\circ}{T}_0}/S(T)^{\nat}}(E^{\bul \leq N,\bul}))$. 
Let 
\begin{align*} 
u_{X_{\bul \leq N,\bul,\os{\circ}{T}_0}/S(T)^{\nat}} \col 
((X_{\bul \leq N,\bul,\os{\circ}{T}_0}/S(T)^{\nat})_{\rm crys},
{\cal O}_{X_{\bul \leq N,\bul,\os{\circ}{T}_0}/S(T)^{\nat}}) 
&\lo ((X_{\bul \leq N,\bul,\os{\circ}{T}_0})_{\rm zar},f^{-1}_{\bul \leq N,\bul}({\cal O}_T))
\tag{1.4.3.3}\label{eqn:uudef}\\
\end{align*}     
be the natural projection. 
Then $u_{X_{\bul \leq N,\os{\circ}{T}_0}/S(T)^{\nat}}\circ \pi_{S(T)^{\nat}{\rm crys}}
=\pi_{\rm zar}\circ u_{X_{\bul \leq N,\bul,\os{\circ}{T}_0}/S(T)^{\nat}}$. 
Hence 
we have the following formula by the log Poincar\'{e} lemma: 
\begin{align*} 
& Ru_{X_{\bul \leq N,\os{\circ}{T}_0}/S(T)^{\nat}*}
(\eps^*_{X_{\bul \leq N,\os{\circ}{T}_0}/S(T)^{\nat}}(E^{\bul \leq N}))\tag{1.4.3.4}\label{eqn:xds} \\
&=
Ru_{X_{\bul \leq N,\os{\circ}{T}_0}/S(T)^{\nat}*}
R\pi_{S(T)^{\nat}{\rm crys}*}
(\eps^*_{X_{\bul \leq N,\bul,\os{\circ}{T}_0}/S(T)^{\nat}}(E^{\bul \leq N,\bul}))\\
&=R\pi_{{\rm zar}*}
Ru_{X_{\bul \leq N,\bul,\os{\circ}{T}_0}/S(T)^{\nat}*}
(\eps^*_{X_{\bul \leq N,\bul,\os{\circ}{T}_0}/S(T)^{\nat}}(E^{\bul \leq N,\bul})) \\
&=R\pi_{{\rm zar}*}
({\cal E}^{\bul \leq N,\bul}
\otimes_{{\cal O}_{{\cal P}^{\rm ex}_{\bul \leq N,\bul}}}
{\Om}^{\bul}_{{\cal P}^{\rm ex}_{\bul \leq N,\bul}/\os{\circ}{T}}). 
\end{align*}  
Since ${\cal E}^{\bul \leq N,\bul}$ is a flat
${\cal O}_{{\mathfrak D}_{\bul \leq N,\bul}}$-module, 
it suffices to prove that the natural morphism 
\begin{align*} 
\theta_{{\cal P}^{\rm ex}_{\bul \leq N,\bul}} \wedge  & 
\col 
{\cal O}_{{\mathfrak D}_{\bul \leq N,\bul}}
\otimes_{{\cal O}_{{\cal P}^{\rm ex}_{\bul \leq N,\bul}}}
{\Om}^i_{{\cal P}^{\rm ex}_{\bul \leq N,\bul}/S(T)^{\nat}}
\lo \\
& \{({\cal O}_{{\mathfrak D}_{\bul \leq N,\bul}}
\otimes_{{\cal O}_{{\cal P}^{\rm ex}_{\bul \leq N,\bul}}}
{\Om}^{i+1}_{{\cal P}^{\rm ex}_{\bul \leq N,\bul}/\os{\circ}{T}}
/P_0   \os{\theta_{{\cal P}^{\rm ex}_{\bul \leq N,\bul}}}{\lo} 
{\cal O}_{{\mathfrak D}_{\bul \leq N,\bul}}
\otimes_{{\cal O}_{{\cal P}^{\rm ex}_{\bul \leq N,\bul}}}
{\Om}^{i+2}_{{\cal P}^{\rm ex}_{\bul \leq N,\bul}/\os{\circ}{T}}/P_1 \\ 
& \os{\theta_{{\cal P}^{\rm ex}_{\bul \leq N,\bul}}}{\lo}\cdots )\}
\end{align*} 
is a quasi-isomorphism. 
Hence it suffices to prove that 
$\theta_{{\cal P}^{\rm ex}_{mn}} \wedge$ 
is a quasi-isomorphism. 
\par 
As in \cite[3.15]{msemi} 
(cf.~\cite[(6.28) (9), (6.29) (1)]{ndw}), 
it suffices to prove that the sequence 
\begin{align*} 
& 0 \lo 
{\rm gr}_0^P
({\cal O}_{{\mathfrak D}_{mn}}
{\otimes}_{{\cal O}_{{\cal P}^{\rm ex}_{mn}}} 
{{\Om}}^{\bul}_{{\cal P}^{\rm ex}_{mn}/\os{\circ}{T}}) 
\os{\theta_{{\cal P}^{\rm ex}_{mn}}\wedge}{\lo}
{\rm gr}_1^P
({\cal O}_{{\mathfrak D}_{mn}}
{\otimes}_{{\cal O}_{{\cal P}^{\rm ex}_{mn}}}
{{\Om}}^{\bul}_{{\cal P}^{\rm ex}_{mn}/\os{\circ}{T}})[1] 
\tag{1.4.3.5}\label{eqn:lgrx}\\
& \os{\theta_{{\cal P}^{\rm ex}_{mn}}\wedge}{\lo}
{\rm gr}_2^P
({\cal O}_{{\mathfrak D}_{mn}}
{\otimes}_{{\cal O}_{{\cal P}_{mn}^{\rm ex}}} 
{{\Om}}^{\bul}_{{\cal P}^{\rm ex}_{mn}/\os{\circ}{T}})[2]  
\os{\theta_{{\cal P}^{\rm ex}_{mn}}\wedge}{\lo} \cdots 
\end{align*} 
is exact. 
By (\ref{lemm:ti}) we have only to prove that 
the following sequence 
\begin{align*} 
0 & \lo P_0({\cal O}_{{\mathfrak D}_{mn}}
{\otimes}_{{\cal O}_{{\cal P}{}^{\rm ex}_{mn}}} 
{\Om}^{\bul}_{{{\cal P}}{}^{\rm ex}_{mn}
/\os{\circ}{T}}) \\ 
&  \lo 
c^{(0)}_{mn *}
({\cal O}_{\os{\circ}{\mathfrak D}{}^{(0)}_{mn}}
{\otimes}_{{\cal O}_{{\cal P}{}^{{\rm ex},(0)}_{mn}}} 
{\Om}^{\bul}_{\os{\circ}{\cal P}{}^{{\rm ex},(0)}_{mn}
/\os{\circ}{T}}
\otimes_{\mab Z}\vp^{(0)}_{\rm zar}
(\os{\circ}{\cal P}{}^{\rm ex}_{mn}/\os{\circ}{T})) \\
&   
\os{\iota^{(0)*}}{\lo} 
c^{(1)}_{mn *}
({\cal O}_{\os{\circ}{\mathfrak D}{}^{(1)}_{mn}}
{\otimes}_{{\cal O}_{{\cal P}{}^{{\rm ex},(1)}_{mn}}}
{\Om}^{\bul}_{
\os{\circ}{\cal P}{}^{{\rm ex},(1)}_{mn}/\os{\circ}{T}}
\otimes_{\mab Z}\vp^{(1)}_{\rm zar}
(\os{\circ}{\cal P}{}^{\rm ex}_{mn}/\os{\circ}{T})) 
\lo \cdots 
\end{align*}
is exact. In (\ref{prop:csrl}) we have already proved this exactness. 
\par 
Next we prove that the isomorphism is independent of 
the choices in (\ref{prop:tefc}).   
\par 
Let $X''_{\bul \leq N,\os{\circ}{T}_0}$ be another disjoint union 
of an affine $N$-truncated simplicial open covering of 
$X_{\bul \leq N,\os{\circ}{T}_0}$. 
Set $X'''_{m,\os{\circ}{T}_0}:=
X'_{m,\os{\circ}{T}_0}\times_{X_{m,\os{\circ}{T}_0}}X''_{m,\os{\circ}{T}_0}$ 
$(0\leq m\leq N)$. 
Then we have the $N$-truncated simplicial log scheme 
$X'''_{\bul \leq N,\os{\circ}{T}_0}$ which is 
the disjoint union of 
the members of an (not necessarily affine) $N$-truncated simplicial open covering 
of $X_{\bul \leq N,\os{\circ}{T}_0}/S_{\os{\circ}{T}_0}$ 
fitting into the following commutative diagram$:$
\begin{equation*}
\begin{CD}
X'''_{\bul \leq N,\os{\circ}{T}_0} @>>> X''_{\bul \leq N,\os{\circ}{T}_0}\\
@VVV @VVV  \\
X'_{\bul \leq N,\os{\circ}{T}_0}@>>> X_{\bul \leq N,\os{\circ}{T}_0}.
\end{CD}
\tag{1.4.3.6}\label{cd:celcxcov}
\end{equation*} 
Set $X_{mn,\os{\circ}{T}_0}:={\rm cosk}^{X_{m,\os{\circ}{T}_0}}_0(X'_{m,\os{\circ}{T}_0})_n$, 
$X'_{mn,\os{\circ}{T}_0}:={\rm cosk}^{X_{m,\os{\circ}{T}_0}}_0(X''_{m,\os{\circ}{T}_0})_n$  
and $X''_{mn,\os{\circ}{T}_0}:={\rm cosk}^{X_{m,\os{\circ}{T}_0}}_0(X'''_{m,\os{\circ}{T}_0})_n$ 
$(0\leq m \leq N, n\in {\mab N})$.  
Let 
$X'_{\bul \leq N,\bul,\os{\circ}{T}_0}\os{\sus}{\lo}
\ol{\cal P}{}'_{\bul \leq N,\bul}$ be 
another immersion in (\ref{eqn:eipxd}).  
Set $\ol{\cal P}{}''_{\bul \leq N,\bul}:=
\ol{\cal P}_{\bul \leq N,\bul}
\times_{\ol{S(T)^{\nat}}}\ol{\cal P}{}'_{\bul \leq N,\bul}$.  
Then we have the following commutative diagram 
\begin{equation*} 
\begin{CD} 
X_{\bul \leq N,\bul,\os{\circ}{T}_0} @>{\sus}>> \ol{\cal P}_{\bul \leq N,\bul} \\ 
@AAA @AAA \\ 
X''_{\bul \leq N,\bul,\os{\circ}{T}_0} @>{\sus}>> 
\ol{\cal P}{}''_{\bul \leq N,\bul} \\ 
@VVV @VVV \\ 
X'_{\bul \leq N,\bul,\os{\circ}{T}_0} @>{\sus}>> 
\ol{\cal P}{}'_{\bul \leq N,\bul}.  
\end{CD} 
\end{equation*} 
Hence we may assume that 
there exists the following commutative diagram 
\begin{equation*} 
\begin{CD} 
X_{\bul \leq N,\bul,\os{\circ}{T}_0} @>{\sus}>> \ol{\cal P}_{\bul \leq N,\bul} \\ 
@VVV @VVV \\ 
X'_{\bul \leq N,\bul,\os{\circ}{T}_0} @>{\sus}>> \ol{\cal P}{}'_{\bul \leq N,\bul}.  
\end{CD} 
\tag{1.4.3.7}\label{cd:sstp} 
\end{equation*}  
Set ${\cal P}'_{\bul \leq N,\bul}:=\ol{\cal P}{}'_{\bul \leq N,\bul}\times_{\ol{S(T)^{\nat}}}S(T)^{\nat}$. 
Let ${\cal P}{}'{}^{\rm ex}_{\!\!\!\bul \leq N,\bul}$  
be the exactification of 
the immersion 
$X'_{\bul \leq N,\bul,\os{\circ}{T}_0} \os{\sus}{\lo}{\cal P}'_{\bul \leq N,\bul}$. 
Let ${\cal E}'{}^{\bul \leq N,\bul}
\otimes_{{\cal O}_{{\cal P}'{}^{\rm ex}_{\!\!\!\bul \leq N,\bul}}}
{\Om}^{\bul}_{{\cal P}'{}^{\rm ex}_{\bul \leq N,\bul}/\os{\circ}{T}}$ 
be an analogous complex 
to 
${\cal E}^{\bul \leq N,\bul}
\otimes_{{\cal O}_{{\cal P}^{\rm ex}_{\bul \leq N,\bul}}}
{\Om}^{\bul}_{{\cal P}^{\rm ex}_{\bul \leq N,\bul}/\os{\circ}{T}}$ 
for ${\cal P}'_{\bul \leq N,\bul}$.   
Then we have the following morphism 
\begin{equation*} 
R\pi_{{\rm zar}*}
(A_{\rm zar}({\cal P}'{}^{\rm ex}_{\! \! \! \bul \leq  N,\bul}/S(T)^{\nat}
,{\cal E}'{}^{\bul \leq N,\bul}),P))
\lo 
R\pi_{{\rm zar}*}
(A_{\rm zar}({\cal P}^{\rm ex}_{\bul \leq N,\bul}/S(T)^{\nat}
,{\cal E}^{\bul \leq N,\bul}),P). 
\tag{1.4.3.8}\label{eqn:flap}
\end{equation*} 
In particular, we have the morphism 
\begin{equation*} 
R\pi_{{\rm zar}*}
(A_{\rm zar}({\cal P}'{}^{\rm ex}_{\! \! \! \bul \leq  N,\bul}/S(T)^{\nat}
,{\cal E}'{}^{\bul \leq N,\bul}))
\lo 
R\pi_{{\rm zar}*}
(A_{\rm zar}({\cal P}^{\rm ex}_{\bul \leq N,\bul}/S(T)^{\nat},{\cal E}^{\bul \leq N,\bul})). 
\tag{1.4.3.9}\label{eqn:flnfap}
\end{equation*} 
This morphism fits into the following commutative diagram 
\begin{equation*} 
\begin{CD} 
 R\pi_{{\rm zar}*}
(A_{\rm zar}({\cal P}'{}^{\rm ex}_{\! \! \! \bul \leq  N,\bul}/S(T)^{\nat}
,{\cal E}'{}^{\bul \leq N,\bul}))
@>>> 
R\pi_{{\rm zar}*}
(A_{\rm zar}({\cal P}^{\rm ex}_{\! \! \! \bul \leq  N,\bul}/S(T)^{\nat}
,{\cal E}^{\bul \leq N,\bul}))\\
@A{R\pi_{{\rm zar}*}(\theta_{{\cal P}'{}^{\rm ex}_{\! \! \! \bul \leq N,\bul}} \wedge)}A{\simeq}A 
@A{\simeq}A{R\pi_{{\rm zar}*}(\theta_{{\cal P}^{\rm ex}_{\!  \bul \leq N,\bul}} \wedge)}A \\
R\pi_{{\rm zar}*}
({\cal E}'{}^{\bul \leq N,\bul}
\otimes_{{\cal O}_{{\cal P}'{}^{\rm ex}_{\bul \leq N,\bul}}}
{\Om}^{\bul}_{{\cal P}'{}^{\rm ex}_{\bul \leq N,\bul}/S(T)^{\nat}})
@=
R\pi_{{\rm zar}*}
({\cal E}^{\bul \leq N,\bul}
\otimes_{{\cal O}_{{\cal P}^{\rm ex}_{\bul \leq N,\bul}}}
{\Om}^{\bul}_{{\cal P}^{\rm ex}_{\bul \leq N,\bul}/S(T)^{\nat}}). 
\end{CD}
\tag{1.4.3.10}\label{cd:fnfap}
\end{equation*} 
This diagram tells us the desired independence 
of the choices in (\ref{prop:tefc}).  
\end{proof}

\par 
Next, by using 
(\ref{lemm:knit}) and (\ref{prop:tefc}),   
we prove that 
$$R\pi_{{\rm zar}*}
((A_{\rm zar}({\cal P}^{\rm ex}_{\bul \leq N,\bul}/S(T)^{\nat},
{\cal E}^{\bul \leq N,\bul}),P))$$ 
depends only on $X_{\bul \leq N,\os{\circ}{T}_0}/(S(T)^{\nat},{\cal J},\del)$:  


\begin{theo}\label{theo:indcr} 
The filtered complex 
$$R\pi_{{\rm zar}*}
((A_{\rm zar}({\cal P}^{\rm ex}_{\bul \leq N,\bul}
/S(T)^{\nat},{\cal E}^{\bul \leq N,\bul}),P))$$
in ${\rm D}^+{\rm F}(f^{-1}_{\bul \leq N}({\cal O}_T))$ 
is independent of the choice of 
an affine $N$-truncated simplicial open covering of 
$X_{\bul \leq N,\os{\circ}{T}_0}$ 
and the choice of an $(N,\infty)$-truncated bisimplicial immersion 
$X_{\bul \leq N,\bul,\os{\circ}{T}_0} \os{\sus}{\lo} 
\ol{\cal P}_{\bul \leq N,\bul}$ over $\ol{S(T)^{\nat}}$.  
\end{theo}
\begin{proof} 
Let the notations be as in the proof of (\ref{prop:tefc}). 
Let us prove that 
the morphism (\ref{eqn:flap}) is a filtered isomorphism. 
\par 
Because there exist 
filtered flasque resolutions   
$$(A_{\rm zar}({\cal P}^{\rm ex}_{\bul \leq N,\bul}/S(T)^{\nat}
,{\cal E}^{\bul \leq N,\bul}),P)
\lo (I^{\bul \leq N,\bul},P)$$ 
and 
$$(A_{\rm zar}({\cal P}'{}^{\rm ex}_{\! \! \! \bul \leq  N,\bul}/S(T)^{\nat}
,{\cal E}'{}^{\bul \leq N,\bul}),P)
\lo (I'{}^{\bul \leq N,\bul},P)$$  
such that 
$$(A_{\rm zar}({\cal P}^{\rm ex}_{m\bul}/S(T)^{\nat}
,{\cal E}^{m\bul})),P)
\lo (I^{m\bul},P) \quad (0\leq \forall m\leq N)$$ 
and 
$$(A_{\rm zar}({\cal P}'{}^{\rm ex}_{\! \! \! m\bul}/S(T)^{\nat},{\cal E}'{}^{m\bul}),P)
\lo (I'{}^{m\bul},P) \quad (0\leq \forall m\leq N)$$  
are filtered flasque resolutions 
(see, e.~g., \cite[(1.5.0.4)]{nh2}), we have only to prove that 
the morphism 
\begin{equation*} 
R\pi_{m,{\rm zar}*}
((A_{\rm zar}({\cal P}'{}^{\rm ex}_{\! \! \! m\bul}/S(T)^{\nat},{\cal E}'{}^{m\bul}),P))
\lo 
R\pi_{m,{\rm zar}*}
((A_{\rm zar}({\cal P}^{\rm ex}_{m\bul}/S(T)^{\nat},{\cal E}^{m\bul}),P))
\tag{1.4.4.1}\label{eqn:flmap}
\end{equation*} 
$(0\leq m\leq N)$ is an isomorphism. 
To prove this, 
we have only to prove that the morphism 
\begin{equation*} 
R\pi_{m,{\rm zar}*}
(A_{\rm zar}({\cal P}'{}^{\rm ex}_{\! \! \! m\bul}/S(T)^{\nat},{\cal E}'{}^{m\bul}))
\lo 
R\pi_{m,{\rm zar}*}(A_{\rm zar}({\cal P}^{\rm ex}_{m\bul}/S(T)^{\nat},{\cal E}^{m\bul}))
\tag{1.4.4.2}\label{eqn:flapwp}
\end{equation*} 
$(0\leq m\leq N)$ is an isomorphism in 
${\rm D}^+(
f^{-1}_m({\cal O}_T))$
and that the morphism 
\begin{equation*} 
R\pi_{m,{\rm zar}*}
(P_kA_{\rm zar}({\cal P}'{}^{\rm ex}_{\! \! \! m\bul}/S(T)^{\nat},
{\cal E}'{}^{m\bul}))
\lo 
R\pi_{m,{\rm zar}*}
(P_kA_{\rm zar}({\cal P}^{\rm ex}_{m\bul}/S(T)^{\nat},{\cal E}^{m\bul}))
\tag{1.4.4.3}\label{eqn:flsap}
\end{equation*} 
$(0\leq m\leq N, k\in {\mab Z})$ is an isomorphism in 
${\rm D}^+(
f^{-1}_m({\cal O}_T))$. 
\par 
By (\ref{eqn:uz}) we see that (\ref{eqn:flapwp}) is an isomorphism. 
Because the question that the morphism (\ref{eqn:flsap}) is an isomorphism is 
local on $X_{m,\os{\circ}{T}_0}$, 
we may assume that $X_{m,\os{\circ}{T}_0}$ is affine, 
in particular, quasi-compact. 
Hence, if $k\gg 0$, then 
\begin{align*}
P_kA_{\rm zar}({\cal P}^{\rm ex}_{m\bul}/S(T)^{\nat},{\cal E}^{m\bul})
=A_{\rm zar}({\cal P}^{\rm ex}_{m\bul}/S(T)^{\nat},{\cal E}^{m\bul}).
\tag{1.4.4.4}\label{ali:sepk} 
\end{align*} 
Consequently 
\begin{align*}
R\pi_{m,{\rm zar}*}(P_kA_{\rm zar}({\cal P}^{\rm ex}_{m\bul}/S(T)^{\nat},{\cal E}^{m\bul}))
=R\pi_{m,{\rm zar}*}(A_{\rm zar}({\cal P}^{\rm ex}_{m\bul}/S(T)^{\nat},{\cal E}^{m\bul}))
\tag{1.4.4.5}\label{ali:pipk} 
\end{align*} 
if $k\gg 0$. 
We have the following formula 
\begin{align*} 
& 
{\rm gr}^P_kR\pi_{m,{\rm zar}*}
(A_{\rm zar}({\cal P}^{\rm ex}_{m\bul}/S(T)^{\nat},{\cal E}^{m\bul}))
\os{\sim}{\lo} 
R\pi_{m,{\rm zar}*}
{\rm gr}^P_kA_{\rm zar}({\cal P}^{\rm ex}_{m\bul}/S(T)^{\nat},{\cal E}^{m\bul})
\tag{1.4.4.6}\label{ali:ruogrvp}\\
{} & \os{\sim}{\lo} 
R\pi_{m,{\rm zar}*}
(\bigoplus_{j\geq \max \{-k,0\}}
(a^{(2j+k)}_{m\bul,T_0*} 
(Ru_{\os{\circ}{X}{}^{(2j+k)}_{m\bul,T_0}/\os{\circ}{T}*}
(E^m_{\os{\circ}{X}{}^{(2j+k)}_{m\bul,T_0}
/\os{\circ}{T}}\otimes_{\mab Z}\vp_{\rm crys}^{(2j+k)}
(\os{\circ}{X}_{m\bul,T_0}/\os{\circ}{T})))\\
&[-2j-k],\nabla)) \\
& =\bigoplus_{j\geq \max \{-k,0\}} 
(a^{(2j+k)}_{m,T_0*}
(Ru_{\os{\circ}{X}{}^{(2j+k)}_{m,T_0}/\os{\circ}{T}*}
(E^m_{\os{\circ}{X}{}^{(2j+k)}_{m,T_0}
/\os{\circ}{T}} \otimes_{\mab Z}\vp_{\rm crys}^{(2j+k)}
(\os{\circ}{X}_{m,T_0}/\os{\circ}{T})))
[-2j-k],\nabla)
\end{align*} 
and the analogous formula 
for $R\pi_{m,{\rm zar}*}((A_{\rm zar}({\cal P}'{}^{\rm ex}_{\! \! \! m\bul}/S(T)^{\nat}
,{\cal E}'{}^{m\bul}),P))$. 
Here we have used \cite[(1.3.4.5)]{nh2} 
(resp.~(\ref{ali:grmpd}))
for the first isomorphism (resp. the second isomorphism); 
we obtain the last equality by the cohomological descent. 
We also obtain the following tautological triangle 
\begin{align*} 
&R\pi_{m,{\rm zar}*}(P_{k-1}A_{\rm zar}({\cal P}^{\rm ex}_{m\bul}/S(T)^{\nat},{\cal E}^{m\bul}))
\lo R\pi_{m,{\rm zar}*}(P_kA_{\rm zar}({\cal P}^{\rm ex}_{m\bul}/S(T)^{\nat},{\cal E}^{m\bul}))
\tag{1.4.4.7}\label{ali:rupprvp}\\
&\lo R\pi_{m,{\rm zar}*}
({\rm gr}^P_kA_{\rm zar}({\cal P}^{\rm ex}_{m\bul}/S(T)^{\nat},{\cal E}^{m\bul}))
\os{+1}{\lo}. 
\end{align*} 
and the analogous triangle for 
$$R\pi_{m,{\rm zar}*}
((A_{\rm zar}({\cal P}'{}^{\rm ex}_{\! \! \! m\bul}/S(T)^{\nat},{\cal E}'{}^{m\bul}),P)).$$ 
Now we see that (\ref{eqn:flsap}) is an isomorphism by 
(\ref{ali:pipk}), (\ref{ali:ruogrvp}) and (\ref{ali:rupprvp}) 
and the descending induction on $k$. 
\par 
We complete the proof of (\ref{theo:indcr}).  
\end{proof}

\begin{rema}\label{rema:asid} 
One can also prove (\ref{theo:indcr}) without using (\ref{prop:tefc}). 
Indeed, because we have only to prove that 
the morphism (\ref{eqn:flmap}) is an isomorphism, 
the problem is local on $X_{m,\os{\circ}{T}_0}$. 
Hence we may assume that $X_{m,\os{\circ}{T}_0}$ is affine, in particular, quasi-compact. 
Hence, if $k\gg 0$, then we may assume that 
$\os{\circ}{X}{}^{(k)}_{m,\os{\circ}{T}_0}=\emptyset$. 
If $k\ll 0$ (resp.~$k\gg 0$), then $2j+k\geq -k\gg 0$ (resp.~$2j+k\geq k\gg 0$). 
Hence the target of the isomorphism (\ref{ali:ruogrvp}) is $0$ in either case. 
Consequently 
\begin{align*} 
P_kR\pi_{m,{\rm zar}*}
(A_{\rm zar}({\cal P}^{\rm ex}_{m\bul}/S(T)^{\nat}
,{\cal E}^{m\bul}))=R\pi_{m,{\rm zar}*}
(A_{\rm zar}({\cal P}^{\rm ex}_{m\bul}/S(T)^{\nat},{\cal E}^{m\bul})).
\tag{1.4.5.1}\label{ali:ased} 
\end{align*}  
Now, by the isomorphism (\ref{ali:ruogrvp}) and the ascending induction on $k$,   
we see that the morphism (\ref{eqn:flmap}) is an isomorphism. 
\par 
We shall use this argument in (\ref{coro:indg}) below. 
\end{rema}

\begin{defi}\label{defi:fdirpd}
We call the filtered direct image 
$R\pi_{{\rm zar}*}
((A_{\rm zar}({\cal P}^{\rm ex}_{\bul \leq N,\bul}/S(T)^{\nat},{\cal E}^{\bul \leq N,\bul}),P))$ 
the {\it zariskian $p$-adic filtered Steenbrink complex} of 
$E^{\bul \leq N}$ for $X_{\bul \leq N,\os{\circ}{T}_0}/(S(T)^{\nat},{\cal J},\del)$. 
We denote it  by
$(A_{\rm zar}(X_{\bul \leq N,\os{\circ}{T}_0}/S(T)^{\nat},E^{\bul \leq N}),P)
\in {\rm D}^+{\rm F}(
f^{-1}({\cal O}_T))$. 
When $E^{\bul \leq N}={\cal O}_{\os{\circ}{X}_{\bul \leq N,T_0}/\os{\circ}{T}}$, we denote 
$(A_{\rm zar}(X_{\bul \leq N,\os{\circ}{T}_0}/S(T)^{\nat},E^{\bul \leq N}),P)$ 
by 
$(A_{\rm zar}(X_{\bul \leq N,\os{\circ}{T}_0}/S(T)^{\nat}),P)$. 
We call $(A_{\rm zar}(X_{\bul \leq N,\os{\circ}{T}_0}/S(T)^{\nat}),P)$ the 
{\it zariskian $p$-adic filtered Steenbrink complex} of 
$X_{\bul \leq N,\os{\circ}{T}_0}/(S(T)^{\nat},{\cal J},\del)$. 
\end{defi}


We restate the formula (\ref{ali:ruogrvp}) as 
the following proposition:  
\begin{prop}\label{prop:graxnp}
For each $0\leq m\leq N$, 
there exists the following canonical isomorphism$:$ 
\begin{align*} 
{\rm gr}^P_kA_{\rm zar}(X_{m,\os{\circ}{T}_0}/S(T)^{\nat},E^m)
\os{\sim}{\lo} \bigoplus_{j\geq \max \{-k,0\}} 
&a^{(2j+k)}_{m,T_0*} 
(Ru_{\os{\circ}{X}{}^{(2j+k)}_{m,T_0}/\os{\circ}{T}*}
(E^m_{\os{\circ}{X}{}^{(2j+k)}_{m,T_0}
/\os{\circ}{T}} \tag{1.4.7.1}\label{ali:ruovp}\\
&\otimes_{\mab Z}\vp_{\rm crys}^{(2j+k)}
(\os{\circ}{X}_{m,T_0}/\os{\circ}{T})))[-2j-k]
\end{align*}
in ${\rm D}^+(f^{-1}({\cal O}_T))$. 
\end{prop}

\begin{rema}\label{rema:trcadwd}
(1) In the case of the trivial coefficient  
$E={\cal O}_{\os{\circ}{X}_{\bul \leq N,T_0}/\os{\circ}{T}}$, 
we need not assume the existence of $\ol{\cal P}_{\bul \leq N,\bul}$ 
in this section; 
we have only to assume the existence of 
${\cal P}_{\bul \leq N,\bul}$ because we need 
only the filtered double complex 
$({\cal O}_{{\mathfrak D}_{\bul \leq N,\bul}}
\otimes_{{\cal O}_{{\cal P}_{\bul \leq N,\bul}}}
{\Om}^{i+j+1}_{{\cal P}_{\bul \leq N,\bul}/\os{\circ}{T}}/P_j,P)_{i,j\in {\mab N}}$ 
for the definition of 
$(A_{\rm zar}(X_{\bul \leq N,\os{\circ}{T}_0}/S(T)^{\nat}),P)$. 
\par 
(2) Let $S$ be as in \S\ref{sec:ldc} or 
let $S$ be a fine log formal scheme 
whose underlying formal scheme is a 
$p$-adic formal ${\cal V}$-scheme 
in the sense of \cite{of}: 
$\os{\circ}{S}$ is a noetherian formal scheme 
over ${\rm Spf}({\cal V})$ 
with the $p$-adic topology 
which is topologically of finite type over 
${\rm Spf}({\cal V})$.   
Let $(T,{\cal J},\del)$ be a log PD ($p$-adic formal) scheme. 
We assume that $\os{\circ}{T}$ is flat over ${\rm Spf}({\cal V})$.
Let $T_0\lo S$ be a morphism of log PD ($p$-adic formal) schemes. 
Let $\pi$ be a nonzero element of 
the maximal ideal of ${\cal V}$.
Assume also that $\pi{\cal O}_T$ has a PD-structure $\del$.    
In a future paper in which 
we discuss a $p$-adic Clemens-Schmid (exact) sequence (cf.~\cite{ctcs}), 
we have to consider 
$(A_{\rm zar}(X_{\bul \leq N,\os{\circ}{T}_0}/S(T)^{\nat},E^{\bul \leq N}),P)$ 
for a more generalized sheaf than $E^{\bul \leq N}$. 
Namely, assume that, for the disjoint union 
$X'_{\bul \leq N,\os{\circ}{T}_0}$ 
of an affine $N$-truncated simplicial open covering of 
$X_{\bul \leq N,\os{\circ}{T}_0}$ and 
an $(N,\infty)$-truncated bisimplicial immersion 
$X_{\bul \leq N,\bul,\os{\circ}{T}_0} \os{\sus}{\lo} \ol{\cal P}_{\bul \leq N,\bul}$ over 
$\ol{S(T)^{\nat}}$ 
($X_{mn,\os{\circ}{T}_0}:={\rm cosk}^{X_{m,\os{\circ}{T}_0}}_0(X'_{m,\os{\circ}{T}_0})_n$ 
$(0\leq m \leq N, n\in {\mab N})$) into a log smooth 
$(N,\infty)$-truncated log scheme over $\ol{S(T)^{\nat}}$, we are given a 
flat quasi-coherent ${\cal O}_{{\mathfrak D}_{\bul \leq N,\bul}}$-module 
or a flat quasi-coherent ${\cal K}_{{\mathfrak D}_{\bul \leq N,\bul}}$-module 
(${\cal K}_{{\mathfrak D}_{\bul \leq N,\bul}}:=
{\cal O}_{{\mathfrak D}_{\bul \leq N,\bul}}\otimes_{\mab Z}{\mab Q}$)
with integrable connection  
$({\cal E}({\mathfrak D}_{\bul \leq N,\bul}),\nabla)$: 
\begin{equation*} 
\nabla \col {\cal E}({\mathfrak D}_{\bul \leq N,\bul})\lo 
{\cal E}({\mathfrak D}_{\bul \leq N,\bul})
\otimes_{{\cal O}_{{\cal P}^{\rm ex}_{\bul \leq N,\bul}}}
\Om^1_{{\cal P}^{\rm ex}_{\bul \leq N,\bul}/\os{\circ}{T}}. 
\end{equation*} 
We endow 
${\cal E}({\mathfrak D}_{\bul \leq N,\bul})
\otimes_{{\cal O}_{{\cal P}^{\rm ex}_{\bul \leq N,\bul}}}
\Om^{\bul}_{{\cal P}^{\rm ex}_{\bul \leq N,\bul}/\os{\circ}{T}}$ 
with the tensor product $P$ of the trivial filtration on 
${\cal E}({\mathfrak D}_{\bul \leq N,\bul})$ and 
the (pre)weight filtration (\ref{eqn:pkdefpw}) 
on $\Om^{\bul}_{{\cal P}^{\rm ex}_{\bul \leq N,\bul}/\os{\circ}{T}}$. 
Then, 
by the same construction as that of 
(\ref{cd:adef}), (\ref{cd:locstbd}) and (\ref{eqn:dblad}), 
we have the filtered double complex 
$(A_{\rm zar}({\cal P}^{\rm ex}_{\bul \leq N,\bul}/S(T)^{\nat},
{\cal E}({\mathfrak D}_{\bul \leq N,\bul}))^{\bul \bul},P)$ 
and the single filtered complex 
$(A_{\rm zar}({\cal P}^{\rm ex}_{\bul \leq N,\bul}/S(T)^{\nat},{\cal E}
({\mathfrak D}_{\bul \leq N,\bul})),P)$  
associated to this double complex. 
We assume that the following three conditions hold: 
\par 
(1) We assume that, for another immersion  
$X'_{\bul \leq N,\bul,\os{\circ}{T}_0} \os{\sus}{\lo} 
{\cal P}{}'_{\bul \leq N,\bul}$ above and 
for the following commutative diagram 
\begin{equation*} 
\begin{CD} 
X_{\bul \leq N,\bul,\os{\circ}{T}_0} @>{\sus}>> {\cal P}_{\bul \leq N,\bul} \\ 
@VVV @VV{g_{\bul \leq N,\bul}}V \\ 
X'_{\bul \leq N,\bul,\os{\circ}{T}_0} @>{\sus}>> 
{\cal P}{}'_{\bul \leq N,\bul},   
\end{CD} 
\tag{1.4.8.1}\label{cd:cop}
\end{equation*} 
there exists a morphism 
\begin{equation*} 
\rho_{g_{\bul \leq N,\bul}} 
\col {\cal E}({\mathfrak D}'_{\bul \leq N,\bul})
\lo g_{\bul \leq N,\bul*}({\cal E}({\mathfrak D}_{\bul \leq N,\bul}))
\tag{1.4.8.2}\label{cd:cngp}
\end{equation*} 
of ${\cal O}_{{\mathfrak D}'_{\bul \leq N,\bul}}$-modules 
or ${\cal K}_{{\mathfrak D}'_{\bul \leq N,\bul}}$-modules 
satisfying the usual cocycle condition 
for the following commutative diagram: 
\begin{equation*} 
\begin{CD} 
X_{\bul \leq N,\bul,\os{\circ}{T}_0} @>{\sus}>> {\cal P}_{\bul \leq N,\bul} \\ 
@VVV @VV{}V \\ 
X'_{\bul \leq N,\bul,\os{\circ}{T}_0} @>{\sus}>> {\cal P}{}'_{\bul \leq N,\bul}\\
@VVV @VV{}V \\ 
X''_{\bul \leq N,\bul,\os{\circ}{T}_0} @>{\sus}>> 
{\cal P}{}''_{\bul \leq N,\bul}.   
\end{CD} 
\tag{1.4.8.3}\label{cd:cctpt}
\end{equation*} 
Here ${\mathfrak D}'_{\bul \leq N,\bul}$ is the analogue of 
${\mathfrak D}_{\bul \leq N,\bul}$ for the immersion 
$X'_{\bul \leq N,\bul,\os{\circ}{T}_0} \os{\sus}{\lo} {\cal P}'_{\bul \leq N,\bul}$. 
\par 
(2) We assume that the following morphism 
\begin{align*} 
R\pi_{{\rm zar}*}(A_{\rm zar}({\cal P}'{}^{\rm ex}_{\bul \leq N,\bul}/S(T)^{\nat},
{\cal E}({\mathfrak D}'_{\bul \leq N,\bul}))\lo 
R\pi_{{\rm zar}*}(A_{\rm zar}({\cal P}^{\rm ex}_{\bul \leq N,\bul}/S(T)^{\nat},
{\cal E}({\mathfrak D}_{\bul \leq N,\bul})) 
\end{align*}
induced by the morphisms 
$g_{\bul \leq N,\bul}\col {\cal P}_{\bul \leq N,\bul}\lo {\cal P}'_{\bul \leq N,\bul}$ 
and $\rho_{g_{\bul \leq N,\bul}}$ is an isomorphism 
in ${\rm D}^+(f^{-1}({\cal O}_T))$. 
\par 
(3) We assume that, for a positive integer $k$, 
for nonnegative integers $m$ and $n$ $(0\leq m \leq N, n\in {\mab N})$), 
the following isomorphism 
\begin{equation*} 
{\rm gr}^P_k({\cal E}({\mathfrak D}'_{mn})
\otimes_{{\cal O}_{{\cal P}'{}^{\rm ex}_{\! \!mn}}}
{\Om}^{\bul}_{{\cal P}'{}^{\rm ex}_{\! \!mn}/\os{\circ}{T}})
\os{\sim}{\lo}  
Rg_{mn*}({\rm gr}^P_k({\cal E}({\mathfrak D}_{mn})
\otimes_{{\cal O}_{{\cal P}^{\rm ex}_{mn}}}
{\Om}^{\bul}_{{\cal P}^{\rm ex}_{mn}/\os{\circ}{T}}))
\tag{1.4.8.4}\label{eqn:ptep}
\end{equation*} 
induced by the morphisms 
$g_{\bul \leq N,\bul}\col {\cal P}_{\bul \leq N,\bul}\lo {\cal P}'_{\bul \leq N,\bul}$ 
and $\rho_{g_{\bul \leq N,\bul}}$ is an isomorphism in ${\rm D}^+(f^{-1}_{mn}({\cal O}_T))$. 
\bigskip
\parno 
Then, by the same proof as that of the well-definedness of 
$(A_{\rm zar}(X_{\bul \leq N,\os{\circ}{T}_0}/S(T)^{\nat},E^{\bul \leq N}),P)$, 
we can prove that the filtered complex 
$R\pi_{{\rm zar}*}
((A_{\rm zar}({\cal P}^{\rm ex}_{\bul \leq N,\bul}/S(T)^{\nat},
{\cal E}({\mathfrak D}_{\bul \leq N,\bul})),P))$ 
is a well-defined filtered complex. 
\end{rema}

\begin{lemm}\label{lemm:flpis} 
Let $N$ be a nonnegative integer. 
Then the following hold$:$
\par 
$(1)$ Let $(U',{\cal K}',\eps')$ and $(U,{\cal K},\eps)$ be fine log PD-schemes 
on which $p$ is locally nilpotent. 
Set $U_0:=\ul{\rm Spec}^{\log}_U({\cal O}_U/{\cal K})$ 
and $U'_0:=\ul{\rm Spec}^{\log}_{U'}({\cal O}_{U'}/{\cal K}')$.  
Assume that $\os{\circ}{U}{}'=\os{\circ}{U}$.  
Let $(U',{\cal K}',\eps')\lo (U,{\cal K},\eps)$ 
be a morphism of log PD-schemes 
such that the underlying morphism $\os{\circ}{U}{}'\lo \os{\circ}{U}$ of schemes 
is equal to ${\rm id}_{\os{\circ}{U}}$. 
Let $Z_{\bul \leq N}$ be a log smooth integral $N$-truncated simplicial log scheme 
over $U_0$. 
Set $Z'_{\bul \leq N}:=Z_{\bul \leq N}\times_{U_0}U'_0$ 
with projection $q\col Z'_{\bul \leq N}\lo Z_{\bul \leq N}$. 
Identify $(Z'_{\bul \leq N})_{\rm zar}$ with $(Z_{\bul \leq N})_{\rm zar}$ 
$($Because $Z\lo U_0$ is integral, 
$\os{\circ}{Z}{}'=\os{\circ}{Z}\times_{\os{\circ}{U}_0}\os{\circ}{U}{}'_0
=\os{\circ}{Z}.)$. 
Let $F^{\bul \leq N}$ be a flat quasi-coherent ${\cal O}_{Z_{\bul \leq N}/U}$-module. 
Then the natural morphism 
\begin{align*} 
Ru_{Z_{\bul \leq N}/U*}(F^{\bul \leq N})\lo 
Ru_{Z'_{\bul \leq N}/U'*}(q^*_{\rm crys}(F^{\bul \leq N}))
\tag{1.4.9.1}\label{ali:nmzqf} 
\end{align*} 
is the identity. 
\par 
$(2)$ Let $(U'_1,{\cal K}'_1,\eps'_1)\lo (U_1,{\cal K}_1,\eps_1)$ 
be a similar morphism to the morphism 
$(U',{\cal K}',\eps')\lo (U,{\cal K},\eps)$ in {\rm (1)} 
fitting into the following commutative diagram 
\begin{equation*} 
\begin{CD} 
(U'_1,{\cal K}'_1,\eps'_1) @>>> (U_1,{\cal K}_1,\eps_1) \\ 
@VVV @VV{}V \\ 
(U',{\cal K}',\eps') @>>> (U,{\cal K},\eps).   
\end{CD} 
\end{equation*}
Set $U_{10}:=\ul{\rm Spec}^{\log}_{U_1}({\cal O}_{U_1}/{\cal K}_1)$ 
and $U'_{10}:=\ul{\rm Spec}^{\log}_{U'_1}({\cal O}_{U'_1}/{\cal K}'_1)$. 
Let $W_{\bul \leq N}/U_{10}$ $($resp.~$W'_{\bul \leq N}/U'_{10})$ 
be a similar log scheme to $Z_{\bul \leq N}/U_0$ 
$($resp.~$Z'_{\bul \leq N}/U'_0)$.  
Let $g\col W_{\bul \leq N}\lo Z_{\bul \leq N}$ be 
a morphism over $U_{10}\lo U_0$. 
Let $g'\col W_{\bul \leq N}\lo Z_{\bul \leq N}$ be 
the base change morphism of $g$ over $U'_{10}\lo U'_0$. 
Let $G^{\bul \leq N}$ be a similar log crystal of ${\cal O}_{W_{\bul \leq N}/U_1}$-modules 
to $F^{\bul \leq N}$ and let 
$F^{\bul \leq N}\lo g_{{\rm crys}*}(G^{\bul \leq N})$ be a morphism of  
${\cal O}_{Z_{\bul \leq N}/U}$-modules.  
Let $r\col W'_{\bul \leq N}\lo W_{\bul \leq N}$ be the projection. 
Then the following diagram is commutative$:$ 
\begin{equation*} 
\begin{CD} 
Ru_{Z_{\bul \leq N}/U*}(F^{\bul \leq N})@= 
Ru_{Z'_{\bul \leq N}/U'*}(q^*_{\rm crys}(F^{\bul \leq N}))\\ 
@VVV @VVV \\
Rg_*Ru_{W_{\bul \leq N}/U_1*}(G^{\bul \leq N})@= 
Rg'_*Ru_{W'_{\bul \leq N}/U'_1*}(r^*_{\rm crys}(G^{\bul \leq N})). 
\end{CD} 
\tag{1.4.9.2}\label{cd:wwggr} 
\end{equation*}  
\end{lemm}
\begin{proof} 
(1): As stated above, 
we identify $(Z'_{\bul \leq N})_{\rm zar}$ with $(Z_{\bul \leq N})_{\rm zar}$. 
The natural morphism 
$F^{\bul \leq N}\lo q_{{\rm crys}*}q^*_{\rm crys}(F^{\bul \leq N})$ 
induces the morphism (\ref{ali:nmzqf}).  
We may assume that $N=0$. Set $Z:=Z_0$ and $F=F^0$. 
The question is local. Hence we may assume that 
$Z$ has a log smooth integral lift ${\cal Z}/U$ by  
the proofs of \cite[(3.14)]{klog1} and \cite[(2.3.14)]{nh2} 
(cf.~the proof of \cite[(4.7)]{ny}). 
Set ${\cal Z}_{U'}:={\cal Z}\times_UU'$ and 
let $Q\col {\cal Z}_{U'}\lo {\cal Z}$ be the projection. 
Because ${\cal Z}\lo U$ is also integral,   
$({\cal Z}_{U'})^{\circ}=\os{\circ}{\cal Z}\times_{\os{\circ}{U}}\os{\circ}{U}{}'=\os{\circ}{\cal Z}$. 
Let $({\cal F},{\nabla})$ be a quasi-coherent ${\cal O}_{\cal Z}$-module 
with integrable connection corresponding to $F$.
Then $Q^*({\cal F})={\cal F}$ and 
$\Om^{\bul}_{{\cal Z}'/U'}=Q^*(\Om^{\bul}_{{\cal Z}/U})=\Om^{\bul}_{{\cal Z}/U}$ 
by the formula in the end of \cite[(1.7)]{klog1}. 
By the log Poincar\'{e} lemma, we obtain the following equality 
\begin{align*} 
Ru_{Z'/U'*}(q^*_{\rm crys}(F))=
Q^*({\cal F})\otimes_{{\cal O}_{{\cal Z}'}}\Om^{\bul}_{{\cal Z}'/U'}
={\cal F}\otimes_{{\cal O}_{\cal Z}}\Om^{\bul}_{{\cal Z}/U}=
Ru_{Z/U*}(F). 
\end{align*} 
(2): (2) is obvious. 
\end{proof} 

\begin{coro}\label{coro:oozu}
Assume that $T$ is restrictively hollow with respective to the morphism $T_0\lo S$. 
Then there exists the following isomorphism 
\begin{align*} 
\theta:=\theta_{X_{\bul \leq N,T_0}/T}\col Ru_{X_{\bul \leq N,T_0}/T*}
(\eps^*_{X_{\bul \leq N,T_0/T}}(E^{\bul \leq N})) 
\os{\sim}{\lo}
A_{\rm zar}(X_{\bul \leq N,\os{\circ}{T}_0}/S(T),E^{\bul \leq N})
\tag{1.4.10.1}\label{eqn:utz} 
\end{align*} 
in ${\rm D}^+(
f^{-1}({\cal O}_T))$.  
\end{coro} 
\begin{proof} 
Because $S(T)^{\nat}=S(T)$, 
we have the following isomorphism by (\ref{ali:nmzqf}) and (\ref{eqn:uz}):  
\begin{align*} 
Ru_{X_{\bul \leq N,T_0}/T*}
(\eps^*_{X_{\bul \leq N,T_0/T}}(E^{\bul \leq N})) &=
Ru_{X_{\bul \leq N,\os{\circ}{T}_0}/S(T)*}
(\eps^*_{X_{\bul \leq N,\os{\circ}{T}_0/S(T)}}(E^{\bul \leq N})) 
\tag{1.4.10.2}\label{eqn:uctsz} \\
&\os{\theta \wedge,\simeq}{\lo} 
A_{\rm zar}(X_{\bul \leq N,\os{\circ}{T}_0}/S(T),E^{\bul \leq N}).
\end{align*}  
\end{proof}

\begin{coro}[{\bf Invariance of the Poincar\'{e} filtration I 
with respect to log structures}]\label{coro:nids} 
Let $q$ be a nonnegative integer. 
Let $(T',{\cal J}',\del')\lo (T,{\cal J},\del)$ be 
a morphism from another log PD-enlargement of $S$. 
Set $T'_0:=\ul{\rm Spec}_{T'}^{\log}({\cal O}_{T'}/{\cal J}')$.   
Assume that $S(T)=S(T')$ and that 
$T$ $($and hence $T')$ is restrictively hollow with respective to the morphism $T_0\lo S$.  
Let $P$ and $P'$ be the induced filtrations on 
$$R^qf_{X_{\bul \leq N,T_0}/T*}(\eps^*_{X_{\bul \leq N,T_0}/T}(E^{\bul \leq N}))
\quad 
{\rm and} \quad  
R^qf_{X_{\bul \leq N,T'_0}/T'*}(\eps^*_{X_{\bul \leq N,T'_0}/T'}(E^{\bul \leq N}))$$ 
by 
$$(A_{\rm zar}(X_{\bul \leq N,\os{\circ}{T}_0}/S(T),E^{\bul \leq N}),P)
\quad 
{\rm and} \quad 
(A_{\rm zar}(X_{\bul \leq N,\os{\circ}{T}{}'_0}/S(T'),E^{\bul \leq N}),P)$$ 
and the isomorphism {\rm (\ref{eqn:utz})}, respectively. 
Then $P=P'$. 
\end{coro} 
\begin{proof} 
(\ref{coro:nids}) immediately follows from the isomorphism (\ref{eqn:utz}) and 
the tautological equality 
$(A_{\rm zar}(X_{\bul \leq N,\os{\circ}{T}_0}/S(T),E^{\bul \leq N}),P)= 
(A_{\rm zar}(X_{\bul \leq N,\os{\circ}{T}{}'_0}/S(T'),E^{\bul \leq N}),P)$. 
\end{proof}

\section{Contravariant functoriality of 
zariskian $p$-adic filtered Steenbrink complexes}\label{sec:fcuc} 
Let the notations be as in the previous section. 
\par 
Let $S'$ be another family of log points.   
Let $(T',{\cal J}',\del')$ be a log PD-enlargement over $S'$. 
Assume that $p$ is locally nilpotent on $\os{\circ}{T}{}'$. 
Let $u\col (S(T)^{\nat},{\cal J},\del) \lo (S'(T')^{\nat},{\cal J}',\del')$ be 
a morphism of fine log schemes. 
Set $T_0:=\ul{\rm Spec}^{\log}_T({\cal O}_T/{\cal J})$ and 
$T'_0:=\ul{\rm Spec}^{\log}_{T'}({\cal O}_{T'}/{\cal J}')$.   
By the definition of $\deg(u)_x$ ((\ref{defi:ddef})), 
we have the following equality: 
\begin{equation*} 
u_{x}^*(\theta_{S'(T')^{\nat},\os{\circ}{u}(x)})= \deg(u)_x\theta_{S(T)^{\nat},x} 
\quad (x\in \os{\circ}{T}). 
\tag{1.5.0.1}\label{eqn:uta}  
\end{equation*}
By (\ref{prop:oot}) (1), $\deg(u)_x\not=0$ for any point $x\in \os{\circ}{T}$. 
Let $X_{\bul \leq N}$ and $Y_{\bul \leq N}$ be 
$N$-truncated simplicial SNCL schemes 
over $S_{\os{\circ}{T}_0}$ and $S'_{\os{\circ}{T}{}'_0}$, respectively. 
Let ${\mathfrak D}(\ol{S'(T')^{\nat}})$ be the log PD-envelope of 
the immersion $S'(T')^{\nat}\os{\sus}{\lo} \ol{S'(T')^{\nat}}$ 
over $(\os{\circ}{T}{}',{\cal J}',\del')$. 
Let 
\begin{equation*} 
\begin{CD} 
X_{\bul \leq N,\os{\circ}{T}_0} 
@>{g_{\bul \leq N}}>> Y_{\bul \leq N,\os{\circ}{T}{}'_0}\\
@VVV @VVV \\ 
S_{\os{\circ}{T}_0} @>>> S'_{\os{\circ}{T}{}'_0} \\ 
@V{\bigcap}VV @VV{\bigcap}V \\ 
S(T)^{\nat} @>{u}>> S'(T')^{\nat}
\end{CD}
\tag{1.5.0.2}\label{eqn:xdxduss}
\end{equation*} 
be a commutative diagram of 
$N$-truncated simplicial SNCL schemes 
over $S_{\os{\circ}{T}_0}$ and $S'_{\os{\circ}{T}{}'_0}$ 
such that $X_{\bul \leq N,\os{\circ}{T}_0}$ and $Y_{\bul \leq N,\os{\circ}{T}{}'_0}$ 
have affine $N$-truncated simplicial open coverings 
$X'_{\bul \leq N,\os{\circ}{T}_0}$ and $Y'_{\bul \leq N,\os{\circ}{T}{}'_0}$ 
of $X_{\bul \leq N,\os{\circ}{T}_0}$ and $Y_{\bul \leq N,\os{\circ}{T}{}'_0}$, 
respectively, 
fitting into the following commutative diagram 
\begin{equation*} 
\begin{CD} 
X'_{\bul \leq N,\os{\circ}{T}_0} @>{g'_{\bul \leq N}}>> Y'_{\bul \leq N,\os{\circ}{T}{}'_0} \\
@VVV @VVV \\ 
X_{\bul \leq N,\os{\circ}{T}_0} @>{g_{\bul \leq N}}>> Y_{\bul \leq N,\os{\circ}{T}{}'_0}. 
\end{CD}
\tag{1.5.0.3}\label{cd:xygxy}
\end{equation*} 
(If $X_{\bul \leq N,\os{\circ}{T}_0}$ and $Y_{\bul \leq N,\os{\circ}{T}{}'_0}$ are split, 
then this condition is satisfied ((\ref{prop:lissmp}) (1)).
We do not assume that $X_{\bul \leq N}$ and 
$Y_{\bul \leq N}$ have affine $N$-truncated simplicial open coverings of 
$X_{\bul \leq N}$ and $Y_{\bul \leq N}$, respectively.) 
\par 
Let $X_{\bul \leq N,\bul,\os{\circ}{T}_0} \os{\sus}{\lo} \ol{\cal P}{}'_{\bul \leq N,\bul}$
and 
$Y_{\bul \leq N,\bul,\os{\circ}{T}{}'_0} \os{\sus}{\lo} \ol{\cal Q}_{\bul \leq N,\bul}$ 
be immersions into $(N,\infty)$-truncated bisimplicial 
log smooth schemes over $\ol{S(T)^{\nat}}$ and $\ol{S'(T')^{\nat}}$, respectively. 
Indeed, these immersions exist by (\ref{eqn:eipxd}). 
Set 
$$\ol{\cal P}_{\bul \leq N,\bul}:=
\ol{\cal P}{}'_{\bul \leq N,\bul}\times_{\ol{S(T)^{\nat}}}
(\ol{\cal Q}_{\bul \leq N,\bul}\times_{\ol{S'(T')^{\nat}}}\ol{S(T)^{\nat}})
=\ol{\cal P}{}'_{\bul \leq N,\bul}\times_{\ol{S'(T')^{\nat}}}
\ol{\cal Q}_{\bul \leq N,\bul}.$$ 
Let 
$\ol{g}_{\bul \leq N,\bul}\col \ol{\cal P}_{\bul \leq N,\bul} 
\lo \ol{\cal Q}_{\bul \leq N,\bul}$
be the second projection. 
Then, as explained in the Introduction or (\ref{rema:uuc}),  
we have the following commutative diagram 
\begin{equation*} 
\begin{CD} 
X_{\bul \leq N,\bul,\os{\circ}{T}_0}
@>{\sus}>> \ol{\cal P}_{\bul \leq N,\bul}\\
@V{g_{\bul \leq N,\bul}}VV 
@VV{\ol{g}_{\bul \leq N,\bul}}V \\ 
Y_{\bul \leq N,\bul,\os{\circ}{T}{}'_0} @>{\sus}>> \ol{\cal Q}_{\bul \leq N,\bul}
\end{CD} 
\tag{1.5.0.4}\label{cd:xpnxp} 
\end{equation*} 
over 
\begin{equation*} 
\begin{CD} 
S_{\os{\circ}{T}_0} @>{\subset}>> \ol{S(T)^{\nat}}\\ 
@VVV @VVV \\ 
S'_{\os{\circ}{T}{}'_0} @>{\subset}>> \ol{S'(T')^{\nat}}. 
\end{CD}
\end{equation*} 
Let $\ol{\mathfrak E}_{\bul \leq N,\bul}$ be the log PD-envelope of  
the immersion 
$Y_{\bul \leq N,\bul,\os{\circ}{T}{}'_0}
\os{\sus}{\lo}\ol{\cal Q}_{\bul \leq N,\bul}$ over 
$(\os{\circ}{T}{}',{\cal J}',\del')$.     
Set ${\mathfrak E}_{\bul \leq N,\bul}:=
\ol{\mathfrak E}_{\bul \leq N,\bul}\times_{{\mathfrak D}(\ol{S'(T')^{\nat}})}S'(T')^{\nat}$. 
By (\ref{cd:xpnxp}) 
we have the following natural morphism  
\begin{equation*} 
\ol{g}{}^{\rm PD}_{\bul \leq N,\bul}\col 
\ol{\mathfrak D}_{\bul \leq N,\bul}\lo \ol{\mathfrak E}_{\bul \leq N,\bul}.   
\tag{1.5.0.5}\label{eqn:gpdf} 
\end{equation*}  
Hence we have the following natural morphism 
\begin{equation*} 
g^{\rm PD}_{\bul \leq N,\bul}\col 
{\mathfrak D}_{\bul \leq N,\bul}\lo {\mathfrak E}_{\bul \leq N,\bul}.   
\tag{1.5.0.6}\label{eqn:gpndf} 
\end{equation*}

\begin{rema}\label{rema:xnad} 
As explained in the Introduction, another argument using 
(\ref{eqn:eipxd}) and (\ref{lemm:obvc}) tells us that 
there exists the following commutative diagram 
\begin{equation*} 
\begin{CD} 
X_{\bul \leq N,\bul,\os{\circ}{T}_0} @>{\sus}>> \ol{\mathfrak D}{}'_{\bul \leq N,\bul}\\
@V{g_{\bul \leq N,\bul}}VV 
@VV{\ol{g}_{\bul \leq N,\bul}}V \\ 
Y_{\bul \leq N,\bul,\os{\circ}{T}{}'_0} 
@>{\sus}>> \ol{\cal Q}{}_{\bul \leq N,\bul}
\end{CD} 
\tag{1.5.1.1}\label{cd:axpnxpd} 
\end{equation*} 
over 
\begin{equation*} 
\begin{CD} 
S_{\os{\circ}{T}} @>{\subset}>> \ol{S(T)^{\nat}}\\ 
@VVV @VVV \\ 
S'_{\os{\circ}{T}{}'} @>{\subset}>> \ol{S'(T')^{\nat}},  
\end{CD}
\end{equation*} 
where 
$\ol{\mathfrak D}{}'_{\bul \leq N,\bul}$ is the log PD-envelope of  
$X_{\bul \leq N,\bul,\os{\circ}{T}_0} \os{\sus}{\lo}\ol{\cal P}{}'_{\bul \leq N,\bul}$ over 
$(\os{\circ}{T},{\cal J},\del)$.    
\end{rema}

Let $E^{\bul \leq N}$ and  $F^{\bul \leq N}$ be flat quasi-coherent crystals of 
${\cal O}_{\os{\circ}{X}_{\bul \leq N,T_0}/\os{\circ}{T}}$-modules  
and ${\cal O}_{\os{\circ}{Y}_{\bul \leq N,T_0'}/\os{\circ}{T}{}'}$-modules, respectively.     
Let 
\begin{align*} 
\os{\circ}{g}{}^*_{\bul \leq N,{\rm crys}}(F^{\bul \leq N})
\lo 
E^{\bul \leq N}
\tag{1.5.1.2}\label{ali:gnfe} 
\end{align*} 
be a morphism of 
${\cal O}_{\os{\circ}{X}_{\bul \leq N,T_0}/\os{\circ}{T}}$-modules. 

\begin{theo}[{\bf Contravariant functoriality I of $A_{\rm zar}$}]\label{theo:funas} 
$(1)$ Assume that 
$\deg(u)_x$ is not divisible by $p$ for any point $x \in \os{\circ}{T}$. 
Then $g_{\bul \leq N}\col X_{\bul \leq N,\os{\circ}{T}{}_0}
\lo Y_{\bul \leq N,\os{\circ}{T}{}'_0}$ induces the following 
well-defined pull-back morphism 
\begin{equation*}  
g_{\bul \leq N}^* \col 
(A_{\rm zar}(Y_{\bul \leq N,\os{\circ}{T}{}'_0}/S'(T')^{\nat},F^{\bul \leq N}),P)
\lo Rg_{\bul \leq N*}((A_{\rm zar}(X_{\bul \leq N,\os{\circ}{T}_0}/S(T)^{\nat}
,E^{\bul \leq N}),P)) 
\tag{1.5.2.1}\label{eqn:fzaxd}
\end{equation*} 
fitting into the following commutative diagram$:$
\begin{equation*} 
\begin{CD}
A_{\rm zar}(Y_{\bul \leq N,\os{\circ}{T}{}'_0}/S'(T')^{\nat},F^{\bul \leq N})
@>{g_{\bul \leq N}^*}>>  \\ 
@A{\theta_{Y_{\bul \leq N,\os{\circ}{T}{}'_0}/S'(T')^{\nat}}\wedge}A{\simeq}A 
\\
Ru_{Y_{\bul \leq N,\os{\circ}{T}{}'_0}/S'(T')^{\nat}*}
(\eps^*_{Y_{\bul \leq N,\os{\circ}{T}{}'_0}/S'(T')^{\nat}}(F^{\bul \leq N}))
@>{g_{\bul \leq N}^*}>>
\end{CD}
\tag{1.5.2.2}\label{cd:pssfccz} 
\end{equation*}
\begin{equation*} 
\begin{CD}
Rg_{\bul \leq N*}(A_{\rm zar}(X_{\bul \leq N,\os{\circ}{T}_0}/S(T)^{\nat}
,E^{\bul \leq N}))\\ 
@A{Rg_{\bul \leq N*}(\theta_{X_{\bul \leq N,\os{\circ}{T}_0}/S(T)^{\nat}} \wedge)}A{\simeq}A\\
Rg_{\bul \leq N*}Ru_{X_{\bul \leq N,\os{\circ}{T}_0}/S(T)^{\nat}*}
(\eps^*_{X_{\bul \leq N,\os{\circ}{T}_0}/S(T)^{\nat}}(E^{\bul \leq N})).
\end{CD} 
\end{equation*}
\par 
$(2)$ Let $S''$ be a family of log points. 
Let $(T'',{\cal J}'',\del'')$ be a log PD-enlargement of $S''$. 
Set  
$T''_0:=\ul{\rm Spec}^{\log}_{T''}({\cal O}_{T''}/{\cal J}'')$. 
Let $v\col (S'(T')^{\nat},{\cal J},\del) \lo (S''(T'')^{\nat},{\cal J}',\del')$ 
and $h_{\bul \leq N}\col Y_{\bul \leq N}\lo Z_{\bul \leq N}$ 
be similar morphisms to $u$ and $g_{\bul \leq N}$, respectively. 
Assume that $\deg(v)_x$ is not divisible by $p$ for any point  $x \in \os{\circ}{T}{}'$. 
Let $G^{\bul \leq N}$ be a flat quasi-coherent crystal of 
${\cal O}_{\os{\circ}{Z}_{\bul \leq N,T''_0}/\os{\circ}{T}{}''}$-modules.    
Let 
\begin{align*} 
\os{\circ}{h}{}^*_{\bul \leq N,{\rm crys}}(G^{\bul \leq N})
\lo 
F^{\bul \leq N}
\tag{1.5.2.3}\label{ali:gznfe} 
\end{align*} 
be a morphism of 
${\cal O}_{\os{\circ}{Y}_{\bul \leq N,T'_0}/\os{\circ}{T}{}'}$-modules. 
Then 
\begin{align*} 
(h_{\bul \leq N}\circ g_{\bul \leq N})^*  =&
Rh_{\bul \leq N*}(g_{\bul \leq N}^*)\circ h_{\bul \leq N}^*    
\col (A_{\rm zar}(Z_{\bul \leq N,\os{\circ}{T}{}''_0}/S''(T'')^{\nat},F^{\bul \leq N}),P)
\tag{1.5.2.4}\label{ali:pdpp} \\ 
& \lo Rh_{\bul \leq N*}Rg_{\bul \leq N*}
(A_{\rm zar}(X_{\bul \leq N,\os{\circ}{T}_0}/S(T)^{\nat},E^{\bul \leq N}),P) \\
& =R(h_{\bul \leq N}\circ g_{\bul \leq N})_*
(A_{\rm zar}(X_{\bul \leq N,\os{\circ}{T}_0}/S(T)^{\nat},E^{\bul \leq N}),P).
\end{align*}  
\par 
$(3)$  
\begin{equation*} 
{\rm id}_{X_{\bul \leq N,\os{\circ}{T}_0}}^*={\rm id} 
\col (A_{\rm zar}(X_{\bul \leq N,\os{\circ}{T}_0}/S(T)^{\nat},E^{\bul \leq N}),P)
\lo (A_{\rm zar}(X_{\bul \leq N,\os{\circ}{T}_0}/S(T)^{\nat},E^{\bul \leq N}),P).  
\tag{1.5.2.5}\label{eqn:fzidd}
\end{equation*} 
\par
$(4)$ Assume that $T$ and $T'$ are restrictively hollow 
with respective to the morphism $T_0\lo S$ and 
$T'_0\lo S'$, respectively. 
Set $X_{\bul \leq N,T_0}:=X_{\bul \leq N}\times_{S}T_0$ and 
$Y_{\bul \leq N,T'_0}:=Y_{\bul \leq N}\times_{S'}T'_0$. 
Let $g_{\bul \leq N,T'T}\col X_{\bul \leq N,T_0} \lo Y_{\bul \leq N,T'_0}$ 
be the base change morphism of $g_{\bul \leq N}$. 
Then the morphism 
{\rm (\ref{eqn:fzaxd})} fits into the following commutative diagram$:$ 
\begin{equation*} 
\begin{CD}
A_{\rm zar}(Y_{\bul \leq N,\os{\circ}{T}{}'_0}/S'(T'),F^{\bul \leq N})
@>{g_{\bul \leq N}^*}>>  
Rg_{\bul \leq N*}(A_{\rm zar}(X_{\bul \leq N,\os{\circ}{T}_0}/S(T)
,E^{\bul \leq N}))\\ 
@A{\theta_{Y_{\bul \leq N,T'_0}/T'}}A{\simeq}A 
@A{Rg_{\bul \leq N*}(\theta_{X_{\bul \leq N,T_0}/T}})A{\simeq}A\\
Ru_{Y_{\bul \leq N,T'_0}/T'*}
(\eps^*_{Y_{\bul \leq N,T'_0}/T'}(F^{\bul \leq N}))
@>{g_{\bul \leq N,T'T}^*}>>Rg_{\bul \leq N*}Ru_{X_{\bul \leq N,T_0}/T*}
(\eps^*_{X_{\bul \leq N,T_0}/T}(E^{\bul \leq N})).
\end{CD}
\tag{1.5.2.6}\label{cd:tmfccz} 
\end{equation*}
Here the morphisms $\theta_{X_{\bul \leq N,T_0}/T}$ and 
$\theta_{Y_{\bul \leq N,T'_0}/T'}$ are the morphisms defined in 
{\rm (\ref{eqn:utz})}.  
\end{theo}
\begin{proof} 
(1): For $F^{\bul \leq N}$, let $(\ol{\cal F}{}^{\bul \leq N,\bul},\nabla)$ 
and $({\cal F}{}^{\bul \leq N,\bul},\nabla)$ 
be the similar objects to 
$(\ol{\cal E}{}^{\bul \leq N,\bul},\nabla)$ and 
$({\cal E}{}^{\bul \leq N,\bul},\nabla)$ in \S\ref{sec:psc}, respectively.  
For simplicity of notation, 
denote $\theta_{{\cal P}^{\rm ex}_{\bul \leq N,\bul}}\in 
{\cal O}_{{\mathfrak D}_{\bul \leq N,\bul}}
\otimes_{{\cal O}_{{\cal P}^{\rm ex}_{\bul \leq N,\bul}}}
{\Om}^1_{{\cal P}^{\rm ex}_{\bul \leq N,\bul}/\os{\circ}{T}}$ 
(resp.~$\theta_{{\cal Q}^{\rm ex}_{\bul \leq N,\bul}}\in 
{\cal O}_{{\mathfrak E}_{\bul \leq N,\bul}}
\otimes_{{\cal O}_{{\cal Q}^{\rm ex}_{\bul \leq N,\bul}}}
{\Om}^1_{{\cal Q}^{\rm ex}_{\bul \leq N,\bul}/\os{\circ}{T}{}'}$) 
simply by $\theta$ (resp.~$\theta'$). 
\par 
First we would like to construct a morphism  
\begin{equation*} 
{\cal F}^{\bul \leq N,\bul}
\otimes_{{\cal O}_{{\cal Q}{}^{\rm ex}_{\bul \leq N,\bul}}}
\Om^{\bul}_{{\cal Q}{}^{\rm ex}_{\bul \leq N,\bul}/\os{\circ}{T}{}'}
\lo 
g^{\rm PD}_{\bul \leq N,\bul*}
({\cal E}^{\bul \leq N,\bul}
\otimes_{{\cal O}_{{\cal P}^{\rm ex}_{\bul \leq N,\bul}}}
\Om^{\bul}_{{\cal P}^{\rm ex}_{\bul \leq N,\bul}/\os{\circ}{T}})
\end{equation*} 
of complexes. 
Because we are given the morphism  
$\os{\circ}{g}{}^*_{\bul \leq N,{\rm crys}}(F^{\bul \leq N})
\lo E^{\bul \leq N}$, we have a morphism 
\begin{align*} 
\ol{g}{}^{{\rm PD}*}_{\bul \leq N,\bul}
\col 
\ol{\cal F}{}^{\bul \leq N,\bul}
\lo 
\ol{g}{}^{{\rm PD}*}_{\bul \leq N,\bul}(\ol{\cal E}{}^{\bul \leq N,\bul})
\end{align*} 
fitting into the following commutative diagram$:$ 
\begin{equation*} 
\begin{CD}  
\ol{\cal F}{}^{\bul \leq N,\bul}
@>{\ol{g}{}^{{\rm PD}*}_{\bul \leq N,\bul}}>>
\ol{g}{}^{\rm PD}_{\bul \leq N,\bul*}(\ol{\cal E}{}^{\bul \leq N,\bul})\\
@VVV @VVV\\
\ol{\cal F}{}^{\bul \leq N,\bul}
\otimes_{{\cal O}_{\ol{\cal Q}{}^{\rm ex}_{\bul \leq N,\bul}}}
\Om^1_{\ol{\cal Q}{}^{\rm ex}_{\bul \leq N,\bul}/\os{\circ}{T}{}'}
@>{\ol{g}{}^{{\rm PD}*}_{\bul \leq N,\bul}}>>
\ol{g}{}^{\rm PD}_{\bul \leq N,\bul*}
(\ol{\cal E}^{\bul \leq N,\bul}
\otimes_{{\cal O}_{\ol{\cal P}{}^{\rm ex}_{\bul \leq N,\bul}}}
\Om^1_{\ol{\cal P}{}^{\rm ex}_{\bul \leq N,\bul}/\os{\circ}{T}})
\end{CD} 
\tag{1.5.2.7}\label{eqn:dolm}
\end{equation*} 
by using (\ref{cd:xpnxp}). Hence we have a morphism 
\begin{align*} 
g^{{\rm PD}*}_{\bul \leq N,\bul}\col {\cal F}^{\bul \leq N,\bul}
\lo g^{\rm PD}_{{\bul \leq N,\bul}*}({\cal E}^{\bul \leq N,\bul})
\tag{1.5.2.8}\label{eqn:defm}
\end{align*} 
fitting into the following commutative diagram$:$ 
\begin{equation*} 
\begin{CD}  
{\cal F}^{\bul \leq N,\bul}
@>{g^{{\rm PD}*}_{\bul \leq N,\bul}}>>
g^{\rm PD}_{\bul \leq N,\bul*}({\cal E}^{\bul \leq N,\bul})\\
@VVV @VVV\\
{\cal F}^{\bul \leq N,\bul}
\otimes_{{\cal O}_{{\cal Q}{}^{\rm ex}_{\bul \leq N,\bul}}}
\Om^1_{{\cal Q}{}^{\rm ex}_{\bul \leq N,\bul}/\os{\circ}{T}{}'}
@>{g^{{\rm PD}*}_{\bul \leq N,\bul}}>>
g^{\rm PD}_{\bul \leq N,\bul*}
({\cal E}^{\bul \leq N,\bul}
\otimes_{{\cal O}_{{\cal P}{}^{\rm ex}_{\bul \leq N,\bul}}}
\Om^1_{{\cal P}{}^{\rm ex}_{\bul \leq N,\bul}/\os{\circ}{T}}).
\end{CD} 
\tag{1.5.2.9}\label{eqn:dgm}
\end{equation*} 
By using (\ref{cd:xpnxp}) again, 
we have the following morphism 
\begin{equation*} 
g^{{\rm PD}*}_{\bul \leq N,\bul}
\col {\cal O}_{{\mathfrak E}_{\bul \leq N,\bul}}
\otimes_{{\cal O}_{{\cal Q}{}^{\rm ex}_{\bul \leq N,\bul}}}
\Om^{\bul}_{{\cal Q}{}^{\rm ex}_{\bul \leq N,\bul}/\os{\circ}{T}{}'}
\lo 
g^{\rm PD}_{\bul \leq N,\bul*}({\cal O}_{{\mathfrak D}_{\bul \leq N,\bul}}
\otimes_{{\cal O}_{{\cal P}^{\rm ex}_{\bul \leq N,\bul}}}
\Om^{\bul}_{{\cal P}^{\rm ex}_{\bul \leq N,\bul}/\os{\circ}{T}}). 
\tag{1.5.2.10}\label{eqn:dgppm}
\end{equation*}  
Set 
$$g^{{\rm PD}*}_{\bul \leq N,\bul}(e\otimes \om):=
g^{{\rm PD}*}_{\bul \leq N,\bul}(e)\otimes g^{{\rm PD}*}_{\bul \leq N,\bul}(\om)
\quad (e\in {\cal F}^{\bul \leq N,\bul},~ 
\om \in \Om^i_{{\cal Q}^{\rm ex}_{\bul \leq N,\bul}/\os{\circ}{T}{}'}~(i\in {\mab N})).$$ 
This $g^{{\rm PD}*}_{\bul \leq N,\bul}$ induces 
the following morphism of complexes: 
\begin{equation*} 
{\cal F}^{\bul \leq N,\bul}
\otimes_{{\cal O}_{{\cal Q}^{\rm ex}_{\bul \leq N,\bul}}}
\Om^{\bul}_{{\cal Q}^{\rm ex}_{\bul \leq N,\bul}/\os{\circ}{T}{}'}
\lo 
g^{\rm PD}_{\bul \leq N,\bul*} 
({\cal E}^{\bul \leq N,\bul}
\otimes_{{\cal O}_{{\cal P}^{\rm ex}_{\bul \leq N,\bul}}}
\Om^{\bul}_{{\cal P}^{\rm ex}_{\bul \leq N,\bul}/\os{\circ}{T}}).  
\tag{1.5.2.11}\label{eqn:dogm}
\end{equation*} 
Because the following diagram
\begin{equation*} 
\begin{CD} 
g^{{\rm PD}*}_{\bul \leq N,\bul} 
({\cal F}^{\bul \leq N,\bul}
\otimes_{{\cal O}_{{\cal Q}^{\rm ex}_{\bul \leq N,\bul}}}
\Om^{\bul}_{{\cal Q}^{\rm ex}_{\bul \leq N,\bul}/\os{\circ}{T}{}'})[1] 
@>{g^{{\rm PD}*}_{\bul \leq N,\bul}}>> 
{\cal E}^{\bul \leq N,\bul}
\otimes_{{\cal O}_{{\cal P}^{\rm ex}_{\bul \leq N,\bul}}}
{\Om}^{\bul}_{{\cal P}^{\rm ex}_{\bul \leq N,\bul}/\os{\circ}{T}}[1] \\   
@A{g^{{\rm PD}*}_{\bul \leq N,\bul} 
(\deg(u)^{-1}\theta'\wedge)}AA @AA{\theta \wedge}A \\ 
g^{{\rm PD}*}_{\bul \leq N,\bul}
({\cal F}^{\bul \leq N,\bul}
\otimes_{{\cal O}_{{\cal Q}^{\rm ex}_{\bul \leq N,\bul}}}
{\Om}^{\bul}_{{\cal Q}^{\rm ex}_{\bul \leq N,\bul}/\os{\circ}{T}{}'})
@>{g^{{\rm PD}*}_{\bul \leq N,\bul}}>> 
{\cal E}^{\bul \leq N,\bul}
\otimes_{{\cal O}_{{\cal P}^{\rm ex}_{\bul \leq N,\bul}}}
{\Om}^{\bul}_{{\cal P}^{\rm ex}_{\bul \leq N,\bul}/\os{\circ}{T}}    
\end{CD} 
\tag{1.5.2.12}\label{cd:fblssod}
\end{equation*} 
is commutative ($\deg(u)^{-1}$ has a meaning by 
the assumption of the $p$-nondivisibility of $\deg(u)$) 
and because 
we have the following commutative diagram for 
$i,j\in {\mab N}$ 
\begin{equation*} 
\begin{CD}
g^{\rm PD}_{\bul \leq N,\bul*}
(P_j({\cal E}^{\bul \leq N,\bul}
\otimes_{{\cal O}_{{\cal P}^{\rm ex}_{\bul \leq N,\bul}}}
{\Om}^{i+j+1}_{{\cal P}^{\rm ex}_{\bul \leq N,\bul}/\os{\circ}{T}{}}))
@>{\subset}>>  \\  
@A{({\rm deg}(u))^{-(j+1)}
g^{{\rm PD}*}_{\bul \leq N,\bul}}AA \\ 
P_j({\cal F}^{\bul \leq N,\bul}
\otimes_{{\cal O}_{{\cal Q}^{\rm ex}_{\bul \leq N,\bul}}}
{\Om}^{i+j+1}_{{\cal Q}^{\rm ex}_{\bul \leq N,\bul}/\os{\circ}{T}{}'})
@>{\subset}>> 
\end{CD}
\tag{1.5.2.13}\label{cd:axfpdwt}
\end{equation*} 
\begin{equation*} 
\begin{CD}
g^{\rm PD}_{\bul \leq N,\bul*}
({\cal E}^{\bul \leq N,\bul}
\otimes_{{\cal O}_{{\cal P}^{\rm ex}_{\bul \leq N,\bul}}}
{\Om}^{i+j+1}_{{\cal P}^{\rm ex}_{\bul \leq N,\bul}/\os{\circ}{T}{}}) \\
@AA{({\rm deg}(u))^{-(j+1)}
g^{{\rm PD}*}_{\bul \leq N,\bul}}A \\
{\cal F}^{\bul \leq N,\bul}
\otimes_{{\cal O}_{{\cal Q}^{\rm ex}_{\bul \leq N,\bul}}}
{\Om}^{i+j+1}_{{\cal Q}^{\rm ex}_{\bul \leq N,\bul}/\os{\circ}{T}{}'}, 
\end{CD} 
\end{equation*} 
we can define the pull-back morphism 
\begin{align*} 
g^*_{\bul \leq N,\bul} \col &
A_{\rm zar}({\cal Q}^{\rm ex}_{\bul \leq N,\bul}/S'(T')^{\nat},{\cal F}^{\bul \leq N,\bul})
\lo 
g^{\rm PD}_{\bul \leq N,\bul*}A_{\rm zar}({\cal P}^{\rm ex}_{\bul \leq N,\bul}/S(T)^{\nat},
{\cal E}^{\bul \leq N,\bul})
\tag{1.5.2.14}\label{eqn:axfdwt}
\end{align*} 
by the following formula 
\begin{align*} 
g^*_{\bul \leq N,\bul}
:= \deg(u)^{-(j+1)}g^{{\rm PD}*}_{\bul \leq N,\bul} 
\col &
A_{\rm zar}({\cal Q}^{\rm ex}_{\bul \leq N,\bul}/S'(T')^{\nat},{\cal F}^{\bul \leq N,\bul})^{ij}
\tag{1.5.2.15}\label{eqn:axdwt}\\
& \lo 
g^{\rm PD}_{\bul \leq N,\bul*}A_{\rm zar}({\cal P}^{\rm ex}_{\bul \leq N,\bul}/S(T)^{\nat},{\cal E}^{\bul \leq N,\bul})^{ij}. 
\end{align*} 
In fact, by (\ref{cd:axfpdwt}), 
we have the following filtered morphism 
\begin{align*} 
g^*_{\bul \leq N,\bul}
:= \deg(u)^{-(j+1)}g^{{\rm PD}*}_{\bul \leq N,\bul} 
\col &
(A_{\rm zar}({\cal Q}^{\rm ex}_{\bul \leq N,\bul}/S'(T')^{\nat},{\cal F}^{\bul \leq N,\bul})^{ij},P) 
\tag{1.5.2.16}\label{eqn:axfit}\\
& \lo 
g^{\rm PD}_{\bul \leq N,\bul*}((A_{\rm zar}({\cal P}^{\rm ex}_{\bul \leq N,\bul}/S(T)^{\nat},
{\cal E}^{\bul \leq N,\bul})^{ij},P)). 
\end{align*} 
Let 
\begin{equation*}
g^*_{\bul \leq N}\col 
(A_{\rm zar}(Y_{\bul \leq N,\os{\circ}{T}{}'_0}/S'(T')^{\nat},F^{\bul \leq N}),P)
\lo Rg_{\bul \leq N *}(A_{\rm zar}(X_{\bul \leq N,\os{\circ}{T}_0}/S(T)^{\nat},E^{\bul \leq N}),P) 
\tag{1.5.2.17}\label{eqn:axnwt}
\end{equation*} 
be the induced morphism by $g^*_{\bul \leq N,\bul}$.  
\par  
Next we check that the morphism (\ref{eqn:axnwt})
is independent of the choice of the diagrams 
(\ref{cd:xygxy}) and (\ref{cd:xpnxp}).  
Let $X''_{\bul \leq N,\os{\circ}{T}_0}$ and $Y''_{\bul \leq N,\os{\circ}{T}{}'_0}$ be 
the disjoint unions of the members of 
affine $N$-truncated simplicial open coverings of 
$X_{\bul \leq N,\os{\circ}{T}_0}$ and $Y_{\bul \leq N,\os{\circ}{T}{}'_0}$, 
respectively, fitting into the following commutative diagram 
(\ref{cd:xygxy}). 
Set $X'''_{m,\os{\circ}{T}_0}:=X'_{m,\os{\circ}{T}_0}
\times_{X_{m,\os{\circ}{T}_0}}X''_{m,\os{\circ}{T}_0}$ 
and $Y'''_{m,\os{\circ}{T}{}'_0}:=Y'_{m,\os{\circ}{T}{}'_0}\times_{Y_{m,\os{\circ}{T}{}'_0}}
Y''_{m,\os{\circ}{T}{}'_0}$ 
$(0\leq m\leq N)$. 
Set $X'_{mn,\os{\circ}{T}_0}:={\rm cosk}^{X_{m,\os{\circ}{T}_0}}_0(X''_{m,\os{\circ}{T}_0})_n$, 
$X''_{mn,\os{\circ}{T}_0}:={\rm cosk}^{X_{m,\os{\circ}{T}_0}}_0(X'''_{m,\os{\circ}{T}_0})_n$, 
$Y'_{mn,\os{\circ}{T}{}'_0}:={\rm cosk}^{Y_{m,\os{\circ}{T}{}'_0}}_0(Y''_{m,\os{\circ}{T}{}'_0})_n$ 
and  
$Y''_{mn,\os{\circ}{T}{}'_0}:={\rm cosk}^{Y_{m,\os{\circ}{T}{}'_0}}_0(Y'''_{m,\os{\circ}{T}{}'_0})_n$
$(0\leq m \leq N, n\in {\mab N})$.  
Then we have the following commutative diagram 
\begin{equation*} 
\begin{CD} 
X'_{\bul \leq N,\bul,\os{\circ}{T}_0}@>{\sus}>> 
\ol{\cal P}{}'^{\rm ex}_{\bul \leq N,\bul} \\ 
@V{g'_{\bul \leq N,\bul}}VV 
@VV{\ol{g}{}'_{\bul \leq N,\bul}}V \\ 
Y'_{\bul \leq N,\bul,\os{\circ}{T}_0} @>{\sus}>> \ol{\cal Q}{}'^{\rm ex}_{\bul \leq N,\bul}, 
\end{CD} 
\tag{1.5.2.18}\label{cd:xqxp} 
\end{equation*} 
where the two horizontal exact immersions above are the exactifications 
of immersions 
$X'_{\bul \leq N,\bul,\os{\circ}{T}_0} \os{\sus}{\lo} \ol{\cal P}{}'_{\bul \leq N,\bul}$
and 
$Y'_{\bul \leq N,\bul,\os{\circ}{T}{}'_0} \os{\sus}{\lo} \ol{\cal Q}{}'_{\bul \leq N,\bul}$ 
into $(N,\infty)$-truncated bisimplicial 
log smooth schemes over $\ol{S(T)^{\nat}}$ and $\ol{S'(T')^{\nat}}$, respectively.  
Set 
$\ol{\cal P}{}''_{\bul \leq N,\bul}
:=\ol{\cal P}_{\bul \leq N,\bul}\times_{\ol{S(T)^{\nat}}}\ol{\cal P}{}'_{\bul \leq N,\bul}$ 
and 
$\ol{\cal Q}{}''_{\bul \leq N,\bul}
:=\ol{\cal Q}_{\bul \leq N,\bul}\times_{\ol{S'(T')^{\nat}}}\ol{\cal Q}{}'_{\bul \leq N,\bul}$.  
Let $\ol{\cal P}{}''{}^{\rm ex}_{\bul \leq N,\bul}$ 
(resp.~$\ol{\cal Q}{}''{}^{\rm ex}_{\bul \leq N,\bul}$) 
be the exactification of 
the diagonal immersion 
$X''_{\bul \leq N,\bul,\os{\circ}{T}_0} \os{\sus}{\lo}
\ol{\cal P}{}''_{\bul \leq N,\bul}$
(resp.~$Y''_{\bul \leq N,\bul,\os{\circ}{T}{}'_0} \os{\sus}{\lo}
\ol{\cal Q}{}''_{\bul \leq N,\bul}$). 
Then we have the following commutative diagram 
\begin{equation*} 
\begin{CD} 
X''_{\bul \leq N,\bul,\os{\circ}{T}_0}@>{\sus}>> 
\ol{\cal P}{}''^{\rm ex}_{\bul \leq N,\bul}\\ 
@V{g''_{\bul \leq N,\bul}}VV 
@VV{\ol{g}{}''_{\bul \leq N,\bul}}V \\ 
Y''_{\bul \leq N,\bul,\os{\circ}{T}{}'_0} @>{\sus}>> \ol{\cal Q}{}''^{\rm ex}_{\bul \leq N,\bul} 
\end{CD} 
\tag{1.5.2.19}\label{cd:xexp} 
\end{equation*} 
over (\ref{cd:xpnxp}) and (\ref{cd:xqxp}). 
The rest of the proof of 
the well-definedness of the morphism (\ref{eqn:axnwt})
is similar to the proof of (\ref{theo:indcr}); we leave the detail of the rest to the reader. 
\par 
By (\ref{prop:tefc}), (\ref{cd:fblssod}) and (\ref{eqn:axdwt}), 
we have the commutative diagram (\ref{cd:pssfccz}).  
\par 
(2): (2) is clear from the construction of the morphism (\ref{eqn:fzaxd}) 
and the relation (\ref{ali:dpvv}).  
\par 
(3): This is obvious. 
\par 
(4): This follows from (1) and (\ref{eqn:utz}).  
\end{proof}

\par 
Next we consider the case where 
$\deg (u)_x$ may be divisible by $p$ for 
some point $x\in \os{\circ}{S}$. 
Let $e_p(x)$ $(x\in  \os{\circ}{S})$ 
be the exponent of $\deg (u)_x$ with respect to $p$: 
$p^{e_p(x)}\vert \vert \deg(u_x)$. 
Then we have a function 
\begin{equation*} 
e_p \col \os{\circ}{S}\owns x \lom e_p(x) \in {\mab N}.   
\tag{1.5.2.20}\label{eqn:expf}
\end{equation*} 

\begin{defi}\label{defi:efu} 
We call $e_p$ the {\it exponent function} of $u$ 
with respect to $p$. 
\end{defi}

Now we assume that 
$\os{\circ}{T}$ and $\os{\circ}{T}{}'$ are $($not necessarily affine$)$ 
$p$-adic formal schemes   
such that $p\in {\cal J}$ and $p\in {\cal J}'$ 
(\cite[7.17, Definition]{bob}). 
Assume that ${\cal O}_T$ is $p$-torsion-free. 
Let 
$f_{\bul}\col {\cal P}_{\bul \leq N,\bul}\lo S(T)^{\nat}$ 
be the structural morphism. 
Furthermore, for any $i,j\in {\mab N}$, 
assume that the induced morphism 
\begin{align*} 
g^{{\rm PD}*}_{\bul \leq N,\bul} 
\col &
{\cal F}^{\bul \leq N,\bul}
\otimes_{{\cal O}_{{\cal Q}^{\rm ex}_{\bul \leq N,\bul}}}
{\Om}^{i+j+1}_{{\cal Q}^{\rm ex}_{\bul \leq N,\bul}/\os{\circ}{T}{}'}/P_j 
\lo 
g^{\rm PD}_{\bul \leq N,\bul*}({\cal E}^{\bul \leq N,\bul}
\otimes_{{\cal O}_{{\cal P}^{\rm ex}_{\bul \leq N,\bul}}}
{\Om}^{i+j+1}_{{\cal P}^{\rm ex}_{\bul \leq N,\bul}/\os{\circ}{T}{}}/P_j) 
\tag{1.5.3.1}\label{eqn:odnl}
\end{align*}
by the following morphism 
\begin{equation*} 
g^{{\rm PD}*}_{\bul \leq N,\bul} 
\col 
{\cal F}^{\bul \leq N,\bul}
\otimes_{{\cal O}_{{\cal Q}^{\rm ex}_{\bul \leq N,\bul}}}
{\Om}^{i+j+1}_{{\cal Q}^{\rm ex}_{\bul \leq N,\bul}/\os{\circ}{T}{}'}
\lo 
g^{\rm PD}_{\bul \leq N,\bul*}
({\cal E}^{\bul \leq N,\bul}
\otimes_{{\cal O}_{{\cal P}^{\rm ex}_{\bul \leq N,\bul}}}
{\Om}^{i+j+1}_{{\cal P}^{\rm ex}_{\bul \leq N,\bul}/\os{\circ}{T}{}}) 
\tag{1.5.3.2}\label{eqn:onfdnl}
\end{equation*}
is divisible by $p^{e_p(j+1)}$, that is, 
$g^{{\rm PD}*}_{\bul \leq N,\bul}$ 
at any point $x\in {\cal P}^{\rm ex}_{\bul \leq N,\bul}$
is divisible by $p^{e_p(f_{\bul}(x))(j+1)}$. 
Because the target of (\ref{eqn:odnl}) is ${\cal O}_T$-flat
by (\ref{coro:flt}), 
the following morphism 
\begin{align*} 
{\rm deg}(u)^{-(j+1)}g^{{\rm PD}*}_{\bul \leq N,\bul} 
\col &
{\cal F}^{\bul \leq N,\bul}
\otimes_{{\cal O}_{{\cal Q}^{\rm ex}_{\bul \leq N,\bul}}}
{\Om}^{i+j+1}_{{\cal Q}^{\rm ex}_{\bul \leq N,\bul}/\os{\circ}{T}{}'}/P_j 
\tag{1.5.3.3}\label{eqn:oddnl}\\
&\lo 
g^{\rm PD}_{\bul \leq N,\bul*}
({\cal E}^{\bul \leq N,\bul}
\otimes_{{\cal O}_{{\cal P}^{\rm ex}_{\bul \leq N,\bul}}}
{\Om}^{i+j+1}_{{\cal P}^{\rm ex}_{\bul \leq N,\bul}/\os{\circ}{T}{}}/P_j)
\end{align*}
for $j\in {\mab N}$ is well-defined.  
By (\ref{coro:flt}) we see that the following morphism 
\begin{align*} 
{\rm deg}(u)^{-(j+1)}g^{{\rm PD}*}_{\bul \leq N,\bul} 
\col &
(P_{2j+k+1}+P_j)({\cal F}^{\bul \leq N,\bul}
\otimes_{{\cal O}_{{\cal Q}^{\rm ex}_{\bul \leq N,\bul}}}
{\Om}^{i+j+1}_{{\cal Q}^{\rm ex}_{\bul \leq N,\bul}/\os{\circ}{T}{}'})/P_j 
\tag{1.5.3.4}\label{eqn:odpdnl}\\
&\lo 
g^{\rm PD}_{\bul \leq N,\bul*}
((P_{2j+k+1}+P_j)({\cal E}^{\bul \leq N,\bul}
\otimes_{{\cal O}_{{\cal P}^{\rm ex}_{\bul \leq N,\bul}}}
{\Om}^{i+j+1}_{{\cal P}^{\rm ex}_{\bul \leq N,\bul}/\os{\circ}{T}{}})/P_j) 
\end{align*}
for $k\in {\mab Z}$ and $j\in {\mab N}$ is also well-defined.  
Note that, if 
$\om \in 
p^nP_k({\cal E}^{\bul \leq N,\bul}
\otimes_{{\cal O}_{{\cal P}^{\rm ex}_{\bul \leq N,\bul}}}
{\Om}^{i}_{{\cal P}^{\rm ex}_{\bul \leq N,\bul}/\os{\circ}{T}{}})$ 
for $n\in {\mab N}$, then 
$p^{-n}\om \in P_k({\cal E}^{\bul \leq N,\bul}
\otimes_{{\cal O}_{{\cal P}^{\rm ex}_{\bul \leq N,\bul}}}
{\Om}^{i}_{{\cal P}^{\rm ex}_{\bul \leq N,\bul}/\os{\circ}{T}{}})$.

\begin{rema}\label{rema:dam}
The condition 
``$g^{{\rm PD}*}_{\bul \leq N,\bul}$ is divisible by $p^{e_p(j+1)}$'' 
does not always hold. 
Indeed, consider the following example. 
\par 
Let $s$ be the log point of a perfect field $\kap$ of characteristic $p>0$. 
Let ${\cal W}$ be the Witt ring of $\kap$ 
and let ${\cal W}(s)$ be the canonical lift of $s$ 
over ${\rm Spf}({\cal W})$.  
Set ${\cal X}:={\mab A}_{{\cal W}}(2,2)$ with natural morphism to ${\cal W}(s)$ 
defined by the diagonal embedding ${\mab N}\os{\sus}{\lo} {\mab N}^2$. 
Let $x,y$ be the canonical coordinates of ${\cal X}$. 
Let $F\col {\cal W}(s)\lo {\cal W}(s)$ be the Frobenius endomorphism. 
Consider an endomorphism ${\mathfrak g}\col {\cal X}\lo {\cal X}$ induced by a morphism 
${\mab N}^2\lo {\mab N}^2$ of monoids defined by 
$(1,0)\lom (1,p-1)$ and $(0,1)\lom (p-1,1)$. 
Then the morphism 
\begin{align*} 
{\mathfrak g}^*\col \Om^1_{{\cal X}/{\cal W}}/P_0\lo 
{\mathfrak g}_*(\Om^1_{{\cal X}/{\cal W}}/P_0)
\end{align*} 
does not satisfy the condition because,   
in the following commutative diagram 
\begin{equation*} 
\begin{CD} 
P_1\Om^1_{{\cal X}/{\cal W}}/P_0@>{{\mathfrak g}^*}>>
{\mathfrak g}_*(P_1\Om^1_{{\cal X}/{\cal W}}/P_0)\\
@V{\rm Res}V{\simeq}V @V{\simeq}V{\rm Res}V \\
{\cal O}_{\os{\circ}{\cal X}{}^{(1)}}@>>>{\cal O}_{\os{\circ}{\cal X}{}^{(1)}},
\end{CD} 
\end{equation*} 
the image of the section $d\log x$ of the upper left
$P_1\Om^1_{{\cal X}/{\cal W}}/P_0$  
in the lower right ${\cal O}_{\os{\circ}{\cal X}{}^{(1)}}$ is equal to 
$(1,p-1)\in {\cal O}_{\cal X}/x\oplus {\cal O}_{\cal X}/y={\cal W}\{y\}\oplus {\cal W}\{x\}$. 
This is not divisible by $p$. 
\end{rema}

\begin{lemm}\label{lemm:divp}  
The divisibility assumption for the morphism 
{\rm (\ref{eqn:odnl})} is independent of the choices of 
the $N$-truncated simplicial affine open coverings of $X_{\bul \leq N,\os{\circ}{T}_0}$ 
and $Y_{\bul \leq N,\os{\circ}{T}{}'_0}$ and the choices of 
the $(N,\infty)$-truncated bisimplicial immersions of 
$X_{\bul \leq N,\bul,\os{\circ}{T}} \os{\sus}{\lo} 
\ol{\cal P}_{\bul \leq N,\bul}$ over $\ol{S(T)^{\nat}}$ and 
$Y_{\bul \leq N,\bul,\os{\circ}{T}{}'_0} 
\os{\sus}{\lo} \ol{\cal Q}_{\bul \leq N,\bul}$ over $\ol{S'(T')^{\nat}}$ 
giving the commutative diagram {\rm (\ref{cd:xpnxp})}. 
\end{lemm} 
\begin{proof} 
Let $X_{\bul \leq N,\bul,\os{\circ}{T}_0}\os{\sus}{\lo} \ol{\cal P}{}'_{\bul \leq N,\bul}$ 
and 
$Y_{\bul \leq N,\bul,\os{\circ}{T}{}'_0}\os{\sus}{\lo} \ol{\cal Q}{}'_{\bul \leq N,\bul}$ 
be other immersions into 
log smooth $(N,\infty)$-truncated bisimplicial schemes over $\ol{S(T)^{\nat}}$ 
and $\ol{S'(T')^{\nat}}$, respectively,  fitting another commutative diagram 
\begin{equation*} 
\begin{CD} 
X_{\bul \leq N,\bul,\os{\circ}{T}_0} 
@>{\subset}>> \ol{\cal P}{}'{}^{\rm ex}_{\! \! \!\bul \leq N,\bul} \\
@VVV @VVV  \\
Y_{\bul \leq N,\bul,\os{\circ}{T}{}'_0} 
@>{\subset}>> \ol{\cal Q}{}'{}^{\rm ex}_{\bul \leq N,\bul}.    
\end{CD}
\end{equation*}  
Then, by considering products as in the proof of (\ref{theo:funas}), 
we may assume that there exists the following commutative diagram 
\begin{equation*} 
\begin{CD} 
X_{\bul \leq N,\bul,\os{\circ}{T}_0} @>{\subset}>> 
\ol{\cal P}{}^{\rm ex}_{\bul \leq N,\bul} @>>> 
\ol{\cal P}{}'{}^{\rm ex}_{\! \! \!\bul \leq N,\bul} \\
@VVV @VVV @VVV \\
Y_{\bul \leq N,\bul,\os{\circ}{T}{}'_0} 
@>{\subset}>> \ol{\cal Q}{}^{\rm ex}_{\bul \leq N,\bul} @>>> 
\ol{\cal Q}{}'{}^{\rm ex}_{\bul \leq N,\bul}.    
\end{CD}
\end{equation*}  
Hence we may assume that $N=0$. 
Set $X:=X_0$, $Y:=Y_0$, $E=E^0$ and $F:=F^0$. 
Because the question is local, we may assume that 
there exists an immersion $X_{\os{\circ}{T}_0}\os{\sus}{\lo} \ol{\cal P}$ 
(resp.~$Y_{\os{\circ}{T}{}'_0}\os{\sus}{\lo} \ol{\cal Q}$) into a log smooth scheme 
over $\ol{S(T)^{\nat}}$  
with morphism $\ol{\cal P}\lo \ol{\cal P}{}'$ 
(resp.~$\ol{\cal Q}\lo \ol{\cal Q}{}'$) of log smooth schemes over $\ol{S(T)^{\nat}}$ 
(resp.~$\ol{S'(T')^{\nat}}$)
such that the composite morphism 
$X_{\os{\circ}{T}_0}\os{\sus}{\lo} \ol{\cal P}\lo \ol{\cal P}{}'$ is also an immersion 
(resp.~$Y_{\os{\circ}{T}{}'_0}\os{\sus}{\lo} \ol{\cal Q}\lo \ol{\cal Q}{}'$).  
Let $\ol{\mathfrak D}$ and $\ol{\mathfrak D}{}'$ 
be the log PD-envelopes of the immersions $X_{\os{\circ}{T}_0}\os{\sus}{\lo} \ol{\cal P}$ 
and $X_{\os{\circ}{T}_0}\os{\sus}{\lo} \ol{\cal P}{}'$ 
over $(\os{\circ}{T},{\cal J},\del)$, respectively.  
We may also assume that 
there exists the following commutative diagram
\begin{equation*} 
\begin{CD} 
X_{\os{\circ}{T}_0} @>{\subset}>> \ol{\cal P}{}^{\rm ex}
@>>>\ol{\cal P}{}'{}^{\rm ex}\\
@VVV @VVV @VVV \\
Y_{\os{\circ}{T}{}'_0} @>{\subset}>> \ol{\cal Q}{}^{\rm ex} @>>> 
\ol{\cal Q}{}'{}^{\rm ex}.   
\end{CD}
\tag{1.5.5.1}\label{cd:xyppq}
\end{equation*} 
Set ${\mathfrak D}:=\ol{\mathfrak D}\times_{{\mathfrak D}(\ol{S(T)^{\nat}})}S(T)^{\nat}$ 
and ${\mathfrak D}':=\ol{\mathfrak D}{}'\times_{{\mathfrak D}(\ol{S(T)^{\nat}})}S(T)^{\nat}$.   
Let $(\ol{\cal E},\ol{\nabla})$ and $(\ol{\cal E}{}',\ol{\nabla}{}')$ be 
an ${\cal O}_{\ol{\mathfrak D}}$-module with integrable connection 
and an ${\cal O}_{\ol{\mathfrak D}{}'}$-module 
with integrable connection obtained by $E$, respectively. 
Set 
${\cal E}:=\ol{\cal E}\otimes_{{\cal O}_{\mathfrak D}(\ol{S(T)^{\nat}})}{\cal O}_T$ 
and ${\cal E}':=\ol{\cal E}{}'\otimes_{{\cal O}_{\mathfrak D}(\ol{S'(T')^{\nat}})}{\cal O}_{T'}$. 
By the local structures of the exact immersions ((\ref{prop:adla})), 
we may assume that 
$\ol{\cal P}{}^{\rm ex}=\ol{\cal P}{}'^{\rm ex}\times_{\ol{S(T)^{\nat}}}{\mab A}_{\ol{S(T)^{\nat}}}^c$ 
for a nonnegative integer $c$. 
Since $E$ is crystal, 
$\ol{\cal E}=\ol{\cal E}{}'\otimes_{{\cal O}_{\ol{S(T)^{\nat}}}}
{\cal O}_{\ol{S(T)^{\nat}}}\langle x_1,\ldots,x_c\rangle$. 
Hence ${\cal E}={\cal E}'\otimes_{{\cal O}_T}
{\cal O}_T\langle x_1,\ldots,x_c\rangle$. 
Because ${\cal O}_T\langle x_1,\ldots,x_c\rangle \otimes_{{\cal O}_T} 
\Om^{i}_{{\mab A}^c_{\os{\circ}{T}}/\os{\circ}{T}} $ 
$(i\in {\mab N})$ is a free ${\cal O}_T$-module, 
the morphism ${\cal O}_T\lo 
{\cal O}_T\langle x_1,\ldots,x_c\rangle \otimes_{{\cal O}_T} 
\Om^{i}_{{\mab A}^c_{\os{\circ}{T}}/\os{\circ}{T}}$ 
is faithfully flat.  
By (\ref{eqn:dopx}) and (\ref{eqn:dpopx}), 
we have the following formula  
\begin{equation*} 
{\cal E}
\otimes_{{\cal O}_{{\cal P}^{\rm ex}}}
\Om^{i+j+1}_{{\cal P}^{\rm ex}/\os{\circ}{T}}/P_j
\simeq 
\bigoplus_{i'+i''=i+j+1}
({\cal E}'\otimes_{{\cal O}_{{\cal P}'{}^{\rm ex}}}
\Om^{i'}_{{\cal P}'{}^{\rm ex}/\os{\circ}{T}}/P_j)
\otimes_{{\cal O}_T}
{\cal O}_T\langle x_1,\ldots, x_c\rangle \otimes_{{\cal O}_T}
\Om^{i''}_{{\mab A}^c_{\os{\circ}{T}}/\os{\circ}{T}}  
\tag{1.5.5.2}\label{eqn:dpdpx} 
\end{equation*}
Because we may assume that the morphism (\ref{eqn:odnl}) 
for the case $N=0$ and $\bul=0$ factors through the morphism 
\begin{align*} 
g^{{\rm PD}*} \col &{\cal F}\otimes_{{\cal O}_{{\cal Q}^{\rm ex}}}
{\Om}^{i+j+1}_{{\cal Q}/\os{\circ}{T}{}'}/P_j \lo 
g^{\rm PD}_{*}({\cal E}'
\otimes_{{\cal O}_{{\cal P}'{}^{\rm ex}}}
{\Om}^{i+j+1}_{{\cal P}'{}^{\rm ex}/\os{\circ}{T}{}}/P_j),  
\end{align*}
the divisibilities in 
${\cal E}\otimes_{{\cal O}_{{\cal P}^{\rm ex}}}
{\Om}^{i+j+1}_{{\cal P}^{\rm ex}/\os{\circ}{T}{}}/P_j$ $(\forall i\in {\mab N})$ and  
${\cal E}'\otimes_{{\cal O}_{{\cal P}'{}^{\rm ex}}}
{\Om}^{i+j+1}_{{\cal P}'{}^{\rm ex}/\os{\circ}{T}{}}/P_j$ 
$(\forall i\in {\mab N})$ are equivalent. 
\end{proof}

\begin{theo}[{\bf Contravariant functoriality II of $A_{\rm zar}$}]\label{theo:funpas} 
$(1)$ Let the notations and the assumptions be as above. 
Then $g_{\bul \leq N}\col X_{\bul \leq N,\os{\circ}{T}{}_0}
\lo Y_{\bul \leq N,\os{\circ}{T}{}'_0}$ induces the following 
well-defined pull-back morphism 
\begin{equation*}  
g_{\bul \leq N}^* \col 
(A_{\rm zar}(Y_{\bul \leq N,\os{\circ}{T}{}'_0}/S'(T')^{\nat},F^{\bul \leq N}),P)
\lo Rg_{\bul \leq N*}
((A_{\rm zar}(X_{\bul \leq N,\os{\circ}{T}_0}/S(T)^{\nat},E^{\bul \leq N}),P))  
\tag{1.5.6.1}\label{eqn:fzapxd}
\end{equation*} 
fitting into the following commutative diagram$:$
\begin{equation*}  
\begin{CD}
A_{\rm zar}(Y_{\bul \leq N,\os{\circ}{T}{}'_0}/S'(T')^{\nat},F^{\bul \leq N})
@>{g_{\bul \leq N}^*}>>  \\ 
@A{\theta_{Y_{\bul \leq N,\os{\circ}{T}{}'_0}/S'(T')^{\nat}}\wedge}A{\simeq}A \\
Ru_{Y_{\bul \leq N,\os{\circ}{T}{}'_0}/S'(T')^{\nat}*}
(\eps^*_{Y_{\bul \leq N,\os{\circ}{T}{}'_0}/S'(T')^{\nat}}(F^{\bul \leq N}))
@>{g_{\bul \leq N}^*}>>
\end{CD}
\tag{1.5.6.2}\label{cd:pssfpccz} 
\end{equation*}
\begin{equation*}  
\begin{CD}
Rg_{\bul \leq N*}(A_{\rm zar}(X_{\bul \leq N,\os{\circ}{T}_0}/S(T)^{\nat},E^{\bul \leq N})) \\ 
@A{Rg_{\bul \leq N*}(\theta_{X_{\bul \leq N,\os{\circ}{T}_0}/S(T)^{\nat}}\wedge)}A{\simeq}A\\
Rg_{\bul \leq N*}Ru_{X_{\bul \leq N,\os{\circ}{T}_0}/S(T)^{\nat}*}
(\eps^*_{X_{\bul \leq N,\os{\circ}{T}_0}/S(T)^{\nat}}(E^{\bul \leq N})). 
\end{CD} 
\end{equation*}
\par 
$(2)$ The similar relation to {\rm (\ref{ali:pdpp})} holds.  
\par 
$(3)$ There exists the similar commutative diagram to {\rm (\ref{cd:tmfccz})}.
\end{theo}
\begin{proof} 
The proof is the same as that of (\ref{theo:funas}). 
\end{proof} 



\par

\par 
Let the notations and the assumptions be 
as in (\ref{theo:funas}) or (\ref{theo:funpas}).  
\par 
Consider the morphism 
\begin{align*} 
{\rm gr}^P_k(g^*_{\bul \leq N}) 
\col & {\rm gr}^P_kA_{\rm zar}
(Y_{\bul \leq N,\os{\circ}{T}{}'_0}/S'(T')^{\nat},F^{\bul \leq N}) 
\lo Rg_{\bul \leq N*}
({\rm gr}^P_kA_{\rm zar}(X_{\bul \leq N,\os{\circ}{T}_0}/S(T)^{\nat},E^{\bul \leq N})). 
\tag{1.5.6.3}\label{ali:gdkcl}
\end{align*}

\par
Fix any integer $0\leq m\leq N$. 
For simplicity of notation, set $X:=X_{m,\os{\circ}{T}_0}$, $Y:=Y_{m,\os{\circ}{T}{}'_0}$, 
$E:=E^m$, $F:=F^m$ and $g:=g_m$. 
The following conditions are 
SNCL versions of 
the conditions \cite[(2.9.2.3), (2.9.2.4)]{nh2} for the open log case.  
\par   
Assume that the following two conditions hold:

\medskip
\parno 
$(1.5.6.4)$: for any smooth component 
$\os{\circ}{X}_{\lam}$ of $\os{\circ}{X}_{T_0}$ over $\os{\circ}{T}_0$, 
there exists a unique smooth component 
$\os{\circ}{Y}_{\mu}$ of 
$\os{\circ}{Y}_{T'_0}$ over $\os{\circ}{T}{}'_0$ such that $g$ 
induces a morphism 
$\os{\circ}{X}_{\lam} \lo \os{\circ}{Y}_{\mu}$. 
(Let $\Lam$ and $M$ be the sets of indices 
of the $\lam$'s and the $\mu$'s, respectively. 
Then we obtain a function 
$\phi \col \Lam \owns \lam \lom \mu \in M$.)
\medskip
\parno
$(1.5.6.5)$: there exist positive integers $e({\lam})$'s  
$(\lam \in \Lam)$ such that 
there exist local equations $x_{\lam}=0$ and 
$y_{\phi(\lam)}=0$ of 
$\os{\circ}{X}_{\lam}$ and $\os{\circ}{Y}_{\phi(\lam)}$, 
respectively,  such that $g^*(y_{\phi(\lam)})=x^{e({\lam})}_{\lam}$.

\begin{prop}\label{prop:xxle} 
Let the assumptions and the notations be as above. 
Set $\Lam(x):=
\{\lam \in \Lam~\vert~x \in \os{\circ}{X}_{\lam}\}$.  
Then $\deg(u)_x=e({\lam})$ for $\lam \in \Lam(x)$.  
In particular, 
$e({\lam})$'s are independent of
the choice of an element of $\Lam(x)$.  
\end{prop}
\begin{proof} 
We have the following commutative diagram 
\begin{equation*} 
\begin{CD} 
M_{X_{\os{\circ}{T}_0},x}/{\cal O}^*_{X_{\os{\circ}{T}_0},x} @<{g^*}<< 
M_{Y_{\os{\circ}{T}{}'_0},\os{\circ}{g}(x)}/{\cal O}^*_{Y_{\os{\circ}{T}{}'_0},\os{\circ}{g}(x)} \\ 
@V{\simeq}VV @AA{\bigcup}A \\ 
{\mab N}^{\Lam(x)}@<<< {\mab N}^{\Lam(x)}\\ 
@A{\rm diag.}AA @AA{\rm diag.}A \\ 
{\mab N}@<<< {\mab N},  
\end{CD} 
\end{equation*} 
where the image of $(1,\ldots,1)$ 
by the middle horizontal morphism 
is $(e({\lam}))_{\lam \in \Lam(x)}$. 
On the other hand, we see that 
the image of $(1,\ldots,1)$ is 
${\rm deg}(u)_x(1,\ldots,1)$ by the definition of 
${\rm deg}(u)$. Hence $e({\lam})=\deg(u)_x$.  
\end{proof} 


\begin{prop}\label{prop:dvok}
Let $n$ be a positive integer. 
Assume that ${\cal J}\subset p^n{\cal O}_T$. 
Assume that $e_p\geq 1$. 
Set $m:=\min\{n,e_p\}$ $(m$ is a function on $\os{\circ}{X}_{T_0})$. 
Then the morphism {\rm (\ref{eqn:odnl})} is divisible by $p^{m(j+1)}$. 
$($As a result, if $e_p\leq n$, then  the morphism {\rm (\ref{eqn:odnl})}
satisfies the divisibility assumption after {\rm (\ref{defi:efu})}.$)$ 
\end{prop}
\begin{proof} 
As in the proof of (\ref{lemm:divp}), we may assume that $N=0$. 
Set $X:=X_0$, $Y:=Y_0$, $E=E^0$ and $F:=F^0$. 
Because the question is local, we may assume that 
there exists the commutative diagram (\ref{cd:xyppq}). 
Furthermore, we may assume that $\ol{\cal Q}{}^{\rm ex}$ 
is a log smooth lift $\ol{\cal Y}$ of 
$Y_{\os{\circ}{T}{}'_0}$ over $\ol{S'(T')^{\nat}}$.
Because the immersion $X_{\os{\circ}{T}_0}\os{\sus}{\lo} \ol{\cal P}{}^{\rm ex}$ 
is exact, we may also assume that 
$\ol{\cal P}{}^{\rm ex}$ is a formal log smooth lift $\ol{\cal X}$ of 
$X_{\os{\circ}{T}_0}$ over $\ol{S(T)^{\nat}}$ ((\ref{prop:adla})). 
Let $\ol{\mathfrak D}$ and $\ol{\mathfrak E}$ 
be the log PD-envelopes of 
the immersions $X_{\os{\circ}{T}_0}\os{\sus}{\lo}\ol{\cal X}$ 
and $Y_{\os{\circ}{T}{}'_0}\os{\sus}{\lo}\ol{\cal Y}$ over 
$(\os{\circ}{T},{\cal J},\del)$ and $(\os{\circ}{T}{}',{\cal J}',\del')$, respectively.  
Set ${\cal X}:=\ol{\cal X}\times_{\ol{S(T)^{\nat}}}S(T)^{\nat}$ and 
${\cal Y}:=\ol{\cal Y}\times_{\ol{S'(T')^{\nat}}}S'(T')^{\nat}$. 
Let $(\ol{\cal E},\ol{\nabla})$ and $(\ol{\cal F},\ol{\nabla})$ 
be the ${\cal O}_{\ol{\mathfrak D}}$-module with integrable connection 
and the ${\cal O}_{\ol{\mathfrak E}}$-module 
with integrable connection corresponding to 
$\eps_{X_{\os{\circ}{T}_0}/\os{\circ}{T}}^*(E)$ and 
$\eps_{Y_{\os{\circ}{T}{}'_0}/\os{\circ}{T}{}'}^*(F)$, respectively. 
Set 
${\cal E}:=\ol{\cal E}\otimes_{{\cal O}_{{\mathfrak D}(\ol{S(T)})}}{\cal O}_{S(T)}$ 
and  
${\cal F}:=\ol{\cal F}\otimes_{{\cal O}_{{\mathfrak D}(\ol{S'(T')})}}{\cal O}_{S'(T')}$.  
The morphism $\ol{\cal X}\lo \ol{\cal Y}$ induces 
the following natural morphism 
\begin{equation*} 
g^{{\rm PD}*}\col {\cal F}
\otimes_{{\cal O}_{\cal Y}}
\Om^{\bul}_{{\cal Y}/\os{\circ}{T}{}'}
\lo 
g^{\rm PD}_*({\cal E}
\otimes_{{\cal O}_{\cal X}}\Om^{\bul}_{{\cal X}/\os{\circ}{T}}) 
\end{equation*}
and we have the following equalities:  
\begin{equation*}
A_{\rm zar}(X_{\os{\circ}{T}_0}/S(T)^{\nat},E) 
= 
s(({\cal E}\otimes_{{\cal O}_{\cal X}}
{\Om}^{i+j+1}_{{\cal X}/\os{\circ}{T}}/P_j)_{i,j \in {\mab N}})
\tag{1.5.8.1}\label{eqn:grdlcal}
\end{equation*}   
and 
\begin{equation*}
A_{\rm zar}(Y_{\os{\circ}{T}{}'_0}/S'(T')^{\nat},F)
= 
s(({\cal F}\otimes_{{\cal O}_{\cal Y}}
{\Om}^{i+j+1}_{{\cal Y}/\os{\circ}{T}{}'}/P_j)_{i,j \in {\mab N}}).
\tag{1.5.8.2}\label{eqn:grycal}
\end{equation*}  
Consider the following morphism:
\begin{align*}
g^{{\rm PD}*}\col  {\cal F}\otimes_{{\cal O}_{\cal Y}}
{\Om}^{i+j+1}_{{\cal Y}/\os{\circ}{T}{}'}/P_j
\lo  g^{\rm PD}_*({\cal E}\otimes_{{\cal O}_{{\cal X}}}
{\Om}^{i+j+1}_{{\cal X}/\os{\circ}{T}}/P_j).  
\tag{1.5.8.3} \\ 
\end{align*} 
Let $\Lam(x)$ be the set in the proof of (\ref{prop:xxle}). 
Let $\os{\circ}{\cal X}_{\lam}$ $(\lam \in \Lam(x))$  
(resp.~$\os{\circ}{\cal Y}_{\phi(\lam)}$) 
be a closed subscheme 
of $\os{\circ}{\cal P}{}^{\rm ex}$ (resp.~$\os{\circ}{\cal Q}{}^{\rm ex}$) 
which is a lift of $\os{\circ}{X}_{\lam}$ 
(resp.~$\os{\circ}{Y}_{\phi(\lam)}$). 
Let $\wt{x}_{\lam}=0$ (resp.~$\wt{y}_{\lam}=0$) be a 
local equation of $\os{\circ}{\cal X}_{\lam}$  
(resp.~$\os{\circ}{\cal Y}_{\phi(\lam)}$). 
Let $\hat{x}_{\lam} \in M_{\cal X}$ 
(resp.~$\hat{y}_{\lam} \in M_{\cal Y}$)  
be the local inverse image of $\wt{x}_{\lam}$ (resp.~$\wt{y}_{\lam}$). 
For $k\geq j$, let $\om =\eta d\log \hat{y}_{\lam_0}\cdots d\log \hat{y}_{\lam_k}$ 
be a local section of
$P_{k+1}({\cal F}\otimes_{{\cal O}_{\cal Y}}
{\Om}^{i+j+1}_{{\cal Y}/\os{\circ}{T}{}'})$, 
where $\eta$ is a local section of 
${\rm Im}({\cal F}\otimes_{{\cal O}_{\cal Y}}
\Om^{i+j-k}_{\os{\circ}{\cal Y}/\os{\circ}{T}{}'}\lo {\cal F}\otimes_{{\cal O}_{\cal Y}}
{\Om}^{i+j-k}_{{\cal Y}/\os{\circ}{T}{}'})$.  
Then, by the conditions (1.5.6.4) and (1.5.6.5) and by (\ref{prop:xxle}),  
$g^{{\rm PD}*}(\hat{y}_i)=\hat{x}_i^{\deg(u)}(1+p^na_i)$ 
for some $a_i\in {\cal O}_{\cal X}$ since ${\cal J}\subset p^n{\cal O}_T$. 
Hence 
\begin{align*}
g^{{\rm PD}*}(d\log \hat{y}_i)&=\deg(u)d\log \hat{x}_i+(1+p^na_i)^{-1}p^nda_i\\
&=p^m((p^{-m}\deg u)d\log \hat{x}_i+p^{n-m}(1+p^na_i)^{-1}da_i).
\end{align*}  
This shows that $g^{{\rm PD}*}(\om)$ is divisible by $p^{m(k+1)}$. 
Since $k\geq j$, $g^{{\rm PD}*}(\om)$ is divisible by $p^{m(j+1)}$. 
This is the desired divisibility. 
\end{proof}

\par 
For a nonnegative integer $k$, 
let $\Lam^{(k)}(\os{\circ}{g})$ be the set of subsets $I$'s of $\Lam$ 
such that $\sharp I= \sharp \phi (I)=k+1$.  
For $\ul{\lam}=\{\lam_0,\ldots, \lam_k\}$ 
$(\ul{\lam}\in \Lam^{(k)}(\os{\circ}{g}))$, 
set 
$\os{\circ}{X}_{\ul{\lam}}:=
\os{\circ}{X}_{\lam_0}\cap \cdots \cap \os{\circ}{X}_{\lam_{k}}$ 
and $\os{\circ}{Y}_{\phi(\ul{\lam})}:=
\os{\circ}{Y}_{\phi(\lam_0)}\cap \cdots 
\cap \os{\circ}{Y}_{\phi(\lam_{k})}$.  
Let 
$\os{\circ}{g}_{\ul{\lam}} \col \os{\circ}{X}_{\ul{\lam}} \lo \os{\circ}{Y}_{\phi(\ul{\lam})}$ 
be the induced morphism by $g$.
\par
Now we change the notation $\ul{\lam}$ for (\ref{prop:grloc}) below. 
For integers $j$ and $k$ such that $j\geq \max \{-k,0\}$, 
set $\ul{\lam}:=\{\lam_0, \ldots, \lam_{2j+k}\}\in 
\Lam^{(2j+k)}(\os{\circ}{g})$. 
Let $a_{\ul{\lam}} \col \os{\circ}{X}_{\ul{\lam}} \os{\sus}{\lo} \os{\circ}{X}_{T_0}$  
and 
$b_{\phi(\ul{\lam})} \col \os{\circ}{Y}_{\phi(\ul{\lam})}\os{\sus}{\lo}  \os{\circ}{Y}_{T'_0}$ 
be the natural closed immersions.  
Let 
$$a_{\ul{\lam}{\rm crys}} 
\col 
((\os{\circ}{X}_{\ul{\lam}}/{\os{\circ}{T}})_{\rm crys},
{\cal O}_{\os{\circ}{X}_{\ul{\lam}}/{\os{\circ}{T}}})
\lo 
((\os{\circ}{X}_{T_0}/{\os{\circ}{T}})_{\rm crys},
{\cal O}_{\os{\circ}{X}_{T_0}/\os{\circ}{T}})$$  
and 
$$b_{\phi(\ul{\lam}){\rm crys}} 
\col 
((\os{\circ}{Y}_{\phi(\ul{\lam})}/{\os{\circ}{T}{}'})_{\rm crys},
{\cal O}_{\os{\circ}{Y}_{\phi(\ul{\lam}),T'_0}/
{\os{\circ}{T}{}'}})
\lo 
((\os{\circ}{Y}_{T'_0}/{\os{\circ}{T}{}'})_{\rm crys},
{\cal O}_{\os{\circ}{Y}_{T'_0}/\os{\circ}{T}{}'})$$  
be the induced morphisms of ringed topoi by 
$a_{\ul{\lam}}$ and $b_{\phi(\ul{\lam})}$, respectively. 
Set $E_{\ul{\lam}}:=a^*_{\ul{\lam}{\rm crys}}(E)$ 
and $F_{\phi(\ul{\lam})}:=b^*_{\phi(\ul{\lam}){\rm crys}}(F)$.  
Let 
$\vp_{\ul{\lam}{\rm crys}}(\os{\circ}{X}_{T_0}/\os{\circ}{T})$ 
(resp.~
$\vp_{\phi(\ul{\lam}){\rm crys}}(\os{\circ}{Y}_{T'_0}/\os{\circ}{T}{}')$) 
be the crystalline orientation sheaf in 
$(\os{\circ}{X}_{\ul{\lam}}/{\os{\circ}{T}})_{\rm crys},$ 
(resp.~$(\os{\circ}{Y}_{\ul{\lam}}/{\os{\circ}{T}{}'})_{\rm crys}$) 
similarly defined in \S\ref{sec:ldc} for the set 
$\{\os{\circ}{X}_{\lam_0}, \ldots, \os{\circ}{X}_{\lam_{2j+k}}\}$ 
(resp.~$\{\os{\circ}{Y}_{\lam_0}, \ldots, \os{\circ}{Y}_{\lam_{2j+k}}\}$). 
Assume that the divisibility condition for the morphism 
(\ref{eqn:odnl}) holds. 
Then consider the following direct factor of the cosimplicial 
degree $m$-part of the morphism (\ref{ali:gdkcl}): 
\begin{align*}
g^*_{\ul{\lam}} \col & 
b_{\phi(\ul{\lam})*}
Ru_{\os{\circ}{Y}_{\phi(\ul{\lam})}/\os{\circ}{T}{}'*}(F_{\phi(\ul{\lam})}\otimes_{\mab Z}
\vp_{\phi(\ul{\lam}){\rm crys}}(\os{\circ}{Y}_{T'_0}
/\os{\circ}{T}{}'))[-2j-k]  
\tag{1.5.8.4}\label{ali:grgm}\\
{} & \lo b_{\phi(\ul{\lam})*}
R\os{\circ}{g}_{\ul{\lam}*}
Ru_{\os{\circ}{X}_{\ul{\lam}}/\os{\circ}{T}*}(E_{\ul{\lam}}
\otimes_{\mab Z}
\vp_{\ul{\lam}{\rm crys}}
(\os{\circ}{X}_{T_0}/\os{\circ}{T}))[-2j-k]. 
\end{align*}

\begin{prop}\label{prop:grloc}
Let the notations 
and the assumptions be as above. 
Then the morphism $g^*_{\ul{\lam}}$ in {\rm (\ref{ali:grgm})} 
is equal to 
$\deg(u)^{j+k}b_{\phi(\ul{\lam})*}
\os{\circ}{g}{}^*_{\ul{\lam}}$ for $j\geq \max\{-k,0\}$. 
\end{prop}
\begin{proof}
(The proof is essentially 
the same as that of \cite[(2.9.3)]{nh2}.) 
As in the proof of (\ref{prop:dvok}), we may assume 
that $X_{\os{\circ}{T}_0}$ and $Y_{\os{\circ}{T}{}'_0}$ have (formal) log smooth lifts 
$\ol{\cal X}$ and $\ol{\cal Y}$
over $\ol{S(T)^{\nat}}$ and $\ol{S'(T')^{\nat}}$, respectively, 
such that there exists a morphism 
$\ol{\cal X}\lo \ol{\cal Y}$ over the morphism $\ol{S(T)^{\nat}}\lo \ol{S'(T')^{\nat}}$
which extends $g\col X_{\os{\circ}{T}_0}\lo Y_{\os{\circ}{T}{}'_0}$. 
Let the notations be as in the proof of (\ref{prop:dvok}).  
Consider the following morphism:
\begin{align*}
\deg(u)^{-(j+1)}g^{{\rm PD}*}\col & 
{\rm gr}^P_{2j+k+1}({\cal F}\otimes_{{\cal O}_{\cal Y}}
{\Om}^{i+j+1}_{{\cal Y}/\os{\circ}{T}{}'})/P_j
\lo g^{\rm PD}_*
({\rm gr}^P_{2j+k+1}({\cal E}\otimes_{{\cal O}_{{\cal X}}}
{\Om}^{i+j+1}_{{\cal X}/\os{\circ}{T}})/P_j).  \tag{1.5.9.1} \\ 
\end{align*} 
\par 
Let $\os{\circ}{\cal X}_{\lam_l}$ 
(resp.~$\os{\circ}{\cal Y}_{\phi(\lam_l)}$) 
be a closed subscheme 
of $\os{\circ}{\cal X}$ (resp.~$\os{\circ}{\cal Y}$) 
which is a lift of $\os{\circ}{X}_{\lam_l}$ 
(resp.~$\os{\circ}{Y}_{\phi(\lam_l)}$). 
Let $\wt{x}_l=0$ (resp.~$\wt{y}_l=0$) be a 
local equation of 
$\os{\circ}{\cal X}_{\lam_l}$  
$(0 \leq l \leq 2j+k)$ 
(resp.~$\os{\circ}{\cal Y}_{\phi(\lam_l)}$). 
Let $\hat{x}_l \in M_{\cal X}$ 
(resp.~$\hat{y}_l \in M_{\cal Y}$)  
be the local inverse image of $\wt{x}_l$ (resp.~$\wt{y}_l$). 
Let $\om =\eta d\log \hat{y}_0\cdots d\log 
\hat{y}_{2j+k}$ be a local section of
$P_{2j+k+1}({\cal F}\otimes_{{\cal O}_{\cal Y}}
{\Om}^{\bul}_{{\cal Y}/\os{\circ}{T}{}'})$ 
$(\eta \in {\rm Im}
({\cal F}\otimes_{{\cal O}_{\cal Y}}
\Om^{\bul-2j-k-1}_{\os{\circ}{\cal Y}/\os{\circ}{T}{}'}
\lo {\cal F}\otimes_{{\cal O}_{\cal Y}}
{\Om}^{\bul-2j-k-1}_{{\cal Y}/\os{\circ}{T}{}'}))$.  
Then, by the assumption (1.5.6.5) and by (\ref{prop:xxle}),  
we have the following equality: 
\begin{equation*} 
\deg(u)^{-(j+1)}g^{{\rm PD}*}(\om)
=\deg(u)^{j+k}
\prod_{l=1}^{2j+k+1}
g^{{\rm PD}*}(\eta)d\log \hat{x}_0 \cdots 
d\log \hat{x}_{2j+k}+\om',
\tag{1.5.9.2}\label{ali:ujg} 
\end{equation*}  
where 
$\om'\in P_{2j+k}({\cal E}\otimes_{{\cal O}_{{\cal X}}}
{\Om}^{\bul}_{{\cal X}/\os{\circ}{T}})$. 
For $\ul{\lam}:=\{\lam_0,\ldots,\lam_{2j+k}\}$ 
$(\lam_l\not=\lam_{l'}$ $(l\not=l'))$, 
set $\os{\circ}{\cal X}_{\ul{\lam}}:=
\os{\circ}{\cal X}_{\lam_0}\cap \cdots \cap 
\os{\circ}{\cal X}_{\lam_{2j+k}}$ 
and 
$\os{\circ}{\cal Y}_{\phi(\ul{\lam})}:=
\os{\circ}{\cal Y}_{\phi(\lam_0)}\cap \cdots \cap 
\os{\circ}{\cal Y}_{\phi(\lam_{2j+k})}$.  
Let 
$${\rm Res}^{\os{\circ}{\cal X}_{\ul{\lam}}}
\col 
{\rm gr}^P_{2j+k+1}({\cal E}
\otimes_{{\cal O}_{\cal X}}{\Om}^{\bul}_{{\cal X}/\os{\circ}{T}}) 
\lo 
{\cal E}_{\ul{\lam}}
\otimes_{{\cal O}_{{\cal X}{}_{\ul{\lam}}}}
\Om^{\bul}_{\os{\circ}{\cal X}{}_{\ul{\lam}}
/\os{\circ}{T}}\otimes_{\mab Z}
\vp_{\ul{\lam}{\rm zar}}
(\os{\circ}{\cal X}/\os{\circ}{T})$$ 
and 
$${\rm Res}^{\os{\circ}{\cal Y}_{\ul{\lam}}}
\col 
{\rm gr}^P_{2j+k+1}({\cal F}\otimes_{{\cal O}_{\cal Y}}
{\Om}^{\bul}_{{\cal Y}/\os{\circ}{T}{}'}) \lo 
{\cal F}_{\phi(\ul{\lam})}
\otimes_{{\cal O}_{{\cal Y}{}_{\phi(\ul{\lam})}}}
\Om^{\bul}_{\os{\circ}{\cal Y}_{\phi(\ul{\lam})}/\os{\circ}{T}{}'}
\otimes_{\mab Z}
\vp_{\phi(\ul{\lam}){\rm zar}}(\os{\circ}{\cal Y}/\os{\circ}{T}{}')$$  
be the residue morphisms with respect to 
$\os{\circ}{\cal X}_{\ul{\lam}}$ 
and $\os{\circ}{\cal Y}_{\phi(\ul{\lam})}$, respectively. 
By (\ref{ali:ujg}) we have the following commutative diagram:
\begin{equation*}
\begin{CD}
{\rm gr}_k^PA_{\rm zar}(Y_{\os{\circ}{T}{}'_0}/S'(T')^{\nat},F)
@>{{\rm gr}_k^P(g^*)}>> \\ 
@V{{\rm Res}^{\os{\circ}{\cal Y}_{\phi(\ul{\lam})}}}VV\\ 
b'_{\ul{\lam}*}
({\cal F}_{\ul{\lam}}
\otimes_{{\cal O}_{\os{\circ}{\cal Y}{}_{\phi(\ul{\lam})}}}
\Om^{\bul}_{\os{\circ}{\cal Y}_{\phi(\ul{\lam})}/\os{\circ}{T}{}'}
\otimes_{\mab Z}
\vp_{\phi(\ul{\lam}){\rm zar}}
(\os{\circ}{\cal Y}/\os{\circ}{T}{}')[-2j-k]
@>{\deg(u)^{j+k}g^{{\rm PD}*}_{\ul{\lam}}}>>  
\end{CD}
\tag{1.5.9.3}\label{cd:grtg}
\end{equation*} 
\begin{equation*}
\begin{CD}
g_{*}A_{\rm zar}(X_{\os{\circ}{T}_0}/S(T)^{\nat},E) \\
@VV{{\rm Res}^{\os{\circ}{\cal X}_{\ul{\lam}}}}V\\ 
b'_{\ul{\lam}*}g_{*}a_{\ul{\lam}*}
({\cal E}_{\ul{\lam}}
\otimes_{{\cal O}_{\os{\circ}{\cal X}_{\ul{\lam}}}}
\Om^{\bul}_{\os{\circ}{\cal X}_{\ul{\lam}}/\os{\circ}{T}}
\otimes_{\mab Z}\vp_{\ul{\lam}{\rm zar}}
(\os{\circ}{\cal X}/\os{\circ}{T})
[-2j-k]. 
\end{CD}
\end{equation*} 
This commutative diagram tells us that (\ref{prop:grloc}) holds. 
\end{proof}

\begin{defi}
Let $v \col {\cal E} \lo {\cal F}$ 
be a morphism of  
$f^{-1}({\cal O}_T)$-modules (resp.~${\cal O}_T$-modules). 
The $D$-{\it twist}(:=degree twist) by $k$  
$$v(-k) \col {\cal E}(-k,u)\lo {\cal F}(-k,u)$$  
of $v$ with respect to $u$ 
is, by definition, the morphism  
$\deg(u)^kv \col {\cal E} \lo {\cal F}$.  
This definition is well-defined for morphisms of objects of 
the derived category 
${\rm D}^+(f^{-1}({\cal O}_T))$ (resp.~${\rm D}^+({\cal O}_T)$). 
\end{defi}

\parno

\begin{coro}\label{coro:fuu}
The morphism 
\begin{align*} 
g^*_m\col 
{\rm gr}^P_kA_{\rm zar}(Y_{m,\os{\circ}{T}{}'_0}/S'(T')^{\nat},F^m)
\lo 
Rg_{m*}({\rm gr}^P_kA_{\rm zar}(X_{m,\os{\circ}{T}_0}/S(T)^{\nat},E^m))
\tag{1.5.11.1}\label{ali:cgrvp} 
\end{align*} 
is equal to 
\begin{align*} 
& \bigoplus_{j\geq \max \{-k,0\}} 
b^{(2j+k)}_{m,T'_0*} 
(Ru_{\os{\circ}{Y}{}^{(2j+k)}_{m,T'_0}/\os{\circ}{T}{}'*}
(F_{\os{\circ}{Y}{}^{(2j+k)}_{m,T'_0}/\os{\circ}{T}{}'} 
\otimes_{\mab Z}\vp_{\rm crys}^{(2j+k)}
(\os{\circ}{Y}_{m,T'_0}/T')))(-j-k,u)
\tag{1.5.11.2}\label{ali:cgrp}\\ 
&[-2j-k]\\ 
& \lo \\
&\bigoplus_{j\geq \max \{-k,0\}} 
\bigoplus_{\ul{\lam}\in \Lam^{(k)}(\os{\circ}{g})}
b_{m,\phi(\ul{\lam})*} 
(Ru_{\os{\circ}{Y}_{m,\phi(\ul{\lam})}/\os{\circ}{T}{}'*}
(F_{\os{\circ}{Y}_{m,\phi(\ul{\lam})}/\os{\circ}{T}{}'} 
\otimes_{\mab Z}\vp_{\phi(\ul{\lam}),{\rm crys}}
(\os{\circ}{Y}_{m,T'_0}/T')))\\
&(-j-k,u)[-2j-k] \\ 
& \os{\sum_{\ul{\lam}\in \Lam^{(k)}(\os{\circ}{g})}
\os{\circ}{g}{}^*_{\ul{\lam}}(-j-k,u)}{\lo} \\
&\bigoplus_{j\geq \max \{-k,0\}} 
\bigoplus_{{\ul{\lam}\in \Lam^{(k)}(\os{\circ}{g})}}
a_{m,\ul{\lam}*} 
(Ru_{\os{\circ}{X}_{m,\ul{\lam}}/\os{\circ}{T}*}
(E_{\os{\circ}{X}_{m,\ul{\lam}}/\os{\circ}{T}} 
\otimes_{\mab Z}\vp_{\ul{\lam},{\rm crys}}
(\os{\circ}{X}_{m,T_0}/T)))\\
&(-j-k,u)[-2j-k] \\ 
& \lo \\
&\bigoplus_{j\geq \max \{-k,0\}} 
a^{(2j+k)}_{m,T_0*} 
(Ru_{\os{\circ}{X}{}^{(2j+k)}_{m,T_0}/\os{\circ}{T}*}
(E_{\os{\circ}{X}{}^{(2j+k)}_{m,T_0}/\os{\circ}{T}} 
\otimes_{\mab Z}\vp_{\rm crys}^{(2j+k)}
(\os{\circ}{X}_{m,T_0}/T)))(-j-k,u)\\
&[-2j-k]. 
\end{align*} 
\end{coro} 
\begin{proof} 
This immediately follows from (\ref{prop:grloc}).  
\end{proof}

\begin{coro}\label{coro:indg}
Let  $h_{\bul \leq N}\col X_{\bul \leq N,\os{\circ}{T}_0}\lo Y_{\bul \leq N,\os{\circ}{T}{}'_0}$ 
be another morphism satisfying the condition {\rm (\ref{cd:xygxy})}, {\rm (1.5.6.4)} and 
{\rm (1.5.6.5)}. 
Assume that $\os{\circ}{g}_{\bul \leq N}=\os{\circ}{h}_{\bul \leq N}$. 
Then, for   for each $0\leq m\leq N$, 
\begin{align*} 
{\cal H}^q(h_m^*)={\cal H}^q(g_m^*) \col &
{\cal H}^q(P_kA_{\rm zar}(Y_{m\os{\circ}{T}{}'_0}/
S'(T')^{\nat},F^{m})) \\
&\lo {\cal H}^q(Rg_{m*}(P_kA_{\rm zar}(X_{m\os{\circ}{T}_0}/
S(T)^{\nat},E^m))) \quad (q\in {\mab N}) 
\end{align*} 
\end{coro}
\begin{proof} 
It suffices to prove that $g_m^*=h_m^*$ for each $0\leq m\leq N$. 
This is a local question on $\os{\circ}{Y}_{m,T'_0}$. 
Hence we may assume that $\os{\circ}{Y}_{m,T'_0}$ is quasi-compact. 
By the argument in (\ref{rema:asid}), it suffices to prove that 
${\rm gr}_k^P(h_m^*)={\rm gr}_k^P(g_m^*)$ $(k\in {\mab Z}, m\in {\mab N})$. 
This follows from (\ref{coro:fuu}). 
\end{proof}

\par  
Let $f \col X_{\bul \leq N,\os{\circ}{T}_0} \lo S(T)^{\nat}$ 
be the structural morphism. 
Let $(S(T)^{\nat})_{\bul \leq N}$ be the 
$N$-truncated constant simplicial log scheme defined by $S(T)^{\nat}$. 
Then $f$ induces the natural morphism 
$f_{\bul \leq N}
\col X_{\bul \leq N,\os{\circ}{T}_0} \lo (S(T)^{\nat})_{\bul \leq N}$. 
For simplicity of notation, denote by 
$(K^{\bul \leq N,\bul},P)$ the filtered complex 
$$Rf_{\bul \leq N *}
((A_{\rm zar}(X_{\bul \leq N,\os{\circ}{T}_0}/S(T)^{\nat},E^{\bul \leq N}),P))
\in 
{\rm D}^+{\rm F}({\cal O}_{T_{\bul \leq N}}).$$  
(The left degree of $K^{\bul \leq N,\bul}$ 
is the truncated cosimplicial degree corresponding 
to the truncated simplicial degree of $(S(T)^{\nat})_{\bul \leq N}$.).   
Then we have the following formula 
\begin{equation*} 
{\bf s}((K^{\bul \leq N,\bul},P))
=Rf_{*}((A_{\rm zar}(X_{\bul \leq N,\os{\circ}{T}_0}/S(T)^{\nat},
E^{\bul \leq N}),P)) 
\in {\rm D}^+{\rm F}({\cal O}_T).  
\end{equation*} 
Let $L$ be the stupid filtration on 
$(K^{\bul \leq N,\bul},P)$ 
with respect to the cosimplicial degree: 
\begin{equation*}
L^m(K^{\bul \leq N,\bul},P)
=\bigoplus_{m'=m}^N(K^{m'\bul},P)
\quad (m\in {\mab N}). 
\tag{1.5.12.1}\label{eqn:lassxs}
\end{equation*}
\par
Let 
$\del(L,P)$ be the diagonal filtration of $L$  and $P$
on ${\bf s}(K^{\bul \leq N,\bul})$
(cf.~\cite[(7.1.6.1), (8.1.22)]{dh3}): 
\begin{align*}
\del(L,P)_k(K^{\bul \leq N,\bul}) & =  
\bigoplus_{m=0}^NP_{k+m}K^{m\bul}. 
\tag{1.5.12.2}\label{ali:ddissfl} 
\end{align*} 
Then we have the following by (\ref{ali:cgrvp}):  
\begin{align*} 
& {\rm gr}^{\del(L,P)}_k{\bf s}(K^{\bul \leq N,\bul})
=\bigoplus_{m=0}^N\bigoplus_{j\geq \max \{-(k+m),0\}} 
Rf_{\os{\circ}{X}{}^{(2j+k+m)}_{m,T_0}/\os{\circ}{T}*}
(E^m_{\os{\circ}{X}{}^{(2j+k+m)}_{m,T_0}
/\os{\circ}{T}}
\tag{1.5.12.3}\label{ali:russrvp}\\
& \otimes_{\mab Z}\vp_{\rm crys}^{(2j+k+m)}
(\os{\circ}{X}_{m,T_0}/\os{\circ}{T})))(-j-k-m,u)[-2j-k-2m]. 
\end{align*} 
By (\ref{ali:russrvp}) 
we have the following spectral sequence 
\begin{align*} 
&E_1^{-k,q+k}=E_1^{-k,q+k}(X_{\bul \leq N,\os{\circ}{T}_0}/S(T)^{\nat}) 
\tag{1.5.12.4}\label{eqn:escssp} \\
&=\bigoplus_{m=0}^N 
\bigoplus_{j\geq \max \{-(k+m),0\}} 
R^{q-2j-k-2m}f_{\os{\circ}{X}{}^{(2j+k+m)}_{m,T_0}
/\os{\circ}{T}*}
(E^m_{\os{\circ}{X}{}^{(2j+k+m)}_{m,T_0}
/\os{\circ}{T}}
\otimes_{\mab Z} \\
&\phantom{R^{q-2j-k-m}f_{(\os{\circ}{X}^{(k)}, 
Z\vert_{\os{\circ}{X}^{(2j+k)}})/S*} 
({\cal O}}
  \vp^{(2j+k+m)}_{\rm crys}(
\os{\circ}{X}_{m,T_0}/\os{\circ}{T}))(-j-k-m,u) \\
&\Lo 
R^qf_{X_{\bul \leq N,\os{\circ}{T}_0}/S(T)^{\nat}*}
(\eps^*_{X_{\bul \leq N,\os{\circ}{T}_0}/S(T)^{\nat}}(E^{\bul \leq N}))  
\quad (q\in {\mab Z}).  
\end{align*}
More generally, for $k\in {\mab Z}$, 
we have the following spectral sequence 
\begin{align*} 
&E_1^{-k',q+k'}=E_1^{-k',q+k'}(X_{\bul \leq N,\os{\circ}{T}_0}/S(T)^{\nat})
\tag{1.5.12.5}\label{eqn:esasp} \\
&=\bigoplus_{m=0}^N 
\bigoplus_{j\geq \max \{-(k'+m),0\}} 
R^{q-2j-k'-2m}f_{\os{\circ}{X}{}^{(2j+k'+m)}_{m,T_0}
/\os{\circ}{T}*}
(E^m_{\os{\circ}{X}{}^{(2j+k'+m)}_{m,T_0}
/\os{\circ}{T}}
\otimes_{\mab Z} \\
&\phantom{R^{q-2j-k'-m}f_{(\os{\circ}{X}^{(k')}, 
Z\vert_{\os{\circ}{X}^{(2j+k')}})/S*} 
({\cal O}}\vp^{(2j+k'+m)}_{\rm crys}(
\os{\circ}{X}_{m,T_0}/\os{\circ}{T}))(-j-k'-m,u) \\
&\Lo 
R^qf_{X_{\bul \leq N,\os{\circ}{T}_0}/S(T)^{\nat}*}
(P_kA_{\rm zar}(X_{\bul \leq N,\os{\circ}{T}_0}/S(T)^{\nat},
E^{\bul \leq N})\quad (k'\leq k, q\in {\mab Z}).  
\end{align*}

\parno
If $T$ is restrictively hollow with respective to the morphism $T_0\lo S$, 
then we have the following spectral sequence 
by (\ref{eqn:escssp}) and (\ref{lemm:flpis}) (1): 
\begin{align*} 
&E_1^{-k,q+k}=E_1^{-k,q+k}(X_{\bul \leq N,T_0}/T) 
\tag{1.5.12.6}\label{eqn:esressp} \\
&=\bigoplus_{m=0}^N 
\bigoplus_{j\geq \max \{-(k+m),0\}} 
R^{q-2j-k-2m}f_{\os{\circ}{X}{}^{(2j+k+m)}_{m,T_0}
/\os{\circ}{T}*}
(E^m_{\os{\circ}{X}{}^{(2j+k+m)}_{m,T_0}
/\os{\circ}{T}}
\otimes_{\mab Z} \\
&\phantom{R^{q-2j-k-m}f_{(\os{\circ}{X}^{(k)}, 
Z\vert_{\os{\circ}{X}^{(2j+k)}})/S*} 
({\cal O}}
  \vp^{(2j+k+m)}_{\rm crys}(
\os{\circ}{X}_{m,T_0}/\os{\circ}{T}))(-j-k-m,u) \\
&\Lo 
R^qf_{X_{\bul \leq N,T_0}/T*}
(\eps^*_{X_{\bul \leq N,T_0}/T}(E^{\bul \leq N}))  
\quad (q\in {\mab Z}).  
\end{align*}

\begin{defi} 
(1) We call the spectral sequence (\ref{eqn:escssp}) (resp.~(\ref{eqn:esressp}))
the {\it Poincar\'{e} spectral sequence} of 
$\eps^*_{X_{\bul \leq N,\os{\circ}{T}_0}/S(T)^{\nat}}(E^{\bul \leq N})$ 
(resp.~$\eps^*_{X_{\bul \leq N,T_0}/T}(E^{\bul \leq N})$).
If $E^{\bul \leq N}={\cal O}_{\os{\circ}{X}_{\bul \leq N,T_0}/\os{\circ}{T}}$ 
and if $p$ is locally nilpotent on $\os{\circ}{T}$,  
then we call the spectral sequence (\ref{eqn:escssp}) 
(resp.~(\ref{eqn:esressp})) the {\it preweight spectral sequence}  of 
$X_{\bul \leq N,\os{\circ}{T}_0}/S(T)^{\nat}$ (resp.~$X_{\bul \leq N,T_0}/T$). 
If $E^{\bul \leq N}={\cal O}_{\os{\circ}{X}_{\bul \leq N,T_0}/\os{\circ}{T}}$ and 
if $\os{\circ}{T}$ is a flat formal ${\mab Z}_p$-scheme,  
then we call the spectral sequence (\ref{eqn:escssp}) 
(resp.~(\ref{eqn:esressp})) the {\it weight spectral sequence} of 
$X_{\bul \leq N,\os{\circ}{T}_0}/S(T)^{\nat}$ (resp.~$X_{\bul \leq N,T_0}/T$).
\par 
(2) We usually denote by 
$P$ (instead of $\del(L,P)$)  
the induced filtration on 
$$R^qf_{X_{\bul \leq N,\os{\circ}{T}_0}/S(T)^{\nat}*}
(\eps^*_{X_{\bul \leq N,\os{\circ}{T}_0}/S(T)^{\nat}}(E^{\bul \leq N}))
\quad 
({\rm resp}.~R^qf_{X_{\bul \leq N,T_0}/T*}
(\eps^*_{X_{\bul \leq N,T_0}/T}(E^{\bul \leq N})))$$ 
by the spectral sequence (\ref{eqn:escssp}) 
(resp.~(\ref{eqn:esressp}))
twisted by $q$ by abuse of notation. 
We call $P$ the {\it Poincar\'{e} filtration} on 
$$R^qf_{X_{\bul \leq N,\os{\circ}{T}_0}/S(T)^{\nat}*}
(\eps^*_{X_{\bul \leq N,\os{\circ}{T}_0}/S(T)^{\nat}}(E^{\bul \leq N}))
\quad 
({\rm resp}.~R^qf_{X_{\bul \leq N,T_0}/T*}
(\eps^*_{X_{\bul \leq N,T_0}/T}(E^{\bul \leq N}))).$$ 
If $E^{\bul \leq N}
={\cal O}_{\os{\circ}{X}_{\bul \leq N,T_0}/\os{\circ}{T}}$ 
and if $p$ is locally nilpotent on $\os{\circ}{T}$,  
then we call $P$ the {\it preweight filtration} on 
$$R^qf_{X_{\bul \leq N,\os{\circ}{T}_0}/S(T)^{\nat}*}
({\cal O}_{X_{\bul \leq N,\os{\circ}{T}_0}/S(T)^{\nat}})\quad 
({\rm resp}.~R^qf_{X_{\bul \leq N,T_0}/T*}
({\cal O}_{X_{\bul \leq N,T_0}/T})).$$  
If $E^{\bul \leq N}
={\cal O}_{\os{\circ}{X}_{\bul \leq N,T_0}/\os{\circ}{T}}$ 
and if $\os{\circ}{T}$ is a flat formal ${\mab Z}_p$-scheme,  
then we call $P$ the {\it weight filtration} on 
$$R^qf_{X_{\bul \leq N,\os{\circ}{T}_0}/S(T)^{\nat}*}
({\cal O}_{X_{\bul \leq N,\os{\circ}{T}_0}/S(T)^{\nat}})\quad 
({\rm resp}.~R^qf_{X_{\bul \leq N,T_0}/T*}
({\cal O}_{X_{\bul \leq N,T_0}/T})).$$  
\end{defi}

\begin{defi}[{\bf  Abrelative Frobenius morphism}]\label{defi:rwd}  
(1) Assume that $\os{\circ}{S}$ is of characteristic $p>0$.  
Let $\os{\circ}{F}_{S}\col \os{\circ}{S} \lo \os{\circ}{S}$ 
be the Frobenius endomorphism of $\os{\circ}{S}$.   
Set $S^{[p]}:=S\times_{\os{\circ}{S},\os{\circ}{F}_{S}}\os{\circ}{S}$. 
Then we have the following natural morphisms  
\begin{align*} 
F_{S/\os{\circ}{S}}\col S\lo S^{[p]}
\end{align*}  
and 
\begin{align*} 
W_{S/\os{\circ}{S}} \col S^{[p]}\lo S.
\end{align*} 
(The underlying morphism of the former morphism is ${\rm id}_{\os{\circ}{S}}$.)
Let $(T,{\cal J},\del)\lo (T',{\cal J}',\del')$ be a morphism of 
$p$-adic formal log PD-enlargements over the morphism $S\lo S^{[p]}$.  
Then we have the following natural morphisms  
\begin{align*} 
S_{\os{\circ}{T}_0}\lo S^{[p]}_{\os{\circ}{T}{}'_0}
\end{align*} 
and 
\begin{align*} 
(S(T)^{\nat},{\cal J},\del)\lo (S^{[p]}(T')^{\nat},{\cal J}',\del').
\end{align*}  
We call the morphisms $S_{\os{\circ}{T}_0}\lo S^{[p]}_{\os{\circ}{T}{}'_0}$ 
and 
$(S(T)^{\nat},{\cal J},\del)\lo (S^{[p]}(T')^{\nat},{\cal J}',\del')$ 
the {\it abrelative Frobenius morphism of base log schemes}  
and the {\it abrelative Frobenius morphism of base log PD-schemes}, respectively. 
(These Frobenius morphisms are essentially absolute 
in the logarithmic structures: these are relative in the scheme structures; 
``abrelative'' is a coined word; 
it means ``absolute and relative'' or ``far from being relative''.) 
In particular, when $(T',{\cal J}',\del')=(T,{\cal J},\del)$ 
with morphism $T_0\lo S$, 
we have the following natural morphisms 
\begin{align*} 
S_{\os{\circ}{T}_0}\lo S^{[p]}_{\os{\circ}{T}_0}
\end{align*} 
and 
\begin{align*} 
(S(T)^{\nat},{\cal J},\del)\lo (S^{[p]}(T)^{\nat},{\cal J},\del)
\end{align*}   
by using a composite morphism 
$T_0\lo S\os{W_{S/\os{\circ}{S}}}{\lo} S^{[p]}$.  
\par 
(2) Let the notations be as in (1).  
Set $X^{[p]}_{\bul \leq N}:=X_{\bul \leq N}\times_SS^{[p]}
=X_{\bul \leq N}\times_{\os{\circ}{S},\os{\circ}{F}_S}\os{\circ}{S}$ 
and 
\begin{align*} 
X^{[p]}_{\bul \leq N,\os{\circ}{T}{}'_0}&:=
X^{[p]}_{\bul \leq N}
\times_{S^{[p]}}S^{[p]}_{\os{\circ}{T}{}'_0}
=X^{[p]}_{\bul \leq N}
\times_{\os{\circ}{S^{[p]}}}\os{\circ}{T}{}'_0. 
\end{align*} 
Then $X^{[p]}_{\bul \leq N,\os{\circ}{T}{}'_0}/S^{[p]}_{\os{\circ}{T}{}'_0}$ 
is an SNCL scheme. 
Let 
$$F^{\rm ar}
_{X_{\bul \leq N,\os{\circ}{T}_0/S_{\os{\circ}{T}_0},S^{[p]}_{\os{\circ}{T}{}'_0}}}\col 
X_{\bul \leq N,\os{\circ}{T}_0}  \lo X^{[p]}_{\bul \leq N,\os{\circ}{T}{}'_0}$$ 
and 
$$F^{\rm ar}_{X_{\bul \leq N,\os{\circ}{T}_0/S(T)^{\nat},S^{[p]}(T')^{\nat}}}
\col 
X_{\bul \leq N,\os{\circ}{T}_0}  \lo X^{[p]}_{\bul \leq N,\os{\circ}{T}{}'_0}$$ 
be the natural morphisms   
over $S_{\os{\circ}{T}_0}\lo (S^{[p]})_{\os{\circ}{T}{}'_0}$ and 
$(S(T)^{\nat},{\cal J},\del)\lo (S^{[p]}(T')^{\nat},{\cal J}',\del')$. 
We call 
$F^{\rm ar}_{X_{\bul \leq N,\os{\circ}{T}_0/S_{\os{\circ}{T}_0},S^{[p]}_{\os{\circ}{T}{}'_0}}}$ 
and 
$F^{\rm ar}_{X_{\bul \leq N,\os{\circ}{T}_0/S(T)^{\nat},S^{[p]}(T')^{\nat}}}$ 
the {\it abrelative Frobenius morphisms}   
of $X_{\bul \leq N,\os{\circ}{T}_0}$ over 
$S_{\os{\circ}{T}_0}\lo S^{[p]}_{\os{\circ}{T}{}'_0}$ 
and $(S(T)^{\nat},{\cal J},\del)\lo (S^{[p]}(T')^{\nat},{\cal J}',\del')$, respectively. 
\par 
Assume that ${\cal O}_T$ is $p$-torsion-free and that ${\cal J}\subset p{\cal O}_T$.  
Let $E^{\bul \leq N}$ and $E'{}^{\bul \leq N}$ be a flat quasi-coherent crystal of 
${\cal O}_{\os{\circ}{X}_{\bul \leq N,T_0}/\os{\circ}{T}}$-modules and 
a flat quasi-coherent crystal of 
${\cal O}_{\os{\circ}{X}{}^{[p]}_{\bul \leq N,T'_0}/\os{\circ}{T}{}'}$-modules, 
respectively.  
Let 
\begin{align*} 
\Phi^{\rm ar}\col 
\os{\circ}{F}{}^{{\rm ar}*}_{X_{\bul \leq N,\os{\circ}{T}_0/
S(T)^{\nat},S^{[p]}(T')^{\nat},{\rm crys}}}(E'{}^{\bul \leq N})
\lo E^{\bul \leq N}
\tag{1.5.14.1}\label{ali:sppts}
\end{align*} 
be a morphism of crystals in 
$(\os{\circ}{X}_{\bul \leq N,T_0}/\os{\circ}{T})_{\rm crys}$.   
Since $\deg (F_{S(T)^{\nat}/S^{[p]}(T')^{\nat}})=p$, 
the divisibility of the morphism (\ref{eqn:odnl}) holds by (\ref{prop:dvok}) 
if ${\cal I}\subset p{\cal O}_T$.  
We call the following induced morphism by $\Phi^{\rm ar}$   
\begin{align*} 
\Phi^{\rm ar} \col &
(A_{\rm zar}(X^{[p]}_{\bul \leq N,\os{\circ}{T}{}'_0}/S^{[p]}(T')^{\nat},E'{}^{\bul \leq N}),P) 
\tag{1.5.14.2}\label{ali:spapts} \\
&\lo 
RF^{\rm ar}_{X_{\bul \leq N,\os{\circ}{T}_0/S(T)^{\nat},S^{[p]}(T')^{\nat}}*}
((A_{\rm zar}(X_{\bul \leq N,\os{\circ}{T}_0}/S(T)^{\nat},E^{\bul \leq N}),P)) 
\end{align*}
the {\it abrelative Frobenius morphism} of 
$$(A_{\rm zar}(X_{\bul \leq N,\os{\circ}{T}_0}/S(T)^{\nat},E^{\bul \leq N}),P)\quad 
{\rm and} \quad  
(A_{\rm zar}(X^{[p]}_{\bul \leq N,\os{\circ}{T}{}'_0}/S^{[p]}(T')^{\nat},E'^{\bul \leq N}),P).$$ 
When $E'{}^{\bul \leq N}={\cal O}_{\os{\circ}{X}{}^{[p]}_{\bul \leq N,T'_0}/\os{\circ}{T}{}'}$, 
we set 
$$(A_{\rm zar}(X^{[p]}_{\bul \leq N,\os{\circ}{T}{}'_0}/S^{[p]}(T')^{\nat}),P) 
:=(A_{\rm zar}(X^{[p]}_{\bul \leq N,\os{\circ}{T}{}'_0}/S^{[p]}(T')^{\nat},E'{}^{\bul \leq N}),P).$$  
Then we have the following {\it  abrelative Frobenius morphism} 
\begin{equation*} 
\Phi^{\rm  ar} \col 
(A_{\rm zar}(X^{[p]}_{\bul \leq N,\os{\circ}{T}{}'_0}/S^{[p]}(T')^{\nat}),P) 
\lo RF^{\rm  ar}_{X_{\bul \leq N,\os{\circ}{T}_0/S(T)^{\nat},S^{[p]}(T')^{\nat}}*}
((A_{\rm zar}(X_{\bul \leq N,\os{\circ}{T}_0}/S(T)^{\nat}),P)) 
\tag{1.5.14.3}
\end{equation*}
of $(A_{\rm zar}(X_{\bul \leq N,\os{\circ}{T}_0}/S(T)^{\nat}),P)$ and 
$(A_{\rm zar}(X^{[p]}_{\bul \leq N,\os{\circ}{T}{}'_0}/S^{[p]}(T')^{\nat}),P)$. 
\end{defi}

\begin{prop}[{\bf Frobenius compatibility I}]\label{prop:fcbar} 
The following diagram is commutative$:$ 
\begin{equation*} 
\begin{CD} 
A_{\rm zar}(X^{[p]}_{\bul \leq N,\os{\circ}{T}{}'_0}/S^{[p]}(T')^{\nat},E'{}^{\bul \leq N})
@>{\Phi^{\rm ar}}>>
 \\
@A{\theta_{X^{[p]}_{\bul \leq N,\os{\circ}{T}_0}/S^{[p]}(T')^{\nat}} \wedge}A{\simeq}A  \\
Ru_{X^{[p]}_{\bul \leq N,\os{\circ}{T}{}'_0}/S^{[p]}(T')^{\nat}*}
(\eps^*_{X^{[p]}_{\bul \leq N,\os{\circ}{T}{}'_0}/S^{[p]}(T')^{\nat}}(E'{}^{\bul \leq N}))
@>{\Phi^{\rm ar}}>>
\end{CD}
\tag{1.5.15.1}\label{cd:rukt}
\end{equation*} 
\begin{equation*} 
\begin{CD} 
RF^{\rm ar}_{X_{\bul \leq N,\os{\circ}{T}_0/S(T)^{\nat},S^{[p]}(T')^{\nat}}*}
(A_{\rm zar}(X_{\bul \leq N,\os{\circ}{T}_0}/S(T)^{\nat},E^{\bul \leq N})) \\
@A{RF^{\rm ar}_{X_{\bul \leq N,
\os{\circ}{T}_0/S(T)^{\nat},S^{[p]}(T')^{\nat}}*}
(\theta_{X_{\bul \leq N,\os{\circ}{T}_0}/S(T)^{\nat}})\wedge}A{\simeq}A 
\\
RF^{\rm ar}_{X_{\bul \leq N,\os{\circ}{T}_0/S(T)^{\nat},S^{[p]}(T')^{\nat}}*}
Ru_{X_{\bul \leq N,\os{\circ}{T}_0}/S(T)^{\nat}*}(E^{\bul \leq N}). 
\end{CD}
\end{equation*} 
This is contravariantly functorial 
for the morphism {\rm (\ref{eqn:xdxduss})} satisfying {\rm (\ref{cd:xygxy})} 
and for the morphism of $F$-crystals 
\begin{equation*} 
\begin{CD} 
\os{\circ}{F}{}^{{\rm ar}*}_{X_{\bul \leq N,\os{\circ}{T}_0/
S(T)^{\nat},S^{[p]}(T')^{\nat},{\rm crys}}}
(E'{}^{\bul \leq N})
@>{\Phi^{\rm ar}}>> E^{\bul \leq N}\\
@AAA @AAA \\
\os{\circ}{g}{}^*_{\bul \leq N}
\os{\circ}{F}{}^{{\rm ar}*}_{Y_{\bul \leq N,\os{\circ}{T}_0
/S(T)^{\nat},S^{[p]}(T')^{\nat},{\rm crys}}}
(F'{}^{\bul \leq N})
@>{\os{\circ}{g}{}^*_{\bul \leq N}(\Phi^{\rm ar})}>> 
\os{\circ}{g}{}^*_{\bul \leq N}(F^{\bul \leq N}), 
\end{CD}
\end{equation*} 
where $F^{\bul \leq N}$ $($resp.~$F'{}^{\bul \leq N})$ is 
a similar quasi-coherent 
${\cal O}_{\os{\circ}{Y}_{\bul \leq N,T_0}/\os{\circ}{T}}$-module 
to $E^{\bul \leq N}$ 
$($resp.~a similar quasi-coherent 
${\cal O}_{\os{\circ}{Y}{}^{[p]}_{\bul \leq N,T_0}/\os{\circ}{T}}$-module to $E'{}^{\bul \leq N})$. 
\end{prop} 
\begin{proof} 
This is a special case of (\ref{theo:funpas}). 
\end{proof} 

\begin{prop}[{\bf Frobenius compatibility II}]\label{prop:nwsfc} 
Assume that $T$ is restrictively hollow with respective to the morphism $T_0\lo S$.  
Assume that the morphism $T'_0\lo S^{[p]}$ factors through the morphism 
$F_{S/{\os{\circ}{S}}}\col S\lo S^{[p]}$.  
Set $X^{\{p\}}_{\bul \leq N}
:=X_{\bul \leq N}\times_{S,F_S}S$ 
and $X^{\{p\}}_{\bul \leq N,\os{\circ}{T}{}'_0}:=
X^{\{p\}}_{\bul \leq N}\times_SS_{\os{\circ}{T}{}'_0}$. 
Let 
$$F^{\rm rel}_{X_{\bul \leq N,\os{\circ}{T}_0/
S_{\os{\circ}{T}_0},S_{\os{\circ}{T}{}'_0}}}
\col X_{\bul \leq N,\os{\circ}{T}_0}  \lo X^{\{p\}}_{\bul \leq N,\os{\circ}{T}{}'_0}$$ 
and  
\begin{align*} 
F^{{\rm rel}*}_{X_{\bul \leq N,T_0/T,T'}}\col 
X_{\bul \leq N,T_0}\lo X_{\bul \leq N,T'_0}
\end{align*} 
be the relative Frobenius morphisms   
over $S_{\os{\circ}{T}_0}\lo S_{\os{\circ}{T}{}'_0}$ and 
$(T,{\cal J},\del)\lo (T',{\cal J},\del)$, respectively. 
Let $E''{}^{\bul \leq N}$ be the pull-back of 
$E'{}^{\bul \leq N}$ to $(\os{\circ}{X}{}^{\{p\}}_{\bul \leq N,T_0}/\os{\circ}{T})_{\rm crys}$. 
Let 
\begin{align*} 
\Phi^{{\rm rel}} \col 
\os{\circ}{F}{}^{{\rm  rel}*}_{X_{\bul \leq N,\os{\circ}{T}_0/T,T',{\rm crys}}}
(E''{}^{\bul \leq N})\lo E^{\bul \leq N}  
\tag{1.5.16.1}\label{ali:spptpps}
\end{align*} 
be the induced morphism by {\rm (\ref{ali:sppts})}. 
Then the following diagram is commutative$:$ 
\begin{equation*} 
\begin{CD} 
A_{\rm zar}(X^{[p]}_{\bul \leq N,\os{\circ}{T}{}'_0}
/S^{[p]}(T')^{\nat},E'{}^{\bul \leq N}) 
@>{\Phi^{\rm  rel}}>>
RF^{{\rm rel}}_{X_{\bul \leq N,\os{\circ}{T}_0/T'T}*}
(A_{\rm zar}(X_{\bul \leq N,\os{\circ}{T}_0}/S(T)^{\nat},E^{\bul \leq N})) \\
@A{\theta_{X^{[p]}_{\bul \leq N,T'_0}/T'} \wedge}A{\simeq}A 
@A{RF^{\rm rel}_{X_{\bul \leq N,\os{\circ}{T}_0/T,T'}*}
(\theta_{X_{\bul \leq N,\os{\circ}{T}_0/T}})\wedge}A{\simeq}A \\
Ru_{X^{\{p\}}_{\bul \leq N,T'_0}/T'*}(\eps^*_{X^{\{p\}}_{\bul \leq N,T'_0}/T'}
(E''{}^{\bul \leq N}))
@>{\Phi^{\rm rel}}>>RF^{\rm rel}_{X_{\bul \leq N,T_0/T,T'}*}
Ru_{X_{\bul \leq N,T_0}/T*}
(\eps^*_{X_{\bul \leq N,T_0}/T}(E^{\bul \leq N})). 
\end{CD}
\tag{1.5.16.2}\label{cd:rekt}
\end{equation*} 
The commutative diagram {\rm (\ref{cd:rekt})} is contravariantly functorial 
for the morphism {\rm (\ref{eqn:xdxduss})} satisfying {\rm (\ref{cd:xygxy})} 
and for the morphism of $F$-crystals 
\begin{equation*} 
\begin{CD} 
\os{\circ}{F}{}^{{\rm rel}*}_{X_{\bul \leq N,\os{\circ}{T}_0/T,T',{\rm crys}}}
(E'{}^{\bul \leq N})
@>{\Phi^{\rm rel}}>> E^{\bul \leq N}\\
@AAA @AAA \\
\os{\circ}{g}{}^*_{\bul \leq N}
\os{\circ}{F}{}^{{\rm rel}*}_{Y_{\bul \leq N,\os{\circ}{T}_0/T,T',{\rm crys}}}
(F'{}^{\bul \leq N})
@>{\os{\circ}{g}{}^*_{\bul \leq N}(\Phi^{\rm rel})}>> 
\os{\circ}{g}{}^*_{\bul \leq N}(F^{\bul \leq N}), 
\end{CD}
\end{equation*} 
where $F^{\bul \leq N}$
$($resp.~$F'{}^{\bul \leq N})$ is 
a similar quasi-coherent 
${\cal O}_{\os{\circ}{Y}_{\bul \leq N,T_0}/\os{\circ}{T}}$-module 
to $E^{\bul \leq N}$ 
$($resp.~${\cal O}_{\os{\circ}{X}{}^{[p]}_{\bul \leq N,T_0}/\os{\circ}{T}}$-module 
to $E'{}^{\bul \leq N})$. 
\end{prop}
\begin{proof} 
Because the morphism 
$T'_0\lo S^{[p]}$ is the composite morphism 
$T'_0\lo S_{\os{\circ}{T}_0}\lo S^{[p]}_{\os{\circ}{T}_0}$, 
$\os{\circ}{X}{}^{[p]}_{\bul \leq N,\os{\circ}{T}{}'_0}
=\os{\circ}{X}{}^{\{p\}}_{\bul \leq N,T'_0}$. 
Consequently  
\begin{align*}
Ru_{X^{[p]}_{\bul \leq N,\os{\circ}{T}_0}/S(T)^{\nat}*}
(\eps^*_{X^{[p]}_{\bul \leq N,\os{\circ}{T}_0}/S(T)^{\nat}}
(E'{}^{\bul \leq N})) 
=
Ru_{X^{\{p\}}_{\bul \leq N,T_0}/T*}(\eps^*_{X^{\{p\}}_{\bul \leq N,T_0}/T}
(E''{}^{\bul \leq N})). 
\end{align*} 
We also have the following equality: 
\begin{align*}
Ru_{X_{\bul \leq N,\os{\circ}{T}_0}/S(T)^{\nat}*}
(\eps^*_{X_{\bul \leq N,\os{\circ}{T}_0}/S(T)^{\nat}}
(E^{\bul \leq N})) 
=
Ru_{X_{\bul \leq N,T_0}/T*}
(\eps^*_{X_{\bul \leq N,T_0}/T}(E^{\bul \leq N})). 
\end{align*} 
Hence we obtain the commutative diagram 
(\ref{cd:rekt}) by (\ref{cd:rukt}) and (\ref{lemm:flpis}). 
\par 
We leave the proof of the functoriality to the reader. 
\end{proof}

\begin{defi}[{\bf Absolute Frobenius endomorphism}]\label{defi:btd}  
Let the notations and the assumptions be as in (\ref{defi:rwd}).  
Let $F_{S} \col S \lo S$ be the Frobenius endomorphism of $S$, 
that is, $\os{\circ}{F}_{S}\col \os{\circ}{S} \lo \os{\circ}{S}$ 
is induced by the $p$-th power endomorphism of ${\cal O}_{S}$ 
and the multiplication by $p$ of the log structure of $S$.   
Let $F_{S_{\os{\circ}{T}_0}} \col S_{\os{\circ}{T}_0} 
\lo S_{\os{\circ}{T}_0}$ be 
the Frobenius endomorphism of $S_{\os{\circ}{T}_0}$.  
Assume that there exists a lift $F_{S(T)^{\nat}} \col S(T)^{\nat}\lo S(T)^{\nat}$ of 
$F_{S_{\os{\circ}{T}_0}}$ 
which gives a PD-morphism 
$F_{S(T)^{\nat}}\col (S(T)^{\nat},{\cal J},\del)\lo (S(T)^{\nat},{\cal J},\del)$.  
Let 
$$F^{\rm abs}_{X_{\bul \leq N,\os{\circ}{T}_0}/S(T)^{\nat}} 
\col X_{\bul \leq N,\os{\circ}{T}_0}  \lo X_{\bul \leq N,\os{\circ}{T}_0}$$ 
be the absolute Frobenius endomorphism over $F_{S(T)^{\nat}}$.   
Let 
\begin{align*} 
\Phi^{\rm abs} \col 
\os{\circ}{F}{}^{{\rm abs}*}_{X_{\bul \leq N,\os{\circ}{T}_0/S(T)^{\nat}},{\rm crys}}
(E^{\bul \leq N})\lo E^{\bul \leq N}
\end{align*} 
be a morphism of crystals in 
$(\os{\circ}{X}_{\bul \leq N,T_0}/\os{\circ}{T})_{\rm crys}$.   
Then the divisibility of the morphism (\ref{eqn:odnl}) 
holds in this situation by (\ref{prop:dvok}) for the case $n=1$ if ${\cal I}\subset p{\cal O}_T$.  
Then we call the induced morphism by $\Phi^{\rm abs}$ and $F_{S(T)^{\nat}}$
\begin{align*} 
\Phi^{\rm abs} \col 
&(A_{\rm zar}(X_{\bul \leq N,\os{\circ}{T}_0}/S(T)^{\nat},E^{\bul \leq N}),P) 
\lo \tag{1.5.17.1}\label{eqn:ebnp}\\
&RF^{\rm abs}_{X_{\bul \leq N,\os{\circ}{T}_0/S(T)^{\nat}}*}
((A_{\rm zar}(X_{\bul \leq N,\os{\circ}{T}_0}/S(T)^{\nat},E^{\bul \leq N}),P)) 
\end{align*}
the {\it absolute Frobenius endomorphism} of 
$(A_{\rm zar}(X_{\bul \leq N,\os{\circ}{T}_0}/S(T)^{\nat},E^{\bul \leq N}),P)$ 
with respect to $F_{S(T)^{\nat}}$.
When $E^{\bul \leq N}={\cal O}_{\os{\circ}{X}_{\bul \leq N,T_0}/\os{\circ}{T}}$, 
we have the following {\it absolute Frobenius endomorphism} 
\begin{equation*} 
\Phi^{\rm abs}{}^* \col 
(A_{\rm zar}(X_{\bul \leq N,\os{\circ}{T}_0}/S(T)^{\nat}),P) 
\lo RF^{\rm abs}_{X_{\bul \leq N,\os{\circ}{T}_0/S(T)^{\nat}}*}
((A_{\rm zar}(X_{\bul \leq N,\os{\circ}{T}_0}/S(T)^{\nat}),P)) 
\tag{1.5.17.2}\label{eqn:abst}
\end{equation*}
of $(A_{\rm zar}(X_{\bul \leq N,\os{\circ}{T}_0}/S(T)^{\nat}),P)$ 
with respect to $F_{S(T)^{\nat}}$.  
\end{defi}  

\begin{rema} 
We leave the formulation of the analogue of 
(\ref{prop:nwsfc}) for the absolute Frobenius endomorphism. 
\end{rema}

\par 
Let the notations be as in (\ref{eqn:esasp}). 
\par 
When $\os{\circ}{S}$ is of characteristic $p>0$, when 
$u\col (S(T)^{\nat},{\cal J},\del)\lo (S'(T')^{\nat},{\cal J}',\del')$ 
is a lift of the abrelative Frobenius morphism 
$F_{S_{\os{\circ}{T}_0}}\col S_{\os{\circ}{T}_0}\lo S^{[p]}_{\os{\circ}{T}{}'_0}$ 
and 
when $g_{\bul \leq N}$ is 
the abrelative Frobenius morphism  
$F^{\rm ar}_{X_{\bul \leq N,\os{\circ}{T}_0/S(T)^{\nat},S^{[p]}(T')^{\nat}}}\col 
X_{\bul \leq N,\os{\circ}{T}_0} \lo X^{[p]}_{\bul \leq N,\os{\circ}{T}_0}$  or 
the absolute Frobenius endomorphism 
of $X_{\bul \leq N,\os{\circ}{T}_0}$, 
we denote $(-j-k-m,u)$ in (\ref{eqn:escssp}), 
(\ref{eqn:esasp}) and (\ref{eqn:esressp}) by $(-j-k-m)$ as usual.

\begin{rema}\label{rema:pob} 
(1) As pointed out before \cite[(11.9)]{ndw}, 
the Frobenius action is not considered 
in the proof of \cite[3.18]{msemi}; 
the action on a similar complex to  
$A_{\rm zar}(X_{\bul \leq N}/S)$ for the case $N=0$ 
in \cite{gkwf} is not either considered in [loc.~cit.]. 
The action (\ref{eqn:axdwt}) is indispensable. 
(Hence all papers before \cite{ndw} considering the Frobenius action 
on the analogues of $A_{\rm zar}(X_{\bul \leq N}/S)$ are incomplete.) 
For example, strictly speaking, 
the $E_2$-degeneration of the spectral sequences in \cite{msemi} 
was not proved since the Frobenius action was not defined in [loc.~cit.] 
and since we cannot use the yoga of weight. 
\par 
(2) 
Let the notations be as in \cite[(4.1)]{gkc}. 
The Frobenius action $\Phi$ on $A^{\bul}$ defined in \cite[(4.1)]{gkc} is mistaken. 
It seems that the operator $F$ in \cite[(9.8.12)]{ndw} is misunderstood by 
the author of \cite{gkc}. The $F\col W_{n+1}A^{ij}_X\lo W_nA^{ij}_X$ 
in \cite[(9.8.12)]{ndw} does not come from a morphism of log schemes, while 
the morphism $F\col A^{\bul}\lo A^{\bul}$ in \cite[(4.1)]{gkc} 
comes from a morphism of log schemes. 
In addition, the proof of the following isomorphism 
\begin{align*} 
{\rm Gr}_kA^{\bul}\simeq \bigoplus_{j\geq {\rm max}\{0,-k\}}
{\mab C}\Om^{\bul}_{Y^{2j+k+1}}[-2j-k](-j-k)
\end{align*} 
has not been given in \cite{gkc}. If one uses the action of 
Frobenius action $\Phi$ on $A^{\bul}$ in \cite[(4.1)]{gkc}, 
the Tate twist $(-j-k)$ will be different. 
Consequently the construction of the weight spectral sequence 
\begin{align*} 
E_1^{-k,i+k}=\simeq \bigoplus_{j\geq {\rm max}\{0,-k\}}
H^{i-2j-k}(Y,{\mab C}\Om^{\bul}_{Y^{2j+k+1}})(-j-k)\Lo H^i_{\rm crys}(Y/T)_{\mab Q}
\end{align*} 
in \cite[(4.1)]{gkc} is not complete. Moreover the phrase ``Passing to $n\to \inf$ and tensoring with ${\mab Q}$'' in the construction of the weight spectral sequence 
is strange because we do not need such processes. 
\end{rema}

\begin{coro}[{\bf Contravariant functoriality 
of the (pre)weight spectral sequence(a generalization of the $p$-adic analogue of 
{\rm \cite[Corollary 2.12]{stwsl}})}]\label{prop:contwt}
Let the notations and the assumption be as above 
and in {\rm (\ref{theo:funas})} or {\rm (\ref{theo:funpas})}.  
Assume that the two conditions 
$(1.5.6.4)$ and $(1.5.6.5)$ hold. 
Let $\os{\circ}{g}{}^{(k)*}_{m}$ be the following  morphism$:$
\begin{align*}  
\os{\circ}{g}{}^{(k)*}_{m}
& :=\sum_{\ul{\lam}\in \Lam^{(k)}(\os{\circ}{g})}
\os{\circ}{g}{}^*_{\ul{\lam}}
\col 
R^qf_{\os{\circ}{Y}{}^{(k)}_{m,T'_0}
/\os{\circ}{T}{}'*}
(F_{\os{\circ}{Y}{}^{(k)}_{m,T'_0}/\os{\circ}{T}{}'}
\otimes_{\mab Z}
\vp^{(k)}_{\rm crys}(\os{\circ}{Y}_{m,T'_0}
/\os{\circ}{T}{}')) \tag{1.5.20.1}\label{eqn:rglm} \\
& \lo 
R^qf_{\os{\circ}{X}{}^{(k)}_{m,T_0}/\os{\circ}{T}{}*}
(E_{\os{\circ}{X}{}^{(k)}_{m,T_0}
/\os{\circ}{T}}
\otimes_{\mab Z}
\vp^{(k)}_{\rm crys}(\os{\circ}{X}_{m,T_0}/T)). 
\end{align*} 
Then there exists the following morphism of 
$($pre$)$weight spectral sequences$:$ 
\begin{align*} 
E_1^{-k,q+k}& = \bigoplus_{m=0}^N
\bigoplus_{j\geq 
\max \{-(k+m),0\}} 
R^{q-2j-k-2m}
f_{\os{\circ}{X}{}^{(2j+k+m)}_{m,T_0}
/\os{\circ}{T}*}
(E^m_{\os{\circ}{X}{}^{(2j+k+m)}_{m,T_0}
/\os{\circ}{T}}
\otimes_{\mab Z} 
\tag{1.5.20.2}\label{eqn:espwfsp} \\
& \vp^{(2j+k+m)}_{\rm crys}
(\os{\circ}{X}_{m,T_0}/\os{\circ}{T}))(-j-k-m,u) 
\Lo 
R^qf_{X_{\bul \leq N,\os{\circ}{T}_0}/S(T)^{\nat}*}
(\eps^*_{X_{\bul \leq N,\os{\circ}{T}_0}/S(T)^{\nat}}(E^{\bul \leq N}))  
\end{align*}
\begin{equation*} 
\begin{CD} 
{ \quad \quad \quad \quad \quad \quad \quad \quad \quad} @. 
{ \quad \quad \quad \quad \quad \quad}\\ 
@A{\bigoplus_{m=0}^N
\bigoplus_{j\geq \max \{-(k+m),0\}}
{\rm deg}(u)^{j+k}\os{\circ}{g}{}^{(2j+k+m)*}_{m}}AA  
@. @. @AA{g^*_{\bul \leq N}}A \\
\end{CD} 
\end{equation*} 
\begin{align*} 
E_1^{-k,q+k}& = \bigoplus_{m\geq 0}
\bigoplus_{j\geq \max \{-(k+m),0\}} 
R^{q-2j-k-2m}
f_{\os{\circ}{Y}^{(2j+k+m)}_{m,T'_0}
/\os{\circ}{T}{}'*}
(F^m_{\os{\circ}{Y}^{(2j+k+m)}_{m,T'_0}
/\os{\circ}{T}{}'}
\otimes_{\mab Z}  \\
& \vp^{(2j+k+m)}_{\rm crys}(\os{\circ}{Y}_{m,T'_0}
/\os{\circ}{T}{}'))(-j-k-m,u) \Lo 
R^qf_{Y_{\bul \leq N,\os{\circ}{T}{}'_0}/S'(T')^{\nat}*}
(\eps^*_{Y_{\bul \leq N,\os{\circ}{T}{}'_0}/S'(T')^{\nat}}(F^{\bul \leq N})). 
\end{align*}  
\end{coro}
\begin{proof} 
This follows from (\ref{theo:funpas}), (\ref{prop:grloc}) 
and from the definition of $\os{\circ}{g}{}^{(l)}_{m}$ 
$(l\in {\mab N})$.  
\end{proof}


\par
Next we describe the boundary morphism 
between the $E_1$-terms of (\ref{eqn:escssp}).
Fix $m\in {\mab N}$ and a decomposition 
$\Del_m:=\{\os{\circ}{X}_{\lam_m}\}_{\lam_m}$ 
of $\os{\circ}{X}_{m,T_0}$ 
by smooth components of $\os{\circ}{X}_{m,T_0}$. 
To describe the boundary morphism, 
we fix a total order on $\Lam$ once and for all 
(cf.~(\ref{eqn:odpaps})). 
Let $\al$ be a nonnegative integer.
Set 
$\ul{\lam}_{m}:=\{\lam_{m0},\ldots, \lam_{m\al}\}$ 
$(\lam_{mi}<\lam_{mj}~(i<j))$,
$\ul{\lam}_{m\bet}:=
\{\lam_{m0},\ldots,  
\widehat{\lam}_{m\bet},  \ldots, \lam_{m\al}\}$ 
$(0 \leq \bet \leq \al)$,
$\os{\circ}{X}_{\ul{\lam}_m}:=
\os{\circ}{X}_{\lam_{m0}} \cap 
\cdots 
\cap \os{\circ}{X}_{\lam_{m\al}}$
and 
$\os{\circ}{X}_{\ul{\lam}_{m\bet}}
:=\os{\circ}{X}_{\lam_{m0}} \cap 
\cdots 
\cap \widehat{\os{\circ}{X}}_{{\lam}_{m\bet}} \cap 
\cdots 
\cap 
\os{\circ}{X}_{\lam_{m\al}}$. 
Here $~~\widehat{}~~$ means the elimination. 
Then $\os{\circ}{X}_{\ul{\lam}_{m}}$ 
is a smooth divisor on 
$\os{\circ}{X}_{\ul{\lam}_{m\bet}}/\os{\circ}{T}_0$.
For a nonnegative integer $k$ and an integer $l$, 
let
\begin{equation*}
(-1)^{\bet}G^{\ul{\lam}_{m\bet}}_{\ul{\lam}_m}
\col 
R^kf_{\os{\circ}{X}_{\ul{\lam}_{m\bet}}
/\os{\circ}{T}*}
(E_{\ul{\lam}_{m\bet}}
\otimes_{\mab Z} 
\vp_{\rm crys}(\os{\circ}{X}_{\ul{\lam}_{m\bet}}/\os{\circ}{T}))(-l,u) 
\tag{1.5.20.3}\label{eqn:egs}
\end{equation*}
\begin{equation*}
\lo 
R^{k+2}f_{\os{\circ}{X}_{\ul{\lam}_m}/\os{\circ}{T}*}(E_{\ul{\lam}_m}
\otimes_{\mab Z} \vp_{\rm crys}(\os{\circ}{X}_{\ul{\lam}_m}/\os{\circ}{T}))(-(l-1))
\end{equation*}
be the obvious sheafied version of the Gysin morphism 
defined in \cite[(2.8.4.5)]{nh2}.
Here 
$$\vp_{\ul{\lam}_{m}{\rm crys}}
(\os{\circ}{X}_{m,T_0}/\os{\circ}{T}) 
\quad \text{and} \quad  
\vp_{\ul{\lam}_{m\bet}{\rm crys}}
(\os{\circ}{X}_{m,T_0}/\os{\circ}{T})$$  
are the crystalline orientation sheaves of
$\os{\circ}{X}_{\ul{\lam}_m}$ and 
$\os{\circ}{X}_{\ul{\lam}_{m\bet}}$ 
in 
$(\os{\circ}{X}_{\ul{\lam}_{m}}
/\os{\circ}{T})_{\rm crys}$ 
and $(\os{\circ}{X}_{\ul{\lam}_{m\bet}}
/\os{\circ}{T})_{\rm crys}$, 
respectively,  
defined similarly in \cite[p.~81, (2.8)]{nh2}.
Set
\begin{align*}
& G_m:=\bigoplus_{j\geq \max \{-(k+m),0\}}
\sum_{\{\lam_{m0},\ldots, \lam_{m,2j+k+m}~ 
\vert ~ \lam_{m\gam} < \lam_{m\al}\; (\gam < \al)\}}
\sum_{\bet=0}^{2j+k+m}(-1)^{\bet}
G_{\ul{\lam}_m}^{\ul{\lam}_{m{\bet}}} \col 
\tag{1.5.20.4}\label{eqn:togsn}\\
&\bigoplus_{j\geq \max \{-(k+m),0\}} 
R^{q-2j-k-2m}
f_{\os{\circ}{X}{}^{(2j+k+m)}_{m,T_0}
/\os{\circ}{T}*}
(E^m_{\os{\circ}{X}{}^{(2j+k+m)}_{m,T_0}
/\os{\circ}{T}}
\otimes_{\mab Z} 
\vp^{(2j+k+m)}_{\rm crys}(
\os{\circ}{X}_{m,T_0}/\os{\circ}{T}))\\ 
&(-j-k-m,u)\lo \\
&\bigoplus_{j\geq \max \{-(k+m)+1,0\}} 
R^{q-2j-k-2m+2}
f_{\os{\circ}{X}{}^{(2j+k+m-1)}_{m,T_0}
/\os{\circ}{T}*}
(E^m_{\os{\circ}{X}{}^{(2j+k+m-1)}_{m,T_0}
/\os{\circ}{T}}
\otimes_{\mab Z} 
\vp^{(2j+k+m-1)}_{\rm crys}(
\os{\circ}{X}_{m,T_0}/\os{\circ}{T}))\\
&(-j-k-m+1,u).
\end{align*}
\par 
Let $\iota_{\ul{\lam}_m}^{\ul{\lam}_{m\bet}} \col 
\os{\circ}{X}_{\ul{\lam}_m} 
\os{\sus}{\lo} 
\os{\circ}{X}_{\ul{\lam}_{m\bet}}$ be 
the natural immersion.  
The morphism $\iota_{\ul{\lam}_m}^{\ul{\lam}_{m\bet}}$ 
induces the morphism  
\begin{equation}
(-1)^{\bet}
\iota_{\ul{\lam}_{m}{\rm crys}}^{\ul{\lam}_{m\bet}*}
\col 
\iota_{\ul{\lam}_m{\rm crys}}^{\ul{\lam}_{m\bet}*}
(E_{\ul{\lam}_{m\bet}}
\otimes_{\mab Z}\vp_{\ul{\lam}_{m\bet}{\rm crys}}
(\os{\circ}{X}_{m,T_0}/\os{\circ}{T})) 
\lo 
E_{\ul{\lam}_m} 
\otimes_{\mab Z}
\vp_{\ul{\lam}_m{\rm crys}}
(\os{\circ}{X}_{m,T_0}/\os{\circ}{T}) 
\tag{1.5.20.5}\label{eqn:defcbd}
\end{equation}
as in \cite[(2.11.1.2)]{nh2} ((\ref{eqn:odpaps})).
Set 
\begin{align*}
& \rho_m:=\bigoplus_{j\geq \max \{-(k+m),0\}}\sum_{\{\lam_{m0},\ldots, \lam_{m,2j+k+m}~ 
\vert ~\lam_{m\gam} \not= \lam_{m\alpha}\; 
(\gam \not= \alpha)\}}
\sum_{\bet=0}^{2j+k+m}(-1)^{\bet}
\iota_{\ul{\lam}_m{\rm crys}}^{\ul{\lam}_{m\bet*}} 
\col 
\tag{1.5.20.6}\label{eqn:rhogsn}\\
&\bigoplus_{j\geq \max \{-(k+m),0\}} 
R^{q-2j-k-2m}
f_{\os{\circ}{X}{}^{(2j+k+m)}_{m,T_0}
/\os{\circ}{T}*}
(E^m_{\os{\circ}{X}{}^{(2j+k+m)}_{m,T_0}
/\os{\circ}{T}}
\otimes_{\mab Z}\vp^{(2j+k+m)}_{\rm crys}
(\os{\circ}{X}_{m,T_0}/\os{\circ}{T}))  \\ 
&(-j-k-m,u)\lo \\ 
& \bigoplus_{j\geq \max \{-(k+m),0\}} 
R^{q-2j-k-2m}
f_{\os{\circ}{X}{}^{(2j+k+m+1)}_{m,T_0}/\os{\circ}{T}}
(E^m_{\os{\circ}{X}{}^{(2j+k+m+1)}_{m,T_0}
/\os{\circ}{T}}
\otimes_{\mab Z}\vp^{(2j+k+m+1)}_{\rm crys}
(\os{\circ}{X}_{m,T_0}/\os{\circ}{T}))\\
& (-j-k-m,u).
\end{align*}

\begin{prop}\label{prop:deccbd} 
The boundary morphism between the $E_1$-terms of 
the spectral sequence {\rm (\ref{eqn:escssp})} 
is given by the following diagram$:$ 
\begin{equation*} 
\tag{1.5.21.1}\label{cd:gsmsd}
\end{equation*} 
\begin{equation*} 
\begin{split} 
{} & \bigoplus_{m=0}^N\bigoplus_{j\geq \max \{-(k+m),0\}} 
R^{q-2j-k-2m}
f_{\os{\circ}{X}{}^{(2j+k+m)}_{m+1,T_0}
/\os{\circ}{T}*}
(E^m_{\os{\circ}{X}{}^{(2j+k+m)}_{m+1,T_0}/\os{\circ}{T}} \\
{} & \phantom{R^{q-2m-k}
f_{(\os{\circ}{X}^{(k+m)}_{m+1}, 
Z\vert_{\os{\circ}{X}^{(k+m)}_m})/S*}
({\cal O}\quad \quad} 
\otimes_{\mab Z}\vp^{(2j+k+m)}_{\rm crys}
(\os{\circ}{X}_{m+1,T_0}/\os{\circ}{T}))(-j-k-m,u). 
\end{split} 
\end{equation*}  
$$\text{\scriptsize
{$\bigoplus_{m=0}^N\sum_{i= 0}^{m+1}(-1)^i\del^i$}}
~\uparrow$$
\begin{equation*} 
\begin{split} 
{} & \bigoplus_{m=0}^N\bigoplus_{j\geq \max \{-(k+m),0\}} 
R^{q-2j-k-2m}
f_{\os{\circ}{X}{}^{(2j+k+m)}_{m,T_0}
/\os{\circ}{T}*}
(E^m_{\os{\circ}{X}{}^{(2j+k+m)}_{m,T_0}
/\os{\circ}{T}} \\
{} & \phantom{R^{q-2m-k}
f_{(\os{\circ}{X}^{(k+m)}_{m}, 
Z\vert_{\os{\circ}{X}^{(k+m)}_m})/S*}
({\cal O}\quad \quad} 
\otimes_{\mab Z}\vp^{(2j+k+m)}_{\rm crys}
(\os{\circ}{X}_{m,T_0}/\os{\circ}{T}))(-j-k-m,u). 
\end{split}  
\end{equation*}  
$$\text{\scriptsize
{${\bigoplus_{m=0}^N(-1)^m
[G_m+\rho_m]}$}}
~\downarrow \quad \quad \quad \quad \quad \quad \quad 
\quad \quad \quad \quad \quad$$
\begin{equation*} 
\begin{split} 
{} & 
\bigoplus_{m=0}^N\bigoplus_{j\geq \max \{-(k+m)+1,0\}} 
R^{q-2j-k-2m+2}
f_{\os{\circ}{X}{}^{(2j+k+m-1)}_{m,T_0}
/\os{\circ}{T}*}
(E^m_{\os{\circ}{X}{}^{(2j+k+m-1)}_{m,T_0}
/\os{\circ}{T}} \\
{} & \phantom{R^{q-2m-k}
f_{(\os{\circ}{X}^{(k+m)}_{m}, 
Z\vert_{\os{\circ}{X}^{(k+m)}_m})/S*}
({\cal O}\quad \quad} 
\otimes_{\mab Z}\vp^{(2j+k+m-1)}_{\rm crys}
(\os{\circ}{X}_{m,T_0}/\os{\circ}{T}))(-j-k-m+1,u). 
\end{split}  
\end{equation*}  
Here $\del^i$ $(0 \leq i \leq m+1)$ is a standard coface morphism.
\end{prop} 
\begin{proof} 
By using (\ref{lemm:ti}), 
this proposition follows from \cite[(2.2.12), (2.2.13)]{nh3} 
and the proof of \cite[(10.1)]{ndw}. 
\end{proof}

\begin{theo}[{\bf Contravariant functoriality III of of $A_{\rm zar}$}]\label{theo:itc} 
Let the notations be as in the beginning of this section.  
Assume that $\os{\circ}{S}$ is a $p$-adic formal scheme. 
Then there exists a morphism 
\begin{align*} 
g^*_{\bul \leq N}  &\col 
(A_{\rm zar}(Y_{\bul \leq N,\os{\circ}{T}{}'_0}/S'(T')^{\nat},
F^{\bul \leq N}),P)\otimes^L_{\mab Z}{\mab Q} \tag{1.5.22.1}\label{eqn:auqd}\\
&\lo
Rg_{{\bul \leq N}*}
((A_{\rm zar}(X_{\bul \leq N,\os{\circ}{T}_0}/S(T)^{\nat},
E^{\bul \leq N}),P)\otimes^L_{\mab Z}{\mab Q})
\end{align*}
of filtered complexes in ${\rm D}^+{\rm F}(f^{-1}_{\bul \leq N}({\cal K}_{T'}))$ fitting into 
the following commutative diagram 
\begin{equation*} 
\begin{CD}
A_{\rm zar}(Y_{\bul \leq N,\os{\circ}{T}{}'_0}/S'(T')^{\nat},F^{\bul \leq N})
\otimes^L_{\mab Z}{\mab Q}
@>{g^*_{\bul \leq N}}>> \\ 
@A{\theta_{Y_{\bul \leq N,\os{\circ}{T}{}'_0}/S'(T')^{\nat}} \wedge}A{\simeq}A 
\\ 
Ru_{Y_{\bul \leq N,T'_0}/T'*}
(\eps^*_{Y_{\bul \leq N,\os{\circ}{T}{}'_0}/S'(T')^{\nat}}(F^{\bul \leq N}))
\otimes^L_{\mab Z}{\mab Q}
@>{g^*_{\bul \leq N}}>> 
\end{CD}
\tag{1.5.22.2}\label{eqn:auqcuc} 
\end{equation*}  
\begin{equation*}
\begin{CD}
Rg_{\bul \leq N*}
(A_{\rm zar}(X_{\bul \leq N,\os{\circ}{T}_0}/S(T)^{\nat},E^{\bul \leq N}))
\otimes^L_{\mab Z}{\mab Q}\\ 
@A{\theta_{X_{\bul \leq N,\os{\circ}{T}_0}/S(T)^{\nat}} \wedge}A{\simeq}A \\
Rg_{\bul \leq N*}
Ru_{X_{\bul \leq N,T_0}/T*}
(\eps^*_{X_{\bul \leq N,\os{\circ}{T}_0}/S(T)^{\nat}}(E^{\bul \leq N}))
\otimes^L_{\mab Z}{\mab Q}. 
\end{CD}
\end{equation*}
This morphism satisfies the similar relation to {\rm (\ref{ali:pdpp})}.  
Assume that $T$ and $T'$ are restrictively hollow 
with respective to the morphism $T_0\lo S$ and $T'_0\lo S'$, respectively.  
Then there exists the similar commutative diagram  to {\rm (\ref{cd:tmfccz})}. 
\end{theo}
\begin{proof} 
Because we take the tensorization $\otimes_{\mab Z}^L{\mab Q}$ in this theorem, 
we need not give care to invert $\deg(u)$. 
Hence the proof of (\ref{theo:funas}) works. 
\end{proof}

\begin{theo}\label{theo:cfsp}
Let the notations be as {\rm (\ref{theo:itc})}.  
Assume that  $(1.5.6.4)$ and $(1.5.6.5)$ hold.   
Then the obvious analogue of {\rm (\ref{coro:fuu})} holds 
for $(A_{\rm zar}(X_{\bul \leq N,\os{\circ}{T}_0}/S(T)^{\nat},E^{\bul \leq N}))
\otimes^L_{\mab Z}{\mab Q}$ and 
$(A_{\rm zar}(Y_{\bul \leq N,\os{\circ}{T}{}'_0}/S'(T')^{\nat},F^{\bul \leq N}))
\otimes^L_{\mab Z}{\mab Q}$. 
\end{theo} 
\begin{proof} 
Because we take the tensorization $\otimes_{\mab Z}^L{\mab Q}$ in this theorem, 
we need not give care to invert $\deg(u)$ as in (\ref{theo:itc}). 
Hence the proof of (\ref{coro:fuu}) works. 
\end{proof}

\section{Filtered base change theorem I}\label{sec:bckf}
In this section we prove 
the filtered base change theorem 
of the zariskian $p$-adic filtered Steenbrink complex 
defined in (\ref{defi:fdirpd}). 
\par
Let the notations be as in \S\ref{sec:psc}. 
Assume that $\os{\circ}{X}_{\bul \leq N}$ is quasi-compact, 
that is, $\os{\circ}{X}_m$ $(0\leq \forall m \leq N)$ is quasi-compact. 
Let $f \col X_{\bul \leq N,\os{\circ}{T}_0} \lo S(T)^{\nat}$ 
be the structural morphism.  
Let $f_T\col X_{\bul \leq N,T_0}\lo T$ be also the structural morphism. 

\begin{prop}\label{prop:bdccd}  
Assume that $\os{\circ}{T}$ is quasi-compact and that $\os{\circ}{f} 
\col \os{\circ}{X}_{\bul \leq N,T_0}\lo \os{\circ}{T}_0$ 
is quasi-compact and quasi-separated. 
Then $Rf_*((A_{\rm zar}(X_{\bul \leq N,\os{\circ}{T}_0}/S(T)^{\nat},E^{\bul \leq N}),P))$ 
is isomorphic to a bounded filtered complex of 
${\cal O}_T$-modules. 
\end{prop}
\begin{proof}
By  \cite[7.6 Theorem]{bob},  
$Rf_{\os{\circ}{X}{}^{(l)}_{m,T_0}/\os{\circ}{T}*}
(E_{\os{\circ}{X}{}^{(l)}_{m,T_0}/\os{\circ}{T}})$ 
$(0\leq m \leq N, l\in {\mab N})$ is bounded. 
Hence 
$Rf_*((A_{\rm zar}(X_{\bul \leq N,\os{\circ}{T}_0}/S(T)^{\nat},E^{\bul \leq N}),P))$ 
is bounded by the spectral sequence (\ref{eqn:esasp}). 
\end{proof}

\begin{theo}[{\bf Log base change theorem of 
$(A_{\rm zar},P)$}]\label{theo:bccange} 
Let the assumptions be as in {\rm (\ref{prop:bdccd})}.  
Let $(T',{\cal J}',\del')$ be another log PD-enlargement over $S$. 
Assume that  ${\cal J}'$ is quasi-coherent. 
Set $T'_0:=\ul{\rm Spec}^{\log}_{T'}({\cal O}_{T'}/{\cal J}')$. 
Let $u\col (S(T')^{\nat},{\cal J}',\del') \lo (S(T)^{\nat},{\cal J},\del)$ be 
a morphism of fine log PD-schemes. 
Let 
$f' \col X_{\bul \leq N,\os{\circ}{T}{}'_0}=
X_{\bul \leq N}\times_{S}S_{\os{\circ}{T}{}'_0} \lo S(T')^{\nat}$ 
be the base change morphism of $f$  
by the morphism $S(T')^{\nat}\lo S(T)^{\nat}$.  
Let $q_{\bul \leq N} \col 
X_{\bul \leq N,\os{\circ}{T}{}'_0} \lo X_{\bul \leq N,\os{\circ}{T}_0}$ 
be the induced morphism by $u$. 
Then there exists 
the following canonical filtered isomorphism
\begin{equation*}
Lu^*Rf_*((A_{\rm zar}(X_{\bul \leq N,\os{\circ}{T}_0}/S(T)^{\nat},E^{\bul \leq N}),P)) 
\os{\sim}{\lo} Rf'_*
((A_{\rm zar}(X_{\bul \leq N,\os{\circ}{T}{}'_0}/S(T')^{\nat},
\os{\circ}{q}{}^{*}_{\bul \leq N,{\rm crys}}(E^{\bul \leq N})),P))
\tag{1.6.2.1}\label{eqn:blucpw}
\end{equation*}
in ${\rm DF}(
f^{-1}_{T'}({\cal O}_{T'}))$. 
\end{theo}
\begin{proof}  
Let the notations be as in \S\ref{sec:psc}. 
Set 
$\ol{\cal P}_{\bul \leq N,\bul,\ol{S(T')^{\nat}}}
:=\ol{\cal P}_{\bul \leq N,\bul}\times_{\ol{S(T)^{\nat}}}\ol{S(T')^{\nat}}$. 
Let   
$\ol{\mathfrak D}{}'_{\bul \leq N,\bul}$ 
be the log PD-envelope of the immersion 
$X_{\bul \leq N,\os{\circ}{T}{}'_0}\os{\sus}{\lo} 
\ol{\cal P}_{\bul \leq N,\bul,\ol{S(T')^{\nat}}}$ 
over $(\os{\circ}{T}{}',{\cal J}',\del')$. 
Then we have the natural morphisms 
$\ol{\cal P}_{\bul \leq N,\bul,\ol{S(T')^{\nat}}}
\lo \ol{\cal P}_{\bul \leq N,\bul}$ 
and 
$\ol{\mathfrak D}{}'_{\bul \leq N,\bul}  \lo \ol{\mathfrak D}_{\bul \leq N,\bul}$.  
We also have the identity morphism 
${\rm id}\col \os{\circ}{q}{}^{*}_{\bul \leq N,{\rm crys}}(E^{\bul \leq N})
\lo \os{\circ}{q}{}^{*}_{\bul \leq N,{\rm crys}}(E^{\bul \leq N})$. 
Hence we have the following natural morphism  
\begin{equation*} 
(A_{\rm zar}(X_{\bul \leq N,\os{\circ}{T}_0}/S(T)^{\nat},E^{\bul \leq N}),P) 
\lo 
Rq_{\bul \leq N*}
((A_{\rm zar}(X_{\bul \leq N,\os{\circ}{T}{}'_0}/S(T')^{\nat},
\os{\circ}{q}{}^{*}_{\bul \leq N,{\rm crys}}(E^{\bul \leq N})),P)). 
\tag{1.6.2.2}\label{eqn:bcxa}
\end{equation*} 
(This is a special case of (\ref{theo:funas}).) 
By applying $Rf_*$ to (\ref{eqn:bcxa}) and using 
the adjoint property of $L$ and $R$ (\cite[(1.2.2)]{nh2}), 
we have the natural morphism (\ref{eqn:blucpw}). 
Here we have used the boundedness in 
(\ref{prop:bdccd}) for the well-definedness of $Lu^*$. 
\par
We  prove that (\ref{eqn:blucpw}) is an isomorphism. 
We may assume 
that $N=0$. Set $X:=X_0$ and $E=E^0$. 
By the filtered cohomological 
descent \cite[(1.5.1) (2)]{nh2}  
and by the same argument as that in 
the proof of \cite[V Proposition 3.5.2]{bb} 
(\cite[7.8 Theorem]{bob}),
we may assume that $\os{\circ}{S}$ is affine 
and that $\os{\circ}{X}$ is an affine scheme over $\os{\circ}{S}$.  
Then $X_{\os{\circ}{T}_0}$ has an SNCL lift ${\cal X}/S(T)^{\nat}$  
by the proofs of \cite[(3.14)]{klog1} and \cite[(2.3.14)]{nh2} 
(cf.~the proof of \cite[(4.7)]{ny}). 
Moreover, by (\ref{lemm:etl}), we can assume that 
${\cal X}$ has a log smooth lift $\ol{\cal X}/\ol{S(T)^{\nat}}$.  
Let $(\ol{\cal E},\ol{\nabla})$ be a quasi-coherent ${\cal O}_{\ol{\cal X}}$-module 
with integrable connection corresponding to 
$\eps^*_{X_{\os{\circ}{T}}/\os{\circ}{T}}(E)$.  
Set ${\cal E}
:=\ol{\cal E}\otimes_{{\cal O}_{{\mathfrak D}(\ol{S(T)})}}{\cal O}_{S(T)}$. 
In this case, 
$(A_{\rm zar}(X_{\bul \leq N,\os{\circ}{T}_0}/S(T)^{\nat},E^{\bul \leq N}),P) 
=(A_{\rm zar}({\cal X}/S(T)^{\nat},{\cal E}),P)$. 
Let ${\mathfrak f} \col {\cal X} \lo S(T)^{\nat}$ be the 
structural morphism.  
For $i,j\in {\mab N}$, 
set ${\mathfrak f}_{*}(P_k):={\mathfrak f}_{*}((P_k+P_j)
A_{\rm zar}({\cal X}/S(T)^{\nat},{\cal E})^{ij}/P_j)$ 
$(i,j\in {\mab N}, k\in {\mab Z})$ 
for simplicity of notation. 
Then we have the following formula 
\begin{equation*}
Rf_*((A_{\rm zar}(X_{\os{\circ}{T}_0}/S(T)^{\nat},E),P)) 
=
s(\cdots, ({\mathfrak f}_{*}
(A_{\rm zar}({\cal X}/S(T)^{\nat},{\cal E})^{ij}),
{\mathfrak f}_{*}(P_{2j+k+1})),\cdots)
\end{equation*}
and we have the same formula for $X_{\os{\circ}{T}{}'_0}/S(T')^{\nat}$. 
By (\ref{coro:flt}), 
${\mathfrak f}_{*}
(A_{\rm zar}({\cal X}/S(T)^{\nat},{\cal E}))$ and 
${\mathfrak f}_{*}(P_{2j+k+1})$ $(\forall k)$  
consist of flat ${\cal O}_T$-modules. 
Therefore the source of (\ref{eqn:blucpw}) is equal to 
$s(\cdots 
(u^*{\mathfrak f}_{*}(A_{\rm zar}({\cal X}/S(T)^{\nat},{\cal E})^{ij}),
u^*{\mathfrak f}_{*}(P_{2j+k+1})) \cdots)$. 
Since  
${\mathfrak f} \col \os{\circ}{\cal X} \lo \os{\circ}{T}$ 
is an affine morphism,
(\ref{theo:bccange}) follows from 
the affine base change theorem as in \cite[7.8 Theorem]{bob}.
\end{proof} 

\par 
Let $\os{\circ}{Y}$ be a smooth scheme over $\os{\circ}{T}$. 
Endow $\os{\circ}{Y}$ with the inverse image of $M_{S_{\os{\circ}{T}}}$ 
and let $Y$ be the resulting log scheme. 
Let ${\cal Y}$ be a log smooth scheme defined in (\ref{coro:connfil}) below. 
Let $D_{{\cal Y}/S(T)^{\nat}}(1)$ be the log PD-envelope of the immersion 
${\cal Y}\os{\sus}{\lo} {\cal Y}\times_{S(T)^{\nat}}{\cal Y}$ over $(S(T)^{\nat},{\cal J},\del)$. 
As in \cite[V]{bb} and \cite[\S7]{bob},  
we have the following two corollaries 
(cf.~\cite[(2.10.5), (2.10.7)]{nh2}) by 
using  (\ref{theo:bccange}) and a fact that 
$p_i \col \os{\circ}{D}_{{\cal Y}/S(T)}(1)\lo \os{\circ}{\cal Y}$ 
$(i=1,2)$ 
is flat (\cite[(6.5)]{klog1}): 

\begin{coro}\label{coro:connfil}
Let $g\col X_{\bul \leq N,\os{\circ}{T}_0} \lo Y$ be an 
$N$-truncated simplicial SNCL scheme 
which has the disjoint union of an 
affine $N$-truncated simplicial open covering of 
$X_{\bul \leq N,\os{\circ}{T}_0}$. 
Assume that $Y$ has a log smooth lift ${\cal Y}$ over $S(T)^{\nat}$. 
Let $q$ be an integer.
Let $g\col X_{\bul \leq N,\os{\circ}{T}_0}\lo {\cal Y}$ 
be the structural morphism. 
Then there exists  a quasi-nilpotent integrable connection 
\begin{align*}
&P_kR^qg_{X_{\bul \leq N,\os{\circ}{T}_0}/{\cal Y}*}
(\eps^*_{X_{\bul \leq N,\os{\circ}{T}_0}/S(T)^{\nat}}(E^{\bul \leq N}))
\os{\nabla_k}{\lo} \tag{1.6.3.1}\\
& P_kR^qg_{X_{\bul \leq N,\os{\circ}{T}_0}/{\cal Y}*}
(\eps^*_{X_{\bul \leq N,\os{\circ}{T}_0}/S(T)^{\nat}}(E^{\bul \leq N}))
{\otimes}_{{\cal O}_{\cal Y}}{\Om}_{{\cal Y}/S(T)^{\nat}}^1
\end{align*}
making the following diagram commutative 
for any two nonnegative integers $k\leq l:$
\begin{equation*}
\begin{CD}
P_kR^qg_{X_{\bul \leq N,\os{\circ}{T}_0}/{\cal Y}*}
(\eps^*_{X_{\bul \leq N,\os{\circ}{T}_0}/S(T)^{\nat}}(E^{\bul \leq N}))
@>{\nabla_k}>> \\
@V{\bigcap}VV  \\
P_lR^qg_{X_{\bul \leq N,\os{\circ}{T}_0}/{\cal Y}*}
(\eps^*_{X_{\bul \leq N,\os{\circ}{T}_0}/S(T)^{\nat}}(E^{\bul \leq N}))
@>{\nabla_l}>>
\end{CD}
\tag{1.6.3.2}
\end{equation*}
\begin{equation*}
\begin{CD}
P_kR^qg_{X_{\bul \leq N,\os{\circ}{T}_0}/{\cal Y}*}
(\eps^*_{X_{\bul \leq N,\os{\circ}{T}_0}/S(T)^{\nat}}(E^{\bul \leq N}))
{\otimes}_{{\cal O}_{\cal Y}}{\Om}_{{\cal Y}/S(T)^{\nat}}^1\\
@V{\bigcap}VV \\ 
P_lR^qg_{X_{\bul \leq N,\os{\circ}{T}_0}/{\cal Y}*}
(\eps^*_{X_{\bul \leq N,\os{\circ}{T}_0}/S(T)^{\nat}}(E^{\bul \leq N})) 
{\otimes}_{{\cal O}_{\cal Y}}{\Om}_{{\cal Y}/S(T)^{\nat}}^1.
\end{CD} 
\end{equation*}
\end{coro}
\begin{proof} 
This follows from (\ref{theo:bccange}) as in \cite[\S7]{bob}.
\end{proof} 

\begin{coro}\label{coro:fctd}
Let the notations and the assumptions be as in $(\ref{prop:bdccd})$. 
Then 
$$Rf_*
(P_kA_{\rm zar}(X_{\bul \leq N,\os{\circ}{T}_0}/S(T)^{\nat},E^{\bul \leq N})) \quad (k \in{\mab N})$$
has finite tor-dimension. 
Moreover, if $\os{\circ}{T}$ is noetherian and 
if $\os{\circ}{f}$ is proper,
then $Rf_*
(P_kA_{\rm zar}(X_{\bul \leq N,\os{\circ}{T}_0}/S(T)^{\nat},E^{\bul \leq N}))$ 
is a perfect complex of ${\cal O}_T$-modules.
\end{coro}

Using \cite[(2.10.10)]{nh2}, 
we have the following corollary 
(cf.~\cite[(2.10.11)]{nh2}): 

\begin{coro}\label{coro:filpcerf}
Let the notations 
and the assumptions be as in $(\ref{coro:fctd})$.
Then the filtered complex 
$Rf_*((A_{\rm zar}(X_{\bul \leq N,\os{\circ}{T}_0}/S(T)^{\nat},E^{\bul \leq N}),P))$ 
is a {\it filtered perfect}
complex of ${\cal O}_T$-modules, that is, 
locally on $T_{\rm zar}$, filteredly quasi-isomorphic to 
a filtered strictly perfect complex {\rm (\cite[(2.10.8)]{nh2})}.
\end{coro}
\begin{proof}
(\ref{coro:filpcerf}) 
immediately follows 
from (\ref{coro:fctd}) and 
\cite[(2.10.10)]{nh2}.
\end{proof}

\section{$p$-adic monodromy operators I and 
modified $P$-filtered log crystalline complexes} 
\label{sec:crcks}
In this section we first define the $p$-adic monodromy operator 
for a truncated simplicial log smooth scheme 
following the method in \cite{hk}. 
To prove the properties of the $p$-adic monodromy operator, 
we try to treat it as a morphism of complexes of 
a derived category as possible. 
This try is different from that in \cite{hk} and \cite{tst}. 
If the truncated simplicial log smooth scheme has 
an affine truncated simplicial open covering, 
then we define the $p$-adic monodromy operator 
following the method in \cite{st1}, \cite{hdw} and \cite{hk}. 
We prove that 
the two definitions of the $p$-adic monodromy operators are equal as in \cite{hk}. 
\par 
First we give the following definition of a
pre-stratification on an object of a derived category:

\begin{defi}[{\bf Pre-stratification on a complex of a derived category}]\label{defi:srdr}
Let $(U,{\cal K},\del)$ be a log PD-scheme 
and let $T\lo U$ be a morphism of log schemes such that 
$\del$ extends to $T$. 
Set $T(1):=T\times_UT$ and $T(2):=T\times_UT\times_UT$. 
Let ${\mathfrak E}$ (resp.~${\mathfrak E}(1)$) 
be the log PD-envelope of 
the diagonal immersion 
$T\os{\sus}{\lo} T(1)$ (resp.~$T\os{\sus}{\lo} T(2)$) 
over $(U,{\cal K},\del)$. 
Let $J$ be the ideal sheaf of the diagonal immersion 
$T\os{\sus}{\lo} T(1)$ and let 
$\Del\col T\os{\sus}{\lo} {\mathfrak E}$ be the canonical immersion. 
Let $q_i\col {\mathfrak E}\lo T$ $(i=1,2)$ be 
the induced morphism 
by the $i$-th projection $T(1)\lo T$. 
\par 
Let ${\cal E}^{\bul}$ be an object of 
the derived category $D^-({\cal O}_T)$ of bounded above  
complexes of ${\cal O}_T$-modules.   
Let $\eps \col Lq_2^*({\cal E}^{\bul})\os{\sim}{\lo} Lq_1^*({\cal E}^{\bul})$ 
be an isomorphism in  $D^-({\cal O}_{\mathfrak E})$. 
We say that $\eps$ is a {\it pre-stratification} on ${\cal E}^{\bul}$ 
if $\eps~{\rm mod}~\ol{J}={\rm id}_{{\cal E}^{\bul}}$. 
Here 
$\eps~{\rm mod}~\ol{J} \col {\cal E}^{\bul}\lo {\cal E}^{\bul}$ 
is the isomorphism $L\Del^*(\eps)\col {\cal E}^{\bul}=
L\Del^*Lq_2^*({\cal E}^{\bul})\os{\sim}{\lo} 
L\Del^*Lq_1^*({\cal E}^{\bul})={\cal E}^{\bul}$.  
\end{defi} 

\begin{lemm}\label{lemm:eit}
Let ${\cal E}^{\bul}_i$ $(i=1,2)$ be an object of $D^-({\cal O}_T)$.     
Let $\eps_i \col Lq_2^*({\cal E}^{\bul}_i)
\os{\sim}{\lo} Lq_1^*({\cal E}^{\bul}_i)$ 
be a pre-stratification on ${\cal E}^{\bul}_i$.  
Then the following composite isomorphism 
\begin{align*} 
Lq_2^*({\cal E}^{\bul}_1\otimes^L_{{\cal O}_{\mathfrak E}}{\cal E}^{\bul}_2)
=Lq_2^*({\cal E}^{\bul}_1)\otimes^L_{{\cal O}_{\mathfrak E}}Lq_2^*({\cal E}^{\bul}_2)
\os{\eps_1\otimes^L \eps_2}{\lo} 
Lq_1^*({\cal E}^{\bul}_1)\otimes^L_{{\cal O}_{\mathfrak E}}Lq_1^*({\cal E}^{\bul}_2)
=Lq_1^*({\cal E}^{\bul}_1\otimes^L_{{\cal O}_{\mathfrak E}}{\cal E}^{\bul}_2)
\tag{1.7.2.1}\label{ali:qqee}
\end{align*} 
is a pre-stratification. 
\end{lemm}
\begin{proof} 
Obvious.
\end{proof}

\begin{lemm}\label{lemm:elq}
Let the notations be as in {\rm (\ref{defi:srdr})}. 
Assume that 
$\os{\circ}{q}_1 \col \os{\circ}{\mathfrak E}\lo \os{\circ}{T}$ induces 
an isomorphism ${\cal O}_{\mathfrak E}\simeq 
\bigoplus_{i\in {\mab N}}(\ol{J}{}^{[i]}/\ol{J}{}^{[i+1]})$ 
as ${\cal O}_T$-algebras. 
Let $\eps \col Lq_2^*({\cal E}^{\bul})\os{\sim}{\lo} Lq_1^*({\cal E}^{\bul})$ be 
a pre-stratification on ${\cal E}^{\bul}$. 
Then there exists a natural morphism   
\begin{align*} 
\nabla \col {\cal E}^{\bul}
\lo {\cal E}^{\bul}\otimes^L_{{\cal O}_T}\Om^1_{T/U}
\tag{1.7.3.1}\label{ali:lotu} 
\end{align*}  
in $D^-({\cal O}_T)$ constructed in the proof below.  
\end{lemm}
\begin{proof}
By using the isomorphism ${\cal O}_{\mathfrak E}\os{\sim}{\lo}
\bigoplus_{i\in {\mab N}}(\ol{J}{}^{[i]}/\ol{J}{}^{[i+1]})$, 
we have the following composite morphism 
\begin{align*} 
\al \col {\cal E}^{\bul}&=
{\cal O}_T\otimes^L_{{\cal O}_T}{\cal E}^{\bul}
\os{q_2^*\otimes {\rm id}_{{\cal E}^{\bul}}}{\lo} 
{\cal O}_{\mathfrak E}\otimes^L_{{\cal O}_T}{\cal E}^{\bul}=
Lq_2^*({\cal E}^{\bul})\os{\eps,\sim}{\lo} Lq_1^*({\cal E}^{\bul})
={\cal E}^{\bul}\otimes^L_{{\cal O}_T}{\cal O}_{\mathfrak E}
\tag{1.7.3.2}\label{ali:eeoe}\\
&=\bigoplus_{i\in {\mab N}}
{\cal E}^{\bul}\otimes^L_{{\cal O}_T}(\ol{J}{}^{[i]}/\ol{J}{}^{[i+1]}). 
\end{align*} 
Since $\ol{J}{}^{[0]}/\ol{J}{}^{[1]}={\cal O}_T$, 
we also have the following composite morphism 
\begin{align*} 
\bet \col {\cal E}^{\bul}=
{\cal E}^{\bul}\otimes^L_{{\cal O}_T}(\ol{J}{}^{[0]}/\ol{J}{}^{[1]})
\os{\subset}{\lo} 
\bigoplus_{i\in {\mab N}}
{\cal E}^{\bul}\otimes^L_{{\cal O}_T}(\ol{J}{}^{[i]}/\ol{J}{}^{[i+1]}). 
\tag{1.7.3.3}\label{ali:oojj1} 
\end{align*} 
Let ${\rm pr}\col 
\bigoplus_{i\in {\mab N}}
{\cal E}^{\bul}\otimes^L_{{\cal O}_T}(\ol{J}{}^{[i]}/\ol{J}{}^{[i+1]})
\lo {\cal E}^{\bul}\otimes^L_{{\cal O}_T}(\ol{J}{}^{[1]}/\ol{J}{}^{[2]})$ 
be the projection. 
Since $\Om^1_{T/U}=\ol{J}{}^{[1]}/\ol{J}{}^{[2]}$, 
we have the following morphism 
\begin{align*} 
\nabla :{\rm pr}\circ(\al-\bet)\col 
{\cal E}^{\bul}\lo 
{\cal E}^{\bul}\otimes^L_{{\cal O}_T}(\ol{J}{}^{[1]}/\ol{J}{}^{[2]})
={\cal E}^{\bul}\otimes^L_{{\cal O}_T}\Om^1_{T/U}.
\end{align*}  
This is the desired morphism. 
\end{proof} 

\begin{defi}\label{defi:dnn}
We call the morphism (\ref{ali:lotu}) the {\it connection} 
associated to the pre-stratification $\eps$. 
\end{defi}

\begin{rema}\label{rema:cat} 
For a local section $f$ of ${\cal O}_T$, consider the following morphisms: 
\begin{align*} 
{\rm id}_{{\cal E}^{\bul}}\otimes df 
\col {\cal E}^{\bul}\owns e\lom e\otimes df \in 
{\cal E}^{\bul} \otimes_{{\cal O}_T}\Om^1_{T/U},  
\end{align*} 
\begin{align*} 
f\nabla  
\col {\cal E}^{\bul}\os{\nabla}{\lo}  
{\cal E}^{\bul} \otimes_{{\cal O}_T}\Om^1_{T/U} 
\os{f\otimes {\rm id}_{\Om^1_{T/U}}}{\lo}  
{\cal E}^{\bul} \otimes_{{\cal O}_T}\Om^1_{T/U} 
\end{align*} 
and 
\begin{align*} 
f\cdot   
\col {\cal E}^{\bul}\os{f \cdot}{\lo} {\cal E}^{\bul}.  
\end{align*} 
Then, by (\ref{lemm:eiat}) below, we have the following ``usual'' formula: 
\begin{align*} 
\nabla(f\cdot)=f\nabla+{\rm id}_{{\cal E}^{\bul}}\otimes df. 
\end{align*} 
\end{rema} 

\begin{defi}\label{defi:ndf}
If there exists an isomorphism 
$R\col \Om^1_{T/U}\os{\sim}{\lo} {\cal O}_T^d$ 
for some $d\in {\mab Z}_{\geq 1}$, 
then 
we have the following composite morphism 
\begin{align*} 
{\cal N} : {\cal E}^{\bul}\lo 
{\cal E}^{\bul}\otimes^L_{{\cal O}_T}\Om^1_{T/U}=
({\cal E}^{\bul})^{\oplus d}.
\tag{1.7.6.1}\label{ali:ebll}
\end{align*}
We call ${\cal N}$ the 
{\it monodromy operator} associated to $\eps$ and $R$. 
\end{defi}

\begin{rema}\label{rema:dmm}
Let $q_{ij}\col {\mathfrak E}(1)\lo {\mathfrak E}$ $(1\leq i< j\leq 3)$ be 
the induced morphism by the $(i,j)$-th projection $T(2)\lo T(1)$. 
Let $\nabla^2$ be the following composite morphism 
\begin{align*} 
{\cal E}^{\bul}\os{\nabla}{\lo} 
{\cal E}^{\bul}\otimes^L_{{\cal O}_T}\Om^1_{T/U}
\os{\nabla \wedge {\rm id}_{\Om^1_{T/U}}+{\rm id}_{{\cal E}^{\bul}}\otimes d}{\lo} 
{\cal E}^{\bul}\otimes^L_{{\cal O}_T}\Om^2_{T/U}.
\end{align*} 
Even if we assume that 
$Lq_{13}^*(\eps)= Lq_{23}^*(\eps)\circ Lq_{12}^*(\eps)$ in 
(\ref{defi:srdr}), I do not know whether $\nabla^2=0$. 
\end{rema} 

\begin{lemm}\label{lemm:eiat}
Let the notations be as in {\rm (\ref{lemm:eit})}. 
Set ${\cal E}^{\bul}:={\cal E}^{\bul}_1\otimes^L_{{\cal O}_T}{\cal E}^{\bul}_2$. 
Let 
$$\tau \col {\cal E}^{\bul}_1\otimes^L_{{\cal O}_T}\Om^1_{T/U}
\otimes^L_{{\cal O}_T}{\cal E}^{\bul}_2\lo 
{\cal E}^{\bul}_1\otimes^L_{{\cal O}_T}{\cal E}^{\bul}_2
\otimes^L_{{\cal O}_T}\Om^1_{T/U}$$ 
be the canonical morphism defined by the morphism 
``$x\otimes \om \otimes y\lom x\otimes y \otimes \om$''. 
Let $\nabla_i$ and $\nabla$ be the connections 
associated to $\eps_i$ and $\eps_1\otimes^L \eps_2$, 
respectively.   
Then 
\begin{align*}
\nabla=
\nabla_1\otimes^L {\rm id}_{{\cal E}^{\bul}_2}+
{\rm id}_{{\cal E}^{\bul}_1}\otimes^L \nabla_2 
\col 
{\cal E}^{\bul}_1\otimes^L_{{\cal O}_T}{\cal E}^{\bul}_2\lo 
{\cal E}^{\bul}_1\otimes^L_{{\cal O}_T}{\cal E}^{\bul}_2
\otimes^L_{{\cal O}_T}\Om^1_{T/U}. 
\tag{1.7.8.1}\label{ali:oelt}
\end{align*}
Here we have identified 
${\cal E}^{\bul}_1\otimes^L_{{\cal O}_T}\Om^1_{T/U}
\otimes^L_{{\cal O}_T}{\cal E}^{\bul}_2$ with 
${\cal E}^{\bul}_1\otimes^L_{{\cal O}_T}{\cal E}^{\bul}_2
\otimes^L_{{\cal O}_T}\Om^1_{T/U}$
by $\tau$. 
\end{lemm}
\begin{proof} 
Let $\al_i$ and $\bet_i$ be the morphisms 
(\ref{ali:eeoe}) and (\ref{ali:oojj1}), respectively for 
$({\cal E}_i^{\bul},\eps_i)$. 
Let $\al_1\cdot \al_2$ be the following composite morphism 
\begin{align*} 
\al_1\cdot \al_2 \col &
{\cal E}^{\bul}_1\otimes^L_{{\cal O}_T}{\cal E}^{\bul}_2
\lo 
Lq_2^*({\cal E}^{\bul}_1)
\otimes^L_{{\cal O}_{\mathfrak E}}Lq_2^*({\cal E}^{\bul}_2)
\os{\eps_1\otimes^L \eps_2,\sim}{\lo} 
Lq_1^*({\cal E}^{\bul}_1)\otimes^L_{{\cal O}_{\mathfrak E}}
Lq_1^*({\cal E}^{\bul}_2)
\tag{1.7.8.2}\label{ali:eeopee}\\
&\os{\sim}{\lo} \bigoplus_{i\in {\mab N}}
{\cal E}^{\bul}_1\otimes^L_{{\cal O}_T}{\cal E}^{\bul}_2
\otimes^L_{{\cal O}_T}{\cal O}_{\mathfrak E}.
\end{align*}
Let $\al$ and $\bet$ be the morphisms 
(\ref{ali:eeoe}) and (\ref{ali:oojj1}), respectively for 
${\cal E}^{\bul}_1\otimes^L_{{\cal O}_T}{\cal E}^{\bul}_2$. 
Because $\al$ is equal to the following morphism 
\begin{align*} 
\al \col {\cal E}^{\bul}_1\otimes^L_{{\cal O}_T}{\cal E}^{\bul}_2\lo &
Lq_2^*({\cal E}^{\bul}_1\otimes^L_{{\cal O}_T}{\cal E}^{\bul}_2)
=Lq_2^*({\cal E}^{\bul}_1)\otimes^L_{{\cal O}_{\mathfrak E}}Lq_2^*({\cal E}^{\bul}_2)
\os{\eps_1\otimes^L \eps_2,\sim}{\lo} 
Lq_1^*({\cal E}^{\bul}_1)\otimes^L_{{\cal O}_{\mathfrak E}}
Lq_1^*({\cal E}^{\bul}_2)
\tag{1.7.8.3}\label{ali:eeoee}\\
&\os{\sim}{\lo} \bigoplus_{i\in {\mab N}}
{\cal E}^{\bul}_1\otimes^L_{{\cal O}_T}{\cal E}^{\bul}_2
\otimes^L_{{\cal O}_T}{\cal O}_{\mathfrak E}, 
\end{align*}
$\al=\al_1\cdot \al_2$. 
Similarly $\bet =\bet_1\cdot \bet_2$. 
Hence we have the following equalities: 
\begin{align*} 
\nabla &={\rm pr}(\al-\bet)
={\rm pr}(\al_1\cdot \al_2-\bet_1\cdot \bet_2) 
\tag{1.7.8.4}\label{ali:nnn}\\
&={\rm pr}((\al_1-\bet_1)\cdot \al_2
+\bet_1\cdot \al_2-\bet_1\cdot \bet_2)
\\
&={\rm pr}((\al_1-\bet_1)\cdot \bet_2)
+{\rm pr}(\bet_1\cdot (\al_2-\bet_2))
=\nabla_1\otimes^L {\rm id}_{{\cal E}^{\bul}_2}+
{\rm id}_{{\cal E}^{\bul}_1}\otimes^L \nabla_2.   
\end{align*} 
This equality is nothing but (\ref{ali:oelt}). 
\end{proof}

\begin{lemm}[{\bf Log base change formula}]\label{lemm:bctn}
Let $h\col (T',{\cal J}',\del') \lo (T,{\cal J},\del)$ be 
a morphism of fine log PD-schemes  on which $p$ is locally nilpotent. 
Set $T_0:=\ul{\rm Spec}^{\log}_T({\cal O}_T/{\cal J})$ and 
$T'_0:=\ul{\rm Spec}^{\log}_{T'}({\cal O}_{T'}/{\cal J}')$. 
Let $h_0\col T'_0\lo T_0$ be the induced morphism by $h$. 
Let $N$ be a nonnegative integer or $\infty$. 
Let $g\col Y_{\bul \leq N}\lo T_0$ be a log smooth integral morphism. 
Assume that $\os{\circ}{g}$ is quasi-compact and quasi-separated.  
Let $g'\col Y'_{\bul \leq N}\lo T'_0$ be 
the base change morphism of $g$ by $h_0$. 
Let $q\col Y'_{\bul \leq N}\lo Y_{\bul \leq N}$ be the projection. 
Let $F^{\bul \leq N}$ be a flat quasi-coherent crystal of 
${\cal O}_{Y_{\bul \leq N}/T}$-modules. 
Then the following natural morphism 
\begin{align*} 
Lh^*Rg_{Y_{\bul \leq N}/T*}(F^{\bul \leq N})\lo 
Rg'_{Y'_{\bul \leq N}/T'*}(q_{\rm crys}^*(F^{\bul \leq N}))
\tag{1.7.9.1}\label{ali:ysen} 
\end{align*} 
is an isomorphism. 
\end{lemm} 
\begin{proof} 
It suffices to prove that the following induced morphism by (\ref{ali:ysen}) 
\begin{align*} 
{\cal H}^q(Lh^*Rg_{Y_{\bul \leq N}/T*}(F^{\bul \leq N}))
\lo {\cal H}^q
(Rg'_{Y'_{\bul \leq N}/T'*}(q_{\rm crys}^*(F^{\bul \leq N}))) 
\quad (q\in {\mab N})
\end{align*} 
is an isomorphism. 
This follows from the following spectral sequences 
\begin{align*} 
E_1^{ij} =R^jg'_{Y'_{i}/T'*}(q_{\rm crys}^*(F^i))
\Lo {\cal H}^{i+j}(Rg'_{Y'_{\bul \leq N}/T'*}
(q_{\rm crys}^*(F^{\bul \leq N}))) \quad (i\leq N)
\end{align*} 
\begin{align*} 
E_1^{ij} ={\cal H}^j(Lh^*Rg_{Y_i/T*}(F^i))
\Lo {\cal H}^{i+j}(Lh^*Rg_{Y_{\bul \leq N}/T*}(F^{\bul \leq N})) 
\quad (i\leq N)
\end{align*}
and Kato's log base change theorem \cite[(6.10)]{klog1}. 
\end{proof}

\par 
Let $(U,{\cal K},\eps)$ be a fine log PD-scheme 
on which $p$ is locally nilpotent. 
Set $U_0:=\ul{\rm Spec}_U^{\log}({\cal O}_U/{\cal K})$. 
Let $Y$ be a fine log scheme over $U_0$ and 
let $g\col Z_{\bul \leq N}\lo Y$ be 
a log smooth integral morphism.  
Assume that $\os{\circ}{g}$ is quasi-compact 
and quasi-separated. 
Let $(Y\os{\sus}{\lo} {\cal Y})=
({\cal Y},{\cal J},\del)$ be an object of $(Y/U)_{\rm crys}$. 
Set ${\cal Y}(1):={\cal Y}\times_U{\cal Y}$ and let 
${\mathfrak D}$ be the log PD-envelope of the diagonal immersion 
${\cal Y}\os{\sus}{\lo}{\cal Y}(1)$ over $(U,{\cal K},\eps)$. 
Assume that $\Om^1_{{\cal Y}/U}$ is 
a locally free ${\cal O}_{\cal Y}$-module of finite rank. 
Assume also that the induced morphism 
$q_i\col \os{\circ}{\mathfrak D}\lo \os{\circ}{\cal Y}$ 
by the $i$-th projection ${\cal Y}(1)\lo {\cal Y}$ $(i=1,2)$ is flat. 
Let $E^{\bul \leq N}$ be a flat quasi-coherent crystal of 
${\cal O}_{Z_{\bul \leq N}/{\cal Y}}$-modules. 
Then, by the log base change theorem (\ref{lemm:bctn}), 
we have the following isomorphisms: 
\begin{align*} 
Lq_2^*Rg_{Z_{\bul \leq N}/{\cal Y}*}(E^{\bul \leq N})
\os{\sim}{\lo}
Rg_{Z_{\bul \leq N}/{\mathfrak D}*}(E^{\bul \leq N})
\os{\sim}{\longleftarrow}
Lq_1^*Rg_{Z_{\bul \leq N}/{\cal Y}*}(E^{\bul \leq N}). 
\end{align*} 
Hence we have the following connection: 
\begin{align*} 
\nabla \col Rg_{Z_{\bul \leq N}/{\cal Y}*}(E^{\bul \leq N})\lo 
Rg_{Z_{\bul \leq N}/{\cal Y}*}(E^{\bul \leq N})
\otimes^L_{{\cal O}_{\cal Y}}\Om^1_{{\cal Y}/U}.  
\tag{1.7.9.2}\label{ali:zycen}
\end{align*} 
In particular, we have the following connection 
\begin{align*} 
\nabla \col R^qg_{Z_{\bul \leq N}/{\cal Y}*}(E^{\bul \leq N})\lo 
R^qg_{Z_{\bul \leq N}/{\cal Y}*}(E^{\bul \leq N})
\otimes_{{\cal O}_{\cal Y}}\Om^1_{{\cal Y}/U} \quad (q\in {\mab N})
\tag{1.7.9.3}\label{ali:zymen}
\end{align*} 
since $\Om^1_{{\cal Y}/U}$ is a locally free ${\cal O}_{\cal Y}$-module. 

\par 
Let $S$, $(T,{\cal J},\del)$, $T_0\lo S$, 
$S_{\os{\circ}{T}_0}$, $S(T)$ and $S(T)^{\nat}$ 
be as in \S\ref{sec:ldc}.  
Set $S(T)^{\nat}(1):=S(T)^{\nat}\times_{\os{\circ}{T}}S(T)^{\nat}$. 
Consider the diagonal immersion 
$\iota\col S(T)^{\nat}\os{\sus}{\lo} S(T)^{\nat}(1)$ over $\os{\circ}{T}$. 
Let $D$ be the log PD-envelope of $\iota$ over 
$(\os{\circ}{T},{\cal J},\del)$. 
Let $J$ be the PD-ideal sheaf of $D$. 
Then ${\cal O}_D=\bigoplus_{i\in \mab N}J^{[i]}$ as an ${\cal O}_T$-module.
Let $\Om^1_{S(T)^{\nat}/\os{\circ}{T}}$ be 
the sheaf of logarithmic differential forms of 
$S(T)^{\nat}/\os{\circ}{T}$ of degree $1$. 
Then $\Om^1_{S(T)^{\nat}/\os{\circ}{T}}=J^{[1]}/J^{[2]}$. 
Consider a local case such that 
$M_{S(T)^{\nat}}/{\cal O}_T^*\simeq {\mab N}$.  
Take a local section $m$ of $M_{S(T)^{\nat}}$ whose image in 
$M_{S(T)^{\nat}}/{\cal O}_T^*\simeq {\mab N}$ is a generator.   
Let $q_i \col D\lo S(T)^{\nat}$ be the induced morphism 
by the projections 
$q_i\col S(T)^{\nat}\times_{\os{\circ}{T}}S(T)^{\nat}\lo S(T)^{\nat}$ $(i=1,2)$. 
Then we have the section 
$\xi=q_2^*(m)q_1^*(m)^{-1}\in M_{D}$. 
(This is a different section from the section in \cite[(3.5)]{hk}.) 
Now we consider the case where $M_{S(T)^{\nat}}/{\cal O}^*_T$ 
is not necessarily isomorphic to ${\mab N}$.)  
Because the section $\xi$ is independent of the choice of $m$, 
the section $\xi$ is a global section of $M_{D}$. 
The local section $d\log m$ is independent of the choice of $m$ 
and it is a global section of $\Om^1_{S(T)^{\nat}/\os{\circ}{T}}$: 
$\Om^1_{S(T)^{\nat}/\os{\circ}{T}}={\cal O}_Td\log m$. 
The morphism 
${\cal O}_T\owns \al \lom \al d\log m\in \Om^1_{S(T)^{\nat}/\os{\circ}{T}}$ 
is independent of the choice of $m$ and this induces 
the isomorphism of ${\cal O}_T$-modules: 
${\cal O}_T\os{\sim}{\lo} \Om^1_{S(T)^{\nat}/\os{\circ}{T}}$. 
Let 
\begin{align*} 
{\rm Res} \col \Om^1_{S(T)^{\nat}/\os{\circ}{T}}\os{\sim}{\lo} {\cal O}_T
\tag{1.7.9.4}\label{ali:stt}
\end{align*} 
be the inverse of this isomorphism. 
\par 
Let $N$ be a nonnegative integer. 
Let $Y_{\bul \leq N}$ be 
an $N$-truncated simplicial log smooth scheme over $S$. 
Set $Y_{\bul \leq N,\os{\circ}{T}_0}:=Y_{\bul \leq N}\times_SS_{\os{\circ}{T}_0}
:=Y_{\bul \leq N}\times_{\os{\circ}{S}}\os{\circ}{T}_0$.  
Let 
$g \col Y_{\bul \leq N,\os{\circ}{T}_0} \lo S_{\os{\circ}{T}_0} \lo S(T)^{\nat}$ 
be the composite structural morphism. 
Let $\ol{F}{}^{\bul \leq N}$ be a flat quasi-coherent crystal of 
${\cal O}_{Y_{\bul \leq N,\os{\circ}{T}_0}/\os{\circ}{T}}$-modules. 
Set 
$F^{\bul \leq N}:=
\eps_{Y_{\bul \leq N,T_0}/S(T)^{\nat}/\os{\circ}{T}}^*(\ol{F}{}^{\bul \leq N})$ 
and 
let $\ol{F}{}^{\bul \leq N}(1)$ be the flat quasi-coherent crystal of 
${\cal O}_{Y_{\bul \leq N,\os{\circ}{T}_0}/D}$-modules obtained by 
$\ol{F}{}^{\bul \leq N}$ by using the structural morphism 
$(D,J,[~])\lo (\os{\circ}{T},{\cal J},\del)$. 
Consider the following two complexes of log crystalline cohomology sheaves: 
\begin{align*} 
{\cal K}^{\bul}:=
Rg_{Y_{\bul \leq N,\os{\circ}{T}_0}/S(T)^{\nat}*}(F^{\bul \leq N}), \quad 
{\cal K}^{\bul}(1)
:=Rg_{Y_{\bul \leq N,\os{\circ}{T}_0}/D*}(\ol{F}{}^{\bul \leq N}(1)). 
\end{align*} 
Then, by the log base change theorem ((\ref{lemm:bctn})), 
we obtain the following composite isomorphism: 
\begin{align*} 
\eps \col Lq_2^*({\cal K}^{\bul})\os{\sim}{\lo} 
{\cal K}^{\bul}(1) \os{\sim}{\longleftarrow} 
Lq_1^*({\cal K}^{\bul}).  
\tag{1.7.9.5}\label{ali:lqk}
\end{align*} 
Hence we have the following connection by (\ref{lemm:elq}): 
\begin{align*} 
\nabla \col 
Rg_{Y_{\bul \leq N,\os{\circ}{T}_0}/S(T)^{\nat}*}(F^{\bul \leq N})
\lo 
Rg_{Y_{\bul \leq N,\os{\circ}{T}_0}/S(T)^{\nat}*}(F^{\bul \leq N})
\otimes^L_{{\cal O}_T}\Om^1_{S(T)^{\nat}/\os{\circ}{T}}. 
\tag{1.7.9.6}\label{ali:lqstk}
\end{align*}

\begin{defi}
We call the morphism (\ref{ali:lqstk}) 
the {\it log crystalline Gauss-Manin connection} of 
$Rg_{Y_{\bul \leq N,\os{\circ}{T}_0}/S(T)^{\nat}*}(F^{\bul \leq N})$. 
\end{defi}

\parno
The morphisms (\ref{ali:lqstk}) and (\ref{ali:stt}) induce 
the following morphism: 
\begin{align*} 
N_{{\rm zar}}:=N_{Y_{\bul \leq N,\os{\circ}{T}_0}/S(T)^{\nat},{\rm zar}}
:={\rm Res}\circ (\nabla)
&\col Rg_{Y_{\bul \leq N,\os{\circ}{T}_0}/S(T)^{\nat}*}(F^{\bul \leq N})
\tag{1.7.10.1}\label{ali:lqnk}\\
&
\lo Rg_{Y_{\bul \leq N,\os{\circ}{T}_0}/S(T)^{\nat}*}(F^{\bul \leq N}). 
\end{align*} 
In the case $N=0$ and $F^0$ is the trivial coefficient, 
then this morphism is equal to the monodromy operator  
defined in \cite[p.~245]{hk}.  
(The notation $(u-1)^{[i]}\otimes {\cal K}$ in \cite[(3.5.1)]{hk} has no sense.)

\begin{defi}\label{defi:pmo}
(1) We call the morphism (\ref{ali:lqnk}) the {\it monodromy operator} of 
$F^{\bul \leq N}$.  
If $\ol{F}{}^{\bul \leq N}=
{\cal O}_{Y_{\bul \leq N,\os{\circ}{T}_0}/\os{\circ}{T}}$, then 
we call the following morphism 
\begin{align*} 
N_{{\rm zar}}:=N_{Y_{\bul \leq N,\os{\circ}{T}_0}/S(T)^{\nat},{\rm zar}}
&\col Rg_{Y_{\bul \leq N,\os{\circ}{T}_0}/S(T)^{\nat}*}
({\cal O}_{Y_{\bul \leq N,\os{\circ}{T}_0}/S(T)^{\nat}})
\tag{1.7.11.1}\label{ali:lqfnk}\\
&\lo Rg_{Y_{\bul \leq N,\os{\circ}{T}_0}/S(T)^{\nat}*}
({\cal O}_{Y_{\bul \leq N,\os{\circ}{T}_0}/S(T)^{\nat}})
\end{align*} 
the {\it monodromy operator} of $Y_{\bul \leq N,\os{\circ}{T}_0}/(S(T)^{\nat},{\cal J},\del)$. 
\par 
(2) 
Let $q$ be a nonnegative integer. Because the morphism 
$\os{\circ}{q}_i\col \os{\circ}{D}\lo (S(T)^{\nat})^{\circ}$ is flat, 
${\cal H}^q(Lq_i^*({\cal K}^{\bul}))=q_i^*{\cal H}^q({\cal K}^{\bul})$. 
Hence we have the following isomorphism by (\ref{ali:lqk}): 
\begin{align*} 
{\cal H}^q(\eps)\col q_2^*{\cal H}^q({\cal K}^{\bul}) \os{\sim}{\lo} 
{\cal H}^q({\cal K}^{\bul}(1)) \os{\sim}{\longleftarrow} 
q_1^*{\cal H}^q({\cal K}^{\bul}). 
\tag{1.7.11.2}\label{ali:qqqhk}
\end{align*} 
By abuse of notation, we denote 
by $N_{{\rm zar}}:=N_{Y_{\bul \leq N,\os{\circ}{T}_0}/S(T)^{\nat},{\rm zar}}$ 
the following induced morphism 
\begin{align*} 
R^qg_{Y_{\bul \leq N,\os{\circ}{T}_0}/S(T)^{\nat}*}(F^{\bul \leq N})
\lo R^qg_{Y_{\bul \leq N,\os{\circ}{T}_0}/S(T)^{\nat}*}(F^{\bul \leq N})
\tag{1.7.11.3}\label{ali:lqfnqk}
\end{align*} 
by (\ref{ali:lqfnk}).  We also call the morphism (\ref{ali:lqfnqk}) 
the {\it monodromy operator} of 
$R^qg_{Y_{\bul \leq N,\os{\circ}{T}_0}/S(T)^{\nat}*}(F^{\bul \leq N})$. 
\end{defi}

\begin{prop}[{\bf cf. \cite[p.~245]{hk}, \cite[(4.3.6)]{tsp}}]\label{prop:thkb}
The following composite morphism 
\begin{align*} 
\theta \col {\cal H}^q({\cal K}^{\bul})& \lo q_2^*{\cal H}^q({\cal K}^{\bul}) \os{\sim}{\lo} 
{\cal H}^q({\cal K}^{\bul}(1)) \os{\sim}{\longleftarrow} 
q_1^*{\cal H}^q({\cal K}^{\bul})=
{\cal O}_T\langle \xi-1\rangle\otimes_{{\cal O}_T}{\cal H}^q({\cal K}^{\bul})\\
&=\oplus_{i\in {\mab N}}{\cal O}_T(\xi-1)^{[i]}\otimes_{{\cal O}_T}{\cal H}^q({\cal K}^{\bul})  
\end{align*} 
is equal to 
\begin{align*} 
\sum_{i\in {\mab N}}(\xi-1)^{[i]}\otimes \prod_{j=0}^{i-1}(N_{\rm zar}-j). 
\tag{1.7.12.1}\label{ali:sinnq}
\end{align*}  
\end{prop} 
\begin{proof}
The proof of this proposition is the same as that of \cite[(4.3.6)]{tsp}.
\par 
Let $D(1)$ be the log PD-envelope of 
the diagonal immersion 
$S(T)^{\nat}\os{\sus}{\lo} 
S(T)^{\nat}\times_{\os{\circ}{T}}S(T)^{\nat}\times_{\os{\circ}{T}}S(T)^{\nat}$ 
over $(\os{\circ}{T},{\cal J},\del)$. 
Then, as usual, we have the ``projections'' 
$q_{ij}\col D(1)\lo D$ $(1\leq i< j\leq 3)$. 
Let $p_i\col D(1)\lo S(T)^{\nat}$ $(1\leq i\leq 3)$ be the $i$-th projection. 
Let ${\cal H}^q(\eps) \col q_2^*({\cal H}^q({\cal K}^{\bul}))\os{\sim}{\lo}  
q_1^*({\cal H}^q({\cal K}^{\bul}))$ be the isomorphism obtained by (\ref{ali:qqqhk}). 
As usual, we have the cocycle condition 
$q_{13}^*({\cal H}^q(\eps))= q_{23}^*({\cal H}^q(\eps))q_{12}^*({\cal H}^q(\eps))$ 
by the log base change formula. 
Let $\theta \col {\cal K}^{\bul} \lo 
{\cal O}_T\langle \xi-1\rangle\otimes_{{\cal O}_T}{\cal K}^{\bul}$ 
be the following composite morphism 
\begin{align*}  
{\cal K}^{\bul} \lo Lq_2^*({\cal K}^{\bul}) \os{\sim}{\lo} 
{\cal K}^{\bul}(1) \os{\sim}{\longleftarrow} 
Lq_1^*({\cal K}^{\bul})=
{\cal O}_T\langle \xi-1\rangle\otimes^L_{{\cal O}_T}{\cal K}^{\bul} 
\end{align*}
and let $\theta_i\col {\cal K}^{\bul} \lo 
{\cal O}_T(\xi-1)^{[i]}\otimes^L_{{\cal O}_T}{\cal K}^{\bul}$ 
be the following composite morphism 
\begin{align*}  
\theta_i\col {\cal K}^{\bul} \os{\theta}{\lo}  
\oplus_{i\in {\mab N}}{\cal O}_T(\xi-1)^{[i]}\otimes^L_{{\cal O}_T}{\cal K}^{\bul} 
\os{{\rm proj}.}{\lo} 
{\cal O}_T(\xi-1)^{[i]}\otimes^L_{{\cal O}_T}{\cal K}^{\bul}. 
\end{align*}
The cocycle condition is equivalent to the following equality: 
\begin{align*} 
\sum_{i\in {\mab N}}(q_{13}^*(\xi)-1)^{[i]}p_1^*\circ {\cal H}^q(\theta_i)=
\sum_{i,j\in {\mab N}}(q_{12}^*(\xi)-1)^{[i]}(q_{23}^*(\xi)-1)^{[j]}
p_1^*\circ {\cal H}^q(\theta_j)\circ {\cal H}^q(\theta_i).  
\tag{1.7.12.2}\label{ali:sinq}
\end{align*} 
Because $q_{13}^*(\xi)=q_{23}^*(\xi)q_{12}^*(\xi)$ and because 
$q_{13}^*(\xi)-1=(q_{23}^*(\xi)-1)(q_{12}^*(\xi)-1)+q_{23}^*(\xi)-1+q_{12}^*(\xi)-1$, 
the morphism containing $(q_{12}^*(\xi)-1)^{[i]}(q_{23}^*(\xi)-1)$ of the left hand side on 
(\ref{ali:sinq}) is 
\begin{align*} 
&(q_{12}^*(\xi)-1)^{[i]}(q_{23}^*(\xi)-1)p_1^*\circ {\cal H}^q(\theta_{i+1})
+(q_{12}^*(\xi)-1)^{[i-1]}(q_{12}^*(\xi)-1)(q_{23}^*(\xi)-1)p_1^*\circ {\cal H}^q(\theta_i)\\
&=(q_{12}^*(\xi)-1)^{[i]}(q_{23}^*(\xi)-1)(p_1^*\circ {\cal H}^q(\theta_{i+1})+ip_1^*\circ {\cal H}^q(\theta_i)).
\end{align*}  
On the other hand, the morphism containing 
$(q_{12}^*(\xi)-1)^{[i]}(q_{23}^*(\xi)-1)$ of the right hand side on 
(\ref{ali:sinq}) is 
$(q_{12}^*(\xi)-1)^{[i]}(q_{23}^*(\xi)-1)p_1^*\circ {\cal H}^q(\theta_1)\circ {\cal H}^q(\theta_i)$.  
Hence 
$p_1^*\circ(i{\cal H}^q(\theta_i)+{\cal H}^q(\theta_{i+1}))=p_1^*\circ {\cal H}^q(\theta_1)\circ 
{\cal H}^q(\theta_{i})$. 
Because the morphism $\os{\circ}{p}_1$ is faithfully flat as in the proof of \cite[(4.3.6)]{tsp}, 
$i{\cal H}^q(\theta_i)+{\cal H}^q(\theta_{i+1})={\cal H}^q(\theta_1)\circ {\cal H}^q(\theta_{i})$. 
Consequently, ${\cal H}^q(\theta_{i+1})=({\cal H}^q(\theta_1)-i)\circ {\cal H}^q(\theta_{i})$. 
This implies (\ref{ali:sinnq}). 
\end{proof}

\begin{rema}\label{rema:lpot}
Let the notations be as in (\ref{prop:thkb}). 
Let $\eps \col Lq_2^*({\cal K}^{\bul})\os{\sim}{\lo}  Lq_1^*({\cal K}^{\bul})$ 
be the isomorphism (\ref{ali:lqk}). 
As usual, we have the cocycle condition 
$Lq_{13}^*(\eps)= Lq_{23}^*(\eps)Lq_{12}^*(\eps)$ 
by the log base change formula. 
The cocycle condition is equivalent to the following equality 
as in the proof of (\ref{prop:thkb}): 
\begin{align*} 
\sum_{i\in {\mab N}}(q_{13}^*(\xi)-1)^{[i]}Lp_1^*\circ \theta_i=
\sum_{i,j\in {\mab N}}(q_{12}^*(\xi)-1)^{[i]}(q_{23}^*(\xi)-1)^{[j]}
Lp_1^*\circ \theta_j\circ \theta_i.  
\tag{1.7.13.1}\label{ali:sirnq}
\end{align*} 
Because $q_{13}^*(\xi)=q_{23}^*(\xi)q_{12}^*(\xi)$,  
we have the following formula: 
\begin{align*}
Lp_1^*\circ(i\theta_i+\theta_{i+1})=Lp_1^*\circ \theta_1\circ \theta_{i}.
\tag{1.7.13.2}\label{ali:sirlpnq}
\end{align*}  
However, I do not know whether the equality 
\begin{align*}
i\theta_i+\theta_{i+1}=\theta_1\circ \theta_{i}
\tag{1.7.13.3}\label{ali:sirlpanq}
\end{align*}  
holds. 
Consequently I do not know whether the equality 
\begin{align*}
\theta =\sum_{i\in {\mab N}}(\xi-1)^{[i]}\otimes^L \prod_{j=0}^{i-1}(N_{\rm zar}-j). 
\tag{1.7.13.4}\label{ali:sirbnq}
\end{align*}  
holds as morphisms 
\begin{align*} 
{\cal K}^{\bul} \lo  
{\cal O}_T\langle \xi-1\rangle\otimes^L_{{\cal O}_T}{\cal K}^{\bul}.  
\end{align*} 
This is nothing but the formula $(**)$ in \cite[p.~236]{hk} holds for 
crystals in the sense of derived category, which Hyodo and Kato has claimed 
that it holds in [loc.~cit.]. However I cannot find the proof of this fact in [loc.~cit.]. 
\end{rema} 

\begin{prop}[{\bf Contravariant functoriality of the monodromy operator}]\label{prop:otj}
Let $S'$, $(T',{\cal J}',\del')$, $T'_0$ and 
$u\col (S(T)^{\nat},{\cal J},\del) \lo (S'(T')^{\nat},{\cal J}',\del')$ be as in 
the beginning of {\rm \S\ref{sec:fcuc}}.   
Let $z\col Y_{\bul \leq N,\os{\circ}{T}_0} 
\lo Z_{\bul \leq N,\os{\circ}{T}{}'_0}$ be a morphism of 
log smooth integral schemes over the morphism 
$S_{\os{\circ}{T}_0}\lo S'_{\os{\circ}{T}{}'_0}$. 
Let $h\col Z_{\bul \leq N,\os{\circ}{T}{}'_0}\lo S'(T')^{\nat}$ 
be the structural morphism. 
Then the following diagram is commutative$:$  
\begin{equation*} 
\begin{CD} 
Ru_*Rg_{Y_{\bul \leq N,\os{\circ}{T}_0}/S(T)^{\nat}*}(F^{\bul \leq N}) 
@>{Ru_*(N_{Y_{\bul \leq N,\os{\circ}{T}_0}/S(T)^{\nat},{\rm zar}})}>>
Ru_*Rg_{Y_{\bul \leq N,\os{\circ}{T}_0}/S(T)^{\nat}*}(F^{\bul \leq N})\\
@A{z^*}AA @AA{z^*}A\\
Rh_{Z_{\bul \leq N,\os{\circ}{T}{}'_0}/S'(T')^{\nat}*}(G^{\bul \leq N}) 
@>{\deg(u)N_{Z_{\bul \leq N,\os{\circ}{T}{}'_0}/S'(T')^{\nat},{\rm zar}}}>>
Rh_{Z_{\bul \leq N,\os{\circ}{T}{}'_0}/S'(T')^{\nat}*}(G^{\bul \leq N}). 
\end{CD}
\tag{1.7.14.1}\label{eqn:ynmumss}
\end{equation*}
\end{prop}
\begin{proof} 
Let ${\cal K}'{}^{\bul}$ and ${\cal K}'{}^{\bul}(1)$ be 
the analogous complexes to 
${\cal K}^{\bul}$ and ${\cal K}^{\bul}(1)$, respectively, for 
$Z_{\bul \leq N,\os{\circ}{T}_0}/S'(T')^{\nat}$. 
Then the following diagram is commutative: 
\begin{equation*} 
\begin{CD} 
Ru_*({\cal K}{}^{\bul})@>>> Ru_*(Lq_2^*({\cal K}{}^{\bul}))@>{\sim}>>
Ru_*({\cal K}{}^{\bul}(1)) @<{\sim}<< Ru_*(Lq_1^*({\cal K}{}^{\bul}))\\
@A{z^*}AA @A{z^*}AA @AA{z^*}A @AA{z^*}A\\
{\cal K}'{}^{\bul}@>>> Lq_2^*({\cal K}'{}^{\bul})@>{\sim}>>
{\cal K}'{}^{\bul}(1) @<{\sim}<< Lq_1^*({\cal K}'{}^{\bul}). 
\end{CD}
\tag{1.7.14.2}\label{ali:lgqk}
\end{equation*} 
Hence we have the following commutative diagram: 
\begin{equation*} 
\begin{CD} 
Ru_*(Rg_{Y_{\bul \leq N,\os{\circ}{T}_0}/S(T)^{\nat}*}(F^{\bul \leq N}))
@>>> 
Ru_*(Rg_{Y_{\bul \leq N,\os{\circ}{T}_0}/S(T)^{\nat}*}(F^{\bul \leq N})
\otimes_{{\cal O}_T}^L\Om^1_{S(T)^{\nat}/\os{\circ}{T}})
\\ 
@A{z^*}AA @AA{z^*}A \\
Rh_{Z_{\bul \leq N,\os{\circ}{T}{}'_0}/S'(T')^{\nat}*}(G^{\bul \leq N})
@>>> 
Rh_{Z_{\bul \leq N,\os{\circ}{T}{}'_0}/S'(T')^{\nat}*}(G^{\bul \leq N})
\otimes^L_{{\cal O}_{T'}}\Om^1_{S'(T')^{\nat}/\os{\circ}{T}{}'}.  
\end{CD}
\tag{1.7.14.3}\label{ali:lssgqk}
\end{equation*} 
Let $m'$ be an analogous section of $M_{S'(T')^{\nat}}$ to $m$. 
Since $u^* \col M_{S'(T')^{\nat}}/{\cal O}_T^*={\mab N}
\lo {\mab N}=M_{S(T)^{\nat}}/{\cal O}_T^*$ 
is the multiplication by $\deg(u)$, 
\begin{align*} 
u^*(d\log m')=\deg(u)d\log m.
\tag{1.7.14.4}\label{ali:lgqmk}
\end{align*} 
By using (\ref{ali:lssgqk}) and (\ref{ali:lgqmk}), 
we obtain the commutative diagram (\ref{eqn:ynmumss}). 
\end{proof}

\begin{rema}\label{rema:acts}
Let the notations be as in (\ref{prop:otj}). 
\par 
(1) Assume that $T_0$ is of characteristic $p>0$ and that 
there exists a lift $u\col S(T)^{\nat}\lo S(T)^{\nat}$ 
of the Frobenius endomorphism of $S_{\os{\circ}{T}_0}$. 
Then, by (\ref{eqn:ynmumss}), 
we have the following relation 
\begin{align*} 
N_{{\rm zar}}z^*=pz^*N_{\rm zar}. 
\tag{1.7.15.1}\label{ali:zgns} 
\end{align*} 
In the case where $N=0$ and $T_0$ is the log point of 
a perfect field of characteristic $p>0$ and $T$ is the canonical lift 
of $T_0$ and $z\col Y_{0,\os{\circ}{T}_0}\lo Y_{0,\os{\circ}{T}_0}$ 
is the Frobenius endomorphism, 
Hyodo and Kato have proved a well-known relation ``$NF=pFN$'' in \cite{hk}. 
(\ref{eqn:ynmumss}) is a generalization of this relation.  
\par 
(2) When $T$ (resp.~$T'$) is restrictively hollow 
with respective to the morphism $T_0\lo S$ (resp.~$T'_0\lo S'$), 
let us denote by $q$ (resp.~$q'$) the natural morphism 
$Y_{\bul \leq N,T_0}\lo Y_{\bul \leq N,\os{\circ}{T}_0}$ 
over $T\lo S(T)$ 
(resp.~$Z_{\bul \leq N,T'_0}\lo Z_{\bul \leq N,\os{\circ}{T}_0}$ 
over $T'\lo S'(T')$). 
Then 
$$Rg_{Y_{\bul \leq N,\os{\circ}{T}_0}/S(T)^{\nat}*}(F^{\bul \leq N}) 
=Rg_{Y_{\bul \leq N,T_0}/T*}(q^*_{\rm crys}(F^{\bul \leq N}))$$
and   
$$Rh_{Z_{\bul \leq N,\os{\circ}{T}{}'_0}/S'(T')^{\nat}*}(G^{\bul \leq N}) 
=Rh_{Z_{\bul \leq N,T_0}/T*}(q'{}^*_{\rm crys}(G^{\bul \leq N}))$$
by (\ref{ali:nmzqf}).  
Hence we have the following zariskian monodromy operators
\begin{equation*}
N_{T,{\rm zar}} \col 
Rg_{Y_{\bul \leq N,T_0}/T*}(q^*_{\rm crys}(F^{\bul \leq N}))  
\lo 
Rg_{Y_{\bul \leq N,T_0}/T*}(q^*_{\rm crys}(F^{\bul \leq N})) 
\tag{1.7.15.2}\label{eqn:montynl}
\end{equation*}
and 
\begin{equation*}
N_{T',{\rm zar}} \col 
Rh_{Z_{\bul \leq N,T'_0}/T'*}(q'{}^*_{\rm crys}(G^{\bul \leq N})) 
\lo 
Rh_{Z_{\bul \leq N,T'_0}/T'*}(q'{}^*_{\rm crys}(G^{\bul \leq N}))  
\tag{1.7.15.3}\label{eqn:mlcmontynl}
\end{equation*}
fitting into the following commutative diagram$:$  
\begin{equation*} 
\begin{CD} 
Ru_*Rg_{Y_{\bul \leq N,T_0}/T*}(F^{\bul \leq N}) 
@>{Ru_*(N_{Y_{\bul \leq N,T_0}/T,{\rm zar}})}>>
Ru_*Rg_{Y_{\bul \leq N,T_0}/T*}(F^{\bul \leq N})\\
@A{z^*}AA @AA{z^*}A\\
Rh_{Z_{\bul \leq N,T'_0}/T'*}(G^{\bul \leq N}) 
@>{\deg(u)N_{Z_{\bul \leq N,T'_0}/T',{\rm zar}}}>>
Rh_{Z_{\bul \leq N,T'_0}/T'*}(G^{\bul \leq N}). 
\end{CD}
\tag{1.7.15.4}\label{eqn:ynmaumss}
\end{equation*}
\end{rema}

\begin{coro}\label{coro:st}
Set $Y_{\bul \leq N,\os{\circ}{T}_0}(1):=
Y_{\bul \leq N,\os{\circ}{T}_0}
\times_{S_{T_0}}Y_{\bul \leq N,\os{\circ}{T}_0}$. 
Let 
$\eps_{12}$ be the following pre-stratification 
\begin{align*} 
&Lq_2^*(Rg_{Y_{\bul \leq N,\os{\circ}{T}_0}/S(T)^{\nat}*}(F^{\bul \leq N}) 
\otimes^L_{{\cal O}_T}
Rg_{Y_{\bul \leq N,\os{\circ}{T}_0}/S(T)^{\nat}*}(G^{\bul \leq N})) \\
&\os{\sim}{\lo} 
Lq_1^*(Rg_{Y_{\bul \leq N,\os{\circ}{T}_0}/S(T)^{\nat}*}(F^{\bul \leq N}) 
\otimes^L_{{\cal O}_T}
Rg_{Y_{\bul \leq N,\os{\circ}{T}_0}/S(T)^{\nat}*}(G^{\bul \leq N})) 
\end{align*} 
obtained by the log K\"{u}nneth formula $(${\rm \cite[(6.12)]{klog1}}$):$  
\begin{align*} 
Rg_{Y_{\bul \leq N,\os{\circ}{T}_0}/S(T)^{\nat}*}(F^{\bul \leq N}) 
\otimes^L_{{\cal O}_T}
Rg_{Y_{\bul \leq N,\os{\circ}{T}_0}/S(T)^{\nat}*}(G^{\bul \leq N}) 
\os{\sim}{\lo} \\
Rg_{Y_{\bul \leq N,\os{\circ}{T}_0}(1)/S(T)^{\nat}*}
(F^{\bul \leq N}\boxtimes G^{\bul \leq N}) 
\end{align*} 
and the following pre-stratification 
\begin{align*} 
Lq_2^*Rg_{Y_{\bul \leq N,\os{\circ}{T}_0}(1)/S(T)^{\nat}*}
(F^{\bul \leq N}\boxtimes G^{\bul \leq N}) 
\os{\sim}{\lo}
Lq_1^*Rg_{Y_{\bul \leq N,\os{\circ}{T}_0}(1)/S(T)^{\nat}*}
(F^{\bul \leq N}\boxtimes G^{\bul \leq N}).  
\end{align*} 
Then 
$\eps_{12}=\eps_1\otimes^L \eps_2$. 
\end{coro}
\begin{proof} 
Let $q_{Y_{\bul \leq N,\os{\circ}{T}_0},i} 
\col Y_{\bul \leq N,\os{\circ}{T}_0}(1)\lo 
Y_{\bul \leq N,\os{\circ}{T}_0}$ be the $i$-th projection. 
We have the following commutative diagram 
\begin{equation*} 
\begin{CD} 
Y_{\bul \leq N,\os{\circ}{T}_0}@<{q_{Y_{\bul \leq N,\os{\circ}{T}_0},1}}<<
Y_{\bul \leq N,\os{\circ}{T}_0}(1)@>{q_{Y_{\bul \leq N,\os{\circ}{T}_0},2}}>> 
Y_{\bul \leq N,\os{\circ}{T}_0}\\
@V{g}VV @VVV @VV{g}V\\
S(T)^{\nat}@=  S(T)^{\nat} @= S(T)^{\nat}\\
@| @V{\bigcap}VV @|\\
S(T)^{\nat}@<{q_1}<< {\mathfrak E}@>{q_2}>> S(T)^{\nat}. 
\end{CD}
\end{equation*} 
Hence we have the following commutative diagram: 
\begin{equation*} 
\begin{CD} 
Lq_1^*Rg_{Y_{\bul \leq N,\os{\circ}{T}_0}(1)/S(T)^{\nat}*}
(F^{\bul \leq N}\boxtimes G^{\bul \leq N}) 
\os{\sim}{\lo}\\
@| \\
Lq_1^*(Rg_{Y_{\bul \leq N,\os{\circ}{T}_0}/S(T)^{\nat}*}(F^{\bul \leq N})
\otimes^L_{{\cal O}_T}
Rg_{Y_{\bul \leq N,\os{\circ}{T}_0}/S(T)^{\nat}*}(G^{\bul \leq N})))
\os{\sim}{\lo} \\
@| \\
Lq_1^*(Rg_{Y_{\bul \leq N,\os{\circ}{T}_0}/S(T)^{\nat}*}(F^{\bul \leq N}))
\otimes^L_{{\cal O}_{\mathfrak E}}
Lq_1^*(Rg_{Y_{\bul \leq N,\os{\circ}{T}_0}/S(T)^{\nat}*}(G^{\bul \leq N}))) 
\os{\sim}{\lo}
\end{CD}
\end{equation*} 
\begin{equation*} 
\begin{CD} 
Rg_{Y_{\bul \leq N,\os{\circ}{T}_0}(1)/{\mathfrak E}*}
(F^{\bul \leq N}\boxtimes G^{\bul \leq N}) \\
@|\\
{\cal O}_{\mathfrak E}\otimes^L_{{\cal O}_T}
(Rg_{Y_{\bul \leq N,\os{\circ}{T}_0}/S(T)^{\nat}*}(F^{\bul \leq N})
\otimes^L_{{\cal O}_T}
Rg_{Y_{\bul \leq N,\os{\circ}{T}_0}/S(T)^{\nat}*}
(G^{\bul \leq N}))\\
@|\\
{\cal O}_{\mathfrak E}\otimes^L_{{\cal O}_T}
Rg_{Y_{\bul \leq N,\os{\circ}{T}_0}/S(T)^{\nat}*}(F^{\bul \leq N})
\otimes^L_{{\cal O}_T}Rg_{Y_{\bul \leq N,\os{\circ}{T}_0}/S(T)^{\nat}*}(G^{\bul \leq N}). 
\end{CD}
\tag{1.7.16.1}\label{cd:ee}
\end{equation*} 
We also have the similar commutative diagram for $i=2$. 
These commutative diagrams tell us that 
the equality $\eps_{12}=\eps_1\otimes^L \eps_2$ holds.


\end{proof}

\begin{coro}\label{coro:fnw}
Let $\nabla_{12}$, $\nabla_1$ and $\nabla_2$ 
be the log crystalline Gauss-Manin connection of
$Rg_{Y_{\bul \leq N,\os{\circ}{T}_0(1)}/S(T)^{\nat}*}
(F^{\bul \leq N}\boxtimes G^{\bul N})$, 
$Rg_{Y_{\bul \leq N,\os{\circ}{T}_0}/S(T)^{\nat}*}(F^{\bul \leq N})$ 
and $Rg_{Y_{\bul \leq N,\os{\circ}{T}_0}/S(T)^{\nat}*}(G^{\bul \leq N})$, 
respectively. 
Then 
\begin{align*}
\nabla_{12}=
\nabla_1\otimes^L {\rm id}+
{\rm id}\otimes^L \nabla_2.  
\tag{1.7.17.1}\label{ali:oet}
\end{align*}
\end{coro} 
\begin{proof} 
This follows from (\ref{coro:st}) and (\ref{lemm:eiat}). 
\end{proof} 

\begin{coro}\label{coro:mdn}
Let $N_{12,{\rm zar}}$, $N_{1,{\rm zar}}$ and $N_{2,{\rm zar}}$ be 
the monodromy operator of 
$$Rg_{Y_{\bul \leq N,\os{\circ}{T}_0(1)}/S(T)^{\nat}*}
(F^{\bul \leq N}\boxtimes G^{\bul N}),~
Rg_{Y_{\bul \leq N,\os{\circ}{T}_0}/S(T)^{\nat}*}(F^{\bul \leq N})~
{\rm and}~
Rg_{Y_{\bul \leq N,\os{\circ}{T}_0}/S(T)^{\nat}*}(G^{\bul \leq N}),$$ 
respectively. 
Then 
\begin{align*}
N_{12,{\rm zar}}=N_{1,{\rm zar}}\otimes^L {\rm id}+
{\rm id}\otimes^L N_{2,{\rm zar}}.  
\tag{1.7.18.1}\label{ali:onet}
\end{align*}
More generally, 
\begin{equation*} 
N^i_{12,{\rm zar}}=\sum_{j=0}^i\binom{i}{j}N^{i-j}_{1,{\rm zar}}\otimes^L N^j_{2,{\rm zar}} 
\quad (i\in {\mab N}).
\tag{1.7.18.2}\label{eqn:bn12dd} 
\end{equation*} 
\end{coro}

\begin{coro}\label{coro:nnn}
Let 
\begin{align*}
\cup & \col 
Rg_{Y_{\bul \leq N,\os{\circ}{T}_0}/S(T)^{\nat}*}
({\cal O}_{Y_{\bul \leq N,\os{\circ}{T}_0}/S(T)^{\nat}})
\otimes^L_{{\cal O}_T}
Rg_{Y_{\bul \leq N,\os{\circ}{T}_0}/S(T)^{\nat}*}
({\cal O}_{Y_{\bul \leq N,\os{\circ}{T}_0}/S(T)^{\nat}})
\tag{1.7.19.1}\label{ali:rgost}\\
&\lo 
Rg_{Y_{\bul \leq N,\os{\circ}{T}_0}/S(T)^{\nat}*}
({\cal O}_{Y_{\bul \leq N,\os{\circ}{T}_0}/S(T)^{\nat}})
\end{align*}
be the cup product. 
Let 
\begin{align*} 
\nabla \col 
Rg_{Y_{\bul \leq N,\os{\circ}{T}_0}/S(T)^{\nat}*}
({\cal O}_{Y_{\bul \leq N,\os{\circ}{T}_0}/S(T)^{\nat}})
\lo 
Rg_{Y_{\bul \leq N,\os{\circ}{T}_0}/S(T)^{\nat}*}
({\cal O}_{Y_{\bul \leq N,\os{\circ}{T}_0}/S(T)^{\nat}})
\otimes^L_{{\cal O}_T}\Om^1_{S(T)^{\nat}/\os{\circ}{T}}
\end{align*} 
be the log crystalline Gauss-Manin connection of 
$Rg_{Y_{\bul \leq N,\os{\circ}{T}_0}/S(T)^{\nat}*}
({\cal O}_{Y_{\bul \leq N,\os{\circ}{T}_0}/S(T)^{\nat}})$. 
Let 
\begin{align*} 
N_{\rm zar}\col Rg_{Y_{\bul \leq N,\os{\circ}{T}_0}/S(T)^{\nat}*}
({\cal O}_{Y_{\bul \leq N,\os{\circ}{T}_0}/S(T)^{\nat}})\lo 
Rg_{Y_{\bul \leq N,\os{\circ}{T}_0}/S(T)^{\nat}*}
({\cal O}_{Y_{\bul \leq N,\os{\circ}{T}_0}/S(T)^{\nat}})
\end{align*} 
be the monodromy operator. 
Then the following hold:
\par 
$(1)$ $\nabla \circ \cup=\nabla \cup {\rm id}+{\rm id}\cup \nabla$. 
\par 
$(2)$ 
$N_{\rm zar} \circ \cup=N_{\rm zar}\cup {\rm id}+{\rm id}\cup N_{\rm zar}$. 
More generally, 
$N^i_{\rm zar} \circ \cup
=\sum_{j=0}^i\binom{i}{j}N^{i-j}_{\rm zar}\cup N^j_{\rm zar}$ $(i\in {\mab N})$. 
\end{coro} 
\begin{proof}
(1): Let $\Del \col 
Y_{\bul \leq N,\os{\circ}{T}_0}\lo 
Y_{\bul \leq N,\os{\circ}{T}_0}(1)$ 
be the diagonal immersion. 
Then the cup product (\ref{ali:rgost}) is equal to 
the following composite morphism 
\begin{align*} 
&Rg_{Y_{\bul \leq N,\os{\circ}{T}_0}/S(T)^{\nat}*}
({\cal O}_{Y_{\bul \leq N,\os{\circ}{T}_0}/S(T)^{\nat}})
\otimes^L_{{\cal O}_T}
Rg_{Y_{\bul \leq N,\os{\circ}{T}_0}/S(T)^{\nat}*}
({\cal O}_{Y_{\bul \leq N,\os{\circ}{T}_0}/S(T)^{\nat}})=\\
&Rg_{Y_{\bul \leq N,\os{\circ}{T}_0}/S(T)^{\nat}*}
({\cal O}_{Y_{\bul \leq N,\os{\circ}{T}_0}(1)/S(T)^{\nat}})
\os{\Del^*}{\lo} 
Rg_{Y_{\bul \leq N,\os{\circ}{T}_0}/S(T)^{\nat}*}
({\cal O}_{Y_{\bul \leq N,\os{\circ}{T}_0}/S(T)^{\nat}}). 
\end{align*} 
Hence (1) follows from (\ref{ali:oet}) and the commutative diagram (\ref{eqn:ynmumss}). 
\par 
(2): (2) follows from (1). 
\end{proof}

\begin{prop}\label{prop:mud}
Let 
\begin{equation*} 
\begin{CD} 
Y_{\bul \leq N,\os{\circ}{T}_0} @>{u_{Y_{\bul \leq N,\os{\circ}{T}_0}}}>> 
Y_{\bul \leq N,\os{\circ}{T}{}_0} \\ 
@VVV @VVV \\ 
S_{\os{\circ}{T}_0} @>{u_0}>> S_{\os{\circ}{T}{}_0} \\ 
@V{\bigcap}VV @VV{\bigcap}V \\ 
S(T)^{\nat} @>{u}>> S(T)^{\nat}
\end{CD}
\tag{1.7.20.1}\label{eqn:ynuuss}
\end{equation*} 
be a commutative diagram of 
$(N$-truncated simplicial$)$ log schemes.  
Then the morphism {\rm (\ref{ali:lqnk})}  
is nothing but a morphism 
\begin{equation*}
N_{S(T)^{\nat},{\rm zar}} \col 
Rg_{Y_{\bul \leq N,\os{\circ}{T}_0}/S(T)^{\nat}*}(F^{\bul \leq N}) 
\lo 
Rg_{Y_{\bul \leq N,\os{\circ}{T}_0}/S(T)^{\nat}*}(F^{\bul \leq N})(-1,u). 
\tag{1.7.20.2}\label{eqn:mlcm1glynl}
\end{equation*}
\end{prop} 
\begin{proof} 
This is a special case of (\ref{prop:otj}).  
\end{proof}

\begin{prop}\label{prop:bcgmm} 
Let $r\col (T',{\cal J}',\del')\lo (T,{\cal J},\del)$ be a morphism from an analogous 
fine log PD-scheme to $(T,{\cal J},\del)$. 
Let us also denote by $r$ the induced morphism 
$r\col (S(T')^{\nat},{\cal J}',\del')\lo (S(T)^{\nat},{\cal J},\del)$ 
by the morphism above.    
Set $Y_{\bul \leq N,\os{\circ}{T}{}'_0}:=Y_{\bul \leq N}\times_{S}S_{T_0'}$. 
Let $p\col Y_{\bul \leq N,\os{\circ}{T}{}'_0}\lo Y_{\bul \leq N,\os{\circ}{T}_0}$ 
be the natural morphism over $\os{\circ}{T}{}'\lo \os{\circ}{T}$. 
Then the following hold$:$ 
\par
$(1)$ The log crystalline Gauss-Manin connection of 
$Rg_{Y_{\bul \leq N,\os{\circ}{T}{}'_0}/S(T')^{\nat}*}(p^*_{{\rm crys}}(F^{\bul \leq N}))$ 
is equal to 
\begin{align*} 
\nabla \otimes^L {\rm id}_{{\cal O}_{T'}} \col 
Rg_{Y_{\bul \leq N,\os{\circ}{T}_0}/S(T)^{\nat}*}(F^{\bul \leq N})\otimes^L_{{\cal O}_T}{\cal O}_{T'}
\lo 
Rg_{Y_{\bul \leq N,\os{\circ}{T}_0}/S(T)^{\nat}*}(F^{\bul \leq N})
\otimes^L_{{\cal O}_T}\Om^1_{S(T)^{\nat}/\os{\circ}{T}}\otimes^L_{{\cal O}_T}{\cal O}_{T'}. 
\tag{1.7.21.1}\label{ali:lqsmgtk}
\end{align*}
\par 
$(2)$ The monodromy operator 
\begin{equation*} 
N_{S(T')^{\nat}{\rm zar}} \col 
Rg_{Y_{\bul \leq N,\os{\circ}{T}{}'_0}/S(T')^{\nat}*}(p^*_{{\rm crys}}(F^{\bul \leq N}))
\lo 
Rg_{Y_{\bul \leq N,\os{\circ}{T}{}'_0}/S(T')^{\nat}*}(p^*_{{\rm crys}}(F^{\bul \leq N}))  
\end{equation*} 
is equal to 
\begin{align*} 
N_{S(T)^{\nat},{\rm zar}}\otimes^L_{{\cal O}_T}{\cal O}_{T'} 
\col & 
Rg_{Y_{\bul \leq N,\os{\circ}{T}_0}/S(T)^{\nat}*}(F^{\bul \leq N})
\otimes^L_{{\cal O}_T}{\cal O}_{T'} 
\tag{1.7.21.2}\label{ali:mimgabv}\\
& \lo Rg_{Y_{\bul \leq N,\os{\circ}{T}_0}/S(T)^{\nat}*}(F^{\bul \leq N})
\otimes^L_{{\cal O}_T}{\cal O}_{T'}.   
\end{align*} 
\end{prop}
\begin{proof}
(1): (1) follows from (\ref{lemm:bctn}) and the construction of $\nabla$. 
\par 
(2): (2) follows from (1). 
\end{proof}

\par 
Now assume that $Y_{\bul \leq N,\os{\circ}{T}_0}$ has the disjoint union 
$Y_{\bul \leq N,\os{\circ}{T}_0}'$
of the member of an affine $N$-truncated simplicial open covering of 
$Y_{\bul \leq N,\os{\circ}{T}_0}$. 
In the following we shall give another definition of the monodromy operator 
under this assumption and we shall prove that 
this is equal to the monodromy operator (\ref{ali:lqnk}). 
This is the $N$-truncated simplicial version of the coincidence of 
the monodromy operators in \cite{hk}, 
though no detailed proof of this fact has been given in [loc.~cit.]; 
to give the proof is a nontrivial work. 
See (\ref{prop:nsg}) below for the proof.

\par 
Let $Y_{\bul \leq N,\bul,\os{\circ}{T}_0}$ be 
the \v{C}ech diagram of $Y'_{\bul \leq N,\os{\circ}{T}_0}$
over $Y_{\bul \leq N,\os{\circ}{T}_0}/S_{\os{\circ}{T}_0}$: 
$Y_{mn,\os{\circ}{T}_0}
={\rm cosk}_0^{Y_{m,\os{\circ}{T}_0}}(Y'_{m,\os{\circ}{T}_0})_n$ 
$(0\leq m\leq N, n\in {\mab N})$. 
Set $Y_{\bul \leq N,T_0}:=Y_{\bul \leq N}\times_{S}T_0$ 
and 
$Y_{\bul \leq N,\bul,T_0}:=Y_{\bul \leq N,\bul}\times_{S}T_0$. 
Let 
$g \col Y_{\bul \leq N,\os{\circ}{T}_0} 
\lo S_{\os{\circ}{T}_0} \os{\sus}{\lo} S(T)^{\nat}$ 
be the composite structural morphisms.  
Let $g_{\bul} \col Y_{\bul \leq N,\bul,\os{\circ}{T}_0} 
\lo S_{\os{\circ}{T}_0} \os{\sus}{\lo} S(T)^{\nat}$ 
and 
$g_{\bul,T} \col Y_{\bul \leq N,\bul,T_0} \lo T_0 
\os{\sus}{\lo} T$  
be also the composite structural morphisms.  
Then $g^{-1}({\cal O}_{S(T)^{\nat}})=g^{-1}_T({\cal O}_T)$ 
and $g^{-1}_{\bul}({\cal O}_{S(T)^{\nat}})=g^{-1}_{\bul,T}({\cal O}_T)$. 
\par 
For $U=T$, $S(T)^{\nat}$ or $\os{\circ}{T}$ 
and $U_0=T_0$ or $\os{\circ}{T}_0$ and 
for $\star=\bul$ or nothing, 
let 
\begin{align*}
\eps_{Y_{\bul \leq N,\star, U_0}/U} \col 
((Y_{\bul \leq N,\star, U_0}/U)_{\rm crys},{\cal O}_{Y_{\bul \leq N,\star, U_0}/U}) 
&\lo 
((\os{\circ}{Y}_{\bul \leq N,\star, U_0}/\os{\circ}{U})_{\rm crys},
{\cal O}_{\os{\circ}{Y}_{\bul \leq N,\star,U_0}/\os{\circ}{U}})
\tag{1.7.21.3}\label{eqn:tqypi}  \\
&=((\os{\circ}{Y}_{\bul \leq N,\star, T_0}/\os{\circ}{T})_{\rm crys},
{\cal O}_{\os{\circ}{Y}_{\bul \leq N,\star,T_0}/\os{\circ}{T}})
\end{align*}   
be the morphism forgetting the log structures of 
$Y_{\bul \leq N,\star, U_0}$ and $U$. 
Let 
\begin{align*} 
\eps_{Y_{\bul \leq N,\star,U_0}/\os{\circ}{T}} \col 
((Y_{\bul \leq N,\star,U_0}/\os{\circ}{T})_{\rm crys},
{\cal O}_{Y_{\bul \leq N,\star,U_0}/\os{\circ}{T}}) 
&\lo 
((\os{\circ}{Y}_{\bul \leq N,\star,U_0}/\os{\circ}{T})_{\rm crys},
{\cal O}_{\os{\circ}{Y}_{\bul \leq N,\star,U_0}
/\os{\circ}{T}})\tag{1.7.21.4}\label{ali:prtsedef}\\
&=
((\os{\circ}{Y}_{\bul \leq N,\star,T_0}/\os{\circ}{T})_{\rm crys},
{\cal O}_{\os{\circ}{Y}_{\bul \leq N,\star,T_0}/\os{\circ}{T}})
\end{align*}   
be the morphism forgetting the log structure of $Y_{\bul \leq N,\star,U_0}$. 
Let 
\begin{equation*}  
\eps_{Y_{\bul \leq N,\star,U_0}/U/\os{\circ}{T}}\col 
((Y_{\bul \leq N,\star,U_0}/U)_{\rm crys},
{\cal O}_{Y_{\bul \leq N,\star,U_0}/U})\lo 
((Y_{\bul \leq N,\star,U_0}/\os{\circ}{T})_{\rm crys},
{\cal O}_{Y_{\bul \leq N,\star,U_0}/\os{\circ}{T}})
\tag{1.7.21.5}\label{eqn:tnspi}  
\end{equation*} 
be the morphism forgetting the log structure of $U$. 
Then 
\begin{equation*}  
\eps_{Y_{\bul \leq N,\star,U_0}/U}=\eps_{Y_{\bul \leq N,\star,U_0}/\os{\circ}{T}}
\circ \eps_{Y_{\bul \leq N,\star,U_0}/U/\os{\circ}{T}}. 
\tag{1.7.21.6}\label{eqn:tnsepi}
\end{equation*} 
\par 
Let 
$Y_{\bul \leq N,\bul,\os{\circ}{T}_0}\os{\sus}{\lo} \ol{{\cal Q}}_{\bul \leq N,\bul}$ 
be an immersion into a log smooth scheme over $\ol{S(T)^{\nat}}$ 
(If each affine member of $Y'_{N,\os{\circ}{T}_0}$ is sufficiently small, 
then this immersion exists by (\ref{lemm:etl1}) and (\ref{prop:ytlft}).). 
Set ${\cal Q}_{\bul \leq N,\bul}
:=\ol{{\cal Q}}_{\bul \leq N,\bul}\times_{\ol{S(T)^{\nat}}}S(T)^{\nat}$.  
Let $\ol{\mathfrak E}_{\bul \leq N,\bul}$ 
be the log PD-envelope of 
the immersion $Y_{\bul \leq N,\bul,\os{\circ}{T}_0} 
\os{\sus}{\lo} \ol{{\cal Q}}_{\bul \leq N,\bul}$ over $(\os{\circ}{T},{\cal J},\del)$. 
Set ${\mathfrak E}_{\bul \leq N,\bul}=\ol{\mathfrak E}_{\bul \leq N,\bul}
\times_{{\mathfrak D}(\ol{S(T)^{\nat}})}S(T)^{\nat}$. 
Let $\theta_{{\cal Q}_{\bul\leq N,\bul}}\in 
{\Om}^1_{{{\cal Q}}_{\bul \leq N,\bul}/\os{\circ}{T}}$ 
be the pull-back of $\theta=``d\log \tau$''$\in  
{\Om}^1_{S(T)^{\nat}/\os{\circ}{T}}$ by the structural morphism 
${\cal Q}_{\bul \leq N,\bul}\lo S(T)^{\nat}$.  
Let $\ol{F}{}^{\bul \leq N}$ be a flat quasi-coherent crystal of 
${\cal O}_{Y_{\bul \leq N,\os{\circ}{T}_0}/\os{\circ}{T}}$-modules. 
Let $\ol{F}{}^{\bul \leq N,\bul}$ be the crystal of 
${\cal O}_{Y_{\bul \leq N,\bul,\os{\circ}{T}_0}/\os{\circ}{T}}$-modules 
obtained by $\ol{F}{}^{\bul \leq N}$. 
Let $(\ol{\cal F}^{\bul \leq N,\bul},\ol{\nabla})$ 
be the quasi-coherent 
${\cal O}_{\ol{\mathfrak E}_{\bul \leq N,\bul}}$-module 
with integrable connection corresponding to $\ol{F}{}^{\bul \leq N,\bul}$. 
Set 
$F^{\bul \leq N}:=
\eps_{Y_{\bul \leq N,\os{\circ}{T}_0}/S(T)^{\nat}/\os{\circ}{T}}^*(\ol{F}{}^{\bul \leq N})$ 
and 
$F^{\bul \leq N,\bul}:=
\eps_{Y_{\bul \leq N,\bul,\os{\circ}{T}_0}/S(T)^{\nat}/\os{\circ}{T}}^*(\ol{F}{}^{\bul \leq N,\bul})$. 
Set 
$({\cal F}^{\bul \leq N,\bul},\nabla)
=
(\ol{\cal F}^{\bul \leq N,\bul},\ol{\nabla})
\otimes_{{\cal O}_{\ol{\mathfrak E}_{\bul \leq N,\bul}}}
{\cal O}_{{\mathfrak E}_{\bul \leq N,\bul}}$.  
Let 
\begin{equation*} 
\nabla \col {\cal F}^{\bul \leq N,\bul}
\lo {\cal F}^{\bul \leq N,\bul}\otimes_{{\cal O}_{{{\cal Q}}_{\bul \leq N,\bul}}}
{\Om}^1_{{\cal Q}_{\bul \leq N,\bul}/\os{\circ}{T}} 
\tag{1.7.21.7}\label{eqn:nidfopltd}
\end{equation*}  
be the induced connection by $\nabla$ ((\ref{eqn:nidpltd})). 
Since ${\Om}^1_{{\cal Q}_{\bul \leq N,\bul}/S(T)^{\nat}}$ 
is a quotient sheaf of 
${\Om}^1_{{\cal Q}_{\bul \leq N,\bul}/\os{\circ}{T}}$, 
we also have the induced connection 
\begin{equation*} 
\nabla_{/S(T)^{\nat}} \col {\cal F}^{\bul \leq N,\bul} 
\lo {\cal F}^{\bul \leq N,\bul}
\otimes_{{\cal O}_{{{\cal Q}}_{\bul \leq N,\bul}}}
\Om^1_{{\cal Q}_{\bul \leq N,\bul}/S(T)^{\nat}}  
\tag{1.7.21.8}\label{eqn:nidfqtopltd}
\end{equation*}   
by $\nabla$.  
The object 
$({\cal F}^{\bul \leq N,\bul},\nabla_{/S(T)^{\nat}})$ 
corresponds to the log crystal $F^{\bul \leq N,\bul}$ 
of ${\cal O}_{Y_{\bul \leq N,\bul,\os{\circ}{T}_0}/S(T)^{\nat}}$-modules 
(\cite[(1.7), (6.4)]{klog1} (cf.~\cite[(2.2.7)]{nh2})).

\begin{prop}\label{prop:mce}
The following sequence 
\begin{align*} 
0 & \lo 
{\cal F}^{\bul \leq N,\bul}
\otimes_{{\cal O}_{{{\cal Q}}_{\bul \leq N,\bul}}}
{\Om}^{\bul}_{{{\cal Q}}_{\bul \leq N,\bul}/S(T)^{\nat}}[-1] 
\os{\theta_{{\cal Q}_{{\bul\leq N,\bul}}}\wedge }{\lo} 
{\cal F}^{\bul \leq N,\bul}
\otimes_{{\cal O}_{{{\cal Q}}_{\bul \leq N,\bul}}}
{\Om}^{\bul}_{{{\cal Q}}_{\bul \leq N,\bul}/\os{\circ}{T}} 
\tag{1.7.22.1}\label{eqn:gsflxd}\\ 
& \lo {\cal F}^{\bul \leq N,\bul}\otimes_{{\cal O}_{{{\cal Q}}_{\bul \leq N,\bul}}}
{\Om}^{\bul}_{{{\cal Q}}_{\bul \leq N,\bul}/S(T)^{\nat}} \lo 0
\end{align*} 
is exact. 
\end{prop}
\begin{proof} 
By (\ref{ali:eom}), 
\begin{align*} 
\nabla (\theta_{{\cal Q}_{{\bul\leq N,\bul}}}\wedge)=
\theta_{{\cal Q}_{{\bul\leq N,\bul}}}\wedge (-\nabla). 
\tag{1.7.22.2}\label{eqn:gasflxd}
\end{align*} 
Hence $\theta_{{\cal Q}_{{\bul\leq N,\bul}}}\wedge $ 
in (\ref{eqn:gsflxd}) is a morphism of complexes. 
It suffices to prove that the following sequence 
\begin{align*} 
0 & \lo {\cal O}_{{\mathfrak E}_{mn}}
\otimes_{{\cal O}_{{{\cal Q}}_{mn}}}
{\Om}^{\bul}_{{{\cal Q}}_{mn}/S(T)^{\nat}}[-1]
\os{\theta_{{\cal Q}_{mn}}\wedge}{\lo} 
{\cal O}_{{\mathfrak E}_{mn}}
\otimes_{{\cal O}_{{{\cal Q}}_{mn}}}
{\Om}^{\bul}_{{{\cal Q}}_{mn}/\os{\circ}{T}} 
\tag{1.7.22.3}\label{eqn:gslxd}\\ 
& \lo {\cal O}_{{\mathfrak E}_{mn}}
\otimes_{{\cal O}_{{{\cal Q}}_{mn}}}
{\Om}^{\bul}_{{{\cal Q}}_{mn}/S(T)^{\nat}} \lo 0
\end{align*} 
is locally split for $0\leq m\leq N$ and  $n\in {\mab N}$. 
Let 
$\ol{p}_{mn} \col \ol{{\cal Q}}_{mn} \lo \ol{S(T)^{\nat}}$ 
and 
$p_{mn} \col {\cal Q}_{mn} \lo S(T)^{\nat}$ 
be the structural morphisms. 
Since $\ol{{\cal Q}}_{mn}\lo \ol{S(T)^{\nat}}$ 
is log smooth, the first fundamental exact sequence 
\begin{equation*} 
0 \lo \ol{p}^*_{mn}({\Om}^1_{\ol{S(T)^{\nat}}/\os{\circ}{T}})\lo 
{\Om}^1_{\ol{{\cal Q}}_{mn}/\os{\circ}{T}}
\lo 
{\Om}^1_{\ol{{\cal Q}}_{mn}/\ol{S(T)^{\nat}}} \lo 0 
\tag{1.7.22.4}\label{eqn:plst}
\end{equation*}
is locally split. 
(The exact sequence 
$f^*(\om^1_{X/S})\lo \om^1_{Y/S}\lo \om^1_{X/Y}\lo 0$ 
in \cite[(3.12)]{klog1} is mistaken: 
the right one is 
$f^*(\om^1_{Y/S})\lo \om^1_{X/S}\lo \om^1_{X/Y}\lo 0$.)
Hence the following sequence 
\begin{equation*} 
0 \lo p^*_{mn}({\Om}^1_{S(T)^{\nat}/\os{\circ}{T}})\lo 
{\Om}^1_{{\cal Q}_{mn}/\os{\circ}{T}}
\lo 
{\Om}^1_{{\cal Q}_{mn}/S(T)^{\nat}} \lo 0 
\tag{1.7.22.5}\label{eqn:plpmnst}
\end{equation*}
is also locally split.  
Hence 
\begin{align*} 
{\Om}^i_{{\cal Q}_{mn}/\os{\circ}{T}}\simeq 
\theta_{{\cal Q}_{mn}}\wedge ({\Om}^{i-1}_{{\cal Q}_{mn}/S(T)^{\nat}}) 
\oplus {\Om}^i_{{\cal Q}_{mn}/S(T)^{\nat}}. 
\end{align*} 
\end{proof}

\begin{rema}\label{rema:welr} 
Let $\kap$ be a perfect field of characteristic $p>0$. 
In the case $T=S$, $\os{\circ}{S}={\rm Spec}(W_n(\kap))$, $N=0$  
and $F^{\bul \leq N}={\cal O}_{Y_{0,\os{\circ}{T}_0}/\os{\circ}{T}}$, 
the exact sequence (\ref{eqn:gsflxd}) is nothing but 
the exact sequence in the last line in \cite[p.~245]{hk}, 
though the exactness of it has not been proved in [loc.~cit].  
\end{rema}

\begin{coro}\label{coro:filt}
If the structural morphism 
$Y_{\bul \leq N,\os{\circ}{T}_0}\lo S_{\os{\circ}{T}_0}$ is integral 
and if  $\ol{F}{}^{\bul \leq N}$ is a flat crystal of 
${\cal O}_{Y_{\bul \leq N,\os{\circ}{T}_0}/\os{\circ}{T}}$-modules, 
then ${\cal F}^{\bul \leq N,\bul}
\otimes_{{\cal O}_{{{\cal Q}}_{\bul \leq N,\bul}}}
{\Om}^i_{{{\cal Q}}_{\bul \leq N,\bul}/\os{\circ}{T}}$ 
$(i\in {\mab N})$ is a sheaf of flat $g_{\bul}^{-1}({\cal O}_T)$-modules. 
\end{coro} 
\begin{proof}  
By \cite[(2.22)]{hk} or (\ref{coro:gmn}), 
${\cal F}^{\bul \leq N,\bul}
\otimes_{{\cal O}_{{{\cal Q}}_{\bul \leq N,\bul}}}
{\Om}^i_{{{\cal Q}}_{\bul \leq N,\bul}/S(T)^{\nat}}$ 
$(i\in {\mab N})$ is a sheaf of flat $g_{\bul}^{-1}({\cal O}_T)$-modules. 
Hence (\ref{coro:filt}) follows from (\ref{prop:mce}). 
\end{proof} 

\begin{coro}\label{coro:ci}
The following sequence 
\begin{align*} 
\cdots &\lo {\cal F}^{\bul \leq N,\bul}
\otimes_{{\cal O}_{{{\cal Q}}_{\bul \leq N,\bul}}}
{\Om}^{\bul}_{{{\cal Q}}_{\bul \leq N,\bul}/\os{\circ}{T}}[-2]  
\os{\theta_{{\cal Q}_{{\bul\leq N,\bul}}}\wedge }{\lo}  
{\cal F}^{\bul \leq N,\bul}
\otimes_{{\cal O}_{{{\cal Q}}_{\bul \leq N,\bul}}}
{\Om}^{\bul}_{{{\cal Q}}_{\bul \leq N,\bul}/\os{\circ}{T}}[-1]
\tag{1.7.25.1}\label{eqn:gsqlxd}\\ 
& \os{\theta_{{\cal Q}_{{\bul\leq N,\bul}}}\wedge }{\lo} 
{\cal F}^{\bul \leq N,\bul}
\otimes_{{\cal O}_{{{\cal Q}}_{\bul \leq N,\bul}}}
{\Om}^{\bul}_{{{\cal Q}}_{\bul \leq N,\bul}/\os{\circ}{T}} 
\lo {\cal F}^{\bul \leq N,\bul}\otimes_{{\cal O}_{{{\cal Q}}_{\bul \leq N,\bul}}}
{\Om}^{\bul}_{{{\cal Q}}_{\bul \leq N,\bul}/S(T)^{\nat}} \lo 0
\end{align*} 
is exact.
\end{coro}
\begin{proof}
This immediately follows from (\ref{eqn:gsflxd}). 
\end{proof}

\par 
Let 
\begin{equation*} 
\pi_{S(T)^{\nat}{\rm crys}} 
\col ((Y_{\bul \leq N,\bul,\os{\circ}{T}_0}/S(T)^{\nat})_{\rm crys},
{\cal O}_{Y_{\bul \leq N,\bul,\os{\circ}{T}_0}/S(T)^{\nat}})   
\lo ((Y_{\bul \leq N,\os{\circ}{T}_0}/S(T)^{\nat})_{\rm crys},
{\cal O}_{Y_{\bul \leq N,\os{\circ}{T}_0}/S(T)^{\nat}})
\tag{1.7.25.2}\label{eqn:tcpi} 
\end{equation*}  
and 
\begin{equation*} 
\pi_{T{\rm crys}} 
\col ((Y_{\bul \leq N,\bul,T_0}/T)_{\rm crys},{\cal O}_{Y_{\bul \leq N,\bul,T_0}/T})   
\lo ((Y_{\bul \leq N,T_0}/T)_{\rm crys},{\cal O}_{Y_{\bul \leq N,T_0}/T})
\tag{1.7.25.3}\label{eqn:tclpi} 
\end{equation*}  
be the natural morphisms of ringed topoi. 
Let 
\begin{equation*} 
\pi_{{}{\rm crys}} 
\col ((Y_{\bul \leq N,\bul,\os{\circ}{T}_0}/\os{\circ}{T})_{\rm crys},
{\cal O}_{Y_{\bul \leq N,\bul,\os{\circ}{T}_0}/\os{\circ}{T}})   
\lo ((Y_{\bul \leq N,\os{\circ}{T}_0}/\os{\circ}{T})_{\rm crys},
{\cal O}_{Y_{\bul \leq N,\os{\circ}{T}_0}/\os{\circ}{T}})
\tag{1.7.25.4}\label{eqn:tcopi} 
\end{equation*} 
and 
\begin{align*} 
\pi_{{\rm zar}} \col &
((Y_{\bul \leq N,\bul,\os{\circ}{T}_0})_{\rm zar},
g^{-1}_{\bul}({\cal O}_T))=
((Y_{\bul \leq N,\bul,T_0})_{\rm zar},g^{-1}_{\bul T}({\cal O}_T))  
\tag{1.7.25.5}\label{eqn:tzar}\\
& \lo ((Y_{\bul \leq N,T_0})_{\rm zar},g^{-1}_T({\cal O}_T))= 
((Y_{\bul \leq N,\os{\circ}{T}_0})_{\rm zar},g^{-1}({\cal O}_T)) 
\end{align*} 
be also the natural morphisms of ringed topoi. 
Then we have 
the following commutative diagrams of ringed topoi: 
\begin{equation*}  
\begin{CD} 
((Y_{\bul \leq N,\bul,T_0}/T)_{\rm crys},
{\cal O}_{Y_{\bul \leq N,\bul,T_0}/T})   
@>{\pi_{T{\rm crys}}}>> ((Y_{\bul \leq N,T_0}/T)_{\rm crys},
{\cal O}_{Y_{\bul \leq N,T_0}/T})  \\ 
@V{\eps_{Y_{\bul \leq N,\bul,T_0}/T/\os{\circ}{T}}}VV 
@VV{\eps_{Y_{\bul \leq N,T_0}/T/\os{\circ}{T}}}V \\ 
((Y_{\bul \leq N,\bul,\os{\circ}{T}_0}/\os{\circ}{T})_{\rm crys},
{\cal O}_{Y_{\bul \leq N,\bul,\os{\circ}{T}_0}/\os{\circ}{T}})  
@>{\pi_{{}{\rm crys}}}>> 
((Y_{\bul \leq N,\os{\circ}{T}_0}/\os{\circ}{T})_{\rm crys},
{\cal O}_{Y_{\bul \leq N,\os{\circ}{T}_0}/\os{\circ}{T}}),       
\end{CD} 
\tag{1.7.25.6}\label{eqn:yntpys}
\end{equation*} 
\begin{equation*}  
\begin{CD} 
((Y_{\bul \leq N,\bul,\os{\circ}{T}_0}/S(T)^{\nat})_{\rm crys},
{\cal O}_{Y_{\bul \leq N,\bul,\os{\circ}{T}_0}/S(T)^{\nat}})   
@>{\pi_{S(T)^{\nat}{\rm crys}}}>> ((Y_{\bul \leq N,\os{\circ}{T}_0}/S(T)^{\nat})_{\rm crys},
{\cal O}_{Y_{\bul \leq N,\os{\circ}{T}_0}/S(T)^{\nat}})  \\ 
@V{\eps_{Y_{\bul \leq N,\bul,\os{\circ}{T}_0}/S(T)^{\nat}/\os{\circ}{T}}}VV 
@VV{\eps_{Y_{\bul \leq N,\os{\circ}{T}_0}/S(T)^{\nat}/\os{\circ}{T}}}V \\ 
((Y_{\bul \leq N,\bul,\os{\circ}{T}_0}/\os{\circ}{T})_{\rm crys},
{\cal O}_{Y_{\bul \leq N,\bul,\os{\circ}{T}_0}/\os{\circ}{T}})  
@>{\pi_{{}{\rm crys}}}>> 
((Y_{\bul \leq N,\os{\circ}{T}_0}/\os{\circ}{T})_{\rm crys},
{\cal O}_{Y_{\bul \leq N,\os{\circ}{T}_0}/\os{\circ}{T}}),      
\end{CD} 
\tag{1.7.25.7}\label{eqn:ynpys}
\end{equation*} 

\begin{equation*}  
\begin{CD} 
((Y_{\bul \leq N,\bul,T_0}/T)_{\rm crys},
{\cal O}_{Y_{\bul \leq N,\bul,T_0}/T})   
@>{\pi_{T{\rm crys}}}>> 
((Y_{\bul \leq N,T_0}/T)_{\rm crys},{\cal O}_{Y_{\bul \leq N,T_0}/T})  
\\ 
@V{u_{Y_{\bul \leq N,\bul,T_0}/T}}VV 
@VV{u_{Y_{\bul \leq N,T_0}/T}}V \\ 
((Y_{\bul \leq N,\bul,T_0})_{\rm zar},g^{-1}_{\bul,T}({\cal O}_T))   
@>{\pi_{{\rm zar}}}>> 
((Y_{\bul \leq N,T_0})_{\rm zar},g^{-1}_{T}({\cal O}_T)),       
\end{CD} 
\tag{1.7.25.8}\label{eqn:ydtiss}
\end{equation*} 
\begin{equation*}  
\begin{CD} 
((Y_{\bul \leq N,\bul,\os{\circ}{T}_0}/S(T)^{\nat})_{\rm crys},
{\cal O}_{Y_{\bul \leq N,\bul,\os{\circ}{T}_0}/S(T)^{\nat}})   
@>{\pi_{S(T)^{\nat}{\rm crys}}}>> 
((Y_{\bul \leq N,\os{\circ}{T}_0}/S(T)^{\nat})_{\rm crys},
{\cal O}_{Y_{\bul \leq N,\os{\circ}{T}_0}/S(T)^{\nat}})  
\\ 
@V{u_{Y_{\bul \leq N,\bul,\os{\circ}{T}_0}/S(T)^{\nat}}}VV 
@VV{u_{Y_{\bul \leq N,\os{\circ}{T}_0}/S(T)^{\nat}}}V \\ 
((Y_{\bul \leq N,\bul,\os{\circ}{T}_0})_{\rm zar},g^{-1}_{\bul}({\cal O}_T))   
@>{\pi_{{\rm zar}}}>> 
((Y_{\bul \leq N,\os{\circ}{T}_0})_{\rm zar},g^{-1}({\cal O}_T)).       
\end{CD} 
\tag{1.7.25.9}\label{eqn:ydiss}
\end{equation*}

\begin{prop}\label{prop:hkt}
The complex $R\pi_{{\rm zar}*}({\cal F}^{\bul \leq N,\bul}
\otimes_{{\cal O}_{{{\cal Q}}_{\bul \leq N,\bul}}}
{\Om}^{\bul}_{{{\cal Q}}_{\bul \leq N,\bul}/\os{\circ}{T}})$ 
is independent of the choice of the disjoint union 
of the member of an affine $N$-truncated simplicial open covering of 
$Y_{\bul \leq N,\os{\circ}{T}}$ and 
an $(N,\infty)$-truncated bisimplicial immersion 
$Y_{\bul \leq N,\bul,\os{\circ}{T}} \os{\sus}{\lo} 
\ol{{\cal Q}}_{\bul \leq N,\bul}$ over $\ol{S(T)^{\nat}}$. 
\end{prop}
\begin{proof} 
By (\ref{eqn:gsflxd}) we have the following triangle 
\begin{align*} 
&R\pi_{{\rm zar}*}({\cal F}^{\bul \leq N,\bul}
\otimes_{{\cal O}_{{{\cal Q}}_{\bul \leq N,\bul}}}
{\Om}^{\bul}_{{{\cal Q}}_{\bul \leq N,\bul}/S(T)^{\nat}})[-1]
\lo 
R\pi_{{\rm zar}*}({\cal F}^{\bul \leq N,\bul}
\otimes_{{\cal O}_{{{\cal Q}}_{\bul \leq N,\bul}}}
{\Om}^{\bul}_{{{\cal Q}}_{\bul \leq N,\bul}/\os{\circ}{T}}) 
\tag{1.7.26.1}\label{eqn:gspiflxd}\\ 
& \lo 
R\pi_{{\rm zar}*}({\cal F}^{\bul \leq N,\bul}
\otimes_{{\cal O}_{{{\cal Q}}_{\bul \leq N,\bul}}}
{\Om}^{\bul}_{{{\cal Q}}_{\bul \leq N,\bul}/S(T)^{\nat}}) \os{+1}{\lo}.  
\end{align*}
Because $({\cal F}^{\bul \leq N,\bul},\nabla_{/S(T)^{\nat}})$ 
corresponds to the log crystal $F^{\bul \leq N,\bul}$ 
of ${\cal O}_{Y_{\bul \leq N,\bul,\os{\circ}{T}_0}/S(T)^{\nat}}$-modules,   
we have the following equality  
\begin{align*} 
Ru_{Y_{\bul \leq N,\bul,\os{\circ}{T}_0}/S(T)^{\nat}*}(F^{\bul \leq N,\bul})
={\cal F}^{\bul \leq N,\bul}
\otimes_{{\cal O}_{{{\cal Q}}_{\bul \leq N,\bul}}}
{\Om}^{\bul}_{{{\cal Q}}_{\bul \leq N,\bul}/S(T)^{\nat}}. 
\tag{1.7.26.2}\label{ali:lpkn}
\end{align*}  
Hence we have the following formula by \cite[(1.7)]{klog1}, 
(\ref{eqn:ydiss}) and the cohomological descent:  
\begin{align*}
R\pi_{{\rm zar}*}({\cal F}^{\bul \leq N,\bul}
\otimes_{{\cal O}_{{{\cal Q}}_{\bul \leq N,\bul}}}
{\Om}^{\bul}_{{{\cal Q}}_{\bul \leq N,\bul}/S(T)^{\nat}})
&= R\pi_{{\rm zar}*}
Ru_{Y_{\bul \leq N,\bul,\os{\circ}{T}_0}/S(T)^{\nat}*}(F^{\bul \leq N,\bul})
\tag{1.7.26.3}\label{ali:ptc}\\ 
&= Ru_{Y_{\bul \leq N,\os{\circ}{T}_0}/S(T)^{\nat}*}
R\pi_{S(T)^{\nat}{\rm crys}*}(F^{\bul \leq N,\bul})\\ 
&= Ru_{Y_{\bul \leq N,\os{\circ}{T}_0}/S(T)^{\nat}*}
R\pi_{S(T)^{\nat}{\rm crys}*}
(\pi^*_{S(T)^{\nat}{\rm crys}}(F^{\bul \leq N}))\\
&= Ru_{Y_{\bul \leq N,\os{\circ}{T}_0}/S(T)^{\nat}*}(F^{\bul \leq N}).
\end{align*} 
Thus the triangle (\ref{eqn:gspiflxd}) is equal to the following triangle: 
\begin{align*} 
&Ru_{Y_{\bul \leq N,\os{\circ}{T}_0}/S(T)^{\nat}*}(F^{\bul \leq N})[-1] 
\lo 
R\pi_{{\rm zar}*}({\cal F}^{\bul \leq N,\bul}
\otimes_{{\cal O}_{{{\cal Q}}_{\bul \leq N,\bul}}}
{\Om}^{\bul}_{{{\cal Q}}_{\bul \leq N,\bul}/\os{\circ}{T}}) 
\tag{1.7.26.4}\label{eqn:glxd}\\ 
& \lo
Ru_{Y_{\bul \leq N,\os{\circ}{T}_0}/S(T)^{\nat}*}(F^{\bul \leq N}) \os{+1}{\lo}.  
\end{align*}
\par 
Assume that we are given another affine $N$-truncated simplicial open covering of 
$Y_{\bul \leq N,\os{\circ}{T}_0}$ 
which gives another $(N,\infty)$-truncated bisimplicial immersion 
$Y'_{\bul \leq N,\bul,\os{\circ}{T}_0} \os{\sus}{\lo} 
\ol{{\cal Q}}{}'_{\bul \leq N,\bul}$ 
into a log smooth $(N,\infty)$-truncated bisimplicial scheme 
over $\ol{S(T)^{\nat}}$.  
Then, by (\ref{prop:lissmp}), we have another such 
$(N,\infty)$-truncated bisimplicial immersion   
$Y''_{\bul \leq N,\bul,\os{\circ}{T}_0} \os{\sus}{\lo} 
\ol{{\cal Q}}{}''_{\bul \leq N,\bul}$ 
fitting into the following commutative diagram over $\ol{S(T)^{\nat}}$: 
\begin{equation*}  
\begin{CD} 
Y_{\bul \leq N,\bul,\os{\circ}{T}_0} 
@>{\sus}>> \ol{{\cal Q}}_{\bul \leq N,\bul}\\ 
@AAA @AAA \\ 
Y''_{\bul \leq N,\bul,\os{\circ}{T}_0} 
@>{\sus}>> \ol{{\cal Q}}{}''_{\bul \leq N,\bul} \\
@VVV @VVV\\ 
Y'_{\bul \leq N,\bul,\os{\circ}{T}_0} 
@>{\sus}>> \ol{{\cal Q}}{}'_{\bul \leq N,\bul}.    
\end{CD} 
\tag{1.7.26.5}\label{eqn:ydss}
\end{equation*} 
By using (\ref{eqn:glxd}), the rest of the proof is a routine work. 
\end{proof}

\begin{defi}\label{defi:wrpt} 
We call 
$$R\pi_{{\rm zar}*}
({\cal F}^{\bul \leq N,\bul}
\otimes_{{\cal O}_{{{\cal Q}}_{\bul \leq N,\bul}}}
{\Om}^{\bul}_{{{\cal Q}}_{\bul \leq N,\bul}/\os{\circ}{T}})$$  
the {\it modified log crystalline complex} of $\ol{F}{}^{\bul \leq N}$ 
on $Y_{\bul \leq N,\os{\circ}{T}_0}/\os{\circ}{T}$. 
By abuse of notation, we denote it by 
$\wt{R}u_{Y_{\bul \leq N,\os{\circ}{T}_0}/\os{\circ}{T}*}(\ol{F}{}^{\bul \leq N})$ 
(this depends on the log scheme $S(T)^{\nat}$ and 
the morphism $Y_{\bul \leq N,\os{\circ}{T}_0}\lo S(T)^{\nat}$ of log schemes).  
When $\ol{F}{}^{\bul \leq N}
={\cal O}_{Y_{\bul \leq N,\os{\circ}{T}_0}/\os{\circ}{T}}$, 
we call it the {\it modified log crystalline complex} of 
$Y_{\bul \leq N,\os{\circ}{T}_0}/\os{\circ}{T}$ 
and denote it by  
$\wt{R}u_{Y_{\bul \leq N,\os{\circ}{T}_0}/\os{\circ}{T}*}
({\cal O}_{Y_{\bul \leq N,\os{\circ}{T}_0}/\os{\circ}{T}})$. 
We denote 
$Rg_*\wt{R}u_{Y_{\bul \leq N,\os{\circ}{T}_0}/\os{\circ}{T}*}(\ol{F}{}^{\bul \leq N})$ 
and ${\cal H}^i
(\wt{R}g_{Y_{\bul \leq N,\os{\circ}{T}_0}/\os{\circ}{T}*}(\ol{F}{}^{\bul \leq N}))$ 
$(i\in {\mab N})$ by 
$$\wt{R}g_{Y_{\bul \leq N,\os{\circ}{T}_0}/\os{\circ}{T}*}(\ol{F}{}^{\bul \leq N})$$ 
and 
$$\wt{R}{}^ig_{Y_{\bul \leq N,\os{\circ}{T}_0}/\os{\circ}{T}*}(\ol{F}{}^{\bul \leq N}),$$ 
respectively. 
\end{defi}

\begin{prop}\label{prop:spbc}
Let $r\col (T',{\cal J}',\del')\lo (T,{\cal J},\del)$ be a morphism from an analogous 
fine log PD-scheme to $(T,{\cal J},\del)$. 
Assume that $\os{\circ}{T}$ is quasi-compact and that $\os{\circ}{f} 
\col \os{\circ}{Y}_{\bul \leq N,T_0}\lo \os{\circ}{T}_0$ 
is quasi-compact and quasi-separated. 
Let us also denote by $r$ the induced morphism 
$r\col (S(T')^{\nat},{\cal J}',\del')\lo (S(T)^{\nat},{\cal J},\del)$ 
by the morphism above.   
Let $q\col Y_{\bul \leq N,\os{\circ}{T}{}'_0}\lo Y_{\bul \leq N,\os{\circ}{T}{}_0}$ 
be the induced morphism over 
$\os{\circ}{T}{}'\lo \os{\circ}{T}$ by $r$. 
Then 
\begin{align*} 
Lr^*(\wt{R}u_{Y_{\bul \leq N,\os{\circ}{T}_0}/\os{\circ}{T}*}(\ol{F}{}^{\bul \leq N}))
=\wt{R}u_{Y_{\bul \leq N,\os{\circ}{T}{}'_0}/\os{\circ}{T'}*}
(q^*_{\rm crys}(\ol{F}{}^{\bul \leq N})). 
\end{align*} 
\end{prop}
\begin{proof} 
It is clear that there exists a canonical morphism 
\begin{align*} 
\wt{R}u_{Y_{\bul \leq N,\os{\circ}{T}_0}/\os{\circ}{T}*}
(\ol{F}{}^{\bul \leq N})
\lo 
Rr_*\wt{R}u_{Y_{\bul \leq N,\os{\circ}{T}{}'_0}/\os{\circ}{T'}*}
(q^*_{\rm crys}(\ol{F}{}^{\bul \leq N}))
\end{align*} 
by the definitions of 
$\wt{R}u_{Y_{\bul \leq N,\os{\circ}{T}_0}/\os{\circ}{T}*}
(\ol{F}{}^{\bul \leq N})$ and 
$\wt{R}u_{Y_{\bul \leq N,\os{\circ}{T}{}'_0}/\os{\circ}{T'}*}
(q^*_{\rm crys}(\ol{F}{}^{\bul \leq N}))$. 
Hence we have the following canonical morphism 
\begin{align*} 
Lr^*(\wt{R}u_{Y_{\bul \leq N,\os{\circ}{T}_0}/\os{\circ}{T}*}
(\ol{F}{}^{\bul \leq N}))
\lo 
\wt{R}u_{Y_{\bul \leq N,\os{\circ}{T}{}'_0}/\os{\circ}{T'}*}
(q^*_{\rm crys}(\ol{F}{}^{\bul \leq N})). 
\end{align*} 
To prove that this is an isomorphism, it suffices to assume that $N=0$ and 
that there exists a lift $\ol{\cal Q}$ of $Y_{0,\os{\circ}{T}_0}$
into a log smooth scheme over $\ol{S(T)^{\nat}}$ by 
the triangle (\ref{eqn:glxd}) and the cohomological descent. 
Hence we can use the simplified notations 
$Y$ for $Y_0$, $F$ for $F^0$, 
${\cal Q}$ for ${\cal Q}_{\bul \leq N,\bul}$ and  
${\cal F}$ for ${\cal F}^{\bul \leq N,\bul}$, respectively.   
Set ${\cal Q}_{\os{\circ}{T}{}'}:={\cal Q}\times_{\os{\circ}{T}}\os{\circ}{T}{}'$. 
Then ${\cal Q}_{\os{\circ}{T}{}'}$ is a lift of $Y_{0,\os{\circ}{T}{}'_0}$. 
Because  
\begin{align*} 
Ru_{Y_{\os{\circ}{T}{}'_0}/\os{\circ}{T}{}'*}(q^*_{\rm crys}(F))=q^*({\cal F})
\otimes_{{\cal O}_{{\cal Q}_{\os{\circ}{T}{}'}}}
\Om^{\bul}_{{\cal Q}_{\os{\circ}{T}{}'}/\os{\circ}{T}{}'}
\end{align*} 
and 
\begin{align*} 
Ru_{Y_{\os{\circ}{T}_0}/\os{\circ}{T}*}(F)=
{\cal F}\otimes_{{\cal O}_{\cal Q}}\Om^{\bul}_{{\cal Q}/\os{\circ}{T}},  
\end{align*} 
the desired equality is obvious since 
${\cal F}\otimes_{{\cal O}_{\cal Q}}
{\Om}^i_{{\cal Q}/\os{\circ}{T}}$ 
$(i\in {\mab N})$ is a sheaf of flat $g^{-1}({\cal O}_T)$-modules 
((\ref{coro:filt})). 
\end{proof}

\begin{prop}[\bf Contravariant functoriality of 
$\wt{R}g_{Y_{\bul \leq N,\os{\circ}{T}_0}/\os{\circ}{T}*}(\ol{F}{}^{\bul \leq N})$]\label{prop:cttu}
$(1)$ Let $u\col (S(T)^{\nat},{\cal J},\del) \lo (S'(T')^{\nat},{\cal J}',\del')$ 
be as in the beginning of {\rm \S\ref{sec:fcuc}}. 
Let 
\begin{equation*} 
\begin{CD} 
Y_{\bul \leq N,\os{\circ}{T}_0} @>{h_{\bul \leq N}}>> Z_{\bul \leq N,\os{\circ}{T}{}'_0}\\ 
@VVV @VVV \\ 
S_{\os{\circ}{T}_0} @>>> S'_{\os{\circ}{T}{}'_0} 
\end{CD}
\tag{1.7.29.1}\label{eqn:xdss}
\end{equation*} 
be a commutative diagram of 
$N$-truncated simplicial log smooth schemes over 
$S_{\os{\circ}{T}_0}$ and $S'_{\os{\circ}{T}{}'_0}$ 
such that $Y_{\bul \leq N,\os{\circ}{T}_0}$ and $Z_{\bul \leq N,\os{\circ}{T}{}'_0}$ 
have the disjoint unions 
$Y'_{\bul \leq N,\os{\circ}{T}_0}$ and $Z'_{\bul \leq N,\os{\circ}{T}{}'_0}$ 
of the members of affine $N$-truncated simplicial open coverings of 
$Y_{\bul \leq N,\os{\circ}{T}_0}$ and $Z_{\bul \leq N,\os{\circ}{T}{}'_0}$, 
respectively, fitting into the following commutative diagram 
\begin{equation*} 
\begin{CD} 
Y'_{\bul \leq N,\os{\circ}{T}_0} @>{h'_{\bul \leq N}}>> 
Z'_{\bul \leq N,\os{\circ}{T}{}'_0} \\
@VVV @VVV \\ 
Y_{\bul \leq N,\os{\circ}{T}_0} @>{h_{\bul \leq N}}>> Z_{\bul \leq N,\os{\circ}{T}{}'_0}. 
\end{CD}
\tag{1.7.29.2}\label{cd:xxy}
\end{equation*} 
Let $\ol{G}{}^{\bul \leq N}$ be 
a flat quasi-coherent crystal of 
${\cal O}_{Z_{\bul \leq N,\os{\circ}{T}{}'_0}/\os{\circ}{T}{}'}$-modules. 
Let 
$$j^{\bul \leq N} \col h^*_{\bul \leq N,{\rm crys}}(\ol{G}{}^{\bul \leq N})\lo 
\ol{F}{}^{\bul \leq N}$$ 
be a morphism of flat quasi-coherent crystals of 
${\cal O}_{Y_{\bul \leq N,\os{\circ}{T}}/\os{\circ}{T}}$-modules. 
Then there exists a natural morphism 
\begin{align*} 
\wt{j}{}^{\bul \leq N} \col 
\wt{R}u_{Z_{\bul \leq N,\os{\circ}{T}{}'_0}/\os{\circ}{T}{}'*}
(\ol{G}{}^{\bul \leq N})
\lo 
Rh_{\bul \leq N*}\wt{R}u_{Y_{\bul \leq N,\os{\circ}{T}_0}/\os{\circ}{T}*}
(\ol{F}{}^{\bul \leq N}), 
\tag{1.7.29.3}\label{ali:xxwy}
\end{align*} 
which is compatible with the compositions of $h_{\bul \leq N*}$'s and 
$j^{\bul \leq N}$'s, 
and $$\wt{({\rm id}_{Y_{\bul \leq N,\os{\circ}{T}_0}})}{}^*
={\rm id}_{\wt{R}u_{Y_{\bul \leq N,\os{\circ}{T}_0}/\os{\circ}{T}*}(\ol{F}{}^{\bul \leq N})}.$$ 
\par 
$(2)$ Set 
$G^{\bul \leq N}:=\eps_{Z_{\bul \leq N,T'_0}/S'(T')^{\nat}/\os{\circ}{T}}^*
(\ol{G}{}^{\bul \leq N})$.  
Then the morphism $\wt{j}{}^{\bul \leq N}$ fits into the following commutative diagram of 
triangles$:$ 
\begin{equation*} 
\begin{CD} 
Rh_{\bul \leq N*}Ru_{Y_{\bul \leq N,\os{\circ}{T}_0}/S(T)^{\nat}*}(F^{\bul \leq N})[-1] 
@>>> Rh_{\bul \leq N*}R\wt{u}_{Y_{\bul \leq N,\os{\circ}{T}_0}/\os{\circ}{T}*}
(\ol{F}{}^{\bul \leq N})
@>>> \\ 
@A{j^{\bul \leq N}}AA @AA{\wt{j}{}^{\bul \leq N}}A \\
Ru_{Z_{\bul \leq N,\os{\circ}{T}{}'_0}/S'(T')^{\nat}*}(G^{\bul \leq N})[-1]
@>>> 
R\wt{u}_{Z_{\bul \leq N,\os{\circ}{T}{}'_0}/\os{\circ}{T}{}'*}(\ol{G}{}^{\bul \leq N})
@>>> 
\end{CD} 
\tag{1.7.29.4}\label{cd:xaxy}
\end{equation*} 
\begin{equation*} 
\begin{CD} 
Rh_{\bul \leq N*}Ru_{Y_{\bul \leq N,\os{\circ}{T}_0}/S(T)^{\nat}*}(F^{\bul \leq N}) 
@>{+1}>> \\
@A{j^{\bul \leq N}}AA \\
Ru_{Z_{\bul \leq N,\os{\circ}{T}{}'_0}/S'(T')^{\nat}*}(G^{\bul \leq N})@>{+1}>>.
\end{CD} 
\end{equation*} 
\end{prop} 
\begin{proof} 
(1): We obtain (1) by the argument before 
(\ref{theo:funas}) and the proof of (\ref{theo:funas}).  
Indeed, the proof of (1) is simpler than that of (\ref{theo:funas}) as follows. 
\par 
We have the following commutative diagram 
\begin{equation*} 
\begin{CD} 
Y_{\bul \leq N,\bul,\os{\circ}{T}_0}
@>{\sus}>> \ol{\cal Q}_{\bul \leq N,\bul}\\
@V{h_{\bul \leq N,\bul}}VV 
@VV{\ol{h}_{\bul \leq N,\bul}}V \\ 
Z_{\bul \leq N,\bul,\os{\circ}{T}{}'_0} @>{\sus}>> \ol{\cal R}_{\bul \leq N,\bul}
\end{CD}  
\end{equation*} 
over 
\begin{equation*} 
\begin{CD} 
S_{\os{\circ}{T}_0} @>{\subset}>> \ol{S(T)^{\nat}}\\ 
@VVV @VVV \\ 
S'_{\os{\circ}{T}{}'_0} @>{\subset}>> \ol{S'(T')^{\nat}}, 
\end{CD}
\end{equation*} 
where the horizontal morphisms 
are immersions into $(N,\infty)$-truncated bisimplicial 
log smooth schemes over $\ol{S(T)^{\nat}}$ and $\ol{S'(T')^{\nat}}$, respectively. 
Set ${\cal R}_{\bul \leq N,\bul}
:=\ol{\cal R}_{\bul \leq N,\bul}\times_{\ol{S'(T')^{\nat}}}S'(T')^{\nat}$. 
Let $\ol{\mathfrak E}_{\bul \leq N,\bul}$ 
and $\ol{\mathfrak F}_{\bul \leq N,\bul}$ be the log PD-envelopes of  
the immersions 
$Y_{\bul \leq N,\bul,\os{\circ}{T}_0}
\os{\sus}{\lo}\ol{\cal Q}_{\bul \leq N,\bul}$ and 
$Z_{\bul \leq N,\bul,\os{\circ}{T}{}'_0}
\os{\sus}{\lo}\ol{\cal R}_{\bul \leq N,\bul}$ over 
$(\os{\circ}{T},{\cal J},\del)$ and $(\os{\circ}{T}{}',{\cal J}',\del')$, respectively.     
Set ${\mathfrak F}_{\bul \leq N,\bul}:=
\ol{\mathfrak F}_{\bul \leq N,\bul}\times_{{\mathfrak D}(\ol{S'(T')^{\nat}})}S'(T')^{\nat}$. 
We have the following natural morphism   
$\ol{h}{}^{\rm PD}_{\bul \leq N,\bul}\col 
\ol{\mathfrak E}_{\bul \leq N,\bul}\lo \ol{\mathfrak F}_{\bul \leq N,\bul}$.   
Hence we have the following natural morphism  
$h^{\rm PD}_{\bul \leq N,\bul}\col 
{\mathfrak E}_{\bul \leq N,\bul}\lo {\mathfrak F}_{\bul \leq N,\bul}$.   
For $\ol{G}^{\bul \leq N}$, let $(\ol{\cal G}{}^{\bul \leq N,\bul},\nabla)$ 
and $({\cal G}{}^{\bul \leq N,\bul},\nabla)$ 
be the similar objects to 
$(\ol{\cal F}{}^{\bul \leq N,\bul},\nabla)$ and 
$({\cal F}{}^{\bul \leq N,\bul},\nabla)$, respectively.  
\par 
As the proof of (\ref{theo:funas}), we can construct the following morphism  
\begin{equation*} 
{\cal G}^{\bul \leq N,\bul}
\otimes_{{\cal O}_{{\cal R}{}^{\rm ex}_{\bul \leq N,\bul}}}
\Om^{\bul}_{{\cal R}{}^{\rm ex}_{\bul \leq N,\bul}/\os{\circ}{T}{}'}
\lo 
h^{\rm PD}_{\bul \leq N,\bul*}
({\cal F}^{\bul \leq N,\bul}
\otimes_{{\cal O}_{{\cal Q}^{\rm ex}_{\bul \leq N,\bul}}}
\Om^{\bul}_{{\cal Q}^{\rm ex}_{\bul \leq N,\bul}/\os{\circ}{T}})
\end{equation*} 
of complexes. 
This morphism induces the morphism (\ref{ali:xxwy}). 
We leave the detail of the rest of the proof to the reader. 
\par 
(2): (2) follows from the definition of $\wt{j}{}^{\bul \leq N}$. 
\end{proof}

\par 
Let 
\begin{equation*} 
{\cal F}^{\bul \leq N,\bul}
\otimes_{{\cal O}_{{{\cal Q}}_{\bul \leq N,\bul}}}
{\Om}^{\bul}_{{{\cal Q}}_{\bul \leq N,\bul}/S(T)^{\nat}} 
\lo 
{\cal F}^{\bul \leq N,\bul}
\otimes_{{\cal O}_{{{\cal Q}}_{\bul \leq N,\bul}}}
{\Om}^{\bul}_{{{\cal Q}}_{\bul \leq N,\bul}/S(T)^{\nat}} 
\tag{1.7.29.5}\label{eqn:gyglxd}
\end{equation*} 
be the boundary morphism of (\ref{eqn:gsflxd}) 
in the derived category 
${\rm D}^+(f^{-1}_{\bul \leq N,\bul}({\cal O}_T))$. 
(We make the convention on the sign of 
the boundary morphism as in \cite[p.~12 (4)]{nh2}.)  
This morphism induces the following morphism by (\ref{ali:ptc}): 
\begin{equation*}
N'_{\rm zar}:=N'_{Y_{\bul \leq N,\os{\circ}{T}_0}/S(T)^{\nat},{\rm zar}} \col 
Ru_{Y_{\bul \leq N,\os{\circ}{T}_0}/S(T)^{\nat}*}(F^{\bul \leq N})\lo 
Ru_{Y_{\bul \leq N,\os{\circ}{T}_0}/S(T)^{\nat}*}(F^{\bul \leq N}).
\tag{1.7.29.6}\label{eqn:nzgslyne}
\end{equation*}

\begin{prop}\label{prop:qftn}
The morphism {\rm (\ref{eqn:nzgslyne})} is well-defined. 
\end{prop} 
\begin{proof} 
The proof of this proposition is the same as that of (\ref{prop:hkt}). 
\end{proof} 

\parno  
By the definitions in (\ref{defi:wrpt}), 
we have the following triangle 
and the long exact sequence: 
\begin{align*} 
& \wt{R}u_{Y_{\bul \leq N,\os{\circ}{T}_0}/\os{\circ}{T}*}(\ol{F}{}^{\bul \leq N}) \lo 
Ru_{Y_{\bul \leq N,\os{\circ}{T}_0}/S(T)^{\nat}*}(F^{\bul \leq N})
\os{N'_{\rm zar}}{\lo} 
Ru_{Y_{\bul \leq N,\os{\circ}{T}_0}/S(T)^{\nat}*}(F^{\bul \leq N}) 
\os{+1}{\lo}\tag{1.7.30.1}\label{ali:eynt}  
\end{align*}  
and 
\begin{align*} 
& \cdots \lo R^{q-1}g_{Y_{\bul \leq N,\os{\circ}{T}_0}/S(T)^{\nat}*}(F^{\bul \leq N})
\lo 
\wt{R}^qg_{Y_{\bul \leq N,\os{\circ}{T}_0}/\os{\circ}{T}*}(\ol{F}{}^{\bul \leq N})
\lo 
\tag{1.7.30.2}\label{ali:lceeynt}\\
& R^qg_{Y_{\bul \leq N,\os{\circ}{T}_0}/S(T)^{\nat}*}(F^{\bul \leq N})
\os{N'_{\rm zar}}{\lo} 
R^qg_{Y_{\bul \leq N,\os{\circ}{T}_0}/S(T)^{\nat}*}(F^{\bul \leq N})
\lo \cdots. 
\end{align*}

\par 
Let the notations be as in (\ref{prop:otj}). 
\par 
Assume that there exist another disjoint union 
$Y'_{\bul \leq N,\os{\circ}{T}_0}$ (resp.~$Z'_{\bul \leq N,\os{\circ}{T}{}'_0}$)
of the members of affine $N$-truncated simplicial open coverings of 
$Y_{\bul \leq N,\os{\circ}{T}_0}$ (resp.~$Z_{\bul \leq N,\os{\circ}{T}_0}$) 
and a morphism 
$z' 
\col Y'_{\bul \leq N,\bul,\os{\circ}{T}_0} \lo Z'_{\bul \leq N,\os{\circ}{T}_0}$ 
fitting into the following commutative diagram: 
\begin{equation*} 
\begin{CD} 
Y'_{\bul \leq N,\os{\circ}{T}_0} @>{z'}>> Z'_{\bul \leq N,\os{\circ}{T}_0} \\ 
@VVV @VVV \\ 
Y_{\bul \leq N,\os{\circ}{T}_0} @>{z}>> Z_{\bul \leq N,\os{\circ}{T}_0}. 
\end{CD}
\tag{1.7.30.3}\label{eqn:ydgguss}
\end{equation*} 
Then the following diagram is commutative$:$  
\begin{equation*} 
\begin{CD} 
Rz_*Ru_{Y_{\bul \leq N,\os{\circ}{T}_0}/S(T)^{\nat}*}(F^{\bul \leq N}) 
@>{Ru_*(N'_{Y_{\bul \leq N,\os{\circ}{T}_0}/S(T)^{\nat},{\rm zar}})}>>
Rz_*Ru_{Y_{\bul \leq N,\os{\circ}{T}_0}/S(T)^{\nat}*}(F^{\bul \leq N})\\
@A{z^*}AA @AA{z^*}A\\
Ru_{Z_{\bul \leq N,\os{\circ}{T}{}'_0}/S'(T')^{\nat}*}(G^{\bul \leq N}) 
@>{\deg(u)N'_{Z_{\bul \leq N,\os{\circ}{T}{}'_0}/S'(T')^{\nat},{\rm zar}}}>>
Ru_{Z_{\bul \leq N,\os{\circ}{T}{}'_0}/S'(T')^{\nat}*}(G^{\bul \leq N}). 
\end{CD}
\tag{1.7.30.4}\label{eqn:yumss}
\end{equation*}
since $u^*(\theta_{{\cal Q}_{\bul \leq N,\bul}})=
\deg(u)\theta_{{\cal R}_{\bul \leq N,\bul}}$.
Here ${\cal R}_{\bul \leq N,\bul}$ is an analogous log scheme 
to ${\cal Q}_{\bul \leq N,\bul}$ for $Z_{\bul \leq N,\bul}$ 
fitting into the following commutative diagram
\begin{equation*} 
\begin{CD} 
Y_{\bul \leq N,\bul,\os{\circ}{T}_0} @>{\subset}>> {\cal Q}_{\bul \leq N,\bul}\\
@VVV @VVV\\
Z_{\bul \leq N,\bul,\os{\circ}{T}_0} @>{\subset}>> {\cal R}_{\bul \leq N,\bul}. 
\end{CD}
\tag{1.7.30.5}\label{eqn:yumyzss}
\end{equation*}
over 
\begin{equation*} 
\begin{CD} 
S_{\os{\circ}{T}_0} @>{\subset}>> S(T)^{\nat}\\
@VVV @VVV\\
S'_{\os{\circ}{T}{}'_0} @>{\subset}>> S'(T')^{\nat}. 
\end{CD}
\end{equation*}
In particular, in the situation (\ref{prop:mud}) with 
the existences of $Y'_{\bul \leq N}$ and $Y''_{\bul \leq N}$ 
fitting into the following commutative diagram 
\begin{equation*} 
\begin{CD} 
Y'_{\bul \leq N,\os{\circ}{T}_0} @>>> Y''_{\bul \leq N,\os{\circ}{T}_0}\\
@VVV @VVV\\
Y_{\bul \leq N,\os{\circ}{T}_0} 
@>{u_{Y_{\bul \leq N,\os{\circ}{T}_0}}}>> Y_{\bul \leq N,\os{\circ}{T}_0},  
\end{CD}
\tag{1.7.30.6}\label{eqn:yumyzs}
\end{equation*}
then the morphism (\ref{eqn:nzgslyne}) 
is nothing but a morphism 
\begin{equation*}
N'_{\rm zar}:=N'_{Y_{\bul \leq N,\os{\circ}{T}_0}/S(T)^{\nat},{\rm zar}} 
\col Ru_{Y_{\bul \leq N,\os{\circ}{T}_0}/S(T)^{\nat}*}(F^{\bul \leq N}) 
\lo 
Ru_{Y_{\bul \leq N,\os{\circ}{T}_0}/S(T)^{\nat}*}(F^{\bul \leq N})(-1,u). 
\tag{1.7.30.7}\label{eqn:mlcepglynl}
\end{equation*}
Here $Y''_{\bul \leq N,\os{\circ}{T}_0}$ is another 
disjoint union of the members of 
affine $N$-truncated simplicial open coverings of 
$Y_{\bul \leq N,\os{\circ}{T}_0}$. 
This morphism induces the following morphism 
\begin{equation*}
N'_{\rm zar}:=N'_{Y_{\bul \leq N,\os{\circ}{T}_0}/S(T)^{\nat},{\rm zar}} 
\col Rg_{Y_{\bul \leq N,\os{\circ}{T}_0}/S(T)^{\nat}*}(F^{\bul \leq N}) 
\lo 
Rg_{Y_{\bul \leq N,\os{\circ}{T}_0}/S(T)^{\nat}*}(F^{\bul \leq N})(-1,u). 
\tag{1.7.30.8}\label{eqn:mgglynl}
\end{equation*}

\begin{defi}\label{defi:ctm}
We call the morphism (\ref{eqn:mlcepglynl}) (resp.~(\ref{eqn:mgglynl})) 
the {\it  zariskian monodromy operator} (resp.~{\it  monodromy operator})
of $F^{\bul \leq N}$. 
When $F^{\bul \leq N}={\cal O}_{Y_{\bul \leq N,\os{\circ}{T}_0}/\os{\circ}{T}}$, 
we call the morphism (\ref{eqn:mlcepglynl}) the 
{\it  zariskian monodromy operator} (resp.~{\it  monodromy operator})
of $Y_{\bul \leq N,\os{\circ}{T}_0}/T$. 
(We shall show that the monodromy operator (\ref{eqn:mgglynl}) is equal to 
the operator (\ref{ali:lqnk}).)
\end{defi}

\begin{prop}\label{prop:nsg}
Let the notations be as after {\rm (\ref{lemm:bctn})}. 
Assume that $Z_{\bul \leq N}$ has the disjoint union $Z_{\bul \leq N}'$
of the member of an affine $N$-truncated simplicial open covering of 
$Z_{\bul \leq N}$.  Then there exists the following connection 
\begin{align*} 
\nabla \col R^qg_{Z_{\bul \leq N}/{\cal Y}*}(E^{\bul \leq N})\lo 
R^qg_{Z_{\bul \leq N}/{\cal Y}*}(E^{\bul \leq N})
\otimes_{{\cal O}_{\cal Y}}\Om^1_{{\cal Y}/U},   
\tag{1.7.32.1}\label{ali:zyen}
\end{align*} 
which is $($essentially$)$ equal to 
the log version of the Gauss-Manin connection in {\rm \cite{ko}} 
in the case $N=0$, ${\cal I}=0$ and ${\cal Y}=Y$. 
\end{prop}
\begin{proof} 
Set $Z_{mn}:={\rm cosk}^{Z_m}_0(Z'_m)_n$ $(0\leq m\leq N,n\in {\mab N})$. 
Let $Z_{\bul \leq N,\bul}\os{\sus}{\lo} {\cal R}_{\bul \leq N,\bul}$ 
be an immersion into a log smooth scheme over ${\cal Y}$.  
Let 
$g_{\bul} \col {\cal R}_{\bul \leq N,\bul}\lo {\cal Y}$ be the structural morphism. 
If each affine member of $Z'_N$ is sufficiently small, 
this immersion indeed exists by (\ref{lemm:etl1}) and (\ref{prop:ytlft}). 
Let ${\mathfrak E}_{\bul \leq N,\bul}$ be the log PD-envelope of 
the immersion $Z_{\bul \leq N,\bul}\os{\sus}{\lo} {\cal R}_{\bul \leq N,\bul}$ 
over $({\cal Y},{\cal J},\del)$. 
Let $E^{\bul \leq N,\bul}$ be the crystal of 
${\cal O}_{Z_{\bul \leq N,\bul}/{\cal Y}}$-modules obtained by $E^{\bul \leq N}$. 
Let $({\cal E}^{\bul \leq N,\bul},\nabla)$ 
be the quasi-coherent 
${\cal O}_{{\mathfrak E}_{\bul \leq N,\bul}}$-module 
with integrable connection corresponding to $E^{\bul \leq N,\bul}$. 
Since ${\cal R}_{\bul \leq N,\bul}/{\cal Y}$ is log smooth, 
the following sequence 
\begin{align*} 
0\lo g^*_{\bul}(\Om^1_{{\cal Y}/U})\lo 
\Om^1_{{\cal R}_{\bul \leq N,\bul}/U}\lo 
\Om^1_{{\cal R}_{\bul \leq N,\bul}/{\cal Y}}\lo 0. 
\end{align*} 
is exact (\cite[(3.12)]{klog1}).  
Consider the following filtration 
\begin{align*} 
{\rm Fil}^i({\cal E}^{\bul \leq N,\bul}
\otimes_{{\cal O}_{{\cal R}_{\bul \leq N,\bul}}}
\Om^{\bul}_{{\cal R}_{\bul \leq N,\bul}/U})
:=&{\rm Im}({\cal E}^{\bul \leq N,\bul}
\otimes_{{\cal O}_{{\cal R}_{\bul \leq N,\bul}}} 
\Om^{\bul-i}_{{\cal R}_{\bul \leq N,\bul}/U}
\otimes_
{{\cal O}_{{\cal R}_{\bul \leq N,\bul}}}g^*_{\bul}(\Om^i_{{\cal Y}/U})
\\
&\lo {\cal E}^{\bul \leq N,\bul}\otimes_
{{\cal O}_{{\cal R}_{\bul \leq N,\bul}}}
\Om^{\bul}_{{\cal R}_{\bul \leq N,\bul}/U}). 
\end{align*} 
Because $\Om^1_{{\cal Y}/U}$ and 
$\Om^1_{{\cal R}_{\bul \leq N,\bul}/{\cal Y}}$ 
are locally free modules of finite rank, 
\begin{align*} 
{\rm gr}^i({\cal E}^{\bul \leq N,\bul}
\otimes_{{\cal O}_{{\cal R}_{\bul \leq N,\bul}}}
\Om^{\bul}_{{\cal R}_{\bul \leq N,\bul}/U})
={\cal E}^{\bul \leq N,\bul}\otimes_
{{\cal O}_{{\cal R}_{\bul \leq N,\bul}}}
\Om^{\bul-i}_{{\cal R}_{\bul \leq N,\bul}/{\cal Y}}
\otimes_
{{\cal O}_{{\cal R}_{\bul \leq N,\bul}}}
g^*_{\bul}(\Om^i_{{\cal Y}/U}). 
\tag{1.7.32.2}\label{ali:pigrr}
\end{align*} 
Let $\pi_{\rm zar} \col (Y_{\bul \leq N,\bul})_{\rm zar}
\lo (Y_{\bul \leq N})_{\rm zar}$ 
be the natural morphism. 
Then we have the following spectral sequence by (\ref{ali:pigrr}): 
\begin{align*} 
E_1^{ij}&=
R^j(g\pi_{\rm zar})_*({\cal E}^{\bul \leq N,\bul}
\otimes_{{\cal O}_{{\cal R}_{\bul \leq N,\bul}}}
\Om^{\bul}_{{\cal R}_{\bul \leq N,\bul}/{\cal Y}})
\otimes_{{\cal O}_{\cal Y}}
\Om^i_{{\cal Y}/U}\tag{1.7.32.3}\label{ali:pierr}\\
&\Lo 
R^{i+j}(g\pi_{\rm zar})_*({\cal E}^{\bul \leq N,\bul}
\otimes_{{\cal O}_{{\cal R}_{\bul \leq N,\bul}}}
\Om^{\bul}_{{\cal R}_{\bul \leq N,\bul}/U}).
\end{align*} 
This spectral sequence is equal to the following spectral sequence: 
\begin{align*} 
E_1^{ij}=
R^jg_{Z_{\bul \leq N}/{\cal Y}*}(E^{\bul \leq N})
\otimes_{{\cal O}_{\cal Y}}
\Om^i_{{\cal Y}/U}
\Lo 
R^{i+j}(g\pi_{\rm zar})_*({\cal E}^{\bul \leq N,\bul}
\otimes_{{\cal O}_{{\cal R}_{\bul \leq N,\bul}}}
\Om^{\bul}_{{\cal R}_{\bul \leq N,\bul}/U}).
\tag{1.7.32.4}\label{ali:piesrr}
\end{align*} 
The boundary morphism 
$d_1^{0j}\col E_1^{0j}\lo E_1^{1j}$ induces the morphism (\ref{ali:zyen}). 
\par 
We claim that this morphism is independent of the choices of 
the disjoint union of the member of an affine $N$-truncated simplicial 
open covering of $Z_{\bul \leq N}$ and the immersion 
$Z_{\bul \leq N,\bul}\os{\sus}{\lo} {\cal R}_{\bul \leq N,\bul}$. 
Because this is a routine work, we leave the proof of the independence 
to the reader.  
\par 
In the case $N=0$, ${\cal I}=0$ and ${\cal Y}=Y$, 
the connection $d_1^{0j}\col E_1^{0j}\lo E_1^{1j}$ 
is $($essentially$)$ 
the log version of the Gauss-Manin connection in {\rm \cite{ko}}. 
\end{proof} 

\begin{coro}\label{coro:inds}
The complex $R(g\pi_{\rm zar})_*({\cal E}^{\bul \leq N,\bul}
\otimes_{{\cal O}_{{\cal R}_{\bul \leq N,\bul}}}
\Om^{\bul}_{{\cal R}_{\bul \leq N,\bul}/U})$ is independent of 
the choices of 
the disjoint union of the member of an affine $N$-truncated simplicial 
open covering of $Z_{\bul \leq N}$ and the immersion 
$Z_{\bul \leq N,\bul}\os{\sus}{\lo} {\cal R}_{\bul \leq N,\bul}$. 
\end{coro}
\begin{proof}
This follows from (\ref{ali:piesrr}). 
\end{proof} 

The following is an $N$-truncated simplicial log version of 
a generalization of \cite[V (3.6.4)]{bb}:

\begin{prop}\label{prop:gmn}
Let the notations be as in {\rm (\ref{prop:nsg})}. 
Then the connections {\rm (\ref{ali:zymen})} and 
{\rm (\ref{ali:zyen})} are equal.  
\end{prop} 
\begin{proof} 
(The proof is an $N$-truncated simplicial log version of 
a generalization of the proof of \cite[V (3.6.4)]{bb}.) 
Let the notations be as in the proof of (\ref{prop:nsg}). 

\par 
Let $j_{{\cal R}_{\bul \leq N,\bul}} 
\col (Z_{\bul \leq N,\bul}/U)_{\rm crys}
\vert_{{\cal R}_{\bul \leq N,\bul}} \lo (Z_{\bul \leq N,\bul}/U)_{\rm crys}$ 
be the localization of log crystalline topoi.

Let 
\begin{align*}
\varphi^* 
\col & \{\text{the category of }
{\cal O}_{{\cal R}_{\bul \leq N,\bul}}\text{-modules}\}  \\
{} & \lo \{\text{the category of log crystals of }
{\cal O}_{Z_{\bul \leq N,\bul}/U}\vert_{{\cal R}_{\bul \leq N,\bul}}\text{-modules}\}
\end{align*}
be a natural functor which is defined as follows as usual: 
for a coherent ${\cal O}_{{\cal R}_{\bul \leq N,\bul}}$-module 
${\cal F}^{\bul \leq N,\bul}$ and for an object $T$ of 
$(Z_{\bul \leq N,\bul}/U)_{\rm crys}\vert_{{\cal R}_{\bul \leq N,\bul}}$ 
with a morphism $T\lo {\cal R}_{\bul \leq N,\bul}$,  
$\varphi^*({\cal F}^{\bul \leq N,\bul}):=
{\cal O}_T\otimes_{{\cal O}_{{\cal R}_{\bul \leq N,\bul}}}{\cal F}^{\bul \leq N,\bul}$.  
Set 
\begin{align*} 
L_{Z_{\bul \leq N,\bul}/U}:=j_{{\cal R}_{\bul \leq N,\bul}}\circ \varphi^* \col 
& \{\text{the category of log crystals of }
{\cal O}_{{\cal R}_{\bul \leq N,\bul}}\text{-modules}\}\\
&\lo  \{\text{the category of log crystals of }
{\cal O}_{Z_{\bul \leq N,\bul}/U}\text{-modules}\}
\end{align*}
as usual (a log linearization functor). 
We can also define the following functor as above: 
\begin{align*} 
L_{Y/U} \col 
& \{\text{the category of }
{\cal O}_{{\cal Y}}\text{-modules}\}\\
&\lo  \{\text{the category of log crystals of }
{\cal O}_{Y/U}\text{-modules}\}. 
\end{align*}
Set 
\begin{align*} 
&{\rm Fil}^iL_{Z_{\bul \leq N,\bul}/U}({\cal E}^{\bul \leq N,\bul}
\otimes_{{\cal O}_{{\cal R}_{\bul \leq N,\bul}}}
\Om^{\bul}_{{\cal R}_{\bul \leq N,\bul}/U}) \\
&:={\rm Im}(L_{Z_{\bul \leq N,\bul}/U}({\cal E}^{\bul \leq N,\bul}
\otimes_{{\cal O}_{{\cal R}_{\bul \leq N,\bul}}}
\Om^{\bul-i}_{{\cal R}_{\bul \leq N,\bul}/U}
\otimes_{{\cal O}_{{\cal R}_{\bul \leq N,\bul}}}g_{\bul}^*(\Om^i_{{\cal Y}/U}))\\
&\lo 
L_{Z_{\bul \leq N,\bul}/U}({\cal E}^{\bul \leq N,\bul}
\otimes_{{\cal O}_{{\cal R}_{\bul \leq N,\bul}}}
\Om^{\bul}_{{\cal R}_{\bul \leq N,\bul}/U})) \quad 
(i\in {\mab Z}). 
\end{align*} 
Then 
\begin{align*} 
&{\rm gr}^i
L_{Z_{\bul \leq N,\bul}/U}({\cal E}^{\bul \leq N,\bul}
\otimes_{{\cal O}_{{\cal R}_{\bul \leq N,\bul}}}
\Om^{\bul}_{{\cal R}_{\bul \leq N,\bul}/U})\\
&=
L_{Z_{\bul \leq N,\bul}/U}({\cal E}^{\bul \leq N,\bul}
\otimes_{{\cal O}_{{\cal R}_{\bul \leq N,\bul}}}
\Om^{\bul-i}_{{\cal R}_{\bul \leq N,\bul}/{\cal Y}}
\otimes_{{\cal O}_{{\cal R}_{\bul \leq N,\bul}}}g_{\bul}^*(\Om^i_{{\cal Y}/U}))
\end{align*} 
because $L_{Z_{\bul \leq N,\bul}/U}$ is an exact functor (\cite[(2.2.4)]{nh2}). 
Hence we have the following spectral sequence 
\begin{align*} 
&E_1^{ij}=R^{j}(u_{Y/U}\pi_{\rm crys}
g_{\bul {\rm crys}})_*L_{Z_{\bul \leq N,\bul}/U}({\cal E}^{\bul \leq N,\bul}
\otimes_{{\cal O}_{{\cal R}_{\bul \leq N,\bul}}}
\Om^{\bul-i}_{{\cal R}_{\bul \leq N,\bul}/{\cal Y}}
\otimes_{{\cal O}_{{\cal R}_{\bul \leq N,\bul}}}g_{\bul}^*(\Om^i_{{\cal Y}/U}))
\tag{1.7.34.1}\label{ali:ygc}\\
&\Lo 
R^{i+j}(u_{Y/U}\pi_{\rm crys}g_{\bul {\rm crys}})_*
L_{Z_{\bul \leq N,\bul}/U}({\cal E}^{\bul \leq N,\bul}
\otimes_{{\cal O}_{{\cal R}_{\bul \leq N,\bul}}}
\Om^{\bul}_{{\cal R}_{\bul \leq N,\bul}/U}). 
\end{align*} 
Because 
\begin{align*} 
E_1^{ij}&=R^{j}(\pi_{\rm zar}g_{\bul}u_{Z_{\bul \leq N,\bul}/U})_*
L_{Z_{\bul \leq N,\bul}/U}({\cal E}^{\bul \leq N,\bul}
\otimes_{{\cal O}_{{\cal R}_{\bul \leq N,\bul}}}
\Om^{\bul-i}_{{\cal R}_{\bul \leq N,\bul}/{\cal Y}}
\otimes_{{\cal O}_{{\cal R}_{\bul \leq N,\bul}}}g_{\bul}^*(\Om^i_{{\cal Y}/U}))\\
&=
R^{j}(\pi_{\rm zar}g_{\bul})_*
({\cal E}^{\bul \leq N,\bul}
\otimes_{{\cal O}_{{\cal R}_{\bul \leq N,\bul}}}
\Om^{\bul-i}_{{\cal R}_{\bul \leq N,\bul}/{\cal Y}}
\otimes_{{\cal O}_{{\cal R}_{\bul \leq N,\bul}}}g_{\bul}^*(\Om^i_{{\cal Y}/U}))\\
&= R^j(\pi_{\rm zar}g_{\bul})_*
({\cal E}^{\bul \leq N,\bul}
\otimes_{{\cal O}_{{\cal R}_{\bul \leq N,\bul}}}
\Om^{\bul-i}_{{\cal R}_{\bul \leq N,\bul}/{\cal Y}})
\otimes_{{\cal O}_{\cal Y}}\Om^i_{{\cal Y}/U},
\end{align*} 
the edge morphism 
$E^{ij}\lo E^{i+1,j}$ of the spectral sequence 
(\ref{ali:ygc}) is equal to that of (\ref{ali:pierr}). 
\par 
Following the proof of \cite[Proposition 3.6.4]{bb}, 
we claim that 
$$Rg_{\bul, {\rm crys}*}
(L_{Z_{\bul \leq N,\bul}/U}({\cal E}^{\bul \leq N,\bul}
\otimes_{{\cal O}_{{\cal R}_{\bul \leq N,\bul}}}
\Om^{\bul-i}_{{\cal R}_{\bul \leq N,\bul}/{\cal Y}}
\otimes_{{\cal O}_{{\cal R}_{\bul \leq N,\bul}}}g_{\bul}^*(\Om^i_{{\cal Y}/U})))$$
is equal to  
$$Rg_{\bul, {\rm crys}*}
(L_{Z_{\bul \leq N,\bul}/U}({\cal E}^{\bul \leq N,\bul}
\otimes_{{\cal O}_{{\cal R}_{\bul \leq N,\bul}}}
\Om^{\bul-i}_{{\cal R}_{\bul \leq N,\bul}/{\cal Y}}))
\otimes_{{\cal O}_{Y/U}}L_{Y/U}(\Om^i_{{\cal Y}/U}).$$ 
First we have the following morphism 
\begin{align*} 
E^{\bul \leq N,\bul}\otimes^L_{{\cal O}_{Z_{\bul \leq N,\bul}/U}}
Lg^*_{\bul,{\rm crys}}(L_{Y/U}(\Om^i_{{\cal Y}/U})) 
\lo & 
E^{\bul \leq N,\bul}\otimes_{{\cal O}_{Z_{\bul \leq N,\bul}/U}}
L_{Z_{\bul \leq N,\bul}/U}(g^*_{\bul}(\Om^i_{{\cal Y}/U}))
\tag{1.7.34.2}\label{ali:ygrc}\\
&=L_{Z_{\bul \leq N,\bul}/U}({\cal E}^{\bul \leq N,\bul}
\otimes_{{\cal O}_{{\cal R}_{\bul \leq N,\bul}/U}}g^*_{\bul}(\Om^i_{{\cal Y}/U}))
\end{align*} 
by using \cite[(2.2.4)]{nh2}.
Hence we have the following morphism 
\begin{align*} 
&E^{\bul \leq N,\bul}
\otimes^L_{{\cal O}_{Z_{\bul \leq N,\bul}/U}}
Lg^*_{\bul,{\rm crys}}(L_{Y/U}(\Om^i_{{\cal Y}/U}))
\lo \\
&L_{Z_{\bul \leq N,\bul}/U}
({\cal E}^{\bul \leq N,\bul}\otimes_{{\cal O}_{{\cal R}_{\bul \leq N,\bul}/U}}
\Om^{\bul}_{{\cal R}_{\bul \leq N,\bul}/{\cal Y}}
\otimes_{{\cal O}_{{\cal R}_{\bul \leq N,\bul}/U}}g^*_{\bul}(\Om^i_{{\cal Y}/U})). 
\end{align*} 
Because $E^{\bul \leq N,\bul}=
L_{Z_{\bul \leq N,\bul}/U}(
{\cal E}^{\bul \leq N,\bul}\otimes_{{\cal O}_{{\cal R}_{\bul \leq N,\bul}/U}}
\Om^{\bul}_{{\cal R}_{\bul \leq N,\bul}/U})$ 
by the log crystalline Poincar\'{e} lemma (\cite[(2.2.8)]{nh2}), 
we have the following morphism
\begin{align*} 
&L_{Z_{\bul \leq N,\bul}/U}(
{\cal E}^{\bul \leq N,\bul}\otimes_{{\cal O}_{{\cal R}_{\bul \leq N,\bul}/U}}
\Om^{\bul}_{{\cal R}_{\bul \leq N,\bul}/U})
\otimes^L_{{\cal O}_{Z_{\bul \leq N,\bul}/U}}
Lg^*_{\bul,{\rm crys}}(L_{Y/U}(\Om^i_{{\cal Y}/U}))
\lo \tag{1.7.34.3}\label{ali:lyybu}\\
&L_{Z_{\bul \leq N,\bul}/U}(
{\cal E}^{\bul \leq N,\bul}\otimes_{{\cal O}_{{\cal R}_{\bul \leq N,\bul}/U}}
\Om^{\bul}_{{\cal R}_{\bul \leq N,\bul}/{\cal Y}}
\otimes_{{\cal O}_{{\cal R}_{\bul \leq N,\bul}/U}}
g^*_{\bul}(\Om^i_{{\cal Y}/U}))
\end{align*}
by (\ref{ali:ygrc}). 
By using the adjunction morphism 
\begin{align*} 
&Lg^*_{\bul,{\rm crys}}Rg_{\bul,{\rm crys}*}
(L_{Z_{\bul \leq N,\bul}/U}({\cal E}^{\bul \leq N,\bul}
\otimes_{{\cal O}_{{\cal R}_{\bul \leq N,\bul}/U}}
\Om^{\bul}_{{\cal R}_{\bul \leq N,\bul}/U}))\\
&\lo 
L_{Z_{\bul \leq N,\bul}/U}(
{\cal E}^{\bul \leq N,\bul}\otimes_{{\cal O}_{{\cal R}_{\bul \leq N,\bul}/U}}
\Om^{\bul}_{{\cal R}_{\bul \leq N,\bul}/U})
\end{align*}
and the morphism (\ref{ali:lyybu}), 
we have the following morphism 
\begin{align*} 
&Lg^*_{\bul,{\rm crys}}Rg_{\bul,{\rm crys}*}(L_{Z_{\bul \leq N,\bul}/U}
({\cal E}^{\bul \leq N,\bul}\otimes_{{\cal O}_{{\cal R}_{\bul \leq N,\bul}/U}}
\Om^{\bul}_{{\cal R}_{\bul \leq N,\bul}/U}))
\otimes^L_{{\cal O}_{Z_{\bul \leq N,\bul}/U}}
Lg^*_{\bul,{\rm crys}}(L_{Y/U}(\Om^i_{{\cal Y}/U}))\tag{1.7.34.4}\label{ali:lynyu}
\\
&\lo L_{Z_{\bul \leq N,\bul}/U}
({\cal E}^{\bul \leq N,\bul}\otimes_{{\cal O}_{{\cal R}_{\bul \leq N,\bul}/U}}
\Om^{\bul}_{{\cal R}_{\bul \leq N,\bul}/{\cal Y}}
\otimes^L_{{\cal O}_{{\cal R}_{\bul \leq N,\bul}/U}}
g^*_{\bul}(\Om^i_{{\cal Y}/U})). 
\end{align*}
The morphism 
give us the following morphism 
\begin{align*} 
&
Rg_{\bul,{\rm crys}*}(L_{Z_{\bul \leq N,\bul}/U}
({\cal E}^{\bul \leq N,\bul}
\otimes_{{\cal O}_{{\cal R}_{\bul \leq N,\bul}/U}}
\Om^{\bul}_{{\cal R}_{\bul \leq N,\bul}/U}))\otimes^L_{{\cal O}_{Y/U}}
L_{Y/U}(\Om^i_{{\cal Y}/U})
\lo \tag{1.7.34.5}\label{ali:lyyu}\\
&Rg_{\bul,{\rm crys}*}(
L_{Z_{\bul \leq N,\bul}/U}
({\cal E}^{\bul \leq N,\bul}\otimes_{{\cal O}_{{\cal R}_{\bul \leq N,\bul}/U}}
\Om^{\bul}_{{\cal R}_{\bul \leq N,\bul}/{\cal Y}}
\otimes_{{\cal O}_{{\cal R}_{\bul \leq N,\bul}}}
g^*_{\bul}(\Om^i_{{\cal Y}/U})). 
\end{align*} 
We claim that (\ref{ali:lyyu}) is an isomorphism. 
It suffices to prove that the value of both hand sides on 
(\ref{ali:lyyu}) at $({\cal Y},{\cal Y},0)$ are equal. 
The value of 
the left hand side on (\ref{ali:lyyu}) at $({\cal Y},{\cal Y},0)$
is equal to 
$$Rg_{\bul*}({\cal E}^{\bul \leq N,\bul}
\otimes_{{\cal O}_{{\cal R}_{\bul \leq N,\bul}/{\cal Y}}}
\Om^{\bul}_{{\cal R}_{\bul \leq N,\bul}/{\cal Y}})\otimes_{{\cal O}_{\cal Y}}
{\cal O}_{\mathfrak D}\otimes_{{\cal O}_{\cal Y}}
\Om^i_{{\cal Y}/U}.$$
On the other hand, the value of 
the right hand side on (\ref{ali:lyyu}) at $({\cal Y},{\cal Y},0)$
is equal to 
\begin{align*} 
&Rg_{\bul*}(
{\cal E}^{\bul \leq N,\bul}\otimes_{{\cal O}_{{\cal R}_{\bul \leq N,\bul}/U}}
\Om^{\bul}_{{\cal R}_{\bul \leq N,\bul}/{\cal Y}}
\otimes_{{\cal O}_{{\cal R}_{\bul \leq N,\bul}/U}}
g^*_{\bul}(\Om^i_{{\cal Y}/U}))\\
&=
Rg_{\bul*}(
{\cal E}^{\bul \leq N,\bul}\otimes_{{\cal O}_{{\cal R}_{\bul \leq N,\bul}/U}}
\Om^{\bul}_{{\cal R}_{\bul \leq N,\bul}/{\cal Y}})
\otimes_{{\cal O}_{\cal Y}}\Om^i_{{\cal Y}/U}.
\end{align*} 
Because they are isomorphic, we see that 
the morphism (\ref{ali:lyyu}) is an isomorphism. 
Thus we have proved that 
\begin{align*} 
E^{ij}_1=
{\cal H}^j(Ru_{Y/S*}(
Rg_{\bul,{\rm crys}*}(E^{\bul \leq N,\bul})
\otimes^L_{{\cal O}_{Y/U}}L_{Y/U}(\Om^i_{{\cal Y}/U}))). 
\end{align*} 
Consider the morphism 
\begin{align*} 
{\rm id}\otimes^L L_{Y/U}(d) \col & 
Rg_{\bul,{\rm crys}*}(E^{\bul \leq N,\bul})
\otimes^L_{{\cal O}_{Z_{\bul \leq N,\bul}/U}}
L_{Y/U}(\Om^i_{{\cal Y}/U})\lo \\
& Rg_{\bul,{\rm crys}*}(E^{\bul \leq N,\bul})
\otimes^L_{{\cal O}_{Z_{\bul \leq N,\bul}/U}}L_{Y/U}(\Om^{i+1}_{{\cal Y}/U}). 
\end{align*}
Then, as in \cite[p.~362]{bb}, we claim that the following diagram is commutative: 
\begin{equation*} 
\begin{CD}
R^{i+j}g_{\bul,{\rm crys}*}
({\rm gr}^iL_{Z_{\bul \leq N,\bul}/U}({\cal E}^{\bul \leq N,\bul}
\otimes_{{\cal O}_{{\cal R}_{\bul \leq N,\bul}}}
\Om^{\bul}_{{\cal R}_{\bul \leq N,\bul}/U}))
@>{\sim}>> \\
@V{d_1^{ij}}VV  \\
R^{i+j+1}g_{\bul,{\rm crys}*}
({\rm gr}^{i+1}L_{Z_{\bul \leq N,\bul}/U}({\cal E}^{\bul \leq N,\bul}
\otimes_{{\cal O}_{{\cal R}_{\bul \leq N,\bul}}}
\Om^{\bul}_{{\cal R}_{\bul \leq N,\bul}/U}))
@>{\sim}>> 
\end{CD}
\tag{1.7.34.6}\label{cd:hdgf}
\end{equation*} 
\begin{equation*} 
\begin{CD}
R^jg_{\bul,{\rm crys}*}(E^{\bul \leq N,\bul})
\otimes_{{\cal O}_{Z_{\bul \leq N,\bul}/U}}
L_{Y/U}(\Om^i_{{\cal Y}/U})\\
@VV{{\rm id}\otimes L_{Y/U}(d)}V\\
R^jg_{\bul,{\rm crys}*}(E^{\bul \leq N,\bul})
\otimes_{{\cal O}_{Z_{\bul \leq N,\bul}/U}}
L_{Y/U}(\Om^{i+1}_{{\cal Y}/U}). 
\end{CD} 
\end{equation*} 
Here the edge morphism $d_1^{ij}$ above is obtained by 
the linearization of the differential 
$d\col \Om^i_{{\cal Y}/U}\lo \Om^{i+1}_{{\cal Y}/U}$. 
Indeed, we have the natural morphism
\begin{align*} 
E^{\bul \leq N,\bul}
\otimes_{{\cal O}_{Z_{\bul \leq N,\bul}/U}}
g^*_{\bul,{\rm crys}}(L_{Y/U}(\Om^{\bul}_{{\cal Y}/U}))
\lo 
L_{Z_{\bul \leq N,\bul}/U}({\cal E}^{\bul \leq N,\bul}
\otimes_{{\cal O}_{{\cal R}_{\bul \leq N,\bul}/U}}
\Om^{\bul}_{{\cal R}_{\bul \leq N,\bul}/U}). 
\end{align*} 
In fact, this is a filtered morphism if one endow 
$E^{\bul \leq N,\bul}
\otimes_{{\cal O}_{Z_{\bul \leq N,\bul}/U}}
g^*_{\bul,{\rm crys}}(L_{Y/U}(\Om^{\bul}_{{\cal Y}/U}))$ 
with 
the induced filtration by the Hodge filtration 
on $\Om^{\bul}_{{\cal Y}/U}$. 
Hence we have the commutative diagram 
(\ref{cd:hdgf}) by the lemma below. 
\end{proof} 

\begin{lemm}[{\bf An obvious log version of \cite[IV Proposition 3.1.4]{bb}}]\label{lemm:lb}
Let $E$ be a crystal of ${\cal O}_{Y/U}$-module and 
let ${\cal E}\lo {\cal E}\otimes_{{\cal O}_{\cal Y}}\Om^1_{{\cal Y}/U}$ 
be the integrable connection corresponding $E$. 
Let ${\cal G}$ be an ${\cal O}_{\cal Y}$-module with the hyper-PD-log stratification 
with respect to $U$. Then the following hold$:$ 
\par 
$(1)$ 
There exists the following canonical isomorphism
\begin{align*} 
L_{Y/U}({\cal E}\otimes_{{\cal O}_{\cal Y}}{\cal G})
\os{\sim}{\lo} E\otimes_{{\cal O}_{Y/U}}L_{Y/U}({\cal G}). 
\tag{1.7.35.1}\label{cd:hdglf}
\end{align*} 
\par 
$(2)$ Let $u$ be a hyper-differential operator from 
${\cal F}\lo {\cal G}$ and let $v$ be the induced 
hyper-differential operator from 
${\cal E}\otimes_{{\cal O}_{\cal Y}}{\cal F}\lo 
{\cal E}\otimes_{{\cal O}_{\cal Y}}{\cal G}$ by $u$. 
Then the following diagram is commutative$:$
 \begin{equation*} 
\begin{CD}
L_{Y/U}({\cal E}\otimes_{{\cal O}_{\cal Y}}{\cal F})
@>{L_{Y/U}(v)}>> L_{Y/U}({\cal E}\otimes_{{\cal O}_{\cal Y}}{\cal G})\\
@V{\simeq}VV @VV{\simeq}V \\
E\otimes_{{\cal O}_{Y/U}}L_{Y/U}({\cal F})
@>{{\rm id}\otimes L_{Y/U}(u)}>>
E\otimes_{{\cal O}_{Y/U}}L_{Y/U}({\cal G}). 
\end{CD}
\tag{1.7.35.2}\label{cd:hdgff}
\end{equation*} 
\end{lemm}
\begin{proof}
The proof is the same as that of \cite[pp.272--273]{bb}. 
\end{proof} 

\begin{prop}\label{prop:theprop}
Let the notations be as above. 
Let 
\begin{equation*}
N_{S(T)^{\nat},{\rm zar}} \col 
Rg_{Y_{\bul \leq N,\os{\circ}{T}_0}/S(T)^{\nat}*}(F^{\bul \leq N}) 
\lo 
Rg_{Y_{\bul \leq N,\os{\circ}{T}_0}/S(T)^{\nat}*}(F^{\bul \leq N})(-1,u)
\tag{1.7.36.1}\label{eqn:mluglynl}
\end{equation*}
be the indued morphism by {\rm (\ref{eqn:mlcepglynl})}.  
Then this morphism is equal to the monodromy operator {\rm (\ref{ali:lqnk})}. 
\end{prop}
\begin{proof} 
This immediately follows from (\ref{prop:gmn}).  
\end{proof} 

We can prove the formula (\ref{eqn:bndd}) for the monodromy operator 
by using the following lemma without using (\ref{prop:theprop}):   

\begin{lemm}\label{lemm:lbg} 
Let $({\cal T},{\cal A})$ be a ringed topos. 
Let $\Om^{\bul}$ and $\wt{\Om}^{\bul}$ be dga's over 
${\cal A}$ with nonnegative degrees.
Let $\theta \in \Gam({\cal T},\wt{\Om}^1)$ be a global section.   
Assume that there exists the following  exact sequence
\begin{equation*} 
0\lo \Om^{\bul}[-1]\os{\theta \wedge}{\lo} \wt{\Om}^{\bul}\lo \Om^{\bul}\lo 0 
\tag{1.7.37.1}\label{eqn:otw} 
\end{equation*}
such that, for each $q\in {\mab N}$, 
\begin{equation*} 
0\lo \Om^{q-1}\os{\theta \wedge}{\lo} \wt{\Om}^q\lo \Om^q\lo 0 
\tag{1.7.37.2}\label{eqn:otyw} 
\end{equation*} 
is locally split. 
Let 
\begin{equation*} 
{\cal N} \col H^q({\cal T},\Om^{\bul})\lo H^q({\cal T},\Om^{\bul})
\tag{1.7.37.3}\label{eqn:bdd} 
\end{equation*} 
be the boundary morphism of {\rm (\ref{eqn:otw})}.  
Let $x$ and $y$ be elements of 
$H^q({\cal T},\Om^{\bul})$ and $H^{q'}({\cal T},\Om^{\bul})$, respectively. 
Then 
\begin{equation*} 
{\cal N}^i(x\cup y)=\sum_{j=0}^i\binom{i}{j}{\cal N}^{i-j}(x)\cup {\cal N}^j(y) 
\quad (i\in {\mab N}).
\tag{1.7.37.4}\label{eqn:bndad} 
\end{equation*}  
\end{lemm}
\begin{proof} 
It suffices to prove (\ref{eqn:bndad}) for the case $i=1$.  
Let $q$ and $q'$ be nonnegative integers. 
Let $U$ be an object of ${\cal T}$ such that the sequence
\begin{equation*} 
0\lo \Om^r(U)[-1]\os{\theta \wedge}{\lo} \wt{\Om}^r(U)\lo \Om^r(U)\lo 0 
\end{equation*}
for $r=q$ and $q'$ is split. 
Let $\iota \col \Om^r(U)\lo \wt{\Om}^r(U)$ 
be a splitting ($r=q$, $q'$). 
Let $x$ and $y$ be elements of ${\rm Ker}(\Om^q(U)\lo \Om^{q+1}(U))$ 
and ${\rm Ker}(\Om^{q'}(U)\lo \Om^{q'+1}(U))$, respectively. 
Let $x'$ and $y'$ be elements of $\Om^q(U)$ and $\Om^{q'}(U)$ 
such that $\theta \wedge x'=\theta \wedge \iota(x')=d\iota(x)$ and 
$\theta \wedge y'=\theta \wedge \iota(y')=d\iota(y)$, respectively. 
Then 
\begin{align*}
d(\iota(x)\wedge \iota(y))& =d\iota(x)\wedge \iota(y)+(-1)^q\iota(x)\wedge d\iota(y)\\
&=\theta \wedge \iota(x')\wedge \iota(y)
+(-1)^q\iota(x)\wedge \theta \wedge \iota(y')\\
&=\theta \wedge \iota(x')\wedge \iota(y)+ \theta \wedge \iota(x)  \wedge \iota(y')\\
& =\theta \wedge (x'\wedge y+ x  \wedge y').
\end{align*}  
Hence the image of $x\wedge y$ by the boundary morphism 
$H^{q+q'}(\Om^{\bul}(U))\lo H^{q+q'}(\Om^{\bul}(U))$ 
is $x'\wedge y+ x  \wedge y'$. 
By using the Godement resolution of 
$\wt{\Om}^{\bul}$ and $\Om^{\bul}$, we obtain (\ref{lemm:lbg}). 
\end{proof}

We can also prove (\ref{prop:bcgmm}) in the case after (\ref{eqn:glxd}): 
\begin{prop}\label{prop:bcmop} 
Let $(T',{\cal J}',\del')\lo (T,{\cal J},\del)$ be as in {\rm (\ref{prop:spbc})}.
Set $Y_{\bul \leq N,\os{\circ}{T}{}'_0}:=Y_{\bul \leq N}\times_{S}S_{T_0'}$. 
Let $p\col Y_{\bul \leq N,\os{\circ}{T}{}'_0}\lo Y_{\bul \leq N,\os{\circ}{T}_0}$ 
be the natural morphism over $\os{\circ}{T}{}'\lo \os{\circ}{T}$. 
Then 
\begin{align*} 
N_{S(T)^{\nat},{\rm zar}}\otimes^L_{{\cal O}_T}{\cal O}_{T'} 
\col & 
Rg_{Y_{\bul \leq N,\os{\circ}{T}_0}/S(T)^{\nat}*}(F^{\bul \leq N})
\otimes^L_{{\cal O}_T}{\cal O}_{T'} 
\tag{1.7.38.1}\label{ali:mibv}\\
& \lo Rg_{Y_{\bul \leq N,\os{\circ}{T}_0}/S(T)^{\nat}*}(F^{\bul \leq N})
\otimes^L_{{\cal O}_T}{\cal O}_{T'}  
\end{align*} 
is equal to 
\begin{equation*} 
N_{S(T')^{\nat}{\rm zar}} \col 
Rg_{Y_{\bul \leq N,\os{\circ}{T}{}'_0}/S(T')^{\nat}*}(p^*_{{\rm crys}}(F^{\bul \leq N}))
\lo 
Rg_{Y_{\bul \leq N,\os{\circ}{T}{}'_0}/S(T')^{\nat}*}(p^*_{{\rm crys}}(F^{\bul \leq N})).  
\end{equation*} 
\end{prop}
\begin{proof}
Because  
${\cal F}^{\bul \leq N,\bul}
\otimes_{{\cal O}_{{{\cal Q}}_{\bul \leq N,\bul}}}
{\Om}^{\bul}_{{{\cal Q}}_{\bul \leq N,\bul}/S(T)^{\nat}}$ 
consists of flat ${\cal O}_T$-modules, 
$(\ref{eqn:gsflxd})\otimes_{{\cal O}_T}{\cal O}_{T'}$ is exact. 
Then we have only to apply the direct image 
$Rg_{T'*}R\pi_{S(T')^{\nat}{\rm zar}*}$ 
to the resulting exact sequence and to use 
Kato's log base change theorem (\cite[(6.10)]{klog1}).  
\end{proof}

\par 
Now we consider the case where 
$Y_{\bul \leq N}/S$ 
is an $N$-truncated simplicial SNCL scheme 
$X_{\bul \leq N}/S$ and 
$F^{\bul \leq N}=
\eps^*_{X_{\bul \leq N,\os{\circ}{T}_0}/\os{\circ}{T}}(E^{\bul \leq N})$, where 
$E^{\bul \leq N}$ is the sheaf in \S\ref{sec:psc}.  
Assume that $X_{\bul \leq N,\os{\circ}{T}_0}$ has 
an affine $N$-truncated simplicial open covering of $X_{\bul \leq N,\os{\circ}{T}_0}$. 
Let $f_{\bul}\col X_{\bul \leq N,\bul,\os{\circ}{T}_0}\lo 
S_{\os{\circ}{T}_0}\os{\sus}{\lo} S(T)^{\nat}$ 
be the composite structural morphism.

\begin{theo-defi}\label{theo:fcp} 
Let ${\cal E}^{\bul \leq N,\bul}$ and ${\cal P}^{\rm ex}_{\bul \leq N,\bul}$ 
be as in {\rm \S\ref{sec:psc}}. Then the following hold$:$
\par 
$(1)$ The filtered complex 
$$R\pi_{{\rm zar}*}(({\cal E}^{\bul \leq N,\bul}
\otimes_{{\cal O}_{{\cal P}^{\rm ex}_{\bul \leq N,\bul}}}
{\Om}^{\bul}_{{\cal P}^{\rm ex}_{\bul \leq N,\bul}/\os{\circ}{T}},P))$$ 
is independent of the choice of an affine $N$-truncated simplicial open covering of 
$X_{\bul \leq N,\os{\circ}{T}_0}$ and 
an $(N,\infty)$-truncated bisimplicial immersion 
$X_{\bul \leq N,\bul,\os{\circ}{T}_0} \os{\sus}{\lo} \ol{\cal P}_{\bul \leq N,\bul}$ over 
$\ol{S(T)^{\nat}}$. 
Set 
\begin{align*} 
(\wt{R}u_{X_{\bul \leq N,\os{\circ}{T}_0}/\os{\circ}{T}*}
(\eps^*_{X_{\bul \leq N,\os{\circ}{T}_0}/\os{\circ}{T}}(E^{\bul \leq N})),P) 
:=R\pi_{{\rm zar}*}(({\cal E}^{\bul \leq N,\bul}
\otimes_{{\cal O}_{{\cal P}^{\rm ex}_{\bul \leq N,\bul}}}
{\Om}^{\bul}_{{\cal P}^{\rm ex}_{\bul \leq N,\bul}/\os{\circ}{T}},P)).
\end{align*}
We call this filtered complex the 
{\it modified $P$-filtered log crystalline complex} of 
$\eps^*_{X_{\bul \leq N,\os{\circ}{T}_0}/\os{\circ}{T}}(E^{\bul \leq N})$
($P$ is the abbreviation of Poincar\'{e}). 
When $E^{\bul \leq N}={\cal O}_{\os{\circ}{X}_{\bul \leq N,T_0}/\os{\circ}{T}}$, 
we call this filtered complex 
the {\it modified $P$-filtered log crystalline complex} of 
$X_{\bul \leq N,\os{\circ}{T}_0}/\os{\circ}{T}$. 
\par 
$(2)$ For each $0\leq m\leq N$, 
there exist the following canonical isomorphisms 
in ${\rm D}^+(f^{-1}_{m,T}({\cal O}_T)):$ 
\begin{align*}
P_0(\wt{R}u_{X_{m,\os{\circ}{T}_0}/\os{\circ}{T}*}
(\eps^*_{X_{m,\os{\circ}{T}_0}/\os{\circ}{T}}(E^m)))
\os{\sim}{\lo} &
{\rm MF}(a^{(0)}_{*}
Ru_{\os{\circ}{X}{}^{(0)}_{m,T_0}/\os{\circ}{T}*}
(E^m_{\os{\circ}{X}{}^{(0)}_{m,T_0}/\os{\circ}{T}}
\otimes_{\mab Z}\vp^{(0)}_{\rm crys}
(\os{\circ}{X}_{m,T_0}/\os{\circ}{T}))\tag{1.7.39.1}\label{ali:p0dept}\\
&\lo 
{\rm MF}(a^{(1)}_{*}
Ru_{\os{\circ}{X}{}^{(1)}_{m,T_0}/\os{\circ}{T}*}
(E^m_{\os{\circ}{X}{}^{(1)}_{m,T_0}/\os{\circ}{T}}
\otimes_{\mab Z}\vp^{(1)}_{\rm crys}
(\os{\circ}{X}_{m,T_0}/\os{\circ}{T}))  \\
&\lo \cdots \\
&\lo 
{\rm MF}(a^{(l)}_{*}
Ru_{\os{\circ}{X}{}^{(l)}_{m,T_0}/\os{\circ}{T}*}
(E^m_{\os{\circ}{X}{}^{(l)}_{m,T_0}/\os{\circ}{T}}
\otimes_{\mab Z}\vp^{(l)}_{\rm crys}(\os{\circ}{X}_{m,T_0}/\os{\circ}{T}))\\
& \lo \cdots
\cdots)\cdots ),  
\end{align*} 
\begin{align*}
&{\rm gr}_k^P(\wt{R}u_{X_{m,\os{\circ}{T}_0}/\os{\circ}{T}*}
(\eps^*_{X_{m,\os{\circ}{T}_0}/\os{\circ}{T}}(E^m))) 
\tag{1.7.39.2}\label{ali:eoppd}\\
&\os{\sim}{\lo}  
a^{(k-1)}_*
Ru_{\os{\circ}{X}{}^{(k-1)}_{m,T_0}/\os{\circ}{T}*}
(E^m_{\os{\circ}{X}{}^{(k-1)}_{m,T_0}/\os{\circ}{T}}
\otimes_{\mab Z}
\vp^{(k-1)}_{\rm crys}(\os{\circ}{X}_{m,T_0}/\os{\circ}{T}))[-k]
\quad (k\in {\mab Z}_{\geq 1}).
\end{align*}
\par 
$(3)$ {\rm {\bf (Contravariant functoriality)}} 
Let $g_{\bul \leq N}$ be as in {\rm (\ref{eqn:xdxduss})} 
satisfying the condition {\rm (\ref{cd:xygxy})}. 
Let $u\col (S(T)^{\nat},{\cal J},\del) \lo (S'(T')^{\nat},{\cal J}',\del')$ 
be as in the beginning of \S\ref{sec:fcuc}. 
Assume that we are given the morphism {\rm (\ref{ali:gnfe})}.  
Then $g_{\bul \leq N}$ induces the following pull-back morphism 
\begin{align*}
\wt{g}_{\bul \leq N}^*\col &
(\wt{R}u_{Y_{\bul \leq N,\os{\circ}{T}{}'_0}/\os{\circ}{T}{}'*}
(\eps^*_{Y_{\bul \leq N,\os{\circ}{T}{}'_0}/\os{\circ}{T}{}'}(F^{\bul \leq N})),P)
\tag{1.7.39.3}\label{ali:ruynt}\\
&\lo 
Rg_{\bul \leq N *}
(\wt{R}u_{X_{\bul \leq N,\os{\circ}{T}_0}/\os{\circ}{T}*}
(\eps^*_{X_{\bul \leq N,\os{\circ}{T}_0}/\os{\circ}{T}}(E^{\bul \leq N})),P). 
\end{align*}
This morphism satisfies the similar transitive relation 
to {\rm (\ref{ali:pdpp})} and the similar relation to {\rm (\ref{eqn:fzidd})}. 
\par 
(4) Let 
$g_{\bul \leq N}\col X_{\bul \leq N}\lo Y_{\bul \leq N}$ and 
$h_{\bul \leq N}\col X_{\bul \leq N}\lo Y_{\bul \leq N}$ be morphisms 
satisfying the conditions (\ref{cd:xygxy}), $(1.5.6.4)$ and $(1.5.6.5)$. 
Assume that $\os{\circ}{g}_{\bul \leq N}=\os{\circ}{h}_{\bul \leq N}$. 
Then 
${\rm gr}^P_k(\wt{g}_{\bul \leq N}^*)={\rm gr}^P_k(\wt{h}_{\bul \leq N}^*)$ $(k\in {\mab Z})$. 
\par 
(5) Let $g_{\bul \leq N}\col X_{\bul \leq N}\lo Y_{\bul \leq N}$ be as in (4). 
Let ${\rm MF}(\os{\circ}{X}_{m,T}/\os{\circ}{T},E^m)$ 
be the right hand side on (\ref{ali:p0dept}). 
Let $g^{(\star)*}_m\col 
{\rm MF}(\os{\circ}{Y}_{m,T'}/\os{\circ}{T}{}',E^m)
\lo Rg_{m*}({\rm MF}(\os{\circ}{X}_{m,T}/\os{\circ}{T},E^m))$ 
be the induced morphism by $\{g^{(l)*}_m\}_{l\in {\mab N}}$, 
where $\os{\circ}{g}{}^{(l)*}_m\col 
a^{(l)}_{*}
Ru_{\os{\circ}{Y}{}^{(l)}_{m,T_0}/\os{\circ}{T}*}
(E^m_{\os{\circ}{X}{}^{(l)}_{m,T_0}/\os{\circ}{T}}
\otimes_{\mab Z}\vp^{(l)}_{\rm crys}(\os{\circ}{X}_{m,T_0}/\os{\circ}{T}))
\lo 
a^{(l)}_{*}
Ru_{\os{\circ}{X}{}^{(l)}_{m,T_0}/\os{\circ}{T}*}
(E^m_{\os{\circ}{X}{}^{(l)}_{m,T_0}/\os{\circ}{T}}
\otimes_{\mab Z}\vp^{(l)}_{\rm crys}(\os{\circ}{X}_{m,T_0}/\os{\circ}{T}))$ 
is the pull-back of 
$\os{\circ}{g}{}^{(l)}_m\col 
\os{\circ}{X}{}^{(l)}_{m,T_0}\lo \os{\circ}{Y}{}^{(l)}_{m,T'_0}$. 
Then $\wt{g}{}^*_m\vert_{P_0}=g^{(\star)*}_m$. 
\end{theo-defi}
\begin{proof} 
(1): By (\ref{prop:hkt}) the filtered complex 
$R\pi_{{\rm zar}*}({\cal E}^{\bul \leq N,\bul}
\otimes_{{\cal O}_{{\cal P}^{\rm ex}_{\bul \leq N,\bul}}}
\Om^{\bul}_{{\cal P}^{\rm ex}_{\bul \leq N,\bul}/\os{\circ}{T}})$
is independent of the choice stated in (1). 
By the same proof as that of (\ref{theo:indcr}), we obtain (1) by (\ref{ali:eoppd}) 
(which we will prove in (2)) and 
the descending induction on the number of the filtration $P$.  
\par 
One can also prove (1) by using only (2) and the ascending induction on the number of 
the filtration $P$. 
\par 
(2): (2) follows from (\ref{ali:pdte}) and (\ref{eqn:eoppd}) as in the proof of (\ref{theo:indcr}). 
\par
(3): The proof of the existence of $\wt{g}{}^*_{\bul \leq N}$ 
is simpler than that of (\ref{theo:funas}). 
We leave the detail to the reader.  
\par 
(4): If $k\geq 1$, then, by a simpler proof than that of (\ref{prop:grloc}), 
we have the following formula  
\begin{align*} 
{\rm gr}_k^P(\wt{g}_{m}^*)=\bigoplus_{\ul{\lam}\in \Lam^{(k)}(\os{\circ}{g})}
\deg(u)^{k}b_{\phi(\ul{\lam})*}
\os{\circ}{g}{}^*_{\ul{\lam}}
\tag{1.7.39.4}\label{ali:kpgmt}
\end{align*}  
under the identification of the isomorphism (\ref{ali:eoppd}). 
We have only to prove that 
${\rm gr}_k^P(\wt{h}_{\bul \leq N}^*)={\rm gr}_k^P(\wt{g}_{\bul \leq N}^*)$ 
$(k\in {\mab N})$. 
If $k\geq 1$, then we have the equality   
${\rm gr}_k^P(\wt{h}_{m}^*)={\rm gr}_k^P(\wt{g}_{m}^*)$ $(0\leq \forall m \leq N)$ by (\ref{ali:kpgmt}).  
The morphism ${\rm gr}_0^P(\wt{g}_{m}^*)$ and ${\rm gr}_0^P(\wt{h}_{m}^*)$ 
is induced by the following morphism 
\begin{align*} 
\bigoplus_{\ul{\lam}\in \Lam^{(k)}(\os{\circ}{g})}&
b_{\phi(\ul{\lam})*}
\os{\circ}{g}{}^*_{\ul{\lam}}\col 
b^{(l)}_{*}
Ru_{\os{\circ}{Y}{}^{(l)}_{m,T_0}/\os{\circ}{T}*}
(E_{\os{\circ}{Y}{}^{(l)}_{m,T_0}/\os{\circ}{T}}
\otimes_{\mab Z}\vp^{(1)}_{\rm crys}
(\os{\circ}{Y}_{m,T_0}/\os{\circ}{T}))\tag{1.7.39.6}\label{ali:kplt}\\
&\lo 
Rg_{m*}a^{(l)}_{*}
Ru_{\os{\circ}{X}{}^{(l)}_{m,T_0}/\os{\circ}{T}*}
(E_{\os{\circ}{X}{}^{(l)}_{m,T_0}/\os{\circ}{T}}
\otimes_{\mab Z}\vp^{(l)}_{\rm crys}
(\os{\circ}{X}_{m,T_0}/\os{\circ}{T}))\quad (l\in {\mab N}). 
\end{align*} 
Hence ${\rm gr}_0^P(\wt{g}_{m}^*)={\rm gr}_0^P(\wt{h}_{m}^*)$. 
We complete the proof of (4). 
\par 
(5): (5) follows from the definition of $\wt{g}{}^*_{\bul \leq N}$. 
\end{proof} 

\begin{coro}\label{coro:pwt}
$(1)$ Set 
\begin{equation*} 
E_1^{-k,q+k} :=
\bigoplus_{\max \{-k,0\}\leq m\leq N}
R^{q-m}f_{*}({\rm gr}^P_{k+m}\wt{R}u_{X_{m,\os{\circ}{T}_0}/\os{\circ}{T}*}
(\eps^*_{X_{m,\os{\circ}{T}_0}/\os{\circ}{T}}(E^m))).
\end{equation*}
Then there exists the following spectral sequence 
\begin{align*} 
E_1^{-k,q+k}\Lo 
\wt{R}{}^qf_{X_{\bul \leq N,\os{\circ}{T}_0}/\os{\circ}{T}*}
(\eps^*_{X_{\bul \leq N,\os{\circ}{T}_0}/\os{\circ}{T}}(E^{\bul \leq N})). 
\tag{1.7.40.1}\label{ali:xnt}
\end{align*} 
Explicitly  
\begin{equation*}
E_1^{-k,q+k}=
\bigoplus_{m=0}^N
R^{q-2m-k}f_{\os{\circ}{X}{}^{(k+m-1)}_{m,T_0}/\os{\circ}{T}*}
(E_{\os{\circ}{X}{}^{(k+m-1)}_{m,T_0}/\os{\circ}{T}}
\otimes_{\mab Z}
\vp^{(k+m-1)}_{\rm crys}(\os{\circ}{X}_{m,T_0}/\os{\circ}{T}))(-k-m,u)
\end{equation*}
for $k>0$, 
\begin{align*}
E_1^{-k,q+k}&=
R^{q+k}f_{*}(P_0(\wt{R}u_{X_{-k,\os{\circ}{T}_0}/\os{\circ}{T}*}
(\eps^*_{X_{-k,\os{\circ}{T}_0}/\os{\circ}{T}}(E^{-k})))\\
&\oplus \bigoplus_{m=-k+1}^N
R^{q-2m-k}f_{\os{\circ}{X}{}^{(k+m-1)}_{m,T_0}/\os{\circ}{T}*}
(E_{\os{\circ}{X}{}^{(k+m-1)}_{m,T_0}/\os{\circ}{T}}
\otimes_{\mab Z}
\vp^{(k+m-1)}_{\rm crys}(\os{\circ}{X}_{m,T_0}/\os{\circ}{T}))(-k-m,u)
\end{align*} 
for $-N\leq k\leq 0$ 
and 
\begin{equation*}
E_1^{-k,q+k}=0
\end{equation*}
for $k<-N$. 
There also exists the following spectral sequence 
\begin{align*}
E_1^{i,q+k-i}&=R^{q+k-i}f_{\os{\circ}{X}{}^{(i)}_{-k,T_0}/\os{\circ}{T}*}
(E_{\os{\circ}{X}{}^{(i)}_{-k,T_0}/\os{\circ}{T}}
\otimes_{\mab Z}
\vp^{(i)}_{\rm crys}(\os{\circ}{X}_{m,T_0}/\os{\circ}{T}))\\
&\Lo 
R^{q+k}f_{*}(P_0(\wt{R}u_{X_{-k,\os{\circ}{T}_0}/\os{\circ}{T}*}
(\eps^*_{X_{-k,\os{\circ}{T}_0}/\os{\circ}{T}}(E^{-k}))).
\end{align*} 
These spectral sequences are contravariantly functorial for the morphism 
$g_{\bul \leq N}$ satisfying {\rm (1.5.6.4)} and {\rm (1.5.6.5)}. 
\par 
$(2)$ 
Let 
$g_{\bul \leq N}\col X_{\bul \leq N}\lo Y_{\bul \leq N}$ and 
$h_{\bul \leq N}\col X_{\bul \leq N}\lo Y_{\bul \leq N}$ be the morphisms 
in {\rm (\ref{theo:fcp}) (4)}. 
Then 
\begin{align*} 
\wt{g}{}^*_{\bul \leq N}=\wt{h}_{\bul \leq N}^*\col  
\wt{R}{}^qf_{Y_{\bul \leq N,\os{\circ}{T}{}'_0}/\os{\circ}{T}{}'*}
(\eps^*_{X_{\bul \leq N,\os{\circ}{T}_0}/\os{\circ}{T}}(F^{\bul \leq N}))
\lo 
u_*\wt{R}{}^qf_{X_{\bul \leq N,\os{\circ}{T}_0}/\os{\circ}{T}*}
(\eps^*_{X_{\bul \leq N,\os{\circ}{T}_0}/\os{\circ}{T}}(E^{\bul \leq N})). 
\end{align*} 
\end{coro} 
\begin{proof}
(1): This follows from (\ref{theo:fcp}). 
\par 
(2): (2) follows from the spectral sequence 
(\ref{ali:xnt}) and (\ref{theo:fcp}) (4). 
\end{proof}

\begin{rema}\label{rema:grph} 
(1) As in (\ref{rema:asid}) we can prove (\ref{theo:fcp}) (1)  by using 
(\ref{ali:p0dept}), (\ref{ali:eoppd}) and the ascending induction on the 
number of the filtration $P$. 
\par 
(2) Even when $E^{\bul \leq N}={\cal O}_{\os{\circ}{X}_{\bul \leq N,T_0}/\os{\circ}{T}}$, 
the graded complex of   
\begin{align*} 
\wt{R}u_{X_{m,\os{\circ}{T}_0}/\os{\circ}{T}*}
({\cal O}_{X_{m,T_0}/\os{\circ}{T}})
:=
\wt{R}u_{X_{m,\os{\circ}{T}_0}/\os{\circ}{T}*}
(\eps^*_{X_{m,\os{\circ}{T}_0}/\os{\circ}{T}}(
{\cal O}_{\os{\circ}{X}_{m,T_0}/\os{\circ}{T}})))
\end{align*} 
by $P$ is ``mixed'' by (\ref{ali:p0dept}). 
Hence we do not use the terminology 
``modified preweight-filtered log crystalline complex'' 
for $(\wt{R}u_{X_{m,\os{\circ}{T}_0}/\os{\circ}{T}*}
({\cal O}_{X_{m,T_0}/\os{\circ}{T}}),P)$. 
\end{rema}

\par 
Assume that $\os{\circ}{S}$ and $\os{\circ}{T}$ are $p$-adic formal schemes. 
Let the notations be as in the beginning of this section. 
Let $e$ be a positive integer. 
Then the following sequence 
\begin{align*} 
0 & \lo {\cal F}^{\bul \leq N,\bul}
\otimes_{{\cal O}_{{{\cal Q}}_{\bul \leq N,\bul}}}
{\Om}^{\bul}_{{{\cal Q}}_{\bul \leq N,\bul}/S(T)^{\nat}}
\otimes_{\mab Z}{\mab Q}[-1] 
\os{e\theta_{{\cal Q}_{{\bul\leq N,\bul}} \wedge}}{\lo}  
{\cal F}^{\bul \leq N,\bul}
\otimes_{{\cal O}_{{{\cal Q}}_{\bul \leq N,\bul}}}
{\Om}^{\bul}_{{{\cal Q}}_{\bul \leq N,\bul}/\os{\circ}{T}}
\otimes_{\mab Z}{\mab Q}
\tag{1.7.41.1}\label{eqn:gnsd}\\ 
& 
\lo {\cal F}^{\bul \leq N,\bul}
\otimes_{{\cal O}_{{{\cal Q}}_{\bul \leq N,\bul}}}
{\Om}^{\bul}_{{{\cal Q}}_{\bul \leq N,\bul}/S(T)^{\nat}}
\otimes_{\mab Z}{\mab Q} \lo 0
\end{align*} 
is exact. 
Let 
\begin{equation*} 
{\cal F}^{\bul \leq N,\bul}
\otimes_{{\cal O}_{{{\cal Q}}_{\bul \leq N,\bul}}}
{\Om}^{\bul}_{{{\cal Q}}_{\bul \leq N,\bul}/S(T)^{\nat}}
\otimes_{\mab Z}{\mab Q}
\lo 
{\cal F}^{\bul \leq N,\bul}
\otimes_{{\cal O}_{{{\cal Q}}_{\bul \leq N,\bul}}}
{\Om}^{\bul}_{{{\cal Q}}_{\bul \leq N,\bul}/S(T)^{\nat}}
\otimes_{\mab Z}{\mab Q} 
\tag{1.7.41.2}\label{eqn:gyssxd}
\end{equation*} 
be the boundary morphism of (\ref{eqn:gnsd}) 
in the derived category 
${\rm D}^+(f^{-1}_{\bul \leq N,\bul}({\cal O}_T))$. 
This morphism induces the following morphism by (\ref{ali:ptc}): 
\begin{equation*}
e^{-1}N_{S(T)^{\nat},{\rm zar}} \col 
Ru_{Y_{\bul \leq N,\os{\circ}{T}_0}/S(T)^{\nat}*}
(F^{\bul \leq N})\otimes_{\mab Z}{\mab Q} \lo 
Ru_{Y_{\bul \leq N,\os{\circ}{T}_0}/S(T)^{\nat}*}(F^{\bul \leq N})
\otimes_{\mab Z}{\mab Q} .
\tag{1.7.41.3}\label{eqn:nyne}
\end{equation*}

\section{Variational $p$-adic monodromy-weight conjecture I 
and variational filtered log $p$-adic hard Lefschetz conjecture I}\label{sec:vpmn} 
In this section we define the $p$-adic (zariskian) quasi-monodromy operator and 
we formulate the variational $p$-adic monodromy-weight conjecture ((\ref{conj:remc})) 
and the variational filtered log $p$-adic hard Lefschetz conjecture ((\ref{conj:lhilc})). 
The former is a generalization of the $p$-adic monodromy-weight conjecture in \cite{msemi}; 
the latter is a $p$-adic generalization of a conjecture in the $l$-adic case 
attributed to K.~Kato in \cite{nlpi}. 
We also define a new filtered complex 
$(B_{\rm zar},P)$ fitting into the following triangle of filtered complexes
\begin{align*} 
(A_{\rm zar},P)[-1]\lo (B_{\rm zar},P)\lo (A_{\rm zar},P)\os{+1}{\lo} 
\end{align*} 
and we prove the contravariant functoriality of 
$(B_{\rm zar},P)$. 
As a corollary of the functoriality, 
we prove the contravariant functoriality of 
the $p$-adic quasi-monodromy operator. 

\par 
Let $A_{\rm zar}({\cal P}^{\rm ex}_{\bul \leq N,\bul}/S(T)^{\nat}, 
{\cal E}^{\bul \leq N,\bul})^{\bul \bul}$ 
be the $(N,\infty)$-truncated bicosimplicial double complex of 
$f^{-1}_{\bul \leq N,\bul}({\cal O}_T)$-modules defined in 
(\ref{cd:adef}) and (\ref{cd:locstbd}). 
Let 
\begin{equation*} 
(A_{\rm zar}({\cal P}^{\rm ex}_{\bul \leq N,\bul}/S(T)^{\nat}, 
{\cal E}^{\bul \leq N,\bul})^{\bul \bul},P)
\lo (I^{\bul \leq N,\bul \bul \bul},P)
\tag{1.8.0.1}\label{eqn:apfi}
\end{equation*} 
be the $(N,\infty)$-truncated bicosimplicial 
filtered Godement resolution, which exist by \cite[(1.1.6)]{nh2}. 
(Note that the topos $(X_{\bul \leq N,\bul,\os{\circ}{T}_0})_{\rm zar}$ has enough points.) 
(A filtered injective resolution of 
$(A_{\rm zar}({\cal P}^{\rm ex}_{\bul \leq N,\bul}/S(T)^{\nat}, 
{\cal E}^{\bul \leq N,\bul})^{\bul \bul},P)$ 
also exists by the argument in \cite[p.~42]{nh2}; 
however we do not use the filtered injective resolution for 
the calculations of ${\rm Ker}(\wt{\nu}_{S(T)^{\nat},{\rm zar}})$, 
${\rm Coker}(\wt{\nu}_{S(T)^{\nat},{\rm zar}})$ and 
${\rm Im}(\wt{\nu}_{S(T)^{\nat},{\rm zar}})$ in a future article, 
where $\wt{\nu}_{S(T)^{\nat},{\rm zar}}$ is the morphism defined in 
(\ref{ali:ntzd}) below.) 
Then we have a resolution 
$(s(A_{\rm zar}({\cal P}^{\rm ex}_{\bul \leq N,\bul}/S(T)^{\nat}, 
{\cal E}^{\bul \leq N,\bul})^{\bul \bul}),P)
\lo  
(s(I^{\bul \leq N,\bul \bul \bul}),P)$,   
where $s$ means the single complex 
of the double complex with respect to the last two degrees. 
Set $(A^{\bul \leq N,\bul \bul},P):=
\pi_{{\rm zar}*}((I^{\bul \leq N,\bul \bul \bul},P))$. 
Then 
$$s((A^{\bul \leq N,\bul \bul},P))
=(A_{\rm zar}(X_{\bul \leq N,\os{\circ}{T}_0}/S(T)^{\nat},E^{\bul \leq N}),P)$$ 
in ${\rm D}^+{\rm F}(f^{-1}_{\bul \leq N}({\cal O}_T))$.  
Consider the following morphism 
\begin{align*} 
&\nu_{S(T)^{\nat},{\rm zar}}({\cal P}^{\rm ex}_{\bul \leq N,\bul}/S(T)^{\nat},
{\cal E}^{\bul \leq N,\bul})^{ij}:={\rm proj} \tag{1.8.0.2}\label{eqn:reslbb} \\ 
& \col 
{\cal E}^{\bul \leq N,\bul}\otimes_{{\cal O}_{{\cal P}^{\rm ex}_{\bul \leq N,\bul}}}
{\Om}^{i+j+1}_{{\cal P}^{\rm ex}_{\bul \leq N,\bul}/\os{\circ}{S}}/P_{j+1}  
\lo {\cal E}^{\bul \leq N,\bul}\otimes_{{\cal O}_{{\cal P}^{\rm ex}_{\bul \leq N,\bul}}}
{\Om}^{i+j+1}_{{\cal P}^{\rm ex}_{\bul \leq N,\bul}/\os{\circ}{S}}/P_{j+2}.   
\end{align*}  
It is easy to check that 
the morphism above actually induces  
a morphism of complexes 
with the boundary morphisms in (\ref{cd:locstbd}).  
Since $(I^{\bul \leq N,\bul \bul \bul},P)$ is 
the filteredly Godement resolution of 
$(s(A_{\rm zar}({\cal P}^{\rm ex}_{\bul \leq N,\bul}/S(T)^{\nat}, 
{\cal E}^{\bul \leq N,\bul})^{\bul \bul}),P)$, 
the morphism (\ref{eqn:reslbb}) 
induces the morphism  
\begin{equation*} 
\wt{\nu}^{ij}_{S(T)^{\nat},{\rm zar}}:={\rm proj}. \col 
A^{\bul \leq N,ij}\lo A^{\bul \leq N,i-1,j+1}
\end{equation*}  
of  sheaves of $f^{-1}_{\bul \leq N}({\cal O}_T)$-modules. 
Set 
\begin{align*} 
\wt{\nu}_{S(T)^{\nat},{\rm zar}}:=
s(\oplus_{i,j\in {\mab N}}\wt{\nu}^{ij}_{S(T)^{\nat},{\rm zar}}).
\tag{1.8.0.3}\label{ali:ntzd}
\end{align*}  
Let 
\begin{equation*} 
\nu_{S(T)^{\nat},{\rm zar}}
\col (A_{\rm zar}(X_{\bul \leq N,\os{\circ}{T}_0}/S(T)^{\nat},E^{\bul \leq N}),P)
\lo 
(A_{\rm zar}(X_{\bul \leq N,\os{\circ}{T}_0}/S(T)^{\nat},E^{\bul \leq N}),P\langle -2\rangle) 
\tag{1.8.0.4}\label{eqn:naxgdxdn}
\end{equation*} 
be a morphism of complexes induced by 
$\{\nu^{ij}_{S(T)^{\nat},{\rm zar}}\}_{i,j \in {\mab N}}$. 
The morphism 
$$\theta_{{\cal P}^{\rm ex}_{\bul \leq N,\bul}} \wedge \col 
(A_{\rm zar}({\cal P}^{\rm ex}_{\bul \leq N,\bul}/S(T)^{\nat}, 
{\cal E}^{\bul \leq N,\bul})^{ij},P)
\lo 
(A_{\rm zar}({\cal P}^{\rm ex}_{\bul \leq N,\bul}/S(T)^{\nat}, 
{\cal E}^{\bul \leq N,\bul})^{i,j+1},P)$$  
in (\ref{cd:locstbd}) induces a morphism 
$$\theta_{{\cal P}^{\rm ex}_{\bul \leq N,\bul}} \wedge 
\col A^{\bul \leq N,\bul,ij}\lo A^{\bul \leq N,\bul,i,j+1}.$$ 
\par 
Let $g_{\bul \leq N}$ be the morphism 
(\ref{eqn:xdxduss}) for the case 
$Y_{\bul \leq N,\os{\circ}{T}{}'_0}= X_{\bul \leq N,\os{\circ}{T}_0}$ 
and $(T',{\cal J}',\del')=(T,{\cal J},\del)$ and $S'=S$ 
satisfying the condition (\ref{cd:xygxy}): 
\begin{equation*} 
\begin{CD} 
X'_{\bul \leq N,\os{\circ}{T}_0} @>{g'_{\bul \leq N}}>> X''_{\bul \leq N,\os{\circ}{T}_0} \\
@VVV @VVV \\ 
X_{\bul \leq N,\os{\circ}{T}_0} @>{g_{\bul \leq N}}>> X_{\bul \leq N,\os{\circ}{T}_0} \\
@VVV @VVV \\ 
S_{\os{\circ}{T}_0} @>>> S_{\os{\circ}{T}{}_0} \\ 
@V{\bigcap}VV @VV{\bigcap}V \\ 
S(T)^{\nat} @>{u}>> S(T)^{\nat}, 
\end{CD}
\tag{1.8.0.5}\label{cd:xgspxy}
\end{equation*}
where $X''_{\bul \leq N,\os{\circ}{T}{}'_0}$ 
is another disjoint union 
of the member of an affine $N$-truncated simplicial open covering of 
$X_{\bul \leq N,\os{\circ}{T}_0}$. 
Assume that $\deg(u)$ is not divisible by $p$ or 
that  $\os{\circ}{T}$ is a $p$-adic formal scheme and 
that the morphism (\ref{eqn:odnl}) is divisible by $p^{e_p(j+1)}$.  
Since the following diagram  
\begin{equation*}
\begin{CD}
A^{\bul \leq N,\bul,ij}
@>{{\rm proj}.}>> A^{\bul \leq N,\bul,i-1,j+1}\\ 
@V{(g^{{\rm PD}*}_{\bul \leq N,\bul})^{(ij)}}VV 
@VV{\deg(u)(g^{{\rm PD}*}_{\bul \leq N,\bul})^{(i-1,j+1)}}V  \\
g^{{\rm PD}*}_{\bul \leq N,\bul*}(A^{\bul \leq N,\bul,ij}) @>{{\rm proj}.}>>
g^{{\rm PD}*}_{\bul \leq N,\bul*}(A^{\bul \leq N,\bul,i-1,j+1}) 
\end{CD}
\tag{1.8.0.6}\label{cd:phipnu}
\end{equation*}
is commutative, the morphism (\ref{eqn:naxgdxdn}) is 
the following morphism 
\begin{equation*} 
\nu_{S(T)^{\nat},{\rm zar}} \col 
(A_{\rm zar}(X_{\bul \leq N,\os{\circ}{T}_0}/S(T)^{\nat}
,E^{\bul \leq N}),P)
\lo 
(A_{\rm zar}(X_{\bul \leq N,\os{\circ}{T}_0}/S(T)^{\nat}
,E^{\bul \leq N}),P\langle -2\rangle)(-1,u). 
\tag{1.8.0.7}\label{eqn:axdxdn}
\end{equation*}

\par  
Set 
\begin{equation*} 
B^{\bul \leq N,\bul ij}
:=A^{\bul \leq N,\bul i-1,j}(-1,u)\oplus A^{\bul \leq N,\bul ij} 
\quad ( i,j\in {\mab N})
\tag{1.8.0.8}\label{eqn:axmdn}
\end{equation*} 
and 
\begin{equation*} 
B^{\bul \leq N,ij}
:=A^{\bul \leq N,i-1,j}(-1,u)\oplus A^{\bul \leq N,ij} 
\quad ( i,j\in {\mab N})
\tag{1.8.0.9}\label{eqn:anxdn}
\end{equation*} 
(cf.~\cite[p.~246]{st1}).
The horizontal boundary morphism
$d' \col B^{\bul \leq N,ij} \lo B^{\bul \leq N,i+1,j}$ 
is, by definition, the induced morphism 
$d''{}^{\bul} \col B^{\bul \leq N,\bul ij} \lo B^{\bul \leq N,\bul, i+1,j}$ 
defined by the following formula: 
\begin{align*} 
d'(\om_1,\om_2)=(\nabla \om_1,-\nabla \om_2) 
\tag{1.8.0.10}\label{ali:sddclxn}
\end{align*} 
and the vertical one 
$d'' \col B^{\bul \leq N,ij} \lo B^{\bul \leq N,i,j+1}$ is the induced morphism 
by a morphism  $d''{}^{\bul} \col B^{\bul \leq N,\bul ij} \lo B^{\bul \leq N,\bul i,j+1}$ 
defined by the following formula: 
\begin{align*} 
d''{}^{\bul}(\om_1,\om_2)=
(-\theta_{{\cal P}^{\rm ex}_{\bul \leq N,\bul}}\wedge \om_1+
\nu_{S(T)^{\nat},{\rm zar}}(\om_2),
\theta_{{\cal P}^{\rm ex}_{\bul \leq N,\bul}}\wedge \om_2).
\tag{1.8.0.11}\label{ali:sddpxn}
\end{align*}  
It is easy to check that 
$B^{\bul \leq N,\bul \bul}$ is 
actually an $N$-truncated cosimplicial double complex. 
Let $B^{\bul \leq N,\bul}$ be the single complex of $B^{\bul \leq N,\bul \bul}$ 
with respect to the last two degrees. 
\par
Let 
\begin{align*} 
\mu_{{X_{\bul \leq N,\os{\circ}{T}_0}/\os{\circ}{T}}} \col 
\wt{R}u_{X_{\bul \leq N,\os{\circ}{T}_0}/\os{\circ}{T}*}
(\eps^*_{X_{\bul \leq N,\os{\circ}{T}_0}/\os{\circ}{T}}(E^{\bul \leq N}))
=&R\pi_{{\rm zar}*}({\cal E}^{\bul \leq N,\bul}
\otimes_{{\cal O}_{{\cal P}^{\rm ex}_{\bul \leq N,\bul}}}
{\Om}^{\bul}_{{\cal P}^{\rm ex}_{\bul \leq N,\bul}/\os{\circ}{T}})\tag{1.8.0.12}\label{ali:sgsclxn}\\
&\lo B^{\bul \leq N,\bul}
\end{align*} 
be a morphism
of complexes induced by the following morphisms
\begin{align*} 
\mu^i_{\bul \leq N,\bul} \col {\cal E}^{\bul \leq N,\bul}
\otimes_{{\cal O}_{{\cal P}^{\rm ex}_{\bul \leq N,\bul}}}
{\Om}^i_{{\cal P}^{\rm ex}_{\bul \leq N,\bul}/\os{\circ}{T}} 
\lo  & ~{\cal E}^{\bul \leq N,\bul}
\otimes_{{\cal O}_{{\cal P}^{\rm ex}_{\bul \leq N,\bul}}}
{\Om}^i_{{\cal P}^{\rm ex}_{\bul \leq N,\bul}/\os{\circ}{T}}/P_0 \oplus  \\
& {\cal E}^{\bul \leq N,\bul}
\otimes_{{\cal O}_{{\cal P}^{\rm ex}_{\bul \leq N,\bul}}}
{\Om}^{i+1}_{{\cal P}^{\rm ex}_{\bul \leq N,\bul}/\os{\circ}{T}}/P_0 
\quad (i\in {\mab N})
\end{align*} 
defined by the following formula 
\begin{align*} 
\mu^i_{\bul \leq N,\bul}(\om):=(\om~{\rm mod}~P_0,
\theta_{{\cal P}^{\rm ex}_{\bul \leq N,\bul}} 
\wedge \om~{\rm mod}~P_0) \quad 
(\om \in 
{\cal E}^{\bul \leq N,\bul}
\otimes_{{\cal O}_{{\cal P}^{\rm ex}_{\bul \leq N,\bul}}}
{\Om}^i_{{\cal P}^{\rm ex}_{\bul \leq N,\bul}/\os{\circ}{T}}).
\end{align*} 
First we consider the case of the trivial coefficient: 
$E^{\bul \leq N}={\cal O}_{\os{\circ}{X}_{\bul \leq N,T_0}/\os{\circ}{T}})$. 

Then we have the following 
morphism of triangles:
\begin{equation*}
\begin{CD}
@>>> 
A_{\rm zar}(X_{\bul \leq N,\os{\circ}{T}_0}/S(T)^{\nat}
,E^{\bul \leq N})[-1] @>{}>>\\
@. @A{(\theta_{X_{\bul \leq N,\os{\circ}{T}_0}/S(T)^{\nat}}\wedge *)[-1]}AA  \\
@>>> 
R\pi_{{\rm zar}*}({\cal E}^{\bul \leq N,\bul}
\otimes_{{\cal O}_{{\cal P}^{\rm ex}_{\bul \leq N,\bul}}}
{\Om}^{\bul}_{{\cal P}^{\rm ex}_{\bul \leq N,\bul}/S(T)^{\nat}})[-1]   
@>{R\pi_{{\rm zar}*}
(\theta_{{\cal P}^{\rm ex}_{\bul \leq N,\bul}} \wedge)}>> 
\end{CD}
\tag{1.8.0.13}\label{cd:sgsctexn}
\end{equation*} 
\begin{equation*}
\begin{CD}
B^{\bul \leq N,\bul} @>>> \\
@A{\mu_{{X_{\bul \leq N,\os{\circ}{T}_0}/\os{\circ}{T}}}}AA \\
R\pi_{{\rm zar}*}({\cal E}^{\bul \leq N,\bul}
\otimes_{{\cal O}_{{\cal P}^{\rm ex}_{\bul \leq N,\bul}}}
{\Om}^{\bul}_{{\cal P}^{\rm ex}_{\bul \leq N,\bul}/\os{\circ}{T}}) @>>> 
\end{CD}
\end{equation*} 
\begin{equation*}
\begin{CD}  
A_{\rm zar}(X_{\bul \leq N,\os{\circ}{T}_0}/S(T)^{\nat}
,E^{\bul \leq N})
@>{+1}>>  \\  
@A{\theta_{X_{\bul \leq N,\os{\circ}{T}_0}/S(T)^{\nat}} \wedge}AA \\ 
R\pi_{{\rm zar}*}({\cal E}^{\bul \leq N,\bul}
\otimes_{{\cal O}_{{\cal P}^{\rm ex}_{\bul \leq N,\bul}}}
{\Om}^{\bul}_{{\cal P}^{\rm ex}_{\bul \leq N,\bul}/S(T)^{\nat}}) 
@>{+1}>>.\\
\end{CD}
\end{equation*} 
This is nothing but the following diagram of triangles: 
\begin{equation*}
\begin{CD}
@>>> 
A_{\rm zar}(X_{\bul \leq N,\os{\circ}{T}_0}/S(T)^{\nat}
,E^{\bul \leq N})[-1] @>{}>>\\
@. @A{(\theta_{X_{\bul \leq N,\os{\circ}{T}_0}/S(T)^{\nat}}\wedge )[-1]}AA  \\
@>>> 
Ru_{X_{\bul \leq N,\os{\circ}{T}_0}/S(T)^{\nat}*}
(\eps^*_{X_{\bul \leq N,\os{\circ}{T}_0}/S(T)^{\nat}}(E^{\bul \leq N}))[-1]   
@>{}>> 
\end{CD}
\tag{1.8.0.14}\label{cd:monodcom}
\end{equation*} 
\begin{equation*}
\begin{CD}
B^{\bul \leq N,\bul} @>>> \\
@A{\mu_{{X_{\bul \leq N,\os{\circ}{T}_0}/\os{\circ}{T}}}}AA \\
\wt{R}u_{X_{\bul \leq N,\os{\circ}{T}_0}/\os{\circ}{T}*}
(\eps^*_{X_{\bul \leq N,\os{\circ}{T}_0}/\os{\circ}{T}}(E^{\bul \leq N})) @>>> 
\end{CD}
\end{equation*} 
\begin{equation*}
\begin{CD}  
A_{\rm zar}(X_{\bul \leq N,\os{\circ}{T}_0}/S(T)^{\nat},E^{\bul \leq N})
@>{+1}>>  \\  
@A{\theta_{X_{\bul \leq N,\os{\circ}{T}_0}/S(T)^{\nat}} \wedge}AA \\ 
Ru_{X_{\bul \leq N,\os{\circ}{T}_0}/S(T)^{\nat}*}
(\eps^*_{X_{\bul \leq N,\os{\circ}{T}_0}/S(T)^{\nat}}(E^{\bul \leq N}))
@>{+1}>>.\\
\end{CD}
\end{equation*}

\parno 
Hence we obtain the following as in \cite[p.~246]{st1} 
and \cite[(11.10)]{ndw}:

\begin{prop}\label{prop:cmzoqm}
The zariskian monodromy operator
\begin{align*} 
N_{S(T)^{\nat},{\rm zar}} \col &
Ru_{X_{\bul \leq N,\os{\circ}{T}_0}/S(T)^{\nat}*}
(\eps^*_{X_{\bul \leq N,\os{\circ}{T}_0}/S(T)^{\nat}}(E^{\bul \leq N}))
\tag{1.8.1.1}\label{eqn:minbv}\\
&\lo  
Ru_{X_{\bul \leq N,\os{\circ}{T}_0}/S(T)^{\nat}*}
(\eps^*_{X_{\bul \leq N,\os{\circ}{T}_0}/S(T)^{\nat}}(E^{\bul \leq N}))(-1,u)
\end{align*}  
is equal to  
$$\nu_{S(T)^{\nat},{\rm zar}} \col 
A_{\rm zar}(X_{\bul \leq N,\os{\circ}{T}_0}/S(T)^{\nat},E^{\bul \leq N}) 
\lo A_{\rm zar}(X_{\bul \leq N,\os{\circ}{T}_0}/S(T)^{\nat},E^{\bul \leq N})(-1,u)$$ 
via the isomorphism {\rm (\ref{eqn:uz})}. 
\end{prop}

\begin{coro}\label{coro:nncilp} 
Assume that $\os{\circ}{X}_{\bul \leq N}$ is quasi-compact. 
Then  the zariskian monodromy operator 
\begin{align*} 
N_{S(T)^{\nat},{\rm zar}} \col & 
Ru_{X_{\bul \leq N,\os{\circ}{T}_0}/S(T)^{\nat}*}
(\eps^*_{X_{\bul \leq N,\os{\circ}{T}_0}/S(T)^{\nat}}(E^{\bul \leq N}))
\tag{1.8.2.1}\label{ali:mdmycoh} \\
&\lo  
Ru_{X_{\bul \leq N,\os{\circ}{T}_0}/S(T)^{\nat}*}
(\eps^*_{X_{\bul \leq N,\os{\circ}{T}_0}/S(T)^{\nat}}(E^{\bul \leq N}))(-1,u)
\end{align*} 
is nilpotent.  
\end{coro}
\begin{proof}
(\ref{coro:nncilp}) immediately follows from 
(\ref{prop:cmzoqm}) since 
$\nu_{S(T)^{\nat},{\rm zar}}$ is nilpotent.
\end{proof}

\begin{coro}\label{prop:cnmorp} 
Let $(T',{\cal J}',\del')\lo (T,{\cal J},\del)$ 
be as in {\rm (\ref{prop:bcmop})}. 
Then 
\begin{align*} 
\nu_{S(T)^{\nat},{\rm zar}}\otimes^L_{{\cal O}_T}{\cal O}_{T'} 
\col & 
 A_{\rm zar}(X_{\bul \leq N,\os{\circ}{T}_0}/S(T)^{\nat}
,E^{\bul \leq N})\otimes^L_{{\cal O}_T}{\cal O}_{T'}  
\tag{1.8.3.1}\label{ali:mcibv}\\
&\lo A_{\rm zar}(X_{\bul \leq N,\os{\circ}{T}_0}/S(T)^{\nat}
,E^{\bul \leq N})(-1,u)
\otimes^L_{{\cal O}_T}{\cal O}_{T'} 
\end{align*} 
is equal to 
\begin{align*} 
\nu_{S'(T')^{\nat},{\rm zar}} \col &
 A_{\rm zar}(X_{\bul \leq N,\os{\circ}{T}{}'_0}/S'(T')^{\nat},\os{\circ}{f'}{}^{*}_{\! \!{\rm crys}}(E^{\bul \leq N})) \\
&\lo A_{\rm zar}(X_{\bul \leq N,\os{\circ}{T}{}'_0}/S'(T')^{\nat},\os{\circ}{f'}{}^{*}_{\! \!{\rm crys}}(E^{\bul \leq N}))(-1,u).  
\end{align*} 
\end{coro}
\begin{proof} 
(\ref{prop:cnmorp}) follows from  
(\ref{prop:cmzoqm}) and (\ref{prop:bcmop}). 
\end{proof} 

\begin{defi}
We call $\nu_{S(T)^{\nat},{\rm zar}}$ the {\it zariskian quasi-monodromy operator} of 
$\eps^*_{X_{\bul \leq N,\os{\circ}{T}_0}/\os{\circ}{T}}(E^{\bul \leq N})$. 
When $E^{\bul \leq N}={\cal O}_{\os{\circ}{X}_{\bul \leq N,T_0}/\os{\circ}{T}}$, 
then we call $\nu_{S(T)^{\nat},{\rm zar}}$ the {\it zariskian quasi-monodromy operator} of 
$X_{\bul \leq N,\os{\circ}{T}_0}/S(T)^{\nat}$.
\end{defi}

\begin{prop}\label{prop:quasimon}
Let $k$ be a positive integer and let $q$ be a nonnegative integer. 
Then 
\begin{align*} 
\nu_{S(T)^{\nat},{\rm zar}}^k \col  & 
R^qf_{*}
({\rm gr}_k^PA_{\rm zar}(X_{\bul \leq N,\os{\circ}{T}_0}/S(T)^{\nat},E^{\bul \leq N})) 
\tag{1.8.5.1}\label{eqn:kpa}\\
& \lo 
R^qf_{*}
({\rm gr}_{-k}^PA_{\rm zar}(X_{\bul \leq N,\os{\circ}{T}_0}/S(T)^{\nat}
,E^{\bul \leq N}))(-k,u)
\end{align*} 
is the identity. 
\end{prop}
\begin{proof}
The proof is the same as that of \cite[(11.7)]{ndw}. 
Here we give only the following easy calculation of 
${\rm gr}_{-k}^P
A_{\rm zar}(X_{\bul \leq N,\os{\circ}{T}_0}/S(T)^{\nat}
,E^{\bul \leq N})$. 
\par 
The complex 
$${\rm gr}_{-k}^PA_{\rm zar}(X_{\bul \leq N,\os{\circ}{T}_0}/S(T)^{\nat}
,E^{\bul \leq N})=
{\rm gr}_{-k}^P
R\pi_{{\rm zar}*}
(A_{\rm zar}({\cal P}^{\rm ex}_{\bul \leq N,\bul}/S(T)^{\nat},
{\cal E}^{\bul \leq N,\bul}))$$  
is equal to 
\begin{align*} 
&R\pi_{{\rm zar}*}({\rm gr}_{-k}^P
A_{\rm zar}({\cal P}^{\rm ex}_{\bul \leq N,\bul}/S(T)^{\nat},
{\cal E}^{\bul \leq N,\bul}))\\
& 
=R\pi_{{\rm zar}*}(\bigoplus_{j\geq 0}
{\rm gr}_{2j-k+1}^P
A_{\rm zar}({\cal P}^{\rm ex}_{\bul \leq N,\bul}/S(T)^{\nat},
{\cal E}^{\bul \leq N,\bul}))^{\bul j}\{-j\})
\end{align*} 
and this complex is equal to 
the following complex:
\begin{align*}
&
R\pi_{{\rm zar}*}[\cdots \os{-\nabla}{\lo} 
\us{(i-k,j+k)}{
{\rm gr}^P_{2(j+k)-k+1}
(A_{\rm zar}({\cal P}^{\rm ex}_{\bul \leq N,\bul}/S(T)^{\nat},
{\cal E}^{\bul \leq N,\bul}))^{i-k, j+k}\{-j-k\}}
\os{-\nabla}{\lo} \cdots ]
\tag{1.8.5.2}\label{eqn:hstjkp} \\ 
&=R\pi_{{\rm zar}*}[\cdots \os{-\nabla}{\lo} 
\us{(i,j)}{
{\rm gr}^P_{2j+k+1}
(A_{\rm zar}({\cal P}^{\rm ex}_{\bul \leq N,\bul}/S(T)^{\nat},
{\cal E}^{\bul \leq N,\bul}))^{ij}\{-j\}}
\os{-\nabla}{\lo} \cdots ] \\
{}& \simeq 
{\rm gr}_k^P
A_{\rm zar}(X_{\bul \leq N,\os{\circ}{T}_0}/S(T)^{\nat}
,E^{\bul \leq N}). 
\end{align*}
\end{proof}


\begin{prop}\label{prop:nucsa} 
The zariskian quasi-monodromy operator 
\begin{equation*} 
\nu_{S(T)^{\nat},{\rm zar}} \col 
A_{\rm zar}(X_{\bul \leq N,\os{\circ}{T}_0}/S(T)^{\nat}
,E^{\bul \leq N})  \lo 
A_{\rm zar}(X_{\bul \leq N,\os{\circ}{T}_0}/S(T)^{\nat}
,E^{\bul \leq N})(-1,u)
\end{equation*}
underlies the following morphism 
of filtered morphism 
\begin{equation*} 
\nu_{S(T)^{\nat},{\rm zar}} \col 
(A_{\rm zar}(X_{\bul \leq N,\os{\circ}{T}_0}/S(T)^{\nat}
,E^{\bul \leq N}),P)
\lo  
(A_{\rm zar}(X_{\bul \leq N,\os{\circ}{T}_0}/S(T)^{\nat}
,E^{\bul \leq N}),
P\langle -2 \rangle)(-1,u), 
\end{equation*}
where $P\langle -2 \rangle$ is a filtration 
defined by $(P\langle -2 \rangle)_k= P_{k-2}$. 
\end{prop}
\begin{proof} 
Obvious. 
\end{proof}

\begin{coro}\label{coro:zmop}
The zariskian monodromy operator {\rm (\ref{ali:mdmycoh})} 
induces the following morphism 
\begin{align*} 
N_{S(T)^{\nat},{\rm zar}} \col & P_kRf_{X_{\bul \leq N,\os{\circ}{T}_0}/S(T)^{\nat}*}
(\eps^*_{X_{\bul \leq N,\os{\circ}{T}_0}/S(T)^{\nat}}(E^{\bul \leq N}))
\lo  \tag{1.8.7.1}\label{eqn:grppd}\\
& P_{k-2}Rf_{X_{\bul \leq N,\os{\circ}{T}_0}/S(T)^{\nat}*}
(\eps^*_{X_{\bul \leq N,\os{\circ}{T}_0}/S(T)^{\nat}}(E^{\bul \leq N}))(-1,u)
\end{align*}
\end{coro}
\begin{proof} 
(\ref{coro:zmop}) immediately follows from 
(\ref{prop:cmzoqm}) and (\ref{prop:nucsa}).
\end{proof}

\par 

\begin{prop}\label{prop:spaqqqq}
Assume that we are given the commutative diagram 
{\rm (\ref{cd:xgspxy})} and a similar commutative diagram 
\begin{equation*} 
\begin{CD} 
Y'_{\bul \leq N,\os{\circ}{T}{}'_0} @>{g'_{Y_{\bul \leq N,\os{\circ}{T}{}'_0}}}>> 
Y''_{\bul \leq N,\os{\circ}{T}{}'_0} \\ 
@VVV @VVV \\ 
Y_{\bul \leq N,\os{\circ}{T}{}'_0} @>{g_{Y_{\bul \leq N,\os{\circ}{T}{}'_0}}}
>> Y_{\bul \leq N,\os{\circ}{T}{}'_0} \\
@VVV @VVV \\ 
S'_{\os{\circ}{T}{}'_0} @>>> S'_{\os{\circ}{T}{}'_0} \\ 
@V{\bigcap}VV @VV{\bigcap}V \\ 
S'(T')^{\nat} @>{u'}>> S'(T')^{\nat}, 
\end{CD}
\tag{1.8.8.1}\label{eqn:yhxuss}
\end{equation*}
to {\rm (\ref{cd:xgspxy})}. 
Let $h_{\bul \leq N}\col X_{\bul \leq N,\os{\circ}{T}_0}\lo Y_{\bul \leq N,\os{\circ}{T}{}'_0}$ 
be a morphism fitting into a morphism from the commutative 
{\rm (\ref{cd:xgspxy})} to {\rm (\ref{eqn:yhxuss})}. 
Let $F^{\bul \leq N}$ be a flat quasi-coherent crystal of 
${\cal O}_{\os{\circ}{Y}_{\bul \leq N,T'_0}/\os{\circ}{T}{}'}$-modules. 
Let 
$$\os{\circ}{h}{}^*_{\bul \leq N,{\rm crys}}(F^{\bul \leq N})\lo E^{\bul \leq N}$$ 
be a morphism of flat quasi-coherent crystals of 
${\cal O}_{Y_{\bul \leq N,\os{\circ}{T}}/\os{\circ}{T}}$-modules. 
Then the following hold$:$
\par 
$(1)$ Let $h_{\bul \leq N}$ be the morphism $g_{\bul \leq N}$
in {\rm (\ref{theo:funas})} or {\rm (\ref{theo:funpas})}. 
Then the following diagram is commutative$:$ 
\begin{equation*} 
\begin{CD} 
Rh_{{\bul \leq N}*}(A_{\rm zar}(X_{\bul \leq N,\os{\circ}{T}_0}/S(T)^{\nat}
,E^{\bul \leq N}),P)@>{\nu_{S(T)^{\nat},{\rm zar}} }>>   \\ 
@A{h^*_{\bul \leq N}}AA \\
(A_{\rm zar}(Y_{\bul \leq N,\os{\circ}{T}{}'_0}/S'(T')^{\nat}
,F^{\bul \leq N}),P)@>{\deg (u) \nu_{S'(T')^{\nat},{\rm zar}} }>>   
\end{CD} 
\tag{1.8.8.2}\label{eqn:aasf}
\end{equation*} 
\begin{equation*} 
\begin{CD} 
Rh_{{\bul \leq N}*}(A_{\rm zar}(X_{\bul \leq N,\os{\circ}{T}_0}/S(T)^{\nat}
,E^{\bul \leq N}),P\langle -2 \rangle)\\
@AA{h^*_{\bul \leq N}}A \\
(A_{\rm zar}(Y_{\bul \leq N,\os{\circ}{T}{}'_0}/S(T)^{\nat},F^{\bul \leq N}),
P\langle -2 \rangle). 
\end{CD} 
\end{equation*} 
\par 
$(2)$ Let the notations be as in {\rm (\ref{theo:itc})} with the replacement 
of $g_{\bul \leq N}$ with $h_{\bul \leq N}$.  
Then the following diagram is commutative$:$ 
\begin{equation*} 
\begin{CD} 
Rh_{{\bul \leq N}*}(A_{\rm zar}(X_{\bul \leq N,\os{\circ}{T}_0}/S(T)^{\nat}
,E^{\bul \leq N}),P)\otimes^L_{\mab Z}{\mab Q}
@>{\nu_{S(T)^{\nat},{\rm zar}} }>>   \\ 
@A{h^*_{\bul \leq N}}AA \\
(A_{\rm zar}(Y_{\bul \leq N,\os{\circ}{T}{}'_0}/S'(T')^{\nat}
,F^{\bul \leq N}),P)\otimes^L_{\mab Z}{\mab Q}
@>{\nu_{S'(T')^{\nat},{\rm zar}} }>>   
\end{CD} 
\tag{1.8.8.3}\label{eqn:aalsf}
\end{equation*} 
\begin{equation*} 
\begin{CD} 
Rh_{{\bul \leq N}*}(A_{\rm zar}(X_{\bul \leq N,\os{\circ}{T}_0}/S(T)^{\nat}
,E^{\bul \leq N}),
P\langle -2 \rangle)(-1,u)\otimes^L_{\mab Z}{\mab Q}\\
@AA{h^*_{\bul \leq N}}A \\
(A_{\rm zar}(Y_{\bul \leq N,\os{\circ}{T}{}'_0}/S(T)^{\nat},F^{\bul \leq N}),
P\langle -2 \rangle)(-1,u')\otimes^L_{\mab Z}{\mab Q}. 
\end{CD} 
\end{equation*} 
\end{prop}
\begin{proof} 
Obvious. 
\end{proof}

Set 
\begin{align*} 
P_kB^{\bul \leq N,ij}
:=P_kA^{\bul \leq N,i-1,j}(-1,u)\oplus P_kA^{\bul \leq N,ij} 
\quad ( i,j\in {\mab N}).
\tag{1.8.8.4}\label{ali:bpa}
\end{align*} 
Then we have a filtered complex 
$(B^{\bul \leq N,\bul},P)$ in ${\rm D}^+{\rm F}(f^{-1}({\cal O}_T))$. 
\begin{prop}\label{prop:bid}
$(1)$ There exists the following sequence of the triangles in 
${\rm D}^+{\rm F}(f^{-1}({\cal O}_T)):$ 
\begin{align*} 
&\lo (A_{\rm zar}(X_{\bul \leq N,\os{\circ}{T}_0}/S(T)^{\nat},E^{\bul \leq N}),P)[-1]
\lo (B^{\bul \leq N,\bul},P)\tag{1.8.9.1}\label{ali:trab} \\
&\lo (A_{\rm zar}(X_{\bul \leq N,\os{\circ}{T}_0}/S(T)^{\nat},E^{\bul \leq N}),P)\os{+1}{\lo}. 
\end{align*} 
\par 
$(2)$ The filtered complex $(B^{\bul \leq N,\bul},P)$ is independent of the choice 
of an affine $N$-truncated simplicial open covering of $X_{\bul \leq N,\os{\circ}{T}_0}$ 
and the choice of an $(N,\infty)$-truncated bisimplicial immersion 
$X_{\bul \leq N,\bul,\os{\circ}{T}_0} \os{\sus}{\lo} 
\ol{\cal P}_{\bul \leq N,\bul}$ over $\ol{S(T)^{\nat}}$.  
\par 
$(3)$ The morphism 
{\rm (\ref{ali:sgsclxn})} is independent of the choices in $(2)$. 
\end{prop} 
\begin{proof} 
(1): By the definition of $(B^{\bul \leq N,\bul},P)$, we obtain the triangle (\ref{ali:trab}).
\par 
(2): Let the notations be as in the proof of (\ref{prop:tefc}). 
We may assume that there exists the commutative diagram (\ref{cd:sstp}).  
Then we have a natural morphism from $(B^{\bul \leq N,\bul},P)$ obtained by 
$\ol{\cal P}{}'_{\bul \leq N,\bul}$ to $(B^{\bul \leq N,\bul},P)$ obtained by 
$\ol{\cal P}_{\bul \leq N,\bul}$. 
Because 
$$(A_{\rm zar}(X_{\bul \leq N,\os{\circ}{T}_0}/S(T)^{\nat},E^{\bul \leq N}),P)$$  
is independent of the choice stated in (2) by (\ref{theo:indcr}), 
so is $(B^{\bul \leq N,\bul},P)$. 
\par 
(3): (3) follows from (2) and the commutative diagram (\ref{cd:monodcom}). 
\end{proof} 

\begin{defi} 
We call $(B^{\bul \leq N,\bul},P)$ 
the {\it extended zariskian $p$-adic filtered Steenbrink complexes} of 
$E^{\bul \leq N}$ for $X_{\bul \leq N,\os{\circ}{T}_0}/(S(T)^{\nat},{\cal J},\del)$. 
We denote it  by
$$(B_{\rm zar}(X_{\bul \leq N,\os{\circ}{T}_0}/S(T)^{\nat},E^{\bul \leq N}),P)
\in {\rm D}^+{\rm F}(
f^{-1}({\cal O}_T)).$$ 
When $E^{\bul \leq N}={\cal O}_{\os{\circ}{X}_{\bul \leq N,T_0}/\os{\circ}{T}}$, we denote 
$(B_{\rm zar}(X_{\bul \leq N,\os{\circ}{T}_0}/S(T)^{\nat},E^{\bul \leq N}),P)$ 
by 
$$(B_{\rm zar}(X_{\bul \leq N,\os{\circ}{T}_0}/S(T)^{\nat}),P).$$
We call $(B_{\rm zar}(X_{\bul \leq N,\os{\circ}{T}_0}/S(T)^{\nat}),P)$ the 
{\it extended zariskian $p$-adic filtered Steenbrink complex} of 
$X_{\bul \leq N,\os{\circ}{T}_0}/(S(T)^{\nat},{\cal J},\del)$. 
\end{defi} 

\begin{theo}[{\bf Contravariant functoriality I of $B_{\rm zar}$}]\label{theo:funcb}
$(1)$ Let the notations be as in {\rm (\ref{theo:funas})} or {\rm (\ref{theo:funpas})}.  
Then $g_{\bul \leq N}\col X_{\bul \leq N,\os{\circ}{T}{}_0}\lo 
Y_{\bul \leq N,\os{\circ}{T}{}'_0}$ 
induces the following well-defined pull-back morphism 
\begin{equation*}  
g_{\bul \leq N}^* \col 
(B_{\rm zar}(Y_{\bul \leq N,\os{\circ}{T}{}'_0}/S'(T')^{\nat},F^{\bul \leq N}),P)
\lo Rg_{\bul \leq N*}((B_{\rm zar}(X_{\bul \leq N,\os{\circ}{T}_0}/S(T)^{\nat}
,E^{\bul \leq N}),P)) 
\tag{1.8.11.1}\label{eqn:fzbaxd}
\end{equation*} 
fitting into the following commutative diagram$:$
\begin{equation*} 
\begin{CD}
B_{\rm zar}(Y_{\bul \leq N,\os{\circ}{T}{}'_0}/S'(T')^{\nat},F^{\bul \leq N})
@>{g_{\bul \leq N}^*}>>  
Rg_{\bul \leq N*}(B_{\rm zar}(X_{\bul \leq N,\os{\circ}{T}_0}/S(T)^{\nat}
,E^{\bul \leq N}))\\ 
@A{\mu_{Y_{\bul \leq N,\os{\circ}{T}{}'_0}/S'(T')^{\nat}}\wedge}A{\simeq}A 
@A{\simeq}A{Rg_{\bul \leq N*}(\mu_{X_{\bul \leq N,\os{\circ}{T}_0}/S(T)^{\nat}}\wedge)}A\\
\wt{R}u_{Y_{\bul \leq N,\os{\circ}{T}{}'_0}/\os{\circ}{T}{}'*}
(\eps^*_{Y_{\bul \leq N,\os{\circ}{T}{}'_0}/\os{\circ}{T}{}'}(F^{\bul \leq N}))
@>{g_{\bul \leq N}^*}>>Rg_{\bul \leq N*}\wt{R}u_{X_{\bul \leq N,\os{\circ}{T}_0}/\os{\circ}{T}*}
(\eps^*_{X_{\bul \leq N,\os{\circ}{T}_0}/\os{\circ}{T}}(E^{\bul \leq N})).
\end{CD}
\tag{1.8.11.2}\label{cd:psbccz} 
\end{equation*}
\par 
$(2)$ Let $h_{\bul \leq N}\col Y_{\bul \leq N,\os{\circ}{T}{}'_0}\lo 
Z_{\bul \leq N,\os{\circ}{T}{}''_0}$ be a morphism over 
$S''(T'')^{\nat}\lo S'(T')^{\nat}$ which is analogous to the morphism $g_{\bul \leq N}$. 
Then 
\begin{align*} 
(h_{\bul \leq N}\circ g_{\bul \leq N})^*  =&
Rh_{\bul \leq N*}(g_{\bul \leq N}^*)\circ h_{\bul \leq N}^*    
\col (B_{\rm zar}(Z_{\bul \leq N,\os{\circ}{T}{}''_0}/S''(T'')^{\nat},F^{\bul \leq N}),P)
\tag{1.8.11.3}\label{ali:pdbpp} \\ 
& \lo Rh_{\bul \leq N*}Rg_{\bul \leq N*}
(B_{\rm zar}(X_{\bul \leq N,\os{\circ}{T}_0}/S(T)^{\nat},E^{\bul \leq N}),P) \\
& =R(h_{\bul \leq N}\circ g_{\bul \leq N})_*
(B_{\rm zar}(X_{\bul \leq N,\os{\circ}{T}_0}/S(T)^{\nat},E^{\bul \leq N}),P).
\end{align*}  
\par 
$(3)$  
\begin{equation*} 
{\rm id}_{X_{\bul \leq N,\os{\circ}{T}_0}}^*={\rm id} 
\col (B_{\rm zar}(X_{\bul \leq N,\os{\circ}{T}_0}/S(T)^{\nat},E^{\bul \leq N}),P)
\lo (B_{\rm zar}(X_{\bul \leq N,\os{\circ}{T}_0}/S(T)^{\nat},E^{\bul \leq N}),P).  
\tag{1.8.11.4}\label{eqn:fzibdd}
\end{equation*} 
\end{theo} 
\begin{proof} 
Let the notations be as in the proof of (\ref{theo:funas}). 
Set 
\begin{align*} 
B_{\rm zar}({\cal P}^{\rm ex}_{\bul \leq N,\bul}/S(T)^{\nat},
{\cal E}^{\bul \leq N,\bul})^{ij}&:=
A_{\rm zar}({\cal P}^{\rm ex}_{\bul \leq N,\bul}/S(T)^{\nat},
{\cal E}^{\bul \leq N,\bul})^{i-1,j}\\
&\oplus 
A_{\rm zar}({\cal P}^{\rm ex}_{\bul \leq N,\bul}/S(T)^{\nat},
{\cal E}^{\bul \leq N,\bul})^{ij}
\end{align*} 
with boundary morphisms (\ref{ali:sddclxn}) and (\ref{ali:sddpxn}). 
Then we have a double complex 
$B_{\rm zar}({\cal P}^{\rm ex}_{\bul \leq N,\bul}/S(T)^{\nat},
{\cal E}^{\bul \leq N,\bul})^{\bul \bul}$. 
Let $B_{\rm zar}({\cal Q}^{\rm ex}_{\bul \leq N,\bul}/S'(T')^{\nat},
{\cal F}^{\bul \leq N,\bul})^{\bul \bul}$ be the similar complex to 
$B_{\rm zar}({\cal P}^{\rm ex}_{\bul \leq N,\bul}/S(T)^{\nat},
{\cal E}^{\bul \leq N,\bul})^{\bul \bul}$. 
We can define the pull-back morphism 
\begin{align*} 
g^*_{\bul \leq N,\bul} \col &
B_{\rm zar}({\cal Q}^{\rm ex}_{\bul \leq N,\bul}/S'(T')^{\nat},{\cal F}^{\bul \leq N,\bul})
\lo g^{\rm PD}_{\bul \leq N,\bul *}
(B_{\rm zar}({\cal P}^{\rm ex}_{\bul \leq N,\bul}/S(T)^{\nat},
{\cal E}^{\bul \leq N,\bul}))
\tag{1.8.11.5}\label{eqn:axfbwt}
\end{align*} 
by the following formula 
\begin{align*} 
g^{*ij}_{\bul \leq N,\bul}
:= (\deg(u)^{-j}g^{{\rm PD}*}_{\bul \leq N,\bul},
\deg(u)^{-(j+1)}g^{{\rm PD}*}_{\bul \leq N,\bul}) 
&\col 
B_{\rm zar}({\cal Q}^{\rm ex}_{\bul \leq N,\bul}/S'(T')^{\nat},{\cal F}^{\bul \leq N,\bul})^{ij}
\tag{1.8.11.6}\label{eqn:axbwt}\\
& \lo 
g^{\rm PD}_{\bul \leq N,\bul *}
B_{\rm zar}({\cal P}^{\rm ex}_{\bul \leq N,\bul}/S(T)^{\nat},{\cal E}^{\bul \leq N,\bul})^{ij}. 
\end{align*} 
This morphism induces the following filtered morphism 
\begin{align*} 
g^{*ij}_{\bul \leq N,\bul}
\col &
(B_{\rm zar}({\cal Q}^{\rm ex}_{\bul \leq N,\bul}/S'(T')^{\nat},{\cal F}^{\bul \leq N,\bul})^{ij},P) 
\lo 
g^{\rm PD}_{\bul \leq N,\bul *}
(B_{\rm zar}({\cal P}^{\rm ex}_{\bul \leq N,\bul}/S(T)^{\nat},
{\cal E}^{\bul \leq N,\bul})^{ij},P).\tag{1.8.11.7}\label{eqn:axbfit}
\end{align*} 
As in the proof of (\ref{theo:funas}), we can check that 
$g^{*ij}_{\bul \leq N,\bul}$ is the following morphism of double complexes: 
\begin{align*} 
g^{*\bul \bul}_{\bul \leq N,\bul}\col 
(B_{\rm zar}({\cal Q}^{\rm ex}_{\bul \leq N,\bul}/S'(T')^{\nat},
{\cal F}^{\bul \leq N,\bul})^{\bul \bul},P) 
\lo 
g^{\rm PD}_{\bul \leq N,\bul *}
(B_{\rm zar}({\cal P}^{\rm ex}_{\bul \leq N,\bul}/S(T)^{\nat},
{\cal E}^{\bul \leq N,\bul})^{\bul \bul},P) 
\end{align*} 
Let 
\begin{align*} 
g^*_{\bul \leq N,\bul}\col 
(B_{\rm zar}({\cal Q}^{\rm ex}_{\bul \leq N,\bul}/S'(T')^{\nat},{\cal F}^{\bul \leq N,\bul}),P) 
\lo 
g^{\rm PD}_{\bul \leq N,\bul *}
(B_{\rm zar}({\cal P}^{\rm ex}_{\bul \leq N,\bul}/S(T)^{\nat},{\cal E}^{\bul \leq N,\bul}),P) 
\end{align*} 
be the morphism of single complexes induced by $g^{*\bul \bul}_{\bul \leq N,\bul}$ 
and 
let 
\begin{equation*}
g^*_{\bul \leq N}\col 
(B_{\rm zar}(Y_{\bul \leq N,\os{\circ}{T}{}'_0}/S'(T')^{\nat},F^{\bul \leq N}),P)
\lo Rg_{\bul \leq N *}
(B_{\rm zar}(X_{\bul \leq N,\os{\circ}{T}_0}/S(T)^{\nat},E^{\bul \leq N}),P) 
\tag{1.8.11.8}\label{eqn:axnbwt}
\end{equation*} 
be the induced morphism by $g^*_{\bul \leq N,\bul}$.  
We can show the well-definedness of $g^*_{\bul \leq N}$ in (\ref{eqn:axnbwt}) as 
in the proof of (\ref{theo:funas}). 
(We leave the detail of the proof of the well-definedness to the reader.)
\end{proof}

\begin{theo}[{\bf Contravariant functoriality II of $B_{\rm zar}$}]\label{theo:funiicb}
Let the notations be as in 
{\rm (\ref{theo:funas})} or {\rm (\ref{theo:funpas})}.
Then there exists a morphism 
\begin{align*} 
g^*_{\bul \leq N}  &\col 
(B_{\rm zar}(Y_{\bul \leq N,\os{\circ}{T}{}'_0}/S'(T')^{\nat},
F^{\bul \leq N}),P)\otimes^L_{\mab Z}{\mab Q}
\tag{1.8.12.1}\label{eqn:aubbqd}\\
&\lo
Rg_{{\bul \leq N}*}
((B_{\rm zar}(X_{\bul \leq N,\os{\circ}{T}_0}/S(T)^{\nat},
E^{\bul \leq N}),P)\otimes^L_{\mab Z}{\mab Q})
\end{align*}
of filtered complexes in ${\rm D}^+{\rm F}(
f^{-1}({\cal K}_{T'}))$ fitting into 
the following commutative diagram 
\begin{equation*} 
\begin{CD}
B_{\rm zar}(Y_{\bul \leq N,\os{\circ}{T}{}'_0}/S'(T')^{\nat},F^{\bul \leq N})\otimes^L_{\mab Z}{\mab Q}
@>{g_{\bul \leq N}^*}>>  \\
@A{\mu_{Y_{\bul \leq N,\os{\circ}{T}{}'_0}/S'(T')^{\nat}}\wedge}A{\simeq}A\\ 
\wt{R}u_{Y_{\bul \leq N,\os{\circ}{T}{}'_0}/\os{\circ}{T}{}'*}
(\eps^*_{Y_{\bul \leq N,\os{\circ}{T}{}'_0}/\os{\circ}{T}{}'}(F^{\bul \leq N}))
\otimes^L_{\mab Z}{\mab Q}
@>{g_{\bul \leq N}^*}>>
\end{CD}
\tag{1.8.12.2}\label{cd:psbcbz} 
\end{equation*}
\begin{equation*} 
\begin{CD}
Rg_{\bul \leq N*}(B_{\rm zar}(X_{\bul \leq N,\os{\circ}{T}_0}/S(T)^{\nat}
,E^{\bul \leq N})\otimes^L_{\mab Z}{\mab Q})\\  
@A{Rg_{\bul \leq N*}(\mu_{X_{\bul \leq N,\os{\circ}{T}_0}/S(T)^{\nat}}\wedge)}A{\simeq}A\\
Rg_{\bul \leq N*}\wt{R}u_{X_{\bul \leq N,\os{\circ}{T}_0}/\os{\circ}{T}*}
(\eps^*_{X_{\bul \leq N,\os{\circ}{T}_0}/\os{\circ}{T}}(E^{\bul \leq N}))
\otimes^L_{\mab Z}{\mab Q}.
\end{CD} 
\end{equation*}
This morphism satisfies the similar relation to {\rm (\ref{ali:pdbpp})}.  
\end{theo} 
\begin{proof} 
The proof is the same as that of (\ref{theo:funcb}). 
\end{proof}

\begin{coro}\label{coro:spfs}
$(1)$ Let the notations be as in {\rm (\ref{theo:funas})} or {\rm (\ref{theo:funpas})}. 
The isomorphisms {\rm (\ref{eqn:uz})} 
for $X_{\bul \leq N,\os{\circ}{T}_0}/S(T)^{\nat}$, $E^{\bul \leq N}$ 
and $Y_{\bul \leq N,\os{\circ}{T}{}'_0}/S'(T')^{\nat}$, $F^{\bul \leq N}$ 
induce a morphism $g_{\bul \leq N}^*$ from {\rm (\ref{eqn:ynmumss})} 
for the case $Z_{\bul \leq N,\os{\circ}{T}{}'_0}=Y_{\bul \leq N,\os{\circ}{T}{}'_0}$ and 
$Y_{\bul \leq N,\os{\circ}{T}_0}=X_{\bul \leq N,\os{\circ}{T}{}_0}$ 
to {\rm (\ref{eqn:aasf})} for $z=g_{\bul \leq N}$.
This morphism satisfies the transitive relation and  
${\rm id}_{X_{\bul \leq N,\os{\circ}{T}_0}/S(T)^{\nat}}^*={\rm id}$. 
\par 
$(2)$ Let the notations be as in {\rm (\ref{theo:itc})}. 
The isomorphisms {\rm (\ref{eqn:uz})} 
for $X_{\bul \leq N,\os{\circ}{T}_0}/S(T)^{\nat}$, $E^{\bul \leq N}$ 
and $Y_{\bul \leq N,\os{\circ}{T}{}'_0}/S'(T')^{\nat}$, $F^{\bul \leq N}$ 
induce a morphism from 
{\rm (\ref{eqn:ynmaumss})}$\otimes^L_{\mab Z}{\mab Q}$ 
to {\rm (\ref{eqn:aalsf})}. 
This morphism satisfies the transitive relation and  
${\rm id}_{X_{\bul \leq N,\os{\circ}{T}_0}/S(T)^{\nat}}^*={\rm id}$. 
\end{coro}

\par 
Now we consider the case $N=0$.  Set $X:=X_0$. 
Assume that $(\os{\circ}{T},{\cal J},\del)$ is a $p$-adic formal PD-scheme.

\begin{conj}[{\bf Variational $p$-adic monodromy-weight conjecture}]\label{conj:rcpmc}
Assume that the structural morphism $\os{\circ}{X}_{T_0} \lo \os{\circ}{T}_0$ is projective. 
Let $q$ be nonnegative integer. 
Then the induced morphism 
\begin{align*} 
\nu^k_{S(T)^{\nat},{\rm zar}} \col & 
{\rm gr}^P_{q+k}R^qf_{X_{\os{\circ}{T}_0}/S(T)^{\nat}*}
({\cal O}_{X_{\os{\circ}{T}_0}/S(T)^{\nat}})
\otimes_{\mab Z}{\mab Q} \tag{1.8.14.1}\label{eqn:grmpd} \\
& \lo 
{\rm gr}^P_{q-k}R^qf_{X_{\os{\circ}{T}_0}/S(T)^{\nat}*}
({\cal O}_{X_{\os{\circ}{T}_0}/S(T)^{\nat}})(-k,u)
\otimes_{\mab Z}{\mab Q} 
\end{align*}
by the quasi-monodromy operator 
is an isomorphism.  
\end{conj}

\begin{rema} 
(1) By (\ref{theo:csoncrdw}) and (\ref{prop:nnuc}) below 
(the comparison between the method of filtered crystalline complexes 
and that of filtered de Rham-Witt complexes), 
we shall see that the conjecture (\ref{conj:rcpmc}) 
is a generalization of Mokrane's 
$p$-adic monodromy-weight conjecture in \cite[(3.27)]{msemi}. 
\par 
(2) Though one can raise a problem when 
the induced morphism 
\begin{align*} 
&\nu^k_{S(T)^{\nat},{\rm zar}} \col 
{\rm gr}^P_{q+k}R^qf_{X_{\os{\circ}{T}_0}/S(T)^{\nat}*}
(\eps^*_{X_{\os{\circ}{T}_0}/S(T)^{\nat}}(E))
\otimes_{\mab Z}{\mab Q}  
\lo \tag{1.8.15.1}\label{eqn:grmppd}\\
&{\rm gr}^P_{q-k}R^qf_{X_{\os{\circ}{T}_0}/S(T)^{\nat}*}
(\eps^*_{X_{\os{\circ}{T}_0}/S(T)^{\nat}}(E))(-k,u)
\otimes_{\mab Z}{\mab Q} 
\end{align*}
by the quasi-monodromy operator 
is an isomorphism, 
I am not familiar with this problem.  
However, see \cite{pvd} and \cite{ccpi}.  
\par 
(3) 
The conjecture (\ref{conj:rcpmc}) has been solved 
in the case $\dim X\leq 2$ by \cite[5.3, \S6]{msemi}, 
\cite[(11.3)]{ndw}, \cite[(8.4)]{nlpi}, the log crystalline analogue of \cite{kaji} 
in the case where $\os{\circ}{T}$ is the formal spectrum of 
the Witt ring of a perfect field of characteristic $p>0$ 
(see (\ref{theo:22}) below for a more general case). 
Assume that $\os{\circ}{T}$ is a $p$-adic formal ${\cal V}$-scheme 
in the sense of \cite{of}, 
that is, a noetherian $p$-adic formal ${\cal V}$-scheme 
which is topologically of finite type over ${\cal V}$, 
where ${\cal V}$ is a complete discrete valuation ring of mixed characteristics 
with perfect residue field.  
By (\ref{rema:nlcfi}), (\ref{prop:convmon}) 
and (\ref{prop:conl}) below and by \cite[(3.17)]{of}, 
we can reduce (\ref{conj:rcpmc}) to the case where 
$S$ is the formal spectrum of the Witt ring of a finite field 
((\ref{theo:nrlp}) below).   
(\ref{conj:rcpmc}) in other interesting cases 
(e.~g.,~\cite{itup}) has been solved. 
In \S\ref{sec:filbo} below we shall give other theorems.  
\par 
(4) If we do not assume the projectivity in (\ref{conj:rcpmc}), 
then there are many examples (rigid analytic reductions of certain elliptic surfaces) 
such that the structural morphisms 
$\os{\circ}{X}_{T_0} \lo \os{\circ}{T}_0$'s are proper 
and such that the morphisms (\ref{eqn:grmpd})'s 
are not isomorphisms.  See \cite{nlpi} for these examples. 
\end{rema}

\par 
We also give a variational filtered log $p$-adic hard Lefschetz conjecture, 
which is a $p$-adic generalization of Kato's conjecture in the $l$-adic case 
in \cite{nlpi}.

\par 
Assume that the structural morphism 
$\os{\circ}{f} \col \os{\circ}{X}_{T_0} \lo \os{\circ}{T}_0$ 
is projective and the relative dimension of 
$\os{\circ}{f}$ is of pure dimension $d$. 
Let $L$ be a relatively ample line bundle on 
$\os{\circ}{X}_{T_0}/\os{\circ}{T}_0$. 
By the same argument as that of  
\cite[\S3]{boi} (see also \cite[p.~130]{nh2}), 
we obtain the log cohomology class 
$\lam=c_{1,{\rm crys},S(T)^{\nat}}(L)$ of $L$
in $R^2f_{X_{\os{\circ}{T}_0}/S(T)^{\nat}*}({\cal O}_{X_{\os{\circ}{T}_0}/S(T)^{\nat}})(1)$. 
\par 
Let us recall this construction. 
More generally, let us define the first log crystalline chern class as follows. 
\par 
Let $(T,{\cal J},\del)$ be a fine log PD-scheme and 
let $T_0$ be the exact closed log subscheme defined by ${\cal J}$. 
Let $Y_{\bul \leq N}$ $(N\in {\mab N})$ be 
a fine $N$-truncated simplicial log scheme over $T_0$. 
Let ${\cal M}_{Y_{\bul \leq N}/T}$ be a sheaf on ${\rm Crys}(Y_{\bul \leq N}/T)$ such that,  
for an object $T'_{\bul \leq N}$ of ${\rm Crys}(Y_{\bul \leq N}/T)$,  
$\Gam(T'_{\bul \leq N},M_{Y_{\bul \leq N}/T}):=\{\Gam(T'_m,M_{T'_m})\}_{m\leq N}$. 
Let $i_{\bul \leq N}\col (Y_{\bul \leq N})_{\rm zar}\lo 
(Y_{\bul \leq N}/T)_{\rm crys}$ be the canonical morphism 
defined by 
$$\Gam(U_{\bul \leq N},i^*(E^{\bul \leq N})):=
E^{\bul \leq N}((U_{\bul \leq N},U_{\bul \leq N},0)).$$ 
Then we have the following exact sequence 
\begin{align*}
0\lo 1+{\cal J}_{Y_{\bul \leq N}/T}\lo M_{Y_{\bul \leq N}/T}\lo 
i_{\bul \leq N*}(M_{Y_{\bul \leq N}})\lo 0. 
\end{align*} 
(This is the log and $N$-truncated simplicial version of 
the exact sequence of \cite[(3.1.2)]{boi}.)
Since $M_{Y_{\bul \leq N}}$ and $M_{Y_{\bul \leq N}/T}$ are sheaves of integral monoids, 
the following sequence is exact: 
\begin{align*}
0\lo 1+{\cal J}_{Y_{\bul \leq N}/T}\lo M^{\rm gp}_{Y_{\bul \leq N}/T}\lo 
i_{\bul \leq N*}(M^{\rm gp}_{Y_{\bul \leq N}})\lo 0. 
\tag{1.8.15.2}\label{ali:mygyp}
\end{align*} 
Hence we have the following composite morphism
\begin{align*} 
c_{1,{\rm crys}}\col M_{Y_{\bul \leq N}}^{\rm gp}\lo 
Ru_{Y_{\bul \leq N}/T*}(1+{\cal J}_{Y_{\bul \leq N}/T})[1]\os{\log}{\lo} 
Ru_{Y_{\bul \leq N}/T*}({\cal J}_{Y_{\bul \leq N}/T})[1]\lo 
Ru_{Y/T*}({\cal O}_{Y_{\bul \leq N}/T})[1]. 
\tag{1.8.15.3}\label{ali:mygp}
\end{align*} 

\begin{defi} 
We call the morphism (\ref{ali:mygp}) 
the {\it first log crystalline chern class map} of $Y_{\bul \leq N}/T$. 
\end{defi}

\par 
Assume that $Y_{\bul \leq N}$ has the disjoint union of the member of 
an affine truncated simplicial open covering of $Y_{\bul \leq N}$. 
Let $Y_{\bul \leq N,\bul}$ be the \v{C}ech diagram of this open covering and 
let $Y_{\bul \leq N,\bul}\os{\sus}{\lo} {\cal Y}_{\bul \leq N,\bul}$ 
be an immersion 
into a log smooth $(N,\infty)$-truncated bisimplicial log scheme over $T$. 
Let ${\mathfrak E}_{\bul \leq N,\bul}$ 
be the log PD-envelope of this immersion over $(T,{\cal J},\del)$. 
Let us consider the following complex 
\begin{align*} 
M^{\rm gp}_{Y_{\bul \leq N,\bul}/T}\os{d\log}{\lo}  
L({\cal O}_{{\mathfrak E}_{\bul \leq N,\bul}}
\otimes_{{\cal O}_{{\cal Y}_{\bul \leq N,\bul}}}
\Om^1_{{\cal Y}_{\bul \leq N,\bul}/T})
\os{L(d)}{\lo} L({\cal O}_{{\mathfrak E}_{\bul \leq N,\bul}}
\otimes_{{\cal O}_{{\cal Y}_{\bul \leq N,\bul}}}
\Om^1_{{\cal Y}_{\bul \leq N,\bul}/T})\lo 
\cdots.   
\end{align*}
Here $L$ means the log linearization functor for 
quasi-coherent ${\cal O}_{{\mathfrak E}_{\bul \leq N,\bul}}$-modules 
(\cite[(2.2.1.2)]{nh2}). 
We denote this complex by 
$L^{\log}({\cal O}_{{\mathfrak E}_{\bul \leq N,\bul}}
\otimes_{{\cal O}_{{\cal Y}_{\bul \leq N,\bul}}}
\Om^{\bul}_{{\cal Y}_{\bul \leq N,\bul}/T})$. 
Set ${\cal K}^{\times}:={\rm Ker}(L^{\log}({\cal O}_{{\mathfrak E}_{\bul \leq N,\bul}}
\otimes_{{\cal O}_{{\cal Y}_{\bul \leq N,\bul}}}
\Om^{\bul}_{{\cal Y}_{\bul \leq N,\bul}/T})\lo i_{\bul \leq N,\bul*}(M^{\rm gp}_{Y_{\bul \leq N,\bul}}))$. 
Then we have the following exact sequence 
\begin{align*} 
0\lo {\cal K}^{\times}\lo L^{\log}({\cal O}_{{\mathfrak E}_{\bul \leq N,\bul}}
\otimes_{{\cal O}_{{\cal Y}_{\bul \leq N,\bul}}}
\Om^{\bul}_{{\cal Y}_{\bul \leq N,\bul}/T})\lo i_{\bul \leq N,\bul*}
(M^{\rm gp}_{Y_{\bul \leq N,\bul}})\lo 0. 
\tag{1.8.16.1}\label{ali:cro}
\end{align*} 
On the other hand we have the following exact sequence 
\begin{align*}
0\lo 1+{\cal J}_{Y_{\bul \leq N,\bul}/T}\lo M^{\rm gp}_{Y_{\bul \leq N,\bul}/T}\lo 
i_{\bul \leq N,\bul*}(M^{\rm gp}_{Y_{\bul \leq N,\bul}})\lo 0. 
\tag{1.8.16.2}\label{ali:mygnyp}
\end{align*} 
These exact sequences fit into the following commutative diagram: 
\begin{equation*} 
\small{
\begin{CD} 
0@>>> 1+{\cal J}_{Y_{\bul \leq N,\bul}/T}@>>> M^{\rm gp}_{Y_{\bul \leq N,\bul}/T}
@>>> 
i_{\bul \leq N,\bul*}(M^{\rm gp}_{Y_{\bul \leq N,\bul}})@>>> 0\\
@. @VVV @V{\bigcap}VV @| @.\\
0@>>> {\cal K}^{\times}@>>> L^{\log}({\cal O}_{{\mathfrak E}_{\bul \leq N,\bul}}
\otimes_{{\cal O}_{{\cal Y}_{\bul \leq N,\bul}}}
\Om^{\bul}_{{\cal Y}_{\bul \leq N,\bul}/T})@>>> 
i_{\bul \leq N,\bul*}(M^{\rm gp}_{Y_{\bul \leq N,\bul}})
@>>> 0.
\end{CD}}
\end{equation*} 
Hence we have the following commutative diagram 
\begin{equation*} 
\begin{CD} 
Ru_{Y_{\bul \leq N,\bul}/T*}(M^{\rm gp}_{Y_{\bul \leq N,\bul}})
@>>> Ru_{Y_{\bul \leq N,\bul}/T*}(1+{\cal J}_{Y_{\bul \leq N,\bul}/T})[1]\\
@VVV @VVV \\
Ru_{Y_{\bul \leq N,\bul}/T*}(M^{\rm gp}_{Y_{\bul \leq N,\bul}})
@>>>Ru_{Y_{\bul \leq N,\bul}/T*}({\cal K}^{\times})[1]. 
\end{CD}
\end{equation*}
Set  ${\cal K}:={\rm Ker}(L({\cal O}_{{\mathfrak E}_{\bul \leq N,\bul}})\lo 
i_{\bul \leq N,\bul*}({\cal O}_{Y_{\bul \leq N,\bul}}))$ and 
$${\cal K}^{\bul}:=({\cal K}\os{L(d)}{\lo}  
L({\cal O}_{{\mathfrak E}_{\bul \leq N,\bul}}
\otimes_{{\cal O}_{{\cal Y}_{\bul \leq N,\bul}}}
\Om^1_{{\cal Y}_{\bul \leq N,\bul}/T})\lo \cdots).$$ 
Then we have the morphism $\log \col  {\cal K}^{\times}\lo {\cal K}^{\bul}$ 
fitting into the following commutative diagram 
\begin{equation*} 
\begin{CD} 
1+{\cal J}_{Y_{\bul \leq N,\bul}/T}
@>{\log}>> {\cal J}_{Y_{\bul \leq N,\bul}/T} @>>> {\cal O}_{Y_{\bul \leq N,\bul}/T}\\
@VVV @VVV @VVV\\ 
{\cal K}^{\times}@>{\log}>>{\cal K}^{\bul}@>>> 
L({\cal O}_{{\mathfrak E}_{\bul \leq N,\bul}}
\otimes_{{\cal O}_{{\cal Y}_{\bul \leq N,\bul}}}
\Om^{\bul}_{{\cal Y}_{\bul \leq N,\bul}/T}). 
\end{CD}
\end{equation*}
Hence we have the following commutative diagram 
\begin{equation*} 
\begin{CD} 
Ru_{Y_{\bul \leq N,\bul}/T*}(M^{\rm gp}_{Y_{\bul \leq N,\bul}})
@>>> Ru_{Y_{\bul \leq N,\bul}/T*}(1+{\cal J}_{Y_{\bul \leq N,\bul}/T})[1]
@>{\log}>>\\
@VVV @VVV \\
Ru_{Y_{\bul \leq N,\bul}/T*}(M^{\rm gp}_{Y_{\bul \leq N,\bul}})[1]
@>>>Ru_{Y_{\bul \leq N,\bul}/T*}({\cal K}^{\times})[1] @>{\log}>>
\end{CD}
\end{equation*}
\begin{equation*} 
\begin{CD} 
Ru_{Y_{\bul \leq N,\bul}/T*}({\cal J}_{Y_{\bul \leq N,\bul}/T})[1] @>>> 
Ru_{Y_{\bul \leq N,\bul}/T*}({\cal O}_{Y_{\bul \leq N,\bul}/T})[1]\\
@VVV @VVV \\
Ru_{Y_{\bul \leq N,\bul}/T*}({\cal K}^{\bul})[1]@>>> 
Ru_{Y_{\bul \leq N,\bul}/T*}(L({\cal O}_{{\mathfrak E}_{\bul \leq N,\bul}}
\otimes_{{\cal O}_{{\cal Y}_{\bul \leq N,\bul}}}
\Om^{\bul}_{{\cal Y}_{\bul \leq N,\bul}/T}))[1]. 
\end{CD}
\end{equation*}
Because 
$Ru_{Y_{\bul \leq N,\bul}/T*}(L({\cal O}_{{\mathfrak E}_{\bul \leq N,\bul}}
\otimes_{{\cal O}_{{\cal Y}_{\bul \leq N,\bul}}}
\Om^{\bul}_{{\cal Y}_{\bul \leq N,\bul}/T}))=
{\cal O}_{{\mathfrak E}_{\bul \leq N,\bul}}
\otimes_{{\cal O}_{{\cal Y}_{\bul \leq N,\bul}}}
\Om^{\bul}_{{\cal Y}_{\bul \leq N,\bul}/T}$, 
we obtain the following commutative diagram 
\begin{equation*} 
\begin{CD} 
M^{\rm gp}_{Y_{\bul \leq N}}
@>{c_{1,{\rm crys}}}>>Ru_{Y_{\bul \leq N}/T*}({\cal O}_{Y_{\bul \leq N}/T})[1]\\
@| @VVV \\
M^{\rm gp}_{Y_{\bul \leq N}}
@>>>
R\pi_{{\rm zar}*}
({\cal O}_{{\mathfrak E}_{\bul \leq N,\bul}}
\otimes_{{\cal O}_{{\cal Y}_{\bul \leq N,\bul}}}
\Om^{\bul}_{{\cal Y}_{\bul \leq N,\bul}/T})[1]. 
\end{CD}
\end{equation*}
Lastly we consider the following complex 
\begin{align*} 
M^{\rm gp}_{{\mathfrak E}_{\bul \leq N,\bul}}\os{d\log}{\lo} 
{\cal O}_{{\mathfrak E}_{\bul \leq N,\bul}}
\otimes_{{\cal O}_{{\cal Y}_{\bul \leq N,\bul}}}
\Om^1_{{\cal Y}_{\bul \leq N,\bul}/T}\os{d}{\lo}
{\cal O}_{{\mathfrak E}_{\bul \leq N,\bul}}
\otimes_{{\cal O}_{{\cal Y}_{\bul \leq N,\bul}}}
\Om^2_{{\cal Y}_{\bul \leq N,\bul}/T} \os{d}{\lo}\cdots, 
\end{align*} 
which we denote by 
${\cal O}_{{\mathfrak E}_{\bul \leq N,\bul}}
\otimes_{{\cal O}_{{\cal Y}_{\bul \leq N,\bul}}}
\Om^{\times}_{{\cal Y}_{\bul \leq N,\bul}/T}$. 
Set 
$${\cal K}^{\times}_{{\cal Y}_{\bul \leq N,\bul}/T}:={\rm Ker}
({\cal O}_{{\mathfrak E}_{\bul \leq N,\bul}}
\otimes_{{\cal O}_{{\cal Y}_{\bul \leq N,\bul}}}
\Om^{\times}_{{\cal Y}_{\bul \leq N,\bul}/T}\lo 
M^{\rm gp}_{Y_{\bul \leq N,\bul}}).$$ 
Set also 
$${\cal J}^{\bul}_{{\cal Y}_{\bul \leq N,\bul}}:={\rm Ker}
({\cal O}_{{\mathfrak E}_{\bul \leq N,\bul}}
\otimes_{{\cal O}_{{\cal Y}_{\bul \leq N,\bul}}}
\Om^{\bul}_{{\cal Y}_{\bul \leq N,\bul}/T}\lo 
{\cal O}_{Y_{\bul \leq N,\bul}}).$$ 
Then we obtain the following morphism as in \cite[(3.2)]{boi}: 
$\log \col {\cal K}^{\times}_{{\cal Y}_{\bul \leq N,\bul}/T}\lo 
{\cal J}^{\bul}_{{\cal Y}_{\bul \leq N,\bul}}$. 
Because following exact sequence 
\begin{align*} 
0\lo {\cal K}^{\times}_{{\cal Y}_{\bul \leq N,\bul}/T}\lo 
{\cal O}_{{\mathfrak E}_{\bul \leq N,\bul}}
\otimes_{{\cal O}_{{\cal Y}_{\bul \leq N,\bul}}}
\Om^{\times}_{{\cal Y}_{\bul \leq N,\bul}/T}\lo 
M^{\rm gp}_{Y_{\bul \leq N,\bul}}\lo 0
\end{align*} 
is exact, 
we obtain the following morphism
\begin{align*} 
M^{\rm gp}_{Y_{\bul \leq N}}\lo 
R\pi_{{\rm zar}*}({\cal K}^{\times}_{{\cal Y}_{\bul \leq N,\bul}/T})[1].
\end{align*} 
Hence we obtain the following composite morphism 
\begin{align*} 
M^{\rm gp}_{Y_{\bul \leq N}}\lo 
R\pi_{{\rm zar}*}({\cal K}^{\times}_{{\cal Y}_{\bul \leq N,\bul}/T})[1]
\os{\log}{\lo} R\pi_{{\rm zar}*}({\cal J}^{\bul}_{{\cal Y}_{\bul \leq N,\bul}})[1]
\lo R\pi_{{\rm zar}*}({\cal O}_{{\mathfrak E}_{\bul \leq N,\bul}}
\otimes_{{\cal O}_{{\cal Y}_{\bul \leq N,\bul}}}
\Om^{\bul}_{{\cal Y}_{\bul \leq N,\bul}/T})[1].
\tag{1.8.16.3}\label{ali:zac}
\end{align*} 
Because the following diagram is commutative
\begin{equation*} 
\begin{CD} 
M^{\rm gp}_{Y_{\bul \leq N}}@>>>Ru_{Y_{\bul \leq N}/T*}
R\pi_{{\rm crys}*}({\cal K}^{\times})[1]@>{\log}>>\\
@| @VVV \\
M^{\rm gp}_{Y_{\bul \leq N}}@>>>
R\pi_{{\rm zar}*}({\cal K}^{\times}_{{\cal Y}_{\bul \leq N,\bul}/T})[1]
@>{\log}>>
\end{CD}
\tag{1.8.16.4}\label{cd:lchnn}
\end{equation*} 
\begin{equation*} 
\begin{CD} 
Ru_{Y_{\bul \leq N}/T*}
R\pi_{{\rm crys}*}({\cal K}^{\bul})[1] @>>>
Ru_{Y_{\bul \leq N,\bul}/T*}(L({\cal O}_{{\mathfrak E}_{\bul \leq N,\bul}}
\otimes_{{\cal O}_{{\cal Y}_{\bul \leq N,\bul}}}
\Om^{\bul}_{{\cal Y}_{\bul \leq N,\bul}/T}))[1]\\
@VVV @VVV \\
R\pi_{{\rm zar}*}({\cal J}^{\bul}_{{\cal Y}_{\bul \leq N,\bul}})[1]
@>>>
R\pi_{{\rm zar}*}({\cal O}_{{\mathfrak E}_{\bul \leq N,\bul}}
\otimes_{{\cal O}_{{\cal Y}_{\bul \leq N,\bul}}}
\Om^{\bul}_{{\cal Y}_{\bul \leq N,\bul}/T})[1]
\end{CD}
\end{equation*}
as in \cite[(3.3.2)]{boi}, we obtain the following: 

\begin{prop}\label{prop:ofc}
The first log crystalline chern class map
$c_{1,{\rm crys}}$ in {\rm (\ref{ali:mygp})} is equal to the morphism 
{\rm (\ref{ali:zac})}
via the identification 
$$Ru_{Y/T*}({\cal O}_{Y_{\bul \leq N}/T})=
R\pi_{{\rm zar}*}({\cal O}_{{\mathfrak E}_{\bul \leq N,\bul}}
\otimes_{{\cal O}_{{\cal Y}_{\bul \leq N,\bul}}}
\Om^{\bul}_{{\cal Y}_{\bul \leq N,\bul}/T}).$$ 
\end{prop}

\begin{conj}[{\bf Variational filtered log $p$-adic hard Lefschetz conjecture}]\label{conj:lhlc}
The following cup product 
\begin{equation*} 
\lam^i_{S(T)^{\nat}} \col 
R^{d-i}f_{X_{\os{\circ}{T}_0}/S(T)^{\nat}*}
({\cal O}_{X_{\os{\circ}{T}_0}/S(T)^{\nat}})
\otimes_{\mab Z}{\mab Q} 
\lo (R^{d+i}f_{X_{\os{\circ}{T}_0}/S(T)^{\nat}*}
({\cal O}_{X_{\os{\circ}{T}_0}/S(T)^{\nat}})
\otimes_{\mab Z}{\mab Q})(i) \quad (i\in {\mab N})
\tag{1.8.18.1}\label{eqn:fcpl} 
\end{equation*}
is an isomorphism. 
In fact, $\lam^i_{S(T)^{\nat}}$ 
is the following isomorphism of filtered sheaves: 
\begin{equation*} 
\lam^i_{S(T)^{\nat}} \col 
(R^{d-i}f_{X_{\os{\circ}{T}_0}/S(T)^{\nat}*}
({\cal O}_{X_{\os{\circ}{T}_0}/S(T)^{\nat}})
\otimes_{\mab Z}{\mab Q},P) 
\os{\sim}{\lo} ((R^{d+i}f_{X_{\os{\circ}{T}_0}/S(T)^{\nat}*}({\cal O}_{X_{\os{\circ}{T}_0}/S(T)^{\nat}})
\otimes_{\mab Z}{\mab Q})(i),P). 
\tag{1.8.18.2}\label{eqn:filqpl} 
\end{equation*}
Here $(i)$ means the Tate twist which has a meaning when 
$S(T)^{\nat}$ has a lift of the Frobenius endomorphism of 
$S(T)^{\nat}~{\rm mod}~p$. 
\end{conj}

\begin{rema} 
(1) In \cite[(9.10)]{nlpi} we have proved that, if the morphism 
(\ref{eqn:fcpl}) in the case $i=1$ is an isomorphism, then  
(\ref{conj:rcpmc}) holds for $q=1$ and $q=2d-1$ 
in the case where $S$ is the log point of a perfect field of characteristic $p> 0$ 
and $T={\cal W}(t)$, where $t$ is a fine log scheme 
whose underlying scheme is the spectrum of a perfect field of characteristic $p> 0$.  
\par
(2) In \cite{kaji} Kajiwara has proved that
the analogous morphism to the morphism (\ref{eqn:fcpl}) in the cases $i=1$ and $i=2d-1$ 
for $l$-adic Kummer \'{e}tale cohomologies is an isomorphism when $S(T)^{\nat}$ is a log point of a field. 
In the Chapter V we shall prove a generalized analogue for log crystalline cohomologies. 
\par 
(3) As in \cite{bbd}, 
one can raise the problem 
when the following morphism  
\begin{equation*} 
\lam^i_{S(T)^{\nat}} \col 
R^{d-i}f_{X_{\os{\circ}{T}_0}/S(T)^{\nat}*}
(\eps^*_{X_{\os{\circ}{T}_0}/S(T)^{\nat}}(E))\otimes_{\mab Z}{\mab Q} 
\lo (R^{d+i}f_{X_{\os{\circ}{T}_0}/S(T)^{\nat}*}
(\eps^*_{X_{\os{\circ}{T}_0}/S(T)^{\nat}}(E))\otimes_{\mab Z}{\mab Q})(i)
\tag{1.8.19.1}\label{eqn:fccpl} 
\end{equation*}
is an isomorphism. 
However 
I cannot formulate the precise statement 
for nontrivial coefficients for which 
(\ref{eqn:fccpl}) should be an isomorphism. 
\end{rema}

\chapter{Weight filtrations and slope filtrations on log crystalline cohomologies
via log de Rham-Witt complexes} 

\parno  
In this chapter we generalize Hyodo-Kato's comparison theorem (\cite{hk}, \cite{ndw})
which tells us that the log de Rham-Witt complex of a log smooth scheme 
over a fine log scheme whose underlying scheme is the spectrum of 
a perfect field of characteristic $p>0$ 
calculate the log crystalline complex of it. 
This generalization is a log version 
and a truncated simplicial version of Etesse's comparison theorem (\cite{et})  
which is a generalization of Bloch-Illusie's comparison theorem (\cite{idw}). 
We construct the $p$-adic filtered Steenbrink complex of 
a simplicial SNCL scheme over the log point of a perfect field of characteristic $p>0$ 
by using the filtered log de Rham-Witt complex of it. 
This is a generalization of the Hyodo-Mokrane-Steenbrink complex defined by Mokrane 
(\cite{msemi}) and corrected by the author (\cite{ndw}).  
We give a comparison theorem between 
the $p$-adic filtered Steenbrink complexes in the Chapters I and II.  
We also define the $p$-adic monodromy operator 
(resp.~the $p$-adic quasi-monodromy operator)
of the log de Rham-Witt complex (resp.~the Hyodo-Mokrane-Steenbrink complex) of 
a simplicial SNCL scheme and give a comparison theorem 
between the $p$-adic monodromy operators and 
the $p$-adic quasi-monodromy operators in the Chapters I and II.

\section{Log de Rham-Witt complexes of crystals I}
\label{sec:ldrwc} 
In this section we construct a 
theory of log de Rham-Witt complexes of 
unit root log $F$-crystals on a log smooth scheme  
of Cartier type over a fine log scheme 
whose underlying scheme is 
the spectrum of a perfect field of characteristic $p>0$. 
This is a log version of Etesse's theory in \cite{et} 
and a generalization of Hyodo-Kato's theory in \cite{hk}.
The proofs of several results in this section 
are the imitations of Etesse's proofs in [loc.~cit.]; 
we give only the proofs which are not the imitations. 
\par 
Let $(S,{\cal I},\gam)$ be a fine log PD-scheme 
such that $p$ is locally nilpotent on $\os{\circ}{S}$ 
or a fine log $p$-adic formal PD-scheme such that 
$p\in {\cal I}$. 
(In this section $S$ is not necessarily a family of log points.)
Set $S(0):=\ul{\rm Spec}^{\log}_S({\cal O}_S/{\cal I})$. 
Let $Y$ be a fine log scheme over $S(0)$. 
Let $g\col Y\lo S(0)$ be the structural morphism. 
By abuse of notation, we denote the composite structural morphism 
$Y\os{g}{\lo} S(0)\os{\sus}{\lo} S$ by $g\col Y\lo S$. 
Let $(U,T,\del)$ be an object of 
the log crystalline site 
${\rm Crys}(Y/(S,{\cal I},\gam))$ of $Y/S$.  
Set ${\cal J}:={\rm Ker}({\cal O}_T\lo {\cal O}_U)$. 
Let ${\mathfrak D}_{T/S}(1)$ be 
the log PD-envelope 
of the diagonal immersion $T\os{\sus}{\lo} T\times_ST$ 
over $(S,{\cal I},\gam)$: 
${\mathfrak D}_{T/S}(1)$ is the usual PD-envelope of 
the immersion $T\os{\sus}{\lo} (T\times_ST)^{\rm ex}$ 
over $(S,{\cal I},\gam)$ endowed with 
the inverse image of the log structure of $T$ by the ``projection'' 
$p_i \col (T\times_ST)^{\rm ex}\lo T$  
$(i=1,2)$ (note that $p_1^*(M_T)=p_2^*(M_T)$).  
Let $\ol{\cal J}$ be the PD-ideal sheaf of 
${\cal O}_{{\mathfrak D}_{T/S}(1)}$. 
By abuse of notation, 
let us denote by $p_i\col {\mathfrak D}_{T/S}(1)\lo T$ 
the induced morphism by the morphism 
$p_i \col (T\times_ST)^{\rm ex}\lo T$.   
Set ${\cal D}^n_{T/S}(1):=
{\cal O}_{{\mathfrak D}_{T/S}(1)}/\ol{\cal J}{}^{[n+1]}$ 
and $\ol{\cal J}_n:=\ol{\cal J}/\ol{\cal J}{}^{[n+1]}$. 
Let ${\cal J}_n$ be the inverse image of 
${\cal J}$ in ${\cal D}^n_{T/S}(1)$. 
Then we have the following commutative diagram of 
split exact sequences:   
\begin{equation*} 
\begin{CD} 
0@>>> \ol{\cal J}_n @>>> {\cal D}^n_{T/S}(1) @>>> {\cal O}_T@>>> 0\\ 
@. @| @A{\bigcup}AA @A{\bigcup}AA \\ 
0@>>> \ol{\cal J}_n @>>> {\cal J}_n @>>> {\cal J} @>>> 0. 
\end{CD} 
\tag{2.0.0.1}\label{eqn:jndts} 
\end{equation*}
The two splittings are given by the two morphisms 
$p_i^*\col {\cal O}_T\lo {\cal D}^n_{T/S}(1)$ $(i=1,2)$. 
Then we have two decompositions 
${\cal J}_n\simeq \ol{\cal J}_n\oplus {\cal J}$ 
by using $p_1^*$ and $p_2^*$. 
Consequently we have two PD-structures 
$\del_1$ and $\del_2$ 
on ${\cal J}_n$ (\cite[3.12 Proposition]{bob}). 
Following \cite[II (1.1.1)$\sim$(1.1.4)]{et}, 
set ${\cal D}^n_{T/S,\del}(1):={\cal D}^n_{T/S}(1)/
\langle (\del_1)_m(x)-(\del_2)_m(x)~\vert~x\in {\cal J}_n, 
m\in {\mab N}\rangle$. 
Let $\ol{\cal J}_n/\langle ~ \rangle$ be the image of 
$\ol{\cal J}_n$ in ${\cal D}^n_{T/S,\del}(1)$. 
Let ${\mathfrak D}^n_{T/S,\del}(1)$ 
be the exact closed log subscheme of 
${\mathfrak D}^n_{T/S}(1)$ whose structure 
sheaf is ${\cal D}^n_{T/S,\del}(1)$. 
Let $\bigoplus_{i\geq 0}
{\Om}^i_{T/S,\del}$ be 
a sheaf of differential graded algebras over 
${\cal O}_S$ which is a quotient of 
$\bigoplus_{i\geq 0}{\Om}^i_{T/S}$ 
divided by the ideal sheaf generated by local sections of the following form 
$d\del_j(a)-\del_{j-1}(a)da$ ($a\in {\cal J}$, $j\geq 1$). 
Then, by (\ref{ali:omyti}) and the calculation in the proof of \cite[II (1.1.4)]{et}, 
we have an equality 
\begin{align*} 
\ol{\cal J}_1/\langle ~ \rangle=\Om^1_{T/S,\del}.  
\tag{2.0.0.2}\label{eqn:joms} 
\end{align*} 
\par 
Let $E$ be a quasi-coherent log crystal of 
${\cal O}_{Y/S}$-modules. 
Because we have the ``projections'' 
$p_1, p_2\col (U,{\mathfrak D}^1_{T/S,\del}(1))\lo (U,T)$, 
we have canonical isomorphisms 
$p_2^*(E_T)\os{\sim}{\lo}
E_{{\mathfrak D}^1_{T/S,\del}(1)}
\os{\sim}{\longleftarrow}p_1^*(E_T)$. 
Using the difference of these two canonical isomorphisms,  
and considering the modulo reduction ``$/\langle~\rangle$'',  
we have the following integrable connection 
\begin{equation*} 
\nabla \col E_T\lo E_T
\otimes_{{\cal O}_T}{\Om}^1_{T/S,\del} 
\tag{2.0.0.3}\label{eqn:intcon}
\end{equation*} 
as usual (cf.~[loc.~cit.]). 
\par 
Let $\kap$ be a perfect field of 
characteristic $p>0$.  
Let $s$ be a fine log scheme 
whose underlying scheme is ${\rm Spec}(\kap)$. 
Let ${\cal W}_n$ (resp.~${\cal W}$) be 
the Witt ring of $\kap$ of length $n\in {\mab Z}_{\geq 1}$ 
(resp.~the Witt ring of $\kap$).  
Let $Y$ be a log smooth scheme of Cartier type over $s$. 
Let  $g\col Y\lo s$ be the structural morphism. 
Let ${\cal W}_n(Y)$ be the canonical lift of $Y$ defined in \cite[(3.1)]{hk}: 
${\cal W}_n(Y)
=(\os{\circ}{Y},M_Y\oplus {\rm Ker}({\cal W}_n({\cal O}_Y)^*\lo {\cal O}_Y^*)
\lo {\cal W}_n({\cal O}_Y))$.   
In particular, we have the canonical lift 
${\cal W}_n(s)$ of $s$ over 
${\cal W}_n$ and hence 
the formal canonical lift ${\cal W}(s)$ over ${\cal W}$.  
We also denote  
the structural morphism $Y\os{g}{\lo} s \os{\sus}{\lo} {\cal W}_n(s)$ 
by $g\col Y\lo {\cal W}_n(s)$. 
The log scheme ${\cal W}_n(Y)$ has 
the canonical PD-structure $[~]$ on $V{\cal W}_{n-1}({\cal O}_Y)$ 
(\cite[0 (1.4.3)]{idw}): 
$[Vx]_m:=p^{m-1}/m!\cdot V(x^m)$ 
$(m\in{\mab Z}_{\geq 1}, x\in {\cal W}_{n-1}({\cal O}_Y))$. 
\par 
Let $E$ be a quasi-coherent log crystal of 
${\cal O}_{Y/{\cal W}(s)}$-modules. 
Let $\iota \col {\cal W}_n(Y)\os{\sus}{\lo} {\cal W}_{n+1}(Y)$ 
be the natural closed immersion.  
Set $E_n:=E_{{\cal W}_n(Y)}$.  
For a quasi-coherent log crystal $E'$ of 
${\cal O}_{Y/{\cal W}_n(s)}$-modules 
$(n\in {\mab Z}_{\geq 1})$, set also 
$E'_r:= E'_{{\cal W}_r(Y)}$ $(1\leq r\leq n)$.  
Since $E$ is a crystal, we have 
the isomorphism $\iota^*(E_{n+1})=E_n$.  
Hence we have the projection $R\col E_{n+1}\lo E_n$. 
We also have the projection $R\col E'_{r+1}\lo E'_r$ 
$(1\leq r\leq n-1)$. 
\par
Let ${\cal W}_n({\cal O}_Y)'$ and ${\cal W}_n({\cal O}_Y)$ be 
the obverse and reverse Witt ring sheaves of $Y$ of length $n$, 
respectively, in the sense of \cite[(7.8)]{ndw}:  
${\cal W}_n({\cal O}_Y)'$ is, by definition,  
the structure sheaf of 
${\cal W}_n(Y)$; on the other hand, 
if there exists an immersion $Y\os{\sus}{\lo} {\cal Q}$ 
into a log smooth scheme over ${\cal W}_n(s)$, 
$${\cal W}_n({\cal O}_Y)=
{\cal H}^0({\cal O}_{{\mathfrak E}\otimes_{{\cal O}_{\cal Q}}}
\Om^*_{{\cal Q}/{\cal W}_n(s)})=
{\cal H}^0(\Om^*_{{\mathfrak E}/{\cal W}_n(s),[~]}),$$ 
where ${\mathfrak E}$ is the log PD-envelope 
of the immersion $Y\os{\sus}{\lo}{\cal Q}$ over 
$({\cal W}_n(s),p{\cal W}_n,[~])$.    
Here we have used (\ref{eqn:fwniwu}) for the last equality. 
We have the Cartier isomorphism 
$C^{-n}\col  {\cal W}_n({\cal O}_Y)'\os{\sim}{\lo} 
{\cal W}_n({\cal O}_Y)$ (\cite[(7.5)]{ndw}).   
If there exists an immersion $Y\os{\sus}{\lo} {\cal Q}$ 
into a log smooth scheme over ${\cal W}_n(s)$,
then $C^{-n}$ is defined by the following morphism  
\begin{align*} 
s_n(0,0)\col {\cal W}_n({\cal O}_Y)'\owns 
(a_0,\ldots, a_{n-1})\lom   
\sum_{i=0}^{n-1}\wt{a}_i^{p^{n-i}}p^i
\in {\cal H}^0(\Om^*_{{\mathfrak E}/{\cal W}_n(s),[~]})=
{\cal W}_n({\cal O}_Y),
\tag{2.0.0.4}\label{ali:owny} 
\end{align*} 
where $\wt{a}_i\in {\cal O}_{\mathfrak E}$ is a lift of $a_i$ (e.~g., \cite[(4.9)]{hk}). 
In [loc.~cit.] Hyodo and Kato have proved that  
the morphism (\ref{ali:owny}) is independent of the choice of the lift 
$\wt{a}_i$. By virtue of this fact, it is easy to check that 
the morphism (\ref{ali:owny}) is independent 
of the choice of the immersion $Y\os{\sus}{\lo} {\cal Q}$ over ${\cal W}_n(s)$. 
The ${\cal W}_n({\cal O}_Y)'$-module $E_n$ becomes a 
${\cal W}_n({\cal O}_Y)$-module via the isomorphism (\ref{ali:owny}). 
\par 
Next we recall the obverse log de Rham-Witt complex 
$({\cal W}_n{\Om}^{\bul}_{Y})''$ in \cite[\S7]{ndw} which 
is a correction of $({\cal W}_n{\Om}^{\bul}_{Y})'$ in \cite[(4.6)]{hk}. 
In this book we denote $({\cal W}_n{\Om}^{\bul}_{Y})''$ in 
\cite[\S7]{ndw} by $({\cal W}_n{\Om}^{\bul}_{Y})'$ 
for simplicity of notation. 
Let $\al \col M_Y \lo {\cal O}_Y$ be the structural morphism.
\par  
For a local section $b\in {\cal O}_Y$, set $V^j(b)=({\us{j~{\rm times}}
{\underbrace{0, \cdots , 0}}},b,0,\ldots,0)\in 
{\cal W}_n({\cal O}_Y)$ $(0 \leq j < n)$. 
The sheaf $({\cal W}_n{\Om}_Y^i)'$ on 
$\os{\circ}{Y}$ is a quotient of 
\begin{equation*}
\{{\cal W}_n({\cal O}_Y)'\otimes_{\mab Z}
\bigwedge^i(M^{\rm gp}_Y/g^{-1}(M^{\rm gp}_s))\}
\oplus 
\{{\cal W}_n({\cal O}_Y)'\otimes_{\mab Z}
\bigwedge^{i-1}(M^{\rm gp}_Y/g^{-1}(M^{\rm gp}_s))\}
\tag{2.0.0.5}\label{eqn:flsodw}
\end{equation*}
divided by a subsheaf ${\cal F}_n$ of ${\mab Z}$-modules 
generated by  the images of local sections of the following type
\begin{equation*}
(V^j(\al(a_1))\otimes(a_1\wedge \cdots 
\wedge a_i), 0)-p^j(0, V^j(\al(a_1))
\otimes (a_2\wedge \cdots \wedge a_i)) 
\tag{2.0.0.6}\label{eqn:frel}
\end{equation*}
$$(a_1, \ldots, a_i \in M_Y,~0 \leq j < n)$$
and a subsheaf ${\cal G}_n$ of ${\mab Z}$-modules 
generated by the images of local sections of the following type 
\begin{equation*}
(0,V^j(\al(a_2)^e)\otimes (a_2\wedge \cdots \wedge a_i)) 
\quad
(a_2, \ldots, a_i\in M_Y,~1 \leq j <n,~ e \in {\mab Z}_{\geq 1})
\tag{2.0.0.7}\label{eqn:nrelp}
\end{equation*}
(in this book we denote ${\cal G}'_n$ in \cite[p.~567]{ndw} 
by ${\cal G}_n$). 
\parno 
We define the boundary morphism $d$ by the following formula
\begin{align*}
&d(b\otimes(a_1 \wedge \cdots \wedge a_i), 
c\otimes(a'_2 \wedge \cdots 
\wedge a'_i))
=(0, b\otimes(a_1 \wedge \cdots \wedge a_i)) 
\tag{2.0.0.8}\label{eqn:dfdbac}\\
&(a_1, \ldots, a_i, a'_2, \cdots, a'_i \in M_Y,~b,c\in {\cal W}_n({\cal O}_{Y'})). 
\end{align*}
It is easy to check that $d$ is well-defined and that $d^2=0$. 
The complex $({\cal W}_n{\Om}^{\bul}_Y)'$ is contravariantly functorial with respect to 
a morphism of log smooth schemes of Cartier type over fine log schemes 
whose underlying schemes are the spectrums of perfect fields of characteristic $p>0$. 
\par 
Let ${\cal W}_n{\Om}^{\bul}_Y$ be the reverse log de Rham-Witt complex 
of $Y/s$ of level $n$ in the sense of \cite[(7.8)]{ndw} (\cite[(4.1)]{hk}): 
\begin{align*} 
{\cal W}_n{\Om}^i_Y:=R^iu_{Y/{\cal W}_n(s)*}({\cal O}_{Y/{\cal W}_n(s)}) 
\quad (i\in {\mab N}).
\tag{2.0.0.9}\label{eqn:lmntwy} 
\end{align*} 
(The sheaf ${\cal W}_n{\Om}^i_Y$ is denoted by $W_n\om^i_Y$ in 
\cite{hdw}, \cite{hk} and \cite{msemi}.)   
In a standard way, we have a boundary morphism 
$d\col {\cal W}_n{\Om}^i_Y\lo {\cal W}_n{\Om}^{i+1}_Y$ 
(\cite[p.~247]{hk}).
As in the case of $({\cal W}_n{\Om}^{\bul}_Y)'$,  
the complex ${\cal W}_n{\Om}^{\bul}_Y$ is also contravariantly functorial. 
There exists the following Cartier isomorphism of complexes of abelian sheaves on $Y$ 
(\cite[(7.5), (7.7)]{ndw}), which is a nontrivial correction of \cite[(4.6)]{hk}: 
\begin{equation*} 
C^{-n} \col ({\cal W}_n{\Om}^i_Y)' \os{\sim}{\lo} 
{\cal W}_n{\Om}^i_Y \quad (i\in {\mab N}).  
\tag{2.0.0.10}\label{eqn:wnymod}
\end{equation*} 
(In \cite[(7.5), (7.7)]{ndw} we have denoted the isomorphism $C^{-n}$ by $s_n$.)
Because we need the definition of $C^{-n}$ in (\ref{prop:nmc}) below, we recall it. 
\par 
As in \cite[(4.9)]{hk} and \cite[p.~589]{ndw}, we define three morphisms as follows. 
\par 
Assume that there exists an immersion $Y\os{\sus}{\lo} {\cal Q}$ 
into a log smooth scheme over ${\cal W}_n(s)$. 
In addition to the morphism $s_n(0,0)$ in (\ref{ali:owny}), 
we define other two morphisms as follows: 
\begin{equation*}
s_n(1,0) \col {\cal W}_n({\cal O}_Y) \owns 
(a_0, \ldots, a_{n-1}) \lom 
\left[\sum_{i=0}^{n-1}\wt{a}_i^{p^{n-i}-1}d\wt{a}_i\right]
\in {\cal H}^1(\Om^*_{{\mathfrak E}/{\cal W}_n(s),[~]}),
\tag{2.0.0.11}\label{eqn:tisn10}
\end{equation*}
\begin{equation*}
d\log  \col M^{\rm gp}_Y/g^{-1}(M^{\rm gp}_s)\owns b 
\lom \left[d\log \wt{b}\right]\in 
{\cal H}^1(\Om^*_{{\mathfrak E}/{\cal W}_n(s),[~]}),
\tag{2.0.0.12}\label{eqn:tdloegb}
\end{equation*}
where $\wt{a}_i\in {\cal O}_{\mathfrak E}$ and 
$\wt{b}\in M_{\mathfrak E}$ are lifts of 
$a_i\in {\cal O}_Y$ and $b\in M_Y\subset M_{{\cal W}_n(Y)}$, respectively 
and $[~]$ means the class. 
Then, by the same proof as that for the well-definedness of $s_n(0,0)$, 
$s_n(1,0)$ and $d\log$ are well-defined (\cite[(4.9)]{hk}). 
Then the following morphism 
\begin{align*} 
&(s_n(0,0)\otimes d\log \wedge \cdots \wedge d\log) \oplus 
(s_n(1,0)\otimes d\log \wedge \cdots \wedge d\log)  \col \\
&\{{\cal W}_n({\cal O}_Y)'\otimes_{\mab Z}
\bigwedge^i(M^{\rm gp}_Y/g^{-1}(M^{\rm gp}_s))\}
\oplus 
\{{\cal W}_n({\cal O}_Y)'\otimes_{\mab Z}
\bigwedge^{i-1}(M^{\rm gp}_Y/g^{-1}(M^{\rm gp}_s))\}
\lo 
{\cal W}_n{\Om}^{\bul}_Y
\end{align*} 
induces the isomorphism (\ref{eqn:wnymod}) (\cite[(7.5)]{ndw}).
\par 
Because ${\cal W}_n{\Om}^i_Y$ $(i\in {\mab N})$ 
has a structure of ${\cal W}_n({\cal O}_Y)$-modules,  
$({\cal W}_n{\Om}^i_Y)' $ 
has a structure of ${\cal W}_n({\cal O}_Y)'$-modules
via the isomorphism (\ref{eqn:wnymod}). 
\par 
The following is a log version of \cite[I (1.9)]{idw} 
(modulo the equality of the de Rham-Witt complexes above 
in the case of the trivial log structure and that in [loc.~cit.]):

\begin{prop}\label{prop:bcdw}
Let $t\lo s$ be a morphism of fine log schemes 
whose underlying schemes are the spectrums of 
perfect fields of characteristic $p>0$. 
Set $\kap_t:=\Gam(t,{\cal O}_t)$ and $Y_t:=Y\times_st$. 
Let $g_t\col Y_t\lo t$ be the base change of $g\col Y\lo s$. 
Let $q\col Y_t \lo Y$ be the projection. 
Let $\star$ be nothing or $\prime$. 
Assume that $g$ is integral. 
Then the canonical morphism 
\begin{equation*} 
q^{-1}(({\cal W}_n{\Om}_Y^i)^{\star})
\otimes_{{\cal W}_n}{\cal W}_n(\kap_t)
\lo ({\cal W}_n{\Om}_{Y_t}^i)^{\star}
\tag{2.1.1.1}\label{eqn:drkp}
\end{equation*}
is an isomorphism. 
In particular, if $\os{\circ}{t}=\os{\circ}{s}$, 
then 
the canonical morphism 
\begin{equation*} 
q^{-1}(({\cal W}_n{\Om}_Y^i)^{\star})
\lo ({\cal W}_n{\Om}_{Y_t}^i)^{\star}
\tag{2.1.1.2}\label{eqn:drtkp}
\end{equation*}
is an isomorphism $($note that $\os{\circ}{Y}_t=\os{\circ}{Y})$. 
\end{prop} 
\begin{proof} 
It suffices to prove (\ref{prop:bcdw}) in the case $\star=$nothing.  
The problem is local on $Y$. Hence we may assume that 
there exists an immersion $Y\os{\sus}{\lo}{\cal Q}$ 
into a log smooth integral scheme over ${\cal W}_n(s)$ 
by the proofs of \cite[(3.14)]{klog1} and \cite[(4.7)]{ny}. 
In this case, 
\begin{equation*} 
{\cal W}_n{\Om}_{Y}^i
={\cal H}^i
({\cal O}_{\mathfrak E}\otimes_{{\cal O}_{\cal Q}}\Om^*_{{\cal Q}/{\cal W}_n(s)}). 
\tag{2.1.1.3}\label{eqn:laey} 
\end{equation*}   
Set ${\cal Q}_{{\cal W}_n(t)}:={\cal Q}\times_{{\cal W}_n(s)}{\cal W}_n(t)$ 
and 
${\mathfrak E}_{{\cal W}_n(t)}:={\mathfrak E}\times_{{\cal W}_n(s)}{\cal W}_n(t)$. 
Because the structural morphism ${\cal Q}\lo {\cal W}_n(s)$ is integral, 
$({\cal Q}_{{\cal W}_n(t)})^{\circ}
=\os{\circ}{\cal Q}\times_{{\cal W}_n(\os{\circ}{s})}{\cal W}_n(\os{\circ}{t})$. 
Because the extension ${\cal W}_n(\kap_t)/{\cal W}_n(\kap)$ is flat, 
${\mathfrak E}_{{\cal W}_n(t)}$ is the log PD-envelope of 
the immersion $Y_t\os{\sus}{\lo} {\cal Q}_{{\cal W}_n(t)}$ 
over $({\cal W}_n(t),p{\cal W}_n(\kap_t),[~])$ by \cite[Proposition 3.21]{bob}. 
Let $\wt{q}\col {\cal Q}_{{\cal W}_n(t)}\lo {\cal Q}$ be the first projection. 
By \cite[(1.7)]{klog1} we have $\Om^i_{{\cal Q}_{{\cal W}_n(t)}/{\cal W}_n(t)}=
\wt{q}{}^*(\Om^i_{{\cal Q}/{\cal W}_n(s)})$. 
Hence we have the following equalities: 
\begin{align*} 
{\cal W}_n{\Om}_{Y_t}^i& ={\cal H}^i
({\cal O}_{{\mathfrak E}_{{\cal W}_n(t)}}
\otimes_{{\cal O}_{{\cal Q}_{{\cal W}_n(t)}}}\Om^*_{{\cal Q}_{{\cal W}_n(t)}/{\cal W}_n(t)}) \\
&=q^{-1}({\cal H}^i({\cal O}_{\mathfrak E}
\otimes_{{\cal O}_{\cal Q}}\Om^*_{{\cal Q}/{\cal W}_n(s)}))\otimes_{{\cal W}_n}
{\cal W}_n(\kap_t)\\
&= q^{-1}(({\cal W}_n{\Om}_Y^i)^{\star})
\otimes_{{\cal W}_n}{\cal W}_n(\kap_t). 
\end{align*} 
\end{proof}

\begin{rema}
In (\ref{coro:lemos}) below we shall generalize  (\ref{prop:bcdw}). 
\end{rema}

\par 
As in \cite[(4.9)]{hk}
and the proof of \cite[(4.19)]{hk} (cf.~\cite[p.~589]{ndw}), 
we see that there exists the following natural morphism 
\begin{equation*} 
{\Om}^{\bul}_{{\cal W}_n(Y)/{\cal W}_n(s), [~]}
\lo {\cal W}_n{\Om}^{\bul}_Y 
\tag{2.1.2.1}\label{eqn:lmwy} 
\end{equation*}  
of dga's over $g^{-1}({\cal O}_{{\cal W}_n(s)})$.  
If there exists an immersion $Y\os{\sus}{\lo}{\cal Q}$ 
into a log smooth scheme over ${\cal W}_n(s)$, the morphism 
(\ref{eqn:lmwy}) is equal to the following morphism: 
\begin{equation*} 
{\Om}^{\bul}_{{\cal W}_n(Y)/{\cal W}_n(s), [~]}
\lo {\cal H}^{\bul}
({\cal O}_{\mathfrak E}\otimes_{{\cal O}_{\cal Q}}\Om^*_{{\cal Q}/{\cal W}_n(s)})
={\cal H}^{\bul}
({\cal O}_{\mathfrak E}\otimes_{{\cal O}_{\cal Q}}\Om^*_{{\cal Q}^{\rm ex}/{\cal W}_n(s)})
={\cal H}^{\bul}(\Om^*_{{\mathfrak E}/{\cal W}_n(s),[~]}). 
\tag{2.1.2.2}\label{eqn:lay} 
\end{equation*}   
This morphism is, by definition, the induced morphism by 
$s_n(0,0)$ and the following two morphisms
\begin{equation*}
s_n(1,0)' \col \Om^1_{{\cal W}_n(Y)/{\cal W}_n(s), [~]} \owns 
d(a_0, \ldots, a_{n-1}) \lom 
\left[\sum_{i=0}^{n-1}\wt{a}_i^{p^{n-i}-1}d\wt{a}_i\right]
\in {\cal H}^1(\Om^*_{{\mathfrak E}/{\cal W}_n(s),[~]}),  
\tag{2.1.2.3}\label{eqn:tise10}
\end{equation*}
\begin{equation*}
\iota  \col \Om^1_{{\cal W}_n(Y)/{\cal W}_n(s), [~]} \owns 
d\log b \lom \left[d\log \wt{b}\right]\in 
{\cal H}^1(\Om^*_{{\mathfrak E}/{\cal W}_n(s),[~]}),
\tag{2.1.2.4}\label{eqn:tdlogb}
\end{equation*}
where $\wt{a}_i\in {\cal O}_{\mathfrak E}$ and 
$\wt{b}\in M_{\mathfrak E}$ are lifts of 
$a_i\in {\cal O}_Y$ and $b\in M_Y\subset M_{{\cal W}_n(Y)}$, respectively.
By the same proof as that for the well-definedness of $s_n(0,0)$ again, 
$s_n(1,0)'$ and $\iota$ are well-defined (\cite[(4.9)]{hk}). 
By the commutative diagram in \cite[(7.18.3)]{ndw},  
this morphism induces a morphism of complexes 
$\Om^{\bul}_{{\cal W}_n(Y)/{\cal W}_n(s), [~]}
\lo {\cal W}_n{\Om}^{\bul}_Y$.

\begin{rema} 
By the definition of the crystalline complex $C_{Y/{\cal W}_n(s)}$ 
in \cite[p.~ 238]{hk}, $C_{Y/{\cal W}_n(s)}$ is a cosimplicial complex
and ${\cal H}^i(C_{Y/{\cal W}_n(s)})\not={\cal W}_n\Om^{\bul}_{Y}$. 
The morphisms $\tau$, $\del$ and $d\log$ in \cite[(4.9)]{hk} are not 
correct morphisms. 
\end{rema} 


\begin{prop}\label{prop:nmc} 
The composite morphism 
\begin{align*} 
{\Om}^{\bul}_{{\cal W}_n(Y)/{\cal W}_n(s), [~]}
\lo {\cal W}_n{\Om}^{\bul}_Y\os{(C^{-n})^{-1}}{\lo} ({\cal W}_n{\Om}^{\bul}_Y)',  
\tag{2.1.4.1}\label{eqn:lnyy} 
\end{align*} 
which is a morphism of dga's over $g^{-1}({\cal O}_{{\cal W}_n(s)})$,   
is characterized by the following$:$ 
\begin{align*} 
&a\lom a \quad (a\in {\cal W}_n({\cal O}_Y)), \quad 
d\log m \lom (1\otimes m,0) \quad (m\in M_Y), \\
&du\lom (0,u\otimes 1) 
\quad (u\in {\cal W}_n({\cal O}_Y)^*). 
\end{align*}  
\end{prop}
\begin{proof} 
The characterization is clear because 
the images of $a\in {\Om}^0_{{\cal W}_n(Y)/{\cal W}_n(s), [~]}$ 
and $a\in ({\cal W}_n{\Om}^0_Y)'$, 
$d\log m\in \Om^1_{{\cal W}_n(Y)/{\cal W}_n(s), [~]}$ 
and $(1\otimes m,0)\in ({\cal W}_n{\Om}^1_Y)'$,  
$du\in \Om^1_{{\cal W}_n(Y)/{\cal W}_n(s), [~]}$ 
and $(0,u\otimes 1)\in ({\cal W}_n{\Om}^1_Y)'$) in 
${\cal W}_n{\Om}^0_Y$ and 
${\cal W}_n{\Om}^1_Y$ are the same, respectively.  
\end{proof}

\par 
In this book, let us denote by $R$ the projection 
$\pi \col {\cal W}_{n+1}{\Om}^{\bul}_Y
\lo {\cal W}_n{\Om}^{\bul}_Y$ in \cite[(4.2)]{hk}.  
In this book,  based on Hyodo's idea in \cite[\S1]{hdw} and our idea in \cite{ndw}, 
we give a different construction of $R$ from $\pi$ in \cite[(4.2)]{hk}. 
\par
Let $A$ be a commutative ring with unit element.
Consider an element $a\in A$.
Let $n$ and $m$ be nonnegative integers.
Let $L$ be an $A/a^{n+m}$-module.
Set $L/a^j:=L/a^jL$ $(0\leq j\leq n+m)$. 
Let $\ul{a}^m\col L/a^n\lo L$ be the induced morphism by 
the morphism $a^m\cdot \col L\lo L$. 
Assume that $\ul{a}^m\col L/a^n\lo L$ is injective. 
(If $a$ is a nonzero divisor of $A$ and if 
$L$ is a flat $A/a^{n+m}$-module, then this injectivity holds 
since the morphism $a^m\col A/a^n \lo A/a^{n+m}$ is injective.) 
Let $f\col L\lo L$ be
an endomorphism of $A$-modules such that
${\rm Im}(f)\subset {\rm Im}(\ul{a}^m)$. Then 
we can define a unique endomorphism $a^{-m}f\col L/a^n\lo L/a^n$ such that
$\ul{a}^m\circ a^{-m}f\circ {\rm pr}.=f$. Here ${\rm pr}.$ is the natural projection
$L\lo L/a^n$. 
\par
Let $n$ be a positive integer and let $i$ be a nonnegative integer. 
To define $R\col {\cal W}_n{\Om}_Y^i\lo {\cal W}_{n-1}{\Om}_Y^i$,
we first assume that there exists an immersion
$Y\os{\sus}{\lo} {\cal Q}_{n+i}$ into a formally log smooth {\it integral} scheme over
${\cal W}_{n+i}(s)$ such that ${\cal Q}_{n+i}$ has a lift 
$\varphi_{n+i} \col {\cal Q}_{n+i}\lo {\cal Q}_{n+i}$
of the Frobenius endomorphism of 
${\cal Q}_{n+i}\times_{{\cal W}_{n+i}(s)}s$.
Let ${\mathfrak E}_{n+i}$ be the log PD-envelope of the immersion
$Y\os{\sus}{\lo} {\cal Q}_{n+i}$ over $({\cal W}_{n+i}(s),p{\cal W}_{n+i},[~])$.
Set ${\cal Q}_j:={\cal Q}_{n+i}\times_{{\cal W}_{n+i}(s)}{\cal W}_j(s)$ and 
${\mathfrak E}_j:={\mathfrak E}_{n+i}\times_{{\cal W}_{n+i}(s)}{\cal W}_j(s)$
$(0\leq j\leq n+i)$.
Consider the following morphism
\begin{align*}
\varphi^*_{n+i}
\col & {\cal O}_{{\mathfrak E}_{n+i}}\otimes_{{\cal O}_{{\cal Q}_{n+i}}}
\Om^i_{{\cal Q}_{n+i}/{\cal W}_{n+i}(s)}\lo
{\cal O}_{{\mathfrak E}_{n+i}}\otimes_{{\cal O}_{{\cal Q}_{n+i}}}
\Om^i_{{\cal Q}_{n+i}/{\cal W}_{n+i}(s)} \quad (i\in {\mab N}).
\end{align*}
Because ${\cal O}_{{\mathfrak E}_{n+i}}\otimes_{{\cal O}_{{\cal Q}_{n+i}}}
\Om^i_{{\cal Q}_{n+i}/{\cal W}_{n+i}(s)}$ is a sheaf of flat
${\cal O}_{{\cal W}_{n+i}(s)}$-modules by (\ref{coro:filt}) 
and because $\varphi^*_{n+i}$ is divisible by $p^i$
(since $\varphi_{n+i}$ is a lift of the  Frobenius endomorphism),
the morphism
\begin{align*}
p^{-(i-1)}\varphi^*_{n+i}
\col & {\cal O}_{{\mathfrak E}_{n+1}}\otimes_{{\cal O}_{{\cal Q}_{n+1}}}
\Om^j_{{\cal Q}_{n+1}/{\cal W}_{n+1}(s)}\lo
{\cal O}_{{\mathfrak E}_{n+1}}\otimes_{{\cal O}_{{\cal Q}_{n+1}}}
\Om^j_{{\cal Q}_{n+1}/{\cal W}_{n+1}(s)} \quad (j\geq i-1)
\end{align*}
is well-defined.
Because the image of the morphism above is contained in
$p{\cal O}_{{\mathfrak E}_{n+1}}\otimes_{{\cal O}_{{\cal Q}_{n+1}}}
\Om^j_{{\cal Q}_{n+1}/{\cal W}_{n+1}(s)}$ for the case $j=i$,
we have the following well-defined morphism (cf.~\cite[Editorial comment (5)]{hdw})
\begin{align*}
p^{-(i-1)}\varphi^*_{n+i} &\col {\cal W}_n{\Om}_{Y}^i
={\cal H}^i({\cal O}_{{\mathfrak E}_n}\otimes_{{\cal O}_{{\cal Q}_n}}
\Om^{\bul}_{{\cal Q}_n/{\cal W}_n(s)})\lo\tag{2.1.4.2}\label{ali:enqw}\\
& 
{\cal H}^i({\cal O}_{{\mathfrak E}_{n+1}}\otimes_{{\cal O}_{{\cal Q}_{n+1}}}
\Om^{\bul}_{{\cal Q}_{n+1}/{\cal W}_{n+1}(s)})
={\cal W}_{n+1}{\Om}_{Y}^i \quad (i\in {\mab N}).
\end{align*}

\begin{rema}\label{rema:infp}
(1) Assume that there exists an immersion 
$Y\os{\sus}{\lo} {\cal Q}$ 
into a formally log smooth {\it integral} scheme over ${\cal W}(s)$ 
such that ${\cal Q}$ has a lift $\varphi \col {\cal Q}\lo {\cal Q}$ 
of the Frobenius endomorphism of 
${\cal Q}\times_{{\cal W}(s)}s$. 
Let ${\mathfrak E}$ be the log PD-envelope of the immersion 
$Y\os{\sus}{\lo} {\cal Q}$ over 
$({\cal W}(s),p{\cal W},[~])$.  
Set ${\cal Q}_n:={\cal Q}\times_{{\cal W}(s)}{\cal W}_n(s)$ and   
${\mathfrak E}_n:={\mathfrak E}\times_{{\cal W}(s)}{\cal W}_n(s)$. 
Consider the following morphism 
\begin{align*} 
\varphi^* \col {\cal O}_{\mathfrak E}
\otimes_{{\cal O}_{\cal Q}}
\Om^i_{{\cal Q}/{\cal W}(s)}
\lo 
{\cal O}_{\mathfrak E}\otimes_{{\cal O}_{\cal Q}}
\Om^i_{{\cal Q}/{\cal W}(s)} \quad (i\in {\mab N}). 
\end{align*} 
Then the morphism 
\begin{align*} 
p^{-(i-1)}\varphi^*
\col {\cal O}_{\mathfrak E}\otimes_{{\cal O}_{\cal Q}}
\Om^j_{{\cal Q}/{\cal W}(s)}
\lo 
{\cal O}_{\mathfrak E}\otimes_{{\cal O}_{\cal Q}}
\Om^j_{{\cal Q}/{\cal W}(s)} \quad (j\geq i-1)
\end{align*} 
is well-defined. 
Hence we have the following morphism for $j\geq i-1$: 
\begin{align*} 
p^{-(i-1)}\varphi^* 
\col 
{\cal O}_{{\mathfrak E}_{n+1}}
\otimes_{{\cal O}_{{\cal Q}_{n+1}}}
\Om^j_{{\cal Q}_{n+1}/{\cal W}_{n+1}(s)}
\lo 
{\cal O}_{{\mathfrak E}_{n+1}}\otimes_{{\cal O}_{{\cal Q}_{n+1}}}
\Om^j_{{\cal Q}_{n+1}/{\cal W}_{n+1}(s)}.  
\end{align*} 
As above, we have the following well-defined morphism 
\begin{align*} 
p^{-(i-1)} \varphi^* \col 
& {\cal W}_n{\Om}_{Y}^i
={\cal H}^i({\cal O}_{{\mathfrak E}_{n}}\otimes_{{\cal O}_{{\cal Q}_n}}
\Om^{\bul}_{{\cal Q}_n/{\cal W}_n(s)})
\tag{2.1.5.1}\label{ali:qnw}\\
& \lo 
{\cal H}^i({\cal O}_{{\mathfrak E}_{n+1}}\otimes_{{\cal O}_{{\cal Q}_{n+1}}}
\Om^{\bul}_{{\cal Q}_{n+1}/{\cal W}_{n+1}(s)})
={\cal W}_{n+1}{\Om}_{Y}^i \quad (i\in {\mab N}). 
\end{align*} 
This morphism is equal to (\ref{ali:enqw}). 
\par 
(2) Let the notations be as in \cite[(6.27)]{ndw}. 
Then we have to assume that the structural morphism 
${\cal Y}\lo ({\rm Spf}(W),W(L))$ is assumed to be integral 
because we need the flatness of the structural morphism 
$\os{\circ}{\cal Y}\lo {\rm Spf}(W)$; the condition of the flatness 
is necessary for the condition \cite[(6.0.2)]{ndw}.
\end{rema}

\begin{prop}\label{prop:nqop} 
$(1)$ The morphism $p^{-(i-1)}\varphi^*$ in {\rm (\ref{ali:enqw})} 
is independent of the choice of 
the immersion $Y\os{\sus}{\lo}{\cal Q}$ and 
the lift $\varphi \col {\cal Q}\lo {\cal Q}$ 
of the Frobenius endomorphism of 
${\cal Q}\times_{{\cal W}(s)}s$.
We set ${\bf p}:=p^{-(i-1)}\varphi^*\col 
{\cal W}_n{\Om}_{Y}^i\lo {\cal W}_{n+1}{\Om}_{Y}^i$. 
\par 
$(2)$ ${\rm Ker}({\bf p}\col {\cal W}_n{\Om}_{Y}^i\lo 
{\cal W}_{n+1}\Om_Y^i)=0$.
\par 
$(3)$ ${\rm Im}({\bf p}\col {\cal W}_n{\Om}_{Y}^i\lo 
{\cal W}_{n+1}{\Om}_{Y}^i)=
{\rm Im}(p\col {\cal W}_{n+1}{\Om}_{Y}^i\lo 
{\cal W}_{n+1}{\Om}_{Y}^i)$. 
\end{prop}
\begin{proof}
(1): (This is a routine work in our formalism because 
we do not consider only a lift but also an immersion.) 
Let $Y\os{\sus}{\lo}{\cal Q}'$ 
be another immersion and another
lift $\varphi' \col {\cal Q}'\lo {\cal Q}'$ 
of the Frobenius endomorphism of ${\cal Q}'\times_{{\cal W}(s)}s$.
Then, considering the products 
${\cal Q}\times_{{\cal W}_n(s)}{\cal Q}'$
and $\varphi \times \varphi'$, we may assume that 
there exists the following commutative diagram  
\begin{equation*} 
\begin{CD} 
Y@>{\subset}>>{\cal Q}@>{\varphi}>>{\cal Q}\\
@| @VVV @VVV \\
Y@>{\subset}>>{\cal Q}'@>{\varphi'}>>{\cal Q}'. 
\end{CD}
\end{equation*} 
Hence we have the following commutative diagram 
\begin{equation*} 
\begin{CD} 
{\cal H}^i({\cal O}_{{\mathfrak E}_n}\otimes_{{\cal O}_{{\cal Q}_n}}
\Om^{\bul}_{{\cal Q}_n/{\cal W}_n(s)})@>{p^{-(i-1)}\varphi^*}>>
{\cal H}^i({\cal O}_{{\mathfrak E}_{n+1}}
\otimes_{{\cal O}_{{\cal Q}_{n+1}}}
\Om^{\bul}_{{\cal Q}_{n+1}/{\cal W}_{n+1}(s)})\\
@A{\simeq}AA @AA{\simeq}A  \\
{\cal H}^i({\cal O}_{{\mathfrak E}'_{n}}\otimes_{{\cal O}_{{\cal Q}'_n}}
\Om^{\bul}_{{\cal Q}'_n/{\cal W}_n(s)})@>{p^{-(i-1)}\varphi'{}^*}>>
{\cal H}^i({\cal O}_{{\mathfrak E}'_{n+1}}
\otimes_{{\cal O}_{{\cal Q}'_{n+1}}}
\Om^{\bul}_{{\cal Q}'_{n+1}/{\cal W}_{n+1}(s)}),
\end{CD}
\end{equation*} 
which shows the desired independence. 
\par 
(2): We may assume that ${\cal Q}$ is a 
formal lift ${\cal Y}/{\cal W}(s)$ of $Y/s$. 
Then 
$(\Om^{\bul},\phi):=(\Om^{\bul}_{{\cal Y}/{\cal W}(s)},\varphi^*)$ 
satisfies the conditions $(6.0.1)\sim (6.0.5)$ in \cite{ndw}
(these conditions are only the abstract versions of 
necessary conditions for \cite[(8.8)]{bob}).  
That is, the following hold:
\medskip
\parno
(2.1.6.1) $\Om^i=0$ for $i<0$.
\medskip
\parno
(2.1.6.2)  $\Om^i$ 
$(\forall i \in {\mab N})$ are sheaves of $p$-torsion-free, $p$-adically 
complete ${\mab Z}_p$-modules.
\medskip
\parno
(2.1.6.3) $\phi(\Om^{i}) \subset 
\{\om \in p^i \Om^i~ \vert~ d\om \in p^{i+1} 
\Om^{i+1}\}$ $(\forall i\in {\mab N})$.
\medskip
\parno
(2.1.6.4) Set $\Om_1^{\bul}:=\Om^{\bul}/p\Om^{\bul}$. 
Then there exists an ${\mab F}_p$-linear isomorphism 
$$C^{-1} \col \Om^i_1\os{\sim}{\lo} {\cal H}^i(\Om^{\bul}_1) \quad 
(\forall i \in {\mab N}).$$
\medskip
\parno
(2.1.6.5) A composite morphism 
(${\rm mod}~p)\circ p^{-i}\phi 
\col \Om^{i} \lo \Om^i \lo \Om^i_1$ factors 
through ${\rm Ker}(d \col\Om^i_1 
\lo \Om^{i+1}_{1})$, and the 
following diagram 
is commutative:
\begin{equation*}
\begin{CD}
\Om^i @>{\mod p}>> 
\Om^i_1\\ 
@V{p^{-i}\phi}VV  @VV{C^{-1}}V \\
\Om^i @>{\mod p}>> 
{\cal H}^i(\Om^{\bul}_1).
\end{CD}
\end{equation*}
Set 
\begin{equation*}
Z^i_n:=\{\om \in \Om^i \vert~ 
d\om \in p^n\Om^{i+1}\}, \quad
B^i_n:=p^n\Om^i+ d\Om^{i-1}, \quad 
{\mathfrak W}_n{\Om}^i= Z^i_n/B^i_n. 
\tag{2.1.6.6}\label{eqn:pnzb}
\end{equation*}
Then ${\cal W}_n{\Om}^i_{Y}= 
{\mathfrak W}_n{\Om}^i$ and 
the morphism 
${\bf p}\col {\cal W}_n{\Om}^i_{Y}\lo 
{\cal W}_{n+1}{\Om}^i_{Y}$ 
is equal to ${\bf p}\col {\mathfrak W}_n{\Om}^i\lo 
{\mathfrak W}_{n+1}{\Om}^i$ in \cite[p.~546]{ndw}. 
Because the latter morphism is injective (\cite[(6.8)]{ndw}), we obtain (2). 
\par 
(3): By \cite[p.~546]{ndw},  
${\rm Im}({\bf p}\col {\mathfrak W}_n{\Om}^i\lo {\mathfrak W}_{n+1}{\Om}^i)
={\rm Im}(p\col {\mathfrak W}_{n+1}{\Om}^i\lo {\mathfrak W}_{n+1}{\Om}^i)$. 
Hence we obtain (3). 
\end{proof}

\begin{rema}\label{rema:csh}
In \cite{hdw} Hyodo has not proved 
the well-definedness of the morphism 
${\bf p} \col {\cal W}_n\Om^i_{Y}\lo {\cal W}_{n+1}\Om^i_{Y}$ 
(The editorial comment (5) in [loc.~cit.] is useless for the proof of the well-definedness).
He has not proved the well-definedness of the projection  
$\pi \col {\cal W}_{n+1}\Om^i_{Y}\lo {\cal W}_n\Om^i_{Y}$ either. 
\end{rema} 

\begin{defi}[{\bf cf.~\cite[(1.3.2)]{hdw}}]\label{defi:hdyd} 
(1) We call the morphism 
${\bf p} \col {\cal W}_n\Om^i_{Y}\lo 
{\cal W}_{n+1}\Om^i_{Y}$ {\it Hyodo's multiplication by} $p$ on 
${\cal W}_{\bul}\Om^i_{Y}$. 
\par 
(2) The morphism $R\col {\cal W}_{n+1}{\Om}^i_{Y}\lo 
{\cal W}_n{\Om}^i_{Y}$ 
is, by definition, the unique morphism fitting into the following commutative diagram: 
\begin{equation*} 
\begin{CD} 
{\cal W}_{n+1}{\Om}^i_{Y}
@>{R}>>{\cal W}_n{\Om}^i_{Y}\\
@V{p}VV @VV{\bf p}V\\ 
{\cal W}_{n+1}{\Om}^i_{Y}
@={\cal W}_{n+1}{\Om}^i_{Y},  
\end{CD}
\end{equation*} 
whose existence is assured by (\ref{prop:nqop}). 
\end{defi}

\begin{rema}\label{rema:sot} 
(1) Our definition is a simpler one than 
the definition of $\pi\col  
{\cal W}_{n+1}\Om^{\bul}_{Y}\lo {\cal W}_n\Om^{\bul}_{Y}$ 
in \cite[(4.2)]{hk} 
because we do not use the crystalline complexes 
(see \cite[(2.25)]{hk} (cf.~\cite[(4.1.6), (4.1.7)]{ols})). 
\par 
(2) Note that the proof of \cite[(2.24)]{hk} is mistaken because 
$\varphi \col Ru^{\log}_{X'/T_n*}({\cal O}_{X'/T_n*}({\cal O}_{X'/T_n})
\lo Ru^{\log}_{X'/T_n*}({\cal O}_{X/T_n*}({\cal O}_{X/T_n})$ in the statement 
of [loc.~cit.] does {\it not} produce the morphism 
$\varphi \col C^q_{X'/T_n}\lo C^q_{X/T_n}$ in [loc.~cit., p.~242]; 
in [loc.~cit.] Hyodo and Kato have not considered 
local lifts of relative Frobenius morphism 
$g\col X\lo X'$ in [loc.~cit., (2.12.1)]. 
\end{rema}

\parno 
By the proof of (\ref{prop:nqop}) and \cite[(6.27)]{ndw},  
the morphism $R$ is equal to the following composite morphism 
in the case where $Y/s$ has a local lift ${\cal Y}/{\cal W}(s)$
(cf.~\cite[(4.2)]{hk}):
\begin{align*}
{\cal W}_{n+1}\wt{\Om}^i_{Y_{\os{\circ}{t}}}= 
Z^i_{n+1}/B^i_{n+1} 
\os{\us{\sim}{p^i}}{\lo} 
p^iZ^i_{n+1}/p^iB^i_{n+1}
& \os{{\rm proj}.}{\lo} 
p^iZ^i_{n+1}/(p^{i+n}Z^i_1+p^{i-1}dZ^{i-1}_1) 
\tag{2.1.9.1}\label{ali:aplqcpj}\\
{} & \os{\us{\sim}{\phi}}{\longleftarrow}
Z^i_n/B^i_n={\cal W}_n\wt{\Om}^i_{Y_{\os{\circ}{t}}}.
\end{align*}
By \cite[(6.5)]{ndw}, 
the morphism $R$ is equal to the 
following composite morphism$:$
\begin{align*}
&{\cal W}_{n+1}\wt{\Om}^i_{Y_{\os{\circ}{t}}}  = 
Z^i_{n+1}/B^i_{n+1}
\os{{\rm proj}.}{\lo} 
Z^i_{n+1}/(p^nZ^i_1+d\Om^{i-1}) 
\tag{2.1.9.2}\label{ali:alocpqj}\\
& 
\os{(p^{-i}\phi)^{-1}}{\os{\sim}{\lo}}
Z^i_n/(p^n\Om^i+pd\Om^{i-1})
\os{{\rm proj}.}{\lo}
Z^i_n/B^i_n={\cal W}_n\wt{\Om}^i_{Y_{\os{\circ}{t}}}.
\end{align*}
In particular, $R$ is surjective.

\par 
Let $R\col ({\cal W}_{n+1}{\Om}^{\bul}_Y)'
\lo ({\cal W}_n{\Om}^{\bul}_Y)'$ be also the projection 
induced by the projection 
${\cal W}_{n+1}({\cal O}_Y)'\lo {\cal W}_n({\cal O}_Y)'$.  
By \cite[(7.1)]{ndw} the two $R$'s are compatible with 
the Cartier isomorphism 
$C^{-m}\col ({\cal W}_m\Om^{\bul}_Y)'\os{\sim}{\lo} {\cal W}_m\Om^{\bul}_Y$ 
$(m=n,n+1)$.  
\par 
In the following we construct a canonical morphism (\ref{eqn:ywcnny}) below. 
This is the first aim in this section. 
\par 
By (\ref{eqn:intcon}) 
we have the following integrable connection 
\begin{equation*} 
E_n\lo E_n\otimes_{{\cal W}_n({\cal O}_Y)'}
{\Om}^1_{{\cal W}_n(Y)/{\cal W}_n(s),[~]}.  
\tag{2.1.9.3}\label{eqn:eefnny} 
\end{equation*} 
This gives the complex 
$E_n\otimes_{{\cal W}_n({\cal O}_Y)'}
{\Om}^{\bul}_{{\cal W}_n(Y)/{\cal W}_n(s),[~]}$ 
as in \cite[II (1.1.5)]{et}.  
By using the morphism (\ref{eqn:lmwy}),  
we have the connection 
\begin{equation*} 
E_n\lo E_n\otimes_{{\cal W}_n({\cal O}_Y)'}{\cal W}_n{\Om}^1_Y
\tag{2.1.9.4}\label{eqn:eewnny} 
\end{equation*} 
and the ``reverse'' log de Rham-Witt complex 
\begin{equation*}
E_n\otimes_{{\cal W}_n({\cal O}_Y)}
{\cal W}_n{\Om}^{\bul}_Y
\tag{2.1.9.5}\label{eqn:ldwy} 
\end{equation*}  
of $E$ of level $n$.  
By using the morphism (\ref{eqn:wnymod}), 
we also have the ``obverse'' log de Rham-Witt complex
\begin{equation*} 
E_n\otimes_{{\cal W}_n({\cal O}_Y)'}
({\cal W}_n{\Om}_Y^{\bul})'
\tag{2.1.9.6}\label{eqn:lmdy} 
\end{equation*} 
of $E$ of level $n$. 
Let $\star$ be nothing or $\prime$. 
We have the following morphism 
\begin{equation*} 
R\col E_{n+1}\otimes_{{\cal W}_{n+1}({\cal O}_Y)^{\star}}
({\cal W}_{n+1}{\Om}^i_Y)^{\star}
\owns x\otimes y \lom R(x)\otimes R(y) \in 
E_n\otimes_{{\cal W}_n({\cal O}_Y)^{\star}}
({\cal W}_n{\Om}^i_Y)^{\star} 
\quad (i\in {\mab N}). 
\tag{2.1.9.7}\label{eqn:reny}
\end{equation*}
Because $R\col {\cal W}_{n+1}{\Om}^{\bul}_Y\lo {\cal W}_n{\Om}^{\bul}_Y$ 
is a morphism of complexes (\cite[(6.8), (6.27)]{ndw}) 
(it is obvious that 
$R\col ({\cal W}_{n+1}{\Om}^{\bul}_Y)'\lo ({\cal W}_n{\Om}^{\bul}_Y)'$ 
is a morphism of complexes), 
the morphism $R$ in (\ref{eqn:reny}) is a morphism of complexes. 
As usual, for $0\leq r\leq n$, 
set 
\begin{align*} 
{\rm Fil}^r(E_n\otimes_{{\cal W}_n({\cal O}_Y)^{\star}}
({\cal W}_n{\Om}^{\bul}_Y)^{\star}) := 
{\rm Ker}(R^{n-r}\col & E_n\otimes_{{\cal W}_n({\cal O}_Y)^{\star}}
({\cal W}_n{\Om}^{\bul}_Y)^{\star} 
\tag{2.1.9.8}\label{ali:fln}\\
&\lo 
E_r\otimes_{{\cal W}_r({\cal O}_Y)^{\star}}
({\cal W}_r{\Om}^{\bul}_Y)^{\star}). 
\end{align*} 
We call ${\rm Fil}:=\{{\rm Fil}^r\}_{r\in {\mab Z}}$ 
the {\it canonical filtration} on  
$E_n\otimes_{{\cal W}_n({\cal O}_Y)^{\star}}({\cal W}_n{\Om}^{\bul}_Y)^{\star}$ 
as in \cite[I (3.1.1)]{idw}. 
The multiplication of $p^r$ $(r\in {\mab N})$ on 
$E_{n+r}\otimes_{{\cal W}_{n+r}({\cal O}_Y)}
{\cal W}_{n+r}{\Om}_Y^{\bul}$ induces the following morphism of complexes
(\cite[(4.5)]{hk}): 
\begin{equation*} 
{\bf p}^r\col 
E_n\otimes_{{\cal W}_n({\cal O}_Y)}
{\cal W}_n{\Om}_Y^{\bul}
\lo 
E_{n+r}\otimes_{{\cal W}_{n+r}({\cal O}_Y)}
{\cal W}_{n+r}{\Om}_Y^{\bul}.  
\tag{2.1.9.9}\label{eqn:preeny} 
\end{equation*} 
The multiplication of $p^r$ $(r\in {\mab N})$ on 
$E_{n+r}\otimes_{{\cal W}_{n+r}({\cal O}_Y)'}
({\cal W}_{n+r}{\Om}_Y^{\bul})'$
also induces the following morphism of complexes 
\begin{equation*} 
{\bf p}^r\col 
E_n\otimes_{{\cal W}_n({\cal O}_Y)'}({\cal W}_n{\Om}_Y^{\bul})'
\lo 
E_{n+r}\otimes_{{\cal W}_{n+r}({\cal O}_Y)'}
({\cal W}_{n+r}{\Om}_Y^{\bul})'.  
\tag{2.1.9.10}\label{eqn:prewny} 
\end{equation*} 
The morphism (\ref{eqn:prewny}) is equal to the morphism 
(\ref{eqn:preeny}) via ${\rm id}_{E_i}\otimes C^{-i}$ $(i=n+r, n)$ 
because $C^{-i}$ is compatible with the two projections. 
If $E$ is a flat quasi-coherent ${\cal O}_{Y/{\cal W}_n(s)}$-modules, 
then the morphism (\ref{eqn:preeny}) and hence the morphism 
(\ref{eqn:prewny}) are injective 
by \cite[(4.5) (1)]{hk} or \cite[(6.8) (2)]{ndw}
(the case $E_n={\cal W}_n({\cal O}_Y)$).

\par 
Now let $N$ be a nonnegative integer or $\infty$. 
Let $Y_{\bul \leq N}$ be 
a log smooth $N$-truncated 
simplicial log scheme of Cartier type over $s$. 
Let $g \col Y_{\bul \leq N} \lo s \os{\sus}{\lo} {\cal W}_n(s)$ be 
the structural morphism.  
Let $E^{\bul \leq N}$ be a flat coherent log crystal of 
${\cal O}_{Y_{\bul \leq N}/{\cal W}_n(s)}$-modules. 
Set $E^{\bul \leq N}_n:=
(E^{\bul \leq N})_{{\cal W}_n(Y_{\bul \leq N})}$. 
Then we have the 
$N$-truncated cosimplicial log de Rham-Witt complex
$E^{\bul \leq N}_n
\otimes_{{\cal W}_n({\cal O}_{Y_{\bul \leq N}})^{\star}}
({\cal W}_n{\Om}^{\bul}_{Y_{\bul \leq N}})^{\star}$. 
Assume that $N\not= \infty$ and that $Y_{\bul \leq N}$ has 
an affine $N$-truncated simplicial open covering of $Y_{\bul \leq N}$. 
As mentioned before, we would like to construct a canonical morphism 
\begin{equation*} 
Ru_{Y_{\bul \leq N}/{\cal W}_n(s)*}
(E^{\bul \leq N}) 
\lo 
E^{\bul \leq N}_n
\otimes_{{\cal W}_n({\cal O}_{Y_{\bul \leq N}})^{\star}}
({\cal W}_n{\Om}^{\bul}_{Y_{\bul \leq N}})^{\star} 
\tag{2.1.9.11}\label{eqn:ywcnny}
\end{equation*} 
in ${\rm D}^+(g^{-1}({\cal W}_n))$ 
and to prove that this is an isomorphism. 
To prove this, we need the following 
which is a variant of a coefficient version of \cite[(4.22) (1)]{nh3}:

\begin{lemm}[{\bf cf.~the proof of \cite[(4.19)]{hk}}]\label{lemm:neh3fi}   
Let $Y$ $($resp.~${\cal Q})$ 
be a log smooth scheme of Cartier type over $s$ 
$($resp.~a log smooth scheme over ${\cal W}_n(s))$.   
Let $Y \os{\sus}{\lo} {\cal Q}$ 
be an immersion over ${\cal W}_n(s)$ and 
let ${\mathfrak E}$ be the log PD-envelope of this immersion 
over $({\cal W}_n(s),p{\cal W}_n,[~])$. 
Let $g \col Y \lo {\cal W}_n(s)$ be the structural morphism.  
Let $E$ be a quasi-coherent crystal of ${\cal O}_{Y/{\cal W}_n(s)}$-modules 
and let $({\cal E},\nabla)$ be the corresponding quasi-coherent 
${\cal O}_{\mathfrak E}$-module with integrable connection to $E$. 
Assume that $\os{\circ}{Y}$ is affine. 
Then the following hold$:$
\par 
$(1)$ There exists a morphism  
${\cal W}_n(Y)\lo {\mathfrak E}$ of log PD-schemes 
over ${\cal W}_n(s)$ such that the composite morphism 
$Y \os{\sus}{\lo} {\cal W}_n(Y)\lo {\mathfrak E}$ 
is the immersion obtained by the given immersion 
$Y\os{\sus}{\lo} {\cal Q}$. 
\par 
$(2)$  
Let ${\cal W}_n(Y)\lo  {\mathfrak E}$ be the morphism in $(1)$. 
Then this morphism induces the following morphism 
\begin{equation*} 
{\cal E}\otimes_{{\cal O}_{\cal Q}}\Om^{\bul}_{{\cal Q}/{\cal W}_n(s)} 
\lo 
E_n\otimes_{{\cal W}_n({\cal O}_Y)}{\cal H}^{\bul}
({\cal O}_{\mathfrak E}\otimes_{{\cal O}_{\cal Q}}
\Om^*_{{\cal Q}/{\cal W}_n(s)})=
E_n\otimes_{{\cal W}_n({\cal O}_Y)}{\cal W}_n\Om^{\bul}_Y 
\tag{2.1.10.1}\label{eqn:fnih3dw} 
\end{equation*} 
in ${\rm C}^+(\os{\circ}{g}{}^{-1}({\cal W}_n))$ 
which will turn out to be a quasi-isomorphism. 
This morphism is functorial in the following sense$:$
for a log smooth log affine scheme $Y_i$ 
of Cartier type over a fine log scheme $s_i$ $(i=1,2)$ 
whose underlying scheme 
is the spectrum of a perfect field of characteristic $p>0$ 
and for an immersion $Y_i\os{\sus}{\lo}{\cal Q}_i$ 
into a log smooth scheme over ${\cal W}_n(s_i)$ 
and for the log PD-envelope ${\mathfrak E}_i$ $(i=1,2)$ of 
the immersion $Y_i \os{\sus}{\lo} {\cal Q}_i$ 
over $({\cal W}_n(s_i),p{\cal W}_n,[~])$ 
and for the following commutative diagram 
\begin{equation*} 
\begin{CD} 
{\cal W}_n(Y_1)@>>>  {\mathfrak E}_1 \\ 
@V{h_n}VV  @VV{h}V \\  
{\cal W}_n(Y_2) @>>> {\mathfrak E}_2
\end{CD} 
\tag{2.1.10.2}\label{cd:fnih2dw} 
\end{equation*} 
of log PD-schemes over the morphism 
${\cal W}_n(s_1) \lo {\cal W}_n(s_2)$ 
such that the composite morphism 
$Y_i\os{\sus}{\lo} {\cal W}_n(Y_i)\lo {\mathfrak E}_i$ $(i=1,2)$ 
is the given immersion and for a morphism $h_{1{\rm crys}}^*(E_2)\lo E_1$ of 
${\cal O}_{Y_1/{\cal W}_n(s_1)}$-modules 
$(E_i${\rm :} a quasi-coherent crystal of 
${\cal O}_{Y_i/{\cal W}_n(s_i)}$-modules$)$, 
the following diagram is commutative$:$ 
\begin{equation*}  
\begin{CD} 
{\cal E}_1{\otimes}_{{\cal O}_{{\cal Q}_1}}\Om^{\bul}_{{\cal Q}_1/{\cal W}_n(s)}
@>>> 
E_{1,n}{\otimes}_{{\cal W}_n({\cal O}_{Y_1})}{\cal H}^{\bul}({\cal O}_{{\mathfrak E}_1} 
{\otimes}_{{\cal O}_{{\cal Q}_1}}
\Om^{*}_{{\cal Q}_1/{\cal W}_n(s)}) \\ 
@A{h^*}AA @AA{h^*_n}A \\ 
h^*({\cal E}_2{\otimes}_{{\cal O}_{{\cal Q}_2}}
\Om^{\bul}_{{\cal Q}_2/{\cal W}_n(s)} )
@>>> 
h_n^*(E_{2,n}{\otimes}_{{\cal W}_n({\cal O}_{Y_2})}
{\cal H}^{\bul}({\cal O}_{{\mathfrak E}_2} 
{\otimes}_{{\cal O}_{{\cal Q}_2}}
\Om^{*}_{{\cal Q}_2/{\cal W}_n(s)} )),   
\end{CD} 
\tag{2.1.10.3}\label{cd:odh3zqph}
\end{equation*} 
where ${\cal E}_i:=E_{{\mathfrak E}_i}$. 
The morphism {\rm (\ref{eqn:fnih3dw})} is compatible 
with the projections. 
\par 
$(3)$ Let $s_i$ and the immersion $Y_i\os{\sus}{\lo}{\cal Q}_i$ 
$(i=1,2)$ be as in $(2)$. 
Let ${\mathfrak E}_i$ be the log PD-envelope  of 
the immersion $Y_i \os{\sus}{\lo} {\cal Q}_i$ 
over $({\cal W}_n(s_i),p{\cal W}_n,[~])$.  
Let ${\cal W}_n(Y_i)\lo {\cal Q}_i$ be the morphism obtained in the proof of $(1)$. 
Assume that there exists the following commutative diagram 
\begin{equation*} 
\begin{CD} 
{\cal W}_n(Y_1)@>>>  {\cal Q}_1 \\ 
@V{h_n}VV  @VV{h}V \\  
{\cal W}_n(Y_2) @>>> {\cal Q}_2, 
\end{CD} 
\tag{2.1.10.4}\label{cd:fnih0dw} 
\end{equation*}  
where $h_n$ is a morphism of log PD-schemes. 
By abuse of notation, denote by $h\col {\mathfrak E}_1\lo {\mathfrak E}_2$ 
the induced morphism by $h\col {\cal Q}_1\lo {\cal Q}_2$. 
Then the diagram 
\begin{equation*} 
\begin{CD} 
{\cal W}_n(Y_1)@>>>  {\mathfrak E}_1 \\ 
@V{h_n}VV  @VV{h}V \\  
{\cal W}_n(Y_2) @>>> {\mathfrak E}_2
\end{CD} 
\tag{2.1.10.5}\label{cd:fnih1dw} 
\end{equation*} 
is commutative. 
\end{lemm} 
\begin{proof} 
(1): Because $\os{\circ}{Y}$ is affine, 
because ${\cal Q}$ is log smooth over ${\cal W}_n(s)$ and 
because the closed immersion $Y\os{\sus}{\lo}{\cal W}_n(Y)$ is nilpotent 
($(V{\cal W}_{n-1}({\cal O}_Y))^n=0$ (\cite[p.~589]{idw})), 
we have a morphism ${\cal W}_n(Y)\lo {\cal Q}$ over ${\cal W}_n(s)$ 
such that the composite morphism $Y\os{\sus}{\lo} {\cal W}_n(Y)\lo {\cal Q}$ 
is the given immersion $Y\os{\sus}{\lo} {\cal Q}$. 
Because the immersion $Y\os{\sus}{\lo}{\cal W}_n(Y)$ is a log PD-immersion
over $({\cal W}_n(s),p{\cal W}_n,[~])$,  
we have a natural morphism 
${\cal W}_n(Y)  \lo {\mathfrak E}$ of log PD-schemes over ${\cal W}_n(s)$ 
by the universality of the log PD-envelope.  
\par 
(2): By (\ref{eqn:fwniwu}) we have the following isomorphism 
\begin{align*} 
{\cal E}\otimes_{{\cal O}_{\cal Q}}
\Om^i_{{\cal Q}/{\cal W}_n(s)}  \os{\sim}{\lo} 
{\cal E}\otimes_{{\cal O}_{\mathfrak E}}
\Om^i_{{\mathfrak E}/{\cal W}_n(s),[~]}  \quad (i\in {\mab N}). 
\tag{2.1.10.6}\label{ali:acdw}
\end{align*} 
The morphism ${\cal W}_n(Y)\lo {\mathfrak E}$ gives 
the following morphism 
\begin{align*} 
{\cal E}\otimes_{{\cal O}_{\mathfrak E}}
\Om^i_{{\mathfrak E}/{\cal W}_n(s),[~]}  
{\lo} E_n\otimes_{{\cal W}_n({\cal O}_Y)'}
\Om^i_{{\cal W}_n(Y)/{\cal W}_n(s),[~]}. 
\tag{2.1.10.7}\label{ali:wsn}
\end{align*} 
By (\ref{eqn:lay}) we have the following morphism  
\begin{align*} 
E_n\otimes_{{\cal W}_n({\cal O}_Y)'}
\Om^i_{{\cal W}_n(Y)/{\cal W}_n(s),[~]}
\lo 
E_n\otimes_{{\cal W}_n({\cal O}_Y)}
{\cal H}^i({\cal O}_{\mathfrak E}\otimes_{{\cal O}_{\cal Q}}
\Om^*_{{\cal Q}/{\cal W}_n(s)}).
\tag{2.1.10.8}\label{ali:sm}
\end{align*}  
By (\ref{ali:acdw}), (\ref{ali:wsn}) and (\ref{ali:sm}), 
we obtain the following composite morphism:  
\begin{align*} 
{\cal E}\otimes_{{\cal O}_{\cal Q}}
\Om^i_{{\cal Q}/{\cal W}_n(s)}  
&\lo E_n\otimes_{{\cal W}_n({\cal O}_Y)'}
\Om^i_{{\cal W}_n(Y)/{\cal W}_n(s),[~]} 
\lo 
E_n\otimes_{{\cal W}_n({\cal O}_Y)}
{\cal H}^i({\cal O}_{\mathfrak E}\otimes_{{\cal O}_{\cal Q}}
\Om^*_{{\cal Q}/{\cal W}_n(s)} ) \tag{2.1.10.9}\label{eqn:fcnidw}.   
\end{align*} 
Because the morphism 
$\Om^{\bul}_{{\mathfrak E}/{\cal W}_n(s),[~]}  \lo 
\Om^{\bul}_{{\cal W}_n(Y)/{\cal W}_n(s),[~]}$ is a morphism of complexes, 
it is easy to check that the morphism 
\begin{align*} 
{\cal E}\otimes_{{\cal O}_{\mathfrak E}}
\Om^{\bul}_{{\mathfrak E}/{\cal W}_n(s),[~]}  
&\lo E_n\otimes_{{\cal W}_n({\cal O}_Y)'}
\Om^{\bul}_{{\cal W}_n(Y)/{\cal W}_n(s),[~]}. 
\tag{2.1.10.10}\label{eqn:fcndadw}
\end{align*} 
is a morphism of complexes. 
Because the target of the second morphism in (\ref{eqn:fcnidw}) 
is equal to 
the degree $i$-part of the complex (\ref{eqn:ldwy}),  
the second morphism is a morphism of complexes. 
As a result, we obtain the following desired composite morphism  
in ${\rm C}^+(\os{\circ}{g}{}^{-1}({\cal W}_n))$:  
\begin{align*} 
{\cal E}\otimes_{{\cal O}_{\cal Q}}
\Om^{\bul}_{{\cal Q}/{\cal W}_n(s)}  
&\lo E_n\otimes_{{\cal W}_n({\cal O}_Y)'}
\Om^{\bul}_{{\cal W}_n(Y)/{\cal W}_n(s),[~]} 
\lo 
E_n\otimes_{{\cal W}_n({\cal O}_Y)}
{\cal H}^{\bul}({\cal O}_{\mathfrak E}\otimes_{{\cal O}_{\cal Q}}
\Om^*_{{\cal Q}/{\cal W}_n(s)} ) \tag{2.1.10.11}\label{ali:fciqdw}.   
\end{align*} 
\par 
The functoriality of the first morphism in (\ref{ali:fciqdw}) is obvious by (\ref{ali:acdw}). 
The functoriality of the second morphism in 
(\ref{ali:fciqdw}) with respect to the diagram 
(\ref{cd:fnih2dw}) 
is clear by the construction in \cite[(4.9)]{hk}.  
\par  
The compatibility of the morphism (\ref{eqn:fnih3dw}) 
with the projections follows from \cite[(7.18)]{ndw}. 
\par 
(3): 
Because the morphism ${\cal W}_n(Y_i) \lo{\mathfrak E}_i$ is a PD-morphism, 
the diagram (\ref{cd:fnih1dw}) is commutative as a diagram of schemes. 
Since the log structure of ${\mathfrak E}_2$ is the pull-back of 
that of ${\cal Q}_2$, we see that the diagram (\ref{cd:fnih1dw}) 
is commutative as a diagram of log schemes.  
\end{proof}

\begin{rema}\label{rema:disd} 
(1) The morphisms ${\cal W}_n(Y)\lo  {\cal Q}$ 
and ${\cal W}_n(Y)\lo  {\mathfrak E}$ in (\ref{lemm:neh3fi}) (3) 
give the first morphism in (\ref{eqn:fcnidw}). 
We can dispense with the isomorphism (\ref{eqn:fwniwu}) 
(cf.~\cite[(4.19)]{hk}) to prove (\ref{lemm:neh3fi}). 
\par 
(2) We can  generalize (\ref{lemm:neh3fi}) (3) as follows. 
\par 
Let $(S_i,{\cal I}_i,\gam_i)$ $(i=1,2)$ be a fine log PD-scheme and 
let $S_{i0}$ be an exact closed subscheme of $S_i$ defined by ${\cal I}_i$. 
Let $Z_i$ be a fine log scheme over $S_{i0}$. 
Let $Z_i\os{\sus}{\lo} {\cal R}_i$ be an immersion over $S_i$ and 
let $Z_i\os{\sus}{\lo} {\mathfrak D}_i$ be a log PD-immersion over $(S_i,{\cal I}_i,\gam_i)$. 
Let ${\mathfrak E}_i$ be the log PD-envelope of the immersion 
$Z_i\os{\sus}{\lo} {\cal R}_i$ over $(S_i,{\cal I}_i,\gam_i)$. 
Assume that there exists the following commutative diagram 
\begin{equation*} 
\begin{CD} 
{\mathfrak D}_1@>>>  {\cal R}_1 \\ 
@VVV  @VVV \\  
{\mathfrak D}_2 @>>> {\cal R}_2,  
\end{CD} 
\tag{2.1.11.1}\label{cd:fniq0dw} 
\end{equation*}  
where the morphism 
${\mathfrak D}_1\lo {\mathfrak D}_2$ is a morphism of 
log PD-schemes. 
Then the induced diagram 
\begin{equation*} 
\begin{CD} 
{\mathfrak D}_1@>>>  {\mathfrak E}_1 \\ 
@VVV  @VVV \\  
{\mathfrak D}_2 @>>> {\mathfrak E}_2
\end{CD} 
\tag{2.1.11.2}\label{cd:fniq1dw} 
\end{equation*} 
is commutative. 
\end{rema}

The following (1) is a log version of Etesse's comparison theorem and 
a generalization of \cite[(7.19)]{ndw}(=a correction of \cite[(4.19)]{hk}).  
(In \cite[(7.12)]{nh3} we have also stated the following  (1)  
in the case of the trivial coefficient.)

\begin{theo}\label{theo:ccrw}
Let $E^{\bul \leq N}$ be a flat coherent log crystal of 
${\cal O}_{Y_{\bul \leq N}/{\cal W}_n(s)}$-modules. 
Assume that $Y_{\bul \leq N}$ has 
an affine $N$-truncated simplicial open covering of $Y_{\bul \leq N}$. 
Let $g \col 
Y_{\bul \leq N} \lo s \os{\sus}{\lo} {\cal W}_n(s)$ be 
the structural morphism. Then the following hold$:$ 
\par 
$(1)$ There exists a canonical isomorphism 
\begin{equation*} 
Ru_{Y_{\bul \leq N}/{\cal W}_n(s)*}
(E^{\bul \leq N}) 
\os{\sim}{\lo} 
E^{\bul \leq N}_n
\otimes_{{\cal W}_n({\cal O}_{Y_{\bul \leq N}})^{\star}}
({\cal W}_n{\Om}^{\bul}_{Y_{\bul \leq N}})^{\star} 
\tag{2.1.12.1}\label{eqn:ywnnny}
\end{equation*} 
in ${\rm D}^+(g^{-1}({\cal W}_n))$.
The isomorphisms $(\ref{eqn:ywnnny})$ for $n$'s are 
compatible with the two projections of both hand sides on $(\ref{eqn:ywnnny})$. 
\par 
$(2)$ Let $s\lo s'$ be a morphism of fine log schemes 
whose underlying schemes are the spectrums of perfect fields 
of characteristic $p>0$. 
The isomorphism $(\ref{eqn:ywnnny})$ is contravariantly functorial 
with respect to a morphism $h\col Y_{\bul \leq N}\lo Z_{\bul \leq N}$ of 
$N$-truncated simplicial log schemes over the morphism 
${\cal W}_n(s)\lo {\cal W}_n(s')$ and 
a morphism 
$h^*_{\rm crys}(F^{\bul \leq N})\lo E^{\bul \leq N}$ of 
${\cal O}_{Y_{\bul \leq N}/{\cal W}_n(s)}$-modules, 
where $Z_{\bul \leq N}$ and $F^{\bul \leq N}$ are similar objects to 
$Y_{\bul \leq N}$ and $E^{\bul \leq N}$, respectively. 
\end{theo}  
\begin{proof} 
(1): It suffices to prove (1) in the case $\star=$nothing. 
By using (\ref{lemm:lisj}) and Tsuzuki's functor as used in \cite[(3.5)]{nh2}, 
we have the morphism (\ref{eqn:ywnnny}) as follows. 
\par 
Let $Y'_{\bul \leq N}$ be the disjoint union of an affine open covering of $Y_{\bul \leq N}$. 
Since $\os{\circ}{Y}{}'_N$ is affine, there exists an immersion  
$Y'_N\os{\sus}{\lo} {\cal Y}'_N$ into a log smooth integral scheme over ${\cal W}_n(s)$.  
In fact, we have a log smooth integral lift ${\cal Y}'_N$ of $Y'_N$ over ${\cal W}_n(s)$ 
(\cite[(3.14)]{klog1}, \cite[(4.7)]{ny}). 
Because the closed immersion $Y'_N\os{\sus}{\lo} {\cal W}_n(Y'_N)$ is nilpotent 
(\cite[0 (1.3.13)]{idw}) and 
because ${\cal Y}'_N$ is log smooth over ${\cal W}_n(s)$,  
there exists a morphism 
${\cal W}_n(Y'_N)\lo {\cal Y}'_N$ over ${\cal W}_n(s)$ 
such that the composite morphism 
$Y'_N\os{\sus}{\lo} {\cal W}_n(Y'_N)\lo {\cal Y}'_N$ 
is the given immersion.   
For a morphism $\gam \col [N] \lo [m]$ in $\Del$, 
let $Y'(\gam)$ be the corresponding morphism 
$Y'_m \lo Y'_N$. 
Then we have the following commutative diagram 
\begin{equation*} 
\begin{CD} 
Y'_m @>{\subset}>>{\cal W}_n(Y'_m) @>>> {\cal Y}{}'_N\\ 
@V{Y'(\gam)}VV @VV{{\cal W}_n(Y'(\gam))}V  @|\\ 
Y'_N @>{\subset}>>{\cal W}_n(Y'_N)  
@>>> {\cal Y}{}'_N.  
\end{CD} 
\tag{2.1.12.2}\label{cd:fubnwn}
\end{equation*} 
By this commutative diagram, 
we have a morphism   
\begin{equation*} 
{\cal W}_n(Y'_{\bul \leq N})\lo \Gam^{{\cal W}_n(s)} _N({\cal Y}'_N)_{\bul \leq N} 
\tag{2.1.12.3}\label{eqn:pwh3ge} 
\end{equation*} 
of $N$-truncated simplicial log schemes, 
where $\Gam^{{\cal W}_n(s)}_N$ is the 
Tsuzuki's functor ((\ref{defi:tkifun})). 
Let 
${\cal W}_n(Y_{\bul \leq N,\bul})$ be the \v{C}ech diagram of 
${\cal W}_n(Y'_{\bul \leq N})$ over ${\cal W}_n(Y_{\bul \leq N})$: 
${\cal W}_n(Y_{ml})
:={\rm cosk}_0^{{\cal W}_n(Y_m)}{\cal W}_n(Y'_m)_l$ 
$(0\leq m \leq N$, $l\in {\mab N})$.  
It is easy to check that 
${\cal W}_n(Y_{\bul \leq N,\bul})$ 
is the canonical lift of $Y_{\bul \leq N,\bul}$.  
Set ${\cal Q}_{ml}
:={\rm cosk}_0^{{\cal W}_n(s)}(\Gam^{{\cal W}_n(s)} _N({\cal Y}'_N)_m)_l$.  
By (\ref{eqn:pwh3ge}) we have the following morphism   
\begin{equation*} 
{\cal W}_n(Y_{\bul \leq N,\bul}) \lo {\cal Q}_{\bul \leq N,\bul}
\tag{2.1.12.4}\label{eqn:ppwh3gep} 
\end{equation*} 
of $(N,\infty)$-truncated bisimplicial log schemes,  
which gives an immersion 
$Y_{\bul \leq N,\bul}\os{\sus}{\lo} {\cal Q}_{\bul \leq N,\bul}$. 
Let ${\mathfrak E}_{\bul \leq N,\bul}$ be the log PD-envelope of 
the immersion $Y_{\bul \leq N,\bul}\os{\sus}{\lo} {\cal Q}_{\bul \leq N,\bul}$ 
over $({\cal W}_n(s),p{\cal W}_n,[~])$. 
Let $E^{\bul \leq N,\bul}$ be the crystal of 
${\cal O}_{Y_{\bul \leq N,\bul}/{\cal W}_n(s)}$-modules 
obtained by $E^{\bul \leq N}$. 
Let $({\cal E}^{\bul \leq N,\bul},\nabla)$ be the corresponding 
${\cal O}_{{\mathfrak E}_{\bul \leq N,\bul}}$-module with integrable connection to $E^{\bul \leq N}$. 
Set $E^{\bul \leq N,\bul}_n:=(E^{\bul \leq N,\bul})_{{\cal W}_n(Y_{\bul \leq N,\bul})}$. 
By (\ref{lemm:neh3fi}) we have the following morphism of complexes: 
\begin{align*} 
{\cal E}^{\bul \leq N,\bul}
{\otimes}_{{\cal O}_{{\cal Q}_{\bul \leq N,\bul}}}
\Om^{\bul}_{{\cal Q}_{\bul \leq N,\bul}/{\cal W}_n(s)}   
\lo E^{\bul \leq N,\bul}_n
\otimes_{{\cal W}_n({\cal O}_{Y_{\bul \leq N,\bul}})}
{\cal H}^{\bul}({\cal O}_{{\mathfrak E}_{\bul \leq N,\bul}}
{\otimes}_{{\cal O}_{{\cal Q}_{\bul \leq N,\bul}}}
\Om^{*}_{{\cal Q}_{\bul \leq N,\bul}/{\cal W}_n(s)} ).  
\tag{2.1.12.5}\label{ali:enoq}
\end{align*} 
Let $g_{\bul}\col Y_{\bul \leq N,\bul}\lo {\cal W}_n(s)$ be 
the structural morphism. 
Let 
\begin{equation*} 
\pi_{\rm zar} \col 
((Y_{\bul \leq N,\bul})_{\rm zar},g^{-1}_{\bul}({\cal W}_n)) \lo 
((Y_{\bul \leq N})_{\rm zar},g^{-1}({\cal W}_n)) 
\tag{2.1.12.6}\label{eqn:ppyzd} 
\end{equation*} 
be the natural morphism of ringed topoi. 
Applying $R\pi_{{\rm zar}*}$ to the morphism (\ref{ali:enoq})
and using \cite[(6.4)]{klog1} or the log Poincar\'{e} lemma (\cite[(2.2.7)]{nh2}) 
and the cohomological descent, 
we obtain the morphism (\ref{eqn:ywnnny}) for the case $\star$=nothing. 
\par 
This morphism is independent of the choices of $Y'_{\bul \leq N}$  
and the immersion $Y'_N\os{\sus}{\lo} {\cal Y}'_N$. 
Indeed, let $Y''_{\bul \leq N}$ be the disjoint union of 
another affine open covering of $Y_{\bul \leq N}$.
By (\ref{prop:lissmp}) (2) we may assume that 
there exists a morphism $Y'_{\bul \leq N}\lo Y''_{\bul \leq N}$ over $Y_{\bul \leq N}$. 
Take an immersion $Y''_N\os{\sus}{\lo} {\cal Y}''_N$ into 
a log smooth scheme over ${\cal W}_n(s)$. 
Then we have a morphism ${\cal W}_n(Y''_N)\os{\sus}{\lo} {\cal Y}''_N$ 
as already explained. 
Then we have the following commutative diagram 
\begin{equation*} 
\begin{CD} 
{\cal W}_n(Y'_N) @>>> {\cal Y}'_N\times_{{\cal W}_n(s)}{\cal Y}''_N \\
@VVV @VVV \\
{\cal W}_n(Y''_N) @>>> {\cal Y}''_N, 
\end{CD}
\tag{2.1.12.7}\label{cd:wyy} 
\end{equation*} 
where the right vertical morphism is the second projection. 
(Note that, even if one obtains the following commutative diagram 
\begin{equation*} 
\begin{CD} 
Y'_N @>>> {\cal Y}'_N \\
@VVV @VVV \\
Y''_N @>>> {\cal Y}''_N, 
\end{CD} 
\end{equation*} 
one cannot necessarily obtain the following commutative diagram 
\begin{equation*} 
\begin{CD} 
{\cal W}_n(Y'_N) @>>> {\cal Y}'_N \\
@VVV @VVV \\
{\cal W}_n(Y''_N) @>>> {\cal Y}''_N 
\end{CD} 
\end{equation*} 
in general.) 
The rest of the proof of the independence of the choices 
is a routine work. 
\par 
Now we have only to prove (1) for the constant 
$N$-truncated simplicial case $Y_{\bul \leq N}=Y$ and 
$E^{\bul \leq N}=E$. 
Then the rest of the proof is the same proof as that of \cite[II (2.1)]{et}. 
However we give the complete proof because 
we would like to give a (slight) simplification 
of a part of the proof of [loc.~cit.].
\par 
The question is local on $Y$; 
we may assume that there exists an immersion 
$Y\os{\sus}{\lo} {\cal Q}$ into a log smooth scheme 
over ${\cal W}_n(s)$. Let ${\mathfrak E}$ be 
the log PD-envelope of this immersion over 
$({\cal W}_n(s),p{\cal W}_n,[~])$. 
Set ${\cal E}:=E_{\mathfrak E}$ and $E_n:=E_{{\cal W}_n(Y)}$. 
By the log Poincar\'{e} lemma, we have the following equality:  
\begin{equation*} 
Ru_{Y/{\cal W}_n(s)*}(E)=
{\cal E}\otimes_{{\cal O}_{\cal Y}}
\Om^{\bul}_{{\cal Y}/{\cal W}_n(s)}. 
\tag{2.1.12.8}\label{eqn:ynse} 
\end{equation*} 
Then the morphism (\ref{eqn:ywnnny}) is equal to 
the morphism (\ref{eqn:fnih3dw}). 
It is obvious that 
the morphism (\ref{eqn:fnih3dw}) is an isomorphism 
in the case $n=1$ since the morphism (\ref{eqn:fnih3dw}) 
is equal to the isomorphism ${\rm id}_{E_1}\otimes C^{-1}$. 
As in \cite[II (2.1)]{et}, consider the $p$-adic filtration $F'$ 
(resp.~the canonical filtration $F:={\rm Fil}$) on 
${\cal E}\otimes_{{\cal O}_{\cal Y}}
\Om^{\bul}_{{\cal Y}/{\cal W}_n(s)}$
(resp.~
$E_n\otimes_{{\cal W}_n({\cal O}_Y)}
{\cal W}_n\Om^{\bul}_Y$). 
Then we have the following commutative diagram 
\begin{equation*} 
\begin{CD}
E_Y\otimes_{{\cal O}_Y}
\Om^{\bul}_{Y/s} @>{{\rm id}_{E_1}\otimes C^{-1},\sim}>> E_1\otimes_{{\cal O}_Y}
\Om^{\bul}_{Y/s}\\
@V{p^m}V{\simeq}V @VV{{\bf p}^m}V \\ 
{\rm gr}^{m}_{F'}({\cal E}\otimes_{{\cal O}_{\cal Y}}
\Om^{\bul}_{{\cal Y}/{\cal W}_n(s)})
@>>> {\rm gr}^{m}_{F}
(E_n\otimes_{{\cal W}_n({\cal O}_Y)}
{\cal W}_n\Om^{\bul}_Y) 
\end{CD}
\tag{2.1.12.9}\label{cd:ysy}
\end{equation*}   
for $0\leq m\leq n-1$ since $E_Y=E_1$. 
Hence it suffices to prove that the injective morphism 
\begin{equation*} 
{\bf p}\col 
E_{n+1}\otimes_{{\cal W}_{n+1}({\cal O}_Y)}
{\rm gr}_F^n({\cal W}\Om^{\bul}_Y)
\lo E_{n+2}\otimes_{{\cal W}_{n+2}({\cal O}_Y)}
{\rm gr}_F^{n+1}({\cal W}\Om^{\bul}_Y)
\tag{2.1.12.10}\label{cd:yn1y}
\end{equation*} 
is a quasi-isomorphism. 
In other words, it suffices to prove that the complex 
\begin{equation*} 
E_{n+2}\otimes_{{\cal W}_{n+2}({\cal O}_Y)}
({\rm gr}_F^{n+1}({\cal W}\Om^{\bul}_Y)
/{\bf p}{\rm gr}_F^n({\cal W}\Om^{\bul}_Y))
\end{equation*}  
is acyclic. Let $\{e_i\}$ be a local basis of $E_{n+2}$. 
Then consider a local section 
$\sum e_i\otimes u_i$ 
$(u_i\in {\rm gr}_F^{n+1}({\cal W}{\Om}^q_Y))$. 
Then $\nabla(\sum e_i\otimes u_i)=  
\sum(\nabla(e_i)\wedge u_i+e_i\otimes du_i)$. 
Express 
$\nabla(e_i)=
\sum_j e_j\otimes \om_{ji}$ 
$(\om_{ji}\in ({\cal W}_{n+2}{\Om}^1_{Y}))$. 
Then 
\begin{align*} 
\nabla(\sum e_i\otimes u_i)=
\sum e_i\otimes (\sum_j\om_{ij}\wedge u_j+du_i).
\end{align*}   
Hence, as in [loc.~cit., II (2.5)], it suffices to prove
the following: 
\parno 
If 
\begin{equation*} 
\sum_j\om_{ij}\wedge u_j+du_i\in 
{\bf p}\,{\rm gr}_F^n({\cal W}{\Om}^{q+1}_{Y}) 
\quad 
(u_i\in ({\rm gr}_F^{n+1}{\cal W}{\Om}^q_Y)), 
\tag{2.1.12.11}\label{eqn:duijcn}
\end{equation*}  
then 
\begin{align*}
u_i\equiv \sum_j\om_{ij}\wedge v_j+dv_i
{\mod {\bf p}}\,
{\rm gr}_F^n({\cal W} {\Om}^q_{Y})
\tag{2.1.12.12}\label{ali:duijcn}
\end{align*} 
for some $v_j\in 
{\rm gr}_F^{n+1}({\cal W}{\Om}^{q-1}_Y)$
([loc.~cit.,~(2.5.7.2), (2.5.7.3)]). 
Because $u_j\in {\rm Fil}^{n+1}{\cal W}_{n+2}\Om^q_Y
=V^{n+1}\Om^q_Y+dV^{n+1}\Om^q_Y$ 
(\cite[I (3.31)]{idw}, \cite[(1.16), p.~258]{lodw}), \cite[(6.15.1)]{ndw}), 
there exists sections $a_j,b_j\in \Om^q_Y$ 
such that 
$u_j=V^{n+1}a_j+dV^{n+1}b_j$. 
Then $F^{n+1}(u_j)=p^{n+1}a_j+db_j=db_j$ and 
$F^{n+1}(du_i)=da_i$. 
Since 
\begin{align*} 
F^{n+1}
({\bf p}\,{\rm gr}_F^n({\cal W}{\Om}^{q+1}_{Y}))
& =
F^{n+1}
({\bf p}\,{\rm Fil}^n({\cal W}_{n+1}{\Om}^{q+1}_{Y}))
\subset 
F^{n+1}({\bf p}({\cal W}_{n+1}{\Om}^{q+1}_{Y}))
=0,
\end{align*} 
we have the following equalities: 
\begin{align*} 
0&=F^{n+1}(\sum_j\om_{ij}\wedge u_j+du_i)
=\sum_j F^{n+1}(\om_{ij})\wedge db_j+da_i
\tag{2.1.12.13}\label{eqn:0om}.
\end{align*} 
Moreover, because $F^{n+1}(\om_{ij})$ is a one form, 
we have the following formula:  
\begin{align*} 
d(F^{n+1}(\om_{ij})\wedge b_j)
& =dF^{n+1}(\om_{ij})\wedge b_j-
F^{n+1}(\om_{ij})\wedge db_j \\
&=p^{n+1}F^{n+1}(d\om_{ij})\wedge b_j-
F^{n+1}(\om_{ij})\wedge db_j\\ 
&=-F^{n+1}(\om_{ij})\wedge db_j.
\end{align*}  
Hence, by (\ref{eqn:0om}), we have 
$d(a_i-\sum_jF^{n+1}(\om_{ij})\wedge b_j)=0$. 
Because 
${\rm Ker}(d\col {\cal W}_n\Om^q_Y
\lo {\cal W}_n\Om^{q+1}_Y)=F^n{\cal W}_{2n}\Om^q_Y$ 
(\cite[I (3.21)]{idw}, \cite[II (1.3)]{ir}, 
\cite[1.3.4]{msemi}, \cite[(6.21.1)]{ndw}), 
there exists a section 
$c_i\in {\cal W}_2\Om^q_Y$
such that 
$a_i-\sum_j(F^{n+1}\om_{ij})b_j=Fc_i$. 
By using the formula $V((Fx)y)=xVy$, 
we have the following: 
\begin{align*} 
u_i &= V^{n+1}(\sum_jF^{n+1}(\om_{ij})\wedge b_j+Fc_i)
+dV^{n+1}b_i  \tag{2.1.12.14}\label{ali:vbnc}\\
& =\sum_j\om_{ij}\wedge V^{n+1}b_j+pV^nc_i+dV^{n+1}b_i \\
& \equiv dV^{n+1}b_i +\sum_j\om_{ij}\wedge V^{n+1}b_j 
\quad {\mod {\bf p}}\,{\rm gr}_F^n{\cal W} {\Om}^q_{Y}, 
\end{align*}
which implies (\ref{ali:duijcn}). 
\par 
The compatibility of the morphism (\ref{eqn:ywnnny})
with the projections follows from \cite[(7.18)]{ndw}.
\par 
(2): 
Let $Z'_N$ be the disjoint union of an affine open covering of $Z_N$. 
Let $Y'_N$ be the disjoint union of an affine open covering of $Y_N$ such that 
there exists a morphism $Y'_N\lo Z'_N$ over the morphism $Y_N\lo Z_N$.  
Let $Y'_N\os{\sus}{\lo} {\cal Y}''_N$ and $Z'_N\os{\sus}{\lo} {\cal Z}'_N$ 
be immersions into log smooth schemes over ${\cal W}_n(s)$ and ${\cal W}_n(s')$, 
respectively.  
Then we have the following morphisms 
${\cal W}_n(Y'_N) \lo {\cal Y}''_N$ and 
${\cal W}_n(Z'_N) \lo {\cal Z}'_N$ such that 
the composite morphisms 
$Y'_N\os{\sus}{\lo} {\cal W}_n(Y'_N) \lo {\cal Y}''_N$
and 
$Z'_N\os{\sus}{\lo} {\cal W}_n(Z'_N) \lo {\cal Z}'_N$ 
are the given immersions. 
Set 
\begin{align*} 
{\cal Y}'_N:= {\cal Y}''_N\times_{{\cal W}_n(s)}
({\cal Z}'_N\times_{{\cal W}_n(s')}{\cal W}_n(s)).
\tag{2.1.12.15}\label{ali:ynws}
\end{align*}  
Then we have the following commutative diagram 
\begin{equation*} 
\begin{CD} 
{\cal W}_n(Y'_N) @>>> {\cal Y}'_N \\
@VVV @VVV \\
{\cal W}_n(Z'_N) @>>> {\cal Z}'_N, 
\end{CD}
\tag{2.1.12.16}\label{cd:wiy} 
\end{equation*} 
where the right vertical morphism is induced by the second projection. 
Let ${\cal W}_n(Z_{\bul \leq N,\bul})$, ${\cal R}_{\bul \leq N,\bul}$ 
and ${\mathfrak F}_{\bul \leq N,\bul}$ be similar objects to 
${\cal W}_n(Y_{\bul \leq N,\bul})$, ${\cal Q}_{\bul \leq N,\bul}$ 
and ${\mathfrak E}_{\bul \leq N,\bul}$, respectively.  
Then we have the following commutative diagram
\begin{equation*} 
\begin{CD}
{\cal W}_n(Y_{\bul \leq N,\bul})@>>> {\cal Q}_{\bul \leq N,\bul}\\
@VVV @VVV \\ 
{\cal W}_n(Z_{\bul \leq N,\bul})@>>> {\cal R}_{\bul \leq N,\bul}.
\end{CD}
\tag{2.1.12.17}\label{cd:kecycd}
\end{equation*} 
Hence, by (\ref{lemm:neh3fi}) (3), 
we have the following commutative diagram
\begin{equation*} 
\begin{CD}
{\cal W}_n(Y_{\bul \leq N,\bul})@>>> {\mathfrak E}_{\bul \leq N,\bul}\\
@VVV @VVV \\ 
{\cal W}_n(Z_{\bul \leq N,\bul})@>>> {\mathfrak F}_{\bul \leq N,\bul}.
\end{CD}
\tag{2.1.12.18}\label{cd:cwefd}
\end{equation*} 
Let ${\cal F}^{\bul \leq N,\bul}$ and $F^{\bul \leq N,\bul}_n$
be similar objects to ${\cal E}^{\bul \leq N,\bul}$ and $E^{\bul \leq N,\bul}_n$, respectively. 
Let $h_{\bul \leq N,\bul}\col Y_{\bul \leq N,\bul}\lo Z_{\bul \leq N,\bul}$ 
be the induced morphism by the morphism 
${\cal W}_n(Y_{\bul \leq N,\bul})\lo
{\cal W}_n(Z_{\bul \leq N,\bul})$. 
By (\ref{lemm:neh3fi}) 
we have the following commutative diagram  
\begin{equation*} 
\begin{CD} 
h_{\bul \leq N,\bul*}({\cal E}^{\bul \leq N,\bul}
{\otimes}_{{\cal O}_{{\cal Q}_{\bul \leq N,\bul}}}
\Om^{\bul}_{{\cal Q}_{\bul \leq N,\bul}/{\cal W}_n(s)})    
@>>>  \\
@AAA \\ 
{\cal F}^{\bul \leq N,\bul}
{\otimes}_{{\cal O}_{{\cal R}_{\bul \leq N,\bul}}}
\Om^{\bul}_{{\cal R}_{\bul \leq N,\bul}/{\cal W}_n(s)}   
@>>>  
\end{CD}
\tag{2.1.12.19}\label{ali:ensoq}
\end{equation*} 
\begin{equation*} 
\begin{CD}
h_{\bul \leq N,\bul*}(E^{\bul \leq N,\bul}_n
\otimes_{{\cal W}_n({\cal O}_{Y_{\bul \leq N,\bul}})}
{\cal H}^{\bul}({\cal O}_{{\mathfrak E}_{\bul \leq N,\bul}}
{\otimes}_{{\cal O}_{{\cal Q}_{\bul \leq N,\bul}}}
\Om^{*}_{{\cal Q}_{\bul \leq N,\bul}/{\cal W}_n(s)}))\\
@AAA \\
F^{\bul \leq N,\bul}_n
\otimes_{{\cal W}_n({\cal O}_{Z_{\bul \leq N,\bul}})}
{\cal H}^{\bul}({\cal O}_{{\mathfrak F}_{\bul \leq N,\bul}}
{\otimes}_{{\cal O}_{{\cal R}_{\bul \leq N,\bul}}}
\Om^{*}_{{\cal R}_{\bul \leq N,\bul}/{\cal W}_n(s)}).  
\end{CD}
\end{equation*}
By this commutative diagram, 
we obtain the desired functoriality. 
\end{proof}

\begin{rema}\label{rema:altn}  
(1) In the proof of (\ref{theo:ccrw}) we have not used 
lifts of the Frobenius endomorphisms of log schemes (cf.~\cite[pp.~601--605]{idw}, 
especially \cite[p.~603]{idw}) nor the log version of the lemma of 
Dwork-Dieudnonn\'{e}-Cartier (\cite[(7.13)]{ndw})
(see (\ref{lemm:ddc}) below). 
In the case of trivial log structures and the constant simplicial case, 
to prove the functoriality of the isomorphism (\ref{eqn:ywnnny}), 
Illusie has used lifts of Frobenii and 
the lemma of Dwork-Dieudnonn\'{e}-Cartier (\cite[p.~605]{idw}). 
In particular, we have simplified the proof in [loc.~cit.] 
by considering the product in (\ref{ali:ynws}).  
\par 
(2) Because the isomorphism in \cite[(4.19)]{hk} has not been shown 
to be independent of the choice of an embedding system, 
I cannot understand the meaning of ``canonical'' 
of the following canonical isomorphism
\begin{align*} 
W_n\om^{\bul}_Y\simeq Ru^{\log}_{Y/(W_n,W_n(L)*}({\cal O}_{Y/W_n})
\tag{2.1.13.1}\label{ali:wnru}
\end{align*} 
in [loc.~cit.]. 
(The statement of \cite[(7.19)]{ndw} and the proof of it have no problem.) 
\par 
(3) In \cite[(4.19)]{hk} Hyodo and Kato have claimed that the isomorphism 
(\ref{ali:wnru}) is compatible with the Frobenii on both hand sides on 
(\ref{ali:wnru}). However they have not proved this compatibility 
because they have not proved neither the functoriality of this isomorphism nor 
a log version of a lemma of Dwork-Dieudnonn\'{e}-Cartier 
(\cite[(7.13)]{ndw}) in [loc.~cit.].  
\par 
(4) The functoriality of the canonical isomorphisms in 
\cite[(2.1)]{et}, \cite[(4.19)]{hk}, \cite[(1.22)]{lodw}, \cite[(7.19)]{ndw} 
and \cite[(4.4.17)]{ols} has not been proved in these articles.  
It should be proved in these articles. 
\par 
(5) Except Berthelot's article \cite[the argument after (1.5)]{bfi}, 
I cannot find a reference in which such a product in (\ref{ali:ynws}) 
and such a commutative diagram in (\ref{cd:wiy})  
are used in the proof of the functoriality of a canonical isomorphism. 
Though we claim the functoriality of 
the canonical filtered isomorphism in \cite[(7.6.1)]{ndw}, 
the proof is not complete:  
the diagram 
\begin{equation*} 
\begin{CD} 
({\cal W}_n(X_{\bul \leq N,\bul}),
{\cal W}_n(D_{\bul \leq N,\bul}\cup Z_{\bul \leq N,\bul}))
@>>> {\cal P}_{\bul \leq N,\bul} \\ 
@VVV @VVV \\ 
({\cal W}_n(Y_{\bul \leq N,\bul}),
{\cal W}_n(E_{\bul \leq N,\bul}\cup W_{\bul \leq N,\bul}))
@>>> {\cal R}_{\bul \leq N,\bul}.  
\end{CD} 
\end{equation*} 
in [loc.~cit.] is not necessarily commutative. 
We can complete the proof by considering a product 
$({\cal X}'_N,{\cal D}'_N\cup {\cal Z}'_N)
\times_{{\cal W}_n}({\cal Y}'_N,{\cal E}'_N\cup {\cal W}'_N)$
as in (\ref{ali:ynws}). Here $({\cal X}'_N,{\cal D}'_N\cup {\cal Z}'_N)$ 
and $({\cal Y}'_N,{\cal E}'_N\cup {\cal W}'_N)$ are log schemes in [loc.~cit.]. 
\end{rema}

The following (see also \cite[(7.16), (7.17)]{ndw}) 
is obtained by a log version of a lemma of Dwork-Dieudnonn\'{e}-Cartier 
(\cite[(7.13)]{ndw}):  we have used this in the proof of \cite[(7.19)]{ndw}. 
(We do not use this lemma in this book.)
\begin{lemm}\label{lemm:ddc} 
Let ${\cal Z}$ be an integral log formal scheme over ${\cal W}(s)$ 
such that ${\cal O}_{\cal Z}$ is $p$-torsion-free.  
Set ${\cal Z}_n:={\cal Z}\otimes_{{\cal W}}{\cal W}_n$ $(n\in {\mab Z}_{\geq 1})$. 
Assume that ${\cal Z}$ has a lift of the Frobenius endomorphism of ${\cal Z}_1$. 
Then there exists a functorial morphism 
$s_{{\cal Z},n}\col {\cal W}_n({\cal Z}_1)\lo {\cal Z}_n$. 
Here we mean by the functoriality a functoriality with respect to morphisms of 
integral log formal schemes which are compatible with lifts of Frobenii 
over a morphism 
${\cal W}(s)\lo {\cal W}(s')$ induced 
by a morphism $s\lo s'$ of fine log schemes 
whose underlying schemes are the spectrums of perfect fields 
of characteristic $p>0$.    
\end{lemm}

Next we prove the log version (\ref{theo:bcn}) below of \cite[I~(1.14)]{idw}: 
the log de Rham-Witt sheaf of a log smooth scheme of Cartier type 
commutes with the pull-back of the log \'{e}tale morphism 
under a certain condition. The proof of this log version 
is different from the proof of [loc.~cit.]. 
To prove it, we need the theorem (\ref{theo:lgetb}).  
This is a log version of \cite[0 (2.2.7)]{idw}.  
The proof of (\ref{theo:lgetb}) 
fills a gap in the explanation for \cite[0 (2.2.7)]{idw}  
as a very special case
(see (\ref{rema:idp}) below). 
To prove (\ref{theo:lgetb}), we need the following: 

\begin{lemm}\label{lemm:isoex}
Let $f\col X\lo Y$ and $g\col Y\lo Z$ be morphisms of fine log schemes. 
Then the following hold$:$
\par 
$(1)$ Assume that $g$ and $g\circ f$ are log \'{e}tale.  
Then $f$ is log \'{e}tale. 
\par 
$(2)$ {\rm (cf.~\cite[IV (3.3.8)]{ob})} 
Assume that the following five conditions hold$:$
\par 
~$({\rm a})$ $f$ is log \'{e}tale.
\par 
~$({\rm b})$ $f$ is weakly purely inseparable over ${\mab F}_p$ 
in the sense of {\rm \cite[III (2.4.2)]{ob}}.  
\par 
~$({\rm c})$ $g\circ f$ is log \'{e}tale or weakly purely inseparable.
\par 
~$({\rm d})$ $g\circ f$ is exact. 
\par 
~$({\rm e})$ the log structure of $Y$ is saturated. 
\par
Then $f$ is an isomorphism. 
\end{lemm}
\begin{proof} 
(1): (1) is well-known. Indeed, by \cite[(3.12)]{klog1} the following sequence 
\begin{align*} 
f^*(\Om^1_{Y/Z}) \lo \Om^1_{X/Z} \lo \Om^1_{X/Y}\lo 0
\end{align*} 
is exact (cf.~the proof of (\ref{prop:mce})). 
By the assumption, $\Om^1_{X/Z}=0$ and $f^*(\Om^1_{Y/Z})=0$.
Obviously the morphism $f^*(\Om^1_{Y/Z}) \lo \Om^1_{X/Z}$ 
is an isomorphism. Hence $f$ is log \'{e}tale by [loc.~cit.] again. 
\par 
(2): (I have heard from C.~Nakayama that \cite[(4.11)]{klog1} is mistaken. 
Hence I do not use [loc.~cit.].)
We prove that $f$ is exact. 
Because $f$ is weakly purely inseparable and log \'{e}tale, 
the pull-back morphism $f^* \col f^*(M_Y^{\rm gp})\lo M_X^{\rm gp}$ 
is an isomorphism. Let $a$ be a local section of $f^*(M_Y^{\rm gp})$ 
such that $f^*(a)\in M_X$. By the condition (c), 
there exists a positive integer $m$ such that 
$f^*(a^m)\in {\rm Im}(f^*g^*(M^{\rm gp}_Z)\lo M^{\rm gp}_X)$. 
Because $g\circ f$ is exact, 
$f^*(a^m)\in {\rm Im}(f^*g^*(M_Z)\lo M_X)$.
In particular, 
$f^*(a^m)\in {\rm Im}(f^*(M_Y)\lo M_X)$. 
Since $M_Y$ is saturated, $a\in {\rm Im}(f^*(M_Y)\lo M_X)$. 
Hence $f$ is exact. This means that $M_X=f^*(M_Y)$. 
Because $f$ is log \'{e}tale, $\os{\circ}{f}$ is \'{e}tale. 
Because $\os{\circ}{f}$ is radicial ({\rm \cite[III (2.4.2)]{ob}}), 
$\os{\circ}{f}$ is an open immersion by \cite[(17.9.1)]{ega4}. 
We complete the proof of (2). 
\end{proof}

\par 
Let $S$ be a fine log scheme of characteristic $p>0$. 
Let $X$ be a fine log smooth scheme of Cartier type over $S$. 
Next we redefine $Z_n\Om^i_{X/S}$ and 
$B_n\Om^i_{X/S}$ from our point of view.  
The way of our definition of these sheaves
are slightly different from that in \cite[0 (2.2.2)]{idw};  
our way is necessary for the proof of (\ref{theo:lgetb}) below. 
If we use the definition in [loc.~cit.], the proof of 
(\ref{theo:lgetb}) will be confusing. 
\par 
Let $F_X\col X\lo X$ be the absolute Frobenius endomorphism of $X$. 
For a nonnegative integer $n$, 
set $X^{(p^n)}:=X\times_{S,F_S^n}S$ and 
let $F^n_{X/S}\col X\lo X^{(p^n)}$ be 
the induced morphism over $S$ by $F_X^n$. 
Here we set $F_S^0:={\rm id}_S$ and $F_X^0:={\rm id}_X$.
Set $Z\Om^i_{X/S}:={\rm Ker}(d\col \Om^i_{X/S}\lo \Om^{i+1}_{X/S})$ 
and $B\Om^i_{X/S}:={\rm Im}(d\col \Om^{i-1}_{X/S}\lo \Om^i_{X/S})$. 
Set also  
\begin{align*} 
Z_1\Om^i_{X/S}: =F_{X/S*}(Z\Om^i_{X/S})\quad 
{\rm and}  \quad
B_1\Om^i_{X/S}:=F_{X/S*}(B\Om^i_{X/S}).
\end{align*}  
It is convenient to set $Z_0\Om^i_{X/S}: =\Om^i_{X/S}$ and 
$B_0\Om^i_{X/S}:=0$. 
Though $Z\Om^i_{X/S}$ and $B\Om^i_{X/S}$ are only abelian sheaves on $X$, 
$Z_1\Om^i_{X/S}$ and $B_1\Om^i_{X/S}$ are naturally 
${\cal O}_{X^{(p)}}$-submodules of 
$F_{X/S*}(\Om^i_{X/S})$. 
By \cite[(4.12) (1)]{klog1} and by using 
a fact that $\os{\circ}{F}_{X/S} 
\col \os{\circ}{X}\lo \os{\circ}{X}{}^{(p)}$ is 
a homeomorphism of topological spaces (\cite[XV Proposition 2 a)]{sga5-2}), 
we have the following Cartier isomorphism 
of ${\cal O}_{X^{(p)}}$-modules: 
\begin{equation*} 
C^{-1}_{X/S} \col \Om^i_{X^{(p)}/S}\os{\sim}{\lo} 
F_{X/S*}({\cal H}^i(\Om^{\bul}_{X/S}))
=Z_1\Om^i_{X/S}/B_1\Om^i_{X/S}.  
\tag{2.1.15.1}\label{eqn:oxps} 
\end{equation*} 
Assume that two ${\cal O}_{X^{(p^n)}}$-modules 
$Z_n\Om^i_{X/S}$ and $B_n\Om^i_{X/S}$ of  
$F^n_{X/S*}(\Om^i_{X/S})$ are defined for a positive integer $n$. 
Applying $F^n_{X^{(p)}/S*}$ to (\ref{eqn:oxps}), 
we have the following isomorphism of ${\cal O}_{X^{(p^{n+1})}}$-modules: 
\begin{equation*} 
F^n_{X^{(p)}/S*}(C^{-1}_{X/S}) \col 
F^n_{X^{(p)}/S*}(\Om^i_{X^{(p)}/S})
\os{\sim}{\lo} 
F^{n+1}_{X/S*}(Z\Om^i_{X/S})
/F^{n+1}_{X/S*}(B\Om^i_{X/S}). 
\tag{2.1.15.2}\label{eqn:oxznps} 
\end{equation*} 
Let $D$ be $Z$ or $B$. 
Then we define $Z_{n+1}\Om^i_{X/S}$ and $B_{n+1}\Om^i_{X/S}$ 
by the following formula: 
\begin{equation*} 
F^n_{X^{(p)}/S*}(C^{-1}_{X/S}) \col 
D_n\Om^i_{X^{(p)}/S}
\os{\sim}{\lo} 
D_{n+1}\Om^i_{X/S}
/F^{n+1}_{X/S*}(B\Om^i_{X/S}) \quad (n\in {\mab N}). 
\tag{2.1.15.3}\label{eqn:ofznps} 
\end{equation*} 
By the definitions of $Z_n\Om^i_{X/S}$ and $B_n\Om^i_{X/S}$ in this book, 
we have the following inclusions 
\begin{align*} 
F_{X^{(p^n)}/S*}(B_{n}\Om^i_{X/S}) \subset B_{n+1}\Om^i_{X/S}\subset 
F^{n+1}_{X/S*}(Z\Om^i_{X/S}) \subset 
F^{n+1}_{X/S*}(\Om^i_{X/S}) 
\tag{2.1.15.4}\label{eqn:ofps} 
\end{align*} 
and 
\begin{align*} 
Z_{n+1}\Om^i_{X/S} \subset F_{X^{(p^n)}/S*}(Z_n\Om^i_{X/S})
\subset F^{n+1}_{X/S*}(Z\Om^i_{X/S})\subset 
F^{n+1}_{X/S*}(\Om^i_{X/S}) . 
\tag{2.1.15.5}\label{eqn:fznp} 
\end{align*} 
Because $F^n_{X/S} \col X\lo X^{(p^n)}$ is 
a homeomorphism of underlying topological spaces, 
we can define two abelian sheaves 
${}_nZ\Om^i_{X/S}$ and ${}_nB\Om^i_{X/S}$ on $X$
by the following formula 
\begin{align*} 
{}_nD\Om^i_{X/S}
:=((F^{n}_{X/S})^{-1})_*(D_n\Om^i_{X/S}). 
\tag{2.1.15.6}\label{eqn:ofznds} 
\end{align*} 
In \cite{idw} (and except this book) 
${}_nZ\Om^i_{X/S}$ and ${}_nB\Om^i_{X/S}$ have been denoted by 
$Z{}_n\Om^i_{X/S}$ and $B{}_n\Om^i_{X/S}$, respectively.  
\par 
By the definition of $Z_n\Om^i_{X/S}$ and $B_n\Om^i_{X/S}$, 
we have the following isomorphism of ${\cal O}_{X^{(p^{n+1})}}$-modules: 
\begin{align*} 
F^n_{X^{(p)}/S*}(C^{-1}_{X/S})\col 
Z_n\Om^i_{X^{(p)}/S}/B_n\Om^i_{X^{(p)}/S}
\os{\sim}{\lo} 
Z_{n+1}\Om^i_{X/S}/B_{n+1}\Om^i_{X/S}. 
\tag{2.1.15.7}\label{eqn:ofznbs} 
\end{align*} 
Iterating the isomorphism above, 
we have the following isomorphism of ${\cal O}_{X^{(p^n)}}$-modules: 
\begin{align*} 
C^{-n}_{X/S}\col \Om^i_{X^{(p^n)}/S}
\os{\sim}{\lo} 
Z_n\Om^i_{X/S}/B_n\Om^i_{X/S}. 
\tag{2.1.15.8}\label{eqn:ocznps} 
\end{align*} 
Consider the following exact sequence 
of ${\cal O}_{X^{(p^{n+1})}}$-modules: 
\begin{align*} 
0\lo Z_1\Om^{i}_{X^{(p^n)}/S}\lo F_{X^{(p^n)}/S*}(\Om^{i}_{X^{(p^n)}/S})
\os{d}{\lo} B_1\Om^{i+1}_{X^{(p^n)}/S}\lo 0
\end{align*} 
(cf.~\cite[p.~520]{idw}).  
Applying $C^{-n}_{X^{(p)}/S}$ to this exact sequence, we have 
the following exact sequence of ${\cal O}_{X^{(p^{n+1})}}$-modules: 
\begin{align*} 
0& \lo Z_{n+1}\Om^{i}_{X^{(p^n)}/S}/F_{X^{(p^n)}/S*}(B_n\Om^{i}_{X/S})
\lo F_{X^{(p^n)}/S*}(Z_n\Om^i_{X^{(p^n)}/S})/F_{X^{(p^n)}/S*}(B_n\Om^{i}_{X/S})
\tag{2.1.15.9}\label{eqn:oclnps} \\
&\os{d}{\lo} B_{n+1}\Om^{i+1}_{X^{(p^n)}/S}/F_{X^{(p^n)}/S*}(B_n\Om^{i+1}_{X/S})\lo 0. 
\end{align*} 
Hence we have the following isomorphism of ${\cal O}_{X^{(p^{n+1})}}$-modules:
\begin{align*} 
F_{X^{(p^n)}/S*}(Z_n\Om^i_{X^{(p^n)}/S})/Z_{n+1}\Om^{i}_{X^{(p^n)}/S}
\os{d,\sim}{\lo} B_{n+1}\Om^{i+1}_{X^{(p^n)}/S}/F_{X^{(p^n)}/S*}(B_n\Om^{i+1}_{X/S}). 
\tag{2.1.15.10}\label{eqn:oacocps} 
\end{align*} 
Set $C^{-n}\Om^{\bul}_{X^{(p^n)}/S}:=Z_n\Om^{\bul}_{X/S}/B_n\Om^{\bul}_{X/S}$. 
By  (\ref{eqn:oacocps})
the inclusion 
$Z_{n+1}\Om^i_{X/S} \subset F_{X^{(p^n)}/S*}(Z_n\Om^i_{X/S})$ 
induces the following isomorphism of 
${\cal O}_{X^{(p^{n+1})}}$-modules: 
\begin{align*} 
Z_{n+1}\Om^{i}_{X/S}/B_{n+1}\Om^i_{X/S}
\os{\sim}{\lo} 
F_{X^{(p^n)}/S*}({\cal H}^i(C^{-n}\Om^{\bul}_{X^{(p^n)}/S}))
\tag{2.1.15.11}\label{eqn:ococps} 
\end{align*}

\begin{theo}\label{theo:lgetb}
Let $S$ be as above. 
Assume that $S$ is saturated. 
Let $g_{}\col Z_{\bul \leq N}\lo Y_{\bul \leq N}$ be 
a log \'{e}tale morphism 
of log smooth $N$-truncated simplicial log schemes of 
Cartier type over $S$. 
Let $n$ be a nonnegative integer. 
Let $g^{(p^n)}_{}\col Z^{(p^n)}_{\bul \leq N}
:=Z_{\bul \leq N}\times_{S,F^n_S}S\lo Y^{(p^n)}_{\bul \leq N}
:=Y_{\bul \leq N}\times_{S,F^n_S}S$ 
be the base change of $g$. 
Assume that $Y_{\bul \leq N}$ is saturated and 
that the structural morphism $Y_{\bul \leq N}\lo S$ and $g_{}$ are saturated. 
Then, for $i\in {\mab N}$, 
\begin{align*} 
g^{(p^n)*}_{}(Z_n\Om^i_{Y_{\bul \leq N}/S})
=Z_n\Om^i_{Z_{\bul \leq N}/S}
\tag{2.1.16.1}\label{ali:zzxn}
\end{align*} 
and 
\begin{align*} 
g^{(p^n)*}_{}(B_n\Om^i_{Y_{\bul \leq N}/S})
=B_n\Om^i_{Z_{\bul \leq N}/S}. 
\tag{2.1.16.2}\label{ali:bbn}
\end{align*} 
\end{theo}
\begin{proof} 
(1): We may assume that $N=0$. 
Denote $Y_0$ and $Z_0$ by $Y$ and $Z$, respectively. 
Because $g$ is log \'{e}tale, $g^{(p^n)*}_{}$ is log \'{e}tale.  
Hence $g^{(p^n)*}_{}(\Om^i_{Y^{(p^n)}/S})=\Om^i_{Z^{(p^n)}/S}$. 
This implies (\ref{ali:zzxn}) in the case $n=0$. 
The equality (\ref{ali:bbn}) is obvious in the case $n=0$. 
The point of the following proof is to prove (\ref{ali:bbn}) in the case $n=1$. 
(Though one may think that the proof of this is easy, it is not so easy.  
We have to be careful for some points of the proof.
See also (\ref{rema:idp}) below in the trivial log case.) 
\par 
We claim that the natural morphism 
\begin{align*} 
Z^{(p^{n})}\lo Z^{(p^{n+1})}\times_{Y^{(p^{n+1})},F_{Y^{(p^n)}/S}}Y^{(p^{n})}
\tag{2.1.16.3}\label{ali:osbxy}
\end{align*} 
obtained by 
the following commutative diagram 
\begin{equation*} 
\begin{CD} 
Z^{(p^{n})}@>{F_{Z^{(p^n)}/S}}>>Z^{(p^{n+1})}\\
@V{g^{(p^n)*}_{}}VV @VV{g^{(p^{n+1})}_{}}V\\
Y^{(p^{n})}@>{F_{Y^{(p^n)}/S}}>>Y^{(p^{n+1})}. 
\end{CD} 
\tag{2.1.16.4}\label{cd:ocsyfxy}
\end{equation*} 
is an isomorphism.  That is, we claim that the diagram 
(\ref{cd:ocsyfxy}) is cartesian.  
\par 
We prove that the morphism (\ref{ali:osbxy}) is log \'{e}tale and 
weakly purely inseparable in turn. 
Indeed, since  
the log \'{e}taleness is stable under base change, 
the morphism 
$Z^{(p^{n+1})}\times_{Y^{(p^{n+1})},F_{Y^{(p^n)}/S}}Y^{(p^{n})}
\lo Y^{(p^{n})}$ is log \'{e}tale. 
Because the morphism 
$Z^{(p^{n})}\lo Y^{(p^{n})}$ is log \'{e}tale 
and this is the composite morphism 
of (\ref{ali:osbxy}) and 
$Z^{(p^n)}\times_{Y^{(p^n)},F_{Y^{(p^n)}/S}}Y^{(p^{n+1})}\lo Y^{(p^{n+1})}$,  
the morphism (\ref{ali:osbxy}) is log \'{e}tale by (\ref{lemm:isoex}) (1). 
\par 
Let $W$ be $Z$ or $Y$. 
Because the weakly pure inseparability is stable under base change 
by \cite[III (2.4.8) 4]{ob} and because $F_S\col S\lo S$ is weakly purely inseparable, 
the base change morphism  
$W^{(p^{n+1})}\lo W^{(p^n)}$ of $F_S$ is weakly purely inseparable. 
By considering the composite morphism 
$W^{(p^n)}\os{F_{W^{(p^n)}/S}}{\lo} W^{(p^{n+1})}\lo W^{(p^n)}$, 
$F_{W^{(p^n)}/S}$ is weakly purely inseparable by \cite[III (2.8.4) 3]{ob}. 
Hence the base change morphism 
$Z^{(p^{n+1})}\times_{Y^{(p^{n+1})},F_{Y^{(p^n)}/S}}Y^{(p^{n})}\lo 
Z^{(p^{n+1})}$ of $F_{Y^{(p^n)}/S}\col Y^{(p^n)}\lo Y^{(p^{n+1})}$ 
is weakly purely inseparable.  
By considering the composite morphism 
$Z^{(p^n)}\os{(\ref{ali:osbxy})}{\lo} 
Z^{(p^{n+1})}\times_{Y^{(p^{n+1})},F_{Y^{(p^n)}/S}}Y^{(p^{n})}\lo Z^{(p^{n+1})}$, 
the morphism (\ref{ali:osbxy}) is weakly purely inseparable.  
\par 
Because the morphism of Cartier type is stable under base change, 
so is the structural morphism $Z^{(p^n)}\lo S$. 
Hence the morphism $F_{Z^{(p^n)}/S}\col Z^{(p^n)}\lo Z^{(p^{n+1})}$ 
is exact. 
\par 
In \cite[II (2.11)]{tsa} Tsuji has proved the following: 
\par 
$(1)$ The composite morphism of two saturated morphisms of integral log schemes is saturated. 
\par 
$(2)$ The saturated morphisms of 
integral log schemes are stable 
under the base change of integral log schemes. 
\par 
$(3)$ Let $h\col `Y\lo `Z$ be a morphism of 
$($fine$)$ integral log schemes. 
Then $h$ is saturated if and only if the base change 
$``Y$ of $`Y$ with respect to a morphism $``Z\lo `Z$ from 
a $($fine$)$ saturated log scheme are saturated. 
\par 
Hence $Y^{(p^n)}$ and 
$Z^{(p^{n+1})}\times_{Y^{(p^{n+1})},F_{Y^{(p^n)}/S}}Y^{(p^n)}$ are saturated. 
By (\ref{lemm:isoex}) (2) and by considering the composite morphism 
$$Z^{(p^n)}\lo Z^{(p^{n+1})}\times_{Y^{(p^{n+1})},F_{Y^{(p^n)}/S}}Y^{(p^n)}
\lo Z^{(p^{n+1})},$$
the morphism (\ref{ali:osbxy}) turns out to be an isomorphism. 
Consequently each of the following commutative diagram is cartesian. 
\begin{equation*} 
\begin{CD} 
Z@>{F_{Z/S}}>>Z^{(p)}@>{F_{Z^{(p)}/S}}>> \cdots 
@>{F_{Z^{(p^{n-1})}/S}}>>Z^{(p^{n})}
@>{F_{Z^{(p^{n})}/S}}>>Z^{(p^{n+1})}\\
@V{g}VV @V{g^{(p)}_{}}VV @. 
@VV{g^{(p^{n})}_{}}V @VV{g^{(p^{n+1})}_{}}V\\
Y@>{F_{Y/S}}>>Y^{(p)}@>{F_{Y^{(p)}/S}}>>\cdots 
@>{F_{Y^{(p^{n-1})}/S}}>>
Y^{(p^{n})}
@>{F_{Y^{(p^{n})}/S}}>>Y^{(p^{n+1})}. 
\end{CD} 
\tag{2.1.16.5}\label{cd:osyfaxy}
\end{equation*} 

\par 
Next we claim that the natural morphism 
\begin{align*}  
&g^{(p^{n+1})*}_{}(F_{Y^{(p^n)}/S*}(\Om^i_{Y^{(p^n)}/S}))
\owns f\otimes \om 
\lom F_{Z^{(p^n)}/S}^*(f) g^{(p^n)*}(\om) \in  
F_{Z^{(p^n)}/S*}(\Om^i_{Z^{(p^n)}/S}) 
\tag{2.1.16.6}\label{ali:osmxy}\\
&(f\in {\cal O}_{Y^{(p^n)}}, \om \in \Om^i_{Y^{(p^n)}/S})
\end{align*} 
is an isomorphism. 
Indeed, because this is a local problem, 
because $\Om^i_{Y^{(p^n)}/S}$ is locally free and because 
$\Om^i_{Z^{(p^n)}/S}=g^{(p^n)*}_{}(\Om^i_{Y^{(p^n)}/S})$,  
it suffices to prove that 
the natural morphism 
\begin{align*} 
{\cal O}_{Z^{(p^{n+1})}}\otimes_{{\cal O}_{Y^{(p^{n+1})}}}
F_{Y^{(p^n)}/S*}({\cal O}_{Y^{(p^{n})}})
\lo 
F_{Z^{(p^n)}/S*}({\cal O}_{Z^{(p^{n})}})
\end{align*}  
is an isomorphism. 
Because (\ref{ali:osbxy}) is an isomorphism, this is an isomorphism. 
Similarly we see that the following natural morphism 
\begin{align*}  
g^{(p^{n+1})*}_{}(F^n_{Y^{(p)}/S*}(\Om^i_{Y^{(p)}/S}))
\lo 
F^n_{Z^{(p)}/S*}(\Om^i_{Z^{(p)}/S})
\tag{2.1.16.7}\label{ali:osmafxy}
\end{align*} 
is an isomorphism. 
More generally, for a locally free ${\cal O}_{Y^{(p)}}$-module ${\cal F}$, 
the following natural morphism 
\begin{align*}  
g^{(p^{n+1})*}_{}(F^n_{Y^{(p)}/S*}({\cal F}))
\lo 
F^n_{Z^{(p)}/S*}(g^{(p)*}_{}({\cal F}))
\tag{2.1.16.8}\label{ali:osmfxy}
\end{align*} 
is an isomorphism. 
\par 
Since the following diagram 
\begin{equation*} 
\begin{CD} 
g^{(p^{n+1})*}_{}(F_{Y^{(p^n)}/S*}{\cal H}^i(\Om^{\bul}_{Y^{(p^n)}/S}))
@>>> F_{Z^{(p^n)}/S*}{\cal H}^i(\Om^{\bul}_{Z^{(p^n)}/S})\\
@A{g^{(p^{n+1})*}_{}(C^{-1}_{Y^{(p^n)}/S})}A{\simeq}A 
@A{\simeq}A{C^{-1}_{Z^{(p^n)}/S}}A \\
g^{(p^{n+1})}_{*}(\Om^i_{Y^{(p^{n+1})}/S})@>{\simeq}>>\Om^i_{Z^{(p^{n+1})}/S}
\end{CD} 
\tag{2.1.16.9}\label{cd:osxy}
\end{equation*} 
is commutative, 
the upper horizontal morphism 
in (\ref{cd:osxy}) is an isomorphism of 
${\cal O}_{Z^{(p^{n+1})}}$-modules: 
\begin{align*} 
g^{(p^{n+1})*}_{}(F_{Y^{(p^n)}/S*}
{\cal H}^i(\Om^{\bul}_{Y^{(p^n)}/S}))
\os{\sim}{\lo} F_{Z^{(p^n)}/S*}({\cal H}^i(\Om^{\bul}_{Z^{(p^n)}/S})). 
\tag{2.1.16.10}\label{ali:oxy}
\end{align*}  
\par 
We prove the equalities (\ref{ali:zzxn}) and (\ref{ali:bbn}) by induction on $n$. 

Consider the complex 
$F_{Y^{(p^n)}/S*}(\Om^{\bul}_{Y^{(p^n)}/S})$ of 
${\cal O}_{Y^{(p^{n+1})}}$-modules.  Note that 
the derivative 
$F_{Y^{(p^n)}/S*}(\Om^i_{Y^{(p^n)}/S})\lo 
F_{Y^{(p^n)}/S*}(\Om^{i+1}_{Y^{(p^n)}/S})$ 
is indeed an ${\cal O}_{Y^{(p^{n+1})}}$-linear morphism. 
Especially consider this complex for the case $n=0$. 
Though the following diagram 
\begin{equation*} 
\begin{CD} 
g^*(\Om^{i-1}_{Y/S})
@>{\sim}>> \Om^{i-1}_{Z/S}\\
@V{{\rm id}_{{\cal O}_Z}\otimes d}VV @VV{d}V \\
g^{*}(\Om^i_{Y/S})@>{\sim}>> \Om^i_{Z/S}
\end{CD} 
\tag{2.1.16.11}\label{cd:aobsxy}
\end{equation*} 
is not commutative at all when $Z\not=Y$ in general, 
the following diagram is commutative 
\begin{equation*} 
\begin{CD} 
g^{(p^n)*}_{}(F^n_{Y/S*}(\Om^{i-1}_{Y/S}))
@>{\sim}>> F^n_{Z/S*}(\Om^{i-1}_{Z/S})\\
@V{g^{(p^n)*}_{}F^n_{Y/S*}(d)}VV @VV{F^n_{Z/S*}(d)}V \\
g^{(p^n)*}_{}(F^n_{Y/S*}(\Om^i_{Y/S}))
@>{\sim}>> F^n_{Z/S*}(\Om^i_{Z/S}) 
\end{CD} 
\tag{2.1.16.12}\label{cd:aosxy}
\end{equation*} 
for $n\in {\mab Z}_{\geq 1}$. 
Because the underlying morphism of 
$g^{(p)}$ is not necessarily flat, 
we need the following argument to prove 
that $g^{(p)*}(B_1\Om^i_{Y/S})=B_1\Om^i_{Z/S}$ ((\ref{ali:bbficn}) below). 
Since $\os{\circ}{F}_{Y/S}\col \os{\circ}{Y}\lo \os{\circ}{Y}{}^{(p)}$ is a homeomorphism, 
the following sequence 
\begin{align*} 
0\lo Z_1\Om^i_{Y/S}\lo F_{Y/S*}(\Om^i_{Y/S})
\os{F_{Y/S*}(d)}{\lo} B_1\Om^{i+1}_{Y/S}\lo 0
\end{align*}   
is exact.  
Because $Z_1\Om^i_{Y/S}/B_1\Om^i_{Y/S}=
F_{Y/S*}{\cal H}^i(Z\Om^i_{Y/S}/B\Om^i_{Y/S})
=\Om^i_{Y^{(p)}/S}$ 
is a locally free ${\cal O}_{Y^{(p)}}$-modules
and 
because $B_1\Om^{i+1}_{Y/S}$ is a locally free 
${\cal O}_{Y^{(p)}}$-modules by \cite[(1.13)]{lodw}, 
we have the following exact sequences of 
${\cal O}_{Z^{(p)}}$-modules: 
\begin{align*} 
0\lo g^{(p)*}_{}(Z_1\Om^i_{Y/S}) \lo g^{(p)*}_{}F_{Y/S*}(\Om^i_{Y/S})
\lo g^{(p)*}_{}(B_1\Om^{i+1}_{Y/S})\lo 0
\tag{2.1.16.13}\label{ali:bbban} 
\end{align*} 
and 
\begin{align*} 
0\lo g^{(p)*}_{}(B_1\Om^{i+1}_{Y/S}) \lo g^{(p)*}_{}(Z_1\Om^{i+1}_{Y/S})
\lo g^{(p)*}_{}(Z_1\Om^i_{Y/S}/B_1\Om^i_{Y/S})\lo 0. 
\tag{2.1.16.14}\label{ali:bbbcn} 
\end{align*} 
Hence 
\begin{align*} 
g^{(p)*}_{}(B_1\Om^{i+1}_{Y/S})
={\rm Im}(g^{(p)*}_{}F_{Y/S*}(d)\col g^{(p)*}_{}F_{Y/S*}(\Om^i_{Y/S})
\lo g^{(p)*}_{}F_{Y/S*}(\Om^{i+1}_{Y/S})).
\tag{2.1.16.15}\label{ali:bbbicn} 
\end{align*}  
Because the morphism (\ref{ali:osmxy}) is an isomorphism 
and because (\ref{cd:aosxy}) is commutative, 
the equality (\ref{ali:bbbicn})  is nothing but the following equality: 
\begin{align*} 
g^{(p)*}_{}(B_1\Om^{i+1}_{Y/S})
&={\rm Im}(F_{Z/S*}(d)\col F_{Z/S*}(\Om^i_{Z/S})\lo F_{Z/S*}(\Om^{i+1}_{Z/S}))
\tag{2.1.16.16}\label{ali:bbficn} \\
&=B_1\Om^{i+1}_{Z/S}.
\end{align*}  
Thus we have proved the desired equality. 
Since  
\begin{align*} 
C^{-1}_{W^{(p)}/S}\col 
Z_n\Om^i_{W^{(p)}/S}
\os{\sim}{\lo} 
Z_{n+1}\Om^i_{W/S}/F^{n}_{W^{(p)}/S*}(B_1\Om^i_{W/S}), 
\end{align*}
we have the following isomorphism for $n\in {\mab N}$: 
\begin{align*} 
&g^{(p^{n+1})*}(C^{-1}_{Y^{(p)}/S})\col 
g^{(p^{n+1})*}(Z_n\Om^i_{Y^{(p)}/S})
\os{\sim}{\lo} 
g^{(p^{n+1})*}_{}(Z_{n+1}\Om^i_{Y/S}/
F^{n}_{Y^{(p)}/S*}(B_1\Om^i_{Y/S})) \\
&\os{\sim}{\lo} 
g^{(p^{n+1})*}_{}(Z_{n+1}\Om^i_{Y/S})/
{\rm Im}(g^{(p^{n+1})*}
(F^{n}_{Y^{(p)}/S*}(B_1\Om^i_{Y/S}))\lo g^{(p^{n+1})*}_{}(Z_{n+1}\Om^i_{Y/S}))\\
&\os{\sim}{\lo} 
g^{(p^{n+1})*}_{}(Z_{n+1}\Om^i_{Y/S})/
F^{n}_{Z^{(p)}/S*}g^{(p)*}_{}(B_1\Om^i_{Y/S})\\
&=g^{(p^{n+1})*}_{}(Z_{n+1}\Om^i_{Y/S})/
F^{n}_{Z^{(p)}/S*}(B_1\Om^{i+1}_{Z/S}). 
\end{align*}
Here we have used the isomorphism (\ref{ali:osmfxy}), 
(\ref{cd:aosxy}) and 
the equality (\ref{ali:bbficn}). 
Induction on $n$ shows that (\ref{ali:zzxn}) holds. 
Similarly induction on $n$ shows that (\ref{ali:bbn}) holds. 
\end{proof} 

\begin{rema}\label{rema:idp}  
(1) Let the notations be as in \cite[0 (2.2.7)]{idw}. 
In particular, let $f\col X'\lo X$ be an \'{e}tale morphism of smooth schemes over 
a scheme $S$ of characteristic $p>0$. 
In [loc.~cit.] Illusie has proved that 
the following morphism 
\begin{align*} 
f^{(p^n)*}(Z_n\Om^i_{X/S}/B_n\Om^i_{X/S})
\lo 
Z_n\Om^i_{X'/S}/B_n\Om^i_{X'/S}
\tag{2.1.17.1}\label{ali:zzbbxn}
\end{align*} 
is an isomorphism. However he has not proved 
the following morphisms 
\begin{align*} 
f^{(p^n)*}(B_n\Om^i_{X/S})\lo B_n\Om^i_{X'/S}
\tag{2.1.17.2}\label{ali:bi}
\end{align*} 
and 
\begin{align*} 
f^{(p^n)*}(Z_n\Om^i_{X/S})\lo Z_n\Om^i_{X'/S}
\tag{2.1.17.3}\label{ali:iz}
\end{align*} 
are isomorphisms. 
I do not think that  
the claim in \cite[p.~320]{idw} that the isomorphism 
$f^*(\Om^i_{X/S}) \os{\sim}{\lo}  \Om^1_{X'/S}$ 
induces the isomorphisms (\ref{ali:bi}) and (\ref{ali:iz}) 
is true 
even if one uses the following commutative diagram 
(cf.~[loc.~cit., 0 (2.1.8)]): 
\begin{equation*} 
\begin{CD} 
\Om^i_{X'{}^{(p)}/S}@>{C_{X'/S},~\sim}>>
F_{X'/S*}{\cal H}^i(\Om^{\bul}_{X'/S})\\
@| @|\\
f^{(p)*}(\Om^i_{X^{(p)}/S})@>{f^{(p)*}(C^{-1}_{X/S}),~\sim}>>
f^{(p)*}F_{X/S*}{\cal H}^i(\Om^{\bul}_{X/S}). 
\end{CD} 
\end{equation*} 
(This commutative diagram is a special case of (\ref{cd:osxy}).) 
I think that the isomorphism $f^*(\Om^i_{X/S}) \os{\sim}{\lo}  \Om^1_{X'/S}$ 
is not directly useful for the proof of the isomorphism
(\ref{ali:bi}) and (\ref{ali:iz}). 
\par 
Let the notations be as in the proof of (\ref{theo:lgetb}). 
The point of the proof of (\ref{theo:lgetb}) is to prove that 
$g^{(p)*}_{}(B_1\Om^{i+1}_{Y/S}) =B_1\Om^{i+1}_{Z/S}$. 
The corresponding statement in the trivial log case 
has not been proved in \cite[p.~520]{idw}. 
To prove it, we need the 
corresponding commutative diagram to (\ref{cd:aosxy}), 
which has not appeared in [loc.~cit.]. 
To prove this, we have only to show that the morphism 
\begin{align*} 
X'{}^{(p^{n})}\lo X'{}^{(p^{n+1})}\times_{Z^{(p^{n+1})},F_{X^{(p^n)}/S}}X^{(p^{n})}
\end{align*} 
is an isomorphism by using \cite[(17.9.1)]{ega4}, 
where $X'$ and $X$ are schemes over 
$S$ in [loc.~cit.].  
\par 
(2) We do not know whether the assumptions  that $Y_{\bul \leq N}$ and $S$ 
are saturated and that the morphisms $Y_{\bul \leq N}$ and 
$g_{}$ is saturated are indispensable for the conclusion of (\ref{theo:lgetb}). 
\end{rema}

\begin{prop}\label{prop:eztb} 
Let $Y_{\bul \leq N}/S$ be 
a log smooth $N$-truncated simplicial log scheme of Cartier type. 
Let $n$ and $i$ be nonnegative integers. Then the following hold$:$ 
\par 
$(1)$  Let $g\col Z_{\bul \leq N}\lo Y_{\bul \leq N}$ be a solid morphism from 
an $N$-truncated simplicial fine log scheme over $S$. 
Then the log scheme $Z_{\bul \leq N}$ is log smooth of Cartier type over $S$. 
\par 
$(2)$ Let ${\cal O}_{Y_{{\bul \leq N},{\rm et}}}$ and 
$(\Om^i_{Y_{\bul \leq N}/S})_{\rm et}$ 
be the \'{e}tale sheaves in the small \'{e}tale topos 
$Y_{{\bul \leq N},\rm et}$ defined by ${\cal O}_{Y_{\bul \leq N}}$ 
and $\Om^i_{Y_{\bul \leq N}/S}$, respectively. 
Then the sheaves $Z_n\Om^i_{Y_{\bul \leq N}/S}$ and $B_n\Om^i_{Y_{\bul \leq N}/S}$ 
extend to \'{e}tale subsheaves of 
$(F_{Y_{\bul \leq N}/S})^n_{\rm et}({\cal O}_{Y_{\bul \leq N}})$-modules 
$(\Om^i_{Y_{\bul \leq N}/S})_{\rm et}$. 
\end{prop} 
\begin{proof} 
(1): 
We may assume that $N=0$. We set $Z:=Z_0$ and $Y:=Y_0$. 
The structural morphism $Z\lo S$ is obviously integral. 
The morphism $F_{Z^{(p)}/S}^*(M_{Z^{(p)}})\lo M_Z$ is equal to 
the morphism 
$F_{Z^{(p)}/S}^*g^{(p)}{}^*(M_{Y^{(p)}})=g^*F_{Y/S}^*(M_{Y^{(p)}})\lo g^*(M_Y)$. 
Since the morphism $F_{Y/s}^*(M_{Y'})\lo M_Y$ is exact, 
we see that the morphism $g^*F_{Y/s}^*(M_{Y'})\lo g^*(M_Y)$ 
is exact. 
\par 
(2): By (1) the morphism $g\col Z\lo S$ is log smooth of Cartier type. 
By the proof of (\ref{theo:lgetb}) 
the morphism $Z^{(p^n)}\lo 
Z^{(p^{n+1})}\times_{Y^{(p^{n+1})},F_{Y^{(p^n)}/S}}Y^{(p^{n})}$ 
is weakly purely inseparable. 
Hence it is an isomorphism by (\ref{lemm:isoex}) (2). 
By (\ref{theo:lgetb}) we obtain the following equalities: 
\begin{align*} 
g^{(p^n)*}_{}(Z_n\Om^i_{Y_{\bul \leq N}/S})
=Z_n\Om^i_{Z_{\bul \leq N}/S}
\end{align*} 
and 
\begin{align*} 
g^{(p^n)*}_{}(B_n\Om^i_{Y_{\bul \leq N}/S})
=B_n\Om^i_{Z_{\bul \leq N}/S}. 
\end{align*} 
These show (\ref{prop:eztb}).  
\end{proof} 

\begin{rema}\label{rema:sr}
It is surprising to me that (\ref{prop:eztb}) (2) 
(even in the trivial log case) has not been proved 
in any other references as far as I know.  
\end{rema}

Let us recall the following which has been proved in \cite[(3.10)]{nlfc}:

\begin{theo}\label{prop:fzg}
Let $Z$ be a log smooth scheme of Cartier type over $s$. 
Let $F \col {\cal W}_n(Z) \lo {\cal W}_n(Z)$ be 
the Frobenius endomorphism of ${\cal W}_n(Z)$. 
Then the following hold$:$
\par 
$(1)$ 
The morphism 
$F^n\col {\cal W}_{n+1}\Om^i_Z\lo {\cal W}_1\Om^i_Z$ 
factors through $Z_n{\cal W}_1\Om^i_Z$ and 
the following sequence 
\begin{align*} 
F_{*}({\cal W}_{n}\Om^i_Z)
\os{V}{\lo} {\cal W}_{n+1}\Om^i_Z
\os{F^n}{\lo} 
Z_n{\cal W}_1\Om^i_Z
\lo 0
\tag{2.1.20.1}\label{ali:mvee}
\end{align*} 
is an exact sequence of 
${\cal W}_{n+1}({\cal O}_Z)$-modules. 
\par 
$(2)$  The morphism 
$F_*(F^{n-1}d)\col {\cal W}_{n+1}\Om^i_Z\lo {\cal W}_1\Om^{i+1}_Z$ 
factors through $B_n{\cal W}_1\Om^{i+1}_Z$ and 
the following sequence 
\begin{align*} 
{\cal W}_{n+1}\Om^i_Z\os{F}{\lo} F_{*}({\cal W}_n\Om^i_Z)
\os{F_*(F^{n-1}d)}{\lo} 
B_n{\cal W}_1\Om^{i+1}_Z
\lo 0
\tag{2.1.20.2}\label{ali:mfee}
\end{align*}
is an exact sequence of 
${\cal W}_n({\cal O}_Z)$-modules. 
\end{theo}

The following is a log version of \cite[I~(1.14)]{idw}; 
the proof of it is different from the proof of [loc.~cit.]: 

\begin{theo}\label{theo:bcn} 
Let $Z_{\bul \leq N}$, $Y_{\bul \leq N}$ and 
$g$ be as in {\rm (\ref{theo:lgetb})}. 
Assume that $S=s$. 
Then, for $i,n\in {\mab N}$, 
\begin{align*} 
g^*({\cal W}_n\Om^i_{Y_{\bul \leq N}/s})={\cal W}_n\Om^i_{Z_{\bul \leq N}/s}. 
\tag{2.1.21.1}\label{ali:zznxn}
\end{align*} 
\end{theo}
\begin{proof} 
By \cite[(1.13)]{lodw} $Z_n\Om^{i}_{Y/s}$ is 
a locally free ${\cal O}_{Y^{(p^n)}}$-module. Hence 
$Z_n{\cal W}_1\Om^{i}_{Y/s}$ is also a locally free ${\cal O}_{Y^{(p^n)}}$-module. 
Because the natural morphism  
${\cal O}_Y\lo {\cal O}_{Y^{(p^n)}}$ is an isomorphism,  
(\ref{theo:bcn}) 
immediately follows from (\ref{theo:lgetb}), (\ref{prop:fzg}) (1) 
and induction on $n$. 
\end{proof}

In (\ref{rema:qok}) (3) below we shall use the argument in the proof of 
the following proposition. 

\begin{coro}\label{coro:lemos}  
Let $Y_{\bul \leq N}$ be a log smooth $N$-truncated simplicial log scheme over $s$. 
Let $g\col Z_{\bul \leq N}\lo Y_{\bul \leq N}$ be a solid morphism from 
an $N$-truncated simplicial fine log scheme over $s$. 
Assume that $\os{\circ}{g}$ is \'{e}tale.  
Then the canonical morphism 
\begin{align*} 
{\cal W}_n({\cal O}_{Z_{\bul \leq N}})^{\star}
\otimes_{{\cal W}_n({\cal O}_{Y_{\bul \leq N}})^{\star}}
({\cal W}_n{\Om}_{Y_{\bul \leq N}}^i)^{\star} 
\lo ({\cal W}_n{\Om}_{Z_{\bul \leq N}}^i)^{\star} \quad (i\in {\mab N}) 
\tag{2.1.22.1}\label{ali:ozs}
\end{align*} 
is an isomorphism. 
\end{coro}
\begin{proof}
We may assume that $N=0$. 
Set $Y:=Y_0$ and $Z:=Z_0$. 
It suffices to prove (\ref{coro:lemos}) in the case $\star=$nothing. 
By (\ref{prop:eztb}) we see that this is a special case of (\ref{ali:zznxn}). 
\end{proof} 


\begin{rema}\label{rema:qc} 
(1) Consider the case $N=0$ for simplicity. 
Let $S$ be a fine log scheme such that $\os{\circ}{S}$ is a perfect scheme of 
characteristic $p>0$. 
Let $h\col Z\lo S$ be a log smooth scheme of Cartier type. 
Let $(Z/{\cal W}_n(S))_{\rm crys}$ be the log crystalline topos 
defined in \cite{klog1} and 
let $u_{Z/{\cal W}_n(S),{\rm et}}\col (Z/{\cal W}_n(S))_{\rm crys}\lo 
\os{\circ}{Z}_{\rm et}$ be the canonical projection. 
Then we obtain the \'{e}tale sheaf 
${\cal W}_n\Om^i_{Z,{\rm et}}
:={\cal H}^i(Ru_{Z/{\cal W}_n(S),{\rm et}*}({\cal O}_{Z/{\cal W}_n(S)}))$ 
$(i\in {\mab N})$.  
This is only an $h^{-1}_{\rm et}({\cal W}_n({\cal O}_{S_{\rm et}}))$-module. 
However, as is well-known, ${\cal W}_n\Om^i_{Z,{\rm et}}$ 
is a ${\cal W}_n({\cal O}_Z)'_{\rm et}$-module, 
where ${\cal W}_n({\cal O}_Z)'_{\rm et}$ is an associated quasi-coherent sheaf 
in $Z_{\rm et}$ to ${\cal W}_n({\cal O}_Z)$. 
However it is not clear that ${\cal W}_n\Om^i_{Z,{\rm et}}$ is indeed a quasi-coherent 
${\cal W}_n({\cal O}_Z)'_{\rm et}$-module. 
By (\ref{coro:lemos}) we see that ${\cal W}_n\Om^i_{Z,{\rm et}}$ is a quasi-coherent 
${\cal W}_n({\cal O}_Z)'_{\rm et}$-module. 
\par 
(2) Though the same statement as (1) has been stated in \cite[(1.17)]{lodw}, 
the proof of \cite[(1.17)]{lodw} is not complete because 
(\ref{prop:eztb}) (2) has not been proved in [loc.~cit.].  
\par 
(3) Let the notations be as in (1). 
The proof of \cite[(2.12.2)]{nh2} is not complete because 
we have used a fact that ${\cal W}_n\Om^i_{Z,{\rm et}}$ is a quasi-coherent 
${\cal W}_n({\cal O}_Z)'_{\rm et}$-module in \cite[(2.12.2.4)]{nh2}. 
If one wants to avoid using the isomorphism (\ref{ali:ozs}) to obtain 
the isomorphism in \cite[(2.12.2.4)]{nh2}, 
then one has to use the Poincar\'{e} residue isomorphism 
for the log de Rham-Witt sheaf of a smooth scheme with an SNCD over $\kap$
and use the comparison theorem 
${\cal W}_n\Om^i_Y\os{\sim}{\lo} {\cal W}_n\Om^i_Y$ of 
{\it ${\cal W}_n({\cal O}_Y)$-modules} 
for a smooth scheme $Y/\kap$ (\cite[p.~142]{ir}), 
where the source (resp. target) of the isomorphism above 
is the de Rham-Witt sheaf defined in \cite{idw} (resp. \cite{ir} and \cite{hk}) 
and use the isomorphism \cite[I (1.14.1)]{idw}. 
\end{rema}

\par 
Let $F_{Y_{\bul \leq N}}\col Y_{\bul \leq N}\lo Y_{\bul \leq N}$ be the 
absolute Frobenius endomorphism(=the $p$-th power endomorphism) of 
$Y_{\bul \leq N}$.

\begin{defi}\label{defi:nfc} 
Let $E^{\bul \leq N}$ be 
a coherent crystal on $Y_{\bul \leq N}/{\cal W}(s)$. 
Let 
$$\Phi_{\rm crys} \col 
F_{Y_{\bul \leq N},{\rm crys}}^*(E^{\bul \leq N})\lo E^{\bul \leq N}$$  
be an ${\cal O}_{Y_{\bul \leq N}/{\cal W}(s)}$-linear morphism. 
We call the pair $(E^{\bul \leq N},\Phi_{\rm crys})$ 
an {\it $F$-crystal} on $Y_{\bul \leq N}/{\cal W}(s)$. 
\end{defi}  

\parno 
Let $E^{\bul \leq N}$ be a flat coherent $F$-crystal 
on $Y_{\bul \leq N}/{\cal W}(s)$.  
Let $F\col {\cal W}_n(Y_{\bul \leq N})
\lo {\cal W}_n(Y_{\bul \leq N})$ 
be the Frobenius endomorphism of ${\cal W}_n(Y_{\bul \leq N})$. 
Set $(E^{\bul \leq N}_n)^{\sig}:=F^*(E^{\bul \leq N}_n)$. 
By the commutative diagram
\begin{equation*} 
\begin{CD} 
Y_{\bul \leq N}@>{F_{Y_{\bul \leq N}}}>>Y_{\bul \leq N}\\
@V{\bigcap}VV @VV{\bigcap}V \\
{\cal W}_n(Y_{\bul \leq N})@>{F}>>{\cal W}_n(Y_{\bul \leq N}), 
\end{CD} 
\tag{2.1.24.1}\label{cd:fyu}
\end{equation*} 
we have the following morphism: 
the ${\cal W}_n({\cal O}_{Y_{\bul \leq N}})$-linear morphism 
$\Phi \col (E^{\bul \leq N}_n)^{\sig}\lo E^{\bul \leq N}_n$. 
As in \cite[II (3.0)]{et} we have an $F^*$-linear morphism 
$$F_{E^{\bul \leq N}} \col 
E^{\bul \leq N}_n\lo (E^{\bul \leq N}_n)^{\sig}\lo E^{\bul \leq N}_n,$$ 
where $E^{\bul \leq N}_n\lo (E^{\bul \leq N}_n)^{\sig}$ 
is the natural morphism $x\lom 1\otimes x$ 
$(x\in E^{\bul \leq N}_n)$. 
Let $F\col E^{\bul \leq N}_{n+1} \lo E^{\bul \leq N}_n$ 
be the composite morphism 
$E^{\bul \leq N}_{n+1}\os{F_{E^{\bul \leq N}}}{\lo} E^{\bul \leq N}_{n+1} 
\os{\rm proj.}{\lo} E^{\bul \leq N}_n$. 
Set 
$$F:=F\otimes F\col E^{\bul \leq N}_{n+1}
\otimes_{{\cal W}_{n+1}({\cal O}_{Y_{\bul \leq N}})^{\star}}
({\cal W}_{n+1}{\Om}^i_{Y_{\bul \leq N}})^{\star}\lo 
E^{\bul \leq N}_n
\otimes_{{\cal W}_n({\cal O}_{Y_{\bul \leq N}})^{\star}}
({\cal W}_n{\Om}^i_{Y_{\bul \leq N}})^{\star}.$$ 
\par 
Now assume that $E^{\bul \leq N}$ is a unit root log $F$-crystal, i.~e., 
$\Phi_{\rm crys}$ is an isomorphism. In particular, 
$\Phi$ is an isomorphism. 
First consider the following morphism 
\begin{equation*} 
\Phi^{-1}\otimes {\rm id} \col 
E^{\bul \leq N}_n
\otimes_{{\cal W}_n({\cal O}_{Y_{\bul \leq N}})^{\star}}
({\cal W}_n\Om^i_{Y_{\bul \leq N}})^{\star}\lo 
(E^{\bul \leq N}_n)^{\sig}
\otimes_{{\cal W}_n({\cal O}_{Y_{\bul \leq N}})^{\star}}
({\cal W}_n{\Om}^i_{Y_{\bul \leq N}})^{\star}. 
\end{equation*} 
Next consider the following isomorphism 
\begin{align*} 
I\col & (E^{\bul \leq N}_n)^{\sig}
\otimes_{{\cal W}_n({\cal O}_{Y_{\bul \leq N}})^{\star}}
({\cal W}_n{\Om}^i_{Y_{\bul \leq N}})^{\star} 
= F^{-1}(E^{\bul \leq N}_n)
\otimes_{F^{-1}({\cal W}_n({\cal O}_{Y_{\bul \leq N}})^{\star})} 
({\cal W}_n{\Om}^i_{Y_{\bul \leq N}})^{\star} \\
&=E^{\bul \leq N}_n
\otimes_{{\cal W}_n({\cal O}_{Y_{\bul \leq N}})^{\star}}
F_*({\cal W}_n{\Om}^i_{Y_{\bul \leq N}})^{\star}. 
\end{align*} 
Consider also the following morphism 
\begin{align*} 
{\rm id}\otimes V \col & 
E^{\bul \leq N}_n
\otimes_{{\cal W}_n({\cal O}_{Y_{\bul \leq N}})^{\star}}
F_*({\cal W}_n{\Om}^i_{Y_{\bul \leq N}})^{\star} \\
& =
(E^{\bul \leq N}_{n+1}
\otimes_{{\cal W}_{n+1}({\cal O}_{Y_{\bul \leq N}})^{\star},R}
{\cal W}_n({\cal O}_{Y_{\bul \leq N}})^{\star})
\otimes_{{\cal W}_n({\cal O}_{Y_{\bul \leq N}})^{\star}} 
F_*({\cal W}_n{\Om}^i_{Y_{\bul \leq N}})^{\star} \\
& \lo 
E^{\bul \leq N}_{n+1}
\otimes_{{\cal W}_{n+1}({\cal O}_{Y_{\bul \leq N}})^{\star}}
({\cal W}_{n+1}{\Om}^i_{Y_{\bul \leq N}})^{\star}. 
\end{align*} 
Here we have used the formula ``$xVy=V(Fxy)$'' 
to define ${\rm id}\otimes V$. 
Let 
\begin{equation*} 
V\col E^{\bul \leq N}_n
\otimes_{{\cal W}_n({\cal O}_{Y_{\bul \leq N}})^{\star}}
({\cal W}_n{\Om}^i_{Y_{\bul \leq N}})^{\star}\lo 
E^{\bul \leq N}_{n+1}
\otimes_{{\cal W}_{n+1}({\cal O}_{Y_{\bul \leq N}})^{\star}}
({\cal W}_{n+1}{\Om}^i_{Y_{\bul \leq N}})^{\star} 
\end{equation*}
be the composite of $\Phi^{-1}\otimes {\rm id}$, 
$I$ and ${\rm id}\otimes V$: 
$V:= ({\rm id}\otimes V)\circ I\circ 
(\Phi^{-1}\otimes {\rm id})$. 
Now we obtain the graded pro-module 
$\{E^{\bul \leq N}_n\otimes_{{\cal W}_n({\cal O}_{Y_{\bul \leq N}})^{\star}}
({\cal W}_n\Om^{\bul}_{Y_{\bul \leq N}})^{\star}\}_{n=1}^{\infty}$ 
over the Raynaud ring ${\mab R}$ of $\kap$. That is, 
$E^{\bul \leq N}_n\otimes_{{\cal W}_n({\cal O}_{Y_{\bul \leq N}})^{\star}}
({\cal W}_n{\Om}^i_{Y_{\bul \leq N}})^{\star}$ is 
a sheaf of ${\cal W}_n({\cal O}_{Y_{\bul \leq N}})^{\star}$-modules  
and we have the following operators 
satisfying the standard relations in \cite[I (1.1)]{ir} 
and \cite[(1.3.1)]{hdw}: 
\begin{align*} 
& F \col E^{\bul \leq N}_{n+1}
\otimes_{{\cal W}_{n+1}({\cal O}_{Y_{\bul \leq N}})^{\star}}
({\cal W}_{n+1}{\Om}_{Y_{\bul \leq N}}^{\bul})^{\star}
\lo 
E^{\bul \leq N}_n\otimes_{{\cal W}_n({\cal O}_{Y_{\bul \leq N}})^{\star}}
({\cal W}_n{\Om}_{Y_{\bul \leq N}}^{\bul})^{\star}, \tag{2.1.24.2}\label{eqn:hyu}\\ 
& V\col E^{\bul \leq N}_n
\otimes_{{\cal W}_n({\cal O}_{Y_{\bul \leq N}})^{\star}} 
({\cal W}_n{\Om}_{Y_{\bul \leq N}}^{\bul})^{\star} \lo 
E^{\bul \leq N}_{n+1}
\otimes_{{\cal W}_{n+1}({\cal O}_{Y_{\bul \leq N}})^{\star}}
({\cal W}_{n+1}{\Om}_{Y_{\bul \leq N}}^{\bul})^{\star}, \\
& \nabla \col E^{\bul \leq N}_n
\otimes_{{\cal W}_n({\cal O}_{Y_{\bul \leq N}})^{\star}}
({\cal W}_n{\Om}_{Y_{\bul \leq N}}^{\bul})^{\star} 
\lo E^{\bul \leq N}_n\otimes_{{\cal W}_n({\cal O}_{Y_{\bul \leq N}})^{\star}}
({\cal W}_n{\Om}_{Y_{\bul \leq N}}^{\bul +1})^{\star}, \\
& {\bf p} \col E^{\bul \leq N}_n
\otimes_{{\cal W}_n({\cal O}_{Y_{\bul \leq N}})^{\star}}
({\cal W}_n{\Om}_{Y_{\bul \leq N}}^{\bul})^{\star} \lo 
E^{\bul \leq N}_{n+1}
\otimes_{{\cal W}_{n+1}({\cal O}_{Y_{\bul \leq N}})^{\star}}
({\cal W}_{n+1}{\Om}_{Y_{\bul \leq N}}^{\bul})^{\star}, \\
& R \col E^{\bul \leq N}_{n+1}
\otimes_{{\cal W}_{n+1}({\cal O}_{Y_{\bul \leq N}})^{\star}}
({\cal W}_{n+1}{\Om}_{Y_{\bul \leq N}}^{\bul})^{\star}
\lo 
E^{\bul \leq N}_n
\otimes_{{\cal W}_n({\cal O}_{Y_{\bul \leq N}})^{\star}}
({\cal W}_n{\Om}_{Y_{\bul \leq N}}^{\bul})^{\star}. 
\end{align*}
(See \cite[II (3.1.3)]{et} for the proof of the relation 
$\nabla F=pF\nabla$; this relation implies the relation 
$F\nabla V=\nabla$ because $E^m_n$ $(0\leq m\leq N)$ 
is generated by $F_{E^m}(E^m_n)$ 
as in \cite[p.86]{et}.)
We also have the following Frobenius morphism: 
\begin{equation*} 
\Phi_n  \col E^{\bul \leq N}_n
\otimes_{{\cal W}_n({\cal O}_{Y_{\bul \leq N}})^{\star}}
({\cal W}_n{\Om}_{Y_{\bul \leq N}}^{\bul})^{\star} \lo 
E^{\bul \leq N}_n
\otimes_{{\cal W}_n({\cal O}_{Y_{\bul \leq N}})^{\star}}
({\cal W}_n{\Om}_{Y_{\bul \leq N}}^{\bul})^{\star}. 
\tag{2.1.24.3}\label{eqn:wlphyu}
\end{equation*}
 
As in \cite[II (3.2.3), (3.3.1), (3.4.1), (4.2.1)]{et}, 
we obtain the following four results$:$ 
\begin{prop} 
\begin{align*} 
& {\rm Fil}^n(E^{\bul \leq N}_{n+r}
\otimes_{{\cal W}_{n+r}({\cal O}_{Y_{\bul \leq N}})^{\star}}
({\cal W}_{n+r}{\Om}^i_{Y_{\bul \leq N}})^{\star})
= \tag{2.1.25.1}\label{eqn:fyau}\\
&V^n(E^{\bul \leq N}_r\otimes_{{\cal W}_r({\cal O}_{Y_{\bul \leq N}})^{\star}}
({\cal W}_r{\Om}^i_{Y_{\bul \leq N}})^{\star})+
\nabla V^n(E^{\bul \leq N}_r\otimes_{{\cal W}_r({\cal O}_{Y_{\bul \leq N}})^{\star}}
({\cal W}_r{\Om}^{i-1}_{Y_{\bul \leq N}})^{\star}).
\end{align*} 
\end{prop} 

\begin{prop}
The following formula holds$:$ 
\begin{equation*} 
F^r({\rm Fil}^n(E^{\bul \leq N}_{n+r}
\otimes_{{\cal W}_{n+r}({\cal O}_{Y_{\bul \leq N}})^{\star}}
({\cal W}_{n+r}{\Om}^i_{Y_{\bul \leq N}})^{\star}))= 
F^r\nabla V^n(E^{\bul \leq N}_r
\otimes_{{\cal W}_r({\cal O}_{Y_{\bul \leq N}})^{\star}}
({\cal W}_r{\Om}^{i-1}_{Y_{\bul \leq N}})^{\star}). 
\tag{2.1.26.1}\label{eqn:frfyu}
\end{equation*} 
Consequently the morphism 
$$F^r \col E^{\bul \leq N}_{n+r}
\otimes_{{\cal W}_{n+r}({\cal O}_{Y_{\bul \leq N}})^{\star}}
({\cal W}_{n+r}{\Om}^i_{Y_{\bul \leq N}})^{\star}\lo 
F^r_*\{F^rE^{\bul \leq N}_n
\otimes_{{\cal W}_n({\cal O}_{Y_{\bul \leq N}})^{\star}}
({\cal W}_n{\Om}^i_{Y_{\bul \leq N}})^{\star}\}$$  
of ${\cal W}_n({\cal O}_{Y_{\bul \leq N}})^{\star}$-modules 
induces the following morphism of 
${\cal W}_n({\cal O}_{Y_{\bul \leq N}})^{\star}$-modules$:$  
\begin{align*} 
\check{F}{}^r \col & 
E^{\bul \leq N}_n
\otimes_{{\cal W}_n({\cal O}_{Y_{\bul \leq N}})^{\star}}
({\cal W}_n{\Om}^i_{Y_{\bul \leq N}})^{\star}\lo 
\tag{2.1.26.2}\label{eqn:fsyu}\\
& F^r_*\{F^r(E^{\bul \leq N}_{n+r}
\otimes_{{\cal W}_{n+r}({\cal O}_{Y_{\bul \leq N}})^{\star}}
({\cal W}_{n+r}{\Om}^i_{Y_{\bul \leq N}})^{\star})
/F^r\nabla V^n(E^{\bul \leq N}_r
\otimes_{{\cal W}_r({\cal O}_{Y_{\bul \leq N}})^{\star}}
({\cal W}_r{\Om}^{i-1}_{Y_{\bul \leq N}})^{\star})\}. 
\end{align*} 
For the case $r=n$, $\check{F}{}^n$ 
induces the following morphism
\begin{align*}
\check{F}{}^n \col 
E^{\bul \leq N}_n
\otimes_{{\cal W}_n({\cal O}_{Y_{\bul \leq N}})^{\star}}
({\cal W}_n{\Om}^i_{Y_{\bul \leq N}})^{\star} \lo 
F^n_*{\cal H}^i(E^{\bul \leq N}_n
\otimes_{{\cal W}_n({\cal O}_{Y_{\bul \leq N}})^{\star}}
({\cal W}_n{\Om}^{\bul}_{Y_{\bul \leq N}})^{\star}). 
\tag{2.1.26.3}\label{eqn:fyru}
\end{align*}
\end{prop}


\begin{theo}\label{theo:cin}  
The morphism $V^r\col F^r_*\{E^{\bul \leq N}_n
\otimes_{{\cal W}_n({\cal O}_{Y_{\bul \leq N}})^{\star}}
({\cal W}_n{\Om}^i_{Y_{\bul \leq N}})^{\star}\} \lo 
E^{\bul \leq N}_{n+r}
\otimes_{{\cal W}_{n+r}({\cal O}_{Y_{\bul \leq N}})^{\star}}
({\cal W}_{n+r}{\Om}^i_{Y_{\bul \leq N}})^{\star}$ 
induces the following morphism of 
${\cal W}_n({\cal O}_{Y_{\bul \leq N}})^{\star}$-modules$:$  
\begin{align*} 
\check{V}{}^r \col & 
F^r_*\{E^{\bul \leq N}_n
\otimes_{{\cal W}_n({\cal O}_{Y_{\bul \leq N}})^{\star}}
({\cal W}_n{\Om}^i_{Y_{\bul \leq N}})^{\star}/
F^r\nabla V^n(E^{\bul \leq N}_r\otimes_{{\cal W}_r({\cal O}_{Y_{\bul \leq N}})^{\star}}
({\cal W}_r{\Om}^{i-1}_{Y_{\bul \leq N}})^{\star})\} 
\tag{2.1.27.1}\label{eqn:fyu}\\
& \lo 
E^{\bul \leq N}_{n+r}
\otimes_{{\cal W}_{n+r}({\cal O}_{Y_{\bul \leq N}})^{\star}}
({\cal W}_{n+r}{\Om}^i_{Y_{\bul \leq N}})^{\star}. 
\end{align*}
\parno 
There exists a generalized Cartier isomorphism 
\begin{align*} 
\check{C}{}^r\col & 
F^r_*\{F^r(E^{\bul \leq N}_{n+r}
\otimes_{{\cal W}_{n+r}({\cal O}_{Y_{\bul \leq N}})^{\star}}
({\cal W}_{n+r}{\Om}^i_{Y_{\bul \leq N}})^{\star})/
F^r\nabla V^n(E^{\bul \leq N}_r
\otimes_{{\cal W}_r({\cal O}_{Y_{\bul \leq N}})^{\star}}
({\cal W}_r{\Om}^{i-1}_{Y_{\bul \leq N}})^{\star})\} \\
& \os{\sim}{\lo}  
E^{\bul \leq N}_n
\otimes_{{\cal W}_n({\cal O}_{Y_{\bul \leq N}})^{\star}}
({\cal W}_n{\Om}^i_{Y_{\bul \leq N}})^{\star} 
\tag{2.1.27.2}\label{eqn:wlfoyu}
\end{align*} 
of ${\cal W}_n({\cal O}_{Y_{\bul \leq N}})^{\star}$-modules,    
which is the inverse of $\check{F}{}^r$. 
The morphism  $\check{C}{}^r$ satisfies a relation 
${\bf p}^r\circ \check{C}{}^r=\check{V}{}^r$. 
In particular, there exist the following isomorphisms 
\begin{align*} 
\check{F}{}^n \col 
E^{\bul \leq N}_n
\otimes_{{\cal W}_n({\cal O}_{Y_{\bul \leq N}})^{\star}}
({\cal W}_n{\Om}^i_{Y_{\bul \leq N}})^{\star}  
\os{\sim}{\lo} 
F^n_*({\cal H}^i(E^{\bul \leq N}_n
\otimes_{{\cal W}_n({\cal O}_{Y_{\bul \leq N}})^{\star}}
({\cal W}_n{\Om}^{\bul}_{Y_{\bul \leq N}})^{\star})) 
\tag{2.1.27.3}\label{eqn:wlfyu}
\end{align*}  
and 
\begin{align*} 
\check{C}{}^n \col 
F^n_*({\cal H}^i(E^{\bul \leq N}_n
\otimes_{{\cal W}_n({\cal O}_{Y_{\bul \leq N}})^{\star}}
({\cal W}_n{\Om}^{\bul}_{Y_{\bul \leq N}})^{\star})) 
\os{\sim}{\lo} 
E^{\bul \leq N}_n
\otimes_{{\cal W}_n({\cal O}_{Y_{\bul \leq N}})^{\star}}
({\cal W}_n{\Om}^i_{Y_{\bul \leq N}})^{\star},   
\tag{2.1.27.4}\label{eqn:wlcyu}
\end{align*}  
which are the inverse of another. 
\end{theo}

\begin{theo}[{\bf Slope spectral sequence}]\label{theo:e1deg}
Set $E^{\bul \leq N}_{\infty}:=\vpl_nE^{\bul \leq N}_n$. 
Then there exists the following spectral sequence 
\begin{equation*} 
E_1^{i,q-i}=H^{q-i}(Y_{\bul \leq N},E^{\bul \leq N}_{\infty}
\otimes_{{\cal W}({\cal O}_{Y_{\bul \leq N}})^{\star}}
({\cal W}{\Om}^i_{Y_{\bul \leq N}})^{\star})
\Lo H^q((Y_{\bul \leq N}/{\cal W}(s))_{\rm crys},E^{\bul \leq N}). 
\tag{2.1.28.1}\label{eqn:wlce1yu}
\end{equation*}
This degenerates at $E_1$ modulo torsion and there exists 
the following canonical isomorphism 
\begin{align*} 
&(H^q((Y_{\bul \leq N}/{\cal W}(s))_{\rm crys},E^{\bul \leq N})
\otimes_{\cal W}K_0)_{[i,i+1)}\tag{2.1.28.2}\label{eqn:wlsd1yu}\\
&=H^{q-i}(Y_{\bul \leq N},E^{\bul \leq N}_{\infty}
\otimes_{{\cal W}({\cal O}_{Y_{\bul \leq N}})^{\star}}
({\cal W}{\Om}^i_{Y_{\bul \leq N}})^{\star})\otimes_{\cal W}K_0. 
\end{align*}
\end{theo}

\par 
The log crystal $E^{\bul \leq N}$ gives a log crystal 
$E^{\bul \leq N}_{Y_{\bul \leq N}/{\cal W}_n(s)}$ 
of ${\cal O}_{Y_{\bul \leq N}/{\cal W}_n(s)}$-modules. 
Set $E^{\bul \leq N}_n:=E^{\bul \leq N}_{{\cal W}_n(Y_{\bul \leq N})}$. 
Assume that $E^{\bul \leq N}$ is a flat coherent unit root log $F$-crystal.  
Then, by (\ref{eqn:wlfyu}) and (\ref{theo:ccrw}), 
\begin{equation*} 
E^{\bul \leq N}_n
\otimes_{{\cal W}_n({\cal O}_{Y_{\bul \leq N}})}{\cal W}_n{\Om}^i_{Y_{\bul \leq N}}=
F^n_*R^iu_{Y_{\bul \leq N}/{\cal W}_n(s)*}(E^{\bul \leq N}_{Y_{\bul \leq N}/{\cal W}_n(s)}).   
\tag{2.1.28.3}\label{eqn:wlyu}
\end{equation*} 
In a standard way,  we can define a boundary morphism 
``$p^{-n}\nabla $'': 
\begin{align*} 
\nabla \col 
R^iu_{Y_{\bul \leq N}/{\cal W}_n(s)*}(E^{\bul \leq N}_{Y_{\bul \leq N}/{\cal W}_n(s)})\lo 
R^{i+1}u_{Y_{\bul \leq N}/{\cal W}_n(s)*}(E^{\bul \leq N}_{Y_{\bul \leq N}/{\cal W}_n(s)}). 
\tag{2.1.28.4}\label{ali:ywns}
\end{align*} 
\par 
Assume that there exists an immersion 
$\iota \col Y_{\bul \leq N} \os{\sus}{\lo} {\cal Q}_{{\bul \leq N},n}$ 
into a log smooth scheme over ${\cal W}_n(s)$. 
Let 
$g_n \col {\cal Q}_{{\bul \leq N},n}\lo {\cal W}_n(s)$ 
be the structural morphism. 
Let $E^{\bul \leq N}$ be as in (\ref{theo:cin}). 
Let ${\mathfrak E}_{{\bul \leq N},n}$ 
be the log PD-envelope of the immersion 
$Y_{\bul \leq N} \os{\subset}{\lo} {\cal Q}_{{\bul \leq N},n}$ 
over $({\cal W}_n(s),p{\cal W}_n,[~])$. 
Let $({\cal E}^{\bul \leq N}_n,\nabla)$ be the corresponding 
${\cal O}_{{\mathfrak E}_{{\bul \leq N},n}}$-module 
with integrable connection to $E^{\bul \leq N}_{Y_{\bul \leq N}/{\cal W}_n(s)}$. 
Then, by (\ref{eqn:wlyu}) and \cite[(6.4)]{klog1} 
(or the log Poincar\'{e} lemma (\cite[(2.2.7)]{nh2})), 
we have the following equality 
\begin{equation*} 
E^{\bul \leq N}_n\otimes_{{\cal W}_n({\cal O}_{Y_{\bul \leq N}})}
{\cal W}_n{\Om}^i_{Y_{\bul \leq N}}=
F^n_*{\cal H}^i({\cal E}^{\bul \leq N}_n{\otimes}_{{\cal O}_{{\cal Q}_{{\bul \leq N},n}}}
{\Om}_{{\cal Q}_{{\bul \leq N},n}/{\cal W}_n(s)}^{\bul}) \quad (i\in {\mab N}). 
\tag{2.1.28.5}\label{eqn:wlh}
\end{equation*}  

\begin{theo}\label{theo:edcisw}
Let $E^{\bul \leq N}$ be a flat coherent unit root log $F$-crystal in 
$(Y_{\bul \leq N}/{\cal W}(s))_{\rm crys}$.   
Consider $E^{\bul \leq N}_{Y_{\bul \leq N}/{\cal W}_n(s)}$ as $E$ in 
{\rm (\ref{theo:ccrw})}. 
Denote ${\cal E}^{\bul \leq N,\bul}$ and ${\cal Q}_{\bul \leq N,\bul}$ 
in the proof of {\rm (\ref{theo:ccrw})} by 
${\cal E}^{\bul \leq N,\bul}_n$ and ${\cal Q}_{\bul \leq N,\bul,n}$, respectively. 
Then the isomorphism in {\rm (\ref{eqn:ywnnny})}
for the case $\star$=nothing is equal to 
the following isomorphism 
\begin{align*} 
R\pi_{{\rm zar}*}({\cal E}^{\bul \leq N,\bul}_n
{\otimes}_{{\cal O}_{{\cal Q}_{\bul \leq N,\bul,n}}}
\Om^{\bul}_{{\cal Q}_{\bul \leq N,\bul,n}/{\cal W}_n(s)})  
\os{\sim}{\lo} 
R\pi_{{\rm zar}*}(F^n_*{\cal H}^{\bul}({\cal E}^{\bul \leq N,\bul}_n
{\otimes}_{{\cal O}_{{\cal Q}_{\bul \leq N,\bul,n}}}
\Om^{*}_{{\cal Q}_{\bul \leq N,\bul,n}/{\cal W}_n(s)})).  
\tag{2.1.29.1}\label{ali:enroq}
\end{align*} 
\end{theo} 
\begin{proof} 
This follows from (\ref{ali:enoq}) and (\ref{eqn:wlh}). 
\end{proof}

\begin{rema}\label{rema:wqn}
Let the notations and the assumptions be in {\rm (\ref{theo:edcisw})}. 
Then, as in \cite{ir}, 
we can consider 
$R^iu_{Y_{\bul \leq N}/{\cal W}_n(s)*}(E^{\bul \leq N}_{Y_{\bul \leq N}/{\cal W}_n(s)})$ 
as the definition of the sheaf 
``$E^{\bul \leq N}_n\otimes_{{\cal W}_n({\cal O}_{Y_{\bul \leq N}})}
{\cal W}_n{\Om}^i_{Y_{\bul \leq N}}$'' 
and $\{R^iu_{Y_{\bul \leq N}/{\cal W}_n(s)*}
(E^{\bul \leq N}_{Y_{\bul \leq N}/{\cal W}_n(s)})\}_{i\in {\mab N}}$ 
becomes a  complex 
$\{R^{\bul}u_{Y_{\bul \leq N}/{\cal W}_n(s)*}(E^{\bul \leq N}_{Y_{\bul \leq N}/{\cal W}_n(s)})\}$. 
Moreover, one can prove that 
$$\{R^{\bul}u_{Y_{\bul \leq N}/{\cal W}_n(s)*}
(E^{\bul \leq N}_{Y_{\bul \leq N}/{\cal W}_n(s)})\}_{n\in {\mab N}}$$ 
has operators $F$, $V$, ${\bf p}$ and $R$ in a standard way. 
(\ref{theo:ccrw}) and (\ref{theo:edcisw}) imply 
that there exists a canonical isomorphism 
\begin{align*} 
E^{\bul \leq N}_n
\otimes_{{\cal W}_n({\cal O}_{Y_{\bul \leq N}})^{\star}}
({\cal W}_n{\Om}^{\bul}_{Y_{\bul \leq N}})^{\star} 
\os{\sim}{\lo} 
R^{\bul}u_{Y_{\bul \leq N}/{\cal W}_n(s)*}(E^{\bul \leq N})
\tag{2.1.30.1}\label{ali:ynw}
\end{align*} 
in ${\rm D}^+(g^{-1}({\cal W}_n))$.
\end{rema}

\section{Log de Rham-Witt complexes of crystals II}\label{sec:ldrwii}
Let the notations be as in the previous section. 
In this section $s$ is assumed to be 
the log point of a perfect field $\kap$ of characteristic $p>0$, 
that is, $s=({\rm Spec}(\kap),{\mab N}\oplus \kap^*\lo \kap)$ 
as in the Introduction. 
Set ${\cal W}_n:=\Gam({\cal W}_n(s),{\cal O}_{{\cal W}_n(s)})$. 
In this section we give a definition of the sheaf 
${\cal W}_n\wt{\Om}^i_Y$ $(i\in {\mab N})$ 
of ${\cal W}_n({\cal O}_Y)$-modules from our point of view. 
The sheaf ${\cal W}_n\wt{\Om}^i_Y$ has been denoted 
by $W_n\wt{\om}^i_Y$ in \cite{hdw} and \cite{msemi}. 
For a quasi-coherent log crystal $\ol{E}$ of 
${\cal O}_{Y/{\cal W}(\os{\circ}{s})}$-modules, 
we consider a complex  $\ol{E}_{n}\otimes_{{\cal W}_{n}({\cal O}_Y)}
{\cal W}_n\wt{{\Om}}_Y^\bul$ $(\ol{E}_n:=\ol{E}_{{\cal W}_n(Y)})$ 
and we give the ``crystalline'' interpretation of this complex. 
In the next section we construct an 
analogous filtered complex to 
$(A_{\rm zar}(X_{\bul \leq N,\os{\circ}{T}_0}/S(T)^{\nat},
E^{\bul \leq N}),P)$ in \S\ref{sec:psc}
by using results in this section in the case 
where $\os{\circ}{T}$ is the spectrum of a perfect field of characteristic $p>0$. 
This analogous filtered complex is a generalization of 
Hyodo-Mokrane-Steenbrink complex in \cite{msemi} and \cite{ndw}. 
\par 
Let $t\lo s$ be a morphism from a fine log scheme
whose underlying scheme is the spectrum of 
a perfect field of characteristic $p>0$. 
Set $s_{\os{\circ}{t}}:=s\times_{\os{\circ}{s}}\os{\circ}{t}$. 
Let $Y$ be a log smooth scheme of Cartier type over $s$. 
Set $Y_{\os{\circ}{t}}:=Y\times_{\os{\circ}{s}}\os{\circ}{t}$ 
and $\kap_t:=\Gam(t,{\cal O}_t)$. 
In this section we consider the following commutative diagram 
\begin{equation*} 
\begin{CD} 
t@>{\subset}>> {\cal W}_n(t)\\
@VVV @VVV \\
s@>{\subset}>>{\cal W}_n(s). 
\end{CD}
\end{equation*} 
We consider 
$S$, $T_0$ and $T$ in \S\ref{sec:ldc} 
as $s$, $t$, ${\cal W}_n(t)$, respectively in this section. 
Because 
$$S(T)=s({\cal W}_n(t))={\cal W}_n(s_{\os{\circ}{t}})$$ 
is hollow, 
${\cal W}_n(s_{\os{\circ}{t}})^{\nat}={\cal W}_n(s_{\os{\circ}{t}})$. 
\par 
Set 
\begin{align*} 
{\cal W}_n\wt{\Om}^i_{Y_{\os{\circ}{t}}}:=
\wt{R}^iu_{Y_{\os{\circ}{t}}/{\cal W}_n(\os{\circ}{t})*}
({\cal O}_{Y_{\os{\circ}{t}}/{\cal W}_n(\os{\circ}{t})}) 
\quad (i\in {\mab N}).
\tag{2.2.0.1}\label{eqn:lmtwy} 
\end{align*} 
(Recall the definition of $\wt{R}^iu_{Y_{\os{\circ}{t}}/{\cal W}_n(\os{\circ}{t})*}
({\cal O}_{Y_{\os{\circ}{t}}/{\cal W}_n(\os{\circ}{t})})$ in (\ref{defi:wrpt}).)

\begin{prop}\label{prop:fuo}  
Let 
\begin{equation*}
\begin{CD} 
t@>>> t'\\
@VVV @VVV \\
s@>>> s'  
\end{CD} 
\tag{2.2.1.1}\label{cd:tts}
\end{equation*} 
be a commutative diagram of fine log schemes whose underlying schemes 
are the spectrums of perfect fields of characteristic $p>0$. 
Assume that $s'$ and $s$ are log points.  
Then the sheaf ${\cal W}_n\wt{\Om}^i_{Y_{\os{\circ}{t}}}$ 
is contravariantly functorial with respect to the following commutative diagram
\begin{equation*} 
\begin{CD} 
Y_{\os{\circ}{t}}@>>> Y'_{\os{\circ}{t}{}'}\\
@VVV @VVV \\
s_{\os{\circ}{t}}@>>> s'_{\os{\circ}{t}{}'}. 
\end{CD} 
\tag{2.2.1.2}\label{cd:yyp}
\end{equation*} 
\end{prop}
\begin{proof} 
This immediately follows from (\ref{prop:cttu}) (1). 
\end{proof} 

\parno 
If there exists an immersion  $Y_{\os{\circ}{t}}\os{\sus}{\lo} {\cal Q}_n$ 
into a log smooth scheme over ${\cal W}_n(s_{\os{\circ}{t}})$,   then 
\begin{align*} 
{\cal W}_n\wt{\Om}^i_{Y_{\os{\circ}{t}}}={\cal H}^i({\cal O}_{{\mathfrak E}_{n}}
\otimes_{{\cal O}_{{\cal Q}_n}}
\Om^{\bul}_{{\cal Q}_n/{\cal W}_n(\os{\circ}{t})})
=
{\cal H}^i({\cal O}_{{\mathfrak E}_{n}}
\otimes_{{\cal O}_{{\cal Q}^{\rm ex}_n}}
\Om^{\bul}_{{\cal Q}^{\rm ex}_n/{\cal W}_n(\os{\circ}{t})}),  
\tag{2.2.1.3}\label{ali:woye}
\end{align*}
where ${\mathfrak E}_n$ is the log PD-envelope of the immersion 
$Y_{\os{\circ}{t}}\os{\sus}{\lo} {\cal Q}_n$ over 
$({\cal W}_n(s_{\os{\circ}{t}}),p{\cal W}_n,[~])$.  
This is equal to the definition 
in \cite[p.~311]{msemi} only in the case where ${\cal Q}_n$ is a lift of 
$Y_{\os{\circ}{t}}$ over ${\cal W}_n(s_{\os{\circ}{t}})$.
Let $Y_{\os{\circ}{t}}\os{\sus}{\lo} {\cal Q}_{2n}$ be an immersion 
into a log smooth scheme over ${\cal W}_{2n}(s)$. 
Set ${\cal Q}_{n}:={\cal Q}_{2n}\times_{{\cal W}_{2n}(s_{\os{\circ}{t}})}{\cal W}_n(s_{\os{\circ}{t}})$. 
Because the following sequence 
\begin{align*} 
0& \lo {\cal O}_{{\mathfrak E}_{n}}
\otimes_{{\cal O}_{{\cal Q}_n}}
\Om^{\bul}_{{\cal Q}_n/{\cal W}_n(\os{\circ}{t})}
\os{p^n}{\lo} 
{\cal O}_{{\mathfrak E}_{2n}}
\otimes_{{\cal O}_{{\cal Q}_{2n}}}
\Om^{\bul}_{{\cal Q}_{2n}/{\cal W}_{2n}(\os{\circ}{t})}
\tag{2.2.1.4}\label{ali:bpnt}\\
&\lo {\cal O}_{{\mathfrak E}_{n}}
\otimes_{{\cal O}_{{\cal Q}_n}}
\Om^{\bul}_{{\cal Q}_n/{\cal W}_n(\os{\circ}{t})}\lo 0
\end{align*} 
is exact ((\ref{coro:filt})), we have the following boundary morphism 
\begin{align*} 
d\col {\cal W}_n\wt{\Om}^i_{Y_{\os{\circ}{t}}}
=
{\cal H}^i({\cal O}_{{\mathfrak E}_{n}}
\otimes_{{\cal O}_{{\cal Q}_n}}
\Om^{\bul}_{{\cal Q}_n/{\cal W}_n(\os{\circ}{t})})
\lo {\cal H}^{i+1}({\cal O}_{{\mathfrak E}_{n}}
\otimes_{{\cal O}_{{\cal Q}_n}}
\Om^{\bul}_{{\cal Q}_n/{\cal W}_n(\os{\circ}{t})})
={\cal W}_n\wt{\Om}^{i+1}_{Y_{\os{\circ}{t}}}. \tag{2.2.1.5}\label{ali:bdft}\\
\end{align*} 
It is easy to prove that $d$ is independent of the choice of the immersion 
$Y_{\os{\circ}{t}}\os{\sus}{\lo} {\cal Q}_{2n}$ as in the proof of (\ref{prop:nop}) below. 
Let $g_n\col {\cal Q}_n\lo {\cal W}_n(s)$ be the structural morphism. 
\par
Let $\tau_n$ be a global section of the log 
structure $M_{{\cal W}_n(s)}$ of ${\cal W}_n(s)$ whose image in 
$\Gamma({\cal W}_n(s),M_{{\cal W}_n(s)}/{\cal O}_{{\cal W}_n(s)}^*)$ is a 
generator. Let $\tau_n$ be also the 
image of $\tau_n$ in $\Gamma({\cal Q}_n,M_{{\cal Q}_n})$. 
Let $\theta_n \in {\cal W}_n\wt{{\Om}}^1_{Y_{\os{\circ}{t}}}
={\cal H}^1({\cal O}_{{\mathfrak E}_{n}}
\otimes_{{\cal O}_{{\cal Q}_n}}
\Om^{\bul}_{{\cal Q}_n/{\cal W}_n(\os{\circ}{t})})$
be the cohomology class of  
$d\log \tau_n\in {\cal O}_{{\mathfrak E}_{n}}
\otimes_{{\cal O}_{{\cal Q}_n}}\Om^1_{{\cal Q}_n/{\cal W}_n(\os{\circ}{t})}$ 
(cf.~\cite[3.4]{msemi}). 

\begin{prop}[{\bf cf.~\cite[(1.4.3)]{hdw}}]\label{prop:plz} 
Set $Y_t:=Y\times_st$. 
Let ${\cal W}_n\Om^{\bul}_{Y_t}$ be the log de Rham-Witt complex of $Y_t/t$.  
Then the following sequence 
\begin{equation*}
\begin{CD}
0 @>>> {\cal W}_n\Om_{Y_t}^{\bul}[-1] 
@>{\theta_n \wedge}>>  {\cal W}_n\wt{\Om}_{Y_{\os{\circ}{t}}}^{\bul}  
@>>> 
{\cal W}_n\Om_{Y_t}^{\bul}@>>> 0 
\end{CD}
\tag{2.2.2.1}\label{cd:lijlmext}
\end{equation*}
is exact. 
\end{prop}
\begin{proof}
Because ${\cal W}_n\Om_{Y_t}^i={\cal W}_n\Om_{Y_{s_{\os{\circ}{t}}}}^i$ by (\ref{prop:bcdw}), 
we may assume that $t=s$.
Because the problem is local,  
we may also assume that there exists an immersion  
$Y\os{\sus}{\lo} {\cal Q}_{2n}$ 
into a log smooth scheme over ${\cal W}_{2n}(s)$. 
Set ${\cal Q}_n:={\cal Q}_{2n}\times_{{\cal W}_{2n}(s)}{\cal W}_n(s)$. 
Then we have the following exact sequence by (\ref{eqn:gsflxd}): 
\begin{align*} 
0&\lo {\cal O}_{{\mathfrak E}_n}\otimes_{{\cal O}_{{\cal Q}_n}}
\Om^{\bul}_{{\cal Q}_n/{\cal W}_n(s)}[-1]
\os{d\log \tau_n\wedge}{\lo} 
{\cal O}_{{\mathfrak E}_n}\otimes_{{\cal O}_{{\cal Q}_n}}
\Om^{\bul}_{{\cal Q}_n/{\cal W}_n(\os{\circ}{s})}\tag{2.2.2.2}\label{eqn:senlst}\\
&\lo 
{\cal O}_{{\mathfrak E}_n}\otimes_{{\cal O}_{{\cal Q}_n}}
\Om^{\bul}_{{\cal Q}_n/{\cal W}_n(s)} 
\lo 0.
\end{align*} 
Consider the following exact sequence 
\begin{align*} 
\cdots 
&\lo {\cal H}^{i-1}
({\cal O}_{{\mathfrak E}_{n}}\otimes_{{\cal O}_{{\cal Q}_n}}
\Om^{\bul}_{{\cal Q}_n/{\cal W}_n(s)})
\os{\theta_n\wedge}{\lo} 
{\cal H}^i({\cal O}_{{\mathfrak E}_{n}}\otimes_{{\cal O}_{{\cal Q}_n}}
\Om^{\bul}_{{\cal Q}_n/{\cal W}_n(\os{\circ}{s})}) \\
&\lo {\cal H}^i({\cal O}_{{\mathfrak E}_{n}}\otimes_{{\cal O}_{{\cal Q}_n}}
\Om^{\bul}_{{\cal Q}_n/{\cal W}_n(s)}) 
\lo \cdots.
\end{align*} 
We have to prove  
the injectivity of the following morphism 
\begin{align*} 
\theta_n \wedge \col {\cal H}^{i-1}
({\cal O}_{{\mathfrak E}_n}\otimes_{{\cal O}_{{\cal Q}_{n}}}
\Om^{\bul}_{{\cal Q}_n/{\cal W}_n(s)})
\lo 
{\cal H}^i({\cal O}_{{\mathfrak E}_n}\otimes_{{\cal O}_{{\cal Q}_n}}
\Om^{\bul}_{{\cal Q}_n/{\cal W}_n(\os{\circ}{s})})
\tag{2.2.2.3}\label{al:wny}
\end{align*}
to prove the injectivity of the morphism 
$\theta_n \wedge \col 
{\cal W}_n\Om_{Y}^{\bul}[-1] \lo  {\cal W}_n\wt \Om_{Y}^{\bul}$. 
Consider the following diagram 
\begin{equation*} 
\begin{CD} 
@. @. {\cal O}_{{\mathfrak E}_{2n}}
\otimes_{{\cal O}_{{\cal Q}_{2n}}}
\Om^{i-1}_{{\cal Q}_{2n}/{\cal W}_{2n}(\os{\circ}{s})}\\ 
@. @. @VV{d}V \\ 
0@>>>{\cal O}_{{\mathfrak E}_{n}}\otimes_{{\cal O}_{{\cal Q}_{n}}}
\Om^{i}_{{\cal Q}_{n}/{\cal W}_n(\os{\circ}{s})} 
@>{p^n}>> 
{\cal O}_{{\mathfrak E}_{2n}}\otimes_{{\cal O}_{{\cal Q}_{2n}}}
\Om^{i}_{{\cal Q}_{2n}/{\cal W}_{2n}(\os{\circ}{s})}. \\ 
\end{CD}
\tag{2.2.2.4}\label{cd:xddefss}
\end{equation*} 
To prove the injectivity, 
we have only to prove the following inclusion holds
\begin{align*} 
d\log \tau_n \wedge (p^n)^{-1}d
({\cal O}_{{\mathfrak E}_{2n}}\otimes_{{\cal O}_{{\cal Q}_{2n}}}
{\Om}^{i-2}_{{\cal Q}_{2n}/{\cal W}_{2n}(\os{\circ}{s})})
 \supset & 
(p^n)^{-1}
d
({\cal O}_{{\mathfrak E}_{2n}}\otimes_{{\cal O}_{{\cal Q}_{2n}}}
{\Om}^{i-1}_{{\cal Q}_{2n}/{\cal W}_{2n}(\os{\circ}{s})}) 
\tag{2.2.2.5}\label{al:hoyw}\\
& \cap 
(d\log \tau_n \wedge
({\cal O}_{{\mathfrak E}_{n}}\otimes_{{\cal O}_{{\cal Q}_{n}}}
{\Om}^{i-1}_{{\cal Q}_{n}/{\cal W}_n(\os{\circ}{s})})).   
\end{align*}
For a local section 
$\om \in {\cal O}_{{\mathfrak E}_{n}}\otimes_{{\cal O}_{{\cal Q}_n}}
{\Om}^{i-1}_{{\cal Q}_n/{\cal W}_n(\os{\circ}{s})}$, 
assume that 
$$d\log \tau_n \wedge \om \in (p^n)^{-1}
d({\cal O}_{{\mathfrak E}_{2n}}\otimes_{{\cal O}_{{\cal Q}_{2n}}}
\Om^{i-1}_{{\cal Q}_{2n}/{\cal W}_{2n}(\os{\circ}{s})}).$$ 
Then, by using a basis of 
${\cal O}_{{\mathfrak E}_{2n}}
\otimes_{{\cal O}_{{\cal Q}_{2n}}}\Om^1_{{\cal Q}/{\cal W}_{2n}(\os{\circ}{s})}$ 
taking $d\log \tau_{2n}$ as a member and 
noting that $d\log \tau_{2n}$ is closed and that $d\tau_{2n}=\tau_{2n}d\log \tau_{2n}=0$, 
it is clear that 
$\om \in (p^n)^{-1}d({\cal O}_{{\mathfrak E}_{2n}}\otimes_{{\cal O}_{{\cal Q}_{2n}}}
\wt{\Om}^{i-2}_{{\cal Q}_{2n}/{\cal W}_{2n}(\os{\circ}{s})})$. 
Hence we obtain the exactness of (\ref{cd:lijlmext}). 
\end{proof}

\begin{rema}\label{rema:qok}
(1) Let the notations be as in \cite{hdw}. 
By expressing a local section of $\wt{\om}^1_{{\cal Y}_n}$ 
by a linear combination of a basis taking $dT/T$ as a member, 
the inclusion  
\begin{align*} 
(d\wt{\om}^{i-1}_{{\cal Y}_n} \cap 
(\wt{\om}^{i-1}_{{\cal Y}_n}
\wedge (dT/T))) \subset d\wt{\om}^{i-2}_{{\cal Y}_n} \wedge dT/T
\tag{2.2.3.1}\label{ali:dyn}
\end{align*} 
in \cite[p.~246]{hdw}
is an obvious formula about which 
the editorial comment \cite[(6), (11A)]{hdw} is unnecessary. 
Because the boundary morphism 
${\cal H}^i(\Om_{{\cal Y}_n/{\cal W}_n(\os{\circ}{s})}) 
\lo 
{\cal H}^{i+1}(\Om_{{\cal Y}_n/{\cal W}_n(\os{\circ}{s})})$ 
is ``$p^{-n}d$'', 
the above formula is not a precisely necessary ingredient 
for the injectivity needed in [loc.~cit.].  
The boundary morphism of ${\cal H}^i(\Om_{{\cal Y}_n/{\cal W}_n(\os{\circ}{s})}) 
\lo 
{\cal H}^{i+1}(\Om_{{\cal Y}_n/{\cal W}_n(\os{\circ}{s})})$ has been misunderstood in 
[loc.~cit.] and the editorial comment for (\ref{ali:dyn}). 
\par 
\par 
(2) In \cite[p.~311]{msemi} it has been claimed that 
the sequence (\ref{cd:lijlmext}) is exact. However the proof of 
this claim (the injectivity of the morphism 
$\theta_n \wedge \col W_n\om^{q-1}_Y\lo W_n\wt{\om}^q_Y$) 
has not been given in [loc.~cit.]. 
\par 
(3) I am not sure that 
the reduction to the local calculations of 
the proof of the well-definedness of 
$W_n\wt{\om}^q_Y$ $(q\in {\mab N})$ 
(and $P_kW_n\wt{\om}^q_Y$ $(k\in {\mab Z})$)
in the end of \cite[\S1]{hdw} (and in \cite[Lemme 3.4]{msemi}) 
are correct since the differentials of 
(log) de Rham complexes are {\it not} linear for the structure sheaves. 
See the proof of (\ref{coro:lemos}) (and \cite[(2.12.2)]{nh2}) 
which can remove this incorrectness. 
\end{rema}
It is known that the projection
$R\col {\cal W}_n\wt \Om_{Y}^{\bul}\lo {\cal W}_{n-1}\wt \Om_{Y}^{\bul}$
is well-defined (\cite[(1.3.2)]{hdw}, cf.~\cite[(4.2)]{hk}, \cite[(6.27)]{ndw}).
However the definition of $\pi(=R)$ in \cite{hdw} is not perfect 
as pointed out in (\ref{rema:csh}) (and (\ref{rema:sot})). 
Moreover, we would like to avoid the calculation of
subquotients of  ${\cal W}_n\wt \Om_{Y}^i$ $(i\in {\mab N})$ for the
calculations of the kernel and the image of
${\bf p}\col {\cal W}_n\wt \Om_{Y}^i\lo
{\cal W}_{n+1}\wt \Om_{Y}^i$. 
In the following, we calculate ${\rm Ker}({\bf p})$ and ${\rm Im}({\bf p})$ 
and construct $R$ in our formalism 
based on Hyodo's idea in \cite[\S1]{hdw} 
as in the previous section;  
we do not use the admissible lift in \cite[(1.1)]{hdw} and \cite[(2.4.3)]{msemi}.

\par
\par
Let $n$ and $i$ be nonnegative integers.
To define $R\col {\cal W}_n\wt \Om_{Y_{\os{\circ}{t}}}^i\lo {\cal W}_{n-1}\wt \Om_{Y_{\os{\circ}{t}}}^i$,
we first assume that there exists an immersion
$Y_{\os{\circ}{t}}\os{\sus}{\lo} {\cal Q}_{n+i}$ into 
a formally log smooth integral scheme over
${\cal W}_{n+i}(s_{\os{\circ}{t}})$
such that ${\cal Q}_{n+i}$ has a lift 
$\varphi_{n+i} \col {\cal Q}_{n+i}\lo {\cal Q}_{n+i}$
of the Frobenius endomorphism of 
${\cal Q}_{n+i}\times_{{\cal W}_{n+i}(s_{\os{\circ}{t}})}s_{\os{\circ}{t}}$.
Let ${\mathfrak E}_{n+i}$ be the log PD-envelope of the immersion
$Y_{\os{\circ}{t}}\os{\sus}{\lo} {\cal Q}_{n+i}$ over $({\cal W}_{n+i}(s_{\os{\circ}{t}}),p{\cal W}_{n+i},[~])$.
Set ${\cal Q}_j:={\cal Q}_{n+i}\times_{{\cal W}_{n+i}(s_{\os{\circ}{t}})}{\cal W}_j(s_{\os{\circ}{t}})$ and 
${\mathfrak E}_j:={\mathfrak E}_{n+i}\times_{{\cal W}_{n+i}(s_{\os{\circ}{t}})}{\cal W}_j(s_{\os{\circ}{t}})$
$(0\leq j\leq n+i)$.
Consider the following morphism
\begin{align*}
\varphi^*_{n+i}
\col & {\cal O}_{{\mathfrak E}_{n+i}}\otimes_{{\cal O}_{{\cal Q}_{n+i}}}
\Om^i_{{\cal Q}_{n+i}/{\cal W}_{n+i}(\os{\circ}{t})}\lo
{\cal O}_{{\mathfrak E}_{n+i}}\otimes_{{\cal O}_{{\cal Q}_{n+i}}}
\Om^i_{{\cal Q}_{n+i}/{\cal W}_{n+i}(\os{\circ}{t})} \quad (i\in {\mab N}).
\end{align*}
%
As in (\ref{ali:enqw}), 
we have the following well-defined morphism (cf.~\cite[Editorial comment (5)]{hdw})
\begin{align*}
p^{-(i-1)}\varphi^*_{n+i} &\col {\cal W}_n\wt{\Om}_{Y_{\os{\circ}{t}}}^i
={\cal H}^i({\cal O}_{{\mathfrak E}_n}\otimes_{{\cal O}_{{\cal Q}_n}}
\Om^{\bul}_{{\cal Q}_n/{\cal W}_n(\os{\circ}{t})})\lo\tag{2.2.3.2}\label{ali:enw}\\
& 
{\cal H}^i({\cal O}_{{\mathfrak E}_{n+1}}\otimes_{{\cal O}_{{\cal Q}_{n+1}}}
\Om^{\bul}_{{\cal Q}_{n+1}/{\cal W}_{n+1}(\os{\circ}{t})})
={\cal W}_{n+1}\wt{\Om}_{Y_{\os{\circ}{t}}}^i \quad (i\in {\mab N}).
\end{align*}

\begin{rema}\label{rema:infap}
Assume that there exists an immersion 
$Y_{\os{\circ}{t}}\os{\sus}{\lo} {\cal Q}$ 
into a formally log smooth integral scheme over ${\cal W}(s_{\os{\circ}{t}})$ 
such that ${\cal Q}$ has a lift $\varphi \col {\cal Q}\lo {\cal Q}$ 
of the Frobenius endomorphism of 
${\cal Q}\times_{{\cal W}(s_{\os{\circ}{t}})}s_{\os{\circ}{t}}$. 
Let ${\mathfrak E}$ be the log PD-envelope of the immersion 
$Y_{\os{\circ}{t}}\os{\sus}{\lo} {\cal Q}$ over 
$({\cal W}(s_{\os{\circ}{t}}),p{\cal W},[~])$.  
Set ${\cal Q}_n:={\cal Q}\times_{{\cal W}(s_{\os{\circ}{t}})}{\cal W}_n(s_{\os{\circ}{t}})$ and   
${\mathfrak E}_n:={\mathfrak E}\times_{{\cal W}(s_{\os{\circ}{t}})}{\cal W}_n(s_{\os{\circ}{t}})$. 
Consider the following morphism 
\begin{align*} 
\varphi^* \col {\cal O}_{\mathfrak E}
\otimes_{{\cal O}_{\cal Q}}
\Om^i_{{\cal Q}/{\cal W}(\os{\circ}{t})}
\lo 
{\cal O}_{\mathfrak E}\otimes_{{\cal O}_{\cal Q}}
\Om^i_{{\cal Q}/{\cal W}(\os{\circ}{t})} \quad (i\in {\mab N}). 
\end{align*} 
As in (\ref{ali:qnw}), 
we have the following well-defined morphism 
\begin{align*} 
p^{-(i-1)} \varphi^* \col 
& {\cal W}_n\wt{\Om}_{Y_{\os{\circ}{t}}}^i
={\cal H}^i({\cal O}_{{\mathfrak E}_{n}}\otimes_{{\cal O}_{{\cal Q}_n}}
\Om^{\bul}_{{\cal Q}_n/{\cal W}_n(\os{\circ}{t})})
\tag{2.2.4.1}\label{ali:nw}\\
& \lo 
{\cal H}^i({\cal O}_{{\mathfrak E}_{n+1}}\otimes_{{\cal O}_{{\cal Q}_{n+1}}}
\Om^{\bul}_{{\cal Q}_{n+1}/{\cal W}_{n+1}(\os{\circ}{t})})
={\cal W}_{n+1}\wt{\Om}_{Y_{\os{\circ}{t}}}^i \quad (i\in {\mab N}). 
\end{align*} 
This morphism is equal to (\ref{ali:enw}). 
\end{rema}

\begin{prop}\label{prop:nop} 
$(1)$ The morphism $p^{-(i-1)}\varphi^*$ is independent of the choice of 
the immersion $Y_{\os{\circ}{t}}\os{\sus}{\lo}{\cal Q}$ and 
the lift $\varphi \col {\cal Q}\lo {\cal Q}$ 
of the Frobenius endomorphism of 
${\cal Q}\times_{{\cal W}(s_{\os{\circ}{t}})}s_{\os{\circ}{t}}$.
We set ${\bf p}:=p^{-(i-1)}\varphi^*\col 
{\cal W}_n\wt{\Om}_{Y_{\os{\circ}{t}}}^i\lo 
{\cal W}_{n+1}\wt{\Om}_{Y_{\os{\circ}{t}}}^i$. 
\par 
$(2)$ ${\rm Ker}({\bf p}\col {\cal W}_n\wt{\Om}_{Y_{\os{\circ}{t}}}^i\lo 
{\cal W}_{n+1}\wt{\Om}_{Y_{\os{\circ}{t}}}^i)=0$.
\par 
$(3)$ ${\rm Im}({\bf p}\col {\cal W}_n\wt{\Om}_{Y_{\os{\circ}{t}}}^i\lo 
{\cal W}_{n+1}\wt{\Om}_{Y_{\os{\circ}{t}}}^i)=
{\rm Im}(p\col {\cal W}_{n+1}\wt{\Om}_{Y_{\os{\circ}{t}}}^i\lo 
{\cal W}_{n+1}\wt{\Om}_{Y_{\os{\circ}{t}}}^i)$. 
\end{prop}
\begin{proof}
The proof is the same as that of (\ref{prop:nqop}). 
\end{proof}

\begin{defi}[{\bf cf.~\cite[(1.3.2)]{hdw}}]\label{defi:hyd} 
The morphism $R\col {\cal W}_{n+1}\wt{\Om}^i_{Y_{\os{\circ}{t}}}\lo 
{\cal W}_n\wt{\Om}^i_{Y_{\os{\circ}{t}}}$ 
is, by definition, the unique morphism fitting into the following commutative diagram: 
\begin{equation*} 
\begin{CD} 
{\cal W}_{n+1}\wt{\Om}^i_{Y_{\os{\circ}{t}}}
@>{R}>>{\cal W}_n\wt{\Om}^i_{Y_{\os{\circ}{t}}}\\
@V{p}VV @VV{\bf p}V\\ 
{\cal W}_{n+1}\wt{\Om}^i_{Y_{\os{\circ}{t}}}
@={\cal W}_{n+1}\wt{\Om}^i_{Y_{\os{\circ}{t}}}. 
\end{CD}
\end{equation*} 
\end{defi}


\parno 
If there exists a formal lift ${\cal Y}/{\cal W}(s_{\os{\circ}{t}})$ of 
$Y_{\os{\circ}{t}}/s_{\os{\circ}{t}}$, 
set $(\Om^{\bul},\phi):=(\Om^{\bul}_{{\cal Y}/{\cal W}(\os{\circ}{t})},\varphi^*)$ 
and 
$$Z^i_n:=\{\om \in \Om^i \vert~ 
d\om \in p^n\Om^{i+1}\}, \quad B^i_n:=p^n\Om^i+ d\Om^{i-1}.$$ 
By the proof of (\ref{prop:nop}) and \cite[(6.27)]{ndw},  
the morphism $R$ is equal to the following composite morphism 
in the case where $Y_{\os{\circ}{t}}/s_{\os{\circ}{t}}$ 
has a local lift ${\cal Y}/{\cal W}(s_{\os{\circ}{t}})$
(cf.~\cite[(4.2)]{hk}):
\begin{align*}
{\cal W}_{n+1}\wt{\Om}^i_{Y_{\os{\circ}{t}}}= 
Z^i_{n+1}/B^i_{n+1} 
\os{\us{\sim}{p^i}}{\lo} 
p^iZ^i_{n+1}/p^iB^i_{n+1}
& \os{{\rm proj}.}{\lo} 
p^iZ^i_{n+1}/(p^{i+n}Z^i_1+p^{i-1}dZ^{i-1}_1) 
\tag{2.2.6.1}\label{ali:aplocpj}\\
{} & \os{\us{\sim}{\phi}}{\longleftarrow}
Z^i_n/B^i_n={\cal W}_n\wt{\Om}^i_{Y_{\os{\circ}{t}}}.
\end{align*}
By \cite[(6.5)]{ndw}, 
the morphism $R$ is equal to the 
following composite morphism$:$
\begin{align*}
&{\cal W}_{n+1}\wt{\Om}^i_{Y_{\os{\circ}{t}}}  = 
Z^i_{n+1}/B^i_{n+1}
\os{{\rm proj}.}{\lo} 
Z^i_{n+1}/(p^nZ^i_1+d\Om^{i-1}) 
\tag{2.2.6.2}\label{ali:alocpj}\\
& 
\os{(p^{-i}\phi)^{-1}}{\os{\sim}{\lo}}
Z^i_n/(p^n\Om^i+pd\Om^{i-1})
\os{{\rm proj}.}{\lo}
Z^i_n/B^i_n={\cal W}_n\wt{\Om}^i_{Y_{\os{\circ}{t}}}.
\end{align*}
In particular, $R$ is surjective. 

\par 
The following is a generalization of \cite[(8.1.2)]{ndw}:

\par
As is well-known,  
$\{{\cal W}_n\wt{{\Om}}^{\bul}_{Y_{\os{\circ}{t}}}\}_{n\in {\mab N}}$ is 
a pro-dga over the Raynaud ring ${\mab R}$ of $\kap_t$, 
that is, there exist operators $F$, $V$, $d$, ${\bf p}$ 
and $R$ on 
$\{{\cal W}_n\wt{{\Om}}^{\bul}_{Y_{\os{\circ}{t}}}\}_{n\in {\mab N}}$
satisfying standard relations.

\begin{prop}\label{prop:t} 
Let $R\col {\cal W}_{n+1}\wt{\Om}^1_{Y}\lo 
{\cal W}_n\wt{\Om}^{\bul}_{Y}$ be the projection. 
Then $R(\theta_{n+1})=\theta_n$. 
\end{prop}
\begin{proof}  
Let the notations be before (\ref{prop:plz}).  
Since $\varphi(\tau_{n+1})=u\tau_{n+1}^p$  
for some $u\in {\cal O}_{{\cal Q}_{n+1}}^*$, 
$\varphi(\theta_{n+1}) \col =p\theta_{n+1}$. 
Hence ${\bf p}(\theta_n)=p\theta_{n+1}$. 
Since ${\bf p}R=p$ and ${\bf p}$ is injective,  
$R(\theta_{n+1})=\theta_n$. 
\end{proof} 
Set $\theta:=\vpl_n\theta_n \in {\cal W}\wt{{\Om}}^1_{Y_{\os{\circ}{t}}}$. 

\begin{coro}\label{coro:exwn}
The following sequence 
\begin{equation*}
\begin{CD}
0 @>>>{\cal W}\Om_{Y_t}^{\bul}[-1] 
@>{\theta \wedge}>>  {\cal W}\wt{\Om}_{Y_{\os{\circ}{t}}}^{\bul} @>>> 
{\cal W}\Om_{Y_t}^{\bul}@>>> 0
\end{CD}
\tag{2.2.8.1}\label{cd:lijplext}
\end{equation*}
is exact. 
\end{coro} 
\begin{proof} 
Because the projection 
$R\col {\cal W}_n\Om_{Y_t}^{\bul}\lo {\cal W}_{n-1}\Om_{Y_t}^{\bul}$ is surjective, 
we obtain the exactness of (\ref{cd:lijplext}) by that of (\ref{cd:lijlmext}).  
\end{proof} 

\par 
As in \cite{nh2}, we can prove that 
the Frobenius endomorphism of $Y$ induces the following morphism 
$$\wt{\Phi}_{\star} \col 
{\cal W}_{\star}\wt{\Om}^i_{Y_{\os{\circ}{t}}} \lo 
{\cal W}_{\star}\wt{\Om}^i_{Y_{\os{\circ}{t}}}  \quad 
(\star=n \in {\mab Z}_{>0} \text{ or nothing}).$$
More generally, 
consider the commutative diagram (\ref{cd:tts}) and 
the following commutative diagram of log smooth schemes of 
Cartier type over $s$ and $s'$:  
\begin{equation*} 
\begin{CD} 
Y_{\os{\circ}{t}} @>{g}>> Y'_{\os{\circ}{t}{}'}\\ 
@VVV @VVV \\ 
s_{\os{\circ}{t}} @>>> s'_{\os{\circ}{t}{}'} \\ 
@V{\bigcap}VV @VV{\bigcap}V \\ 
{\cal W}_n(s_{\os{\circ}{t}}) @>{u_n}>> {\cal W}_n(s'_{\os{\circ}{t}{}'}) \\
@V{\bigcap}VV @VV{\bigcap}V \\ 
{\cal W}(s_{\os{\circ}{t}}) @>{u}>> {\cal W}(s'_{\os{\circ}{t}{}'}). 
\end{CD}
\tag{2.2.8.2}\label{eqn:xdxdrss}
\end{equation*} 
By (\ref{prop:fuo}) we have the following morphism 
\begin{align*} 
\wt{g}{}^*_{\star} \col 
g^*({\cal W}_{\star}\wt{\Om}^i_{Y'_{\os{\circ}{t}{}'}}) 
\lo {\cal W}_{\star}\wt{\Om}^i_{Y_{\os{\circ}{t}}} \quad 
(\star=n \in {\mab Z}_{>0} \text{ or nothing}). 
\tag{2.2.8.3}\label{eqn:xttrss}
\end{align*}
The sequences (\ref{cd:lijlmext}) 
and (\ref{cd:lijplext}) are equal to the following exact sequences 
\begin{align*}
0 \lo {\cal W}_n\Om^{\bul}_{Y_t}(-1,u_n)[-1] \os{\theta_n \wedge}{\lo} 
{\cal W}_n\wt{\Om}^{\bul}_{Y_{\os{\circ}{t}}}\lo {\cal W}_n\Om^{\bul}_{Y_t}\lo 0
\tag{2.2.8.4}\label{eqn:lijfmext}
\end{align*}
and 
\begin{align*}
0 \lo {\cal W}\Om^{\bul}_{Y_t}(-1,u)[-1] 
\os{\theta \wedge}{\lo} 
{\cal W}\wt{\Om}^{\bul}_{Y_{\os{\circ}{t}}}
\lo {\cal W}\Om^{\bul}_{Y_t}\lo 0,  
\tag{2.2.8.5}\label{eqn:lijpfext}
\end{align*}
respectively. 
\par
Let $\al \col M_{Y_{\os{\circ}{t}}} \lo {\cal O}_{Y_{\os{\circ}{t}}}$ be 
the structural morphism. 
Consider two abelian 
subsheaves $\wt{\cal F}_n$ 
and $\wt{\cal G}_n$ in
\begin{equation*}
({\cal W}_n({\cal O}_{Y_{\os{\circ}{t}}})'\otimes_{\mab Z}
\bigwedge^i(M^{\rm gp}_{Y_{\os{\circ}{t}}}/\al^{-1}(\kap^*)))\oplus 
({\cal W}_n({\cal O}_{Y_{\os{\circ}{t}}})'\otimes_{\mab Z}
\bigwedge^{i-1}(M^{\rm gp}_{Y_{\os{\circ}{t}}}/\al^{-1}(\kap^*))): 
\tag{2.2.8.6}\label{eqn:bgtlwlam}
\end{equation*}
the sheaf $\wt{\cal F}_n$ is, by definition, 
generated by the images of local sections of the type (\ref{eqn:frel}) 
and $\wt{\cal G}_n$ is, 
by definition, generated by the images of local sections of 
the type (\ref{eqn:nrelp}). 
Set 
\begin{equation*}
({\cal W}_n\wt{\Om}{}^i_{Y_{\os{\circ}{t}}})':=
\{({\cal W}_n({\cal O}_{Y_{\os{\circ}{t}}})'\otimes_{\mab Z}
\bigwedge^i(M^{\rm gp}_{Y_{\os{\circ}{t}}}/\al^{-1}(\kap^*)))
\oplus 
({\cal W}_n({\cal O}_{Y_{\os{\circ}{t}}})'\otimes_{\mab Z}
\bigwedge^{i-1}(M^{\rm gp}_{Y_{\os{\circ}{t}}}/\al^{-1}(\kap^*)))\}
/(\wt{\cal F}_n+\wt{\cal G}_n). 
\tag{2.2.8.7}\label{eqn:obtwlm}
\end{equation*}
Here we have considered  
$({\cal W}_n\wt{\Om}{}^i_{Y_{\os{\circ}{t}}})'$ only as 
an abelian sheaf on $Y_{\os{\circ}{t},\rm zar}$.
We define a boundary morphism 
$d\col ({\cal W}_n\wt{\Om}{}^i_{Y_{\os{\circ}{t}}})' \lo 
({\cal W}_n\wt{\Om}{}^{i+1}_{Y_{\os{\circ}{t}}})'$ by 
the following formula 
\begin{align*}
&d(b\otimes(a_1 \wedge \cdots \wedge a_i), 
c\otimes(a'_2 \wedge \cdots \wedge a'_i))=
(0, b\otimes(a_1 \wedge \cdots \wedge a_i)) 
\tag{2.2.8.8}\label{eqn:dfbac}\\
&(a_1, \ldots, a_i, a'_2, \cdots, a'_i \in M_{Y_{\os{\circ}{t}}},
~b,c\in {\cal W}_n({\cal O}_{Y_{\os{\circ}{t}}})). 
\end{align*}
It is easy to check that $d$ is well-defined and $d^2=0$. 
Thus we have a complex 
$({\cal W}_n\wt{{\Om}}^{\bul}_{Y_{\os{\circ}{t}}})'$ 
of abelian sheaves. Set 
$({\cal W}\wt{{\Om}}^{\bul}_{Y_{\os{\circ}{t}}})':=
\vpl_n({\cal W}_n\wt{{\Om}}^{\bul}_{Y_{\os{\circ}{t}}})'$.
The Frobenius endomorphism of $Y$ induces the following endomorphism
\begin{equation*} 
\wt{\Phi}{}^*_{\star} \col 
({\cal W}_{\star}\wt{{\Om}}^i_{Y_{\os{\circ}{t}}})' \lo 
({\cal W}_{\star}\wt{{\Om}}^i_{Y_{\os{\circ}{t}}})' \quad 
(\star=n \in {\mab Z}_{>0} 
\text{ or nothing}). 
\end{equation*} 
More generally, we have the morphism 
\begin{equation*} 
\wt{g}{}^*_{\star} \col 
g^*({\cal W}_{\star}\wt{{\Om}}^i_{Y'_{\os{\circ}{t}{}'}})' \lo 
({\cal W}_{\star}\wt{{\Om}}^i_{Y_{\os{\circ}{t}}})' \quad 
(\star=n \in {\mab Z}_{>0} 
\text{ or nothing}) 
\end{equation*} 
for the commutative diagrams (\ref{cd:tts}) and (\ref{eqn:xdxdrss}). 
Let $\theta'_{\star}:=(d\log \tau)_{\star}$ be 
the image of $1\otimes \tau \in 
{\cal W}_{\star}({\cal O}_{Y_{\os{\circ}{t}}})'\otimes M_{Y_{\os{\circ}{t}}}^{\rm gp}$ in 
$({\cal W}_{\star}\wt{{\Om}}^1_{Y_{\os{\circ}{t}}})'$.
Then $\theta'_{\star}$ is independent of the choice of $\tau$.

\par 
Recall the following theorem: 

\begin{theo}[{\bf \cite[(11.1)]{ndw}}]\label{theo:tiobre}
$(1)$ There exists a canonical 
isomorphism 
\begin{equation*} 
\wt{C}^{-n} \col 
({\cal W}_n\wt{{\Om}}^{\bul}_{Y_{\os{\circ}{t}}})' \os{\sim}{\lo} 
{\cal W}_n\wt{{\Om}}^{\bul}_{Y_{\os{\circ}{t}}} \quad (n\in {\mab Z}_{\geq 1})
\tag{2.2.9.1}\label{eqn:snisot}
\end{equation*}
of complexes of abelian sheaves on $Y$ which makes 
the following two diagrams $($of exact sequences$)$
commutative$:$
\begin{equation*}
\begin{CD}
0 @>>> ({\cal W}_n{\Om}_{Y_t}^{\bul})'(-1,u)[-1] 
@>{\theta'_n  \wedge}>>  
({\cal W}_n\wt{{\Om}}^{\bul}_{Y_{\os{\circ}{t}}})' @>>> 
({\cal W}_n{\Om}_{Y_{t}}^{\bul})'@>>> 0\\ 
@. @V{C^{-n},~\simeq}VV  
@V{\wt{C}{}^{-n},~\simeq}VV @V{C^{-n},~\simeq}VV \\
0 @>>> {\cal W}_n{\Om}_{Y_t}^{\bul}(-1,u)[-1] 
@>{\theta_n \wedge}>>  
{\cal W}_n\wt{{\Om}}^{\bul}_{Y_{\os{\circ}{t}}} @>>> 
{\cal W}_n{{\Om}}_{Y_t}^{\bul}@>>> 0,
\end{CD}
\tag{2.2.9.2;$n$}\label{cd:orexttlam}
\end{equation*}
\begin{equation*}
\begin{CD}
({\cal W}_{n+1}\wt{{\Om}}^{\bul}_{Y_{\os{\circ}{t}}})'  
@>{\os{\wt{C}{}^{-(n+1)}}{\sim}}>> 
{\cal W}_{n+1}\wt{{\Om}}^{\bul}_{Y_{\os{\circ}{t}}}\\ 
@V{{\rm proj}.}VV  @VV{R}V \\
({\cal W}_n\wt{{\Om}}^{\bul}_{Y_{\os{\circ}{t}}})'  
@>{\os{\wt{C}{}^{-n}}{\sim}}>> 
{\cal W}_n\wt{{\Om}}^{\bul}_{Y_{\os{\circ}{t}}}.
\end{CD}
\tag{2.2.9.3}\label{cd:pjpitlm}
\end{equation*}
\par
$(2)$ The two projections 
$R \col 
({\cal W}_{n+1}\wt{{\Om}}^{\bul}_{Y_{\os{\circ}{t}}})' \lo
({\cal W}_n\wt{{\Om}}^{\bul}_{Y_{\os{\circ}{t}}})'$, 
$R \col 
(W_{n+1}{{\Om}}_{Y_{\os{\circ}{t}}}^{\bul})' \lo 
({\cal W}_n{{\Om}}_{Y_{\os{\circ}{t}}}^{\bul})'$ 
and the other two projections
$R \col {\cal W}_{n+1}\wt{{\Om}}^{\bul}_{Y_{\os{\circ}{t}}} \lo
{\cal W}_n\wt{{\Om}}^{\bul}_{Y_{\os{\circ}{t}}}$, 
$R \col W_{n+1}{{\Om}}_{Y_{t}}^{\bul} \lo 
{\cal W}_n{{\Om}}_{Y_{t}}^{\bul}$ induce 
a morphism from $(2.2.9.2;n+1)$ to 
{\rm (\ref{cd:orexttlam})}.
\par
$(3)$ There exists the following 
commutative diagram of exact sequences$:$
\begin{equation*}
\begin{CD}
0 @>>> ({\cal W}{\Om}_{Y_t}^{\bul})'(-1)[-1] 
@>{\theta' \wedge}>>  
({\cal W}\wt{{\Om}}^{\bul}_{Y_{\os{\circ}{t}}})' @>>> 
({\cal W}{\Om}_{Y_t}^{\bul})'@>>> 0\\ 
@. @V{\vpl_nC^{-n},~\simeq}VV  
@V{\vpl_n\wt{C}{}^{-n},~\simeq}VV 
@V{\vpl_nC^{-n},~\simeq}VV \\
0 @>>> {\cal W}{\Om}_{Y_t}^{\bul}(-1)[-1] 
@>{\theta \wedge}>>  
{\cal W}\wt{{\Om}}^{\bul}_{Y_{\os{\circ}{t}}} @>>> 
{\cal W}{\Om}_{Y_t}^{\bul}@>>> 0. 
\end{CD}
\tag{2.2.9.4}\label{cd:limwlmext}
\end{equation*}
\end{theo}

\begin{coro}\label{coro:wtf}
$(1)$ The sheaves ${\cal W}\wt{\Om}_{Y_{\os{\circ}{t}}}^i$ 
and $({\cal W}\wt{\Om}_{Y_{\os{\circ}{t}}}^i)'$ $(i\in {\mab N})$ are 
sheaves of flat ${\cal O}_{{\cal W}(\os{\circ}{t})}$-modules.  
\par 
$(2)$ Let the notations be as in {\rm (\ref{prop:bcdw})}. 
Let $\star$ be nothing or $\prime$. 
Then the canonical morphism 
\begin{equation*} 
q^{-1}(({\cal W}_n\wt{\Om}_Y^i)^{\star})
\otimes_{{\cal W}_n}{\cal W}_n(\kap_t)
\lo ({\cal W}_n\wt{\Om}_{Y_{\os{\circ}{t}}}^i)^{\star}
\tag{2.2.10.1}\label{eqn:dkp}
\end{equation*}
is an isomorphism. 
\par 
$(3)$ Let $g\col Z\lo Y$ be a solid morphism from a fine log scheme over $s$. 
Assume that $\os{\circ}{g}$ is \'{e}tale.  Then the canonical morphism 
\begin{align*} 
{\cal W}_n({\cal O}_{Z})^{\star}
\otimes_{{\cal W}_n({\cal O}_{Y})^{\star}}
({\cal W}_n\wt{\Om}_{Y_{\os{\circ}{t}}}^i)^{\star} 
\lo ({\cal W}_n\wt{\Om}_{Z_{\os{\circ}{t}}}^i)^{\star} \quad (i\in {\mab Z})
\tag{2.2.10.2}\label{ali:owtzs}
\end{align*} 
is an isomorphism. 
\end{coro} 
\begin{proof} 
(1): By \cite[(4.5) (1)]{hk} or \cite[(6.8) (2)]{ndw}, 
${\cal W}\Om_{Y_{\os{\circ}{t}}}^i$ $(i\in {\mab N})$ 
is a sheaf of flat ${\cal W}(\kap_t)$-modules.  
Hence (1) follows from (\ref{cd:lijplext}) and  (\ref{cd:limwlmext}). 
((1) also immediately follows from \cite[(6.8) (3)]{ndw}.) 
\par 
(2): (2) follows from (\ref{prop:plz}), (\ref{prop:bcdw}) and (\ref{theo:tiobre}).
\par 
(3): It suffices to prove (3) in the case $\star=\prime$. 
Because $\os{\circ}{g}$ is \'{e}tale, so is 
${\cal W}_n(\os{\circ}{Z})\lo {\cal W}_n(\os{\circ}{Y})$ (\cite[0 (1.5.8)]{idw}). 
Hence we have the following exact sequence: 
\begin{align*} 
0 &\lo {\cal W}_n({\cal O}_{Z})\otimes_{{\cal W}_n({\cal O}_{Y})}
({\cal W}_n{{\Om}}_{Y}^i)'(-1,u)[-1] 
\os{\theta'_n  \wedge}{\lo} 
 {\cal W}_n({\cal O}_{Z})\otimes_{{\cal W}_n({\cal O}_{Y})}
({\cal W}_n\wt{{\Om}}^i_{Y})' \\
&\lo {\cal W}_n({\cal O}_{Z})\otimes_{{\cal W}_n({\cal O}_{Y})}
({\cal W}_n{{\Om}}_{Y}^i)'\lo  0. 
\end{align*} 
By (\ref{ali:ozs}) and (\ref{theo:tiobre}), 
we see that the morphism (\ref{ali:owtzs}) is an isomorphism. 
\end{proof} 

\begin{coro}\label{coro:nyt}
The following formula holds$:$
\begin{align*} 
{\cal W}_n\wt{\Om}^i_{Y_{\os{\circ}{t}}}
={\cal W}_n(\kap_t)\otimes_{{\cal W}_n}{\cal W}_n\wt{\Om}^i_Y
 \quad (i\in {\mab N}).
\tag{2.2.11.1}\label{ali:lmoiwy} 
\end{align*} 
\end{coro} 
\begin{proof} 
This follows from (\ref{prop:spbc}) and (\ref{coro:wtf}) (1). 
\end{proof}

\par 
Let $\ol{E}$ be a quasi-coherent log crystal of 
${\cal O}_{Y_{\os{\circ}{t}}/{\cal W}(\os{\circ}{t})}$-modules. 
Set $\ol{E}_{n}:=\ol{E}_{{\cal W}_n(Y_{\os{\circ}{t}})}$. 
By (\ref{eqn:intcon}) 
we have the following integrable connection 
\begin{equation*} 
\ol{E}_{n}\lo \ol{E}_{n}\otimes_{{\cal W}_n({\cal O}_{Y_t})'}
\Om^1_{{\cal W}_n({Y_{\os{\circ}{t}}})
/{\cal W}_n(\os{\circ}{t}),[~]},  
\tag{2.2.11.2}\label{eqn:eettnny}  
\end{equation*} 
which gives the complex 
$\ol{E}_{n}\otimes_{{\cal W}_n({\cal O}_{Y_t})'}
{\Om}^{\bul}_{{\cal W}_n({Y_t})/{\cal W}_n(\os{\circ}{t}),[~]}$ 
as in \cite[II (1.1.5)]{et}.    
As in \cite[pp.~251--252]{hk}, the proof of \cite[(4.19)]{hk} and \cite[p.~589]{ndw},
we can define the following morphism of complexes 
by the same formula as that in [loc.~cit.] 
in spite of the non-log smoothness of $Y_{\os{\circ}{t}}$ over $\os{\circ}{t}$:  
\begin{align*} 
\Om^{\bul}_{{\cal W}_n(Y_{\os{\circ}{t}})/{\cal W}_n(\os{\circ}{t}),[~]}
\lo  {\cal W}_n\wt{\Om}^{\bul}_{Y_{\os{\circ}{t}}}.
\tag{2.2.11.3}\label{eqn:lttmawy}
\end{align*} 
Indeed, if there exists an immersion $Y_{\os{\circ}{t}}\os{\sus}{\lo}{\cal Q}$ 
into a log smooth integral scheme over ${\cal W}_n(s_{\os{\circ}{t}})$, 
the morphism (\ref{eqn:lttmawy}) is equal to the following morphism: 
\begin{align*} 
\Om^{\bul}_{{\cal W}_n(Y_{\os{\circ}{t}})/{\cal W}_n(\os{\circ}{t}),[~]}
\lo 
{\cal H}^{\bul}({\cal O}_{\mathfrak E}\otimes_{{\cal O}_{\cal Q}}
\Om^*_{{\cal Q}/{\cal W}_n(\os{\circ}{t})}),  
\tag{2.2.11.4}\label{eqn:lttmwy}
\end{align*} 
where ${\mathfrak E}$ is the log PD-envelope of the immersion 
$Y_{\os{\circ}{t}}\os{\sus}{\lo}{\cal Q}$ over $({\cal W}_n(s_{\os{\circ}{t}}),p{\cal W}_n,[~])$. 
We define the morphism (\ref{eqn:lttmwy}) by the similar morphism to 
(\ref{ali:owny}), (\ref{eqn:tisn10}) and (\ref{eqn:tdlogb}). 
We note that 
a local section $d\log m$ $(m\in M_{Y_{\os{\circ}{t}}})$ 
is mapped to the cohomology class of 
$d\log \wt{m}$ by the morphism (\ref{eqn:lttmwy}), 
where $\wt{m}\in M_{\cal Q}$ is any lift of $m$. 
This morphism is compatible with the 
${\cal W}_n({\cal O}_{Y_{\os{\circ}{t}}})'$-module structure and 
the ${\cal W}_n({\cal O}_{Y_{\os{\circ}{t}}})$-module structure for each morphism 
${\Om}^i_{{\cal W}_n({Y_{\os{\circ}{t}}})/{\cal W}_n(\os{\circ}{t}), [~]}
\lo {\cal W}_n\wt{{\Om}}^i_{Y_{\os{\circ}{t}}}$ $(i\in {\mab Z})$. 
Hence we have the log de Rham-Witt complex 
$\ol{E}_{n}\otimes_{{\cal W}_n({\cal O}_{Y_{\os{\circ}{t}}})}
{\cal W}_n\wt{{\Om}}^{\bul}_{Y_{\os{\circ}{t}}}$. 
By (\ref{eqn:snisot}) we have the log de Rham-Witt complex 
$\ol{E}_{n}\otimes_{{\cal W}_n({\cal O}_{Y_{\os{\circ}{t}}})'}
({\cal W}_n\wt{{\Om}}^{\bul}_{Y_{\os{\circ}{t}}})'$. 
Let $\star$ be nothing or $\prime$. 
We have the following morphism 
\begin{equation*} 
R\col \ol{E}_{n+1}\otimes_{{\cal W}_{n+1}({\cal O}_{Y_{\os{\circ}{t}}})^{\star}}
({\cal W}_{n+1}\wt{{\Om}}^i_{Y_{\os{\circ}{t}}})^{\star}
\owns x\otimes y \lom R(x)\otimes R(y) \in 
\ol{E}_{n}\otimes_{{\cal W}_n({\cal O}_{Y_{\os{\circ}{t}}})^{\star}}
({\cal W}_n\wt{{\Om}}^i_{Y_{\os{\circ}{t}}})^{\star} 
\quad (i\in {\mab N}). 
\tag{2.2.11.5}\label{eqn:rettny}
\end{equation*}
As in the case for 
$E_n\otimes_{{\cal W}_n({\cal O}_{Y_{\os{\circ}{t}}})}{\cal W}_n\Om^{\bul}_{Y_{\os{\circ}{t}}}$ 
in the previous section, 
we can define the {\it canonical filtration} ${\rm Fil}$ 
on $\ol{E}_{n}
\otimes_{{\cal W}_n({\cal O}_{Y_{\os{\circ}{t}}})^{\star}}
({\cal W}_n\wt{{\Om}}^i_{Y_{\os{\circ}{t}}})^{\star}$ $(0\leq r\leq n)$ 
and the morphisms   
\begin{equation*} 
{\bf p}^r\col 
\ol{E}_{n}\otimes_{{\cal W}_n({\cal O}_{Y_{\os{\circ}{t}}})'}
({\cal W}_n\wt{{\Om}}_{Y_{\os{\circ}{t}}}^i)'
\lo 
\ol{E}_{n+r}\otimes_{{\cal W}_{n+r}({\cal O}_{Y_{\os{\circ}{t}}})'}
({\cal W}_{n+r}\wt{{\Om}}_{Y_{\os{\circ}{t}}}^i)'   
\tag{2.2.11.6}\label{eqn:pretwny} 
\end{equation*} 
and 
\begin{equation*} 
{\bf p}^r\col 
\ol{E}_{n}\otimes_{{\cal W}_n({\cal O}_{Y_{\os{\circ}{t}}})}{\cal W}_n\wt{{\Om}}_{Y_{\os{\circ}{t}}}^i
\lo 
\ol{E}_{n+r}
\otimes_{{\cal W}_{n+r}({\cal O}_{Y_{\os{\circ}{t}}})}
{\cal W}_{n+r}\wt{{\Om}}_{Y_{\os{\circ}{t}}}^i.  
\tag{2.2.11.7}\label{eqn:preteny} 
\end{equation*} 
If $\ol{E}$ is flat and coherent, then 
the two ${\bf p}^r$'s are injective by (\ref{prop:nop}) (2). 
\par 
Set $E:=\eps_{Y_t/{\cal W}_n(s_{\os{\circ}{t}})/{\cal W}_n(\os{\circ}{t})}^*(\ol{E})$ 
and $E_{n}:=E_{{\cal W}_n(Y_t)}$. 
By (\ref{eqn:eefnny}) we have the following integrable connection 
\begin{equation*} 
E_n\lo E_n
\otimes_{{\cal W}_n({\cal O}_{Y_t})'}
{\Om}^1_{{\cal W}_n(Y_t)/{\cal W}_n(t),[~]}.   
\tag{2.2.11.8}\label{eqn:eenny} 
\end{equation*}  
Then we have the following log de Rham complexes by 
(\ref{eqn:ldwy}) and (\ref{eqn:lmdy}):  
\begin{equation*} 
E_n\otimes_{{\cal W}_n({\cal O}_{Y_t})}
{\cal W}_n{\Om}^{\bul}_{Y_t} \quad 
{\rm and} \quad  
E_n\otimes_{{\cal W}_n({\cal O}_{Y_t})'}
({\cal W}_n{\Om}^{\bul}_{Y_t})'.
\tag{2.2.11.9}\label{eqn:eenytny} 
\end{equation*}  

\begin{lemm}\label{lemm:ntetfi}   
Let $Y$ $($resp.~$\ol{\cal Q})$ 
be a log smooth scheme of Cartier type over $s$ 
$($resp.~a log smooth integral scheme over $\ol{{\cal W}_n(s_{\os{\circ}{t}})})$.   
Let $Y_{\os{\circ}{t}} \os{\sus}{\lo} \ol{\cal Q}$ 
be an immersion over $s_{\os{\circ}{t}}\lo \ol{{\cal W}_n(s_{\os{\circ}{t}})}$ and 
let $\ol{\mathfrak E}$ be the log PD-envelope of this immersion 
over $({\cal W}_n(\os{\circ}{t}),p{\cal W}_n,[~])$. 
Set ${\cal Q}:=\ol{\cal Q}\times_{\ol{{\cal W}_n(s_{\os{\circ}{t}})}}{\cal W}_n(s_{\os{\circ}{t}})$. 
Let $g \col Y_{\os{\circ}{t}} \lo {\cal W}_n(s_{\os{\circ}{t}})$ be the structural morphism.  
Let $\ol{E}$ be a quasi-coherent log crystal of 
${\cal O}_{Y_{\os{\circ}{t}}/{\cal W}_n(\os{\circ}{t})}$-modules 
and let $(\ol{\cal E},\ol{\nabla})$ be the corresponding quasi-coherent 
${\cal O}_{\ol{\mathfrak E}}$-module with integrable connection. 
Set 
${\mathfrak E}
:=\ol{\mathfrak E}\times_{{\mathfrak D}(\ol{{\cal W}_n(s_{\os{\circ}{t}})})}{\cal W}_n(s_{\os{\circ}{t}})$ 
and $({\cal E},\nabla):=(\ol{\cal E},\ol{\nabla})\otimes_{{\cal O}_{\ol{\mathfrak E}}}
{\cal O}_{\mathfrak E}$. 
Assume that $\os{\circ}{Y}$ is affine. 
Then the following hold$:$
\par 
$(1)$ There exists a morphism 
${\cal W}_n(Y_{\os{\circ}{t}})\lo \ol{\mathfrak E}$ over $\ol{{\cal W}_n(s_{\os{\circ}{t}})}$ 
such that the composite morphism 
$Y_{\os{\circ}{t}}\os{\sus}{\lo} {\cal W}_n(Y_{\os{\circ}{t}})\lo \ol{\mathfrak E}$ 
is the immersion obtained by the given immersion 
$Y_{\os{\circ}{t}}\os{\sus}{\lo} \ol{\cal Q}$.
\par 
$(2)$ 
Let ${\cal W}_n(Y_{\os{\circ}{t}})\lo  \ol{\mathfrak E}$ be the morphism 
over $\ol{{\cal W}_n(s_{\os{\circ}{t}})}$ in $(1)$. 
Then this morphism induces the following morphism of complexes$:$  
\begin{equation*} 
{\cal E}\otimes_{{\cal O}_{\cal Q}}
\Om^{\bul}_{{\cal Q}/{\cal W}_n(\os{\circ}{t})} 
\lo 
\ol{E}_{n}\otimes_{{\cal W}_n({\cal O}_{Y_{\os{\circ}{t}}})}
{\cal H}^{\bul}({\cal O}_{\mathfrak E}\otimes_{{\cal O}_{\cal Q}}
\Om^*_{{\cal Q}/{\cal W}_n(\os{\circ}{t})}).  
\tag{2.2.12.1}\label{eqn:fnithtdw} 
\end{equation*} 
This morphism is functorial in the following sense$:$
\par 
For the following commutative diagram  
\begin{equation*} 
\begin{CD} 
t_1 @>>> s_1 \\ 
@VVV @VVV \\
t_2 @>>> s_2
\end{CD} 
\tag{2.2.12.2}\label{cd:imtpol}
\end{equation*} 
of fine log schemes whose underlying schemes 
are the spectrums of perfect fields of characteristic $p>0$ 
and for 
a log smooth affine scheme $Y_{i,\os{\circ}{t}_i}$ 
of Cartier type over a log point $s_{i,\os{\circ}{t}_i}$ $(i=1,2)$ 
whose underlying scheme 
is the spectrum of a perfect field of characteristic $p>0$ 
and for a morphism ${\cal W}_n(Y_{i,{\os{\circ}{t}_i}})\lo \ol{\cal Q}_i$ 
to a log smooth integral scheme $\ol{\cal Q}_i$ over $\ol{{\cal W}_n(s_{i,\os{\circ}{t}_i})}$ 
such that the composite morphism 
$Y_{i,\os{\circ}{t}_i} \os{\sus}{\lo} {\cal W}_n(Y_{i,\os{\circ}{t}_i})\lo \ol{\cal Q}_i$ 
is an immersion 
and for the log PD-envelope $\ol{\mathfrak E}_i$ $(i=1,2)$ of 
the immersion $Y_{i,\os{\circ}{t}_i} \os{\sus}{\lo} \ol{\cal Q}_i$ 
over $(\ol{{\cal W}_n(s_{i,\os{\circ}{t}_i})},p{\cal O}_{\ol{{\cal W}_n(s_{i,\os{\circ}{t}_i})}},[~])$ 
and 
for the following commutative diagram 
\begin{equation*} 
\begin{CD} 
{\cal W}_n(Y_{1,\os{\circ}{t}_1}) @>>> \ol{\mathfrak E}_1 \\ 
@V{h_n}VV @VV{\ol{h}}V \\ 
{\cal W}_n(Y_{2,\os{\circ}{t}_2}) @>>> \ol{\mathfrak E}_2
\end{CD} 
\tag{2.2.12.3}\label{cd:fndw} 
\end{equation*} 
of log schemes over the commutative diagram 
$($cf.~{\rm (\ref{prop:bar}))} 
\begin{equation*} 
\begin{CD} 
{\cal W}_n(s_{1,\os{\circ}{t}_1}) @>>> \ol{{\cal W}_n(s_{1,\os{\circ}{t}_1})} \\ 
@VVV @VVV \\ 
{\cal W}_n(s_{2,\os{\circ}{t}_2}) @>>> \ol{{\cal W}_n(s_{2,\os{\circ}{t}_2})}
\end{CD} 
\tag{2.2.12.4}\label{cd:impol}
\end{equation*}  
such that the composite morphism 
$Y_{i,\os{\circ}{t}_i}\os{\sus}{\lo} {\cal W}_n(Y_{i,\os{\circ}{t}}) \lo \ol{\mathfrak E}_i$ 
is the immersion obtained by the immersion 
$Y_{i,\os{\circ}{t}_i}\os{\sus}{\lo} \ol{\cal Q}_i$ 
and for a morphism $h_{1{\rm crys}}^*(E_2)\lo E_1$ 
$(E_i$ $(i=1,2):$ a quasi-coherent crystal of 
${\cal O}_{Y_{i,\os{\circ}{t}_i}/{\cal W}_n(\os{\circ}{t}_i)}$-modules$)$ of 
${\cal O}_{Y_{1,\os{\circ}{t}_1}/{\cal W}_n(\os{\circ}{t}_1)}$-modules,  
the following diagram is commutative$:$ 
\begin{equation*}  
\begin{CD} 
{\cal E}_1{\otimes}_{{\cal O}_{{\cal Q}_1}}
\Om^{\bul}_{{\cal Q}_1/{\cal W}_n(\os{\circ}{t}_1)}@>>> 
\ol{E}_{1,n}{\otimes}_{{\cal W}_n({\cal O}_{Y_{1,\os{\circ}{t}_1}})}
{\cal H}^{\bul}({\cal O}_{{\mathfrak E}_1} 
{\otimes}_{{\cal O}_{{\cal Q}_1}}
\Om^{*}_{{\cal Q}_1/{\cal W}_n(\os{\circ}{t}_1)}) \\ 
@A{h^{{\rm PD}*}}AA @AA{h^*_n}A \\ 
h^*({\cal E}_2{\otimes}_{{\cal O}_{{\cal Q}_2}}
\Om^{\bul}_{{\cal Q}_2/{\cal W}_n(\os{\circ}{t}_2)})
@>>> 
h_n^*(\ol{E}_{2,n}{\otimes}_{{\cal W}_n({\cal O}_{Y_{2,\os{\circ}{t}_2}})}
{\cal H}^{\bul}({\cal O}_{{\mathfrak E}_2} 
{\otimes}_{{\cal O}_{{\cal Q}_2}}
\Om^{*}_{{\cal Q}_2/{\cal W}_n(\os{\circ}{t}_2)})),   
\end{CD} 
\tag{2.2.12.5}\label{cd:odhtqph}
\end{equation*} 
where ${\cal E}_i:=\ol{E}_{\ol{\mathfrak E}_i}
\otimes_{{\cal O}_{\mathfrak D}(\ol{{\cal W}_n(s_{i,\os{\circ}{t}_i})})}
{\cal O}_{{\cal W}_n(s_{i,\os{\circ}{t}_i})}$. 
The morphism {\rm (\ref{eqn:fnithtdw})} is compatible 
with the projections.  
\end{lemm} 
\begin{proof} 
(1): The proof is the same as that of (\ref{lemm:neh3fi}) (1). 
\par 
(2): By (\ref{eqn:fwuou}) we have the following isomorphism 
\begin{align*} 
{\cal E}\otimes_{{\cal O}_{{\cal Q}}}
\Om^i_{{\cal Q}/{\cal W}_n(\os{\circ}{t})}  \os{\sim}{\lo} 
{\cal E}\otimes_{{\cal O}_{{\mathfrak E}}}
\Om^i_{{\mathfrak E}/{\cal W}_n(\os{\circ}{t}),[~]}  \quad (i\in {\mab N}). 
\tag{2.2.12.6}\label{ali:actdw}
\end{align*} 
The morphism 
${\cal W}_n(Y_{\os{\circ}{t}})
\lo \ol{\mathfrak E}$ gives the following morphism 
\begin{align*} 
\ol{\cal E}\otimes_{{\cal O}_{\ol{\mathfrak E}}}
\Om^i_{\ol{\mathfrak E}/{\cal W}_n(\os{\circ}{t}),[~]}  
{\lo} \ol{E}_n\otimes_{{\cal W}_n({\cal O}_{Y_{\os{\circ}{t}}})'}
\Om^i_{{\cal W}_n(Y_{\os{\circ}{t}})/{\cal W}_n(\os{\circ}{t}),[~]}. 
\tag{2.2.12.7}\label{ali:wstn}
\end{align*} 
Because the morphism 
${\cal W}_n(Y_{\os{\circ}{t}})\lo \ol{\mathfrak E}$ factors through 
the morphism ${\cal W}_n(Y_{\os{\circ}{t}})\lo {\mathfrak E}$, 
we have the following morphism 
\begin{align*} 
{\cal E}\otimes_{{\cal O}_{\mathfrak E}}
\Om^i_{{\mathfrak E}/{\cal W}_n(\os{\circ}{t}),[~]}  
{\lo} \ol{E}_n\otimes_{{\cal W}_n({\cal O}_{Y_{\os{\circ}{t}}})'}
\Om^i_{{\cal W}_n(Y_{\os{\circ}{t}})/{\cal W}_n(\os{\circ}{t}),[~]}. 
\tag{2.2.12.8}\label{ali:wtsn}
\end{align*} 
By using the morphism (\ref{eqn:lttmwy}), 
we have the following morphism  
\begin{align*} 
\ol{E}_n\otimes_{{\cal W}_n({\cal O}_{Y_{\os{\circ}{t}}})'}
\Om^i_{{\cal W}_n(Y_{\os{\circ}{t}})/{\cal W}_n(\os{\circ}{t}),[~]}
\lo 
\ol{E}_n\otimes_{{\cal W}_n({\cal O}_{Y_{\os{\circ}{t}}})}
{\cal H}^i({\cal O}_{\mathfrak E}\otimes_{{\cal O}_{\cal Q}}
\Om^*_{{\cal Q}/{\cal W}_n(\os{\circ}{t})}).
\tag{2.2.12.9}\label{ali:ttm}
\end{align*}  
By (\ref{ali:actdw}), (\ref{ali:wtsn}) and (\ref{ali:ttm}), 
we obtain the following composite morphism:  
\begin{align*} 
{\cal E}\otimes_{{\cal O}_{\cal Q}}
\Om^i_{{\cal Q}/{\cal W}_n(\os{\circ}{t})}  
&\lo \ol{E}_n\otimes_{{\cal W}_n({\cal O}_{Y_{\os{\circ}{t}}})'}
\Om^i_{{\cal W}_n(Y_{\os{\circ}{t}})/{\cal W}_n(\os{\circ}{t}),[~]} 
\lo 
\ol{E}_n\otimes_{{\cal W}_n({\cal O}_{Y_{\os{\circ}{t}}})}
{\cal H}^i({\cal O}_{\mathfrak E}\otimes_{{\cal O}_{\cal Q}}
\Om^*_{{\cal Q}/{\cal W}_n(\os{\circ}{t})} ) \tag{2.2.12.10}\label{eqn:fctdw}.   
\end{align*} 
As in the proof of (\ref{theo:ccrw}), 
it is easy to check that the morphism (\ref{eqn:fctdw}) induces 
the morphism of complexes.
This is the desired morphism (\ref{eqn:fnithtdw}).    
\par
As in the proof of (\ref{theo:ccrw}),  
the functoriality of the morphism (\ref{eqn:fctdw}) is clear.  
\par  
The compatibility of the morphism (\ref{eqn:fnithtdw}) 
with the projections follows from \cite[(7.18)]{ndw}. 
\end{proof}


\par 
Now let $N$ be a nonnegative integer or $\infty$. 
Let $Y_{\bul \leq N}$ be 
a log smooth $N$-truncated 
simplicial log scheme of Cartier type over $s$. 
Set $Y_{\bul \leq N,t}:=Y_{\bul \leq N}\times_st$ 
and $Y_{\bul \leq N,\os{\circ}{t}}:=Y_{\bul \leq N}\times_{\os{\circ}{s}}\os{\circ}{t}$.  
Let $\ol{E}{}^{\bul \leq N}$ be a flat coherent log crystal of 
${\cal O}_{Y_{\bul \leq N,\os{\circ}{t}}/{\cal W}_n(\os{\circ}{t})}$-modules. 
Set $\ol{E}{}^{\bul \leq N}_n:=
\ol{E}{}^{\bul \leq N}_{{\cal W}_n(Y_{\bul \leq N,\os{\circ}{t}})}$. 
Then we have the 
$N$-truncated cosimplicial log de Rham-Witt complex
$\ol{E}{}^{\bul \leq N}_n
\otimes_{{\cal W}_n({\cal O}_{Y_{\bul \leq N,\os{\circ}{t}}})^{\star}}
({\cal W}_n\wt{{\Om}}^{\bul}_{Y_{\bul \leq N,\os{\circ}{t}}})^{\star}$. 
Set $E^{\bul \leq N}:=
\eps^*_{Y_{\bul \leq N,\os{\circ}{t}}/{\cal W}_n(s_{\os{\circ}{t}})/{\cal W}_n(\os{\circ}{t})}
(\ol{E}{}^{\bul \leq N})$. 
Let $g_t \col Y_{\bul \leq N,t} \lo t \os{\sus}{\lo} {\cal W}_n(t)$ 
be the structural morphism.  

\par 
By (\ref{cd:orexttlam}) and (\ref{cd:limwlmext}), 
we have the following commutative diagrams: 
\begin{equation*}
\begin{CD}
0 @>>> 
E^{\bul \leq N}_{n}\otimes_{{\cal W}_n({\cal O}_{Y_{\bul \leq N,\os{\circ}{t}}})'}
({\cal W}_n{\Om}^{\bul}_{Y_{\bul \leq N,\os{\circ}{t}}})'(-1,u)[-1] 
@>{{\rm id}\otimes (\theta'_n  \wedge)}>>  \\ 
@. @V{{\rm id}\otimes C^{-n}}V{\simeq}V  \\
0 @>>> E^{\bul \leq N}_{n}\otimes_{{\cal W}_n({\cal O}_{{Y_{\bul \leq N,\os{\circ}{t}}}})}
{\cal W}_n{\Om}^{\bul}_{Y_{\bul \leq N,\os{\circ}{t}}}(-1,u)[-1] 
@>{{\rm id}\otimes (\theta_n \wedge)}>>  
\end{CD}
\tag{2.2.12.11;$n$}\label{cd:oreeam}
\end{equation*}
\begin{equation*}
\begin{CD}
\ol{E}{}^{\bul \leq N}_{n}\otimes_{{\cal W}_n({\cal O}_{{Y_{\bul \leq N,\os{\circ}{t}}}})'}
({\cal W}_n\wt{{\Om}}^{\bul}_{Y_{\bul \leq N,\os{\circ}{t}}})' @>>> 
E^{\bul \leq N}_{n}\otimes_{{\cal W}_n({\cal O}_{{Y_{\bul \leq N,\os{\circ}{t}}}})'}
({\cal W}_n{\Om}^{\bul}_{Y_{\bul \leq N,\os{\circ}{t}}})'@>>> 0\\ 
@V{{\rm id}\otimes \wt{C}{}^{-n}}V{\simeq}V @V{{\rm id}\otimes C^{-n}}V{\simeq}V\\
\ol{E}{}^{\bul \leq N}_{n}\otimes_{{\cal W}_n({\cal O}_{{Y_{\bul \leq N,\os{\circ}{t}}}})}
{\cal W}_n\wt{\Om}^{\bul}_{Y_{\bul \leq N,\os{\circ}{t}}}
@>>> E^{\bul \leq N}_{n}\otimes_{{\cal W}_n({\cal O}_{{Y_{\bul \leq N,\os{\circ}{t}}}})}
{\cal W}_n{\Om}^{\bul}_{Y_{\bul \leq N,\os{\circ}{t}}}@>>> 0, 
\end{CD}
\end{equation*}
\begin{equation*}
\begin{CD}
0 @>>> 
\vpl_n(E^{\bul \leq N}_{n})\otimes_{{\cal W}({\cal O}_{{Y_{\bul \leq N,\os{\circ}{t}}}})'}
({\cal W}{\Om}^{\bul}_{Y_{\bul \leq N,\os{\circ}{t}}})'(-1,u)[-1] 
   \\ 
@. @V{{\rm id}\otimes \vpl_n(C^{-n})}V{\simeq}V   \\
0 @>>> 
\vpl_n(E^{\bul \leq N}_{n})\otimes_{{\cal W}({\cal O}_{Y_{\bul \leq N,\os{\circ}{t}}})}
{\cal W}{\Om}^{\bul}_{Y_{\bul \leq N,\os{\circ}{t}}}(-1,u)[-1] 
\end{CD}
\tag{2.2.12.12}\label{cd:liemext}
\end{equation*}
\begin{equation*}
\begin{CD}
@>{{\rm id}\otimes (\theta' \wedge)}>>
\vpl_n(\ol{E}{}^{\bul \leq N}_{n})\otimes_{{\cal W}({\cal O}_{Y_{\bul \leq N,\os{\circ}{t}}})'}
({\cal W}\wt{\Om}{}^{\bul}_{Y_{\bul \leq N,\os{\circ}{t}}})'@>>> 
\\
@. @V{{\rm id}\otimes \vpl_n(\wt{C}{}^{-n})}V{\simeq}V \\
@>{{\rm id}\otimes (\theta \wedge)}>> 
\vpl_n(\ol{E}{}^{\bul \leq N}_{n})\otimes_{{\cal W}({\cal O}_{Y_{\bul \leq N,\os{\circ}{t}}})}
{\cal W}\wt{\Om}^{\bul}_{Y_{\bul \leq N,\os{\circ}{t}}}
@>>> 
\end{CD}
\end{equation*}
\begin{equation*}
\begin{CD}
\vpl_n(E^{\bul \leq N}_{n})\otimes_{{\cal W}({\cal O}_{Y_{\bul \leq N,\os{\circ}{t}}})'}
({\cal W}{\Om}^{\bul}_{Y_{\bul \leq N,\os{\circ}{t}}})'@>>> 0\\
@V{{\rm id}\otimes \vpl_n(C^{-n})}V{\simeq}V\\
\vpl_n(E^{\bul \leq N}_{n})\otimes_{{\cal W}({\cal O}_{Y_{\bul \leq N,\os{\circ}{t}}})}
{\cal W}{\Om}^{\bul}_{Y_{\bul \leq N,\os{\circ}{t}}}@>>> 0.
\end{CD}
\end{equation*}
The horizontal lines in (\ref{cd:oreeam}) and (\ref{cd:liemext}) 
are exact.

\begin{prop}[\bf Contravariant functoriality]\label{prop:ctwtu}
$(1)$ Let $u_n$ be as in {\rm (\ref{eqn:xdxdrss})}. 
Let 
\begin{equation*} 
\begin{CD} 
Y_{\bul \leq N,\os{\circ}{t}} @>{g_{\bul \leq N}}>> Z_{\bul \leq N,\os{\circ}{t}{}'}\\ 
@VVV @VVV \\ 
s_{\os{\circ}{t}} @>>> s'_{\os{\circ}{t}{}'} 
\end{CD}
\tag{2.2.13.1}\label{eqn:xwdss}
\end{equation*} 
be a commutative diagram of 
$N$-truncated simplicial log smooth schemes over 
$s_{\os{\circ}{t}}$ and $s'_{\os{\circ}{t}{}'}$.  
Let $\ol{F}{}^{\bul \leq N}$ be 
a flat quasi-coherent crystal of 
${\cal O}_{Z_{\bul \leq N,\os{\circ}{t}{}'}/{\cal W}_n(\os{\circ}{t}{}')}$-modules. 
Let 
$$h \col g^*_{\bul \leq N,{\rm crys}}(\ol{F}{}^{\bul \leq N})\lo \ol{E}{}^{\bul \leq N}$$ 
be a morphism of flat quasi-coherent crystals of 
${\cal O}_{Y_{\bul \leq N,\os{\circ}{t}}/{\cal W}_n(\os{\circ}{t})}$-modules. 
Then there exists a natural morphism 
\begin{align*} 
\wt{g}{}^*_{\bul \leq N} \col 
\ol{F}{}^{\bul \leq N}_{n}\otimes_{{\cal W}_n({\cal O}_{{Z_{\bul \leq N,\os{\circ}{t}}}})}
({\cal W}_n\wt{\Om}^{\bul}_{Z_{\bul \leq N,\os{\circ}{t}}})^{\star}
\lo 
Rg_{\bul \leq N*}(\ol{E}{}^{\bul \leq N}_{n}\otimes_{{\cal W}_n({\cal O}_{{Y_{\bul \leq N,\os{\circ}{t}}}})}
({\cal W}_n\wt{\Om}^{\bul}_{Y_{\bul \leq N,\os{\circ}{t}}})^{\star}), 
\tag{2.2.13.2}\label{ali:eozy}
\end{align*} 
which is compatible with the composition of $g_{\bul \leq N}$'s and 
$$\wt{({\rm id}_{Y_{\bul \leq N,\os{\circ}{T}_0}})}{}^*
={\rm id}_{\ol{E}{}^{\bul \leq N}_{n}\otimes_{{\cal W}_n({\cal O}_{{Y_{\bul \leq N,\os{\circ}{t}}}})}
{\cal W}_n\wt{\Om}^{\bul}_{Y_{\bul \leq N,\os{\circ}{t}}}}.$$ 
\par 
$(2)$ Set 
$F^{\bul \leq N}:=\eps_{Z_{\bul \leq N,T'_0}/
{\cal W}_n(s_{\os{\circ}{t}})/{\cal W}_n(\os{\circ}{t})}^*(\ol{F}{}^{\bul \leq N})$.  
Then the morphism $\wt{g}{}^*_{\bul \leq N}$ fits into the following commutative diagram of 
triangles$:$ 
\begin{equation*} 
\begin{CD} 
Rg_{\bul \leq N*}
(E^{\bul \leq N}_{n}\otimes_{{\cal W}_n({\cal O}_{{Y_{\bul \leq N,t}}})}
({\cal W}_n\Om^{\bul}_{Y_{\bul \leq N,t}})^{\star})[-1] 
@>>> 
\\ 
@A{g^*_{\bul \leq N}}AA \\
F^{\bul \leq N}_{n}
\otimes_{{\cal W}_n({\cal O}_{{Z_{\bul \leq N,t'}}})}
({\cal W}_n\Om^{\bul}_{Z_{\bul \leq N,t'}})^{\star}[-1] 
@>>> 
\end{CD} 
\tag{2.2.13.3}\label{cd:xawxy}
\end{equation*} 
\begin{equation*} 
\begin{CD} 
Rg_{\bul \leq N*}
(\ol{E}{}^{\bul \leq N}_{n}
\otimes_{{\cal W}_n({\cal O}_{{Y_{\bul \leq N,\os{\circ}{t}}}})}
({\cal W}_n\wt{\Om}^{\bul}_{Y_{\bul \leq N,\os{\circ}{t}}})^{\star}) @>>>\\ 
@AA{\wt{g}{}*_{\bul \leq N}}A \\
\ol{F}{}^{\bul \leq N}_{n}
\otimes_{{\cal W}_n({\cal O}_{{Z_{\bul \leq N,\os{\circ}{t}{}'}}})}
({\cal W}_n\wt{\Om}^{\bul}_{Z_{\bul \leq N,\os{\circ}{t}{}'}})^{\star}
@>>>
\end{CD} 
\end{equation*} 
\begin{equation*} 
\begin{CD} 
Rg_{\bul \leq N*}(E^{\bul \leq N}_{n}
\otimes_{{\cal W}_n({\cal O}_{{Y_{\bul \leq N,t}}})}
({\cal W}_n\Om^{\bul}_{Y_{\bul \leq N,t}})^{\star})@>{+1}>> \\
@A{g^*_{\bul \leq N}}AA \\
F^{\bul \leq N}_{n}
\otimes_{{\cal W}_n({\cal O}_{{Z_{\bul \leq N,\os{\circ}{t}{}'}}})}
({\cal W}_n{\Om}^{\bul}_{Z_{\bul \leq N,\os{\circ}{t}{}'}})^{\star}@>{+1}>>. 
\end{CD} 
\end{equation*} 
\end{prop} 
\begin{proof} 
(1): (1) follows from (\ref{prop:fuo}).   
\par 
(2): (2) follows from (\ref{prop:fuo}) and (\ref{prop:plz}). 
\end{proof}

\begin{theo}\label{theo:ccttrw}
Let $\ol{E}{}^{\bul \leq N}$ be a flat coherent log crystal of 
${\cal O}_{Y_{\bul \leq N,\os{\circ}{t}}/{\cal W}_n(\os{\circ}{s})}$-modules. 
Assume that $Y_{\bul \leq N,\os{\circ}{t}}$ has 
an affine $N$-truncated simplicial open covering. 
Then the following hold$:$ 
\par 
$(1)$ There exists a canonical isomorphism 
\begin{equation*} 
\wt{R}u_{Y_{\bul \leq N,\os{\circ}{t}}/{\cal W}_n(\os{\circ}{t})*}(\ol{E}{}^{\bul \leq N})
\os{\sim}{\lo} 
\ol{E}{}^{\bul \leq N}_{n}
\otimes_{{\cal W}_n({\cal O}_{Y_{\bul \leq N,\os{\circ}{t}}})^{\star}}
({\cal W}_n\wt{{\Om}}^{\bul}_{Y_{\bul \leq N,\os{\circ}{t}}})^{\star} 
\tag{2.2.14.1}\label{eqn:ywnttnny}
\end{equation*} 
in ${\rm D}^+(g^{-1}({\cal W}_n))$.
The isomorphisms $(\ref{eqn:ywnttnny})$ for $n$'s are 
compatible with the two projections of 
both hand sides on $(\ref{eqn:ywnttnny})$. 
The isomorphism $(\ref{eqn:ywnttnny})$ 
fits into the following commutative diagram 
of triangles$:$
\begin{equation*}
\begin{CD}
@>>> 
Ru_{Y_{\bul \leq N,t}/{\cal W}_n(t)*}({E}{}^{\bul \leq N})(-1,u_n)[-1] 
@>{}>>   \\ 
@. @V{\simeq}VV   \\
@>>> E^{\bul \leq N}_{n}\otimes_{{\cal W}_n({\cal O}_{Y_{\bul \leq N,t}})^{\star}}
({\cal W}_n{\Om}^{\bul}_{Y_{\bul \leq N,t}})^{\star}(-1,u_n)[-1] 
@>{{\rm id}\otimes ((\theta_n)^{\star} \wedge)}>>  
\end{CD}
\tag{2.2.14.2}\label{cd:oreedsam}
\end{equation*}
\begin{equation*}
\begin{CD}
\wt{R}u_{Y_{\bul \leq N,\os{\circ}{t}}/{\cal W}_n(\os{\circ}{t})*}(\ol{E}{}^{\bul \leq N})
@>>> Ru_{Y_{\bul \leq N,t}/{\cal W}_n(t)*}(E^{\bul \leq N}) @>{+1}>> \\ 
@V{\simeq}VV @V{\simeq}VV \\
\ol{E}{}^{\bul \leq N}_{n}\otimes_{{\cal W}_n({\cal O}_{Y_{\bul \leq N,t}})^{\star}}
({\cal W}_n\wt{\Om}^{\bul}_{Y_{\bul \leq N,t}})^{\star}
@>>> E^{\bul \leq N}_{n}\otimes_{{\cal W}_n({\cal O}_{Y_{\bul \leq N,t}})^{\star}}
({\cal W}_n{\Om}^{\bul}_{Y_{\bul \leq N,t}})^{\star}@>{+1}>> . 
\end{CD}
\end{equation*}
The commutative diagram {\rm (\ref{cd:oreedsam})} 
is compatible with $n$'s. Here 
the isomorphism 
$$Ru_{Y_{\bul \leq N,t}/{\cal W}_n(t)*}({E}{}^{\bul \leq N})
\os{\sim}{\lo} 
E^{\bul \leq N}_{n}\otimes_{{\cal W}_n({\cal O}_{Y_{\bul \leq N,t}})^{\star}}
({\cal W}_n{\Om}^{\bul}_{Y_{\bul \leq N,t}})^{\star}$$ 
in the commutative diagram above is the isomorphism {\rm (\ref{eqn:ywnnny})}. 
\par 
$(2)$ Let $u_n$ be as in {\rm (\ref{eqn:xdxdrss})}.  
The isomorphism $(\ref{eqn:ywnttnny})$ is functorial 
with respect to a morphism 
$g_{\bul \leq N}\col 
Y_{\bul \leq N,\os{\circ}{t}}\lo Z_{\bul \leq N,\os{\circ}{t}{}'}$ of 
$N$-truncated simplicial log schemes over the morphism 
$u_n\col {\cal W}_n(s_{\os{\circ}{t}})\lo {\cal W}_n(s'_{\os{\circ}{t}{}'})$ 
satisfying the condition {\rm (\ref{cd:xygxy})} 
and a morphism 
$g^*_{{\bul \leq N},{\rm crys}}(\ol{F}{}^{\bul \leq N})\lo \ol{E}{}^{\bul \leq N}$ 
of ${\cal O}_{Y_{\bul \leq N},\os{\circ}{t}
/{\cal W}_n(s_{\os{\circ}{t}})}$-modules, 
where $Z_{\bul \leq N}$ and $\ol{F}{}^{\bul \leq N}$ are similar objects to 
$Y_{\bul \leq N}$ and $\ol{E}{}^{\bul \leq N}$, respectively. 
\par 
$(3)$ Let the notations be as in $(2)$. 
The morphisms {\rm (\ref{ali:xxwy})}, {\rm (\ref{ali:eozy})} 
and the following morphisms 
$$g^*_{\bul \leq N}\col 
F^{\bul \leq N}_{n}\otimes_{{\cal W}_n({\cal O}_{Z_{\bul \leq N,t'}})^{\star}}
({\cal W}_n{\Om}^{\bul}_{Z_{\bul \leq N,t'}})^{\star}\lo 
Rg_{\bul \leq N*}(E^{\bul \leq N}_{n}\otimes_{{\cal W}_n({\cal O}_{Y_{\bul \leq N,t}})^{\star}}
({\cal W}_n{\Om}^{\bul}_{Y_{\bul \leq N,t}})^{\star}),$$
$$g^*_{\bul \leq N}\col 
\ol{F}{}^{\bul \leq N}_{n}\otimes_{{\cal W}_n({\cal O}_{{Z_{\bul \leq N,t'}}})}
({\cal W}_n\wt{\Om}{}^{\bul}_{Z_{\bul \leq N,\os{\circ}{t}{}'}})^{\star}\lo 
Rg_{\bul \leq N*}
(\ol{E}{}^{\bul \leq N}_{n}\otimes_{{\cal W}_n({\cal O}_{{Y_{\bul \leq N,t}}})}
({\cal W}_n\wt{\Om}{}^{\bul}_{Y_{\bul \leq N,\os{\circ}{t}}})^{\star})$$ 
induce a morphism from the commutative diagram 
{\rm (\ref{cd:oreedsam})} for 
$Z_{\bul \leq N,\os{\circ}{t}{}'}$ and $\ol{F}{}^{\bul \leq N}$ 
to the commutative diagram 
$Rg_{\bul \leq N*}(${\rm (\ref{cd:oreedsam}))}.
\end{theo}  
\begin{proof}
(1): By (\ref{cd:oreeam}) 
it suffices to prove (\ref{theo:ccttrw}) in the case $\star=$nothing. 
By using (\ref{lemm:lisj}) and Tsuzuki's functor $\Gam$, 
we have the canonical morphism (\ref{eqn:ywnttnny}) 
as in the proof of (\ref{theo:ccrw}). 
Indeed, let $Y'_{N,\os{\circ}{t}}$ be the disjoint union of a log affine open covering of 
$Y_{N,\os{\circ}{t}}$ 
such that there exists an immersion $Y'_{N,\os{\circ}{t}}\os{\sus}{\lo} \ol{\cal Y}{}'_N$ 
into a log smooth integral scheme over $\ol{{\cal W}_n(s_{\os{\circ}{t}})}$. 
Then there exists a morphism 
${\cal W}_n(Y'_{N,\os{\circ}{t}})\lo \ol{\cal Y}{}'_N$ over $\ol{{\cal W}_n(s_{\os{\circ}{t}})}$ 
such that the composite morphism 
$Y'_{N,\os{\circ}{t}}\os{\sus}{\lo} {\cal W}_n(Y'_{N,\os{\circ}{t}})\lo \ol{\cal Y}{}'_N$ 
is the given immersion.   
As in the proof of (\ref{theo:ccrw}), 
we have the following natural morphism 
\begin{equation*} 
{\cal W}_n(Y'_{\bul \leq N,\os{\circ}{t}})\lo 
{\Gam}{}^{\ol{{\cal W}_n(s_{\os{\circ}{t}})}}_N(\ol{\cal Y}{}'_N)_{\bul \leq N} 
\tag{2.2.14.3}\label{eqn:pwholge} 
\end{equation*} 
of $N$-truncated simplicial log schemes. 
Let 
${\cal W}_n(Y_{\bul \leq N,\bul,\os{\circ}{t}})$ be the \v{C}ech diagram of 
${\cal W}_n(Y'_{\bul \leq N,\os{\circ}{t}})$.  
Set $\ol{\cal Q}_{mn}
:={\rm cosk}_0^{\ol{{\cal W}_n(s_{\os{\circ}{t}})}}
(\Gam^{\ol{{\cal W}_n(s_{\os{\circ}{t}})}}_N(\ol{\cal Y}{}'_N)_m)_n$ 
$(0\leq m \leq N$, $n\in {\mab N})$. 
By (\ref{eqn:pwholge}) we have the following morphism 
\begin{equation*} 
{\cal W}_n(Y_{\bul \leq N,\bul,\os{\circ}{t}})\lo \ol{\cal Q}_{\bul \leq N,\bul}
\tag{2.2.14.4}\label{eqn:ppwht3gep} 
\end{equation*} 
of $(N,\infty)$-truncated bisimplicial log schemes. 
Set ${\cal Q}_{\bul \leq N,\bul}:=\ol{\cal Q}_{\bul \leq N,\bul}
\times_{\ol{{\cal W}_n(s)}}{\cal W}_n(s)$. 
Let $\ol{\mathfrak E}_{\bul \leq N,\bul}$ be the log PD-envelope of 
the immersion 
$Y_{\bul \leq N,\bul,\os{\circ}{t}}\os{\sus}{\lo} \ol{\cal Q}_{\bul \leq N,\bul}$ 
over $({\cal W}_n(\os{\circ}{t}),p{\cal W}_n,[~])$. 
Let $\ol{E}{}^{\bul \leq N,\bul}$ be the crystal of 
${\cal O}_{Y_{\bul \leq N,\bul}/{\cal W}_n(\os{\circ}{t})}$-modules 
obtained by $\ol{E}{}^{\bul \leq N}$. 
Let $(\ol{\cal E}{}^{\bul \leq N,\bul},\nabla)$ be the corresponding 
${\cal O}_{\ol{\mathfrak E}_{\bul \leq N,\bul}}$-module with integrable connection 
to $\ol{E}{}^{\bul \leq N,\bul}$. 
Set ${\mathfrak E}_{\bul \leq N,\bul}:=\ol{\mathfrak E}_{\bul \leq N,\bul}
\times_{{\mathfrak D}(\ol{{\cal W}_n(s)})}{\cal W}_n(s)$. 
Then ${\mathfrak E}_{\bul \leq N,\bul}$ is the log PD-envelope of 
the immersion ${\cal Y}_{\bul \leq N,\bul,\os{\circ}{t}}\os{\sus}{\lo} {\cal Q}_{\bul \leq N,\bul}$ 
over $({\cal W}_n(s_{\os{\circ}{t}}),p{\cal W}_n,[~])$. 
Set 
$({\cal E}^{\bul \leq N,\bul},\nabla):=
(\ol{\cal E}{}^{\bul \leq N,\bul},\nabla)\otimes_{{\cal O}_{\ol{\mathfrak E}_{\bul \leq N,\bul}}}
{\cal O}_{{\mathfrak E}_{\bul \leq N,\bul}}$: 
$$\nabla \col {\cal E}^{\bul \leq N,\bul}\lo {\cal E}^{\bul \leq N,\bul}
\otimes_{{\cal O}_{{\cal Q}_{\bul \leq N,\bul}}}\Om^1_{{\cal Q}_{\bul \leq N,\bul}
/{\cal W}_n(\os{\circ}{t})}.$$ 
The connection $\nabla$ induces the following connection: 
$$\nabla_{/{\cal W}_n(t)} \col {\cal E}^{\bul \leq N,\bul}\lo {\cal E}^{\bul \leq N,\bul}
\otimes_{{\cal O}_{{\cal Q}_{\bul \leq N,\bul}}}
\Om^1_{{\cal Q}_{\bul \leq N,\bul}/{\cal W}_n(t)}.$$ 
Set $\ol{E}{}^{\bul \leq N,\bul}_n
:=(\ol{E}{}^{\bul \leq N,\bul})_{{\cal W}_n(Y_{\bul \leq N,\bul,\os{\circ}{t}})}$. 
Then we have natural morphisms  
${\cal W}_n(Y_{\bul \leq N,\bul,\os{\circ}{t}})\lo \ol{\mathfrak E}_{\bul \leq N,\bul}$ 
and 
${\cal W}_n(Y_{\bul \leq N,\bul,\os{\circ}{t}})\lo {\mathfrak E}_{\bul \leq N,\bul}$. 
By (\ref{lemm:ntetfi})
we have the following morphism of complexes: 
\begin{align*} 
&{\cal E}^{\bul \leq N,\bul}
{\otimes}_{{\cal O}_{{\cal Q}_{\bul \leq N,\bul}}}
\Om^{\bul}_{{\cal Q}_{\bul \leq N,\bul}/{\cal W}_n(\os{\circ}{t})}   
\lo \tag{2.2.14.5}\label{ali:enotq}\\
&\ol{E}{}^{\bul \leq N,\bul}_{n}
\otimes_{{\cal W}_n({\cal O}_{Y_{\bul \leq N,\bul}})}
{\cal H}^{\bul}({\cal O}_{{\mathfrak E}_{\bul \leq N,\bul}}
{\otimes}_{{\cal O}_{{\cal Q}_{\bul \leq N,\bul}}}
\Om^{*}_{{\cal Q}_{\bul \leq N,\bul}/{\cal W}_n(\os{\circ}{t})} ).  
\end{align*} 
Applying $R\pi_{{\rm zar}*}$ ((\ref{eqn:tzar}))
to this morphism 
and using the log Poincar\'{e} lemma and the cohomological descent, 
we obtain the morphism (\ref{eqn:ywnttnny}) for the case $\star$=nothing. 
As in the proof of (\ref{theo:ccrw}) (1), we see that this morphism is independent 
of the choices of $Y'_{N,\os{\circ}{t}}$ and the immersion $Y'_{N,\os{\circ}{t}}\os{\sus}{\lo} \ol{\cal Y}{}'_N$. 
\par 
By the definition of the morphisms (\ref{ali:enoq}) and (\ref{ali:enotq}), 
we have the following commutative diagram with exact rows: 
\begin{equation*}
\begin{CD}  
0@>>> {\cal E}^{\bul \leq N,\bul}{\otimes}_{{\cal O}_{{\cal Q}_{\bul \leq N,\bul}}}
\Om^{\bul}_{{\cal Q}_{\bul \leq N,\bul}/{\cal W}_n(t)}[-1]   \\
@. @VVV \\
0@>>> E^{\bul \leq N,\bul}_{n}
\otimes_{{\cal W}_n({\cal O}_{Y_{\bul \leq N,\bul,\os{\circ}{t}}})}
{\cal H}^{\bul}({\cal O}_{{\mathfrak E}_{\bul \leq N,\bul}}
{\otimes}_{{\cal O}_{{\cal Q}_{\bul \leq N,\bul}}}
\Om^{*}_{{\cal Q}_{\bul \leq N,\bul}/{\cal W}_n(t)})[-1]  
\end{CD} 
\tag{2.2.14.6}\label{ali:enntq}
\end{equation*} 
\begin{equation*}
\begin{CD}  
@>{d\log \tau \wedge}>> {\cal E}^{\bul \leq N,\bul}
{\otimes}_{{\cal O}_{{\cal Q}_{\bul \leq N,\bul}}}
\Om^{\bul}_{{\cal Q}_{\bul \leq N,\bul}/{\cal W}_n(\os{\circ}{t})}   \\
@. @VVV \\
@>{d\log \tau \wedge}>> \ol{E}{}^{\bul \leq N,\bul}_{n}
\otimes_{{\cal W}_n({\cal O}_{Y_{\bul \leq N,\bul,\os{\circ}{t}}})}
{\cal H}^{\bul}({\cal O}_{{\mathfrak E}_{\bul \leq N,\bul}}
{\otimes}_{{\cal O}_{{\cal Q}_{\bul \leq N,\bul}}}
\Om^{*}_{{\cal Q}_{\bul \leq N,\bul}/{\cal W}_n(\os{\circ}{t})} )  
\end{CD}  
\end{equation*} 
\begin{equation*}
\begin{CD}  
@>>> {\cal E}^{\bul \leq N,\bul}
{\otimes}_{{\cal O}_{{\cal Q}_{\bul \leq N,\bul}}}
\Om^{\bul}_{{\cal Q}_{\bul \leq N,\bul}/{\cal W}_n(t)} @>>>0  \\
@. @VVV \\
@>>> E^{\bul \leq N,\bul}_{n}
\otimes_{{\cal W}_n({\cal O}_{Y_{\bul \leq N,\bul,\os{\circ}{t}}})}
{\cal H}^{\bul}({\cal O}_{{\mathfrak E}_{\bul \leq N,\bul}}
{\otimes}_{{\cal O}_{{\cal Q}_{\bul \leq N,\bul}}}
\Om^{*}_{{\cal Q}_{\bul \leq N,\bul}/{\cal W}_n(t)})  @>>> 0.
\end{CD}  
\end{equation*} 
Let 
\begin{equation*} 
\pi_{\rm zar} \col 
((Y_{\bul \leq N,\bul,t})_{\rm zar},g^{-1}_{\bul}({\cal W}_n)) \lo 
((Y_{\bul \leq N,t})_{\rm zar},g^{-1}({\cal W}_n)) 
\tag{2.2.14.7}\label{eqn:ppygzd} 
\end{equation*} 
be the natural morphism of ringed topoi. 
Applying $R\pi_{{\rm zar}*}$ to the diagram (\ref{ali:enntq}) and using 
the isomorphism (\ref{eqn:ywnnny}), 
we obtain the diagram (\ref{cd:oreedsam}). 
\par 
The compatibility of the morphism (\ref{eqn:ywnttnny})
with the projections follows from the proof of \cite[(7.18)]{ndw} 
((\ref{lemm:ntetfi})).
\par  
Now the compatibility of the commutative diagram (\ref{cd:oreedsam}) with $n$'s is clear. 
\par 
(2):  
Let $Z'_{N,\os{\circ}{t}{}'}$ be the disjoint union of a log affine open covering of 
$Z_{N,\os{\circ}{t}{}'}$ 
such that there exists an immersion $Z'_{N,\os{\circ}{t}{}'}\os{\sus}{\lo} \ol{\cal Z}{}'_N$ 
into a log smooth scheme over $\ol{{\cal W}_n(s'_{\os{\circ}{t}{}'})}$. 
Let $Y'_{N,\os{\circ}{t}}\os{\sus}{\lo} 
\ol{\cal Y}{}''_{N}$ and $Z'_{N,\os{\circ}{t}{}'}\os{\sus}{\lo} 
\ol{\cal Z}{}'_{N}$ 
be immersions into log smooth schemes over $\ol{{\cal W}_n(s_{\os{\circ}{t}})}$ and 
$\ol{{\cal W}(s'_{\os{\circ}{t}{}'})}$, respectively.  Then we have the morphisms 
${\cal W}_n(Y'_{N,\os{\circ}{t}}) \lo \ol{\cal Y}{}''_{N}$ and 
${\cal W}_n(Z'_{N,\os{\circ}{t}{}'}) \lo \ol{\cal Z}{}'_{N}$ 
such that the composite morphisms 
$Y_{N,\os{\circ}{t}}\os{\sus}{\lo} {\cal W}_n(Y'_{N,\os{\circ}{t}}) \lo \ol{\cal Y}{}''_{N}$
and 
$Z_{N,\os{\circ}{t}{}'}\os{\sus}{\lo} {\cal W}_n(Z'_{N,\os{\circ}{t}}) \lo \ol{\cal Z}{}'_{N}$ 
are the given immersions. 
Set 
\begin{align*} 
\ol{\cal Y}{}'_{N}
:= \ol{\cal Y}{}''_{N}\times_{{\cal W}_n(s_{\os{\circ}{t}})}
(\ol{\cal Z}{}'_{N}\times_{\ol{\cal W}_n(s'_{\os{\circ}{t}{}'})}{\cal W}_n(s_{\os{\circ}{t}})).
\tag{2.2.14.8}\label{ali:ynwws}
\end{align*}   
Then we have the following commutative diagram 
\begin{equation*} 
\begin{CD} 
{\cal W}_n(Y'_{N,\os{\circ}{t}}) @>>> \ol{\cal Y}{}'_N \\
@VVV @VVV \\
{\cal W}_n(Z'_{N,\os{\circ}{t}{}'}) @>>> \ol{\cal Z}{}'_N, 
\end{CD}
\tag{2.2.14.9}\label{cd:wziy} 
\end{equation*} 
where the right vertical morphism is induced by the second projection. 
By using (\ref{cd:odhtqph}), 
the rest of the proof is almost the same as that of (\ref{theo:ccrw}) (2). 
We leave the detailed proof to the reader. 
\par 
(3):  We leave the detailed proof to the reader. 
\end{proof}

\par  
Let $F\col {\cal W}(Y_{\bul \leq N,\os{\circ}{t}})\lo {\cal W}(Y_{\bul \leq N,\os{\circ}{t}})$ 
be the Frobenius endomorphism. 
By replacing ${\cal W}(s)$ in (\ref{defi:nfc}) with ${\cal W}(\os{\circ}{t})$, 
we obtain the notion of a coherent log $F$-crystal 
on $Y_{\bul \leq N,\os{\circ}{t}}/{\cal W}(\os{\circ}{t})$ and a unit root log $F$-crystal on 
$Y_{\bul \leq N,\os{\circ}{t}}/{\cal W}(\os{\circ}{t})$.  
Let $\ol{E}{}^{\bul \leq N}$ be a flat coherent log $F$-crystal on 
$Y_{\bul \leq N,\os{\circ}{t}}/{\cal W}_n(\os{\circ}{t})$. 
Set $\ol{E}{}^{\bul \leq N}_{n}:=\ol{E}{}^{\bul \leq N}_{{\cal W}_n(Y_{\bul \leq N,\os{\circ}{t}})}$ and $(\ol{E}{}^{\bul \leq N}_n)^{\sig}
:=F^*(\ol{E}{}^{\bul \leq N}_{n})$.  
Then we have the relative Frobenius morphism of 
$\ol{E}{}^{\bul \leq N}_n$: 
the ${\cal O}_{{\cal W}_n(Y_{\bul \leq N,\os{\circ}{t}})}$-linear endomorphism 
$\Phi \col (\ol{E}{}^{\bul \leq N}_n)^{\sig}\lo \ol{E}{}^{\bul \leq N}_n$. 
As in \cite[II (3.0)]{et}, 
we have the $F^*$-linear morphism 
$F \col \ol{E}{}^{\bul \leq N}_n\lo \ol{E}{}^{\bul \leq N}_n$. 
\par 
Now assume that $\ol{E}{}^{\bul \leq N}$ is a unit root log $F$-crystal. 
Then $\Phi \col (\ol{E}{}^{\bul \leq N}_n)^{\sig}\lo \ol{E}{}^{\bul \leq N}_n$ is an isomorphism. 
As in \S\ref{sec:ldrwc}, 
we obtain the graded pro-module 
$\{\ol{E}{}^{\bul \leq N}_{n}\otimes_{{\cal W}_n({\cal O}_{Y_{\bul \leq N,\os{\circ}{t}}})^{\star}}
({\cal W}_n\wt{{\Om}}^{\bul}_{Y_{\bul \leq N,\os{\circ}{t}}})^{\star}\}_{n\in {\mab N}}$ 
over the Raynaud ring ${\mab R}$ of $\kap_t$;  
$\ol{E}{}^{\bul \leq N}_{n}
\otimes_{{\cal W}_n({\cal O}_{Y_{\bul \leq N,\os{\circ}{t}}})^{\star}}
({\cal W}_n\wt{{\Om}}^i_{Y_{\bul \leq N,\os{\circ}{t}}})^{\star}$ is 
a ${\cal W}_n({\cal O}_{Y_t})^{\star}$-module and 
we have the following operators 
satisfying the standard relations in \cite[I (1.1)]{ir} 
and \cite[(1.3.1)]{hdw}: 
\begin{align*} 
& F \col \ol{E}{}^{\bul \leq N}_{n+1}
\otimes_{{\cal W}_{n+1}({\cal O}_{Y_{\bul \leq N,\os{\circ}{t}}})^{\star}}
({\cal W}_{n+1}\wt{{\Om}}_{Y_{\bul \leq N,\os{\circ}{t}}}^{\bul})^{\star} 
\lo 
\ol{E}{}^{\bul \leq N}_{n}\otimes_{{\cal W}_n({\cal O}_{Y_{\bul \leq N,\os{\circ}{t}}})^{\star}}
({\cal W}_n\wt{{\Om}}_{Y_{\bul \leq N,\os{\circ}{t}}}^{\bul})^{\star}, \\ 
& V\col \ol{E}{}^{\bul \leq N}_{n}
\otimes_{{\cal W}_n({\cal O}_{Y_{\bul \leq N,\os{\circ}{t}}})^{\star}} 
({\cal W}_n\wt{{\Om}}_{Y_{\bul \leq N,\os{\circ}{t}}}^{\bul})^{\star} \lo 
\ol{E}{}^{\bul \leq N}_{n+1}
\otimes_{{\cal W}_{n+1}({\cal O}_{Y_{\bul \leq N,\os{\circ}{t}}})^{\star}}
({\cal W}_{n+1}\wt{{\Om}}_{Y_{\bul \leq N,\os{\circ}{t}}}^{\bul})^{\star}, \\
& \nabla\col \ol{E}{}^{\bul \leq N}_{n}
\otimes_{{\cal W}_n({\cal O}_{Y_{\bul \leq N,\os{\circ}{t}}})^{\star}}
({\cal W}_n\wt{{\Om}}_{Y_{\bul \leq N,\os{\circ}{t}}}^{\bul})^{\star} 
\lo  \ol{E}{}^{\bul \leq N}_{n}
\otimes_{{\cal W}_n({\cal O}_{Y_{\bul \leq N,\os{\circ}{t}}})^{\star}}
({\cal W}_n\wt{{\Om}}_{Y_{\bul \leq N,\os{\circ}{t}}}^{\bul +1})^{\star}, \\
& {\bf p} \col 
\ol{E}{}^{\bul \leq N}_{n}\otimes_{{\cal W}_n({\cal O}_{Y_{\bul \leq N,\os{\circ}{t}}})^{\star}}
({\cal W}_n\wt{{\Om}}_{Y_{\bul \leq N,\os{\circ}{t}}}^{\bul})^{\star} \lo 
\ol{E}{}^{\bul \leq N}_{n+1}
\otimes_{{\cal W}_{n+1}({\cal O}_{Y_{\bul \leq N,\os{\circ}{t}}})^{\star}}
({\cal W}_{n+1}\wt{{\Om}}_{Y_{\bul \leq N,\os{\circ}{t}}}^{\bul})^{\star}, \\
& R \col \ol{E}{}^{\bul \leq N}_{n+1}
\otimes_{{\cal W}_{n+1}({\cal O}_{Y_{\bul \leq N,\os{\circ}{t}}})^{\star}}
({\cal W}_{n+1}\wt{{\Om}}_{Y_{\bul \leq N,\os{\circ}{t}}}^{\bul})^{\star}
\lo \ol{E}{}^{\bul \leq N}_{n}
\otimes_{{\cal W}_n({\cal O}_{Y_{\bul \leq N,\os{\circ}{t}}})^{\star}}
({\cal W}_n\wt{{\Om}}_{Y_{\bul \leq N,\os{\circ}{t}}}^{\bul})^{\star}. 
\end{align*} 
As in \cite[(3.2.3), (3.3.1)]{et}, the following holds$:$ 

\begin{prop}\label{prop:fnet} 
\begin{align*} 
&{\rm Fil}^n(\ol{E}{}^{\bul \leq N}_{n+r}
\otimes_{{\cal W}_{n+r}({\cal O}_{Y_{\bul \leq N,\os{\circ}{t}}})^{\star}}
({\cal W}_{n+r}\wt{{\Om}}^i_{Y_{\bul \leq N,\os{\circ}{t}}})^{\star})
\tag{2.2.15.1}\label{eqn:foilu} \\
&=V^n(\ol{E}{}^{\bul \leq N}_{r}
\otimes_{{\cal W}_r({\cal O}_{Y_{\bul \leq N,\os{\circ}{t}}})^{\star}}
({\cal W}_r\wt{{\Om}}^i_{Y_{\bul \leq N,\os{\circ}{t}}})^{\star})+
\nabla V^n(\ol{E}{}^{\bul \leq N}_{r}
\otimes_{{\cal W}_r({\cal O}_{Y_{\bul \leq N,\os{\circ}{t}}})^{\star}}
({\cal W}_r\wt{{\Om}}^{i-1}_{Y_{\bul \leq N,\os{\circ}{t}}})^{\star}).
\end{align*} 
\end{prop} 
\begin{proof} 
The proof is the same as that of \cite[II (3.2.3)]{et}. 
\end{proof} 

\begin{prop}\label{prop:tff}
The following formula holds$:$ 
\begin{align*} 
&F^r({\rm Fil}^n(\ol{E}{}^{\bul \leq N}_{n+r}
\otimes_{{\cal W}_{n+r}({\cal O}_{Y_{\bul \leq N,\os{\circ}{t}}})^{\star}}
({\cal W}_{n+r}\wt{{\Om}}^i_{Y_{\bul \leq N,\os{\circ}{t}}})^{\star})) 
\tag{2.2.16.1}\label{eqn:wlsfou}\\
&= F^r\nabla V^n(\ol{E}{}^{\bul \leq N}_{n+r}
\otimes_{{\cal W}_{n+r}({\cal O}_{Y_{\bul \leq N,\os{\circ}{t}}})^{\star}}
({\cal W}_{n+r}\wt{{\Om}}^{i-1}_{Y_{\bul \leq N,\os{\circ}{t}}})^{\star}). 
\end{align*} 
Consequently the morphism 
$$F^r \col \ol{E}{}^{\bul \leq N}_{n+r}
\otimes_{{\cal W}_{n+r}({\cal O}_{Y_{\bul \leq N,\os{\circ}{t}}})^{\star}}
({\cal W}_{n+r}\wt{{\Om}}^i_{Y_{\bul \leq N,\os{\circ}{t}}})^{\star}\lo 
F^r_*\{\ol{E}{}^{\bul \leq N}_{n}
\otimes_{{\cal W}_n({\cal O}_{Y_{\bul \leq N,\os{\circ}{t}}})^{\star}}
({\cal W}_n\wt{{\Om}}^i_{Y_{\bul \leq N,\os{\circ}{t}}})^{\star}\}$$  
induces the following morphism 
\begin{align*} 
\check{F}{}^r \col & 
\ol{E}{}^{\bul \leq N}_{n}
\otimes_{{\cal W}_n({\cal O}_{Y_{\bul \leq N,\os{\circ}{t}}})^{\star}}
({\cal W}_n\wt{{\Om}}^i_{Y_{\bul \leq N,\os{\circ}{t}}})^{\star}
\lo \tag{2.2.16.2}\label{eqn:wlsfwu}\\
& 
F^r_*\{F^r(\ol{E}{}^{\bul \leq N}_{n+r}
\otimes_{{\cal W}_{n+r}({\cal O}_{Y_{\bul \leq N,\os{\circ}{t}}})^{\star}}
({\cal W}_{n+r}\wt{{\Om}}^i_{Y_{\bul \leq N,\os{\circ}{t}}})^{\star})/ 
F^r\nabla V^n(\ol{E}{}^{\bul \leq N}_{r}
\otimes_{{\cal W}_r({\cal O}_{Y_{\bul \leq N,\os{\circ}{t}}})^{\star}}
({\cal W}_n\wt{{\Om}}^{i-1}_{Y_{\bul \leq N,\os{\circ}{t}}})^{\star})\}. 
\end{align*} 
For the case $r=n$, $\check{F}{}^n$ 
induces the surjection 
\begin{align*}
\check{F}{}^n \col 
\ol{E}{}^{\bul \leq N}_{n}
\otimes_{{\cal W}_n({\cal O}_{Y_{\bul \leq N,\os{\circ}{t}}})^{\star}}
({\cal W}_n\wt{{\Om}}^i_{Y_{\bul \leq N,\os{\circ}{t}}})^{\star} \lo 
F^n_*{\cal H}^i(\ol{E}{}^{\bul \leq N}_{n}\otimes_{{\cal W}_r({\cal O}_{Y_{\bul \leq N,\os{\circ}{t}}})^{\star}}
({\cal W}_n\wt{{\Om}}^{\bul}_{Y_{\bul \leq N,\os{\circ}{t}}})^{\star}). 
\tag{2.2.16.3}\label{eqn:wlyfyu}
\end{align*}
\end{prop}
\begin{proof} 
(\ref{prop:tff}) follows from (\ref{prop:fnet}).  
\end{proof} 

\par

\begin{theo}\label{theo:ctin}  
The morphism 
\begin{align*} 
V^r\col F^r_*\{\ol{E}{}^{\bul \leq N}_{n}
\otimes_{{\cal W}_n({\cal O}_{Y_{\bul \leq N,\os{\circ}{t}}})^{\star}}
({\cal W}_n\wt{{\Om}}^i_{Y_{\bul \leq N,\os{\circ}{t}}})^{\star}\} \lo 
\ol{E}{}^{\bul \leq N}_{n+r}
\otimes_{{\cal W}_{n+r}({\cal O}_{Y_{\bul \leq N,\os{\circ}{t}}})^{\star}}
({\cal W}_{n+r}\wt{{\Om}}^i_{Y_{\bul \leq N,\os{\circ}{t}}})^{\star} 
\end{align*} 
induces the following morphism 
\begin{align*} 
\check{V}{}^r \col & 
F^r_*\{\ol{E}{}^{\bul \leq N}_{n}
\otimes_{{\cal W}_n({\cal O}_{Y_{\bul \leq N,\os{\circ}{t}}})^{\star}}
({\cal W}_n\wt{{\Om}}^i_{Y_{\bul \leq N,\os{\circ}{t}}})^{\star}/
F^r\nabla V^n(\ol{E}{}^{\bul \leq N}_{n}
\otimes_{{\cal W}_n({\cal O}_{Y_{\bul \leq N,\os{\circ}{t}}})^{\star}}
({\cal W}_n\wt{{\Om}}^{i-1}_{Y_{\bul \leq N,\os{\circ}{t}}})^{\star})\} 
\tag{2.2.17.1}\label{eqn:wlsiwu}\\
& \lo 
\ol{E}{}^{\bul \leq N}_{n+r}
\otimes_{{\cal W}_{n+r}({\cal O}_{Y_{\bul \leq N,\os{\circ}{t}}})^{\star}}
({\cal W}_{n+r}\wt{{\Om}}^i_{Y_{\bul \leq N,\os{\circ}{t}}})^{\star}. 
\end{align*}
There exists a generalized Cartier isomorphism 
\begin{align*} 
\check{C}{}^r\col & 
F^r_*\{F^r(\ol{E}{}^{\bul \leq N}_{n+r}
\otimes_{{\cal W}_{n+r}({\cal O}_{Y_{\bul \leq N,\os{\circ}{t}}})^{\star}}
({\cal W}_{n+r}\wt{{\Om}}^i_{Y_{\bul \leq N,\os{\circ}{t}}})^{\star})/ \tag{2.2.17.2}\label{eqn:wlsfyu}\\
&F^r\nabla V^n( 
\ol{E}{}^{\bul \leq N}_{r}\otimes_{{\cal W}_r({\cal O}_{Y_{\bul \leq N,\os{\circ}{t}}})^{\star}}
({\cal W}_r\wt{{\Om}}^{i-1}_{Y_{\bul \leq N,\os{\circ}{t}}})^{\star})\} \\
& \os{\sim}{\lo} 
\ol{E}{}^{\bul \leq N}_{n}
\otimes_{{\cal W}_n({\cal O}_{Y_{\bul \leq N,\os{\circ}{t}}})^{\star}}
({\cal W}_n\wt{{\Om}}^i_{Y_{\bul \leq N,\os{\circ}{t}}})^{\star} 
\end{align*} 
which is the inverse of $\check{F}{}^r$. 
The morphism  $\check{C}{}^r$ satisfies a relation 
${\bf p}^r\circ \check{C}{}^r=\check{V}{}^r$. 
In particular, there exist the following isomorphisms 
\begin{align*} 
\check{F}{}^n \col 
\ol{E}{}^{\bul \leq N}_{n}
\otimes_{{\cal W}_n({\cal O}_{Y_{\bul \leq N,\os{\circ}{t}}})^{\star}}
({\cal W}_n\wt{{\Om}}^i_{Y_{\bul \leq N,\os{\circ}{t}}})^{\star}  
\os{\sim}{\lo} 
F^n_*({\cal H}^i(\ol{E}{}^{\bul \leq N}_{n}
\otimes_{{\cal W}_n({\cal O}_{Y_{\bul \leq N,\os{\circ}{t}}})^{\star}}
({\cal W}_n\wt{{\Om}}^{\bul}_{Y_{\bul \leq N,\os{\circ}{t}}})^{\star})) 
\tag{2.2.17.3}\label{eqn:wltfyu}
\end{align*}  
and 
\begin{align*} 
\check{C}{}^n \col 
F^n_*({\cal H}^i(\ol{E}{}^{\bul \leq N}_{n}
\otimes_{{\cal W}_n({\cal O}_{Y_{\bul \leq N,\os{\circ}{t}}})^{\star}}
({\cal W}_n\wt{{\Om}}^{\bul}_{Y_{\bul \leq N,\os{\circ}{t}}})^{\star})) 
\os{\sim}{\lo} 
\ol{E}{}^{\bul \leq N}_{n}\otimes_{{\cal W}_n({\cal O}_{Y_{\bul \leq N,\os{\circ}{t}}})^{\star}}
({\cal W}_n\wt{{\Om}}^i_{Y_{\bul \leq N,\os{\circ}{t}}})^{\star},   
\tag{2.2.17.4}\label{eqn:wlctyu}
\end{align*}  
which are the inverse of another. 
\end{theo} 


\par  
Then, by (\ref{theo:ccttrw}) and (\ref{eqn:wltfyu}), 
\begin{equation*} 
\ol{E}{}^{\bul \leq N}_{n}\otimes_{{\cal W}_n({\cal O}_{Y_{\bul \leq N,\os{\circ}{t}}})^{\star}}
({\cal W}_n\wt{{\Om}}^i_{Y_{\bul \leq N,\os{\circ}{t}}})^{\star}
=
\wt{R}^iu_{Y_{\bul \leq N,\os{\circ}{t}}/{\cal W}_n(\os{\circ}{t})*}(\ol{E}{}^{\bul \leq N}).   
\tag{2.2.17.5}\label{eqn:wltyu}
\end{equation*} 
In a standard way,  we can define a boundary morphism 
``$p^{-n}\nabla $'': 
\begin{align*} 
\nabla \col 
\wt{R}^iu_{Y_{\bul \leq N,\os{\circ}{t}}/{\cal W}_n(\os{\circ}{t})*}(\ol{E}{}^{\bul \leq N})
\lo 
\wt{R}^{i+1}u_{Y_{\bul \leq N,\os{\circ}{t}}/{\cal W}_n(\os{\circ}{t})*}(\ol{E}{}^{\bul \leq N}). 
\tag{2.2.17.6}\label{ali:ytwns}
\end{align*}

\par 
Assume that there exists an immersion 
$\iota \col Y_{\bul \leq N,\os{\circ}{t}} \os{\sus}{\lo} \ol{\cal Q}_{{\bul \leq N},n}$ 
into a log smooth scheme over 
$\ol{{\cal W}_n(s_{\os{\circ}{t}})}$. 
Let $\ol{\mathfrak E}_{{\bul \leq N},n}$ 
be the log PD-envelope of the immersion 
$Y_{\bul \leq N,\os{\circ}{t}} \os{\subset}{\lo} \ol{\cal Q}_{{\bul \leq N},n}$ 
over $({\cal W}_n(\os{\circ}{t}),p{\cal W}_n,[~])$. 
Set ${\mathfrak E}_{{\bul \leq N},n}
:=\ol{\mathfrak E}_{{\bul \leq N},n}
\times_{{\mathfrak D}(\ol{{\cal W}_n(s_{\os{\circ}{t}}))}}
{\cal W}_n(s_{\os{\circ}{t}})$.  
Let $(\ol{\cal E}{}^{\bul \leq N},{\nabla})$ be the corresponding 
${\cal O}_{\ol{\mathfrak E}_{{\bul \leq N},n}}$-modules with integrable 
connection to $\ol{E}{}^{\bul \leq N}$. 
Set $({\cal E}^{\bul \leq N},\nabla):=(\ol{\cal E}{}^{\bul \leq N},{\nabla})
\otimes_{{\cal O}_{\ol{\mathfrak E}_{{\bul \leq N},n}}}
{\cal O}_{{\mathfrak E}_{{\bul \leq N},n}}$. 
Then, by (\ref{eqn:wltyu}), 
we have the following equality for $i\in {\mab N}$:
\begin{align*} 
\ol{E}{}^{\bul \leq N}_{n}\otimes_{{\cal W}_n({\cal O}_{Y_{\bul \leq N,\os{\circ}{t}}})^{\star}}
({\cal W}_n\wt{{\Om}}^i_{Y_{\bul \leq N,\os{\circ}{t}}})^{\star}&=
{\cal H}^i({\cal E}^{\bul \leq N}{\otimes}_{{\cal O}_{{\cal Q}_{n}}}
\Om_{{\cal Q}_{n}/{\cal W}_n(\os{\circ}{t})}^{\bul})
\tag{2.2.17.7}\label{ali:wlth}\\
&=
{\cal H}^i(
\ol{E}{}^{\bul \leq N}_{n}{\otimes}_{{\cal W}_n({\cal O}_{Y_{\bul \leq N,\os{\circ}{t}}})}
{\cal W}_n\wt{\Om}_{Y_{\bul \leq N,\os{\circ}{t}}}^{\bul}) 
\quad (i\in {\mab N}). 
\end{align*}  

\begin{theo}\label{theo:edctisw}
Let the notations be as in {\rm (\ref{theo:ccttrw})} and 
the proof of {\rm (\ref{theo:ccttrw})}. 
Assume that $\ol{E}{}^{\bul \leq N}$ is a flat coherent unit root log $F$-crystal.  
Then the isomorphism in {\rm (\ref{theo:ccttrw}) (1)} 
for the case $\star$=nothing is equal to 
the following isomorphism 
\begin{align*} 
&R\pi_{{\rm zar}*}({\cal E}^{\bul \leq N,\bul}
{\otimes}_{{\cal O}_{{\cal Q}_{\bul \leq N,\bul}}}
\Om^{\bul}_{{\cal Q}_{\bul \leq N,\bul}/{\cal W}_n(\os{\circ}{t})})  
\os{\sim}{\lo} R\pi_{{\rm zar}*}(
{\cal H}^{\bul}({\cal E}^{\bul \leq N,\bul}
{\otimes}_{{\cal O}_{{\cal Q}_{\bul \leq N,\bul}}}
\Om^{*}_{{\cal Q}_{\bul \leq N,\bul}/{\cal W}_n(\os{\circ}{t})})).  
\tag{2.2.18.1}\label{ali:entroq}
\end{align*} 
\end{theo} 
\begin{proof} 
This follows from (\ref{ali:enotq}), (\ref{eqn:ywnttnny}) 
and (\ref{ali:wlth}). 
\end{proof}

\begin{rema}\label{rema:wqtn}
Let the notations and the assumptions be in {\rm (\ref{theo:edctisw})}. 
Then, as in (\ref{rema:wqn}),
we can consider 
$\wt{R}{}^iu_{Y_{\bul \leq N,\os{\circ}{t}}/{\cal W}_n(\os{\circ}{t})*}(\ol{E}{}^{\bul \leq N})$   
as the definition of the sheaf 
``$\ol{E}{}^{\bul \leq N}_{n}
\otimes_{{\cal W}_n({\cal O}_{Y_{\bul \leq N,\os{\circ}{t}}})}
{\cal W}_n\wt{{\Om}}^i_{Y_{\bul \leq N,\os{\circ}{t}}}$''  
and $\{\wt{R}{}^iu_{Y_{\bul \leq N,\os{\circ}{t}}/{\cal W}_n(\os{\circ}{t})*}(\ol{E}{}^{\bul \leq N})\}_{i\in {\mab N}}$ 
becomes a complex 
$\{\wt{R}{}^{\bul}u_{Y_{\bul \leq N,\os{\circ}{t}}/{\cal W}_n(\os{\circ}{t})*}
(\ol{E}{}^{\bul \leq N})\}$. 
Moreover, one can prove that 
$\{\wt{R}{}^{\bul}u_{Y_{\bul \leq N,\os{\circ}{t}}/{\cal W}_n(\os{\circ}{t})*}
(\ol{E}{}^{\bul \leq N})\}_{n\in {\mab N}}$ 
has operators $F$, $V$, ${\bf p}$ and $R$ in a standard way. 
(\ref{theo:ccttrw}) and (\ref{theo:edctisw}) imply 
that there exists a canonical isomorphism 
\begin{align*} 
\ol{E}{}^{\bul \leq N}_{n}
\otimes_{{\cal W}_n({\cal O}_{Y_{\bul \leq N,\os{\circ}{t}}})}
{\cal W}_n\wt{\Om}{}^{\bul}_{Y_{\bul \leq N,\os{\circ}{t}}}
\os{\sim}{\lo} 
\wt{R}{}^{\bul}u_{Y_{\bul \leq N,\os{\circ}{t}}/{\cal W}_n(\os{\circ}{t})*}
(\ol{E}{}^{\bul \leq N})
\tag{2.2.19.1}\label{ali:ytnw}
\end{align*} 
in ${\rm D}^+(g^{-1}({\cal O}_{{\cal W}_n(\os{\circ}{t})}))$.
\end{rema}

\section{Zariskian $p$-adic filtered Steenbrink complexes  
via log de Rham-Witt complexes}\label{sec:flgdw}
Let $\kap$ be a perfect field of 
characteristic $p>0$. 
Let $s$ be the log point 
$({\rm Spec}(\kap), {\mab N}\oplus \kap^*\lo \kap)$. 
Let ${\cal W}$ (resp.~${\cal W}_n$) be the Witt ring of $\kap$ 
(resp.~the Witt ring of $\kap$ of length $n \in {\mab Z}_{>0}$). 
Let $t\lo s$ be a morphism from a fine log scheme
whose underlying scheme is the spectrum of 
a perfect field of characteristic $p>0$. 
Set $\kap_t:=\Gam(t,{\cal O}_t)$. 
Let ${\cal W}(t)$ be the formal canonical lift of $t$ 
over ${\rm Spf}({\cal W}(\kap_t))$.  
Let $N$ be a nonnegative integer or $\infty$. 
Let $X_{\bul \leq N}$ be 
an $N$-truncated simplicial SNCL scheme over $s$.  
Set $X_{\bul \leq N,\os{\circ}{t}}:=X_{\bul \leq N}\times_{\os{\circ}{s}}\os{\circ}{t}$  
and $X_{\bul \leq N,t}:=X_{\bul \leq N}\times_st$.  
Let 
$f\col X_{\bul \leq N,\os{\circ}{t}}\lo {\cal W}_n(s_{\os{\circ}{t}})$ 
and $f_{{\cal W}_n(t)}\col X_{\bul \leq N,t}\lo {\cal W}_n(t)$
be the structural morphisms.  
Let $E$ be a flat coherent crystal of 
${\cal O}_{\os{\circ}{X}_{\bul \leq N,t}/{\cal W}_n(\os{\circ}{t})}$-modules.   
In this section we construct a filtered complex 
$({\cal W}_nA_{X_{\bul \leq N,\os{\circ}{t}}}(E^{\bul \leq N}),P)$
of $f^{-1}({\cal W}_n(\kap_t))$-modules 
for a flat coherent log crystal $E^{\bul \leq N}$ of 
${\cal O}_{\os{\circ}{X}_{\bul \leq N,t}/{\cal W}_n(\os{\circ}{t})}$-modules. 
In the case where 
$E^{\bul \leq N}=
{\cal O}_{\os{\circ}{X}_{\bul \leq N,t}/{\cal W}_n(\os{\circ}{t})}$, 
we denote the filtered complex 
$({\cal W}_nA_{X_{\bul \leq N,\os{\circ}{t}}}
({\cal O}_{\os{\circ}{X}_{\bul \leq N,t}/
{\cal W}_n(\os{\circ}{t})}),P)$ by 
$({\cal W}_nA_{X_{\bul \leq N},\os{\circ}{t}},P)$. 
This is the $N$-truncated cosimplicial 
filtered version of the Hyodo-Mokrane-Steenbrink complex  
constructed by Mokrane (\cite[3.8]{msemi}) and corrected by the author of this book.  
Because there are numerous incomplete and mistaken parts in 
\cite{msemi}, see \cite{ndw}, (\ref{prop:plz}),  
(\ref{rema:qok}) (2), (3), (\ref{rema:po}) (2) below and (\ref{prop:prm}) below.   
In this section, if $N\not= \infty$ and if $X_{\bul \leq N}$ 
has an affine $N$-truncated simplicial open covering, 
then we prove that 
there exists a canonical isomorphism 
\begin{equation*} 
(A_{\rm zar}(X_{\bul \leq N,\os{\circ}{t}}/{\cal W}_n(s_{\os{\circ}{t}}),E^{\bul \leq N}),P)
\os{\sim}{\lo} 
({\cal W}_nA_{X_{\bul \leq N,\os{\circ}{t}}}(E^{\bul \leq N}),P). 
\tag{2.3.0.1}\label{eqn:canzad}
\end{equation*}  
Here note that ${\cal W}_n(s_{\os{\circ}{t}})^{\nat}={\cal W}_n(s_{\os{\circ}{t}})$. 
\par 
Let $X$ be an SNCL scheme over $s$. 
First we give a new definition 
of $P_k{\cal W}_n\wt{{\Om}}{}^i_{X_{\os{\circ}{t}}}$ $(i\in {\mab N})$ in \cite[(3.5)]{msemi} 
as follows. 
\par 
Set 
\begin{align*} 
P_k{\cal W}_n\wt{{\Om}}{}^i_{X_{\os{\circ}{t}}}:=
P_k\wt{R}{}^iu_{X_{\os{\circ}{t}}/{\cal W}_n(\os{\circ}{t})*}
({\cal O}_{X_{\os{\circ}{t}}/{\cal W}_n(\os{\circ}{t})}):=
{\cal H}^i(P_k\wt{R}u_{X_{\os{\circ}{t}}/{\cal W}_n(\os{\circ}{t})*}
({\cal O}_{X_{\os{\circ}{t}}/{\cal W}_n(\os{\circ}{t})})). 
\tag{2.3.0.2}\label{ali:omxt} 
\end{align*} 
(Recall the definition (\ref{theo:fcp}).) 
Then we have a natural morphism 
\begin{align*} 
P_k{\cal W}_n\wt{{\Om}}{}^i_{X_{\os{\circ}{t}}}\lo 
{\cal W}_n\wt{{\Om}}{}^i_{X_{\os{\circ}{t}}}. 
\tag{2.3.0.3}\label{ali:omnxt} 
\end{align*} 

Moreover $P_k{\cal W}_n\wt{{\Om}}{}^i_{X_{\os{\circ}{t}}}$ is naturally a sheaf of 
${\cal W}_n({\cal O}_{X_t})= 
\wt{R}{}^0u_{X_{\os{\circ}{t}}/{\cal W}_n(\os{\circ}{t})*}
({\cal O}_{X_t/{\cal W}_n(\os{\circ}{t})})$-modules. 
(This equality follows from (\ref{eqn:wltyu}).)

\begin{prop}\label{prop:fupo}  
$(1)$ The sheaf $P_k{\cal W}_n\wt{\Om}^i_{X_{\os{\circ}{t}}}$ is contravariantly functorial 
with respect to the following commutative diagram  
\begin{equation*} 
\begin{CD} 
X@>>> X'\\
@VVV @VVV \\
s_{\os{\circ}{t}}@>>> s'_{\os{\circ}{t}{}'} 
\end{CD} 
\tag{2.3.1.1}\label{cd:xxp}
\end{equation*} 
and {\rm (\ref{cd:tts})},  
where $X'$, $s'$ and $t'$ are similar objects to $X$, $s$ and $t$,   
respectively. 
\par 
$(2)$ The differential $d\col {\cal W}_n\wt{{\Om}}{}^i_{X_{\os{\circ}{t}}}
\lo {\cal W}_n\wt{{\Om}}{}^{i+1}_{X_{\os{\circ}{t}}}$
induces a morphism 
$d\col P_k{\cal W}_n\wt{{\Om}}{}^i_{X_{\os{\circ}{t}}}
\lo P_k{\cal W}_n\wt{{\Om}}{}^{i+1}_{X_{\os{\circ}{t}}}$ of complexes. 
\end{prop}
\begin{proof} 
(1): (1) immediately follows from (\ref{prop:cttu}) (1). 
\par 
(2): Assume that there exists an immersion $X_{\os{\circ}{t}} \os{\sus}{\lo} {\cal P}_{2n}$ 
into a log smooth integral scheme over ${\cal W}_{2n}(s_{\os{\circ}{t}})$. 
Set 
${\cal P}_n:=
{\cal P}_{2n}\times_{{\cal W}_{2n}(s_{\os{\circ}{t}})}{\cal W}_n(s_{\os{\circ}{t}})$. 
Let ${\mathfrak D}_m$ $(m=2n,n)$ be the log PD-envelope of 
the immersion $X_{\os{\circ}{t}} \os{\sus}{\lo} {\cal P}_m$ 
over $({\cal W}_m(s_{\os{\circ}{t}}),p{\cal W}_m,[~])$.  
Then 
\begin{align*} 
P_k{\cal W}_m\wt{{\Om}}{}^i_{X_{\os{\circ}{t}}}=
{\cal H}^i(P_k({\cal O}_{{\mathfrak D}_m}
\otimes_{{\cal O}_{{\cal P}^{\rm ex}_m}}
{\Om}^{\bul}_{{\cal P}^{\rm ex}_m/{\cal W}_m(\os{\circ}{t})})). 
\tag{2.3.1.2}\label{ali:omwxt} 
\end{align*}  
((\ref{ali:omwxt}) is a generalized description 
of $P_k{\cal W}_m\wt{{\Om}}{}^i_{X_{\os{\circ}{t}}}$ in \cite[(3.5)]{msemi}.) 
By (\ref{coro:flt}) we have the following exact sequence 
\begin{align*} 
0& \lo P_k({\cal O}_{{\mathfrak D}_{n}}
\otimes_{{\cal O}_{{\cal P}^{\rm ex}_n}}
\Om^{\bul}_{{\cal P}^{\rm ex}_n/{\cal W}_n(\os{\circ}{t})})
\os{p^n}{\lo} 
P_k({\cal O}_{{\mathfrak D}_{2n}}
\otimes_{{\cal O}_{{\cal P}^{\rm ex}_{2n}}}
\Om^{\bul}_{{\cal P}^{\rm ex}_{2n}/{\cal W}_{2n}(\os{\circ}{t})})
\tag{2.3.1.3}\label{ali:bopnt}\\
&\lo P_k({\cal O}_{{\mathfrak D}_{n}}
\otimes_{{\cal O}_{{\cal P}^{\rm ex}_n}}
\Om^{\bul}_{{\cal P}^{\rm ex}_n/{\cal W}_n(\os{\circ}{t})})\lo 0. 
\end{align*} 
Because the differential 
\begin{align*} 
d\col {\cal O}_{{\mathfrak D}_{2n}}
\otimes_{{\cal O}_{{\cal P}^{\rm ex}_{2n}}}
\Om^i_{{\cal P}^{\rm ex}_{2n}/{\cal W}_{2n}(\os{\circ}{t})}
\lo 
{\cal O}_{{\mathfrak D}_{2n}}
\otimes_{{\cal O}_{{\cal P}^{\rm ex}_{2n}}}
\Om^{i+1}_{{\cal P}^{\rm ex}_{2n}/{\cal W}_{2n}(\os{\circ}{t})}\quad (i\in {\mab N})
\tag{2.3.1.4}\label{ali:bopdt}
\end{align*}
preserves $P$, 
we have the following boundary morphism of the exact sequence 
associated to (\ref{ali:bopnt}): 
\begin{align*} 
d:=``p^{-n}d''\col P_k{\cal W}_n\wt{{\Om}}{}^i_{X_{\os{\circ}{t}}}&=
{\cal H}^i(P_k({\cal O}_{{\mathfrak D}_{n}}
\otimes_{{\cal O}_{{\cal P}^{\rm ex}_n}}
\Om^{\bul}_{{\cal P}^{\rm ex}_n/{\cal W}_n(\os{\circ}{t})}))
\tag{2.3.1.5}\label{ali:bopant}\\
& \lo 
{\cal H}^{i+1}(P_k({\cal O}_{{\mathfrak D}_{n}}
\otimes_{{\cal O}_{{\cal P}^{\rm ex}_n}}
\Om^{\bul}_{{\cal P}^{\rm ex}_n/{\cal W}_n(\os{\circ}{t})}))
=P_k{\cal W}_n\wt{{\Om}}{}^{i+1}_{X_{\os{\circ}{t}}}. 
\end{align*} 
This proves (2). 
\end{proof}

\begin{rema}\label{rema:po}
(1) In \cite[(3.5)]{msemi}
the sheaf ${\cal W}_n\wt{\Om}^i_{X_{\os{\circ}{t}}}$ 
(=$P_kW_n\wt{\om}^i_{X_{\os{\circ}{t}}}$ in [loc.~cit.])  
was defined by only an admissible lift of $X_{\os{\circ}{t}}$ over a polynomial ring 
in one variable over ${\cal W}_n(\kap_t)$. 
We allow an immersion from $X_{\os{\circ}{t}}$ into 
a log smooth integral scheme ${\cal P}$ 
over ${\cal W}_n(s_{\os{\circ}{t}})$ 
for the definition of ${\cal W}_n\wt{\Om}^i_{X_{\os{\circ}{t}}}$ 
because the immersion is stable under the product of two immersions 
$X_{\os{\circ}{t}}\os{\sus}{\lo} {\cal P}$'s. 
\par 
(2) As pointed out in (\ref{rema:qok}) (3), 
I am not sure that the proof of 
\cite[Lemme 3.4]{msemi} needed for the well-definedness of 
$P_k{\cal W}_n\wt{{\Om}}{}^i_{X}$ is perfect. 
\end{rema}

\par 
Let the notations be as in the proof of (\ref{prop:fupo}) (2). 
By the definition of (\ref{ali:bopant}), the following diagram is commutative: 
\begin{equation*} 
\begin{CD}
P_k{\cal W}_n\wt{{\Om}}{}^i_{X_{\os{\circ}{t}}}@>{d}>>
P_k{\cal W}_n\wt{{\Om}}{}^{i+1}_{X_{\os{\circ}{t}}}\\
@VVV @VVV \\
P_{k+1}{\cal W}_n\wt{{\Om}}{}^i_{X_{\os{\circ}{t}}}@>{d}>>
P_{k+1}{\cal W}_n\wt{{\Om}}{}^{i+1}_{X_{\os{\circ}{t}}}. 
\end{CD}
\end{equation*} 
We have not yet claimed that the vertical morphisms above are injective. 
\par 
Though the following Poincar\'{e} residue isomorphism is different from 
the Poincar\'{e} residue isomorphism in \cite{msemi}, 
the following proposition is essentially the same as  
\cite[Corollaire 3.7]{msemi}. 

\begin{prop}[{\bf cf.~\cite[Lemme 3.3, Lemme 3.4]{msemi}}]
\label{prop:prm} 
Let $X$ be an SNCL scheme over $s$. Then 
there exists 
the following Poincar\'{e} residue morphism 
\begin{align*} 
{\rm Res} \col & 
P_k{\cal W}_n\wt{{\Om}}{}^{\bul}_{X_{\os{\circ}{t}}}
\lo a^{(k-1)}_{*}
({\cal W}_n\Om^{\bul}_{\os{\circ}{X}{}^{(k-1)}_t}
\otimes_{\mab Z}\vp^{(k-1)}_{\rm zar}(\os{\circ}{X}_t/\os{\circ}{t}))[-k] 
\quad (k\in {\mab Z}_{\geq 1}) \tag{2.3.3.1}\label{eqn:retsd}
\end{align*} 
fitting into the following exact sequence 
\begin{align*} 
0 \lo & 
P_{k-1}{\cal W}_n\wt{{\Om}}{}^{\bul}_{X_{\os{\circ}{t}}}
\lo P_k{\cal W}_n\wt{{\Om}}{}^{\bul}_{X_{\os{\circ}{t}}}
\os{\rm Res}{\lo}  
a^{(k-1)}_{*}
({\cal W}_n\Om^{\bul}_{\os{\circ}{X}{}^{(k-1)}_t}
\otimes_{\mab Z}\vp^{(k-1)}_{\rm zar}(\os{\circ}{X}_t/\os{\circ}{t}))[-k] 
\lo 0. \tag{2.3.3.2}\label{eqn:ptwa}
\end{align*} 
\end{prop} 
\begin{proof}  
Take an open log subscheme of $X_{\os{\circ}{t}}$ such that 
there exists an immersion from the open log subscheme 
into a log smooth integral scheme over ${\cal W}_n(s_{\os{\circ}{t}})$. 
By abuse of notation, we denote this open log subscheme by $X_{\os{\circ}{t}}$. 
Let $X_{\os{\circ}{t}} \os{\sus}{\lo} {\cal P}_n$ be an immersion 
into a log smooth integral scheme over ${\cal W}_n(s_{\os{\circ}{t}})$. 
Let ${\mathfrak D}_n$ be the log PD-envelope of 
this immersion over $({\cal W}_n(s_{\os{\circ}{t}}),p{\cal W}_n,[~])$. 
Let $\os{\circ}{\mathfrak D}{}^{(k-1)}_n$ be the PD-envelope of 
the immersion $\os{\circ}{X}{}^{(k-1)}_t \os{\sus}{\lo} \os{\circ}{\cal P}{}^{(k-1)}_n$ 
over $({\cal W}_n(\os{\circ}{t}),p{\cal W}_n,[~])$.  
Let $c^{(k-1)}\col \os{\circ}{\mathfrak D}{}^{(k-1)}_n\lo \os{\circ}{\mathfrak D}_n$ 
be the natural morphism. 
We define the morphism (\ref{eqn:retsd}) as the following morphism: 
\begin{align*} 
&{\cal H}^i(P_k({\cal O}_{{\mathfrak D}_{n}}
\otimes_{{\cal O}_{{\cal P}^{\rm ex}_n}}
{\Om}^{\bul}_{{\cal P}^{\rm ex}_n/{\cal W}_n(\os{\circ}{t})}))
\tag{2.3.3.3}\label{ali:poe}\\
& \os{{\cal H}^i({\rm Res})}{\lo}  
c^{(k-1)}_*
{\cal H}^i({\cal O}_{\os{\circ}{\mathfrak D}{}^{(k-1)}_{n}}
\otimes_{{\cal O}_{\os{\circ}{\cal P}{}^{{\rm ex},(k-1)}_n}}
\Om^{\bul}_{\os{\circ}{\cal P}{}^{{\rm ex},(k-1)}_n/{\cal W}_n}
\otimes_{\mab Z}
\vp^{(k-1)}_{\rm zar}({\cal P}^{\rm ex}_n/{\cal W}_n(s))[-k]). 
\end{align*} 
Here the morphism ${\rm Res}$ is the morphism (\ref{eqn:mprrn}) in the special case. 
This morphism is independent of the choice of 
the immersion $X_{\os{\circ}{t}}\os{\sus}{\lo}{\cal P}_n$ 
over ${\cal W}_n(s_{\os{\circ}{t}})$. 
Indeed, let $X_{\os{\circ}{t}}\os{\sus}{\lo}{\cal P}'_n$ be another immersion  
into a log smooth integral scheme over ${\cal W}_n(s_{\os{\circ}{t}})$. 
Then, by considering the product ${\cal P}_n\times_{s_{\os{\circ}{t}}}{\cal P}'_n$,  
we may assume that we have 
the following commutative diagram over ${\cal W}_n(s_{\os{\circ}{t}})$:
\begin{equation*} 
\begin{CD} 
X_{\os{\circ}{t}}@>>> {\cal P}_n \\
@| @VVV \\
X_{\os{\circ}{t}}@>>> {\cal P}'_n.
\end{CD}
\tag{2.3.3.4}\label{cd:xpppp} 
\end{equation*} 
By this commutative diagram and 
the contravariant functoriality of the Poincar\'{e} residue morphism 
((\ref{prop:rescos})), the morphism (\ref{ali:poe}) turns out to 
be the well-defined morphism (\ref{eqn:retsd}). 
\par 
By (\ref{ali:poe}) we see that  the morphism (\ref{eqn:retsd}) kills 
the image of 
$P_{k-1}{\cal W}_n\wt{{\Om}}{}^i_{X_{\os{\circ}{t}}}$ in 
${\cal W}_n\wt{{\Om}}{}^i_{X_{\os{\circ}{t}}}$. 
By (\ref{prop:injf}) (2) and (\ref{prop:perkl}), 
we have the following exact sequence 
\begin{align*} 
& 0 \lo P_{k-1}({\cal O}_{{\mathfrak D}_{n}}
\otimes_{{\cal O}_{{\cal P}^{\rm ex}_n}}
\Om^{\bul}_{{\cal P}^{\rm ex}_n/{\cal W}_n(\os{\circ}{t})} )
\lo
P_k({\cal O}_{{\mathfrak D}_{n}}
\otimes_{{\cal O}_{{\cal P}^{\rm ex}_n}}
{\Om}^{\bul}_{{\cal P}^{\rm ex}_n/{\cal W}_n(\os{\circ}{t})})
\tag{2.3.3.5}\label{ali:powe}\\
& \os{\rm Res}{\lo}  
b^{(k-1)}_*
({\cal O}_{\os{\circ}{\mathfrak D}{}^{(k-1)}_{n}}
\otimes_{{\cal O}_{\os{\circ}{\cal P}{}^{{\rm ex},(k-1)}_n}}
\Om^{\bul}_{\os{\circ}{\cal P}{}^{{\rm ex},(k-1)}_n
/{\cal W}_n}
\otimes_{\mab Z}
\vp^{(k-1)}_{\rm zar}({\cal P}^{\rm ex}_n/{\cal W}_n(\os{\circ}{t})))[-k]
\lo 0. 
\end{align*} 
Because the exactness of the sequence (\ref{eqn:ptwa}) is a local question, 
we may assume that 
there exists the cartesian diagram (\ref{cd:pwtp}) 
for ${\cal P}={\cal P}_n$.  
By (\ref{coro:wtf}) (3) 
we may assume that there exists a local retraction 
of the closed immersion 
$\os{\circ}{\cal P}{}^{{\rm ex},(k-1)}
\os{\sus}{\lo} \os{\circ}{\cal P}$ 
as in \cite[1.2]{msemi}. 
As in [loc.~cit.], 
we can prove that the long exact sequence 
associated to (\ref{ali:powe}) is decomposed into 
the following exact sequence 
\begin{align*} 
0 \lo & 
P_{k-1}({\cal W}_n\wt{{\Om}}{}^i_{X_{\os{\circ}{t}}})  
\lo P_k({\cal W}_n\wt{{\Om}}{}^i_{X_{\os{\circ}{t}}})  
\os{\rm Res}{\lo} \tag{2.3.3.6}\label{eqn:pwpa}\\ 
&b^{(k-1)}_*({\cal H}^{i-k}(
{\cal O}_{\os{\circ}{\mathfrak D}{}^{(k-1)}_{n}}
\otimes_{{\cal O}_{\os{\circ}{\cal P}{}^{{\rm ex},(k-1)}_n}}
\Om^{\bul}_{\os{\circ}{\cal P}{}^{{\rm ex},(k-1)}_n/{\cal W}_n}
\otimes_{\mab Z}
\vp^{(k-1)}_{\rm zar}({\cal P}^{\rm ex}_n/{\cal W}_n(\os{\circ}{t}))))\lo 0.   
\end{align*} 
It is easy to check that ${\rm Res}$ is a morphism of complexes. 
Here note the signs of the derivatives of de Rham(-Witt) complexes. 
We complete the proof of (\ref{prop:prm}). 
\end{proof}

\begin{rema}\label{rema:lcrdc}  
Note that we have used (\ref{coro:wtf}) (3) for the proof of 
(\ref{prop:prm}), which has not been proved in \cite{msemi}.
\end{rema}

\begin{coro}[{\bf cf.~\cite[Lemme 3.3]{msemi}}]\label{coro:cord}
The natural morphism 
\begin{align*} 
P_k({\cal W}_n\wt{{\Om}}{}^i_{X_{\os{\circ}{t}}})\lo 
{\cal W}_n\wt{{\Om}}{}^i_{X_{\os{\circ}{t}}} \quad (k\in {\mab N})
\tag{2.3.5.1}\label{ali:nmo} 
\end{align*} 
is injective. 
Consequently we have a filtered complex 
$({\cal W}_n\wt{{\Om}}{}^{\bul}_{X_{\os{\circ}{t}}},P)$ of 
${\rm C}^+(f^{-1}({\cal W}_n)$-modules, 
where $f\col X_{\os{\circ}{t}}\lo {\cal W}_n(s_{\os{\circ}{t}})$ 
is the structural morphism. 
\end{coro} 

\begin{prop}\label{prop:repmaf}
There exists the following canonical functorial isomorphism 
\begin{align*}
P_0({\cal W}_n\wt{{\Om}}{}^i_{X_{\os{\circ}{t}}})
\os{\sim}{\lo} 
&{\cal H}^i[{\rm MF}(a^{(0)}_*
({\cal W}_n\Om^{\bul}_{\os{\circ}{X}{}^{(0)}_t}
\otimes_{\mab Z}\vp^{(0)}_{\rm zar}
(\os{\circ}{X}_{t_0}/\os{\circ}{t}))
\tag{2.3.6.1}\label{ali:pdete}\\
&\lo 
{\rm MF}(a^{(1)}_*({\cal W}_n\Om^{\bul}_{\os{\circ}{X}{}^{(1)}_t}
\otimes_{\mab Z}\vp^{(1)}_{\rm zar}
(\os{\circ}{X}_{t_0}/\os{\circ}{t}))  \\
&\lo \cdots \\
&\lo 
{\rm MF}(a^{(m)}_*({\cal W}_n\Om^{\bul}_{\os{\circ}{X}{}^{(m)}_t}
\otimes_{\mab Z}\vp^{(m)}_{\rm zar}
(\os{\circ}{X}_{t_0}/\os{\circ}{t}))\\
& \lo \cdots
\cdots)\cdots )]. 
\end{align*}
\end{prop} 
\begin{proof}
This follows from (\ref{ali:pdte}) and 
the well-known equality 
${\cal W}_n\Om^{\bul}_{Y}
=Ru_{Y/{\cal W}_n}({\cal O}_{Y/{\cal W}_n})$ 
for a proper smooth scheme $Y/\kap$ 
(\cite[(7.19)]{ndw}(=a correction of \cite[(4.19)]{hk}), cf.~\cite[II (1.4)]{idw}). 
\end{proof}

\begin{coro}\label{coro:bcptdw}
Let the notations be as in {\rm (\ref{prop:bcdw})}. 
Then the canonical morphism 
\begin{equation*} 
q^{-1}(P_k{\cal W}_n\wt{\Om}_X^i)
\otimes_{{\cal W}_n}{\cal W}_n(\kap_t)
\lo P_k{\cal W}_n\wt{\Om}_{X_{\os{\circ}{t}}}^i
\tag{2.3.7.1}\label{eqn:drckp}
\end{equation*}
is an isomorphism. 
\end{coro} 
\begin{proof}
The same proof as that of (\ref{prop:bcdw}) also works.
\end{proof}

\begin{prop}\label{prop:ovp}
The morphisms 
$F \col  {\cal W}_{n+1}\wt{{\Om}}_{X_{\os{\circ}{t}}}^{\bul}
\lo {\cal W}_n\wt{{\Om}}_{X_{\os{\circ}{t}}}^{\bul}$
and 
$V \col  {\cal W}_n\wt{{\Om}}_{X_{\os{\circ}{t}}}^{\bul}
\lo {\cal W}_{n+1}\wt{{\Om}}_{X_{\os{\circ}{t}}}^{\bul}$
preserve $P$'s. 
\end{prop}
\begin{proof}
This is obvious.  
\end{proof}

Let us recall the following 

\begin{prop}[{\rm {\bf {\cite[(8.4) (1), (2)]{ndw}}}}]\label{prop:pua}  
$(1)$
The morphism 
${\bf p} \col  {\cal W}_n\wt{{\Om}}_{X_{\os{\circ}{t}}}^{\bul}
\lo {\cal W}_{n+1}\wt{{\Om}}_{X_{\os{\circ}{t}}}^{\bul}$
preserves $P$'s. 
\par 
$(2)$ The projection 
$R \col  {\cal W}_{n+1}\wt{{\Om}}_{X_{\os{\circ}{t}}}^{\bul}
\lo {\cal W}_n\wt{{\Om}}_{X_{\os{\circ}{t}}}^{\bul}$
preserves $P$'s. 
\end{prop}
\begin{proof}(We give the proof of (1) 
because we do not use the admissible lift explicitly in this proof, 
which was used in \cite[(8.4) (1)]{ndw}.)
\par 
(1): Let $F_{X_{\os{\circ}{t}}} 
\col X_{\os{\circ}{t}}\lo X_{\os{\circ}{t}}$ 
be the Frobenius endomorphism of $X_{\os{\circ}{t}}$. 
By (\ref{ali:ruynt}) we have the following pull-back morphism 
\begin{align*}
{\cal H}^i(\wt{F}{}^*_{X_{\os{\circ}{t}}}) 
\col 
{\cal H}^i
(P_k\wt{R}u_{X_{\os{\circ}{t}}/{\cal W}_{n+i}(\os{\circ}{t})*}
({\cal O}_{X_{\os{\circ}{t}}/{\cal W}_{n+i}(\os{\circ}{t})}))
\lo 
{\cal H}^i
(P_k\wt{R}u_{X_{\os{\circ}{t}}/{\cal W}_{n+i}(\os{\circ}{t})*}
({\cal O}_{X_{\os{\circ}{t}}/{\cal W}_{n+i}(\os{\circ}{t})})). 
\tag{2.3.9.1}\label{ali:ruypnt}\\
\end{align*}
Because the problem is local, 
we may assume that there exists an immersion 
$X_{\os{\circ}{t}} \os{\sus}{\lo} {\cal P}_{n+i}$ 
into a log smooth integral scheme over ${\cal W}_{n+i}(s_{\os{\circ}{t}})$ 
such that there exists an endomorphism  
$\varphi_{n+i}\col {\cal P}_{n+i}\lo {\cal P}_{n+i}$ which is a lift of  
the Frobenius endomorphism of 
${\cal P}_{n+i}\times_{{\cal W}_{n+i}(s_{\os{\circ}{t}})}s_{\os{\circ}{t}}$.
Let ${\mathfrak D}_{n+i}$ be the log PD-envelope of 
the immersion $X_{\os{\circ}{t}} \os{\sus}{\lo} {\cal P}_{n+i}$ 
over $({\cal W}_{n+i}(s_{\os{\circ}{t}}),p{\cal W}_{n+i},[~])$. 
The morphism (\ref{ali:ruypnt}) is nothing but the following morphism 
\begin{align*} 
\varphi^*_{n+i}\col &
{\cal H}^i(P_k({\cal O}_{{\mathfrak D}_{n+i}}
\otimes_{{\cal O}_{{\cal P}^{\rm ex}_{n+i}}}
\Om^{\bul}_{{\cal P}^{\rm ex}_{n+i}/{\cal W}_{n+i}(\os{\circ}{t})}))
\lo 
{\cal H}^i(P_k({\cal O}_{{\mathfrak D}_{n+i}}
\otimes_{{\cal O}_{{\cal P}^{\rm ex}_{n+i}}}
\Om^{\bul}_{{\cal P}^{\rm ex}_{n+i}/{\cal W}_{n+i}(\os{\circ}{t})})). 
\end{align*} 
By the argument before (\ref{rema:infp}), 
the restriction of the morphism ${\bf p}$ to $P_k{\cal W}_n\wt{\Om}_{X_t}^{\bul}$
is equal to the following morphism 
\begin{align*} 
p^{-(i-1)}\varphi^*_{n+i}\col &
{\cal H}^i(P_k({\cal O}_{{\mathfrak D}_{n}}
\otimes_{{\cal O}_{{\cal P}^{\rm ex}_{n}}}
\Om^{\bul}_{{\cal P}^{\rm ex}_{n}/{\cal W}_{n}(\os{\circ}{t})}))
\lo 
{\cal H}^i(P_k({\cal O}_{{\mathfrak D}_{n+1}}
\otimes_{{\cal O}_{{\cal P}^{\rm ex}_{n+1}}}
\Om^{\bul}_{{\cal P}^{\rm ex}_{n+1}/{\cal W}_{n+1}(\os{\circ}{t})}))\\ 
&\os{\sus}{\lo} 
{\cal H}^i({\cal O}_{{\mathfrak D}_{n+1}}
\otimes_{{\cal O}_{{\cal P}^{\rm ex}_{n+1}}}
\Om^{\bul}_{{\cal P}^{\rm ex}_{n+1}/{\cal W}_{n+1}(\os{\circ}{t})}). 
\end{align*}  
This proves (1). 
\par 
(2): See the proof of \cite[(8.4) (2)]{ndw}. 
\end{proof} 


\begin{prop}\label{prop:aco}  
The morphism ${\rm Res}$ in {\rm (\ref{eqn:retsd})} 
is compatible with $F$, $V$, $d$, ${\bf p}$ and $R$ on 
$\{P_k{\cal W}_m\wt{{\Om}}{}^i_{X_{\os{\circ}{t}}}\}_{m=1}^n$.
\end{prop} 
\begin{proof} 
In the proof of (\ref{prop:prm}),   
we have already mentioned that ${\rm Res}$ is compatible with $d$. 
In \cite[(8.4.1), (8.4.2), (8.4.3)]{ndw} we have proved that 
the morphism (\ref{eqn:retsd}) is compatible with ${\bf p}$ and $R$ on 
$P_k{\cal W}_n\wt{{\Om}}{}^i_{X_{\os{\circ}{t}}}$.  
The compatibility of $F$ and $V$ is obvious. 
\end{proof}

\begin{prop}\label{prop:sprj} 
Let $k$ be a nonnegative integer and let $n$ be a positive integer. 
Then the projection  
\begin{align*} 
R\col P_k{\cal W}_{n+1}\wt{\Om}_{X_{\os{\circ}{t}}}^{\bul}
\lo 
P_k{\cal W}_n\wt{\Om}_{X_{\os{\circ}{t}}}^{\bul}
\tag{2.3.11.1}\label{ali:wno}
\end{align*}  
is surjective. 
\end{prop} 
\begin{proof}
Let $i$ be a nonnegative integer. 
Because the problem is local, we may assume that there exists a positive integer 
$K$ such that $P_K{\cal W}_j\wt{\Om}_{X_{\os{\circ}{t}}}^i={\cal W}_j\wt{\Om}_{X_{\os{\circ}{t}}}^i$ 
$(j=n,n+1)$. By \cite[(6.28) (8)]{ndw} the projection 
$R\col {\cal W}_{n+1}\wt{\Om}_{X_{\os{\circ}{t}}}^i\lo {\cal W}_n\wt{\Om}_{X_{\os{\circ}{t}}}^i$ is surjective. 
(This surjectivity also follows from (\ref{cd:pjpitlm}) immediately.)
Consider the following commutative diagram of exact sequences 
\begin{equation*}
\begin{CD}
0 @>>> P_{k-1}{\cal W}_{n+1}\wt{\Om}_{X_{\os{\circ}{t}}}^i @>>>  
P_k{\cal W}_{n+1}\wt{\Om}_{X_{\os{\circ}{t}}}^i @>{\rm Res}>> 
{\cal W}_{n+1}\Om_{\os{\circ}{X}{}^{(k-1)}_t}^{i-k}@>>> 0\\ 
@. @V{R}VV  
@V{R}VV @V{R}VV \\
0 @>>> P_{k-1}{\cal W}_n\wt{\Om}_{X_{\os{\circ}{t}}}^i @>>>  
P_k{\cal W}_n\wt{\Om}_{X_{\os{\circ}{t}}}
@>{\rm Res}>> 
{\cal W}_n\Om_{\os{\circ}{X}{}^{(k-1)}_t}^{i-k}@>>> 0. 
\end{CD}
\tag{2.3.11.2}\label{cd:orexgram}
\end{equation*}
By the snake lemma and the descending induction, 
we have only to prove that 
the morphism 
\begin{align*} 
\{{\rm Ker}(P_k{\cal W}_{n+1}\wt{\Om}_{X_{\os{\circ}{t}}}^i\lo 
P_k{\cal W}_n\wt{\Om}_{X_{\os{\circ}{t}}}^i)\} \lo 
{\rm Fil}^n({\cal W}_{n+1}\Om_{\os{\circ}{X}{}^{(k-1)}_t}^{i-k})
=V^n\Om_{\os{\circ}{X}{}^{(k-1)}_t}^{i-k}
+dV^n\Om_{\os{\circ}{X}{}^{(k-1)}_t}^{i-1-k}
\end{align*} 
is surjective. 
Consider a subsheaf 
$V^nP_k\wt{\Om}_{X_{\os{\circ}{t}}}^i+dV^nP_k\wt{\Om}_{X_{\os{\circ}{t}}}^{i-1}$ in 
${\rm Fil}^n({\cal W}_{n+1}\wt{\Om}_{X_{\os{\circ}{t}}}^i)$. 
Because ${\rm Res}$ is compatible with $V$ and $d$, 
the upper ${\rm Res}$ in (\ref{cd:orexgram}) 
induces a surjective morphism 
$V^nP_k\wt{\Om}_{X_{\os{\circ}{t}}}^i+dV^nP_k\wt{\Om}_{X_{\os{\circ}{t}}}^{i-1}\lo 
V^n\Om_{\os{\circ}{X}{}^{(k-1)}_t}^{i-k}
+dV^n\Om_{\os{\circ}{X}{}^{(k-1)}_t}^{i-1-k}$. 
Because $V$ preserves $P$, 
$V^nP_k\wt{\Om}_{X_{\os{\circ}{t}}}^i+dV^nP_k\wt{\Om}_{X_{\os{\circ}{t}}}^{i-1}$ 
is contained in ${\rm Ker}(P_k{\cal W}_{n+1}\wt{\Om}_{X_{\os{\circ}{t}}}^i\lo 
P_k{\cal W}_n\wt{\Om}_{X_{\os{\circ}{t}}}^i)$. 
Thus we have proved the desired surjectivity.   
\end{proof}

Let $E$ be a flat coherent log crystal of 
${\cal O}_{\os{\circ}{X}_t/{\cal W}_n(\os{\circ}{t})}$-modules.  
Set  
\begin{align*} 
& P_k((\eps^*_{X_{\os{\circ}{t}}
/{\cal W}_n(\os{\circ}{t})}(E))_{{\cal W}_n(X_{\os{\circ}{t}})}
\otimes_{{\cal W}_n({\cal O}_{X_{\os{\circ}{t}}})}
{\cal W}_n\wt{{\Om}}{}^i_{X_{\os{\circ}{t}}}) 
\tag{2.3.11.3}\label{eqn:ply}\\
& :=
(\eps^*_{X_{{\os{\circ}{t}}}/{\cal W}_n(\os{\circ}{t})}(E))_{{\cal W}_n(X_{\os{\circ}{t}})}
\otimes_{{\cal W}_n({\cal O}_{X_{\os{\circ}{t}}})}
P_k{\cal W}_n\wt{{\Om}}{}^i_{X_{\os{\circ}{t}}} 
\quad (i,k\in {\mab Z}).  
\end{align*} 
Because $E$ is flat, 
$$P_k((\eps^*_{X_{\os{\circ}{t}}/
{\cal W}_n(\os{\circ}{t})}(E))_{{\cal W}_n(X_{\os{\circ}{t}})}
\otimes_{{\cal W}_n({\cal O}_{X_{\os{\circ}{t}}})}
{\cal W}_n\wt{{\Om}}{}^i_{X_{\os{\circ}{t}}})$$ 
is a subsheaf of 
$(\eps^*_{X_{\os{\circ}{t}}/{\cal W}_n(\os{\circ}{t})}(E))_{{\cal W}_n(X_{\os{\circ}{t}})}
\otimes_{{\cal W}_n({\cal O}_{X_{\os{\circ}{t}}})}{\cal W}_n\wt{{\Om}}{}^i_{X_{\os{\circ}{t}}}$.

\begin{prop}\label{prop:lbsl} 
$(1)$ The connection 
\begin{align*} 
\nabla \col & (\eps^*_{X_{\os{\circ}{t}}/{\cal W}_n(\os{\circ}{t})}(E))_{
{\cal W}_n(X_{\os{\circ}{t}})}
\otimes_{{\cal W}_n({\cal O}_{X_{\os{\circ}{t}}})}
{\cal W}_n\wt{{\Om}}{}^i_{X_{\os{\circ}{t}}}\lo 
(\eps^*_{X_{\os{\circ}{t}}/{\cal W}_n(\os{\circ}{t})}(E))_{{\cal W}_n(X_{\os{\circ}{t}})}
\otimes_{{\cal W}_n({\cal O}_{X_{\os{\circ}{t}}})}
{\cal W}_n\wt{{\Om}}{}^{i+1}_{X_{\os{\circ}{t}}} 
\end{align*}  
preserves $P$'s. 
\par 
$(2)$ 
Let $E$ be a flat coherent log crystal of 
${\cal O}_{\os{\circ}{X}_t/{\cal W}_{n+1}(\os{\circ}{t})}$-modules.  
By abuse of notation, denote by the same symbol $E$ 
the log crystal of ${\cal O}_{\os{\circ}{X}_t/{\cal W}_{n}(\os{\circ}{t})}$-modules 
obtained by $E$. 
Then the following morphisms 
\begin{align*} 
{\bf p} \col & 
(\eps^*_{X_{\os{\circ}{t}}/{\cal W}_n(\os{\circ}{t})}(E))_{
{\cal W}_n(X_{\os{\circ}{t}})}\otimes_{{\cal W}_n({\cal O}_{X_{\os{\circ}{t}}})}
{\cal W}_n\wt{{\Om}}_{X_{\os{\circ}{t}}}^{\bul}
\lo \\ 
& 
(\eps^*_{X_{\os{\circ}{t}}/{\cal W}_{n+1}(\os{\circ}{t})}(E))_{
{\cal W}_{n+1}(X_{\os{\circ}{t}})}\otimes_{{\cal W}_{n+1}
({\cal O}_{X_{\os{\circ}{t}}})}
{\cal W}_{n+1}\wt{{\Om}}_{X_{\os{\circ}{t}}}^{\bul}, 
\end{align*} 
\begin{align*} 
R \col & 
(\eps^*_{X_{\os{\circ}{t}}/{\cal W}_{n+1}(\os{\circ}{t})}(E))_{
{\cal W}_{n+1}(X_{\os{\circ}{t}})}
\otimes_{{\cal W}_{n+1}({\cal O}_{X_{\os{\circ}{t}}})}
{\cal W}_{n+1}\wt{{\Om}}{}^{\bul}_{X_{\os{\circ}{t}}}
\lo \\
& (\eps^*_{X_{\os{\circ}{t}}/{\cal W}_n(\os{\circ}{t})}(E))_{
{\cal W}_n(X_{\os{\circ}{t}})}
\otimes_{{\cal W}_n({\cal O}_{X_{\os{\circ}{t}}})}
{\cal W}_n\wt{{\Om}}{}^{\bul}_{X_{\os{\circ}{t}}} 
\end{align*} 
and 
\begin{align*} 
F \col & 
(\eps^*_{X_{\os{\circ}{t}}/{\cal W}_{n+1}(\os{\circ}{t})}(E))_{
{\cal W}_{n+1}(X_{\os{\circ}{t}})}
\otimes_{{\cal W}_{n+1}({\cal O}_{X_{\os{\circ}{t}}})}
{\cal W}_{n+1}\wt{{\Om}}{}^{\bul}_{X_{\os{\circ}{t}}}
\lo \\
& (\eps^*_{X_{\os{\circ}{t}}/{\cal W}_n(\os{\circ}{t})}(E))_{
{\cal W}_n(X_{\os{\circ}{t}})}
\otimes_{{\cal W}_n({\cal O}_{X_{\os{\circ}{t}}})}
{\cal W}_n\wt{{\Om}}{}^{\bul}_{X_{\os{\circ}{t}}} 
\end{align*} 
preserve $P$'s. 
\end{prop} 
\begin{proof} 
(1): The problem is local. Let 
$X_{\os{\circ}{t}}\os{\sus}{\lo} \ol{\cal P}$ be an immersion  
into a log smooth integral scheme over $\ol{{\cal W}_n(s_{\os{\circ}{t}})}$. 
Let $\ol{\mathfrak D}$ be the log PD-envelope of this immersion 
over $({\cal W}_n(\os{\circ}{t}),p{\cal W}_n,[~])$. 
Set ${\cal P}:=\ol{\cal P}\times_{\ol{{\cal W}_n(s)}}{\cal W}_n(s)$ 
and ${\mathfrak D}:=\ol{\mathfrak D}\times_{{\mathfrak D}(\ol{{\cal W}_n(s)})}
{\cal W}_n(s)$.  
Set $\ol{E}_{n}:=E_{{\cal W}_n(\os{\circ}{X}_t)}$. 
Because $E$ is a crystal of 
${\cal O}_{\os{\circ}{X}_t/{\cal W}_n(\os{\circ}{t})}$-modules, 
$\ol{E}_{n}=
(\eps^*_{X_{\os{\circ}{t}}/{\cal W}_n(\os{\circ}{t})}(E))_{{\cal W}_n(X_{\os{\circ}{t}})}$. 
By (\ref{eqn:intcon}) (for the trivial logarithmic case), 
we have the following integral connection:  
\begin{equation*} 
\nabla \col \ol{E}_{n}\lo \ol{E}_{n}\otimes_{{\cal W}_n({\cal O}_{X_t})'}
\Om^1_{{\cal W}_n({\os{\circ}{X}_t})/{\cal W}_n(\os{\circ}{t}),[~]}.   
\tag{2.3.12.1}\label{eqn:eexnny}  
\end{equation*} 
Using this connection, 
we have the complex 
$\ol{E}_{n}\otimes_{{\cal W}_n({\cal O}_{X_t})'}
\Om^{\bul}_{{\cal W}_n({\os{\circ}{X}_t})/{\cal W}_n(\os{\circ}{t}),[~]}$. 
We define the following two morphisms: 
\begin{equation*}
\wt{s}_n(0,0) \col {\cal W}_n({\cal O}_{X_{\os{\circ}{t}}}) 
\owns (a_0, \ldots, a_{n-1}) \lom 
\left[\sum_{i=0}^{n-1}\wt{a}_i^{p^{n-i}}p^i\right]
\in 
{\cal H}^0({\cal O}_{\ol{\mathfrak D}}\otimes_{{\cal O}_{\ol{\cal P}{}^{\rm ex}}}
\Om^*_{\os{\circ}{\ol{\cal P}}{}^{\rm ex}/{\cal W}_n(\os{\circ}{t})}),
\tag{2.3.12.2}\label{eqn:tis0}
\end{equation*}
\begin{equation*}
\wt{s}_n(1,0)' \col {\cal W}_n({\cal O}_{X_{\os{\circ}{t}}}) \owns 
d(a_0, \ldots, a_{n-1}) \lom 
\left[\sum_{i=0}^{n-1}\wt{a}_i^{p^{n-i}-1}d\wt{a}_i\right]
\in {\cal H}^1({\cal O}_{\ol{\mathfrak D}}\otimes_{{\cal O}_{\ol{\cal P}{}^{\rm ex}}}
\Om^*_{\os{\circ}{\ol{\cal P}}{}^{\rm ex}/{\cal W}_n(\os{\circ}{t})}),
\tag{2.3.12.3}\label{eqn:tis10}
\end{equation*}
where $\wt{a}_i\in {\cal O}_{\ol{\mathfrak D}}$ 
is a lift of $a_i\in {\cal O}_{X_{\os{\circ}{t}}}$.
Then, by the same proof as that of \cite[(4.9)]{hk}, 
$\wt{s}_n(0,0)$ and $\wt{s}_n(1,0)'$ are  well-defined. 
In a standard way, we have a complex 
${\cal H}^{\bul}({\cal O}_{\ol{\mathfrak D}}\otimes_{{\cal O}_{\ol{\cal P}{}^{\rm ex}}}
\Om^*_{\os{\circ}{\ol{\cal P}}{}^{\rm ex}/{\cal W}_n(\os{\circ}{t})})$. 
The morphisms $\wt{s}_n(0,0)$ and $\wt{s}_n(1,0)'$ induces 
the following morphism of complexes as in \S\ref{sec:ldrwc}:
\begin{align*} 
\Om^{\bul}_{{\cal W}_n({\os{\circ}{X}_t})/{\cal W}_n(\os{\circ}{t}),[~]}
\lo 
{\cal H}^{\bul}({\cal O}_{\ol{\mathfrak D}}\otimes_{{\cal O}_{\ol{\cal P}{}^{\rm ex}}}
\Om^*_{\os{\circ}{\ol{\cal P}}{}^{\rm ex}/{\cal W}_n(\os{\circ}{t})}).
\tag{2.3.12.4}\label{ali:dopc} 
\end{align*}
By (\ref{eqn:pops}) we have the following natural morphism
\begin{align*} 
{\cal H}^{\bul}({\cal O}_{\ol{\mathfrak D}}\otimes_{{\cal O}_{\ol{\cal P}{}^{\rm ex}}}
\Om^*_{\os{\circ}{\ol{\cal P}}{}^{\rm ex}/{\cal W}_n(\os{\circ}{t})})
\lo 
{\cal H}^{\bul}(P_0({\cal O}_{{\mathfrak D}}\otimes_{{\cal O}_{{\cal P}{}^{\rm ex}}}
\Om^*_{{\cal P}^{\rm ex}/{\cal W}_n(\os{\circ}{t})})).
\tag{2.3.12.5}\label{ali:ddpc} 
\end{align*}
Hence we have the following composite morphism 
\begin{align*} 
\Om^{\bul}_{{\cal W}_n({\os{\circ}{X}})/{\cal W}_n(\os{\circ}{t}),[~]}
\lo 
{\cal H}^{\bul}(P_0({\cal O}_{{\mathfrak D}}\otimes_{{\cal O}_{{\cal P}{}^{\rm ex}}}
\Om^*_{{\cal P}^{\rm ex}/{\cal W}_n(\os{\circ}{t})}))=P_0{\cal W}_n\wt{\Om}{}^{\bul}_{X_{\os{\circ}{t}}}.
\tag{2.3.12.6}\label{ali:dipc} 
\end{align*}
This morphism induces the following morphism 
\begin{align*} 
\Om^{\bul}_{{\cal W}_n({\os{\circ}{X}_t})/{\cal W}_n(\os{\circ}{t}),[~]}
\lo P_0{\cal W}_n\wt{\Om}{}^{\bul}_{X_{\os{\circ}{t}}}.
\tag{2.3.12.7}\label{ali:dinc} 
\end{align*}
Thus we have the complex  
$\ol{E}_{n}\otimes_{{\cal W}_n({\cal O}_{X})}
P_0{\cal W}_n\wt{\Om}{}^{\bul}_{X_{\os{\circ}{t}}}$. 
Because $P_k{\cal W}_n\wt{{\Om}}{}^i_{X_{\os{\circ}{t}}}$ is 
stable under the actions of the operator $d$ on 
$P_k{\cal W}_n\wt{{\Om}}{}^i_{X_{\os{\circ}{t}}}$ ((\ref{ali:bopant})), 
$P_k((\eps^*_{X_{\os{\circ}{t}}/{\cal W}_n(\os{\circ}{t})}(
E_{\os{\circ}{X}_{t}/{\cal W}_n(\os{\circ}{t})}))_{{\cal W}_n(X_{\os{\circ}{t}})}
\otimes_{{\cal W}_n({\cal O}_X)}{\cal W}_n\wt{{\Om}}{}^i_X)$ is also stable.  
\par 
(2): 
Because $P_k{\cal W}_m\wt{{\Om}}{}^i_{X_{\os{\circ}{t}}}$ $(m=n,n+1)$ is 
stable under the actions of the operators ${\bf p}$ and $R$ ((\ref{prop:pua})), 
$P_k((\eps^*_{X_{\os{\circ}{t}}/{\cal W}_m(\os{\circ}{t})}(E))_{{\cal W}_m(X_{\os{\circ}{t}})}
\otimes_{{\cal W}_m({\cal O}_{X_{\os{\circ}{t}}})}{\cal W}_m\wt{{\Om}}{}^i_{X_{\os{\circ}{t}}})$  
is also stable.  By the local description of $F$, 
it is obvious that 
$P_k((\eps^*_{X_{\os{\circ}{t}}/{\cal W}_m(\os{\circ}{t})}(E))_{{\cal W}_m(X_{\os{\circ}{t}})}
\otimes_{{\cal W}_m({\cal O}_{X_{\os{\circ}{t}}})}{\cal W}_m\wt{{\Om}}{}^i_{X_{\os{\circ}{t}}})$ $(m=n,n+1)$ is 
stable under the action of the operator $F$. 
\end{proof}

\begin{prop}\label{prop:nip} 
$(1)$ Let $E$ be a flat coherent log crystal of 
${\cal O}_{\os{\circ}{X}_t/{\cal W}_n(\os{\circ}{t})}$-modules.  
Let 
\begin{align*} 
a^{(k-1)}_{\rm crys} \col 
((\os{\circ}{X}{}^{(k-1)}_t/{\cal W}_n(\os{\circ}{t}))_{\rm crys},
{\cal O}_{\os{\circ}{X}{}^{(k-1)}_t/{\cal W}_n(\os{\circ}{t})})
\lo 
((\os{\circ}{X}_t/{\cal W}_n(\os{\circ}{t}))_{\rm crys},
{\cal O}_{\os{\circ}{X}_t/{\cal W}_n(\os{\circ}{t})}) 
\quad (k\in {\mab Z}_{\geq 1})
\end{align*}
be the natural morphism of ringed topoi. 
Let $E_{\os{\circ}{X}{}^{(k-1)}_t/{\cal W}_n(\os{\circ}{t})}$ be the crystal of 
${\cal O}_{\os{\circ}{X}{}^{(k-1)}_t/{\cal W}_n(\os{\circ}{t})}$-modules 
obtained by $E$.  
Then there exists 
the following Poincar\'{e} residue morphism for 
$k\in {\mab Z}_{\geq 1}$ 
\begin{align*} 
{\rm Res} \col & 
P_k((\eps^*_{X_{\os{\circ}{t}}/{\cal W}_n(\os{\circ}{t}){\rm crys}}(E))_{{\cal W}_n(X_{\os{\circ}{t}})}
\otimes_{{\cal W}_n({\cal O}_{X_{\os{\circ}{t}}})}
{\cal W}_n\wt{{\Om}}{}^i_{X_{\os{\circ}{t}}})
\lo \tag{2.3.13.1}\label{eqn:resd} \\
& a^{(k-1)}_*((E_{\os{\circ}{X}{}^{(k-1)}_t/
{\cal W}_n(\os{\circ}{t})})_{{\cal W}_n(\os{\circ}{X}{}^{(k-1)}_t)})
\otimes_{{\cal W}_n({\cal O}_{\os{\circ}{X}{}^{(k-1)}_t})}
{\cal W}_n\Om^{i-k}_{\os{\circ}{X}{}^{(k-1)}_t}
\otimes_{\mab Z}\vp^{(k-1)}_{\rm zar}(\os{\circ}{X}_t/\os{\circ}{t}) 
\end{align*} 
fitting into the following exact sequence 
\begin{align*} 
0 \lo & 
P_{k-1}((\eps^*_{X_{\os{\circ}{t}}/{\cal W}_n(\os{\circ}{t})}(E))_{{\cal W}_n(X_{\os{\circ}{t}})}
\otimes_{{\cal W}_n({\cal O}_{X_{\os{\circ}{t}}})}
{\cal W}_n\wt{{\Om}}{}^i_{X_{\os{\circ}{t}}})
\lo \tag{2.3.13.2}\label{eqn:pwa}\\
& P_k((\eps^*_{X_{\os{\circ}{t}}/{\cal W}_n(\os{\circ}{t})}(E))_{{\cal W}_n(X_{\os{\circ}{t}})}
\otimes_{{\cal W}_n({\cal O}_{X_{\os{\circ}{t}}})}
{\cal W}_n\wt{{\Om}}{}^i_{X_{\os{\circ}{t}}})
\os{\rm Res}{\lo} \\
& 
a^{(k-1)}_*
((E_{\os{\circ}{X}{}^{(k-1)}_t/{\cal W}_n(\os{\circ}{t})})_{{\cal W}_n(\os{\circ}{X}{}^{(k-1)}_t)}
\otimes_{{\cal W}_n({\cal O}_{\os{\circ}{X}{}^{(k-1)}_t})}
{\cal W}_n\Om^{i-k}_{\os{\circ}{X}{}^{(k-1)}_t}
\otimes_{\mab Z}\vp^{(k-1)}_{\rm zar}(\os{\circ}{X}_t/\os{\circ}{t}) 
\lo 0.   
\end{align*} 
\par 
$(2)$ The Poincar\'{e} residue morphism {\rm (\ref{eqn:resd})} 
is compatible with $F$, $d$, ${\bf p}$ and $R$. 
\end{prop}
\begin{proof}
(1): (1) follows from (\ref{eqn:ply}) and (\ref{prop:prm}).  
\par 
(2): (2) follows from (\ref{prop:aco}).  
\end{proof}

\par 
Let $E$ be a flat coherent log crystal of 
${\cal O}_{\os{\circ}{X}_t/{\cal W}(\os{\circ}{t})}$-modules.   
By abuse of notation, denote by the same symbol $E$ 
the log crystal of 
${\cal O}_{\os{\circ}{X}_t/{\cal W}_{n}(\os{\circ}{t})}$-modules 
obtained by $E$. 
Then, by (\ref{prop:lbsl}) (2), we have the projective system 
\begin{equation*} 
\{((\eps^*_{X_t/{\cal W}_n(\os{\circ}{t})}(E))_{{\cal W}_n(X_{\os{\circ}{t}})}
\otimes_{{\cal W}_n({\cal O}_{X_{\os{\circ}{t}}})}
{\cal W}_n\wt{{\Om}}{}^{\bul}_{X_{\os{\circ}{t}}},P)\}_{n=1}^{\infty}
\end{equation*}  
of filtered complexes in $X_{t,{\rm zar}}$. 
Set 
\begin{align*} 
&((\eps^*_{X_{\os{\circ}{t}}/{\cal W}(\os{\circ}{t})}(E))_{{\cal W}(X_{\os{\circ}{t}})}
\otimes_{{\cal W}({\cal O}_{X_{\os{\circ}{t}}})}
{\cal W}\wt{{\Om}}{}^{\bul}_X,P)\\
&:=  
\vpl_n((\eps^*_{X_{\os{\circ}{t}}/{\cal W}_n(\os{\circ}{t})}(E))_{
{\cal W}_n(X_{\os{\circ}{t}})}\otimes_{{\cal W}_n({\cal O}_{X_{\os{\circ}{t}}})}
{\cal W}_n\wt{{\Om}}{}^{\bul}_{X_{\os{\circ}{t}}},P)\\
&=  
(\vpl_n(\eps^*_{X_{\os{\circ}{t}}/
{\cal W}_n(\os{\circ}{t})}(E))_{
{\cal W}_n(X_{\os{\circ}{t}})})
\otimes_{{\cal W}({\cal O}_{X_{\os{\circ}{t}}})}
{\cal W}\wt{{\Om}}{}^{\bul}_{X_{\os{\circ}{t}}},P). 
\end{align*} 
By (\ref{ali:wno}) and (\ref{eqn:pwa}) we have 
the following exact sequence 
\begin{align*} 
0 \lo & 
P_{k-1}((\eps^*_{X_{\os{\circ}{t}}/{\cal W}(\os{\circ}{t})}(E))_{{\cal W}(X_{\os{\circ}{t}})}
\otimes_{{\cal W}({\cal O}_{X_{\os{\circ}{t}}})}{\cal W}\wt{{\Om}}{}^i_{X_{\os{\circ}{t}}})
\lo \tag{2.3.13.3}\label{eqn:pwia}\\
& P_k((\eps^*_{X_{\os{\circ}{t}}/{\cal W}(\os{\circ}{t})}(E))_{{\cal W}(X_{\os{\circ}{t}})}
\otimes_{{\cal W}({\cal O}_{X_{\os{\circ}{t}}})}
{\cal W}\wt{{\Om}}{}^i_{X_{\os{\circ}{t}}})
\os{\rm Res}{\lo} \\
& 
a^{(k-1)}_*
(E_{\os{\circ}{X}{}^{(k-1)}_t/{\cal W}(\os{\circ}{t})})_{{\cal W}(\os{\circ}{X}{}^{(k-1)}_t)}
\otimes_{{\cal W}({\cal O}_{\os{\circ}{X}{}^{(k-1)}_t})}
{\cal W}\Om^{i-k}_{\os{\circ}{X}{}^{(k-1)}_t}
\otimes_{\mab Z}\vp^{(k-1)}_{\rm zar}(\os{\circ}{X}_t/\os{\circ}{t}) 
\lo 0.   
\end{align*}

\begin{defi}\label{defi:cthi}  
Let $\star$ be noting or a positive integer $n$. 
Let $E$ be a flat coherent log crystal of 
${\cal O}_{\os{\circ}{X}_t/{\cal W}_{\star}(\os{\circ}{t})}$-modules.  
We call the isomorphism 
\begin{align*}  
& {\rm Res} \col 
{\rm gr}^P_k
((\eps^*_{X_{\os{\circ}{t}}/{\cal W}_{\star}(\os{\circ}{t})}(E))_{{\cal W}_{\star}(X_{\os{\circ}{t}})}
\otimes_{{\cal W}_{\star}({\cal O}_{X_{\os{\circ}{t}}})}
{\cal W}_{\star}\wt{{\Om}}{}^{\bul}_
{X_{\os{\circ}{t}}}) 
\os{\sim}{\lo} \tag{2.3.14.1}\label{eqn:pris}\\
& 
a^{(k)}_{{\cal W}_{\star}(\os{\circ}{t})*} 
((E_{\os{\circ}{X}{}^{(k)}_t/{\cal W}_{\star}
(\os{\circ}{t})})_{{\cal W}_{\star}(\os{\circ}{X}{}^{(k)}_{t})})
\otimes_{{\cal W}_{\star}(\os{\circ}{X}{}^{(k)}_t)}
{\cal W}_{\star}\Om^{\bul}_{\os{\circ}{X}{}^{(k)}_t}
\otimes_{\mab Z}\vp^{(k)}_{\rm zar}(X_{\os{\circ}{t}}/{\os{\circ}{t}}))[-k](-k) 
\end{align*} 
the {\it Poincar\'{e} residue isomorphism 
of the log de Rham-Witt complex} 
of 
$\eps^*_{X_{\os{\circ}{t}}/{\cal W}_{\star}(\os{\circ}{t})}(E)$. 
\end{defi} 

\begin{defi}\label{defi:ene} 
Let $E$ be a flat coherent log crystal of 
${\cal O}_{\os{\circ}{X}_t/{\cal W}(\os{\circ}{t})}$-modules.   
Let $E_{\os{\circ}{X}_t/{\cal W}_{n}(\os{\circ}{t})}$ be the log crystal of 
${\cal O}_{\os{\circ}{X}_t/{\cal W}_{n}(\os{\circ}{t})}$-modules 
obtained by $E$. 
Assume that $E_{\os{\circ}{X}_t/{\cal W}_{n}(\os{\circ}{t})}$ 
is a unit root log $F$-crystal. 
Then we have the following isomorphism 
\begin{align*} 
\Phi^{-1}\otimes {\rm id} \col &
(\eps^*_{X_{\os{\circ}{t}}/{\cal W}_n(\os{\circ}{t})}(E))_{{\cal W}_n(X_{\os{\circ}{t}})}
\otimes_{{\cal W}_n({\cal O}_{X_{\os{\circ}{t}}})}
{\cal W}_n\wt{\Om}^{\bul}_{X_{\os{\circ}{t}}}
\os{\sim}{\lo}  \\
&
((\eps^*_{X_{\os{\circ}{t}}/{\cal W}_n(\os{\circ}{t})}(E))_{{\cal W}_n(X_{\os{\circ}{t}})})^{\sig}
\otimes_{{\cal W}_n({\cal O}_{X_{\os{\circ}{t}}})}
{\cal W}_n\wt{\Om}^{\bul}_{X_{\os{\circ}{t}}}. 
\end{align*} 
Next consider the following isomorphism 
\begin{align*} 
I\col & 
((\eps^*_{X_{\os{\circ}{t}}/{\cal W}_n(\os{\circ}{t})}(E))_{{\cal W}_n(X_{\os{\circ}{t}})})^{\sig}
\otimes_{{\cal W}_n({\cal O}_{X_{\os{\circ}{t}}})}
{\cal W}_n\wt{\Om}^{\bul}_{X_{\os{\circ}{t}}} \\
&= 
F^{-1}((\eps^*_{X_{\os{\circ}{t}}/{\cal W}_n(\os{\circ}{t})}(E))_{{\cal W}_n(X_{\os{\circ}{t}})})
\otimes_{F^{-1}({\cal W}_n({\cal O}_{X_{\os{\circ}{t}}}))} 
{\cal W}_n\wt{\Om}^{\bul}_{X_{\os{\circ}{t}}}  \\
&=(\eps^*_{X_{\os{\circ}{t}}/{\cal W}_n(\os{\circ}{t})}(E))
_{{\cal W}_n(X_{\os{\circ}{t}}))}
\otimes_{{\cal W}_n({\cal O}_{X_{\os{\circ}{t}}}))}
F_*({\cal W}_n{\Om}^i_{X_{\os{\circ}{t}}}). 
\end{align*} 
Consider also the following morphism 
\begin{align*} 
{\rm id}\otimes V \col & 
(\eps^*_{X_{\os{\circ}{t}}/{\cal W}_n(\os{\circ}{t})}(E))_{{\cal W}_n(X_{\os{\circ}{t}})}
\otimes_{{\cal W}_n({\cal O}_{X_{\os{\circ}{t}}})}
F_*({\cal W}_n{\Om}^i_{X_{\os{\circ}{t}}}) \\
& =
((\eps^*_{X_{\os{\circ}{t}}/{\cal W}_{n+1}(\os{\circ}{t})}(E))_{{\cal W}_{n+1}(X_{\os{\circ}{t}})}
\otimes_{{\cal W}_{n+1}({\cal O}_{Y_{\bul \leq N}})^{\star},R}
{\cal W}_n({\cal O}_{X_{\os{\circ}{t}}}))
\otimes_{{\cal W}_n({\cal O}_{X_{\os{\circ}{t}}})} 
F_*({\cal W}_n{\Om}^i_{X_{\os{\circ}{t}}}) \\
& \lo 
(\eps^*_{X_{\os{\circ}{t}}/{\cal W}_{n+1}(\os{\circ}{t})}(E))_{{\cal W}_{n+1}(X_{\os{\circ}{t}})}
\otimes_{{\cal W}_{n+1}({\cal O}_{X_{\os{\circ}{t}}})}
({\cal W}_{n+1}{\Om}^i_{X_{\os{\circ}{t}}}).  
\end{align*} 
(Here we have used the formula ``$xVy=V(Fxy)$'' 
to define ${\rm id}\otimes V$.) 
Let 
\begin{align*} 
V\col &(\eps^*_{X_{\os{\circ}{t}}/{\cal W}_n(\os{\circ}{t})}(E))_{{\cal W}_n(X_{\os{\circ}{t}})}
\otimes_{{\cal W}_n({\cal O}_{X_{\os{\circ}{t}}})}
{\cal W}_n{\Om}^{\bul}_{X_{\os{\circ}{t}}}\\
&\lo 
(\eps^*_{X_{\os{\circ}{t}}/{\cal W}_{n+1}(\os{\circ}{t})}(E))_{{\cal W}_{n+1}(X_{\os{\circ}{t}})}
\otimes_{{\cal W}_{n+1}({\cal O}_{X_{\os{\circ}{t}}})}
{\cal W}_{n+1}{\Om}^{\bul}_{X_{\os{\circ}{t}}}
\end{align*}
be the composite of $\Phi^{-1}\otimes {\rm id}$, 
$I$ and ${\rm id}\otimes V$: 
$V:= ({\rm id}\otimes V)\circ I\circ 
(\Phi^{-1}\otimes {\rm id})$.  
\end{defi} 

\begin{prop}\label{prop:lbasl} 
Let the notations be as in {\rm (\ref{prop:lbsl})}.  
Assume that $E$ is a unit root $F$-crystal on 
$\os{\circ}{X}_t/{\cal W}_n(\os{\circ}{t})$. 
Then the following hold$:$
\par 
$(1)$ 
The following operator
\begin{align*} 
V \col & 
(\eps^*_{X_{\os{\circ}{t}}/{\cal W}_n(\os{\circ}{t})}(E))_{{\cal W}_n(X_{\os{\circ}{t}})}
\otimes_{{\cal W}_n({\cal O}_{X_{\os{\circ}{t}}})}
{\cal W}_n\wt{{\Om}}{}^i_{X_{\os{\circ}{t}}} \lo \\
&(\eps^*_{X_{\os{\circ}{t}}/{\cal W}_{n+1}(\os{\circ}{t})}(E))_{
{\cal W}_{n+1}(X_{\os{\circ}{t}})}
\otimes_{{\cal W}_{n+1}({\cal O}_{X_{\os{\circ}{t}}})}
{\cal W}_{n+1}\wt{{\Om}}{}^i_{X_{\os{\circ}{t}}}
\tag{2.3.16.1}\label{ali:vesx}
\end{align*} 
preserves $P$'s. 
\par 
$(2)$ The Poincar\'{e} residue isomorphism {\rm (\ref{eqn:resd})} 
is compatible with $V$. 
\end{prop}
\begin{proof} 
(1): Set 
$E_{n}
:=(\eps^*_{X_{\os{\circ}{t}}/{\cal W}_n(\os{\circ}{t})}(E))_{{\cal W}_n(X_{\os{\circ}{t}})}$. 
By (\ref{prop:aco}) the morphism 
$V\col {\cal W}_n\wt{{\Om}}{}^i_{X_{\os{\circ}{t}}}\lo 
{\cal W}_{n+1}\wt{{\Om}}{}^i_{X_{\os{\circ}{t}}}$
preserve $P$'s. 
Hence the following morphism 
\begin{align*} 
{\rm id}\otimes V \col & 
E_n\otimes_{{\cal W}_n({\cal O}_{X_{\os{\circ}{t}}})^{\star}}
F_*({\cal W}_n{\Om}^i_{X_{\os{\circ}{t}}})^{\star} \\
& =
(E_{n+1}
\otimes_{{\cal W}_{n+1}({\cal O}_{X_{\os{\circ}{t}}})^{\star},R}
{\cal W}_n({\cal O}_{X_{\os{\circ}{t}}})^{\star})
\otimes_{{\cal W}_n({\cal O}_{X_{\os{\circ}{t}}})^{\star}} 
F_*({\cal W}_n{\Om}^i_{X_{\os{\circ}{t}}})^{\star} \\
& \lo 
E_{n+1}
\otimes_{{\cal W}_{n+1}({\cal O}_{X_{\os{\circ}{t}}})^{\star}}
({\cal W}_{n+1}{\Om}^i_{X_{\os{\circ}{t}}})^{\star}
\end{align*} 
preserves $P$'s. 
Furthermore it is clear that the morphisms 
\begin{align*} 
\Phi^{-1}\otimes {\rm id} \col & 
E_n\otimes_{{\cal W}_n({\cal O}_{X_{\os{\circ}{t}}})}
({\cal W}_n\wt{\Om}{}^i_{X_{\os{\circ}{t}}})^{\star} \lo (E_n)^{\sig}
\otimes_{{\cal W}_n({\cal O}_{X_{\os{\circ}{t}}})}
({\cal W}_n\wt{\Om}^i_{X_{\os{\circ}{t}}})
\end{align*} 
and 
\begin{align*} 
I\col  (E_n)^{\sig}
\otimes_{{\cal W}_n({\cal O}_{X_{\os{\circ}{t}}})}
({\cal W}_n\wt{\Om}^i_{X_{\os{\circ}{t}}})
&= F^{-1}(E_n)\otimes_{F^{-1}({\cal W}_n({\cal O}_{X_{\os{\circ}{t}}}))} 
({\cal W}_n\wt{\Om}^i_{X_{\os{\circ}{t}}}) \\
&=E_n\otimes_{{\cal W}_n({\cal O}_{X_{\os{\circ}{t}}})^{\star}}
F_*({\cal W}_n\wt{\Om}^i_{X_{\os{\circ}{t}}})^{\star}. 
\end{align*} 
preserve $P$'s. 
Because the morphism (\ref{ali:vesx}) is the composite morphism of 
$\Phi^{-1}\otimes {\rm id}$, $I$ and 
${\rm id}\otimes V$, 
the morphism (\ref{ali:vesx}) preserves $P$'s. 
\par 
(2): (2) follows from (\ref{prop:aco}). 
\end{proof}


\begin{prop}\label{prop:mfunc} 
Let $X'$ be an SNCL scheme over a log point 
$s'$ whose underlying scheme is 
the spectrum of a perfect field of characteristic $p>0$. 
For the commutative diagram {\rm (\ref{cd:tts})}  
and a commutative diagram 
\begin{equation*} 
\begin{CD} 
X_{\os{\circ}{t}} @>{g}>> X'_{\os{\circ}{t}{}'} \\ 
@VVV @VVV \\ 
s_{\os{\circ}{t}} @>{u}>> s'_{\os{\circ}{t}{}'} 
\end{CD} 
\tag{2.3.17.1}\label{cd:funct} 
\end{equation*} 
and a morphism of 
$E'\lo \os{\circ}{g}_{{\rm crys}*}(E)$, 
where $E'$ $($resp.~$E)$ is a flat coherent 
crystal of ${\cal O}_{\os{\circ}{X}{}'_{t'}/
{\cal W}_n(\os{\circ}{t}{}')}$-modules 
$($resp.~a flat coherent 
crystal of 
${\cal O}_{\os{\circ}{X}_t/{\cal W}_n(\os{\circ}{t}{})}$-modules$)$,   
there exists the following natural pull-back morphism 
\begin{align*} 
g^* \col & 
(\eps^*_{X'_{\os{\circ}{t}{}'}/{\cal W}_n(\os{\circ}{t}{}')}(E))_{
{\cal W}_n(X'_{\os{\circ}{t}{}'})}
\otimes_{{\cal W}_n({\cal O}_{X'_{\os{\circ}{t}{}'}})}
{\cal W}_n\wt{{\Om}}{}^{\bul}_{X'_{\os{\circ}{t}{}'}},P) 
\lo \tag{2.3.17.2}\label{eqn:wwlx} \\
& g_*(((\eps^*_{X_{\os{\circ}{t}}/{\cal W}_n(\os{\circ}{t})}(E))_{
{\cal W}_n(X_{\os{\circ}{t}})}\otimes_{{\cal W}_n({\cal O}_{X_{\os{\circ}{t}}})}
{\cal W}_n\wt{{\Om}}{}^{\bul}_{X_{\os{\circ}{t}}},P)),     
\end{align*} 
which is compatible with projections.  
\end{prop} 
\begin{proof} 
Because the proof the existence of $g^*$ is standard, we omit the proof. 
Because we have the following commutative diagrams 
\begin{equation*} 
\begin{CD} 
(\eps^*_{X'_{\os{\circ}{t}{}'}/{\cal W}_n(\os{\circ}{t}{}')}(E))_{
{\cal W}_n(X'_{\os{\circ}{t}{}'})}
\otimes_{{\cal W}_n({\cal O}_{X'_{\os{\circ}{t}{}'}})}
{\cal W}_n\wt{{\Om}}{}^{\bul}_{X'_{\os{\circ}{t}{}'}},P) 
@>{g^*}>> 
 \\ 
@V{{\bf p}}V{\bigcap}V \\
(\eps^*_{X'_{\os{\circ}{t}{}'}/{\cal W}_{n+1}(\os{\circ}{t}{}')}(E))_{
{\cal W}_{n+1}(X'_{\os{\circ}{t}{}'})}
\otimes_{{\cal W}_{n+1}({\cal O}_{X'_{\os{\circ}{t}{}'}})}
{\cal W}_{n+1}\wt{{\Om}}{}^{\bul}_{X'_{\os{\circ}{t}{}'}},P) 
@>{g^*}>> 
\end{CD} 
\end{equation*} 
\begin{equation*} 
\begin{CD} 
g_*((\eps^*_{X_{\os{\circ}{t}}/{\cal W}_n(\os{\circ}{t})}(E))_{
{\cal W}_n(X_{\os{\circ}{t}})}
\otimes_{{\cal W}_n({\cal O}_{X_{\os{\circ}{t}}})}
{\cal W}_n\wt{{\Om}}{}^{\bul}_{X_{\os{\circ}{t}}},P))\\
@V{\bigcap}V{\bf p}V \\ 
g_*((\eps^*_{X_{\os{\circ}{t}}/
{\cal W}_{n+1}(\os{\circ}{t})}(E))_{{\cal W}_{n+1}(X_{\os{\circ}{t}})}
\otimes_{{\cal W}_{n+1}({\cal O}_{X_{\os{\circ}{t}}})}
{\cal W}_{n+1}\wt{{\Om}}{}^{\bul}_{X_{\os{\circ}{t}}},P)) 
\end{CD} 
\end{equation*} 
and 
\begin{equation*} 
\begin{CD} 
(\eps^*_{X'_{\os{\circ}{t}{}'}/{\cal W}_n(\os{\circ}{t}{}')}(E))_{
{\cal W}_n(X'_{\os{\circ}{t}{}'})}
\otimes_{{\cal W}_n({\cal O}_{X'_{\os{\circ}{t}{}'}})}
{\cal W}_n\wt{{\Om}}{}^{\bul}_{X'_{\os{\circ}{t}{}'}},P) 
@>{g^*}>> 
 \\ 
@V{p}VV \\
(\eps^*_{X'_{\os{\circ}{t}{}'}/{\cal W}_{n+1}(\os{\circ}{t}{}')}(E))_{
{\cal W}_{n+1}(X'_{\os{\circ}{t}{}'})}
\otimes_{{\cal W}_{n+1}({\cal O}_{X'_{\os{\circ}{t}{}'}})}
{\cal W}_{n+1}\wt{{\Om}}{}^{\bul}_{X'_{\os{\circ}{t}{}'}},P) 
@>{g^*}>> 
\end{CD} 
\end{equation*} 
\begin{equation*} 
\begin{CD} 
g_*((\eps^*_{X_{\os{\circ}{t}}/{\cal W}_n(\os{\circ}{t})}(E))_{
{\cal W}_n(X_{\os{\circ}{t}})}
\otimes_{{\cal W}_n({\cal O}_{X_{\os{\circ}{t}}})}
{\cal W}_n\wt{{\Om}}{}^{\bul}_{X_{\os{\circ}{t}}},P))\\
@VV{p}V \\ 
g_*((\eps^*_{X_{\os{\circ}{t}}/
{\cal W}_{n+1}(\os{\circ}{t})}(E))_{{\cal W}_{n+1}(X_{\os{\circ}{t}})}
\otimes_{{\cal W}_{n+1}({\cal O}_{X_{\os{\circ}{t}}})}
{\cal W}_{n+1}\wt{{\Om}}{}^{\bul}_{X_{\os{\circ}{t}}},P))
\end{CD} 
\end{equation*} 
and because $p={\bf p}\circ R$, 
$g^*\circ R=R\circ g^*$. 
\end{proof}

\par 
In \cite[3.4]{msemi} Mokrane has 
defined a well-defined class 
$\theta_n \in {\cal W}_n\wt{{\Om}}_{X_{\os{\circ}{t}}}^1$
and he has obtained the following morphism 
\begin{equation*} 
\theta_n \wedge \col 
{\cal W}_n\wt{{\Om}}_{X_{\os{\circ}{t}}}^{j}
\lo {\cal W}_n\wt{{\Om}}_{X_{\os{\circ}{t}}}^{j+1}. 
\end{equation*} 
In \cite[(8.1.2)]{ndw} 
we have proved that the following diagram  
\begin{equation*} 
\begin{CD} 
{\cal W}_{n+1}\wt{{\Om}}_{X_{\os{\circ}{t}}}^{j+1} @>{R}>> 
{\cal W}_n\wt{{\Om}}_{X_{\os{\circ}{t}}}^{j+1} \\ 
@A{{\theta}_{n+1}\wedge}AA @AA{{\theta}_n\wedge}A \\ 
{\cal W}_{n+1}\wt{{\Om}}_{X_{\os{\circ}{t}}}^j 
@>{R}>> {\cal W}_n\wt{{\Om}}_{X_{\os{\circ}{t}}}^j 
\end{CD} 
\tag{2.3.17.3}\label{cd:futc} 
\end{equation*} 
is commutative. 

\par 
Let $E$ be a flat coherent log crystal of 
${\cal O}_{\os{\circ}{X}_t/{\cal W}_{n+1}(\os{\circ}{t})}$-modules. 
It is straightforward to generalize  
(\ref{cd:futc}) to the following commutative diagram by using 
the fact that the isomorphism 
${\cal W}_n({\cal O}_{X_{\os{\circ}{t}}})'\os{\sim}{\lo} 
{\cal W}_n({\cal O}_{X_{\os{\circ}{t}}})$ is compatible with 
projections: 
\begin{equation*} 
\begin{CD} 
(\eps^*_{X_{\os{\circ}{t}}/{\cal W}_{n+1}(\os{\circ}{t})}(E))_{
{\cal W}_{n+1}(X_{\os{\circ}{t}})}
\otimes_{{\cal W}_{n+1}({\cal O}_{X_{\os{\circ}{t}}})}
{\cal W}_{n+1}\wt{{\Om}}_{X_{\os{\circ}{t}}}^{j+1} 
@>{R:=R\otimes R}>>  \\ 
@A{{\theta}_{n+1}\wedge}AA \\
(\eps^*_{X_{\os{\circ}{t}}/{\cal W}_{n+1}(\os{\circ}{t})}
(E))_{
{\cal W}_{n+1}(X_{\os{\circ}{t}})}
\otimes_{{\cal W}_{n+1}({\cal O}_{X_{\os{\circ}{t}}})}
{\cal W}_{n+1}\wt{{\Om}}_{X_{\os{\circ}{t}}}^j 
@>{R:=R\otimes R}>> 
\end{CD} 
\tag{2.3.17.4}\label{cd:futnxc} 
\end{equation*} 
\begin{equation*} 
\begin{CD} 
(\eps^*_{X_{\os{\circ}{t}}/{\cal W}_n(\os{\circ}{t})}
(E))_{
{\cal W}_n(X_{\os{\circ}{t}})}
\otimes_{{\cal W}_n({\cal O}_{X_{\os{\circ}{t}}})}
{\cal W}_n\wt{{\Om}}_{X_{\os{\circ}{t}}}^{j+1}  \\
@AA{{\theta}_n\wedge}A \\
(\eps^*_{X_{\os{\circ}{t}}/{\cal W}_n(\os{\circ}{t})}
(E))_{
{\cal W}_n(X_{\os{\circ}{t}})}
\otimes_{{\cal W}_n({\cal O}_{X_{\os{\circ}{t}}})}
{\cal W}_n\wt{{\Om}}_{X_{\os{\circ}{t}}}^j. 
\end{CD}  
\end{equation*} 
For the commutative diagrams (\ref{cd:tts}) and (\ref{cd:funct}), 
we have the following commutative diagram 
\begin{equation*} 
\begin{CD} 
(\eps^*_{X'_{\os{\circ}{t}{}'}/{\cal W}_n(\os{\circ}{t'})}(E'))_{
{\cal W}_n(X'_{\os{\circ}{t}{}'})}
\otimes_{{\cal W}_n({\cal O}_{X'_{\os{\circ}{t}{}'}})}
{\cal W}_n\wt{{\Om}}_{X'_{\os{\circ}{t}{}'}}^{j+1}  @>{g^*}>> \\ 
@A{{\theta}'_n\wedge}AA  \\ 
(\eps^*_{X'_{\os{\circ}{t}{}'}/{\cal W}_n(\os{\circ}{t}{}')}(E'))_{
{\cal W}_n(X'_{\os{\circ}{t}{}'})}\otimes_{{\cal W}_n({\cal O}_{X'_{\os{\circ}{t}{}'}})}
{\cal W}_n\wt{{\Om}}_{X'_{\os{\circ}{t}{}'}}^j  @>{g^*}>> 
\end{CD} 
\tag{2.3.17.5}\label{cd:futct} 
\end{equation*} 
\begin{equation*} 
\begin{CD} 
g_*((\eps^*_{X_{\os{\circ}{t}}/{\cal W}_n(\os{\circ}{t})}(E))_{
{\cal W}_n(X_{\os{\circ}{t}})}
\otimes_{{\cal W}_n({\cal O}_{X_{\os{\circ}{t}}})}
{\cal W}_n\wt{{\Om}}_{X_{\os{\circ}{t}}}^{j+1}) \\ 
@AA{{\rm deg}(u){\theta}_n\wedge}A \\
g_*((\eps^*_{X_{\os{\circ}{t}}/{\cal W}_n(\os{\circ}{t})}(E))_{
{\cal W}_n(X_{\os{\circ}{t}})}
\otimes_{{\cal W}_n({\cal O}_{X_{\os{\circ}{t}}})}
{\cal W}_n\wt{{\Om}}_{X_{\os{\circ}{t}}}^j). 
\end{CD}  
\end{equation*} 
Here $\theta'_n\in {\cal W}_n\wt{{\Om}}_{X'_{\os{\circ}{t}{}'}}^1$ 
is the analogous differential form for the log point $s'_{\os{\circ}{t}{}'}$ 
to $\theta'_n\in {\cal W}_n\wt{{\Om}}_{X'_{\os{\circ}{t}{}'}}^1$ 
for the log point $s_{\os{\circ}{t}}$.  
\par 
Let $\ol{E}$ be a flat coherent log crystal of 
${\cal O}_{X_{\os{\circ}{t}}/{\cal W}_n(\os{\circ}{t})}$-modules. 
Set $\ol{E}_n:=\ol{E}_{{\cal W}_n(X_{\os{\circ}{t}})}$. 
Let 
$X_{\os{\circ}{t}}\os{\sus}{\lo} \ol{\cal P}_n$ be an immersion  
into a log smooth integral scheme over $\ol{{\cal W}_n(s_{\os{\circ}{t}})}$. 
Let $\ol{\mathfrak D}$ be the log PD-envelope of this immersion 
over $({\cal W}_n(\os{\circ}{t}),p{\cal W}_n,[~])$. 
Set ${\cal P}_n:=\ol{\cal P}_n\times_{\ol{{\cal W}_n(s)}}{\cal W}_n(s)$ 
and ${\mathfrak D}:=\ol{\mathfrak D}\times_{{\mathfrak D}(\ol{{\cal W}_n(s)})}
{\cal W}_n(s)$.  
Let $(\ol{\cal E},\ol{\nabla})$ be the 
coherent ${\cal O}_{\ol{\mathfrak D}_n}$-module with integrable connection  
obtained by $\ol{E}$. 
Set $({\cal E},\nabla):=(\ol{\cal E},\ol{\nabla})\otimes_{{\cal O}_{\ol{\mathfrak D}}}
{\cal O}_{\mathfrak D}$. 
In (\ref{eqn:fctdw}) we have constructed the following composite morphism$:$
\begin{align*} 
{\cal E}
\otimes_{{\cal O}_{{\cal P}_n}}
{{\Om}}{}^{\bul}_{{\cal P}_n/{\cal W}_n(\os{\circ}{t})} 
& \lo   
\ol{E}_{n}\otimes_{{\cal W}_n({\cal O}_Y)'}
\Om^{\bul}_{{\cal W}_n(X)/{\cal W}_n(\os{\circ}{t}),[~]} 
\tag{2.3.17.6}\label{eqn:fnidw} \\ 
&\lo \ol{E}_{n}
\otimes_{{\cal W}_n({\cal O}_{X_{\os{\circ}{t}}})}
{\cal W}_n\wt{{\Om}}^{\bul}_{X_{\os{\circ}{t}}}
\end{align*} 
in ${\rm C}^+(f^{-1}({\cal W}_n(\kap_t)))$.  
In fact, this is a quasi-isomorphism by (\ref{theo:ccttrw}). 
Let $E$ be a flat coherent crystal of 
${\cal O}_{\os{\circ}{X}_t/{\cal W}_n(\os{\circ}{t})}$-modules. 
Consider the case $\ol{E}:=\eps^*_{X_{\os{\circ}{t}}/{\cal W}_n(\os{\circ}{t})}(E)$. 

\begin{defi} 
Set 
\begin{align*} 
P_k(\ol{E}_{n}\otimes_{{\cal W}_n({\cal O}_{X_{\os{\circ}{t}}})'}
\Om^{\bul}_{{\cal W}_n(X_{\os{\circ}{t}})/{\cal W}_n(\os{\circ}{t}),[~]}) 
:=
\ol{E}_{n}\otimes_{{\cal W}_n({\cal O}_{X_{\os{\circ}{t}}})'}
P_k\Om^{\bul}_{{\cal W}_n(X_{\os{\circ}{t}})/{\cal W}_n(\os{\circ}{t}),[~]},  
\end{align*} 
where $\{P_k\Om^{\bul}_{{\cal W}_n(X_{\os{\circ}{t}})/{\cal W}_n(\os{\circ}{t}),[~]}\}_{k\in {\mab Z}}$ 
is the induced filtration by the filtration $P$ on 
$\Om^{\bul}_{{\cal W}_n(X_{\os{\circ}{t}})/{\cal W}_n(\os{\circ}{t})}$ ((\ref{eqn:pkdefpw})). 
\end{defi}

\parno 
It is clear that the first morphism in (\ref{eqn:fnidw}) 
induces the following filtered morphism 
\begin{align*} 
({\cal E}
\otimes_{{\cal O}_{{\cal P}_n}}
{{\Om}}{}^{\bul}_{{\cal P}_n/{\cal W}_n(\os{\circ}{t})},P) 
& \lo   
(\ol{E}_{n}\otimes_{{\cal W}_n({\cal O}_{X_{\os{\circ}{t}}})'}
\Om^{\bul}_{{\cal W}_n(X_{\os{\circ}{t}})/{\cal W}_n(\os{\circ}{t}),[~]},P).  
\tag{2.3.18.1}\label{eqn:fnfdw} 
\end{align*} 
By the note after (\ref{eqn:lttmwy}), the second morphism in (\ref{eqn:fnidw}) 
induces the following filtered morphism 
\begin{align*} 
(\ol{E}_{n}\otimes_{{\cal W}_n({\cal O}_{X_{\os{\circ}{t}}})'}
\Om^{\bul}_{{\cal W}_n(X)/{\cal W}_n(\os{\circ}{t}),[~]},P)  
&\lo (\ol{E}_n\otimes_{{\cal W}_n({\cal O}_{X_{\os{\circ}{t}}})}
{\cal W}_n\wt{{\Om}}^{\bul}_{X_{\os{\circ}{t}}},P).
\tag{2.3.18.2}\label{eqn:fnkdw} 
\end{align*} 
Hence the morphism (\ref{eqn:fnidw})  
induces the following filtered morphism 
\begin{align*} 
({\cal E}
\otimes_{{\cal O}_{{\cal P}_n}}
{\Om}^{\bul}_{{\cal P}_n/{\cal W}_n(\os{\circ}{t})},P)  
\lo &
(\ol{E}_n\otimes_{{\cal W}_n({\cal O}_{X_{\os{\circ}{t}}})}
{\cal W}_n\wt{{\Om}}{}^{\bul}_{X_{\os{\circ}{t}}},P)
\tag{2.3.18.3}\label{eqn:fbidw} \\
\end{align*} 
in ${\rm C}^+{\rm F}(f^{-1}({\cal W}_n))$. 
The morphism (\ref{eqn:fbidw}) 
is a filtered quasi-isomorphism by (\ref{eqn:eoppd}), (\ref{eqn:pwa})  
and (\ref{eqn:ywnnny}) 
(or Etesse's comparison theorem \cite[II (2.1)]{et}
(=(\ref{theo:ccrw}) for the case of the trivial log structure).

\par 
For the proof of (\ref{theo:csoncrdw}) below, 
we need the following, 
which is a variant of a coefficient version of 
an SNCL version of \cite[(4.22) (2)]{nh3}:  

\begin{lemm}\label{lemm:nefi}   
Assume that $\os{\circ}{X}$ is affine.  
Let $({\cal E},\nabla)$ be as above. 
Then the following hold$:$ 
\par 
$(1)$ Assume that we are given the commutative diagram {\rm (\ref{cd:imtpol})}.
The morphism {\rm (\ref{eqn:fbidw})} 
is functorial in the following sense$:$
for an affine SNCL scheme $X_i$ over $s_i$ $(i=1,2)$, 
for a morphism ${\cal W}_n(X_{i,\os{\circ}{t}_i})\lo \ol{\cal P}_{i,n}$ 
to a log smooth integral scheme over $\ol{{\cal W}_n(s_{i,\os{\circ}{t}_i})}$ 
such that the composite morphism 
$X_{i,\os{\circ}{t}_i} \os{\sus}{\lo} {\cal W}_n(X_{i,\os{\circ}{t}_i})\lo \ol{\cal P}_{i,n}$ 
is an immersion 
and for the log PD-envelope $\ol{\mathfrak D}_{i,n}$ $(i=1,2)$ of 
the immersion $X_{i,\os{\circ}{t}_i} \os{\sus}{\lo} \ol{\cal P}_{i,n}$ 
over $({\cal W}_n(s_{i,\os{\circ}{t}_i}),p{\cal W}_n,[~])$ 
and for the following commutative diagram 
\begin{equation*} 
\begin{CD} 
{\cal W}_n(X_{1,\os{\circ}{t}_1}) @>>> \ol{\mathfrak D}_{1,n}  \\ 
@V{g_n}VV @VV{\ol{g}}V \\ 
{\cal W}_n(X_{2,\os{\circ}{t}_2}) @>>>\ol{\mathfrak D}_{2,n}
\end{CD} 
\tag{2.3.19.1}\label{cd:yyppd}
\end{equation*} 
of log schemes over the following commutative diagram 
\begin{equation*} 
\begin{CD} 
{\cal W}_n(s_{1,\os{\circ}{t}_1}) @>>> \ol{{\cal W}_n(s_{1,\os{\circ}{t}_1})} \\ 
@VVV @VVV \\ 
{\cal W}_n(s_{2,\os{\circ}{t}_2}) @>>> \ol{{\cal W}_n(s_{2,\os{\circ}{t}_2})}
\end{CD} 
\tag{2.3.19.2}\label{cd:imsol}
\end{equation*}  
and 
for the log PD-envelope ${\mathfrak D}_{i,n}$ 
$(i=1,2)$ of the immersion $X_i \os{\sus}{\lo} {\cal P}_{i,n}$ 
$({\cal P}_{i,n}:=
\ol{\cal P}_{i,n}\times_{\ol{{\cal W}_n(s_{i,\os{\circ}{t}_i})}}{\cal W}_n(s_{i,\os{\circ}{t}_i})$ 
over $({\cal W}_n(s_{i,\os{\circ}{t}_i}),
p{\cal O}_{{\cal W}_n(s_{i,\os{\circ}{t_i}})},[~])$ 
and for 
a flat coherent crystal $E_i$ of 
${\cal O}_{\os{\circ}{X}_{i,t_i}/{\cal W}_n(\os{\circ}{t}_i)}$-modules 
with a morphism 
$\os{\circ}{g}{}^*_{1{\rm crys}}(E_2)\lo E_1$  
and the corresponding 
integrable connection $(\ol{\cal E}_i,\nabla_i)$ to 
$\eps^*_{X_{i,\os{\circ}{t}_i}/{\cal W}_n(\os{\circ}{t}_i)}(E_i)$, 
the following diagram is commutative$:$
\begin{equation*}  
\begin{CD} 
({\cal E}_1{\otimes}_{{\cal O}_{{\cal P}^{\rm ex}_{1,n}}}
\Om^{\bul}_{{\cal P}^{\rm ex}_{1,n}/{\cal W}_n(\os{\circ}{t}_1)},P)
@>{\sim}>> ({\cal E}_1{\otimes}_{{\cal O}_{{\cal P}^{\rm ex}_{1,n}}}
{\cal H}^{\bul}({\Om}^{*}_{{\cal P}^{\rm ex}_{1,n}/{\cal W}_n(\os{\circ}{t}_1)}),P) \\ 
@A{g^*}AA @AA{g^*_n}A \\ 
(g^*({\cal E}_2
{\otimes}_{{\cal O}_{{\cal P}^{\rm ex}_{2,n}}}
{\Om}^{\bul}_{{\cal P}^{\rm ex}_{2,n}/{\cal W}_n(\os{\circ}{t}_2)}),P)
@>{\sim}>> 
(g^*({\cal E}_2{\otimes}_{{\cal O}_{{\cal P}^{\rm ex}_{2,n}}}
{\cal H}^{\bul}({\Om}^{*}_{{\cal P}^{\rm ex}_{2,n}/{\cal W}_n(\os{\circ}{t}_2)})),P),   
\end{CD} 
\tag{2.3.19.3}\label{cd:odzqph}
\end{equation*} 
where ${\cal E}_i:=
\ol{\cal E}_i\otimes_{{\cal O}_{{\mathfrak D}(\ol{{\cal W}_n(s_{i,\os{\circ}{t}_i})})}}
{\cal O}_{{\cal W}_n(s_{i,\os{\circ}{t}_i})}$ 
and $g\col {\mathfrak D}_{1,n}\lo {\mathfrak D}_{2,n}$ is the induced morphism 
by $\ol{g}$. 
In the $p$-adic case over 
${\cal W}(s):=({\rm Spf}({\cal W}),
{\mab N}\oplus {\cal W}^*)$ and the case where 
$E$ is a coherent flat crystal of ${\cal O}_{\os{\circ}{X}_t/{\cal W}(\os{\circ}{t})}$-modules, 
the morphism {\rm (\ref{eqn:fbidw})} is compatible 
with the projections. 
\par 
$(2)$ The following diagram is commutative$:$
\begin{equation*} 
\begin{CD} 
{\cal E}
\otimes_{{\cal O}_{{\cal P}^{\rm ex}_n}}
{{\Om}}{}^{i+1}_{{\cal P}^{\rm ex}_n/{\cal W}_n(\os{\circ}{t})} 
@>{\sim}>>
(\eps^*_{X_{\os{\circ}{t}}/{\cal W}_n(\os{\circ}{t})}
(E))_{
{\cal W}_n(X_{\os{\circ}{t}})}\otimes_{{\cal W}_n({\cal O}_{X_{\os{\circ}{t}}})}
{\cal W}_n\wt{{\Om}}{}^{i+1}_{X_{\os{\circ}{t}}}\\ 
@A{\theta_{{\cal P}^{\rm ex}_n}\wedge}AA 
@AA{\theta_n\wedge}A \\
{\cal E}\otimes_{{\cal O}_{{\cal P}^{\rm ex}_n}}
{\Om}{}^i_{{\cal P}^{\rm ex}_n/{\cal W}_n(\os{\circ}{t})}  
@>{\sim}>>
(\eps^*_{X_{\os{\circ}{t}}/{\cal W}_n(\os{\circ}{t})}
(E))_{{\cal W}_n(X_{\os{\circ}{t}})}
\otimes_{{\cal W}_n({\cal O}_{X_{\os{\circ}{t}}})}
{\cal W}_n\wt{{\Om}}{}^i_{X_{\os{\circ}{t}}}
\end{CD} 
\tag{2.3.19.4}\label{cd:fedw} 
\end{equation*} 
The commutative diagram {\rm (\ref{cd:fedw})} 
is compatible with projections. 
\end{lemm} 
\begin{proof} 
(1): Let $g_{1{\rm crys}} 
\col (X_{1,\os{\circ}{t}_1}/{\cal W}_n(\os{\circ}{t}_1))_{\rm crys}\lo 
(X_{2,\os{\circ}{t}_2}/{\cal W}_n(\os{\circ}{t}_2))_{\rm crys}$ be 
the induced morphism 
by $g_1\col X_{1,\os{\circ}{t}_1}\lo X_{2,\os{\circ}{t}_2}$.  
Because we are given a morphism 
$\os{\circ}{g}{}^*_{1{\rm crys}}(E_2)\lo E_1$, 
we have a morphism 
$g^*_{1{\rm crys}} \eps^*_{X_{2,\os{\circ}{t}_2}/{\cal W}_n(\os{\circ}{t}_2)}(E_2)\lo 
\eps^*_{X_{1,\os{\circ}{t}_1}/{\cal W}_n(\os{\circ}{t}_1)}(E_1)$.   
The rest of the proof is the same as that of (\ref{lemm:ntetfi}).  
\par 
(2): (2) follows immediately from the construction of 
the morphism (\ref{eqn:fbidw}) and by the note after (\ref{eqn:lttmwy}). 
The compatibility of (\ref{cd:fedw}) with respect to the projections follows 
from (\ref{cd:futc}). 
\end{proof}

\par  
Now let us come back to the beginning of this section. 

\begin{theo}[{\bf Comparison theorem}]\label{theo:csssh}
$(1)$ Assume that $N\not=\infty$. 
Assume that $X_{\bul \leq N,\os{\circ}{t}}$ has 
an affine $N$-truncated simplicial open covering. 
Then there exists a canonical filtered isomorphism 
\begin{equation*} 
(\wt{R}u_{X_{\bul\leq N,\os{\circ}{t}}/{\cal W}_n(\os{\circ}{t})*}(\ol{E}{}^{\bul \leq N}),P) 
\os{\sim}{\lo} 
(\ol{E}{}^{\bul \leq N}_{n}
\otimes_{{\cal W}_n({\cal O}_{X_{\bul \leq N,\os{\circ}{t}}})}
{\cal W}_n\wt{{\Om}}^{\bul}_{X_{\bul \leq N,\os{\circ}{t}}},P) 
\tag{2.3.20.1}\label{eqn:ywony}
\end{equation*} 
in ${\rm D}^+{\rm F}(f^{-1}_{\bul \leq N}({\cal W}_n(\kap_t)))$.
The isomorphisms $(\ref{eqn:ywony})$ for $n$'s are 
compatible with two projections of 
both hand sides on $(\ref{eqn:ywony})$. 
\par 
$(2)$ 
Let the notations be as in {\rm (\ref{cd:tts})}.  
The isomorphism {\rm (\ref{eqn:ywony})} is contravariantly functorial 
with respect to a morphism 
$g_{\bul \leq N}\col X_{\bul \leq N}\lo Y_{\bul \leq N}$ of 
$N$-truncated simplicial SNCL schemes over the morphism 
${\cal W}_n(s_{\os{\circ}{t}})\lo {\cal W}_n(s'_{\os{\circ}{t}{}'})$ 
satisfying the condition {\rm (\ref{cd:xygxy})} 
for the case $S=s$, $S'=s'$, $T={\cal W}_n(t)$ and $T'={\cal W}_n(t')$
and 
$\os{\circ}{g}{}^*_{\rm crys}(F^{\bul \leq N})\lo E^{\bul \leq N}$, 
where $Y_{\bul \leq N}$ and $F^{\bul \leq N}$ are similar objects to 
$X_{\bul \leq N}$ and $E^{\bul \leq N}$, respectively. 
$($The morphism $g_{\bul \leq N}$ is the following morphism 
fitting into the following commutative diagram
\begin{equation*} 
\begin{CD} 
X_{\bul \leq N,\os{\circ}{t}_0} @>{g_{\bul \leq N}}>> Y_{\bul \leq N,\os{\circ}{t}{}'_0}\\
@VVV @VVV \\ 
s_{\os{\circ}{t}} @>>> s'_{\os{\circ}{t}{}'} \\ 
@V{\bigcap}VV @VV{\bigcap}V \\ 
{\cal W}_n(s_{\os{\circ}{t}}) @>{u}>> {\cal W}_n(s'_{\os{\circ}{t}{}'})
\end{CD}
\tag{2.3.20.2}\label{eqn:xdxdwuss}
\end{equation*} 
of $N$-truncated simplicial SNCL schemes 
over $s_{\os{\circ}{t}}$ and $s'_{\os{\circ}{t}{}'}$ 
such that $X_{\bul \leq N,\os{\circ}{t}}$ and $Y_{\bul \leq N,\os{\circ}{t}{}'}$ 
have the disjoint unions $X'_{\bul \leq N,\os{\circ}{t}}$ and $Y'_{\bul \leq N,\os{\circ}{t}{}'}$ 
of the members of affine $N$-truncated simplicial open coverings of 
$X_{\bul \leq N,\os{\circ}{t}}$ and $Y_{\bul \leq N,\os{\circ}{t}{}'}$, 
respectively, such that there exists a morphism 
$g'_{\bul \leq N}\col X'_{\bul \leq N,\os{\circ}{t}} \lo 
Y'_{\bul \leq N,\os{\circ}{t}{}'}$ 
fitting into the following commutative diagram 
\begin{equation*} 
\begin{CD} 
X'_{\bul \leq N,\os{\circ}{t}} @>{g'_{\bul \leq N}}>> Y'_{\bul \leq N,\os{\circ}{t}{}'} \\
@VVV @VVV \\ 
X_{\bul \leq N,\os{\circ}{t}} @>{g_{\bul \leq N}}>> Y_{\bul \leq N,\os{\circ}{t}{}'}.) 
\end{CD}
\tag{2.3.20.3}\label{cd:xwygxy}
\end{equation*} 
\end{theo} 
\begin{proof} 
(1): Let the notations be as in the proof of (\ref{theo:ccttrw}) 
for the case $Y_{\bul \leq N,\os{\circ}{t}}=X_{\bul \leq N,\os{\circ}{t}}$. 
In this proof  
we denote ${\cal Q}_{\bul \leq N,\bul}$ in (\ref{theo:ccttrw}) 
by ${\cal P}_{\bul \leq N,\bul}$. 
By (\ref{eqn:fbidw}) and (\ref{lemm:nefi}) (and the argument before (\ref{lemm:nefi})), 
we have the following morphism 
of filtered complexes: 
\begin{align*} 
({\cal E}^{\bul \leq N,\bul}
{\otimes}_{{\cal O}_{{\cal P}^{\rm ex}_{\bul \leq N,\bul}}}
{\Om}^{\bul}_{{\cal P}^{\rm ex}_{\bul \leq N,\bul}/{\cal W}_n(\os{\circ}{t})},P)  
& \lo (E^{\bul \leq N,\bul}_{n} 
{\otimes}_{{\cal O}_{{\cal P}^{\rm ex}_{\bul \leq N,\bul}}}
{\cal H}^{\bul}({\Om}^{*}_{{\cal P}^{\rm ex}_{\bul \leq N,\bul}/{\cal W}_n(\os{\circ}{t})}),P) 
\tag{2.3.20.4}\label{ali:pno}\\ 
&=(E^{\bul \leq N,\bul}_{n} \otimes_{{\cal W}_n{(\cal O}_{X_{\bul \leq N,\bul,\os{\circ}{t}}})}
{\cal W}_n\wt{{\Om}}^{\bul}_{X_{\bul \leq N,\bul,\os{\circ}{t}}},P)  
\end{align*}  
in ${\rm C}^+{\rm F}
(f^{-1}_{\bul \leq N,\bul}({\cal W}_n(\kap_t)))$. 
As in the proof of (\ref{theo:ccttrw}), we obtain the filtered isomorphism 
(\ref{eqn:ywony}) by the argument before (\ref{lemm:nefi}). 
\par 
(2): By using (\ref{cd:odzqph}), 
the proof is almost the same as that of (\ref{theo:ccrw}) (2). 
We leave the detailed proof to the reader. 
\end{proof}

\par 
Let $\star$ be a positive integer $n$ or nothing. 
Let $E^{\bul \leq N}$ be a flat coherent log crystal of 
${\cal O}_{\os{\circ}{X}_{\bul \leq N,t}/{\cal W}_{\star}(\os{\circ}{t})}$-modules. 
Next we give the definition of a zariskian filtered complex 
$({\cal W}_{\star}A_{X_{\bul \leq N,\os{\circ}{t}}}(E^{\bul \leq N}),P)$ as follows. 
(This is a generalization of the zariskian $p$-adic filtered 
Hyodo-Mokrane-Steenbrink complex 
in \cite[3.15]{msemi} (see also \cite[(9.9)]{ndw}).) 
\par 
Set 
\begin{align*} 
&{\cal W}_{\star}A_{X_{\bul \leq N,\os{\circ}{t}}}(E^{\bul \leq N})^{ij} := 
((\eps^*_{X_{\bul\leq N,\os{\circ}{t}}/{\cal W}_{\star}(\os{\circ}{t})}(E^{\bul \leq N}))_
{{\cal W}_{\star}(X_{\bul \leq N,\os{\circ}{t}})}
\otimes_{{\cal W}_{\star}({\cal O}_{X_{\bul \leq N,\os{\circ}{t}}})}
{\cal W}_{\star}\wt{{\Om}}_{X_{\bul \leq N,\os{\circ}{t}}}^{i+j+1})/ 
\tag{2.3.20.5}\label{eqn:wlpj} \\
& P_j(
(\eps^*_{X_{\bul \leq N,\os{\circ}{t}}/{\cal W}_{\star}(\os{\circ}{t})}(E^{\bul \leq N}))_
{{\cal W}_{\star}(X_{\bul \leq N,\os{\circ}{t}})}
\otimes_{{\cal W}_{\star}({\cal O}_{X_{\bul \leq N,\os{\circ}{t}}})}
{\cal W}_{\star}\wt{{\Om}}_{X_{\bul \leq N,\os{\circ}{t}}}^{i+j+1}) 
\quad (i,j \in {\mab N}).  
\end{align*}  
As in \cite[(5.17)]{fup} and (\ref{cd:lccbd}),  
we consider 
the following boundary morphisms of the double complexes 
${\cal W}_{\star}A_{X_{\bul \leq N,\os{\circ}{t}}}(E^{\bul \leq N})$: 
\begin{equation*}
\begin{CD}
{\cal W}_{\star}A_{X_{\bul \leq N,\os{\circ}{t}}}(E^{\bul \leq N})^{i,j+1}  @.  \\ 
@A{\theta_{\star}\wedge}AA  @. \\
{\cal W}_{\star}A_{X_{\bul \leq N,\os{\circ}{t}}}(E^{\bul \leq N})^{ij}
@>{-\nabla}>> 
{\cal W}_{\star}A_{X_{\bul \leq N,\os{\circ}{t}}}(E^{\bul \leq N})^{i+1,j}.\\
\end{CD}
\tag{2.3.20.6}\label{cd:lcwnbd} 
\end{equation*}   
Then ${\cal W}_{\star}A_{X_{\bul \leq N,\os{\circ}{t}}}(E^{\bul \leq N})^{\bul \bul}$ 
becomes 
an $N$-truncated cosimplicial double complex 
of $f^{-1}_{\bul \leq N}({\cal W}_{\star}(\kap_t))$-modules. 
Let ${\cal W}_{\star}A_{X_{\bul \leq N,\os{\circ}{t}}}(E^{\bul \leq N})$ 
be the single complex of  
${\cal W}_{\star}A_{X_{\bul \leq N,\os{\circ}{t}}}(E^{\bul \leq N})^{\bul \bul}$. 
We endow ${\cal W}_{\star}A_{X_{\bul \leq N,\os{\circ}{t}}}(E^{\bul \leq N})$ 
with a filtration $P$ as follows: 
\begin{align*} 
P_k{\cal W}_{\star}A_{X_{\bul \leq N,\os{\circ}{t}}}(E^{\bul \leq N})^{ij}&:= 
{\rm Im} (P_{2j+k+1}
(((\eps^*_{X_{\bul \leq N,\os{\circ}{t}}/{\cal W}_{\star}(\os{\circ}{t})}
(E^{\bul \leq N}))_{{\cal W}_{\star}(X_{\bul \leq N,\os{\circ}{t}})}
\tag{2.3.20.7}\label{eqn:dbwlad}\\
& \otimes_{{\cal W}_{\star}({\cal O}_{X_{\bul \leq N,\os{\circ}{t}}})}
{\cal W}_{\star}\wt{{\Om}}^{i+j+1}_{X_{\bul \leq N,\os{\circ}{t}}}) \lo
{\cal W}_{\star}A_{X_{\bul \leq N,\os{\circ}{t}}}(E^{\bul \leq N})^{ij}). 
\end{align*} 
Let $({\cal W}_{\star}A_{X_{\bul \leq N,\os{\circ}{t}}}(E^{\bul \leq N}),P)$ 
be the filtered single complex of the filtered double complex 
$({\cal W}_{\star}A_{X_{\bul \leq N,\os{\circ}{t}}}(E^{\bul \leq N})^{\bul \bul},P)$.  

\begin{defi}\label{defi:fcwaxp}
We call the filtered complex 
$({\cal W}_{\star}A_{X_{\bul \leq N,\os{\circ}{t}}}(E^{\bul \leq N}),P)$ 
the {\it filtered Hyodo-Mokrane-Steenbrink complex} of   
$E^{\bul \leq N}$. 
In the case where $E^{\bul \leq N}
={\cal O}_{\os{\circ}{X}_{\bul \leq N,t}/{\cal W}_{\star}(\os{\circ}{t})}$, 
we denote $({\cal W}_{\star}A_{X_{\bul \leq N,\os{\circ}{t}}}(E^{\bul \leq N}),P)$ 
by $({\cal W}_{\star}A_{X_{\bul \leq N,\os{\circ}{t}}},P)$ and we call it 
the {\it filtered Hyodo-Mokrane-Steenbrink complex} of   
$X_{\bul \leq N,\os{\circ}{t}}/\os{\circ}{t}$.  
(In the case $N=0$ and $t=s$, in \cite{msemi}, 
Mokrane has called $({\cal W}_{\star}A_{X_{0,\os{\circ}{t}}}(E^0),P)$ 
the filtered Steenbrink-Hyodo complex.) 
\end{defi}
 
\begin{prop}\label{prop:qidw}
The following morphism 
\begin{align*} 
\theta_{\star} \wedge \col 
(\eps^*_{X_{\bul \leq N,\os{\circ}{t}}/{\cal W}_{\star}(\os{\circ}{t})}
(E^{\bul \leq N}))_{{\cal W}_{\star}(X_{\bul \leq N,\os{\circ}{t}})}
\otimes_{{\cal W}_{\star}({\cal O}_{X_{\bul \leq N,\os{\circ}{t}}})}
{\cal W}_{\star}{\Om}^{\bul}_{X_{\bul \leq N,\os{\circ}{t}}}
\lo {\cal W}_{\star}A_{X_{\bul \leq N,\os{\circ}{t}}}(E^{\bul \leq N})  
\tag{2.3.22.1}\label{eqn:wlxa} 
\end{align*} 
is a quasi-isomorphism.  
\end{prop} 
\begin{proof} 
It suffices to prove that  
the following morphism 
\begin{align*} 
\theta_{\star} \wedge \col 
(\eps^*_{X_{m,\os{\circ}{t}}/{\cal W}_{\star}(\os{\circ}{t})}
(E^{m}))_{{\cal W}_{\star}(X_{m,\os{\circ}{t}})}
\otimes_{{\cal W}_{\star}({\cal O}_{X_{m,\os{\circ}{t}}})}
{\cal W}_{\star}{\Om}^i_{X_{m,\os{\circ}{t}}}
\lo {\cal W}_{\star}A_{X_{m,\os{\circ}{t}}}(E^m)^{i\bul}  
\tag{2.3.22.2}\label{eqn:wlmxa} 
\end{align*} 
is a quasi-isomorphism for $0\leq m\leq N$.  
We may assume that $E^{m}={\cal O}_{\os{\circ}{X}_{m,t}/{\cal W}(\os{\circ}{t})}$. 
By \cite[3.15]{msemi} and \cite[(6.28) (9), (6.29) (1)]{ndw} 
(see also the proof of (\ref{prop:tefc})), 
the following morphism
\begin{align*} 
\theta_{\star} \wedge \col 
{\cal W}_{\star}{\Om}^i_{X_{m,\os{\circ}{t}}}
\lo {\cal W}_{\star}A_{X_{m,\os{\circ}{t}}}^{i \bul}
\end{align*} 
is a quasi-isomorphism. 
Hence the morphism (\ref{eqn:wlmxa}) is a quasi-isomorphism.  
\end{proof}


\begin{theo}[{\bf Contravariant functoriality}]\label{theo:ctrw} 
$(1)$ Let the notations be as in the beginning of {\rm (\ref{theo:csssh})}. 
Assume that 
$\deg(u)_x$ is not divisible by $p$ for any point  
$x \in \os{\circ}{t}$. 
Then $g_{\bul \leq N}$ induces the following 
well-defined pull-back morphism 
\begin{align*}  
g_{\bul \leq N}^* \col &
({\cal W}_nA_{Y_{\bul\leq N,\os{\circ}{t}{}'}}(F^{\bul \leq N}),P) 
\lo 
Rg_{\bul \leq N*}(({\cal W}_nA_{X_{\bul\leq N,\os{\circ}{t}}}(E^{\bul \leq N}),P))
\tag{2.3.23.1}\label{eqn:fzdwd}
\end{align*} 
fitting into the following commutative diagram$:$
\begin{equation*} 
\begin{CD}
{\cal W}_nA_{Y_{\bul\leq N,{\os{\circ}{t}{}'}}}(F^{\bul \leq N})
@>{g_{\bul \leq N}^*}>> \\ 
@A{\theta_n \wedge}A{\simeq}A \\
Ru_{Y_{\bul \leq N,\os{\circ}{t}{}'}/{\cal W}_n(s'_{\os{\circ}{t}{}'})*}
(\eps^*_{Y_{\bul \leq N,\os{\circ}{t}{}'}/
{\cal W}_n(s'_{\os{\circ}{t}{}'})}(F^{\bul \leq N}))@>{g_{\bul \leq N}^*}>>
\end{CD}
\tag{2.3.23.2}\label{cd:psccz} 
\end{equation*} 
\begin{equation*} 
\begin{CD}
Rg_{\bul \leq N*}({\cal W}_nA_{X_{\bul\leq N,\os{\circ}{t}}}(E^{\bul \leq N}))\\ 
@A{Rg_{\bul \leq N*}(\theta_n \wedge)}A{\simeq}A \\
Rg_{\bul \leq N*}
Ru_{X_{\bul \leq N,\os{\circ}{t}}/{\cal W}_n(s_{\os{\circ}{t}})*}
(\eps^*_{X_{\bul \leq N,\os{\circ}{t}}/{\cal W}_n(s_{\os{\circ}{t}})}
(E^{\bul \leq N})).
\end{CD}
\end{equation*}
\par 
$(2)$ Let the notations be as in $(1)$. 
Let $s''$, $t''$, $S''={\cal W}_n(s'')$ and $T''={\cal W}_n(t'')$ 
be similar objects to $s'$, $t'$, $S'={\cal W}_n(s')$ and $T'={\cal W}_n(t')$, 
respectively.  
Let $h_{\bul \leq N}\col Y_{\bul \leq N,\os{\circ}{t}{}'}\lo 
Z_{\bul \leq N,\os{\circ}{t}{}''}$ be a similar morphism 
over $s'_{\os{\circ}{t}{}'}\lo s''_{\os{\circ}{t}{}''}$ 
to $g_{\bul \leq N}\col X_{\bul \leq N,\os{\circ}{t}}\lo Y_{\bul \leq N,\os{\circ}{t}{}'}$.  
Let 
\begin{align*} 
\os{\circ}{h}{}^*_{\bul \leq N,{\rm crys}}(G^{\bul \leq N})\lo F^{\bul \leq N} 
\tag{2.3.23.3}\label{ali:rssntpm}
\end{align*} 
be as in {\rm (\ref{ali:gnfe})}. Let $v\col {\cal W}_n(s'_{\os{\circ}{t}{}'})\lo 
{\cal W}_n(s''_{\os{\circ}{t}{}''})$ be a similar morphism to 
$u\col {\cal W}_n(s_{\os{\circ}{t}{}})\lo {\cal W}_n(s'_{\os{\circ}{t}{}'})$.  
Assume that 
$\deg(v)_y$ is not divisible by $p$ for any point  $y \in \os{\circ}{t}{}'$. 
Then 
\begin{align*} 
(h_{\bul \leq N}\circ g_{\bul \leq N})^* =
Rh_{\bul \leq N*}(g_{\bul \leq N}^*)\circ h_{\bul \leq N}^*     
\col &({\cal W}_nA_{Z_{\bul\leq N,\os{\circ}{t}{}''}}(G^{\bul \leq N})),P) \lo 
\tag{2.3.23.4}\label{ali:pwdpp} \\ 
& Rh_{\bul \leq N*}
Rg_{\bul \leq N*}({\cal W}_nA_{X_{\bul\leq N,\os{\circ}{t}}}(E^{\bul \leq N}),P)  \\
& =R(h_{\bul \leq N}\circ 
g_{\bul \leq N})_*({\cal W}_nA_{X_{\bul\leq N,\os{\circ}{t}}}(E^{\bul \leq N}),P).
\end{align*}  
\par 
$(3)$ \begin{equation*} 
{\rm id}_{X_{\bul \leq N,\os{\circ}{t}}}^*={\rm id} 
\col ({\cal W}_nA_{X_{\bul\leq N,\os{\circ}{t}}}(E^{\bul \leq N}),P)
\lo ({\cal W}_nA_{X_{\bul\leq N,\os{\circ}{t}}}(E^{\bul \leq N}),P). 
\tag{2.3.23.5}\label{eqn:fzidxd}
\end{equation*} 
\par 
$(4)$ 
Let the notations be as in {\rm (2)} without assuming the 
nondivisibility condition, 
where $S={\cal W}(s)$, $S'={\cal W}(s')$, 
$T={\cal W}(t)$ and $T={\cal W}(t')$.  
Let $E^{\bul \leq N}$ and  $F^{\bul \leq N}$ 
be flat quasi-coherent crystals of 
${\cal O}_{\os{\circ}{X}_{\bul \leq N}/{\cal W}(\os{\circ}{t})}$-modules  
and 
${\cal O}_{\os{\circ}{Y}_{\bul \leq N}/{\cal W}(\os{\circ}{t}{}')}$-modules, 
respectively.     
Let 
\begin{align*} 
\os{\circ}{g}{}^*_{\bul \leq N,{\rm crys}}(F^{\bul \leq N})
\lo 
E^{\bul \leq N} 
\end{align*} 
be a morphism of 
${\cal O}_{\os{\circ}{X}_{\bul \leq N,t}/{\cal W}(\os{\circ}{t})}$-modules. 
Set 
\begin{align*} 
&(\eps^*_{X_{\bul \leq N,\os{\circ}{t}}/{\cal W}(\os{\circ}{t})}
(E^{\bul \leq N}))_{{\cal W}(X_{\bul \leq N,\os{\circ}{t}})}
\otimes_{{\cal W}({\cal O}_{X_{\bul \leq N,\os{\circ}{t}})}}
{\cal W}\wt{\Om}_{X_{\bul \leq N,\os{\circ}{t}}}^j \\
&:=\vpl_n((\eps^*_{X_{\bul \leq N,\os{\circ}{t}}/{\cal W}_n(\os{\circ}{t})}
(E^{\bul \leq N}))_{{\cal W}_n(X_{\bul \leq N,\os{\circ}{t}})}
\otimes_{{\cal W}_n({\cal O}_{X_{\bul \leq N,\os{\circ}{t}})}}
{\cal W}_n\wt{\Om}_{X_{\bul \leq N,\os{\circ}{t}}}^j)\\
&=
\vpl_n((\eps^*_{X_{\bul \leq N,\os{\circ}{t}}/{\cal W}_n(\os{\circ}{t})}
(E^{\bul \leq N}))_{{\cal W}_n(X_{\bul \leq N,\os{\circ}{t}})})
\otimes_{{\cal W}({\cal O}_{X_{\bul \leq N,\os{\circ}{t}})}}
{\cal W}\wt{\Om}_{X_{\bul \leq N,\os{\circ}{t}}}^j
\end{align*} 
and let 
\begin{align*} 
&(\eps^*_{Y_{\bul \leq N,\os{\circ}{t}{}'}/{\cal W}(\os{\circ}{t}{}')}
(F^{\bul \leq N}))_{{\cal W}(Y_{\bul \leq N,\os{\circ}{t}{}'})}
\otimes_{{\cal W}({\cal O}_{Y_{\bul \leq N,\os{\circ}{t}{}'})}}
{\cal W}\wt{\Om}_{Y_{\bul \leq N,\os{\circ}{t}{}'}}^j 
\end{align*} 
be the similar sheaf to 
$(\eps^*_{X_{\bul \leq N,\os{\circ}{t}}/{\cal W}(\os{\circ}{t})}
(E^{\bul \leq N}))_{{\cal W}(X_{\bul \leq N,\os{\circ}{t}})}
\otimes_{{\cal W}({\cal O}_{X_{\bul \leq N,\os{\circ}{t}})}}
{\cal W}\wt{\Om}_{X_{\bul \leq N,\os{\circ}{t}}}^j$. 
If the morphism 
\begin{align*} 
&g^*_{\bul \leq N} \col 
(\eps^*_{Y_{\bul \leq N,\os{\circ}{t}{}'}/{\cal W}(\os{\circ}{t}{}')}
(F^{\bul \leq N}))_{{\cal W}(Y_{\bul \leq N,\os{\circ}{t}{}'})}
\otimes_{{\cal W}({\cal O}_{Y_{\bul \leq N,\os{\circ}{t}{}'})}}
{\cal W}\wt{\Om}_{Y_{\bul \leq N,\os{\circ}{t}{}'}}^{i+j+1}/P_j 
\tag{2.3.23.6}\label{ali:odnl} \\
&\lo g_{\bul \leq N*}
((\eps^*_{X_{\bul \leq N,\os{\circ}{t}}/{\cal W}(\os{\circ}{t})}
(E^{\bul \leq N}))_{{\cal W}(X_{\bul \leq N,\os{\circ}{t}})}
\otimes_{{\cal W}({\cal O}_{X_{\bul \leq N,\os{\circ}{t}})}}
{\cal W}\wt{\Om}_{X_{\bul \leq N,\os{\circ}{t}}}^{i+j+1})/P_j 
\end{align*} 
is divisible by $p^{e_p(j+1)}$ for any $i,j\in {\mab N}$, 
$($Because  the sheaf ${\cal W}\wt{{\Om}}_{X_{\bul \leq N,\os{\circ}{t}}}^{i+j+1}/P_j$ is a 
sheaf of flat ${\cal W}(\os{\circ}{t})$-modules,  
the divisibility above is well-defined.$)$, 
then the conclusion of $(1)$ holds. 
The similar relation to {\rm (\ref{ali:pwdpp})} for the $g_{\bul \leq N}^*$ above holds. 
\par 
$(5)$ Let the notations be as in $(4)$ without assuming the 
divisibility condition. 
Then $g_{\bul \leq N}$ induces the following 
well-defined pull-back morphism 
\begin{align*}  
g_{\bul \leq N}^* \col &
({\cal W}A_{Y_{\bul\leq N,\os{\circ}{t}{}'}}(F^{\bul \leq N})\otimes_{\mab Z}{\mab Q},P) 
\lo 
Rg_{\bul \leq N*}(({\cal W}A_{X_{\bul\leq N,\os{\circ}{t}}}
(E^{\bul \leq N})\otimes_{\mab Z}{\mab Q},P))
\tag{2.3.23.7}\label{eqn:fzpwd}
\end{align*} 
fitting into the following commutative diagram$:$
\begin{equation*} 
\begin{CD}
{\cal W}A_{Y_{\bul\leq N,{\os{\circ}{t}{}'}}}(F^{\bul \leq N})\otimes_{\mab Z}{\mab Q}
@>{g_{\bul \leq N}^*}>> \\ 
@A{\theta \wedge}A{\simeq}A \\
Ru_{Y_{\bul \leq N,\os{\circ}{t}{}'}/{\cal W}(s'_{\os{\circ}{t}{}'})*}
(\eps^*_{Y_{\bul \leq N,\os{\circ}{t}{}'}/
{\cal W}(s'_{\os{\circ}{t}{}'})}(F^{\bul \leq N}))\otimes_{\mab Z}{\mab Q}
@>{g_{\bul \leq N}^*}>>
\end{CD}
\tag{2.3.23.8}\label{cd:pscpz} 
\end{equation*} 
\begin{equation*} 
\begin{CD}
Rg_{\bul \leq N*}({\cal W}A_{X_{\bul\leq N,\os{\circ}{t}}}(E^{\bul \leq N}))
\otimes_{\mab Z}{\mab Q}\\ 
@A{Rg_{\bul \leq N*}(\theta \wedge)}A{\simeq}A \\
Rg_{\bul \leq N*}
Ru_{X_{\bul \leq N,\os{\circ}{t}}/{\cal W}(s_{\os{\circ}{t}})*}
(\eps^*_{X_{\bul \leq N,\os{\circ}{t}}/{\cal W}(s_{\os{\circ}{t}})}
(E^{\bul \leq N}))\otimes_{\mab Z}{\mab Q}.
\end{CD}
\end{equation*}
The similar relation to {\rm (\ref{ali:pwdpp})} for the $g_{\bul \leq N}^*$ above 
holds. 
\end{theo} 
\begin{proof} 
The proof of this theorem is easier than the proofs of {\rm (\ref{theo:funas})}, 
(\ref{theo:funpas}) and (\ref{theo:itc}). We leave the detail of the proof to the reader. 
\end{proof}

\begin{prop}\label{prop:dvwok}
Assume that the two conditions $(1.5.6.4)$ and $(1.5.6.5)$ hold 
in the case $S=s$, $S'=s'$, where $S={\cal W}(s)$, $S'={\cal W}(s')$, 
$T={\cal W}(t)$ and $T={\cal W}(t')$.
Assume also that $e_p\leq 1$. 
Then the morphism {\rm (\ref{ali:odnl})} is divisible by $p^{j+1}$. 
\end{prop}
\begin{proof} 
The proof is the same as that of (\ref{prop:dvok}). 
\end{proof} 

Let the assumptions be as in (\ref{prop:dvwok}).  
Consider the morphism 
\begin{align*} 
{\rm gr}^P_k(g^*_{\bul \leq N}) 
\col & {\rm gr}^P_k{\cal W}_{\star}A_{Y_{\bul \leq N,\os{\circ}{t}{}'}}
(F^{\bul \leq N}) \lo 
& Rg_{\bul \leq N*}
{\rm gr}^P_k{\cal W}_{\star}A_{X_{\bul \leq N,\os{\circ}{t}}}
(E^{\bul \leq N}).
\tag{2.3.24.1}\label{ali:gadkcl} 
\end{align*}
\par
Fix an integer $0\leq m\leq N$. 
Let the notations be before (\ref{ali:grgm}). 
Consider the following direct factor of the cosimplicial 
degree $m$-part of the morphism (\ref{ali:gadkcl}): 
\begin{align*}
g^* \col & 
b_{\phi(\ul{\lam}),t'*}
Ru_{\os{\circ}{Y}_{\phi(\ul{\lam}),t'}/{\cal W}_{\star}(\os{\circ}{t}{}')*}
(F_{\phi(\ul{\lam}),{\cal W}_{\star}(\os{\circ}{t}{}')}
\otimes_{\mab Z}
\vp_{\phi(\ul{\lam}){\rm crys}}(\os{\circ}{Y}_{t'}
/{\cal W}_{\star}(\os{\circ}{t}{}')))[-2j-k]
\tag{2.3.24.2}\label{ali:grgwm}\\
{} & \lo b_{\phi(\ul{\lam}),t'*}
Rg_{\ul{\lam},{\cal W}_{\star}(t'){\cal W}_{\star}(t)*}
Ru_{\os{\circ}{X}_{\ul{\lam},t}/{\cal W}_{\star}(\os{\circ}{t})*}
(E_{\ul{\lam}}
\otimes_{\mab Z}
\vp_{\ul{\lam}{\rm crys}}
(\os{\circ}{X}_t/{\cal W}_{\star}(\os{\circ}{t})))[-2j-k]. 
\end{align*}

\begin{prop}\label{prop:graloc}
Let the notations 
and the assumptions be as above. 
Then the morphism 
$g^*$ in {\rm (\ref{ali:gadkcl})} 
is equal to 
$\bigoplus_{j\geq \max\{-k,0\}}
\deg(u)^{j+k}b_{\phi(\ul{\lam}),t'*}
\os{\circ}{g}{}^*_{\ul{\lam}}$. 
\end{prop}
\begin{proof}
The proof is only a suitable modification of the proof of 
(\ref{prop:grloc}). We leave the detail of the proof to the reader.  
\end{proof}

\begin{coro}\label{coro:reis}
There exists the following isomorphism$:$ 
\begin{align*}
&{\rm gr}^P_k{\cal W}_{\star}A_{X_{m,\os{\circ}{t}}}(E^m)  
\tag{2.3.26.1}\label{ali:grawt}\\ 
& \os{\sim}{\lo} \bigoplus_{j\geq \max \{-k,0\}} 
a^{(2j+k)}_{m*}
(E_{\os{\circ}{X}{}^{(2j+k)}_{m,t}})_{
{\cal W}_n(\os{\circ}{X}{}^{(2j+k)}_{m,t})}
\otimes_{{\cal O}_{\os{\circ}{X}{}^{(2j+k)}_{m,t}}}
{\cal W}_n\Om^{\bul}_{\os{\circ}{X}{}^{(2j+k)}_{m,t}}\otimes_{\mab Z}\\
&\phantom{=\bigoplus_{\ul{t}\geq 0}R^{q-2\ul{t}_r-k}
f_{(D_{\ul{t}}^{(\ul{t}_r+k)},Z}}
\vp^{(2j+k)}_{\rm zar}(\os{\circ}{X}_{m,t}/{\cal W}_n(\os{\circ}{t})))(-j-k,u)[-2j-k]. \\ 
\end{align*} 
\end{coro} 

\par
The following is a main result in this section:

\begin{theo}[{\bf Comparison theorem I}]\label{theo:csoncrdw}
Let $N$ be a nonnegative integer and let $n$ be a positive integer. 
Let $E^{\bul \leq N}$ be a flat coherent log crystal of 
${\cal O}_{\os{\circ}{X}_{\bul \leq N,t}/{\cal W}_n(\os{\circ}{t})}$-modules.  
Assume that $X_{\bul \leq N,\os{\circ}{t}}$ 
has an affine $N$-truncated simplicial open covering. 
Then the following hold$:$ 
\par 
$(1)$ In ${\rm D}^+{\rm F}(f^{-1}({\cal W}_n(\kap_t)))$
there exists the following canonical isomorphism$:$
\begin{equation*}
(A_{\rm zar}(X_{\bul \leq N,\os{\circ}{t}}/{\cal W}_n(s_{\os{\circ}{t}}),E^{\bul \leq N}),P)
\os{\sim}{\lo}
({\cal W}_nA_{X_{\bul \leq N,\os{\circ}{t}}}(E^{\bul \leq N}),P).
\tag{2.3.27.1}\label{eqn:brildw}
\end{equation*}
The isomorphisms $(\ref{eqn:brildw})$ for $n$'s are 
compatible with two projections of 
both hand sides on $(\ref{eqn:brildw})$. 
\par 
$(2)$ Let the notations be as in {\rm (\ref{theo:ctrw}) (1)} or $(4)$.   
Then the isomorphism $(\ref{eqn:brildw})$ is contravariantly functorial 
for the commutative diagrams 
{\rm (\ref{eqn:xdxdwuss})} and {\rm (\ref{cd:xwygxy})}.   
\par 
$(3)$ The isomorphism 
{\rm (\ref{eqn:brildw})} forgetting the filtrations   
fits into the following commutative diagram$:$  
\begin{equation*} 
\begin{CD} 
A_{\rm zar}(X_{\bul \leq N,\os{\circ}{t}}/{\cal W}_{\star}(s_{\os{\circ}{t}}),E^{\bul \leq N})
@>{(\ref{eqn:brildw}),~\sim}>> \\
@A{\theta_{\star} \wedge}A{\simeq}A \\
Ru_{X_{\bul \leq N,t}/{\cal W}_{\star}(t)*}
(\eps^*_{X_{\bul \leq N,\os{\circ}{t}}/{\cal W}_{\star}(s_{\os{\circ}{t}})}(E^{\bul \leq N})) 
@>{(\ref{eqn:ywnnny}), \sim}>> 
\end{CD} 
\tag{2.3.27.2}\label{cd:axwl}
\end{equation*} 
\begin{equation*} 
\begin{CD} 
{\cal W}_{\star}A_{X_{\bul \leq N,\os{\circ}{t}}}(E^{\bul \leq N}) \\ 
@A{\theta_{\star}\wedge}A{\simeq}A \\
(\eps^*_{X_{\bul \leq N,\os{\circ}{t}/{\cal W}_{\star}(s_{\os{\circ}{t}})}}
(E^{\bul \leq N}))_{{\cal W}_{\star}(X_{\bul \leq N,\os{\circ}{t}})}
\otimes_{{\cal W}_{\star}({\cal O}_{X_{\bul \leq N,\os{\circ}{t}})}}
{\cal W}_{\star}{\Om}^{\bul}_{X_{\bul \leq N,\os{\circ}{t}}}, 
\end{CD} 
\end{equation*} 
where the upper vertical morphism 
$\theta_{\star} \wedge$ is the morphism $(\ref{eqn:uz})$ 
for the case $S={\cal W}_{\star}(s)$ and $T={\cal W}_{\star}(t)$. 
\par 
$(4)$ The commutative diagram {\rm (\ref{cd:axwl})} is contravariantly functorial. 
\end{theo}
\begin{proof} 
Let the notations be as in \S\ref{sec:psc}.  
\par 
(1): First we construct the morphism (\ref{eqn:brildw}). 
Let $X'_{\bul \leq N,\os{\circ}{t}}$ be the disjoint union of 
an affine $N$-truncated simplicial open covering of 
$X_{\bul \leq N,\os{\circ}{t}}$. 
Let $X'_{N,\os{\circ}{t}} \os{\sus}{\lo} \ol{\cal P}{}'_N$ be 
an immersion into a log smooth integral scheme over 
$\ol{{\cal W}_n(s_{\os{\circ}{t}})}$.
Since $X'_{N,\os{\circ}{t}}$ is affine, there exists a morphism 
${\cal W}_n(X'_{N,\os{\circ}{t}})\lo \ol{\cal P}{}'_N$ 
over $\ol{{\cal W}_n(s_{\os{\circ}{t}})}$ 
such that the composite morphism 
$X'_{N,\os{\circ}{t}}\os{\sus}{\lo} 
{\cal W}_n(X'_{N,\os{\circ}{t}})\lo \ol{\cal P}{}'_N$ 
is the closed immersion above.  
As in the proof of (\ref{theo:ccrw}),  
we have a natural morphism 
\begin{equation*} 
{\cal W}_n(X'_{\bul \leq N,\os{\circ}{t}})\lo 
\Gam^{\ol{{\cal W}_n(s_{\os{\circ}{t}})}}_N(\ol{\cal P}{}'_N)_{\bul \leq N} 
\tag{2.3.27.3}\label{eqn:pflwge} 
\end{equation*} 
of $N$-truncated simplicial log schemes. 
Let 
${\cal W}_n(X_{\bul \leq N,\bul,\os{\circ}{t}})$ 
be the \v{C}ech diagram of 
${\cal W}_n(X'_{\bul \leq N,\os{\circ}{t}})$ over 
${\cal W}_n(X_{\bul \leq N,\os{\circ}{t}})$.  
It is easy to check that 
${\cal W}_n(X_{\bul \leq N,\bul,\os{\circ}{t}})$ 
is the canonical lift of $X_{\bul \leq N,\bul,\os{\circ}{t}}$.  
Set $\ol{\cal P}_{mn}
:={\rm cosk}_0^{\ol{{\cal W}_n(s_{\os{\circ}{t}})}}
(\Gam^{\ol{{\cal W}_n(s_{\os{\circ}{t}})}}_N(\ol{\cal P}{}'_N)_m)_n$.  
As in (\ref{eqn:eixd}) we have the following morphism 
\begin{equation*} 
{\cal W}_n(X_{\bul \leq N,\bul,\os{\circ}{t}}) 
\lo \ol{\cal P}_{\bul \leq N,\bul}
\tag{2.3.27.4}\label{eqn:ppwgep} 
\end{equation*} 
of $(N,\infty)$-truncated bisimplicial log schemes 
such that the composite morphism  
\begin{equation*} 
X_{\bul \leq N,\bul,\os{\circ}{t}} \os{\sus}{\lo} 
{\cal W}_n(X_{\bul \leq N,\bul,\os{\circ}{t}}) \lo 
\ol{\cal P}_{\bul \leq N,\bul} 
\end{equation*} 
is the exact immersion in 
(\ref{eqn:eipxd}) in the case $S=s$ and $T={\cal W}_n(t)$. 
Set ${\cal P}_{\bul \leq N,\bul}:=
\ol{\cal P}_{\bul \leq N,\bul}
\times_{\ol{{\cal W}_n(s_{\os{\circ}{t}})}}{\cal W}_n(s_{\os{\circ}{t}})$. 
Let ${\cal P}^{\rm ex}_{\bul \leq N,\bul}$ 
be the exactification of the immersion 
$X_{\bul \leq N,\bul,\os{\circ}{t}} \os{\sus}{\lo} {\cal P}_{\bul \leq N,\bul}$. 
Let $\ol{\mathfrak D}_{\bul \leq N,\bul}$ 
be the log PD-envelope of 
the immersion $X_{\bul \leq N,\bul,\os{\circ}{t}}\os{\sus}{\lo}
\ol{\cal P}_{\bul \leq N,\bul}$ over 
$({\cal W}_n(\os{\circ}{t}),p{\cal W}_n,[~])$. 
Set ${\mathfrak D}_{\bul \leq N,\bul}
:=\ol{\mathfrak D}_{\bul \leq N,\bul}
\times_{{\mathfrak D}(\ol{{\cal W}_n(s_{\os{\circ}{t}})})}{\cal W}_n(s_{\os{\circ}{t}})$.  
Let $E^{\bul \leq N,\bul}$ be a crystal of 
${\cal O}_{\os{\circ}{X}_{\bul \leq N,\bul,t}
/{\cal W}_n(\os{\circ}{t})}$-modules obtained by 
$E^{\bul \leq N}$. 
Let $(\ol{\cal E}{}^{\bul \leq N,\bul},\ol{\nabla})$ 
be the coherent 
${\cal O}_{\ol{\mathfrak D}_{\bul \leq N,\bul}}$-module 
with integrable connection 
obtained by 
$\eps^*_{X_{\bul \leq N,\bul,\os{\circ}{t}}/{\cal W}_n(\os{\circ}{t})}
(E^{\bul \leq N,\bul})$. 
Set $({\cal E}^{\bul \leq N,\bul},\nabla)
:=(\ol{\cal E}{}^{\bul \leq N,\bul},\ol{\nabla})
\otimes_{{\cal O}_{{\mathfrak D}(\ol{{\cal W}_n(s_{\os{\circ}{t}})})}}
{\cal W}_n(\kap_t)$.  
By using the morphism (\ref{ali:pno}) and (\ref{lemm:nefi}) (2), 
we have the following commutative diagrams: 
\begin{equation*} 
\begin{CD} 
({\cal E}^{\bul \leq N,\bul}
{\otimes}_{{\cal O}_{{\cal P}^{\rm ex}_{\bul \leq N,\bul}}}
{\Om}^{i+j+2}_{{\cal P}^{\rm ex}_{\bul \leq N,\bul}
/{\cal W}_n(\os{\circ}{t})}/P_{j+1},P)  
@>>> \\
@A{d}AA  \\
({\cal E}^{\bul \leq N,\bul}
{\otimes}_{{\cal O}_{{\cal P}^{\rm ex}_{\bul \leq N,\bul}}}
{\Om}^{i+j+1}_{{\cal P}^{\rm ex}_{\bul \leq N,\bul}
/{\cal W}_n(\os{\circ}{t})}/P_{j+1},P)  
@>>> 
\end{CD}
\tag{2.3.27.5}\label{ali:pnho}
\end{equation*}   
\begin{equation*} 
\begin{CD}
(E^{\bul \leq N,\bul}_{n} 
\otimes_{{\cal W}_n{(\cal O}_{X_{\bul \leq N,\bul,\os{\circ}{t}}})}
{\cal W}_n\wt{{\Om}}^{i+j+2}_{X_{\bul \leq N,\bul,\os{\circ}{t}}}/P_{j+1},P) \\
@AA{d}A \\
(E^{\bul \leq N,\bul}_{n} 
\otimes_{{\cal W}_n{(\cal O}_{X_{\bul \leq N,\bul,\os{\circ}{t}}})}
{\cal W}_n\wt{{\Om}}^{i+j+1}_{X_{\bul \leq N,\bul,\os{\circ}{t}}}/P_{j+1},P),   
\end{CD}
\end{equation*} 
\begin{equation*} 
\begin{CD} 
({\cal E}^{\bul \leq N,\bul}
{\otimes}_{{\cal O}_{{\cal P}^{\rm ex}_{\bul \leq N,\bul}}}
{\Om}^{i+j+2}_{{\cal P}^{\rm ex}_{\bul \leq N,\bul}
/{\cal W}_n(\os{\circ}{t})}/P_{j+2},P)  
@>>> \\
@A{\theta_{{\cal P}^{\rm ex}_{\bul \leq N,\bul}}\wedge }AA \\
({\cal E}^{\bul \leq N,\bul}
{\otimes}_{{\cal O}_{{\cal P}^{\rm ex}_{\bul \leq N,\bul}}}
{\Om}^{i+j+1}_{{\cal P}^{\rm ex}_{\bul \leq N,\bul}
/{\cal W}_n(\os{\circ}{t})}/P_{j+1},P)  
@>>>   
\end{CD}
\tag{2.3.27.6}\label{ali:pnvo}
\end{equation*}   
\begin{equation*} 
\begin{CD} 
(E^{\bul \leq N,\bul}_{n} 
\otimes_{{\cal W}_n{(\cal O}_{X_{\bul \leq N,\bul,\os{\circ}{t}}})}
{\cal W}_n\wt{{\Om}}^{i+j+2}_{X_{\bul \leq N,\bul,\os{\circ}{t}}}/P_{j+2},P) \\
@AA{\theta_n \wedge}A \\ 
(E^{\bul \leq N,\bul}_{n} 
\otimes_{{\cal W}_n{(\cal O}_{X_{\bul \leq N,\bul,\os{\circ}{t}}})}
{\cal W}_n\wt{{\Om}}^{i+j+1}_{X_{\bul \leq N,\bul,\os{\circ}{t}}}/P_{j+1},P). 
\end{CD} 
\end{equation*}  
Hence we have the following filtered morphism in ${\rm C}^+{\rm F}
(f^{-1}_{\bul \leq N,\bul}({\cal W}_n(\kap_t)))$: 
\begin{equation*} 
(A_{\rm zar}({\cal P}^{\rm ex}_{\bul \leq N,\bul,\os{\circ}{t}}
/{\cal W}_n(s_{\os{\circ}{t}}),{\cal E}^{\bul \leq N,\bul}),P) 
\lo ({\cal W}_nA_{X_{\bul \leq N,\bul,\os{\circ}{t}}}(E^{\bul \leq N}),P).
\tag{2.3.27.7}\label{eqn:cadmp}
\end{equation*} 
Applying $R\pi_{{\rm zar}*}$ to (\ref{eqn:cadmp}), we have 
the filtered morphism (\ref{eqn:brildw}). 
\par 
To prove that (\ref{eqn:brildw}) is a filtered isomorphism 
in ${\rm D}^+{\rm F}
(f^{-1}_{\bul \leq N}({\cal W}_n(\kap_t)))$, 
we have only to prove that 
the following morphism 
\begin{equation*}
(A_{\rm zar}(X_{m,\os{\circ}{t}}/{\cal W}_n(s_{\os{\circ}{t}}),E^m),P) 
\lo
({\cal W}_nA_{X_{m,\os{\circ}{t}}}(E^m),P) 
\quad (0\leq m\leq N)
\tag{2.3.27.8}\label{eqn:fixtdw}
\end{equation*} 
is a filtered isomorphism. 
Because ${\rm gr}^P_k((\ref{eqn:fixtdw}))$ 
is the following canonical isomorphism 
\begin{align*}
\bigoplus_{j\geq \max \{-k,0\}} &
a^{(2j+k)}_{m*}
Ru_{\os{\circ}{X}{}^{(2j+k)}_{m,t}/{\cal W}_n(\os{\circ}{t})*}
(E_{\os{\circ}{X}{}^{(2j+k)}_{m,t}
/{\cal W}_n(\os{\circ}{t})}
\otimes_{\mab Z}\vp^{(2j+k)}_{\rm crys}(\os{\circ}{X}_{m,t}
/{\cal W}_n(\os{\circ}{t})))[-2j-k]  \tag{2.3.27.9}\label{ali:grwt}\\
 \os{\sim}{\lo} \bigoplus_{j\geq \max \{-k,0\}} &
(a^{(2j+k)}_{m*}
(E_{\os{\circ}{X}{}^{(2j+k)}_{m,t}})_{
{\cal W}_n(\os{\circ}{X}{}^{(2j+k)}_{m,t})}
\otimes_{{\cal O}_{\os{\circ}{X}{}^{(2j+k)}_{m,t}}}
{\cal W}_n\Om^{\bul}_{\os{\circ}{X}{}^{(2j+k)}_{m,t}}
\otimes_{\mab Z}\\
&\vp^{(2j+k)}_{\rm zar}(\os{\circ}{X}_{m,t}/{\cal W}_n(\os{\circ}{t})))[-2j-k], 
\end{align*} 
by (\ref{ali:cgrp}) and (\ref{ali:grawt}),
the morphism (\ref{eqn:fixtdw}) is 
a filtered isomorphism by the comparison theorem in \cite[II (2.1)]{et}.  
\par 
Next we note that the isomorphism 
(\ref{eqn:brildw}) is independent of 
the choices of an open covering of 
$X_{N,\os{\circ}{t}}$ and 
a closed immersion 
$X'_{N,\os{\circ}{t}} \os{\sus}{\lo} \ol{\cal P}{}'_N$.  
Indeed, 
assume that we are given another such closed immersion 
$X''_{N,\os{\circ}{t}} \os{\sus}{\lo} \ol{\cal P}{}''_N$ over $\ol{{\cal W}_n(s_{\os{\circ}{t}})}$, 
where $X''_{N,\os{\circ}{t}}$ is the disjoint union of an affine open covering of 
$X_{N,\os{\circ}{t}}$ and $\ol{\cal P}{}''_N$ is a log smooth integral scheme over 
$\ol{{\cal W}_n(s_{\os{\circ}{t}})}$.
Then there exists another disjoint union $X'''_{N,\os{\circ}{t}}$ 
and a closed immersion 
$X'''_{N,\os{\circ}{t}} \os{\sus}{\lo} \ol{\cal P}{}'''_N$ 
into a log smooth integral scheme over $s\lo \ol{{\cal W}_n(s)}$  
fitting into the following commutative diagram 
\begin{equation*} 
\begin{CD} 
X'_{N,\os{\circ}{t}} @>{\sus}>> \ol{\cal P}{}'_N \\
@AAA @AAA \\
X'''_{N,\os{\circ}{t}} @>{\sus}>>\ol{\cal P}{}'''_N \\
@VVV @VVV \\
X''_{N,\os{\circ}{t}} @>{\sus}>>\ol{\cal P}{}''_N. 
\end{CD} 
\tag{2.3.27.10}\label{cd:xpppdpppz}
\end{equation*} 
Now, by using the argument for 
the construction of the morphism (\ref{eqn:brildw}), 
it is a routine work to 
prove that the filtered isomorphism (\ref{eqn:brildw}) is 
independent of the choices of 
an open covering of $X_{N,\os{\circ}{t}}$ 
and the closed immersion 
$X'_{N,\os{\circ}{t}} \os{\sus}{\lo} \ol{\cal P}{}'_N$. 
\par 
By a well-known argument (e.~g., \cite[(7.1), (7.19)]{ndw}), 
the isomorphism (\ref{eqn:brildw}) 
is compatible with the projections.  
\par 
(2): The proof of (2) is only a modification of the proof of 
(\ref{theo:ccrw}) (2) by using the commutative diagrams (\ref{ali:pnho}) and 
(\ref{ali:pnvo}) and by the definitions of the pull-back $g^*_{\bul \leq N}$ 
on both hand sides on (\ref{eqn:brildw}). 
We leave the proof of the contravariant functoriality of (\ref{eqn:wlxa}) to the reader.
\par 
(3): Using the commutative diagrams (\ref{ali:pnho}) and 
(\ref{ali:pnvo}), we immediately obtain the commutative diagram 
(\ref{cd:axwl}). 
\par 
(4):  As in (2), we leave the proof to the reader. 
\end{proof}


\par 
We restate (\ref{ali:grwt}) for our memory. 

\begin{coro} 
For each $0\leq m\leq N$, 
the following diagram is commutative$:$ 
\begin{equation*}
\begin{CD} 
{\rm gr}^P_kA_{\rm zar}(X_{m,t}/{\cal W}_n(s_{\os{\circ}{t}}),E^{\bul \leq N}))
@>{(\ref{eqn:brildw}),~\sim}>> 
{\rm gr}^P_k{\cal W}_nA_{X_{m,\os{\circ}{t}}}(E^{\bul \leq N}))\\
@VVV @VVV\\
\end{CD} 
\tag{2.3.28.1}\label{cd:pkag}
\end{equation*}
\begin{equation*}
\begin{split}
&\bigoplus_{j\geq \max \{-k,0\}} 
a^{(2j+k)}_{m*}
Ru_{\os{\circ}{X}{}^{(2j+k)}_{m,t}/{\cal W}_n(\os{\circ}{t})*}
(E_{\os{\circ}{X}{}^{(2j+k)}_{m,t}
/{\cal W}_n(\os{\circ}{t})}
\otimes_{\mab Z} \\
&\phantom{=\bigoplus_{\ul{t}\geq 0}R^{q-2\ul{t}_r-k}
f_{(D_{\ul{t}}^{(\ul{t}_r+k)},Z\vert_{D^{(\ul{t}_r+k)}}}}
\vp^{(2j+k)}_{\rm crys}(\os{\circ}{X}_{m,t}
/{\cal W}_n(\os{\circ}{t})))(-j-k,u)[-2j-k] \\
& \os{(\ref{ali:grwt}),\sim}{\lo} \bigoplus_{j\geq \max \{-k,0\}} 
a^{(2j+k)}_{m*}
(E_{\os{\circ}{X}{}^{(2j+k)}_{m,t}/{\cal W}_n(\os{\circ}{t})})_{
{\cal W}_n(\os{\circ}{X}{}^{(2j+k)}_{m,t})}
\otimes_{{\cal O}_{\os{\circ}{X}{}^{(2j+k)}_{m,t}}}
{\cal W}_n\Om^{\bul}_{\os{\circ}{X}{}^{(2j+k)}_{m,t}}\otimes_{\mab Z}\\ 
&\phantom{=\bigoplus_{\ul{t}\geq 0}R^{q-2\ul{t}_r-k}
f_{(D_{\ul{t}}^{(\ul{t}_r+k)},Z\vert_{D^{(\ul{t}_r+k)}}}}
\vp^{(2j+k)}_{\rm crys}(\os{\circ}{X}_{m,t}
/{\cal W}_n(\os{\circ}{t})))(-j-k,u)[-2j-k]. \\ 
&
\end{split} 
\end{equation*} 
\end{coro} 
\begin{proof} 
This follows from (\ref{ali:grwt}), (\ref{ali:cgrp}) and (\ref{ali:grawt}). 
\end{proof}


\par 
Let $\star$ be a positive integer $n$ or nothing. 
Let $E^{\bul \leq N}$  
be a flat quasi-coherent crystal of 
${\cal O}_{\os{\circ}{X}_{\bul \leq N,t}/{\cal W}_{\star}(\os{\circ}{t})}$-modules.   
Let $f_{\bul \leq N}\col X_{\bul \leq N,\os{\circ}{t}}\lo {\cal W}_{\star}(s_{\os{\circ}{t}})$ and 
$f_{m}\col X_{m,\os{\circ}{t}}\lo {\cal W}_{\star}(s_{\os{\circ}{t}})$ $(0\leq m \leq N)$ 
be the structural morphisms. 
Let $f_{\bul \leq N}\col X_{\bul \leq N,\os{\circ}{t}}\lo 
{\cal W}_{\star}(\os{\circ}{t})_{\bul \leq N}$
be the induced morphisms.  
Let $L$ be the stupid filtration on 
$Rf_{{\bul \leq N}*}
({\cal W}_{\star}A_{X_{\bul \leq N,\os{\circ}{t}}}(E^{\bul \leq N}))$
with respect to the cosimplicial degrees: 
\begin{equation*}
L^m(Rf_{\bul \leq N*}
({\cal W}_{\star}A_{X_{\bul \leq N,\os{\circ}{t}}}(E^{\bul \leq N})))
= \bigoplus_{m'\geq m}
Rf_{m'*}({\cal W}_{\star}A_{X_{m',\os{\circ}{t}}}(E^{m'}))
\quad (m\in {\mab N}). 
\tag{2.3.28.2}\label{eqn:ltcwp}
\end{equation*}
\par 
Let 
$\del(L,P)$ 
be the diagonal filtration of $L$  and $P$ on $Rf_{*}
({\cal W}_{\star}A_{X_{\bul \leq N,\os{\circ}{t}}}(E^{\bul \leq N}))$
(cf.~\cite[(7.1.6.1), (8.1.22)]{dh3}): 
\begin{align*}
\del(L,P)_k(Rf_{*}{\cal W}_{\star}A_{X_{\bul \leq N,\os{\circ}{t}}}(E^{\bul \leq N}))=  
\bigoplus_{m=0}^N
P_{k+m}Rf_{m*}{\cal W}_{\star}A_{X_{m,\os{\circ}{t}}}(E^m). 
\tag{2.3.28.3}\label{eqn:ddifl} 
\end{align*} 
(Strictly speaking, we have to take the representative of 
$Rf_{*}{\cal W}_{\star}A_{X_{\bul \leq N,\os{\circ}{t}}}(E^{\bul \leq N})$.)
Then we easily have the following equality:  
\begin{align*} 
& {\rm gr}^{\del(L,P)}_k
Rf_{*}({\cal W}_{\star}A_{X_{\bul \leq N,\os{\circ}{t}}}(E^{\bul \leq N})) 
= \tag{2.3.28.4}\label{ali:ruwagrvp}\\
& \bigoplus_{m=0}^N\bigoplus_{j\geq \max \{-(k+m),0\}} 
Rf_{\os{\circ}{X}{}^{(2j+k+m)}_{m*}}
((E^m_{\os{\circ}{X}{}^{(2j+k+m)}_{m,t}/{\cal W}_{\star}(\os{\circ}{t})})_{{\cal W}_{\star}(\os{\circ}{X}{}^{(2j+k+m)}_{m,t})}
\otimes_{{\cal W}_{\star}({\cal O}_{\os{\circ}{X}{}^{(2j+k+m)}_{m,t}})}
 \\
&{\cal W}_{\star}\Om^{\bul}_{\os{\circ}{X}{}^{(2j+k+m)}_{m,t}}
\otimes_{\mab Z}\vp_{\rm zar}^{(2j+k+m)}
(\os{\circ}{X}_{m,t}/\os{\circ}{t}))(-j-k-m,u)[-2j-k-2m]. 
\end{align*} 
By (\ref{ali:ruwagrvp}), (\ref{eqn:wlxa}) and (\ref{eqn:ywnnny}), 
we have the following spectral sequence 
$($cf.~{\rm \cite[3.23]{msemi}, 
\cite[(9.9.4)]{ndw})}${\rm :}$
\begin{align*} 
E_1^{-k,q+k}
= & \bigoplus_{m=0}^N
\us{j\geq {\rm max}\{-(k+m),0\}}
\bigoplus
H^{q-2j-k-2m}(\os{\circ}{X}{}^{(2j+k+m)}_{m,t}, 
(E^m_{\os{\circ}{X}{}^{(2j+k+m)}_{m,t}/{\cal W}_{\star}(\os{\circ}{t})})
_{{\cal W}_{\star}(\os{\circ}{X}{}^{(2j+k+m)}_{m,t})}
\tag{2.3.28.5}\label{eqn:trhkwsp}\\
&
\otimes_{{\cal W}_{\star}({\cal O}_{\os{\circ}{X}{}^{(2j+k+m)}_{m,t}})}
{\cal W}_{\star}\Om^{\bul}_{\os{\circ}{X}{}^{(2j+k+m)}_{m,t}}
\otimes_{\mab Z}
\vp^{(2j+k+m)}_{\rm zar}(\os{\circ}{X}_{m,t}/\os{\circ}{t})) \\
& (-j-k-m,u)\Lo 
H^q((X_{\bul \leq N,\os{\circ}{t}}/{\cal W}_{\star}(\os{\circ}{t}))_{\rm crys},
\eps^*_{X_{\bul \leq N,\os{\circ}{t}}/{\cal W}_{\star}(\os{\circ}{t})}(E^{\bul \leq N})).
\end{align*}
Here we use the following natural identification  
\begin{align*}
&H^{*}(\os{\circ}{X}{}^{(e)}_{m,t}, (E^m_{\os{\circ}{X}{}^{(e)}_{m,t}/{\cal W}_{\star}(\os{\circ}{t})})
_{{\cal W}_{\star}(\os{\circ}{X}{}^{(e)}_{m,t})}
\otimes_{{\cal W}_{\star}({\cal O}_{\os{\circ}{X}{}^{(e)}_{m}})}
{\cal W}_{\star}\Om^{\bul}_{\os{\circ}{X}{}^{(e)}_m},-d))\tag{2.3.28.6}\label{eqn:fidpmbd}\\
&=
H^{*}(\os{\circ}{X}{}^{(e)}_{m,t},
(E^m_{\os{\circ}{X}{}^{(e)}_{m,t}/{\cal W}_{\star}(\os{\circ}{t})})
_{{\cal W}_{\star}(\os{\circ}{X}{}^{(e)}_{m,t})}
\otimes_{{\cal W}_{\star}({\cal O}_{\os{\circ}{X}{}^{(e)}_{m,t}})}
{\cal W}_{\star}\Om^{\bul}_{\os{\circ}{X}{}^{(e)}_{m,t}}) \\ 
&=H^{*}((\os{\circ}{X}{}^{(e)}_{m,t}/{\cal W}_{\star}(\os{\circ}{t}))_{\rm crys},
E^m_{\os{\circ}{X}{}^{(e)}_{m,t}/{\cal W}_{\star}(\os{\circ}{t})}) 
\quad (e\in {\mab N}). 
\end{align*} 

\begin{defi}[{\bf Abrelative Frobenius morphism}]\label{defi:rawd} 
(1) Let $F_{{\cal W}(s)} \col {\cal W}(s) \lo {\cal W}(s)$ 
(resp.~$F_s\col s \lo s$) be 
the Frobenius endomorphism of ${\cal W}(s)$ (resp.~$s$).   
Set $s^{[p]}:=s\times_{\os{\circ}{s},\os{\circ}{F}_{s}}\os{\circ}{s}$. 
Then we have natural morphisms 
\begin{align*} 
F_{s/\os{\circ}{s}}\col s\lo s^{[p]} \quad {\rm and} \quad  
W_{s/\os{\circ}{s}} \col s^{[p]}\lo s.
\end{align*} 
(The underlying morphism of the former morphism is ${\rm id}_{\os{\circ}{s}}$.) 
Let 
\begin{equation*}
\begin{CD} 
t@>>> t'\\
@VVV @VVV \\
s@>>> s^{[p]}
\end{CD} 
\tag{2.3.29.1}\label{cd:ttps}
\end{equation*} 
be a commutative diagram of fine log schemes whose underlying schemes 
are the spectrums of perfect fields of characteristic $p>0$. 
Then we have the following natural morphisms  
\begin{align*} 
s_{\os{\circ}{t}}\lo s^{[p]}_{\os{\circ}{t}{}'}\quad 
{\rm and}  \quad {\cal W}(s_{\os{\circ}{t}})\lo {\cal W}(s^{[p]}_{\os{\circ}{t}{}'}).
\end{align*}  
We call the morphisms $s_{\os{\circ}{t}}\lo s^{[p]}_{\os{\circ}{t}{}'}$ 
and 
${\cal W}(s_{\os{\circ}{t}})\lo {\cal W}(s^{[p]}_{\os{\circ}{t}{}'})$ 
the {\it abrelative Frobenius morphism of base log schemes}  
and the {\it abrelative Frobenius morphism of base log PD-schemes}, respectively. 
In particular, when $t'=t$, we have natural morphisms 
\begin{align*} 
s_{\os{\circ}{t}}\lo s^{[p]}_{\os{\circ}{t}} \quad 
{\rm and} \quad  
{\cal W}(s_{\os{\circ}{t}})\lo s^{[p]}({\cal W}(t))
\end{align*}   
by using a composite morphism 
$t\lo s\os{W_{s/\os{\circ}{s}}}{\lo} s^{[p]}$.  
\par 
(2) Let the notations be as in (1).  
Set $X^{[p]}_{\bul \leq N}:=X_{\bul \leq N}\times_ss^{[p]}
=X_{\bul \leq N}\times_{\os{\circ}{s},\os{\circ}{F}_s}\os{\circ}{s}$ 
and 
\begin{align*} 
X^{[p]}_{\bul \leq N,\os{\circ}{t}{}'}&:=X^{[p]}_{\bul \leq N}
\times_{s^{[p]}}s^{[p]}_{\os{\circ}{t}{}'}=X^{[p]}_{\bul \leq N}
\times_{\os{\circ}{s^{[p]}}}\os{\circ}{t}{}'. 
\end{align*} 
Then $X^{[p]}_{\bul \leq N,\os{\circ}{t}{}'}/s^{[p]}_{\os{\circ}{t}{}'}$ 
is an SNCL scheme. 
Let 
$$F^{\rm ar}_{X_{\bul \leq N,\os{\circ}{t}/s_{\os{\circ}{t}},s^{[p]}_{\os{\circ}{t}{}'}}}\col 
X_{\bul \leq N,\os{\circ}{t}}  \lo X^{[p]}_{\bul \leq N,\os{\circ}{t}{}'}$$ 
and 
$$F^{\rm ar}_{X_{\bul \leq N,\os{\circ}{t}/{\cal W}(s_{\os{\circ}{t}}),
{\cal W}(s^{[p]}_{\os{\circ}{t}{}'})}}\col 
X_{\bul \leq N,\os{\circ}{t}}  \lo X^{[p]}_{\bul \leq N,\os{\circ}{t}{}'}$$ 
be the abrelative Frobenius morphism  
over $s_{\os{\circ}{t}}\lo s^{[p]}_{\os{\circ}{t}{}'}$ and 
${\cal W}(s_{\os{\circ}{t}})\lo {\cal W}(s^{[p]}_{\os{\circ}{t}{}'})$.  
Let $E^{\bul \leq N}$ and $E'{}^{\bul \leq N}$ be 
a flat quasi-coherent crystal of 
${\cal O}_{\os{\circ}{X}_{\bul \leq N,t}/{\cal W}(\os{\circ}{t})}$-modules and 
a flat quasi-coherent crystal of 
${\cal O}_{\os{\circ}{X}_{\bul \leq N,t'}/{\cal W}(\os{\circ}{t}{}')}$-modules, 
respectively.  
Let 
\begin{align*} 
\Phi^{\rm ar} \col 
\os{\circ}{F}{}^{{\rm ar}*}_{X_{\bul \leq N,\os{\circ}{t}
/{\cal W}(s_{\os{\circ}{t}}),{\cal W}(s^{[p]}_{\os{\circ}{t}{}'}),{\rm crys}}}(E'{}^{\bul \leq N})
\lo E^{\bul \leq N}
\tag{2.3.29.2}\label{ali:spwpts}
\end{align*} 
be a morphism of crystals in 
$(\os{\circ}{X}_{\bul \leq N,t}/{\cal W}(\os{\circ}{t}))_{\rm crys}$.   
Since $\deg (F_{s/\os{\circ}{s}})=p$, 
the divisibility of the morphism (\ref{ali:odnl}) holds by (\ref{prop:dvwok}).  
We call the following induced morphism by $\Phi^{\rm ar}$   
\begin{align*} 
\Phi^{\rm ar} \col &
({\cal W}_{\star}A_{X^{[p]}_{\bul \leq N,\os{\circ}{t}{}'}}(E'{}^{\bul \leq N}),P) 
\lo 
RF^{\rm ar}_{X_{\bul \leq N,\os{\circ}{t}/{\cal W}(s_{\os{\circ}{t}}),
{\cal W}(s^{[p]}_{\os{\circ}{t}{}'}),{\rm crys}}*}
(({\cal W}_{\star}A_{X_{\bul \leq N,\os{\circ}{t}{}}}(E^{\bul \leq N}),P)) 
\tag{2.3.29.3}\label{ali:spwswpts} 
\end{align*}
the {\it abrelative Frobenius morphism} of 
$$({\cal W}_{\star}A_{X_{\bul \leq N,\os{\circ}{t}{}}}(E^{\bul \leq N}),P)
\quad 
{\rm and} \quad  
({\cal W}_{\star}A_{X^{[p]}_{\bul \leq N,\os{\circ}{t}{}'}}(E'{}^{\bul \leq N}),P).$$ 
When $E'{}^{\bul \leq N}={\cal O}_{\os{\circ}{X}{}'_{\bul \leq N,t'}/{\cal W}(\os{\circ}{t}{}')}$, 
we set 
$$({\cal W}_{\star}A_{X^{[p]}_{\bul \leq N,\os{\circ}{t}{}'}},P) 
:=({\cal W}_{\star}A_{X^{[p]}_{\bul \leq N,\os{\circ}{t}{}'}}(E'{}^{\bul \leq N}),P).$$  
Then we have the following {\it abrelative Frobenius morphism} 
\begin{equation*} 
\Phi^{\rm ar} \col 
({\cal W}_{\star}A_{X^{[p]}_{\bul \leq N,\os{\circ}{t}{}'}},P) 
\lo RF^{\rm ar}_{X_{\bul \leq N,\os{\circ}{t}/{\cal W}(s_{\os{\circ}{t}}),
{\cal W}(s^{[p]}_{\os{\circ}{t}{}'}),{\rm crys}}*}
(({\cal W}_{\star}A_{X_{\bul \leq N,\os{\circ}{t}}},P)) 
\tag{2.3.29.4}
\end{equation*}
of $({\cal W}_{\star}A_{X_{\bul \leq N,\os{\circ}{t}{}}},P)$ and 
$({\cal W}_{\star}A_{X^{[p]}_{\bul \leq N,\os{\circ}{t}{}'}},P)$. 
\end{defi}

\begin{prop}[{\bf Frobenius compatibility I}]\label{prop:fcar} 
The following diagram is commutative$:$ 
\begin{equation*} 
\begin{CD} 
{\cal W}_{\star}A_{X^{[p]}_{\bul \leq N,\os{\circ}{t}{}'}}(E'{}^{\bul \leq N})
@>{\Phi^{\rm ar}}>>
 \\
@A{\theta_{\star} \wedge}A{\simeq}A  \\
Ru_{X^{[p]}_{\bul \leq N,\os{\circ}{t}{}'}/{\cal W}_{\star}(s^{[p]}_{\os{\circ}{t}{}'})*}
(\eps^*_{X^{[p]}_{\bul \leq N,\os{\circ}{t}{}'}/{\cal W}_{\star}(s^{[p]}_{\os{\circ}{t}{}'})}
(E'{}^{\bul \leq N}))
@>{\Phi^{\rm ar}}>>
\end{CD}
\tag{2.3.30.1}\label{cd:rwukt}
\end{equation*} 
\begin{equation*} 
\begin{CD} 
RF^{\rm ar}_{X_{\bul \leq N,\os{\circ}{t}
/{\cal W}(s_{\os{\circ}{t}}),{\cal W}(s^{[p]}_{\os{\circ}{t}{}'})}*}
({\cal W}_{\star}A_{X_{\bul \leq N,\os{\circ}{t}{}}}(E^{\bul \leq N})) \\
@A{RF^{\rm ar}_{X_{\bul \leq N,\os{\circ}{t}/{\cal W}(s_{\os{\circ}{t}}),
{\cal W}(s^{[p]}_{\os{\circ}{t}{}'})}*}(\theta_{\star})}A{\simeq}A \\
RF^{\rm ar}_{X_{\bul \leq N,\os{\circ}{t}/{\cal W}_{\star}(s_{\os{\circ}{t}})}*}
Ru_{X_{\bul \leq N,\os{\circ}{t}}/{\cal W}_{\star}(s_{\os{\circ}{t}})*}
(\eps^*_{X_{\bul \leq N,\os{\circ}{t}}/{\cal W}_{\star}(s_{\os{\circ}{t}})}(E^{\bul \leq N})). 
\end{CD}
\end{equation*} 
The commutative diagram {\rm (\ref{cd:rwukt})} is contravariantly functorial 
for the morphism {\rm (\ref{eqn:xdxdwuss})} satisfying {\rm (\ref{cd:xwygxy})} 
and for the morphism of $F$-crystals 
\begin{equation*} 
\begin{CD} 
\os{\circ}{F}{}^{{\rm ar}*}_{X_{\bul \leq N,\os{\circ}{t}
/{\cal W}(s_{\os{\circ}{t}}),{\cal W}(s^{[p]}_{\os{\circ}{t}{}'}),{\rm crys}}}
(E'{}^{\bul \leq N})
@>{\Phi^{\rm ar}}>> E^{\bul \leq N}\\
@AAA @AAA \\
\os{\circ}{g}{}^*_{\bul \leq N}
\os{\circ}{F}{}^{{\rm ar}*}_{Y_{\bul \leq N,\os{\circ}{t}
/{\cal W}(s_{\os{\circ}{t}}),{\cal W}(s^{[p]}_{\os{\circ}{t}{}'}),{\rm crys}}}
(F'{}^{\bul \leq N})
@>{\os{\circ}{g}{}^*_{\bul \leq N}(\Phi^{\rm ar})}>> 
\os{\circ}{g}{}^*_{\bul \leq N}(F^{\bul \leq N}), 
\end{CD}
\end{equation*} 
where $F^{\bul \leq N}$ is a similar quasi-coherent 
${\cal O}_{\os{\circ}{Y}_{\bul \leq N,t}/{\cal W}_n(\os{\circ}{t})}$-module 
to $E^{\bul \leq N}$. 
\end{prop} 
\begin{proof} 
This is a special case of (\ref{theo:ctrw}). 
\end{proof} 

\begin{prop}[{\bf Frobenius compatibility II}]\label{prop:nwlsfc} 
Assume that the morphism $t'\lo s^{[p]}$ factors through the morphism 
$F_{s/{\os{\circ}{s}}}\col s\lo s^{[p]}$. 
Set $X^{\{p\}}_{\bul \leq N}:=X_{\bul \leq N}\times_{s,F_s}s$ 
and $X^{\{p\}}_{\bul \leq N,\os{\circ}{t}{}'}:=
X^{\{p\}}_{\bul \leq N}\times_ss_{\os{\circ}{t}{}'}$. 
Let 
$$F^{{\rm rel}}_{X_{\bul \leq N,\os{\circ}{t}/s_{\os{\circ}{t}},s_{\os{\circ}{t}{}'}}}
\col X_{\bul \leq N,\os{\circ}{t}}  \lo X^{\{p\}}_{\bul \leq N,\os{\circ}{t}{}'}$$ 
and 
$$F^{\rm rel}_{X_{\bul \leq N,\os{\circ}{t}/s(T'),S(T)}}\col 
X_{\bul \leq N,\os{\circ}{t}}  \lo X^{\{p\}}_{\bul \leq N,\os{\circ}{t}{}'}$$ 
be the relative Frobenius morphism  
over $s_{\os{\circ}{t}}\lo s_{\os{\circ}{t}{}'}$ and 
$({\cal W}_{\star}(s_{\os{\circ}{t}}),p{\cal W}(\kap_t),[~])\lo 
({\cal W}_{\star}(s_{\os{\circ}{t}{}'}),p{\cal W}(\kap_{t'}),[~])$. 
Let $E''{}^{\bul \leq N}$ be the pull-back of 
$E'{}^{\bul \leq N}$ to $(\os{\circ}{X}{}^{\{p\}}_{\bul \leq N,t}/\os{\circ}{t})_{\rm crys}$. 
Let 
\begin{align*} 
\Phi^{{\rm rel}} \col 
\os{\circ}{F}{}^{{\rm  rel}*}_{X_{\bul \leq N,\os{\circ}{t}/
{\cal W}_{\star}(s_{\os{\circ}{t}}),{\cal W}_{\star}(s_{\os{\circ}{t}{}'}),{\rm crys}}}
(E''{}^{\bul \leq N})\lo E^{\bul \leq N}  
\tag{2.3.31.1}\label{ali:swptpps}
\end{align*} 
be the induced morphism by {\rm (\ref{ali:spwpts})}. 
Then the following diagram
\begin{equation*} 
\begin{CD} 
{\cal W}_{\star}A_{X^{[p]}_{\bul \leq N,\os{\circ}{t}{}'}}(E'{}^{\bul \leq N})
@>{\Phi^{\rm  rel}}>> \\
@A{\theta_{\star} \wedge}A{\simeq}A  \\
Ru_{X^{\{p\}}_{\bul \leq N,\os{\circ}{t}{}'}/
{\cal W}_{\star}(s_{\os{\circ}{t}{}'})*}
(\eps^*_{X^{\{p\}}_{\bul \leq N,\os{\circ}{t}{}'}/{\cal W}_{\star}(s_{\os{\circ}{t}{}'})}
(E''{}^{\bul \leq N}))
@>{\Phi^{\rm rel}}>>
\end{CD}
\tag{2.3.31.2}\label{cd:rekwtw}
\end{equation*} 
\begin{equation*} 
\begin{CD} 
RF^{{\rm rel}}_{X_{\bul \leq N,\os{\circ}{t}/{\cal W}_{\star}(s_{\os{\circ}{t}}),
{\cal W}_{\star}(s_{\os{\circ}{t}{}'})}*}
({\cal W}_{\star}A_{X_{\bul \leq N,\os{\circ}{t}{}'}}(E'{}^{\bul \leq N}))\\
@A{\simeq}A{RF^{\rm rel}_{X_{\bul \leq N,\os{\circ}{t}
/{\cal W}_{\star}(s_{\os{\circ}{t}}),{\cal W}_{\star}(s_{\os{\circ}{t}{}'})}*}
(\theta_{\star})\wedge}A\\
RF^{\rm rel}_{X_{\bul \leq N,\os{\circ}{t}/{\cal W}_{\star}(s_{\os{\circ}{t}}),
{\cal W}_{\star}(s_{\os{\circ}{t}{}'})}*}
Ru_{X_{\bul \leq N,t}/{\cal W}_{\star}(s_{\os{\circ}{t}})*}
(\eps^*_{X_{\bul \leq N,\os{\circ}{t}}/{\cal W}_{\star}(s_{\os{\circ}{t}})}(E^{\bul \leq N})) 
\end{CD} 
\end{equation*} 
is commutative. The commutative diagram {\rm (\ref{cd:rekwtw})} 
is contravariantly functorial 
for the morphism {\rm (\ref{eqn:xdxdwuss})} satisfying {\rm (\ref{cd:xwygxy})}
and for the morphism of $F$-crystals 
\begin{equation*} 
\begin{CD} 
\os{\circ}{F}{}^{{\rm rel}*}_{X_{\bul \leq N,\os{\circ}{t}
/{\cal W}_{\star}(s_{\os{\circ}{t}}),
{\cal W}_{\star}(s_{\os{\circ}{t}{}'}),{\rm crys}}}(E''{}^{\bul \leq N})
@>{\Phi^{\rm rel}}>> E^{\bul \leq N}\\
@AAA @AAA \\
\os{\circ}{g}{}^*_{\bul \leq N}
\os{\circ}{F}{}^{{\rm rel}*}_{Y_{\bul \leq N,\os{\circ}{t}/
{\cal W}_{\star}(s_{\os{\circ}{t}}),
{\cal W}_{\star}(s_{\os{\circ}{t}{}'}),{\rm crys}}}
(F'{}^{\bul \leq N})
@>{\os{\circ}{g}{}^*_{\bul \leq N}(\Phi^{\rm rel})}>> 
\os{\circ}{g}{}^*_{\bul \leq N}(F^{\bul \leq N}), 
\end{CD}
\end{equation*} 
where $F^{\bul \leq N}$ is a similar quasi-coherent 
${\cal O}_{\os{\circ}{Y}_{\bul \leq N,t}/\os{\circ}{T}}$-module 
to $E^{\bul \leq N}$. 
\end{prop}
\begin{proof} 
Because the morphism 
$t\lo s^{[p]}$ is the composite morphism 
$t\lo s_{\os{\circ}{t}}\lo s^{[p]}_{\os{\circ}{t}}$, 
$X^{[p]}_{\bul \leq N}\times_{s^{[p]},F^{\rm ar}_s}s\times_{\os{\circ}{s}}{\os{\circ}{t}}
=X^{\{p\}}_{\bul \leq N,\os{\circ}{t}}$. 
Consequently  
\begin{align*}
Ru_{X^{[p]}_{\bul \leq N,\os{\circ}{t}}/{\cal W}_{\star}(s_{\os{\circ}{t}})*}
(\eps^*_{X^{[p]}_{\bul \leq N,\os{\circ}{t}}/{\cal W}_{\star}(s_{\os{\circ}{t}})}
(E'{}^{\bul \leq N})) 
=
Ru_{X^{\{p\}}_{\bul \leq N,\os{\circ}{t}}/{\cal W}_{\star}(s_{\os{\circ}{t}})*}
(\eps^*_{X^{\{p\}}_{\bul \leq N,\os{\circ}{t}}/{\cal W}_{\star}(s_{\os{\circ}{t}})}
(E''{}^{\bul \leq N})). 
\end{align*} 
Hence we obtain the commutative diagram 
(\ref{cd:rekwtw}) by (\ref{cd:rwukt}) and (\ref{lemm:flpis}). 
\par 
We leave the proof of the functoriality to the reader. 
\end{proof}

\begin{defi}[{\bf Absolute Frobenius endomorphism}]\label{defi:bwtd}  
Let the notations and the assumptions be as in (\ref{defi:rwd}).  
Let $F_{s_{\os{\circ}{t}}} \col s_{\os{\circ}{t}} \lo s_{\os{\circ}{t}}$ be 
the Frobenius endomorphism of $s_{\os{\circ}{t}}$.  
Note that there exists a natural Frobenius endomorphism 
$$F_{{\cal W}_{\star}(s_{\os{\circ}{t}})} \col  {\cal W}_{\star}(s_{\os{\circ}{t}})\lo {\cal W}_{\star}(s_{\os{\circ}{t}}),$$ 
which is a lift of  
$F_{s_{\os{\circ}{t}}}$.  This gives a PD-morphism 
$$F_{{\cal W}_{\star}(s_{\os{\circ}{t}})}\col ({\cal W}_{\star}(s_{\os{\circ}{t}}),p{\cal W}(\kap_t),[~])\lo 
({\cal W}_{\star}(s_{\os{\circ}{t}}),p{\cal W}_{\star}(\kap_t),[~])$$ 
of $p$-adic formal family of log points.  
Let 
$$F^{\rm abs}_{X_{\bul \leq N,\os{\circ}{t}_0}/{\cal W}_{\star}(s_{\os{\circ}{t}})} 
\col X_{\bul \leq N,\os{\circ}{t}}  \lo X_{\bul \leq N,\os{\circ}{t}}$$ 
be the absolute Frobenius endomorphism over $F_{{\cal W}_{\star}(s_{\os{\circ}{t}})}$.   
Let 
\begin{align*} 
\Phi^{\rm abs} \col 
\os{\circ}{F}{}^{{\rm abs}*}_{X_{\bul \leq N,\os{\circ}{t}/{\cal W}_{\star}(s_{\os{\circ}{t}})},{\rm crys}}
(E^{\bul \leq N})\lo E^{\bul \leq N}
\end{align*} 
be a morphism of crystals in 
$(\os{\circ}{X}_{\bul \leq N,t}/{\cal W}_{\star}(\os{\circ}{t}))_{\rm crys}$.   
Then the divisibility of the morphism (\ref{ali:odnl}) 
holds in this situation with  $e_p=1$.  
Then we call the induced morphism by $\Phi^{\rm abs}$ and $F_{{\cal W}_{\star}(s_{\os{\circ}{t}})}$
\begin{equation*} 
\Phi^{\rm abs} \col 
({\cal W}_{\star}A_{X_{\bul \leq N,\os{\circ}{t}}}(E^{\bul \leq N}),P) 
\lo RF^{\rm abs}_{X_{\bul \leq N,\os{\circ}{t}/{\cal W}(s_{\os{\circ}{t}})}*}
(({\cal W}_{\star}A_{X_{\bul \leq N,\os{\circ}{t}}}(E^{\bul \leq N}),P)) 
\tag{2.3.32.1}\label{eqn:ebwnp}
\end{equation*}
the {\it absolute Frobenius endomorphism} of 
$({\cal W}A_{X_{\bul \leq N,\os{\circ}{t}}}(E^{\bul \leq N}),P)$.
When $E^{\bul \leq N}={\cal O}_{\os{\circ}{X}_{\bul \leq N,t}/{\cal W}(\os{\circ}{t})}$, 
we have the following {\it absolute Frobenius endomorphism} 
\begin{equation*} 
\Phi^{\rm abs}{}^* \col 
({\cal W}A_{X_{\bul \leq N,\os{\circ}{t}}}(E^{\bul \leq N}),P) 
\lo RF^{\rm abs}_{X_{\bul \leq N,T_0/T}*}
(({\cal W}A_{X_{\bul \leq N,\os{\circ}{t}}}(E^{\bul \leq N}),P)) 
\tag{2.3.32.2}\label{eqn:abwst}
\end{equation*}
of $({\cal W}A_{X_{\bul \leq N,\os{\circ}{t}}}(E^{\bul \leq N}),P)$.  
\end{defi}  

\begin{rema} 
We leave the formulation of the analogue of 
(\ref{prop:fcar}) for the absolute Frobenius endomorphism. 
\end{rema}

\par 
Let the notations be as in (\ref{defi:bwtd}).  
Then the morphism $\Phi^{\rm abs}$ in (\ref{eqn:ebwnp}) 
is the singlization of the following morphism of filtered double complexes: 
\begin{align*} 
\Phi^{{\rm abs}(\bul *)} \col 
({\cal W}A_{X_{\bul\leq N,\os{\circ}{t}}}(E^{\bul \leq N})^{\bul *},P) 
\lo {\cal W}A_{X_{\bul\leq N,\os{\circ}{t}}}(E^{\bul \leq N})^{\bul *}. 
\tag{2.3.33.1}\label{ali:wtc} 
\end{align*}  
This morphism induces the following morphism 
\begin{align*} 
\Phi^{{\rm abs}(\bul *)} \col 
({\cal W}_nA_{X_{\bul\leq N,\os{\circ}{t}}}(E^{\bul \leq N})^{\bul *},P) 
\lo ({\cal W}_nA_{X_{\bul\leq N,\os{\circ}{t}}}(E^{\bul \leq N})^{\bul *},P). 
\tag{2.3.33.2}\label{ali:wtxtc} 
\end{align*} 
On the other hand, 
we have the following composite morphism 
\begin{align*} 
F\col (E^{\bul \leq N})_{{\cal W}_{n+1}(X_{\bul \leq N,\os{\circ}{t}})}
\os{F_{E^{\bul \leq N}}}{\lo} 
(E^{\bul \leq N})_{{\cal W}_{n+1}(X_{\bul \leq N,\os{\circ}{t}})}
\os{\rm proj}{\lo} 
(E^{\bul \leq N})_{{\cal W}_n(X_{\bul \leq N,\os{\circ}{t}})}.
\end{align*} 
Consequently we have the following morphism 
\begin{align*} 
F:=F\otimes F \col 
({\cal W}_{n+1}A_{X_{\bul\leq N,\os{\circ}{t}}}(E^{\bul \leq N})^{\bul *},P) 
\lo ({\cal W}_nA_{X_{\bul\leq N,\os{\circ}{t}}}(E^{\bul \leq N})^{\bul *},P). 
\tag{2.3.33.3}\label{ali:wtnnc} 
\end{align*}  
By \cite[(9.9)]{ndw} the morphisms $\Phi^{{\rm abs}(\bul *)}$ 
and $F$ fit into the following commutative diagram$:$
\begin{equation*}
\begin{CD}
{\cal W}_{n+1}A_{X_{\bul\leq N,\os{\circ}{t}}}(E^{\bul \leq N})^{\bul *} 
@>{R}>> 
{\cal W}_nA_{X_{\bul\leq N,\os{\circ}{t}}}(E^{\bul \leq N})^{\bul *} \\
@V{p^{\bul}F}VV 
@VV{\Phi^{{\rm abs}(\bul *)}_{{\cal W}_n(t)}}V   \\
{\cal W}_nA_{X_{\bul\leq N,\os{\circ}{t}}}(E^{\bul \leq N})^{\bul *} @= 
{\cal W}_nA_{X_{\bul\leq N,\os{\circ}{t}}}(E^{\bul \leq N})^{\bul *}.
\end{CD}
\tag{2.3.33.4}\label{cd:wpiphipfa}
\end{equation*}

\par
Let $i$ be a fixed nonnegative integer. 
Let $N$ be a nonnegative integer. 
We conclude this section by constructing the 
Poincar\'{e} spectral sequence of 
${\cal W}_{\star}{\Om}^i_{X_{\bul \leq N},\os{\circ}{t}}$ 
and describing the boundary morphisms 
between the $E_1$-terms 
of  the spectral sequences.
\par 
The following is 
the $N$-truncated cosimplicial 
and sheafied version of 
the preweight spectral sequence of 
log Hodge-Witt cohomologies in 
\cite[(4.1)]{ndw}:

\begin{theo}\label{theo:hwwt} 
Let $u\col {\cal W}_{\star}(s_{\os{\circ}{t}})\lo {\cal W}_{\star}(s_{\os{\circ}{t}})$ 
be a morphism of log schemes. Then the following hold$:$
\par 
$(1)$ Let $E^{\bul \leq N}$ be a flat quasi-coherent crystal of 
${\cal O}_{\os{\circ}{X}_{\bul \leq N,t}/{\cal W}_{\star}(\os{\circ}{t})}$-modules.
Then there exists the following spectral sequences for 
$q\in {\mab Z}:$ 
\begin{align*}
& E_1^{-k,q+k}
=\bigoplus_{m=0}^N
\bigoplus_{j\geq \max \{-(k+m),0\}} 
H^{q-i-j-m}(\os{\circ}{X}{}^{(2j+k+m)}_{m,t}, 
(E^m_{\os{\circ}{X}{}^{(2j+k+m)}_{m,t}/
{\cal W}_{\star}(\os{\circ}{t})})_{
{\cal W}_{\star}(\os{\circ}{X}{}^{(2j+k+m)}_{m,t})}
\tag{2.3.34.1}\label{eqn:esplasp}\\
& 
\otimes_{{\cal W}_{\star}
({\cal O}_{\os{\circ}{X}{}^{(2j+k+m)}_{m,t}})}
{\cal W}_{\star}\Om^{i-j-k-m}_{\os{\circ}{X}{}^{(2j+k+m)}_{m,t}}
\otimes_{\mab Z}
\vp^{(2j+k+m)}_{\rm zar}(\os{\circ}{X}_{m,t}/\os{\circ}{t}))(-j-k-m,u) 
\\ 
&  
\Lo 
H^{q-i}(X_{\bul \leq N,t},
(\eps^*_{X_{\bul \leq N,t}/{\cal W}_{\star}(t)}
(E^{\bul \leq N}))_{{\cal W}_{\star}(X_{\bul \leq N,t})}
\otimes_{{\cal W}_{\star}({\cal O}_{X_{\bul \leq N,t}})}
{\cal W}_{\star}{\Om}^i_{X_{\bul \leq N,t}}). 
\end{align*}
\par 
$(2)$ Let $E^{\bul \leq N}$ be a flat quasi-coherent crystal of 
${\cal O}_{\os{\circ}{X}_{\bul \leq N,t}/{\cal W}_{\star}(\os{\circ}{t})}$-modules.
Then there exists the following  spectral sequence$:$
\begin{align*}
& E_1^{-k,q+k}
=\bigoplus_{m\geq 0}^N
\bigoplus_{j\geq \max \{-(k+m),0\}} 
H^{q-i-j-m}(\os{\circ}{X}{}^{(2j+k+m)}_{m,t},
(E^m_{\os{\circ}{X}{}^{(2j+k+m)}_{m,t}/{\cal W}(\os{\circ}{t})})_
{{\cal W}_{\star}(\os{\circ}{X}{}^{(2j+k+m)}_{m,t})} \tag{2.3.34.2}\label{eqn:espmuasp}\\
& \otimes_{{\cal W}_{\star}
({\cal O}_{\os{\circ}{X}{}^{(2j+k+m)}_{m,t}})}
{\cal W}\Om^{i-j-k-m}_{\os{\circ}{X}{}^{(2j+k+m)}_{m,t}}
\otimes_{\mab Z}
\vp^{(2j+k+m)}_{\rm zar}(\os{\circ}{X}_{m,t}/\os{\circ}{t}))(-j-k-m,u)\\ 
& \Lo H^{q-i}(X_{\bul \leq N,t},
(\eps^*_{X_{\bul \leq N,t}/{\cal W}(t)}(E))_{{\cal W}(X_{\bul \leq N,t})}
\otimes_{{\cal W}({\cal O}_{X_{\bul \leq N,t}})}
{\cal W}{\Om}^i_{X_{\bul \leq N,t}}) 
\quad (q\in {\mab Z}). 
\end{align*}
\end{theo}
\begin{proof} 
By (\ref{prop:bcdw}),  
${\cal W}_n\Om_{Y_t}^i={\cal W}_n\Om_{Y_{s_{\os{\circ}{t}}}}^i$.
Hence we may assume that $t=s_{\os{\circ}{t}}$. 
By \cite[3.15]{msemi} and \cite[(6.28) (9), (6.29) (1)]{ndw}, 
the following  sequence
\begin{equation*}
0\lo {\cal W}_{\star}{{\Om}}_{X_{\bul \leq N,\os{\circ}{t}}}^i 
\os{\theta_{\star} \wedge}{\lo} 
{\cal W}_{\star}A_{X_{\bul \leq N,\os{\circ}{t}}}^{i0} 
\os{\theta_{\star} \wedge }{\lo} 
{\cal W}_{\star}A_{X_{\bul \leq N,\os{\circ}{t}}}^{i1}
\os{\theta_{\star} \wedge}{\lo} 
{\cal W}_{\star}A_{X_{\bul \leq N,\os{\circ}{t}}}^{i2}
\os{\theta_{\star} \wedge}{\lo} \cdots 
\tag{2.3.34.3}\label{eqn:evdc}
\end{equation*}
is exact. 
Hence 
the following  sequence
\begin{align*}
& 0\lo 
(\eps^*_{X_{\bul \leq N,\os{\circ}{t}}/{\cal W}_{\star}(t)}(E^{\bul \leq N})
)_{{\cal W}_{\star}(X_{\bul \leq N,\os{\circ}{t}})}
\otimes_{{\cal W}_{\star}({\cal O}_{X_{\bul \leq N,\os{\circ}{t}}})}
{\cal W}_{\star}{{\Om}}_{X_{\bul \leq N,\os{\circ}{t}}}^i 
\os{\theta_{\star} \wedge}{\lo} \tag{2.3.34.4}\label{eqn:eevdc}\\
& 
{\cal W}_{\star}A_{X_{\bul \leq N,\os{\circ}{t}}}^{i0}
(E^{\bul \leq N})
\os{\theta_{\star} \wedge }{\lo} 
{\cal W}_{\star}A_{X_{\bul \leq N,\os{\circ}{t}}}^{i1}
(E^{\bul \leq N}) \os{\theta_{\star} \wedge}{\lo} 
{\cal W}_{\star}A_{X_{\bul \leq N,\os{\circ}{t}}}^{i2}
(E^{\bul \leq N})
\os{\theta_{\star} \wedge}{\lo} \cdots 
\end{align*}
is exact. 
Let us consider the following single complex 
\begin{align*}
& {\cal W}_{\star}A_{X_{\bul \leq N,\os{\circ}{t}}}^{i\bul}
(E^{\bul \leq N})
:=(\cdots 
\os{\theta_{\star} \wedge}{\lo} 
\tag{2.3.34.5}\label{eqn:vst} \\
& (\eps^*_{X_{\bul \leq N,\os{\circ}{t}}/{\cal W}_{\star}(\os{\circ}{t})}
(E^{\bul \leq N}))_{{\cal W}
(X_{\bul \leq N,\os{\circ}{t}})}
\otimes_{{\cal W}_{\star}({\cal O}_{X_{\bul \leq N,\os{\circ}{t}}})}
{\cal W}_{\star}\wt{{\Om}}_{X_{\bul \leq N,\os{\circ}{t}}}^{i+j+1}
/P_j{\cal W}_{\star}\wt{{\Om}}_{X_{\bul \leq N,\os{\circ}{t}}}^{i+j+1} 
\os{\theta_{\star} \wedge}{\lo} 
\cdots)_{j\geq 0},
\end{align*}
and let us define a Poincar\'{e} filtration 
$P=\{P_k\}_{k\in {\mab Z}}$ on 
${\cal W}_{\star}A_{X_{\bul \leq N,\os{\circ}{t}}}^{i\bul}(E^{\bul \leq N})$ 
as follows: 
\begin{align*}
&P_k{\cal W}_{\star}A_{X_{\bul \leq N,\os{\circ}{t}}}^{i\bul}
(E^{\bul \leq N})
= (\cdots \os{\theta_{\star} \wedge}{\lo}
\tag{2.3.34.6}\label{eqn:vstpw}\\
& (P_{2j+k+1}+P_j)
((\eps^*_{X_{\bul \leq N,\os{\circ}{t}}/{\cal W}_{\star}(t)}
(E^{\bul \leq N}))_{{\cal W}_{\star}(X_{\bul \leq N,\os{\circ}{t}})}
\otimes_{{\cal W}_{\star}({\cal O}_{X_{\bul \leq N,\os{\circ}{t}})}}
{\cal W}_{\star}\wt{{\Om}}_{X_{\bul \leq N,\os{\circ}{t}}}^{i+j+1})/P_j  \\
& \os{\theta_{\star} \wedge}{\lo} \cdots)_{j\geq 0}. 
\end{align*} 
By (\ref{eqn:pwa})
we have the following equalities: 
\begin{align*}
&{\rm gr}^{\del(L,P)}_k
Rf_{*}{\cal W}_{\star}A_{X_{\bul \leq N,\os{\circ}{t}}}^{i\bul}
(E^{\bul \leq N}) 
= \bigoplus_{m=0}^N
\bigoplus_{j\geq \max \{-(k+m),0\}} 
\tag{2.3.34.7}\label{eqn:mkres}\\
&{\rm gr}_{2j+k+m+1}^P
Rf_{*}((\eps^*_{X_{m,\os{\circ}{t}}/{\cal W}_{\star}(t)}(E^m))
_{{\cal W}_{\star}(X_{m,\os{\circ}{t}})}
\otimes_{{\cal W}_{\star}({\cal O}_{X_{m,\os{\circ}{t}})}}
{\cal W}_{\star}\wt{{\Om}}_{X_{m,\os{\circ}{t}}}^{i+j+1}\{-j\},-\nabla) \\
&= \bigoplus_{m\geq 0}
\bigoplus_{j\geq \max \{-(k+m),0\}} 
Rf_{*}(\eps^*_{X_{m,t}/{\cal W}_{\star}(\os{\circ}{t})}
(E^{\bul \leq N})_{{\cal W}_{\star}(X_{m,t})}
\otimes_{{\cal W}_{\star}({\cal O}_{X_{m,t}})} 
{\cal W}_{\star}\Om_{\os{\circ}{X}{}^{(2j+k+m)}_{m,t}}^{i-j-k-m}\{-j\},-\nabla). 
\end{align*}
The compatibility of (\ref{eqn:mkres}) with $u$ 
can be proved as in \cite[(9.9)]{ndw}.  
Thus we obtain (\ref{eqn:esplasp}).
\par
By a sheafied version of \cite[(8.6)]{ndw} (2), 
we obtain the complex 
${\cal W}A_{X_{\bul \leq N,\os{\circ}{t}}}^{i\bul}
(E^{\bul \leq N})$ 
by taking the projective limit with respect to 
the projection 
$R \col
{\cal W}_{n+1}A_{X_{\bul \leq N,\os{\circ}{t}}}^{i\bul}
(E^{\bul \leq N})  
\lo {\cal W}_{\star}A_{X_{\bul \leq N,\os{\circ}{t}}}^{i\bul}
(E^{\bul \leq N})$ 
$(n \in {\mab Z}_{>0})$. 
By  \cite[(8.6) (5)]{ndw}  and by the proof of [loc.~cit., (9.9)] 
(cf.~[loc.~cit., (9.11)]), we obtain (\ref{eqn:espmuasp}).
\end{proof}

\begin{defi}\label{defi:pw}
We call (\ref{eqn:esplasp}) the {\it Poincar\'{e} spectral sequence} of 
$$(\eps^*_{X_{\bul \leq N,\os{\circ}{t}}/{\cal W}_{\star}(t)}
(E^{\bul \leq N}))_{{\cal W}_{\star}(X_{\bul \leq N,\os{\circ}{t}})}
\otimes_{{\cal W}_{\star}({\cal O}_{X_{\bul \leq N,\os{\circ}{t}})}}
{\cal W}_{\star}{\Om}_{X_{\bul \leq N,\os{\circ}{t}}}^i.$$   
When 
$E^{\bul \leq N}={\cal O}_{\os{\circ}{X}_{\bul \leq N,t}/{\cal W}_{\star}(\os{\circ}{t})}$,  
we call (\ref{eqn:esplasp}) the 
{\it preweight spectral sequence} of 
$X_{\bul \leq N,\os{\circ}{t}}/{\cal W}_{\star}(t)$.   
We call (\ref{eqn:espmuasp}) the {\it Poincar\'{e} spectral sequence} of 
$$(\eps^*_{X_{\bul \leq N,\os{\circ}{t}}/{\cal W}(t)}(E^{\bul \leq N}))
_{{\cal W}(X_{\bul \leq N,\os{\circ}{t}})}
\otimes_{{\cal W}({\cal O}_{X_{\bul \leq N,\os{\circ}{t}}})}
{\cal W}{\Om}_{X_{\bul \leq N,\os{\circ}{t}}}^i.$$ 
When $E^{\bul \leq N}={\cal O}_{\os{\circ}{X}_{\bul \leq N,t}/{\cal W}(\os{\circ}{t})}$,  
we call (\ref{eqn:espmuasp}) 
the {\it weight spectral sequence} 
of ${\cal W}{\Om}_{X_{\bul \leq N,\os{\circ}{t}}}^i$. 
Let $\star$ be a positive integer $n$ or nothing. 
Let $E^{\bul \leq N}$ be a flat quasi-coherent crystal of 
${\cal O}_{\os{\circ}{X}_{\bul \leq N}/{\cal W}_{\star}(\os{\circ}{s})}$-modules. 
We define the {\it filtration} 
$P$ on 
$$H^{q-i}(X_{\bul \leq N,\os{\circ}{t}},
(\eps^*_{X_{\bul \leq N,\os{\circ}{t}}/{\cal W}_{\star}(t)}
(E^{\bul \leq N}))_{{\cal W}_{\star}(X_{\bul \leq N,\os{\circ}{t}})}
\otimes_{{\cal W}_{\star}({\cal O}_{X_{\bul \leq N,\os{\circ}{t}}})}
{\cal W}_{\star}{\Om}_{X_{\bul \leq N,\os{\circ}{t}}}^i)$$  
as follows: 
\begin{align*} 
& {\rm gr}_k^PH^{q-i}(X_{\bul \leq N,\os{\circ}{t}},
(\eps^*_{X_{\bul \leq N,\os{\circ}{t}}
/{\cal W}_{\star}(\os{\circ}{t})}
(E^{\bul \leq N}))_{{\cal W}_{\star}(X_{\bul \leq N,\os{\circ}{t}})}
\otimes_{{\cal W}_{\star}({\cal O}_{X_{\bul \leq N,\os{\circ}{t}}})}
{\cal W}_{\star}{\Om}_{X_{\bul \leq N,\os{\circ}{t}}}^i) \\ 
&=
E_{\infty}^{h-k, k}(
(\eps^*_{X_{\bul \leq N,\os{\circ}{t}}/{\cal W}_{\star}(\os{\circ}{t})}
(E^{\bul \leq N}))_{{\cal W}_{\star}(X_{\bul \leq N,\os{\circ}{t}})}
\otimes_{{\cal W}_{\star}({\cal O}_{X_{\bul \leq N,\os{\circ}{t}}})}
{\cal W}_{\star}{\Om}_{X_{\bul \leq N,\os{\circ}{t}}}^i).
\end{align*} 
\end{defi}

As in (\ref{prop:deccbd}), we have the following: 

\begin{prop}\label{prop:deschwbd} 
Let $\star$ be a positive integer $n$ or nothing.  
Then the boundary morphism between the $E_1$-terms of 
the spectral sequences {\rm (\ref{eqn:esplasp})} and 
{\rm (\ref{eqn:espmuasp})} are given by the following diagram$:$ 
\begin{align*} 
& \bigoplus_{m\geq 0}
\bigoplus_{j\geq \max \{-(k+m),0\}} 
H^{q-i-j-m}(\os{\circ}{X}{}^{(2j+k+m)}_{m+1,t},
(E^m_{\os{\circ}{X}{}^{(2j+k+m)}_{m+1,t}/{\cal W}_{\star}(\os{\circ}{t})}
\tag{2.3.36.1}\label{cd:gbhwsd}\\
&\otimes_{
{\cal W}_{\star}({\cal O}_{\os{\circ}{X}{}^{(2j+k+m)}_{m+1,t}})}
{\cal W}_{\star}\Om^{i-j-k-m}_{
\os{\circ}{X}{}^{(2j+k+m)}_{m+1,t}}
\otimes_{\mab Z}
\vp^{(2j+k+m)}_{\rm zar}
(\os{\circ}{X}_{m+1,t}/\os{\circ}{t}))(-j-k-m,u) 
\end{align*}  
$$\text{\scriptsize
{$\sum_{i= 0}^{t_j+1}(-1)^i\del^i_{j}$}}
~\uparrow~(1 \leq j \leq r)$$
\begin{align*} 
&\bigoplus_{m\geq 0}
\bigoplus_{j\geq \max \{-(k+m),0\}} 
H^{q-i-j-m}(\os{\circ}{X}{}^{(2j+k+m)}_{m,t},
(E^m_{\os{\circ}{X}{}^{(2j+k+m)}_{m,t}/{\cal W}_{\star}}\\
&\otimes_{
{\cal W}_{\star}({\cal O}_{\os{\circ}{X}{}^{(2j+k+m)}_{m,t}})}
{\cal W}_{\star}\Om^{i-j-k-m}_{\os{\circ}{X}{}^{(2j+k+m)}_{m,t}}
\otimes_{\mab Z}
\vp^{(2j+k+m)}_{\rm zar}
(\os{\circ}{X}_{m,t}/\os{\circ}{t}))(-j-k-m,u) 
\end{align*}  
$$\text{\scriptsize
{$(-1)^m{\sum_{j\geq {\rm max}\{-(k+m),0\}}
[G_m+\rho_m]}$}}
~\downarrow$$
\begin{align*} 
&\bigoplus_{m\geq 0}
\bigoplus_{j\geq \max \{-(k+m),0\}} 
H^{q-i-j-m+1}(\os{\circ}{X}{}^{(2j+k+m)}_{m,t},
(E^m_{\os{\circ}{X}{}^{(2j+k+m)}_{m,t}/{\cal W}_{\star}(\os{\circ}{t})}\\
&\otimes_{
{\cal W}_{\star}({\cal O}_{\os{\circ}{X}{}^{(2j+k+m)}_{m,t}})}
{\cal W}_{\star}\Om^{i-j-k-m+1}_{\os{\circ}{X}{}^{(2j+k+m)}_{m,t}}
\otimes_{\mab Z}
\vp^{(2j+k+m)}_{\rm zar}
(\os{\circ}{X}_{m,t}/\os{\circ}{t}))(-j-k-m+1,u).  
\end{align*}  
Here $G_m$ is the \v{C}ech-Gysin morphism in Hodge-Witt cohomologies $($cf.~{\rm \cite[(4.4.2)]{ndw}}$)$  
which is defined as in {\rm (\ref{eqn:togsn})} 
and $\rho_m$ is 
the morphism obtained by the pull-back of the closed immersions in Hodge-Witt cohomologies 
which is defined as in {\rm (\ref{eqn:rhogsn})}.  
\end{prop}

\begin{prop}\label{prop:pwt}
Set 
\begin{equation*} 
E_1^{-k,q+k} :=
\bigoplus_{\max \{-k,0\}\leq m\leq N}
H^{q-m}(X_m,{\rm gr}^P_{k+m}
(\eps^*_{X_{m,\os{\circ}{t}}/{\cal W}_{\star}(\os{\circ}{t})}
(E^m))_{{\cal W}_{\star}(X_{m,\os{\circ}{t}})}
\otimes_{{\cal W}_{\star}({\cal O}_{X_{m,\os{\circ}{t}}})}
{\cal W}\wt{\Om}{}^{\bul}_{X_{m,\os{\circ}{t}_0}})).\\ 
\end{equation*}
Then there exists the following spectral sequence 
\begin{align*} 
E_1^{-k,q+k}\Lo 
H^q(X_{\bul \leq N,\os{\circ}{t}_0},
(\eps^*_{X_{\bul \leq N,\os{\circ}{t}_0}/{\cal W}_{\star}(\os{\circ}{t})}
(E^{\bul \leq N}))_{{\cal W}_{\star}(X_{\bul \leq N,\os{\circ}{t}})}
\otimes_{{\cal W}_{\star}({\cal O}_{X_{\bul \leq N,\os{\circ}{t}}})}
{\cal W}\wt{\Om}{}^{\bul}_{X_{\bul \leq N,\os{\circ}{t}_0}}). 
\tag{2.3.37.1}\label{ali:xnat}
\end{align*} 
Explicitly 
\begin{align*}
E_1^{-k,q+k}=
\bigoplus_{m=0}^N
H^{q-2m-k}&((\os{\circ}{X}{}^{(k+m-1)}_{m,t_0}/{\cal W}_{\star}(\os{\circ}{t}))_{\rm crys}, 
E^m_{\os{\circ}{X}{}^{(k+m-1)}_{m,t_0}/{\cal W}_{\star}(\os{\circ}{t})}\\
&\otimes_{\mab Z}
\vp^{(k+m-1)}_{\rm crys}(\os{\circ}{X}_{m,t_0}/{\cal W}_{\star}(\os{\circ}{t})))(-k-m,u)
\end{align*}
for $k>0$, 
\begin{align*}
E_1^{-k,q+k}&=
H^{q+k}(X_{-k,\os{\circ}{t}},
P_0(\eps^*_{X_{-k,\os{\circ}{t}}/{\cal W}_{\star}(\os{\circ}{t})}
(E^{-k})_{{\cal W}_{\star}(X_{-k,\os{\circ}{t}})}
\otimes_{{\cal W}_{\star}({\cal O}_{X_{-k,\os{\circ}{t}}})}
{\cal W}\wt{\Om}{}^{\bul}_{X_{-k,\os{\circ}{t}_0}}))\\
&\oplus \bigoplus_{m=-k+1}^N
H^{q-2m-k}((\os{\circ}{X}{}^{(k+m-1)}_{m,t_0}/{\cal W}_{\star}(\os{\circ}{t}))_{\rm crys},\\
&(E^m_{\os{\circ}{X}{}^{(k+m-1)}_{m,t_0}/\os{\circ}{t}}
\otimes_{\mab Z}
\vp^{(k+m-1)}_{\rm zar}(\os{\circ}{X}_{m,t_0}/\os{\circ}{t}))(-k-m,u)
\end{align*} 
for $-N\leq k\leq 0$ 
and 
\begin{equation*}
E_1^{-k,q+k}=0
\end{equation*}
for $k<-N$. 
There also exists the following spectral sequence 
\begin{align*}
E_1^{i,q+k-i}&=H^{q+k-i}((\os{\circ}{X}{}^{(i)}_{-k,t_0}
/{\cal W}_{\star}(\os{\circ}{t}))_{\rm crys},
E^{-k}_{\os{\circ}{X}{}^{(i)}_{-k,t_0}/{\cal W}_{\star}(\os{\circ}{t})}
\otimes_{\mab Z}
\vp^{(i)}_{\rm crys}(\os{\circ}{X}_{m,t_0}/\os{\circ}{t}))\\
&\Lo 
H^{q+k}(X_{-k,\os{\circ}{t}},
P_0(\eps^*_{X_{-k,\os{\circ}{t}}/{\cal W}_{\star}(\os{\circ}{t})}
(E^{-k})_{{\cal W}_{\star}(X_{-k,\os{\circ}{t}})}
\otimes_{{\cal W}_{\star}({\cal O}_{X_{-k,\os{\circ}{t}}})}
{\cal W}\wt{\Om}{}^{\bul}_{X_{-k,\os{\circ}{t}_0}})).
\end{align*} 
\end{prop} 
These spectral sequences are contravariantly functorial for the morphism 
$g_{\bul \leq N}$ satisfying {\rm (1.5.6.4)} and {\rm (1.5.6.5)} 
in the case $S=s$ and $S'=s'$. 
\begin{proof}
This follows from (\ref{ali:pdete}) and (\ref{eqn:pwa}). 
\end{proof}

\begin{prop}\label{prop:ihw}
Let $i$ be a nonnegative integer. 
Set 
\begin{equation*} 
E_1^{-k,q+k} :=
\bigoplus_{\max \{-k,0\}\leq m\leq N}
H^{q-m}(X_m,{\rm gr}^P_{k+m}
(\eps^*_{X_{m,\os{\circ}{t}}/{\cal W}_{\star}(\os{\circ}{t})}
(E^m))_{{\cal W}_{\star}(X_{m,\os{\circ}{t}})}
\otimes_{{\cal W}_{\star}({\cal O}_{X_{m,\os{\circ}{t}}})}
{\cal W}\wt{\Om}{}^i_{X_{m,\os{\circ}{t}_0}}))\\ 
\end{equation*}
Then there exists the following spectral sequence 
\begin{align*} 
E_1^{-k,q+k}\Lo 
H^q(X_{\bul \leq N,\os{\circ}{t}_0},
(\eps^*_{X_{\bul \leq N,\os{\circ}{t}_0}/{\cal W}_{\star}(\os{\circ}{t})}
(E^{\bul \leq N}))_{{\cal W}_{\star}(X_{\bul \leq N,\os{\circ}{t}})}
\otimes_{{\cal W}_{\star}({\cal O}_{X_{\bul \leq N,\os{\circ}{t}}})}
{\cal W}\wt{\Om}{}^i_{X_{\bul \leq N,\os{\circ}{t}_0}}). 
\tag{2.3.38.1}\label{ali:xwnt}
\end{align*} 
Here 
\begin{align*}
E_1^{-k,q+k}=
\bigoplus_{m=0}^N
H^{q-2m-k}&(\os{\circ}{X}{}^{(k+m-1)}_{m,t_0}, 
E^m_{\os{\circ}{X}{}^{(k+m-1)}_{m,t_0}/{\cal W}_{\star}(\os{\circ}{t})}
\otimes_{{\cal W}_{\star}({\cal O}_{\os{\circ}{X}{}^{(k+m-1)}_{m,\os{\circ}{t}}})}\\
&{\cal W}{\Om}{}^{i-m-k}_{\os{\circ}{X}{}^{(k+m-1)}_{m,\os{\circ}{t}_0}}
\otimes_{\mab Z}
\vp^{(k+m-1)}_{\rm zar}(\os{\circ}{X}_{m,t_0}/{\cal W}_{\star}(\os{\circ}{t})))(-k-m,u)
\end{align*}
for $k>0$, 
\begin{align*}
E_1^{-k,q+k}&=
H^{q+k}(X_{-k,\os{\circ}{t}},
P_0(\eps^*_{X_{-k,\os{\circ}{t}}/{\cal W}_{\star}(\os{\circ}{t})}
(E^{-k})_{{\cal W}_{\star}(X_{-k,\os{\circ}{t}})}
\otimes_{{\cal W}_{\star}({\cal O}_{X_{-k,\os{\circ}{t}}})}
{\cal W}\wt{\Om}{}^{\bul}_{X_{-k,\os{\circ}{t}_0}}))\\
&\oplus \bigoplus_{m=-k+1}^N
H^{q-2m-k}(\os{\circ}{X}{}^{(k+m-1)}_{m,t_0}, 
(E_{\os{\circ}{X}{}^{(k+m-1)}_{m,t_0}/\os{\circ}{t}}
\otimes_{{\cal W}_{\star}({\cal O}_{\os{\circ}{X}{}^{(k+m-1)}_{m,\os{\circ}{t}}})}
{\cal W}{\Om}{}^{i-m-k}_{\os{\circ}{X}{}^{(k+m-1)}_{m,\os{\circ}{t}_0}}\\
&\otimes_{\mab Z}
\vp^{(k+m-1)}_{\rm zar}(\os{\circ}{X}_{m,t_0}/\os{\circ}{t}))(-k-m,u)
\end{align*} 
for $-N\leq k\leq 0$ 
and 
\begin{equation*}
E_1^{-k,q+k}=0
\end{equation*}
for $k<-N$. 
There also exists the following spectral sequence 
\begin{align*}
E_1^{j,q+k-j}&=H^{q+k-j}((\os{\circ}{X}{}^{(j)}_{-k,t_0}
/{\cal W}_{\star}(\os{\circ}{t}))_{\rm crys},
E_{\os{\circ}{X}{}^{(j)}_{-k,t_0}/{\cal W}_{\star}(\os{\circ}{t})}
\otimes_{\mab Z}
\vp^{(i)}_{\rm crys}(\os{\circ}{X}_{m,t_0}/\os{\circ}{t}))\\
&\Lo 
H^{q+k}(X_{-k,\os{\circ}{t}},
P_0(\eps^*_{X_{-k,\os{\circ}{t}}/{\cal W}_{\star}(\os{\circ}{t})}
(E^{-k})_{{\cal W}_{\star}(X_{-k,\os{\circ}{t}})}
\otimes_{{\cal W}_{\star}({\cal O}_{X_{-k,\os{\circ}{t}}})}
{\cal W}\wt{\Om}{}^i_{X_{-k,\os{\circ}{t}_0}})).
\end{align*} 
These spectral sequences are contravariantly functorial for the morphism 
$g_{\bul \leq N}$ satisfying {\rm (1.5.6.4)} and {\rm (1.5.6.5)} 
in the case $S=s$ and $S'=s'$. 
\end{prop}
\begin{proof} 
This follows from (\ref{ali:pdete}) and (\ref{eqn:pwa}). 
\end{proof}

\section{$p$-adic monodromy operators via log de Rham-Witt complexes}\label{sec:pmr}
In this section we give 
generalizations of 
results in \cite{hdw}, \cite{hk}, \cite{msemi} and \cite{ndw} about 
$p$-adic monodromy operators for flat coherent log crystals. 
The results in this section are easy corollaries of results in previous sections. 
\par 
Let $\kap$, ${\cal W}_n$, $s$, $t$, ${\cal W}_n(s)$ 
and ${\cal W}(s)$, ${\cal W}_n(\os{\circ}{t})$ and ${\cal W}(\os{\circ}{t})$ 
be as in the previous section. 
Let $Y_{\bul \leq N}/s$ be 
a log smooth integral $N$-truncated simplicial log scheme over $s$ 
with structural morphism $g \col Y_{\bul \leq N} \lo s$. 
Let us also denote by $g$ the structural morphism 
$Y_{\bul \leq N}\lo s\lo {\cal W}_n(s)$. 
Assume that $Y_{\bul \leq N}$ has the disjoint union 
$Y_{\bul \leq N}'$
of the member of an affine $N$-truncated simplicial open covering of 
$Y_{\bul \leq N}$. 
Let $Y_{\bul \leq N,\bul}$ be the \v{C}ech diagram of 
$Y_{\bul \leq N}'$ over $Y_{\bul \leq N}$. 
Let $\ol{F}{}^{\bul \leq N}$ be a flat coherent log crystal of 
${\cal O}_{Y_{\bul \leq N,\os{\circ}{t}}/{\cal W}_n(\os{\circ}{t})}$-modules. 
Set 
\begin{align*} 
F^{\bul \leq N}:=
\eps^*_{Y_{\bul \leq N,t}/{\cal W}_n(t)/{\cal W}_n(\os{\circ}{t})}(\ol{F}{}^{\bul \leq N}).
\tag{2.4.0.1}\label{ali:fdf} 
\end{align*}  
Set 
\begin{align*} 
\ol{F}{}^{\bul \leq N}_n:=(\ol{F}{}^{\bul \leq N})_{{\cal W}_n(Y_{\bul \leq N,\os{\circ}{t}})} 
\tag{2.4.0.2}\label{ali:qofdf} 
\end{align*} 
and 
\begin{align*} 
F^{\bul \leq N}_{n}:=(F^{\bul \leq N})_{{\cal W}_n(Y_{\bul \leq N,t})}.
\tag{2.4.0.3}\label{ali:qntdf} 
\end{align*} 
Then $\ol{F}{}^{\bul \leq N}_{n}=F^{\bul \leq N}_{n}$. 
Let 
\begin{align*} 
\ol{\nabla}^{\bul \leq N,\bul}
\col 
\ol{F}{}^{\bul \leq N}_{n}\lo \ol{F}{}^{\bul \leq N}_{n}
\otimes_{{\cal W}_n({\cal O}_{Y_{\bul \leq N,\os{\circ}{t}}})}
{\cal W}_n\wt{{\Om}}^1_{Y_{\bul \leq N,\os{\circ}{t}}}
\tag{2.4.0.4}\label{eqn:nidfwtopltd}
\end{align*}  
and 
\begin{align*} 
\nabla^{\bul \leq N,\bul}
\col 
F^{\bul \leq N}_n
\lo F^{\bul \leq N}_n\otimes_{{\cal W}_n({\cal O}_{Y_{\bul \leq N,t}})}
{\cal W}_n{\Om}^1_{Y_{\bul \leq N,t}}
\tag{2.4.0.5}\label{eqn:nidfwopltd}
\end{align*}  
be the natural connections obtained by the log crystals $\ol{F}{}^{\bul \leq N}$
and $F^{\bul \leq N}$, respectively (\S\ref{sec:ldrwc}). 
\par 
Let $i$ be a nonnegative integer. 
By (\ref{cd:oreeam}) we have the following exact sequences:   
\begin{align*} 
0 &\lo F^{\bul \leq N}_n\otimes_{{\cal W}_n({\cal O}_{Y_{\bul \leq N,t}})'}
({\cal W}_n\Om^{\bul \geq i-1}_{Y_{\bul \leq N,t}})'\{-(i-1)\}(-1,u)[-1] 
 \tag{2.4.0.6}\label{eqn:wby}\\ 
&\os{(1\otimes \theta'_n) \wedge}{\lo} 
\ol{F}{}^{\bul \leq N}_{n}\otimes_{{\cal W}_n({\cal O}_{Y_{\bul \leq N,\os{\circ}{t}}})'} 
({\cal W}_n\wt{\Om}^{\bul \geq i}_{Y_{\bul \leq N,\os{\circ}{t}}})'\{-i\} \\
&\lo F^{\bul \leq N}_{n}\otimes_{{\cal W}_n({\cal O}_{Y_{\bul \leq N,t}})'}
({\cal W}_n\Om^{\bul \geq i}_{Y_{\bul \leq N,t}})'\{-i\} \lo 0,    
\end{align*}
\begin{align*} 
0 &\lo F^{\bul \leq N}_{n}
\otimes_{{\cal W}_n({\cal O}_{Y_{\bul \leq N,t}})}
{\cal W}_n\Om^{\bul \geq i-1}_{Y_{\bul \leq N,t}}\{-(i-1)\}(-1,u)[-1] \tag{2.4.0.7}\label{eqn:wnlby}\\ 
&\os{(1\otimes\theta_n) \wedge}{\lo}  
\ol{F}{}^{\bul \leq N}_{n}\otimes_{{\cal W}_n({\cal O}_{Y_{\bul \leq N,\os{\circ}{t}}})} 
{\cal W}_n\wt{\Om}^{\bul \geq i}_{Y_{\bul \leq N,\os{\circ}{t}}}\{-i\} \\
& \lo 
F^{\bul \leq N}_n\otimes_{{\cal W}_n({\cal O}_{Y_{\bul \leq N,t}})}
{\cal W}_n\Om^{\bul \geq i}_{Y_{\bul \leq N,t}}\{-i\} \lo 0.   
\end{align*}
\parno 
By (\ref{eqn:wby}) we have the following boundary morphism 
\begin{align*} 
(N^{\geq i}_{{\rm dRW},n})' \col  &F^{\bul \leq N}_{n}
\otimes_{{\cal W}_n({\cal O}_{Y_{\bul \leq N,t}})'}
({\cal W}_n{{\Om}}^{\bul \geq i}_{Y_{\bul \leq N,t}})'\{-i\} \tag{2.4.0.8}\label{ali:ndw}\\
& \lo 
F^{\bul \leq N}_{n}\otimes_{{\cal W}_n({\cal O}_{Y_{\bul \leq N,\os{\circ}{t}}})'}
({\cal W}_n\Om^{\bul \geq i-1}_{Y_{\bul \leq N,t}})'\{-(i-1)\}(-1,u). 
\end{align*}
By (\ref{eqn:wnlby}) we have the following boundary morphism 
\begin{align*} 
N^{\geq i}_{{\rm dRW},n} &\col F^{\bul \leq N}_{n}
\otimes_{{\cal W}_n({\cal O}_{Y_{\bul \leq N,t}})}
{\cal W}_n\Om_{Y_{\bul \leq N,t}}^{\bul \geq i}\{-i\} 
\tag{2.4.0.9}\label{eqn:ndrwnp}\\ 
& \lo F^{\bul \leq N}_{n}
\otimes_{{\cal W}_n({\cal O}_{Y_{\bul \leq N,t}})}
{\cal W}_n\Om_{Y_{\bul \leq N,t}}^{\bul \geq i-1}\{-(i-1)\}(-1,u). 
\end{align*}
\parno 
Set 
$F^{\bul \leq N}_{\infty}:=\vpl_n(F^{\bul \leq N}_{n})$. 
By (\ref{eqn:wby}) we have the following boundary morphisms: 
\begin{align*} 
(N^{\geq i}_{{\rm dRW}})' \col  &F^{\bul \leq N}_{\infty}
\otimes_{{\cal W}({\cal O}_{Y_{\bul \leq N,t}})}
({\cal W}{{\Om}}^{\bul \geq i}_{Y_{\bul \leq N,t}})'\{-i\} \lo 
F^{\bul \leq N}_{\infty}
\otimes_{{\cal W}({\cal O}_{Y_{\bul \leq N,t}})}
({\cal W}\Om^{\bul \geq i-1}_{Y_{\bul \leq N,t}})'\{-(i-1)\}. 
\tag{2.4.0.10}\label{ali:nldw}
\end{align*}
By (\ref{eqn:wnlby}) we have the following boundary morphism 
\begin{align*} 
N^{\geq i}_{{\rm dRW}} \col F^{\bul \leq N}_{\infty}
\otimes_{{\cal W}({\cal O}_{Y_{\bul \leq N,t}})}
{\cal W}\Om_{Y_{\bul \leq N,t}}^{\bul \geq i}\{-i\} \lo 
F^{\bul \leq N}_{\infty}
\otimes_{{\cal W}({\cal O}_{Y_{\bul \leq N,t}})}
{\cal W}\Om_{Y_{\bul \leq N,t}}^{\bul \geq i-1}\{-(i-1)\}(-1,u). 
\tag{2.4.0.11}\label{eqn:nlwnp}
\end{align*}
For a flat coherent crystal $\ol{F}{}^{\bul \leq N}$ of 
${\cal O}_{Y_{\bul \leq N,\os{\circ}{t}}/{\cal W}_{\star}(\os{\circ}{t})}$-modules, 
via the identifications 
\begin{align*} 
F^{\bul \leq N}_{n}
\otimes_{{\cal W}_n({\cal O}_{Y_{\bul \leq N,t}})'}
({\cal W}_n{\Om}_{Y_{\bul \leq N,t}}^{\bul \geq i})'\{-i\} 
\os{\sim}{\lo} 
F^{\bul \leq N}_{n}
\otimes_{{\cal W}_n({\cal O}_{Y_{\bul \leq N,t}})}
{\cal W}_n{\Om}_{Y_{\bul \leq N,t}}^{\bul \geq i}\{-i\},  
\end{align*}
\begin{align*} 
F^{\bul \leq N}_{\infty}
\otimes_{{\cal W}({\cal O}_{Y_{\bul \leq N,t}})'}
({\cal W}{\Om}_{Y_{\bul \leq N,t}}^{\bul \geq i})'\{-i\}
\os{\sim}{\lo} 
F^{\bul \leq N}_{\infty}
\otimes_{{\cal W}({\cal O}_{Y_{\bul \leq N,t}})}
{\cal W}{\Om}_{Y_{\bul \leq N,t}}^{\bul \geq i}\{-i\},  
\end{align*}
$(N^{\geq i}_{{\rm dRW},\star})'=N^{\geq i}_{{\rm dRW},\star}$ 
$( \star=n$ or nothing$)$.

\begin{defi} 
We call $(N^{\geq i}_{{\rm dRW},n})'$ and $N^{\geq i}_{{\rm dRW},n}$ the 
{\it obverse log de Rham-Witt  monodromy operator} of $\ol{F}{}^{\bul \leq N}_{n}$ 
and the  
{\it reverse log de Rham-Witt  monodromy operator} of $\ol{F}{}^{\bul \leq N}_{n}$, 
respectively. 
We call $(N^{\geq i}_{{\rm dRW}})'$ and $N^{\geq i}_{{\rm dRW}}$ the 
{\it obverse log de Rham-Witt  monodromy operator} of 
$F^{\bul \leq N}_{\infty}$ 
and the  
{\it reverse log de Rham-Witt  monodromy operator} of 
$F^{\bul \leq N}_{\infty}$, 
respectively. 
When $i=0$, we denote 
$(N^{\geq i}_{{\rm dRW},{\star}})'$ and $(N^{\geq i}_{{\rm dRW},{\star}})$ 
$(\star=n$ or nothing) simply by 
$(N_{{\rm dRW},{\star}})'$ and $(N_{{\rm dRW},{\star}})$, respectively. 
\end{defi}

\begin{prop}\label{prop:sld} 
Assume that $F^{\bul \leq N}$ is a unit root log $F$-crystal. 
In this proposition, let $s'$, $t'$, $Y'_{\bul \leq N}$ and $F'{}^{\bul \leq N}$ be 
similar objects to $s$, $t$, $Y^{\bul \leq N}$ and $F^{\bul \leq N}$, respectively. 
Let 
\begin{equation*} 
\begin{CD} 
Y_{\bul \leq N,\os{\circ}{t}}@>{g}>>Y'_{\bul \leq N,\os{\circ}{t}{}'}\\
@VVV @VVV \\
s_{\os{\circ}{t}}@>{u_s}>>s'_{\os{\circ}{t}{}'} \\
@V{\bigcap}VV @VV{\bigcap}V \\
{\cal W}(s_{\os{\circ}{t}})@>{u}>>{\cal W}(s'_{\os{\circ}{t}{}'}) 
\end{CD}
\tag{2.4.2.1}\label{eqn:gwwsp}
\end{equation*} 
be a commutative diagram. 
Let $\star$ be nothing or $\prime$. 
Then the monodromy operator 
\begin{align*} 
(N^{\geq i}_{\rm dRW})^{\star}
\col &
H^{q-i}(Y_{\bul \leq N,t},F^{\bul \leq N}_{\infty}
\otimes_{{\cal W}({\cal O}_{Y_{\bul \leq N,t}})^{\star}}
({\cal W}{\Om}_{Y_{\bul \leq N,t}}^{\bul \geq i})^{\star})
\lo \tag{2.4.2.2}\label{ali:igeq}\\
& H^{q-i+1}(Y_{\bul \leq N,t},F^{\bul \leq N}_{\infty}
\otimes_{{\cal W}({\cal O}_{Y_{\bul \leq N,t}})^{\star}}
({\cal W}{\Om}_{Y_{\bul \leq N,t}}^{\bul \geq i-1})^{\star})(-1,u)
\end{align*} 
induces the following morphism 
\begin{align*} 
(N^{q-i,i}_{\rm dRW})^{\star}
\col &
H^{q-i}(Y_{\bul \leq N,t},F^{\bul \leq N}_{\infty}
\otimes_{{\cal W}({\cal O}_{Y_{\bul \leq N,t}})^{\star}}
({\cal W}{\Om}_{Y_{\bul \leq N,t}}^i)^{\star})\otimes_{\cal W}K_0
\lo \tag{2.4.2.3}\label{ali:igeqoi}\\
&H^{q-i+1}(Y_{\bul \leq N,t},F^{\bul \leq N}_{\infty}
\otimes_{{\cal W}({\cal O}_{Y_{\bul \leq N,t}})^{\star}}
({\cal W}{\Om}_{Y_{\bul \leq N,t}}^{i-1})^{\star})(-1,u)
\otimes_{\cal W}K_0.  
\end{align*} 
\end{prop} 
\begin{proof} 
The following commutative diagram 
\begin{equation*} 
\begin{CD} 
H^{q-i-1}(Y_{\bul \leq N,t},F^{\bul \leq N}_{\infty}
\otimes_{{\cal W}({\cal O}_{Y_{\bul \leq N,t}})^{\star}}
({\cal W}{\Om}_{Y_{\bul \leq N,t}}^{\bul \geq i+1})^{\star})
@>{(N^{\geq i+1}_{\rm dRW})^{\star}}>>\\
@VVV \\
H^{q-i}(Y_{\bul \leq N,t},F^{\bul \leq N}_{\infty}
\otimes_{{\cal W}({\cal O}_{Y_{\bul \leq N,t}})^{\star}}
({\cal W}{\Om}_{Y_{\bul \leq N,t}}^{\bul \geq i})^{\star})
@>{(N^{\geq i}_{\rm dRW})^{\star}}>> \\
\end{CD}
\tag{2.4.2.4}\label{eqn:gwdrwsp}
\end{equation*} 
\begin{equation*} 
\begin{CD} 
H^{q-i}(Y_{\bul \leq N,t},F^{\bul \leq N}_{{\infty}}
\otimes_{{\cal W}({\cal O}_{Y_{\bul \leq N,t}})^{\star}}
({\cal W}{\Om}_{Y_{\bul \leq N,t}}^{\bul \geq i})^{\star})
(-1,u) \\
@VVV \\
H^{q-i+1}(Y_{\bul \leq N,t},
F^{\bul \leq N}_{\infty}\otimes_{{\cal W}({\cal O}_{Y_{\bul \leq N,t}})^{\star}}
({\cal W}{\Om}_{Y_{\bul \leq N,t}}^{\bul \geq i-1})^{\star})
(-1,u)
\end{CD} 
\end{equation*} 
and (\ref{theo:e1deg}) show (\ref{prop:sld}). 
\end{proof}

\parno  
Let $\ol{F}{}^{\bul \leq N}$ be a flat coherent crystal of 
${\cal O}_{Y_{\bul \leq N,\os{\circ}{t}}/{\cal W}_n(\os{\circ}{t})}$-modules.  
By (\ref{eqn:ywnnny}) (cf.~\cite[(7.19)]{ndw}(=a correction of \cite[(4.19)]{hk})), 
we have the following canonical isomorphism 
\begin{align*}
&Ru_{Y_{\bul \leq N,t}/{\cal W}_n(t)*}(F^{\bul \leq N})
\os{\sim}{\lo}
F^{\bul \leq N}_{n}\otimes_{{\cal W}_n({\cal O}_{Y_{\bul \leq N,t}})^{\star}}
({\cal W}_n\Om^{\bul}_{Y_{\bul \leq N},t})^{\star}. 
\tag{2.4.2.5}\label{eqn:uxwns}
\end{align*}

By (\ref{cd:oreedsam}) we obtain the following: 
\begin{prop}\label{prop:tfdic} 
The following diagram is commutative$:$
\begin{equation*} 
\begin{CD}
Ru_{Y_{\bul \leq N,t}/{\cal W}_n(t)*}(F^{\bul \leq N})
@>{N_{\rm zar}}>> Ru_{Y_{\bul \leq N,t}/{\cal W}_n(t)*}(F^{\bul \leq N})\\
@V{\simeq}VV @VV{\simeq}V \\ 
F^{\bul \leq N}_{n}
\otimes_{{\cal W}_n({\cal O}_{Y_{\bul \leq N,t}})^{\star}}
({\cal W}_n\Om^{\bul}_{Y_{\bul \leq N,t}})^{\star}
@>{(N_{{\rm dRW},n})^{\star}}>> F^{\bul \leq N}_{n}
\otimes_{{\cal W}_n({\cal O}_{Y_{\bul \leq N,t}})^{\star}}
({\cal W}_n\Om^{\bul}_{Y_{\bul \leq N,t}})^{\star}. 
\end{CD}
\tag{2.4.3.1}\label{eqn:uywns}
\end{equation*}
The commutative diagram {\rm (\ref{eqn:uywns})} 
is compatible with the two projections of  both hand sides on {\rm (\ref{eqn:uywns})}. 
\end{prop}

\par 
Let the notation be as in the beginning of this section. 
Let $e$ be a positive integer. 
Let $\ol{F}{}^{\bul \leq N}$ be a flat coherent crystal  of 
${\cal O}_{Y_{\bul \leq N,\os{\circ}{t}}/{\cal W}_{\star}(\os{\circ}{t})}$-modules.  
By (\ref{cd:liemext}) we have the following exact sequences:   
\begin{align*} 
0 &\lo 
F^{\bul \leq N}_{\infty}
\otimes_{{\cal W}({\cal O}_{Y_{\bul \leq N,t}})'}
({\cal W}\Om_{Y_{\bul \leq N,t}}^{\bul \geq i-1})'(-1,u)
\{-(i-1)\}[-1]\otimes_{\mab Z}{\mab Q} 
\os{(1\otimes e\theta) \wedge}{\lo}  
\tag{2.4.3.2}\label{eqn:wfnby}\\ 
&F^{\bul \leq N}_{\infty}
\otimes_{{\cal W}({\cal O}_{Y_{\bul \leq N,\os{\circ}{t}}})'}
({\cal W}\wt{\Om}^{\bul \geq i}_{Y_{\bul \leq N,\os{\circ}{t}}})' \{-i\}
\otimes_{\mab Z}{\mab Q} \\
&\lo  F^{\bul \leq N}_{\infty}\otimes_{{\cal W}({\cal O}_{Y_{\bul \leq N,t}})'}
({\cal W}\Om^{\bul \geq i}_{Y_{\bul \leq N,t}})'\{-i\}\otimes_{\mab Z}{\mab Q}  
\lo 0,    
\end{align*}
\begin{align*} 
0 &\lo F^{\bul \leq N}_{\infty}
\otimes_{{\cal W}({\cal O}_{Y_{\bul \leq N,t}})}
{\cal W}\Om_{Y_{\bul \leq N,t}}^{\bul \geq i-1}(-1,u)\{-(i-1)\}[-1]
\otimes_{\mab Z}{\mab Q}  
\os{(1\otimes e\theta)\wedge}{\lo}  
\tag{2.4.3.3}\label{eqn:wenlby}\\ 
& F^{\bul \leq N}_{\infty}
\otimes_{{\cal W}({\cal O}_{Y_{\bul \leq N,\os{\circ}{t}}})}
{\cal W}\wt{\Om}^{\bul \geq i}_{Y_{\bul \leq N,\os{\circ}{t}}}\{-i\}
\otimes_{\mab Z}{\mab Q} \\
&\lo 
F^{\bul \leq N}_{\infty}\otimes_{{\cal W}({\cal O}_{Y_{\bul \leq N,t}})}
{\cal W}\Om_{Y_{\bul \leq N,t}}^{\bul \geq i}\{-i\}\otimes_{\mab Z}{\mab Q} 
\lo 0.   
\end{align*}
Let 
\begin{align*} 
e^{-1}(N_{{\rm dRW}})' \col  & F^{\bul \leq N}_{\infty}
\otimes_{{\cal W}({\cal O}_{Y_{\bul \leq N,t}})}
({\cal W}{{\Om}}^{\bul \geq i}_{Y_{\bul \leq N,t}})'\{-i\}
\otimes_{\mab Z}{\mab Q}  
\tag{2.4.3.4}\label{ali:nldewdw}\\
&\lo F^{\bul \leq N}_{\infty}
\otimes_{{\cal W}({\cal O}_{Y_{\bul \leq N,\os{\circ}{t}}})'}
({\cal W}\Om^{\bul \geq i-1}_{Y_{\bul \leq N,t}})'
\{-(i-1)\}(-1,u)\otimes_{\mab Z}{\mab Q}  
\end{align*}
be the boundary morphism of (\ref{eqn:wfnby}).  
Let 
\begin{align*} 
e^{-1}N_{{\rm dRW}} \col & F^{\bul \leq N}_{\star}
\otimes_{{\cal W}({\cal O}_{Y_{\bul \leq N,t}})}
{\cal W}\Om_{Y_{\bul \leq N,t}}^{\bul \geq i}\{-i\}
\otimes_{\mab Z}{\mab Q}  
\tag{2.4.3.5}\label{eqn:nlwdewnp} \\
& \lo F^{\bul \leq N}_{\star}
\otimes_{{\cal W}({\cal O}_{Y_{\bul \leq N,\os{\circ}{t}}})}
{\cal W}\Om_{Y_{\bul \leq N,t}}^{\bul \geq i-1}\{-(i-1)\}(-1,u)
\otimes_{\mab Z}{\mab Q}  
\end{align*}
be the boundary morphism of (\ref{eqn:wenlby}). 
Via the identification 
\begin{align*} 
F^{\bul \leq N}_{\infty}
\otimes_{{\cal W}({\cal O}_{Y_{\bul \leq N,t}})'}
({\cal W}{\Om}_{Y_{\bul \leq N,t}}^{\bul \geq j})'
\os{\sim}{\lo} 
F^{\bul \leq N}_{\infty}
\otimes_{{\cal W}({\cal O}_{Y_{\bul \leq N,t}})}
{\cal W}{\Om}_{Y_{\bul \leq N,t}}^{\bul \geq j},  
\quad (j=i,i-1)
\end{align*}
we obtain an equality 
$e^{-1}(N_{\rm dRW})'=e^{-1}N_{\rm dRW}$.

\begin{prop}[{\bf Contravariant functoriality of 
the monodromy operator}]\label{prop:ctrwmon} 
Let the notations be as in {\rm (\ref{prop:cttu})}. 
Assume that $S=s$, $S'=s'$, $T={\cal W}(t)$ and $T={\cal W}(t')$.
Let $g_{\bul \leq N}\col Y_{\bul \leq N,\os{\circ}{t}}\lo Z_{\bul \leq N,\os{\circ}{t}{}'}$ 
be a morphism over the morphism $u\col s_{\os{\circ}{t}}\lo s'_{\os{\circ}{t}{}'}$. 
Let $\star$ be nothing  or $\prime$ and 
let $\sharp=n\in {\mab Z}_{\geq 1}$ or $\infty$. 
Then the following diagrams are commutative$:$
\begin{equation*} 
\begin{CD} 
g_*(F^{\bul \leq N}_{\sharp}
\otimes_{{\cal W}_{\sharp}({\cal O}_{Y_{\bul \leq N,t}})^{\star}}
({\cal W}_{\sharp}\Om_{Y_{\bul \leq N,t}}^{\bul \geq i})^{\star}\{-i\}) 
@>{(N^{\geq i}_{\rm dRW})^{\star}}>>\\
@A{g^*_{\bul \leq N}}AA \\
G^{\bul \leq N}_{\sharp}\otimes_{{\cal W}_{\sharp}({\cal O}_{Z_{\bul \leq N,t}})^{\star}}
({\cal W}_{\sharp}\Om_{Z_{\bul \leq N,t}}^{\bul \geq i})^{\star}\{-i\}
@>{\deg (u)(N^{\geq i}_{\rm dRW})^{\star}}>>
\end{CD}
\tag{2.4.4.1}\label{eqn:yaluwss}
\end{equation*}
\begin{equation*} 
\begin{CD} 
g_*(F^{\bul \leq N}_{\sharp}
\otimes_{{\cal W}_{\sharp}({\cal O}_{Y_{\bul \leq N,t}})^{\star}}
({\cal W}_{\sharp}\Om_{Y_{\bul \leq N,t}}^{\bul \geq i-1})^{\star})\{-i\}[1](-1,u)\\
@AA{g^*_{\bul \leq N}}A\\
G^{\bul \leq N}_{\sharp}
\otimes_{{\cal W}_{\sharp}({\cal O}_{Z_{\bul \leq N,t}})^{\star}}
({\cal W}_{\sharp}\Om_{Z_{\bul \leq N,t}}^{\bul \geq i-1})^{\star}\{-i\}[1](-1,u), 
\end{CD}
\end{equation*}
\begin{equation*} 
\begin{CD} 
g_*(F^{\bul \leq N}_{\sharp}
\otimes_{{\cal W}_{\sharp}({\cal O}_{Y_{\bul \leq N,t}})^{\star}}
({\cal W}_{\sharp}\Om_{Y_{\bul \leq N,t}}^{\bul \geq i})^{\star}\{-i\})
\otimes_{\mab Z}{\mab Q}   
@>{e^{-1}(N^{\geq i}_{\rm dRW})^{\star}}>>\\
@A{g^*_{\bul \leq N}}AA \\
G^{\bul \leq N}_{\sharp}\otimes_{{\cal W}_{\sharp}({\cal O}_{Z_{\bul \leq N,t}})^{\star}}
({\cal W}_{\sharp}\Om_{Z_{\bul \leq N,t}}^{\bul \geq i})^{\star}\{-i\}\otimes_{\mab Z}{\mab Q}  
@>{e^{-1}(N^{\geq i}_{\rm dRW})^{\star}}>>
\end{CD}
\tag{2.4.4.2}\label{eqn:yadws}
\end{equation*}
\begin{equation*} 
\begin{CD} 
g_*(F^{\bul \leq N}_{\sharp}\otimes_{{\cal W}_{\sharp}({\cal O}_{Y_{\bul \leq N,t}})^{\star}}
({\cal W}_{\sharp}\Om_{Y_{\bul \leq N,t}}^{\bul \geq i-1})^{\star})\{-i\}[1](-1,u)
\otimes_{\mab Z}{\mab Q}  \\
@AA{g^*_{\bul \leq N}}A\\
G^{\bul \leq N}_{\sharp}
\otimes_{{\cal W}_{\sharp}({\cal O}_{Z_{\bul \leq N,t}})^{\star}}
({\cal W}_{\sharp}\Om_{Z_{\bul \leq N,t}}^{\bul \geq i-1})^{\star}\{-i\}[1](-1,u)
\otimes_{\mab Z}{\mab Q}. 
\end{CD}
\end{equation*}
\end{prop}
\begin{proof} 
This immediately follows from the formula $g^*(\theta')=\deg(u)\theta$, 
where $\theta'$ is an analogous form to $\theta$ for $s'_{\os{\circ}{t}{}'}$.
\end{proof}

\section{$p$-adic monodromy-weight conjecture 
and filtered log $p$-adic hard Lefschetz conjecture 
via log de Rham-Witt complexes}\label{sec:pdwmn} 
In this section we show the variational $p$-adic monodromy-weight conjecture  
(\ref{eqn:grmpd}) in a special case 
is a generalization of Mokrane $p$-adic monodromy-weight conjecture in \cite{msemi}. 
We also formulate the variational filtered log $p$-adic hard Lefschetz conjecture 
via log de Rham-Witt complexes. This is equivalent to (\ref{eqn:filqpl}) in a special case. 
We also define a filtered complex 
$({\cal W}_nB_{\rm zar},P)$ fitting into the following triangle of filtered complexes
\begin{align*} 
({\cal W}_nA_{\rm zar},P)[-1]\lo ({\cal W}_nB_{\rm zar},P)\lo ({\cal W}_nA_{\rm zar},P)\os{+1}{\lo} 
\end{align*} 
and we prove the contravariant functoriality of 
$({\cal W}_nB_{\rm zar},P)$. 
(Though Mokrane has essentially defined the complex 
${\cal W}_nB_{\rm zar}$ in [loc.~cit.], we have had to correct it 
about the signs of the boundary morphisms of 
${\cal W}_nB_{\rm zar}$ in \cite{msemi}. See \cite[(11.9)]{ndw} for details.)  
As a corollary of the functoriality, 
we prove the contravariant functoriality of 
the $p$-adic quasi-monodromy operator.

\par 
Let the notations be as in \S\ref{sec:flgdw}.

\par 
Let $N$ be a nonnegative integer or $\infty$. 
Let $X_{\bul \leq N}$ be 
an $N$-truncated simplicial SNCL scheme over $s$. 
Let $E^{\bul \leq N}$ be a flat coherent log crystal of 
${\cal O}_{\os{\circ}{X}_{\bul \leq N,t}/{\cal W}_{\star}(\os{\circ}{t})}$-modules.

Let 
\begin{equation*}
\nu_{{\rm dRW},\star} \col 
{\cal W}_{\star}A_{X_{\bul \leq N,\os{\circ}{t}}}(E^{\bul \leq N}) 
\lo {\cal W}_{\star}A_{X_{\bul \leq N,\os{\circ}{t}}}(E^{\bul \leq N})
\tag{2.5.0.1}\label{eqn:wadnu}
\end{equation*}
be the induced morphism by the following morphism 
$${\rm proj}. 
\col {\cal W}_{\star}A_{X_{\bul \leq N,\os{\circ}{t}}}(E^{\bul \leq N})^{ij}
\lo {\cal W}_{\star}A_{X_{\bul \leq N,\os{\circ}{t}}}(E^{\bul \leq N})^{i-1,j+1} 
\quad (i,j\in {\mab Z}).
$$
It is easy to check that $\nu_{{\rm dRW}}$ is actually a morphism of 
complexes. 

\begin{defi}\label{defi:dwdf} 
We call $\nu_{\rm dRW}$ 
the {\it reverse log de Rham-Witt quasi-monodromy operator} 
of $E^{\bul \leq N}$. 
When $E^{\bul \leq N}=
{\cal O}_{\os{\circ}{X}_{\bul \leq N,t}/{\cal W}_{\star}(\os{\circ}{t})}$, 
we call $\nu_{{\rm dRW},\star}$
the {\it reverse log de Rham-Witt quasi-monodromy operator} of 
$X_{\bul \leq N,\os{\circ}{t}}/{\cal W}_{\star}(s_{\os{\circ}{t}})$.
\end{defi}

\par 
Assume that $N<\infty$. 
Let $g_{\bul \leq N}$ be the morphism 
(\ref{eqn:xdxduss}) for the case 
$Y_{\bul \leq N,\os{\circ}{T}{}'_0}= X_{\bul \leq N,\os{\circ}{t}}$,  
$(T',{\cal J}',\del')=(T,{\cal J},\del)=({\cal W}(t),p{\cal W}(\kap_t),[~])$ 
and $S'=S=s$ 
satisfying the condition (\ref{cd:xygxy}): 
\begin{equation*} 
\begin{CD} 
X'_{\bul \leq N,\os{\circ}{t}} @>{g'_{\bul \leq N}}>> X''_{\bul \leq N,\os{\circ}{t}} \\
@VVV @VVV \\ 
X_{\bul \leq N,\os{\circ}{t}} @>{g_{\bul \leq N}}>> X_{\bul \leq N,\os{\circ}{t}},  
\end{CD}
\tag{2.5.1.1}\label{cd:xgwwy}
\end{equation*}
where $X''_{\bul \leq N,\os{\circ}{t}}$ is another disjoint union 
of the member of an affine $N$-truncated simplicial open covering of 
$X_{\bul \leq N,\os{\circ}{t}}$. 
Assume that $\deg(u)$ is not divisible by $p$ or 
that $\star$ is nothing and 
that the morphism (\ref{ali:odnl}) is divisible by $p^{e_p(j+1)}$.  
By a similar diagram 
for ${\cal W}_{\star}A_{X_{\bul \leq N,\os{\circ}{t}}}(E^{\bul \leq N})^{\bul \bul}$ 
to (\ref{cd:phipnu}), the morphism (\ref{eqn:wadnu}) is 
the following morphism 
\begin{equation*} 
\nu_{{\rm dRW},{\star}} \col 
({\cal W}_{\star}A_{X_{\bul \leq N,\os{\circ}{t}}}(E^{\bul \leq N}),P)
\lo 
({\cal W}_{\star}A_{X_{\bul \leq N,\os{\circ}{t}}}(E^{\bul \leq N}),P\langle -2\rangle)
(-1,u). 
\tag{2.5.1.2}\label{eqn:axwdxdn}
\end{equation*}

\par 
The following follows from the proof of 
the comparison theorem (\ref{theo:csoncrdw}) and 
from the definition of $\nu_{\rm zar}$ and 
$\nu_{\rm dRW}$: 

\begin{prop}\label{prop:nnuc} 
The following diagram is commutative$:$ 
\begin{equation*}
\begin{CD}
A_{\rm zar}(X_{\bul \leq N,\os{\circ}{t}}/{\cal W}_{\star}(\os{\circ}{t}),E^{\bul \leq N})
@>{\nu_{\rm zar}}>> 
A_{\rm zar}(X_{\bul \leq N,\os{\circ}{t}}/{\cal W}_{\star}(\os{\circ}{t}),E^{\bul \leq N})(-1,u)\\
@V{\simeq}VV   
@VV{\simeq}V \\
{\cal W}_{\star}A_{X_{\bul \leq N,\os{\circ}{t}}}(E^{\bul \leq N})
@>{\nu_{{\rm dRW},\star}}>> 
{\cal W}_{\star}A_{X_{\bul \leq N,\os{\circ}{t}}}(E^{\bul \leq N})(-1,u).
\end{CD}
\tag{2.5.2.1}\label{cd:nuedwqln}
\end{equation*} 
\end{prop}

\par
Let us give an obvious generalization 
the double complex $W_{\star}B_X^{\bul \bul}$ defined in 
\cite[\S11]{ndw} (cf.~\cite[p.~246]{st1}, \cite[p.~318]{msemi}). 
The $(i,j)$-component ${\cal W}_{\star}B_{X_{\bul \leq N,\os{\circ}{t}}}(E^{\bul \leq N})^{ij}$ 
$(i,j\in {\mab Z}_{\geq 0})$ is defined by the following formula: 
\begin{align*} 
{\cal W}_{\star}B_{X_{\bul \leq N,\os{\circ}{t}}}(E^{\bul \leq N})^{ij}:=
{\cal W}_{\star}A_{X_{\bul \leq N,\os{\circ}{t}}}(E^{\bul \leq N})^{i-1,j}(-1,u)
\oplus {\cal W}_{\star}A_{X_{\bul \leq N,\os{\circ}{t}}}(E^{\bul \leq N})^{ij}.
\tag{2.5.2.2}\label{ali:sgbxn}
\end{align*}  
The horizontal boundary morphism
$d' \col {\cal W}_{\star}B_{X_{\bul \leq N,\os{\circ}{t}}}(E^{\bul \leq N})^{ij} \lo 
{\cal W}_{\star}B_{X_{\bul \leq N,\os{\circ}{t}}}(E^{\bul \leq N})^{i+1,j}$ 
is, by definition,
$$d'(\om_1, \om_2)=(\nabla \om_1,
-\nabla\om_2) $$
and the vertical one 
$d'' \col {\cal W}_{\star}B_{X_{\bul \leq N,\os{\circ}{t}}}(E^{\bul \leq N})^{ij} \lo 
{\cal W}_{\star}B_{X_{\bul \leq N,\os{\circ}{t}}}(E^{\bul \leq N})^{i,j+1}$ is 
$$d''(\om_1, \om_2)=
(-\theta_{\star}\wedge \om_1+\nu_{\rm dRW}(\om_2),\theta_{\star}\wedge \om_2).$$ 
Here we have omitted a notation 
${\rm mod}~P_{j+1}$ in the 
definition of $d''$ for simplicity.
It is easy to check that 
${\cal W}_{\star}B_{X_{\bul \leq N,\os{\circ}{t}}}(E^{\bul \leq N})^{\bul \bul}$ 
is indeed a double complex. Let 
${\cal W}_{\star}B_{X_{\bul \leq N,\os{\circ}{t}}}(E^{\bul \leq N})$ be the single complex of 
${\cal W}_{\star}B_{X_{\bul \leq N,\os{\circ}{t}}}(E^{\bul \leq N})^{\bul \bul}$. 
\par
Let 
\begin{align*} 
\mu_{\star} \col E^{\bul \leq N}_{\star}
\otimes_{{\cal W}_{\star}({\cal O}_{X_{\bul \leq N,\os{\circ}{t}}})}
{\cal W}_{\star}\wt{{\Om}}^{\bul}_{X_{\bul \leq N,\os{\circ}{t}}} 
\lo {\cal W}_{\star}B_{X_{\bul \leq N,\os{\circ}{t}}}(E^{\bul \leq N})
\end{align*} 
be a morphism of complexes defined by 
$$\mu_{\star}(\om):=(\om~{\rm mod}~P_0,
\theta_{\star} \wedge \om~{\rm mod}~P_0) \quad 
( \om \in  E^{\bul \leq N}_{\star}
\otimes_{{\cal W}_{\star}({\cal O}_{X_{\bul \leq N,\os{\circ}{t}}})}
{\cal W}_{\star}\wt{{\Om}}^{\bul}_{X_{\bul \leq N,\os{\circ}{t}}}).$$
Then we have the following commutative diagram of exact sequences:
\begin{equation*}
\begin{CD}
0 @>>> {\cal W}_{\star}A_{X_{\bul \leq N,\os{\circ}{t}}}(E^{\bul \leq N})(-1,u)[-1]   
 \\
@. @A{(\theta_{\star}\wedge *)[-1] }AA   \\
0 @>>> (\eps^*_{X_{\bul \leq N,\os{\circ}{t}}/{\cal W}_{\star}(\os{\circ}{t})}
(E^{\bul \leq N}))_{{\cal W}_{\star}(X_{\bul \leq N,\os{\circ}{t}})}
\otimes_{{\cal W}_{\star}({\cal O}_{X_{\bul \leq N,\os{\circ}{t}}})}
{\cal W}_{\star}{\Om}^{\bul}_{X_{\bul \leq N,\os{\circ}{t}}}(-1,u)[-1]   \\
\end{CD}
\tag{2.5.2.3}\label{cd:sgstexn}
\end{equation*}
\begin{equation*}
\begin{CD}
@>{}>>{\cal W}_{\star}B_{X_{\bul \leq N,\os{\circ}{t}}}(E^{\bul \leq N}) @>{}>>  \\
@. @A{\mu_{\star}}AA\\
@>{\theta_{\star} \wedge}>>E^{\bul \leq N}_{\star}
\otimes_{{\cal W}_{\star}({\cal O}_{X_{\bul \leq N,\os{\circ}{t}}})}
{\cal W}_{\star}\wt{{\Om}}^{\bul}_{X_{\bul \leq N,\os{\circ}{t}}}@>{}>>  
\end{CD} 
\end{equation*}
\begin{equation*}
\begin{CD}
{\cal W}_{\star}A_{X_{\bul \leq N,\os{\circ}{t}}}(E^{\bul \leq N})
@>>> 0 \\
@A{\theta_{\star} \wedge}AA\\
(\eps^*_{X_{\bul \leq N,\os{\circ}{t}}/{\cal W}_{\star}(\os{\circ}{t})}
(E^{\bul \leq N}))_{{\cal W}_{\star}(X_{\bul \leq N,\os{\circ}{t}})}
\otimes_{{\cal W}_{\star}({\cal O}_{X_{\bul \leq N,\os{\circ}{t}}})}
{\cal W}_{\star}{\Om}^{\bul}_{X_{\bul \leq N,\os{\circ}{t}}}@>>> 0.  
\end{CD} 
\end{equation*}
(By (\ref{eqn:wby}), the lower sequence is exact.) 

As in \cite[(11.10)]{ndw} we can easily prove the following by (\ref{cd:sgstexn}):

\begin{prop}\label{prop:qneqn}
Let $\star$ be a positive integer $n$ or nothing. 
Then the following diagram is commutative$:$
\begin{equation*}
\begin{CD}
{\cal W}_{\star}A_{X_{\bul \leq N,\os{\circ}{t}}}(E^{\bul \leq N})
@>{\nu_{{\rm dRW},\star}}>>  \\
@A{\theta_{\star}\wedge,~ \simeq}AA  \\ 
(\eps^*_{X_{\bul \leq N,t}/{\cal W}_{\star}(t)}
(E^{\bul \leq N}))_{{\cal W}_{\star}({\cal O}_{X_{\bul \leq N,t}})}
\otimes_{{\cal W}_{\star}({\cal O}_{X_{\bul \leq N,t}})}
{\cal W}_{\star}{{\Om}}^{\bul}_{X_{\bul \leq N,t}}
@>{N_{{\rm dRW},\star}}>> 
\end{CD}
\tag{2.5.3.1}\label{cd:nueqln}
\end{equation*} 
\begin{equation*}
\begin{CD}
{\cal W}_{\star}A_{X_{\bul \leq N,\os{\circ}{t}}}(E^{\bul \leq N})(-1,u) \\
@AA{\theta_{\star} \wedge, ~\simeq}A \\
(\eps^*_{X_{\bul \leq N,t}/{\cal W}_{\star}(t)}
(E^{\bul \leq N}))_{{\cal W}_{\star}(X_{\bul \leq N,t})}
\otimes_{{\cal W}_{\star}({\cal O}_{X_{\bul \leq N,t}})}
{\cal W}_{\star}{{\Om}}^{\bul}_{X_{\bul \leq N,t}}(-1,u).
\end{CD} 
\end{equation*} 
\end{prop}

\par


In the following we only state analogous results  for 
$({\cal W}_{\star}A_{X_{\bul \leq N,\os{\circ}{t}}}(E^{\bul \leq N}),P)$ 
to results $(A_{\rm zar}(X_{\bul \leq N,\os{\circ}{T}_0}/S(T)^{\nat},E^{\bul \leq N})),P)$ 
in \S\ref{sec:vpmn}. We omit the proofs of them because we have only to 
make suitable modifications 
of proofs in \S\ref{sec:vpmn}.

\begin{prop}\label{prop:quaimon}
Let $k$ be a positive integer and let $q$ be a nonnegative integer. 
Then 
\begin{align*} 
\nu_{{\rm dRW},\star}^k \col  & 
R^qf_{*}
({\rm gr}_k^P{\cal W}_{\star}A_{X_{\bul \leq N,\os{\circ}{t}}}(E^{\bul \leq N}))
\lo 
R^qf_{*}
({\rm gr}_{-k}^P{\cal W}_{\star}A_{X_{\bul \leq N,\os{\circ}{t}}}(E^{\bul \leq N}))(-k,u)
\tag{2.5.4.1}\label{eqn:kwpa}
\end{align*} 
is the identity. 
\end{prop}

\begin{prop}\label{prop:nucswa} 
The quasi-monodromy operator 
\begin{equation*} 
\nu_{{\rm dRW},\star} \col 
{\cal W}_{\star}A_{X_{\bul \leq N,\os{\circ}{t}}}(E^{\bul \leq N})  \lo 
{\cal W}_{\star}A_{X_{\bul \leq N,\os{\circ}{t}}}(E^{\bul \leq N})(-1,u)
\end{equation*}
underlies the following morphism 
of filtered morphism 
\begin{equation*} 
\nu_{{\rm dRW},\star} \col 
({\cal W}_{\star}A_{X_{\bul \leq N,\os{\circ}{t}}}(E^{\bul \leq N}),P)
\lo  
({\cal W}_{\star}A_{X_{\bul \leq N,\os{\circ}{t}}}(E^{\bul \leq N}),P\langle -2 \rangle)(-1,u), 
\end{equation*}
where $P\langle -2 \rangle$ is a filtration 
defined by $(P\langle -2 \rangle)_k= P_{k-2}$. 
\end{prop}

\begin{coro}\label{coro:zmwop}
The quasi-monodromy operator {\rm (\ref{defi:dwdf})} 
induces the following morphism 
\begin{align*} 
\nu_{{\rm dRW}} \col & P_kRf_{X_{\bul \leq N,\os{\circ}{t}}/{\cal W}(s_{\os{\circ}{t}})*}
(\eps^*_{X_{\bul \leq N,\os{\circ}{t}}/{\cal W}(s_{\os{\circ}{t}})}(E^{\bul \leq N}))
\lo  \tag{2.5.6.1}\label{eqn:grpwpd}\\
& P_{k-2}Rf_{X_{\bul \leq N,\os{\circ}{t}}/{\cal W}(s_{\os{\circ}{t}})*}
(\eps^*_{X_{\bul \leq N,\os{\circ}{t}}/{\cal W}(s_{\os{\circ}{t}})}(E^{\bul \leq N}))(-1,u). 
\end{align*}
\end{coro}

\begin{prop}\label{prop:spqwqqqq}
Let the notations be as in {\rm (\ref{theo:itc})}. 
Assume that $S=s$, $S=s'$, $T={\cal W}_{\star}(t)$ and $T'={\cal W}_{\star}(t')$. 
Then the following hold$:$
\par 
$(1)$ The following diagram is commutative$:$ 
\begin{equation*} 
\begin{CD} 
Rg_{{\bul \leq N}*}(({\cal W}_{\star}A_{X_{\bul \leq N,\os{\circ}{t}}}(E^{\bul \leq N}),P))
@>{\nu_{{\rm dRW},\star}}>> 
Rg_{{\bul \leq N}*}(({\cal W}_{\star}A_{X_{\bul \leq N,\os{\circ}{t}}}(E^{\bul \leq N}),
P\langle -2 \rangle))\\ 
@A{g^*_{\bul \leq N}}AA @AA{g^*_{\bul \leq N}}A \\
({\cal W}_{\star}A_{Y_{\bul \leq N,\os{\circ}{t}{}'}}(F^{\bul \leq N}),P)
@>{\deg (u)\nu_{{\rm dRW},\star}}>>
({\cal W}_{\star}A_{Y_{\bul \leq N,\os{\circ}{t}{}'}}(F^{\bul \leq N}),P\langle -2 \rangle). 
\end{CD} 
\tag{2.5.7.1}\label{eqn:aawsf}
\end{equation*} 
\par 
$(2)$ Let the notations be as in {\rm (\ref{theo:itc})}.  
Then the following diagram is commutative$:$ 
\begin{equation*} 
\begin{CD} 
Rg_{{\bul \leq N}*}
(({\cal W}_{\star}A_{X_{\bul \leq N,\os{\circ}{t}}}(E^{\bul \leq N}),P))
\otimes^L_{\mab Z}{\mab Q}
@>{\nu_{{\rm dRW},\star}}>>   \\ 
@A{g^*_{\bul \leq N}}AA \\
({\cal W}_{\star}A_{Y_{\bul \leq N,\os{\circ}{t}{}'}}(F^{\bul \leq N}),P)
\otimes^L_{\mab Z}{\mab Q}
@>{\deg (u) \nu_{{\rm dRW},\star}}>>   
\end{CD} 
\tag{2.5.7.2}\label{eqn:aalwsf}
\end{equation*} 
\begin{equation*} 
\begin{CD} 
Rg_{{\bul \leq N}*}
(({\cal W}_{\star}A_{X_{\bul \leq N,\os{\circ}{t}}}(E^{\bul \leq N}),
P\langle -2 \rangle)(-1,u))
\otimes^L_{\mab Z}{\mab Q}\\
@AA{g^*_{\bul \leq N}}A \\
({\cal W}_{\star}A_{Y_{\bul \leq N,\os{\circ}{t}{}'}}(F^{\bul \leq N}),
P\langle -2 \rangle)(-1,u')\otimes^L_{\mab Z}{\mab Q}. 
\end{CD} 
\end{equation*} 
\end{prop}

Set 
\begin{align*} 
P_k{\cal W}_{\star}B_{X_{\bul \leq N,\os{\circ}{t}}}(E^{\bul \leq N})^{ij}
:=&P_k{\cal W}_{\star}A_{X_{\bul \leq N,\os{\circ}{t}}}(E^{\bul \leq N})^{i-1,j}(-1,u)
\oplus P_k{\cal W}_{\star}A_{X_{\bul \leq N,\os{\circ}{t}}}(E^{\bul \leq N})^{ij} 
\tag{2.5.7.3}\label{eqn:aabbsf}\\
&(i,j\in {\mab N}).
\end{align*} 

\begin{prop}\label{prop:bwid}
$(1)$ There exists the following sequence of the triangles in 
${\rm C}^+{\rm F}(f^{-1}_{\bul \leq N}({\cal O}_T)):$ 
\begin{align*} 
&\lo ({\cal W}_{\star}A_{X_{\bul \leq N,\os{\circ}{t}}}(E^{\bul \leq N}),P)[-1]
\lo ({\cal W}_{\star}B_{X_{\bul \leq N,\os{\circ}{t}}}(E^{\bul \leq N}),P)
\tag{2.5.8.1}\label{eqn:aabebsf}\\
&\lo ({\cal W}_{\star}A_{X_{\bul \leq N,\os{\circ}{t}}}(E^{\bul \leq N}),P)\os{+1}{\lo}. 
\end{align*} 
\end{prop} 

\begin{defi} 
We call $({\cal W}_{\star}B_{X_{\bul \leq N,\os{\circ}{t}}}(E^{\bul \leq N}),P)$ 
the {\it extended zariskian $p$-adic filtered Steenbrink complexes} of 
$E^{\bul \leq N}$ for $X_{\bul \leq N,\os{\circ}{t}}/{\cal W}(s_{\os{\circ}{t}})$. 
When $E^{\bul \leq N}={\cal O}_{\os{\circ}{X}_{\bul \leq N,t}/\os{\circ}{t}}$, 
we denote it by 
$({\cal W}_{\star}B_{X_{\bul \leq N,\os{\circ}{t}}},P)$ and call this 
the {\it extended zariskian $p$-adic filtered Steenbrink complex} of 
$X_{\bul \leq N,\os{\circ}{t}}/{\cal W}(s_{\os{\circ}{t}})$. 
\end{defi} 

\begin{theo}[{\bf Contravariant functoriality I of ${\cal W}_{\star}B_{\rm zar}$}]
\label{theo:funwcb}
Let the notations and the assumptions be as in {\rm (\ref{theo:ctrw})}. 
Then the following hold$:$ 
\par 
$(1)$ 
The morphism 
$g_{\bul \leq N}\col X_{\bul \leq N,\os{\circ}{t}}\lo Y_{\bul \leq N,\os{\circ}{t}{}'}$ 
induces the following pull-back morphism 
\begin{equation*}  
g_{\bul \leq N}^* \col 
({\cal W}_{\star}B_{Y_{\bul \leq N,\os{\circ}{t}{}'}}(F^{\bul \leq N}),P)
\lo Rg_{\bul \leq N*}(({\cal W}_{\star}B_{X_{\bul \leq N,\os{\circ}{t}}}(E^{\bul \leq N}),P)) 
\tag{2.5.10.1}\label{eqn:fzwaxd}
\end{equation*} 
fitting into the following commutative diagram$:$
\begin{equation*} 
\begin{CD}
{\cal W}_{\star}B_{Y_{\bul \leq N,\os{\circ}{t}{}'}}(F^{\bul \leq N})
@>{g_{\bul \leq N}^*}>>  
Rg_{\bul \leq N*}({\cal W}_{\star}B_{X_{\bul \leq N,\os{\circ}{t}}}(E^{\bul \leq N}))\\ 
@A{\mu_{\star}\wedge}A{\simeq}A 
@A{Rg_{\bul \leq N*}(\mu_{\star})}A{\simeq}A\\
F^{\bul \leq N}_{\star}
\otimes_{{\cal W}_{\star}({\cal O}_{Y_{\bul \leq N,\os{\circ}{t}{}'}})}
{\cal W}_{\star}\wt{{\Om}}^{\bul}_{Y_{\bul \leq N,\os{\circ}{t}{}'}}
@>{g_{\bul \leq N}^*}>>Rg_{\bul \leq N*}(E^{\bul \leq N}_{\star}
\otimes_{{\cal W}_{\star}({\cal O}_{X_{\bul \leq N,\os{\circ}{t}}})}
{\cal W}_{\star}\wt{{\Om}}^{\bul}_{X_{\bul \leq N,\os{\circ}{t}}}).
\end{CD}
\tag{2.5.10.2}\label{cd:psbwcz} 
\end{equation*}
\par 
$(2)$ Let $h_{\bul \leq N}\col Y_{\bul \leq N,\os{\circ}{t}{}'}\lo 
Z_{\bul \leq N,\os{\circ}{t}{}''}$ be an analogous morphism to 
$g_{\bul \leq N}\col X_{\bul \leq N,\os{\circ}{t}}\lo 
Y_{\bul \leq N,\os{\circ}{t}{}'}$. 
\begin{align*} 
(h_{\bul \leq N}\circ g_{\bul \leq N})^*  =&
Rh_{\bul \leq N*}(g_{\bul \leq N}^*)\circ h_{\bul \leq N}^*    
\col ({\cal W}_{\star}B_{Z_{\bul \leq N,\os{\circ}{t}{}''}}(G^{\bul \leq N}),P)
\tag{2.5.10.3}\label{ali:pwwp} \\ 
& \lo Rh_{\bul \leq N*}Rg_{\bul \leq N*}
({\cal W}_{\star}B_{X_{\bul \leq N,\os{\circ}{t}}}(E^{\bul \leq N}),P) \\
& =R(h_{\bul \leq N}\circ g_{\bul \leq N})_*
({\cal W}_{\star}B_{X_{\bul \leq N,\os{\circ}{t}}}(E^{\bul \leq N}),P).
\end{align*}  
\par 
$(3)$  
\begin{equation*} 
{\rm id}_{X_{\bul \leq N,\os{\circ}{T}_0}}^*={\rm id} 
\col ({\cal W}_{\star}B_{X_{\bul \leq N,\os{\circ}{t}}}(E^{\bul \leq N}),P)
\lo ({\cal W}_{\star}B_{X_{\bul \leq N,\os{\circ}{t}}}(E^{\bul \leq N}),P).  
\tag{2.5.10.4}\label{eqn:fzwdd}
\end{equation*} 
\end{theo} 

\begin{theo}[{\bf Contravariant functoriality II of ${\cal W}B_{\rm zar}$}]\label{theo:fwuniicb}
Let the notations be as in {\rm (\ref{theo:itc})}.  
Assume that $S=s$, $S=s'$, $T={\cal W}(t)$ and $T'={\cal W}(t')$. 
Then there exists a morphism 
\begin{align*} 
g^*_{\bul \leq N}  \col 
({\cal W}B_{Y_{\bul \leq N,\os{\circ}{t}{}'}}(F^{\bul \leq N}),P)\otimes^L_{\mab Z}{\mab Q}
\lo
Rg_{{\bul \leq N}*}
(({\cal W}B_{X_{\bul \leq N,\os{\circ}{t}}}(E^{\bul \leq N}),P)\otimes^L_{\mab Z}{\mab Q})
\tag{2.5.11.1}\label{eqn:awqd}
\end{align*}
of filtered complexes in 
${\rm D}^+{\rm F}(f^{-1}_{\bul \leq N}({\cal W}(\kap_{t'})\otimes_{\mab Z}{\mab Q})$ fitting into 
the following commutative diagram 
\begin{equation*} 
\begin{CD}
{\cal W}B_{Y_{\bul \leq N,\os{\circ}{t}{}'}}(F^{\bul \leq N})\otimes^L_{\mab Z}{\mab Q}
@>{g_{\bul \leq N}^*}>>  
Rg_{\bul \leq N*}({\cal W}B_{X_{\bul \leq N,\os{\circ}{t}}}(E^{\bul \leq N})
\otimes^L_{\mab Z}{\mab Q})\\ 
@A{\mu \wedge}A{\simeq}A 
@A{Rg_{\bul \leq N*}(\mu)}A{\simeq}A\\
F^{\bul \leq N}_{\infty}
\otimes_{{\cal W}({\cal O}_{Y_{\bul \leq N,\os{\circ}{t}{}'}})}
{\cal W}\wt{{\Om}}^{\bul}_{Y_{\bul \leq N,\os{\circ}{t}{}'}}
\otimes^L_{\mab Z}{\mab Q}
@>{g_{\bul \leq N}^*}>>Rg_{\bul \leq N*}(E^{\bul \leq N}_{\infty}
\otimes_{{\cal W}({\cal O}_{X_{\bul \leq N,\os{\circ}{t}}})}
{\cal W}\wt{{\Om}}^{\bul}_{X_{\bul \leq N,\os{\circ}{t}}})
\otimes^L_{\mab Z}{\mab Q}.
\end{CD}
\tag{2.5.11.2}\label{cd:pswbz} 
\end{equation*}
This morphism satisfies the similar relation to {\rm (\ref{ali:pwwp})}.   
\end{theo}

\begin{coro}\label{coro:swpfs}
$(1)$ Let the notations be as in {\rm (\ref{theo:funwcb}) (1)}. 
The isomorphisms {\rm (\ref{eqn:wlxa})} 
for $X_{\bul \leq N,\os{\circ}{t}}/{\cal W}(s_{\os{\circ}{t}{}})$, $E^{\bul \leq N}$ 
and $Y_{\bul \leq N,\os{\circ}{t}{}'}/{\cal W}(s'_{\os{\circ}{t}{}'})$, $F^{\bul \leq N}$ 
induce a morphism $g_{\bul \leq N}^*$ from {\rm (\ref{eqn:yaluwss})} 
for the case $i=0$, $S=s$, $S'=s'$, $T={\cal W}(t)$ and $T'={\cal W}(t')$ 
to {\rm (\ref{eqn:aawsf})} forgetting the filtrations $P$'s.
This morphism satisfies the transitive relation and  
${\rm id}_{X_{\bul \leq N,\os{\circ}{t}}/{\cal W}(s_{\os{\circ}{t}})}^*={\rm id}$. 
\par 
$(2)$ Let the notations be as in {\rm (\ref{theo:itc})}. 
The isomorphisms {\rm (\ref{eqn:wlxa})} 
for $X_{\bul \leq N,\os{\circ}{T}_0}/S(T)^{\nat}$, $E^{\bul \leq N}$ 
and $Y_{\bul \leq N,\os{\circ}{T}{}'_0}/S'(T')^{\nat}$, $F^{\bul \leq N}$ 
induce a morphism from {\rm (\ref{eqn:yadws})} 
to {\rm (\ref{eqn:aalwsf})}. 
This morphism satisfies the transitive relation and  
${\rm id}_{X_{\bul \leq N,\os{\circ}{T}_0}/S(T)^{\nat}}^*={\rm id}$. 
\end{coro} 

\begin{theo}[{\bf Comparison theorem}]\label{theo:cmpb}
Let the notations be as in {\rm (\ref{theo:csoncrdw})}. 
Then the following hold$:$ 
\par 
$(1)$ In ${\rm D}^+{\rm F}(f^{-1}({\cal W}_n(\kap_t)))$
there exists the following canonical isomorphism$:$
\begin{equation*}
(B_{\rm zar}(X_{\bul \leq N,\os{\circ}{t}}/{\cal W}_n(s_{\os{\circ}{t}}),E^{\bul \leq N}),P)
\os{\sim}{\lo}
({\cal W}_nB_{X_{\bul \leq N,\os{\circ}{t}}}(E^{\bul \leq N}),P).
\tag{2.5.13.1}\label{eqn:brbdw}
\end{equation*}
The isomorphisms $(\ref{eqn:brbdw})$ for $n$'s are 
compatible with two projections of 
both hand sides on $(\ref{eqn:brbdw})$. The isomorphism $(\ref{eqn:brbdw})$ 
is contravariantly functorial. 
\par 
$(2)$ The isomorphism 
{\rm (\ref{eqn:brbdw})} forgetting the filtrations   
fits into the following commutative diagram$:$  
\begin{equation*} 
\begin{CD} 
B_{\rm zar}(X_{\bul \leq N,\os{\circ}{t}}/{\cal W}_{\star}(s_{\os{\circ}{t}}),E^{\bul \leq N})
@>{(\ref{eqn:brbdw}),~\sim}>> \\
@A{\mu_{\star} \wedge}A{\simeq}A \\
\wt{R}u_{X_{\bul \leq N,t}/{\cal W}_{\star}(t)*}
(\eps^*_{X_{\bul \leq N,\os{\circ}{t}}/{\cal W}_{\star}(\os{\circ}{t})}(E^{\bul \leq N})) 
@>{(\ref{eqn:ywnttnny}), \sim}>> 
\end{CD} 
\tag{2.5.13.2}\label{cd:axbwbl}
\end{equation*} 
\begin{equation*} 
\begin{CD} 
{\cal W}_{\star}B_{X_{\bul \leq N,\os{\circ}{t}}}(E^{\bul \leq N}) \\ 
@A{\mu_{\star}}A{\simeq}A \\
(\eps^*_{X_{\bul \leq N,\os{\circ}{t}/{\cal W}_{\star}(\os{\circ}{t})}}
(E^{\bul \leq N}))_{{\cal W}_{\star}(X_{\bul \leq N,\os{\circ}{t}})}
\otimes_{{\cal W}_{\star}({\cal O}_{X_{\bul \leq N,\os{\circ}{t}})}}
{\cal W}_{\star}\wt{\Om}{}^{\bul}_{X_{\bul \leq N,\os{\circ}{t}}}. 
\end{CD} 
\end{equation*} 
\end{theo} 
\begin{proof} 
(1): By the definitions of 
$(B_{\rm zar}(X_{\bul \leq N,\os{\circ}{t}}/{\cal W}_{\star}(s_{\os{\circ}{t}}),E^{\bul \leq N}),P)$ 
((\ref{eqn:anxdn}), (\ref{ali:bpa})) 
and 
$({\cal W}_nB_{X_{\bul \leq N,\os{\circ}{t}}}(E^{\bul \leq N}),P)$
((\ref{ali:sgbxn}), (\ref{eqn:aabbsf})) and by the proof of (\ref{theo:csssh}), 
we have the canonical morphism (\ref{eqn:brbdw}). 
By (\ref{ali:trab}), (\ref{eqn:aabebsf}) and 
(\ref{eqn:brildw}), this morphism is an isomorphism.  
The desired contravariant functoriality follows from (\ref{theo:csoncrdw}) (2). 
\par 
(2): The problem is local. Then (2) follows from 
the constructions of (\ref{eqn:ywnttnny}) and (\ref{eqn:brbdw}). 
\end{proof} 

\par 
Assume that $\os{\circ}{S}$ is a $p$-adic formal scheme. 
Now we consider the case $N=0$.  Set $X:=X_0$. 
In \cite[(3.24)]{msemi} 
Mokrane has already conjectured the $p$-adic monodromy-weight conjecture 
for a projective strict semistable family over a complete discrete valuation ring of 
mixed characteristics. As in (\ref{conj:rcpmc}) we formulate the conjecture for 
a projective SNCL scheme over the log point of a perfect field of characteristic $p>0$ 
as follows (cf.~\cite[(3.28) (4)]{msemi}): 

\begin{conj}[{\bf $p$-adic monodromy-weight conjecture 
via log de Rham-Witt complexes}]\label{conj:rcwmc}
Assume that $\os{\circ}{X}_{t_0} \lo \os{\circ}{t}_0$ is projective. 
Let $q$ be nonnegative integer. 
Then the induced morphism 
\begin{align*} 
\nu^k_{{\rm dRW}} \col & 
{\rm gr}^P_{q+k}H^q_{\rm crys}(X_{t_0}/{\cal W}(t))
\otimes_{\mab Z}{\mab Q} \lo 
{\rm gr}^P_{q-k}H^q_{\rm crys}(X_{t_0}/{\cal W}(t))(-k,u)
\otimes_{\mab Z}{\mab Q} 
\tag{2.5.14.1}\label{eqn:grwpd} 
\end{align*}
by the de Rham-Witt quasi-monodromy operator 
is an isomorphism.  
\end{conj}

\begin{prop}\label{prop:dwqm}
{\rm (\ref{conj:rcwmc})} is equivalent to {\rm (\ref{conj:rcpmc})} 
in the case $S=s$ and $T={\cal W}(t)$.  
\end{prop}
\begin{proof} 
This follows from (\ref{prop:tefc}), (\ref{eqn:wlxa}) and (\ref{cd:nuedwqln}).  
\end{proof}

\par 
Assume that the structural morphism 
$\os{\circ}{f} \col \os{\circ}{X}_{t} \lo \os{\circ}{t}$ 
is projective and the relative dimension of 
$\os{\circ}{f}$ is of pure dimension $d$. 
Let $L$ be a relatively ample line bundle on 
$\os{\circ}{X}_{t}/\os{\circ}{t}$. 
Because the definition of ${\cal W}_n\Om^{\bul}_{Y}$ in this book is 
different from that of ${\cal W}_n\Om^{\bul}_{Y}$ in \cite[I (3.27), (3.29.2)]{idw} 
even in the case of trivial log structures, we need the following even in this case: 

\begin{lemm}\label{lemm:wny}
Let ${\cal Y}/{\cal W}(t)$ be a log smooth integral lift of $Y/t$ over ${\cal W}(t)$. 
Let $M_{\cal Y}$ and $M_Y$ be the log structures of ${\cal Y}$ and $Y$, respectively. 
Then the morphism 
\begin{align*} 
M_{\cal Y}\owns m\lom [d\log m]\in 
{\cal H}^1(\Om^{\bul}_{{\cal Y}/{\cal W}(t)})
\end{align*} 
factors through the projection $M_{\cal Y}\lo M_Y$. 
Here $[~]$ means the cohomology class. 
\end{lemm}  
\begin{proof}
(cf.~the proof of \cite[(6.27)]{ndw})
Let $a$ be a local section of ${\cal O}_{\cal Y}$. 
Consider the local section $m(1+pa)$ of $M_{\cal Y}$. 
Then 
\begin{align*} 
[d\log m(1+pa)]&=[d\log m+\dfrac{pda}{1+pa}]=[d\log m+
(p\sum_{n=1}^{\infty}(-1)^{n-1}(pa)^{n-1}da)]\tag{2.5.16.1}\label{ali:aoy}\\
&=
[d\log m+pd(\sum_{n=1}^{\infty}(-1)^{n-1}\dfrac{p^{n-1}}{n}a^n)=[d\log m].
\end{align*}  
This proves (\ref{lemm:wny}). 
\end{proof} 
By (\ref{lemm:wny}) we can define the following morphisms 
(cf.~\cite[I (3.27), (3.29.2)]{idw}): 
\begin{align*} 
c_{1,{\rm dRW}}
\col M_Y\lo {\cal W}_n\Om^1_{Y} \quad (n\in {\mab Z}_{\geq 1}~{\rm or}~
{\rm nothing}), 
\end{align*} 
\begin{align*} 
c_{1,{\rm dRW}}\col M_Y\lo {\cal W}_n\Om^{\bul}_{Y} 
\quad (n\in {\mab Z}_{\geq 1}~{\rm or}~
{\rm nothing}). 
\end{align*}

\begin{prop}\label{prop:ccc}
Let $Y/t$ be a log smooth integral scheme. 
Let $n$ be a positive integer or nothing. 
Then the following diagram 
\begin{equation*} 
\begin{CD} 
M_Y[-1] @>{c_{1,{\rm dRW}}}>> {\cal W}_n\Om^{\bul}_{Y}\\
@| @AA{\simeq}A \\
M_Y[-1]
@>{c_{1,{\rm crys}}}>> Ru_{Y/{\cal W}_n(t)*}({\cal O}_{Y/{\cal W}_n(t)})
\end{CD}
\end{equation*} 
is commutative. 
If $n$ is a positive integer, then 
this commutative diagram is compatible with $n$'s and 
the projections ${\cal W}_{n+1}\Om^{\bul}_{Y}\lo {\cal W}_n\Om^{\bul}_{Y}$ 
and $Ru_{Y/{\cal W}_{n+1}(t)*}({\cal O}_{Y/{\cal W}_{n+1}(t)})\lo 
Ru_{Y/{\cal W}_n(t)*}({\cal O}_{Y/{\cal W}_n(t)})$. 
\end{prop} 
\begin{proof} 
This is a local problem. 
Hence we may assume that there exists a lift 
${\cal Y}/{\cal W}(t)$ of $Y/t$ over ${\cal W}(t)$. 
Set ${\cal Y}_n:={\cal Y}\otimes_{\cal W}{\cal W}_n$. 
Hence we have only to prove that 
the following diagram 
\begin{equation*} 
\begin{CD} 
M_Y[-1] @>{c_{1,{\rm dRW}}}>> 
{\cal W}_n\Om^{\bul}_{Y}\\
@| @AA{\simeq}A \\
M_Y[-1]@>{c_{1,{\rm dR}}}>> \Om^{\bul}_{{\cal Y}_n/{\cal W}_n(t)}
\end{CD}
\end{equation*} 
is commutative (cf.~\cite[(2.3)]{beill}). 
Indeed, it suffices to prove that 
the following diagram 
\begin{equation*} 
\begin{CD} 
M_Y @>{c_{1,{\rm dRW}}}>> 
{\cal H}^1(\Om^{\bul}_{{\cal Y}_n/{\cal W}_n(t)})\\
@| @AAA \\
M_Y@>{c_{1,{\rm dR}}}>> \Om^1_{{\cal Y}_n/{\cal W}_n(t)}
\end{CD}
\end{equation*} 
is commutative.  
Let $m$ be a local section of $M_{\cal Y}$. 
Let $\Phi_{\cal W}\col {\cal W}(t)\lo {\cal W}(t)$ be 
the canonical lift of the Frobenius endomorphism of $t$. 
We may assume that there exists a lift 
$\Phi \col {\cal Y}\lo {\cal Y}$ of the Frobenius endomorphism 
of $Y$ over $\Phi_{\cal W}$.  
Express $\Phi^*(m)=m^p(1+pa)$ for some $a\in {\cal O}_{\cal Y}$. 
Then the composite of the following morphism 
\begin{align*} 
M_{\cal Y}\owns m\lom p^{-1}d\log \Phi^*(m)&=[d\log m+\dfrac{da}{1+pa}]=
[d\log m+d(\sum_{n=1}^{\infty}(-1)^{n-1}\dfrac{p^{n-1}}{n}a^n)]
\tag{2.5.17.1}\label{ali:aopy}\\
&=[d\log m]
\in {\cal H}^1(\Om^{\bul}_{{\cal Y}/{\cal W}(t)})
\end{align*}  
(cf.~(\ref{ali:aoy}))
and the projection ${\cal H}^1(\Om^{\bul}_{{\cal Y}/{\cal W}(t)}) \lo 
{\cal H}^1(\Om^{\bul}_{{\cal Y}_n/{\cal W}_n(t)})$ 
induces a composite morphism 
$M_Y\lo \Om^1_{{\cal Y}_n/{\cal W}_n(t)}\lo 
{\cal H}^1(\Om^{\bul}_{{\cal Y}_n/{\cal W}_n(t)})$ by (\ref{lemm:wny}). 
This composite morphism is equal to the composite morphism 
$M_Y\os{c_{1,{\rm dR}}}{\lo} \Om^1_{{\cal Y}_n/{\cal W}_n(t)}
\lo {\cal H}^1(\Om^{\bul}_{{\cal Y}_n/{\cal W}_n(t)})$. 
This morphism is nothing but the morphism 
$c_{1,{\rm dRW}}\col M_Y \lo  
{\cal H}^1(\Om^{\bul}_{{\cal Y}_n/{\cal W}_n(t)})$. 
\end{proof} 

\begin{coro}\label{coro:fdc}
The following diagram is commutative: 
\begin{equation*} 
\begin{CD} 
{\rm Pic}(\os{\circ}{X}_t)@>c_{1,{\rm dRW}}>> H^2(X_t,{\cal W}\Om^{\bul}_{X_t})(1)\\
@| @AA{\simeq}A\\
{\rm Pic}(\os{\circ}{X}_t)@>{c_{1,{\rm crys}}}>> H^2_{\rm crys}(X_t/{\cal W}(t))(1). 
\end{CD} 
\tag{2.5.18.1}\label{cd:drcrpic}
\end{equation*} 
\end{coro}
\begin{proof} 
This follows from (\ref{prop:ccc}). 
\end{proof}

We obtain the log cohomology class 
$\lam_{\rm dRW}=c_{1,{\rm dRW}}(L)$ of $L$
in $H^2(X_t,{\cal W}\Om^{\bul}_{X_t})(1)$.

\begin{conj}[{\bf Filtered log $p$-adic hard Lefschetz conjecture via log de Rham-Witt complexes}]
\label{conj:lhwlc}
The following cup product 
\begin{equation*} 
\lam^i_{\rm dRW} \col 
H^{d-i}_{\rm crys}(X_t/{\cal W}(t))\otimes_{\mab Z}{\mab Q}
\lo H^{d+i}_{\rm crys}(X_t/{\cal W}(t))(i)
\otimes_{\mab Z}{\mab Q} \quad (i\in {\mab N})
\tag{2.5.19.1}\label{eqn:fcwpl} 
\end{equation*}
is an isomorphism. 
In fact, $\lam^i_{\rm dRW}$ 
is the following isomorphism of filtered sheaves: 
\begin{equation*} 
\lam^i_{\rm dRW} \col 
(H^{d-i}_{\rm crys}(X_{t}/{\cal W}(t))\otimes_{\mab Z}{\mab Q},P) 
\os{\sim}{\lo} (H^{d+i}(X_t/{\cal W}(t))(i)
\otimes_{\mab Z}{\mab Q},P). 
\tag{2.5.19.2}\label{eqn:fiwqpl} 
\end{equation*}
\end{conj}

\begin{prop} 
The conjecture {\rm (\ref{conj:lhwlc})} is equivalent to 
{\rm (\ref{conj:lhlc})} in the case  $S=s$ and $T={\cal W}(t)$.   
\end{prop}
\begin{proof} 
This follows from the commutative diagram (\ref{cd:drcrpic}). 
\end{proof}

\chapter{Weight filtrations on log isocrystalline cohomology sheaves} 

\parno 
In this chapter we give the definition of 
a truncated simplicial base change of SNCL schemes 
and the definition of the admissible immersion. 
These notions give us an essentially important framework in this book.  
By assuming the existence of an admissible immersion, 
we construct a filtered complex 
which is a generalization of the zariskian $p$-adic filtered Steenbrink complex 
(modulo torsion) defined in the Chapter I. 
To construct this filtered complex, 
we need to calculate the graded complex of 
the filtered complex of the base change above.  
We also give the definitions of a family of 
successive split truncated simplicial SNCL schemes and  
a log smooth split truncated simplicial log scheme associated to this family.  
This log smooth split truncated simplicial log scheme 
gives us a typical example of a truncated simplicial base change 
of SNCL schemes and an admissible immersion.

\section{Truncated simplicial base changes of SNCL schemes}\label{sec:bcsncl}
In this section we give the definition of 
the truncated simplicial base change of SNCL schemes 
and the definition of the admissible immersion. 
These definitions have come to me by my effort to 
find a framework to treat a truncated simplicial SNCL scheme 
defined in \S\ref{sec:tdiai} and  
families of successive split truncated simplicial SNCL scheme 
defined in the next section at the same time. 
\par 
Let $\{Y_i\}_{i=1}^l$ be a set of 
finitely many formal SNC schemes over a formal scheme $T$. 
(Note that a (log) scheme can be considered 
as a formal (log) scheme by a trivial ideal sheaf.) 
Set 
\begin{equation*} 
X:= \coprod_{i=1}^lY_i \quad (l\in {\mab Z}_{\geq 1}).
\tag{3.1.0.1}\label{eqn:xyi}
\end{equation*}  
\medskip  
\parno 
In this case,  we have the following equality:  
\begin{align*} 
X^{(k)}=\coprod_{i=1}^lY^{(k)}_i \quad (k\in {\mab N}).
\tag{3.1.0.2}\label{ali:yst} 
\end{align*}  
Let $\{Y_{\lam(i)}\}_{\lam(i)}$ be a decomposition of $Y_i$ by 
smooth components of $Y_i/T$ $(0\leq i\leq l)$. 
The set $\{Y_{\lam(i)}\}_{i,\lam(i)}$
defines an orientation sheaf $\vp^{(k)}_{\rm zar}(X/T)$ 
$(k\in {\mab N})$ in $X^{(k)}_{\rm zar}$ 
(cf.~\cite[(3.1.4)]{dh2} and \cite[(2.2.18)]{nh2}).  
If $T$ is a closed subscheme of $T'$ defined by 
a quasi-coherent nil-ideal sheaf ${\cal J}$ 
which has a 
PD-structure $\del$, 
then $\vp^{(k)}_{\rm zar}(X/T)$ 
extends to an abelian sheaf 
$\vp^{(k)}_{\rm crys}(X/T')$ in 
$(X^{(k)}/(T',{\cal J},\del))_{\rm crys}$. 

\begin{defi}\label{defi:snb} 
(1) Let $X$ be a fine log formal scheme over a formal family $S$ of log points. 
We say that $X/S$ is a {\it base change of formal $($proper or projective$)$ SNCL schemes} if 
there exist finitely many morphisms $S\lo S_i$ $(1\leq i\leq l)$ 
to formal families of log points and finitely many  
formal (proper or projective) SNCL schemes $Y_i/S_i$'s  
such that there exists an isomorphism 
\begin{align*} 
X\os{\sim}{\lo} \coprod_{i=1}^l(Y_i\times_{S_i}S).
\tag{3.1.1.1}\label{ali:yxst} 
\end{align*}  
Moreover, if there exists a morphism $S_i\lo S_0$ 
to a formal family of log points for $1\leq \forall i\leq l$, 
then we say that $X/S$ is a 
{\it base change of formal $($proper or projective$)$ SNCL schemes augmented to} $S_0$. 
\par 
(2) Let $N$ be a nonnegative integer or $\infty$.  
Let $X_{\bul \leq N}$ be an $N$-truncated simplicial 
fine log formal scheme over a formal family $S$ of log points. 
We say that $X_{\bul \leq N}/S$ is an  
{\it $N$-truncated simplicial base change of formal $($proper or projective$)$ SNCL schemes} if 
$X_m/S$ is a base change of formal SNCL schemes for $0\leq \forall m \leq N$. 
Moreover, if $X_m/S$ $(0\leq \forall m\leq N)$ is 
a base change of formal $($proper or projective$)$ SNCL schemes augmented to 
a formal family $S_0$ of log points, 
then we say that $X_{\bul \leq N}/S$ is an 
{\it $N$-truncated simplicial base change of formal $($proper or projective$)$ SNCL schemes 
augmented to} $S_0$. 
\end{defi}

\begin{rema}\label{rema:ibs}
(1) Let the notations be as in (\ref{defi:snb}) (1). 
Since the morphism $Y_i\lo S_i$ is integral by \cite[(4.4)]{klog1}, 
we have an isomorphism 
$\os{\circ}{X}\simeq 
\coprod_{i=1}^l(\os{\circ}{Y}_i\times_{\os{\circ}{S}_i}\os{\circ}{S})$.
Consequently $\os{\circ}{X}$ satisfies the condition (\ref{eqn:xyi}).   
\par
(2) Let the notations be as in (\ref{defi:snb}) (2). 
For a morphism $S'\lo S$ of formal families of log points,  
the fiber product $X_{\bul \leq N}\times_SS'$ is 
an $N$-truncated simplicial base change of formal SNCL schemes. 
\par 
(3) Let the notations be as in (\ref{defi:snb}) (1). 
Let $X\os{\sim}{\lo} \coprod_{i=1}^kX_i$ be the disjoint union
of the connected components of $\os{\circ}{X}$. 
Then each $X_i$ is a base change of SNCL schemes.  
\end{rema}


\begin{defi}\label{defi:bca}
(1) Let $S$ be a formal family of log points 
(over a formal family $S_0$ of log points). 
Let $X/S$ be a base change of SNCL schemes 
(augmented to a formal family $S_0$ of log points). 
Let $(T,{\cal J},\del)$ be a formal log PD-enlargement of $S$. 
Let $X'_{\os{\circ}{T}_0}$ be the disjoint union of the member of 
an affine open covering of $X_{\os{\circ}{T}_0}$. 
Let $X_{\os{\circ}{T}_0,\bul}$ be the \v{C}ech diagram of $X'_{\os{\circ}{T}_0}$ 
over $X_{\os{\circ}{T}_0}$.    
Let $\iota\col X_{\os{\circ}{T}_0,\bul}\os{\sus}{\lo} {\cal P}_{\bul}$ be 
an immersion into a formally log smooth simplicial 
log formal scheme over $S(T)^{\nat}$ (augmented to $S_0(T)^{\nat}$). 
We say that $\iota$ is {\it admissible} if the following conditions are satisfied: 
\par 
(a) There exist finitely many morphisms $S\lo S_i$ $(1\leq i\leq l)$ 
to families of log points and finitely many SNCL schemes $Y_i/S_{i,\os{\circ}{T}_0}$ 
such that there exists an isomorphism 
$X_{\os{\circ}{T}_0}\os{\sim}{\lo} \coprod_{i=1}^l(Y_i\times_{S_{i,\os{\circ}{T}_0}}
S_{\os{\circ}{T}_0})$. Here $S_{i,\os{\circ}{T}_0}
:=S_i\times_{\os{\circ}{S}_i}{\os{\circ}{T}_0}$. 
\par 
(b) There exists the disjoint union $Y'_{i}$ of the member of 
an affine open covering of $Y_{i}$ such that 
$X_{\os{\circ}{T}_0,\bul} 
\os{\sim}{\lo} \coprod_{i=1}^l
(Y_{i\bul}\times_{S_{i,\os{\circ}{T}_0}}S_{\os{\circ}{T}_0})$ 
over the isomorphism in (a). 
Here $Y_{i\bul}$ is the \v{C}ech diagram of $Y'_i$ over $Y_i$. 
\par 
(c) There exist finitely many
formally log smooth simplicial log formal schemes ${\cal Q}_{i\bul}$ 
over $S_i(T)^{\nat}$ such that 
there exists an isomorphism ${\cal P}_{\bul}\os{\sim}{\lo} 
\coprod_{i=1}^l({\cal Q}_{i\bul}\times_{S_i(T)^{\nat}}S(T)^{\nat})$ 
fitting into the following commutative diagrams
\begin{equation*} 
\begin{CD} 
Y_{i\bul}@>{\subset}>> {\cal Q}_{i\bul}\\
@VVV @VVV\\
S_{i,\os{\circ}{T}_0}@>{\subset}>> S_i(T)^{\nat},  
\end{CD} 
\tag{3.1.3.1}\label{cd:yiqi}
\end{equation*}
\begin{equation*} 
\begin{CD} 
X_{\bul,\os{\circ}{T}_0}@>{\subset}>> {\cal P}_{\bul}\\
@V{\simeq}VV @VV{\simeq}V\\
\coprod_{i=1}^l(Y_{i\bul}\times_{S_{i,\os{\circ}{T}_0}}S_{\os{\circ}{T}_0})
@>{\subset}>>
\coprod_{i=1}^l({\cal Q}_{i\bul}\times_{S_{i}(T)^{\nat}}S(T)^{\nat}). 
\end{CD} 
\tag{3.1.3.2}\label{cd:xpqst}
\end{equation*} 
\par 
(2) Let $S$, $S_0$ and $(T,{\cal J},\del)$ be as in (1). 
Let $X/S$ be a base change of SNCL schemes. 
Let $X'_{\os{\circ}{T}_0}$ and $X_{\os{\circ}{T}_0,\bul}$ be as in (1). 
Let $\ol{\iota}\col X_{\os{\circ}{T}_0,\bul}\os{\sus}{\lo} \ol{\cal P}_{\bul}$ be 
an immersion into a formally log smooth simplicial log formal scheme over $\ol{S(T)^{\nat}}$ 
(augmented to $\ol{S_0(T)^{\nat}}$). 
We say that $\ol{\iota}$ is {\it admissible} 
if the conditions (a) and (b) hold in (1) and if the following condition instead of 
(c) is satisfied:  
\par 
$(\ol{{\rm c}})$ 
There exist finitely many formally log smooth simplicial log formal schemes  
$\ol{\cal Q}_{i\bul}$ 
over $\ol{S_i(T)^{\nat}}$ such that 
there exists an isomorphism $\ol{\cal P}_{\bul}\os{\sim}{\lo} 
\coprod_{i=1}^l(\ol{\cal Q}_{i\bul}\times_{\ol{S_i(T)^{\nat}}}\ol{S(T)^{\nat}})$ 
fitting into the following commutative diagrams
\begin{equation*} 
\begin{CD} 
Y_{i\bul}@>{\subset}>> \ol{\cal Q}_{i\bul}\\
@VVV @VVV\\
S_{i,\os{\circ}{T}_0}@>{\subset}>> \ol{S_i(T)^{\nat}}, 
\end{CD} 
\end{equation*}
\begin{equation*} 
\begin{CD} 
X_{\os{\circ}{T}_0,\bul}@>{\subset}>> \ol{\cal P}_{\bul}\\
@V{\simeq}VV @VV{\simeq}V\\
\coprod_{i=1}^l(Y_{i\bul}\times_{S_{i,\os{\circ}{T}_0}}S_{\os{\circ}{T}_0})
@>{\subset}>>\coprod_{i=1}^l(\ol{\cal Q}_{i\bul}\times_{\ol{S_i(T)^{\nat}}}
\ol{S(T)^{\nat}}). 
\end{CD} 
\end{equation*} 
\par 
(3)  Let $S$, $S_0$ and $(T,{\cal J},\del)$ be as in (1). 
Let $N$ be a nonnegative integer or $\infty$.  
Let $X_{\bul \leq N}$ be an $N$-truncated simplicial base change of 
formal SNCL schemes over $S$ (augmented to $S_0)$. 
Assume that $X_{\bul \leq N,\os{\circ}{T}_0}$ has the disjoint union 
$X'_{\bul \leq N,\os{\circ}{T}_0}$ of 
the member of an affine $N$-truncated simplicial open covering of 
$X_{\bul \leq N,\os{\circ}{T}_0}$. 
Let $X_{\bul \leq N,\bul,\os{\circ}{T}_0}$ be the \v{C}ech diagram of 
$X'_{\bul \leq N,\os{\circ}{T}_0}$. 
Let $\iota \col X_{\bul \leq N,\bul,\os{\circ}{T}_0}\os{\sus}{\lo} {\cal P}_{\bul \leq N,\bul}$ 
be an immersion into a formally log smooth $(N,\infty)$-truncated 
bisimplicial log formal scheme over $S(T)^{\nat}$. 
Then we say that $\iota$ is {\it admissible} if 
the immersion $X_{m\bul,\os{\circ}{T}_0} 
\os{\sus}{\lo} {\cal P}_{m\bul}$ over $S(T)^{\nat}$ $(0\leq \forall m\leq N)$ 
is admissible. 
\par 
(4) Let $S$ and $S_0$ be as in (2). 
Let $N$ be a nonnegative integer or $\infty$.  
Let $X_{\bul \leq N}$ be an $N$-truncated simplicial base change of 
formal SNCL schemes over $S$ 
(augmented to $S_0$). 
Assume that $X_{\bul \leq N,\os{\circ}{T}_0}$ has the disjoint union 
$X'_{\bul \leq N,\os{\circ}{T}_0}$ of 
the member of an affine $N$-truncated simplicial open covering of 
$X_{\bul \leq N,\os{\circ}{T}_0}$. 
Let $X_{\bul \leq N,\bul,\os{\circ}{T}_0}$ be the \v{C}ech diagram of 
$X'_{\bul \leq N,\os{\circ}{T}_0}$. 
Let $\ol{\iota} \col X_{\bul \leq N,\bul,\os{\circ}{T}_0}\os{\sus}{\lo} \ol{\cal P}_{\bul \leq N,\bul}$ 
be an immersion into a formally log smooth 
$(N,\infty)$-truncated bisimplicial log formal scheme over $\ol{S(T)^{\nat}}$. 
Then we say that $\ol{\iota}$ is {\it admissible} if 
the immersion $X_{m\bul} \os{\sus}{\lo} \ol{\cal P}_{m\bul}$ over $\ol{S(T)^{\nat}}$ 
$(0\leq \forall m\leq N)$ is admissible. 
\end{defi}  

\begin{exem}\label{exem:xpm} 
(1) Let the notations be as in the Chapter I. 
Let $X_{\bul \leq N,\os{\circ}{T}_0}/S_{\os{\circ}{T}_0}$ be 
an $N$-truncated simplicial formal SNCL scheme 
which has the disjoint union of 
an affine $N$-truncated simplicial open covering of 
$X_{\bul \leq N,\os{\circ}{T}_0}$. 
Then, by (\ref{prop:xbn}), we have an admissible immersion 
$\ol{\iota} \col X_{\bul \leq N,\bul,\os{\circ}{T}_0}\os{\sus}{\lo} \ol{\cal P}_{\bul \leq N,\bul}$ 
over $\ol{S(T)^{\nat}}$. 
\par 
(2) Let the notations be as in (\ref{defi:bca}) (3). 
Let $S'\lo S$ be a morphism of formal families of log points. 
Let $(T',{\cal J}',\del')\lo (T,{\cal J},\del)$ be a morphism of log formal 
PD-enlargements over $S'\lo S$. 
Then the base change morphism 
$X_{\bul \leq N,\bul,S'_{\os{\circ}{T}{}'_0}}:=
X_{\bul \leq N,\bul,\os{\circ}{T}_0}\times_{S_{\os{\circ}{T}{}_0}}
S'_{\os{\circ}{T}{}'_0}
\os{\sus}{\lo} {\cal P}_{\bul \leq N,\bul,S'(T')^{\nat}}$ 
is an admissible immersion. 
\par 
(3) Let the notations be as in (2) and (\ref{defi:bca}) (4). 
Then the base change morphism 
$X_{\bul \leq N,\bul,S'_{\os{\circ}{T}{}'_0}}\os{\sus}{\lo} 
\ol{\cal P}_{\bul \leq N,\bul,\ol{S'(T')^{\nat}}}$ 
is an admissible immersion. 
\par 
(4) In (\ref{theo:thenad}) below we give an important example of 
an admissible immersion. 
\end{exem}

\par 
Let the notations be as in (\ref{defi:bca}) (4). 
We need the following notations for later sections. 
Consider the immersion 
$X_{m\bul} \os{\sus}{\lo} \ol{\cal P}_{m\bul}$ over $\ol{S(T)^{\nat}}$. 
By the definition of the admissibility, 
there exist a positive integer $l(m)$ and finitely many fine log formal schemes $S_{i(m)}$ $(1\leq i(m)\leq l(m))$, 
there exist finitely many morphisms $S\lo S_i$ $(1\leq i(m)\leq l(m))$ 
to families of log points and finitely many SNCL schemes $Y_{i(m)}/S_{i(m)}$ 
such that there exists an isomorphism 
\begin{align*} 
X_m\os{\sim}{\lo} \coprod_{i(m)=1}^{l(m)}(Y_{i(m)}
\times_{S_{i(m),\os{\circ}{T}_0}}S_{\os{\circ}{T}_0})
\tag{3.1.4.1}\label{ali:xmym}
\end{align*}  
and 
there exists the disjoint union $Y'_{i(m)}$ of the member of 
an affine truncated simplicial open covering of $Y_{i(m)}$ 
such that
\begin{align*} 
X_{m\bul,\os{\circ}{T}_0} \os{\sim}{\lo} 
\coprod_{i=1}^{l(m)}(Y_{i(m)\bul}\times_{S_{i(m),\os{\circ}{T}_0}}S_{\os{\circ}{T}_0})
\tag{3.1.4.2}\label{ali:xmbym}
\end{align*}  
over the isomorphism in (\ref{ali:xmym})
($Y_{i(m)\bul}$ is the \v{C}ech diagram of $Y'_{i(m)}$ over $Y_{i(m)}$)  
and we have the following commutative diagram 
\begin{equation*} 
\begin{CD} 
X_{m\bul,\os{\circ}{T}_0}@>{\subset}>> \ol{\cal P}_{m\bul}\\
@V{\simeq}VV @VV{\simeq}V\\
\coprod_{i(m)=1}^{l(m)}(Y_{i(m)\bul}\times_{S_{i(m),\os{\circ}{T}_0}}
S_{\os{\circ}{T}_0})@>{\subset}>>
\coprod_{i(m)=1}^{l(m)}(\ol{\cal Q}_{i(m)\bul}\times_{\ol{S_{i(m)}(T)^{\nat}}}\ol{S(T)^{\nat}}),  
\end{CD} 
\tag{3.1.4.3}\label{cd:yqst}
\end{equation*} 
where $\ol{\cal Q}_{i(m)\bul}$ is a formally log smooth simplicial formal log scheme over 
$\ol{S_{i(m)}(T)^{\nat}}$. 
Hence, by (\ref{prop:xpls}), we have  the following commutative diagram 
\begin{equation*} 
\begin{CD} 
X_{m\bul}@>{\subset}>> \ol{\cal P}{}^{\rm ex}_{m\bul}\\
@| @|\\
\coprod_{i(m)=1}^{l(m)}(Y_{i(m)\bul}\times_{S_{i(m),\os{\circ}{T}_0}}S_{\os{\circ}{T}_0})
@>{\subset}>>
\coprod_{i(m)=1}^{l(m)}(\ol{\cal Q}{}^{\rm ex}_{i(m)\bul}
\times_{\ol{S_{i(m)}(T)^{\nat}}}\ol{S(T)^{\nat}}).  
\end{CD} 
\tag{3.1.4.4}\label{cd:yqexst}
\end{equation*} 
Set ${\cal P}_{m\bul}:=\ol{\cal P}_{m\bul}\times_{\ol{S(T)^{\nat}}}S(T)^{\nat}$ and 
${\cal Q}_{i(m)\bul}:=\ol{\cal Q}_{i(m)\bul}\times_{\ol{S_{i(m)}(T)^{\nat}}}S_{i(m)}(T)^{\nat}$. 
By (\ref{prop:nexeo}) (2),  
$\ol{\cal Q}{}^{\rm ex}_{i(m)\bul}$ is a formal strict semistable log scheme 
over $\ol{S_{i(m)}(T)^{\nat}}$. Hence we have the following: 

\begin{prop}\label{prop:xxemn}  
$(1)$ For any $0\leq m\leq N$, 
${\cal P}^{\rm ex}_{m\bul}$ is 
a simplicial base change of formal SNCL schemes over $S(T)^{\nat}$. 
\par 
$(2)$ For any $0\leq m\leq N$,  
$\ol{\cal P}{}^{\rm ex}_{m\bul}$ a simplicial base changes of 
formal strict semistable log schemes over $\ol{S(T)^{\nat}}$.
\end{prop} 

\parno 
We have an equality 
$\os{\circ}{\cal P}{}^{{\rm ex},(k)}_{mn}
=\coprod_{i(m)=1}^{l(m)}\os{\circ}{\cal Q}{}^{{\rm ex},(k)}_{i(m)n}$
for nonnegative integers $k$ and $n$. 
In fact, by (\ref{lemm:knit}) (1), we have the simplicial formal scheme 
$\os{\circ}{\cal P}{}^{{\rm ex},(k)}_{m\bul}
=\coprod_{i(m)=1}^{l(m)}\os{\circ}{\cal Q}{}^{{\rm ex},(k)}_{i(m)\bul}$.

\begin{lemm}\label{lemm:mnl}
Let $m_i$ $(i=1,2)$ and $n_i$ $(i=1,2)$ be positive integers. 
Let $\Del \col {\mab N}\lo {\mab N}^n$ be the diagonal immersion. 
Let $m\cdot \col {\mab N}\lo {\mab N}$ be the multiplication by $m$. 
If there exists an isomorphism 
${\mab N}^{n_1}\oplus_{\Del,{\mab N},m_1\cdot}{\mab N}
\os{\sim}{\lo} 
{\mab N}^{n_2}\oplus_{\Del,{\mab N},m_2\cdot}{\mab N}$ of monoids 
fitting into the following commutative diagram
\begin{equation*} 
\begin{CD} 
{\mab N}^{n_1}\oplus_{\Del,{\mab N},m_1\cdot}{\mab N}
@>{\sim}>>{\mab N}^{n_2}\oplus_{\Del,{\mab N},m_2\cdot}{\mab N}\\
@V{{\rm 2nd~proj}.}VV @VV{{\rm 2nd~proj}.}V\\
{\mab N}@={\mab N}, 
\end{CD}
\tag{3.1.6.1}\label{cd:p}
\end{equation*} 
then $n_1=n_2$ and $m_1=m_2$. 
\end{lemm}
\begin{proof} 
Since $({\mab N}^{n_i}\oplus_{\Del,{\mab N},m_i\cdot}{\mab N})^{\rm gp}
={\mab Z}^{n_i}$, $n_1=n_2$. 
Set $n:=n_1$. 
Let $\al \col 
{\mab N}^{n}\oplus_{\Del,{\mab N},m_1\cdot}{\mab N}
\lo 
{\mab N}^{n}\oplus_{\Del,{\mab N},m_2\cdot}{\mab N}$ 
be an isomorphism of monoids. We may assume that $m_1\geq m_2$. 
Because ${\rm Aut}({\mab N}^n)={\mathfrak S}_n$ (\cite[p.~47]{nh3}), 
we may assume that 
$\al((e_i,0))=(e_i,l_i)$ for some $l_i\in {\mab N}$ and $1\leq \forall i\leq n$. 
Hence $\al((0,m_1))=\sum_{i=1}^n\al((e_i,0))
=\sum_{i=1}^n(e_i,l_i)=(0,m_2+\sum_{i=1}^nl_i)$. 
By (\ref{cd:p}), $m_1=m_2+\sum_{i=1}^nl_i$. 
Since $m_1\leq m_2$, $m_1=m_2$. 
\end{proof}

\begin{prop}\label{prop:sxst}
Two admissible immersions of 
$X_{\bul \leq N,\os{\circ}{T}_0}/S_{\os{\circ}{T}_0}$ 
over $S(T)^{\nat}$ $($resp.~$\ol{S(T)^{\nat}})$ 
are covered by an admissible immersions 
over $S(T)^{\nat}$ $($resp.~$\ol{S(T)^{\nat}})$. 
\end{prop}
\begin{proof} 
We prove (\ref{prop:sxst}) only in the case of admissible immersions over $\ol{S(T)^{\nat}}$. 
\par 
Let the notations be before (\ref{prop:sxst}). 
Let $X''_{\bul \leq N,\os{\circ}{T}_0}$ be another disjoint union of 
the member of an affine $N$-truncated simplicial open covering of 
$X_{\bul \leq N,\os{\circ}{T}_0}$. 
Let $X'_{\bul \leq N,\bul,\os{\circ}{T}_0}$ be the \v{C}ech diagram of 
$X''_{\bul \leq N,\os{\circ}{T}_0}$ over $S_{\os{\circ}{T}_0}$. 
Let $\ol{\iota}{}' \col X'_{\bul \leq N,\bul,\os{\circ}{T}_0}\os{\sus}{\lo} \ol{\cal P}{}'_{\bul \leq N,\bul}$ 
be another admissible immersion over $\ol{S(T)^{\nat}}$. 
Fix $0\leq m\leq N$. 
Then we have $l(m)'$, $Z_{i(m)}$, $Z'_{i(m)}$, $Z_{i(m)\bul}$, $S'_{i(m),\os{\circ}{T}_0}$ 
$\ol{\cal R}^{\rm ex}_{i(m)\bul}$ and $\ol{S'_{i(m)}(T)^{\nat}}$ for 
$\ol{\cal P}{}'_{m\bul}$ 
which are similar to $l(m)$, $Y_{i(m)}$, $Y'_{i(m)}$, $Y_{i(m)\bul}$, $S_{i(m),\os{\circ}{T}_0}$ 
$\ol{\cal Q}{}^{\rm ex}_{i(m)\bul}$ and $\ol{S_{i(m)}(T)^{\nat}}$, respectively 
and we have the following commutative diagram 
\begin{equation*} 
\begin{CD} 
X'_{m\bul,\os{\circ}{T}_0}@>{\subset}>> \ol{\cal P}{}'_{m\bul}\\
@V{\simeq}VV @VV{\simeq}V\\
\coprod_{i(m)=1}^{l(m)'}(Z_{i(m)\bul}\times_{S'_{i(m),\os{\circ}{T}_0}}
S_{\os{\circ}{T}_0})@>{\subset}>>
\coprod_{i(m)=1}^{l(m)'}(\ol{\cal R}_{i(m)\bul}\times_{\ol{S'_{i(m)}(T)^{\nat}}}\ol{S(T)^{\nat}}).   
\end{CD} 
\tag{3.1.7.1}\label{cd:ypqst}
\end{equation*} 
By decomposing $X'_{m,\os{\circ}{T}_0}$ and $X''_{m,\os{\circ}{T}_0}$ 
into connected components and 
taking large $l(m)$ and $l(m)'$, we may assume that $l(m)=l'(m)$ and 
$Y_{i(m)}\times_{S_{i(m),\os{\circ}{T}_0}}S_{\os{\circ}{T}_0}
=Z_{i(m)}\times_{S'_{i(m),\os{\circ}{T}_0}}S_{\os{\circ}{T}_0}$. 
Since the underlying schemes of 
$Y_{i(m)}\times_{S_{i(m),\os{\circ}{T}_0}}S_{\os{\circ}{T}_0}$ and 
$Z_{i(m)}\times_{S'_{i(m),\os{\circ}{T}_0}}S_{\os{\circ}{T}_0}$ are the same 
and $Y_{i(m)}$ and $Z_{i(m)}$ are SNCL schemes over 
$S_{i(m),\os{\circ}{T}_0}$ and $S'_{i(m),\os{\circ}{T}_0}$, 
we have an equality $S_{i(m),\os{\circ}{T}_0}=S'_{i(m),\os{\circ}{T}_0}$ by 
(\ref{lemm:mnl}) (since the underlying schemes of 
$S_{i(m),\os{\circ}{T}_0}$ and $S'_{i(m),\os{\circ}{T}_0}$ are $\os{\circ}{T}_0$). 
Set $X'''_{m,\os{\circ}{T}_0}:=X'_{m,\os{\circ}{T}_0}
\times_{X_{m,\os{\circ}{T}_0}}X''_{m,\os{\circ}{T}_0}$ $(0\leq m\leq N)$. 
Refining $X'''_{m,\os{\circ}{T}_0}$
(and hence refining the open coverings of $Y_{i(m)}$ and $Z_{i(m)}$), 
we may assume that $Y'_{i(m)}=Z'_{i(m)}$. 
In this way, we have another \v{C}ech diagram 
$X''_{\bul \leq N,\bul,\os{\circ}{T}_0}$ with morphisms 
$X''_{\bul \leq N,\bul,\os{\circ}{T}_0}\lo X_{\bul \leq N,\bul,\os{\circ}{T}_0}$ 
and $X''_{\bul \leq N,\bul,\os{\circ}{T}_0}\lo X'_{\bul \leq N,\bul,\os{\circ}{T}_0}$.  
Consider the following composite immersion 
$$X''_{\bul \leq N,\bul,\os{\circ}{T}_0}\os{\sus}{\lo} 
X_{\bul \leq N,\bul,\os{\circ}{T}_0}\times_{S_{\os{\circ}{T}_0}}
X'_{\bul \leq N,\bul,\os{\circ}{T}_0}\os{\sus}{\lo} 
\ol{\cal P}_{\bul \leq N,\bul}\times_{\ol{S(T)^{\nat}}}\ol{\cal P}{}'_{\bul \leq N,\bul}.$$ 
Then we have the following commutative diagram 
\begin{equation*} 
\begin{CD} 
X_{m,\bul,\os{\circ}{T}_0}@>{\subset}>>\ol{\cal P}_{m\bul}\\
@AAA @AAA \\
X''_{m,\bul,\os{\circ}{T}_0}@>{\subset}>>
(\ol{\cal P}_{m\bul}\times_{\ol{S(T)}}\ol{\cal P}{}'_{m\bul})^{\rm ex}\\
@VVV @VVV \\
X'_{m,\bul,\os{\circ}{T}_0}@>{\subset}>>\ol{\cal P}{}'_{m\bul}. 
\end{CD} 
\tag{3.1.7.2}\label{cd:ppppxmmm}
\end{equation*} 
By the argument above, 
 we have the following commutative diagram 
\begin{equation*} 
\begin{CD} 
X''_{m\bul,\os{\circ}{T}_0}@>{\subset}>> 
(\ol{\cal P}{}_{m\bul}\times_{\ol{S(T)^{\nat}}}\ol{\cal P}{}'_{m\bul})^{\rm ex}\\
@V{\simeq}VV @VV{\simeq}V\\
\coprod_{i(m)=1}^{l(m)}(Y_{i(m)\bul}\times_{S_{i(m),\os{\circ}{T}_0}}
S_{\os{\circ}{T}_0})@>{\subset}>>
\coprod_{i(m)=1}^{l(m)}((\ol{\cal Q}_{i(m)\bul}\times_{\ol{S_{i(m)}(T)^{\nat}}}
\ol{\cal R}_{i(m)\bul})^{\rm ex}\times_{\ol{S_{i(m)}(T)^{\nat}}}\ol{S(T)^{\nat}}). 
\end{CD} 
\tag{3.1.7.3}\label{cd:yppdqst}
\end{equation*} 
Hence the middle immersion 
in (\ref{cd:ppppxmmm}) is an admissible immersion. 
\end{proof}

\section{Families of successive 
split truncated simplicial SNCL schemes and 
bisimplicial (exact) immersions II}\label{sec:beii} 
In this section we give the definition of a family of 
successive split truncated simplicial SNCL schemes (see (\ref{defi:xsnndf}) below). 
We obtain a log smooth split truncated simplicial log scheme associated to this family. 
This is a typical example of a truncated simplicial base change of SNCL schemes 
in the previous section. 
Though our results in this section are not the imitations of du Bois'   
construction in \cite{db}, we are influenced by du Bois' article. 
The main result in this section is to show that there exists an admissible immersion 
for the truncated simplicial log scheme above. 
\par  
Let $N$ be a nonnegative integer.  
Let $p$ be a prime number. 
Let $S_0$ be a (($p$-adic) formal) family of log points.

\begin{defi}\label{defi:fld} 
Let 
\begin{equation*} 
S_N \lo S_{N-1} \lo \cdots \lo S_0
\tag{3.2.1.1}\label{eqn:stcnniseq}
\end{equation*} 
be a sequence of (($p$-adic) formal) families of log points.  
Let $s_N$ be a point of $\os{\circ}{S}_N$ and 
let $s_0$ be the image of $s_N$ in $S_0$. 
If the natural morphism
$M_{S_0,s_0}/{\cal O}^*_{S_0,s_0} \lo M_{S_N,s_N}
/{\cal O}^*_{S_N,s_N}$ is an isomorphism 
for any point $s_N\in \os{\circ}{S}_N$, 
we call (\ref{eqn:stcnniseq}) a 
{\it sequence of nonramified families of log points} 
and, otherwise, we call (\ref{eqn:stcnniseq}) a 
{\it sequence of ramified families of log points}.   
For simplicity of notation, 
we denote the sequence (\ref{eqn:stcnniseq}) by $\{S_m\}_{m=0}^N$. 
We define a morphism 
\begin{align*} 
\{S'_{m'}\}_{m=0}^{N'}\lo \{S_m\}_{m=0}^N
\tag{3.2.1.2}\label{eqn:stcs0s}
\end{align*} 
of sequences of (($p$-adic) formal) families of log points as follows. 
\par 
If $N'\not= N$, then the set of the morphisms is empty and 
if $N'=N$, then the morphism is, by definition, 
the following commutative diagrams: 
\begin{align*} 
\begin{CD} 
S'_{N} @>>> S'_{N-1} @>>> \cdots @>>> S'_0\\
@VVV @VVV @. @VVV\\
S_N @>>> S_{N-1} @>>> \cdots @>>> S_0. 
\end{CD} 
\tag{3.2.1.3}\label{cd:snsn0}
\end{align*}
\end{defi} 

\parno  
For simplicity of notation, 
denote the sequence (\ref{eqn:stcnniseq}) 
by $\{S_m\}_{m=0}^N$. 
For a fine $m$-truncated simplicial log scheme 
$Y_{\bul \leq m}$ over $S_m$,  
let $s^{m-1}_i \col Y_{m-1} \lo Y_m$ 
$(0\leq m\leq N, 0\leq i \leq m-1)$ be 
the degeneracy morphism 
before (\ref{lemm:lisj}).   
Let $\{X(m)_{\bul \leq m}/S_m\}_{m=0}^N$ 
be a family of log smooth split truncated simplicial log schemes 
such that the intersection of the complements of 
$s^{m-1}_i(X(m)_{m-1})$ $(0 \leq i \leq m-1)$ 
in $X(m)_m$ 
is an SNCL scheme $N_m$ over $S_m$ $(0\leq m\leq N)$.

\begin{defi}\label{defi:xsnndf}
(1) We say that 
$\{X(m)_{\bul \leq m}/S_m\}_{m=0}^N$ 
is a {\it family of 
successive split truncated simplicial SNCL schemes 
with respect to the sequence} 
(\ref{eqn:stcnniseq}) 
if 
\begin{equation*} 
X(m)_{\bul \leq m-1}=
X(m-1)_{\bul \leq m-1}
\times_{S_{m-1}}S_m
\tag{3.2.2.1}\label{eqn:xm12sm} 
\end{equation*}  
and 
\begin{equation*} 
X(m)_m=\coprod_{0\leq l \leq m}\coprod_{[m] 
\twoheadrightarrow [l]}(N_l\times_{S_l}S_m)  
\tag{3.2.2.2}\label{eqn:xmm12sm}  
\end{equation*} 
for all $0\leq m\leq N$, 
where the subscripts $[m]\twoheadrightarrow [l]$'s 
run through the surjective nondecreasing morphisms $[m]\lo [l]$'s. 
\par 
(2) Set $X_{\bul \leq N}:=X(N)_{\bul \leq N}$ 
and $S:=S_N$. 
We call $X_{\bul \leq N}/S$ 
the {\it associated $($log smooth$)$ split $N$-truncated simplicial log 
scheme} to the family 
$\{X(m)_{\bul \leq m}/S_m\}_{m=0}^N$.  
\end{defi} 

\begin{prop}\label{prop:xmbc}
The log scheme $X(m)_{\bul \leq m}/S_m$ $(0\leq m\leq N)$ is an 
$m$-truncated simplicial base change of 
SNCL schemes augmented to $S_0$. 
\end{prop} 
\begin{proof} 
By the formulas (\ref{eqn:xm12sm}) and (\ref{eqn:xmm12sm}), 
we see that 
\begin{equation*} 
X(m)_k=\coprod_{0\leq l \leq k}
\coprod_{[k]\twoheadrightarrow [l]}
(N_l\times_{S_l}S_m) \quad (0\leq k\leq m).  
\tag{3.2.3.1}\label{eqn:xmnlm} 
\end{equation*}  
This implies (\ref{prop:xmbc}). 
\end{proof}  

\begin{rema}\label{rema:knm}  
(1) Let the notations be as in (\ref{defi:xsnndf}). 
When $S_m=S_0$ for all $0\leq m\leq N$,  
the log scheme $X_{\bul \leq N}/S$ associated to $\{X(m)_{\bul \leq m}/S_m\}_{m=0}^N$ 
is nothing but a split $N$-truncated simplicial SNCL scheme 
$X_{\bul \leq N}/S$.  
\par 
(2) Set $S_N^m:=S_m\times_{\os{\circ}{S}_m}\os{\circ}{S}_N$. 
Then we have a sequence 
\begin{equation*} 
S=S_N \lo S^{N-1}_N \lo \cdots \lo S^0_N
\tag{3.2.4.1}\label{eqn:stnneq}
\end{equation*} 
of (($p$-adic) formal) families of log points and 
$N_m\times_{S_m}S_N^m$ is an SNCL scheme over $S_N^m$.  
Set 
\begin{equation*} 
X_N(m)_m=\coprod_{0\leq l \leq m}\coprod_{[m] 
\twoheadrightarrow [l]}(N_l\times_{S_l}S_N^m)  
\tag{3.2.4.2}\label{eqn:xmm2sm}  
\end{equation*} 
for all $0\leq m\leq N$. Then we have a new family 
$\{X_N(m)_{\bul \leq m}\}_{m=0}^N$ of 
successive split truncated simplicial SNCL schemes 
with respect to the sequence (\ref{eqn:stnneq}).  
Moreover $X_N(N)_{\bul \leq N}=X(N)_{\bul \leq N}$. 
Because we are only interested in $X(N)_{\bul \leq N}$, 
we can assume that $\os{\circ}{S}_m=\os{\circ}{S}_N$ for $0\leq \forall m\leq N$ 
without loss of generality. 
\par 
(3) The sequence (\ref{eqn:stnneq}) is equivalent to the following two data: 
\par
(a) the (($p$-adic) formal) log scheme $S_N$, 
\par
(b) the finite increasing filtration of log structures of submonoids of $M_{S_N}$  
each of which is generated by one element modulo ${\cal O}_{S_N}^*$
\begin{align*} 
M_0\subset M_1\subset M_2 \subset \cdots \subset M_{N-1}\subset M_{S_N}. 
\end{align*}
\end{rema} 

\par 
Let the notations be as in (\ref{defi:xsnndf}). 
Then, for $k\in {\mab N}$ and $0\leq m\leq N$,  
\begin{equation*} 
\os{\circ}{X}{}^{(k)}_m=
\coprod_{0\leq l \leq m}
\coprod_{[m]\twoheadrightarrow [l]}
(\os{\circ}{N}{}^{(k)}_l
\times_{\os{\circ}{S}_l}\os{\circ}{S}).  
\tag{3.2.4.3}\label{eqn:dfm}
\end{equation*} 

\par 
Let the notations be as in (\ref{eqn:stcnniseq}). 
Let $T_0$ be a fine log scheme over $S$.  
Then we have the following  sequence of fine log schemes: 
\begin{align*} 
S_{\os{\circ}{T}_0}\lo S_{N-1,\os{\circ}{T}_0}\lo \cdots \lo S_{0,\os{\circ}{T}_0}.
\tag{3.2.4.4}\label{ali:stnttt} 
\end{align*} 
Set $X(m)_{\bul \leq m,\os{\circ}{T}_0}:=
X(m)_{\bul \leq m}\times_{S_m}S_{m,\os{\circ}{T}_0}$. 
Then we have a family 
$\{X(m)_{\bul \leq m,\os{\circ}{T}_0}/S_{m,\os{\circ}{T}_0}\}_{m=0}^N$ of 
successive split truncated simplicial SNCL schemes 
with respect to the sequence (\ref{ali:stnttt}).  


\begin{defi}\label{defi:fso}  
Let  
\begin{equation*} 
S_N(=S) \lo S_{N-1} \lo \cdots \lo S_0
\tag{3.2.5.1}\label{eqn:stcsseq}
\end{equation*} 
be as in (\ref{eqn:stcnniseq}). 
Let $(T,{\cal J},\del)$ be a fine log (($p$-adic) formal) PD-enlargement of $S$. 
Set $T_0:=\ul{\rm Spec}^{\log}_T({\cal O}_T/{\cal J})$. 
Then we have the sequence (\ref{ali:stnttt}) and 
the following sequence of  
fine log (($p$-adic) formal) families of log points:  
\begin{equation*} 
(S(T)^{\nat},{\cal J},\del) :=(S_N(T)^{\nat},{\cal J},\del)
\lo (S_{N-1}(T)^{\nat},{\cal J},\del) \lo \cdots \lo (S_0(T)^{\nat},{\cal J},\del).   
\tag{3.2.5.2}\label{eqn:stcseq}
\end{equation*} 
If the sequence  (\ref{eqn:stcsseq}) is nonramified (resp.~ramified), 
then we call (\ref{eqn:stcseq}) a 
{\it nonramified $($resp.~ramified$)$ sequence of 
$((p$-adic$)$ formal$)$ PD-families of log points}.   
\end{defi}

\par 
In the following,  we assume that the transitive morphism 
$\os{\circ}{S}_m \lo \os{\circ}{S}_{m-1}$ $(0\leq m\leq N)$ 
is affine. 
Let 
$\{X(m)_{\bul \leq m}/S_m\}_{m=0}^N$ 
be a family of successive 
split truncated simplicial SNCL schemes 
with respect to the sequence (\ref{eqn:stcnniseq}).   
Let us construct admissible immersions  
of $X_{\bul \leq N,\os{\circ}{T}_0}/S(T)^{\nat}$ and 
$X_{\bul \leq N,\os{\circ}{T}_0}/\ol{S(T)^{\nat}}$.   
For this purpose, we reexamine the construction of 
$G_{\bul \leq N,\bul}$ after (\ref{prop:ytlft}).   
\par 
First we reconstruct the disjoint union $X(m)'_{\bul \leq m}$
of the member of an affine $m$-truncated simplicial open covering of 
$X(m)_{\bul \leq m}$ $(0\leq m\leq N)$ (cf.~\cite[(6.1)]{nh3}). 
\par 
Let $N'_0$ be the disjoint union of 
log affine open (formal) subschemes which cover
$N_0$ and whose images in $S_0$ are contained in 
log affine open (formal) subschemes of $S_0$. 
Then $N'_0\times_{S_0}S_l$ 
$(0\leq l\leq N)$ is the disjoint union of 
the member of an affine open covering of 
$N_0\times_{S_0}S_l$ 
because the morphism  
$\os{\circ}{S}_l\lo \os{\circ}{S}_0$ is affine. 
Let $m$ be a positive integer.
Assume that we are given $X(m-1)'_{\bul \leq m-1,\os{\circ}{T}_0}$. 
We construct $X(m)'_{\bul \leq m,\os{\circ}{T}_0}$ as follows.  
First we set 
$X(m)'_{\bul \leq m-1,\os{\circ}{T}_0}:=
X(m-1)'_{\bul \leq m-1,\os{\circ}{T}_0}\times_{S_{m-1}}S_m$.
Hence 
\begin{equation*} 
X(m)'_{k,\os{\circ}{T}_0}=X(k)'_{k,\os{\circ}{T}_0}\times_{S_k}S_m 
\quad (0\leq k\leq m-1). 
\tag{3.2.5.3}\label{eqn:skm}
\end{equation*}  
For a nonnegative integer $r$, 
let $\Del_{[r]}$ 
$(r\in {\mab N})$ 
be the category of objects of $\Del$ augmented to $[r]$. 
For a nonnegative integer $l$, 
let $\Del_{l[r]}$ 
be a full subcategory of $\Del_{[r]}$ 
whose objects are $[q] \lo [r]$'s 
such that $q\leq l$. 
The inverse limit of a finite inverse 
system exists in the category of fine log (formal) schemes 
(over a fine log (formal) scheme) (\cite[(1.6), (2.8)]{klog1}). 
Recall the following explicit description of 
the coskelton 
(e.~g., \cite[${\rm V}^{\rm bis}$ (3.0.1.2)]{sga4-2}):
\begin{equation*}
{\rm cosk}^{S_{m-1,\os{\circ}{T}_0}}_{m-1}
(X(m-1)_{\bul \leq m-1,\os{\circ}{T}_0})_r
={\vpl}_{[q]\in{\Del}_{(m-1)[r]}}X(m-1)_{q,\os{\circ}{T}_0} \quad (l, r\in {\mab N}),
\tag{3.2.5.4}\label{eqn:cokynr}
\end{equation*}
where the projective limits are taken in 
the category of log (formal) schemes over 
$S_{m-1,\os{\circ}{T}_0}$. 
The log scheme 
${\rm cosk}^{S_{m-1,\os{\circ}{T}_0}}_{m-1}
(X(m-1)'_{\bul \leq m-1,\os{\circ}{T}_0})_m$ 
is the disjoint union of 
the members of a log open covering  of 
${\rm cosk}^{S_{m-1,\os{\circ}{T}_0}}_{m-1}
(X(m-1)_{\bul \leq m-1,\os{\circ}{T}_0})_m$. 
Set $N_{m,\os{\circ}{T}_0}:=N_m\times_{S_m}S_{m,\os{\circ}{T}_0}$. 
Consider the natural composite morphism 
$N_{m,\os{\circ}{T}_0} \os{\sus}{\lo} X(m)_{m,\os{\circ}{T}_0} 
\lo 
{\rm cosk}^{S_{m,\os{\circ}{T}_0}}_{m-1}
(X(m)_{\bul \leq m-1,\os{\circ}{T}_0})_m= 
{\rm cosk}^{S_{m-1,\os{\circ}{T}_0}}_{m-1}
(X(m-1)_{\bul \leq m-1,\os{\circ}{T}_0})_m
\times_{S_{m-1,\os{\circ}{T}_0}}S_{m,\os{\circ}{T}_0}$ 
and a log affine open covering of $N_{m,\os{\circ}{T}_0}$
which refines the inverse image of the open covering of 
${\rm cosk}^{S_{m-1,\os{\circ}{T}_0}}_{m-1}(X(m-1)_{\bul \leq m-1,\os{\circ}{T}_0})_m
\times_{S_{m-1,\os{\circ}{T}_0}}S_{m,\os{\circ}{T}_0}$. 
Then the image of each member of the log affine open (formal) 
subschemes of $N_{m,\os{\circ}{T}_0}$ in $S_{m,\os{\circ}{T}_0}$ 
is contained in a log affine open (formal) subscheme of $S_{m,\os{\circ}{T}_0}$. 
Let $N'_{m,\os{\circ}{T}_0}$ be the disjoint union of the members of this open covering.
Then we have the following 
commutative diagram:
\begin{equation*}
\begin{CD}
N_{m,\os{\circ}{T}_0}' @>>> 
{\rm cosk}^{S_{m,\os{\circ}{T}_0}}_{m-1}
(X(m)'_{\bul \leq m-1,\os{\circ}{T}_0})_m \\
@VVV @VVV  \\ 
N_{m,\os{\circ}{T}_0} @>>> {\rm cosk}^{S_{m,\os{\circ}{T}_0}}_{m-1}
(X(m)_{\bul \leq m-1,\os{\circ}{T}_0})_m.
\end{CD}
\tag{3.2.5.5}\label{cd:spcyny}
\end{equation*}
Set 
\begin{align*} 
X(m)'_{m,\os{\circ}{T}_0}=\coprod_{0\leq l \leq m}\coprod_{[m] 
\twoheadrightarrow [l]}(N_{l,\os{\circ}{T}_0}'\times_{S_{l,\os{\circ}{T}_0}}S_{m,\os{\circ}{T}_0}). 
\tag{3.2.5.6}\label{cd:spcnxy}
\end{align*}  
Since $X(m)_{\bul \leq m,\os{\circ}{T}_0}$ is split, $X(m)'_{\bul \leq m,\os{\circ}{T}_0}$ 
is a desired log scheme over $S_{m,\os{\circ}{T}_0}$, and we have a natural morphism 
$X(m)'_{\bul \leq m,\os{\circ}{T}_0} \lo X(m)_{\bul \leq m,\os{\circ}{T}_0}$ 
of simplicial log (formal) schemes over 
$S_{m,\os{\circ}{T}_0}$ by 
(\ref{cd:spcyny}) and 
\cite[${\rm V}^{\rm bis}$ (5.1.3)]{sga4-2}. 
\par
Next we reexamine the construction of $\Gam'_{\bul}$ 
in (\ref{prop:ytlft}). 
\par  
Set $X'_{\bul \leq N,\os{\circ}{T}_0}:=X(N)'_{\bul \leq N,\os{\circ}{T}_0}$.  
By (\ref{cd:spcnxy}),  
$X'_{N,\os{\circ}{T}_0}=\coprod_{0\leq m \leq N}
\coprod_{[N] \twoheadrightarrow [m]}
(N'_{m,\os{\circ}{T}_0}\times_{S_{m,\os{\circ}{T}_0}}S_{\os{\circ}{T}_0})$.  
Let $N'_{m,\os{\circ}{T}_0} \os{\sus}{\lo} \ol{\cal N}{}'_m$ 
$(0\leq m \leq N)$ be an immersion into a log smooth (($p$-adic) formal) 
scheme  over $\ol{S_m(T)^{\nat}}$. In fact, by (\ref{eqn:xdpelda}), 
we can assume the existence of a (formal) SNCL lift $\ol{\cal N}{}'_m$ of $N'_m$ over 
$\ol{S_m(T)^{\nat}}$ if each of the log affine open (formal) subscheme of 
$N'_m$ $(0\leq m\leq N)$ is small.   
Set ${\cal N}'_m:=\ol{\cal N}{}'_m\times_{\ol{S_m(T)^{\nat}}}S_m(T)^{\nat}$. 
Set $`{\cal X}_N:=\coprod_{0\leq m \leq N}
\coprod_{[N] \twoheadrightarrow [m]}
({\cal N}_m'\times_{S_m(T)^{\nat}}S(T)^{\nat})$ and 
$`\ol{\cal X}_N:=\coprod_{0\leq m \leq N}
\coprod_{[N] \twoheadrightarrow [m]}
(\ol{\cal N}{}'_m\times_{\ol{S_m(T)^{\nat}}}\ol{S(T)^{\nat}})$. 
Set ${\cal P}'_m:=
\prod_{{\rm Hom}_{\Del}([N],[m])}`{\cal X}_N$ 
and 
$\ol{\cal P}{}'_m:=
\prod_{{\rm Hom}_{\Del}([N],[m])}`\ol{\cal X}_N$ 
$(m\in {\mab N})$. 
Then we have the following natural composite immersion  
\begin{align*} 
X'_{m,\os{\circ}{T}_0} & \os{{\rm diag}., \sus}{\lo}
\prod_{{\rm Hom}_{\Del}([N],[m])}X'_{m,\os{\circ}{T}_0}
\os{\sus}{\lo} 
\prod_{{\rm Hom}_{\Del}([N],[m])}X'_{N,\os{\circ}{T}_0} \tag{3.2.5.7}\label{ali:xgpmn}\\
&\os{\sus}{\lo} 
\prod_{{\rm Hom}_{\Del}([N],[m])}`{\cal X}_N= 
{\cal P}'_m \os{\sus}{\lo} \ol{\cal P}{}'_m.
\end{align*}
Set $X_{mn,\os{\circ}{T}_0}
:={\rm cosk}_0^{X_{m,\os{\circ}{T}_0}}(X'_{m,\os{\circ}{T}_0})_n$, 
${\cal P}_{mn}:={\rm cosk}_0^{S(T)^{\nat}}({\cal P}'_m)_n$ 
and $\ol{\cal P}_{mn}:={\rm cosk}_0^{\ol{S(T)^{\nat}}}(\ol{\cal P}{}'_m)_n$. 
Then we have the following natural immersions  
\begin{equation*} 
X_{mn,\os{\circ}{T}_0} \os{\sus}{\lo} {\cal P}_{mn} \os{\sus}{\lo} \ol{\cal P}_{mn}.
\tag{3.2.5.8}\label{eqn:xgmn}
\end{equation*} 
In conclusion, we obtain the desired immersions 
\begin{align*} 
X_{\bul \leq N,\bul,\os{\circ}{T}_0}\os{\sus}{\lo} {\cal P}_{\bul \leq N,\bul} 
\os{\sus}{\lo}  \ol{\cal P}_{\bul \leq N,\bul}
\tag{3.2.5.9}\label{ali:xgsmn}
\end{align*} 
over $S_{\os{\circ}{T}_0}\os{\sus}{\lo} S(T)^{\nat}$ and 
over $S_{\os{\circ}{T}_0}\os{\sus}{\lo} \ol{S(T)^{\nat}}$.  
Because 
\begin{align*} 
X'_{m,\os{\circ}{T}_0}=\coprod_{0\leq l \leq m}
\coprod_{[m] \twoheadrightarrow [l]}
(N_{l,\os{\circ}{T}_0}'\times_{S_{l,\os{\circ}{T}_0}}S_{\os{\circ}{T}_0}),
\end{align*} 
\begin{align*} 
&(\underset{n~{\rm pieces}}{\underbrace{
N'_{l,\os{\circ}{T}_0}\times_{X_{l,\os{\circ}{T}_0}}
\cdots  \times_{X_{l,\os{\circ}{T}_0}}N'_{l,\os{\circ}{T}_0}}})
\times_{S_{l,\os{\circ}{T}_0}}S_{\os{\circ}{T}_0}=
(\underset{n~{\rm pieces}}{\underbrace{
N'_{l,\os{\circ}{T}_0}\times_{N_{l,\os{\circ}{T}_0}}\cdots 
\times_{N_{l,\os{\circ}{T}_0}}N'_{l,\os{\circ}{T}_0}}})\times_{S_{l,\os{\circ}{T}_0}}S_{\os{\circ}{T}_0} \\
&\quad  ([m] \twoheadrightarrow [l]) 
\end{align*} 
is a component of $X_{mn,\os{\circ}{T}_0}$, that is, 
there exists an open and closed log (formal) subscheme 
$V$ of $X_{mn,\os{\circ}{T}_0}$ such that 
$((N'_{l,\os{\circ}{T}_0}
\times_{N_{l,\os{\circ}{T}_0}}\cdots  \times_{N_{l,\os{\circ}{T}_0}}N'_{l,\os{\circ}{T}_0})
\times_{S_{l,\os{\circ}{T}_0}}S_{\os{\circ}{T}_0})\coprod V=X_{mn,\os{\circ}{T}_0}$.

\begin{prop}\label{prop:keyp2} 
Extract $N'_l\times_{S_{l,\os{\circ}{T}_0}}S_{\os{\circ}{T}_0}$ 
in the source of the immersion {\rm (\ref{ali:xgpmn})}. 
Then there exists the following commutative diagram$:$
\begin{equation*}  
\begin{CD}  
X'_{m,\os{\circ}{T}_0}@>{\subset}>>{\cal P}'_m@>{\subset}>>\ol{\cal P}{}'_m\\
@A{\bigcup}AA @AA{\bigcup}A @AA{\bigcup}A\\
N'_{l,\os{\circ}{T}_0}\times_{S_{l,\os{\circ}{T}_0}}S_{\os{\circ}{T}_0}
@>{\subset}>>{\cal L}'_l\times_{S_l(T)^{\nat}}S(T)^{\nat} 
@>{\subset}>>\ol{\cal L}{}'_l\times_{\ol{S_l(T)^{\nat}}}\ol{S(T)^{\nat}}, 
\end{CD} 
\tag{3.2.6.1}\label{cd:xppm} 
\end{equation*}  
where ${\cal L}'_l$ and $\ol{\cal L}{}'_l$ are 
$((p$-adically$)$ formally$)$ log smooth schemes over $S_l(T)^{\nat}$ 
and $\ol{S_l(T)^{\nat}}$, respectively. 
Consequently there exists the following commutative diagram$:$
\begin{equation*}  
\begin{CD}  
X_{mn,\os{\circ}{T}_0}@>{\subset}>>{\cal P}_{mn}\\
@A{\bigcup}AA @AA{\bigcup}A \\
(N'_{l,\os{\circ}{T}_0}\times_{X_{l,\os{\circ}{T}_0}}\cdots 
\times_{X_{l,\os{\circ}{T}_0}}N'_{l,\os{\circ}{T}_0})
\times_{S_{l,\os{\circ}{T}_0}}S_{\os{\circ}{T}_0}
@>{\subset}>>({\cal L}'_l\times_{S_l(T)^{\nat}}\cdots 
\times_{S_l(T)^{\nat}}{\cal L}'_l)\times_{S_l(T)^{\nat}}S(T)^{\nat} 
\end{CD} 
\tag{3.2.6.2}\label{cd:xpsm} 
\end{equation*}   
\begin{equation*}  
\begin{CD}  
@>{\subset}>>\ol{\cal P}{}'_{mn}\\
@. @AA{\bigcup}A \\
@>{\subset}>>(\ol{\cal L}{}'_l\times_{\ol{S_l(T)^{\nat}}}\cdots 
\times_{\ol{S_l(T)^{\nat}}}\ol{\cal L}{}'_l)\times_{\ol{S_l(T)^{\nat}}}\ol{S(T)^{\nat}}. 
\end{CD} 
\end{equation*} 
\end{prop}
\begin{proof} 
Let $0\leq l \leq m \leq N$ be nonnegative integers. 
Consider a surjection $[m] \lo [l]$ in $\Del$. 
For an element $\gam \in {\rm Hom}_{\Del}([N],[m])$, 
we claim that the composite morphism 
$N'_{l,\os{\circ}{T}_0} \os{\sus}{\lo} X'_{m,\os{\circ}{T}_0} 
\os{X'(\gam)}{\lo} X'_{N,\os{\circ}{T}_0}$ 
$(\gam \in {\rm Hom}_{\Del}([N],[m]))$
factors through an open and closed 
log subscheme of $X'_{N,\os{\circ}{T}_0}$ which is the base change to $S_{l,\os{\circ}{T}_0}$ 
of the disjoint union of finitely many SNCL schemes 
over $S_{l',\os{\circ}{T}_0}$ for various $l'\leq l$. 
To prove this,  
let us investigate the morphism 
$X'(\gam) \col X'_{m,\os{\circ}{T}_0} \lo X'_{N,\os{\circ}{T}_0}$. 
Note that $\gam$ is the composite of standard degeneracy maps and face maps 
(\cite[(8.1.2)]{weib}).  
First consider the degeneracy morphism 
$s^{l-1}_i \col {X'_{l-1,\os{\circ}{T}_0}} \lo X'_{l,\os{\circ}{T}_0}$ 
$(l\in {\mab Z}_{>0}, 0\leq i \leq l-1)$ 
corresponding to the standard degeneracy map 
$\partial^i_l \col [l] \lo [l-1]$ before (\ref{lemm:lisj}). 
Then, by 
the proof of 
\cite[${\rm V}^{\rm bis}$ (5.1.3)]{sga4-2}, 
\begin{align*} 
s^{l-1}_i \col {X'_{l-1,\os{\circ}{T}_0}}=
\coprod_{0\leq k \leq l-1}\coprod_{[l-1] 
\twoheadrightarrow [k]}
(N_{k,\os{\circ}{T}_0}'\times_{S_{k,\os{\circ}{T}_0}}S_{\os{\circ}{T}_0}) 
\lo {X'_{l,\os{\circ}{T}_0}}=\coprod_{0\leq k \leq l}\coprod_{[l] 
\twoheadrightarrow [k]}(N_{k,\os{\circ}{T}_0}'
\times_{S_{k,\os{\circ}{T}_0}}S_{\os{\circ}{T}_0}) 
\end{align*} 
is induced by the identity ${\rm id} \col N_{k,\os{\circ}{T}_0}' \lo N_{k,\os{\circ}{T}_0}'$ 
corresponding to the components with indexes 
$[l-1] \twoheadrightarrow [k]$ 
and 
$[l] \os{\partial^i_l}{\lo} 
[l-1] \twoheadrightarrow [k]$. 
Hence $s^{l-1}_i\vert_{N'_{k,\os{\circ}{T}_0}\times_{S_{k,\os{\circ}{T}_0}}S_{\os{\circ}{T}_0}}$ 
$([l-1]\twoheadrightarrow [k])$ 
is nothing but the inclusion morphism 
\begin{equation*} 
N'_{k,\os{\circ}{T}_0}\times_{S_{k,\os{\circ}{T}_0}}S_{\os{\circ}{T}_0} 
\os{\sus}{\lo} 
\coprod_{0\leq k \leq l}
\coprod_{[l]\twoheadrightarrow [k]}
(N_{k,\os{\circ}{T}_0}'\times_{S_{k,\os{\circ}{T}_0}}S_{\os{\circ}{T}_0})
=X'_{l,\os{\circ}{T}_0}  
\tag{3.2.6.3}\label{eqn:nklxl} 
\end{equation*} 
corresponding to the index $[l]\twoheadrightarrow [k]$. 
We have checked that our claim holds for $s^{l-1}_i$. 
\par 
Next consider the standard face morphism 
$\del'^l_i \col X'_{l,\os{\circ}{T}_0} \lo X'_{l-1,\os{\circ}{T}_0}$ $(0\leq i \leq l)$ 
corresponding to the morphism 
$\sig^l_i \col [l-1] \lo [l]$ defined by 
the formulas $\sig^l_i(j)=j$ $(0\leq j < i)$ 
and $\sig^l_i(j)=j+1$ $(i\leq j \leq l)$. 
Let 
$\del''^l_i \col 
{\rm cosk}_{l-1}^{S_{\os{\circ}{T}_0}}(X'_{\bul \leq l-1,\os{\circ}{T}_0})_l 
\lo {\rm cosk}_{l-1}^{S_{\os{\circ}{T}_0}}
(X'_{\bul \leq l-1,\os{\circ}{T}_0})_{l-1}= X'_{l-1,\os{\circ}{T}_0}$ 
be also the standard face morphism. 
Then, by the proof of 
\cite[${\rm V}^{\rm bis}$ (5.1.3)]{sga4-2}, 
the face morphism 
\begin{align*} 
\del'{}^l_i \col X'_{l,\os{\circ}{T}_0}=
\coprod_{0\leq k \leq l}
\coprod_{[l]\twoheadrightarrow [k]}
(N_{k,\os{\circ}{T}_0}'\times_{S_{k,\os{\circ}{T}_0}}S_{\os{\circ}{T}_0})
\lo X'_{l-1,\os{\circ}{T}_0}(=\coprod_{0\leq k \leq l-1}
\coprod_{[l-1] \twoheadrightarrow [k]}
(N_{k,\os{\circ}{T}_0}'\times_{S_{k,\os{\circ}{T}_0}}S_{\os{\circ}{T}_0})) 
\end{align*} 
is induced by the following composite morphisms 
for a surjective non-decreasing morphism 
$[l]\twoheadrightarrow [k]$: 
\begin{align*} 
N_{k,\os{\circ}{T}_0}'\times_{S_{k,\os{\circ}{T}_0}}S_{\os{\circ}{T}_0} 
& \os{\sus}{\lo} X'_{k,\os{\circ}{T}_0}
={\rm cosk}_{l-1}^{S_{\os{\circ}{T}_0}}(X'_{\bul \leq l-1,\os{\circ}{T}_0})_k  
\tag{3.2.6.4}\label{eqn:nxkl}\\
& \lo 
{\rm cosk}_{l-1}^{S_{\os{\circ}{T}_0}}
(X'_{\bul \leq l-1,\os{\circ}{T}_0})_l 
\os{\del''^l_i}{\lo} X'_{l-1,\os{\circ}{T}_0}\quad (0 \leq k<l)  
\end{align*} 
\parno 
and   
\begin{equation*} 
N_{l,\os{\circ}{T}_0}'\times_{S_{l,\os{\circ}{T}_0}}S_{\os{\circ}{T}_0} \lo 
{\rm cosk}_{l-1}^{S_{\os{\circ}{T}_0}}(X'_{\bul \leq l-1,\os{\circ}{T}_0})_l 
\os{\del''^l_i}{\lo} X'_{l-1,\os{\circ}{T}_0}.
\tag{3.2.6.5}\label{eqn:nxlkl} 
\end{equation*} 
Because the claim holds for the morphism (\ref{eqn:nxlkl}),  
we have only to consider 
the morphism (\ref{eqn:nxkl}) to prove that the claim holds. 
First consider the case where the composite map  
$[l-1]\lo [l] \twoheadrightarrow [k]$ is surjective. 
Then the morphism 
$X'_{k,\os{\circ}{T}_0}\lo X'_{l-1,\os{\circ}{T}_0}$ obtained by (\ref{eqn:nxkl}) 
is the composite map of degeneracy maps 
(\cite[(8.1.2)]{weib}). 
By the argument in the previous paragraph, we see that 
the morphism  (\ref{eqn:nxkl}) is the composite of inclusion morphisms.  
Next consider the case where the composite map is not surjective. 
Then this map is the composite map 
of a surjective map $[l-1]\lo [k-1]$ and 
the map $\sig_j^k\col [k-1]\lo [k]$ for some $j$ 
$(0\leq j\leq k)$. 
Hence the morphism (\ref{eqn:nxkl}) is the following composite morphism 
\begin{equation*} 
N_{k,\os{\circ}{T}_0}'\times_{S_{k,\os{\circ}{T}_0}}S_{\os{\circ}{T}_0} 
\lo {\rm cosk}_{k-1}^{S_{\os{\circ}{T}_0}}
(X'_{\bul \leq k-1,\os{\circ}{T}_0})_k \os{\del''^k_j}{\lo} X'_{k-1,\os{\circ}{T}_0}
\os{\sus}{\lo} X'_{l-1,\os{\circ}{T}_0},   
\tag{3.2.6.6}\label{eqn:nxxlmn}
\end{equation*}
where the morphism $X'_{k-1,\os{\circ}{T}_0}\os{\sus}{\lo} X'_{l-1,\os{\circ}{T}_0}$ 
is the already mentioned inclusion morphism. 
Now we have proved that the claim holds. 
\par 
Set $\Del'_m :=\prod_{{\rm Hom}_{\Del}([N],[m])}X'_{N,\os{\circ}{T}_0}$ 
$(m\in {\mab N})$ and 
$H_{mn}:={\rm cosk}_0^{S_{m,\os{\circ}{T}_0}}(\Del'_m)_n$. 
Then we have the natural immersions  
$X'_{m,\os{\circ}{T}_0} \os{\sus}{\lo} \Del'_m$ 
and $X_{mn,\os{\circ}{T}_0} \os{\sus}{\lo} H_{mn}$. 
By the argument in the previous paragraph, 
we see that the restriction of 
the immersion $X'_{m,\os{\circ}{T}_0}\os{\sus}{\lo} \Del'_m$ 
to $N'_l\times_{S_{l,\os{\circ}{T}_0}}S_{\os{\circ}{T}_0}$ 
corresponding to the index $[m]\twoheadrightarrow [l]$ is mapped into 
$L_{l,\os{\circ}{T}_0}\times_{S_{l,\os{\circ}{T}_0}}S_{\os{\circ}{T}_0}$, 
where $L_{l,\os{\circ}{T}_0}$ is a log smooth scheme over $S_{l,\os{\circ}{T}_0}$. 
Now we have only to notice that the immersions 
(\ref{ali:xgpmn}) and (\ref{eqn:xgmn}) 
are the following composite immersions, respectively: 
\begin{equation*} 
X'_{m,\os{\circ}{T}_0} \os{\sus}{\lo}  \Del'_m \os{\sus}{\lo} \ol{\cal P}{}'_m, 
\tag{3.2.6.7}\label{eqn:xgdmn} 
\end{equation*} 
\begin{equation*} 
X_{mn,\os{\circ}{T}_0} \os{\sus}{\lo} H_{mn} \os{\sus}{\lo} \ol{\cal P}_{mn}.
\tag{3.2.6.8}\label{eqn:xghmn} 
\end{equation*}
We complete the proof. 
\end{proof}

\parno 
Let ${\cal P}^{\rm ex}_{mn}$ and 
$\ol{\cal P}{}^{\rm ex}_{mn}$ be the exactifications of 
the immersions $X_{mn} \os{\sus}{\lo} {\cal P}_{mn}$ 
and $X_{mn} \os{\sus}{\lo} \ol{\cal P}_{mn}$, respectively. 
Then we have the following $(N,\infty)$-truncated bisimplicial immersions: 
\begin{equation*} 
X_{\bul \leq N,\bul,\os{\circ}{T}_0} \os{\sus}{\lo} 
{\cal P}^{\rm ex}_{\bul \leq N,\bul} 
\os{\sus}{\lo} 
\ol{\cal P}{}^{\rm ex}_{\bul \leq N,\bul}. 
\tag{3.2.6.9}\label{eqn:xxnbmne}
\end{equation*}

\begin{theo}[{\bf Existence of admissible immersions}]\label{theo:thenad} 
The immersions 
$X_{\bul \leq N,\bul,\os{\circ}{T}_0} \os{\sus}{\lo} 
{\cal P}^{\rm ex}_{\bul \leq N,\bul}$ and  
$X_{\bul \leq N,\bul,\os{\circ}{T}_0} 
\os{\sus}{\lo} \ol{\cal P}{}^{\rm ex}_{\bul \leq N,\bul}$ 
are admissible over $S_{\os{\circ}{T}_0}$ and $\ol{S(T)^{\nat}}$ 
augmented to $S_{0,\os{\circ}{T}_0}$ and $\ol{S_0(T)^{\nat}}$, respectively.  
\end{theo} 
\begin{proof} 
Extract $(N'_{l,\os{\circ}{T}_0}\times_{N_{l,\os{\circ}{T}_0}}
\cdots  \times_{N_{l,\os{\circ}{T}_0}}N'_{l,\os{\circ}{T}_0})
\times_{S_{l,\os{\circ}{T}_0}}S_{\os{\circ}{T}_0}$ 
in the source of the immersion (\ref{eqn:xgmn}). 
By (\ref{prop:keyp2}) this is mapped into 
$(\ol{\cal L}{}'_l\times_{\ol{S_l(T)^{\nat}}}\cdots 
\times_{\ol{S_l(T)^{\nat}}}\ol{\cal L}{}'_l)\times_{\ol{S_l(T)^{\nat}}}\ol{S(T)^{\nat}}$, 
where $\ol{\cal L}{}'_l$ is a log smooth scheme over $\ol{S_l(T)^{\nat}}$. 
By (\ref{prop:xpls}) we have 
$((\ol{\cal L}{}'_l\times_{\ol{S_l(T)^{\nat}}} \cdots 
\times_{\ol{S_l(T)^{\nat}}}\ol{\cal L}{}'_l)
\times_{\ol{S_l(T)^{\nat}}}\ol{S(T)^{\nat}})^{\rm ex}
=(\ol{\cal L}{}'_l\times_{\ol{S_l(T)^{\nat}}} \cdots 
\times_{\ol{S_l(T)^{\nat}}}\ol{\cal L}{}'_l)^{\rm ex}
\times_{\ol{S_l(T)^{\nat}}}\ol{S(T)^{\nat}}$. 
By (\ref{prop:nexeo}) we see that $(\ol{\cal L}{}'_l\times_{\ol{S_l(T)^{\nat}}} 
\cdots \times_{\ol{S_l(T)^{\nat}}}\ol{\cal L}{}'_l)^{\rm ex}$
is a (formal) strict semistable log scheme over $\ol{S_l(T)^{\nat}}$.  
\par 
Fix $0\leq m\leq N$. 
Then we have the following decompositions   
according to the decomposition 
$X'_m=
\coprod_{l=0}^m
\coprod_{[m] \twoheadrightarrow [l]}(N'_{l,\os{\circ}{T}_0}
\times_{S_{l,\os{\circ}{T}_0}}S_{\os{\circ}{T}_0})$: 
\begin{equation*} 
{\cal P}'^{\rm ex}_m=\coprod_{l=0}^m
\coprod_{[m] \twoheadrightarrow [l]}
({\cal N}(m,l)\times_{S_l(T)^{\nat}}S(T)^{\nat}) 
\tag{3.2.7.1}\label{eqn:pexm}
\end{equation*}  
and 
\begin{equation*} 
\ol{\cal P}{}'{}^{\rm ex}_{m}=\coprod_{l=0}^m
\coprod_{[m] \twoheadrightarrow [l]}
(\ol{\cal N}(m,l)\times_{\ol{S_l(T)^{\nat}}}\ol{S(T)^{\nat}}), 
\tag{3.2.7.2}
\end{equation*}  
where ${\cal N}(m,l)$ and $\ol{\cal N}(m,l)$ 
are a (formal) SNCL scheme over $S_{l,\os{\circ}{T}_0}$ and 
a (formal) strict semistable log scheme over $\ol{S_l(T)^{\nat}}$, respectively.  
They fit into the following commutative diagram, respectively: 
\begin{equation*} 
\begin{CD}
X'_{m,\os{\circ}{T}_0} @>{\subset}>> {\cal P}'{}^{\rm ex}_{m} 
@>{\subset}>> \ol{\cal P}{}'{}^{\rm ex}_{m}\\
@A{\bigcup}AA @A{\bigcup}AA @AA{\bigcup}A \\
N'_{l,\os{\circ}{T}_0}\times_{S_{l,\os{\circ}{T}_0}}S_{\os{\circ}{T}_0} 
@>{\subset}>> 
{\cal N}(m,l)\times_{S_l(T)^{\nat}}S(T)^{\nat} @>{\subset}>> 
\ol{\cal N}(m,l)\times_{\ol{S_l(T)^{\nat}}}\ol{S(T)^{\nat}}. 
\end{CD}
\tag{3.2.7.3}\label{cd:xpmn} 
\end{equation*}  
Consequently we have the simplicial log scheme  
${\cal N}(m,l)_{\bul}$ over $S(T)^{\nat}$ and $\ol{\cal N}(m,l)_{\bul}$ over $\ol{S(T)^{\nat}}$ 
and the following commutative diagram:   
\begin{equation*} 
\begin{CD}
X_{m\bul,\os{\circ}{T}_0} @>{\subset}>> {\cal P}{}^{\rm ex}_{m\bul} 
@>{\subset}>> \ol{\cal P}{}^{\rm ex}_{m\bul}\\
@A{\bigcup}AA @A{\bigcup}AA @AA{\bigcup}A \\
N_{l\bul,\os{\circ}{T}_0}\times_{S_{l,\os{\circ}{T}_0}}S_{\os{\circ}{T}_0}
@>{\subset}>> 
{\cal N}(m,l)_{\bul}\times_{S_l(T)^{\nat}}S(T)^{\nat} 
@>{\subset}>> \ol{\cal N}(m,l)_{\bul}\times_{\ol{S_l(T)^{\nat}}}\ol{S(T)^{\nat}}. 
\end{CD}
\tag{3.2.7.4}\label{cd:xpmsn} 
\end{equation*}  
This commutative diagram shows the admissibility of 
the immersions 
$X_{\bul \leq N,\bul} \os{\sus}{\lo} 
{\cal P}^{\rm ex}_{\bul \leq N,\bul}$ and  
$X_{\bul \leq N,\bul} \os{\sus}{\lo} \ol{\cal P}{}^{\rm ex}_{\bul \leq N,\bul}$.  
\end{proof}

\par 

Let 
$\{S_{m1}\}_{m=0}^N\lo  \{S_m\}_{m=0}^N$
be a morphism of sequences of families of log points. 
Then we have a morphism $\{S_{m1}\}_{m=0}^N\lo \{S_{m,\os{\circ}{S}_{m1}}\}_{m=0}^N$. 
Assume that the morphism $T_0\lo S_N$ factors through a morphism $T_0\lo S_{N1}$. 
Then we have the following two morphisms  
$$\{S_{m1}(T)\}_{m=0}^N\lo \{S_{m,\os{\circ}{S}_{m1}}(T)\}_{m=0}^N$$ 
and 
$$\{\ol{S_{m1}(T)^{\nat}}\}_{m=0}^N\lo \{\ol{S_{m,\os{\circ}{S}_{m1}}(T)^{\nat}}\}_{m=0}^N.$$ 
Since $N_m$ is an SNCL scheme over $S_m$, 
$N_{m,S_{m,\os{\circ}{S}_{m1},\os{\circ}{T}_0}}:=
N_m\times_{S_m}S_{m,\os{\circ}{S}_{m1},\os{\circ}{T}_0}$ 
is an SNCL scheme over $S_{m,\os{\circ}{S}_{m1},\os{\circ}{T}_0}$. 
Using $\{N_{m,S_{m,\os{\circ}{S}_{m1},\os{\circ}{T}_0}}\}_{m=0}^N$, 
we have a family 
$$\{X(m)_{\bul \leq m,S_{m,\os{\circ}{S}_{m1},\os{\circ}{T}_0}}
/S_{m,\os{\circ}{S}_{m1},\os{\circ}{T}_0}\}_{m=0}^N$$ 
of successive 
split truncated simplicial SNCL schemes 
with respect to the sequence (\ref{eqn:stcnniseq}). 
By the argument in this section, 
we obtain an admissible immersion 
\begin{align*} 
X_{\bul \leq N,\bul,\os{\circ}{T}_0}\os{\sus}{\lo} \ol{\cal Q}_{\bul \leq N,\bul} 
\tag{3.2.7.5}\label{ali:sstq} 
\end{align*}
over $\ol{S_{\os{\circ}{S}_1}(T)^{\nat}}$. 
Consider the fiber product 
$\ol{\cal R}_{\bul \leq N,\bul}:=
\ol{\cal P}_{\bul \leq N,\bul}\times_{\ol{S(T)^{\nat}}}
(\ol{\cal Q}_{\bul \leq N,\bul}\times_{\ol{S_{\os{\circ}{S}_1}(T)^{\nat}}}\ol{S(T)^{\nat}})$. 
Then we have the following commutative diagram 
\begin{equation*} 
\begin{CD} 
X_{\bul \leq N,\bul,\os{\circ}{T}_0}@>>> \ol{\cal R}_{\bul \leq N,\bul}\\
@VVV @VVV \\
X_{\bul \leq N,\bul,S_{\os{\circ}{S}_1},\os{\circ}{T}_0}
@>>> \ol{\cal Q}_{\bul \leq N,\bul}
\end{CD}
\end{equation*}  
over $\ol{S(T)^{\nat}}= \ol{S_{\os{\circ}{S}_1}(T)^{\nat}}$.

\section{Iso-zariskian $p$-adic filtered Steenbrink complexes}\label{sec:pgensc} 
Let $S\lo S_0$ be a morphism of $p$-adic formal PD-families of log points.  
Let $X_{\bul \leq N}/S$ be 
an $N$-truncated simplicial base change of SNCL schemes augmented to $S_0$. 
Let $(T,{\cal J},\del)$ be a fine log $p$-adic formal PD-enlargements of $S$. 
Set $T_0:=\ul{\rm Spec}_T^{\log}({\cal O}_T/{\cal J})$.  
Let the notations be as in (\ref{defi:bca}) (4).  
Let $E^{\bul \leq N}$ be a flat quasi-coherent 
${\cal O}_{\os{\circ}{X}_{\bul \leq N,T_0}/\os{\circ}{T}}$-module 
(without any other assumptions!). 
In this section we generalize results in \S\ref{sec:psc} modulo torsion 
for an admissible immersion of 
the augmented truncated simplicial base change of SNCL schemes 
over a $p$-adic formal PD-family of log points. 
That is, we construct a filtered complex 
$$(A_{{\rm zar},{\mab Q}}
(X_{\bul \leq N,\os{\circ}{T}_0}/S(T)^{\nat},E^{\bul \leq N}),P),$$ 
which is a generalization of  
$(A_{\rm zar}(X_{\bul \leq N,\os{\circ}{T}_0}/S(T)^{\nat},E^{\bul \leq N}),P)
\otimes^L_{\mab Z}{\mab Q}$ 
defined in (\ref{defi:fdirpd}), 
and we calculate 
${\rm gr}_k^P
A_{{\rm zar},{\mab Q}}(X_{\bul \leq N,\os{\circ}{T}_0}/S(T)^{\nat},E^{\bul \leq N})$. 
\par
As already mentioned, assume that there exists an admissible immersion 
\begin{align*} 
\iota \col X_{\bul \leq N,\bul,\os{\circ}{T}_0}\os{\sus}{\lo} \ol{\cal P}_{\bul \leq N,\bul}
\tag{3.3.0.1}\label{alm:adstef} 
\end{align*} 
into a formally log smooth $(N,\infty)$-truncated bisimplicial log formal scheme over $\ol{S(T)^{\nat}}$. 
(By (\ref{theo:thenad}) we have the admissible immersion 
for a family of successive split truncated simplicial SNCL schemes.) 
By the definition of the admissible immersion, 
we have the commutative diagram (\ref{cd:yqst}). 
\par 
For an exact closed point $s$ of $S$, 
let $s_0$ be the image (with the exact log structure) of $s$ in $S_0$ 
by the structural morphism $t\col S \lo S_0$. 
The morphism $t$ induces a map  
$e_s \col {\mab N}\simeq (M_{S_0}/{\cal O}^*_{S_0})_{s_0}
\lo (M_S/{\cal O}^*_S)_s\simeq {\mab N}$ 
and a locally constant morphism  
$e\col t^{-1}(M_{S_0}/{\cal O}^*_{S_0})\lo M_S/{\cal O}^*_S$.  
We can regard $e$ as a locally constant function on $\os{\circ}{S}$ 
and then regard it as a locally constant function on a formal scheme 
over $\os{\circ}{S}$. 
Let $\tau$ be a local section of $M_S$ which gives a local basis of 
$M_S/{\cal O}_S^*$. 
Let $\theta_{{\cal P}^{\rm ex}_{\bul \leq N,\bul}}$ be the pull-back of 
$d\log \tau$ by the composite morphism 
$\os{\circ}{\cal P}{}^{\rm ex}_{\bul \leq N,\bul}\lo \os{\circ}{T}\lo \os{\circ}{S}$. 
As in \S\ref{sec:psc}, 
we define $f_{\bul \leq N}\col X_{\bul \leq N,\os{\circ}{T}_0}\lo S(T)^{\nat}$, 
$f_{m}\col X_{m,\os{\circ}{T}_0}\lo S(T)^{\nat}$ 
$(0\leq m\leq N)$,  
and   
$a_{m\bul,T_0} \col \os{\circ}{X}{}^{(k)}_{m\bul,T_0} 
\lo \os{\circ}{X}_{m\bul,T_0}$ $(0\leq m\leq N,k\in {\mab N})$.  
As in \S\ref{sec:psc} again, we define 
$\ol{\mathfrak D}_{\bul \leq N,\bul}$, ${\mathfrak D}_{\bul \leq N,\bul}$, 
${\mathfrak D}_{\bul \leq N,\bul,T}:=
{\mathfrak D}_{\bul \leq N}\times_{S(T)^{\nat}}T$ 
(when $T$ is restrictively hollow),  
${\cal E}^{\bul \leq N,\bul}$ by using the admissible immersion 
(\ref{alm:adstef}).  
Because ${\Om}^{\bul}_{{\cal P}^{\rm ex}_{\bul \leq N,\bul}/\os{\circ}{T}}$ 
has the filtration $P$ ((\ref{eqn:pkdefpw})), 
we have the induced filtration $P$ on 
${\cal E}^{\bul \leq N,\bul}
\otimes_{{\cal O}_{{\cal P}^{\rm ex}_{\bul \leq N,\bul}}}
{\Om}^{\bul}_{{\cal P}^{\rm ex}_{\bul \leq N,\bul}
/\os{\circ}{T}}\otimes_{\mab Z}{\mab Q}$. 
Set  
\begin{align*} 
A_{{\rm zar},{\mab Q}}({\cal P}^{\rm ex}_{\bul \leq N,\bul}/S(T)^{\nat},
{\cal E}^{\bul \leq N,\bul})^{ij} 
& :=({\cal E}^{\bul \leq N,\bul}
\otimes_{{\cal O}_{{\cal P}^{\rm ex}_{\bul \leq N,\bul}}}
{\Om}^{i+j+1}_{{\cal P}^{\rm ex}_{\bul \leq N,\bul}
/\os{\circ}{T}}\otimes_{\mab Z}{\mab Q})/P_j 
\quad (i,j \in {\mab N}) 
\tag{3.3.0.2}\label{cd:adsef} 
\end{align*}   
as in (\ref{cd:adef}). 
We consider the following boundary morphisms of 
the double complex: 
\begin{equation*}
\begin{CD}
A_{{\rm zar},{\mab Q}}({\cal P}^{\rm ex}_{\bul \leq N,\bul}/S(T)^{\nat},
{\cal E}^{\bul \leq N,\bul})^{i,j+1} @.  \\ 
@A{e\theta_{{\cal P}^{\rm ex}_{\bul \leq N,\bul,T}}\wedge}AA 
@. \\
A_{{\rm zar},{\mab Q}}({\cal P}^{\rm ex}_{\bul \leq N,\bul}/S(T)^{\nat},
{\cal E}^{\bul \leq N,\bul})^{ij} 
@>{-\nabla}>> A_{{\rm zar},{\mab Q}}({\cal P}^{\rm ex}_{\bul \leq N,\bul}/S(T)^{\nat},
{\cal E}^{\bul \leq N,\bul})^{i+1,j}\\
\end{CD}
\tag{3.3.0.3}\label{cd:locstsbd} 
\end{equation*}  
as in (\ref{cd:locstbd}). 
Set 
\begin{equation*} 
A_{{\rm zar},{\mab Q}}({\cal P}^{\rm ex}_{\bul \leq N,\bul}/S(T)^{\nat},
{\cal E}^{\bul \leq N,\bul})
:=
s(A_{{\rm zar},{\mab Q}}({\cal P}^{\rm ex}_{\bul \leq N,\bul}/S(T)^{\nat},
{\cal E}^{\bul \leq N,\bul})^{\bul \bul}). 
\end{equation*} 
As in (\ref{eqn:dblad}), we have the filtrations $P$'s on 
$$A_{{\rm zar},{\mab Q}}({\cal P}^{\rm ex}_{\bul \leq N,\bul}/S(T)^{\nat},
{\cal E}^{\bul \leq N,\bul})^{\bul \bul}$$ 
and 
$$A_{{\rm zar},{\mab Q}}({\cal P}^{\rm ex}_{\bul \leq N,\bul}/S(T)^{\nat},
{\cal E}^{\bul \leq N,\bul}).$$   

\par 
Let 
\begin{equation*} 
\pi \col X_{\bul \leq N,\bul,\os{\circ}{T}_0} \lo X_{\bul \leq N,\os{\circ}{T}_0}
\tag{3.3.0.4}\label{eqn:pznpd} 
\end{equation*} 
and 
\begin{equation*} 
\pi_m \col X_{m\bul,\os{\circ}{T}_0} \lo X_{m,\os{\circ}{T}_0} \quad (0\leq m\leq N)
\tag{3.3.0.5}\label{eqn:pzmxd} 
\end{equation*} 
be the natural morphisms of 
$((N,)\infty)$-truncated (bi)simplicial log schemes over $T_0$. 
Then we have the following morphism of ringed topoi: 
\begin{equation*} 
\pi_{{\rm zar}} \col 
((X_{\bul \leq N,\bul,\os{\circ}{T}_0})_{\rm zar}, 
f^{-1}_{\bul \leq N,\bul}({\cal O}_T)
\otimes_{\mab Z}{\mab Q}) 
\lo 
((X_{\bul \leq N,\os{\circ}{T}_0})_{\rm zar}, 
f^{-1}_{\bul \leq N}({\cal O}_T)\otimes_{\mab Z}{\mab Q}),  
\tag{3.3.0.6}\label{eqn:pzrpd} 
\end{equation*} 
\begin{equation*} 
\pi_{m,{\rm zar}} \col 
((X_{m\bul,\os{\circ}{T}_0})_{\rm zar},
f^{-1}_{m\bul}({\cal O}_{T})\otimes_{\mab Z}{\mab Q}) 
\lo 
((X_{m,\os{\circ}{T}_0})_{\rm zar},
f^{-1}_{m}({\cal O}_T)\otimes_{\mab Z}{\mab Q}).  
\tag{3.3.0.7}\label{eqn:pzrmpd} 
\end{equation*}

\begin{defi}\label{defi:ndfr} 
Let $v\col T'\lo T$ be a morphism of fine log formal schemes. 
We say that $v$ {\it preserves the essential ranks of log structures} 
if, for any point $x$ of $T'$, 
the pull-back morphism 
\begin{align*} 
(M^{\rm gp}_T/{\cal O}_T^*)_{v(x)}\otimes_{\mab Z}{\mab Q}\lo 
(M^{\rm gp}_{T'}/{\cal O}_{T'}^*)_x\otimes_{\mab Z}{\mab Q}
\end{align*} 
is an isomorphism. 
\end{defi}

\begin{prop}\label{prop:flie}
Let the notations be as in {\rm (\ref{defi:ndfr})}.  
Then the natural morphism 
\begin{align*} 
v^*({\Om}^1_{T/\os{\circ}{T}{}})\otimes_{\mab Z}{\mab Q}
\lo 
{\Om}^1_{T'/\os{\circ}{T}{}'}
\otimes_{\mab Z}{\mab Q} 
\tag{3.3.2.1}\label{ali:mtzaq}
\end{align*} 
is an isomorphism. 
\end{prop}
\begin{proof} 
Since 
$$1\otimes u=u^{-1}(u\otimes u)=0 \quad (u\in {\cal O}_T^*)$$  
in 
$\Om^1_{T/\os{\circ}{T}}=
({\cal O}_T\otimes_{\mab Z}M_T^{\rm gp})
/(\alpha(m)\otimes m~\vert~m\in M_T)$, 
we have the following equalities: 
\begin{align*} 
v^{-1}(\Om^1_{T/\os{\circ}{T}})\otimes_{\mab Z}{\mab Q}&=
v^{-1}({\cal O}_T)\otimes_{\mab Z}v^{-1}(M_T^{\rm gp})
\otimes_{\mab Z}{\mab Q}/
(\alpha(m)\otimes  m~\vert~m\in M_T) \\
&=
v^{-1}({\cal O}_T)\otimes_{\mab Z}v^{-1}(M_T^{\rm gp}/{\cal O}_T^*)
\otimes_{\mab Z}{\mab Q}/
(\alpha(m)\otimes m~\vert~m\in v^{-1}(M_T))\\
&=
v^{-1}({\cal O}_T)\otimes_{\mab Z}(M_{T'}^{\rm gp}/{\cal O}_{T'}^*)
\otimes_{\mab Z}{\mab Q}/
(\alpha(m)\otimes m~\vert~m\in M_{T'})\\
&=
v^{-1}({\cal O}_T)\otimes_{\mab Z}M_{T'}^{\rm gp}
\otimes_{\mab Z}{\mab Q}/
(\alpha(m)\otimes m~\vert~m\in M_{T'}). 
\end{align*} 
Hence we obtain the following equalities: 
\begin{align*} 
v^*(\Om^1_{T/\os{\circ}{T}})\otimes_{\mab Z}{\mab Q}
&={\cal O}_{T'}\otimes_{v^{-1}({\cal O}_T)}
v^{-1}(\Om^1_{T/\os{\circ}{T}})\otimes_{\mab Z}{\mab Q} \\
&={\cal O}_{T'}\otimes_{\mab Z}M_{T'}^{\rm gp}\otimes_{\mab Z}{\mab Q}/
(\alpha(m)\otimes m~\vert~m\in M_{T'})\\
&=\Om^1_{T'/\os{\circ}{T}{}'}\otimes_{\mab Z}{\mab Q}. 
\end{align*} 
\end{proof}

\begin{exem}\label{exem:ssm}  
Let $S'\lo S$ be a morphism of 
($p$-adic formal) family of log points. 
Then it is clear that the morphism preserves 
the essential ranks of log structures. 
\end{exem}

Though the proof of the following lemma is easy, 
this is one of key lemmas for the calculation of 
${\rm gr}^P_kA_{{\rm zar},{\mab Q}}({\cal P}^{\rm ex}_{\bul \leq N,\bul}/S(T)^{\nat},
{\cal E}^{\bul \leq N,\bul})$. 

\begin{lemm}\label{lemm:nqgr} 
Let $v\col T'\lo T$ be as in {\rm (\ref{defi:ndfr})}.  
Let $Y$ be a fine log formal scheme over $T$. 
Set $Y':=Y\times_TT'$. 
Let $q\col Y'\lo Y$ be the first projection. 
Let $i$ be a nonnegative integer. 
Assume that $Y$ is formally log smooth over $T$. 
Assume also that 
the natural morphism 
\begin{align*} 
v^*({\Om}^1_{T/\os{\circ}{T}{}})
\otimes_{\mab Z}{\mab Q}\lo 
{\Om}^1_{T'/\os{\circ}{T}{}'}
\otimes_{\mab Z}{\mab Q} 
\tag{3.3.4.1}\label{ali:mtzq}
\end{align*} 
is an isomorphism. 
Then the following hold$:$ 
\par 
$(1)$ The natural morphism 
\begin{align*} 
q^*({\Om}^i_{Y/\os{\circ}{T}{}})\otimes_{\mab Z}{\mab Q}
\lo 
{\Om}^i_{Y'/\os{\circ}{T}{}'}
\otimes_{\mab Z}{\mab Q} 
\tag{3.3.4.2}\label{ali:mtyzq}
\end{align*} 
is an isomorphism. 
\par 
$(2)$ 
Let ${\cal E}$ be a quasi-coherent ${\cal O}_Y$-module. 
Assume that 
${\rm gr}^P_k({\cal E}\otimes_{{\cal O}_Y}
{\Om}^i_{Y/\os{\circ}{T}})$ 
is a flat ${\cal O}_T$-module for any $k\in{\mab Z}$. 
Then there exists the following commutative diagram of exact rows$:$
\begin{equation*} 
\begin{CD} 
0 @>>> {\cal O}_{T'}
\otimes_{{\cal O}_T}P_{k-1}
({\cal E}\otimes_{{\cal O}_{Y}}{\Om}^i_{Y/\os{\circ}{T}})
\otimes_{\mab Z}{\mab Q} 
@>>> {\cal O}_{T'}
\otimes_{{\cal O}_T}P_k
({\cal E}\otimes_{{\cal O}_{Y}}{\Om}^i_{Y/\os{\circ}{T}})
\otimes_{\mab Z}{\mab Q} \\
@. @V{\simeq}VV @V{\simeq}VV \\
0 @>>> P_{k-1}
(q^*({\cal E})
\otimes_{{\cal O}_{Y'}}{\Om}^i_{Y'/\os{\circ}{T}{}'})\otimes_{\mab Z}{\mab Q} 
@>>> P_k
(q^*({\cal E})
\otimes_{{\cal O}_{Y'}}{\Om}^i_{Y'/\os{\circ}{T}{}'})  
\otimes_{\mab Z}{\mab Q}
\end{CD} 
\tag{3.3.4.3}\label{eqn:oppt}
\end{equation*} 
\begin{equation*} 
\begin{CD} 
@>>>{\cal O}_{T'}
\otimes_{{\cal O}_T}{\rm gr}^P_k
({\cal E}\otimes_{{\cal O}_Y}
{\Om}^i_{Y/\os{\circ}{T}})\otimes_{\mab Z}{\mab Q} 
@>>> 0 \\ 
@. @V{\simeq}VV @. \\
@>>> {\rm gr}^P_k
(q^*({\cal E})\otimes_{{\cal O}_{Y'}}
{\Om}^i_{Y'/\os{\circ}{T}{}'})\otimes_{\mab Z}{\mab Q} 
@>>> 0.
\end{CD} 
\end{equation*} 
\end{lemm}
\begin{proof} 
(1): Let $g\col Y\lo T$ and $g'\col Y'\lo T'$ be 
the structural morphisms. 
Consider the following commutative diagram 
\begin{equation*} 
\begin{CD} 
0 @>>>
{\cal O}_{T'}\otimes_{{\cal O}_T}
g^*({\Om}^1_{T/\os{\circ}{T}{}})\otimes_{\mab Z}{\mab Q}
@>>> {\cal O}_{T'}\otimes_{{\cal O}_T}
{\Om}^1_{Y/\os{\circ}{T}{}}\otimes_{\mab Z}{\mab Q}\\
@. @V{\simeq}VV @V{}VV \\
0 @>>> g'{}^*({\Om}^1_{T'/\os{\circ}{T}{}'})
\otimes_{\mab Z}{\mab Q}  
@>>> {\Om}^1_{Y'/\os{\circ}{T}{}'}\otimes_{\mab Z}{\mab Q} 
\end{CD} 
\end{equation*} 
\begin{equation*} 
\begin{CD} 
@>>>{\cal O}_{T'}\otimes_{{\cal O}_T}
{\Om}^1_{Y/T}
\otimes_{\mab Z}{\mab Q}
@>>> 0 \\ 
@. @V{}VV @. \\
@>>> {\Om}^1_{Y'/T'}\otimes_{\mab Z}{\mab Q} 
@>>> 0.
\end{CD} 
\end{equation*} 
Because $Y$ (resp.~$Y'$) is formally log smooth over $T$ (resp.~$T'$), 
the upper horizontal sequence (resp.~the lower horizontal sequence) 
is exact by \cite[(3.12)]{klog1}. 
Because the right vertical arrow is an isomorphism 
by [loc.~cit., (1.7)], so is the middle vertical arrow. 
Taking the $i$-times wedge product of the middle vertical arrow, 
we obtain the isomorphism (\ref{ali:mtyzq}).   
\par 
(2): We have the following natural commutative diagram 
\begin{equation*} 
\begin{CD} 
{\cal O}_{T'}\otimes_{{\cal O}_T}
P_k({\cal E}
\otimes_{{\cal O}_Y}{\Om}^i_{Y/\os{\circ}{T}})\otimes_{\mab Z}{\mab Q}  
@>>> P_k(q^*({\cal E})\otimes_{{\cal O}_{Y'}}
{\Om}^i_{Y'/\os{\circ}{T}{}'})\otimes_{\mab Z}{\mab Q} \\ 
@VVV @VV{\bigcap}V \\ 
{\cal O}_{T'}\otimes_{{\cal O}_T}
{\cal E}\otimes_{{\cal O}_Y}
{\Om}^i_{Y/\os{\circ}{T}}\otimes_{\mab Z}{\mab Q}  
@>>> q^*({\cal E})\otimes_{{\cal O}_{Y'}}
{\Om}^i_{Y'/\os{\circ}{T}{}'}\otimes_{\mab Z}{\mab Q}. 
\end{CD} 
\tag{3.3.4.4}\label{eqn:opt}
\end{equation*} 
By (1) the lower horizontal morphism 
${\cal O}_{T'}\otimes_{{\cal O}_T}{\Om}^i_{Y/\os{\circ}{T}}\otimes_{\mab Z}{\mab Q}  
\lo {\Om}^i_{Y'/\os{\circ}{T}{}'}\otimes_{\mab Z}{\mab Q}$ in (\ref{eqn:opt}) 
for the case ${\cal E}={\cal O}_Y$ 
is an isomorphism and hence so is 
the lower horizontal morphism for any ${\cal E}$.  
It is easy to check that the upper horizontal morphism 
in (\ref{eqn:opt}) is surjective. 
By the assumption and the descending induction on $k$, 
the left vertical morphism is injective. 
Hence the upper horizontal morphism is an isomorphism. 
\end{proof}

\par  
Let $E^{m\bul}$ be 
the flat quasi-coherent ${\cal O}_{\os{\circ}{X}_{m\bul,T_0}/\os{\circ}{T}}$-modules 
obtained by $E^m$. 
Let $E_{i(m)}^{\bul}$ be the restriction of $E^{m\bul}$ to 
$(\os{\circ}{Y}_{i(m)\bul}/\os{\circ}{T})_{\rm crys}$ 
by the isomorphism (\ref{ali:xmbym}): 
\begin{align*} 
X_{m\bul,\os{\circ}{T}_0}\os{\sim}{\lo} 
\coprod_{i(m)=1}^{l(m)}(Y_{i(m)\bul}\times_{S_{i(m),\os{\circ}{T}_0}}S_{\os{\circ}{T}_0}). 
\tag{3.3.4.5}\label{eqn:oxypt}
\end{align*}  
Set $Y_{i(m)\bul,S_{\os{\circ}{T}_0}}:=
Y_{i(m)\bul}\times_{S_{i(m),\os{\circ}{T}_0}}S_{\os{\circ}{T}_0}$ 
and 
$\ol{\cal Q}_{i(m)\bul,\ol{S(T)^{\nat}}}
:=\ol{\cal Q}_{i(m)\bul}\times_{\ol{S_{i(m)}(T)^{\nat}}}\ol{S(T)^{\nat}}$ 
in (\ref{cd:yqst}). 
Let $\ol{\mathfrak D}_{i(m)\bul}$ be the log PD-envelope of 
the immersion $Y_{i(m)\bul,S_{\os{\circ}{T}_0}}\os{\sus}{\lo} \ol{\cal Q}_{i(m)\bul,\ol{S(T)^{\nat}}}$
over $(\os{\circ}{T},{\cal J},\del)$ defined in (\ref{cd:yqst}). 
Set ${\mathfrak D}_{i(m)\bul}:=\ol{\mathfrak D}_{i(m)\bul}
\times_{{\mathfrak D}(\ol{S(T)^{\nat}})}S(T)^{\nat}$. 
Let $(\ol{\cal E}{}^{\bul}_{i(m)},\ol{\nabla})$ 
be the quasi-coherent 
${\cal O}_{\ol{\mathfrak D}_{i(m)\bul}}$-module 
with integrable connection corresponding to 
$\eps^*_{Y_{i(m)\bul,S_{\os{\circ}{T}_0}}/\os{\circ}{T}}(E_{i(m)}^{\bul})$. 
Set $({\cal E}^{\bul}_{i(m)},\nabla):=
(\ol{\cal E}{}^{\bul}_{i(m)},\ol{\nabla})
\otimes_{{\cal O}_{\ol{\mathfrak D}_{i(m)\bul}}}
{\cal O}_{{\mathfrak D}_{i(m)\bul}}$. 
Let $u_{i(m)}\col S_{i(m)}\lo S_0$ be the structural morphism.
Set $d_{i(m)}:=\deg(u_{i(m)})$. 
\par 
By the functoriality of the Poincar\'{e} residue isomorphism ((\ref{prop:rescos})), 
we obtain the following isomorphism: 
\begin{align*} 
{\rm Res}\otimes{\rm id}_{\mab Q}\col {\rm gr}_{k}^P({\cal E}^{\bul}_{i(m)}
\otimes_{{\cal O}_{{\cal Q}^{\rm ex}_{i(m)\bul}}}
{\Om}^{\bul}_{{\cal Q}^{\rm ex}_{i(m)\bul}/\os{\circ}{T}}\otimes_{\mab Z}{\mab Q})
&\os{\sim}{\lo} 
{\cal E}^{\bul}_{i(m)}
\otimes_{{\cal O}_{{\cal Q}^{\rm ex}_{i(m)\bul}}}
{\Om}^{\bul -k}_{{\cal Q}^{{\rm ex},(k-1)}_{i(m)\bul}/\os{\circ}{T}}
\otimes_{\mab Z}{\mab Q}. \tag{3.3.4.6}\label{ali:opst}
\end{align*} 
Hence, by (\ref{lemm:nqgr}) we have the following composite isomorphism 
\begin{align*} 
{\rm Res}_{\mab Q}\col & 
{\rm gr}_k^P({\cal E}^{\bul}_{i(m)}
\otimes_{{\cal O}_{{\cal Q}^{\rm ex}_{i(m)\bul,S(T)^{\nat}}}}
\Om^{\bul}_{{\cal Q}^{\rm ex}_{i(m)\bul,S(T)^{\nat}}/\os{\circ}{T}}
\otimes_{\mab Z}{\mab Q})\tag{3.3.4.7}\label{ali:opisst}\\ 
&={\rm gr}_k^P({\cal E}^{\bul}_{i(m)}
\otimes_{{\cal O}_{{\cal Q}^{\rm ex}_{i(m)\bul},S(T)^{\nat}}}
{\cal O}_{S(T)^{\nat}}\otimes_{{\cal O}_{S_{i(m)}(T)^{\nat}}}
{\Om}^{\bul}_{{\cal Q}^{\rm ex}_{i(m)\bul}/\os{\circ}{T}}
\otimes_{\mab Z}{\mab Q})\\
&={\cal E}^{\bul}_{i(m)}
\otimes_{{\cal O}_{{\cal Q}^{\rm ex}_{i(m)\bul}}}
{\rm gr}_k^P(\Om^{\bul}_{{\cal Q}^{\rm ex}_{i(m)\bul}/\os{\circ}{T}}
\otimes_{\mab Z}{\mab Q})\\
&\os{\sim}{\lo} 
{\cal E}^{\bul}_{i(m)}\otimes_{{\cal O}_{{\cal Q}^{\rm ex}_{i(m)\bul}}}
\Om^{\bul -k}_{{\cal Q}^{{\rm ex},(k-1)}_{i(m)\bul}/\os{\circ}{T}}
\otimes_{\mab Z}{\mab Q}. 
\end{align*} 

\begin{coro}\label{coro:grsem} 
Set ${\cal K}_T:={\cal O}_T\otimes_{\mab Z}{\mab Q}$. 
Then the complexes 
${\cal O}_{{\mathfrak D}_{mn}}\otimes_{{\cal O}_{{\cal P}^{\rm ex}_{mn}}}
{\Om}^{\bul}_{{\cal P}^{\rm ex}_{mn}/\os{\circ}{T}}\otimes_{\mab Z}{\mab Q}$,  
${\rm gr}^P_k
({\cal O}_{{\mathfrak D}_{mn}}\otimes_{{\cal O}_{{\cal P}^{\rm ex}_{mn}}}
{\Om}^{\bul}_{{\cal P}^{\rm ex}_{mn}/\os{\circ}{T}})\otimes_{\mab Z}{\mab Q}$ 
$(k\in {\mab N})$ 
and 
$P_k
({\cal O}_{{\mathfrak D}_{mn}}
\otimes_{{\cal O}_{{\cal P}^{\rm ex}_{mn}}}
{\Om}^{\bul}_{{\cal P}^{\rm ex}_{mn}/\os{\circ}{T}})\otimes_{\mab Z}{\mab Q}$ 
consist of locally free ${\cal K}_T$-modules.   
\end{coro}
\begin{proof} 
This immediately follows from (\ref{ali:opisst}) and (\ref{coro:flt}). 
\end{proof}

\begin{theo}\label{theo:tesufc}
There exists the following isomorphism 
\begin{align*} 
&\theta:=\theta_{X_{\bul \leq N,\os{\circ}{T}_0}/S(T)^{\nat}/S_0(T)^{\nat}}:=
eR\pi_{{\rm zar}*}(\theta_{{\cal P}^{\rm ex}_{\bul \leq N,\bul}})\wedge 
\col \tag{3.3.6.1}\label{ali:usz} \\
&
Ru_{X_{\bul \leq N,\os{\circ}{T}_0}/S(T)^{\nat}*}
(\eps^*_{X_{\bul \leq N,T_0}/T}(E^{\bul \leq N}))
\otimes_{\mab Z}{\mab Q}
\os{\sim}{\lo} 
R\pi_{{\rm zar}*}
(A_{{\rm zar},{\mab Q}}({\cal P}^{\rm ex}_{\bul \leq N,\bul}/S(T)^{\nat},
{\cal E}^{\bul \leq N,\bul})). 
\end{align*} 
This isomorphism is independent of the choice of 
an affine $N$-truncated simplicial open covering of 
$X_{\bul \leq N,\os{\circ}{T}_0}$ 
and the choice of an admissible immersion 
$X_{\bul \leq N,\bul,\os{\circ}{T}_0} \os{\sus}{\lo} 
\ol{\cal P}_{\bul \leq N,\bul}$ over $\ol{S(T)^{\nat}}$.
In particular, the complex 
$R\pi_{{\rm zar}*}
(A_{{\rm zar},{\mab Q}}({\cal P}^{\rm ex}_{\bul \leq N,\bul}/S(T)^{\nat},
{\cal E}^{\bul \leq N,\bul}))$ 
is independent of the choices above. 
\end{theo}
\begin{proof}
Because $\iota \col X_{m\bul}\os{\sus}{\lo} \ol{\cal P}_{m\bul}$ 
is admissible, there exists the commutative diagram (\ref{cd:yqst}). 
By considering each component of the lower immersion in (\ref{cd:yqst}) 
and by using (\ref{ali:opisst}), 
the same argument in the proof as that of (\ref{prop:tefc}) shows that 
(\ref{ali:usz}) is an isomorphism. By using (\ref{prop:sxst}), 
the same argument above also shows the independence. 
\end{proof}

\begin{lemm}\label{lemm:axp}  
Let $E^m_{\os{\circ}{X}{}^{(2j+k)}_{m,T_0}/\os{\circ}{T}}$ be the 
${\cal O}_{\os{\circ}{X}{}^{(2j+k)}_{m,T_0}/\os{\circ}{T}}$-module 
obtained by $E^m$ $(0\leq m\leq N)$. 
Then there exists the following isomorphism$:$
\begin{align*} 
& {\rm gr}^P_kR\pi_{m,{\rm zar}*}
(A_{{\rm zar},{\mab Q}}({\cal P}^{\rm ex}_{m\bul}/S(T)^{\nat},{\cal E}^{m\bul}))
\os{\sim}{\lo}\tag{3.3.7.1}\label{ali:xnsnp}\\
&\bigoplus_{j\geq \max \{-k,0\}} 
a^{(2j+k)}_{m*} 
Ru_{\os{\circ}{X}{}^{(2j+k)}_{m,T_0}
/\os{\circ}{T}*}
(E_{\os{\circ}{X}{}^{(2j+k)}_{m,T_0}/\os{\circ}{T}}
\otimes_{\mab Z}
\vp_{\rm crys}^{(2j+k)}(\os{\circ}{X}_{m,T_0}/\os{\circ}{T}))
\otimes_{\mab Z}{\mab Q})[-2j-k]. 
\end{align*}  
\end{lemm} 
\begin{proof}  
Let $a^{(k)}_{\bul}\col 
\os{\circ}{X}{}^{(k)}_{m\bul,T_0}\lo \os{\circ}{X}_{m\bul,T_0}$ be 
the natural morphism of schemes. 
Let $g_{i(m)}\col Y_{i(m),\os{\circ}{T}_0}\lo S_{i(m)}(T)^{\nat}$ 
and $g_{i(m)\bul}\col Y_{i(m)\bul,\os{\circ}{T}_0}\lo S_{i(m)}(T)^{\nat}$ 
be the structural morphisms. 
Let 
\begin{equation*} 
\pi_{i(m),{\rm zar}} \col 
((Y_{i(m)\bul,\os{\circ}{T}_0})_{\rm zar},g^{-1}_{i(m)\bul}({\cal O}_T)\otimes_{\mab Z}{\mab Q}) 
\lo 
((Y_{i(m),\os{\circ}{T}_0})_{\rm zar},
g^{-1}_{i(m)}({\cal O}_T)\otimes_{\mab Z}{\mab Q})  
\tag{3.3.7.2}\label{eqn:pzrjpd} 
\end{equation*} 
be the natural morphisms of ringed topoi. 
By (\ref{ali:opisst}) we obtain the following isomorphism: 
\begin{align*} 
& {\rm gr}^P_kR\pi_{m,{\rm zar}*}
A_{{\rm zar},{\mab Q}}({\cal P}^{\rm ex}_{m\bul}/S(T)^{\nat},{\cal E}^{m\bul})
\tag{3.3.7.3}\label{ali:omnp}\\
& \os{\sim}{\lo} R\pi_{m,{\rm zar}*}
{\rm gr}^P_k
A_{{\rm zar},{\mab Q}}({\cal P}^{\rm ex}_{m\bul}/S(T)^{\nat},{\cal E}^{m\bul})
\\
& \os{\sim}{\lo} \bigoplus_{j\geq \max \{-k,0\}}
R\pi_{m,{\rm zar}*} (\cdots \lo {\rm gr}_{2j+k+1}^P
({\cal E}^{m\bul}
\otimes_{{\cal O}_{{\cal P}^{\rm ex}_{m\bul}}}
{\Om}^{\bul +j+1}_{{\cal P}^{\rm ex}_{m\bul}/\os{\circ}{T}}\otimes_{\mab Z}{\mab Q})\{-j\},-\nabla)  \\
&~~~~~~~\lo\cdots) \\ 
& \os{\sim}{\lo} 
\bigoplus_{i(m)=1}^{l(m)}\bigoplus_{j\geq \max \{-k,0\}}
(R\pi_{i(m),{\rm zar}*}(\cdots \lo {\rm gr}_{2j+k+1}^P
({\cal E}^{\bul}_{i(m)}
\otimes_{{\cal O}_{{\cal Q}^{\rm ex}_{i(m)\bul,S(T)^{\nat}}}}
\Om^{\bul +j+1}_{{\cal Q}^{\rm ex}_{i(m)\bul,S(T)^{\nat}}/\os{\circ}{T}}\\
&~~~~~~~\otimes_{\mab Z}{\mab Q})\{-j\},-\nabla)\lo \cdots) \\ 
&\os{{\rm (\ref{ali:opisst})},\simeq}{\lo}
\bigoplus_{i(m)=1}^{l(m)}\bigoplus_{j\geq \max \{-k,0\}}(R\pi_{i(m),{\rm zar}*}
(\cdots \lo ({\cal E}^{\bul}_{i(m)}\otimes_{{\cal O}_{{\cal Q}^{\rm ex}_{i(m)\bul}}}
\Om^{\bul}_{{\cal Q}^{{\rm ex},(2j+k)}_{i(m)\bul}/\os{\circ}{T}}\\
&~~~~~~~\otimes_{\mab Z}{\mab Q}[-2j-k]\lo \cdots)   \os{\sim}{\lo} 
\\
&
\bigoplus_{j\geq \max \{-k,0\}}
R\pi_{m,{\rm zar}*}
a^{(2j+k)}_{m\bul *}
(Ru_{\os{\circ}{X}{}^{(2j+k)}_{m\bul,\os{\circ}{T}_0/\os{\circ}{T}}*}
(E_{\os{\circ}{X}{}^{(2j+k)}_{m\bul,\os{\circ}{T}_0/\os{\circ}{T}}}\otimes_{\mab Z} 
\vp_{\rm crys}^{(2j+k)}(\os{\circ}{X}_{m\bul,\os{\circ}{T}}/\os{\circ}{T}))
\otimes_{\mab Z}{\mab Q})[-2j-k]\\
& =\bigoplus_{j\geq \max \{-k,0\}} 
a^{(2j+k)}_{m*} 
(Ru_{\os{\circ}{X}{}^{(2j+k)}_{m,T_0}/\os{\circ}{T}*}
(E_{\os{\circ}{X}{}^{(2j+k)}_{m,T_0}/\os{\circ}{T}}
\otimes_{\mab Z}\vp_{\rm crys}^{(2j+k)}(\os{\circ}{X}_{m,T_0}/\os{\circ}{T}))
\otimes_{\mab Z}{\mab Q})[-2j-k].  
\end{align*} 
Here we have used \cite[(1.3.4)]{nh2} 
to obtain the first isomorphism.
\end{proof}

\begin{theo}\label{theo:ingendcr} 
The filtered complex 
$R\pi_{{\rm zar}*}
((A_{{\rm zar},{\mab Q}}({\cal P}^{\rm ex}_{\bul \leq N,\bul}/S(T)^{\nat},
{\cal E}^{\bul \leq N,\bul}),P))$ 
is independent of the choice of 
the disjoint union $X_{\bul \leq N,\os{\circ}{T}_0}'$ 
and the choice of an immersion 
$X_{\bul \leq N,\bul,\os{\circ}{T}_0} \os{\sus}{\lo} \ol{\cal P}_{\bul \leq N,\bul}$ 
over $\ol{S(T)^{\nat}}$. 
\end{theo}
\begin{proof} 
Assume that we are given two admissible immersions  
\begin{equation*} 
X_{\bul \leq N,\bul,\os{\circ}{T}_0} \os{\sus}{\lo} \ol{\cal P}_{\bul \leq N,\bul}, \quad 
X'_{\bul \leq N,\bul,\os{\circ}{T}_0} \os{\sus}{\lo} \ol{\cal Q}_{\bul \leq N,\bul} 
\tag{3.3.8.1}\label{eqn:twosn}
\end{equation*}
into $(N,\infty)$-truncated bisimplicial log schemes over $\ol{S(T)^{\nat}}$.  
By (\ref{prop:sxst}) we may assume that 
there exists the following commutative diagram  
over $\ol{S(T)^{\nat}}$: 
\begin{equation*} 
\begin{CD} 
X_{\bul \leq N,\bul,\os{\circ}{T}_0}
@>{\sus}>>\ol{\cal P}_{\bul \leq N,\bul} \\
@VVV @VVV \\
X'_{\bul \leq N,\bul,\os{\circ}{T}_0} @>{\sus}>> \ol{\cal Q}_{\bul \leq N,\bul}.  
\end{CD} 
\tag{3.3.8.2}\label{cd:pcspx}
\end{equation*}
The commutative diagram 
(\ref{cd:pcspx})
gives us the following natural morphism 
\begin{equation*} 
R\pi_{{\rm zar}*}
(A_{{\rm zar},{\mab Q}}({\cal Q}^{\rm ex}_{\bul \leq N,\bul}/S(T)^{\nat},
{\cal E}^{\bul \leq N,\bul}),P)
\lo 
R\pi_{{\rm zar}*}
((A_{{\rm zar},{\mab Q}}({\cal P}^{\rm ex}_{\bul \leq N,\bul}/S(T)^{\nat},
{\cal E}^{\bul \leq N,\bul}),P)). 
\tag{3.3.8.3}\label{cd:pcsfpx} 
\end{equation*} 
To prove that the morphism (\ref{cd:pcsfpx}) is 
a filtered isomorphism, it suffices to prove that 
the following morphism 
\begin{equation*} 
R\pi_{m,{\rm zar}*}
((A_{{\rm zar},{\mab Q}}({\cal Q}^{\rm ex}_{m\bul}/S(T)^{\nat},{\cal E}^{m\bul})),P)) 
\lo 
R\pi_{m,{\rm zar}*}
((A_{{\rm zar},{\mab Q}}({\cal P}^{\rm ex}_{m\bul}/S(T)^{\nat},{\cal E}^{m\bul}),P))  
\tag{3.3.8.4}\label{cd:pmqp} 
\end{equation*} 
is a filtered isomorphism. 
By the isomorphism (\ref{ali:xnsnp}), 
this morphism is a filtered isomorphism 
because we have the following equality 
for $0\leq m\leq N$: 
\begin{align*} 
{\rm gr}^P_kR\pi_{m,T,{\rm zar}*}
(A_{{\rm zar},{\mab Q}}({\cal Q}^{\rm ex}_{m\bul}/S(T)^{\nat},{\cal E}^{m\bul})) 
={\rm gr}^P_kR\pi_{m,T,{\rm zar}*}
(A_{{\rm zar},{\mab Q}}({\cal P}^{\rm ex}_{m\bul}/S(T)^{\nat},{\cal E}^{m\bul})). 
\end{align*}  
\par 
We finish the proof of (\ref{theo:ingendcr}).  
\end{proof}

\begin{defi}\label{defi:wcthd}
We call the filtered direct image 
$$R\pi_{{\rm zar}*}
((A_{{\rm zar},{\mab Q}}({\cal P}^{\rm ex}_{\bul \leq N,\bul}/S(T)^{\nat},
{\cal E}^{\bul \leq N,\bul}),P))$$ 
the {\it iso-zariskian $p$-adic filtered  Steenbrink complex} of $E^{\bul \leq N}$ 
with respect to 
the morphism $(T,{\cal J},\del)\lo (S,{\cal I},\gam)$.   
We denote it by 
\begin{align*}  
(A_{{\rm zar},{\mab Q}}(X_{\bul \leq N,\os{\circ}{T}_0}/S(T)^{\nat}, 
E^{\bul \leq N}),P).
\tag{3.3.9.1}\label{ali:aqxnb}
\end{align*}   
When we have to emphasize the morphism $T_0\lo S$, 
we denote it  by 
\begin{align*} 
(A_{{\rm zar},{\mab Q}}((X_{\bul \leq N}/S(T)^{\nat})_{T_0\lo S},E^{\bul \leq N}),P).
\tag{3.3.9.2}\label{ali:aqexnb}
\end{align*}  
When $E^{\bul \leq N}
={\cal O}_{\os{\circ}{X}_{\bul \leq N,T_0}/\os{\circ}{T}_{\bul \leq N}}$, 
we call $(A_{{\rm zar},{\mab Q}}(X_{\bul \leq N,\os{\circ}{T}_0}/S(T)^{\nat}, 
E^{\bul \leq N}),P)$  
the {\it iso-zariskian $p$-adic filtered  Steenbrink complex} of 
$X_{\bul \leq N,\os{\circ}{T}_0}/(S(T)^{\nat},{\cal J},\del)$. 
We denote it by  
$(A_{{\rm zar},{\mab Q}}(X_{\bul \leq N,\os{\circ}{T}_0}/S(T)^{\nat}),P)$. 
\end{defi}

\begin{coro}\label{coro:oonzu}
Assume that $T$ is restrictively hollow with respective to the morphism $T_0\lo S$. 
Then there exists the following isomorphism 
\begin{align*} 
\theta:=\theta_{X_{\bul \leq N,T_0}/T}\col Ru_{X_{\bul \leq N,T_0}/T*}
(\eps^*_{X_{\bul \leq N,T_0/T}}(E^{\bul \leq N})) 
\os{\sim}{\lo}
A_{{\rm zar},{\mab Q}}(X_{\bul \leq N,\os{\circ}{T}_0}/S(T),E^{\bul \leq N})
\tag{3.3.10.1}\label{eqn:uqntz} 
\end{align*} 
in ${\rm D}^+(f^{-1}_{\bul \leq N}({\cal O}_T))$.  
\end{coro} 
\begin{proof} 
The proof is the same as that of (\ref{coro:oozu}). 
\end{proof}

The following is nothing but (\ref{ali:xnsnp}). 

\begin{prop}\label{prop:wht}
\begin{align*} 
{\rm gr}^P_kA_{{\rm zar},{\mab Q}}(X_{m,\os{\circ}{T}_0}/S(T)^{\nat},E^{\bul \leq N})
= \bigoplus_{j\geq \max \{-k,0\}} &
(a^{(2j+k)}_{m,T_0*} 
Ru_{\os{\circ}{X}{}^{(2j+k)}_{m,T_0}
/\os{\circ}{T}*}
(E_{\os{\circ}{X}{}^{(2j+k)}_{m,T_0}/\os{\circ}{T}}\otimes_{\mab Z}
\tag{3.3.11.1}\label{ali:jkgr}\\
& 
\vp_{\rm crys}^{(2j+k)}(\os{\circ}{X}_{m,T_0}/\os{\circ}{T}))
\otimes_{\mab Z}{\mab Q}, (-1)^{j+1}\nabla)[-2j-k]. 
\end{align*} 
\end{prop}

\begin{rema} 
By using the arguments in this section, 
we can simplify the proof of du Bois'  result in (\cite{db}). 
Especially we do not need the semistable reduction theorem 
and resolution of singularities; 
de Jong's semistable reduction theorem is enough. 
(We leave the detail to the reader.) 
\end{rema}


\begin{lemm}\label{lemm:cop}
Let $Y_{\bul \leq N}$ be an SNCL scheme over $S_0$. 
Set $X_{\bul \leq N}:=Y_{\bul \leq N}\times_{S_0}S$.  
Then there exists a canonical isomorphism 
\begin{align*} 
Ru_{Y_{\bul \leq N,\os{\circ}{T}_0}/S_0(T)^{\nat}*}
(\eps^*_{Y_{\bul \leq N,\os{\circ}{T}_0}/S_0(T)^{\nat}}(E^{\bul \leq N}))
\otimes^L_{\mab Z}{\mab Q}
\os{\sim}{\lo} 
Ru_{X_{\bul \leq N,\os{\circ}{T}_0}/S(T)^{\nat}*}
(\eps^*_{X_{\bul \leq N,\os{\circ}{T}_0}/S(T)^{\nat}}
(E^{\bul \leq N}))\otimes^L_{\mab Z}{\mab Q}.
\tag{3.3.13.1}\label{ali:auzrxt} 
\end{align*} 
\end{lemm}
\begin{proof} 
We may assume that $\bul \leq N$ is equal to $m$ for some $0\leq m\leq N$. 
In this case (\ref{lemm:cop}) 
immediately follows from the log Poincar\'{e} lemma \cite[(6.4)]{klog1}. 
\end{proof}

\begin{theo}[{\bf Comparison theorem}]\label{theo:ctbfd}
Let $Y_{\bul \leq N}$ be an SNCL scheme over $S_0$. 
Let $S\lo S_0$ be a morphism of family of log points  
and set $X_{\bul \leq N}:=Y_{\bul \leq N}\times_{S_0}S$.  
Assume that $Y_{\bul \leq N,\os{\circ}{T}_0}$ has the disjoint union 
of an affine simplicial open covering of $Y_{\bul \leq N,\os{\circ}{T}_0}$. 
Then there exists a canonical filtered isomorphism 
\begin{align*} 
(A_{\rm zar}(Y_{\bul \leq N,\os{\circ}{T}_0}/S_0(T)^{\nat},E^{\bul \leq N}),P)
\otimes^L_{\mab Z}{\mab Q}
\os{\sim}{\lo} 
(A_{{\rm zar},{\mab Q}}(X_{\bul \leq N,\os{\circ}{T}_0}/S(T)^{\nat},E^{\bul \leq N}),P)
\tag{3.3.14.1}\label{ali:azrxt} 
\end{align*} 
fitting into the following commutative diagram 
\begin{equation*} 
\begin{CD} 
A_{\rm zar}(Y_{\bul \leq N,\os{\circ}{T}_0}/S_0(T)^{\nat},E^{\bul \leq N})
\otimes^L_{\mab Z}{\mab Q}
@>{\sim}>>  \\
@A{\theta_{Y_{\bul \leq N,\os{\circ}{T}_0}/S_0(T)^{\nat}}}A{\simeq}A\\
Ru_{Y_{\bul \leq N,\os{\circ}{T}_0}/S_0(T)^{\nat}*}
(\eps^*_{Y_{\bul \leq N,\os{\circ}{T}_0}/S_0(T)^{\nat}}(E^{\bul \leq N}))
\otimes^L_{\mab Z}{\mab Q}
@= 
\end{CD} 
\tag{3.3.14.2}\label{eqqn:azerxt} 
\end{equation*} 
\begin{equation*} 
\begin{CD} 
A_{{\rm zar},{\mab Q}}(X_{\bul \leq N,\os{\circ}{T}_0}/S(T)^{\nat},E^{\bul \leq N})\\
@A{\theta_{X_{\bul \leq N,\os{\circ}{T}_0}/S(T)^{\nat}/S_0(T)^{\nat}}}A{\simeq}A\\
Ru_{X_{\bul \leq N,\os{\circ}{T}_0}/S(T)^{\nat}*}
(\eps^*_{X_{\bul \leq N,\os{\circ}{T}_0}/S(T)^{\nat}}
(E^{\bul \leq N}))\otimes^L_{\mab Z}{\mab Q}.
\end{CD}  
\end{equation*}
\end{theo}
\begin{proof} 
Let $Y_{\bul \leq N,\bul,\os{\circ}{T}_0}$ be as in \S\ref{sec:tdiai} and 
let $Y_{\bul \leq N,\bul,\os{\circ}{T}_0}\os{\sus}{\lo}\ol{\cal P}_{\bul \leq N,\bul}$ 
be an immersion into a formally log smooth $(N,\infty)$-truncated bisimplicial 
log formal scheme over $\ol{S_0(T)^{\nat}}$. 
Then the immersion 
$Y_{\bul \leq N,\bul,\os{\circ}{T}_0}\times_{S_{0,\os{\circ}{T}_0}}S_{\os{\circ}{T}_0}
\os{\sus}{\lo}\ol{\cal P}_{\bul \leq N,\bul}\times_{\ol{S_0(T)^{\nat}}}\ol{S(T)^{\nat}}$ 
is admissible. 
Let $\ol{\mathfrak D}_{\bul \leq N,\bul}$ be the log PD-envelope of the immersion 
$Y_{\bul \leq N,\bul,\os{\circ}{T}_0}\times_{S_{0,\os{\circ}{T}_0}}S_{\os{\circ}{T}_0}
\os{\sus}{\lo}\ol{\cal P}_{\bul \leq N,\bul}$ over $(\os{\circ}{T},{\cal J},\del)$. 
Then $\ol{\mathfrak D}_{\bul \leq N,\bul}\times_{{\mathfrak D}(\ol{S_0(T)^{\nat}})}
{\mathfrak D}(\ol{S(T)^{\nat}})$ is the log PD-envelope of 
this admissible immersion above over $(\os{\circ}{T},{\cal J},\del)$. 
(Note that the underlying morphism of schemes of the morphism
${\mathfrak D}(\ol{S_0(T)^{\nat}})\lo {\mathfrak D}(\ol{S(T)^{\nat}})$ 
is the identity.)
Set ${\mathfrak D}_{\bul \leq N,\bul}:=\ol{\mathfrak D}_{\bul \leq N,\bul}
\times_{{\mathfrak D}(\ol{S_0(T)^{\nat}})}S_0(T)^{\nat}$.
Then 
\begin{align*} 
\ol{\mathfrak D}_{\bul \leq N,\bul}\times_{{\mathfrak D}(\ol{S_0(T)^{\nat}})}
{\mathfrak D}(\ol{S(T)^{\nat}})\times_{{\mathfrak D}(\ol{S(T)^{\nat}})}S(T)^{\nat}
&=
\ol{\mathfrak D}_{\bul \leq N,\bul}\times_{{\mathfrak D}(\ol{S_0(T)^{\nat}})}S(T)^{\nat}\\
&={\mathfrak D}_{\bul \leq N,\bul}\times_{S_0(T)^{\nat}}S(T)^{\nat}. 
\end{align*} 
Hence the underlying formal scheme of 
$\ol{\mathfrak D}_{\bul \leq N,\bul}\times_{{\mathfrak D}(\ol{S_0(T)^{\nat}})}
{\mathfrak D}(\ol{S(T)^{\nat}})\times_{{\mathfrak D}(\ol{S(T)^{\nat}})}S(T)^{\nat}$ 
is equal to the underlying formal scheme of ${\mathfrak D}_{\bul \leq N,\bul}$. 
In particular, 
\begin{align*} 
{\cal O}_{\ol{\mathfrak D}_{\bul \leq N,\bul}\times_{{\mathfrak D}(\ol{S_0(T)^{\nat}})}
{\mathfrak D}(\ol{S(T)^{\nat}})\times_{{\mathfrak D}(\ol{S(T)^{\nat}})}S(T)^{\nat}}
={\cal O}_{{\mathfrak D}_{\bul \leq N,\bul}}.
\end{align*}  
Let $({\cal E}^{\bul \leq N,\bul},\nabla)$ be the 
${\cal O}_{{\mathfrak D}_{\bul \leq N,\bul}}$-module with integrable connection 
obtained by $E^{\bul \leq N}$. 
This is equal to the 
${\cal O}_{\ol{\mathfrak D}_{\bul \leq N,\bul}\times_{{\mathfrak D}(\ol{S_0(T)^{\nat}})}
{\mathfrak D}(\ol{S(T)^{\nat}})\times_{{\mathfrak D}(\ol{S(T)^{\nat}})}S(T)^{\nat}}$-module  
with integrable connection 
obtained by $E^{\bul \leq N}$. 
Then, by (\ref{exem:ssm}) and (\ref{lemm:nqgr}),  we have the following equalities: 
\begin{align*} 
&P_kA_{{\rm zar},{\mab Q}}({\cal P}^{\rm ex}_{\bul \leq N,\bul}/S(T)^{\nat},
{\cal E}^{\bul \leq N,\bul})^{ij}\tag{3.3.14.3}\label{ali:idtt}\\
&=P_{2j+k+1}({\cal E}^{\bul \leq N,\bul}
\otimes_{{\cal O}_{{\cal P}^{\rm ex}_{\bul \leq N,\bul,S(T)^{\nat}}}}
{\Om}^{i+j+1}_{{\cal P}^{\rm ex}_{\bul \leq N,\bul,S(T)^{\nat}}
/\os{\circ}{T}}\otimes_{\mab Z}{\mab Q})+P_j/P_j \\
&=
P_{2j+k+1}({\cal E}^{\bul \leq N,\bul}
\otimes_{{\cal O}_{{\cal P}^{\rm ex}_{\bul \leq N,\bul}}}
{\Om}^{i+j+1}_{{\cal P}^{\rm ex}_{\bul \leq N,\bul}
/\os{\circ}{T}}\otimes_{\mab Z}{\mab Q})+P_j/P_j\\
&=
P_kA_{\rm zar}({\cal P}^{\rm ex}_{\bul \leq N,\bul}/S(T)^{\nat},
{\cal E}^{\bul \leq N,\bul})^{ij}\otimes_{\mab Z}{\mab Q}
\end{align*}   
for $k\in {\mab Z}$ or $k=\infty$. 
We also have the following commutative diagram: 
\begin{equation*} 
\begin{CD} 
A_{{\rm zar},{\mab Q}}({\cal P}^{\rm ex}_{\bul \leq N,\bul}/S(T)^{\nat},
{\cal E}^{\bul \leq N,\bul})^{i,j+1}
@= A_{\rm zar}({\cal P}^{\rm ex}_{\bul \leq N,\bul}/S_0(T)^{\nat},
{\cal E}^{\bul \leq N,\bul})^{i,j+1}\otimes_{\mab Z}{\mab Q}\\
@A{e\theta_{{\cal P}^{\rm ex}_{\bul \leq N,\bul}}\wedge}A{\simeq}A 
@A{\simeq}A{\theta_{{\cal P}^{\rm ex}_{\bul \leq N,\bul}}\wedge}A \\
A_{{\rm zar},{\mab Q}}({\cal P}^{\rm ex}_{\bul \leq N,\bul}/S(T)^{\nat},
{\cal E}^{\bul \leq N,\bul})^{ij}
@= A_{\rm zar}({\cal P}^{\rm ex}_{\bul \leq N,\bul}/S_0(T)^{\nat},
{\cal E}^{\bul \leq N,\bul})^{ij}\otimes_{\mab Z}{\mab Q}.
\end{CD} 
\tag{3.3.14.4}\label{eqqn:aittxt} 
\end{equation*} 
Here note that the definitions of $\theta_{{\cal P}^{\rm ex}_{\bul \leq N,\bul}}$'s 
in the vertical morphisms in (\ref{eqqn:aittxt}) are different.  
Hence the equality (\ref{ali:idtt}) induces the isomorphism (\ref{ali:azrxt}). 
We also have the following commutative diagram: 
\begin{equation*} 
\begin{CD} 
A_{{\rm zar},{\mab Q}}({\cal P}^{\rm ex}_{\bul \leq N,\bul}/S(T)^{\nat},
{\cal E}^{\bul \leq N,\bul})^{ij}
@= A_{\rm zar}({\cal P}^{\rm ex}_{\bul \leq N,\bul}/S_0(T)^{\nat},
{\cal E}^{\bul \leq N,\bul})^{ij}\otimes_{\mab Z}{\mab Q}\\
@A{e\theta_{{\cal P}^{\rm ex}_{\bul \leq N,\bul}}\wedge}A{\simeq}A 
@A{\simeq}A{\theta_{{\cal P}^{\rm ex}_{\bul \leq N,\bul}}\wedge}A \\
{\cal E}^{\bul \leq N,\bul}
\otimes_{{\cal O}_{{\cal P}^{\rm ex}_{\bul \leq N,\bul}}}
{\Om}^{i+j+1}_{{\cal P}^{\rm ex}_{\bul \leq N,\bul}/S(T)^{\nat}}
\otimes_{\mab Z}{\mab Q}
@= {\cal E}^{\bul \leq N,\bul}
\otimes_{{\cal O}_{{\cal P}^{\rm ex}_{\bul \leq N,\bul}}}
{\Om}^{i+j+1}_{{\cal P}^{\rm ex}_{\bul \leq N,\bul}/S_0(T)^{\nat}}
\otimes_{\mab Z}{\mab Q}.
\end{CD} 
\tag{3.3.14.5}\label{eqqn:aizerxt} 
\end{equation*} 
This gives the commutative diagram (\ref{eqqn:azerxt}).  
\par 
We leave the reader to the proof that the isomorphism 
(\ref{ali:azrxt}) is independent of the open covering 
of $Y_{\bul \leq N,\os{\circ}{T}_0}$ and the choice of the admissible immersion of 
$Y_{\bul \leq N,\os{\circ}{T}_0}$ because the proof of this is the same as that 
of (\ref{prop:tefc}). 
\end{proof} 


\section{Contravariant functoriality of 
iso-zariskian $p$-adic filtered Steenbrink complexes}\label{sec:pfnsc} 
Let the notations be as in the previous section. 
Let $S'$ be a $p$-adic formal PD-family of log points.   
Let $S\lo S'$ be a morphism of 
$p$-adic formal families of log points over a morphism 
$u_0\col S_0\lo S'_0$ of $p$-adic formal families of log points. 
Let $(T,{\cal J},\del)\lo (T',{\cal J}',\del')$ be a morphism of 
log $p$-adic enlargements over the morphism $S\lo S'$. 
Let $u\col (S(T)^{\nat},{\cal J},\del) \lo (S'(T')^{\nat},{\cal J}',\del')$ 
be the induced morphism by this morphism. 
Let 
\begin{equation*} 
\begin{CD} 
X_{\bul \leq N,\os{\circ}{T}_0} @>{g_{\bul \leq N}}>> Y_{\bul \leq N,\os{\circ}{T}{}'_0}\\ 
@VVV @VVV \\ 
S_{\os{\circ}{T}_0} @>>> S'_{\os{\circ}{T}{}'_0} \\ 
@V{\bigcap}VV @VV{\bigcap}V \\ 
S(T)^{\nat} @>{u}>> S'(T')^{\nat}
\end{CD}
\tag{3.4.0.1}\label{eqn:xdbquss}
\end{equation*} 
be a commutative diagram of $N$-truncated simplicial base changes of SNCL schemes 
over $S_{\os{\circ}{T}_0}$ and $S'_{\os{\circ}{T}{}'_0}$ 
augmented to $S_{0,\os{\circ}{T}_0}$ and $S'_{0,\os{\circ}{T}{}'_0}$, respectively, 
such that $X_{\bul \leq N,\os{\circ}{T}_0}$ and $Y_{\bul \leq N,\os{\circ}{T}{}'_0}$ 
have the disjoint unions 
$X'_{\bul \leq N,\os{\circ}{T}_0}$ and $Y'_{\bul \leq N,\os{\circ}{T}{}'_0}$ 
of the members of affine $N$-truncated simplicial open coverings of 
$X_{\bul \leq N,\os{\circ}{T}_0}$ and $Y_{\bul \leq N,\os{\circ}{T}{}'_0}$, respectively, 
fitting into the following commutative diagram 
\begin{equation*} 
\begin{CD} 
X'_{\bul \leq N,\os{\circ}{T}_0} @>{g'_{\bul \leq N}}>> Y'_{\bul \leq N,\os{\circ}{T}{}'_0} \\
@VVV @VVV \\ 
X_{\bul \leq N,\os{\circ}{T}_0} @>{g_{\bul \leq N}}>> Y_{\bul \leq N,\os{\circ}{T}{}'_0}. 
\end{CD}
\tag{3.4.0.2}\label{cd:xdjqbxy}
\end{equation*} 
As in the previous section, 
we assume that there exist admissible immersions 
\begin{align*} 
X_{\bul \leq N,\bul,\os{\circ}{T}_0} \os{\sus}{\lo} \ol{\cal P}{}'_{\bul \leq N,\bul}
\tag{3.4.0.3}\label{cd:adxpn}
\end{align*} 
and 
\begin{align*} 
Y_{\bul \leq N,\bul,\os{\circ}{T}{}'_0} \os{\sus}{\lo} \ol{\cal Q}_{\bul \leq N,\bul}
\tag{3.4.0.4}\label{cd:adypn}
\end{align*} 
into $(N,\infty)$-truncated bisimplicial 
formally log smooth schemes over $\ol{S(T)^{\nat}}$ augmented to $\ol{S_0(T)^{\nat}}$ 
and over $\ol{S'(T')^{\nat}}$ augmented to $\ol{S'_0(T')^{\nat}}$, respectively. 
Set 
\begin{align*} 
\ol{\cal P}_{\bul \leq N,\bul}:=\ol{\cal P}{}'_{\bul \leq N,\bul}\times_{\ol{S(T)^{\nat}}}
(\ol{\cal Q}_{\bul \leq N,\bul}\times_{\ol{S'(T')^{\nat}}}\ol{S(T)^{\nat}}). 
\end{align*} 
Then we have the following commutative diagram 
\begin{equation*}
\begin{CD} 
X_{\bul \leq N,\bul,\os{\circ}{T}_0} 
@>{\sus}>> \ol{\cal P}_{\bul \leq N,\bul} \\ 
@V{g_{\bul \leq N,\bul}}VV 
@VV{}V \\ 
Y_{\bul \leq N,\bul,\os{\circ}{T}{}'_0} @>{\sus}>> \ol{\cal Q}_{\bul \leq N,\bul}
\end{CD}
\tag{3.4.0.5}\label{cd:necmd}
\end{equation*} 
over the morphism $\ol{S(T)^{\nat}}\lo \ol{S'(T')^{\nat}}$. 
Consequently we have the following commutative diagram 
\begin{equation*}
\begin{CD} 
X_{\bul \leq N,\bul,\os{\circ}{T}_0} 
@>{\sus}>> \ol{\cal P}{}^{\rm ex}_{\bul \leq N,\bul} \\ 
@V{g_{\bul \leq N,\bul}}VV 
@VV{}V \\ 
Y_{\bul \leq N,\bul,\os{\circ}{T}{}'_0} @>{\sus}>> \ol{\cal Q}{}^{\rm ex}_{\bul \leq N,\bul} 
\end{CD}
\tag{3.4.0.6}\label{cd:necemd}
\end{equation*} 
over the morphism $\ol{S(T)^{\nat}}\lo \ol{S'(T')^{\nat}}$. 
Let $\ol{\mathfrak D}_{\bul \leq N,\bul}$ and 
$\ol{\mathfrak E}_{\bul \leq N,\bul}$ be the log PD-envelopes of  
$X_{\bul \leq N,\bul,\os{\circ}{T}_0} \os{\sus}{\lo}\ol{\cal P}_{\bul \leq N,\bul}$ over 
$(\os{\circ}{T},{\cal J},\del)$ and 
$Y_{\bul \leq N,\bul,\os{\circ}{T}{}'_0} \os{\sus}{\lo}\ol{\cal Q}_{\bul \leq N,\bul}$ over 
$(\os{\circ}{T}{}',{\cal J}',\del')$, respectively.      
Set  
${\mathfrak E}_{\bul \leq N,\bul}:=
\ol{\mathfrak E}_{\bul \leq N,\bul}\times_{{\mathfrak D}(\ol{S'(T')^{\nat}})}S'(T')^{\nat}$. 
Then we have the natural morphism  
\begin{equation*} 
\ol{g}{}^{\rm PD}_{\bul \leq N,\bul}\col 
\ol{\mathfrak D}_{\bul \leq N,\bul}\lo \ol{\mathfrak E}_{\bul \leq N,\bul}.   
\tag{3.4.0.7}\label{eqn:gpidf} 
\end{equation*}  
Hence we have the natural morphism 
\begin{equation*} 
g^{\rm PD}_{\bul \leq N,\bul}\col 
{\mathfrak D}_{\bul \leq N,\bul}\lo {\mathfrak E}_{\bul \leq N,\bul}.   
\tag{3.4.0.8}\label{eqn:gpnbdf} 
\end{equation*}  
Let $E^{\bul \leq N}$ and  $F^{\bul \leq N}$ be flat quasi-coherent crystals of 
${\cal O}_{\os{\circ}{X}_{\bul \leq N,T_0}/\os{\circ}{T}}$-modules  
and 
${\cal O}_{\os{\circ}{Y}_{\bul \leq N,T'_0}/\os{\circ}{T}{}'}$-modules,  
respectively.     
Let 
\begin{align*} 
\os{\circ}{g}{}^*_{\bul \leq N,{\rm crys}}(F^{\bul \leq N})
\lo 
E^{\bul \leq N}
\tag{3.4.0.9}\label{ali:gnafe} 
\end{align*} 
be a morphism of 
${\cal O}_{\os{\circ}{X}_{\bul \leq N,T_0}/\os{\circ}{T}}$-modules. 
Let $T'_0$, $Y_{\bul \leq N,\os{\circ}{T}{}'_0}$ and  ${\mathfrak E}_{\bul \leq N,\bul}$ 
be the similar objects to 
$T_0$, $X_{\bul \leq N,\os{\circ}{T}_0}$, ${\mathfrak D}_{\bul \leq N,\bul}$, 
in \S\ref{sec:pgensc}, respectively. 

\begin{theo}[{\bf Contravariant functoriality}]\label{theo:fugennas} 
$(1)$ The morphism  $g_{\bul \leq N}$ induces the following 
well-defined pull-back morphism 
\begin{align*}  
g_{\bul \leq N}^* \col & 
(A_{{\rm zar},{\mab Q}}(Y_{\bul \leq N,\os{\circ}{T}{}'_0}/S'(T')^{\nat},F^{\bul \leq N}),P)
\lo Rg_{\bul \leq N*}((A_{{\rm zar},{\mab Q}}
(X_{\bul \leq N,\os{\circ}{T}_0}/S(T)^{\nat},E^{\bul \leq N}),P))  
\tag{3.4.1.1}\label{eqn:fzgenaxd}
\end{align*} 
fitting into the following commutative diagram$:$
\begin{equation*} 
\begin{CD}
A_{{\rm zar},{\mab Q}}(Y_{\bul \leq N,\os{\circ}{T}{}'_0}/S'(T')^{\nat},F^{\bul \leq N})
@>{g_{\bul \leq N}^*}>>  \\ 
@A{\theta_{Y_{\bul \leq N,\os{\circ}{T}{}'_0}/S'(T')^{\nat}/S'_0(T')^{\nat}} \wedge}A{\simeq}A \\
Ru_{Y_{\bul \leq N,\os{\circ}{T}{}'_0}/S'(T')^{\nat}*}
(\eps^*_{Y_{\bul \leq N,\os{\circ}{T}{}'_0}/S'(T')^{\nat}}(F^{\bul \leq N}))\otimes_{\mab Z}{\mab Q}
@>{g_{\bul \leq N,}^*}>>
\end{CD}
\tag{3.4.1.2}\label{cd:pssfgenccz} 
\end{equation*}
\begin{equation*} 
\begin{CD}
Rg_{\bul \leq N*}
(A_{{\rm zar},{\mab Q}}(X_{\bul \leq N,\os{\circ}{T}_0}/S(T)^{\nat},E^{\bul \leq N}))\\ 
@A{Rg_{\bul \leq N*}
(\theta_{X_{\bul \leq N,\os{\circ}{T}_0}/S(T)^{\nat}/S_0(T)^{\nat}}\wedge)}A{\simeq}A \\
Rg_{\bul \leq N*}Ru_{X_{\bul \leq N,\os{\circ}{T}_0}/S(T)^{\nat}*}
(\eps^*_{X_{\bul \leq N,\os{\circ}{T}_0}/S(T)^{\nat}}(E^{\bul \leq N}))\otimes_{\mab Z}{\mab Q}.
\end{CD}
\end{equation*}
\par 
$(2)$ Let $v\col S'\lo S''$, $v_0\col S'_0\lo S''_0$ and 
$h_{\bul \leq N}\col Y_{\bul \leq N}\lo Z_{\bul \leq N}$ 
be similar morphisms to $u$, $u_0$ and $g_{\bul \leq N}$, respectively. 
Let $(T',{\cal J}',\del') \lo (T'',{\cal J}'',\del'')$ be a morphism of 
log $p$-adic PD-enlargements over $v$. 
Assume that there exists the disjoint union 
$Z'_{\bul \leq N,\os{\circ}{T}{}''_0}$ of the member of an affine simplicial covering 
$Z_{\bul \leq N,\os{\circ}{T}{}''_0}$ and that there exists an admissible immersion 
$Z_{\bul \leq N,\bul,\os{\circ}{T}{}''_0} \os{\sus}{\lo} \ol{\cal R}_{\bul \leq N,\bul}$ 
into a formally log smooth $(N,\infty)$-truncated bisimplicial log formal scheme over $\ol{S''(T'')^{\nat}}$. 
Let $G^{\bul \leq N}$ be a flat quasi-coherent crystal of 
${\cal O}_{\os{\circ}{Z}_{\bul \leq N,T''_0}/\os{\circ}{T}{}''}$-modules.    
Let 
\begin{align*} 
\os{\circ}{h}{}^*_{\bul \leq N,{\rm crys}}(G^{\bul \leq N})
\lo 
F^{\bul \leq N}
\tag{3.4.1.3}\label{ali:gznbfe} 
\end{align*} 
be a morphism of 
${\cal O}_{\os{\circ}{Y}_{\bul \leq N,T'_0}/\os{\circ}{T}{}'}$-modules. 
Then 
\begin{align*} 
(h_{\bul \leq N}\circ g_{\bul \leq N})^* =&
Rh_{\bul \leq N*}(g_{\bul \leq N}^*)\circ 
h_{\bul \leq N}^*    \col 
(A_{{\rm zar},{\mab Q}}(Z_{\bul \leq N,\os{\circ}{T}{}''_0}/S''(T'')^{\nat},G^{\bul \leq N}),P)
\tag{3.4.1.4}\label{ali:pdgenpp} \\ 
&\lo Rh_{\bul \leq N*}Rg_{\bul \leq N*}((A_{{\rm zar},{\mab Q}}(X_{\bul \leq N,\os{\circ}{T}_0}/S(T)^{\nat}
,E^{\bul \leq N}),P)) \\
& =R(h_{\bul \leq N}\circ g_{\bul \leq N})_*
((A_{{\rm zar},{\mab Q}}(X_{\bul \leq N,\os{\circ}{T}_0}/S(T)^{\nat},E^{\bul \leq N}),P)).
\end{align*}  
$(3)$  
\begin{equation*} 
{\rm id}_{X_{\bul \leq N,T_0}}^*={\rm id} 
\col (A_{{\rm zar},{\mab Q}}(X_{\bul \leq N,\os{\circ}{T}_0}/S(T)^{\nat},E^{\bul \leq N}),P)
\lo (A_{{\rm zar},{\mab Q}}(X_{\bul \leq N,\os{\circ}{T}_0}/S(T)^{\nat},E^{\bul \leq N}),P).  
\tag{3.4.1.5}\label{eqn:fziadd}
\end{equation*} 
\end{theo}
\begin{proof} 
We have only to replace $\deg (u)$ in the proof of (\ref{theo:funas}) 
by $\deg (u_0)$. 
\end{proof}

The following is the contravariant functoriality stated in the Introduction. 

\begin{coro}\label{coro:uip} 
Let the notations be as in {\rm (\ref{theo:ctbfd})}. 
Denote $Y_{\bul \leq N,\os{\circ}{T}_0}$ and $X_{\bul \leq N,\os{\circ}{T}_0}$ by 
$Y_{1,\bul \leq N,\os{\circ}{T}_0}$ and $X_{1,\bul \leq N,\os{\circ}{T}_0}$ in this Corollary. 
Let $g_{\bul \leq N}\col Y_{1,\bul \leq N,\os{\circ}{T}_0}\lo  Y_{2,\bul \leq N,\os{\circ}{T}_0}$ 
be a morphism of SNCL schemes over $S_{0,\os{\circ}{T}_0}$, where 
$Y_{2,\bul \leq N,\os{\circ}{T}_0}/S_{0,\os{\circ}{T}_0}$ is a similar SNCL scheme 
to  $Y_{1,\bul \leq N,\os{\circ}{T}_0}/S_{0,\os{\circ}{T}_0}$. 
Set $X_{2,\bul \leq N,\os{\circ}{T}_0}:=Y_{2,\bul \leq N,\os{\circ}{T}_0}
\times_{S_0,\os{\circ}{T}_0}S_{\os{\circ}{T}_0}$ and 
$g^*_{\bul \leq N,S_{\os{\circ}{T}_0}}
:=g_{\bul \leq N}^*\times_{S_{0,\os{\circ}{T}_0}}S_{\os{\circ}{T}_0}$. 
Then the canonical filtered isomorphism {\rm (\ref{ali:azrxt})} 
fits into the following commutative diagram$:$ 
\begin{equation*} 
\begin{CD} 
(A_{\rm zar}(Y_{1,\bul \leq N,\os{\circ}{T}_0}/S_0(T)^{\nat},E^{\bul \leq N}),P)
\otimes^L_{\mab Z}{\mab Q}@>{\sim}>>  
(A_{{\rm zar},{\mab Q}}(X_{1,\bul \leq N,\os{\circ}{T}_0}/S(T)^{\nat},E^{\bul \leq N}),P) \\
@A{g_{\bul \leq N}^*}AA @AA{g^*_{\bul \leq N,S_{\os{\circ}{T}_0}}}A \\
(A_{\rm zar}(Y_{2,\bul \leq N,\os{\circ}{T}_0}/S_0(T)^{\nat},E^{\bul \leq N}),P)
\otimes^L_{\mab Z}{\mab Q}
@>{\sim}>>
(A_{{\rm zar},{\mab Q}}(X_{2,\bul \leq N,\os{\circ}{T}_0}/S(T)^{\nat},E^{\bul \leq N}),P). 
\end{CD} 
\tag{3.4.2.1}\label{ali:sim} 
\end{equation*} 
The morphisms $g_{\bul \leq N}^*$'s and $g_{\bul \leq N,S_{\os{\circ}{T}_0}}^*$'s 
have the transitive relation for the commutative diagram {\rm (\ref{ali:sim})}. 
\end{coro}
\begin{proof} 
This immediately follows from the construction of $g_{\bul \leq N}^*$ in 
(\ref{theo:fugennas}). 
\end{proof} 


\begin{rema}\label{rema:cqa}
(\ref{coro:uip}) tells us that the filtered complex 
$A_{\rm zar}(Y_{i,\bul \leq N,\os{\circ}{T}_0}/S_0(T)^{\nat},E^{\bul \leq N}),P)
\otimes^L_{\mab Z}{\mab Q}$ $(i=1,2)$ has the contravariant functoriality 
for a morphism $h_{\bul \leq N} \col 
Y_{1,\bul \leq N,\os{\circ}{T}_0}\lo Y_{2,\bul \leq N,\os{\circ}{T}_0}$ 
over $S_{\os{\circ}{T}_0}$ not only over $S_{0,\os{\circ}{T}_0}$. 
\end{rema}

\begin{defi}[{\bf Abrelative Frobenius morphism}]\label{defi:qrwd}  
Let $S$, $S^{[p]}$, $F_{S/\os{\circ}{S}}\col S\lo S^{[p]}$ and 
$W_{S/\os{\circ}{S}} \col S^{[p]}\lo S$, $(T,{\cal J},\del)\lo (T',{\cal J}',\del')$ 
be as in (\ref{defi:rwd}).    
Set $X^{[p]}_{\bul \leq N}:=X_{\bul \leq N}\times_SS^{[p]}
=X_{\bul \leq N}\times_{\os{\circ}{S},\os{\circ}{F}_S}\os{\circ}{S}$ 
and $X^{[p]}_{\bul \leq N,\os{\circ}{T}{}'_0}:=
X^{[p]}_{\bul \leq N}
\times_{S^{[p]}}(S^{[p]})_{\os{\circ}{T}{}'_0}$. 
Let 
$$F^{\rm ar}_{X_{\bul \leq N,\os{\circ}{T}_0/
S^{[p]}_{\os{\circ}{T}{}'_0},S_{\os{\circ}{T}_0}}}\col 
X_{\bul \leq N,\os{\circ}{T}_0}  \lo X^{[p]}_{\bul \leq N,\os{\circ}{T}{}'_0}$$ 
and 
$$F^{\rm ar}_{X_{\bul \leq N,\os{\circ}{T}_0/S^{[p]}(T')^{\nat},S(T)^{\nat}}}\col 
X_{\bul \leq N,\os{\circ}{T}_0}  \lo X^{[p]}_{\bul \leq N,\os{\circ}{T}{}'_0}$$ 
be the relative Frobenius morphism  
over $S_{\os{\circ}{T}_0}\lo (S^{[p]})_{\os{\circ}{T}{}'_0}$ and 
$(S(T)^{\nat},{\cal J},\del)\lo (S^{[p]}(T')^{\nat},{\cal J}',\del')$.  
Let $E^{\bul \leq N}$ and $E'{}^{\bul \leq N}$ be a flat quasi-coherent crystal of 
${\cal O}_{\os{\circ}{X}_{\bul \leq N,T_0}/\os{\circ}{T}}$-modules and 
a flat quasi-coherent crystal of 
${\cal O}_{\os{\circ}{X}{}^{[p]}_{\bul \leq N,T_0}/\os{\circ}{T}}$-modules, 
respectively.  
Let 
\begin{align*} 
\Phi^{\rm ar} \col 
\os{\circ}{F}{}^{{\rm ar}*}_{X_{\bul \leq N,\os{\circ}{T}_0/S(T)^{\nat},
S^{[p]}(T')^{\nat},{\rm crys}}}(E'{}^{\bul \leq N})
\lo E^{\bul \leq N}
\tag{3.4.4.1}\label{ali:qsppts}
\end{align*} 
be a morphism of crystals of 
${\cal O}_{\os{\circ}{X}_{\bul \leq N,T_0}/\os{\circ}{T}}$-modules.   
We call the following induced morphism by $\Phi^{\rm ar}$   
\begin{align*} 
\Phi^{\rm ar} \col &
(A_{{\rm zar},{\mab Q}}(X^{[p]}_{\bul \leq N,\os{\circ}{T}{}'_0}/
S^{[p]}(T')^{\nat},E'{}^{\bul \leq N}),P) 
\tag{3.4.4.2}\label{ali:sqpapts} \\
&\lo 
RF^{\rm ar}_{X_{\bul \leq N,\os{\circ}{T}_0/S(T)^{\nat},S^{[p]}(T')^{\nat}}*}
((A_{{\rm zar},{\mab Q}}(X_{\bul \leq N,\os{\circ}{T}_0}/S(T)^{\nat},E^{\bul \leq N}),P)) 
\end{align*}
the {\it abrelative Frobenius morphism} of 
$$(A_{{\rm zar},{\mab Q}}(X_{\bul \leq N,\os{\circ}{T}_0}/S(T)^{\nat},E^{\bul \leq N}),P)
\quad {\rm and} \quad  
(A_{{\rm zar},{\mab Q}}(X^{[p]}_{\bul \leq N,\os{\circ}{T}{}'_0}/S^{[p]}(T')^{\nat},
E'^{\bul \leq N}),P).$$ 
When $E'{}^{\bul \leq N}={\cal O}_{\os{\circ}{X}{}'_{\bul \leq N,T'_0}/\os{\circ}{T}{}'}$, 
we denote 
$(A_{{\rm zar},{\mab Q}}(X^{[p]}_{\bul \leq N,\os{\circ}{T}{}'_0}/S^{[p]}(T')^{\nat},E'{}^{\bul \leq N}),P)$ by 
$(A_{{\rm zar},{\mab Q}}(X^{[p]}_{\bul \leq N,\os{\circ}{T}{}'_0}/S^{[p]}(T')^{\nat}),P)$.   
Then we have the following {\it abrelative Frobenius morphism} 
\begin{equation*} 
\Phi^{\rm ar} \col 
(A_{{\rm zar},{\mab Q}}(X^{[p]}_{\bul \leq N,\os{\circ}{T}{}'_0}/S^{[p]}(T')^{\nat}),P) 
\lo RF^{\rm ar}_{X_{\bul \leq N,\os{\circ}{T}_0/S(T)^{\nat},S^{[p]}(T')^{\nat}}*}
((A_{{\rm zar},{\mab Q}}(X_{\bul \leq N,\os{\circ}{T}_0}/S(T)^{\nat}),P)) 
\tag{3.4.4.3}
\end{equation*}
of $(A_{{\rm zar},{\mab Q}}(X_{\bul \leq N,\os{\circ}{T}_0}/S(T)^{\nat}),P)$ and 
$(A_{{\rm zar},{\mab Q}}(X^{[p]}_{\bul \leq N,\os{\circ}{T}{}'_0}/S^{[p]}(T')^{\nat}),P)$. 
\end{defi}

\begin{prop}[{\bf Frobenius compatibility I}]\label{prop:afcar} 
The following diagram is commutative$:$ 
\begin{equation*} 
\begin{CD} 
A_{{\rm zar},{\mab Q}}
(X^{[p]}_{\bul \leq N,\os{\circ}{T}{}'_0}/S^{[p]}(T')^{\nat},E'{}^{\bul \leq N})
@>{\Phi^{\rm ar}}>>
 \\
@A{\theta_{X^{[p]}_{\bul \leq N,\os{\circ}{T}_0}/S^{[p]}(T')^{\nat}/
S^{[p]}_0(T)^{\nat}} \wedge}A{\simeq}A  \\
Ru_{X^{[p]}_{\bul \leq N,\os{\circ}{T}{}'_0}/S^{[p]}(T')*}
(\eps^*_{X^{[p]}_{\bul \leq N,\os{\circ}{T}{}'_0}/S^{[p]}(T')}(E'{}^{\bul \leq N}))
\otimes_{\mab Z}^L{\mab Q}
@>{\Phi^{\rm ar}}>>
\end{CD}
\tag{3.4.5.1}\label{cd:ruakt}
\end{equation*} 
\begin{equation*} 
\begin{CD} 
RF^{\rm ar}_{X_{\bul \leq N,\os{\circ}{T}_0/T}*}
(A_{{\rm zar},{\mab Q}}(X_{\bul \leq N,\os{\circ}{T}_0}/S(T)^{\nat},E^{\bul \leq N}) \\
@A{\simeq}A{RF^{\rm rel}_{X_{\bul \leq N,\os{\circ}{T}_0/T}*}(
\theta_{X_{\bul \leq N,\os{\circ}{T}_0}/S(T)^{\nat}/S_0(T)^{\nat}})\wedge}A \\
RF^{\rm ar}_{X_{\bul \leq N,\os{\circ}{T}_0/S^{[p]}(T')^{\nat},S(T)^{\nat}}*}
Ru_{X_{\bul \leq N,\os{\circ}{T}_0}/S(T)^{\nat}*}(E^{\bul \leq N})
\otimes_{\mab Z}^L{\mab Q}. 
\end{CD}
\end{equation*} 
The commutative diagram {\rm (\ref{cd:ruakt})} is contravariantly functorial 
for the morphism {\rm (\ref{eqn:xdbquss})} satisfying {\rm (\ref{cd:xdjqbxy})} 
and for the morphism of $F$-crystals 
\begin{equation*} 
\begin{CD} 
\os{\circ}{F}{}^{{\rm ar}*}_{X_{\bul \leq N,\os{\circ}{T}_0/S^{[p]}(T')^{\nat},S(T)^{\nat},{\rm crys}}}
(E'{}^{\bul \leq N})
@>{\Phi^{\rm ar}}>> E^{\bul \leq N}\\
@AAA @AAA \\
\os{\circ}{g}{}^*_{\bul \leq N}
\os{\circ}{F}{}^{{\rm ar}*}_{Y_{\bul \leq N,\os{\circ}{T}_0/S^{[p]}(T')^{\nat},S(T)^{\nat},{\rm crys}}}
(F'{}^{\bul \leq N})
@>{\os{\circ}{g}{}^*_{\bul \leq N}(\Phi^{\rm ar})}>> 
\os{\circ}{g}{}^*_{\bul \leq N}(F^{\bul \leq N}), 
\end{CD}
\end{equation*} 
where $F^{\bul \leq N}$ and $F'{}^{\bul \leq N}$ are similar quasi-coherent 
${\cal O}_{\os{\circ}{Y}_{\bul \leq N,T'_0}/\os{\circ}{T}{}'}$-modules  
to $E^{\bul \leq N}$ and 
${\cal O}_{\os{\circ}{Y}{}^{[p]}_{\bul \leq N,T'_0}/\os{\circ}{T}{}'}$-modules 
to $E'{}^{\bul \leq N}$, 
respectively.  
\end{prop} 
\begin{proof} 
This is a special case of (\ref{theo:fugennas}). 
\end{proof} 

\begin{prop}[{\bf Frobenius compatibility II}]\label{prop:nwsfcs} 
Assume that $T$ is restrictively hollow with respective to the morphism $T_0\lo S$.  
Assume that the morphism $T'_0\lo S^{[p]}$ factors through the morphism 
$F_{S/{\os{\circ}{S}}}\col S\lo S^{[p]}$.  
Set $X^{\{p\}}_{\bul \leq N}
:=X_{\bul \leq N}\times_{S,F_S}S$ 
and $X^{\{p\}}_{\bul \leq N,\os{\circ}{T}{}'_0}:=
X^{\{p\}}_{\bul \leq N}\times_{S,F_S}S\times_SS_{\os{\circ}{T}{}'_0}$. 
Let 
$$F^{{\rm rel}}_{X_{\bul \leq N,\os{\circ}{T}_0/S_{\os{\circ}{T}{}'},S_{\os{\circ}{T}_0}}}
\col X_{\bul \leq N,\os{\circ}{T}_0}  \lo X^{\{p\}}_{\bul \leq N,\os{\circ}{T}{}'_0}$$ 
and 
$$F^{\rm rel}_{X_{\bul \leq N,\os{\circ}{T}_0/S(T'),S(T)}}\col 
X_{\bul \leq N,\os{\circ}{T}_0}  \lo X^{\{p\}}_{\bul \leq N,\os{\circ}{T}{}'_0}$$ 
be the relative Frobenius morphisms  
over $S_{\os{\circ}{T}_0}\lo S_{\os{\circ}{T}{}'_0}$ and 
$(S(T)^{\nat},{\cal J},\del)\lo (S(T')^{\nat},{\cal J},\del)$, respectively. 
Let $E''{}^{\bul \leq N}$ be the pull-back of 
$E'{}^{\bul \leq N}$ to $(\os{\circ}{X}{}^{\{p\}}_{\bul \leq N,T_0}/\os{\circ}{T})_{\rm crys}$. 
Let 
\begin{align*} 
\Phi^{{\rm rel}} \col 
\os{\circ}{F}{}^{{\rm  rel}*}_{X_{\bul \leq N,\os{\circ}{T}_0/S(T'),S(T),{\rm crys}}}
(E''{}^{\bul \leq N})\lo E^{\bul \leq N}  
\tag{3.4.6.1}\label{ali:spphtpps}
\end{align*} 
be the induced morphism by {\rm (\ref{ali:qsppts})}. 
Then the following diagram is commutative$:$ 
\begin{equation*} 
\begin{CD} 
A_{\rm zar}(X^{[p]}_{\bul \leq N,\os{\circ}{T}{}'_0}
/S^{[p]}(T')^{\nat},E'{}^{\bul \leq N}) 
@>{\Phi^{\rm  rel}}>>
RF^{{\rm rel}}_{X_{\bul \leq N,\os{\circ}{T}_0/T'T}*}
(A_{\rm zar}(X_{\bul \leq N,\os{\circ}{T}_0}/S(T)^{\nat},E^{\bul \leq N})) \\
@A{\theta \wedge}A{\simeq}A 
@A{RF^{\rm rel}_{X_{\bul \leq N,\os{\circ}{T}_0/T'T}*}(\theta)\wedge}A{\simeq}A \\
Ru_{X^{\{p\}}_{\bul \leq N,T'_0}/T'*}(\eps^*_{X^{\{p\}}_{\bul \leq N,T'_0}/T'}
(E''{}^{\bul \leq N}))
@>{\Phi^{\rm rel}}>>RF^{\rm rel}_{X_{\bul \leq N,\os{\circ}{T}_0/T',T}*}
Ru_{X_{\bul \leq N,T_0}/T*}
(\eps^*_{X_{\bul \leq N,T_0}/T}(E^{\bul \leq N})). 
\end{CD}
\tag{3.4.6.2}\label{cd:reqqkt}
\end{equation*} 
The commutative diagram {\rm (\ref{cd:reqqkt})} is contravariantly functorial 
for the morphism {\rm (\ref{eqn:xdbquss})} satisfying {\rm (\ref{cd:xdjqbxy})}
and for the morphism of $F$-crystals 
\begin{equation*} 
\begin{CD} 
\os{\circ}{F}{}^{{\rm rel}*}_{X_{\bul \leq N,\os{\circ}{T}_0/S(T'),S(T),{\rm crys}}}
(E''{}^{\bul \leq N})
@>{\Phi^{\rm rel}}>> E^{\bul \leq N}\\
@AAA @AAA \\
\os{\circ}{g}{}^*_{\bul \leq N}
\os{\circ}{F}{}^{{\rm rel}*}_{Y_{\bul \leq N,\os{\circ}{T}_0/S(T')^{\nat},S(T)^{\nat},{\rm crys}}}
(F''{}^{\bul \leq N})
@>{\os{\circ}{g}{}^*_{\bul \leq N}(\Phi^{\rm rel})}>> 
\os{\circ}{g}{}^*_{\bul \leq N}(F^{\bul \leq N}), 
\end{CD}
\end{equation*} 
where $F^{\bul \leq N}$ is a similar quasi-coherent 
${\cal O}_{\os{\circ}{Y}_{\bul \leq N,T_0}/\os{\circ}{T}}$-module 
to $E^{\bul \leq N}$. 
\end{prop}
\begin{proof} 
Because the morphism 
$T_0\lo S^{[p]}$ is the composite morphism 
$T_0\lo S_{\os{\circ}{T}_0}\lo S^{[p]}_{\os{\circ}{T}_0}$, 
$X^{[p]}_{\bul \leq N,T_0}=X^{\{p\}}_{\bul \leq N,T_0}$. 
Consequently  
\begin{align*}
Ru_{X^{[p]}_{\bul \leq N,T_0}/T*}(\eps^*_{X^{[p]}_{\bul \leq N,T_0}/T}
(E'{}^{\bul \leq N})) 
=
Ru_{X^{\{p\}}_{\bul \leq N,T_0}/T*}(\eps^*_{X^{\{p\}}_{\bul \leq N,T_0}/T}
(E''{}^{\bul \leq N})). 
\end{align*} 
Hence we obtain the commutative diagram 
(\ref{cd:reqqkt}) by (\ref{cd:ruakt}) and (\ref{lemm:flpis}). 
\par 
We leave the proof of the functoriality to the reader. 
\end{proof}

\begin{defi}[{\bf Absolute Frobenius endomorphism}]\label{defi:btqd}  
Let the notations be as in (\ref{prop:nwsfcs}).  
Let 
$$F^{\rm abs}_{X_{\bul \leq N,\os{\circ}{T}_0}/S(T)^{\nat}} 
\col X_{\bul \leq N,\os{\circ}{T}_0}  \lo X_{\bul \leq N,\os{\circ}{T}_0}$$ 
be the absolute Frobenius endomorphism over $F_{S(T)^{\nat}}$.   
Let 
\begin{align*} 
\Phi^{\rm abs} \col 
\os{\circ}{F}{}^{{\rm abs}*}_{X_{\bul \leq N,\os{\circ}{T}_0/S(T)^{\nat}},{\rm crys}}
(E^{\bul \leq N})\lo E^{\bul \leq N}
\end{align*} 
be a morphism of crystals in 
$(\os{\circ}{X}_{\bul \leq N,T_0}/\os{\circ}{T})_{\rm crys}$.   
Then we call the induced morphism by $\Phi^{\rm abs}$ and 
$F^{\rm abs}_{X_{\bul \leq N,\os{\circ}{T}_0}/S(T)^{\nat}}$
\begin{equation*} 
\Phi^{\rm abs} \col 
(A_{{\rm zar},{\mab Q}}(X_{\bul \leq N,\os{\circ}{T}_0}/S(T)^{\nat},E^{\bul \leq N}),P) 
\lo RF^{\rm abs}_{X_{\bul \leq N,\os{\circ}{T}_0/S(T)^{\nat}}*}
((A_{\rm zar}(X_{\bul \leq N,\os{\circ}{T}_0}/S(T)^{\nat},E^{\bul \leq N}),P)) 
\tag{3.4.7.1}\label{eqn:eqbqnp}
\end{equation*}
the {\it absolute Frobenius endomorphism} of 
$(A_{{\rm zar},{\mab Q}}(X_{\bul \leq N,\os{\circ}{T}_0}/S(T)^{\nat},E^{\bul \leq N}),P)$ 
with respect to $F_S$ and $F_{S(T)^{\nat}}$.
When $E^{\bul \leq N}={\cal O}_{\os{\circ}{X}_{\bul \leq N,T_0}/\os{\circ}{T}}$, 
we have the following {\it absolute Frobenius endomorphism} 
\begin{equation*} 
\Phi^{\rm abs}{}^* \col 
(A_{\rm zar}(X_{\bul \leq N,\os{\circ}{T}_0}/S(T)^{\nat}),P) 
\lo RF^{\rm abs}_{X_{\bul \leq N,T_0/T}*}
((A_{\rm zar}(X_{\bul \leq N,\os{\circ}{T}_0}/S(T)^{\nat}),P)) 
\tag{3.4.7.2}\label{eqn:abqst}
\end{equation*}
of $(A_{\rm zar}(X_{\bul \leq N,\os{\circ}{T}_0}/S(T)^{\nat}),P)$ 
with respect to $F_S$ and $F_{S(T)^{\nat}}$.  
\end{defi}  

\begin{rema} 
We leave the formulation of the analogue of 
(\ref{prop:nwsfcs}) for the absolute Frobenius endomorphism. 
\end{rema}

\par 
Set $X:=X_m$ and $Y:=Y_m$. 
Fix $m$ $(0\leq m\leq N)$.  
As in (1.5.6.4) and (1.5.6.5), 
assume that the following two conditions hold:
\medskip
\parno 
$(3.4.8.1)$: for any smooth component 
$\os{\circ}{X}_{\lam}$ of the SNC scheme $\os{\circ}{X}_{T_0}$ over $\os{\circ}{T}_0$, 
there exists a unique smooth component 
$\os{\circ}{Y}_{\mu}$ of $\os{\circ}{Y}_{T'_0}$ over $\os{\circ}{T}{}'_0$ 
such that $g$ induces a morphism 
$g_{\lam}\col \os{\circ}{X}_{\lam} \lo \os{\circ}{Y}_{\mu}$. 
(Let $\Lam$ and $M$ be the sets of indices 
of the $\lam$'s and the $\mu$'s, respectively. 
We obtain a function 
$\phi \col \Lam \owns \lam \lom \mu \in M$.)
\medskip
\par
\medskip 
\parno
$(3.4.8.2)$: 
there exist positive integers $e({\lam})$'s  
$(\lam \in \Lam)$ such that 
there exist local equations $x_{\lam}=0$ and 
$y_{\phi(\lam)}=0$ of $\os{\circ}{X}_{\lam}$ and $\os{\circ}{Y}_{\phi(\lam)}$, 
respectively,  
such that $g^*(y_{\phi(\lam)})=x^{e({\lam})}_{\lam}$.  
\medskip
\par 
Consider the case $m=N$. 
For $X:=X_N$ and $Y:=Y_N$, we use the same notations 
$\Lam^{(k)}(g)$ ($k\in {\mab Z}_{>0}$), $\ul{\lam}=\{\lam_0,\ldots, \lam_k\}$ 
$(\ul{\lam}\in \Lam^{(k)}(g))$, $\os{\circ}{X}_{\ul{\lam}}$, 
$\os{\circ}{Y}_{\phi(\ul{\lam})}$, 
$\os{\circ}{g}_{\ul{\lam}} 
\col \os{\circ}{X}_{\ul{\lam}} \lo \os{\circ}{Y}_{\phi(\ul{\lam})}$ 
$(\ul{\lam}\in \Lam^{(k)}(g))$, 
$a_{\ul{\lam}}
\col \os{\circ}{X}_{\ul{\lam}}
\lo \os{\circ}{X}_{T_0}$ and 
$b_{\phi(\ul{\lam})}\col \os{\circ}{Y}_{\ul{\lam}}\lo \os{\circ}{Y}_{T'_0}$ 
as the notations after (\ref{prop:dvok}). 
\par
For integers $j$ and $k$ such that $j\geq \max \{-k,0\}$, 
set $\ul{\lam}:=\{\lam_0, \ldots, \lam_{2j+k}\}\in 
\Lam^{(2j+k)}(g)$. 
Let 
$$a_{\ul{\lam}{\rm crys}} 
\col 
((\os{\circ}{X}_{\ul{\lam}}/{\os{\circ}{T}})_{\rm crys},
{\cal O}_{\os{\circ}{X}_{\ul{\lam},T_0}/{\os{\circ}{T}}})
\lo 
((\os{\circ}{X}_{T_0}/{\os{\circ}{T}})_{\rm crys},
{\cal O}_{\os{\circ}{X}_{T_0}/\os{\circ}{T}})$$  
and 
$$b_{\phi(\ul{\lam}){\rm crys}} 
\col 
((\os{\circ}{Y}_{\phi(\ul{\lam}),T'_0}/{\os{\circ}{T}{}'})_{\rm crys},
{\cal O}_{\os{\circ}{Y}_{\phi(\ul{\lam}),T'_0}/{\os{\circ}{T}{}'}})
\lo 
((\os{\circ}{Y}_{T'_0}/{\os{\circ}{T}{}'})_{\rm crys},
{\cal O}_{\os{\circ}{Y}_{T'_0}/\os{\circ}{T}{}'})$$  
be the induced morphisms of ringed topoi by 
$a_{\ul{\lam}}$ 
and $b_{\phi(\ul{\lam})}$, respectively. 
Set $E_{\ul{\lam}}:=a^*_{\ul{\lam}{\rm crys}}(E)$ 
and $F_{\phi(\ul{\lam})}:=b^*_{\phi(\ul{\lam}){\rm crys}}(F)$.  
Let 
$\vp_{\ul{\lam}{\rm crys}}(\os{\circ}{X}_{T_0}/\os{\circ}{T})$ 
(resp.~
$\vp_{\phi(\ul{\lam}){\rm crys}}(\os{\circ}{Y}_{T'_0}/\os{\circ}{T}{}')$) 
be the crystalline orientation sheaf in 
$(\os{\circ}{X}_{\ul{\lam}}/{\os{\circ}{T}})_{\rm crys},$ 
(resp.~$(\os{\circ}{Y}_{\ul{\lam}}/{\os{\circ}{T}{}'})_{\rm crys}$) 
defined similarly in \S\ref{sec:ldc} for the set 
$\{\os{\circ}{X}_{\lam_0}, \ldots, \os{\circ}{X}_{\lam_{2j+k}}\}$ 
(resp.~$\{\os{\circ}{Y}_{\lam_0}, \ldots, \os{\circ}{Y}_{\lam_{2j+k}}\}$). 
\par 
Let $Z$ be $X$ or $Y$. Let $\ul{\mu}$ be $\ul{\lam}$ or $\phi(\ul{\lam})$. 
Let $U$ be $T$ or $T'$. 
Let $G$ be an ${\cal O}_{\os{\circ}{Z}_{U_0}/\os{\circ}{U}}$-module. 
Let $G_{\os{\circ}{Z}_{\ul{\mu}}/\os{\circ}{U}}$ be 
an ${\cal O}_{\os{\circ}{Z}_{\ul{\mu}}/\os{\circ}{U}}$-module 
which is the pull-back of $G$.   
Consider the morphism 
\begin{align*} 
&{\rm gr}^P_k(g^*_{\bul \leq N})
\col {\rm gr}^P_kA_{{\rm zar},{\mab Q}}
(Y_{\bul \leq N,\os{\circ}{T}{}'_0}/S'(T')^{\nat},F^{\bul \leq N})
\lo 
Rg_{\bul \leq N*}({\rm gr}^P_k
A_{{\rm zar},{\mab Q}}(X_{\bul \leq N,\os{\circ}{T}_0}/S(T)^{\nat},
E^{\bul \leq N})) \tag{3.4.8.3}\label{ali:gdgskcl}. 
\end{align*}
Consider the following direct factor of the cosimplicial degree 
$m$-part of the morphism (\ref{ali:gdgskcl}): 
\begin{align*}
g^{*}_{\ul{\lam}} \col 
& 
b_{\phi(\ul{\lam})*}
Ru_{\os{\circ}{Y}_{\phi(\ul{\lam})}/T'*}
(F_{\os{\circ}{Y}_{\phi(\ul{\lam})}/\os{\circ}{T}{}'}
\otimes_{\mab Z}
\vp_{\phi(\ul{\lam}){\rm crys}}
(\os{\circ}{Y}_{T'}/\os{\circ}{T}{}'_0))
[-2j-k]  \tag{3.4.8.4}\label{ali:grgsgm}\\
{} & \lo b_{\phi(\ul{\lam})*}
Rg_{\ul{\lam}*}Ru_{\os{\circ}{X}_{\ul{\lam}}/\os{\circ}{T}*}
(E_{\os{\circ}{X}_{\ul{\lam}}
/\os{\circ}{T}}
\otimes_{\mab Z}
\vp_{\ul{\lam}{\rm crys}}
(\os{\circ}{X}_{T_0}/\os{\circ}{T}))[-2j-k]. 
\end{align*}

\begin{prop}\label{prop:grlgoc}
Let the notations and the assumptions be as above.
Let 
$$g_{\ul{\lam}}\col \os{\circ}{X}_{\ul{\lam}}
\lo \os{\circ}{Y}_{\phi(\ul{\lam})}$$ 
be the induced morphism by $g$.
Then the morphism 
$g^*_{\ul{\lam}}$ in 
{\rm (\ref{ali:grgsgm})}  is equal to 
$$\prod_{j=1}^k\sum_{i(N)=1}^{l(N)}d_{i(N)}^{j+k}
b_{\phi(\ul{\lam}){\rm crys}*}g^*_{\ul{\lam}}$$ 
for $k >0$. 
$($Recall the notation $l(N)$ in {\rm (\ref{eqn:oxypt})}$)$.
\end{prop}
\begin{proof}
We use the similar notations to those in the proof of (\ref{prop:grloc}). 
We have the following commutative diagram 
by the proof of (\ref{prop:grloc}) and 
by using the residue isomorphism (\ref{ali:opisst}): 
\begin{equation*}
\begin{CD}
{\rm gr}_k^PA_{{\rm zar},{\mab Q}}(Y_{\os{\circ}{T}{}'_0}/S'(T')^{\nat},F) 
@>{{\rm gr}_k^P(g^*)}>> \\ 
@V{{\rm Res}^{\os{\circ}{\cal Y}_{\phi(\ul{\lam})}}}VV\\ 
b'_{\ul{\lam},T'*}
({\cal F}_{\ul{\lam},T'}
\otimes_{{\cal O}_{\os{\circ}{\cal Y}{}_{\phi(\ul{\lam})}}}
\Om^{\bul}_{\os{\circ}{\cal Y}_{\phi(\ul{\lam})}/\os{\circ}{T}{}'}
\otimes_{\mab Z}
\vp_{\phi(\ul{\lam}){\rm zar}}
(\os{\circ}{\cal Y}/\os{\circ}{T}{}')_{\mab Q}[-2j-k]
@>{d_{i(N)}^{j+k}g^{{\rm PD}*}_{\ul{\lam},T'T}}>>  \\
\end{CD}
\tag{3.4.9.1}\label{cd:gratg}
\end{equation*} 
\begin{equation*}
\begin{CD}
g_{*}A_{{\rm zar},{\mab Q}}(X_{\os{\circ}{T}{}_0}/S(T)^{\nat},E) \\
@VV{{\rm Res}^{\os{\circ}{\cal X}_{\ul{\lam}}}}V\\ 
b'_{\ul{\lam}*}g_{*}a_{\ul{\lam}*}
({\cal E}_{\ul{\lam}}
\otimes_{{\cal O}_{\os{\circ}{\cal X}_{\ul{\lam}}}}
\Om^{\bul}_{\os{\circ}{\cal X}_{\ul{\lam}}/\os{\circ}{T}}
\otimes_{\mab Z}\vp_{\ul{\lam}{\rm zar}}
(\os{\circ}{\cal X}/\os{\circ}{T})_{\mab Q}[-2j-k]. 
\end{CD}
\end{equation*} 
We complete the proof.  
\end{proof}

\par 
Let $f\col X_{\bul \leq N,\os{\circ}{T}{}_0} \lo S(T)^{\nat}$ 
be the structural morphism. 
Let $(S(T)^{\nat})_{\bul \leq N}$ be the 
$N$-truncated constant simplicial log scheme defined by $S(T)^{\nat}$. 
The morphism $f$ induces the natural morphism 
$f_{\bul \leq N}
\col X_{\bul \leq N,\os{\circ}{T}_0}\lo (S(T)^{\nat})_{\bul \leq N}$.  
Then we have the $N$-truncated cosimplicial filtered complex 
\begin{equation*} 
Rf_{{\bul \leq N}*}
((A_{{\rm zar},{\mab Q}}(X_{\bul \leq N,\os{\circ}{T}_0}/S(T)^{\nat},
E^{\bul \leq N}),P)) 
\in {\rm D}^+{\rm F}
({\cal O}_{T_{\bul \leq N}}\otimes_{\mab Z}{\mab Q}).  
\end{equation*} 
Let ${\bf s} \col {\rm D}^+{\rm F}
({\cal O}_{T_{\bul \leq N}}\otimes_{\mab Z}{\mab Q}) 
\lo {\rm D}^+{\rm F}({\cal O}_T\otimes_{\mab Z}{\mab Q})$ 
be the single complex functor (see \cite[(2.0.1)]{nh3} 
for the convention on the signs of boundary morphisms). 
Then we have the following formula 
in ${\rm D}^+{\rm F}({\cal O}_T\otimes_{\mab Z}{\mab Q})$:
\begin{align*} 
& Rf_*
((A_{{\rm zar},{\mab Q}}(X_{\bul \leq N,\os{\circ}{T}_0}/S(T)^{\nat},
E^{\bul \leq N}),P) 
\os{\sim}{\lo}  {\bf s}Rf_{{\bul \leq N}*}
((A_{{\rm zar},{\mab Q}}(X_{\bul \leq N,\os{\circ}{T}_0}/S(T)^{\nat},
E^{\bul \leq N}),P)).  
\end{align*} 
Let $L$ be the stupid filtration on 
$$Rf_{{\bul \leq N}*}
(A_{{\rm zar},{\mab Q}}(X_{\bul \leq N,\os{\circ}{T}_0}/S(T)^{\nat},
E^{\bul \leq N}))$$  
with respect to the cosimplicial degree:  
\begin{align*}
&L^mRf_{{\bul \leq N}*}
(A_{{\rm zar},{\mab Q}}(X_{\bul \leq N,\os{\circ}{T}_0}/S(T)^{\nat},
E^{\bul \leq N})) 
=\bigoplus_{m'\geq m}
Rf_{m',T*}
(A_{{\rm zar},{\mab Q}}(X_{m',\os{\circ}{T}_0}/S(T)^{\nat},E^{m'})). 
\tag{3.4.9.2}\label{eqn:lasxs}
\end{align*}
\parno
Let 
$\del(L,P)$ 
be the diagonal filtration of $L$  and $P$
on 
$$Rf_{{\bul \leq N}*}
(A_{{\rm zar},{\mab Q}}(X_{\bul \leq N,\os{\circ}{T}_0}/S(T)^{\nat},
E^{\bul \leq N})) $$ 
(cf.~\cite[(7.1.6.1), (8.1.22)]{dh3}): 
\begin{align*}
&\del(L,P)_k
Rf_{\bul \leq N,T*} 
(A_{{\rm zar},{\mab Q}}(X_{\bul \leq N,\os{\circ}{T}_0}/S(T)^{\nat},
E^{\bul \leq N})) 
\tag{3.4.9.3}\label{ali:ddisfl} \\
& =  
\bigoplus_{m \geq 0}
P_{k+m}Rf_{m,T*} 
(A_{{\rm zar},{\mab Q}}(X_{m,\os{\circ}{T}_0}/S(T)^{\nat}
,E^m)). 
\end{align*} 
Then we have the following by (\ref{ali:xnsnp}):  
\begin{align*} 
& {\rm gr}^{\del(L,P)}_k
Rf_{{\bul \leq N}*}
(A_{{\rm zar},{\mab Q}}(X_{\bul \leq N,\os{\circ}{T}_0}/S(T)^{\nat},
E^{\bul \leq N}))
= \tag{3.4.9.4}\label{ali:rutgrvp}\\
& \bigoplus_{m\geq 0}\bigoplus_{j\geq \max \{-(k+m),0\}} 
a^{(2j+k+m)}_{m*} 
(Ru_{\os{\circ}{X}{}^{(2j+k+m)}_{m,T}/\os{\circ}{T}*}
(E_{\os{\circ}{X}{}^{(2j+k+m)}_{m,T_0}/\os{\circ}{T}} 
\otimes_{\mab Z}\vp_{\rm crys}^{(2j+k+m)}
(\os{\circ}{X}_{m,T_0}/\os{\circ}{T}))) \\
& (-j-k-m,u)[-2j-k-2m].  
\end{align*} 
Set $f_{X_{\bul \leq N,T_0}/T}:=f\circ u_{X_{\bul \leq N,T_0}/T}$. 
By (\ref{ali:usz}) and (\ref{ali:rutgrvp}),   
we have the following spectral sequence, 
which is a generalization of (\ref{eqn:getpsp}): 
\begin{align*} 
& E_1^{-k,q+k} = \bigoplus_{m=0}^N
\bigoplus_{j\geq \max \{-(k+m),0\}} 
R^{q-2j-k-2m}
f_{\os{\circ}{X}{}^{(2j+k+m)}_{m,T}/\os{\circ}{T}*}
(E_{\os{\circ}{X}{}^{(2j+k+m)}_{m,T_0}/\os{\circ}{T}}
\otimes_{\mab Z} 
\tag{3.4.9.5}\label{eqn:espssp} \\
& \vp^{(2j+k+m)}_{\rm crys}
(\os{\circ}{X}_{m,T_0}/\os{\circ}{T}))(-j-k-m,u)_{\mab Q} 
\Lo 
R^qf_{X_{\bul \leq N,\os{\circ}{T}_0}/S(T)^{\nat}*}
(\eps^*_{X_{\bul \leq N,\os{\circ}{T}_0}/S(T)^{\nat}}
(E^{\bul \leq N}))_{\mab Q} \\
&\quad (q\in {\mab Z}).  
\end{align*}  
Here $(-j-k-m,u)$ means the twist 
``$\bigoplus_{i(m)=1}^{l(m)}d_{i(m)}^{j+k}\cdot$''. 
More generally, for $k\in {\mab Z}$, 
we have the following spectral sequence 
\begin{align*} 
& E_1^{-k',q+k'} = \bigoplus_{m=0}^N
\bigoplus_{j\geq \max \{-(k'+m),0\}} 
R^{q-2j-k'-2m}
f_{\os{\circ}{X}{}^{(2j+k'+m)}_{m,T}/\os{\circ}{T}*}
(E_{\os{\circ}{X}{}^{(2j+k'+m)}_{m,T_0}/\os{\circ}{T}}
\otimes_{\mab Z}\tag{3.4.9.6}\label{eqn:ekpssp}\\
&\vp^{(2j+k'+m)}_{\rm crys}
(\os{\circ}{X}_{m,T_0}/\os{\circ}{T}))(-j-k'-m,u)_{\mab Q}  \\
& \Lo 
R^qf_{X_{\bul \leq N,T_0}/T*}
(P_kA_{{\rm zar},{\mab Q}}(X_{\bul \leq N,\os{\circ}{T}_0}/S(T)^{\nat},
E^{\bul \leq N})\quad (k'\leq k).   
\end{align*}  
When $T$ is restrictively hollow with respective to the morphism $T_0\lo S$, 
then we have the following spectral sequence 
by (\ref{eqn:espssp}) and  (\ref{lemm:flpis}): 

\begin{align*} 
& E_1^{-k,q+k} = \bigoplus_{m=0}^N
\bigoplus_{j\geq \max \{-(k+m),0\}} 
R^{q-2j-k-2m}
f_{\os{\circ}{X}{}^{(2j+k+m)}_{m,T}/\os{\circ}{T}*}
(E_{\os{\circ}{X}{}^{(2j+k+m)}_{m,T_0}/\os{\circ}{T}}
\otimes_{\mab Z} 
\tag{3.4.9.7}\label{eqn:emtttsp} \\
& \vp^{(2j+k+m)}_{\rm crys}
(\os{\circ}{X}_{m,T_0}/\os{\circ}{T}))(-j-k-m,u)_{\mab Q} 
\Lo 
R^qf_{X_{\bul \leq N,\os{\circ}{T}_0}/T*}
(\eps^*_{X_{\bul \leq N,T_0}/T}(E^{\bul \leq N}))_{\mab Q} 
\quad (q\in {\mab Z}).  
\end{align*}  

\par 
When $\os{\circ}{S}$ is of characteristic $p>0$, when 
$u\col (S(T)^{\nat},{\cal J},\del)\lo (S(T)^{\nat},{\cal J},\del)$ is a lift of $F_{S_{\os{\circ}{T}_0}}$, 
when $g_{\bul \leq N}$ is the absolute Frobenius endomorphism 
of $X_{\bul \leq N,\os{\circ}{T}_0}$ or 
the relative Frobenius morphism of 
$X_{\bul \leq N,\os{\circ}{T}_0} \lo X^{[p]}_{\bul \leq N,\os{\circ}{T}_0}$ and when 
$E^{\bul \leq N}={\cal O}_{\os{\circ}{X}_{\bul \leq N,T_0}/\os{\circ}{T}}$, 
we denote $(-j-k-m,u)$ by $(-j-k-m)$ as usual.

\begin{defi} 
(1) We call the spectral sequence (\ref{eqn:espssp}) 
(resp.~(\ref{eqn:emtttsp}))
the {\it Poincar\'{e} spectral sequence} of 
$\eps^*_{X_{\bul \leq N,\os{\circ}{T}_0}/S(T)^{\nat}}(E^{\bul \leq N})$ 
on $X_{\bul \leq N,\os{\circ}{T}_0}/S(T)^{\nat}$ (resp.~
$\eps^*_{X_{\bul \leq N,T_0}/T}(E^{\bul \leq N})$ on 
$X_{\bul \leq N,T_0}/T$). 
When 
$E^{\bul \leq N}={\cal O}_{\os{\circ}{X}_{\bul \leq N,T_0}/\os{\circ}{T}}$, 
we call (\ref{eqn:espssp}) (resp.~(\ref{eqn:emtttsp})) 
the {\it weight spectral sequence} of 
$X_{\bul \leq N,\os{\circ}{T}_0}/S(T)^{\nat}$ 
(resp.~$X_{\bul \leq N,\os{\circ}{T}_0}/T$).   
\par 
(2) We usually denote by 
$P$ (instead of $\del(L,P)$)  
the induced filtration on 
$$R^qf_{X_{\bul \leq N,T_0}/T*}
(\eps^*_{X_{\bul \leq N,T_0}/T}(E^{\bul \leq N}))_{\mab Q}$$
by the spectral sequence (\ref{eqn:espssp})
by abuse of notation. We call $P$ 
the {\it Poincar\'{e} filtration} on 
$R^qf_{X_{\bul \leq N,T_0}/T*}
(\eps^*_{X_{\bul \leq N,T_0}/T}
(E^{\bul \leq N}))_{\mab Q}$. 
If $E^{\bul \leq N}={\cal O}_{\os{\circ}{X}_{\bul \leq N,T_0}/\os{\circ}{T}}$, 
then we call $P$ the {\it weight filtration} on 
$R^qf_{X_{\bul \leq N,\os{\circ}{T}_0}/S(T)^{\nat}*}
({\cal O}_{X_{\bul \leq N,\os{\circ}{T}_0}/S(T)^{\nat}})_{\mab Q}$ 
(resp.~$R^qf_{X_{\bul \leq N,T_0}/T*}({\cal O}_{X_{\bul \leq N,T_0}/T})_{\mab Q}$).
\end{defi}

\par


\begin{prop}[{\bf Contravariant functoriality}]\label{prop:cosuntwt}
Let $g_{\bul \leq N}$ be as in {\rm (\ref{prop:grlgoc})}. 
Let $\os{\circ}{g}{}^{(k)*}_{T'T}$ be the following  morphism$:$
\begin{align*}  
\os{\circ}{g}{}^{(k)*}_{T'T}
& :=\sum_{\ul{\lam}\in \Lam^{(k)}(g)}
\os{\circ}{g}{}^*_{\ul{\lam},T'T}
\col 
R^qf_{\os{\circ}{Y}{}^{(k)}_{m,T'_0}
/\os{\circ}{T}{}'*}
(F_{\os{\circ}{Y}{}^{(k)}_{m,T'_0}/\os{\circ}{T}{}'}
\otimes_{\mab Z}
\vp^{(k)}_{\rm crys}(\os{\circ}{Y}_{m,T'_0}
/\os{\circ}{T}{}')) \tag{3.4.11.1}\label{eqn:rgnlm} \\
& \lo 
R^qf_{\os{\circ}{X}{}^{(k)}_{m,T_0}/\os{\circ}{T}{}*}
(E_{\os{\circ}{X}{}^{(k)}_{m,T_0}
/\os{\circ}{T}}
\otimes_{\mab Z}
\vp^{(k)}_{\rm crys}(\os{\circ}{X}_{m,T_0}/T)). 
\end{align*} 
Then there exists the following morphism of 
the following spectral sequences: 
\begin{align*} 
E_1^{-k,q+k}& = \bigoplus_{m=0}^N
\bigoplus_{j\geq \max \{-(k+m),0\}} 
R^{q-2j-k-2m}
f_{\os{\circ}{X}{}^{(2j+k+m)}_{m,T}/\os{\circ}{T}*}
(E_{\os{\circ}{X}{}^{(2j+k+m)}_{m,T_0}/\os{\circ}{T}}
\otimes_{\mab Z} 
\tag{3.4.11.2}\label{eqn:espsfsp} \\
& \vp^{(2j+k+m)}_{\rm crys}
(\os{\circ}{X}_{m,T_0}
/\os{\circ}{T}))(-j-k-m,u) \Lo  
R^qf_{X_{\bul \leq N,\os{\circ}{T}_0}/S(T)^{\nat}*}
(\eps^*_{X_{\bul \leq N,\os{\circ}{T}_0}/S(T)^{\nat}}(E^{\bul \leq N}))
\end{align*}
\begin{equation*} 
\begin{CD} 
{\quad \quad \quad \quad \quad \quad \quad \quad \quad} 
@. 
{\quad \quad \quad \quad \quad \quad}\\ 
@A{\bigoplus_{m\geq 0}
\bigoplus_{j\geq \max \{-(k+m),0\}}\os{\circ}{g}{}^{(2j+k+m)*}}AA  
@. @. @AA{g^*_{\bul \leq N}}A \\
\end{CD} 
\end{equation*} 
\begin{align*} 
E_1^{-k,q+k}& 
= \bigoplus_{m=0}^N
\bigoplus_{j\geq \max \{-(k+m),0\}} 
R^{q-2j-k-2m}
f_{\os{\circ}{Y}^{(2j+k+m)}_{m,T'_0}/\os{\circ}{T}{}'}
(F_{\os{\circ}{Y}^{(2j+k+m)}_{m,T'_0/\os{\circ}{T}{}'}}
\otimes_{\mab Z}  \\
& \vp^{(2j+k+m)}_{\rm crys}(\os{\circ}{Y}_{m,T'_0}/\os{\circ}{T}{}'))(-j-k-m,u) 
\Lo 
R^qf_{Y_{\bul \leq N,\os{\circ}{T}{}'_0}/S'(T')^{\nat}}
(\eps^*_{Y_{\bul \leq N,\os{\circ}{T}{}'_0}/S'(T')^{\nat}}
(F^{\bul \leq N})). 
\end{align*}
\end{prop}
\begin{proof} 
The proof is the same as that of (\ref{prop:contwt}). 
\end{proof}

\par
Fix a total order on $\Lam$ once and for all. 
Let 
\begin{align*}
G_m 
\col 
&\bigoplus_{j\geq \max \{-(k+m),0\}}  
R^{q-2j-k-2m}
f_{\os{\circ}{X}{}^{(2j+k+m)}_{m,T_0}
/\os{\circ}{T}*}
(E_{\os{\circ}{X}{}^{(2j+k+m)}_{m,T_0}
/\os{\circ}{T}} 
\tag{3.4.11.3}\label{eqn:tossn}\\ 
& 
\otimes_{\mab Z}
\vp^{(2j+k+m)}_{\rm crys}(\os{\circ}{X}_{m,T_0}
/\os{\circ}{T}))(-j-k-m,u)_{\mab Q}\\ 
& \lo \bigoplus_{j\geq \max \{-(k+m)+1,0\}} 
R^{q-2j-k-2m+2}f_{\os{\circ}{X}{}^{(2j+k+m-1)}_{m,T_0}
/\os{\circ}{T}*}
(E_{\os{\circ}{X}{}^{(2j+k+m-1)}_{m,T_0}/\os{\circ}{T}}
\otimes_{\mab Z}\\
& \vp^{(2j+k+m-1)}_{\rm crys}
(\os{\circ}{X}_{m,T_0}/\os{\circ}{T}))
(-j-k-m+1,u)_{\mab Q}
\end{align*}
be the obvious analogue of (\ref{eqn:togsn}). 

\par 
Let 
\begin{align*}
\rho_m
\col 
&\bigoplus_{j\geq \max \{-(k+m),0\}} 
R^{q-2j-k-2m}
f_{\os{\circ}{X}{}^{(2j+k+m)}_{m,T_0}
/\os{\circ}{T}*}
(E_{\os{\circ}{X}{}^{(2j+k+m)}_{m,T_0}
/\os{\circ}{T}}
\otimes_{\mab Z} 
\tag{3.4.11.4}\label{eqn:rhsgsn}\\
& 
\vp^{(2j+k+m)}_{\rm crys}
(\os{\circ}{X}_{m,T_0}/\os{\circ}{T}))
(-j-k-m,u)_{\mab Q}  \\ 
&\lo \bigoplus_{j\geq \max \{-(k+m),0\}} 
R^{q-2j-k-2m}f_{\os{\circ}{X}{}^{(2j+k+m+1)}_{m,T_0}
/\os{\circ}{T}*}
(E_{\os{\circ}{X}{}^{(2j+k+m+1)}_{m,T_0}
/\os{\circ}{T}}
\otimes_{\mab Z} \\ 
&\vp^{(2j+k+m+1)}_{\rm crys}
(\os{\circ}{X}_{m,T_0}/\os{\circ}{T}))(-j-k-m,u)_{\mab Q}.
\end{align*}
be the obvious analogue of (\ref{eqn:rhogsn}). 
(However note that, 
on each $\os{\circ}{N}{}^{(2j+k+m)}_m$, 
the $\rho_m$ in (\ref{eqn:rhsgsn}) 
is defined by the replacement of 
$(-1)^{\bet}
\iota_{\ul{\lam}_m{\rm crys}}^{\ul{\lam}_{m\bet*}}$ by 
$\bigoplus_{i(m)}^{l(m)}d_{i(m)}(-1)^{\bet}
\iota_{\ul{\lam}_m{\rm crys}}^{\ul{\lam}_{m\bet*}}$ 
in the morphism (\ref{eqn:rhogsn}).  

\begin{prop}\label{prop:exbd}
Let the notations be as above. Then 
the obvious analogue of {\rm (\ref{prop:deccbd})} 
holds. 
\end{prop}
\begin{proof}
The proof is the same as that of (\ref{prop:deccbd}).  
\end{proof}

\parno
\section{Filtered base change theorem II}\label{sec:bcakf}
Let the notations be as in the previous section. 
In this section we prove the filtered base change theorem of 
$(A_{{\rm zar},{\mab Q}}(X_{\bul \leq N,\os{\circ}{T}_0}/S(T)^{\nat},E^{\bul \leq N}),P)$. 
This is a nontrivial generalization of (\ref{theo:bccange}) 
modulo torsion even in the situation in \S\ref{sec:bckf}
because we allow a ramified base change 
morphism of log points in this section, 
while we have not allowed in (\ref{theo:bccange}).  

\begin{prop}\label{prop:bbtdc}  
Assume that $\os{\circ}{T}$ is quasi-compact and that $\os{\circ}{f} 
\col \os{\circ}{X}_{\bul \leq N,T_0}\lo \os{\circ}{T}_0$ 
is quasi-compact and quasi-separated. 
Then 
$Rf_*((A_{{\rm zar},{\mab Q}}(X_{\bul \leq N,\os{\circ}{T}_0}/S(T)^{\nat},E^{\bul \leq N}),P))$ 
is isomorphic to a bounded filtered complex of ${\cal O}_T$-modules. 
\end{prop}
\begin{proof}
By using the spectral sequence (\ref{eqn:ekpssp}) instead of 
(\ref{eqn:escssp}), the proof is the same as that of 
(\ref{prop:bdccd}).   
\end{proof}

\begin{theo}[{\bf Log base change theorem of 
$(A_{{\rm zar},{\mab Q}},P)$}]\label{theo:bcfqa}
Let the assumptions be as in {\rm (\ref{prop:bbtdc})}.  
Let $S'\lo S$ be a morphism of $p$-adic formal families of log points. 
Let $u\col (T',{\cal J}',\del') \lo (T,{\cal J},\del)$ 
be a morphism of log PD-enlargements over the morphism $S'\lo S$. 
Set $T'_0:=\ul{\rm Spec}^{\log}_{T'}({\cal O}_{T'}/{\cal J}')$  
and $Y_{\bul \leq N}:=X_{\bul \leq N}\times_SS'$. 
Let 
$f' \col Y_{\bul \leq N,\os{\circ}{T}{}'_0}=
Y_{\bul \leq N}\times_{S'}S'_{\os{\circ}{T}{}'_0} \lo S'(T')^{\nat}$ 
be the base change morphism of $f$  
by the morphism $S'(T')^{\nat}\lo S(T)^{\nat}$.  
Let $g \col Y_{\bul \leq N,\os{\circ}{T}{}'_0} \lo X_{\bul \leq N,\os{\circ}{T}_0}$ 
be the induced morphism by $u$. 
Then there exists 
the following canonical filtered isomorphism
\begin{equation*}
Lu^*Rf_*((A_{{\rm zar},{\mab Q}}(X_{\bul \leq N,\os{\circ}{T}_0}/S(T)^{\nat},
E^{\bul \leq N}),P)) 
\os{\sim}{\lo} Rf'_*
((A_{{\rm zar},{\mab Q}}(Y_{\bul \leq N,\os{\circ}{T}{}'_0}/S'(T')^{\nat},
\os{\circ}{g}{}^*(E^{\bul \leq N})),P))
\tag{3.5.2.1}\label{eqn:blqxnpw}
\end{equation*}
in ${\rm DF}(f^{-1}_{\bul \leq N,T'}({\cal K}_{T'}))$, where 
${\cal K}_{T'}:={\cal O}_{T'}\otimes_{\mab Z}{\mab Q}$.  
\end{theo}
\begin{proof}  
Note that $X_{\bul \leq N,\os{\circ}{T}{}'_0}/S'_{\os{\circ}{T}{}'_0}$ 
is an $N$-truncated simplicial base change of SNCL schemes 
((\ref{rema:ibs}))
and that the base change morphism 
$Y_{\bul \leq N,\bul,\os{\circ}{T}{}'_0}\os{\sus}{\lo} \ol{\cal P}_{\bul \leq N,\bul,\ol{S'(T')^{\nat}}}$ 
is an admissible immersion ((\ref{exem:xpm})). 
We have the natural morphisms 
$\ol{\cal P}_{\bul \leq N,\bul,\ol{S'(T')^{\nat}}}
:=\ol{\cal P}_{\bul \leq N,\bul}\times_{\ol{S(T)^{\nat}}}\ol{S(T')^{\nat}}
\lo \ol{\cal P}_{\bul \leq N,\bul}$ 
and 
$\ol{\mathfrak D}_{\bul \leq N,\bul,\ol{S(T')^{\nat}}} 
:=\ol{\mathfrak D}_{\bul \leq N,\bul}\times_{\ol{S(T)^{\nat}}}\ol{S(T')^{\nat}}
\lo \ol{\mathfrak D}_{\bul \leq N,\bul}$.  
Hence we have the following natural morphism  
\begin{equation*} 
(A_{{\rm zar},{\mab Q}}
(X_{\bul \leq N,\os{\circ}{T}_0}/S(T)^{\nat},E^{\bul \leq N}),P) 
\lo 
Rg_*((A_{{\rm zar},{\mab Q}}(X'_{\bul \leq N,\os{\circ}{T}{}'_0}/S'(T')^{\nat},
\os{\circ}{g}{}^*(E^{\bul \leq N})),P)). 
\tag{3.5.2.2}\label{eqn:bcxega}
\end{equation*} 
By applying $Rf_*$ to (\ref{eqn:bcxega}) and using 
the adjoint property of $L$ and $R$ (\cite[(1.2.2)]{nh2}), 
we have the natural morphism (\ref{eqn:blqxnpw}). 
Here we have used the boundedness in 
(\ref{prop:bbtdc}) for the well-definedness of $Lu^*$. 
By using (\ref{coro:grsem}), 
the rest of the proof is the same as that of (\ref{theo:bccange}). 
\end{proof} 

\begin{rema}\label{rema:gn7}
The analogue of (\ref{coro:connfil}) for the filtered complex 
$$(A_{{\rm zar},{\mab Q}}(X_{\bul \leq N,\os{\circ}{T}_0}/S(T)^{\nat},
E^{\bul \leq N}),P)$$ 
holds. 
\end{rema}

\begin{coro}\label{coro:fcqtd}
Let the notations and the assumptions be as in $(\ref{prop:bbtdc})$. 
Then 
$$Rf_*(P_kA_{{\rm zar},{\mab Q}}(X_{\bul \leq N,\os{\circ}{T}_0}/S(T)^{\nat},
E^{\bul \leq N})) \quad (k \in{\mab N})$$
has finite tor-dimension. 
Moreover, if $\os{\circ}{T}$ is noetherian and 
if $\os{\circ}{f}$ is proper,
then $Rf_*(P_k(A_{{\rm zar},{\mab Q}}(X_{\bul \leq N,\os{\circ}{T}_0}/S(T)^{\nat},
E^{\bul \leq N})))$ 
is a perfect complex of ${\cal K}_T$-modules.
\end{coro}

Using \cite[(2.10.10)]{nh2}, 
we have the following corollary 
(cf.~\cite[(2.10.11)]{nh2}): 

\begin{coro}\label{coro:filqerf}
Let the notations 
and the assumptions be as in $(\ref{coro:fcqtd})$.
Then the filtered complex 
$Rf_*((P_kA_{{\rm zar},{\mab Q}}(X_{\bul \leq N,\os{\circ}{T}_0}/S(T)^{\nat},
E^{\bul \leq N}),P))$ 
is a {\it filtered perfect}
complex of ${\cal K}_T$-modules, that is, 
locally on $T_{\rm zar}$, filteredly quasi-isomorphic to 
a filtered strictly perfect complex {\rm (\cite[(2.10.8)]{nh2})}.
\end{coro}
\begin{proof}
(\ref{coro:filqerf}) immediately follows from (\ref{coro:fcqtd}) and 
\cite[(2.10.10)]{nh2}.
\end{proof}

\begin{rema}\label{rema:iff}
(\ref{coro:filqerf}) will be a key ingredient 
for the proof of the log convergence 
of the weight filtration on the log isocrystalline cohomology sheaves.  
See \S\ref{sec:lcw} below, especially the proof of 
(\ref{lemm:pnlcfi}) below.
\end{rema}

\section{$p$-adic monodromy operators II and  
$P$-filtered log isocrystalline complexes}\label{sec:mndomp}
Let the notations be as in \S\ref{sec:pfnsc}.  
Assume that $X_{\bul \leq N,\os{\circ}{T}_0}$ has 
an affine $N$-truncated simplicial open covering of $X_{\bul \leq N,\os{\circ}{T}_0}$ 
and that there exists an admissible immersion 
$X_{\bul \leq N,\bul,\os{\circ}{T}_0}\os{\sus}{\lo} \ol{\cal P}_{\bul \leq N,\bul}$ over $\ol{S(T)^{\nat}}$. 
Let $f_{\bul \leq N,\bul}\col 
X_{\bul \leq N,\bul,\os{\circ}{T}_0}\lo S_{\os{\circ}{T}_0}\os{\sus}{\lo} S(T)^{\nat}$ 
be the composite structural morphism.   
In this section we give analogous results to those in \S\ref{sec:vpmn} for 
$X_{\bul \leq N}/S$ in this section. Because the proofs for the 
analogous results to those in \S\ref{sec:vpmn} are 
the same as the proofs in \S\ref{sec:vpmn}, we omit the proofs. 
\par 
First we define the $p$-adic monodromy operator and 
the $p$-adic quasi-monodromy operator
for $X_{\bul \leq N}/S$, a log $p$-adic PD-enlargement $(T,{\cal J},\del)$ of $S$ 
and a flat quasi-coherent log crystal $\ol{F}{}^{\bul \leq N}$ of 
${\cal O}_{X_{\bul \leq N,T_0}/\os{\circ}{T}}$-modules. 
\par  
Let $\ol{F}{}^{\bul \leq N}$ be a flat quasi-coherent log crystal of 
${\cal O}_{X_{\bul \leq N,\os{\circ}{T}_0}/\os{\circ}{T}}$-modules. 
Set $F^{\bul \leq N}:=
\eps^*_{X_{\bul \leq N,\os{\circ}{T}_0/S(T)^{\nat}}}(\ol{F}{}^{\bul \leq N})$. 
Let $({\cal F}^{\bul \leq N,\bul},\nabla)$ be before in (\ref{eqn:nidfopltd}) 
for $\ol{F}{}^{\bul \leq N}$ in the case $Y_{\bul \leq N}=X_{\bul \leq N}$. 
Then the following sequence 
\begin{align*} 
0 & \lo {\cal F}^{\bul \leq N,\bul}
\otimes_{{\cal O}_{{\cal P}^{\rm ex}_{\bul \leq N,\bul}}}
{\Om}^{\bul}_{{\cal P}^{\rm ex}_{\bul \leq N,\bul}/S(T)^{\nat}}
\otimes_{\mab Z}{\mab Q}[-1] 
\os{e\theta_{{\cal P}^{\rm ex}_{\bul\leq N,\bul} \wedge}}{\lo} 
\tag{3.6.0.1}\label{eqn:gnpsd}\\ 
&{\cal F}^{\bul \leq N,\bul}
\otimes_{{\cal O}_{{\cal P}^{\rm ex}_{\bul \leq N,\bul}}}
{\Om}^{\bul}_{{\cal P}^{\rm ex}_{\bul \leq N,\bul}/\os{\circ}{T}}
\otimes_{\mab Z}{\mab Q} 
\lo  {\cal F}^{\bul \leq N,\bul}
\otimes_{{\cal O}_{{\cal P}^{\rm ex}_{\bul \leq N,\bul}}}
{\Om}^{\bul}_{{\cal P}^{\rm ex}_{\bul \leq N,\bul}/S(T)^{\nat}}
\otimes_{\mab Z}{\mab Q} \lo 0
\end{align*} 
is exact. 
Let 
\begin{equation*} 
{\cal F}^{\bul \leq N,\bul}
\otimes_{{\cal O}_{{\cal P}^{\rm ex}_{\bul \leq N,\bul}}}
{\Om}^{\bul}_{{\cal P}^{\rm ex}_{\bul \leq N,\bul}/S(T)^{\nat}}
\otimes_{\mab Z}{\mab Q}
\lo 
{\cal F}^{\bul \leq N,\bul}
\otimes_{{\cal O}_{{\cal P}^{\rm ex}_{\bul \leq N,\bul}}}
{\Om}^{\bul}_{{\cal P}^{\rm ex}_{\bul \leq N,\bul}/S(T)^{\nat}}
\otimes_{\mab Z}{\mab Q} 
\tag{3.6.0.2}\label{eqn:gypsxd}
\end{equation*} 
be the boundary morphism of (\ref{eqn:gnpsd}) 
in the derived category 
${\rm D}^+(f^{-1}_{\bul \leq N,\bul}({\cal O}_T))$. 
(Recall that $e$ is the degree function of the morphism 
$S\lo S_0$.)
This morphism induces the following morphism by (\ref{ali:ptc}): 
\begin{equation*}
N_{S(T)^{\nat},{\rm zar}} \col Ru_{X_{\bul \leq N,\os{\circ}{T}_0}/S(T)^{\nat}*}
(F^{\bul \leq N})\otimes^L_{\mab Z}{\mab Q} \lo 
Ru_{X_{\bul \leq N,\os{\circ}{T}_0}/S(T)^{\nat}*}(F^{\bul \leq N})
\otimes^L_{\mab Z}{\mab Q}.
\tag{3.6.0.3}\label{eqn:nypne}
\end{equation*}

\begin{defi}\label{defi:mdop}
We call the morphism (\ref{eqn:nypne}) the 
{\it  zariskian monodromy operator} of $\ol{F}{}^{\bul \leq N}$ (or 
$F^{\bul \leq N}$ by abuse of terminology). 
When $\ol{F}^{\bul \leq N}={\cal O}_{X_{\bul \leq N,\os{\circ}{T}_0}/\os{\circ}{T}}$, 
we call the morphism (\ref{eqn:nypne}) the 
{\it  zariskian monodromy operator} 
of $X_{\bul \leq N,\os{\circ}{T}_0}/S(T)^{\nat}$ for simplicity. 
\end{defi}

The following holds as in (\ref{prop:otj})

\begin{prop}[{\bf Contravariant functoriality of 
the $p$-adic monodromy operator}]\label{prop:ctrbmon} 
Assume that we are given a commutative diagram {\rm (\ref{eqn:ydgguss})}. 
Then the following diagram is commutative$:$
\begin{equation*} 
\begin{CD} 
Rg_*Ru_{Y_{\bul \leq N,\os{\circ}{T}_0}/S(T)^{\nat}*}(F^{\bul \leq N})
\otimes^L_{\mab Z}{\mab Q} 
@>{{\rm deg}(u)N_{S(T)^{\nat},{\rm zar}}}>>
Rg_*Ru_{Y_{\bul \leq N,\os{\circ}{T}_0}/S(T)^{\nat}*}(F^{\bul \leq N})
\otimes^L_{\mab Z}{\mab Q}(-1,u)\\
@A{z^*}AA @AA{z^*}A\\
Ru_{Z_{\bul \leq N,\os{\circ}{T}_0}/S(T)^{\nat}*}(G^{\bul \leq N})
\otimes^L_{\mab Z}{\mab Q} 
@>{N_{S(T)^{\nat},{\rm zar}}}>>
Ru_{Z_{\bul \leq N,\os{\circ}{T}_0}/S(T)^{\nat}*}(G^{\bul \leq N})(-1,u')
\otimes^L_{\mab Z}{\mab Q}
\end{CD}
\tag{3.6.2.1}\label{eqn:ydhdagguss}
\end{equation*}
\end{prop}
\begin{proof} 
This is only the formal version of (\ref{cd:xaxy}). 
\end{proof}


\par
In the rest of this section, we give analogues of 
the results in \S\ref{sec:vpmn} for 
$$(A_{{\rm zar},{\mab Q}}
(X_{\bul \leq N,\os{\circ}{T}_0}/S(T)^{\nat},E^{\bul \leq N}),P).$$ 
Because the proof of these results are the same as those of 
the results in \S\ref{sec:vpmn}, we omit the proof. 

\par 
Let $E^{\bul \leq N}$ be as in \S\ref{sec:pgensc}.  
Let 
\begin{equation*} 
(A_{{\rm zar},{\mab Q}}({\cal P}^{\rm ex}_{\bul \leq N,\bul}/S(T)^{\nat}, 
{\cal E}^{\bul \leq N,\bul})^{\bul \bul},P)
\lo (I^{\bul \leq N,\bul \bul \bul},P)
\tag{3.6.2.2}\label{eqn:aipfi}
\end{equation*} 
be the $(N,\infty)$-truncated bicosimplicial filtered Godement resolution. 
Then we have the following resolution 
\begin{align*} 
(s(A_{{\rm zar},{\mab Q}}({\cal P}^{\rm ex}_{\bul \leq N,\bul}/S(T)^{\nat}, 
{\cal E}^{\bul \leq N,\bul})^{\bul \bul}),P)
\lo  
(s(I^{\bul \leq N,\bul \bul \bul}),P).
\end{align*}  
Set $(A^{\bul \leq N,\bul \bul},P):=
\pi_{{\rm zar}*}((I^{\bul \leq N,\bul \bul \bul},P))$. 
Then 
$$s((A^{\bul \leq N,\bul \bul},P))
=(A_{{\rm zar},{\mab Q}}(X_{\bul \leq N,\os{\circ}{T}_0}/S(T)^{\nat},
E^{\bul \leq N}),P)$$ 
in ${\rm D}^+{\rm F}(f^{-1}_{\bul \leq N}({\cal O}_T))$.  
Consider the following morphism 
\begin{align*} 
& \nu_{\rm zar}({\cal P}^{\rm ex}_{\bul \leq N,\bul}/S(T)^{\nat},
{\cal E}^{\bul \leq N,\bul})^{ij}:={\rm proj}\col \tag{3.6.2.3}\label{eqn:reislbb} \\
&({\cal E}^{\bul \leq N,\bul}
\otimes_{{\cal O}_{{\cal P}^{\rm ex}_{\bul \leq N,\bul}}}
{\Om}^{i+j+1}_{{\cal P}^{\rm ex}_{\bul \leq N,\bul}
/\os{\circ}{T}})\otimes_{\mab Z}{\mab Q}
/P_{j+1}  
\lo({\cal E}^{\bul \leq N,\bul}
\otimes_{{\cal O}_{{\cal P}^{\rm ex}_{\bul \leq N,\bul}}}
{\Om}^{i+j+1}_{{\cal P}^{\rm ex}_{\bul \leq N,\bul}
/\os{\circ}{T}})\otimes_{\mab Z}{\mab Q}/P_{j+2}.   
\end{align*}  
This morphism induces a morphism of complexes 
with the same boundary morphisms as those in (\ref{cd:locstbd}).  
The morphism (\ref{eqn:reislbb}) induces the following morphism 
\begin{equation*} 
\wt{\nu}^{ij}_{{\rm zar}}:={\rm proj}. \col 
A^{\bul \leq N,ij}\lo A^{\bul \leq N,i-1,j+1}
\end{equation*}  
of  sheaves of $f^{-1}_{\bul \leq N}({\cal O}_T)$-modules. 
Set 
\begin{align*} 
\wt{\nu}_{{\rm zar}}:=
s(\oplus_{i,j\in {\mab N}}\wt{\nu}^{ij}_{{\rm zar}}).
\tag{3.6.2.4}\label{ali:nitzd}
\end{align*}  
Let 
\begin{equation*} 
\nu_{S(T)^{\nat},{\rm zar}} \col 
(A_{{\rm zar},{\mab Q}}(X_{\bul \leq N,\os{\circ}{T}_0}/S(T)^{\nat},
E^{\bul \leq N}),P)
\lo 
(A_{{\rm zar},{\mab Q}}(X_{\bul \leq N,\os{\circ}{T}_0}/S(T)^{\nat},
E^{\bul \leq N}),P\langle -2\rangle) 
\tag{3.6.2.5}\label{eqn:naixgdxdn}
\end{equation*} 
be a morphism of complexes induced by 
$\{\nu^{ij}_{{\rm zar}}\}_{i,j \in {\mab N}}$. 
The morphism 
$$\theta_{{\cal P}^{\rm ex}_{\bul \leq N,\bul}} \wedge \col 
(A_{{\rm zar},{\mab Q}}({\cal P}^{\rm ex}_{\bul \leq N,\bul}/S(T)^{\nat}, 
{\cal E}^{\bul \leq N,\bul})^{ij},P)
\lo 
(A_{{\rm zar},{\mab Q}}({\cal P}^{\rm ex}_{\bul \leq N,\bul}/S(T)^{\nat}, 
{\cal E}^{\bul \leq N,\bul})^{i,j+1},P)$$  
induces a morphism 
$$\theta_{{\cal P}^{\rm ex}_{\bul \leq N,\bul}} \wedge 
\col A^{\bul \leq N,\bul ij}\lo 
A^{\bul \leq N,\bul i,j+1}.$$ 
\par 
Let $g_{\bul \leq N}$ be the morphism 
in the beginning of \S\ref{sec:pfnsc} for the case 
$Y_{\bul \leq N,\os{\circ}{T}{}'_0}= X_{\bul \leq N,\os{\circ}{T}_0}$ 
and $(T',{\cal J}',\del')=(T,{\cal J},\del)$ and $S'=S$ 
fitting into the following commutative diagram: 
\begin{equation*} 
\begin{CD} 
X'_{\bul \leq N,\os{\circ}{T}_0} @>{g'_{\bul \leq N}}>> X''_{\bul \leq N,\os{\circ}{T}_0} \\
@VVV @VVV \\ 
X_{\bul \leq N,\os{\circ}{T}_0} @>{g_{\bul \leq N}}>> X_{\bul \leq N,\os{\circ}{T}_0}\\
@VVV @VVV \\ 
S_{\os{\circ}{T}_0} @>>> S_{\os{\circ}{T}{}_0} \\ 
@V{\bigcap}VV @VV{\bigcap}V \\ 
S(T)^{\nat} @>{u}>> S(T)^{\nat},  
\end{CD}
\tag{3.6.2.6}\label{cd:xgsqpxy}
\end{equation*}
where $X''_{\bul \leq N,\os{\circ}{T}{}'_0}$ 
is another disjoint union 
of the member of an affine $N$-truncated simplicial open covering of 
$X_{\bul \leq N,\os{\circ}{T}_0}$. 
The morphism (\ref{eqn:naixgdxdn}) is 
the following morphism 
\begin{align*} 
\nu_{S(T)^{\nat},{\rm zar}} \col &
(A_{{\rm zar},{\mab Q}}(X_{\bul \leq N,\os{\circ}{T}_0}/S(T)^{\nat}
,E^{\bul \leq N}),P) \tag{3.6.2.7}\label{eqn:axdaxdn}\\
& \lo 
(A_{{\rm zar},{\mab Q}}(X_{\bul \leq N,\os{\circ}{T}_0}/S(T)^{\nat}
,E^{\bul \leq N}),P\langle -2\rangle)(-1,u). 
\end{align*}

As in (\ref{prop:cmzoqm}) we have the following:


\begin{prop}\label{prop:cmzaoqm}
The zariskian monodromy operator
\begin{align*} 
N_{S(T)^{\nat},{\rm zar}} \col &
Ru_{X_{\bul \leq N,\os{\circ}{T}_0}/S(T)^{\nat}*}
(\eps^*_{X_{\bul \leq N,\os{\circ}{T}_0}/S(T)^{\nat}}
(E^{\bul \leq N}))\otimes_{\mab Z}^L{\mab Q} 
\tag{3.6.3.1}\label{eqn:minabv}\\
& \lo  
Ru_{X_{\bul \leq N,\os{\circ}{T}_0}/S(T)^{\nat}*}
(\eps^*_{X_{\bul \leq N,\os{\circ}{T}_0}/S(T)^{\nat}}
(E^{\bul \leq N}))\otimes_{\mab Z}^L{\mab Q}(-1,u)
\end{align*}  
is equal to  
$$\nu_{{\rm zar}} \col 
A_{{\rm zar},{\mab Q}}(X_{\bul \leq N,\os{\circ}{T}_0}/S(T)^{\nat},E^{\bul \leq N}) 
\lo A_{{\rm zar},{\mab Q}}(X_{\bul \leq N,\os{\circ}{T}_0}/S(T)^{\nat},E^{\bul \leq N})(-1,u)$$ 
via the canonical isomorphism {\rm (\ref{ali:usz})}.  
\end{prop}

\begin{coro}\label{coro:nancilp} 
Assume that $\os{\circ}{X}_{\bul \leq N}$ is quasi-compact. 
Then the zariskian monodromy operator 
\begin{align*} 
N_{S(T)^{\nat},{\rm zar}} \col &Ru_{X_{\bul \leq N,T_0}/T*}
(\eps^*_{X_{\bul \leq N,T_0}/T}(E^{\bul \leq N}))\otimes_{\mab Z}^L{\mab Q}
\tag{3.6.4.1}\label{ali:madmycoh}\\
&\lo  
Ru_{X_{\bul \leq N,T_0}/T*}
(\eps^*_{X_{\bul \leq N,T_0}/T}(E^{\bul \leq N}))\otimes_{\mab Z}^L{\mab Q}(-1,u)
\end{align*} 
is nilpotent.
\end{coro}

\begin{coro}\label{prop:canmorp} 
Let $(T',{\cal J}',\del')\lo (T,{\cal J},\del)$ 
be as in {\rm (\ref{prop:bcmop})}. 
Then 
\begin{align*} 
\nu_{S(T)^{\nat},{\rm zar}}\otimes^L_{{\cal O}_T}{\cal O}_{T'} 
\col & 
A_{{\rm zar},{\mab Q}}(X_{\bul \leq N,\os{\circ}{T}_0}/S(T)^{\nat},E^{\bul \leq N})
\otimes^L_{{\cal O}_T}{\cal O}_{T'}  
\tag{3.6.5.1}\label{ali:macibv}\\
&\lo A_{{\rm zar},{\mab Q}}(X_{\bul \leq N,\os{\circ}{T}_0}/S(T)^{\nat},E^{\bul \leq N})(-1,u)
\otimes^L_{{\cal O}_T}{\cal O}_{T'} 
\end{align*} 
is equal to 
\begin{equation*} 
\nu_{S'(T')^{\nat},{\rm zar}} \col 
 A_{{\rm zar},{\mab Q}}(X_{\bul \leq N,\os{\circ}{T}{}'_0}/S'(T')^{\nat},E^{\bul \leq N})
\lo A_{{\rm zar},{\mab Q}}(X_{\bul \leq N,\os{\circ}{T}{}'_0}/S'(T')^{\nat},
E^{\bul \leq N})(-1,u).  
\end{equation*} 
\end{coro}

\begin{defi}
We call $\nu_{S(T)^{\nat},{\rm zar}}$ the {\it zariskian quasi-monodromy operator} of 
the sheaf $\eps^*_{X_{\bul \leq N,\os{\circ}{T}_0}/\os{\circ}{T}}(E^{\bul \leq N})$. 
When $E^{\bul \leq N}={\cal O}_{\os{\circ}{X}_{\bul \leq N,T_0}/\os{\circ}{T}}$, 
then we call $\nu_{S(T)^{\nat},{\rm zar}}$ the {\it zariskian quasi-monodromy operator} of 
$X_{\bul \leq N,\os{\circ}{T}_0}/S(T)^{\nat}$.
\end{defi}

\begin{prop}\label{prop:qauasimon}
Let $k$ be a positive integer and let $q$ be a nonnegative integer. 
\begin{align*} 
\nu_{S(T)^{\nat},{\rm zar}}^k \col  & 
R^qf_{T*}({\rm gr}_k^PA_{{\rm zar},{\mab Q}}(X_{\bul \leq N,\os{\circ}{T}_0}/S(T)^{\nat},E^{\bul \leq N})) 
\tag{3.6.7.1}\label{eqn:kapa}\\
& \lo 
R^qf_{T*}
({\rm gr}_{-k}^PA_{{\rm zar},{\mab Q}}(X_{\bul \leq N,\os{\circ}{T}_0}/S(T)^{\nat},E^{\bul \leq N}))(-k,u)
\end{align*} 
is the identity. 
\end{prop}


\begin{prop}\label{prop:nucsaa} 
The quasi-monodromy operator 
\begin{equation*} 
\nu_{S(T)^{\nat},{\rm zar}} \col 
A_{{\rm zar},{\mab Q}}(X_{\bul \leq N,\os{\circ}{T}_0}/S(T)^{\nat},E^{\bul \leq N})  \lo 
A_{{\rm zar},{\mab Q}}(X_{\bul \leq N,\os{\circ}{T}_0}/S(T)^{\nat},E^{\bul \leq N})(-1,u)
\end{equation*}
underlies the following morphism 
of filtered morphism 
\begin{align*} 
\nu_{S(T)^{\nat},{\rm zar}} \col &
(A_{{\rm zar},{\mab Q}}(X_{\bul \leq N,\os{\circ}{T}_0}/S(T)^{\nat},E^{\bul \leq N}),P)\lo  \\
&(A_{{\rm zar},{\mab Q}}(X_{\bul \leq N,\os{\circ}{T}_0}/S(T)^{\nat},E^{\bul \leq N}),
P\langle -2 \rangle)(-1,u), 
\end{align*}
where $P\langle -2 \rangle$ is a filtration 
defined by $(P\langle -2 \rangle)_k= P_{k-2}$. 
\end{prop}

\begin{coro}\label{coro:zamop}
The zariskian monodromy operator {\rm (\ref{ali:madmycoh})} 
induces the following morphism 
\begin{align*} 
N_{S(T)^{\nat},{\rm zar}}\col & 
P_kRu_{X_{\bul \leq N,\os{\circ}{T}_0}/S(T)^{\nat}*}
(\eps^*_{X_{\bul \leq N,\os{\circ}{T}_0}/S(T)^{\nat}}
(E^{\bul \leq N}))_{\mab Q}
\lo  \tag{3.6.9.1}\label{eqn:garppd}\\
& P_{k-2}Ru_{X_{\bul \leq N,\os{\circ}{T}_0}/S(T)^{\nat}*}
(\eps^*_{X_{\bul \leq N,\os{\circ}{T}_0}/S(T)^{\nat}}
(E^{\bul \leq N}))_{\mab Q}(-1,u). 
\end{align*}
\end{coro}

\par 

\par  
Set 
$$B^{\bul \leq N,\bul ij}
:=A^{\bul \leq N,\bul i-1,j}(-1,u)\oplus A^{\bul \leq N,\bul ij} 
\quad ( i,j\in {\mab N})$$
and 
$$B^{\bul \leq N,ij}
:=A^{\bul \leq N,i-1,j}(-1,u)\oplus A^{\bul \leq N,ij} 
\quad ( i,j\in {\mab N}).$$
The horizontal boundary morphism
$d' \col B^{\bul \leq N,ij} \lo B^{\bul \leq N,i+1,j}$ 
is, by definition, the induced morphism 
$d''{}^{\bul} \col B^{\bul \leq N,\bul ij} \lo B^{\bul \leq N,\bul, i+1,j}$ 
defined by the following formula 
\begin{align*} 
d'(\om_1,\om_2)=(\nabla \om_1,-\nabla\om_2) 
\tag{3.6.9.2}\label{ali:sddclxan}
\end{align*} 
and the vertical one 
$d'' \col B^{\bul \leq N,ij} \lo B^{\bul \leq N,i,j+1}$ is the induced morphism 
by a morphism  $d''{}^{\bul} \col B^{\bul \leq N,\bul ij} \lo B^{\bul \leq N,\bul i,j+1}$ 
defined by the following formula: 
\begin{align*} 
d''{}^{\bul}(\om_1,\om_2)=
(-e\theta_{{\cal P}^{\rm ex}_{\bul \leq N,\bul}}\wedge \om_1+
\nu_{S(T)^{\nat},{\rm zar}}(\om_2),e\theta_{{\cal P}^{\rm ex}_{\bul \leq N,\bul}}\wedge \om_2).
\tag{3.6.9.3}\label{ali:sddpxan}
\end{align*}  
It is easy to check that 
$B^{\bul \leq N,\bul \bul}$ is 
actually an $N$-truncated cosimplicial double complex. 
Let $B^{\bul \leq N,\bul}$ be the single complex of $B^{\bul \leq N,\bul \bul}$ 
with respect to the last two degrees. 
\par
Let 
\begin{align*} 
\mu_{{X_{\bul \leq N,\os{\circ}{T}_0}/\os{\circ}{T}}} \col 
\wt{R}u_{X_{\bul \leq N,\os{\circ}{T}_0}/\os{\circ}{T}*}
(\eps^*_{X_{\bul \leq N,\os{\circ}{T}_0}/\os{\circ}{T}}(E^{\bul \leq N}))_{\mab Q}
=&R\pi_{{\rm zar}*}({\cal E}^{\bul \leq N,\bul}
\otimes_{{\cal O}_{{\cal P}^{\rm ex}_{\bul \leq N,\bul}}}
{\Om}^{\bul}_{{\cal P}^{\rm ex}_{\bul \leq N,\bul}/\os{\circ}{T}}\otimes_{\mab Z}{\mab Q})
\tag{3.6.9.4}\label{ali:sgsclxan}\\
&\lo B^{\bul \leq N,\bul}
\end{align*} 
be a morphism
of complexes induced by the following morphisms
\begin{align*} 
\mu^i_{\bul \leq N,\bul} \col {\cal E}^{\bul \leq N,\bul}
\otimes_{{\cal O}_{{\cal P}^{\rm ex}_{\bul \leq N,\bul}}}
{\Om}^i_{{\cal P}^{\rm ex}_{\bul \leq N,\bul}/\os{\circ}{T}}\otimes_{\mab Z}{\mab Q} 
\lo  & ~{\cal E}^{\bul \leq N,\bul}
\otimes_{{\cal O}_{{\cal P}^{\rm ex}_{\bul \leq N,\bul}}}
{\Om}^i_{{\cal P}^{\rm ex}_{\bul \leq N,\bul}/\os{\circ}{T}}\otimes_{\mab Z}{\mab Q}/P_0 \oplus  \\
& {\cal E}^{\bul \leq N,\bul}
\otimes_{{\cal O}_{{\cal P}^{\rm ex}_{\bul \leq N,\bul}}}
{\Om}^{i+1}_{{\cal P}^{\rm ex}_{\bul \leq N,\bul}/\os{\circ}{T}}\otimes_{\mab Z}{\mab Q}/P_0 
\quad (i\in {\mab N})
\end{align*} 
defined by the following formula 
\begin{align*} 
\mu^i_{\bul \leq N,\bul}(\om):=(\om~{\rm mod}~P_0,
e\theta_{{\cal P}^{\rm ex}_{\bul \leq N,\bul}} 
\wedge \om~{\rm mod}~P_0) \quad 
(\om \in 
{\cal E}^{\bul \leq N,\bul}
\otimes_{{\cal O}_{{\cal P}^{\rm ex}_{\bul \leq N,\bul}}}
{\Om}^i_{{\cal P}^{\rm ex}_{\bul \leq N,\bul}/\os{\circ}{T}}\otimes_{\mab Z}{\mab Q}).
\end{align*} 
Then we have the following 
morphism of triangles:
\begin{equation*}
\begin{CD}
@>>> 
A_{{\rm zar},{\mab Q}}(X_{\bul \leq N,\os{\circ}{T}_0}/S(T)^{\nat}
,E^{\bul \leq N})[-1] @>{}>>\\
@. @A{(\theta_{X_{\bul \leq N,\os{\circ}{T}_0}/S(T)^{\nat}/S_0(T)^{\nat}}\wedge *)[-1]}AA  \\
@>>> 
R\pi_{{\rm zar}*}({\cal E}^{\bul \leq N,\bul}
\otimes_{{\cal O}_{{\cal P}^{\rm ex}_{\bul \leq N,\bul}}}
{\Om}^{\bul}_{{\cal P}^{\rm ex}_{\bul \leq N,\bul}/S(T)^{\nat}}\otimes_{\mab Z}{\mab Q})[-1]   
@>{R\pi_{{\rm zar}*}
(e\theta_{{\cal P}^{\rm ex}_{\bul \leq N,\bul}} \wedge)}>> 
\end{CD}
\tag{3.6.9.5}\label{cd:sgsctexan}
\end{equation*} 
\begin{equation*}
\begin{CD}
B^{\bul \leq N,\bul} @>>> \\
@A{\mu_{{X_{\bul \leq N,\os{\circ}{T}_0}/\os{\circ}{T}}}}AA \\
R\pi_{{\rm zar}*}({\cal E}^{\bul \leq N,\bul}
\otimes_{{\cal O}_{{\cal P}^{\rm ex}_{\bul \leq N,\bul}}}
{\Om}^{\bul}_{{\cal P}^{\rm ex}_{\bul \leq N,\bul}/\os{\circ}{T}}\otimes_{\mab Z}{\mab Q}) 
@>>> 
\end{CD}
\end{equation*} 
\begin{equation*}
\begin{CD}  
A_{{\rm zar},{\mab Q}}(X_{\bul \leq N,\os{\circ}{T}_0}/S(T)^{\nat},E^{\bul \leq N})
@>{+1}>>  \\  
@A{\theta_{X_{\bul \leq N,\os{\circ}{T}_0}/S(T)^{\nat}/S_0(T)^{\nat}} \wedge}AA \\ 
R\pi_{{\rm zar}*}({\cal E}^{\bul \leq N,\bul}
\otimes_{{\cal O}_{{\cal P}^{\rm ex}_{\bul \leq N,\bul}}}
{\Om}^{\bul}_{{\cal P}^{\rm ex}_{\bul \leq N,\bul}/S(T)^{\nat}}
\otimes_{\mab Z}{\mab Q}) 
@>{+1}>>.\\
\end{CD}
\end{equation*} 
This is nothing but the following diagram of triangles: 
\begin{equation*}
\begin{CD}
@>>> 
A_{{\rm zar},{\mab Q}}(X_{\bul \leq N,\os{\circ}{T}_0}/S(T)^{\nat}
,E^{\bul \leq N})[-1] @>{}>>\\
@. @A{(\theta_{X_{\bul \leq N,\os{\circ}{T}_0}/S(T)^{\nat}/S_0(T)^{\nat}}\wedge )[-1]}AA  \\
@>>> 
Ru_{X_{\bul \leq N,\os{\circ}{T}_0}/S(T)^{\nat}*}
(\eps^*_{X_{\bul \leq N,\os{\circ}{T}_0}/S(T)^{\nat}}(E^{\bul \leq N}))_{\mab Q}[-1]   
@>{}>> 
\end{CD}
\tag{3.6.9.6}\label{cd:monodam}
\end{equation*} 
\begin{equation*}
\begin{CD}
B^{\bul \leq N,\bul} @>>> \\
@A{\mu_{{X_{\bul \leq N,\os{\circ}{T}_0}/\os{\circ}{T}}}}AA \\
\wt{R}u_{X_{\bul \leq N,\os{\circ}{T}_0}/\os{\circ}{T}*}
(\eps^*_{X_{\bul \leq N,\os{\circ}{T}_0}/\os{\circ}{T}}(E^{\bul \leq N}))_{\mab Q} @>>> 
\end{CD}
\end{equation*} 
\begin{equation*}
\begin{CD}  
A_{{\rm zar},{\mab Q}}(X_{\bul \leq N,\os{\circ}{T}_0}/S(T)^{\nat},E^{\bul \leq N})
@>{+1}>>  \\  
@A{\theta_{X_{\bul \leq N,\os{\circ}{T}_0}/S(T)^{\nat}/S_0(T)^{\nat}} \wedge}AA \\ 
Ru_{X_{\bul \leq N,\os{\circ}{T}_0}/S(T)^{\nat}*}
(\eps^*_{X_{\bul \leq N,\os{\circ}{T}_0}/S(T)^{\nat}}(E^{\bul \leq N}))_{\mab Q}
@>{+1}>>.\\
\end{CD}
\end{equation*}

\begin{prop}\label{prop:spabqq}
Let $u\col (S(T)^{\nat},{\cal J},\del) \lo 
(S(T)^{\nat},{\cal J},\del)$ and  
$g_{\bul \leq N}\col 
X_{\bul \leq N,\os{\circ}{T}_0} \lo Y_{\bul \leq N,\os{\circ}{T}{}'_0}$ 
be as in the beginning of {\rm \S\ref{sec:pfnsc}}.   
In particular, assume that we are given the commutative diagrams 
{\rm (\ref{eqn:xdbquss})} and {\rm (\ref{cd:xdjqbxy})}.  
Let $u'\col (S'(T')^{\nat},{\cal J}',\del') \lo (S'(T')^{\nat},{\cal J}',\del')$ 
be an analogous morphism to $u$. 
Assume that we are given the commutative diagram 
{\rm (\ref{cd:xgsqpxy})} and a similar commutative diagram 
\begin{equation*} 
\begin{CD} 
Y'_{\bul \leq N,\os{\circ}{T}{}'_0} @>{g'_{Y_{\bul \leq N,\os{\circ}{T}{}'_0}}}>> 
Y''_{\bul \leq N,\os{\circ}{T}{}'_0} \\ 
@VVV @VVV \\ 
Y_{\bul \leq N,\os{\circ}{T}{}'_0} @>{g_{Y_{\bul \leq N,\os{\circ}{T}{}'_0}}}
>> Y_{\bul \leq N,\os{\circ}{T}{}'_0}  \\ 
@VVV @VVV \\ 
S'_{\os{\circ}{T}{}'_0} @>>> S'_{\os{\circ}{T}{}'_0} \\ 
@V{\bigcap}VV @VV{\bigcap}V \\ 
S'(T')^{\nat} @>{u'}>> S'(T')^{\nat},  
\end{CD}
\tag{3.6.10.1}\label{eqn:yhxubqsas}
\end{equation*}
to {\rm (\ref{cd:xgsqpxy})}. 
Denote 
$g_{\bul \leq N}$ in {\rm (\ref{cd:xgsqpxy})} by 
$g_{X_{\bul \leq N,\os{\circ}{T}_0}}$. 
Assume that 
$g_{\bul \leq N}\col X_{\bul \leq N,\os{\circ}{T}_0}\lo Y_{\bul \leq N,\os{\circ}{T}{}'_0}$ 
is a morphism from the commutative 
{\rm (\ref{cd:xdjqbxy})} to {\rm (\ref{eqn:yhxubqsas})}. 
Let 
$\os{\circ}{g}{}^*_{\bul \leq N,{\rm crys}}(F^{\bul \leq N})
\lo E^{\bul \leq N}$ be a morphism of 
${\cal O}_{\os{\circ}{X}_{\bul \leq N,T_0}/\os{\circ}{T}}$-modules in {\rm (\ref{ali:gnafe})}. 
Let $F^{\bul \leq N}$ be a flat quasi-coherent crystal of 
${\cal O}_{\os{\circ}{Y}_{\bul \leq N,T'_0}/\os{\circ}{T}{}'}$-modules. 
Let 
$$\os{\circ}{g}{}^*_{\bul \leq N,{\rm crys}}(F^{\bul \leq N})\lo E^{\bul \leq N}$$ 
be a morphism of flat quasi-coherent crystals of 
${\cal O}_{Y_{\bul \leq N,\os{\circ}{T}}/\os{\circ}{T}}$-modules. 
Then the following 
diagram is commutative$:$ 
\begin{equation*} 
\begin{CD} 
Rg_{{\bul \leq N}*}(A_{{\rm zar},{\mab Q}}(X_{\bul \leq N,\os{\circ}{T}_0}/S(T)^{\nat}
,E^{\bul \leq N}),P)@>{\nu_{S(T)^{\nat},{\rm zar}} }>>   \\ 
@A{g^*_{\bul \leq N}}AA \\
(A_{{\rm zar},{\mab Q}}(Y_{\bul \leq N,\os{\circ}{T}{}'_0}/S'(T')^{\nat}
,F^{\bul \leq N}),P)@>{\nu_{S'(T')^{\nat},{\rm zar}} }>>   
\end{CD} 
\tag{3.6.10.2}\label{eqn:aasabf}
\end{equation*} 
\begin{equation*} 
\begin{CD} 
Rg_{{\bul \leq N}*}(A_{{\rm zar},{\mab Q}}(X_{\bul \leq N,\os{\circ}{T}_0}/S(T)^{\nat}
,E^{\bul \leq N}),P\langle -2 \rangle)(-1,u)\\
@AA{g^*_{\bul \leq N}}A \\
(A_{{\rm zar},{\mab Q}}(Y_{\bul \leq N,\os{\circ}{T}{}'_0}/S(T)^{\nat},F^{\bul \leq N}),
P\langle -2 \rangle)(-1,u'). 
\end{CD} 
\end{equation*} 
\end{prop}

Set 
$$P_kB^{\bul \leq N,\bul ij}
:=P_kA^{\bul \leq N,\bul i-1,j}(-1,u)\oplus P_kA^{\bul \leq N,\bul ij} 
\quad ( i,j\in {\mab N})$$
and 
$$P_kB^{\bul \leq N,ij}
:=P_kA^{\bul \leq N,i-1,j}(-1,u)\oplus P_kA^{\bul \leq N,ij} 
\quad ( i,j\in {\mab N}).$$
Then we have a filtered complex 
$(B^{\bul \leq N,\bul},P)$ in ${\rm D}^+{\rm F}(f^{-1}_{\bul \leq N}({\cal O}_T))$. 
\begin{prop}\label{prop:baid}
$(1)$ There exists the following sequence of the triangles in 
${\rm D}^+{\rm F}(f^{-1}_{\bul \leq N}({\cal O}_T)):$ 
\begin{align*} 
&\lo (A_{{\rm zar},{\mab Q}}(X_{\bul \leq N,\os{\circ}{T}_0}/S(T)^{\nat},E^{\bul \leq N}),P)[-1]
\lo (B^{\bul \leq N,\bul},P)\\
&\lo (A_{{\rm zar},{\mab Q}}(X_{\bul \leq N,\os{\circ}{T}_0}/S(T)^{\nat},E^{\bul \leq N}),P)\os{+1}{\lo}. 
\end{align*} 
\par 
$(2)$ The filtered complex $(B^{\bul \leq N,\bul},P)$ is independent of the choice 
of an affine $N$-truncated simplicial open covering of $X_{\bul \leq N,\os{\circ}{T}_0}$ 
and the choice of an $(N,\infty)$-truncated bisimplicial immersion 
$X_{\bul \leq N,\bul,\os{\circ}{T}_0} \os{\sus}{\lo} 
\ol{\cal P}_{\bul \leq N,\bul}$ over $\ol{S(T)^{\nat}}$.  
\par 
$(3)$ The morphism 
{\rm (\ref{ali:sgsclxan})} is independent of the choices in $(2)$. 
\end{prop} 

\begin{defi} 
We call $(B^{\bul \leq N,\bul},P)$ 
the {\it extended zariskian $p$-adic filtered Steenbrink complexes} of 
$E^{\bul \leq N}$ for $X_{\bul \leq N,\os{\circ}{T}_0}/(S(T)^{\nat},{\cal J},\del)$. 
We denote it  by
$$(B_{{\rm zar},{\mab Q}}(X_{\bul \leq N,\os{\circ}{T}_0}/S(T)^{\nat},E^{\bul \leq N}),P)
\in {\rm D}^+{\rm F}(f^{-1}_{\bul \leq N}({\cal O}_T)).$$ 
When $E^{\bul \leq N}={\cal O}_{\os{\circ}{X}_{\bul \leq N,T_0}/\os{\circ}{T}}$, we denote 
$(B_{{\rm zar},{\mab Q}}(X_{\bul \leq N,\os{\circ}{T}_0}/S(T)^{\nat},E^{\bul \leq N}),P)$ 
by 
$$(B_{{\rm zar},{\mab Q}}(X_{\bul \leq N,\os{\circ}{T}_0}/S(T)^{\nat}),P).$$
We call $(B_{{\rm zar},{\mab Q}}(X_{\bul \leq N,\os{\circ}{T}_0}/S(T)^{\nat}),P)$ the 
{\it extended zariskian $p$-adic filtered Steenbrink complex} of 
$X_{\bul \leq N,\os{\circ}{T}_0}/(S(T)^{\nat},{\cal J},\del)$. 
\end{defi} 

\begin{theo}[{\bf Contravariant functoriality of $B_{\rm zar}$}]\label{theo:funacb}
$(1)$ Let 
$g_{\bul \leq N}\col 
X_{\bul \leq N,\os{\circ}{T}{}_0}\lo Y_{\bul \leq N,\os{\circ}{T}{}'_0}$ be 
the morphism in the beginning of {\rm \S\ref{sec:pfnsc}}. 
Then $g_{\bul \leq N}$ induces the following 
well-defined pull-back morphism 
\begin{equation*}  
g_{\bul \leq N}^* \col 
(B_{{\rm zar},{\mab Q}}
(Y_{\bul \leq N,\os{\circ}{T}{}'_0}/S'(T')^{\nat},F^{\bul \leq N}),P)
\lo Rg_{\bul \leq N*}
((B_{{\rm zar},{\mab Q}}(X_{\bul \leq N,\os{\circ}{T}_0}/S(T)^{\nat}
,E^{\bul \leq N}),P)) 
\tag{3.6.13.1}\label{eqn:fzbabbad}
\end{equation*} 
fitting into the following commutative diagram$:$
\begin{equation*} 
\begin{CD}
B_{{\rm zar},{\mab Q}}(Y_{\bul \leq N,\os{\circ}{T}{}'_0}/S'(T')^{\nat},F^{\bul \leq N})
@>{g_{\bul \leq N}^*}>>  
Rg_{\bul \leq N*}(B_{{\rm zar},{\mab Q}}(X_{\bul \leq N,\os{\circ}{T}_0}/S(T)^{\nat}
,E^{\bul \leq N}))\\ 
@A{\mu_{Y_{\bul \leq N,\os{\circ}{T}{}'_0}/S'(T')^{\nat}}\wedge}A{\simeq}A 
@A{Rg_{\bul \leq N*}(\mu_{X_{\bul \leq N,\os{\circ}{T}_0}/S(T)^{\nat}}\wedge)}A{\simeq}A\\
\wt{R}u_{Y_{\bul \leq N,\os{\circ}{T}{}'_0}/\os{\circ}{T}{}'*}
(\eps^*_{Y_{\bul \leq N,\os{\circ}{T}{}'_0}/\os{\circ}{T}{}'}(F^{\bul \leq N}))_{\mab Q}
@>{g_{\bul \leq N}^*}>>Rg_{\bul \leq N*}\wt{R}u_{X_{\bul \leq N,\os{\circ}{T}_0}/\os{\circ}{T}*}
(\eps^*_{X_{\bul \leq N,\os{\circ}{T}_0}/\os{\circ}{T}}(E^{\bul \leq N}))_{\mab Q}.
\end{CD}
\tag{3.6.13.2}\label{cd:psbcabcz} 
\end{equation*}
\par 
$(2)$ 
\begin{align*} 
(h_{\bul \leq N}\circ g_{\bul \leq N})^*  =&
Rh_{\bul \leq N*}(g_{\bul \leq N}^*)\circ h_{\bul \leq N}^*    
\col (B_{{\rm zar},{\mab Q}}(Z_{\bul \leq N,\os{\circ}{T}{}''_0}/S''(T'')^{\nat},F^{\bul \leq N}),P)
\tag{3.6.13.3}\label{ali:pdbbap} \\ 
& \lo Rh_{\bul \leq N*}Rg_{\bul \leq N*}
(B_{{\rm zar},{\mab Q}}(X_{\bul \leq N,\os{\circ}{T}_0}/S(T)^{\nat},E^{\bul \leq N}),P) \\
& =R(h_{\bul \leq N}\circ g_{\bul \leq N})_*
(B_{{\rm zar},{\mab Q}}(X_{\bul \leq N,\os{\circ}{T}_0}/S(T)^{\nat},E^{\bul \leq N}),P).
\end{align*}  
\par 
$(3)$  
\begin{equation*} 
{\rm id}_{X_{\bul \leq N,\os{\circ}{T}_0}}^*={\rm id} 
\col (B_{{\rm zar},{\mab Q}}(X_{\bul \leq N,\os{\circ}{T}_0}/S(T)^{\nat},E^{\bul \leq N}),P)
\lo (B_{{\rm zar},{\mab Q}}(X_{\bul \leq N,\os{\circ}{T}_0}/S(T)^{\nat},E^{\bul \leq N}),P).  
\tag{3.6.13.4}\label{eqn:fzibdabd}
\end{equation*} 
\end{theo}

\begin{coro}\label{coro:sqpfs}
The isomorphisms {\rm (\ref{ali:usz})} and {\rm (\ref{ali:sgsclxan})} 
for $X_{\bul \leq N,\os{\circ}{T}_0}/S(T)^{\nat}$, $E^{\bul \leq N}$ 
and $Y_{\bul \leq N,\os{\circ}{T}{}'_0}/S'(T')^{\nat}$, $F^{\bul \leq N}$ 
induce a morphism $g_{\bul \leq N}^*$ from {\rm (\ref{eqn:ydhdagguss})} 
for the case $Y_{\bul \leq N}=X_{\bul \leq N}$ and $Z_{\bul \leq N}=Y_{\bul \leq N}$
there to {\rm (\ref{eqn:aasabf})}.
This morphism satisfies the transitive relation and  
${\rm id}_{X_{\bul \leq N,\os{\circ}{T}_0}/S(T)^{\nat}}^*={\rm id}$. 
\end{coro}

\begin{theo-defi}\label{theo:facp} 
The following hold$:$
\par 
$(1)$ The filtered complex 
$$R\pi_{{\rm zar}*}(({\cal E}^{\bul \leq N,\bul}
\otimes_{{\cal O}_{{\cal P}^{\rm ex}_{\bul \leq N,\bul}}}
{\Om}^{\bul}_{{\cal P}^{\rm ex}_{\bul \leq N,\bul}/\os{\circ}{T}}\otimes_{\mab Z}{\mab Q},P))$$ 
is independent of the choice of the disjoint union 
of the member of an affine $N$-truncated simplicial open covering of 
$X_{\bul \leq N,\os{\circ}{T}_0}$ and 
an $(N,\infty)$-truncated bisimplicial immersion 
$X_{\bul \leq N,\bul,\os{\circ}{T}_0} \os{\sus}{\lo} \ol{\cal P}_{\bul \leq N,\bul}$ over 
$\ol{S(T)^{\nat}}$. 
Set 
\begin{align*} 
(\wt{R}u_{X_{\bul \leq N,\os{\circ}{T}_0}/\os{\circ}{T}*}
(\eps^*_{X_{\bul \leq N,\os{\circ}{T}_0}/\os{\circ}{T}}
(E^{\bul \leq N}))_{\mab Q} ,P)
:=R\pi_{{\rm zar}*}(({\cal E}^{\bul \leq N,\bul}
\otimes_{{\cal O}_{{\cal P}^{\rm ex}_{\bul \leq N,\bul}}}
{\Om}^{\bul}_{{\cal P}^{\rm ex}_{\bul \leq N,\bul}/\os{\circ}{T}}\otimes_{\mab Z}{\mab Q},P)).
\end{align*}
We call this filtered complex the 
{\it modified $P$-filtered log isocrystalline complex} of 
$\eps^*_{X_{\bul \leq N,\os{\circ}{T}_0}/\os{\circ}{T}}(E^{\bul \leq N})$. 
When $E^{\bul \leq N}={\cal O}_{\os{\circ}{X}_{\bul \leq N,T_0}/\os{\circ}{T}}$, 
we call this filtered complex 
the {\it modified $P$-filtered log isocrystalline complex} of 
$X_{\bul \leq N,\os{\circ}{T}_0}/\os{\circ}{T}$. 
\par 
$(2)$ For each $0\leq m\leq N$, there exist the following canonical isomorphisms 
in ${\rm D}^+(f^{-1}_{m,T}({\cal O}_T)):$ 
\begin{align*}
P_0(\wt{R}u_{X_{m,\os{\circ}{T}_0}/\os{\circ}{T}*}
(\eps^*_{X_{m,\os{\circ}{T}_0}/\os{\circ}{T}}(E^m)))_{\mab Q}
\os{\sim}{\lo} &
{\rm MF}(a^{(0)}_{*}
Ru_{\os{\circ}{X}{}^{(0)}_{m,T_0}/\os{\circ}{T}*}
(E^m_{\os{\circ}{X}{}^{(0)}_{m,T_0}/\os{\circ}{T}}
\otimes_{\mab Z}\vp^{(0)}_{\rm crys}
(\os{\circ}{X}_{m,T_0}/\os{\circ}{T}))_{\mab Q}\tag{3.6.15.1}\label{ali:pa0dept}\\
&\lo 
{\rm MF}(a^{(1)}_{*}
Ru_{\os{\circ}{X}{}^{(1)}_{m,T_0}/\os{\circ}{T}*}
(E^m_{\os{\circ}{X}{}^{(1)}_{m,T_0}/\os{\circ}{T}}
\otimes_{\mab Z}\vp^{(1)}_{\rm crys}
(\os{\circ}{X}_{m,T_0}/\os{\circ}{T}))_{\mab Q}  \\
&\lo \cdots \\
&\lo 
{\rm MF}(a^{(l)}_{*}
Ru_{\os{\circ}{X}{}^{(l)}_{m,T_0}/\os{\circ}{T}*}
(E^m_{\os{\circ}{X}{}^{(l)}_{m,T_0}/\os{\circ}{T}}
\otimes_{\mab Z}\vp^{(l)}_{\rm crys}(\os{\circ}{X}_{m,T_0}/\os{\circ}{T}))_{\mab Q}\\
& \lo \cdots
\cdots)\cdots ),  
\end{align*} 
\begin{align*}
&{\rm gr}_k^P(\wt{R}u_{X_{m,\os{\circ}{T}_0}/\os{\circ}{T}*}
(\eps^*_{X_{m,\os{\circ}{T}_0}/\os{\circ}{T}}(E^m))_{\mab Q}) 
\tag{3.6.15.2}\label{ali:eaoppd}\\
&\os{\sim}{\lo}  
a^{(k-1)}_*
Ru_{\os{\circ}{X}{}^{(k-1)}_{m,T_0}/\os{\circ}{T}*}
(E^m_{\os{\circ}{X}{}^{(k-1)}_{m,T_0}/\os{\circ}{T}}
\otimes_{\mab Z}
\vp^{(k-1)}_{\rm crys}(\os{\circ}{X}_{m,T_0}/\os{\circ}{T}))_{\mab Q}  
\quad (k\in {\mab Z}_{\geq 1}).
\end{align*}
\par 
$(3)$ {\bf {\rm (Contravariant functoriality)}} 
Let $g_{\bul \leq N}$ be the morphism in the beginning of \S\ref{sec:pfnsc}. 
Then $g_{\bul \leq N}$ induces the following pull-back morphism 
\begin{align*}
g_{\bul \leq N}^*\col &
(\wt{R}u_{Y_{\bul \leq N,\os{\circ}{T}{}'_0}/\os{\circ}{T}{}'*}
(\eps^*_{Y_{\bul \leq N,\os{\circ}{T}{}'_0}/\os{\circ}{T}{}'}(F^{\bul \leq N}))_{\mab Q},P)
\tag{3.6.15.3}\label{ali:rauynt}\\
&\lo 
Rg_{\bul \leq N *}
((\wt{R}u_{X_{\bul \leq N,\os{\circ}{T}_0}/\os{\circ}{T}*}
(\eps^*_{X_{\bul \leq N,\os{\circ}{T}_0}/\os{\circ}{T}}(E^{\bul \leq N}))_{\mab Q},P)). 
\end{align*}
\par 
(4) Assume that, for any $0\leq m\leq N$, the conditions $(3.4.8.1)$ and $(3.4.8.2)$ hold. 
Then the morphism $g_{\bul \leq N}^*$ in $(3)$ depends only on 
$\os{\circ}{g}_{\bul \leq N}$, $\os{\circ}{T}{}'\lo \os{\circ}{S}{}'$
and $\os{\circ}{T}\lo \os{\circ}{S}$. 
\end{theo-defi}
\begin{proof} 
The proof is the same as that of (\ref{theo:fcp}). 
\end{proof}

\chapter{Weight filtrations and slope filtrations on log isocrystalline cohomologies} 

\parno 
In this chapter we construct a filtered complex, 
which is a generalization of the filtered Hyodo-Mokrane-Steenbrink complex 
defined in the chapter II modulo torsion. 
To construct this filtered complex, 
we need to calculate the graded complex of the Poincar\'{e} filtration on 
the filtered complex which is the log de Rham-Witt version of the 
filtered complex in (\ref{theo:facp}). 

\section{Iso-zariskian $p$-adic filtered Steenbrink complexes  
via log de Rham-Witt complexes}\label{sec:flgsiidw}
Let $\kap$, ${\cal W}$, $s$ and ${\cal W}_n(s)$ 
be as in \S\ref{sec:flgdw}. 
Let $\kap_0$, ${\cal W}_0$, $s_0$ and ${\cal W}_n(s_0)$ 
be an analogous objects to $\kap$, ${\cal W}$, $s$ and ${\cal W}_n(s)$, 
respectively. 
In this section we denote the sequence (\ref{eqn:stcnniseq}) 
in the case of log points of perfect fields of characteristic $p$ by   
\begin{equation*} 
s=s_N \lo s_{N-1} \lo \cdots \lo s_0. 
\tag{4.1.0.1}\label{eqn:sttsseq}
\end{equation*} 
Let $g\col t\lo s$ be the morphism of log schemes in \S\ref{sec:flgdw}.  
Set $\kap_t:=\Gam(t,{\cal O}_t)$. 
Let $N$ be a nonnegative integer or $\infty$. 
Let $X_{\bul \leq N}/s$ be 
the $N$-truncated simplicial base change of SNCL schemes. 
Let $E^{\bul \leq N}$ be a unit root crystal of 
${\cal O}_{\os{\circ}{X}_{\bul \leq N,t}/{\cal W}(\os{\circ}{t})}$-modules 
in \S\ref{sec:pgensc} for the situation above. 
\par 
In this section we construct a weight filtration 
on ${\cal W}\wt{\Om}^{\bul}_{X_{\bul \leq N,\os{\circ}{t}}}$ 
and a certain filtered complex 
$({\cal W}A_{X_{\bul \leq N,\os{\circ}{t}},{\mab Q}}(E^{\bul \leq N}),P)$, 
which is a generalization of the filtered complex  
$({\cal W}A_{X_{\bul \leq N,\os{\circ}{t}}}(E^{\bul \leq N})
\otimes_{\mab Z}{\mab Q},P)$, where 
$({\cal W}A_{X_{\bul \leq N,\os{\circ}{t}}}(E^{\bul \leq N}),P)$
is a filtered complex defined in \S\ref{sec:flgdw}.     
In this section 
we prove that there exists a canonical isomorphism 
\begin{align*} 
(A_{{\rm zar},{\mab Q}}(X_{\bul \leq N,\os{\circ}{T}_0}/S(T)^{\nat},E^{\bul \leq N}),P)
\os{\sim}{\lo} ({\cal W}A_{X_{\bul \leq N,\os{\circ}{t}},{\mab Q}}(E^{\bul \leq N}),P)
\tag{4.1.0.2}\label{ali:aqxts} 
\end{align*}  
in the case where $S=s$ and $T={\cal W}(t)$.   
\par 
First we consider the case $N=0$.

\begin{prop}\label{prop:bcpw}
Let $s'\lo s$ be a morphism of log points of 
perfect fields of characteristic $p>0$. 
Set $\kap':=\Gam(s',{\cal O}_{s'})$ and let ${\cal W}'$ be the Witt ring of $\kap'$. 
Assume that the morphism $t\lo s$ factors through this morphism. 
Let $Y$ be a fine log smooth integral scheme over $s$. 
Set  $Y':=Y\times_ss'$.  
Then the canonical morphism 
\begin{align*} 
{\cal W}(\kap_{t})\otimes_{{\cal W}}
{\cal W}\wt{\Om}^i_{Y}\otimes_{\mab Z}{\mab Q}\lo 
{\cal W}(\kap_{t})\otimes_{{\cal W}'}{\cal W}\wt{\Om}^i_{Y'}\otimes_{\mab Z}{\mab Q}
\tag{4.1.1.1}\label{cd:lijyzqext}
\end{align*} 
is an isomorphism. 
\end{prop}
\begin{proof} 
It suffices to prove that 
the canonical morphism 
\begin{align*} 
{\cal W}'\otimes_{{\cal W}}
{\cal W}\wt{\Om}^i_{Y}\otimes_{\mab Z}{\mab Q}\lo 
{\cal W}\wt{\Om}^i_{Y'}\otimes_{\mab Z}{\mab Q}
\end{align*} 
is an isomorphism. 
Set $\theta=d\log \tau$ and $\theta'=d\log \tau'$, 
where $\tau$ and $\tau'$  are global sections of $M_s$ and $M_{s'}$ 
whose images in $M_s/{\cal O}_s^*$ and $M_{s'}/{\cal O}_{s'}^*$ are generators. 
Then the morphism 
${\mab N}=M_s/{\cal O}_s^*\lo M_{s'}/{\cal O}_{s'}^*={\mab N}$ is 
a multiplication for a positive integer $e'$. 
Then we have the following commutative diagram of exact sequences by (\ref{prop:plz}):  
\begin{equation*}
\begin{CD}
0 @>>> {\cal W}\Om_{Y'}^{\bul}[-1] 
@>{\theta' \wedge}>>  {\cal W}\wt \Om_{Y'}^{\bul} @>>> \\ 
@. @A{e'\times}AA @AAA \\
0 @>>> {\cal W}'\otimes_{{\cal W}}{\cal W}\Om_{Y}^{\bul}[-1] 
@>{\theta \wedge}>>  {\cal W}'\otimes_{{\cal W}}{\cal W}\wt \Om_{Y}^{\bul} 
@>>> 
\end{CD}
\tag{4.1.1.2}\label{cd:lijwqext}
\end{equation*}
\begin{equation*}
\begin{CD}
{\cal W}\Om_{Y'}^{\bul}@>>> 0\\
@AAA \\
{\cal W}'\otimes_{\cal W}{\cal W}\Om_{Y}^{\bul}@>>> 0. 
\end{CD} 
\end{equation*}
By (\ref{prop:bcdw}) the morphism 
${\cal W}'\otimes_{\cal W}{\cal W}\Om_{Y}^{\bul}\lo {\cal W}\Om_{Y'}^{\bul}$ 
is an isomorphism. 
Hence, by (\ref{cd:lijwqext}), we see that 
(\ref{cd:lijyzqext}) is an isomorphism.  
\end{proof} 

\par 
Let $X/s$ be a base change of SNCL schemes defined in (\ref{defi:snb}) (1). 
Let the notations be as in (\ref{cd:yiqi}) and (\ref{cd:xpqst}) 
in the case $N=0$. 
Set 
\begin{align*} 
P_k{\cal W}\wt{\Om}{}^i_{X_{\os{\circ}{t}},{\mab Q}}:=
R\pi_{{\rm zar}*}
({\cal H}^i(P_k\Om^{\bul}_{{\cal P}^{\rm ex}_{\bul}/{\cal W}}
\otimes_{\mab Z}{\mab Q})).  
\end{align*}  
It is clear that 
\begin{align*} 
{\cal W}\wt{\Om}{}^i_{X_{\os{\circ}{t}},{\mab Q}}
:=P_{\infty}{\cal W}\wt{\Om}{}^i_{X_{\os{\circ}{t}},{\mab Q}}
={\cal W}\wt{\Om}{}^i_{X_{\os{\circ}{t}}}\otimes_{\mab Z}{\mab Q}, 
\tag{4.1.1.3}\label{ali:qpwo} 
\end{align*} 
where ${\cal W}\wt{\Om}{}^i_X$ is the sheaf in (\ref{eqn:lmtwy}) for $Y=X$.

\begin{prop}\label{prop:grbrs}
There exists the following residue isomorphism
\begin{align*} 
{\rm Res} \col & 
P_k{\cal W}\wt{{\Om}}{}^{\bul}_{X_{\os{\circ}{t}},{\mab Q}}
\lo a^{(k-1)}_{*}
({\cal W}\Om^{\bul}_{\os{\circ}{X}{}^{(k-1)}_t}
\otimes_{\mab Z}\vp^{(k-1)}_{\rm zar}(\os{\circ}{X}_t/\os{\circ}{t}))[-k]
\otimes_{\mab Z}{\mab Q} 
\quad (k\in {\mab Z}_{\geq 1}) \tag{4.1.2.1}\label{ali:retbsd}
\end{align*} 
fitting into the following exact sequence 
\begin{align*} 
0 \lo & 
P_{k-1}{\cal W}\wt{{\Om}}{}^{\bul}_{X_{\os{\circ}{t}},{\mab Q}}
\lo P_k{\cal W}\wt{{\Om}}{}^{\bul}_{X_{\os{\circ}{t}},{\mab Q}}
\os{\rm Res}{\lo}  
a^{(k-1)}_{*}
({\cal W}\Om^{\bul}_{\os{\circ}{X}{}^{(k-1)}_t}
\otimes_{\mab Z}\vp^{(k-1)}_{\rm zar}(\os{\circ}{X}_t/\os{\circ}{t}))[-k] 
\otimes_{\mab Z}{\mab Q} 
\lo 0. \tag{4.1.2.2}\label{eqn:ptbwa}
\end{align*} 
\end{prop}
\begin{proof} 
Set $\kap_i:=\Gam(s_i,{\cal O}_{s_i})$. 
Let ${\cal Q}^{\rm ex}_{i\bul}$ $(1\leq i\leq l)$ 
be the log scheme in (\ref{defi:bca}). 
By the functoriality of the Poincar\'{e} residue isomorphism ((\ref{prop:rescos})),  
we obtain the following exact sequence: 
\begin{align*} 
0\lo P_{k-1}\Om^{\bul}_{{\cal Q}^{\rm ex}_{i\bul}/{\cal W}(\os{\circ}{s}_i)}
\otimes_{\mab Z}{\mab Q}\lo 
P_k\Om^{\bul}_{{\cal Q}^{\rm ex}_{i\bul}/{\cal W}(\os{\circ}{s}_i)}
\otimes_{\mab Z}{\mab Q} 
\lo \Om^{\bul}_{{\cal Q}^{{\rm ex},(k-1)}_{i\bul}/{\cal W}(\os{\circ}{s}_i)}[-k]
\otimes_{\mab Z}{\mab Q}\lo 0.  
 \tag{4.1.2.3}\label{ali:obgst}   
\end{align*}
As in \cite[1.2]{msemi} we may assume that this sequence is split. 
Hence the following sequence is also split: 
\begin{align*} 
0&\lo \bigoplus_{i=1}^l{\cal W}(\kap_t)\otimes_{{\cal W}(\kap_i)}
P_{k-1}({\cal E}^{\bul}
\otimes_{{\cal O}_{{\cal Q}^{\rm ex}_{i\bul}}}
{\Om}^{\bul}_{{\cal Q}^{\rm ex}_{i\bul}/{\cal W}(\os{\circ}{s}_i)})
\otimes_{\mab Z}{\mab Q}\tag{4.1.2.4}\label{ali:obtgst}  \\
&\lo 
\bigoplus_{i=1}^l{\cal W}(\kap_t)\otimes_{{\cal W}(\kap_i)}
P_k({\cal E}^{\bul}
\otimes_{{\cal O}_{{\cal Q}^{\rm ex}_{i\bul}}}
{\Om}^{\bul}_{{\cal Q}^{\rm ex}_{i\bul}/{\cal W}(\os{\circ}{s}_i)})
\otimes_{\mab Z}{\mab Q}   \\
&\lo \bigoplus_{i=1}^l{\cal W}(\kap_t)\otimes_{{\cal W}(\kap_i)}{\cal E}^{\bul}
\otimes_{{\cal O}_{{\cal Q}^{\rm ex}_{i\bul}}}
{\Om}^{\bul}_{{\cal Q}^{{\rm ex},(k-1)}_{i\bul}/{\cal W}(\os{\circ}{s}_i)}[-k]
\otimes_{\mab Z}{\mab Q}\lo 0.  
\end{align*}
By (\ref{prop:bcpw}) this is isomorphic to the following split exact sequence
\begin{align*}
0 &\lo 
P_{k-1}
\Om^{\bul}_{{\cal Q}^{\rm ex}_{i\bul}\times_{{\cal W}(\os{\circ}{s})}
{\cal W}(\os{\circ}{t})/{\cal W}(\kap_t)}\otimes_{\mab Z}{\mab Q}
\lo 
P_k\Om^{\bul}_{{\cal Q}^{\rm ex}_{i\bul}\times_{{\cal W}(\os{\circ}{s})}
{\cal W}(\os{\circ}{t})/{\cal W}(\kap_t)}\otimes_{\mab Z}{\mab Q}\\
&\lo 
\Om^{\bul}_{{{\cal Q}^{{\rm ex},(k-1)}_{i\bul}\times_{{\cal W}(\os{\circ}{s})}
{\cal W}(\os{\circ}{t})}/{\cal W}(\kap_t)}[-k]
\otimes_{\mab Z}{\mab Q}\lo 0.  
\end{align*}
Taking the cohomologies of the exact sequence above, 
we obtain the exact sequence (\ref{eqn:ptbwa}). 
\end{proof}

\begin{coro}\label{coro:sqe}
The sheaf $P_k{\cal W}\wt{\Om}{}^i_{X_{\os{\circ}{t}},{\mab Q}}$ 
$(k\in {\mab N})$ is well-defined. 
\end{coro}
\begin{proof} 
(\ref{coro:sqe}) follows from (\ref{ali:qpwo}), 
(\ref{prop:grbrs}) and the descending induction. 
\end{proof}

\begin{prop}\label{prop:p0}
There exists the following canonical functorial isomorphism
\begin{align*}
P_0({\cal W}_n\wt{{\Om}}{}^i_{X_{\os{\circ}{t}},{\mab Q}})
\os{\sim}{\lo} 
&{\cal H}^i[{\rm MF}(a^{(0)}_*
({\cal W}_n\Om^{\bul}_{\os{\circ}{X}{}^{(0)}_t,{\mab Q}}
\otimes_{\mab Z}\vp^{(0)}_{\rm zar}
(\os{\circ}{X}_{t_0}/\os{\circ}{t}))
\tag{4.1.4.1}\label{ali:pdae}\\
&\lo 
{\rm MF}(a^{(1)}_*({\cal W}_n\Om^{\bul}_{\os{\circ}{X}{}^{(1)}_t,{\mab Q}}
\otimes_{\mab Z}\vp^{(1)}_{\rm zar}
(\os{\circ}{X}_{t_0}/\os{\circ}{t}))  \\
&\lo \cdots \\
&\lo 
{\rm MF}(a^{(m)}_*({\cal W}_n\Om^{\bul}_{\os{\circ}{X}{}^{(m)}_t,{\mab Q}}
\otimes_{\mab Z}\vp^{(m)}_{\rm zar}
(\os{\circ}{X}_{t_0}/\os{\circ}{t}))\\
& \lo \cdots
\cdots)\cdots )]. 
\end{align*}
\end{prop}  
\begin{proof} 
This follows from (\ref{ali:pa0dept}) and \cite[II, (1.4)]{idw}.
\end{proof}

\par 
We have the 
$N$-truncated cosimplicial log de Rham-Witt complex 
$(E^{\bul \leq N}_{\infty}\otimes_{{\cal W}({\cal O}_{X_{\bul \leq N,\os{\circ}{t}}})}
{\cal W}\wt{{\Om}}^{\bul}_{X_{\bul \leq N,\os{\circ}{t}}},P)$ 
associated to $X_{\bul \leq N,\os{\circ}{t}}/{\cal W}(s_{\os{\circ}{t}})$ and 
$E^{\bul \leq N}$. 
Let $\theta_{\bul \leq N} \in {\cal W}\wt{\Om}_{X_{\bul \leq N,\os{\circ}{t}}}^1$
be the 1-form essentially defined in \cite[3.4]{msemi} by using 
$\theta_{{\cal P}^{\rm ex}_{\bul \leq N,\bul}}$ in \S\ref{sec:pgensc}. 
Set 
\begin{align*} 
\theta_{X_{\bul \leq N,\os{\circ}{t}}/{\cal W}(s_{\os{\circ}{t}})/{\cal W}(s_{0,\os{\circ}{t}})}
:=e\theta_{\bul \leq N}
\end{align*} 
and  
\begin{align*} 
{\cal W}A_{X_{\bul \leq N,\os{\circ}{t}},{\mab Q}}(E^{\bul \leq N})^{ij}
:=&
E^{\bul \leq N}_{\infty}
\otimes_{{\cal W}({\cal O}_{X_{\bul \leq N,\os{\circ}{t}}})}
{\cal W}\wt{\Om}^{i+j+1}_{X_{\bul \leq N,\os{\circ}{t}}}
\otimes_{\mab Z}{\mab Q}
\tag{4.1.4.2}\label{eqn:wxlpj} \\
& /P_j(E^{\bul \leq N}_{\infty}\otimes_{{\cal W}({\cal O}_{X_{\bul \leq N,\os{\circ}{t}}})}
{\cal W}\wt{\Om}^{i+j+1}_{X_{\bul \leq N,\os{\circ}{t}}})
\otimes_{\mab Z}{\mab Q}
\quad (i,j \in {\mab N}).  
\end{align*} 
We consider the following boundary morphisms of the double complex
\begin{equation*}
\begin{CD}
{\cal W}A_{X_{\bul \leq N,\os{\circ}{t}},{\mab Q}}(E^{\bul \leq N})^{i,j+1}@.  \\ 
@AA{\theta_{X_{\bul \leq N,\os{\circ}{t}}/{\cal W}(s_{\os{\circ}{t}})/{\cal W}(s_{0,\os{\circ}{t}})}\wedge}A  @. \\
{\cal W}A_{X_{\bul \leq N,\os{\circ}{t}},{\mab Q}}(E^{\bul \leq N})^{ij}@>{-\nabla}>> 
{\cal W}A_{X_{\bul \leq N,\os{\circ}{t}},{\mab Q}}(E^{\bul \leq N})^{i+1,j}.\\
\end{CD}
\tag{4.1.4.3}\label{cd:lcwnsbd} 
\end{equation*}   
Then 
${\cal W}A_{X_{\bul \leq N,\os{\circ}{t}},{\mab Q}}(E^{\bul \leq N})^{\bul \bul}$
becomes 
an $N$-truncated cosimplicial double complex 
of $f^{-1}_{\bul \leq N}(K_0(\kap_t))$-modules, where 
$K_0(\kap_t):={\rm Frac}({\cal W}(\kap_t))$. 
Let 
${\cal W}A_{X_{\bul \leq N,\os{\circ}{t}},{\mab Q}}(E^{\bul \leq N})$ 
be the single complex of  
${\cal W}A_{X_{\bul \leq N,\os{\circ}{t}},{\mab Q}}(E^{\bul \leq N})^{\bul \bul}$.

\begin{prop}\label{prop:mnb}
The following morphism
\begin{equation*} 
\theta_{X_{\bul \leq N,\os{\circ}{t}}/{\cal W}(s_{\os{\circ}{t}})/{\cal W}(s_{0,\os{\circ}{t}})} \wedge \col 
E^{\bul \leq N}_{\infty}\otimes_{{\cal W}({\cal O}_{X_{\bul \leq N,t}})}
{\cal W}{\Om}_{X_{\bul \leq N,t}}^i
\otimes_{{\mab Z}}{\mab Q}
\lo 
{\cal W}A_{X_{\bul \leq N,\os{\circ}{t}},{\mab Q}}(E^{\bul \leq N})^{i \bul}
\tag{4.1.5.1}\label{eqn:wltqtpxa}
\end{equation*} 
is a quasi-isomorphism. 
Hence the following morphism  
\begin{equation*} 
\theta_{X_{\bul \leq N,\os{\circ}{t}}/
{\cal W}(s_{\os{\circ}{t}})/{\cal W}(s_{0,\os{\circ}{t}})} \wedge \col 
E^{\bul \leq N}_{\infty}\otimes_{{\cal W}({\cal O}_{X_{\bul \leq N,t}})}
{\cal W}{\Om}^{\bul}_{X_{\bul \leq N,t}}
\otimes_{{\mab Z}}{\mab Q}
\lo 
{\cal W}A_{X_{\bul \leq N,\os{\circ}{t}},{\mab Q}}(E^{\bul \leq N})
\tag{4.1.5.2}\label{eqn:wltntpxa}
\end{equation*} 
is a quasi-isomorphism. 
\end{prop}
\begin{proof} 
By using (\ref{prop:bcpw}), the same proof as that of 
\cite[3.15.1]{msemi} and \cite[(6.29)]{ndw} works 
(cf.~the proof of (\ref{prop:tefc})). 
\end{proof} 

\par 
We endow ${\cal W}A_{X_{\bul \leq N,\os{\circ}{t}},{\mab Q}}(E^{\bul \leq N})$ 
with a filtration $P$ as follows: 
\begin{align*} 
P_k
{\cal W}A_{X_{\bul \leq N,\os{\circ}{t}},{\mab Q}}(E^{\bul \leq N})^{ij} :=
{\rm Im}&(P_{2j+k+1}
(E^{\bul \leq N}\otimes_{{\cal W}({\cal O}_{X_{\bul \leq N,\os{\circ}{t}}})}
{\cal W}\wt{{\Om}}^{i+j+1}_{X_{\bul \leq N,\os{\circ}{t}}})\otimes_{{\mab Z}}{\mab Q}
\tag{4.1.5.3}\label{eqn:dblnad}\\
&\lo 
{\cal W}A_{X_{\bul \leq N,\os{\circ}{t}},{\mab Q}}(E^{\bul \leq N})^{ij}). 
\end{align*} 
Let $({\cal W}A_{X_{\bul \leq N,\os{\circ}{t}},{\mab Q}}(E^{\bul \leq N}),P)$ 
be the filtered single complex of the filtered double complex 
$({\cal W}A_{X_{\bul \leq N,\os{\circ}{t}},{\mab Q}}(E^{\bul \leq N})^{\bul \bul},P)$.  

\begin{defi}\label{defi:fcwnaxp}
We call the filtered complex 
$({\cal W}A_{X_{\bul \leq N,\os{\circ}{t}},{\mab Q}}(E^{\bul \leq N}),P)$ 
the {\it filtered Hyodo-Mokrane-Steenbrink complex} of   
$X_{\bul \leq N,\os{\circ}{t}}/s_{\os{\circ}{t}}$ and $E^{\bul \leq N}$.  
\end{defi}

\par 
Recall that $\kap_t=\Gam(t,{\cal O}_t)$. 
As in (\ref{lemm:axp}), 
we have the following formula by using (\ref{eqn:ptbwa}):
\begin{align*} 
{\rm gr}^P_k
{\cal W}A_{X_{m,t},{\mab Q}}(E^{\bul \leq N})
&\os{\sim}{\lo}\bigoplus_{j\geq \max \{-k,0\}} 
a^{(2j+k)}_{m,t*} 
(Ru_{\os{\circ}{X}{}^{(2j+k)}_{m,t}/{\cal W}(\kap_t)*}
(E_{\os{\circ}{X}{}^{(2j+k)}_{m,t}
/{\cal W}(\kap_t)} \tag{4.1.6.1}\label{ali:xnwxnqsnp}\\
& \otimes_{\mab Z}
\vp_{\rm crys}^{(2j+k)}
(\os{\circ}{X}_{m,t}/{\cal W}(\kap_t))))
\otimes_{\mab Z}{\mab Q}[-2j-k]. 
\end{align*}  
Hence, as in \S\ref{sec:flgdw}, 
we obtain the following exact sequence by (\ref{prop:mnb}): 
\begin{align*} 
E_1^{-k,q+k}
= &\bigoplus_{m=0}^N
\us{j\geq {\rm max}\{-(k+m),0\}}
{\bigoplus}
H^{q-2j-k-2m}
((\os{\circ}{X}{}^{(2j+k+m)}_{m,t}/{\cal W}(\kap_t))_{\rm crys}, 
E_{\os{\circ}{X}{}^{(2j+k+m)}_{m,t}/{\cal W}(\kap_t)} 
\tag{4.1.6.2}\label{eqn:trhskwsp}\\
& 
\otimes_{\mab Z}
\vp^{(2j+k+m)}_{\rm crys}(\os{\circ}{X}_{m,t}/{\cal W}_n))
(-j-k-m,u)_{\mab Q}  
\end{align*}
$$\Lo 
H^q
((X_{\bul \leq N,t}/{\cal W}(t))_{\rm crys},
\eps^*_{X_{\bul \leq N,t}/{\cal W}(t)}(E^{\bul \leq N}))_{\mab Q}
\quad ( n\in{\mab Z}_{>0}).
$$

\par
The following is the main result in this section:

\begin{theo}[{\bf Comparison theorem I}]\label{theo:cqsoncrdw}
Let $N$ be a nonnegative integer.  
Then the following hold$:$ 
\par 
$(1)$ In ${\rm D}^+{\rm F}(f^{-1}_{\bul \leq N}(K_0(\kap_t)))$ 
there exists the following 
canonical isomorphism$:$
\begin{align*}
(A_{{\rm zar},{\mab Q}}
(X_{\bul \leq N,\os{\circ}{t}}/{\cal W}(s_{\os{\circ}{t}}),E^{\bul \leq N}),P) 
\os{\sim}{\lo} ({\cal W}A_{X_{\bul \leq N,\os{\circ}{t}},{\mab Q}}(E^{\bul \leq N}),P).
\tag{4.1.7.1}\label{eqn:brntildw}
\end{align*} 
Let $t'\lo s'$ be an analogous morphism to $t\lo s$. 
The isomorphism $(\ref{eqn:brntildw})$ is 
contravariantly functorial for the commutative diagrams for 
{\rm (\ref{eqn:xdbquss})} and {\rm (\ref{cd:xdjqbxy})} for the case 
$S={\cal W}(s)$ snd $S'={\cal W}(s')$ 
and 
the following commutative diagram 
\begin{equation*} 
\begin{CD} 
t@>>>t' \\ 
@VVV @VVV \\ 
s@>>>s', 
\end{CD} 
\end{equation*} 
where we assume 
the existence of the immersion {\rm (\ref{cd:adypn})}. 
\par 
$(2)$ The isomorphism 
{\rm (\ref{eqn:brntildw})} forgetting the filtrations   
fits into the following commutative diagram$:$  
\begin{equation*} 
\begin{CD} 
A_{{\rm zar},{\mab Q}}
(X_{\bul \leq N,\os{\circ}{t}}/{\cal W}(s_{\os{\circ}{t}}),E^{\bul \leq N})
@>{(\ref{eqn:wltntpxa}),~\simeq}>> \\
@A{\theta_{X_{\bul \leq N,\os{\circ}{t}}/{\cal W}(s_{\os{\circ}{t}})/{\cal W}(s_{0,\os{\circ}{t}})} \wedge}A{\simeq}A\\
Ru_{X_{\bul \leq N,t}/{\cal W}(t)*}
(\eps^*_{X_{\bul \leq N,t}/{\cal W}(t)}(E^{\bul \leq N}))
\otimes_{\mab Z}{\mab Q}  @>{(\ref{eqn:ywnnny}),~\sim}>> 
\end{CD} 
\tag{4.1.7.2}\label{cd:axswl}
\end{equation*} 
\begin{equation*} 
\begin{CD} 
{\cal W}A_{X_{\bul \leq N,\os{\circ}{t}},{\mab Q}}(E^{\bul \leq N}) \\
@A{\simeq}A{\theta_{X_{\bul \leq N,\os{\circ}{t}}/{\cal W}(s_{\os{\circ}{t}})/{\cal W}(s_{0,\os{\circ}{t}})} \wedge}A\\
(\eps^*_{X_{\bul \leq N,t}/{\cal W}(t)}(E^{\bul \leq N}))_{{\cal W}(t)}
\otimes_{{\cal W}({\cal O}_{X_{\bul \leq N,t})}}
{\cal W}{\Om}^{\bul}_{X_{\bul \leq N,t}}\otimes_{\mab Z}{\mab Q}. 
\end{CD} 
\end{equation*} 
\end{theo}
\begin{proof} 
In this proof we construct only the morphism (\ref{eqn:brntildw}) 
because the rest of the proof is almost the same as that of 
(\ref{theo:csoncrdw}). 
\par 
Replace the immersion 
$X'_{N,\os{\circ}{t}}\os{\sus}{\lo} \ol{\cal P}{}'_N$ in 
(\ref{theo:csoncrdw}) with an immersion 
$X'_{N,\os{\circ}{t}} \os{\sus}{\lo} \ol{\cal P}{}'_N$ 
into a formally log smooth scheme over $\ol{{\cal W}(s_{\os{\circ}{t}})}$.
Then we have immersions 
$X_{\bul \leq N,\bul,\os{\circ}{t}} \os{\sus}{\lo} \ol{\cal P}_{\bul \leq N,\bul}$ 
and 
$X_{\bul \leq N,\bul,\os{\circ}{t}} \os{\sus}{\lo} {\cal P}_{\bul \leq N,\bul}$ 
into formally log smooth $(N,\infty)$-truncated bisimplicial 
formal log schemes over $\ol{{\cal W}(s_{\os{\circ}{t}})}$
and ${\cal W}(s_{\os{\circ}{t}})$, respectively. 
Let $\ol{\mathfrak D}_{\bul \leq N,\bul}$ 
be the $p$-adic completion of the log PD-envelope of 
the immersion $X_{\bul \leq N,\bul,\os{\circ}{t}}\os{\sus}{\lo}
\ol{\cal P}_{\bul \leq N,\bul}$ over 
$({\cal W}(\os{\circ}{t}),p{\cal W},[~])$. 
Let ${\cal P}^{\rm ex}_{\bul \leq N,\bul}$ 
be the exactification of the immersion 
$X_{\bul \leq N,\bul,\os{\circ}{t}} \os{\sus}{\lo} {\cal P}_{\bul \leq N,\bul}$.  
Set ${\cal P}_{\bul \leq N,\bul,n}:={\cal P}_{\bul \leq N,\bul}\mod p^n$. 
Set ${\mathfrak D}_{\bul \leq N,\bul}
:=\ol{\mathfrak D}_{\bul \leq N,\bul}
\times_{{\mathfrak D}(\ol{{\cal W}(s_{\os{\circ}{t}})})}{\cal W}(s_{\os{\circ}{t}})$.  
Let $E^{\bul \leq N,\bul}$ be a crystal of 
${\cal O}_{\os{\circ}{X}_{\bul \leq N,\bul,t}/{\cal W}(\os{\circ}{t})}$-modules 
obtained by $E^{\bul \leq N}$. 
Let $(\ol{\cal E}{}^{\bul \leq N,\bul},\ol{\nabla})$ 
be the coherent 
${\cal O}_{\ol{\mathfrak D}_{\bul \leq N,\bul}}$-module 
with integrable connection 
obtained by 
$\eps^*_{X_{\bul \leq N,\bul,\os{\circ}{t}}/{\cal W}(\os{\circ}{t})}
(E^{\bul \leq N,\bul})$. 
Set $({\cal E}^{\bul \leq N,\bul},\nabla)
:=(\ol{\cal E}{}^{\bul \leq N,\bul},\ol{\nabla})
\otimes_{{\cal O}_{{\mathfrak D}(\ol{{\cal W}(s_{\os{\circ}{t}})})}}{\cal W}$.  
Because we have the following commutative diagrams
\begin{equation*} 
\begin{CD} 
({\cal E}^{\bul \leq N,\bul}
{\otimes}_{{\cal O}_{{\cal P}^{\rm ex}_{\bul \leq N,\bul}}}
{\Om}^{i+j+2}_{{\cal P}^{\rm ex}_{\bul \leq N,\bul,n}
/{\cal W}_n(\os{\circ}{t})}/P_{j+1},P)  
@>>>  \\
@A{-\nabla}AA  \\
({\cal E}^{\bul \leq N,\bul}
{\otimes}_{{\cal O}_{{\cal P}^{\rm ex}_{\bul \leq N,\bul}}}
{\Om}^{i+j+1}_{{\cal P}^{\rm ex}_{\bul \leq N,\bul,n}
/{\cal W}_n(\os{\circ}{t})}/P_{j+1},P)  
@>>>    
\end{CD}
\tag{4.1.7.3}\label{ali:pnfnho}
\end{equation*}   
\begin{equation*} 
\begin{CD}
(E^{\bul \leq N,\bul}_n
\otimes_{{\cal W}_n{(\cal O}_{X_{\bul \leq N,\bul,\os{\circ}{t}}})}
{\cal W}_n\wt{{\Om}}^{i+j+2}_{X_{\bul \leq N,\bul,\os{\circ}{t}}}/P_{j+1},P)\\
@AA{-\nabla}A\\
(E^{\bul \leq N,\bul}_n 
\otimes_{{\cal W}_n{(\cal O}_{X_{\bul \leq N,\bul,\os{\circ}{t}}})}
{\cal W}_n\wt{{\Om}}^{i+j+1}_{X_{\bul \leq N,\bul,\os{\circ}{t}}}/P_{j+1},P),
\end{CD}
\end{equation*} 
\begin{equation*} 
\begin{CD} 
({\cal E}^{\bul \leq N,\bul}
{\otimes}_{{\cal O}_{{\cal P}^{\rm ex}_{\bul \leq N,\bul}}}
{\Om}^{i+j+2}_{{\cal P}^{\rm ex}_{\bul \leq N,\bul,n}
/{\cal W}_n(\os{\circ}{t})}/P_{j+2},P)  
@>>>  \\
@A{e\theta_{{\cal P}^{\rm ex}_{\bul \leq N,\bul}}\wedge }AA \\
({\cal E}^{\bul \leq N,\bul}
{\otimes}_{{\cal O}_{{\cal P}^{\rm ex}_{\bul \leq N,\bul}}}
{\Om}^{i+j+1}_{{\cal P}^{\rm ex}_{\bul \leq N,\bul,n}
/{\cal W}_n(\os{\circ}{t})}/P_{j+1},P) @>>>  
\end{CD}
\tag{4.1.7.4}\label{ali:pnflvo}
\end{equation*} 
\begin{equation*} 
\begin{CD}  
(E^{\bul \leq N,\bul}_n
\otimes_{{\cal W}_n{(\cal O}_{X_{\bul \leq N,\bul,\os{\circ}{t}}})}
{\cal W}_n\wt{{\Om}}^{i+j+2}_{X_{\bul \leq N,\bul,\os{\circ}{t}}}/P_{j+2},P) \\
@AA{e\theta_{{\cal P}^{\rm ex}_{\bul \leq N,\bul}} \wedge}A \\ 
(E^{\bul \leq N,\bul}_n 
\otimes_{{\cal W}_n{(\cal O}_{X_{\bul \leq N,\bul,\os{\circ}{t}}})}
{\cal W}_n\wt{{\Om}}^{i+j+1}_{X_{\bul \leq N,\bul,\os{\circ}{t}}}/P_{j+1},P)   
\end{CD}
\end{equation*} 
as in the proof of (\ref{theo:csoncrdw})   
and because these commutative diagrams are compatible with projections, 
we have the following filtered morphism 
\begin{equation*} 
(A_{{\rm zar},{\mab Q}}({\cal P}^{\rm ex}_{\bul \leq N,\bul,\os{\circ}{t}}/{\cal W}(t),E^{\bul \leq N}),P) 
\lo ({\cal W}A_{X_{\bul \leq N,\bul,\os{\circ}{t}},{\mab Q}}(E^{\bul \leq N}),P).
\tag{4.1.7.5}\label{eqn:caflmp}
\end{equation*} 
Applying $R\pi_{{\rm zar}*}$ to (\ref{eqn:caflmp}), 
we have the filtered morphism (\ref{eqn:brntildw}). 
The rest of the proof is the same as that of (\ref{theo:csoncrdw}). 
\end{proof}

\begin{theo}[{\bf Contravariant functoriality}]\label{theo:cwtrw} 
$(1)$ Let the notations be as in {\rm (\ref{theo:fugennas})}, 
where $S={\cal W}(s)$, $S'={\cal W}(s')$, 
$T={\cal W}(t)$ and $T'={\cal W}(t')$ 
$(s'$ is the log point of a perfect field of characteristic $p$ and 
$t$ and $t'$ are fine log schemes whose underlying schemes are 
the spectrums of perfect fields of characteristic $p)$. 
Let $\os{\circ}{g}{}^*_{{\rm crys}}(E^{\bul \leq N})\lo F^{\bul \leq N}$
be the morphism in {\rm (\ref{ali:gnafe})} for the case above.  
Then $g_{\bul \leq N}$ 
induces the following 
well-defined pull-back morphism 
\begin{align*}  
g^*_{\bul \leq N} \col &
({\cal W}A_{Y_{\bul\leq N,\os{\circ}{t}{}'},{\mab Q}}(F^{\bul \leq N}),P) 
&\lo Rg_{\bul \leq N*}
(({\cal W}A_{X_{\bul\leq N,\os{\circ}{t}},{\mab Q}}(E^{\bul \leq N}),P))
\tag{4.1.8.1}\label{eqn:fzntwd}\\
\end{align*} 
fitting into the following commutative diagram$:$
\begin{equation*} 
\begin{CD}
{\cal W}A_{Y_{\bul\leq N,\os{\circ}{t}{}'},{\mab Q}}(F^{\bul \leq N})
@>{g_{\bul \leq N}^*}>> \\ 
@A{\theta_{Y_{\bul \leq N,\os{\circ}{t}{}'}/
{\cal W}(s'_{\os{\circ}{t}{}'})/{\cal W}(s'_{0,\os{\circ}{t}{}'})} \wedge}A{\simeq}A \\
F^{\bul \leq N}_{\infty}
\otimes_{{\cal W}({\cal O}_{Y_{\bul \leq N,t}'})}
{\cal W}{\Om}_{Y_{\bul \leq N,t'}}^{\bul}
\otimes_{{\mab Z}}{\mab Q}
@>{g_{\bul \leq N}^*}>>
\end{CD}
\tag{4.1.8.2}\label{cd:pscatcz} 
\end{equation*} 
\begin{equation*} 
\begin{CD}
Rg_{\bul \leq N*}
({\cal W}A_{X_{\bul\leq N,\os{\circ}{t}},{\mab Q}}(E^{\bul \leq N}))\\ 
@A{\simeq}A{Rg_{\bul \leq N*}
(\theta_{X_{\bul \leq N,\os{\circ}{t}}/{\cal W}(s_{\os{\circ}{t}})/{\cal W}(s_{0,\os{\circ}{t}})} \wedge)}A \\
Rg_{\bul \leq N*}
((E^{\bul \leq N}_{\infty}
\otimes_{{\cal W}({\cal O}_{X_{\bul \leq N,t}})}
{\cal W}{\Om}_{X_{\bul \leq N,t}}^{\bul})\otimes_{{\mab Z}}{\mab Q}).
\end{CD}
\end{equation*}
\par 
$(2)$ Let the notations be as in $(1)$. 
Let $s''$, $t''$, $S''={\cal W}(s'')$ and $T''={\cal W}(t'')$ 
be similar objects to $s'$, $t'$, $S'={\cal W}(s')$ and $T'={\cal W}(t')$, 
respectively.  Let $v\col s'_{\os{\circ}{t}{}'}\lo s''_{\os{\circ}{t}{}''}$ 
be a similar morphism to $u$. 
Let $h_{\bul \leq N}\col Y_{\bul \leq N,\os{\circ}{t}{}'}\lo Z_{\bul \leq N,\os{\circ}{t}{}''}$ be 
a similar morphism over $v\col s'\lo s''$ 
to $g_{\bul \leq N}\col X_{\bul \leq N,\os{\circ}{t}}\lo Y_{\bul \leq N,\os{\circ}{t}{}'}$.  
Let 
\begin{align*} 
\os{\circ}{h}{}^*_{\bul \leq N,{\rm crys}}(G^{\bul \leq N})\lo F^{\bul \leq N} 
\tag{4.1.8.3}\label{ali:rsswtpm}
\end{align*} 
be as in {\rm (\ref{ali:gnafe})}.  
Then 
\begin{align*} 
& (h_{\bul \leq N}\circ g_{\bul \leq N})^* =
Rh_{\bul \leq N*}(g_{\bul \leq N}^*)\circ 
h_{\bul \leq N}^* \col ({\cal W}A_{Z_{\bul\leq N,\os{\circ}{t}{}''}}(G^{\bul \leq N}),P) \lo 
\tag{4.1.8.4}\label{ali:piwdpp} \\ 
& Rh_{\bul \leq N*}
Rg_{\bul \leq N*}({\cal W}A_{X_{\bul\leq N,\os{\circ}{t}}}(E^{\bul \leq N}),P)  =R(h_{\bul \leq N}\circ g_{\bul \leq N})_*({\cal W}A_{X_{\bul\leq N,\os{\circ}{t}}}(E^{\bul \leq N}),P).
\end{align*}  
\par 
$(3)$ \begin{equation*} 
{\rm id}_{X_{\bul \leq N,\os{\circ}{t}}}^*={\rm id} 
\col ({\cal W}A_{X_{\bul\leq N,\os{\circ}{t}}}(E^{\bul \leq N}),P)
\lo ({\cal W}A_{X_{\bul\leq N,\os{\circ}{t}}}(E^{\bul \leq N}),P). 
\tag{4.1.8.5}\label{eqn:fziwdxd}
\end{equation*} 
\end{theo} 
\begin{proof} 
It is much easier than the proofs of {\rm (\ref{theo:funas})} and 
(\ref{theo:funpas}). We leave the detail of the proof to the reader. 
\end{proof}

\par
Let $i$ be a fixed nonnegative integer. 
Let $N$ be a nonnegative integer. 
We conclude this section by constructing a 
preweight spectral sequence of 
$E^{\bul \leq N}
\otimes_{{\cal W}({\cal O}_{X_{\bul \leq N,\os{\circ}{t}}})}
{\cal W}{\Om}^i_{X_{\bul \leq N,\os{\circ}{t}}} 
\otimes_{\mab Z}{\mab Q}$ 
and describing the boundary morphisms 
between the $E_1$-terms 
of  the spectral sequence.
\par 
The proof of the following is 
the same as that of (\ref{theo:hwwt}):  

\begin{theo}\label{theo:hbwwt} 
There exists the following spectral sequence$:$ 
\begin{align*}
&E_1^{-k,q+k}=\bigoplus_{m=0}^N
\bigoplus_{j\geq \max \{-(k+m),0\}} 
H^{q-i-j-m}(\os{\circ}{X}{}^{(2j+k+m)}_{m,t},
(E^m_{\os{\circ}{X}{}^{(2j+k+m)}_{m,t}/{\cal W}(t)})_
{{\cal W}(\os{\circ}{X}{}^{(2j+k+m)}_{m,t})} 
\tag{4.1.9.1}\label{eqn:sesplasp} \\
&\otimes_{{\cal W}({\cal O}_{\os{\circ}{X}{}^{(2j+k+m)}_{m}})}
{\cal W}\Om^{i-j-k-m}_{\os{\circ}{X}{}^{(2j+k+m)}_{m}}
\otimes_{\mab Z}{\mab Q}
\otimes_{\mab Z} 
\vp^{(2j+k+m)}_{\rm zar}(\os{\circ}{X}_{m,t}/\os{\circ}{t})))(-j-k-m,u)  \\
& \Lo 
H^{q-i}(X_{\bul \leq N,t},E^{\bul \leq N}_{\infty}
\otimes_{{\cal W}({\cal O}_{X_{\bul \leq N,t}})}
{\cal W}{\Om}^i_{X_{\bul \leq N,t}} \otimes_{\mab Z}{\mab Q})
\quad (q\in {\mab Z}). 
\end{align*}
\end{theo}

\begin{defi}\label{defi:spw}
We call (\ref{eqn:sesplasp}) the 
{\it weight spectral sequence} of 
$E^{\bul \leq N}_{\infty}
\otimes_{{\cal W}({\cal O}_{X_{\bul \leq N,t}})}
{\cal W}{\Om}^i_{X_{\bul \leq N,t}}
\otimes_{\mab Z}{\mab Q}$. 
We define the {\it weight filtration} 
$P$  on $H^{q-i}(X_{\bul \leq N,t},
E^{\bul \leq N}_{\infty}
\otimes_{{\cal W}({\cal O}_{X_{\bul \leq N,t}})}
{\cal W}{\Om}^i_{X_{\bul \leq N,t}})$ 
as follows: 
$${\rm gr}_k^PH^{q-i}(X_{\bul \leq N,t},
E^{\bul \leq N}_{\infty}
\otimes_{{\cal W}({\cal O}_{X_{\bul \leq N,t}})}
{\cal W}{\Om}^i_{X_{\bul \leq N,t}}
\otimes_{\mab Z}{\mab Q}) =
E_{\infty}^{q-k,k}.$$
\end{defi}

Fix a total order on $\Lam$ once and for all. 
As in (\ref{prop:deccbd}), we have the following: 

\begin{prop}\label{prop:dsschwbd} 
The boundary morphism between the $E_1$-terms of 
the spectral sequence {\rm (\ref{eqn:sesplasp})} 
is given by the following diagram$:$ 
\begin{equation*} 
\tag{4.1.11.1}\label{cd:gbshwsd}
\end{equation*} 
\begin{equation*} 
\begin{split} 
{} & \bigoplus_{m=0}^N
\bigoplus_{j\geq \max \{-(k+m),0\}} 
H^{q-i-j-m}(\os{\circ}{X}{}^{(2j+k+m)}_{m+1,t},
(E^{m+1}_{\os{\circ}{X}{}^{(2j+k+m)}_{m+1,t}/{\cal W}(t)})_
{{\cal W}(\os{\circ}{X}{}^{(2j+k+m)}_{m+1,t})} 
 \\ 
{} & 
\otimes_{{\cal W}({\cal O}_{\os{\circ}{X}{}^{(2j+k+m)}_{m+1,t}})}
{\cal W}\Om^{i-j-k-m}_{\os{\circ}{X}{}^{(2j+k+m)}_{m+1,t}}))
\otimes_{\mab Z}
\vp^{(2j+k+m)}_{\rm zar}
(\os{\circ}{X}_{m+1,t}/\os{\circ}{t}))
(-j-k-m,u)_{\mab Q}
\end{split} 
\end{equation*}  
$$\text{\scriptsize
{$\sum_{i= 0}^{m+1}(-1)^i\del^i_{j}$}}
~\uparrow~(1 \leq j \leq r)$$
\begin{equation*} 
\begin{split} 
{} & \bigoplus_{m=0}^N
\bigoplus_{j\geq \max \{-(k+m),0\}} 
H^{q-i-j-m}(\os{\circ}{X}{}^{(2j+k+m)}_m,
(E^m_{\os{\circ}{X}{}^{(2j+k+m)}_{m,t}/{\cal W}(t)})_
{{\cal W}(\os{\circ}{X}{}^{(2j+k+m)}_{m,t})}  \\ 
{} & 
\otimes_{{\cal W}({\cal O}_{\os{\circ}{X}{}^{(2j+k+m)}_{m}})}
{\cal W}\Om^{i-j-k-m}_{\os{\circ}{X}{}^{(2j+k+m)}_{m}}
\otimes_{\mab Z}
\vp^{(2j+k+m)}_{\rm zar}
(\os{\circ}{X}_{m,t}/\os{\circ}{t}))
(-j-k-m,u)_{\mab Q}
\end{split} 
\end{equation*} 
$$\text{\scriptsize
{${\sum_{j\geq {\rm max}\{-(k+m),0\}}
(-1)^m[G_m+\rho_m]}$}}
~\downarrow \quad \quad \quad \quad \quad \quad \quad 
\quad \quad \quad \quad \quad$$
\begin{equation*} 
\begin{split} 
{} & \bigoplus_{m=0}^N
\bigoplus_{j\geq \max \{-(k+m)+1,0\}} 
H^{q-i-j-m+1}(\os{\circ}{X}{}^{(2j+k+m-1)}_{m,t},
(E^m_{\os{\circ}{X}{}^{(2j+k+m-1)}_{m,t}/{\cal W}(t)})_
{{\cal W}(\os{\circ}{X}{}^{(2j+k+m-1)}_{m,t})}  \\ 
{} & 
\otimes_{{\cal W}({\cal O}_{\os{\circ}{X}{}^{(2j+k+m-1)}_{m,t}})}
{\cal W}\Om^{i-j-k-m+1}_{\os{\circ}{X}{}^{(2j+k+m-1)}_{m,t}}
\otimes_{\mab Z}
\vp^{(2j+k+m-1)}_{\rm zar}
(\os{\circ}{X}_{m,t}/\os{\circ}{t}))
(-j-k-m+1,u)_{\mab Q}. 
\end{split} 
\end{equation*}  
Here $G_m$ is the \v{C}ech-Gysin morphism 
in Hodge-Witt cohomologies {\rm \cite[(4.4.2)]{ndw}} 
which is defined analogously as in {\rm (\ref{eqn:togsn})} 
and $\rho_m$ is 
the morphism obtained by the pull-back of the closed immersions in Hodge-Witt cohomologies 
which is defined analogously as in {\rm (\ref{eqn:rhsgsn})}.  
\end{prop}

\par 
By (\ref{eqn:wenlby}) we have the following  exact sequence  
\begin{align*} 
0 \lo &E^{\bul \leq N}_{\infty}\otimes_{{\cal W}({\cal O}_{X_{\bul \leq N,t})}}
{\cal W}\Om^{\bul}_{X_{\bul \leq N,t}}
\otimes_{\mab Z}{\mab Q}(-1,u)[-1]  
\os{\theta_{X_{\bul \leq N,\os{\circ}{t}}
/{\cal W}(s_{\os{\circ}{t}})/{\cal W}(s_{0,\os{\circ}{t}})} \wedge}{\lo} 
\tag{4.1.11.2}\label{eqn:wnltby} \\
&E^{\bul \leq N}_{\infty}\otimes_{{\cal W}({\cal O}_{X_{\bul \leq N,t})}}
{\cal W}\wt{{\Om}}^{\bul}_{X_{\bul \leq N,\os{\circ}{t}}} 
\otimes_{\mab Z}{\mab Q}
\lo E^{\bul \leq N}_{\infty}\otimes_{{\cal W}({\cal O}_{X_{\bul \leq N,t})}}
{\cal W}\Om^{\bul}_{X_{\bul \leq N,t}}
\otimes_{\mab Z}{\mab Q}
\lo 0.  
\end{align*}
\par
Using the extension (\ref{eqn:wnltby}), 
we have the following boundary morphism 
\begin{equation*} 
e^{-1}N_{\rm dRW} \col 
E^{\bul \leq N}_{\infty}\otimes_{{\cal W}({\cal O}_{X_{\bul \leq N,t})}}
{\cal W}\Om^{\bul}_{X_{\bul \leq N,t}}
\otimes_{\mab Z}{\mab Q}
\lo 
E^{\bul \leq N}_{\infty}\otimes_{{\cal W}({\cal O}_{X_{\bul \leq N,t})}}
{\cal W}\Om^{\bul}_{X_{\bul \leq N,t}}
\otimes_{\mab Z}{\mab Q}. 
\end{equation*}

\begin{defi} 
We call $e^{-1}N_{\rm dRW}$ the {\it reverse de Rham-Witt  monodromy 
operator} of $X_{\bul \leq N,t}/t$. 
\end{defi}

\par 
Let 
\begin{equation*}
\nu_{\rm dRW} \col 
{\cal W}A_{X_{\bul \leq N,\os{\circ}{t}},{\mab Q}}(E^{\bul \leq N})
\lo {\cal W}A_{X_{\bul \leq N,\os{\circ}{t}},{\mab Q}}(E^{\bul \leq N})(-1,u)
\tag{4.1.12.1}\label{eqn:wadnsu}
\end{equation*}
be the induced morphism 
by morphisms 
$${\rm proj}. 
\col {\cal W}A_{X_{\bul \leq N,\os{\circ}{t}},{\mab Q}}(E^{\bul \leq N})^{ij}
\lo {\cal W}A_{X_{\bul \leq N,\os{\circ}{t}},{\mab Q}}(E^{\bul \leq N})^{i-1,j+1}.$$
Note that $\nu_{\rm dRW}$ is actually a morphism of 
complexes. 

\begin{defi}\label{defi:rrw} 
We call $\nu_{\rm dRW}$ the {\it reverse de Rham-Witt 
quasi-monodromy operator}.
\end{defi}

We omit the proofs of the following propositions and corollaries 
because we have only to make suitable modifications 
of proofs in \S\ref{sec:vpmn}.

\begin{prop}\label{prop:e1n}
The following diagram is commutative$:$
\begin{equation*}
\begin{CD}
Rg_{\bul \leq N*}(E^{\bul \leq N}_{\infty}\otimes_{{\cal W}({\cal O}_{X_{\bul \leq N,t})}}
{\cal W}\Om^{\bul}_{X_{\bul \leq N,t}}\otimes_{\mab Z}{\mab Q})
@>{e^{-1}N_{\rm dRW}}>> 
Rg_{\bul \leq N*}(E^{\bul \leq N}_{\infty}\otimes_{{\cal W}({\cal O}_{X_{\bul \leq N,t})}}
{\cal W}\Om^{\bul}_{X_{\bul \leq N,t}}\otimes_{\mab Z}{\mab Q})\\
@A{g_{\bul \leq N}^*}AA @AA{g_{\bul \leq N}^*}A  \\
F^{\bul \leq N}_{\infty}\otimes_{{\cal W}({\cal O}_{Y_{\bul \leq N,t'})}}
{\cal W}\Om^{\bul}_{X_{\bul \leq N,t}}\otimes_{\mab Z}{\mab Q}
@>{e^{-1}N_{\rm dRW}}>> 
F^{\bul \leq N}_{\infty}\otimes_{{\cal W}({\cal O}_{Y_{\bul \leq N,t'})}}
{\cal W}\Om^{\bul}_{Y_{\bul \leq N,t'}}\otimes_{\mab Z}{\mab Q}. 
\end{CD}
\tag{4.1.14.1}\label{cd:nuaseqln}
\end{equation*} 
\end{prop}

\begin{prop}\label{prop:qaaimon}
Let $k$ be a positive integer and let $q$ be a nonnegative integer. 
Then 
\begin{align*} 
\nu_{{\rm dRW}}^k \col  & 
R^qf_{*}
({\rm gr}_k^P{\cal W}A_{X_{\bul \leq N,\os{\circ}{t}},{\mab Q}}(E^{\bul \leq N}))
\lo 
R^qf_{*}
({\rm gr}_{-k}^P{\cal W}A_{X_{\bul \leq N,\os{\circ}{t}},{\mab Q}}(E^{\bul \leq N}))(-k,u)
\tag{4.1.15.1}\label{eqn:kawpa}
\end{align*} 
is the identity. 
\end{prop}

\begin{prop}\label{prop:nucaswa} 
The quasi-monodromy operator 
\begin{equation*} 
\nu_{\rm dRW} \col 
{\cal W}A_{X_{\bul \leq N,\os{\circ}{t}},{\mab Q}}(E^{\bul \leq N})  \lo 
{\cal W}A_{X_{\bul \leq N,\os{\circ}{t}},{\mab Q}}(E^{\bul \leq N})(-1,u)
\end{equation*}
underlies the following morphism 
of filtered morphism 
\begin{equation*} 
\nu_{{\rm dRW}} \col 
({\cal W}A_{X_{\bul \leq N,\os{\circ}{t}},{\mab Q}}(E^{\bul \leq N}),P)
\lo  
({\cal W}A_{X_{\bul \leq N,\os{\circ}{t}},{\mab Q}}(E^{\bul \leq N}),P\langle -2 \rangle)(-1,u), 
\end{equation*}
where $P\langle -2 \rangle$ is a filtration 
defined by $(P\langle -2 \rangle)_k= P_{k-2}$. 
\end{prop}

\begin{coro}\label{coro:zmawop}
The quasi-monodromy operator {\rm (\ref{defi:rrw})} 
induces the following morphism 
\begin{align*} 
\nu_{{\rm dRW}} \col & P_kRf_{X_{\bul \leq N,\os{\circ}{t}}/{\cal W}(s_{\os{\circ}{t}})*}
(\eps^*_{X_{\bul \leq N,\os{\circ}{t}}/{\cal W}(s_{\os{\circ}{t}})}(E^{\bul \leq N}))
\otimes_{\mab Z}{\mab Q}
\lo  \tag{4.1.17.1}\label{eqn:grpwapd}\\
& P_{k-2}Rf_{X_{\bul \leq N,\os{\circ}{t}}/{\cal W}(s_{\os{\circ}{t}})*}
(\eps^*_{X_{\bul \leq N,\os{\circ}{t}}/{\cal W}(s_{\os{\circ}{t}})}(E^{\bul \leq N}))(-1,u)
\otimes_{\mab Z}{\mab Q}
\end{align*}
\end{coro}

As in \cite[$(10.6.1;\star)$]{ndw} we can easily 
prove the following:

\begin{prop}\label{prop:qnseqn}
There exists the following 
commutative diagram$:$
\begin{equation*}
\begin{CD}
A_{{\rm zar},{\mab Q}}
(X_{\bul \leq N,\os{\circ}{t}}/{\cal W}(s_{\os{\circ}{t}}),E^{\bul \leq N})
@>{\nu_{\rm zar}}>> \\
@A{\theta_{X_{\bul \leq N,\os{\circ}{t}}
/{\cal W}(s_{\os{\circ}{t}})/{\cal W}(s_{0,\os{\circ}{t}})}
 \wedge}A{\simeq}A   \\
E^{\bul \leq N}\otimes_{{\cal W}({\cal O}_{X_{\bul \leq N,t})}}
{\cal W}\Om^{\bul}_{X_{\bul \leq N,t}}\otimes_{\mab Z}{\mab Q}
@>{e^{-1}N_{\rm dRW}}>> 
\end{CD}
\tag{4.1.18.1}\label{cd:nuseqln}
\end{equation*} 
\begin{equation*}
\begin{CD}
A_{{\rm zar},{\mab Q}}
(X_{\bul \leq N,\os{\circ}{t}}/{\cal W}(s_{\os{\circ}{t}}),E^{\bul \leq N})(-1,u) \\
@A{\simeq}A{\theta_{X_{\bul \leq N,\os{\circ}{t}}
/{\cal W}(s_{\os{\circ}{t}})/{\cal W}(s_{0,\os{\circ}{t}})}\wedge}A \\
E^{\bul \leq N}_t\otimes_{{\cal W}({\cal O}_{X_{\bul \leq N,t})}}
{\cal W}\Om^{\bul}_{X_{\bul \leq N,t}}\otimes_{\mab Z}{\mab Q}(-1,u).
\end{CD} 
\end{equation*} 
\end{prop}

The following is obvious 
by the definition of $\nu_{\rm zar}$ and 
$\nu_{\rm dRW}$: 

\begin{prop}  
\begin{equation*}
\begin{CD}
A_{{\rm zar},{\mab Q}}
(X_{\bul \leq N,\os{\circ}{t}}/{\cal W}(s_{\os{\circ}{t}}),E^{\bul \leq N})
@>{\nu_{\rm zar}}>> A_{{\rm zar},{\mab Q}}
(X_{\bul \leq N,\os{\circ}{t}}/{\cal W}(s_{\os{\circ}{t}}),E^{\bul \leq N})
(-1,u)\\
@V{\simeq}VV  @VV{\simeq}V  \\
{\cal W}A_{X_{\bul \leq N,\os{\circ}{t}},{\mab Q}}(E^{\bul \leq N})
@>{\nu_{\rm dRW}}>> 
{\cal W}A_{X_{\bul \leq N,\os{\circ}{t}},{\mab Q}}(E^{\bul \leq N}))(-1,u).
\end{CD}
\tag{4.1.19.1}\label{cd:nuesdwqln}
\end{equation*} 
\end{prop}

\par 
We also omit the proof of the following

\begin{prop}\label{prop:spqaqqqq} 
Let the notations be as in {\rm (\ref{theo:cwtrw})}.  
Then the following diagram is commutative$:$ 
\begin{equation*} 
\begin{CD} 
Rg_{{\bul \leq N}*}(({\cal W}_{\star}A_{X_{\bul \leq N,\os{\circ}{t}},{\mab Q}}(E^{\bul \leq N}),P))
@>{\nu_{{\rm dRW},\star}}>> 
Rg_{{\bul \leq N}*}(({\cal W}_{\star}A_{X_{\bul \leq N,\os{\circ}{t}},{\mab Q}}(E^{\bul \leq N}),
P\langle -2 \rangle))(-1,u)  \\ 
@A{g^*_{\bul \leq N}}AA @AA{g^*_{\bul \leq N}}A \\
({\cal W}_{\star}A_{Y_{\bul \leq N,\os{\circ}{t}{}'},{\mab Q}}(F^{\bul \leq N}),P)
@>{\nu_{{\rm dRW},\star}}>>
({\cal W}_{\star}A_{Y_{\bul \leq N,\os{\circ}{t}{}'},{\mab Q}}(F^{\bul \leq N}),P\langle -2 \rangle)(-1,u'). 
\end{CD} 
\tag{4.1.20.1}\label{eqn:aawasf}
\end{equation*} 
\end{prop}

\par 
Let us define the double complex 
${\cal W}B_{X_{\bul \leq N,\os{\circ}{t}},{\mab Q}}(E^{\bul \leq N})^{\bul \bul}$ 
as follows.
\par  
The $(i,j)$-component ${\cal W}B_{X_{\bul \leq N,\os{\circ}{t}},{\mab Q}}
(E^{\bul \leq N})^{ij}$ 
$(i,j\in {\mab Z}_{\geq 0})$ is defined by the similar formula to (\ref{ali:sgbxn}). 
The horizontal boundary morphism
$d' \col {\cal W}B_{X_{\bul \leq N,\os{\circ}{t}},{\mab Q}}(E^{\bul \leq N})^{ij} \lo 
{\cal W}B_{X_{\bul \leq N,\os{\circ}{t}},{\mab Q}}(E^{\bul \leq N})^{i+1,j}$ 
is, by definition,
$$d'(\om_1, \om_2)=(\nabla \om_1,-\nabla\om_2) $$
and the vertical one 
$d'' \col {\cal W}B_{X_{\bul \leq N,\os{\circ}{t}},{\mab Q}}(E^{\bul \leq N})^{ij} \lo 
{\cal W}B_{X_{\bul \leq N,\os{\circ}{t}},{\mab Q}}(E^{\bul \leq N})^{i,j+1}$ is 
$$d''(\om_1, \om_2)=
(-\theta_{X_{\bul \leq N,\os{\circ}{t}}/{\cal W}(s_{\os{\circ}{t}})/
{\cal W}(s_{0,\os{\circ}{t}})}\wedge \om_1+\nu_{\rm dRW}(\om_2),
\theta_{X_{\bul \leq N,\os{\circ}{t}}/{\cal W}(s_{\os{\circ}{t}})/
{\cal W}(s_{0,\os{\circ}{t}})}\wedge \om_2).$$ 
Let 
${\cal W}B_{X_{\bul \leq N,\os{\circ}{t}},{\mab Q}}(E^{\bul \leq N})$ be 
the single complex of 
${\cal W}B_{X_{\bul \leq N,\os{\circ}{t}},{\mab Q}}(E^{\bul \leq N})^{\bul \bul}$. 
\par
Let 
\begin{align*} 
\mu_{{X_{\bul \leq N,\os{\circ}{T}_0}/\os{\circ}{T}}} \col E^{\bul \leq N}_{\infty}
\otimes_{{\cal W}({\cal O}_{X_{\bul \leq N,\os{\circ}{t}}})}
{\cal W}\wt{{\Om}}^{\bul}_{X_{\bul \leq N,\os{\circ}{t}}}\otimes_{\mab Z}{\mab Q}
\lo {\cal W}B_{X_{\bul \leq N,\os{\circ}{t}},{\mab Q}}(E^{\bul \leq N})
\end{align*} 
be a morphism of complexes defined by 
\begin{align*} 
\mu_{{X_{\bul \leq N,\os{\circ}{T}_0}/\os{\circ}{T}}}(\om):=&(\om~{\rm mod}~P_0,
\theta_{X_{\bul \leq N,\os{\circ}{t}}
/{\cal W}(s_{\os{\circ}{t}})/{\cal W}(s_{0,\os{\circ}{t}})} 
\wedge \om~{\rm mod}~P_0) \\
&(\om \in  E^{\bul \leq N}_{\infty}
\otimes_{{\cal W}({\cal O}_{X_{\bul \leq N,\os{\circ}{t}}})}
{\cal W}\wt{{\Om}}^{\bul}_{X_{\bul \leq N,\os{\circ}{t}}}\otimes_{\mab Z}{\mab Q}).
\end{align*} 
Then we have the following commutative diagram of exact sequences:
\begin{equation*}
\begin{CD}
0 @>>> {\cal W}A_{X_{\bul \leq N,\os{\circ}{t}},{\mab Q}}(E^{\bul \leq N})(-1,u)[-1]   
 \\
@. @A{(\theta_{X_{\bul \leq N,\os{\circ}{t}}/{\cal W}(s_{\os{\circ}{t}})/
{\cal W}(s_{0,\os{\circ}{t}})}\wedge *)[-1] }AA   \\
0 @>>> (\eps^*_{X_{\bul \leq N,\os{\circ}{t}}/{\cal W}(\os{\circ}{t})}
(E^{\bul \leq N}))_{{\cal W}(X_{\bul \leq N,\os{\circ}{t}})}
\otimes_{{\cal W}({\cal O}_{X_{\bul \leq N,\os{\circ}{t}}})}
{\cal W}{\Om}^{\bul}_{X_{\bul \leq N,\os{\circ}{t}}}\otimes_{\mab Z}{\mab Q}(-1,u)[-1]   \\
\end{CD}
\tag{4.1.20.2}\label{cd:sgsttn}
\end{equation*}
\begin{equation*}
\begin{CD}
@>{}>>{\cal W}B_{X_{\bul \leq N,\os{\circ}{t}},{\mab Q}}(E^{\bul \leq N}) @>{}>>  \\
@. @A{\mu_{{X_{\bul \leq N,\os{\circ}{T}_0}/\os{\circ}{T}}}}AA\\
@>{\theta_{X_{\bul \leq N,\os{\circ}{t}}/{\cal W}(s_{\os{\circ}{t}})/
{\cal W}(s_{0,\os{\circ}{t}})}\wedge}>>\ol{E}{}^{\bul \leq N}
\otimes_{{\cal W}({\cal O}_{X_{\bul \leq N,\os{\circ}{t}}})}
{\cal W}\wt{{\Om}}^{\bul}_{X_{\bul \leq N,\os{\circ}{t}}}
\otimes_{\mab Z}{\mab Q}@>{}>>  
\end{CD} 
\end{equation*}
\begin{equation*}
\begin{CD}
{\cal W}A_{X_{\bul \leq N,\os{\circ}{t}},{\mab Q}}(E^{\bul \leq N})
@>>> 0 \\
@A{\theta_{X_{\bul \leq N,\os{\circ}{t}}/{\cal W}(s_{\os{\circ}{t}})/
{\cal W}(s_{0,\os{\circ}{t}})}\wedge}AA\\
(\eps^*_{X_{\bul \leq N,\os{\circ}{t}}/{\cal W}(\os{\circ}{t})}
(E^{\bul \leq N}))_{{\cal W}(X_{\bul \leq N,\os{\circ}{t}})}
\otimes_{{\cal W}({\cal O}_{X_{\bul \leq N,\os{\circ}{t}}})}
{\cal W}{\Om}^{\bul}_{X_{\bul \leq N,\os{\circ}{t}}}\otimes_{\mab Z}{\mab Q}@>>> 0.  
\end{CD} 
\end{equation*}
(By (\ref{eqn:wby}) the lower sequence is exact.) 
\par 
Set 
\begin{align*} 
P_k{\cal W}B_{X_{\bul \leq N,\os{\circ}{t}},{\mab Q}}(E^{\bul \leq N})^{ij}
:=&P_k{\cal W}A_{X_{\bul \leq N,\os{\circ}{t}},{\mab Q}}(E^{\bul \leq N})^{i-1,j}(-1,u)
\oplus P_k{\cal W}A_{X_{\bul \leq N,\os{\circ}{t}},{\mab Q}}(E^{\bul \leq N})^{ij} \\
& (i,j\in {\mab N}).
\end{align*} 
\begin{prop}\label{prop:bqid}
$(1)$ There exists the following sequence of the triangles in 
${\rm C}^+{\rm F}(f^{-1}_{\bul \leq N}({\cal K}_T)):$ 
\begin{align*} 
&\lo ({\cal W}A_{X_{\bul \leq N,\os{\circ}{t}},{\mab Q}}(E^{\bul \leq N}),P)[-1]
\lo ({\cal W}B_{X_{\bul \leq N,\os{\circ}{t}},{\mab Q}}(E^{\bul \leq N}),P)\\
&\lo ({\cal W}A_{X_{\bul \leq N,\os{\circ}{t}},{\mab Q}}(E^{\bul \leq N}),P)\os{+1}{\lo}. 
\end{align*} 
\end{prop} 

\begin{defi} 
We call $({\cal W}B_{X_{\bul \leq N,\os{\circ}{t}},{\mab Q}}(E^{\bul \leq N}),P)$ 
the {\it extended zariskian $p$-adic filtered Steenbrink complexes} of 
$E^{\bul \leq N}$ for $X_{\bul \leq N,\os{\circ}{t}}/{\cal W}(s_{\os{\circ}{t}})$. 
When $E^{\bul \leq N}={\cal O}_{\os{\circ}{X}_{\bul \leq N,t}/\os{\circ}{t}}$, 
we denote it by 
$({\cal W}B_{X_{\bul \leq N,\os{\circ}{t}},{\mab Q}},P)$ and call this 
the {\it extended zariskian $p$-adic filtered Steenbrink complex} of 
$X_{\bul \leq N,\os{\circ}{t}}/{\cal W}(s_{\os{\circ}{t}})$. 
\end{defi} 

\begin{theo}[{\bf Contravariant functoriality of 
${\cal W}B_{{\rm zar},{\mab Q}}$}]
\label{theo:fuqncb}
$(1)$ Let the notations be as in {\rm \S\ref{sec:pfnsc}}.  
Assume that $S=s$, $S'=s'$, $T={\cal W}(t)$ and $T'={\cal W}(t')$. 
Then $g_{\bul \leq N}\col X_{\bul \leq N,\os{\circ}{t}}\lo Y_{\bul \leq N,\os{\circ}{t}{}'}$ 
induces the following 
well-defined pull-back morphism 
\begin{equation*}  
g_{\bul \leq N}^* \col 
({\cal W}B_{Y_{\bul \leq N,\os{\circ}{t}{}'},{\mab Q}}(F^{\bul \leq N})_{\mab Q},P)
\lo Rg_{\bul \leq N*}(({\cal W}B_{X_{\bul \leq N,\os{\circ}{t}},{\mab Q}}
(E^{\bul \leq N})_{\mab Q},P)) 
\tag{4.1.23.1}\label{eqn:fzqaxd}
\end{equation*} 
fitting into the following commutative diagram$:$
\begin{equation*} 
\begin{CD}
{\cal W}B_{Y_{\bul \leq N,\os{\circ}{t}{}',{\mab Q}},{\mab Q}}(F^{\bul \leq N})
@>{g_{\bul \leq N}^*}>>  
Rg_{\bul \leq N*}({\cal W}B_{X_{\bul \leq N,\os{\circ}{t}},{\mab Q}}(E^{\bul \leq N}))\\ 
@A{\mu_{Y_{\bul \leq N,\os{\circ}{t}}/{\cal W}(\os{\circ}{t})}\wedge}A{\simeq}A 
@A{Rg_{\bul \leq N*}(\mu_{{X_{\bul \leq N,\os{\circ}{T}_0}/\os{\circ}{T}}})}A{\simeq}A\\
\ol{F}{}^{\bul \leq N}_{\infty}
\otimes_{{\cal W}_{\star}({\cal O}_{Y_{\bul \leq N,\os{\circ}{t}{}'}})}
{\cal W}_{\star}\wt{{\Om}}^{\bul}_{Y_{\bul \leq N,\os{\circ}{t}{}'}}
\otimes_{\mab Z}{\mab Q}
@>{g_{\bul \leq N}^*}>>Rg_{\bul \leq N*}(\ol{E}{}^{\bul \leq N}_{\infty}
\otimes_{{\cal W}_{\star}({\cal O}_{X_{\bul \leq N,\os{\circ}{t}}})}
{\cal W}_{\star}\wt{{\Om}}^{\bul}_{X_{\bul \leq N,\os{\circ}{t}}}\otimes_{\mab Z}{\mab Q}).
\end{CD}
\tag{4.1.23.2}\label{cd:psqwcz} 
\end{equation*}
\par 
$(2)$ 
\begin{align*} 
(h_{\bul \leq N}\circ g_{\bul \leq N})^*  =&
Rh_{\bul \leq N*}(g_{\bul \leq N}^*)\circ h_{\bul \leq N}^*    
\col ({\cal W}_{\star}B_{Y_{\bul \leq N,\os{\circ}{t}{}'},{\mab Q}}(F^{\bul \leq N}),P)
\tag{4.1.23.3}\label{ali:pwqwp} \\ 
& \lo Rh_{\bul \leq N*}Rg_{\bul \leq N*}
({\cal W}_{\star}B_{X_{\bul \leq N,\os{\circ}{t}},{\mab Q}}(E^{\bul \leq N}),P) \\
& =R(h_{\bul \leq N}\circ g_{\bul \leq N})_*
({\cal W}_{\star}B_{X_{\bul \leq N,\os{\circ}{t}},{\mab Q}}(E^{\bul \leq N}),P).
\end{align*}  
\par 
$(3)$  
\begin{equation*} 
{\rm id}_{X_{\bul \leq N,\os{\circ}{T}_0}}^*={\rm id} 
\col ({\cal W}_{\star}B_{X_{\bul \leq N,\os{\circ}{t}},{\mab Q}}(E^{\bul \leq N}),P)
\lo ({\cal W}_{\star}B_{X_{\bul \leq N,\os{\circ}{t}},{\mab Q}}(E^{\bul \leq N}),P).  
\tag{4.1.23.4}\label{eqn:fzwqdd}
\end{equation*} 
\end{theo} 

\begin{coro}\label{coro:spqfs}
Let the notations be as in {\rm (\ref{theo:cwtrw})}.
The isomorphisms {\rm (\ref{eqn:wltntpxa})} 
for $X_{\bul \leq N,\os{\circ}{t}}/{\cal W}(s_{\os{\circ}{t}})$, $E^{\bul \leq N}$ 
and $Y_{\bul \leq N,\os{\circ}{t}{}'}/{\cal W}(s'_{\os{\circ}{t}{}'})$, $F^{\bul \leq N}$ 
induce a morphism $g_{\bul \leq N}^*$ from {\rm (\ref{cd:nuaseqln})} 
to {\rm (\ref{eqn:aawasf})}.
This morphism satisfies the transitive relation and  
${\rm id}_{X_{\bul \leq N,\os{\circ}{t}}/{\cal W}(s_{\os{\circ}{t}})}^*={\rm id}$. 
\end{coro}

\begin{theo}[{\bf Comparison theorem}]\label{theo:cmapb}
Let the notations be as in {\rm (\ref{theo:funacb})} and after 
{\rm (\ref{prop:spqaqqqq})}. 
Then the following hold$:$ 
\par 
$(1)$ In ${\rm D}^+{\rm F}(f^{-1}({\cal W}(\kap_t)))$
there exists the following canonical isomorphism$:$
\begin{equation*}
(B_{{\rm zar},{\mab Q}}(X_{\bul \leq N,\os{\circ}{t}}
/{\cal W}(s_{\os{\circ}{t}}),E^{\bul \leq N}),P)
\os{\sim}{\lo}
({\cal W}B_{X_{\bul \leq N,\os{\circ}{t}},{\mab Q}}(E^{\bul \leq N}),P).
\tag{4.1.25.1}\label{eqn:brbadaw}
\end{equation*} 
The isomorphism $(\ref{eqn:brbadaw})$ 
is contravariantly functorial. 
\par 
$(2)$ The isomorphism 
{\rm (\ref{eqn:brbadaw})} forgetting the filtrations   
fits into the following commutative diagram$:$  
\begin{equation*} 
\begin{CD} 
B_{{\rm zar},{\mab Q}}
(X_{\bul \leq N,\os{\circ}{t}}/{\cal W}_{\star}(s_{\os{\circ}{t}}),E^{\bul \leq N})
@>{(\ref{eqn:brbadaw}),~\sim}>> \\
@A{\mu_{X_{\bul \leq N,\os{\circ}{t}}/{\cal W}(\os{\circ}{t})} \wedge}A{\simeq}A \\
\wt{R}u_{X_{\bul \leq N,t}/{\cal W}(t)*}
(\eps^*_{X_{\bul \leq N,\os{\circ}{t}}/{\cal W}(\os{\circ}{t})}(E^{\bul \leq N}))
\otimes_{\mab Z}{\mab Q} 
@>{\vpl_n(\ref{eqn:ywnttnny}),~\sim}>> 
\end{CD} 
\tag{4.1.25.2}\label{cd:axbwabl}
\end{equation*} 
\begin{equation*} 
\begin{CD} 
{\cal W}B_{X_{\bul \leq N,\os{\circ}{t}},{\mab Q}}(E^{\bul \leq N}) \\ 
@A{\mu_{X_{\bul \leq N,\os{\circ}{t}}/{\cal W}(\os{\circ}{t})}}A{\simeq}A \\
(\eps^*_{X_{\bul \leq N,\os{\circ}{t}/{\cal W}(\os{\circ}{t})}}
(E^{\bul \leq N}))_{{\cal W}_{\star}(X_{\bul \leq N,\os{\circ}{t}})}
\otimes_{{\cal W}_{\star}({\cal O}_{X_{\bul \leq N,\os{\circ}{t}})}}
{\cal W}\wt{\Om}{}^{\bul}_{X_{\bul \leq N,\os{\circ}{t}}}\otimes_{\mab Z}{\mab Q} . 
\end{CD} 
\end{equation*} 
\end{theo}

\chapter{Results for weight filtrations on log isocrystalline cohomology sheaves}

\parno 
In this chapter we prove various important properties of  
the weight filtration on the log isocrystalline cohomology sheaf of 
a truncated simplicial base change of proper SNCL schemes 
with admissible immersions. 
Especially we prove the following:  
\par 
(1) the log convergence of the weight filtration on 
the log isocrystalline cohomology sheaf of 
a truncated simplicial base change of proper SNCL schemes 
with admissible immersions,  
\par 
(2) the infinitesimal deformation invariance of the weight filtration on  
the log isocrystalline cohomology sheaf with weight filtration of 
the truncated simplicial base change in (1), 
\par 
(3) the degeneration at $E_2$ of the weight spectral sequence
of the log isocrystalline cohomology sheaf in (1),  
\par 
(4): the strict compatibility of the weight filtration in (1) 
with respect to the pull-back of a morphism of 
truncated simplicial base changes of proper SNCL schemes, 
\par 
(5) the monodromy-weight conjecture for a proper strict semistable family 
over a complete discrete valuation ring of equal positive characteristic 
and for a projective SNCL family  
over a complete discrete valuation ring of mixed characteristics, 
\par 
(6) the log hard Lefschetz conjecture for a projective strict semistable family 
over a complete discrete valuation ring of equal positive characteristic 
and for a projective SNCL family 
over a complete discrete valuation ring of mixed characteristics.

\section{Log convergent $F^{\infty}$-isospans}\label{sec:cfi}
In this section, as in \cite{ollc}, we work over a fine log formal scheme 
whose underlying formal scheme 
is a $p$-adic formal scheme in the sense of \cite[\S1]{of} 
(we review the definition of this formal scheme soon) 
and develop a filtered version of theory of log convergent $F^{\infty}$-isospans  in \cite[\S5]{ollc}, 
including theory of filtered log convergent $F$-isocrystals. 
In \cite{of} in the case of the trivial log structures, 
Ogus has used an equivalence between 
the category of convergent $F$-isocrystals 
and the category of $p$-adically convergent $F$-isocrystals 
for proving the convergence of the 
isocrystalline cohomology sheaf of a proper smooth scheme in characteristic $p>0$. 
In this section we prove a natural equivalence between 
the category of filtered log convergent $F^{\infty}$-isospans  
and the category of filtered log $p$-adically convergent $F^{\infty}$-isospans 
under certain mild assumptions.  
(Strictly speaking, we prove the two equivalences of 
the two categories of ``solid'' and ``not necessarily solid'' filtered log convergent $F^{\infty}$-isospans  
and two categories of ``solid'' and ``not necessarily solid'' 
filtered log $p$-adically convergent $F^{\infty}$-isospans: 
see (\ref{prop:fisoeq}), (\ref{prop:ge}) and (\ref{theo:zbfs}) below in detail.)  
In the next section, we shall use these equivalences  
for proving the log convergence of the weight filtration 
on the log isocrystalline cohomology sheaf of 
a truncated simplicial base change of proper SNCL schemes in characteristic $p>0$. 
\par 
See also \cite[2 (e)]{fao} for $F$-isocrystals in the case where 
fine log schemes are smooth schemes with SNCD's. 
\par  
Let $\kap$, ${\cal W}$ and ${\cal W}(s)$ be as in \S\ref{sec:flgdw}. 
Let ${\cal V}$ be a complete discrete valuation ring of mixed characteristics 
with perfect residue field $\kap_{}$ of characteristic $p>0$.  
Let $F \col {\cal W}(s)\lo {\cal W}(s)$ be the Frobenius endomorphism of ${\cal W}(s)$.  
Set $K_{0}:={\rm Frac}({\cal W}_{})$ 
and $K_{}:={\rm Frac}({\cal V}_{})$. 
\par 
Let $B$ be a fine log $p$-adic formal scheme 
whose underlying formal scheme is ${\rm Spf}({\cal V})$. 
Let $Z$ be a fine log $p$-adic formal scheme over $B$
whose underlying formal scheme is a $p$-adic formal ${\cal V}$-scheme 
in the sense of \cite[\S1]{of} (e.~g., ${\cal V}/p$-schemes), that is, 
a noetherian formal scheme over ${\rm Spf}({\cal V})$ 
with the $p$-adic topology which is topologically of finite type over ${\rm Spf}({\cal V})$.  
In this book, we call this $Z$ a {\it fine log $p$-adic formal $B$-scheme} for short. 
When $B=({\rm Spf}({\cal V}),{\cal V}^*)$, 
we call the $Z$ a {\it fine log $p$-adic formal ${\cal V}$-scheme}. 
\par 
We give the definition of a (solid) log ($p$-adic) enlargement as follows. 
This is a log version of Ogus' ($p$-adic) enlargement in \cite{of}. 
See also \cite[(1.1)]{oc}, \cite[\S3]{ollc}, \cite[Definition 2.1.1, 2.1.9]{s2}, 
\cite[I (2.2)]{s3} and \cite[\S2]{nhw}. 
\par 
Let $T$ be a fine log $p$-adic formal $B$-scheme. 
Assume that $\os{\circ}{T}$ is flat over $\os{\circ}{B}$. 
Set $T_1:=\ul{{\rm Spec}}^{\log}_T({\cal O}_T/p{\cal O}_T)$
and $T_0:=\ul{{\rm Spec}}^{\log}_T(({\cal O}_T/p{\cal O}_T)_{\rm red})$. 
We say that $T=(T,z_i)$ $(i=1$ (resp.~$0$))  
(strictly speaking, $T/B:=(T/B, z_i)$)
is a log $p$-adic enlargement (resp.~log enlargement) of $Z/B$ 
if $z_1\col T_1\lo Z$ (resp.~$z_0\col T_0\lo Z$) is a morphism 
of fine log formal schemes over $B$ 
fitting into the following commutative diagram 
\begin{equation*} 
\begin{CD} 
T_i @>{\subset}>> T \\
@V{z_i}VV @VVV \\
Z @>>> B. 
\end{CD} 
\end{equation*}  
For a log $p$-adic enlargement $T=(T,z_1)$,  
endow $T$ with the canonical PD-structure on $p{\cal O}_T$. 
If the morphism $z_1\col T_1\lo Z$ (resp.~$z_0\col T_0\lo Z$)  is solid, 
we call $T=(T,z_i)$ $(i=1$ (resp.~$0$))  
a {\it solid} log $p$-adic enlargement (resp.~solid log enlargement) of $Z/B$.  
Here note that the solidness of $z_i$ has not been assumed 
in \cite[Definition 2.1.1, 2.1.9]{s2} and \cite[I (2.1), (2.2)]{s3}, 
while this has been assumed in \cite[\S3]{ollc}. 
Note also that $Z$ is assumed to be saturated in \cite[\S3]{ollc}; 
we do not assume that $Z$ is saturated. 
Because the log formal scheme $Z$ will be a family of log points over ${\cal V}$ 
in our main applications in the next section, there is no obstacle in this book even if 
one assumes that $Z$ is saturated. 
We define a morphism of (solid) log ($p$-adic) enlargements of $Z/B$ 
in a standard way. 
Let ${\rm Enl}_p(Z/B)$ and ${\rm Enl}(Z/B)$ 
be the category of log $p$-adic enlargements of $Z/B$ 
and the category of log enlargements of $Z/B$, respectively; 
let ${\rm Enl}_p^{\rm sld}(Z/B)$ and ${\rm Enl}^{\rm sld}(Z/B)$ 
be the category of solid log $p$-adic enlargements of $Z/B$ 
and the category of solid log enlargements of $Z/B$, respectively 
(``sld'' is the abbreviation of ``solid'').  
In this section, let $\star$ be $p$ or nothing 
and let $\sq$ be sld or nothing. 
Then ${\rm Enl}_{\star}^{\rm sld}(Z/B)$ is a full subcategory of 
${\rm Enl}_{\star}(Z/B)$. 
\par 
Let $(T,z_1)$ be an object of ${\rm Enl}_p^{\sq}(Z/B)$.  
Then $(T,z_1)$ naturally defines an object of  
${\rm Enl}^{\sq}(Z/B)$ and consequently we 
have a functor 
\begin{align*} 
{\rm Enl}_p^{\sq}(Z/B)\lo  {\rm Enl}^{\sq}(Z/B).
\tag{5.0.0.1}\label{ali:enlzb}
\end{align*} 
When $B=({\rm Spf}({\cal V}),{\cal V}^*)$, we denote 
${\rm Enl}^{\sq}_{\star}(Z/B)$ by 
${\rm Enl}^{\sq}_{\star}(Z/{\cal V})$.  
We also have the following two natural functors 
\begin{align*} 
{\rm Enl}^{\sq}_{\star}(Z/B)\lo 
{\rm Enl}^{\sq}_{\star}(Z/\os{\circ}{B})={\rm Enl}^{\sq}_{\star}(Z/{\cal V}) 
\lo {\rm Enl}_{\star}(\os{\circ}{Z}/\os{\circ}{B}).
\tag{5.0.0.2}\label{ali:wedl}
\end{align*}

\begin{defi}\label{defi:nsd} 
(1)  Let $T=(T,z_i)$ ($i=1$ (resp.~$i=0$)) be 
an object of ${\rm Enl}^{\sq}_p(Z/B)$ (resp.~${\rm Enl}^{\sq}(Z/B)$). 
Let $\{(T_j,z_{ij})\}_j$ be a family of objects of ${\rm Enl}^{\sq}_{\star}(Z/B)$ 
with morphism $(T_j,z_{ij})\lo (T,z_{i})$. 
We say that $\{(T_j,z_{ij})\}_j$ is a {\it covering family} if $\{T_j\}_j$ is 
a Zariski open covering of $T$. 
By using the covering family, endow  
${\rm Enl}^{\sq}_{\star}(Z/B)$ with a natural topology in a standard way. 
The isostructure sheaf ${\cal K}^{\sq}_{Z/B}$ on ${\rm Enl}^{\sq}_{\star}(Z/B)$ is defined by 
the following formula: 
$\Gam((T,z_i),{\cal K}^{\sq}_{Z/B}):={\cal K}(T):=
\Gam(T,{\cal O}_T)\otimes_{\mab Z}{\mab Q}$. 
\par 
(2) As usual, a sheaf on ${\rm Enl}^{\sq}_{\star}(Z/B)$ is equivalent to the following data: 
a Zariski sheaf $E_T$ for $T\in {\rm Enl}^{\sq}_{\star}(Z/B)$ and a morphism $\rho_u\col 
u^{-1}(E_{T'})\lo E_T$ of Zariski sheaves on $E_T$ for a morphism 
$u\col T\lo T'$ in ${\rm Enl}^{\sq}_{\star}(Z/B)$, which is compatible 
with respect to a composite morphism $T\lo T'\lo T''$ in ${\rm Enl}^{\sq}_{\star}(Z/B)$ 
and which is an isomorphism in the case where $u$ is an open immersion.  
\par 
(3) Let $Z'/B'$ be a similar fine log $p$-adic formal scheme to $Z/B$ 
and let $g \col Z'\lo Z$ be a morphism over $B'\lo B$ such that 
the morphism ${\cal V}\lo \Gam(B',{\cal O}_{B'})=:{\cal V}'$ is a finite extension of 
complete discrete valuation rings. 
Let $E$ be a sheaf on ${\rm Enl}^{\sq}_{\star}(Z'/B')$. 
For an object $(T,z_i)$ of ${\rm Enl}^{\sq}_{\star}(Z/B)$,  
we can define a sheaf $g^{-1}((T,z_i))$ on  ${\rm Enl}^{\sq}_{\star}(Z'/B')$ 
in a standard way (e.~g., ~\cite[5.6]{bob}): 
\begin{align*} 
\Gam((T',z'_i),g^{-1}((T,z_i))):={\rm Hom}((T',z'_i),(T,z_i)), 
\tag{5.1.1.1}\label{ali:ephoffp} 
\end{align*}
where ${\rm Hom}((T',z'_i),(T,z_i))$ is, by definition, 
the set of morphisms $T'\lo T$ over $B'\lo B$ such that the induced morphism  
$T'_i\lo T_i$ fitting into the following commutative diagram 
\begin{equation*} 
\begin{CD} 
T'_i @>>> T_i\\
@V{z'_i}VV @VV{z_i}V \\
Z'@>{g}>> Z. 
\end{CD} 
\tag{5.1.1.2}\label{cd:ztiz} 
\end{equation*} 
Then the {\it push-forward} $g_*(E)$ of $E$ by $g$ is defined as follows: 
\begin{align*} 
g_*(E)((T,z_i)):={\rm Hom}(g^{-1}((T,z_i)),E) \quad (i=1,0),  
\tag{5.1.1.3}\label{ali:eplffp} 
\end{align*} 
where ${\rm Hom}$ means a morphism of sheaves on 
${\rm Enl}^{\sq}_p(Z/B)$ (resp.~${\rm Enl}^{\sq}(Z/B))$.   
\end{defi} 

\par 
\begin{defi}\label{defi:eenl} 
Let $E$ be a sheaf on ${\rm Enl}^{\sq}_{\star}(Z/B)$. 
We can define a sheaf $g^{-1}(E)$ on 
${\rm Enl}^{\sq}_{\star}(Z'/B')$ in a usual way as follows. 
For an object $(T',z'_i)\in {\rm Enl}^{\sq}_{\star}(Z'/B')$, 
consider the following presheaf 
\begin{align*} 
(g^{\bul}(E))((T',z'_i)):=\vil E((T,z_i))
\tag{5.1.2.1}\label{ali:epfp} 
\end{align*} 
on ${\rm Enl}^{\sq}_{\star}(Z'/B')$, where $(T,z_i)$'s 
are objects of ${\rm Enl}^{\sq}_{\star}(Z/B)$ with morphism 
$T'\lo T$ over $B'\lo B$ 
fitting into the commutative diagram (\ref{cd:ztiz}).  
Let $g^{-1}(E)$ be the sheafification of $g^{\bul}(E)$. 
\par 
We can easily prove the adjoint property for $g_*$ and $g^{-1}$. 
\par 
Let $E$ be a sheaf of ${\cal K}^{\sq}_{Z/B}$-modules. 
The {\it pull-back} $g^*(E)$ of $E$ by $g$ is defined as follows: 
\begin{align*} 
g^*(E):={\cal K}^{\sq}_{Z'/B'}\otimes_{g^{-1}({\cal K}^{\sq}_{Z/B})}g^{-1}(E).
\tag{5.1.2.2}\label{ali:epgsfp} 
\end{align*}   
\end{defi}

\par 
As in the nonfiltered case in \cite[Definition (2.7)]{of}, 
we define a filtered (solid) log $p$-adically convergent isocrystal on 
${\rm Enl}^{\sq}_p(Z/B)$ and 
a filtered (solid) log convergent isocrystal on ${\rm Enl}^{\sq}(Z/B)$ as follows:

\begin{defi}\label{defi:lci} 
Let $T=(T,z_i)$ ($i=1$ (resp.~$i=0$)) be 
an object of ${\rm Enl}^{\sq}_p(Z/B)$ (resp.~${\rm Enl}^{\sq}(Z/B)$). 
We call a family $(E,P):=\{(E_T,P_T)\}_T$ 
a {\it filtered $($solid$)$ log $p$-adically convergent isocrystal} on ${\rm Enl}^{\sq}_p(Z/B)$ 
(resp.~{\it filtered $($solid$)$ log convergent isocrystal} on 
${\rm Enl}^{\sq}(Z/B)$) if $(E,P)_T:=(E_T,P_T)$ is a coherent filtered 
${\cal K}_T$-module, that is, 
$E_T$ and $(P_T)_k(E_T)$ for $\forall k\in {\mab Z}$ 
are coherent ${\cal K}_T$-modules 
and for a morphism $g\col T'\lo T$ 
in ${\rm Enl}^{\sq}_p(Z/B)$ (resp.~${\rm Enl}^{\sq}(Z/B)$), 
there exists an isomorphism 
$$\rho_g\col g^*((E,P)_T):=(g^*(E_T), \{g^*(P_{T,k}(E_T))\}_{k\in {\mab Z}})
\os{\sim}{\lo} (E,P)_{T'}$$ 
of filtered ${\cal K}_{T'}$-modules  
which satisfies the usual cocycle condition and $\rho_{{\rm id}_T}={\rm id}_{(E,P)_T}$. 
\par 
A {\it morphism} $\al \col (E,P)\lo (F,Q)$ 
{\it of filtered $($solid$)$ log $p$-adically convergent isocrystals}  
on ${\rm Enl}^{\sq}_p(Z/B)$ (resp.~{\it filtered $($solid$)$ log convergent isocrystals} on 
${\rm Enl}^{\sq}(Z/B)$) 
is a collection $\al:=\{\al_T\col (E,P)_T\lo (F,Q)_T\}$, 
where $\al_T$'s for various $T$'s 
are morphisms of filtered ${\cal K}_T$-modules which are compatible with 
$\rho_g$'s for all $g$'s above.      
\par 
Let ${\rm IsocF}^{\sq}_p(Z/B)$ (resp.~${\rm IsocF}^{\sq}(Z/B)$) denote 
the category of filtered (solid) log $p$-adically convergent isocrystals on 
${\rm Enl}^{\sq}_{\star}(Z/B)$ 
(resp.~the category of filtered (solid) log convergent isocrystals on 
${\rm Enl}^{\sq}(Z/B)$). 
The category ${\rm IsocF}^{\sq}_{\star}(Z/B)$ is an additive category. 
If the filtration is trivial, then we denote 
${\rm IsocF}^{\sq}_{\star}(Z/B)$ by ${\rm Isoc}^{\sq}_{\star}(Z/B)$.  
\end{defi}   


\parno
By (\ref{ali:wedl}) 
we have the following natural functors 
\begin{align*} 
{\rm IsocF}_{\star}(\os{\circ}{Z}/\os{\circ}{B})
\lo {\rm IsocF}_{\star}^{\sq}(Z/\os{\circ}{B})
\lo {\rm IsocF}_{\star}^{\sq}(Z/B). 
\tag{5.1.3.1}\label{ali:czbff}
\end{align*} 
\par 

The calculation of the pull-back in (\ref{defi:eenl}) for the non-solid case is easy 
(cf.~\cite[(2.7.1)]{of}): 

\begin{prop}\label{prop:pce} 
Let the notations be as in {\rm (\ref{defi:nsd}) (3)} 
and {\rm (\ref{defi:eenl})}. 
Then the following hold$:$
\par 
$(1)$ 
Let $(T',z'_i)$ be an object of ${\rm Enl}_{\star}(Z'/B')$. 
Then 
\begin{align*} 
\{g^*((E,P))\}_{(T',z'_i)}=\{g^{-1}((E,P))\}_{(T',z'_i)}=(E,P)_{(T',g\circ z'_i)}. 
\tag{5.1.4.1}\label{ali:epnsefp} 
\end{align*}  
\par 
$(2)$ 
Let 
\begin{align*} 
\iota^* \col {\rm IsocF}(Z/B) \lo {\rm IsocF}_p(Z/B)
\tag{5.1.4.2}\label{ali:rstzb} 
\end{align*} 
and 
\begin{align*} 
\iota'{}^* \col {\rm IsocF}(Z'/B') \lo {\rm IsocF}_p(Z'/B')
\tag{5.1.4.3}\label{ali:rstzbzb} 
\end{align*} 
be the restriction functors. 
Then the following diagram 
\begin{equation*} 
\begin{CD} 
{\rm IsocF}(Z'/B')@>{\iota'{}^*}>> {\rm IsocF}_p(Z'/B') \\
@A{g^*}AA @AA{g^*}A \\
{\rm IsocF}(Z/B)@>{\iota^*}>> {\rm IsocF}_p(Z/B) 
\end{CD}
\tag{5.1.4.4}\label{ali:rstzfb} 
\end{equation*}
is commutative. 
\end{prop}
\begin{proof} 
(1): It is obvious that, if 
the following diagram 
\begin{equation*} 
\begin{CD} 
T'_i @>>> T_i\\
@V{z'_i}VV @VV{z_i}V \\
Z'@>{g}>> Z. 
\end{CD} 
\end{equation*} 
is commutative for an object $(T,z_i)\in {\rm Enl}(Z/B)$, 
then the following diagram 
\begin{equation*} 
\begin{CD} 
T'_i @>>> T_i\\
@V{g\circ z'_i}VV @VV{z_i}V \\
Z@= Z. 
\end{CD} 
\end{equation*}
is commutative. 
Hence $\{g^{-1}((E,P))\}_{(T',z_i)}=(E,P)_{(T',g\circ z'_i)}$. 
Because $(T',g\circ z'_i)$ is an object of ${\rm Enl}_{\star}(Z/B)$ for a not 
necessarily solid morphism $g\col Z\lo Z'$,  we obtain 
the equality (\ref{ali:epnsefp}). 
\par
(2): By using (1), we immediately obtain (2). 
\end{proof} 

\begin{rema}\label{rema:exep} 
(1) Let $E$ be an object of ${\rm IsocF}_{\star}^{\rm sld}(Z'/B')$. 
If the morphism $g\col Z'\lo Z$ is solid, then 
\begin{align*} 
\{g^*((E,P))\}_{(T',z'_i)}=\{g^{-1}((E,P))\}_{(T',z'_i)}=(E,P)_{(T',g\circ z'_i)} 
\tag{5.1.5.1}\label{ali:epgefp} 
\end{align*}  
because $(T',g\circ z'_i)$ is an object of ${\rm Enl}^{\rm sld}_{\star}(Z/B)$. 
For a general $g$, see (\ref{coro:npc}) (1) below for the calculation of $g^*$ 
under a mild assumption on $T'$. 
\par 
(2) Let the notations be as in (\ref{defi:eenl}). 
When $E$ is equal to a representable object 
$(T',z'_i)$ on ${\rm Enl}_{\star}(Z'/B')$, 
the following equality does {\it no}t necessarily hold 
on ${\rm Enl}(Z/B)$: 
\begin{align*} 
g_*(E)=(T',g\circ z'_i). 
\tag{5.1.5.1}\label{ali:epgzfp} 
\end{align*}  
Indeed, if this holds, then, 
for an object $(T'',z''_i)$ of ${\rm Enl}_{\star}(Z'/B)$, 
the following two commutative diagrams 
\begin{equation*} 
\begin{CD} 
T@>>> T'\\
@A{\bigcup}AA @AA{\bigcup}A \\
T_i @>>> T'_i\\
@V{z_i}VV @VV{g\circ z'_i}V \\
Z@= Z 
\end{CD} 
\end{equation*}
and 
\begin{equation*} 
\begin{CD} 
T''@>>> T\\
@A{\bigcup}AA @AA{\bigcup}A \\
T''_i @>>> T_i\\
@V{z''_i}VV @VV{z_i}V \\
Z'@>{g}>> Z
\end{CD} 
\end{equation*}
are given, then we have to obtain the following commutative diagram: 
\begin{equation*} 
\begin{CD} 
T''@>>> T'\\
@A{\bigcup}AA @AA{\bigcup}A \\
T''_i @>>> T'_i\\
@V{z''_i}VV @VV{z'_i}V \\
Z'@= Z'. 
\end{CD} 
\end{equation*}
However the two given diagrams imply only the following commutative diagram: 
\begin{equation*} 
\begin{CD} 
T''@>>> T@>>>T'\\
@A{\bigcup}AA  @A{\bigcup}AA @AA{\bigcup}A \\
T''_i @>>> T_i @>>> T'_i\\
@V{z''_i}VV @V{z_i}VV@VV{g\circ z'_i}V \\
Z'@>{g}>> Z@=Z. 
\end{CD} 
\end{equation*}
Hence our sheaf $g_*((T',z'_i))$ is different from Ogus' notation 
``$g_*((T',z'_i)):=(T',g\circ z'_i)$'' in \cite[p.~784]{of} (in the trivial logarithmic case). 
If we use his notation, then we have the following strange formula by 
(\ref{ali:epgefp}): 
\begin{align*} 
{\rm Hom}((T',z'_i),g^*(E))={\rm Hom}(g_*((T',z'_i)),E), 
\tag{5.1.6.2}\label{ali:epegp} 
\end{align*}  
where ${\rm Hom}$ on the left hand side (resp.~the right hand side) above 
is the morphisms in the category of sheaves on ${\rm Enl}^{\sq}_{\star}(Z'/B')$ 
(resp.~${\rm Enl}^{\sq}_{\star}(Z/B)$).  
For this reason, we do not use his notation. 
\end{rema}

\begin{exem} 
(1) Consider a case where $B=({\rm Spf}({\cal V}),{\cal V}^*)$ 
and $S$ is a $p$-adic formal ${\cal V}$-scheme with trivial log structure 
and $X$ is a proper SNC scheme over $S$ with trivial log structure.  
Let $g\col X\lo S$ be the structural morphism. 
Let $k$ and $q$ be nonnegative integers. 
For a $p$-adic enlargement $T$ of $S/B$, 
set $X^{(k)}_{T_1}:=X^{(k)}\times_ST_1$ $(k\in {\mab Z}_{\geq 0})$ and 
let ${\rm Coh}({\cal K}_T)$ be the category of coherent ${\cal K}_T$-modules. 
Then 
\begin{align*} 
{\rm Enl}_p(S/{\cal V})\owns T\lom 
R^qg_{X^{(k)}_{T_1}/T*}({\cal O}_{X^{(k)}_{T_1}/T})_{\mab Q}
\in {\rm Coh}({\cal K}_T)
\tag{5.1.6.1}
\end{align*} 
is an object of ${\rm Isoc}_p(S/{\cal V})$ by \cite[(3.1)]{of}, 
which we denote by $R^qg_*({\cal O}_{X^{(k)}/K})$ as in [loc.~cit.].  
In fact, $R^qg_*({\cal O}_{X^{(k)}/K})$ prolongs naturally 
to an object of ${\rm Isoc}(S/{\cal V})$ by [loc.~cit., (3.7)]. 
\par 
(2) Consider a case where $B=({\rm Spf}({\cal V}),{\cal V}^*)$ 
and $S$ is a $p$-adic formal family of log points such that $\os{\circ}{S}$ is 
a $p$-adic formal ${\cal V}$-scheme. 
Let $g\col X\lo S$ be a proper SNCL scheme. 
Then 
\begin{align*} 
{\rm Enl}_p(\os{\circ}{S}/{\cal V})\owns \os{\circ}{T}\lom 
R^qg_{\os{\circ}{X}{}^{(k)}_{T_1}/\os{\circ}{T}*}
({\cal O}_{\os{\circ}{X}{}^{(k)}_{T_1}/\os{\circ}{T}})_{\mab Q}
\in {\rm Coh}({\cal K}_T)
\tag{5.1.6.2}
\end{align*} 
is an object of ${\rm Isoc}_p(\os{\circ}{S}/{\cal V})$ by (1). 
Hence this is an object of ${\rm Isoc}^{\sq}_p(S/{\cal V})$,  
which we denote by $R^qg_*({\cal O}_{\os{\circ}{X}{}^{(k)}/K})$.  
In fact, $R^qg_*({\cal O}_{\os{\circ}{X}{}^{(k)}/K})$ prolongs naturally 
to an object of ${\rm Isoc}^{\sq}(\os{\circ}{S}/{\cal V})$ 
and ${\rm Isoc}^{\sq}(S/{\cal V})$ by (1). 
\end{exem} 

\par 
Let $Z/B$ be as in the beginning of this section. 
Let $(T,z_1)$ (resp.~$(T,z_0)$) be an object of ${\rm Enl}_p(Z/B)$ 
(resp.~${\rm Enl}(Z/B)$). 
From now on, we mean this sentence by the following 
sentence ``Let $(T,z_i)$ $(i=0,1)$ be an object of ${\rm Enl}_{\star}(Z/B)$.''

\begin{defi}\label{dfm} 
Let the notations be as in the beginning of this section 
(we do not assume that $z_i$ $(i=0,1)$ is solid).   
Set $Z_{\os{\circ}{T}_i}:=Z\times_{\os{\circ}{Z}}\os{\circ}{T}_i$.
Let $Z(T)$ be a (not necessarily fine) log formal scheme 
whose underlying formal scheme is $\os{\circ}{T}$ and whose log structure 
is the inverse image of the image of the following morphism 
\begin{align*} 
M_{Z_{\os{\circ}{T}_i}}/{\cal O}_{Z_{\os{\circ}{T}_i}}^*\lo 
M_{T_i}/{\cal O}_{T_i}^*
\tag{5.1.7.1}\label{ali:ztio}
\end{align*} 
by the composite morphism 
$M_T\lo M_T/{\cal O}_T^*\os{\sim}{\lo} M_{T_i}/{\cal O}_{T_i}^*$. 
The log formal scheme $Z(T)$ is independent of $i$. 
The morphism $z_1\col T_1\lo Z$ (resp.~$z_0\col T_0\lo Z$) 
induces a morphism 
$Z(z_1)\col Z(T)_1\lo Z$ (resp.~$Z(z_0)\col Z(T)_0\lo Z$).  
When we have to clarify the morphism $z_i$,  
we denote $Z(T)$ and $Z_{\os{\circ}{T}_i}$ by 
$Z_{z_i}(T)$ and $Z_{(\os{\circ}{T}_i,z_i)}$, respectively.  
(By the definition of $Z(T)_i$, the log structure of $Z(T)_i$ is the inverse image of 
$Z$ by $z_i$.)
\end{defi} 
By the definition of $Z(T)$, we have a natural morphism $T\lo Z(T)$ over $B$.

\begin{exem}\label{exem:oot}  
Assume that $Z$ is a formal family of log points. 
Then the morphism (\ref{ali:ztio}) is injective by (\ref{prop:oot}). 
\end{exem} 

\par 
We have essentially proved the following (2) in \cite[Lemma 2.3.1]{nh2}: 

\begin{prop}\label{prop:nlctt} 
Assume that the morphism {\rm (\ref{ali:ztio})} is injective. 
Then the following hold$:$ 
\par 
$(1)$ $Z(T)_i=Z_{\os{\circ}{T}_i}$ $(i=0,1)$. 
\par 
$(2)$ 
Assume that, for any point $z\in \os{\circ}{Z}$, there exists a local chart 
$P\lo M_Z$ around $z$ such that $P^{\rm gp}$ is $p$-torsion-free. 
Then $M_{Z(T)}$ is fine. 
\par 
$(3)$ Let $(T',z'_i)\lo (T,z_i)$ be a morphism in 
${\rm Enl}^{\sq}_{\star}(Z/B)$. 
Assume that the morphism {\rm (\ref{ali:ztio})} for $(T',z'_i)$ is injective. 
Then $Z(T')=Z(T)\times_{\os{\circ}{T}}\os{\circ}{T}{}'$. 
\end{prop} 
\begin{proof} 
(1): Since the log structure of $Z(T)_i$ is equal to 
${\rm Im}(M_{Z_{\os{\circ}{T}_i}}\lo M_{T_i})$, (1) is obvious. 
\par 
(2): Let $t\in \os{\circ}{Z}_{\os{\circ}{T}_i}=\os{\circ}{T}_i$ 
be a point above $z$. 
The sheaf $M_{Z_{\os{\circ}{T}_i}}$ of monoids on $Z_{\os{\circ}{T}_i}$ has a local chart $P\lo M_{Z_{\os{\circ}{T}_i}}$ around $t$. 
Set $Z(T)_n:=Z(T)\mod p^n$ $(n\in {\mab Z}_{>0})$. 
Set ${\cal I}_n:={\rm Ker}({\cal O}_{T_n}\lo {\cal O}_{Z_{\os{\circ}{T}_i}})$. 
Consider the following exact sequence 
\begin{align*} 
1\lo 1+{\cal I}_n\lo M_{Z(T)_n}\lo M_{Z_{\os{\circ}{T}_i}}\lo 1. 
\end{align*}
Since $P^{\rm gp}$ is $p$-torsion-free and $1+{\cal I}_n$ is $p$-torsion, 
${\cal E}{\it xt}^1(P^{\rm gp},1+{\cal I}_n)=0$. 
Hence the morphism $P\lo M_{Z_{\os{\circ}{T}_i}}$ lifts to a morphism 
$P\lo M_{Z(T)_n}$. As in \cite[Lemma 2.3.1]{nh2}, 
we can easily check that this is a local chart of $M_{Z(T)_n}$. 
By the similar reasoning, we see that 
the natural morphism 
${\rm Hom}(P^{\rm gp},M_{Z(T)_{n+1}})\lo {\rm Hom}(P^{\rm gp},M_{Z(T)_n})$ 
is surjective. 
Hence we have a local chart $P\lo M_{Z(T)}$. 
This shows (2). 
\par 
(3): The proof is straightforward. 
\end{proof}

\begin{coro}\label{coro:tlfs}
If $M_Z$ is fs, then $M_{Z(T)}$ is fs. 
\end{coro} 
\begin{proof} 
(cf.~the proof of ~\cite[(7.1)]{fkato})
We have only to prove that $M_{Z(T)}/{\cal O}_{Z(T)}^*$ is finitely generated. 
Since $M_Z$ is saturated, there exists a local chart $P\lo M_Z$  
such that $P^{\rm gp}$ is torsion-free. 
Indeed, if $x^n=1$ $(x\in M_Z^{\rm gp})$, then $x\in {\cal O}_Z^*\subset M_Z$. 
Hence the image of the torsion part 
$P^{\rm gp}_{\rm tor}$ of $P^{\rm gp}$ in $M_Z^{\rm gp}$ 
is contained in ${\cal O}_Z^*$. 
Let $P^{\rm gp}_{\rm free}$ be the free part of $P^{\rm gp}$. 
Let $a_1,\ldots, a_k$ be a system of generators of $P$. 
Then $a_i=b_ic_i$, 
where $a_i=(b_i,c_i)\in P^{\rm gp}_{\rm free}\oplus 
P^{\rm gp}_{\rm tor}=P^{\rm gp}$.  
Set $P':=\langle b_1,\ldots,b_k\rangle$. 
Since $a_i^n=b_i^n$ for a large integer $n$ and since $M_Z$ is saturated, 
the image of $b_i$ in $M_Z^{\rm gp}$ belongs to $M_Z$. 
Since the image of 
$a_i$ is equal to the image of $b_i$ modulo 
a torsion element of $M_Z^{\rm gp}$, 
the morphism $P'\lo M_Z$ is indeed a chart of $M_Z$. 
Hence $M_{Z(T)}$ is fine by (\ref{prop:nlctt}) (2).  
In fact, $M_{Z(T)}$ is fs. 
\end{proof}

\begin{prop}\label{prop:bsch}    
Assume that there exists a local chart in {\rm (\ref{prop:nlctt}) (2)}.   
Then the following hold$:$
\par 
$(1)$ Assume that the morphism {\rm (\ref{ali:ztio})} is injective. 
Then the log formal scheme $Z(T)=(Z(T),Z(z_i))$ $(i=1,0)$  
is a log $(p$-adic$)$ enlargement of $Z/B$.  
\par 
$(2)$ Assume that the morphism {\rm (\ref{ali:ztio})} is injective for any object 
$T\in {\rm Enl}_{\star}(Z/B)$. 
Then the log formal scheme $Z(T)$ defines the following functor 
\begin{align*} 
Z(?)\col {\rm Enl}_{\star}(Z/B)\owns (T,z_i)\lom (Z(T),Z(z_i))\in {\rm Enl}^{\rm sld}_{\star}(Z/B). 
\tag{5.1.11.1}\label{ali:tzi}
\end{align*} 
The restriction $Z(?)\vert_{{\rm Enl}_{\star}^{\rm sld}(Z/B)}$ of the functor 
$Z(?)$ is equal to ${\rm id}_{{\rm Enl}_{\star}^{\rm sld}(Z/B)}$. 
\par 
$(3)$ Let $\iota \col {\rm Enl}_{\star}^{\rm sld}(Z/B)\lo {\rm Enl}_{\star}(Z/B)$ 
be the inclusion functor. Then the functor $Z(?)$ is the left adjoint functor of $\iota$.   
\par 
$(4)$ Assume that the morphism {\rm (\ref{ali:ztio})} is injective for any object 
$T\in {\rm Enl}_{\star}(Z/B)$. 
Then the following restriction functor is an equivalence of categories: 
\begin{align*} 
R\col {\rm IsocF}_{\star}(Z/B)\os{\sim}{\lo} {\rm IsocF}^{\rm sld}_{\star}(Z/B). 
\tag{5.1.11.2}\label{ali:tzzi}
\end{align*} 
\end{prop}
\begin{proof} 
(1): By (\ref{prop:nlctt}) (1) and the definition of $Z(T)$, the morphism 
$Z_{\os{\circ}{T}_i}=Z(T)_i \lo Z(T)$ is an exact closed immersion. 
Let $p\col T\lo B$ be the structural morphism. 
(1) is clear because 
$M_{Z(T)}$ contains ${\rm Im}(p^*(M_B)\lo M_T)$. 
\par 
(2): (2) is clear. 
\par 
(3): (3) follows from (2). 
\par 
(4):  It is clear that $R$ is functorial. 
Assume that we are given an object $(E,P)$ of ${\rm IsocF}^{\rm sld}_{\star}(Z/B)$. 
For an object $T$ of ${\rm Enl}_{\star}(Z/B)$, 
set $(E'_T,P'_T):=(E_{Z(T)},P_{Z(T)})$. 
Then we see that $(E',P'):=\{(E'_T,P'_T)\}_{T\in {\rm Enl}_{\star}(Z/B)}$ 
is an object of ${\rm IsocF}_{\star}(Z/B)$ by (\ref{prop:nlctt}) (3). 
Obviously $R((E',P'))=(E,P)$. 
Hence $R$ is essentially surjective. 
\par 
Let $(E,P)$ be an object of ${\rm IsocF}_{\star}(Z/B)$. 
Because we have a natural morphism 
$(T,z_i)\lo (Z(T),Z(z_i))$ in ${\rm Enl}_{\star}(Z/B)$ and 
$\os{\circ}{Z(T)}=\os{\circ}{T}$, 
$(E,P)_{(T,z_i)}=(E,P)_{(Z(T),Z(z_i))}$. 
Using this formula, we easily see that $R$ is fully faithful. 
\par 
Because $R$ is essentially surjective and fully faithful, 
$R$ is an equivalence of categories. 
\end{proof}

\begin{prop}\label{prop:ips}
Let $Z/B$, $Z'/B'$ and $T'= (T',z'_i)$ and $T=(T,z_i)$ be as in {\rm (\ref{defi:nsd}) (3)}.  
Assume that the morphism $z_i\col T_i\lo Z$ is solid. 
Consider $Z(T')$ by the use of the composite morphism $g\circ z'_i\col T'_i\lo Z$.  
$($We can consider $Z(T')$ as a log $p$-adic formal $B$-scheme.$)$
Assume that the following morphism 
\begin{align*} 
M_{Z_{\os{\circ}{T}{}'_i}}/{\cal O}_{Z_{\os{\circ}{T}{}'_i}}^*\lo 
M_{T'_i}/{\cal O}_{T'_i}^*
\tag{5.1.12.1}\label{ali:ztoio}
\end{align*} 
is injective. 
Then the morphism $T'\lo T$ over $B'\lo B$ 
decomposes uniquely into the following composite 
morphism $T'\lo Z(T')\lo T$ over $B'\lo B;$ 
the commutative diagram {\rm (\ref{cd:ztiz})} 
decomposes uniquely into the following commutative diagram$:$ 
\begin{equation*} 
\begin{CD} 
T'_i @>>> Z(T')_i @>>> T_i\\
@V{z'_i}VV @V{Z(g\circ z'_i)}VV @VV{z_i}V \\
Z'@>{g}>>Z@= Z. 
\end{CD} 
\tag{5.1.12.2}\label{cd:ztviz} 
\end{equation*}
$($Note that $(Z(T'),Z(g\circ z'_i))$ is an object of ${\rm Enl}^{\rm sld}_{\star}(Z/B)$.$)$    
\end{prop}
\begin{proof} 
Since $Z(T)=T$, this immediately follows from (\ref{prop:bsch}) (2). 
\end{proof} 


\begin{coro}\label{coro:npc} 
Let the notations be as in {\rm (\ref{defi:nsd}) (3)} and 
let the assumptions be as in {\rm (\ref{prop:ips})}. 
Then the following hold$:$ 
\par 
$(1)$ 
Let $g^*\col {\rm IsocF}^{\rm sld}_{\star}(Z/B)\lo {\rm IsocF}^{\rm sld}_{\star}(Z'/B')$ 
be the pull-back morphism. 
Then 
\begin{align*} 
g^*(E)_{T'}=E_{Z(T')}.
\tag{5.1.13.1}\label{ali:rstzsfb} 
\end{align*} 
\par 
$(2)$ 
Let 
\begin{align*} 
\iota^* \col {\rm IsocF}^{\rm sld}(Z/B) \lo {\rm IsocF}_p^{\rm sld}(Z/B)
\tag{5.1.13.2}\label{ali:rsstzb} 
\end{align*} 
and 
\begin{align*} 
\iota'{}^* \col {\rm IsocF}^{\rm sld}(Z'/B') \lo {\rm IsocF}_p^{\rm sld}(Z'/B')
\tag{5.1.13.3}\label{ali:rstsbzb} 
\end{align*} 
be the restriction functors. 
Then the following diagram 
\begin{equation*} 
\begin{CD} 
{\rm IsocF}^{\rm sld}(Z'/B')@>{\iota'{}^*}>> {\rm IsocF}_p^{\rm sld}(Z'/B') \\
@A{g^*}AA @AA{g^*}A \\
{\rm IsocF}^{\rm sld}(Z/B)@>{\iota^*}>> {\rm IsocF}_p^{\rm sld}(Z/B) 
\end{CD}
\tag{5.1.13.4}\label{cd:rstzsfb} 
\end{equation*}
is commutative. 
\end{coro}
\begin{proof} 
(1): By (\ref{ali:epfp}) and (\ref{prop:ips}) and by noting 
that $(Z(T'),Z(g\circ z'_i))$ is an object of ${\rm Enl}^{\rm sld}_{\star}(Z/B)$, 
$g^{\bul}(E)_{T'}=E_{Z(T')}$. Hence  
$g^{-1}(E)_{T'}=E_{Z(T')}$ by (\ref{ali:epfp}). 
In particular, 
$g^{-1}({\cal K}_{Z/B})_{T'}={\cal K}_{Z(T')}={\cal K}_{T'}$ 
(since $\os{\circ}{Z(T')}=\os{\circ}{T}{}'$). 
Hence $g^*(E)_{T'}=E_{Z(T')}$. 
\par 
(2): Let $T'$ be an object of ${\rm Enl}_p^{\rm sld}(Z'/B')$ and 
let $(E,P)$ be an object of ${\rm IsocF}^{\rm sld}(Z/B)$. 
We omit to write the filtration $P$ of $(E,P)$. 
Then $(\iota^*g^*(E))_{T'}=(g^*(E))_{T'}=E_{Z(T')}$. 
On the other hand, 
$(g^*\iota^*(E))_{T'}=(\iota^*(E))_{Z(T')}=E_{Z(T')}$.
\par 
Another proof: (2) immediately follows from (\ref{ali:tzzi}). 
\end{proof} 

The following is a filtered $F^{\infty}$-isospan version 
of an $F^{\infty}$-span of crystals 
in \cite[Definition 15]{ollc}. See also [loc.~cit,. \S5] for a different formulation 
in a more restricted case. 

\begin{defi}\label{defi:finf}
Assume that ${\cal V}={\cal W}$ and that $B$ is the canonical lift of $B\mod p$. 
Then $B$ has a canonical lift $F_B\col B\lo B$ 
of the Frobenius endomorphism of $B~{\rm mod}~p$. 
Assume also that $Z$ is of characteristic $p>0$.  
Let $F_Z\col Z\lo Z$ be the Frobenius endomorphism of $Z$ over $F_B\col B\lo B$.   
Let $n$ be a positive integer.   
We call a family $\{((E_n,P_n),\Phi_n)\}_{n=0}^{\inf}$ of pairs 
a {\it filtered $($solid$)$ log $p$-adically convergent $F^{\inf}$-isospan} 
(resp.~a {\it filtered $($solid$)$ log convergent $F^{\inf}$-isospan}) if 
$(E_n,P_n)\in {\rm IsocF}^{\sq}_p(Z/B)$,  
(resp.~$(E_n,P_n)\in {\rm IsocF}^{\sq}(Z/B)$) 
and $\Phi_n$ is the following isomorphism 
\begin{align*} 
\Phi_n \col F^*_Z((E_{n+1},P_{n+1})) \os{\sim}{\lo} (E_n,P_n)
\quad (n\geq 0)
\tag{5.1.14.1}\label{ali:islep} 
\end{align*}  
in ${\rm IsocF}^{\sq}_p(Z/B)$ 
(resp.~$(E_n,P_n)\in {\rm IsocF}^{\sq}(Z/B)$). 
We define a {\it morphism of filtered $($solid$)$ log ($p$-adically) 
convergent $F^{\inf}$-isospans} in an obvious way. 
We denote the category of filtered (solid) 
log $p$-adically convergent $F^{\infty}$-isospans 
(resp.~the category of filtered (solid) log convergent $F^{\infty}$-isospans)
by $F^{\infty}{\textrm -}{\rm IsosF}^{\sq}_p(Z/B)$ and 
$F^{\infty}{\textrm -}{\rm IsosF}^{\sq}(Z/B)$, 
respectively. When the filtration is trivial, we denote 
$F^{\infty}{\textrm -}{\rm IsosF}^{\sq}_p(Z/B)$ and 
$F^{\infty}{\textrm -}{\rm IsosF}^{\sq}(Z/B)$ by 
$F^{\infty}{\textrm -}{\rm Isos}^{\sq}_p(Z/B)$ and 
$F^{\infty}{\textrm -}{\rm Isos}^{\sq}(Z/B)$, respectively. 
The category $F^{\infty}{\textrm -}{\rm IsosF}^{\sq}_{\star}(Z/B)$ is an additive category. 
\par 
As usual,  by considering the case $(E_n,P_n)= (E_{n-1},P_{n-1})$ $(\forall n\geq 1)$, 
we have the notions of a {\it filtered $($solid$)$ log $p$-adically convergent $F$-isocrystal} 
(resp.~a {\it filtered $($solid$)$ log convergent $F$-isocrystal}) 
and a {\it morphism of filtered $($solid$)$ log $p$-adically convergent $F$-isocrystals} 
(resp.~a {\it morphism of filtered $($solid$)$ log convergent $F$-isocrystals}). 
In this case, we denote 
$$F^{\infty}{\textrm -}{\rm IsosF}^{\sq}_{\star}(Z/B)~{\rm and}~ 
F^{\infty}{\textrm -}{\rm IsosF}^{\sq}_{\star}(Z/B)$$
by 
$$F{\textrm -}{\rm IsocF}^{\sq}_{\star}(Z/B)~{\rm and}~
F{\textrm -}{\rm IsocF}^{\sq}_{\star}(Z/B),$$
respectively. 
\end{defi}

\begin{prop}\label{prop:fisoeq} 
Assume that the morphism {\rm (\ref{ali:ztio})} is injective and 
that there exists a local chart in {\rm (\ref{prop:nlctt}) (2)}.   
Then the equivalence {\rm (\ref{ali:tzzi})} 
induces the following equivalence of categories$:$
\begin{align*} 
R=\{R_n\}_{n=0}^{\inf} \col F^{\infty}{\textrm -}{\rm IsosF}_{\star}(Z/B) 
\os{\sim}{\lo} F^{\infty}{\textrm -}{\rm IsosF}^{\rm sld}_{\star}(Z/B). 
\tag{5.1.15.1}\label{ali:fsriszb} 
\end{align*} 
\end{prop}
\begin{proof} 
Let $E$ be an object of ${\rm Isoc}_{\star}(Z/B)$.  
Let $R$ be the functor (\ref{ali:tzzi}). 
We would like to prove that $R(F_Z^*(E))=F_Z^*(R(E))$. 
Let $(T,z_i)$ be an object of ${\rm Enl}^{\rm sld}_{\star}(Z/B)$ 
with structural morphism $z_i\col T_i\lo Z$. 
By abuse of notation, 
set $T_i^{[p]}:=T_i\times_{\os{\circ}{T}_i,\os{\circ}{F}_{T_i}}\os{\circ}{T}_i$.
Let $z^{[1]}_i\col T_i^{[p]}\lo Z$ be the composite morphism 
$T_i^{[p]}\os{\rm proj}{\lo}T_i\os{z_i}{\lo} Z$. 
By using the composite morphism $F_Z\circ z_i$, 
we have a fine log scheme $Z_{F_Z\circ z_i}(T)$. 
Set $T^{[p]}:=(Z_{F_Z\circ z_i}(T),Z_{F_Z\circ z_i}(T)\lo Z)$. 
We claim that there exists a canonical isomorphism 
$Z_{F_Z\circ z_i}(T)_i\os{\sim}{\lo}T^{[p]}_i$. 
Indeed, the underlying scheme of $Z_{F_Z\circ z_i}(T)_i$ 
is equal to $\os{\circ}{T}_i=(T^{[p]}_i)^{\circ}$. 
We claim that the log structure 
$M_{Z_{F_Z\circ z_i}(T)_i}$ is the associated log structure to 
the composite morphism 
${\rm Im}(p\times \col M_{T_i}\lo M_{T_i})\os{\subset}{\lo} 
M_{T_i}\lo {\cal O}_{T_i}$. 
Indeed, by (\ref{prop:nlctt}) (1), 
$Z_{F_Z\circ z_i}(T)_i=Z_{(\os{\circ}{T}_i,F_Z\circ z_i)}$. 
By the following obvious commutative diagram 
\begin{equation*} 
\begin{CD} 
M_{T_i}/{\cal O}_{T_i}^*@<{p\times}<< M_{T_i}/{\cal O}_{T_i}^*\\
@| @| \\
M_{Z_{(\os{\circ}{T}_i,z_i)}}/{\cal O}_{T_i}^*
@<{p\times}<< M_{Z_{(\os{\circ}{T}_i,z_i)}}/{\cal O}_{T_i}^*, 
\end{CD} 
\end{equation*} 
we see that the claim holds. 
Next we claim that this log structure is isomorphic to $M_{T^{[p]}_i}$ 
(cf.~\cite[pp.~196--198]{ollc}). 
The abrelative Frobenius morphism 
$F_{T_i/\os{\circ}{T}_i}\col T_i \lo T^{[p]}_i$ 
induces the following commutative diagram 
\begin{equation*} 
\begin{CD} 
1@>>> {\cal O}_{T_i}^*@>>> M_{T_i}@>>> M_{T_i}/{\cal O}_{T_i}^*@>>> 1\\
@. @| @AAA @AA{p\times }A \\
1@>>> {\cal O}_{T^{[p]}_i}^*
@>>> M_{T^{[p]}_i}@>>> M_{T_i}/{\cal O}_{T_i}^*@>>> 1. 
\end{CD}
\end{equation*} 
Here we have used the following equalities 
$M_{T^{[p]}_i}/{\cal O}_{T^{[p]}_i}^*={\rm proj}^{-1}(M_{T_i}/{\cal O}_{T_i}^*)
=M_{T_i}/{\cal O}_{T_i}^*$. 
By the assumption of the local chart,  
we see that the right vertical morphism induces an 
isomorphism $M_{T^{[p]}_i}/{\cal O}_{T^{[p]}_i}^* \os{\sim}{\lo}
{\rm Im}(p\times \col M_{T_i}/{\cal O}_{T_i}^*\lo M_{T_i}/{\cal O}_{T_i}^*)$. 
Hence we see that the last claim also holds. 
In conclusion, we have a canonical isomorphism 
$Z_{F_Z\circ z_i}(T)_i\os{\sim}{\lo}T^{[p]}_i$. 
\par 
Let $(T',z'_i)$ be an object of ${\rm Enl}_{\star}^{\rm sld}(Z/B)$ 
such that the induced morphism   
$T_i\lo T'_i$ fits into the following commutative diagram 
\begin{equation*} 
\begin{CD} 
T_i @>>> T'_i\\
@V{z_i}VV @VV{z'_i}V \\
Z@>{F_Z}>> Z. 
\end{CD} 
\tag{5.1.15.2}\label{cd:ztfiz} 
\end{equation*} 
Since $\os{\circ}{T}_i= (T_i^{[p]})^{\circ}$ and $\os{\circ}{T}= (T^{[p]})^{\circ}$, 
it is clear that the morphism $(T,z_i)\lo (T',z'_i)$ factors through 
$(T^{[p]},z^{[1]}_i)$. 
Hence, by the definition of $F_Z^*$ in ${\rm Enl}^{\rm sld}_{\star}(Z/B)$,  
$(F_Z^*(R(E)))_{(T,z_i)}=E_{(T^{[p]},z^{[1]}_i)}$. 
On the other hand, we have the following equalities: 
\begin{align*} 
R(F_Z^*(E))_{(T,z_i)}=(F_Z^*(E))_{(T,z_i)}=E_{(T,F_Z\circ z_i)}. 
\end{align*} 
Since we have the following natural commutative diagram 
\begin{equation*} 
\begin{CD}
T@>>>T^{[p]} \\
@A{\bigcup}AA @AA{\bigcup}A \\
T_i @>{F_{T_i/\os{\circ}{T}_i}}>> T^{[p]}_i \\
@V{z_i}VV @VV{z^{[1]}_i}V \\
Z@>{F_Z}>>Z, 
\end{CD}
\end{equation*} 
we have a natural morphism $(T,F_Z\circ z_i)\lo (T^{[p]},z^{[1]}_i)$. 
Since $E$ is an isocrystal on $Z/B$ and ${\cal O}_T={\cal O}_{T^{[p]}}$, 
$E_{(T,F_Z\circ z_i)}=E_{(T^{[p]},z^{[1]}_i)}$. 
Hence 
\begin{align*} 
F_Z^*(R(E))=R(F_Z^*(E)).
\tag{5.1.15.3}\label{ali:fzre}
\end{align*}  
Let the notations be as in (\ref{defi:finf}). 
Using (\ref{ali:fzre}),   
we have the following isomorphism 
\begin{align*} 
F^*_ZR((E_{n+1},P_{n+1}))=RF^*_Z((E_{n+1},P_{n+1}))\os{\sim}{\lo} R((E_{n},P_n)).
\tag{5.1.15.4}\label{ali:fzzre}
\end{align*}
Hence the functor 
$R \col {\rm IsocF}_{\star}(Z/B) \lo {\rm IsocF}^{\rm sld}_{\star}(Z/B)$ 
induces the functor (\ref{ali:fsriszb}).  
Since the functor 
$R \col {\rm IsocF}_{\star}(Z/B) \lo {\rm IsocF}^{\rm sld}_{\star}(Z/B)$ 
is an equivalence of categories ((\ref{ali:tzzi})), 
the functor (\ref{ali:fsriszb}) is also an equivalence of categories.
\end{proof}

\par 
The following lemma plays an important role in the proof of 
(\ref{prop:ge}) and (\ref{theo:zbfs}) below. 

\begin{lemm}\label{lemm:fianf}  
Let ${\cal V}$, $B$ and $Z$ be as in {\rm (\ref{defi:finf})}. 
Let $z\col U(0)\lo Z$ be a solid morphism of log schemes. 
Let $\iota\col U(0)\os{\sus}{\lo}U$ be a nilpotent exact immersion of 
fine log schemes of characteristic $p>0$. 
Let $m$ be a positive integer. 
Set $U^{[p^m]}:=U\times_{\os{\circ}{U},\os{\circ}{F}{}^m_U}\os{\circ}{U}$ 
and $U(0)^{[p^m]}:=U(0)\times_{\os{\circ}{U}(0),\os{\circ}{F}{}^m_{U(0)}}\os{\circ}{U}(0)$.  
Let $\iota^{[p^m]}\col U(0)^{[p^m]}\lo U^{[p^m]}$ be 
the natural exact immersion induced by $\iota$. 
Consider the following commutative diagram
\begin{equation*} 
\begin{CD} 
U(0)^{[p^{m}]}@>{\iota^{[p^m]}}>> U^{[p^{m}]}\\
@V{\rm proj.}VV @VV{\rm proj.}V \\
U(0)@>{\iota}>>U.
\end{CD} 
\tag{5.1.16.1}\label{ali:fumd}
\end{equation*} 
Then the following hold$:$
\par 
$(1)$ Let $U(0)\lo U(0)^{[p^m]}$ be the induced morphism 
by the $m$-th iterate of the Frobenius endomorphism 
$F_{U(0)}\col U(0)\lo U(0)$ of $U(0)$. 
Let $F_{m,0,U}\col U(0)\lo U(0)^{[p^m]}\os{\sus}{\lo} U^{[p^m]}$ 
be the natural composite morphism. 
Then there exist a positive integer $m$ and a morphism 
$\rho^{[m]}\col U^{[p^m]}\lo U(0)$ such that 
$\rho^{[m]} \circ F_{m,0,U}=F^m_{U(0)}$ 
and 
$F_{m,0,U} \circ \rho^{[m]}=F^m_{U(0)^{[p^m]}}$. 
\par 
$(2)$ Let $m$ be a positive integer in $(1)$. 
Set $B^{[p^m]}:=B\times_{\os{\circ}{B},\os{\circ}{F}{}^m_B}\os{\circ}{B}$. 
Let 
\begin{align*} 
F_m:=F_{U^{[p^{m}]}/(U^{[p^{m}]})^{\circ}}\col U^{[p^{m}]}
\lo U^{[p^{m+1}]}=(U^{[p^{m}]})^{[p]}
\tag{5.1.16.2}\label{ali:fmd}
\end{align*} 
be the abrelative Frobenius morphism of $U^{[p^{m}]}/(U^{[p^{m}]})^{\circ}$ 
over the natural morphism $B^{[p^{m}]}\lo B^{[p^{m+1}]}$. 
Then the following diagram 
\begin{equation*} 
\begin{CD} 
U^{[p^{m}]}@>{F_m}>> U^{[p^{m+1}]}\\
@V{\rho^{[m]}}VV @VV{\rho^{[m+1]}}V \\
U(0)@>{F_{U(0)}}>>U(0)
\end{CD}
\tag{5.1.16.3}\label{cd:uumr}
\end{equation*} 
is commutative. 
\par 
$(3)$ {\rm (\cite[p.~203]{ollc})} 
Let $m$ be a positive integer in $(1)$.  
Then there exists a morphism 
$\rho^{(m)}\col U\lo U(0)$ such that 
the composite morphisms 
$$\rho^{(m)}\circ \iota \col U(0)\lo U(0)$$ 
and 
$$\iota \circ \rho^{(m)} \col U\lo U(0)\lo U$$  
are $F^m_{U(0)}$ and $F^m_U$, respectively. 
\par 
$(4)$ Let the notations be as in $(3)$. 
Let $F_U\col U\lo U$ be the Frobenius endomorphism of $U$. 
Then the following diagram 
\begin{equation*} 
\begin{CD} 
U@>{F_U}>> U\\
@V{\rho^{(m+1)}}VV @VV{\rho^{(m)}}V \\
U(0)@=U(0)
\end{CD}
\tag{5.1.16.4}\label{cd:uum0r}
\end{equation*} 
is commutative. 
\end{lemm}
\begin{proof} 
(1): Because the exact closed immersion 
$\os{\circ}{\iota}\col \os{\circ}{U}(0) \os{\sus}{\lo} \os{\circ}{U}$ is nilpotent, 
we see that there exists a morphism 
$\os{\circ}{\rho}{}^{(m)}\col \os{\circ}{U}\lo \os{\circ}{U}(0)$ 
such that $\os{\circ}{\rho}{}^{(m)} \circ \os{\circ}{\iota}=\os{\circ}{F}{}^{m}_{U(0)}$ 
and $\os{\circ}{\iota}\circ \os{\circ}{\rho}{}^{(m)}=\os{\circ}{F}{}^{m}_{U}$ 
for a large integer $m$ (\cite[(2.1)]{boi}). 
Hence 
\begin{align*} 
U^{[p^m]}=U\times_{\os{\circ}{U},\os{\circ}{F}{}^m_U}\os{\circ}{U}
=U\times_{\os{\circ}{U},\iota}\os{\circ}{U}(0)
\times_{\os{\circ}{U}(0),\os{\circ}{\rho}{}^{(m)}_U}\os{\circ}{U}
=U(0)\times_{\os{\circ}{U}(0),\os{\circ}{\rho}{}^{(m)}}\os{\circ}{U}.
\tag{5.1.16.5}\label{ali:u0m0r}
\end{align*}  
Let us denote 
the projection $U^{[p^m]}
=U(0)\times_{\os{\circ}{U}(0),\os{\circ}{\rho}{}^{(m)}}\os{\circ}{U}\lo U(0)$ by 
$\rho^{[m]}$.
Because the following diagrams 
\begin{equation*} 
\begin{CD} 
U(0)@>{F^m_{U(0)}}>> U(0)@= U(0)\\
@V{\bigcap}VV @. @| \\
U@>{F_{m,0,U}}>>
U^{[p^m]}=U(0)\times_{\os{\circ}{U}(0),\os{\circ}{\rho}{}^{(m)}}\os{\circ}{U}
@>{\rho^{[m]}}>> U(0) \\
@|  @. @VV{\bigcap}V\\
U@>{F^m_{U}}>> U@= U
\end{CD}
\tag{5.1.16.6}\label{cd:uur}
\end{equation*} 
and 
\begin{equation*} 
\begin{CD} 
U^{[p^m]}=U(0)\times_{\os{\circ}{U}(0),\os{\circ}{\rho}{}^{(m)}}\os{\circ}{U}
@>{\rho^{[m]}}>> U(0) \\
@|   @VV{\bigcap}V\\
U\times_{\os{\circ}{U},\os{\circ}{F}{}^{m}_U}\os{\circ}{U}
@>{{\rm proj}.}>> U 
\end{CD}
\tag{5.1.16.7}\label{cd:0uur}
\end{equation*} 
are commutative, we obtain (1). 
\par 
(2): It is clear that the following diagram 
\begin{equation*} 
\begin{CD} 
(U^{[p^{m}]})^{\circ}=\os{\circ}{U}@= (U^{[p^{m+1}]})^{\circ}=\os{\circ}{U}\\
@V{\os{\circ}{\rho}{}^{[m]}}VV @VV{\os{\circ}{\rho}{}^{[m+1]}}V \\
\os{\circ}{U}(0)@>{F_{\os{\circ}{U}(0)}}>>\os{\circ}{U}(0)
\end{CD}
\tag{5.1.16.8}\label{cd:unf}
\end{equation*} 
is commutative. 
Because the morphism $F^*_m\col 
F^{-1}_m(M_{U^{[p^{m+1}]}}/{\cal O}_{U^{[p^{m+1}]}}^*)
=M_U/{\cal O}_U^*
\lo M_U/{\cal O}_U^*=M_{U^{[p^{m}]}}/{\cal O}_{U^{[p^{m}]}}^*$
is the multiplication by $p$, 
it is clear that the diagram (\ref{cd:uumr}) is commutative. 
\par 
(3): Let $\os{\circ}{\rho}{}^{[m]}\col \os{\circ}{U}\lo \os{\circ}{U}(0)$ be 
the underlying morphism in the proof of (1).  
The morphism $\os{\circ}{\rho}{}^{(m)}$ 
induces a morphism ${\cal O}^*_{U(0)}\lo {\cal O}^*_U$. 
Since $M_{U(0)}=M_U\oplus_{{\cal O}^*_U}{\cal O}^*_{U(0)}$, 
we see that there exists a morphism 
$\rho^{(m)}\col U\lo U(0)$ such that 
$\rho^{(m)} \circ \iota=F^m_{U(0)}$ 
and $\iota \circ \rho^{(m)} =F^m_U$. 
\par 
(4): This is obvious by the construction of $\rho^{(m)}$. 
\end{proof}

\par
The following is a filtered log $F^{\infty}$-isospan version 
of \cite[(2.18)]{of}:

\begin{prop}\label{prop:ge}
The restriction functor {\rm (\ref{ali:rstzb})} induces the following restriction functor$:$ 
\begin{align*} 
\iota^*:=\{\iota^*_n\}_{n=0}^{\infty} \col F^{\infty}{\textrm -}{\rm IsosF}(Z/B) 
\lo F^{\infty}{\textrm -}{\rm IsosF}_p(Z/B). 
\tag{5.1.17.1}\label{ali:fspkzzb} 
\end{align*} 
This is an equivalence of categories. 
\end{prop}
\begin{proof} 
By using (\ref{cd:rstzsfb}) for the case where $g$ is equal to 
the Frobenius endomorphism $F_Z$, we see that 
the functor $\{\iota^*_n\}_{n=0}^{\infty}$ is compatible with 
the structures of $F^{\infty}$-isospans. 
\par 
Following the method in the proof of \cite[(2.18)]{of} and \cite[Remark 16]{ollc},  
we construct a functor 
\begin{align*} 
Q=\{Q_n\}_{n=0}^{\inf} \col F^{\infty}{\textrm -}{\rm IsosF}_p(Z/B)
\lo F^{\infty}{\textrm -}{\rm IsosF}(Z/B),  
\tag{5.1.17.2}\label{ali:fsqipnszb} 
\end{align*} 
which will turn out to be a quasi-inverse of $\iota^*$. 
Only for simplicity of notation, we prove (\ref{prop:ge}) for the trivially filtered case.  
\par 
Let $\{(E_n,\Phi_n)\}_{n=0}^{\infty}$ be an object of 
$F^{\infty}{\textrm -}{\rm Isos}_p(Z/B)$. 
Let $n$ be a fixed nonnegative integer. 
Let $(T,z)$ be an object of ${\rm Enl}(Z/B)$;  
$z$ is the structural morphism $z\col T_{0}\lo Z$. 
Let $\iota \col T_0\os{\sus}{\lo} T_1$ be the natural exact closed immersion 
and let $F_{T_i}\col T_i\lo T_i$ be the Frobenius endomorphism of $T_i$ $(i=1,0)$. 
Apply {\rm (\ref{lemm:fianf})} (3) for the case $U(0)=T_0$ and $U=T_1$. 
Then there exist a nonnegative integer $m$ and a morphism 
$\rho^{(m)}\col T_1\lo T_0$ such that 
$\rho^{(m)} \circ \iota=F^m_{T_0}$ 
and $\iota \circ \rho^{(m)} =F^m_{T_1}$. 
Set $z^{(m)}:=z\circ \rho^{(m)}\col T_1\lo T_0\lo Z$. 
Then $(z^{(m)})^{(m')}=z^{(m+m')}$. 
Here we have considered the object 
$(T,z^{(m)})\in {\rm Enl}_p(Z/B)$ as an object of ${\rm Enl}(Z/B)$. 
By using the closed immersion  
$\iota \col T_0\os{\sus}{\lo} T_1$,  
we have the following commutative diagram 
\begin{equation*} 
\begin{CD} 
T_0 @>{\iota}>>T_1@>{\rho^{(m)}}>>T_0\\
@V{z}VV  @V{z^{(m)}}VV @V{z}VV \\
Z @>{F^m_Z}>> Z @=Z.   
\end{CD}
\tag{5.1.17.3}\label{cd:tnnspm} 
\end{equation*} 
Set  
\begin{align*} 
(Q_n(E_n))_{(T,z)}:=(E_{m+n})_{(T,z^{(m)})}.
\tag{5.1.17.4}\label{ali:nsamtz}
\end{align*} 
\par 
We claim that $(E_{m+n})_{(T,z^{(m)})}$ is independent of the choice of $m$. 
By (\ref{cd:uum0r}) we have the following equalities: 
\begin{align*} 
F_{Z}\circ z^{(m)}=F_{Z}\circ z\circ \rho^{(m)}=
z\circ F_{T_0}\circ \rho^{(m)}=z\circ \rho^{(m+1)}
=z^{(m+1)}.
\end{align*}  
Hence we have the following equality by (\ref{ali:epnsefp}): 
\begin{align*} 
(F_Z^*(E_{m+1+n})))_{(T,z^{(m)})}=(E_{m+1+n})_{(T,z^{(m+1)})}.
\tag{5.1.17.5}\label{cd:tntnspzm} 
\end{align*}  
By using the isomorphism $\Phi_{m+n}$, 
we have the following canonical isomorphism 
\begin{align*} 
(E_{m+n})_{(T,z^{(m)})}\os{\Phi_{m+n},\sim}{\longleftarrow}
(F_Z^*(E_{m+1+n}))_{(T,z^{(m)})}
=(E_{m+1+n})_{(T,z^{(m+1)})}.
\tag{5.1.17.6}\label{ali:tntidnspzm} 
\end{align*}
As a result, the sheaf $(Q_n(E_n))_{(T,z)}$ in (\ref{ali:nsamtz}) 
is independent of the choice of $m$. 
It is straightforward to check that $Q_n(E_n)$ is indeed 
an object of ${\rm Isoc}(Z/B)$ 
(because the diagram (\ref{cd:tnnspm}) is functorial)
and that the construction of $Q_n$ is functorial for $\{(E_n,\Phi_n)\}_{n=0}^{\infty}$.
We claim that 
\begin{align*} 
Q_{n+1}(F_Z^*(E_{n+1}))=F_Z^*(Q_{n+1}(E_{n+1})).
\tag{5.1.17.7}\label{ali:nmnsfntz}
\end{align*} 
Indeed, $Q_{n+1}(F_Z^*(E_{n+1}))_{(T,z)}=Q_n(E_{n})_{(T,z)}
=(E_{m+n})_{(T,z^{(m)})}$. By (\ref{lemm:fianf}) (2) and the argument above,  
we have the following equalities: 
\begin{align*} 
(F_Z^*(Q_{n+1}(E_{n+1})))_{(T,z)}&\os{{\rm (\ref{cd:tntnspzm})}}{=}(Q_n(E_{n+1}))_{(T,z^{(1)})}=(E_{m+n+1})_{(T,(z^{(1)})^{(m)})}
\\
&=(E_{m+n+1})_{(T,z^{(m+1)})}\os{{\rm (\ref{ali:tntidnspzm})}}{=}(E_{m+n})_{(T,z^{(m)})}.
\end{align*} 
Hence we have proved that the claim holds. 
Consequently we have 
the following canonical isomorphism 
\begin{align*} 
Q_n(\Phi_n) \col F_Z^*(Q_{n+1}(E_{n+1}))
=Q_{n+1}(F_Z^*(E_{n+1}))\os{\sim}{\lo} Q_n(E_n). 
\tag{5.1.17.8}\label{ali:nmnspntz}
\end{align*} 
\par 
Let $\{(E_n,\Phi_n)\}_{n=0}^{\infty}$ be an object of 
$F^{\infty}{\textrm -}{\rm Isos}_p(Z/B)$ and 
let $(T,z_1)$ be an object of ${\rm Enl}_p(Z/B)$. 
Then 
$$\iota^*_n(Q_n(E_n))_{(T,z_1)}=(Q_n(E_n))_{(T,z_1)}=(E_{n})_{(T,z_1)}.$$ 
Hence $\iota^*_n(Q_n(E_n))=E_n$. 
Consequently we see that $\iota^*$ is essentially surjective. 
\par 
The rest we have to prove is the full faithfulness of the functor $\iota^*$ 
(cf.~the former part of the proof of \cite[(2.18)]{of}). 
\par 
Let $(T,z)$ be an object of ${\rm Enl}(Z/B)$. 
Let the notations be as above. 
Since we have the natural morphism  
\begin{align*} 
(T,z^{(m)}\circ \iota)\lo (T,z^{(m)})
\tag{5.1.17.9}\label{ali:eaofm}
\end{align*} 
in ${\rm Enl}(Z/B)$ and $E$ is an isocrystal, 
\begin{align*} 
(F_Z^{m*}(E))_{(T,z)}=E_{(T,z^{(m)}\circ \iota)}
=E_{(T,z^{(m)})}=(\iota^*(E))_{(T,z^{(m)})}
\tag{5.1.17.10}\label{ali:eofm}
\end{align*}  
by the commutative diagram (\ref{cd:tnnspm}). 
Let $h=\{h_n\}_{n=0}^{\inf}\col E_1=(E_{1n})_{n=0}^{\inf}\lo (E_{2n})_{n=0}^{\inf}=E_2$ 
be a morphism in $F^{\inf}{\textrm -}{\rm Isos}(Z/B)$. 
Then we have the following commutative diagram 
\begin{equation*} 
\begin{CD} 
(\iota^*_{m+n}(E_{1,m+n}))_{(T,z^{(m)})}
@=(F_Z^{m*}(E_{1,m+n}))_{(T,z)}@>{\simeq}>>(E_{1n})_{(T,z)}\\
@VV{(\iota^*_{m+n}(h_{m+n}))_{(T,z^{(m)})}}V @VV{(F_Z^{m*}(h_{m+n}))_{(T,z)}}V@VV{(h_n)_{(T,z)}}V \\
(\iota^*_{m+n}(E_{2,m+n}))_{(T,z^{(m)})}@=(F_Z^{m*}(E_{2,m+n}))_{(T,z)}@>{\simeq}>>(E_{2n})_{(T,z)}. 
\end{CD}
\tag{5.1.17.11}\label{cd:efm}
\end{equation*} 
Hence the map 
\begin{align*} 
{\rm Hom}_{F^{\infty}{\textrm -}{\rm IsosF}(Z/B)}(E_1, E_2)
\owns h\lom \iota^*h\in 
{\rm Hom}_{F^{\infty}{\textrm -}{\rm IsosF}_p(Z/B)}(\iota^*(E_1),\iota^*(E_2))
\tag{5.1.17.12}\label{ali:pzb}  
\end{align*}
is injective. 
\par 
Conversely, let $g=\{g_n\}_{n=0}^{\inf}$ be 
an element of ${\rm Hom}(\iota^*(E_1),\iota^*(E_2))$.  
Let 
$$(h_n)_{(T,z)}\col (E_{1n})_{(T,z)}\lo (E_{2n})_{(T,z)}$$ 
be a unique morphism  
making into the following diagram commutative: 
\begin{equation*} 
\begin{CD} 
(\iota^*_{m+n}(E_{1,m+n}))_{(T,z^{(m)})}@=(F_Z^{m*}(E_{1,m+n}))_{(T,z)}@>{\simeq}>>(E_{1n})_{(T,z)}\\
@VV{(g_{m+n})_{(T,z^{(m)})}}V @. @VV{(h_n)_{(T,z)}}V \\
(\iota^*_{m+n}(E_{2,m+n}))_{(T,z^{(m)})}@=(F_Z^{m*}(E_{2,m+n}))_{(T,z)}@>{\simeq}>>(E_{2n})_{(T,z)}. 
\end{CD}
\tag{5.1.17.13;$n$}\label{cd:efem}
\end{equation*} 
Since $g$ is a morphism of $F^{\inf}$-isospans, we see 
that $h$ is also a morphism of  $F^{\inf}$-isospans 
by considering the diagram $F_Z^*((5.1.17.13;n+1))$. 
To show that the map (\ref{ali:pzb})  
is surjective, it suffices to show that 
$(g_{m+n})_{(T,z^{(m)})}
=(\iota^*_{m+n}(h_{m+n}))_{(T,z^{(m)})}=(h_{m+n})_{(T,z^{(m)})}$ 
for a given morphism 
$g=\{g_n\}_{n=0}^{\inf}\in 
{\rm Hom}_{F^{\infty}{\textrm -}{\rm IsosF}_p(Z/B)}(\iota^*(E_1),\iota^*(E_2))$. 
This equality follows from the following commutative diagram 
\begin{equation*} 
\begin{CD} 
(\iota^*_{m+n}(E_{1,m+n}))_{(T,z^{(m)}_1)}@
=(F_Z^{m*}(\iota^*_{m+n}(E_{1,m+n})))_{(T,z_1)}
@>{\simeq}>>(\iota^*_n(E_{1,n}))_{(T,z_1)}\\
@V{(g_{m+n})_{(T,z^{(m)}_1)}}VV  
@VV{((F_Z^m)^*(g_{m+n}))_{(T,z_1)}}V  @VV{(g_{n})_{(T,z_1)}}V \\
(\iota^*_{m+n}(E_{2,m+n}))_{(T,z^{(m)}_1)}@
=(F_Z^{m*}(\iota^*_{m+n}(E_{2,m+n})))_{(T,z_1)}
@>{\simeq}>>(\iota^*_n(E_{2,n}))_{(T,z_1)}. 
\end{CD}
\tag{5.1.17.14}\label{cd:effmem}
\end{equation*} 
\par 
We complete the proof of (\ref{prop:ge}). 
\end{proof}

\begin{rema}\label{rema:rmo}
Note that the proof of the formula (\ref{ali:eofm}) is different from the 
proof in \cite[(2.18)]{of}. 
It is useless to denote $(T,z^{(m)}\circ \iota)$ by 
$\iota_*((T,z^{(m)}))$ as in the proof in [loc.~cit.]; 
we have to use a fact that there exists the morphism (\ref{ali:eaofm}) and 
the assumption that $E$ is an isocrystal. 
\end{rema} 

\begin{coro}\label{coro:ndr} 
Assume that the morphism {\rm (\ref{ali:ztio})} is injective for any object 
$T\in {\rm Enl}_{\star}(Z/B)$. 
Then the restriction functor 
\begin{align*} 
\iota^*:=\{\iota_n\}_{n=0}^{\infty} \col F^{\infty}{\textrm -}{\rm IsosF}^{\rm sld}(Z/B) 
\lo F^{\infty}{\textrm -}{\rm IsosF}^{\rm sld}_p(Z/B)
\tag{5.1.19.1}\label{ali:fspkizb} 
\end{align*} 
is an equivalence of categories. 
\end{coro} 
\begin{proof} 
(\ref{coro:ndr}) follows from (\ref{prop:ge}) and (\ref{prop:fisoeq}). 
\end{proof}




In fact we can prove the following without using 
the assumption in (\ref{coro:ndr}); 
the following proof is much more complicated than 
that of (\ref{prop:ge}) because there does not 
necessarily exist a morphism 
$T_1^{[p^{m}]}\lo Z^{[p^m]}$ 
making the resulting two squares
in (\ref{cd:tnacpm}) commutative.

\begin{theo}\label{theo:zbfs} 
Let ${\cal V}$, $B$ and $Z$ be as in {\rm (\ref{defi:finf})}. 
Assume that $M^{\rm gp}_Z/{\cal O}_Z^*$ is $p$-torsion-free. 
Then the restriction functor 
\begin{align*} 
\iota^*=\{\iota^*_n\}_{n=0}^{\inf} \col F^{\infty}{\textrm -}{\rm IsosF}^{\rm sld}(Z/B) 
\lo F^{\infty}{\textrm -}{\rm IsosF}^{\rm sld}_p(Z/B)
\tag{5.1.20.1}\label{ali:fspzb} 
\end{align*} 
is an equivalence of categories. 
\end{theo}
\begin{proof} 
Only for simplicity of notation, we prove (\ref{theo:zbfs}) 
in the case of the trivial filtration. 
We need to take care of the delicate proof below. 
\par 
Recall the morphism 
$F_m \col Z^{[p^{m}]}\lo Z^{[p^{m+1}]}=(Z^{[p^{m}]})^{[p]}$ $(m\geq 0)$ ((\ref{ali:fmd})).  
Set 
\begin{align*} 
F_{m,0}:=F_{m-1} \circ \cdots \circ F_0\col Z\lo Z^{[p^{m}]}. 
\end{align*} 
\par 
We construct a functor 
\begin{align*} 
Q=\{Q_n\}_{n=0}^{\inf} \col F^{\infty}{\textrm -}{\rm IsosF}^{\rm sld}_p(Z/B)
\lo F^{\infty}{\textrm -}{\rm IsosF}^{\rm sld}(Z/B),  
\tag{5.1.20.2}\label{ali:fsqipzb} 
\end{align*} 
which will turn out to be a quasi-inverse of $\iota^*$. 
Let $\{(E_n,\Phi_n)\}_{n=0}^{\infty}$ be an object of 
$F^{\infty}{\textrm -}{\rm Isos}^{\rm sld}_p(Z/B)$. 
Let $n$ be a fixed nonnegative integer. 
Let $(T,z)$ be an object of ${\rm Enl}^{\rm sld}(Z/B)$; 
$z$ is the structural morphism $z\col T_{0}\lo Z$. 
Let $\iota \col T_{0}\os{\sus}{\lo} T_{1}$ be the natural exact closed immersion 
and let $F_{T_i}\col T_{i} \lo T_{i}$ be the Frobenius endomorphism of $T_{i}$ $(i=1,0)$. 
Apply (\ref{lemm:fianf}) (1) for the case $U(0)=T_0$ and $U=T_1$.   
$($Note that the exact closed immersion 
$\os{\circ}{\iota}\col \os{\circ}{T}_{0} \os{\sus}{\lo} \os{\circ}{T}_{1}$ is nilpotent.$)$ 
Then we have the following commutative diagram 
\begin{equation*} 
\begin{CD} 
T_0@>>>  T_1^{[p^{m}]}@>{\rho^{[m]}}>>T_0\\
@V{z}VV  @. @VV{z}V \\
Z @>{F_{m,0}}>> Z^{[p^m]}@>>>Z,  
\end{CD}
\tag{5.1.20.3}\label{cd:tnacpm} 
\end{equation*} 
where the upper horizontal composite morphism is $F^m_{T_0}$. 
Note that there does not necessarily exist a middle vertical morphism 
making the resulting new two diagrams commutative 
since $T_1^{[p^{m}]}=T_0\times_{\os{\circ}{T}_0,\os{\circ}{\rho}{}^{[m]}}\os{\circ}{T}_1$ 
and $Z^{[p^m]}=Z\times_{\os{\circ}{Z},\os{\circ}{F}{}^m_Z}\os{\circ}{Z}$ and 
since there does not necessarily exist a morphism $\os{\circ}{T}_1\lo \os{\circ}{Z}$. 
Set $z^{[m]}:=z\circ \rho^{[m]}\col T^{[p^m]}_1\lo T_0\lo Z$ 
for a large integer $m$. 
Let $U_0^m$ be a log scheme whose underlying scheme 
$\os{\circ}{T}_0$ and whose log structure is associated to the morphism 
${\rm Im}(M_{T_1^{[p^{m}]}}\lo M_{T_0})\lo {\cal O}_{T_0}$. 
Because the natural morphism 
$p\col M_Z/{\cal O}_Z^*\lo M_Z/{\cal O}_Z^*$ 
is injective, the natural morphism 
$M_{T^{[p^m]}_1}/{\cal O}_{T_1}^*
\os{``p\times{\textrm '}{\textrm '}}{\lo} M_{T_0}/{\cal O}_{T_0}^*$ 
is injective. Using the natural morphism $T_0\lo U_0^m$, 
we have a formal log scheme $T^m:=U_0^m(T)$ 
with a natural closed immersion $U_0^m\os{\sus}{\lo} T^m$ 
such that 
$(T^m)^{\circ}=\os{\circ}{T}$ and 
$M_{T^m}/{\cal O}_T^*\os{\sim}{\lo} M_{U_0^m}/{\cal O}_{T_0}^*$.  
Set $T^m_1:=T^m~{\rm mod}~p$. 
By the definition of $T^m$, 
we have the following commutative diagram 
\begin{equation*} 
\begin{CD} 
M_{T_0}/{\cal O}_{T_0}^*@<{\supset}<< M_{T^{[p^m]}_1}/{\cal O}_{T_1}^*\\
@A{\simeq}AA @VV{\simeq}V \\
M_T/{\cal O}_T^* @<{\supset}<< M_{T^m}/{\cal O}^*_T. 
\end{CD}
\tag{5.1.20.4}\label{cd:tnpmzm} 
\end{equation*} 
Hence we have the following commutative diagram 
\begin{equation*} 
\begin{CD} 
T@>>> T^m @<{\supset}<<T^m_1\\
@A{\bigcup}AA @A{\bigcup}AA @VV{\simeq}V \\
T_0 @>>> U_0^m @>{\subset}>>T_{1}^{[p^{m}]}\\
@V{z}VV @. @VV{z^{[m]}}V\\
Z @>{F_{m,0}}>> Z^{[p^m]}@>>> Z. 
\end{CD}
\tag{5.1.20.5}\label{cd:tnpzm} 
\end{equation*} 
By abuse of notation, denote the right vertical composite morphism by $z^{[m]}$ again. 
Hence we have an object $(T^m,z^{[m]})$ of 
${\rm Enl}^{\rm sld}_p(Z/B)$ and  
we can consider $(E_{m+n})_{(T^m,z^{[m]})}$. 
Set  
\begin{align*} 
(Q_n(E_n))_{(T,z)}:=(E_{m+n})_{(T^m,z^{[m]})}.
\tag{5.1.20.6}\label{ali:nmtz}
\end{align*} 
\par 
We claim that $(E_{m+n})_{(T^m,z^{[m]})}$ is independent of the choice of $m$. 
Consider the following diagram 
\begin{equation*} 
\begin{CD} 
\os{\circ}{T}_1@>{F_{\os{\circ}{T}_1}}>>\os{\circ}{T}_1 
@>{\os{\circ}{\rho}{}^{[m]}}>>\os{\circ}{T}_0\\
@. @. @VV{z}V \\
@. @. \os{\circ}{Z}. 
\end{CD} 
\end{equation*} 
By (\ref{cd:unf}) the upper horizontal composite morphism is 
equal to $\os{\circ}{\rho}{}^{[m+1]}$. 
Hence we have a natural solid morphism 
$T^{[p^{m+1}]}_1\lo T^{[p^{m}]}_1$ over 
$F_{\os{\circ}{T}_1}\col \os{\circ}{T}_1\lo \os{\circ}{T}_1$. 
We also have the natural morphism 
$T^{[p^{m}]}_1\lo T^{[p^{m+1}]}_1$ such that the composite morphism 
$T^{[p^{m}]}_1\lo T^{[p^{m+1}]}_1\lo T^{[p^{m}]}_1$ is equal to $F_{T^{[p^{m}]}_1}$. 
By the definition of $(T^m,z^{[m]})$ and by the commutative diagram 
(\ref{cd:uumr}), 
we have the following commutative diagram 
\begin{equation*} 
\begin{CD} 
T_1^{m}@>>>T_1^{m+1}\\
@V{\simeq}VV @VV{\simeq}V \\
T_1^{[p^m]}@>{F_{T_1^{[p^m]}/(T_1^{[p^m]})^{\circ}}}>> T_1^{[p^{m+1}]}\\
@V{z^{[m]}}VV @VV{z^{[m+1]}}V \\
Z @>{F_Z}>> Z. 
\end{CD}
\tag{5.1.20.7}\label{cd:tntbtpzm} 
\end{equation*} 
Let $(T',z')$ be an object of ${\rm Enl}^{\rm sld}_p(Z/B)$ such that 
there exists a morphism $(T^m,z^{[m]})\lo (T',z')$ of solid log enlargements 
over $F_Z\col Z \lo Z$. 
We claim that $(T^{m+1},z^{[m+1]})$ is the initial object for the $(T',z')$'s.  
It is clear that the underlying enlargement of $(T^{m+1},z^{[m+1]})$
is the initial object of the underlying enlargement of $(T',z')$'s 
since 
$(T_1^{[p^m]})^{\circ}=\os{\circ}{T}_1=(T_1^{[p^{m+1]}})^{\circ}$.  
Moreover, by using the solidness of the morphism $T_0'\lo Z$,
we see that the morphism $M_{T'_0}\lo M_{T^m_1}\lo M_{T^m}$ factors uniquely 
through the morphism  $M_{T^{m+1}}\lo M_{T^{m}}$. 
Hence  $(T^{m+1},z^{[m+1]})$ is the initial object for the $(T',z')$'s.  
Consequently we have the following equality: 
\begin{align*} 
(F_Z^*(E_{m+1+n}))_{(T^m,z^{[m]})}
=(E_{m+1+n})_{(T^{m+1},z^{[m+1]})}.
\tag{5.1.20.8}\label{cd:tnt1pzm} 
\end{align*}  
By using the isomorphism $\Phi_{m+n}$, 
we have the following canonical isomorphism 
\begin{align*} 
(E_{m+n})_{(T^m,z^{[m]})}\os{\Phi_{m+n},\sim}{\longleftarrow}
(F_Z^*(E_{m+1+n})))_{(T^m,z^{[m]})}
=(E_{m+1+n})_{(T^{m+1},z^{[m+1]})}.
\tag{5.1.20.9}\label{cd:tntidpzm} 
\end{align*}
As a result, the sheaf $(Q_n(E_n))_{(T,z)}$ in (\ref{ali:nmtz}) 
is independent of the choice of $m$. 
It is straightforward to check that $Q_n(E_n)$ is indeed 
an object of ${\rm Isoc}(Z/B)$ 
and that the construction of $Q_n$ is functorial.
We claim that 
\begin{align*} 
Q_{n+1}(F_Z^*(E_{n+1}))=F_Z^*(Q_{n+1}(E_{n+1})).
\tag{5.1.20.10}\label{ali:nmfntz}
\end{align*} 
Indeed, $Q_{n+1}(F_Z^*(E_{n+1}))_{(T,z)}\os{\sim}{\lo}Q_n(E_n)_{(T,z)}
=(E_{m+n})_{(T^m,z^{[m]})}$. By the argument above,  
we have the following equalities: 
\begin{align*} 
(F_Z^*(Q_{n+1}(E_{n+1})))_{(T,z)}&=(Q_{n+1}(E_{n+1}))_{(T^1,z^{(1)})}
=(E_{m+n+1})_{(T^{m+1},z^{[m+1]})}\\
&=(E_{m+n})_{(T^m,z^{[m]})}.
\end{align*} 
Hence we have the following canonical isomorphism 
\begin{align*} 
Q_n(\Phi_n) \col F_Z^*(Q_{n+1}(E_{n+1}))=Q_{n+1}(F_Z^*(E_{n+1}))\os{\sim}{\lo} 
Q_n(E_{n}). 
\tag{5.1.20.11}\label{ali:nmpntz}
\end{align*}
Set 
\begin{align*} 
Q(\{(E_n,\Phi_n)\}_{n=0}^{\infty}):=\{(Q_n(E_n),Q_n(\Phi_n))\}_{n=0}^{\infty}. 
\tag{5.1.20.12}\label{ali:nmitz}
\end{align*} 
\par 
Let $\{(E_n,\Phi_n)\}_{n=0}^{\infty}$ be an object of 
$F^{\infty}{\textrm -}{\rm Isos}_p(Z/B)$ and 
let $(T,z)$ be an object of ${\rm Enl}^{\rm sld}_p(Z/B)$. 
Then 
$$\iota^*_n(Q_n(E_n))_{(T,z)}=(Q_n(E_n))_{(T,z)}=(E_{n})_{(T,z)}.$$ 
Hence $\iota^*_n(Q_n(E_n))=E_n$. 
Consequently $\iota^*$ is essentially surjective. 
\par 
The rest we have to prove is the full faithfulness of the functor $\iota^*$. 
(The following is a log version of a part of the proof of \cite[(2.18)]{of}.)
\par 
Let $(T,z)$ be an object of ${\rm Enl}^{\rm sld}(Z/B)$. 
Let the notations be as above. 
Set $T_0^{[p^m]}:=T_0\times_{\os{\circ}{T}_0,\os{\circ}{F}{}^m_{T_0}}\os{\circ}{T}_0$. 
Then we have a natural exact closed immersion 
$\iota_m \col T_0^{[p^m]}\os{\sus}{\lo}T_1^{[p^m]}$ and 
the following commutative diagram 
\begin{equation*} 
\begin{CD} 
T_0@>>>T_0^{[p^m]}@>{\iota_m}>> T_1^{[p^m]} @>{\rho^{[m]}}>>T_0\\
@V{z}VV @VVV @. @VV{z}V\\
Z@>>>Z^{[p^m]}@= Z^{[p^m]} @>>>Z, 
\end{CD}
\end{equation*}
where the upper horizontal morphism (resp.~the lower horizontal morphism) 
is $F^m_{T_0}$ (resp.~$F^m_Z$). 
Then, by the argument for the well-definedness of $Q$, we see that 
$(F^m_Z)_*((T,z))=(T^m,T^m_0)$, where 
$T^m_0:={\rm Spec}^{\log}_{T^m}(({\cal O}_{T^m}/p)_{\rm red})\simeq T^{[p^m]}_0$.  
Hence we have the following equalities:
\begin{align*} 
(F^m_Z)_*((T,z))=(T^m,z^{[m]}\circ \iota_m)
\tag{5.1.20.13}\label{ali:zmb} 
\end{align*}
Consequently 
$(F_Z^{m*}(E_n))_{(T,z)}=(\iota^*_n(E_n))_{(T^m,z^{[m]})}$ as in (\ref{ali:eaofm}). 
Let $h:=\{h_n\}_{n=0}^{\inf} 
\col E_1=(E_{1n})_{n=0}^{\inf}\lo (E_{2n})_{n=0}^{\inf}=E_2$ 
be a morphism in $F^{\inf}{\textrm -}{\rm Isos}(Z/B)$. 
Then we have the following commutative diagram: 
\begin{equation*} 
\begin{CD} 
(i^*_{m+n}(E_{1,m+n}))_{(T^m,z^{[m]})}@=(F_Z^{m*}(E_{1,m+n}))_{(T,z)}@>{\simeq}>>(E_{1n})_{(T,z)}\\
@V{(i^*_{m+n}(h_{m+n}))_{(T^m,z^{[m]})}}VV @V{(F_Z^{m*}(h_n))_{(T,z)}}VV@VV{h_{n,(T,z)}}V \\
(i^*_{m+n}(E_{2,m+n}))_{(T^m,z^{[m]})}@=(F_Z^{m*}(E_{2,m+n}))_{(T,z)}@>{\simeq}>>(E_{2n})_{(T,z)}. 
\end{CD}
\end{equation*} 
Hence the map 
\begin{align*} 
{\rm Hom}_{F^{\infty}{\textrm -}{\rm IsosF}^{\rm sld}(Z/B)}(E_1, E_2)
\owns h\lom \iota^*h\in 
{\rm Hom}_{F^{\infty}{\textrm -}{\rm IsosF}^{\rm sld}_p(Z/B)}(\iota^*(E_1),\iota^*(E_2))
\tag{5.1.20.14}\label{ali:pmzb}  
\end{align*}
is injective. 
\par
By replacing $z^{(m)}$ in the proof of (\ref{prop:ge}) with $z^{[m]}$, we see that 
the map (\ref{ali:pmzb}) is surjective.
\par 
We complete the proof of (\ref{theo:zbfs}).  
\end{proof}


The following is a filtered variant of 
the $F^{\infty}$-isospan in \cite[\S5]{ollc} and (\ref{defi:finf}).  
Though we can dispense with this in this book, 
it may be of independent interest. 

\begin{defi}\label{defi:fisb}  
(1) Let ${\cal V}$, $B$ and $Z$ be as in (\ref{defi:finf}). 
Let $n$ be a nonnegative integer. 
Set $B^{[p^n]}:=B\times_{\os{\circ}{B},\os{\circ}{F}{}^n_B}\os{\circ}{B}$ 
and set $Z^{[p^n]}:=Z\times_{\os{\circ}{Z},\os{\circ}{F}{}^n_Z}\os{\circ}{Z}$. 
Let 
\begin{align*} 
F_n:=F_{Z^{[p^{n}]}/(Z^{[p^{n}]})^{\circ}}\col Z^{[p^{n}]}
\lo Z^{[p^{n+1}]}=(Z^{[p^{n}]})^{[p]}
\end{align*} 
be the abrelative Frobenius morphism of $Z^{[p^{n}]}$ over 
the natural morphism $B^{[p^{n}]}\lo B^{[p^{n+1}]}$. 
We call a sequence $\{((E_n,P_n),\Phi_n)\}_{n=0}^{\inf}$ of pairs 
a {\it filtered log $p$-adically convergent abrelative $($solid$)$-$F^{\inf}$-isospan} 
(resp.~a {\it filtered log convergent abrelative $($solid$)$-$F^{\inf}$-isospan}) if 
$(E_n,P_n)\in {\rm IsocF}^{\sq}_p(Z^{[p^n]}/B)$,  
(resp.~$(E_n,P_n)\in {\rm IsocF}^{\sq}(Z^{[p^n]}/B)$) 
and $\Phi_n$ is the following isomorphism 
\begin{align*} 
\Phi_n \col F^*_n((E_{n+1},P_{n+1})) \os{\sim}{\lo} (E_{n},P_{n})
\quad (n\geq 0)
\tag{5.1.21.1}\label{ali:issslep} 
\end{align*}  
in ${\rm IsocF}^{\sq}_p(Z^{[p^n]}/B)$ 
(resp.~$(E,P)\in {\rm IsocF}^{\sq}(Z^{[p^n]}/B)$). 
We define a {\it morphism} of filtered log ($p$-adically) 
convergent abrelative (solid)$F^{\inf}$-isospans in an obvious way. 
We denote the category of 
filtered log $p$-adically convergent abrelative $F^{\infty}$-isospans 
(resp.~the category of 
filtered log convergent abrelative (solid)-$F^{\infty}$-isospans)
by ${}_{\rm ar}F^{\infty}{\textrm -}{\rm IsosF}^{\sq}_p(Z/B)$ and 
${}_{\rm ar}F^{\infty}{\textrm -}{\rm IsosF}^{\sq}(Z/B)$, 
respectively. 
The category ${}_{\rm ar}F^{\infty}{\textrm -}{\rm IsosF}^{\sq}_{\star}(Z/B)$  
is an additive category. 
\par 
When the filtration is trivial, we denote 
${}_{\rm ar}F^{\infty}{\textrm -}{\rm IsosF}^{\sq}_p(Z/B)$ and 
${}_{\rm ar}F^{\infty}{\textrm -}{\rm IsosF}^{\sq}(Z/B)$ by 
${}_{\rm ar}F^{\infty}{\textrm -}{\rm Isos}^{\sq}_p(Z/B)$ and 
${}_{\rm ar}F^{\infty}{\textrm -}{\rm Isos}^{\sq}(Z/B)$, respectively. 
\par 
(2) Let the notations be as in (1). 
Let ${}_{\rm ar}F^{\infty}{\textrm -}{\rm IsocF}^{\sq}_{\star}(Z/B)$ 
be a full subcategory of 
${}_{\rm ar}F^{\infty}{\textrm -}{\rm IsosF}^{\sq}_{\star}(Z/B)$ 
whose objects are $\{((E_n,P_n),\Phi_n)\}_{n=0}^{\inf}$'s 
such that $E_n=V^*_{0,n}(E_0)$. 
\end{defi}

\begin{prop}\label{prop:fff}
There exists the following functorial functor 
\begin{align*} 
F^{\infty}{\textrm -}{\rm IsosF}^{\sq}_{\star}(Z/B)
\lo {}_{\rm ar}F^{\infty}{\textrm -}{\rm IsosF}^{\sq}_{\star}(Z/B). 
\tag{5.1.22.1}\label{ali:fnfe}
\end{align*}
Here the functoriality means the contravariant functoriality for $Z/B$'s. 
In particular, there exists the following functorial functor 
\begin{align*} 
F{\textrm -}{\rm IsocF}^{\sq}_{\star}(Z/B)
\lo {}_{\rm ar}F^{\infty}{\textrm -}{\rm IsocF}^{\sq}_{\star}(Z/B). 
\tag{5.1.22.2}\label{ali:ffe}
\end{align*}
The latter functor is faithful. 
\end{prop} 
\begin{proof} 
Let $(E,P)=\{(E_n,P_n),\Phi_n\}_{n=0}^{\inf}$ be an object of 
$F^{\infty}{\textrm -}{\rm IsosF}^{\sq}_{\star}(Z/B)$. 
We omit to write the filtration $P$. 
The Frobenius endomorphism $F_{Z^{[p^n]}}\col Z^{[p^n]}\lo Z^{[p^n]}$ 
of $Z^{[p^n]}$ factors naturally through a morphism 
$V_n\col Z^{[p^{n+1}]}\lo Z^{[p^n]}$.  
Set $V_{0,n}:=V_0\circ \cdots \circ V_{n-1}\col Z^{[p^n]}\lo Z$. 
Consider the following diagram  
$$Z^{[p^n]}\os{F_n}{\lo} Z^{[p^{n+1}]}\os{V_n}{\lo} Z^{[p^n]}\os{V_{0,n}}{\lo} Z.$$ 
Set $E'_n:=V_{0,n}^*(E_n)\in {\rm Isoc}^{\sq}_{\star}(Z^{[p^n]}/B^{[p^n]})$. 
Then we have the following equalities and the following isomorphism: 
\begin{align*} 
\Phi'_n\col F_n^*(E'_{n+1})&=F_n^*V^*_{0,n+1}(E_{n+1})=F_n^*V^*_nV^*_{0,n}(E_{n+1})
=(V_n\circ F_n)^*V^*_{0,n}(E_{n+1})\\
&=F_{Z^{[p^n]}}^*V^*_{0,n}(E_{n+1})=V^*_{0,n}F_Z^*(E_{n+1})
\os{V^*_{0,n}(\Phi_n),\simeq}{\lo} V^*_{0,n}(E_n)=E'_n. 
\end{align*} 
Hence we obtain an object 
$\{(E'_n,\Phi'_n)\}_{n=0}^{\inf}$ of 
${}_{\rm ar}F^{\infty}{\textrm -}{\rm IsosF}^{\sq}_{\star}(Z/B)$.
Let 
\begin{align*} 
h=\{h_n\}_{n=0}^{\inf}\col \{(E'_{1n},\Phi'_{1n})\}_{n=0}^{\inf}\lo 
\{(E'_{2n},\Phi'_{2n})\}_{n=0}^{\inf}
\tag{5.1.22.3}\label{ali:e12}
\end{align*} 
be a morphism in 
$F^{\infty}{\textrm -}{\rm IsosF}^{\sq}_{\star}(Z/B)$.
Then we obtain a morphism $\{V_{0,n}^*(h_n)\}_{n=0}^{\inf}$ in 
${}_{\rm ar}F^{\infty}{\textrm -}{\rm IsosF}^{\sq}_{\star}(Z/B)$.
By the construction, the functoriality is clear. 
\par 
The faithfulness 
of the functor (\ref{ali:ffe}) is obvious. 
\end{proof}

\par 
Let us define the action of a morphism in 
the positive part $W^+_{\rm crys}(K)$ 
of the homomorphisms in the crystalline Weil category 
in \cite[p.~190]{boi} and \cite[p.~797]{of} 
in our situation.
\par 
Because I think that the crystalline Weil category itself and 
objects and homomorphisms in  the crystalline Weil category
are confusing in \cite{of} 
(e.~g.~, $W_{\rm cris}(K)$ in [p.~796], 
$\psi =(a,b,d)\in W_{\rm crys}(K)$ and 
$K'\in {\rm Ob}W_{\rm crys}(K)$ in \cite[p.~797]{of}), 
we have to make a clarification to distinguish them. 
\par 
For a finite extension of $\kap'$ of $\kap$, 
set $K_0(\kap'):={\rm Frac}({\cal W}(\kap'))$. 
(In [loc.~cit.] $K_0(\kap')$ has been denoted by 
$K(\kap')$.)
For a morphism $c\col K'\lo K''$ of finite extensions of $K$, 
let $K_0(c)\col K_0(\kap')\lo K_0(\kap'')$ be the induced morphism 
of fields by $c$, where $\kap'$ and $\kap''$ are 
residue fields of the integer rings ${\cal V}':={\cal O}_{K'}$ 
and ${\cal V}'':={\cal O}_{K''}$ of $K'$ and $K''$, respectively.  
Let $K$ be the fraction field of ${\cal V}$. 
Let ${\rm CrW}(K)$ be the crystalline Weil category 
defined by Ogus (\cite[p.~190]{boi}, \cite[p.~796]{of}): 
an object of ${\rm CrW}(K)$ is a finite extension $K'$ of $K$; 
a morphism $K'\lo K''$ in ${\rm CrW}(K)$ is a triple $(b,a,d)$, 
where $d$ is an integer, $a\col K'\lo K''$ is a morphism of $K$-algebras 
and $b\col K'\lo K''$ is a morphism of fields such that  
``$b=a\circ F^d$'' and ``$b=F^d \circ a$'', i.~e.,  
$b\vert_{K_0(\kap')}=K_0(a)\circ F^d$ and $b\vert_{K_0(\kap')}=F^d\circ K_0(a)$, 
where the first $F^d$ (resp.~the second $F^d$) 
is the $d$-th power of the Frobenius endomorphism $K_0(\kap')$ 
(resp.~$K_0(\kap'')$). 
We denote $(d,a,b)$ in [loc.~cit.] by $(b,a,d)$ 
because we shall consider ``the action of $b$'' in (\ref{ali:pozz}) below.
Set $\psi=(b,a,d)$ and $d:={\rm deg}(\psi)$. 
For two homomorphisms $\psi$ and $\psi'$ in ${\rm CrW}(K)$, 
The composite 
$${\rm Hom}_{{\rm CrW}(K)}(K',K'')\times
{\rm Hom}_{{\rm CrW}(K)}(K'',K''')\lo {\rm Hom}_{{\rm CrW}(K)}(K',K''')$$ 
is defined in an evident way.   
Let $W_{\rm crys}(K)$ be the set of $\psi$'s above. 
Let $W^+_{\rm crys}(K)$ be a subset of $W_{\rm crys}(K)$
whose elements are those of $W_{\rm crys}(K)$ and whose 
morphisms are $\psi$'s such that ${\rm deg}(\psi)\geq 0$ ([loc.~cit., p.~797]). 
\par 
Take $\psi=(b,a,d)\in W^+_{\rm crys}(K)$. 
Then $a\col K'\lo K''$ induces a morphism 
$a_{\cal V}\col {\cal V}'\lo {\cal V}''$. 
Let $B$ and $Z$ be as in the beginning of this section. 
Set $Z':=Z\otimes_{\cal V}{\cal V}'$, $B':=B\otimes_{\cal V}{\cal V}'$,  
$Z'':=Z'\otimes_{{\cal V}',a_{\cal V}}{\cal V}''$ and 
$B'':=B\otimes_{{\cal V}',a_{\cal V}}{\cal V}''$. 
Let $a_Z\col  Z''\lo Z'$ be the induced base change morphism by $a_{\cal V}$. 
Assume that $\os{\circ}{Z}$ is a scheme over $\kap$. 
Set 
$$\psi_0:=F^d_{Z'}\circ a_{Z} \col Z''\lo Z'.$$ 
For a family $(E,P)=\{(E_n,P_n)\}_{n=0}^{\inf}$ of objects of 
${\rm IsosF}^{\sq}(Z/B)$, 
set $(E',P'):=\{(E'_n,P'_n)\}_{n=0}^{\inf}$ 
(resp.~$(E'',P''):=\{(E''_n,P''_n)\}_{n=0}^{\inf}$), 
where $(E'_n,P'_n)\in  {\rm IsocF}^{\sq}(Z'/B')$ 
(resp.~$(E''_n,P''_n)\in {\rm IsocF}^{\sq}(Z''/B'')$) 
is the pull-back of $(E_n,P_n)$. 
Set 
\begin{align*} 
\psi((E'_{n+d},P'_{n+d})):=\psi_0^*((E'_{n+d},P'_{n+d}))\in 
{\rm IsocF}^{\sq}(Z''/B'').
\tag{5.1.22.4}\label{ali:pozz}
\end{align*} 

\begin{defi}\label{deif:convfi} 
An action of $W^+_{\rm crys}(K)$ on a family 
$(E,P)=\{(E_n,P_n)\}_{n=0}^{\inf}$ of an object of ${\rm IsosF}^{\sq}(Z/B)$ 
is, by definition, for an element $\psi=(b,a,d) \in W^+_{\rm crys}(K)$, 
to give an isomorphism 
\begin{align*} 
\Phi_{\psi}\col \psi((E'_{n+d},P'_{n+d})) \os{\sim}{\lo} (E''_n,P''_n) 
\quad (n\geq 0)
\tag{5.1.23.1}\label{ali:nmgwtz} 
\end{align*} 
satisfying the following two conditions: 
\par 
(a) If $\psi=(a,a,0)\col K'\lo K''$, then $\Phi_{\psi}$ is 
the following canonical isomorphism 
\begin{align*} 
\psi((E'_n,P'_n))=a_Z^*((E'_n,P'_n))\os{\sim}{\lo} (E''_n,P''_n).
\end{align*} 
\par 
(b)  
If $\psi \col K'\lo K''$ and $\psi'=(b',a',d') \col K''\lo K'''$ 
are morphisms in $W^+_{\rm crys}(K)$, 
then $\Phi_{\psi'}\circ \psi_0'{\!}^*(\Phi_{\psi})=\Phi_{\psi' \circ \psi}\col 
\psi'\psi((E'_{n+d+d'},P'_{n+d+d'})) 
\os{\sim}{\lo} (E'''_n,P'''_n)$. 
Here ${\cal V}''':={\cal O}_{K'''}$, $Z''':=Z''\otimes_{{\cal V}'',a'_{\cal V}}{\cal V}'''$, 
$\psi_0':=F^{d'}_{Z''}\circ a'_{Z'} \col Z'''\lo Z''$ 
and $(E'''_n,P'''_n)\in {\rm IsocF}^{\sq}(Z'''/{\cal V}''')$ is the pull-back of $(E''_n,P''_n)$. 
\end{defi} 

\begin{defi}
A filtered $F^{\inf}$-isospan is a pair  
$(\{(E_n,P_n)\}^{\inf}_{n=0},\Phi)$, where $(E_n,P_n) \in {\rm IsosF}^{\sq}(Z/B)$ 
and $\Phi$ is an action of $W^+_{\rm crys}(K)$ defined in (\ref{deif:convfi}).
A morphism of $F^{\inf}$-isospans is a morphism in 
${\rm IsosF}^{\sq}(Z/B)$ which is compatible with the action of 
$W^+_{\rm crys}(K)$. 
By using these notions, we obtain a category 
$F^{\infty}{\textrm -}{\rm IsosF}(Z/B)$ of 
filtered $F^{\inf}$-isospans on $Z/B$. 
\end{defi}

\par 
As in \cite[p.~798]{of}, we obtain the following: 

\begin{prop}\label{prop:wofth}
Let ${\cal V}$, $B$ and $Z$ be as in {\rm (\ref{defi:finf})}. 
Then a filtered $F^{\inf}$-isospan in the sense of 
{\rm (\ref{deif:convfi})} is equivalent to that in the sense of {\rm (\ref{defi:finf})}. 
\end{prop}
\begin{proof}
The explanation in [loc.~cit.] works in our situation. 
We give the proof of (\ref{prop:wofth}) for the completeness of this book. 
\par 
First assume that we are given an action of $W_{\rm crys}^+(K)$ on 
$(E,P)=\{(E_n,P_n)\}_{n=0}^{\inf}$ of an object of ${\rm IsosF}^{\sq}(Z/B)$.  
Set $\phi:=(F_K,{\rm id}_{K},1)\in W_{\rm crys}^+(K)$, 
where $F_K$ is the Frobenius endomorphism of $K$. 
Then $\phi((E_{n+1},P_{n+1}))=
F^*_{Z}((E_{n+1},P_{n+1}))$ by (\ref{ali:pozz}). 
Hence the action (\ref{ali:nmgwtz}) gives the datum (\ref{ali:issslep}).
Conversely we are given the datum (\ref{ali:issslep}) and a morphism  
$\psi=(b,a,d)$ in $W_{\rm crys}^+(K)$. 
Let $\iota'\col K_0(\kap')\os{\sus}{\lo} K'$ and 
$\iota''\col K_0(\kap'')\os{\sus}{\lo} K''$ be the natural inclusions. 
Set $\iota':=(\iota',\iota',0)$, 
$\iota'':=(\iota'',\iota'',0)$ and $\al:=(K_0(a),K_0(a),0)$ 
in $W^+_{\rm crys}(K)$. 
Then $\psi \circ \iota'=\iota''\circ \al \circ \phi^d$ as in [loc.~cit]. 
Hence, by the conditions (a) and (b), 
the isomorphism 
$\Phi_{\psi}\col \psi((E'_{n+d},P'_{n+d})) \os{\sim}{\lo} (E''_n,P''_n)$  
is uniquely determined. 
\end{proof}

\begin{lemm}\label{lemm:esf}
Assume that, for any object 
$E\in {\rm Isoc}^{\sq}(Z/B)$ and any object $T\in {\rm Enl}(Z/B)$, 
$E_T$ is a flat ${\cal K}_T$-module. Then ${\rm Isoc}^{\sq}(Z/B)$ has an internal hom. 
\end{lemm} 
\begin{proof} 
Let $E$ and $F$ be objects of ${\rm Isoc}(Z/B)$. 
Since $\os{\circ}{T}$ is noetherian, 
${\cal H}{\it om}_{{\cal K}_T}(E_T,F_T)$ is locally free 
${\cal K}_T$-modules. Hence (\ref{lemm:esf}) is obvious. 
\end{proof}

As in \cite[(2.21)]{of}, we obtain the following: 
\begin{prop}\label{prop:trecd}
Let the assumption be as in {\rm (\ref{lemm:esf})}. 
Let ${\cal V}'/{\cal V}$ be a totally ramified extension of 
complete discrete valuation rings of mixed characteristics.
Set $B':=B\hat{\otimes}_{\cal V}{\cal V}'$. 
Then the restriction functor 
$$R:=\{R_n\}_{n=0}^{\inf} \col 
F^{\infty}{\textrm -}{\rm IsosF}^{\rm sld}(Z/B)\lo 
F^{\infty}{\textrm -}{\rm IsosF}^{\rm sld}(Z/B')$$ 
is an equivalence of categories.   
\end{prop}
\begin{proof} 
The explanation in [loc.~cit.] works in our situation. 
We give the proof of (\ref{prop:trecd}) for the completeness of this book. 
\par 
We omit to write the filtration $P$. 
Let $(T,z)$ be an object of ${\rm Enl}^{\rm sld}(Z/B)$. 
Let $(T',z')$ be an object of ${\rm Enl}^{\rm sld}(Z/B')$ 
obtained by the base change morphism 
with respect to the morphism 
${\rm Spf}({\cal V}')\lo {\rm Spf}({\cal V})$. 
Let $q\col T'\lo T$ be the projection. 
Then $R_n(E_n)_{(T',z')}=q^*((E_n)_{(T,z)})$. 
Obviously the natural morphism 
$(E_n)_{(T,z)}\lo q_*q^*((E_n)_{(T,z)})$ is injective. 
Using this injectivity, we see that the functor $R_n$ is faithful.
(We do not need the existence of the internal Hom in ${\rm Isoc}(Z/B)$ 
which was used in [loc.~cit.] in the trivial log case.) 
Hence we may assume that the extension $K'/K$ is Galois 
by taking the Galois closure of $K'/K$ as in [loc.~cit.] 
and we assume this 
until the end of this proof.  
Let $b$ be an element of ${\rm Gal}(K'/K)$. 
Let $i \col K\os{\sus}{\lo} K'$ be the inclusion. 
Set $\bet :=(b,{\rm id}_{K'},0)\in W^+_{\rm crys}(K')$ 
and $\iota:=(i,i,0)\in W^+_{\rm crys}(K)$.  
Then $b$ induces a natural morphism $b_{B'}\col B'\lo B'$. 
Let $(T,z)$ and $(T',z')$ be as above. 
Since $K'/K$ is totally ramified, we have the following commutative diagram: 
\begin{equation*} 
\begin{CD} 
Z@>{\bet_0}>> Z \\
@VVV @VVV \\ 
B'@>{b_{B'}}>> B'. 
\end{CD} 
\end{equation*} 
We also have the solid morphism 
$b_{T'}\col T'\lo T'$ over $b_{B'} \col B'\lo B'$. 
Hence $\bet(E'_n)=b_{T'}^*(E'_n)$
by (\ref{ali:epgefp}). 
Since $\bet \circ \iota=\iota$ in $W^+_{\rm crys}(K)$, 
we have the equality 
$\Phi_{\bet}\circ \bet^*_0(\Phi_{\iota})=\Phi_{\iota}\col \bet(E'_n) 
\os{\sim}{\lo} E'_n$. 
Consequently 
$$(\Phi_{\bet})_{T'}\circ b^*_{T'}((\Phi_{\iota})_{T'})=(\Phi_{\iota})_{T'}\col 
b^*_{T'}((E'_n)_{T'}) \os{\sim}{\lo} (E'_n)_{T'}.$$ 
This is nothing but a Galois descent data on $(E'_n)_{T'}$. 
By Galois descent, we see that $R_n$ is fully faithful. 
Furthermore, $\{\Phi_{\bet}\}_{\bet \in {\rm Gal}(K'/K)}$ 
defines a descent data for each object of 
${\rm IsocF}(Z'/B')$ and each $T'\lo T$ for any object 
$(T',z')\in {\rm Enl}^{\sq}(Z/B')$. 
Hence $R$ is essentially surjective. 
Consequently $R$ is an equivalence of categories. 
\end{proof}

\par 
Now assume that $Z$ is hollow. Then the image of 
any local section of 
$M_{Z_{\os{\circ}{T}_i}}\setminus {\cal O}_{Z_{\os{\circ}{T}_i}}^*$ 
in ${\cal O}_{Z_{\os{\circ}{T}_i}}$ is zero.
Hence the image of any local section of 
$M_{Z(T)}\setminus {\cal O}_T^*$ in ${\cal O}_T$ is 
topologically nilpotent. 

\begin{lemm}\label{lemm:ijk}
Assume that $Z$ is hollow. Then the following hold$:$
\par 
$(1)$ The morphism {\rm (\ref{ali:ztio})} is injective. 
\par 
$(2)$ For a morphism 
$g\col (T',z'_i)\lo (T,z_i)$ in ${\rm Enl}_{\star}(Z/B)$ $(i=0,1)$, 
the induced pull-back morphism 
$g^*\col g^*(M_{Z(T)})\lo M_{Z(T')}$ induces a morphism 
$g^*\col 
g^*(M_{Z(T)}\setminus {\cal O}_{T}^*)\lo 
M_{Z(T')}\setminus {\cal O}_{T'}^*$. 
\end{lemm}
\begin{proof} 
(1): Obvious. 
\par 
(2): 
Assume that there exists a local section 
$m\in M_{Z(T)}\setminus {\cal O}_{T}^*$ such that $g^*(m)\in {\cal O}_{T'}^*$. 
Let $\ol{m}$ be the image of $m$ in 
$M_{Z(T)}/{\cal O}_T^*=M_{Z_{\os{\circ}{T}_i}}/{\cal O}_{T_i}^*
=z_i^{-1}(M_Z/{\cal O}_Z^*)$. 
Then $\ol{m}$ is not a local unit section. 
Since $M_Z$ is hollow, the image of $\ol{m}$ in ${\cal O}_{T_i}$ is zero. 
This is a contradiction by the following commutative diagram 
\begin{equation*} 
\begin{CD} 
g^{-1}(M_{Z_{\os{\circ}{T}_i}})@>>> M_{Z_{\os{\circ}{T}{}'_i}} \\
@VVV @VVV \\
g^{-1}({\cal O}_{T_i})@>>> {\cal O}_{T'_i}. 
\end{CD} 
\end{equation*}
\end{proof}

\begin{defi}\label{defi:whlp} 
We say that $(T,z_1)$, $(T,z_0)$ or $T$ is {\it restrictively hollow with respect to} $Z$ 
(and the log structure of $T$ is {\it restrictively hollow with respect to} $Z$) if $Z(T)$ is hollow. 
\end{defi} 

\parno 
It is obvious that, if the log structure of $T$ is hollow,  
then the log structure of $T$ is restrictively hollow with respect to $Z$. 
We define a morphism of (restrictively) hollow log ($p$-adic) enlargements in an obvious way. 
Let ${\rm hEnl}^{\sq}_p(Z/B)$,  ${\rm rhEnl}^{\sq}_p(Z/B)$, 
${\rm hEnl}^{\sq}(Z/B)$ and ${\rm rhEnl}^{\sq}(Z/B)$ be the category of 
hollow log $p$-adic enlargements of $Z/B$, 
the category of restrictively hollow (solid) log $p$-adic enlargements of $Z/B$, 
the category of hollow (solid) log enlargements of $Z/B$ and 
the category of restrictively hollow (solid) log enlargements of $Z/B$, respectively 
(In \cite{ollc} Ogus has already considered 
categories ${\rm hEnl}^{\sq}_{\star}(Z/B)$ $(\star=p$ or nothing).). 
Then we have the following commutative diagram: 
\begin{equation*}
\begin{CD}  
{\rm hEnl}^{\sq}_p(Z/B) @>{\subset}>> {\rm rhEnl}^{\sq}_p(Z/B) @>{\subset}>>  {\rm Enl}^{\sq}_p(Z/B)\\
@VVV @VVV @VVV \\ 
{\rm hEnl}^{\sq}(Z/B) @>{\subset}>> {\rm rhEnl}^{\sq}(Z/B) @>{\subset}>> {\rm Enl}^{\sq}(Z/B). 
\end{CD}
\tag{5.1.29.1}\label{cd:rzbb} 
\end{equation*} 
Here the horizontal $\os{\subset}{\lo}$'s mean that the sources 
are full subcategories of the targets and 
the vertical arrows mean natural functors. 


\begin{prop}[{\bf cf.~\cite[Proposition 23]{ollc}}]\label{prop:nupf} 
Let $\pi\col Z\lo B$ be the structural morphism. 
Assume that $B=({\rm Spf}({\cal V}),{\cal V}^*)$ and that $Z$ is hollow. 
Assume that $M_Z$ is split, i.~e., 
there exists a subsheaf $N$ of monoids on $Z$ 
with unit local section $e$ such that  $M_Z=N\oplus {\cal O}_Z^*$ 
with structural morphism $N\owns n \lom 0\in {\cal O}_Z$ $(n\not= e)$.  
Let $(\os{\circ}{T},\os{\circ}{z}_1)$ $($resp.~$(\os{\circ}{T},\os{\circ}{z}_0))$ 
be an object of ${\rm Enl}_p(\os{\circ}{Z}/{\cal V})$ 
$($resp.~${\rm Enl}(\os{\circ}{Z}/{\cal V}))$.  
Then there exists the following functors  
\begin{align*}
{\rm Enl}_p(\os{\circ}{Z}/{\cal V})\owns (\os{\circ}{T},\os{\circ}{z}_1)\lom 
(Z(\os{\circ}{T}),z_1)\in {\rm Enl}_p(Z/{\cal V})
\tag{5.1.30.1}\label{ali:zbtz} 
\end{align*}  
and 
\begin{align*}
{\rm Enl}(\os{\circ}{Z}/{\cal V})\owns (\os{\circ}{T},\os{\circ}{z}_0)
\lom (Z(\os{\circ}{T}),z_0)\in {\rm Enl}(Z/{\cal V}), 
\tag{5.1.30.2}\label{ali:zpzbtz} 
\end{align*}  
where the log $p$-adic formal scheme $Z(\os{\circ}{T})/B$ 
and the morphisms $z_1$, $z_0$ will be constructed in the proof below.
Moreover one can consider 
${\rm Enl}_p(\os{\circ}{Z}/{\cal V})$ and ${\rm Enl}(\os{\circ}{Z}/{\cal V})$ 
as full subcategories of ${\rm Enl}_p(Z/{\cal V})$ and ${\rm Enl}(Z/{\cal V})$ by 
the functors {\rm (\ref{ali:zbtz})} and {\rm (\ref{ali:zpzbtz})}, respectively. 
Consequently the restriction of an object of 
${\rm IsocF}_{\star}(Z/{\cal V})$ to ${\rm Enl}_{\star}(\os{\circ}{Z}/{\cal V})$ 
gives an object of 
${\rm IsocF}_{\star}(\os{\circ}{Z}/{\cal V})$.
\end{prop} 
\begin{proof} 
Set $Z_{\os{\circ}{T}_i}:=Z\times_{\os{\circ}{Z}}\os{\circ}{T}_i$ 
$(i=1)$ (resp.~$(i=0)$) with projection $z_i \col Z_{\os{\circ}{T}_i}\lo Z$.  
Since $\os{\circ}{T}$ is homeomorphic to $Z_{\os{\circ}{T}_i}$, 
we can consider $z_i^{-1}(N)$ as a sheaf of monoids on $\os{\circ}{T}$. 
Endow $\os{\circ}{T}$ with the log structure 
$z_i^{-1}(N)\oplus {\cal O}_T^*$ 
with structural morphism $z_i^{-1}(N)\owns n \lom 0\in {\cal O}_{\os{\circ}{T}}$ 
$(n\not= e)$. 
Let $Z(\os{\circ}{T})$ be the resulting log formal scheme. 
This is hollow and fine. 
Because $M_B={\cal O}_B^*$, 
the structural morphism $p\col \os{\circ}{T}\lo \os{\circ}{B}$ induces 
a morphism $Z(\os{\circ}{T})\lo B$ of log formal schemes. 
Consequently we have the following commutative diagram 
\begin{equation*} 
\begin{CD} 
Z_{\os{\circ}{T}_i}@>{\subset}>> Z({\os{\circ}{T}}) \\
@V{z_i}VV @VVV \\ 
Z@>>> B. 
\end{CD} 
\end{equation*} 
and we obtain a solid hollow log ($p$-adic) enlargement $(Z(\os{\circ}{T}),z_i)$ of $Z/B$. 
(The log formal scheme $Z(\os{\circ}{T})/B$ is independent of the choice of $i$.) 
It is easy to check that 
the morphisms (\ref{ali:zbtz}) and (\ref{ali:zpzbtz}) are functorial. 
Conversely, a morphism 
$(Z(\os{\circ}{T}),z_i)\lo (Z(\os{\circ}{T}{}'),z_i)$ 
in ${\rm Enl}_{\star}(Z/B)$ induces a morphism 
$(\os{\circ}{T},z_i)\lo (\os{\circ}{T}{}',z_i)$ in 
${\rm Enl}_{\star}(\os{\circ}{Z}/\os{\circ}{B})$. 
Hence, using the morphisms (\ref{ali:zbtz}) and (\ref{ali:zpzbtz}), 
we can consider 
${\rm Enl}_p(\os{\circ}{Z}/\os{\circ}{B})$ and ${\rm Enl}(\os{\circ}{Z}/\os{\circ}{B})$ 
as full subcategories of ${\rm Enl}_p(Z/B)$ and ${\rm Enl}(Z/B)$, respectively. 
\end{proof} 



\parno 
One may think that the notation $Z(\os{\circ}{T})$ is misleading. 
However the following tells us that we have a reasoning for this notation:

\begin{prop}\label{prop:csds} 
Let the assumptions be as in {\rm (\ref{prop:nupf})}.  
Let $(T,z_1)$ $($resp.~$(T,z_0))$ 
be an object of ${\rm rhEnl}^{\sq}_p(Z/{\cal V})$ $($resp.~${\rm rhEnl}^{\sq}(Z/{\cal V}))$.  
Assume that $M_T$ is split. 
$($This is automatically satisfied when $\sq=$sld.
Then $(Z(T),Z(z_i))\simeq (Z(\os{\circ}{T}),z_i)$. 
$($This isomorphism depends on the splitting of $M_T)$. 
\end{prop} 
\begin{proof} 
Since $M_T$ is split, there exists a subsheaf $L_T$ of monoids of $M_T$  
such that the composite isomorphism 
$L_T\os{\sus}{\lo}M_T\lo M_T/{\cal O}_T^*$ is an isomorphism. 
Let $N_T$ be the inverse image $z^{-1}_i((N\oplus {\cal O}_T^*)/{\cal O}_T^*)$ 
by the composite morphism
$L_T\os{\sim}{\lo} M_T/{\cal O}_T^*=M_{T_i}/{\cal O}_{T_i}^* 
\supset z^{-1}_i(M_Z/{\cal O}_Z^*)$. 
Obviously we have the isomorphism $N_T\os{\sim}{\lo}z^{-1}_i(N)$.  
By (\ref{lemm:ijk}) (1), the morphism (\ref{ali:ztio}) is injective. 
Hence $M_{Z(T)}=N_T\oplus {\cal O}_T^*$.  
Since $(T,z_i)\in {\rm rhEnl}^{\sq}_{\star}(Z/B)$, 
the image of any non-identity section of $N_{T}$ is zero in ${\cal O}_T$. 
Hence $(Z(T),z_i)\simeq (Z(\os{\circ}{T}),z_i)$ over ${\cal V}$. 
\end{proof} 


\begin{prop}\label{prop:plaec}
Let $\pi$ be a uniformizer of ${\cal V}$ and 
set $B_n:=B\mod \pi^n$ $(1\leq n\leq \infty)$ $(B_{\infty}:=B)$. 
Assume that $Z$ is a log formal scheme over $B_n$ and 
that $Z$ has a hollow lift $\wt{Z}$ over $B$, 
which is formally log smooth over $B$. 
Then ${\rm rhEnl}^{\sq}_p(Z/B)={\rm Enl}^{\sq}_p(Z/B)$ and 
${\rm rhEnl}^{\sq}(Z/B)={\rm Enl}^{\sq}(Z/B)$. 
\end{prop} 
\begin{proof} 
Let $T=(T,z_1)$ (resp.~$T=(T,z_0)$) be an object of ${\rm Enl}_p(Z/B)$ 
(resp.~${\rm Enl}(Z/B)$), respectively.  
We would like to prove that the image of any section of 
$M_{Z(T)}\setminus {\cal O}_T^*$ in ${\cal O}_T$ is zero. 
Hence we may assume that $T$ is affine. 
Because $\wt{Z}$ is formally log smooth over $B$, 
we have a morphism $T\lo \wt{Z}$ over $B$ extending the morphism 
$z_i\col T_i\lo Z$ $(i=1,0)$ over $B$. 
Then $M_{Z(T)}$ is nothing but the inverse image of 
the log structure of $\wt{Z}$. 
Hence we see that the image of any section of 
$M_{Z(T)}\setminus {\cal O}_T^*$ in ${\cal O}_T$ is zero. 
\end{proof}

\par 
Henceforth assume that the log structure of $Z$ is constant in the sense of 
\cite[Definition 4]{ollc}, that is, the sheaf $M_Z/{\cal O}_Z^*$ is locally constant. 
\par 
For a (solid) log ($p$-adic) enlargement $T$ of ${\rm Enl}^{\sq}_{\star}(Z/B)$, 
we say that $T$ is {\it constant} if $M_T/{\cal O}_T^*$ is locally constant as above. 
Let ${\rm CEnl}^{\sq}_{\star}(Z/B)$ be the full subcategory of 
${\rm Enl}^{\sq}_{\star}(Z/B)$ whose objects are constant 
($p$-adic) enlargements of $Z/B$.  
The functor (\ref{ali:tzi}) factors through 
${\rm CEnl}^{\sq}_{\star}(Z/B)\os{\sus}{\lo}  {\rm Enl}^{\sq}_{\star}(Z/B)$ 
and we have the following functor 
\begin{align*} 
Z(?)\col {\rm Enl}^{\sq}_{\star}(Z/B)\owns T\lom Z(T)\in  {\rm CEnl}^{\sq}_{\star}(Z/B). 
\tag{5.1.32.1}\label{cd:crabb} 
\end{align*} 
The functor $Z(?)$ induces the following equivalence of categories: 
\begin{align*} 
Z(?)\col  {\rm Enl}^{\rm sld}(Z/B) \os{\sim}{\lo} {\rm CEnl}^{\rm sld}(Z/B). 
\end{align*} 
Set 
\begin{align*} 
{\rm ChEnl}^{\sq}_{\star}(Z/B):=
{\rm hEnl}^{\sq}_{\star}(Z/B)\cap {\rm CEnl}^{\sq}_{\star}(Z/B). 
\tag{5.1.32.2}\label{cd:rebb} 
\end{align*}  
By \cite[Definition 4]{ollc} we have the following functor 
\begin{align*} 
(?)^{\nat}\col {\rm CEnl}_{\star}(Z/B)\lo {\rm ChEnl}_{\star}(Z/B). 
\tag{5.1.32.3}\label{cd:rahbb} 
\end{align*}
This functor is the left adjoint functor of  
the natural inclusion functor 
${\rm ChEnl}^{\sq}_{\star}(Z/B)\os{\sus}{\lo}{\rm CEnl}^{\sq}_{\star}(Z/B)$. 
By noting that $M_{Z(T)}/{\cal O}_{Z(T)}^*=M_Z/{\cal O}_Z^*$ 
($T\in {\rm Enl}_{\star}(Z/B)$) is also locally constant and by (\ref{lemm:ijk}) (2), 
we have the following composite functor 
\begin{align*} 
Z(?)^{\nat}:=\nat\circ Z(?)\col {\rm Enl}^{\sq}_{\star}(Z/B)\os{(\ref{cd:crabb} )}{\lo} 
{\rm CEnl}^{\rm sld}_{\star}(Z/B)\os{\nat}{\lo}{\rm ChEnl}^{\rm sld}_{\star}(Z/B). 
\tag{5.1.32.4}\label{cd:raebb} 
\end{align*}

\begin{defi} 
(1) Let $N$ be a nonnegative integer or $\infty$.  
Let $Y_{\bul \leq N}\lo Z$ be a morphism over $B$ from 
an $N$-truncated simplicial fine log $p$-adic formal $B$-scheme to 
a fine log $p$-adic formal $B$-scheme.  
Let ${\cal C}_p$ be a subcategory of  ${\rm Enl}_p(Z/B)$.  
For each object $T=(T,z_1)$ of ${\cal C}_p$, assume that 
we are given an $N$-truncated cosimplicial filtered crystal $(F^{\bul \leq N}(T),P(T))$ 
of filtered coherent ${\cal O}_{Y_{\bul \leq N,T_1}/T}$-modules 
such that for a morphism $g\col T'\lo T$ in ${\cal C}_p$, 
there exists an isomorphism 
$\rho(g)\col g^*((F^{\bul \leq N}(T),P(T)))\os{\sim}{\lo} (F^{\bul \leq N}(T'),P(T'))$ 
of filtered ${\cal O}_{{Y_{\bul \leq N,T'_1}/T'}}$-modules  
which satisfies the usual cocycle condition 
and $\rho({\rm id}_T)={\rm id}_{(F^{\bul \leq N}(T),P(T))}$. 
Then we call the family 
$$(F^{\bul \leq N},P):=\{(F^{\bul \leq N}(T),P(T))\}_{T\in {\cal C}_p}
:=\{(F^{\bul \leq N}((T,z_1)),P((T,z_1)))\}_{(T,z_1)\in {\cal C}_p}
$$  
a {\it filtered coherent crystal} of ${\cal O}_{\{Y_{\bul \leq N,T_1}/T\}_{T\in {\cal C}_p}}$-modules. 
We define a {\it morphism} of filtered coherent crystals of 
${\cal O}_{\{Y_{\bul \leq N,T_1}/T\}_{T\in {\cal C}_p}}$-modules in an obvious way. 
\par 
(2) Let the notations be as in (1) and (\ref{defi:nsd}) (3).  
Let $Y'_{\bul \leq N}/Z'$ be a similar object to $Y_{\bul \leq N}/Z$. 
Let $h\col Y'_{\bul \leq N}\lo Y_{\bul \leq N}$ be a morphism over $g\col Z'\lo Z$.  
For an object $(T',z'_1)\in {\cal C}'_p$, let 
$h(T')\col Y'_{\bul \leq N,T'_1}=Y'_{\bul \leq N}\times_{Z',z_1}T'_1
\lo Y_{\bul \leq N}\times_{Z',g\circ z'_1}T'_1$ 
be the base change morphism of $h$ over $T'$. 
Set 
$$h^*((F^{\bul \leq N},P)):=\{h(T')^*_{\rm crys}((F^{\bul \leq N}((T',g\circ z'_1)),
P((T',g\circ z'_1))))\}_{(T',z'_1)\in {\cal C}'_p}. 
$$  
We call $h^*((F^{\bul \leq N},P))$ the {\it pull-back} of 
$(F^{\bul \leq N},P)$ by $h$. 
\par 
(3) Let the notations be as in (1). 
For each object $T=(T,z_1)$ of ${\cal C}_p$, assume that 
we are given an $N$-truncated cosimplicial filtered crystal 
$(E^{\bul \leq N}(\os{\circ}{T}),P(\os{\circ}{T}))$ 
of filtered coherent ${\cal O}_{\os{\circ}{Y}_{\bul \leq N,T_1}/\os{\circ}{T}}$-modules 
such that for a morphism $g\col T'\lo T$ in ${\cal C}_p$, 
there exists an isomorphism 
$\rho(g)\col g^*((E^{\bul \leq N}(\os{\circ}{T}),P(\os{\circ}{T})))
\os{\sim}{\lo} (E^{\bul \leq N}(\os{\circ}{T}{}'),P(\os{\circ}{T}{}'))$ 
of filtered ${\cal O}_{{\os{\circ}{Y}_{\bul \leq N,T'_1}/\os{\circ}{T}{}'}}$-modules  
which satisfies the usual cocycle condition and 
$\rho({\rm id}_T)={\rm id}_{(E^{\bul \leq N}(\os{\circ}{T}),P(\os{\circ}{T}))}$. 
Then we call the family 
$$(E^{\bul \leq N},P):=\{(E^{\bul \leq N}(\os{\circ}{T}),P(\os{\circ}{T}))\}_{T\in {\cal C}_p}
:=\{(E^{\bul \leq N}((\os{\circ}{T},z_1)),P((\os{\circ}{T},z_1)))\}_{(T,z_1)\in {\cal C}_p}
$$  
a {\it filtered coherent crystal} of 
${\cal O}_{\{\os{\circ}{Y}_{\bul \leq N,T_1}/\os{\circ}{T}\}_{T\in {\cal C}_p}}$-modules. 
We define a {\it morphism} of filtered coherent crystals of 
${\cal O}_{\{\os{\circ}{Y}_{\bul \leq N,T_1}/\os{\circ}{T}\}_{T\in {\cal C}_p}}$-modules 
in an obvious way. 
\par 
If we are given $(E^{\bul \leq N},P)$ in (3), 
we obtain a coherent crystal of ${\cal O}_{\{Y_{\bul \leq N,T_1}/T\}_{T\in {\cal C}_p}}$-modules 
in the following way:
$$(F^{\bul \leq N},P):=
\{(\eps^*_{Y_{\bul \leq N,T_1}/T}(E^{\bul \leq N}(\os{\circ}{T})),P(\os{\circ}{T}))\}_{T\in {\cal C}_p}.$$
\end{defi}

\begin{rema}\label{rema:fb} 
Assume that $Z$ has a $p$-adically formally log smooth lift ${\cal Z}$ over $B$ 
and that ${\cal Z}\in {\cal C}_p$.  
Then it is obvious that we have the following functor: 
\begin{align*} 
\{{\rm filtered~coherent}~{\cal O}_{\{Y_{\bul \leq N,T_1}/T\}_{T\in {\cal C}_p}}{\textrm -}{\rm modules}\}
\lo  \{{\rm filtered~coherent}~{\cal O}_{Y_{\bul \leq N}/{\cal Z}}{\textrm -}{\rm modules}\}.  
\end{align*} 
\end{rema}

Lastly let us recall Ogus' definition of an $F^{\infty}$-span and his theorem in \cite{ollc}.   

\begin{defi}[{\bf \cite[Definition 15]{ollc}}]\label{defi:fspan}
Assume that $\os{\circ}{Z}$ is of characteristic $p>0$. 
An $F^{\infty}$-span on $Z/{\cal W}$ is 
a sequence  $(E,\Phi,V)=(\{E_n\}_{n=0}^{\inf},\{\Phi_n\}_{n=0}^{\inf},\{V_n\}_{n=0}^{\inf})$ 
of triples, where $E_n$ is a crystal of ${\cal O}_{Z/{\cal W}}$-modules and  
$\Phi_n\col F^*_{Z/{\cal W}}(E_{n+1})\lo E_n$ and 
$V_n\col E_n\lo F^*_{Z/{\cal W}}(E_{n+1})$ are morphisms of ${\cal O}_{Z/{\cal W}}$-modules 
such that $\Phi_n\circ V_n$ and $V_n\circ \Phi_n$ are multiplications by 
$p^m$ for a fixed some $m\in {\mab Z}_{\geq 0}$ (independent of $n$). 
\par 
We say that $E$ is flat (resp.~coherent)  
if $E_n$ $(n\in {\mab Z}_{\geq 0})$ is a flat ${\cal O}_{Z/{\cal W}}$-modules 
(resp.~a coherent ${\cal O}_{Z/{\cal W}}$-modules). 
\end{defi}

\par 
The following Ogus' theorem is a log version and a coefficient version of his theorem 
\cite[(3.7)]{of}: 

\begin{theo}[{\bf \cite[Theorem 4]{ollc}}]\label{theo:cvnl}  
Let $Z$ be an fs log scheme over $\kap$ with constant log structure, i.~e., 
$M_Z/{\cal O}_Z^*$ is locally constant. 
Let $Y\lo Z$ be a proper log smooth morphism 
of fs log schemes over $\kap$ of Cartier type. 
Let $E=\{E_n\}_{n=1}^{\infty}$ be an $F^{\infty}$-span on $Y/{\rm Spf}({\cal W})$. 
Then there exists a convergent $F^{\infty}$-isospan $E^q:=\{E^q_n\}_{n=0}^{\infty}$ 
of flat ${\cal K}_{Z/{\cal W}}$-modules 
such that $E^q_{n,T}=R^qf_{Y_{T_1}/T*}(p_{T{\rm crys}}^*(E_n))_{\mab Q}$ 
for any object $T\in {\rm Enl}_p(Z/{\cal W})$. 
Here $p_T\col Y_{T_1}\lo Y$ is the first projection over $T\lo {\rm Spf}({\cal W})$. 
\end{theo}

In this book we need only the following in the trivial log case. 

\begin{theo}\label{theo:otf}
Endow ${\cal V}$ with the canonical PD-structure $[~]$ on $p{\cal V}$. 
Let $Z$ be a $p$-adic formal ${\cal V}$-scheme and 
let $f\col Y\lo Z_1$ be a proper smooth morphism. 
Set ${\cal E}_p:={\rm Enl}_p(Z/{\cal V})$ and ${\cal E}:={\rm Enl}(Z/{\cal V})$.  
Let $E:=\{E(T)\}_{T\in {\cal E}_p}$ be a 
flat coherent crystals of 
$\{{\cal O}_{Y_{T_1}/T}\}_{T\in {\cal E}_p}$-modules. 
Then the following hold$:$
\par 
$(1)$ $($A coefficient version of {\rm \cite[(3.1)]{of}}$)$ 
There exists an object $R^qf_*(E_{K})$ of ${\rm Isoc}_p(Z/{\cal V})$ 
such that 
$R^qf_*(E_{K})_{T}=
R^qf_{Y_{T_1}/T*}(E(T))\otimes_{\cal V}K$ for any object $T$ of ${\cal E}_p$. 
\par 
$(2)$ $($A coefficient version of {\rm \cite[(3.7)]{of}}$)$ 
Let us consider $R^qf_*(E_K)$ as an object of ${\rm Isoc}_p(Z/{\cal W})$ 
by using the natural morphism ${\rm Spf}({\cal V})\lo {\rm Spf}({\cal W})$ 
as in {\rm \cite[(3.6)]{of}}. 
For an object $T\in {\rm Enl}_p(Z/{\cal W})$, 
set $Y'_{T_1}:=(Y\times_{Z_1,F_{Z_1}}Z_1)\times_{Z_1}T_1$  
and 
let $W_{Y_{T_1}/T}\col Y'_{T_1}\lo Y_{T_1}$ be the natural projection over $T$. 
Let $F_{Y/Z_1}\col Y\lo Y'$ be the relative Frobenius morphism over $Z$ and 
let $F_{Y_{T_1}/T}\col Y_{T_1}\lo Y'_{T_1}$ be the base change morphism of 
$F_{Y/Z_1}$ over $T$.  
Let $\Phi:=\{\Phi(T)\}_{T\in {\cal E}_p}$ 
be a morphism of $\{{\cal O}_{Y_{T_1}/T}\}_{T\in {\cal E}_p}$-modules, where 
$$\Phi(T) \col F_{Y_{T_1}/T{\rm crys}}^*W_{Y_{T_1}/T{\rm crys}}^*(E(T))\lo E(T)$$  
is a morphism of $F$-crystals of ${\cal O}_{Y_{T_1}/T}$-modules. 
Assume that the induced morphism 
\begin{align*} 
R^qf_{Y'_{T_1}/T*}(W^*_{Y_{T_1}/T{\rm crys}}(E(T)))\otimes_{\cal V}K\lo  
R^qf_{Y_{T_1}/T*}(E(T))\otimes_{\cal V}K 
\tag{5.1.37.1}\label{ali:wep}
\end{align*}
by $\Phi(T)$ is an isomorphism for any object $T$ of ${\cal E}_p$. 
Then $R^qf_*(E_K)$ in $(1)$ prolongs naturally to a convergent $F$-isocrystal on 
${\cal E}$ as explained in {\rm \cite[(3.6), (3.7)]{of}}. 
\end{theo} 
\begin{proof} 
(1): We have only to imitate the proof of \cite[(3.1)]{of}. 
Indeed, let $g\col T'\lo T$ be a morphism in ${\rm Enl}_p(Z/{\cal V})$. 
We may assume that $Z$, $T$ and $T'$ are affine: 
$T':={\rm Spf}(C')$ and $T:={\rm Spf}(C)$. 
Since $g^*_{\rm crys}(E(T))=E(T')$, we have 
\begin{align*} 
C'\otimes^L_{C}R\Gam(X_{T_1}/T,E(T))=R\Gam(X_{T'_1}/T',E(T'))
\end{align*} 
by the base change of crystalline cohomologies ([loc.~cit., (3.3)]). 
The rest of the proof is the same as that of [loc.~cit., (3.1)]. 
\par 
(2): We have only to imitate the proof of \cite[(3.7)]{of}. 
Indeed, as in [loc.~cit.], we may assume that ${\cal V}={\cal W}$. 
Let $T$ be an object of ${\cal E}_p$. 
Then 
\begin{align*} 
(F_{Z_1}^*R^qf_*(E_K))_T=
R^qf_{Y'_{T_1}/T*}(W_{Y_{T_1}/T{\rm crys}}^*(E(T))\otimes_{\cal V}K. 
\tag{5.1.37.2}\label{ali:fzc}
\end{align*}
By the assumption (\ref{ali:wep}), this is isomorphic to 
$R^qf_{Y_{T_1}/T*}(E(T))\otimes_{\cal V}K$.   
Hence the Frobenius morphism  
\begin{align*} 
(F_{Z_1}^*R^qf_*(E_K))_T\lo R^qf_*(E_K)_T
\tag{5.1.37.3}\label{ali:frzc}
\end{align*}
is an isomorphism. 
The rest of the proof is the same as that of \cite[(3.7)]{of} by using 
\cite[(2.18)]{of}.
\end{proof}

\section{Log convergences of the weight filtrations on the log isocrystalline cohomology sheaves}\label{sec:lcw} 
In the section we prove the log convergence of the weight filtration 
on the log isocrystalline cohomology sheaf of 
a truncated simplicial base change of proper SNCL schemes in characteristic $p>0$ 
by using results in previous sections. 
We also reformulate the monodromy-weight conjecture and the log hard Lefschetz conjecture 
in terms of log convergent isocrystals. 
\par 
Let the notations be as in the previous section. 
In this section, assume that $B=({\rm Spf}({\cal V}),{\cal V}^*)$.  
\par 
Let $N$ be a nonnegative integer.
Let $S$ be a $p$-adic formal family of log points over $B$ 
such that $\os{\circ}{S}$ is a ${\cal V}/p$-scheme. 
Let $X_{\bul \leq N}/S$ be an $N$-truncated simplicial base change of SNCL schemes. 
Assume that $\os{\circ}{X}_{\bul \leq N}/\os{\circ}{S}$ is proper.

\par 
Let $n$ be a nonnegative integer. 
Let $T$ be an object of ${\rm Enl}^{\sq}_{\star}(S^{[p^n]}/{\cal V})$ 
$(\star=p$ or nothing, $\sq=$sld or nothing).  
Then the hollow out $S^{[p^n]}(T)^{\nat}$ of $S^{[p^n]}(T)$ is a formal family of log points. 
Let $T_i \lo S^{[p^n]}$ $(i=0,1)$ be the structural morphism. 
Set $X^{[p^n]}_{\bul \leq N,\os{\circ}{T}_i}:=X_{\bul \leq N}\times_SS^{[p^n]}_{\os{\circ}{T}_i} 
=X^{[p^n]}_{\bul \leq N}\times_{\os{\circ}{S}}\os{\circ}{T}_i$, 
where $X^{[p^n]}_{\bul \leq N}:=X_{\bul \leq N}\times_SS^{[p^n]}$. 
Let 
$f^{[p^n]}_{\os{\circ}{T}} \col X^{[p^n]}_{\bul \leq N,\os{\circ}{T}_i}\lo S^{[p^n]}(T)^{\nat}$ 
be the structural morphism. 
(Note that this notation is different from 
$f$ in \S\ref{sec:psc} since we add the symbol $\os{\circ}{T}$ to $f^{[p^n]}$.) 
\par 
Set ${\cal E}^{\sq}_{p,n}:={\rm Enl}^{\sq}_p(S^{[p^n]}/{\cal V})$. 
Assume that we are given a flat coherent crystal 
$E^{\bul \leq N}_n=\{E^{\bul \leq N}_n(\os{\circ}{T})\}_{T\in {\cal E}^{\sq}_{p,n}}$ of 
${\cal O}_{\{\os{\circ}{X}{}^{[p^n]}_{\bul \leq N,T_1}
/\os{\circ}{T}\}_{T\in {\cal E}^{\sq}_{p,n}}}$-modules.  
Let $T$ be an object of ${\cal E}^{\sq}_{p,n}$.
Let $z\col T_1\lo S^{[p^n]}$ be the structural morphism. 
Because $S^{[p^n]}(T)^{\nat}$  is a $p$-adic formal family of log points,  
we have an $N$-truncated simplicial base change of proper SNCL schemes  
$X^{[p^n]}_{\bul \leq N,\os{\circ}{T}_1}/S^{[p^n]}_{\os{\circ}{T}_1}$ 
with an exact PD-closed immersion 
$S^{[p^n]}_{\os{\circ}{T}_1} \os{\sus}{\lo} S^{[p^n]}(T)^{\nat}$. 
Assume that, for each object $T$ of ${\cal E}^{\sq}_{p,n}$, 
there exists a disjoint union 
$X^{[p^n]}{}'_{\! \! \! \! \! \! \! \! \! \! \bul \leq N,\os{\circ}{T}_1}$ of 
the member of an affine $N$-truncated simplicial open covering of 
$X^{[p^n]}_{\bul \leq N,\os{\circ}{T}_1}$ and 
that there exists an admissible immersion 
\begin{align*} 
X^{[p^n]}_{\bul \leq N,\bul,\os{\circ}{T}_1}
\os{\sus}{\lo} \ol{\cal P}_{\bul \leq N,\bul}(\ol{S^{[p^n]}(T)^{\nat}})
\tag{5.2.0.1}\label{eqn:audtds}
\end{align*} 
over $\ol{S^{[p^n]}(T)^{\nat}}$. 
\parno 
Then we obtain the $p$-adic iso-zariskian filtered Steenbrink complex
\begin{align*} 
(A_{{\rm zar},{\mab Q}}(X^{[p^n]}_{\bul \leq N,\os{\circ}{T}_1}/
S^{[p^n]}(T)^{\nat},E^{\bul \leq N}(\os{\circ}{T})),P)
\in {\rm D}^+{\rm F}(
f^{-1}_{\os{\circ}{T}}({\cal K}_T))
\tag{5.2.0.2}\label{eqn:auxtds}
\end{align*} 
for each $T\in {\cal E}^{\sq}_{p,n}$. 
The admissible immersion (\ref{eqn:audtds}) exists in the following two cases: 
\par 
(1) 
The case where 
$X_{\bul \leq N}$ is an $N$-truncated simplicial SNCL scheme over $S$ 
(\S\ref{sec:psc}). 
\par 
(2) The case of the successive truncated simplicial case ((\ref{theo:thenad})).

\begin{rema}\label{rema:tep} 
Let $T$ be an object of ${\cal E}^{\sq}_{p,n}$.  
If there exists an admissible immersion 
\begin{align*} 
X^{[p^n]}_{\bul \leq N,\bul,\os{\circ}{T}_1}\os{\sus}{\lo} 
\ol{\cal P}_{\bul \leq N,\bul}(\ol{S^{[p^n]}(T)^{\nat}})
\tag{5.2.1.1}\label{eqn:autdtds}
\end{align*} 
over $\ol{S^{[p^n]}(T)^{\nat}}$, then 
there exists an admissible immersion 
$$X^{[p^n]}_{\bul \leq N,\bul, \os{\circ}{T}_1}\os{\sus}{\lo} 
\ol{\cal P}_{\bul \leq N,\bul}(\ol{S^{[p^n]}(T)^{\nat}})
\times_{\ol{S^{[p^n]}(T)^{\nat}}}\ol{S^{[p^n]}(T')^{\nat}}$$  
for any morphism $T'\lo T$ in ${\cal E}^{\sq}_{p,n}$. 
\end{rema}

\begin{prop}\label{prop:tptt} 
Let ${\mathfrak g}\col T'\lo T$ be a morphism in ${\cal E}^{\sq}_{p,n}$.  
Then the induced morphism 
$g_{\bul \leq N} \col X^{[p^n]}_{\bul \leq N,\os{\circ}{T}{}'_1}\lo X^{[p^n]}_{\bul \leq N,\os{\circ}{T}_1}$ 
by ${\mathfrak g}$
gives us the following natural morphism 
\begin{align*} 
g^*_{\bul \leq N}\col &
(A_{{\rm zar},{\mab Q}}(X^{[p^n]}_{\bul \leq N,\os{\circ}{T}_1}/
S^{[p^n]}(T)^{\nat},E^{\bul \leq N}(\os{\circ}{T})),P)
\tag{5.2.2.1}\label{eqn:auqeds}\\
& \lo
Rg_{\bul \leq N*}
((A_{{\rm zar},{\mab Q}}(X^{[p^n]}_{\bul \leq N,\os{\circ}{T}{}'_1}/S^{[p^n]}(T')^{\nat},
E^{\bul \leq N}(\os{\circ}{T}{}')),P))
\end{align*}
of filtered complexes in 
${\rm D}^+{\rm F}(
f_{\os{\circ}{T}}^{-1}({\cal K}_T))$ 
fitting into the following commutative diagram$:$ 
\begin{equation*} 
\begin{CD} 
A_{{\rm zar},{\mab Q}}(X^{[p^n]}_{\bul \leq N,\os{\circ}{T}_1}/
S^{[p^n]}(T)^{\nat},E^{\bul \leq N}(\os{\circ}{T}))
@>{g^*_{\bul \leq N}}>> \\
@A{\simeq}AA \\
Ru_{X^{[p^n]}_{\bul \leq N,\os{\circ}{T}_1}/
S^{[p^n]}(T)^{\nat}*}(\eps_{X^{[p^n]}_{\bul \leq N,\os{\circ}{T}_1}/
S^{[p^n]}(T)^{\nat}}(E^{\bul \leq N}(\os{\circ}{T})))\otimes_{\mab Z}{\mab Q}
@>{g^*_{\bul \leq N}}>>
\end{CD} 
\end{equation*} 
\begin{equation*} 
\begin{CD} 
Rg_{\bul \leq N*}
((A_{{\rm zar},{\mab Q}}(X^{[p^n]}_{\bul \leq N,\os{\circ}{T}{}'_1}/S^{[p^n]}(T')^{\nat},
E^{\bul \leq N}(\os{\circ}{T}{}')),P))\\ 
@AA{\simeq}A \\
Rg_{\bul \leq N*}(Ru_{X^{[p^n]}_{\bul \leq N,\os{\circ}{T}{}'_1}/
S^{[p^n]}(T')^{\nat}*}(\eps^*_{X^{[p^n]}_{\bul \leq N,\os{\circ}{T}{}'_1}/
S^{[p^n]}(T')^{\nat}}(E^{\bul \leq N}(\os{\circ}{T}{}'))\otimes_{\mab Z}{\mab Q}).
\end{CD} 
\end{equation*} 
For a similar morphism ${\mathfrak h}\col T''\lo T'$ to ${\mathfrak g}$ and 
a similar morphism 
$h_{\bul \leq N}\col 
X_{\bul \leq N,\os{\circ}{T}{}''_1}\lo X_{\bul \leq N,\os{\circ}{T}{}'_1}$ to $g_{\bul \leq N}$, 
the following relation 
\begin{align*} 
(h_{\bul \leq N}\circ g_{\bul \leq N})^*=
Rh_{\bul \leq N*}(g^*_{\bul \leq N})
\circ h^*_{\bul \leq N}
\end{align*} 
holds.  
\begin{equation*} 
{\rm id}_{X^{[p^n]}_{\bul \leq N,\os{\circ}{T}_1}}^*={\rm id} 
_{(A_{{\rm zar},{\mab Q}}(X^{[p^n]}_{\bul \leq N,\os{\circ}{T}_0}/
S(T)^{\nat},E^{\bul \leq N}(\os{\circ}{T}),P)}.  
\end{equation*} 
\end{prop} 
\begin{proof} 
This immediately follows from 
the contravariant functoriality in (\ref{theo:fugennas}). 
\end{proof} 

\parno 
The morphism (\ref{eqn:auqeds}) induces the following morphism 
\begin{align*} 
{\mathfrak g}^*\col &
R
f^{[p^n]}_{\os{\circ}{T}*}((A_{{\rm zar},{\mab Q}}(X^{[p^n]}_{\bul \leq N,\os{\circ}{T}_1}
/S^{[p^n]}(T)^{\nat},E^{\bul \leq N}(\os{\circ}{T})),P))
\tag{5.2.2.2}\label{eqn:auaeds}\\
&\lo
R{\mathfrak g}_*R
f^{[p^n]}_{\os{\circ}{T}{}'*}
((A_{{\rm zar},{\mab Q}}(X^{[p^n]}_{\bul \leq N,\os{\circ}{T}{}'_1}/S^{[p^n]}(T')^{\nat},
E^{\bul \leq N}(\os{\circ}{T}{}')),P))
\end{align*}
of filtered complexes in ${\rm D}^+{\rm F}({\cal K}_T)$.

\par  
The following is a key lemma for (\ref{theo:pwfec}) below.

\begin{lemm}\label{lemm:pnlcfi}
Assume that $M_S$ is split. Let $k$ be a nonnegative integer or $\infty$. 
Then there exists an object
\begin{align*} 
(R^q
f^{[p^n]}_*
(P_kA_{{\rm zar},{\mab Q}}(X^{[p^n]}_{\bul \leq N}/K,E^{\bul \leq N})),P)
\tag{5.2.3.1}\label{ali:pkalkf} 
\end{align*} 
of ${\rm IsocF}_p(\os{\circ}{S}/{\cal V})$ 
such that 
\begin{align*} 
&(R^qf^{[p^n]}_*
(P_kA_{{\rm zar},{\mab Q}}(X^{[p^n]}_{\bul \leq N}/K,E^{\bul \leq N})),P)_{\os{\circ}{T}}=
\tag{5.2.3.2}\label{ali:pkakf} \\
&(R^qf^{[p^n]}_{\os{\circ}{T}*}
(P_kA_{{\rm zar},{\mab Q}}(X^{[p^n]}_{\bul \leq N,\os{\circ}{T}_1}/S^{[p^n]}({\os{\circ}{T})},
E^{\bul \leq N}(\os{\circ}{T}))),P)
\end{align*} 
for any object $\os{\circ}{T}$ in ${\rm Enl}_p(\os{\circ}{S}/{\cal V})$. 
$($Note that $(S^{[p^n]})^{\circ}=\os{\circ}{S}$ and 
recall the log formal scheme $S^{[p^n]}({\os{\circ}{T})}$ 
in {\rm (\ref{prop:nupf})}.$)$
In particular, there exists an object
\begin{align*} 
(R^q
f^{[p^n]}_{*}
(\eps^*_{X^{[p^n]}_{\bul \leq N}/K}(E^{\bul \leq N}_K)),P)
\tag{5.2.3.3}\label{ali:pklkf} 
\end{align*} 
of ${\rm IsocF}_p(\os{\circ}{S}/{\cal V})$ 
such that 
\begin{align*} 
(R^q
f^{[p^n]}_{*}
(\eps^*_{X^{[p^n]}_{\bul \leq N}/K}(E^{\bul \leq N}_K)),P)_{\os{\circ}{T}}=
(R^qf^{[p^n]}_{X^{[p^n]}_{\bul \leq N,\os{\circ}{T}_1/S^{[p^n]}(\os{\circ}{T})*}}
(\eps^*_{X_{\bul \leq N,\os{\circ}{T}_1/S^{[p^n]}(\os{\circ}{T})}}
(E^{\bul \leq N}(\os{\circ}{T})))_{\mab Q},P)
\tag{5.2.3.4}\label{ali:pkkf} 
\end{align*} 
for any object $\os{\circ}{T}$ in ${\rm Enl}_p(\os{\circ}{S}/{\cal V})$. 
\end{lemm} 
\begin{proof} 
Set ${\cal E}_p:={\rm Enl}_p(\os{\circ}{S}/{\cal V})$ and 
let $\os{\circ}{T}$ be an object of ${\cal E}_p$. 
Set 
\begin{align*} 
(E,P)
:=\{(E(\os{\circ}{T}),P(\os{\circ}{T}))\}_{\os{\circ}{T}\in {\cal E}_p}:=
(R^q
f^{[p^n]}_{\os{\circ}{T}*}
(P_kA_{{\rm zar},{\mab Q}}
(X^{[p^n]}_{\bul \leq N,\os{\circ}{T}_1}/S^{[p^n]}(\os{\circ}{T}),
E^{\bul \leq N}(\os{\circ}{T})),P).
\end{align*}  
We would like to prove that 
$\{(E(\os{\circ}{T}),P(\os{\circ}{T}))\}_{\os{\circ}{T}\in {\cal E}_p}$ 
defines an object of ${\rm IsocF}_p(\os{\circ}{S}/{\cal V})$. 
Let $g\col \os{\circ}{T}{}'\lo \os{\circ}{T}$ be 
a morphism in ${\cal E}_p$. 
We may assume that $\os{\circ}{T}$ and $\os{\circ}{T}{}'$ are affine formal schemes. 
Let $z\col \os{\circ}{T}_1\lo \os{\circ}{S}=(S^{[p^n]})^{\circ}$ be the structural morphism. 
\par 
Because we need a more general situation in (\ref{theo:pwfec}) below, 
we consider a morphism $g\col T'\lo T$ of log formal schemes 
in ${\rm Enl}^{\sq}_p(S^{[p^n]}/{\cal V})$ more generally. 
Since $T',T\in {\rm Enl}^{\sq}_p(S^{[p^n]}/{\cal V})$, 
we obtain the $p$-adic iso-zariskian filtered  Steenbrink complexes
$$(A_{{\rm zar},{\mab Q}}
(X^{[p^n]}_{\bul \leq N,\os{\circ}{T}_1}/S^{[p^n]}(T)^{\nat},E^{\bul \leq N}(\os{\circ}{T})),P)$$
and
$$(A_{{\rm zar},{\mab Q}}
(X^{[p^n]}_{\bul \leq N,\os{\circ}{T}{}'_1}/S^{[p^n]}(T')^{\nat},E^{\bul \leq N}(\os{\circ}{T}{}')),P).$$  
Let $h_{\bul \leq N}\col X_{\bul \leq N,\os{\circ}{T}{}'_1}\lo X_{\bul \leq N,\os{\circ}{T}_1}$ 
be the natural morphism over $g\col S_{\os{\circ}{T}{}'}\lo S_{\os{\circ}{T}}$. 
By the contravariant functoriality ((\ref{theo:fugennas})), 
we have the following natural morphism 
\begin{align*} 
&(A_{{\rm zar},{\mab Q}}(X^{[p^n]}_{\bul \leq N,\os{\circ}{T}_1}/
S^{[p^n]}(T)^{\nat},E^{\bul \leq N}(\os{\circ}{T})),P)
\lo \\
&Rh_{\bul \leq N*}((A_{{\rm zar},{\mab Q}}
(X^{[p^n]}_{\bul \leq N,\os{\circ}{T}{}'_1}/
S^{[p^n]}(T')^{\nat},E^{\bul \leq N}(\os{\circ}{T}{}')),P)).  
\end{align*} 
Applying this morphism to 
$Rf^{[p^n]}_{\os{\circ}{T}*}$, 
we obtain the following morphism 
\begin{align*} 
&Rf^{[p^n]}_{\os{\circ}{T}*}
(A_{{\rm zar},{\mab Q}}(X^{[p^n]}_{\bul \leq N,\os{\circ}{T}_1}/
S^{[p^n]}(T)^{\nat},E^{\bul \leq N}(\os{\circ}{T})),P)
\\
&\lo Rf^{[p^n]}_{\os{\circ}{T}*}
Rh_{\bul \leq N*}((A_{{\rm zar},{\mab Q}}
(X^{[p^n]}_{\bul \leq N,\os{\circ}{T}{}'_1}/
S^{[p^n]}(T')^{\nat},E^{\bul \leq N}(\os{\circ}{T}{}')),P))\\
&=Rg_*Rf^{[p^n]}_{\os{\circ}{T}{}'*}((A_{{\rm zar},{\mab Q}}
(X^{[p^n]}_{\bul \leq N,\os{\circ}{T}{}'_1}/
S^{[p^n]}(T')^{\nat},E^{\bul \leq N}(\os{\circ}{T}{}')),P)).  
\end{align*} 
By (\ref{prop:bbtdc}), 
$Rf^{[p^n]}_{\os{\circ}{T}*}
(A_{{\rm zar},{\mab Q}}(X^{[p^n]}_{\bul \leq N,\os{\circ}{T}_1}/
S^{[p^n]}(T)^{\nat},E^{\bul \leq N}(\os{\circ}{T})),P)$ is filtered bounded above.
By the adjunction of $Rg_*$ and $Lg^*$ (\cite[(1.2.2)]{nh2}), 
we have the following morphism 
\begin{align*} 
&Lg^*Rf^{[p^n]}_{\os{\circ}{T}*}
(A_{{\rm zar},{\mab Q}}(X^{[p^n]}_{\bul \leq N,\os{\circ}{T}_1}/
S^{[p^n]}(T)^{\nat},E^{\bul \leq N}(\os{\circ}{T})),P)
\tag{5.2.3.5}\label{ali:gfeij}\\
&\lo Rf^{[p^n]}_{\os{\circ}{T}{}'*}((A_{{\rm zar},{\mab Q}}
(X^{[p^n]}_{\bul \leq N,\os{\circ}{T}{}'_1}/
S^{[p^n]}(T')^{\nat},E^{\bul \leq N}(\os{\circ}{T}{}')),P)).  
\end{align*} 
We claim that the following equality holds:
\begin{align*} 
&{\cal H}^q(Lg^*Rf^{[p^n]}_{\os{\circ}{T}*}
(P_kA_{{\rm zar},{\mab Q}}(X^{[p^n]}_{\bul \leq N,\os{\circ}{T}_1}/
S^{[p^n]}(T)^{\nat},E^{\bul \leq N}(\os{\circ}{T}))))\tag{5.2.3.6}\label{ali:gfeiaj}\\
&=g^*R^qf^{[p^n]}_{\os{\circ}{T}*}
(P_kA_{{\rm zar},{\mab Q}}(X^{[p^n]}_{\bul \leq N,\os{\circ}{T}_1}/
S^{[p^n]}(T)^{\nat},E^{\bul \leq N}(\os{\circ}{T}))) \quad (k\in {\mab Z}).
\end{align*} 
Indeed, by (\ref{theo:bcfqa}), 
$Rf^{[p^n]}_{\os{\circ}{T}*}
(P_kA_{{\rm zar},{\mab Q}}(X^{[p^n]}_{\bul \leq N,\os{\circ}{T}_1}/
S^{[p^n]}(T)^{\nat},E^{\bul \leq N}(\os{\circ}{T}))))$ is 
a strictly perfect complex of ${\cal K}_T$-modules. 
Let $L^{\bul}:=(L^{\bul},d^{\bul})$ be a representative of this complex 
such that $L^{\bul}$ is bounded and such that 
each $L^i$ $(i\in {\mab Z})$ is a locally free ${\cal K}_T$-module. 
Then the source of (\ref{ali:gfeij}) is represented by $g^*(L^{\bul})$. 
Set ${\cal Z}^q:={\rm Ker}(d^q)$, ${\cal B}^q:={\rm Im}(d^q)$ 
and ${\cal H}^q:={\cal Z}^q/{\cal B}^q$ as usual. Then, 
by (\ref{eqn:ekpssp}), (\ref{theo:otf}) (1) and \cite[(2.9)]{of}, 
$${\cal H}^q:=R^qf^{[p^n]}_{\os{\circ}{T}*}
(P_kA_{{\rm zar},{\mab Q}}(X^{[p^n]}_{\bul \leq N,\os{\circ}{T}_1}/
S^{[p^n]}(T)^{\nat},E^{\bul \leq N}(\os{\circ}{T}))))$$ 
is a flat ${\cal K}_T$-module. 
By using this flatness, 
it is easy to check that 
$g^*({\cal H}^q)={\cal H}^q(g^*(L^{\bul}))$. 
This is nothing but the equality (\ref{ali:gfeiaj}). 
\par 
By (\ref{ali:gfeij}) and (\ref{ali:gfeiaj}), we have the following morphism 
\begin{align*}
&g^*R^qf^{[p^n]}_{\os{\circ}{T}*}
(P_kA_{{\rm zar},{\mab Q}}(X^{[p^n]}_{\bul \leq N,\os{\circ}{T}_1}/
S^{[p^n]}(T)^{\nat},E^{\bul \leq N}(\os{\circ}{T}))))
\tag{5.2.3.7}\label{ali:gfegsiaj}\\
&\lo R^qf^{[p^n]}_{\os{\circ}{T}{}'*}((A_{{\rm zar},{\mab Q}}
(X^{[p^n]}_{\bul \leq N,\os{\circ}{T}{}'_1}/
S^{[p^n]}(T')^{\nat},E^{\bul \leq N}(\os{\circ}{T}{}'))).  
\end{align*} 
Let ${\rm SS}_k(S^{[p^n]}(T)^{\nat})$ and 
${\rm SS}_k(S^{[p^n]}(T')^{\nat})$ be the  spectral sequences 
of 
$$R^qf^{[p^n]}_{\os{\circ}{T}*}
(P_kA_{{\rm zar},{\mab Q}}(X^{[p^n]}_{\bul \leq N,\os{\circ}{T}_1}
/S^{[p^n]}(T)^{\nat},E^{\bul \leq N}(\os{\circ}{T})))$$ 
and 
$$R^qf^{[p^n]}_{\os{\circ}{T}{}'*}
(P_kA_{{\rm zar},{\mab Q}}(X^{[p^n]}_{\bul \leq N,T'_1}/
S^{[p^n]}(T')^{\nat},E^{\bul \leq N}(\os{\circ}{T}{}'))),$$  
respectively ((\ref{eqn:esasp})). 
Let $E_1^{\bul \bul}(S^{[p^n]}(T)^{\nat})_k$ and 
$E_1^{\bul \bul}(S^{[p^n]}(T')^{\nat})_k$ be the $E_1$-terms of the spectral sequences of  
${\rm SS}_k(S^{[p^n]}(T)^{\nat})$ and ${\rm SS}_k(S^{[p^n]}(T')^{\nat})$, respectively. 
Set 
$E_1^{\bul \bul}(S^{[p^n]}(T)^{\nat}):=E_1^{\bul \bul}(S^{[p^n]}(T)^{\nat})_{\infty}$
and $E_1^{\bul \bul}(S^{[p^n]}(T')^{\nat}):=E_1^{\bul \bul}(S^{[p^n]}(T')^{\nat})_{\infty}$. 
\par 
Now consider the morphism (\ref{ali:gfegsiaj}) for 
the case $S^{[p^n]}(T)=S^{[p^n]}(\os{\circ}{T})$ 
and $S^{[p^n]}(T')=S^{[p^n]}(\os{\circ}{T}{}')$:  
\begin{align*}
&{\cal H}^q(Lg^*Rf^{[p^n]}_{\os{\circ}{T}*}
(P_kA_{{\rm zar},{\mab Q}}(X^{[p^n]}_{\bul \leq N,\os{\circ}{T}_1}/
S^{[p^n]}(\os{\circ}{T})^{\nat},E^{\bul \leq N}(\os{\circ}{T}))))
\tag{5.2.3.8}\label{ali:gfgriaj}\\
&\lo R^qf^{[p^n]}_{\os{\circ}{T}{}'*}(P_kA_{{\rm zar},{\mab Q}}
(X^{[p^n]}_{\bul \leq N,\os{\circ}{T}{}'_1}/
S^{[p^n]}(\os{\circ}{T}{}')^{\nat},E^{\bul \leq N}(\os{\circ}{T}{}'))).   
\end{align*} 
By (\ref{theo:otf}) (1) we see that the natural morphism 
$g^*E_1^{\bul \bul}(S^{[p^n]}(\os{\circ}{T}))\lo 
E_1^{\bul \bul}(S^{[p^n]}(\os{\circ}{T}{}'))$ 
is an isomorphism. 
Moreover $E_1^{ij}(S^{[p^n]}(\os{\circ}{T}))$ $(i,j\in {\mab Z})$ 
is a flat coherent ${\cal K}_T$-module. 
Hence the natural morphism (\ref{ali:gfgriaj}) is equal to the following isomorphism 
\begin{align*} 
& g^*R^qf^{[p^n]}_{\os{\circ}{T}*}
(P_kA_{{\rm zar},{\mab Q}}(X^{[p^n]}_{\bul \leq N,\os{\circ}{T}_1}
/S^{[p^n]}(\os{\circ}{T}),E^{\bul \leq N}(\os{\circ}{T})))) 
\os{\sim}{\lo} \tag{5.2.3.9}\label{ali:gpapfx}\\
& R^qf^{[p^n]}_{\os{\circ}{T}{}'*}
(P_kA_{{\rm zar},{\mab Q}}
(X^{[p^n]}_{\bul \leq N,\os{\circ}{T}{}'_1}
/S^{[p^n]}(\os{\circ}{T}{}'),E^{\bul \leq N}(\os{\circ}{T}{}')) 
\end{align*} 
is an isomorphism. 
Because 
\begin{align*} 
& P_{q+k'}R^qf^{[p^n]}_{\os{\circ}{T}*}
(P_kA_{{\rm zar},{\mab Q}}(X^{[p^n]}_{\bul \leq N,\os{\circ}{T}_1}/S^{[p^n]}(\os{\circ}{T}), 
E^{\bul \leq N}(\os{\circ}{T}))) \\
&=
{\rm Im}(R^q
f^{[p^n]}_{\os{\circ}{T}*}
(P_{k'}A_{{\rm zar},{\mab Q}}(X^{[p^n]}_{\bul \leq N,\os{\circ}{T}_1}
/S^{[p^n]}(\os{\circ}{T}),E^{\bul \leq N}(\os{\circ}{T})))\\
&\quad \quad \quad 
\lo R^q
f^{[p^n]}_{\os{\circ}{T}*}
(P_kA_{{\rm zar},{\mab Q}}
(X^{[p^n]}_{\bul \leq N,\os{\circ}{T}_1}
/S^{[p^n]}(\os{\circ}{T}),E^{\bul \leq N}(\os{\circ}{T})))) \quad (k'\leq k), 
\end{align*} 
the following natural morphism 
\begin{align*} 
& g^*(P_{q+k'}R^q
f^{[p^n]}_{\os{\circ}{T}*}
(P_kA_{{\rm zar},{\mab Q}}
(X^{[p^n]}_{\bul \leq N,\os{\circ}{T}_1}/S^{[p^n]}(\os{\circ}{T}), 
E^{\bul \leq N}(\os{\circ}{T}))) \lo \tag{5.2.3.9}\label{ali:gppkpfx}\\
& 
P_{q+k'}
R^qf^{[p^n]}_{\os{\circ}{T}{}'*}
(P_kA_{{\rm zar},{\mab Q}}
(X^{[p^n]}_{\bul \leq N,\os{\circ}{T}{}'_1}/S^{[p^n]}(T')^{\nat}, 
E^{\bul \leq N}(\os{\circ}{T}{}')) 
\end{align*} 
is an isomorphism. 
Consequently the natural morphism 
\begin{align*} 
& g^*(R^q f^{[p^n]}_{\os{\circ}{T}*}
(P_kA_{{\rm zar},{\mab Q}}(X^{[p^n]}_{\bul \leq N,\os{\circ}{T}_1}/S^{[p^n]}({\os{\circ}{T})},
E^{\bul \leq N}(\os{\circ}{T}))),P) 
\os{\sim}{\lo} \tag{5.2.3.10}\label{ali:gppfpfx}\\
& (R^qf^{[p^n]}_{\os{\circ}{T}{}'*}
(P_kA_{{\rm zar},{\mab Q}}(X^{[p^n]}_{\bul \leq N,\os{\circ}{T}{}'_1}/S^{[p^n]}({\os{\circ}{T}{}')},
E^{\bul \leq N}(\os{\circ}{T}{}'))),P)
\end{align*} 
is an isomorphism.  
\end{proof}

\par 
Let $A$ be a commutative ring with unit element. 
Recall that we have said that a filtered $A$-module $(M,P)$ is filteredly flat 
if $M$ and $M/P_kM$ $(\forall k\in {\mab Z})$ are flat $A$-modules 
(\cite[(1.1.14)]{nh2}).

\begin{coro}\label{coro:flft}
For a hollow log $p$-adic enlargement $T$ of $S^{[p^n]}/{\cal V}$,    
the filtered sheaf
\begin{align*}  
(R^qf^{[p^n]}_{\os{\circ}{T}*}
(P_kA_{{\rm zar},{\mab Q}}(X^{[p^n]}_{\bul \leq N,\os{\circ}{T}_1}/S^{[p^n]}(T)^{\nat},
E^{\bul \leq N}(\os{\circ}{T}))),P) 
\tag{5.2.4.1}\label{ali:pce}
\end{align*} 
is a filteredly flat ${\cal K}_T$-modules. 
In particular, the filtered sheaf 
\begin{align*}  
(R^qf^{[p^n]}_{\os{\circ}{T}*}
(\eps^*_{X^{[p^n]}_{\bul \leq N,\os{\circ}{T}_1}
/S^{[p^n]}(T)^{\nat}}E^{\bul \leq N}(\os{\circ}{T})),P)
\tag{5.2.4.2}\label{ali:pcte}
\end{align*} 
is a filteredly flat ${\cal K}_T$-module. 
\end{coro} 
\begin{proof}
Since the question is local on $S^{[p^n]}(T)^{\nat}$, we may assume that 
$M_{S^{[p^n]}(T)^{\nat}_i}/{\cal O}_{S^{[p^n]}(T)^{\nat}_i}={\mab N}$. 
Since $M_{S^{[p^n]}(T)^{\nat}}/{\cal O}_{S^{[p^n]}(T)^{\nat}}^*
=M_{S^{[p^n]}(T)^{\nat}_i}/{\cal O}_{S^{[p^n]}(T)^{\nat}_i}$, 
the extension $1\lo {\cal O}_{S^{[p^n]}(T)^{\nat}}^*\lo 
M_{S^{[p^n]}(T)^{\nat}} \lo M_{S^{[p^n]}(T)^{\nat}}/{\cal O}^*_{S^{[p^n]}(T)^{\nat}}\lo 1$ 
is split. 
Because $T$ is hollow, we may assume that  
$T=S({\os{\circ}{T}})$ (this is the argument in the proof of \cite[Corollary 25]{ollc}). 
By (\ref{lemm:pnlcfi}) and \cite[(2.9)]{of},  
$$P_{q+k'}R^qf^{[p^n]}_{\os{\circ}{T}*}
(P_kA_{{\rm zar},{\mab Q}}
(X^{[p^n]}_{\bul \leq N,\os{\circ}{T}_1}/S^{[p^n]}(T)^{\nat},E^{\bul \leq N}(\os{\circ}{T})))$$ 
are flat ${\cal K}_T$-modules.  
\end{proof}

\begin{theo}[{\bf Log $p$-adic convergence of the weight filtration}]
\label{theo:pwfec} 
Let $k$, $q$ be two nonnegative integers. 
Then there exists a unique object
\begin{align*} 
(R^qf^{[p^n]}_{*}
(P_kA_{{\rm zar},{\mab Q}}(X^{[p^n]}_{\bul \leq N}/K,E^{\bul \leq N}))^{\sq},P)
\tag{5.2.5.1}\label{ali:pcpkxee}
\end{align*} 
of ${\rm IsocF}^{\sq}_p(S^{[p^n]}/{\cal V})$ 
such that  
\begin{align*} 
&(R^qf^{[p^n]}_{*}
(P_kA_{{\rm zar},{\mab Q}}(X^{[p^n]}_{\bul \leq N}/K,E^{\bul \leq N}))^{\sq},P)_{T}
\tag{5.2.5.2}\label{ali:pkaskf} \\
&=
(R^qf^{[p^n]}_{\os{\circ}{T}*}
(P_kA_{{\rm zar},{\mab Q}}
(X^{[p^n]}_{\bul \leq N,\os{\circ}{T}_1}/S^{[p^n]}(T)^{\nat},E^{\bul \leq N}(\os{\circ}{T}))),P)
\end{align*} 
for any object $T$ of ${\rm Enl}^{\sq}_p(S^{[p^n]}/{\cal V})$.  
In particular, there exists a unique object
\begin{align*} 
(R^qf^{[p^n]}_{*}
(\eps^*_{X^{[p^n]}_{\bul \leq N}/K}(E^{\bul \leq N}_K))^{\nat,\sq},P)
\tag{5.2.5.3}\label{ali:pcee}
\end{align*}  
of ${\rm IsocF}^{\sq}_p(S^{[p^n]}/{\cal V})$ 
such that  
\begin{align*} 
(R^q f^{[p^n]}_{*}
(\eps^*_{X_{\bul \leq N}/K}(E^{\bul \leq N}_K))^{\nat,\sq},P)_T=
(R^qf^{[p^n]}_{X^{[p^n]}_{\bul \leq N,\os{\circ}{T}_1}/S^{[p^n]}(T)^{\nat}*}
(\eps^*_{X^{[p^n]}_{\bul \leq N,T_1}/S^{[p^n]}(T)^{\nat}}
(E^{\bul \leq N}(\os{\circ}{T}))),P)
\tag{5.2.5.4}\label{ali:apce}
\end{align*} 
for any object $T$ of ${\rm Enl}^{\sq}_p(S^{[p^n]}/{\cal V})$. 
\end{theo} 
\begin{proof} 
Let $g\col T'\lo T$ be a morphism in ${\rm Enl}_p(S^{[p^n]}/{\cal V})$. 
We may assume that $\os{\circ}{T}$ and $\os{\circ}{T}{}'$ are affine formal schemes. 
Let $z\col T_1\lo S^{[p^n]}$ be the structural morphism. 
We may assume that $S^{[p^n]}$ has a global chart ${\mab N}\lo {\cal O}_{S^{[p^n]}}$. 
Because the question is local on $T$, 
we may assume that $M_{S^{[p^n]}(T)^{\nat}}$ is split. 
By the commutative diagram 
\begin{equation*} 
\begin{CD} 
M_{S^{[p^n]}(T')^{\nat}}@>>> M_{S^{[p^n]}(T')^{\nat}}/{\cal O}_{T'}^*={\mab N}\\
@AAA @AAA \\
g^{-1}(M_{S^{[p^n]}(T)^{\nat}})@>>> g^{-1}(M_{S^{[p^n]}(T)^{\nat}}/{\cal O}_{T}^*)={\mab N}
\end{CD}
\end{equation*} 
a splitting of $M_{S^{[p^n]}(T)^{\nat}}$ also gives a splitting of $M_{S^{[p^n]}(T')^{\nat}}$. 
Hence, by (\ref{prop:csds}), we may assume that  
$S^{[p^n]}(T)^{\nat}$ and $S^{[p^n]}(T')^{\nat}$ are equal to $S^{[p^n]}({\os{\circ}{T}})$ and 
$S^{[p^n]}({\os{\circ}{T}{}'})$, respectively. 
By (\ref{lemm:pnlcfi}) we obtain (\ref{theo:pwfec}). 
\end{proof}


\parno 
The following is a filtered version of \cite[(3.5)]{of}. 
To consider the category ${\rm IsocF}^{\rm sld}_p(S^{[p^n]}/{\cal V})$ 
(not ${\rm IsocF}_p(S^{[p^n]}/{\cal V})$) is important. 

\begin{lemm}\label{lemm:nmr}
Let ${\cal V}'$ be a finite extension of complete discrete valuation ring of ${\cal V}$. 
Let $h\col S'\lo S$ be a morphism of $p$-adic formal families of log points over 
${\rm Spec}({\cal V}')\lo {\rm Spec}({\cal V})$. 
Let $g\col Y_{\bul \leq N}\lo X_{\bul \leq N}$ be 
a morphism of $N$-truncated simplicial base changes 
of SNCL schemes over $S'\lo S$.
Assume that 
the admissible immersion {\rm (\ref{eqn:audtds})} for $Y_{\bul \leq N}/S'$ 
and for any object of ${\rm Enl}^{\rm sld}_p((S')^{[p^n]}/{\cal V}')$ exists. 
Let $f'\col Y_{\bul \leq N}\lo S'$ be the structural morphism. 
Then there exists a natural morphism 
\begin{align*} 
g^*\col &
g^*((R^qf^{[p^n]}_{*}
(P_kA_{{\rm zar},{\mab Q}}(X^{[p^n]}_{\bul \leq N}/K,E^{\bul \leq N}))^{\rm sld},P))\\
&\lo (R^qf'{}^{[p^n]}_{*}
(P_kA_{{\rm zar},{\mab Q}}
(Y^{[p^n]}_{\bul \leq N}/K',\os{\circ}{g}{}^*(E^{\bul \leq N})))^{\rm sld},P). 
\end{align*}
in ${\rm IsocF}^{\rm sld}_p((S')^{[p^n]}/{\cal V}')$.  
If $Y_{\bul \leq N}=X_{\bul \leq N}\times_SS'
$, 
then this morphism is an isomorphism. 
\end{lemm}
\begin{proof} 
Note that the log formal scheme $S^{[p^n]}_{\os{\circ}{S}{}'{}^{[p^n]}}$ 
has the universality for the log formal schemes 
$U$ with morphisms $S'{}^{[p^n]}\lo U$ and $U\lo S^{[p^n]}$ such that 
the composite morphism $S'\lo U\lo S$ is equal to the morphism $S'{}^{[p^n]}\lo S^{[p^n]}$ 
and the morphism $U\lo S^{[p^n]}$ is solid, i.e., 
$U=S^{[p^n]}\times_{\os{\circ}{S^{[p^n]}}}\os{\circ}{U}$. 
Now (\ref{lemm:nmr}) follows from (\ref{ali:rstzsfb}) and 
(\ref{theo:pwfec}). 
Indeed, for an object $T'$ of 
${\rm Enl}^{\rm sld}_p((S')^{[p^n]}/{\cal V}')$, 
\begin{align*} 
[g^*((R^qf^{[p^n]}_{*}
(P_kA_{{\rm zar},{\mab Q}}(X^{[p^n]}_{\bul \leq N}/K,E^{\bul \leq N}))^{\rm sld},P))]_{T'}
=
(R^qf^{[p^n]}_{*}(A_{{\rm zar},{\mab Q}}(X^{[p^n]}_{\bul \leq N,\os{\circ}{T}{}'}/
S(T')^{\nat},E^{\bul \leq N})(\os{\circ}{T}{}'),P),  
\end{align*} 
where $S(T')^{\nat}$ is defined by 
the composite morphism $T'_1\lo S'_1\lo S_1$. 
\end{proof} 

\begin{rema}\label{rema:ud}
As in \cite[(3.6)]{of}, 
$(R^qf^{[p^n]}_{*}
(P_kA_{{\rm zar},{\mab Q}}(X^{[p^n]}_{\bul \leq N}/K,E^{\bul \leq N})),P)$ 
in ${\rm IsocF}^{\sq}_p(S^{[p^n]}/{\cal V})$ descends to 
the object 
$(R^qf^{[p^n]}_{*}
(P_kA_{{\rm zar},{\mab Q}}(X^{[p^n]}_{\bul \leq N}/K_0,E^{\bul \leq N})),P)$
of ${\rm IsocF}^{\sq}_p(S^{[p^n]}/{\cal W})$.  
\end{rema}

\par 
To obtains the main result (\ref{theo:pwfaec}) below in this section, 
we assume the existence of the immersion 
(\ref{eqn:audtds}) only for the case $n=0$ and $1$. 
\par 

\begin{theo}[{\bf Log convergence of the weight filtration}]\label{theo:pwfaec}
Let the notations and the assumptions be as above.  
Consider the morphisms $X^{[p^m]}_{\bul \leq N}\lo S^{[p^m]}$ 
over ${\rm Spf}({\cal V})$ $(m=0,1)$ 
as a morphism $X^{[p^m]}_{\bul \leq N}\lo S^{[p^m]}$ 
over ${\rm Spf}({\cal W})$, respectively. 
Set ${\cal E}^{\sq}_{\star,{\cal W}}:={\rm Enl}^{\sq}_{\star}(S/{\cal W})$.  
Let 
$F^{\rm ar}_{X_{\bul \leq N}/S,S^{[p]}/
{\cal W}} \col X_{\bul \leq N}\lo X^{[p]}_{\bul \leq N}$ 
be the abrelative Frobenius morphism 
over the morphism $S\lo S^{[p]}$ over ${\rm Spf}({\cal W})$.  
Let $W_{X_{\bul \leq N}/S^{[p]}, S/{\cal W}}
\col X^{[p]}_{\bul \leq N}\lo X_{\bul \leq N}$ 
be the projection. 
Let $E^{\bul \leq N}:=\{E^{\bul \leq N}_n\}_{n=0}^{\inf}$ be 
a sequence of flat coherent 
$\{{\cal O}_{\os{\circ}{X}{}_{\bul \leq N,\os{\circ}{T}_1}
/\os{\circ}{T}}\}_{T\in {\cal E}^{\sq}_{p,{\cal W}}}$-modules 
with a morphism 
\begin{align*} 
\Phi^{\bul \leq N}_n\col 
F_{\os{\circ}{X}{}_{\bul \leq N}}^* 
(E^{\bul \leq N}_{n+1})\lo E^{\bul \leq N}_n
\tag{5.2.8.1}\label{ali:pchkee} 
\end{align*} 
of $\{{\cal O}_{\os{\circ}{X}{}_{\bul \leq N,\os{\circ}{T}_1}
/\os{\circ}{T}}\}_{T\in {\cal E}^{\sq}_{p,{\cal W}}}$-modules.  
Let $\os{\circ}{W}{}^{(l)}_{m,T}
\col (\os{\circ}{X}{}^{[p]}_{m,T_1})^{(l)}\lo \os{\circ}{X}{}^{(l)}_{m,T_1}$ 
$(l\in {\mab N}, 0\leq m\leq N)$ 
be also the projection over $\os{\circ}{T}$. 
Assume that, for any $0\leq m\leq N$, $l\geq 0$ and $n\geq 0$, 
the morphism 
\begin{align*} 
R^qf_{(\os{\circ}{X}{}^{[p]}_{m,T_1})^{(l)}/\os{\circ}{T}*}
(\os{\circ}{W}{}^{(l)*}_{m,T,{\rm crys}}
(E^m_{n+1}(\os{\circ}{T})_{(\os{\circ}{X}{}^{[p]}_{m,T_1})^{(l)}/\os{\circ}{T}}))
_{\mab Q}\lo  
R^qf_{\os{\circ}{X}{}^{(l)}_{T_1}/\os{\circ}{T}*}
(E^m_n(\os{\circ}{T})_{\os{\circ}{X}{}^{(l)}_{T_1}/\os{\circ}{T}})_{\mab Q}
\end{align*}
induced by $\Phi^m_n$ is an isomorphism for any object $T$ of 
${\cal E}^{\sq}_{p,{\cal W}}$.  
Then there exists an object
\begin{align*} 
\{((R^q f_*
(P_kA_{{\rm zar},{\mab Q}}
(X_{\bul \leq N}/K,E^{\bul \leq N}_n))^{\sq},P),\Phi_n)\}_{n=0}^{\inf}
\tag{5.2.8.2}\label{ali:pcxpkee} 
\end{align*}  
of $F^{\inf}{\textrm -}{\rm IsosF}^{\sq}(S/{\cal V})$ such that 
\begin{align*} 
(R^q f_*(P_kA_{{\rm zar},{\mab Q}}(X_{\bul \leq N}/K,E^{\bul \leq N}_n))^{\sq}_T,P)=
(R^q f_*(P_kA_{{\rm zar},{\mab Q}}(X_{\bul \leq N,T_1}/S(T)^{\nat},
E^{\bul \leq N}_n(\os{\circ}{T}))),P)
\tag{5.2.8.3}\label{ali:pcpnle}
\end{align*} 
for any object $T$ of ${\rm Enl}^{\sq}_p(S/{\cal V})$. 
In particular, there exists an object
\begin{align*} 
\{(R^qf_{*}
(\eps^*_{X_{\bul \leq N}/K}(E^{\bul \leq N}_{n,K}))^{\nat,\sq},P),
\Phi_n\}_{n=0}^{\inf}
\tag{5.2.8.4}\label{ali:pcxakee} 
\end{align*}  
of $F^{\infty}{\textrm -}{\rm IsosF}^{\sq}(S/{\cal V})$ such that 
\begin{align*} 
\{(R^qf_*
(\eps^*_{X_{\bul \leq N}/K}(E^{\bul \leq N}_{n,K}))^{\nat,\sq}_T,P)\}_{n=0}^{\inf}=
(R^qf_{X_{\bul \leq N,\os{\circ}{T}_1}/S(T)^{\nat}*}
(\eps^*_{X_{\bul \leq N,\os{\circ}{T}_1}/S(T)^{\nat}}
(E^{\bul \leq N}_n(\os{\circ}{T})))_{\mab Q},P) 
\tag{5.2.8.5}\label{ali:pcxnle}
\end{align*} 
for any object $T$ of ${\rm Enl}^{\sq}_p(S/{\cal V})$. 
\end{theo}
\begin{proof}
By (\ref{ali:fsriszb}) we may assume that $\sq={\rm sld}$. 
(Note that the morphism {\rm (\ref{ali:ztio})} is injective for $Z=S$.)
By (\ref{rema:ud}) we may assume that ${\cal V}={\cal W}$ and we denote 
${\cal E}^{\rm sld}_{p,{\cal W}}$ simply by ${\cal E}^{\rm sld}_p$.  
We may assume that $p{\cal O}_S=0$. 
Let $F_S\col S\lo S$ be the Frobenius endomorphism of $S$ 
over the Frobenius endomorphism 
${\rm Spf}({\cal W})\lo {\rm Spf}({\cal W})$ of ${\rm Spf}({\cal W})$. 
We have to prove that there exists a filtered isomorphism 
\begin{align*} 
\Psi_n &\col 
F_S^*((R^q f_*(P_kA_{{\rm zar},{\mab Q}}
(X_{\bul \leq N}/K,E^{\bul \leq N}_n))^{\rm sld},P))
\tag{5.2.8.6}\label{ali:pckae}\\
& \os{\sim}{\lo} 
(R^q f_*(P_kA_{{\rm zar},{\mab Q}}(X_{\bul \leq N}/K,E^{\bul \leq N}_n))^{\rm sld},P). 
\end{align*} 
By the proof of (\ref{lemm:pnlcfi}), we may forget the filtration $P$'s 
in (\ref{ali:pckae}). 
\par 
Let $T=(T,z)$ be an object of ${\cal E}^{\rm sld}_p$; 
$z\col T_1\lo S$ is the structural morphism. 
The abrelative Frobenius morphism 
\begin{align*} 
F^{\rm ar}_{X_{\bul \leq N,1}/S_1,S^{[p]}_1}
\col X_{\bul \leq N,1} \lo X^{[p]}_{\bul \leq N,1}
\end{align*} 
over $S_{\os{\circ}{T}_1}\lo S^{[p]}_{\os{\circ}{T}_1}$
induces the abrelative Frobenius morphism   
\begin{align*} 
F^{\rm ar}_{X_{\bul \leq N,\os{\circ}{T}_1}/S(T)^{\nat},S^{[p]}(T)^{\nat}} 
\col X_{\bul \leq N,T_1} \lo X^{[p]}_{\bul \leq N,T_1}
\tag{5.2.8.7}\label{ali:tzxfs} 
\end{align*} 
over 
$S(T)^{\nat}\lo S^{[p]}(T)^{\nat}$. 
Then, by (\ref{coro:npc})  (cf.~(\ref{exem:oot}), (\ref{prop:bsch})) 
and (\ref{lemm:nmr}), we obtain the following equalities:   
\begin{align*} 
&(F_S^*(R^q f_*(P_kA_{{\rm zar},{\mab Q}}
(X_{\bul \leq N}/K,E^{\bul \leq N}_{n+1}))^{\rm sld}))_T
=(R^q f_*(P_kA_{{\rm zar},{\mab Q}}
(X_{\bul \leq N}/K,E^{\bul \leq N}_{n+1}))^{\rm sld}))_{S^{[p]}(T)}\\
&=R^qf^{[p]}_*
(A_{{\rm zar},{\mab Q}}(X^{[p]}_{\bul \leq N,\os{\circ}{T}_1}/S^{[p]}(T)^{\nat},
\os{\circ}{W}{}^*_{X_{\bul \leq N,\os{\circ}{T}_1}/\os{\circ}{T}}
(E^{\bul \leq N}_{n+1}(\os{\circ}{T})))). 
\end{align*} 
The morphism (\ref{ali:pchkee}) induces the following morphism 
\begin{align*} 
\Phi^{\bul \leq N}_n\col 
\os{\circ}{F}{}^{{\rm ar}*}_{X_{\bul \leq N}/S,S^{[p]}/{\cal W}} 
(\os{\circ}{W}{}^*_{X_{\bul \leq N}/{\cal W}}
(E^{\bul \leq N}_{n+1}))\lo E^{\bul \leq N}_n. 
\tag{5.2.8.8}\label{ali:pcaaee} 
\end{align*} 
Hence we have the following morphism 
\begin{align*} 
(\Phi^{\bul \leq N}_n)_T\col 
\os{\circ}{F}{}^*_{X_{\bul \leq N,\os{\circ}{T}_1}/S(T),S^{[p]}(T)}
(\os{\circ}{W}{}^*_{X_{\bul \leq N,\os{\circ}{T}_1}/\os{\circ}{T}}
(E^{\bul \leq N}_{n+1}(\os{\circ}{T})))\lo E^{\bul \leq N}_n(\os{\circ}{T})
\tag{5.2.8.9}\label{ali:pcaee} 
\end{align*} 
of ${\cal O}_{\os{\circ}{X}{}_{\bul \leq N,\os{\circ}{T}_1}/\os{\circ}{T}}$-modules.  
The morphisms (\ref{ali:tzxfs}) and (\ref{ali:pcaee}) induce the following morphism 
\begin{align*} 
\Phi^{\rm ar}_T \col &
(A_{{\rm zar},{\mab Q}}
(X^{[p]}_{\bul \leq N,\os{\circ}{T}{}_1}/S^{[p]}(T)^{\nat},
\os{\circ}{W}{}^*_{X_{\bul \leq N,\os{\circ}{T}_1}/\os{\circ}{T}}
(E^{\bul \leq N}_{n+1}(\os{\circ}{T})),P)
\tag{5.2.8.10}\label{ali:sqpabpts} \\
&\lo 
RF^{\rm ar}_{X_{\bul \leq N,\os{\circ}{T}_1/S^{[p]}(T)^{\nat},S(T)^{\nat}}*}
((A_{{\rm zar},{\mab Q}}(X_{\bul \leq N,\os{\circ}{T}_1}/S(T)^{\nat},
E^{\bul \leq N}_n(\os{\circ}{T})),P)).  
\end{align*}  
Hence we have the following morphism 
\begin{align*} 
\Phi^{\rm ar}_T\col &
R^qf^{[p]}_{\os{\circ}{T}*}(P_kA_{{\rm zar},{\mab Q}}
(X^{[p]}_{\bul \leq N,\os{\circ}{T}_1}/S^{[p]}(T)^{\nat}, 
\os{\circ}{W}{}^*_{X_{\bul \leq N,\os{\circ}{T}_1}/\os{\circ}{T}}
(E^{\bul \leq N}_{n+1}(\os{\circ}{T}))))
\tag{5.2.8.11}\label{ali:sarchpts}\\
&\lo 
R^qf_{\os{\circ}{T}*}(P_kA_{{\rm zar},{\mab Q}}
(X_{\bul \leq N,\os{\circ}{T}_1}/S(T)^{\nat},E^{\bul \leq N}_n(\os{\circ}{T}))). 
\end{align*} 
\par 
Let ${\rm SS}_k
(X_{\bul \leq N,\os{\circ}{T}_1}/S(T)^{\nat},E^{\bul \leq N}_n(\os{\circ}{T}))$ 
and 
${\rm SS}_k
(X^{[p]}_{\bul \leq N,\os{\circ}{T}_1}/S^{[p]}(T)^{\nat},
\os{\circ}{W}{}^*_{X_{\bul \leq N,\os{\circ}{T}_1}/\os{\circ}{T}}
(E^{\bul \leq N}_{n+1}(\os{\circ}{T})))$  
be the spectral sequences 
of 
$$R^qf_{\os{\circ}{T}*}(P_kA_{{\rm zar},{\mab Q}}
(X_{\bul \leq N,\os{\circ}{T}_1}/S(T)^{\nat},E^{\bul \leq N}_n(\os{\circ}{T})))$$
and 
$$R^qf^{[p]}_{\os{\circ}{T}*}(P_kA_{{\rm zar},{\mab Q}}
(X^{[p]}_{\bul \leq N,\os{\circ}{T}_1}/S^{[p]}(T)^{\nat}, 
\os{\circ}{W}{}^*_{X_{\bul \leq N,\os{\circ}{T}_1}/\os{\circ}{T}}
(E^{\bul \leq N}_{n+1}(\os{\circ}{T})))),$$ 
respectively  ((\ref{eqn:ekpssp})). 
Then
we have the pull-back morphism 
$${\rm SS}_k(X^{[p]}_{\bul \leq N,\os{\circ}{T}_1}/S^{[p]}(T)^{\nat},
\os{\circ}{W}{}^*_{X_{\bul \leq N,\os{\circ}{T}_1}/\os{\circ}{T}}
(E^{\bul \leq N}_{n+1}(\os{\circ}{T})))\lo 
{\rm SS}_k(X_{\bul \leq N,\os{\circ}{T}_1}/S(T)^{\nat},E^{\bul \leq N}_n(\os{\circ}{T})).$$ 
Let  
$E_{1,n}^{\bul \bul}$ 
and
$E_{1,n+1}^{\bul \bul}$ 
be the $E_1$-terms of the spectral sequences of 
$${\rm SS}_k(X_{\bul \leq N,\os{\circ}{T}_1}/S(T)^{\nat},E^{\bul \leq N}_n(\os{\circ}{T}))$$
and  
$${\rm SS}_k(X^{[p]}_{\bul \leq N,\os{\circ}{T}_1}/S^{[p]}(T)^{\nat},
\os{\circ}{W}{}^*_{X_{\bul \leq N,\os{\circ}{T}_1}/\os{\circ}{T}}
(E^{\bul \leq N}_{n+1}(\os{\circ}{T}))),$$ 
respectively. 
Then, by the assumption on the induced morphism 
\begin{align*} 
E_{1,n+1}^{\bul \bul}\lo E_{1,n}^{\bul \bul} 
\tag{5.2.8.10}\label{ali:epkpk} 
\end{align*} 
by the abrelative Frobenius morphism, 
the morphism (\ref{ali:epkpk}) is an isomorphism. 
Hence the induced morphism 
\begin{align*} 
&R^qf^{[p]}_{\bul \leq N,\os{\circ}{T}*}
(P_kA_{{\rm zar},{\mab Q}}(X^{[p]}_{\bul \leq N,\os{\circ}{T}_1}
/S^{[p]}(T)^{\nat},E^{\bul \leq N}_{n+1}(\os{\circ}{T})))\\
& \lo 
R^qf_{\bul \leq N,\os{\circ}{T}*}
(P_kA_{{\rm zar},{\mab Q}}(X_{\bul \leq N,\os{\circ}{T}_1}
/S(T)^{\nat},E^{\bul \leq N}_n(\os{\circ}{T}))) 
\end{align*} 
by the abrelative Frobenius morphism 
is an isomorphism. 
Consequently, the induced morphism
\begin{align*} 
\Phi^*_n\col &
(R^qf^{[p]}_{\bul \leq N,\os{\circ}{T}*}(P_kA_{{\rm zar},{\mab Q}}
(X^{[p]}_{\bul \leq N,\os{\circ}{T}_1}/S^{[p]}(T)^{\nat},
E^{\bul \leq N}_{n+1}(\os{\circ}{T}))),P)
\lo \tag{5.2.8.11}\label{ali:xpteqp}\\
&(R^qf_{\bul \leq N,\os{\circ}{T}*}(P_kA_{{\rm zar},{\mab Q}}
(X_{\bul \leq N,\os{\circ}{T}_1}/S(T)^{\nat},
E^{\bul \leq N}_n(\os{\circ}{T}))),P)
\end{align*} 
by the abrelative Frobenius morphism is an isomorphism. 
Because $M^{\rm gp}_S/{\cal O}_S^*$ is $p$-torsion-free and 
because  
we have already constructed an object of 
$F^{\infty}{\textrm -}{\rm IsosF}^{\rm sld}_p(S/{\cal W})$, 
we obtain an object  
$F^{\infty}{\textrm -}{\rm IsosF}^{\rm sld}(Z/{\cal W})$ by (\ref{theo:zbfs}). 
\end{proof}

\begin{coro}\label{coro:filc} 
Let the notations and the assumption be as in {\rm (\ref{theo:pwfaec})}. 
Then there exists a filtered $F^{\infty}$-isospan 
\begin{align*} 
\{(R^qf_{*}
(\eps^*_{X_{\bul \leq N}/K}(E^{\bul \leq N}_{n,K}))^{\sq},P),
\Phi_n\}_{n=0}^{\inf}
\tag{5.2.9.1}\label{ali:pcesee} 
\end{align*}  
on ${\rm rhEnl}^{\sq}(S/{\cal V})$ such that 
\begin{align*} 
\{(R^qf_*
(\eps^*_{X_{\bul \leq N}/K}(E^{\bul \leq N}_{n,K}))^{\sq}_T,P)\}_{n=0}^{\inf}=
(R^qf_{X_{\bul \leq N,\os{\circ}{T}_1}/S(T)*}
(\eps^*_{X_{\bul \leq N,\os{\circ}{T}_1}/S(T)}
(E^{\bul \leq N}_n(\os{\circ}{T})))_{\mab Q},P) 
\tag{5.2.9.2}\label{ali:pqple}
\end{align*} 
for any object $T$ of 
${\rm Enl}^{\sq}_{\star}(S/{\cal V})$. 
\end{coro} 
\begin{proof} 
Let $T$ be an object of 
${\rm rhEnl}^{\sq}_p(S/{\cal V})$. 
Then $S(T)=S(T)^{\nat}$. 
Hence this corollary follows from (\ref{theo:pwfaec}).  
\end{proof}

\begin{rema}\label{rema:nlcfi} 
Let the notations and the assumption be as in {\rm (\ref{theo:pwfaec})}. 
As in (\ref{lemm:pnlcfi}), we see that the filtered log convergent $F$-isocrystal 
$$\{(R^qf_*(P_kA_{{\rm zar},{\mab Q}}(X_{\bul \leq N}/K,E^{\bul \leq N}_n)),
P)\}_{n=0}^{\inf}$$  
on ${\rm Enl}^{\sq}_{\star}(S/{\cal V})$ also 
gives a filtered convergent $F^{\inf}$-isospan 
on ${\rm Enl}_{\star}(\os{\circ}{S}/{\cal V})$ locally on $\os{\circ}{S}$. 
\end{rema}

\begin{coro}\label{coro:fenlt} 
Assume that $E^{\bul \leq N}_n=E^{\bul \leq N}_0$ $(n\in {\mab N})$. 
Set $E^{\bul \leq N}:=E^{\bul \leq N}_0$. 
For any object $T$ of ${\rm Enl}^{\sq}(S/{\cal V})$,    
$$(R^qf_{T*}(P_kA_{{\rm zar},{\mab Q}}(X_{\bul \leq N,\os{\circ}{T}_1}/S(T)^{\nat},
E^{\bul \leq N}(\os{\circ}{T}))),P)$$ 
is a filteredly flat ${\cal K}_T$-module. 
In particular, the filtered sheaf 
\begin{align*}  
(R^qf_{T*}(\eps^*_{X_{\bul \leq N,\os{\circ}{T}_1}/S(T)^{\nat}}
(E^{\bul \leq N}(\os{\circ}{T})))_{\mab Q},P) 
\tag{5.2.11.1}\label{ali:bpce}
\end{align*} 
is a filteredly flat ${\cal K}_T$-module. 
\end{coro} 
\begin{proof} 
This follows from (\ref{coro:flft}), (\ref{theo:pwfaec}), (\ref{prop:csds}) 
and \cite[(2.9)]{of}. 
\end{proof}

\begin{exem}\label{exam:ofl} 
Let the notations be as before (\ref{theo:pwfaec}) (1). 
Then there exists an object
$$(R^qf_*(P_kA_{{\rm zar},{\mab Q}}(X_{\bul \leq N}/K))^{\sq},P)$$  
of $F{\textrm -}{\rm IsosF}^{\sq}(S/{\cal V})$
such that 
\begin{align*} 
(R^qf_*(P_kA_{{\rm zar},{\mab Q}}(X_{\bul \leq N}/K))^{\sq},P)_T
=(R^qf_*(P_kA_{{\rm zar},{\mab Q}}(X_{\bul \leq N,\os{\circ}{T}_1}/S(T)^{\nat}),P)
\tag{5.2.12.1}\label{ali:pcoe}
\end{align*} 
for any object $T$ of ${\rm Enl}^{\sq}_p(S/{\cal V})$.   
Indeed, the assumption in (\ref{theo:pwfaec})  
is satisfied by the base change of \cite[(1.3)]{boi} 
(cf.~the proof of \cite[(3.7)]{of}).
In particular,  there exists an object 
\begin{align*} 
(R^qf_*({\cal O}_{X_{\bul \leq N}/K})^{\nat,\sq},P)
\tag{5.2.12.2}\label{ali:pcxnee} 
\end{align*}  
of $F{\textrm -}{\rm IsosF}^{\sq}(S/{\cal V})$ such that 
\begin{align*} 
(R^qf_*({\cal O}_{X_{\bul \leq N}/K})^{{\nat},{\sq}},P)_T=
(R^qf_{X_{\bul \leq N,\os{\circ}{T}_1}/S(T)^{\nat}*}
({\cal O}_{X_{\bul \leq N,\os{\circ}{T}_1}/S(T)^{\nat}})_{\mab Q},P) 
\tag{5.2.12.3}\label{ali:pctoe}
\end{align*} 
for any object $T$ of ${\rm Enl}_p^{\sq}(S/{\cal V})$. 
We also have a filtered $F^{\infty}$-isospan 
\begin{align*} 
(R^qf_*({\cal O}_{X_{\bul \leq N}/K})^{\sq},P)
\tag{5.2.12.4}\label{ali:pckee} 
\end{align*}  
on ${\rm rhEnl}^{\sq}(S/{\cal V})$ such that 
\begin{align*} 
(R^qf_*({\cal O}_{X_{\bul \leq N}/K})^{\sq},P)_T=
(R^qf_{X_{\bul \leq N,\os{\circ}{T}_1}/S(T)*}
({\cal O}_{X_{\bul \leq N,\os{\circ}{T}_1}/S(T)})_{\mab Q},P) 
\tag{5.2.12.5}\label{ali:pctkoe}
\end{align*} 
for any object $T$ of ${\rm rhEnl}_p^{\sq}(S/{\cal V})$. 
\end{exem}

\bigskip
\parno
{\bf (2) Functoriality}
\bigskip
\begin{prop}\label{prop:fcuccv}
$(1)$ 
Let $g$ be as in {\rm (\ref{theo:itc})}. 
Let the notations and the assumption be as in 
{\rm (\ref{theo:pwfaec})}. 
Then the log $p$-adically convergent isocrystal 
$P_kR^qf_{X_{\bul \leq N}/K*}
(\eps^*_{X_{\bul \leq N}/K}
(E^{\bul \leq N}_{n,K}))$ 
is contravariantly functorial. 
\end{prop} 
\begin{proof} 
This is obvious. 
\end{proof}

\begin{rema} 
We leave the reader to the formulations of the log convergent 
$F$-isocrystal versions of (\ref{prop:fcuccv}) (cf.~(\ref{prop:spcnvuc}) below). 
\end{rema}

\begin{prop}\label{prop:spcnvuc}
Let the notations and the assumption be as in 
{\rm (\ref{theo:pwfaec})}. Assume that $E^{\bul \leq N}_n=
E^{\bul \leq N}_0$ $(n\in {\mab N})$. 
Set $E^{\bul \leq N}:=E^{\bul \leq N}_0$. 
Let ${\cal V}'/{\cal V}$ be a finite extension.  
Let $S'\lo S$ be a morphism of log $p$-adic formal families of log points 
over $Spf({\cal V}')\lo Spf({\cal V})$. 
Let $T'$ and $T$ be log $(p$-adic$)$ enlargements of $S'$ and $S$, respectively. 
Let $T'\lo T$ be a morphism of log $(p$-adic$)$ enlargements over 
$S'\lo S$. Let $u\col S'(T')^{\nat}\lo S(T)^{\nat}$ be the induced morphism.  
Let $Y_{\bul \leq N}$ be a log scheme over $S'$ 
which is similar to $X_{\bul \leq N}$ over $S$. 
Let $F^{\bul \leq N}$ be 
a similar $F$-isocrystal of 
${\cal O}_{\{\os{\circ}{Y}_{\bul \leq N,T'_1}
/\os{\circ}{T}{}'\}_{T'\in {\rm Enl}^{\sq}_p(S'/{\cal V}')}}$-modules to 
$E^{\bul \leq N}$. 
Let $k$ and $q$ be nonnegative integers. 
Let 
$R^qf_{\os{\circ}{X}{}^{(k)}_m/K*}
(E^m_{\os{\circ}{X}{}^{(k)}_m/K} 
\otimes_{\mab Z}\varpi^{(k)}
(\os{\circ}{X}_m/K))$ be an object of $F{\textrm -}{\rm Isoc}^{\sq}_{\star}(S/{\cal V})$  
such that
$$R^qf_{\os{\circ}{X}{}^{(k)}_m/K*}
(E^m_{\os{\circ}{X}{}^{(k)}_m/K} 
\otimes_{\mab Z}\varpi^{(k)}
(\os{\circ}{X}_m/K))_T= 
R^qf_{\os{\circ}{X}{}^{(k)}_{m,T_1}
/\os{\circ}{T}}
(E^m(T)_{\os{\circ}{X}{}^{(k)}_{m,T_1}
/\os{\circ}{T}} \otimes_{\mab Z}\varpi^{(k)}_{\rm crys}(\os{\circ}{X}_{m,T_1}/\os{\circ}{T}))
_{\mab Q}$$  
for any object $T$ of ${\rm Enl}^{\sq}_p(S/{\cal V})$. 
Then 
there exists the following spectral sequence 
in $F{\textrm -}{\rm Isoc}^{\sq}(S/{\cal V}):$
\begin{equation*} 
\begin{split} 
{} & 
E_1^{-k,q+k}=
\bigoplus_{m\geq 0}\bigoplus_{j\geq \max \{-(k+m),0\}} 
R^{q-2j-k-2m}
f_{\os{\circ}{X}{}^{(2j+k+m)}_m/K*}
(E^m_{\os{\circ}{X}{}^{(2j+k+m)}_m/K} \\ 
{} & \phantom{R^{q-2j-k-m}f_{(\os{\circ}{X}^{(k)}, 
Z\vert_{\os{\circ}{X}^{(2j+k)}})/K*} 
({\cal O}}\otimes_{\mab Z}\varpi^{(2j+k+m)}(\os{\circ}{X}_m/K))(-j-k-m,u) \\
&\Lo 
R^qf_{X_{\bul \leq N}/K*}
(\eps^*_{X_{\bul \leq N}/K}
(E^{\bul \leq N}_K))
\quad (q\in {\mab Z}). 
\end{split} 
\tag{5.2.15.1}\label{eqn:getcpsp}
\end{equation*}  
\end{prop}
\begin{proof} 
This follows from (\ref{prop:grlgoc}). 
\end{proof}

\begin{defi}
We call (\ref{eqn:getcpsp}) 
the {\it Poincar\'{e} spectral sequence} of 
$$R^qf_{X_{\bul \leq N}/K*}
(\eps^*_{X_{\bul \leq N}/K}(E^{\bul \leq N}_K))$$ 
in ${\rm Isoc}_p^{\sq}(S/{\cal V})$ and $F{\textrm -}{\rm Isoc}^{\sq}(S/{\cal V})$, 
respectively.
\end{defi}


\bigskip
\parno
{\bf (3) Monodromy and the cup product of a line bundle}
\bigskip
\parno
Let the notations and the assumptions be as in {\rm (\ref{theo:pwfaec})}. 
Assume that $E^{\bul \leq N}_n=
E^{\bul \leq N}_0$ $(n\in {\mab N})$. 
Set $E^{\bul \leq N}:=E^{\bul \leq N}_0$.

\begin{prop}\label{prop:convmon} 
$(1)$ There exists the monodromy operator 
\begin{align*} 
N_{\rm zar}=\nu_{\rm zar} \col & 
(R^qf_{X_{\bul \leq N}/K*}
(\eps^*_{X_{\bul \leq N}/K}
(E^{\bul \leq N}_K)),P)
\lo \\
&(R^qf_{X_{\bul \leq N}/K*}
(\eps^*_{X_{\bul \leq N}/K}
(E^{\bul \leq N}_K)),
P\langle -2 \rangle)(-1,u) 
\end{align*} 
in ${\rm IsocF}_p^{\sq}(S/{\cal V})$. 
\par 
$(2)$ 
There exists the monodromy operator  
\begin{align*} 
N_{\rm zar}=\nu_{\rm zar} \col & 
(R^qf_{X_{\bul \leq N}/S*}
(\eps^*_{X_{\bul \leq N}/K}
(E^{\bul \leq N}_K)),P)
\lo \\
&(R^qf_{X_{\bul \leq N}/S*}
(\eps^*_{X_{\bul \leq N}/K}
(E^{\bul \leq N}_K)),
P\langle -2 \rangle)(-1,u) 
\end{align*} 
in $F{\textrm -}{\rm IsocF}^{\sq}(S/{\cal V})$.
\end{prop}
\begin{proof} 
By (\ref{prop:bcmop}) and  (\ref{prop:canmorp}), 
$N_{\rm zar}$ or $\nu_{\rm zar}$ are compatible with the base change morphism of 
sequences of log $p$-adic formal schemes over 
families of log points. 
Hence we obtain (\ref{prop:convmon}). 
\end{proof}


\parno 
The log ($p$-adically) convergent version of 
(\ref{conj:rcpmc}) is the following:  

\begin{conj}[{\bf Log convergent monodromy-weight conjecture}]
\label{conj:ccpwmc}
Assume that $N=0$. 
Let $q$ be nonnegative integer. 
Then the induced morphism 
\begin{align*} 
\nu^k \col 
{\rm gr}^P_{q+k}R^qf_*({\cal O}_{X/K})
\otimes_{\mab Z}{\mab Q}  
\lo 
{\rm gr}^P_{q-k}R^qf_*({\cal O}_{X/K})(-k,u)
\otimes_{\mab Z}{\mab Q} 
\tag{5.2.18.1}\label{eqn:grmmppd}
\end{align*}
is an isomorphism in ${\rm Isoc}^{\sq}(S/{\cal V})$.  
\end{conj}

\parno 
Obviously the conjecture {\rm (\ref{conj:rcpmc})} implies {\rm (\ref{conj:ccpwmc})}. 


\begin{prop}\label{prop:cxst}  
If the conjecture {\rm (\ref{conj:rcpmc})} for 
an object $(T,p{\cal O}_T,[~])$ 
$(T\in {\rm Enl}^{\sq}_p(S/{\cal V}))$ holds and if 
$T'\lo T$ is a morphism in ${\rm Enl}^{\sq}_p(S/{\cal V})$, 
then the conjecture {\rm (\ref{conj:rcpmc})} for $(T',p{\cal O}_{T'},[~])$ holds. 
\end{prop} 
\begin{proof}  
We obtain the following by the filtered base change theorem 
in log crystalline cohomologies ((\ref{theo:bccange})) 
and the flatness of 
${\cal H}^q(P_kRf_{X_{\os{\circ}{T}{}'}/S(T')^{\nat}*}
({\cal O}_{X_{\os{\circ}{T}{}'}/S(T')^{\nat}})_{\mab Q})$ 
((\ref{coro:fenlt})): 
\begin{align*} 
&P_kR^qf_{X_{\os{\circ}{T}{}'_1}/S(T')^{\nat}*}
({\cal O}_{X_{\os{\circ}{T}{}'_1}/S(T')^{\nat}})_{\mab Q} 
\tag{5.2.19.1}\label{ali:offhfrgm}\\
&={\rm Im}(H^q(P_kRf_{X_{\os{\circ}{T}{}'_1}/S(T')^{\nat}*}
({\cal O}_{X_{\os{\circ}{T}{}'_1}/S(T')^{\nat}})_{\mab Q})
\lo R^qf_{X_{\os{\circ}{T}{}'_1}/S(T')^{\nat}*}
({\cal O}_{X_{\os{\circ}{T}{}'_1}/S(T')^{\nat}})_{\mab Q})
\\
& ={\rm Im}
(H^q({\cal O}_{T'}\otimes^L_{{\cal O}_T}
P_kRf_{X_{\os{\circ}{T}_1}/S(T)^{\nat}*}
({\cal O}_{X_{\os{\circ}{T}_1}/S(T)^{\nat}})_{\mab Q})\\
&\lo H^q({\cal O}_{T'}\otimes^L_{{\cal O}_T}Rf_{X_{\os{\circ}{T}_1}/S(T)^{\nat}*}
({\cal O}_{X_{\os{\circ}{T}_1}/S(T)^{\nat}})_{\mab Q}))
\\
& ={\rm Im}
(H^q({\cal K}_{T'}\otimes^L_{{\cal K}_T}
P_kRf_{X_{\os{\circ}{T}_1}/S(T)^{\nat}*}
({\cal O}_{X_{\os{\circ}{T}_1}/S(T)^{\nat}})_{\mab Q})\\
&\lo H^q({\cal K}_{T'}\otimes^L_{{\cal K}_T}Rf_{X_{\os{\circ}{T}_1}/S(T)^{\nat}*}
({\cal O}_{X_{\os{\circ}{T}_1}/S(T)^{\nat}})_{\mab Q}))\\
&=
{\rm Im}
(H^q({\cal K}_{T'}\otimes_{{\cal K}_T}
P_kRf_{X_{\os{\circ}{T}_1}/S(T)^{\nat}*}
({\cal O}_{X_{\os{\circ}{T}_1}/S(T)^{\nat}})_{\mab Q})\\
&\lo H^q({\cal K}_{T'}\otimes_{{\cal K}_T}Rf_{X_{\os{\circ}{T}_1}/S(T)^{\nat}*}
({\cal O}_{X_{\os{\circ}{T}_1}/S(T)^{\nat}})_{\mab Q}))\\
& =
{\cal K}_{T'}\otimes_{{\cal K}_T}
P_kR^qf_{X_{\os{\circ}{T}_1}/S(T)^{\nat}*}
({\cal O}_{X_{\os{\circ}{T}_1}/S(T)^{\nat}})_{\mab Q}. 
\end{align*} 
Because $(P_kR^qf_{X_{\os{\circ}{T}_1}/S(T)^{\nat}*}
({\cal O}_{X_{\os{\circ}{T}_1}/S(T)^{\nat}})_{\mab Q},P)$ is a filteredly flat 
${\cal K}_T$-module ((\ref{coro:fenlt})), 
\begin{align*} 
{\rm gr}^P_k(R^qf_{X_{\os{\circ}{T}_1}/S(T)^{\nat}*}
({\cal O}_{X_{\os{\circ}{T}_1}/S(T)^{\nat}})_{\mab Q})
\otimes_{{\cal K}_T}{\cal K}_{T'}=
{\rm gr}^P_k(R^qf_{X_{\os{\circ}{T}{}'_1}/S(T')^{\nat}*}
({\cal O}_{X_{\os{\circ}{T}{}'_1}/S(T')^{\nat}})_{\mab Q}).
\tag{5.2.19.2}\label{ali:oxttq} 
\end{align*}  
Hence, if 
the morphism 
\begin{align*} 
\nu^k_{\rm zar} \col 
{\rm gr}^P_{q+k}R^qf_{X_{\os{\circ}{T}_1}/S(T)^{\nat}*}
({\cal O}_{X_{\os{\circ}{T}_1}/S(T)^{\nat}})_{\mab Q}  
\lo 
{\rm gr}^P_{q-k}R^qf_{X_{\os{\circ}{T}_1}/S(T)^{\nat}*}
({\cal O}_{X_{\os{\circ}{T}_1}/S(T)^{\nat}})(-k,u)_{\mab Q} 
\tag{5.2.19.3}\label{ali:inftctxsx}
\end{align*}
is an isomorphism, 
then the morphism $(\ref{ali:inftctxsx})\otimes_{{\cal K}_T}{\cal K}_{T'}$ 
is an isomorphism, which is equal to the morphism 
\begin{align*} 
\nu^k_{\rm zar} \col 
{\rm gr}^P_{q+k}R^qf_{X_{\os{\circ}{T}{}'_1}/S(T')^{\nat}*}
({\cal O}_{X_{\os{\circ}{T}{}'_1}/S(T')^{\nat}})_{\mab Q} 
\lo 
{\rm gr}^P_{q-k}R^qf_{X_{\os{\circ}{T}{}'_1}/S(T')^{\nat}*}
({\cal O}_{X_{\os{\circ}{T}{}'_1}/S(T')^{\nat}})(-k,u)_{\mab Q}. 
\tag{5.2.19.4}\label{ali:incbcsx}
\end{align*}
\end{proof} 

\begin{prop}\label{prop:conl}
Assume that $N=0$. 
Then the morphism {\rm (\ref{eqn:fcpl})} 
prolongs to a morphism in $F{\textrm -}{\rm Isoc}^{\sq}(S/{\cal V})$.  
\end{prop}
\begin{proof} 
This is clear by the construction of the morphism (\ref{eqn:fcpl}). 
\end{proof}

\parno 
By (\ref{prop:conl}) the log convergent version of (\ref{conj:lhlc}) is the following:

\begin{conj}[{\bf Filtered log convergent log hard Lefschetz conjecture}]
\label{conj:clhvlc} 
Assume that $N=0$.  
Then the following cup product 
\begin{equation*} 
\eta^i \col R^{d-i}f_*({\cal O}_{X/K})\lo R^{d+i}f_*({\cal O}_{X/K})(i) \quad (i\in {\mab N})
\tag{5.2.21.1}\label{eqn:fcicpl} 
\end{equation*}
is an isomorphism in $F{\textrm -}{\rm Isoc}^{\sq}(S/{\cal V})$.   
In fact, $\eta^i$ is the following isomorphism of filtered sheaves: 
\begin{equation*} 
\eta^i \col 
(R^{d-i}f_*({\cal O}_{X/K}),P) 
\os{\sim}{\lo} (R^{d+i}f_*({\cal O}_{X/K})(i),P). 
\tag{5.2.21.2}\label{eqn:filqiopl} 
\end{equation*}
Here $P_k(R^{d+i}f_*({\cal O}_{X/K})(i)):=P_{k+2i}R^{d+i}f_*({\cal O}_{X/K})$. 
\end{conj}

\begin{prop}\label{prop:oxlp} 
If the conjecture {\rm (\ref{conj:lhlc})} for 
an object $(T,p{\cal O}_T,[~])$ 
$(T\in {\rm Enl}^{\sq}_p(S/{\cal V}))$ holds and if 
$T'\lo T$ is a morphism in ${\rm Enl}^{\sq}_p(S/{\cal V})$, 
then the conjecture {\rm (\ref{conj:lhlc})} for $(T',p{\cal O}_{T'},[~])$ holds. 
\end{prop} 
\begin{proof} 
Because the proof of (\ref{prop:oxlp}) is easier than that of (\ref{prop:cxst}), 
we leave the proof to the reader.    
\end{proof}

\bigskip
\parno
{\bf (4) Duality of weight spectral sequences.}
\bigskip
\parno

\begin{prop}\label{prop:e2tmc} 
The $E_2$-terms of the spectral sequence 
$(\ref{eqn:espssp})$ $($and hence $(\ref{eqn:escssp}))$ 
for the case $E^{\bul \leq N}
={\cal O}_{\os{\circ}{X}_{\bul \leq N}/\os{\circ}{S}}$ 
prolongs to an object of $F{\textrm -}{\rm Isoc}(\os{\circ}{S}/{\cal V})$.  
In particular, the $E_2$-terms are flat sheaves of 
${\cal O}_T\otimes_{\mab Z}{\mab Q}$-modules.
\end{prop} 
\begin{proof} 
By \cite[(2.18), (3.1), (3.7)]{of}, the $E_1$-terms of the spectral sequence 
$(\ref{eqn:espssp})$ for the case 
$E^{\bul \leq N}
={\cal O}_{\os{\circ}{X}_{\bul \leq N}/\os{\circ}{S}}$ 
prolongs to convergent $F$-isocrystals on $\os{\circ}{S}/{\cal V}$.  
By \cite[(3.13)]{of} we see that 
the Gysin morphism is a morphism of convergent $F$-isocrystals. 
Hence we see that 
the $E_2$-terms are 
convergent $F$-isocrystals on $\os{\circ}{S}/{\cal V}$
by (\ref{prop:exbd}) (and hence (\ref{prop:deccbd})) and by \cite[(2.10)]{of}. 
\end{proof}

\par
Let us recall the following, 
which is the correction of the duality in \cite[4.15]{msemi}:

\begin{prop}[{\bf \cite[(10.5)]{ndw}}]\label{prop:dual} 
Let $s$ be as in {\rm \S\ref{sec:flgdw}} 
and let $X$ be a proper SNCL scheme over $s$. 
Assume that $\os{\circ}{X}$ is of pure dimension $d$. 
Then the following hold$:$
\par
$(1)$ Let $n$ be a positive integer. 
Let $\{E^{\bul \bul}_{r,n}\}_{r \geq 1}$ 
be the $E_r$-terms of the preweight spectral sequence 
{\rm (\ref{eqn:trhkwsp})} for the case $\star=n$, $N=0$ and 
$E^0={\cal O}_{\os{\circ}{X}/{\cal W}_n(\os{\circ}{s})}$. 
Then the Poincar\'{e} duality pairing 
\begin{equation*}
\langle~,~\rangle \col 
E_{1,n}^{-k,2d-h-k}\otimes_{{\cal W}_n}
E_{1,n}^{k,h+k}  
\lo {\cal W}_n(-d) 
\tag{5.2.24.1}\label{eqn:duale1}
\end{equation*}
induces the following perfect pairing 
\begin{equation*}
\langle~,~\rangle \col 
E_{2,n}^{-k,2d-h-k} 
\otimes_{{\cal W}_n}
E_{2,n}^{k,h+k} 
\lo {\cal W}_n(-d).  
\tag{5.2.24.2}\label{eqn:duale2}
\end{equation*}
\par
$(2)$ Holds the analogue of $(1)$ 
for the weight spectral sequence 
$(\ref{eqn:trhkwsp})\otimes_{\cal W}K_0$ 
for the case $\star=$nothing, $N=0$ and 
$E^0={\cal O}_{\os{\circ}{X}/{\cal W}(\os{\circ}{s})}$. 
\end{prop}

\begin{rema}\label{rema:ctrdu}
(1) The statement \cite[4.15]{msemi} is mistaken. 
See \cite[(10.6)]{ndw} for details. 
\par 
(2) Let the notations be as in (\ref{prop:dual}). 
Then, by the same proof as that of \cite[(10.5)]{ndw}, 
we obtain the analogous duality for the weight spectral sequence 
(\ref{eqn:escssp}) directly without using [loc.~cit] nor (\ref{theo:csoncrdw}). 
\end{rema}

\begin{prop}\label{prop:due} 
Let $T$ be a $p$-adic formal family of log points over ${\cal V}$. 
Let $Y/T_1$ be a proper SNCL scheme. 
Assume that $Y^{(0)}$ is projective over $\os{\circ}{T}$. 
Let $\{E^{\bul \bul}_{r}\}_{r \geq 1}$ 
be the $E_r$-terms of the $($pre$)$weight spectral sequence 
{\rm (\ref{eqn:escssp})} for the case $\star=$nothing, $N=0$ and 
$E^0={\cal O}_{\os{\circ}{Y}/\os{\circ}{T}}$. Then $\{E^{\bul \bul}_{1}\}$ and 
$\{E^{\bul \bul}_{2}\}$ prolong to objects 
$E^{\bul \bul}_1(Y/K)$ and $E^{\bul \bul}_2(Y/K)$ of 
$F{\textrm -}{\rm Isoc}(\os{\circ}{T}/{\cal V})$ 
such that 
the Poincar\'{e} duality perfect pairing $(${\rm \cite[(3.12)]{of}}$)$
\begin{equation*}
\langle~,~\rangle \col 
E_1^{-k,2d-h-k}(Y/K)\otimes_{{\cal O}_{T/K}}E_1^{k,h+k}(Y/K)  
\lo {\cal O}_{T/K}(-d) 
\tag{5.2.26.1}\label{eqn:dua11}
\end{equation*}
induces the following perfect pairing 
\begin{equation*}
\langle~,~\rangle \col 
E_2^{-k,2d-h-k}(Y/K)
\otimes_{{\cal O}_{T/K}}E_2^{k,h+k}(Y/K) \lo {\cal O}_{T/K}(-d).  
\tag{5.2.26.2}\label{eqn:dua22}
\end{equation*}
\end{prop}
\begin{proof} 
This follows from \cite[(3.12), (3.13)]{of}, (\ref{prop:deccbd}) 
and the proof of \cite[(10.5)]{ndw}
\end{proof}

\begin{theo}[{\bf Filtered log Berthelot-Ogus isomorphism}]\label{theo:bofis}
Let the notations and the assumptions be as in {\rm (\ref{theo:pwfaec})}.  
Let $T$ be an object of ${\rm Enl}^{\sq}(S/{\cal V})$. 
Let $T_0\lo S$ be the structural morphism. 
Let $f_0 \col X_{\bul \leq N,\os{\circ}{T}_0}\lo T_0$ 
be the structural morphism.  
If there exists an SNCL lift 
$f_1 \col (X_{\bul \leq N})_1 \lo T_1$ of $f_0$, 
then there exists the following canonical 
filtered isomorphism
\begin{equation*} 
(R^qf_*({\cal O}_{X_{\bul \leq N}/K})^{\nat}_T,P)
\os{\sim}{\lo}
(R^qf_{X_{\bul \leq N,\os{\circ}{T}_1}/S(T)^{\nat}*} 
({\cal O}_{X_{\bul \leq N,\os{\circ}{T}_1}/S(T)^{\nat}}),P).
\tag{5.2.27.1}\label{eqn:xntp}
\end{equation*} 
\end{theo}
\begin{proof} 
This follows from the proof of \cite[(3.8)]{of} and that of (\ref{theo:pwfaec}). 
\end{proof} 

\section{Infinitesimal deformation invariance of iso-zarisikian 
$p$-adic filtered Steenbrink complexes}\label{sec:infdi}  
In this section we prove the infinitesimal deformation invariance 
of the pull-back of a morphism of truncated simplicial base changes of 
SNCL schemes with admissible immersions 
on iso-zariskian $p$-adic filtered Steenbrink complexes. 
As a corollary we obtain the infinitesimal deformation invariance 
of the pull-back of the morphism on log isocrystalline cohomologies with weight filtrations.  
These are nontrivial filtered log versions of the infinitesimal deformation invariances  
for isocrystalline cohomological complexes and isocrystalline cohomologies in \cite{boi}. 
As in [loc.~cit.], to prove the invariance, we use Dwork's idea for enlarging 
the radius of convergence of log $F$-isocrystals by the use of the relative Frobenius.  
Strictly speaking, in our case, 
we use the base change of the truncated simplicial base changes 
by the iteration of the {\it abrelative} Frobenius morphism of the base scheme. 
The notion of the truncated simplicial base change of SNCL schemes 
and admissible immersions defined in \S\ref{sec:bcsncl} 
gives us an appropriate framework. 
\par 
For a complex $K$ in a filtered derived category (\cite{nh2}), 
denote $K\otimes_{\mab Z}^L{\mab Q}$ by  $K_{\mab Q}$ 
for simplicity of notation.

\begin{lemm}\label{lemm:invlem} 
Let $\star$ be nothing  or $\prime$. 
Let $T^{\star}$ be a fine log $p$-adic formal scheme.  
Assume that 
$(\os{\circ}{T}{}^{\star},{\cal J}^{\star},\del^{\star})$ 
is a $p$-adic formal PD-scheme in the sense of {\rm \cite[(7.17)]{bob}}. 
Set $T^{\star}_0:=\ul{\rm Spec}^{\log}_{T^{\star}}({\cal O}_{T^{\star}}/{\cal J}^{\star})$ 
$($note that $T^{\star}_0$ is of characteristic $p)$. 
Let $f^{\star} \col Y^{\star}_{\bul \leq N} \lo T^{\star}_0$ 
be a log smooth $N$-truncated simplicial log scheme over $T^{\star}_0$. 
Set $T^{\star}_0:=\ul{\rm Spec}_{T^{\star}}^{\log}({\cal O}_{T^{\star}}/{\cal J}^{\star})$.  
Assume that $\os{\circ}{Y}{}^{\star}_m$ $(0\leq m\leq N)$ 
and $\os{\circ}{T}{}^{\star}_0$ are quasi-compact. 
Let $\iota^{\star} \col T^{\star}_0(0) \os{\subset}{\lo} T^{\star}_0$  
be an exact closed nilpotent immersion. 
Set $Y^{\star}_{\bul \leq N}(0)
:=Y^{\star}_{\bul \leq N}\times_{T^{\star}_0}T^{\star}_0(0)$.   
Let 
\begin{align*} 
g_{0\bul \leq N}\col Y_{\bul \leq N}(0) \lo Y'_{\bul \leq N}(0)
\tag{5.3.1.1}\label{eqn:lntyvn}
\end{align*}   
be a morphism of log schemes over $T\lo T'$. 
Identify $(Y^{\star}_{\bul \leq N})_{\rm zar}$ with 
$(Y^{\star}_{\bul \leq N}(0))_{\rm zar}$. 
Then the following hold$:$ 
\par 
$(1)$ There exists a canonical morphism
\begin{align*}
g^*_{0\bul \leq N} &: 
Ru_{Y'_{\bul \leq N}/T'*}({\cal O}_{Y'_{\bul \leq N}/T'})_{\mab Q}
\lo Rg_{0\bul \leq N*}(Ru_{Y_{\bul \leq N}/T*}({\cal O}_{Y_{\bul \leq N}/T})_{\mab Q}). 
\tag{5.3.1.2}\label{eqn:ldeyfinvn}
\end{align*}  
\par 
$(2)$ Let $g'_{0\bul \leq N}\col Y'_{\bul \leq N}(0)\lo Y''_{\bul \leq N}(0)$ be 
a similar morphism to $g_{0\bul \leq N}$. 
Then 
\begin{align*} 
(g'_{0\bul \leq N}\circ g_{0\bul \leq N})^*
=Rg'_{0\bul \leq N*}(g^*_{0\bul \leq N})\circ g'{}^*_{\!\!0\bul \leq N}. 
\tag{5.3.1.3}\label{eqn:ldggnvn}
\end{align*}
\par 
$(3)$ 
\begin{align*}
{\rm id}_{Y_{\bul \leq N}(0)}^*=
{\rm id}_{Ru_{Y_{\bul \leq N}/T*}({\cal O}_{Y_{\bul \leq N}/T})_{\mab Q}}.
\tag{5.3.1.4}\label{eqn:dlbyyyvn}
\end{align*}
\par 
$(4)$ 
If $g_{0\bul \leq N}$ has a lift $g_{\bul \leq N} \col Y_{\bul \leq N}\lo Y'_{\bul \leq N}$ 
over $T\lo T'$, then $g^*_{0\bul \leq N}$ is equal to 
$$g^*_{\bul \leq N,{\rm crys}}\otimes^L{\rm id}_{\mab Q}
\col Ru_{Y'_{\bul \leq N}/T'*}({\cal O}_{Y'_{\bul \leq N}/T'})_{\mab Q}
\lo Rg_{\bul \leq N*}(Ru_{Y_{\bul \leq N}/T*}({\cal O}_{Y_{\bul \leq N}/T})_{\mab Q}).$$
\end{lemm} 
\begin{proof}
(By using \cite[(2.24)]{hk} (cf.~(\ref{rema:sot}) (2)), 
the proof of this lemma is only an imitation of \cite[(2.1)]{boi}.) 
\par
(1): Because 
$\os{\circ}{\iota}{}^{\star}\col 
\os{\circ}{T}{}^{\star}_0(0) \os{\sus}{\lo} \os{\circ}{T}{}^{\star}_0$ is nilpotent, 
we can apply {\rm (\ref{lemm:fianf})} (3) for the case $U(0)=T^{\star}_0(0)$ and $U=T^{\star}_0$. 
Hence there exists a morphism 
$\rho^{(n)}{}^{\star}\col T^{\star}_0\lo T^{\star}_0(0)$ 
such that $\rho^{(n)}{}^{\star} \circ \iota^{\star}=F^n_{T^{\star}_0(0)}$ 
and 
$\iota^{\star}\circ \rho^{(n)}{}^{\star} =F^n_{T^{\star}_0}$ for some positive integer $n$. 
Hence 
\begin{align*} 
Y^{\star}{}^{\{p^n\}}_{\! \!\!\bul \leq N}
:=Y^{\star}_{\bul \leq N}
\times_{T^{\star}_0,F^n_{T^{\star}_0}}T^{\star}_0
=Y^{\star}_{\bul \leq N}(0)
\times_{T^{\star}_0(0),\rho^{(n)}{}^{\star}}T^{\star}_0.
\tag{5.3.1.5}\label{eqn:ldyayvn}
\end{align*}   
Because we are given the morphism 
$g_{0\bul \leq N} \col Y_{\bul \leq N}(0)\lo Y'_{\bul \leq N}(0)$, 
we have the following base change morphism  
\begin{align*} 
g^{\{p^n\}}_{0\bul \leq N}\col Y^{\{p^n\}}_{\bul \leq N}\lo Y'^{\{p^n\}}_{\bul \leq N}
\tag{5.3.1.6}\label{eqn:ldyuyvn}
\end{align*}  
of $g_{0\bul \leq N}$ over $T\lo T'$ by (\ref{eqn:ldyayvn}). 
Consequently we have the pull-back morphism by the functoriality: 
\begin{align*}  
g^{\{p^n\}*}_{0\bul \leq N}\col Ru_{Y'^{\{p^n\}}_{\bul \leq N}/T'*}
({\cal O}_{Y'^{\{p^n\}}_{\bul \leq N}/T'})\lo 
Rg^{\{p^n\}}_{0\bul \leq N*}
(Ru_{Y^{\{p^n\}}_{\bul \leq N}/T*}({\cal O}_{Y^{\{p^n\}}_{\bul \leq N}/T})). 
\tag{5.3.1.7}\label{ali:yxylpn}
\end{align*}  
By \cite[(2.24)]{hk} (cf.~(\ref{rema:sot}) (2)) 
the iteration of the relative Frobenius morphism
$$F^n_{Y^{\star}_{\bul \leq N}/T^{\star}}\col 
Y^{\star}_{\bul \leq N}\lo Y^{\star}{}^{\{p^n\}}_{\! \! \!\bul \leq N}$$ 
over $T^{\star}$ 
induces an isomorphism
\begin{align*} 
F^{n*}_{Y^{\star}_{\bul \leq N}/T^{\star}} 
\col 
& Ru_{Y^{\star}{}^{\{p^n\}}_{\! \! \!\bul \leq N}/T^{\star}*}
({\cal O}_{Y^{\star}{}^{\{p^n\}}_{\! \! \!\bul \leq N}/T^{\star}})_{\mab Q} 
\os{\sim}{\lo} 
Ru_{Y^{\star}_{\bul \leq N}/T*}
({\cal O}_{Y^{\star}_{\bul \leq N}/T^{\star}})_{\mab Q}. 
\tag{5.3.1.8}\label{ali:yxylqn}
\end{align*} 
Here, as in [loc.~cit.], we identify 
$(Y^{\star}_{\bul \leq N})_{\rm zar}$ with 
$(Y^{\star}{}^{\{p^n\}}_{\! \! \!\bul \leq N})_{\rm zar}$. 
Consequently we obtain the following diagram: 
\begin{equation*} 
\begin{CD} 
Ru_{Y'{}^{\{p^n\}}_{\! \! \!\bul \leq N}/T'*}
({\cal O}_{Y'{}^{\{p^n\}}_{\! \! \!\bul \leq N}/T'})_{\mab Q}
@>{g^{\{p^n\}*}_{0\bul \leq N}}>> Ru_{Y^{\{p^n\}}_{\bul \leq N}/T*}
({\cal O}_{Y^{\{p^n\}}_{\bul \leq N}/T})_{\mab Q}\\
@V{F^{n*}_{Y'_{\bul \leq N}/T'}}V{\simeq}V 
@V{\simeq}V{F^{n*}_{Y_{\bul \leq N}/T}}V \\
Ru_{Y'_{\bul \leq N}/T'*}({\cal O}_{Y'_{\bul \leq N}/T'})_{\mab Q} 
@.Ru_{Y_{\bul \leq N}/T*}({\cal O}_{Y_{\bul \leq N}/T})_{\mab Q}. 
\end{CD} 
\tag{5.3.1.9}\label{cd:xuywuy}
\end{equation*}
The morphism $g^*_{0\bul \leq N}$ in (\ref{eqn:ldeyfinvn}),  
is, by definition, the morphism 
$$Ru_{Y^{\{p^n\}}_{\bul \leq N}/T*}
({\cal O}_{Y^{\{p^n\}}_{\bul \leq N}/T})_{\mab Q}\lo 
Ru_{Y_{\bul \leq N}/T*}({\cal O}_{Y_{\bul \leq N}/T})_{\mab Q}$$  
making the resulting diagram in (\ref{cd:xuywuy}) commutative. 
This morphism is independent of the choice of $n$ by 
the transitive relation of the pull-back morphism of log crystalline complexes. 
\par 
(2) By the transitive relation of the pull-back morphism again, 
$g^*_{0\bul \leq N}$ and $g'{}^*_{0\bul \leq N}$ are compatible 
with the composition of $g_{0\bul \leq N}$ and $g'_{0\bul \leq N}$. 
\par 
(3) 
The formula (\ref{eqn:dlbyyyvn}) immediately follows 
from the definition of $g^*_{0\bul \leq N}$ in (\ref{cd:xuywuy}). 
\par 
(4) 
Assume that $g_{0\bul \leq N}$ has a lift 
$g_{\bul \leq N} \col Y_{\bul \leq N} \lo Y'_{\bul \leq N}$ 
over $T\lo T'$. 
Then we have the relation 
$g^*_{0\bul \leq N}=g^*_{\bul \leq N,\rm crys}\otimes^L{\rm id}_{\mab Q}$ 
by (\ref{cd:xuywuy}) and the following commutative diagram 
\begin{equation*} 
\begin{CD} 
Y_{\bul \leq N}@>{g_{\bul \leq N}}>> Y'_{\bul \leq N}\\
@V{F^n_{Y_{\bul \leq N}/T'}}VV 
@VV{F^n_{Y'_{\bul \leq N}/T}}V \\
Y^{\{p^n\}}_{\bul \leq N}@>{g^{\{p^n\}}_{0\bul \leq N}}>> Y'{}^{\{p^n\}}_{\bul \leq N}. 
\end{CD} 
\end{equation*}
\end{proof}

\begin{coro}\label{coro:finvce}
If $T'=T$ and $Y_{\bul \leq N}(0)=Y'_{\bul \leq N}(0)$, 
then 
\begin{align*} 
Ru_{Y_{\bul \leq N}/T*}({\cal O}_{Y_{\bul \leq N}/T})_{\mab Q}
= Ru_{Y'_{\bul \leq N}/T'*}({\cal O}_{Y'_{\bul \leq N}/T'})_{\mab Q}. 
\tag{5.3.2.1}\label{ali:xdunz}
\end{align*} 
\end{coro} 

\parno 
By (\ref{lemm:flpis}) we obtain the following: 

\begin{coro}\label{coro:otf}
Let $(T^{\star},{\cal J}^{\star},\del^{\star})$, $T^{\star}_0$ and $T^{\star}_0(0)$, 
$Y^{\star}_{\bul \leq N}$, $Y^{\star}_{\bul \leq N}(0)$ and $g_{0\bul \leq N}$ 
be as in {\rm (\ref{lemm:invlem})}. 
Let 
\begin{equation*} 
\begin{CD}
(U,{\cal K},\eps)@>>> (U',{\cal K}',\eps')\\
@VVV @VVV \\
(T,{\cal J},\del)@>>> (T',{\cal J}',\del')
\end{CD}
\end{equation*}
be a commutative diagram of fine log $p$-adic formal PD-schemes. 
Set $U^{\star}_0:=\ul{\rm Spec}^{\log}_{U^{\star}}({\cal O}_{U^{\star}}/{\cal K}^{\star})$ 
and $U^{\star}_0(0):=U^{\star}_0\times_{T^{\star}_0}T^{\star}_0(0)$.
Set $Y^{\star}_{\bul \leq N,U^{\star}_0}
:=Y_{\bul \leq N}\times_{T^{\star}_0}U^{\star}_0$ 
and 
$Y^{\star}_{\bul \leq N,U^{\star}_0}(0):=
Y_{\bul \leq N}(0)\times_{T^{\star}_0(0)}U^{\star}_0(0)$.  
Let 
\begin{align*} 
g_{0\bul \leq N,U^{\star}_0(0)} \col Y_{\bul \leq N,U_0}(0) \lo Y'_{\bul \leq N,U'_0}(0)
\tag{5.3.3.1}\label{eqn:lnyvun}
\end{align*}   
be the base change morphism of $g_{0\bul \leq N}$. 
Assume that $\os{\circ}{U}{}^{\star}=\os{\circ}{T}{}^{\star}$. 
Then the following diagram is commutative$:$  
\begin{equation*} 
\begin{CD}
Ru_{Y'_{\bul \leq N,U'_0}/U'*}({\cal O}_{Y'_{\bul \leq N,U'_0}/U'})_{\mab Q}
@>{g^*_{0\bul \leq N,U^{\star}_0(0)}}>> 
Rg_{0\bul \leq N,U_0(0)*}(Ru_{Y_{\bul \leq N,U_0}/U*}({\cal O}_{Y_{\bul \leq N,U_0}/U})_{\mab Q}) \\
@| @|\\
Ru_{Y'_{\bul \leq N}/T'*}({\cal O}_{Y'_{\bul \leq N,T'_0}/T'})_{\mab Q}
@>{g^*_{0\bul \leq N}}>> 
Rg_{0\bul \leq N*}(Ru_{Y_{\bul \leq N}/T*}({\cal O}_{Y_{\bul \leq N}/T})_{\mab Q}). 
\end{CD} 
\tag{5.3.3.2}\label{cd:ouq}
\end{equation*} 
\end{coro}
\begin{proof}
Obvious.  
\end{proof}

\par 
The lemma (\ref{lemm:invlem}) itself is not sufficient for obtaining the main result 
(\ref{theo:definv}) below in this section because we have to consider 
the abrelative Frobenius instead of the usual relative Frobenius to obtain it. 
We need the following corollary: 

\begin{lemm}\label{lemm:som}
Let $(T^{\star},{\cal J}^{\star},\del^{\star})$, $T^{\star}_0$ and  
$\iota^{\star}$ be as in {\rm (\ref{lemm:invlem})}.  
Let $S^{\star}$ be a family of log points. 
Assume that $S^{\star}$ is of characteristic $p>0$. 
Set $S^{{\star}[p^n]}:=S^{\star}\times_{\os{\circ}{S}{}^{\star},\os{\circ}{F}{}^n_{S^{\star}}}
\os{\circ}{S}{}^{\star}$ {\rm ((\ref{defi:rwd}))}. 
Assume that $(T^{\star},{\cal J}^{\star},\del^{\star})$ is 
a log $p$-adic formal PD-enlargement of $S^{\star}$.  
Let $g^{\star} \col Z^{\star}_{\bul \leq N} \lo S^{\star}(T)^{\star}_0$ 
be a log smooth $N$-truncated simplicial log scheme over $S^{\star}(T)^{\star}_0$.   
Assume that $\os{\circ}{Z}{}^{\star}_m$ $(0\leq m\leq N)$ 
and $\os{\circ}{T}{}^{\star}_0$ are quasi-compact.  
Set $Z^{\star}_{\bul \leq N}(0)
:=Z^{\star}_{\bul \leq N}\times_{S^{\star}_{T^{\star}_0}}S^{\star}_{T^{\star}_0(0)}$.   
Let 
\begin{align*} 
g_{0\bul \leq N}\col Z_{\bul \leq N}(0) \lo Z'_{\bul \leq N}(0)
\tag{5.3.4.1}\label{eqn:lntzvn}
\end{align*}   
be a morphism of log schemes over $S(T)^{\nat}\lo S'(T')^{\nat}$. 
Identify $(Z^{\star}_{\bul \leq N})_{\rm zar}$ with 
$(Z^{\star}_{\bul \leq N}(0))_{\rm zar}$. 
Then the following hold$:$ 
\par 
$(1)$ There exists a canonical morphism
\begin{align*}
g^*_{0\bul \leq N} &: 
Ru_{Z'_{\bul \leq N}/S'(T')^{\nat}*}({\cal O}_{Z'_{\bul \leq N}/S'(T')^{\nat}})_{\mab Q}
\lo Rg_{0\bul \leq N*}(Ru_{Z_{\bul \leq N}/T*}({\cal O}_{Z_{\bul \leq N}/S(T)^{\nat}})_{\mab Q}). 
\tag{5.3.4.2}\label{eqn:ldzznvn}
\end{align*}  
\par 
$(2)$ Let $g'_{0\bul \leq N}\col Y'_{\bul \leq N}(0)\lo Y''_{\bul \leq N}(0)$ be 
a similar morphism to $g_{0\bul \leq N}$. 
Then 
\begin{align*} 
(g'_{0\bul \leq N}\circ g_{0\bul \leq N})^*
=Rg'_{0\bul \leq N*}(g^*_{0\bul \leq N})\circ g'{}^*_{\!\!0\bul \leq N}. 
\tag{5.3.4.3}\label{eqn:ldznvn}
\end{align*}
\par 
$(3)$ 
\begin{align*}
{\rm id}_{Z_{\bul \leq N}(0)}^*=
{\rm id}_{Ru_{Z_{\bul \leq N}/S(T)^{\nat}*}({\cal O}_{Z_{\bul \leq N}/S(T)^{\nat}})_{\mab Q}}.
\tag{5.3.4.4}\label{eqn:dlbzzyvn}
\end{align*}
\par 
$(4)$ 
If $g_{0\bul \leq N}$ has a lift $g_{\bul \leq N} \col Z_{\bul \leq N}\lo Z'_{\bul \leq N}$ 
over $S(T)^{\nat}\lo S'(T')^{\nat}$, then $g^*_{0\bul \leq N}$ is equal to 
$$g^*_{\bul \leq N,{\rm crys}}\otimes^L{\rm id}_{\mab Q}
\col Ru_{Z'_{\bul \leq N}/S'(T')^{\nat}*}({\cal O}_{Z'_{\bul \leq N}/S'(T')^{\nat}})_{\mab Q}
\lo Rg_{\bul \leq N*}(Ru_{Z_{\bul \leq N}/S(T)^{\nat}*}
({\cal O}_{Z_{\bul \leq N}/S(T)^{\nat}})_{\mab Q}).$$
\end{lemm} 
\begin{proof} 
(1), (2), (3), (4): The key ingredient of the proof 
of (\ref{lemm:invlem}) is the isomorphism (\ref{ali:yxylqn}). 
We can obtain this isomorphism by \cite[(2.24)]{hk} in which 
the theory of gauges in \cite{bob} has been used. 
Because the local lifted description of the induced morphism 
\begin{align*} 
F^{{\rm ar}*}_{Z^{\star}_{\bul \leq N}/S^{\star}(T^{\star})^{\nat},
S^{\star [p]}(T^{\star})^{\nat}} &
\col 
Ru_{Y^{\star}{}^{[p]}_{\! \! \!\bul \leq N}/S^{\star}{}^{[p]}(T^{\star})^{\nat}*}
({\cal O}_{Y^{\star}{}^{[p]}_{\! \! \!\bul \leq N}/S^{[p]}(T^{\star})^{\nat}})
\os{\sim}{\lo} 
Ru_{Y^{\star}_{\bul \leq N}/S^{\star}(T^{\star})*}
({\cal O}_{Y^{\star}_{\bul \leq N}/S^{\star}(T^{\star})})
\tag{5.3.4.5}\label{ali:kmcfd}\\
\end{align*} 
by the abrelative Frobenius morphism
$F^{\rm ar}_{Z^{\star}_{\bul \leq N}/S^{\star}(T^{\star})^{\nat},S^{\star}{}^{[p]}(T)^{\nat}} 
\col 
Z^{\star}_{\bul \leq N}\lo Z^{\star}{}^{[p]}_{\! \! \!\bul \leq N}$ 
over $S^{\star}(T^{\star})^{\nat}\lo S^{\star}{}^{[p]}(T^{\star})^{\nat}$ 
satisfies the axiom of (6.0.1)$\sim$(6.0.5) in \cite{ndw}, 
we can use the theory of gauge in this situation. 
Hence we can prove that the morphism 
(\ref{eqn:ldzznvn}) is an isomorphism. 
The rest of the proof of (\ref{lemm:som}) is the same as that of 
(\ref{lemm:invlem}); in the rest we give only the construction of 
the morphism (\ref{eqn:ldzznvn}) as follows.  
\par 
Because 
$\os{\circ}{\iota}{}^{\star}\col 
\os{\circ}{T}{}^{\star}_0(0) \os{\sus}{\lo} \os{\circ}{T}{}^{\star}_0$ is nilpotent, 
there exists a morphism 
$\os{\circ}{\rho}{}^{[n]}{}^{\star}\col \os{\circ}{T}{}^{\star}_0\lo \os{\circ}{T}{}^{\star}_0(0)$ 
such that $\rho^{(n)}{}^{\star} \circ \os{\circ}{\iota}{}^{\star}=F^n_{\os{\circ}{T}{}^{\star}_0(0)}$ 
and 
$\os{\circ}{\iota}{}^{\star}\circ \rho^{[n]}{}^{\star} =\os{\circ}{F}{}^n_{T^{\star}_0}$ 
for some positive integer $n$. 
Hence 
\begin{align*} 
Z^{\star}{}^{[p^n]}_{\! \!\!\bul \leq N}
:=Z^{\star}_{\bul \leq N}
\times_{\os{\circ}{T}{}^{\star}_0,\os{\circ}{F}{}^n_{T^{\star}_0}}\os{\circ}{T}{}^{\star}_0
=Z^{\star}_{\bul \leq N}(0)
\times_{\os{\circ}{T}{}^{\star}_0(0),\os{\circ}{\rho}{}^{[n]}{}^{\star}}\os{\circ}{T}{}^{\star}_0.
\tag{5.3.4.6}\label{eqn:ldypayvn}
\end{align*}   
Because we are given the morphism 
$g_{0\bul \leq N} \col Z_{\bul \leq N}(0)\lo Z'_{\bul \leq N}(0)$, 
we have the following base change morphism  
\begin{align*} 
g^{[p^n]}_{0\bul \leq N}\col Z^{[p^n]}_{\bul \leq N}\lo Z'^{[p^n]}_{\bul \leq N}
\tag{5.3.4.7}\label{eqn:ldvn}
\end{align*}  
of $g_{0\bul \leq N}$ over $S(T)^{\nat}\lo S'(T')^{\nat}$ by (\ref{eqn:ldyayvn}). 
Consequently we have the pull-back morphism by the functoriality: 
\begin{align*}  
g^{[p^n]*}_{0\bul \leq N}\col Ru_{Z'^{[p^n]}_{\bul \leq N}/S'(T')^{\nat}}
({\cal O}_{Z'^{[p^n]}_{\bul \leq N}/S'(T')^{\nat}})\lo 
Rg^{[p^n]}_{0\bul \leq N*}
(Ru_{Z^{[p^n]}_{\bul \leq N}/S(T)^{\nat}*}({\cal O}_{Z^{[p^n]}_{\bul \leq N}/S(T)^{\nat}})). 
\tag{5.3.4.8}\label{ali:yxaylnn}
\end{align*}  
The morphism (\ref{eqn:ldzznvn}) is, by definition, 
the morphism 
\begin{align*} 
Ru_{Z'_{\bul \leq N}/S'(T')^{\nat}*}({\cal O}_{Z'_{\bul \leq N}/S'(T')^{\nat}})_{\mab Q} 
\lo Ru_{Z_{\bul \leq N}/S(T)^{\nat}*}({\cal O}_{Z_{\bul \leq N}/S(T)^{\nat}})_{\mab Q}.
\end{align*} 
making the following diagram commutative:  
\begin{equation*} 
\begin{CD} 
Ru_{Z'{}^{[p^n]}_{\! \! \!\bul \leq N}/S'(T')^{\nat}*}
({\cal O}_{Z'{}^{[p^n]}_{\! \! \!\bul \leq N}/S'(T')^{\nat}})_{\mab Q}
@>{g^{[p^n]*}_{0\bul \leq N}}>> Ru_{Z^{[p^n]}_{\bul \leq N}/S(T)^{\nat}*}
({\cal O}_{Z^{[p^n]}_{\bul \leq N}/S(T)^{\nat}})_{\mab Q}\\
@V{F^{{\rm ar},n*}_{Z'_{\bul \leq N}/S'(T')^{\nat}{\rm crys}}}V{\simeq}V 
@V{\simeq}V{F^{{\rm ar},n*}_{Y_{\bul \leq N}/T{\rm crys}}}V \\
Ru_{Z'_{\bul \leq N}/S'(T')^{\nat}*}({\cal O}_{Z'_{\bul \leq N}/S'(T')^{\nat}})_{\mab Q} 
@.Ru_{Z_{\bul \leq N}/S(T)^{\nat}*}({\cal O}_{Z_{\bul \leq N}/S(T)^{\nat}})_{\mab Q}. 
\end{CD} 
\tag{5.3.4.9}\label{cd:xuybwuy}
\end{equation*}
We leave the detail of the rest of the proof to the reader. 
\end{proof} 

\begin{coro}\label{coro:dis} 
If $S'(T')^{\nat}=S(T)^{\nat}$ and $Z_{\bul \leq N}(0)=Z'_{\bul \leq N}(0)$, 
then 
\begin{align*} 
Ru_{Z_{\bul \leq N}/S(T)^{\nat}*}({\cal O}_{Z_{\bul \leq N}/S(T)^{\nat}})_{\mab Q}
= Ru_{Z'_{\bul \leq N}/S'(T')^{\nat}*}({\cal O}_{Z'_{\bul \leq N}/S'(T')^{\nat}})_{\mab Q}. 
\tag{5.3.5.1}\label{ali:xdsznz}
\end{align*}  
\end{coro}

\par 
The following is a main result in this section.  
To assume the existence of the admissible immersion (\ref{ali:tass}) 
is a point in this result. 

\begin{theo}[{\bf Infinitesimal deformation invariance of filtered iso-zarisikian 
$p$-adic Steenbrink complexes with respect to the pull-back of a morphism}]
\label{theo:definv}
Let $\star$ be nothing  or $\prime$. 
Let $n$ be a positive integer.  
Let $S^{\star}$ be a family of log points. 
Assume that $S^{\star}$ is of characteristic $p>0$. 
Set $S^{{\star}[p^n]}:=S^{\star}\times_{\os{\circ}{S}{}^{\star},\os{\circ}{F}{}^n_{S^{\star}}}
\os{\circ}{S}{}^{\star}$ {\rm ((\ref{defi:rwd}))}. 
Let $(T^{\star},{\cal J}^{\star},\del^{\star})$ be 
a log $p$-adic formal PD-enlargement of $S^{\star}$. 
Set $T^{\star}_0:=\ul{\rm Spec}^{\log}_{T^{\star}}({\cal O}_{T^{\star}}/{\cal J}^{\star})$. 
Let $f^{\star} \col X^{\star}_{\bul \leq N} \lo S^{\star}$ 
be an $N$-truncated simplicial base change of SNCL schemes over $S^{\star}$. 
Assume that $\os{\circ}{X}{}^{\star}_m$ $(0\leq m\leq N)$,  
$\os{\circ}{S}{}^{\star}$ and $\os{\circ}{T}{}^{\star}_0$ are quasi-compact. 
Let $\iota^{\star} \col T^{\star}_0(0) \os{\subset}{\lo} T^{\star}_0$  
be an exact closed nilpotent immersion.  
Set $S^{\star}_{\os{\circ}{T}{}^{\star}_{0}}
:=S^{\star}\times_{\os{\circ}{S}{}^{\star}}\os{\circ}{T}{}^{\star}_0$ 
and 
$S^{\star}_{\os{\circ}{T}{}^{\star}_0(0)}:=S^{\star}
\times_{\os{\circ}{S}{}^{\star}}\os{\circ}{T}{}^{\star}_0(0)$ 
and $X^{\star}_{\bul \leq N,\os{\circ}{T}_0}:=
X^{\star}_{\bul \leq N}\times_{S^{\star}}S^{\star}_{\os{\circ}{T}{}^{\star}_0}$.  
$X^{\star}_{\bul \leq N,\bul,\os{\circ}{T}{}^{\star}_0}(0)
:=X^{\star}_{\bul \leq N,\bul}
\times_{S^{\star}}S^{\star}_{\os{\circ}{T}{}^{\star}_0(0)}$. 
Assume that $X^{\star}_{\bul \leq N,\os{\circ}{T}_0}$ has the disjoint union of 
the members of an affine simplicial open covering of $X^{\star}_{\bul \leq N,\os{\circ}{T}_0}$. 
Let $X^{\star}_{\bul \leq N,\bul,\os{\circ}{T}_0}$ be the \v{C}ech diagram of 
this affine simplicial open covering. 
Set $X^{\star}_{\bul \leq N,\bul,\os{\circ}{T}{}^{\star}_0}(0)
:=X^{\star}_{\bul \leq N,\bul}\times_{S^{\star}}S^{\star}_{\os{\circ}{T}{}^{\star}_0(0)}$. 
Assume that there exists an admissible immersion 
\begin{align*} 
X^{\star}_{\bul \leq N,\bul,\os{\circ}{T}_0}\os{\sus}{\lo} 
{\cal P}^{\star}_{\bul \leq N,\bul}
\tag{5.3.6.1}\label{ali:taxss}
\end{align*}
over $S^{\star}(T)^{\nat}$. 
Set $X^{\star}{}^{[p^n]}_{\! \! \! \bul \leq N}:=
X^{\star}_{\bul \leq N}\times_{S^{\star}}S^{{\star}[p^n]}$ 
and $X^{\star}{}^{[p^n]}_{\! \! \! \bul \leq N,\os{\circ}{T}_0}:=
X^{\star}{}^{[p^n]}
\times_{S^{{\star}[p^n]}}S^{{\star}[p^n]}_{\os{\circ}{T}_0}$. 
$($Note that the underlying schemes of 
$X^{\star}{}^{[p^n]}_{\! \! \! \bul \leq N,\os{\circ}{T}_0}$
and $X^{\star}_{\bul \leq N,\os{\circ}{T}_0}$ are the same and 
that we have the log scheme 
$S^{{\star}[p^n]}_{\os{\circ}{T}_0}$ by using the composite morphism 
$T_0\lo S\lo S^{[p^n]}.)$
Let $X^{\star}{}^{[p^n]}_{\! \! \! \bul \leq N,\bul,\os{\circ}{T}_0}$ 
be the \v{C}ech diagram of the 
disjoint union of the members of 
an affine simplicial open covering of $X^{\star}{}^{[p^n]}_{\! \! \! \bul \leq N,\os{\circ}{T}_0}$ 
obtained by that of $X^{\star}_{\bul \leq N,\os{\circ}{T}_0}$. 
Let $n$ be a positive integer such that the pull-back morphism 
$F^{n*}_{T^{\star}_0}\col 
{\cal O}_{T^{\star}_0}\lo {\cal O}_{T^{\star}_0}$ 
kills ${\rm Ker}({\cal O}_{T^{\star}_0}\lo {\cal O}_{T^{\star}_0})$.
Assume that there exists an admissible immersion 
\begin{align*} 
X^{{\star}[p^n]}_{\bul \leq N,\bul,\os{\circ}{T}_0}
\os{\sus}{\lo} {\cal Q}^{\star}_{\bul \leq N,\bul}
\tag{5.3.6.2}\label{ali:tass}
\end{align*}
over $S^{\star[p^n]}(T^{\star})$.   
Here $S^{\star[p^n]}(T^{\star})$ is defined by 
the following composite morphism 
$T_1\lo S^{\star}\lo S^{\star[p]}$, where 
the morphism $S^{\star}\lo S^{\star[p]}$ is 
the abrelative Frobenius morphism of $S^{\star}$. 
Let 
\begin{align*} 
g_{0\bul \leq N} \col X_{\bul \leq N,\os{\circ}{T}_0}(0) 
\lo X'_{\bul \leq N,\os{\circ}{T}{}'_0}(0)
\tag{5.3.6.3}\label{eqn:ldehtvn}
\end{align*}   
be a morphism of log schemes over $S(T)^{\nat}\lo S'(T')^{\nat}$. 
Then the following hold$:$ 
\par 
$(1)$ There exist a canonical filtered morphism
\begin{align*}
g^*_{0\bul \leq N} &: 
(A_{{\rm zar},{\mab Q}}(X'_{\bul \leq N,\os{\circ}{T}{}'_0}/S'(T')^{\nat}),P) 
\lo 
Rg_{0\bul \leq N*}((A_{{\rm zar},{\mab Q}}(X_{\bul \leq N,\os{\circ}{T}_0}/S(T)^{\nat}),P)). 
\tag{5.3.6.4}\label{eqn:ldefinvn}
\end{align*}  
and a canonical isomorphism 
\begin{align*} 
Ru_{X^{\star}_{\bul \leq N,\os{\circ}{T}{}^{\star}_0}/S^{\star}(T^{\star})^{\nat}*}
({\cal O}_{X^{\star}_{\bul \leq N,\os{\circ}{T}{}^{\star}_0}/S^{\star}(T^{\star})^{\nat}})_{\mab Q} 
\os{\sim}{\lo} A_{{\rm zar},{\mab Q}}(X^{\star}_{\bul \leq N,\os{\circ}{T}{}^{\star}_0}/
S^{\star}(T^{\star})^{\nat})
\tag{5.3.6.5}\label{cd:ldttvn}
\end{align*} 
fitting into the following commutative diagram
\begin{equation*} 
\begin{CD}
A_{{\rm zar},{\mab Q}}(X'_{\bul \leq N,\os{\circ}{T}{}'_0}/S'(T')^{\nat})@>{g^*_{0\bul \leq N}}>>
Rg_{0\bul \leq N*}(A_{{\rm zar},{\mab Q}}(X_{\bul \leq N,\os{\circ}{T}_0}/S(T)^{\nat}))\\
@A{\simeq}AA @AA{\simeq}A \\
Ru_{X'_{\bul \leq N,\os{\circ}{T}{}'_0}/S'(T')^{\nat}*}
({\cal O}_{X'_{\bul \leq N,\os{\circ}{T}{}'_0}/S'(T')^{\nat}})_{\mab Q} 
@>{g^*_{0\bul \leq N}}>>
Rg_{0\bul \leq N*}Ru_{X_{\bul \leq N,\os{\circ}{T}_0}/S(T)^{\nat}*}
({\cal O}_{X_{\bul \leq N,\os{\circ}{T}_0}/S(T)^{\nat}})_{\mab Q}. 
\end{CD}
\tag{5.3.6.6}\label{cd:ldtvn}
\end{equation*} 
\par 
$(2)$ Let $g'_{0\bul \leq N}\col X'_{\bul \leq N,\os{\circ}{T}{}'_0}(0)\lo 
X''_{\bul \leq N,\os{\circ}{T}{}''_0}(0)$ be 
a similar morphism to $g_{0\bul \leq N}$. 
Let $n'$ be a similar positive integer to $n$. 
Assume that $X''_{\bul \leq N,\os{\circ}{T}{}''_0}$ has the disjoint union of 
the members of an affine simplicial open covering of 
$X''_{\bul \leq N,\os{\circ}{T}{}''_0}$. 
Assume also that there exist admissible immersions 
\begin{align*} 
X''_{\bul \leq N,\bul,\os{\circ}{T}_0}
\os{\sus}{\lo} {\cal P}''_{\bul \leq N,\bul}
\tag{5.3.6.7}\label{ali:ppnb}
\end{align*}
over $S''(T'')^{\nat}$ and  
\begin{align*} 
X''{}^{[p^n]}_{\bul \leq N,\bul,\os{\circ}{T}{}''_0}
\os{\sus}{\lo} {\cal Q}''_{\bul \leq N,\bul}
\tag{5.3.6.8}\label{ali:ppnq}
\end{align*}
over $S''{}^{[p^n]}(T'')^{\nat}$. 
Then 
\begin{align*} 
(g'_{0\bul \leq N}\circ g_{0\bul \leq N})^*
=Rg'_{0\bul \leq N*}(g^*_{0\bul \leq N})\circ g'{}^*_{\! \!0\bul \leq N}. 
\tag{5.3.6.9}\label{eqn:ldfilnvn}
\end{align*}
\par 
$(3)$ 
\begin{align*}
{\rm id}^*_{X_{\bul \leq N,T^{\star}_0}(0)}=
{\rm id}_{(A_{{\rm zar},{\mab Q}}(X_{\bul \leq N,\os{\circ}{T}_0}/S(T)^{\nat})),P)}.
\tag{5.3.6.10}\label{eqn:ldeoxnvn}
\end{align*}
\par 
$(4)$ If $g_{0\bul \leq N}$ has a lift $g_{\bul \leq N} \col X_{\bul \leq N,\os{\circ}{T}_0}
\lo X'_{\bul \leq N,\os{\circ}{T}{}'_0}$ over 
$S_{\os{\circ}{T}_0}\lo S'_{\os{\circ}{T}{}'_0}$, 
then $g^*_{0\bul \leq N}$ is equal to the induced morphism by $g^*_{\bul \leq N}$ 
in $(\ref{eqn:fzgenaxd})$ 
in $E^{\bul \leq N}={\cal O}_{X_{\bul \leq N,\os{\circ}{T}{}_0}/S(T)^{\nat}}$
and 
$F^{\bul \leq N}={\cal O}_{Y_{\bul \leq N,\os{\circ}{T}{}'_0}/S'(T')^{\nat}}$.
\end{theo}
\begin{proof}
(1): Because the exact closed immersion 
$\os{\circ}{\iota}\col \os{\circ}{T}{}^{\star}_0(0) \os{\sus}{\lo} 
\os{\circ}{T}{}^{\star}_0$ is nilpotent, 
we see that there exists a morphism 
$\os{\circ}{\rho}{}^{(n)\star}\col \os{\circ}{T}{}^{\star}_0\lo \os{\circ}{T}{}^{\star}_0(0)$ 
such that $\os{\circ}{\rho}{}^{(n)\star} \circ \os{\circ}{\iota}=\os{\circ}{F}{}^n_{T^{\star}_0(0)}$ 
and $\os{\circ}{\iota}\circ \os{\circ}{\rho}{}^{(n)\star}
=\os{\circ}{F}{}^{n}_{T^{\star}_0}$. 
Hence 
\begin{align*} 
X^{\star}{}^{[p^n]}_{\! \! \! \bul \leq N,\os{\circ}{T}{}^{\star}_0}
&=X^{\star}_{\bul \leq N}\times_{\os{\circ}{S}{}^{\star},\os{\circ}{F}{}^n_{S^{\star}}}
\os{\circ}{S}{}^{\star}\times_{\os{\circ}{S}{}^{\star}}\os{\circ}{T}{}^{\star}_0
=X^{\star}_{\bul \leq N,\os{\circ}{T}{}^{\star}_0}
\times_{\os{\circ}{T}{}^{\star}_0,\os{\circ}{F}{}^n_{T^{\star}_0}}
\os{\circ}{T}{}^{\star}_0
=X^{\star}_{\bul \leq N,\os{\circ}{T}{}^{\star}_0}
\times_{\os{\circ}{T}{}^{\star}_0,\os{\circ}{\iota}\circ \os{\circ}{\rho}{}^{(n)\star}}
\os{\circ}{T}{}^{\star}_0\tag{5.3.6.11}\label{eqn:ldet0vn}\\
&=
X^{\star}_{\bul \leq N,\os{\circ}{T}{}^{\star}_0}(0)
\times_{\os{\circ}{T}{}^{\star}_0(0),\os{\circ}{\rho}{}^{(n)\star}}
\os{\circ}{T}{}^{\star}_0
\end{align*}
(cf.~(\ref{ali:u0m0r})). 
\par 
By the assumption (\ref{ali:tass}) we obtain the admissible immersion 
$X^{\star}{}^{[p^n]}_{\! \! \! \bul \leq N,\bul,\os{\circ}{T}{}^{\star}_0}
\os{\sus}{\lo} {\cal Q}^{\star}_{\bul \leq N,\bul}$
over $S^{\star[p^n]}(T)^{\nat}$. 
Consequently we obtain the filtered complex 
$$(A_{{\rm zar},{\mab Q}}
(X^{\star}{}^{[p^n]}_{\bul \leq N,\os{\circ}{T}_0}/S^{\star}{}^{[p^n]}(T^{\star})^{\nat}),P)$$   
((\ref{defi:wcthd})).  
Because we are given the morphism 
$g_{0\bul \leq N} \col X_{\bul \leq N,\os{\circ}{T}_0}(0)
\lo X'_{\bul \leq N,\os{\circ}{T}{}'_0}(0)$, 
we have the following base change morphism 
\begin{align*} 
g^{[p^n]}_{0\bul \leq N}
\col 
X^{[p^n]}_{\bul \leq N,\os{\circ}{T}_0}=
X_{\bul \leq N,\os{\circ}{T}_0}(0)
\times_{\os{\circ}{T}_0(0),\os{\circ}{\rho}{}^{(n)}}\os{\circ}{T}_0
\lo X'_{\bul \leq N,\os{\circ}{T}{}'_0}(0)
\times_{\os{\circ}{T}{}'_0(0),\os{\circ}{\rho}{}^{(n)'}}\os{\circ}{T}{}'_0
=X'{}^{[p^n]}_{\!\!\!\bul \leq N,\os{\circ}{T}{}'_0}. 
\end{align*} 
Hence 
we have the pull-back morphism by the functoriality of 
$(A_{{\rm zar},{\mab Q}},P)$ ((\ref{theo:fugennas})): 
\begin{align*}  
(A_{{\rm zar},{\mab Q}}(X'{}^{[p^n]}_{\! \! \!\bul \leq N,\os{\circ}{T}{}'_0}/S'(T')^{\nat}),P) \lo 
Rg^{[p^n]}_{0\bul \leq N*}
((A_{{\rm zar},{\mab Q}}(X^{[p^n]}_{\bul \leq N,\os{\circ}{T}_0}/S(T)^{\nat}),P)). 
\tag{5.3.6.12}\label{ali:yxlpn}
\end{align*}  
The abrelative Frobenius morphism 
$$F^{{\rm ar},n}_{X^{\star}_{\bul \leq N,\os{\circ}{T}{}^{\star}_0}/S^{\star}(T^{\star})^{\nat}} 
\col X^{\star}_{\bul \leq N,\os{\circ}{T}{}^{\star}_0}
\lo X^{\star [p^n]}_{\bul \leq N,\os{\circ}{T}{}^{\star}_0}$$ 
over the morphism $S^{\star}(T^{\star})^{\nat}\lo S^{\star[p^n]}(T^{\star})^{\nat}$
induces an isomorphism
\begin{align*} 
F^{{\rm ar},n*}_{X^{\star}_{\bul \leq N,
\os{\circ}{T}{}^{\star}_0}/S^{\star}(T^{\star})^{\nat}} 
\col &
(A_{{\rm zar},{\mab Q}}
(X^{\star}{}^{[p^n]}_{\! \! \!\bul \leq N,\os{\circ}{T}{}^{\star}_0}/S^{\star}(T^{\star})^{\nat}),P) 
\os{\sim}{\lo} \\
&(A_{{\rm zar},{\mab Q}}(X^{\star}_{\bul \leq N,\os{\circ}{T}{}^{\star}_0}/
S^{\star}(T^{\star})^{\nat}),P)
\quad (k\in {\mab Z}). 
\end{align*} 
by (\ref{prop:grlgoc}) and the proof of (\ref{theo:pwfaec}) and (\ref{exam:ofl})
because the base change over $\os{\circ}{T}{}^{\star}$ of 
the abrelative Frobenius morphism induces an isomorphism of  
the classical isocrystalline complex of a smooth scheme  
(\cite[(2.24)]{hk}, cf.~\cite[(1.3)]{boi}).
Here, as in \cite[(2.24)]{hk}, we identity 
$(X^{\star}{}^{[p^n]}_{\! \! \!\bul \leq N,\os{\circ}{T}{}^{\star}_0})_{\rm zar}$ 
with $(X^{\star}_{\bul \leq N,\os{\circ}{T}{}^{\star}_0})_{\rm zar}$ 
via the canonical equivalence.  
Consequently we obtain the following diagram: 
\begin{equation*} 
\begin{CD} 
(A_{{\rm zar},{\mab Q}}(X'{}^{[p^n]}_{\bul \leq N,\os{\circ}{T}{}'_0}/S'(T')^{\nat}),P)
@>{g^{[p^n]*}_{0\bul \leq N}}>> 
Rg^{[p^n]}_{0\bul \leq N*}((A_{{\rm zar},{\mab Q}}(X^{[p^n]}_{\bul \leq N,\os{\circ}{T}_0}/S(T)^{\nat}),P))\\
@V{F^{{\rm ar},n*}_{X'{}_{\bul \leq N,\os{\circ}{T}{}'_0}/S'(T')^{\nat}}}V{\simeq}V 
@V{\simeq}V{F^{{\rm ar},n*}_{X_{\bul \leq N,\os{\circ}{T}_0}/S(T)^{\nat}}}V \\
(A_{{\rm zar},{\mab Q}}(X'_{\bul \leq N,\os{\circ}{T}{}'_0}/S'(T')^{\nat}),P)
@.Rg_{0\bul \leq N*}(A_{{\rm zar},{\mab Q}}(X_{\bul \leq N,\os{\circ}{T}_0}/S(T)^{\nat}),P). 
\end{CD} 
\tag{5.3.6.13}\label{cd:xy}
\end{equation*}
This diagram gives us the morphism $g^*_{0\bul \leq N}$ in (\ref{eqn:ldefinvn}). 
This morphism is independent of the choice of $n$ by 
the transitive relation of the pull-back morphism ((\ref{ali:pdgenpp})). 
\par 
By (\ref{cd:pssfgenccz})  we obtain the following commutative diagram: 
\begin{equation*} 
\begin{CD}
A_{{\rm zar},{\mab Q}}(X'{}^{[p^n]}_{\bul \leq N,\os{\circ}{T}{}'_0}/S'(T')^{\nat})
@>{g^{[p^n]*}_{0\bul \leq N}}>>  \\ 
@A{\theta_{X'{}^{[p^n]}_{\bul \leq N,\os{\circ}{T}{}'_0}/S'(T')^{\nat}
/S'_0(T')^{\nat}} \wedge}A{\simeq}A \\
Ru_{X'{}^{[p^n]}_{\bul \leq N,\os{\circ}{T}{}'_0}/S'(T')^{\nat}*}
({\cal O}_{X'{}^{[p^n]}_{\bul \leq N,\os{\circ}{T}{}'_0}/S'(T')^{\nat}})_{\mab Q}
@>{g^{[p^n]*}_{0\bul \leq N}}>>
\end{CD}
\tag{5.3.6.14}\label{cd:psssnccz}
\end{equation*}
\begin{equation*} 
\begin{CD}
Rg^{[p^n]}_{0\bul \leq N*}(A_{{\rm zar},{\mab Q}}(X^{[p^n]}_{\bul \leq N,\os{\circ}{T}_0}/S(T)^{\nat}))\\ 
@A{Rg^{[p^n]}_{0\bul \leq N*}
(\theta_{X^{[p^n]}_{\bul \leq N,\os{\circ}{T}_0}/S(T)^{\nat}/S_0(T)^{\nat}}\wedge)}A{\simeq}A \\
Rg^{[p^n]}_{0\bul \leq N*}Ru_{X^{[p^n]}_{\bul \leq N,\os{\circ}{T}_0}/S(T)^{\nat}*}
({\cal O}_{X^{[p^n]}_{\bul \leq N,\os{\circ}{T}_0}/S(T)^{\nat}})_{\mab Q}.
\end{CD}
\end{equation*}
By (\ref{cd:xuybwuy}), (\ref{cd:xy}) and (\ref{cd:psssnccz}) 
we obtain the commutative diagram 
(\ref{cd:ldtvn}). 
\par 
(2): By the transitive relation of the pull-back morphism again, 
$g^*_{0\bul \leq N}$ and $g'{}^*_{\!\!0\bul \leq N}$ 
are compatible with the composition of 
$g_{0\bul \leq N}$ and $g'_{0\bul \leq N}$. 
\par 
(3): The formula (\ref{eqn:ldeoxnvn}) immediately follows from 
the definition of $g^*_{0\bul \leq N}$. 
\par 
(4): Assume that $g_{0\bul \leq N}$ has the lift $g_{\bul \leq N}$ in the statement of this theorem. 
Then we have the relation $g^*_{0\bul \leq N}=g^*_{{\bul \leq N}}$ by 
the following commutative diagram, the transitivity (\ref{ali:pdgenpp}) and 
the diagram (\ref{cd:xy}): 
\begin{equation*} 
\begin{CD} 
X_{\bul \leq N,\os{\circ}{T}_0}@>{g_{\bul \leq N}}>>X'_{\bul \leq N,\os{\circ}{T}{}'_0}\\
@V{F^{{\rm ar},n}_{X_{\bul \leq N,\os{\circ}{T}_0}/S(T)^{\nat}}}VV 
@VV{F^{{\rm ar},n}_{X'{}_{\bul \leq N,\os{\circ}{T}{}'_0}/S'(T')^{\nat}}}V \\
X^{[p^n]}_{\bul \leq N,\os{\circ}{T}_0}@>{g_{0\bul \leq N}^{[p^n]}}>>X^{'[p^n]}_{\bul \leq N,\os{\circ}{T}{}'_0}. 
\end{CD} 
\end{equation*}
\end{proof}

\begin{rema} 
(1) By the existence of the admissible immersion (\ref{ali:sstq}), 
we obtain the admissible immersion (\ref{ali:tass}) 
in the successive split truncated  simplicial SNCL scheme. 
It is clear that 
we also obtain the admissible immersion (\ref{ali:tass}) for 
$X_{\bul \leq N}/S$ in \S\ref{sec:psc}. 
\par 
(2) The reader should note that $X^{\{p^n\}}_{\bul \leq N,\bul,\os{\circ}{T}_0}=
X^{[p^n]}_{\bul \leq N,\bul,\os{\circ}{T}_0}\times_{S^{[p^n]}_{\os{\circ}{T}_0}}
S_{\os{\circ}{T}_0}$ does not necessarily have 
the admissible immersion even if 
$X^{[p^n]}_{\bul \leq N,\bul,\os{\circ}{T}_0}$ 
has the admissible immersion (\ref{ali:tass}). 
Indeed, consider the case where $X$ is an SNCL scheme over $S$ 
and $S=T_0$ and 
we are given an immersion $X\os{\sus}{\lo} {\cal P}$ into an SNCL scheme 
over $S(T)^{\nat}$. Then we have 
an immersion $X^{[p]}\os{\sus}{\lo} {\cal P}\times_{S(T)^{\nat}}S^{[p]}(T)^{\nat}$ 
into an SNCL scheme over $S^{[p]}(T)^{\nat}$. 
However, because there does not necessarily exist a morphism 
$T'\lo T$ fitting into the following commutative diagram 
\begin{equation*} 
\begin{CD} 
T'@>>> T\\
@A{\bigcup}AA @AA{\bigcup}A\\
T'_0@>>> T_0\\
@VVV @|\\
S@>{F_S}>>S, 
\end{CD} 
\end{equation*} 
we do not necessarily have the admissible immersion 
of 
$X^{\{p\}}_{\bul \leq N,\bul,\os{\circ}{T}{}'_0}$ 
over $S(T')^{\nat}$. 
Here $(T',{\cal J}',\del')$ is a log $p$-adic formal PD-enlargement of $S$ 
and $T'_0:=\ul{\rm Spec}^{\log}_{T'}({\cal O}_{T'}/{\cal J}')$. 
\end{rema}

\begin{coro}[{\bf Infinitesimal deformation invariance of filtered iso-zarisikian 
Steenbrink complexes}]\label{coro:finvcae}
If $T'=T$ and $X_{\bul \leq N,\os{\circ}{T}_0}(0)=X'_{\bul \leq N,\os{\circ}{T}_0}(0)$, then 
\begin{align*} 
(A_{{\rm zar},{\mab Q}}(X_{\bul \leq N,\os{\circ}{T}_0}/S(T)^{\nat}),P)
= (A_{{\rm zar},{\mab Q}}(X'_{\bul \leq N,\os{\circ}{T}_0}/S(T)^{\nat}),P). 
\tag{5.3.8.1}\label{ali:xdnz}
\end{align*} 
\end{coro}

\begin{coro}[{\bf Infinitesimal deformation invariance of log isocrystalline 
cohomologies with weight filtrations}]\label{coro:finvliae}
Let the notations be as in {\rm (\ref{coro:finvcae})}. 
Then 
\begin{align*} 
(R^qf_{X_{\bul \leq N,\os{\circ}{T}_0}/S(T)^{\nat}*}
({\cal O}_{X_{\bul \leq N,\os{\circ}{T}_0}/S(T)^{\nat}})_{\mab Q},P)&=
(R^qf_{X'_{\bul \leq N,\os{\circ}{T}_0}/S(T)^{\nat}*}
({\cal O}_{X'_{\bul \leq N,\os{\circ}{T}_0}/S(T)^{\nat}})_{\mab Q},P).
\tag{5.3.9.1}\label{ali:wfoa}\\
\end{align*}  
Moreover, if $T$ is restrictively hollow, then 
\begin{align*} 
(R^qf_{X_{\bul \leq N,T_0}/T*}
({\cal O}_{X_{\bul \leq N,T_0}/T})_{\mab Q},P)&=
(R^qf_{X'_{\bul \leq N,T_0}/T*}({\cal O}_{X'_{\bul \leq N,T_0}/T})_{\mab Q},P).
\tag{5.3.9.2}\label{ali:wfota}\\
\end{align*}  
\end{coro}


\section{The $E_2$-degeneration of 
the $p$-adic weight spectral sequence and the strict compatibility}
\label{sec:filbo}
Let the notations and the assumptions be as in the previous section. 
Assume that $\os{\circ}{X}_{\bul \leq N}$ is proper over $\os{\circ}{S}$ in this section. 
In this section we prove two important theorems: 
the $E_2$-degeneration of the $p$-adic weight spectral sequence of 
$X_{\bul \leq N,\os{\circ}{T}_1}/S(T)^{\nat}$ 
and the strict compatibility of the weight filtration 
with respect to the pull-back of a morphism of 
$N$-truncated simplicial base change of 
proper SNCL schemes over $S_{\os{\circ}{T}_0}$. 
To prove the $E_2$-degeneration, we need 
the infinitesimal deformation invariance of isocrystalline cohomologies 
with weight filtrations in the previous section ((\ref{ali:wfoa})). 
To prove the strict compatibility, we use the log convergence 
of the weight filtration ((\ref{theo:pwfaec})).  
The strict compatibility will play a key role in 
the definition of the limit of the weight filtration on the infinitesimal cohomology 
in \S\ref{sec:pff} below.

\begin{theo}[{\bf $E_2$-degeneration I}]\label{theo:e2dam}  
Let $s$ be the log point of a perfect field of characteristic $p>0$. 
The spectral sequence $(\ref{eqn:espssp})$ for the case $S=s$ 
and $E^{\bul \leq N}={\cal O}_{\os{\circ}{X}_{\bul \leq N,T_0}/\os{\circ}{T}}$ 
degenerates at $E_2$.  
\end{theo} 
\begin{proof} 
Let ${\cal W}$ be the Witt ring of $\kap:=\Gam(s,{\cal O}_s)$ 
and $K_0$ the fraction field of ${\cal W}$. 
Let $(T,p{\cal O}_T,[~])$ be an object of ${\rm Enl}^{\sq}(S/{\cal W})$. 
As in the proof of \cite[(3.1)]{of}, 
we have the following equalities (cf.~the proof of (\ref{prop:cxst})):    
\begin{align*} 
&R^qf_{X_{\bul \leq N,\os{\circ}{T}_0}/S(T)^{\nat}*}
({\cal O}_{X_{\bul \leq N,\os{\circ}{T}_0}/S(T)^{\nat}})_{\mab Q} 
={\cal H}^q(Rf_{X_{\bul \leq N,\os{\circ}{T}_0}/S(T)^{\nat}*}
({\cal O}_{X_{\bul \leq N,\os{\circ}{T}_0}/S(T)^{\nat}})_{\mab Q})
\tag{5.4.1.1}\label{ali:obrgm}\\
&={\cal H}^q({\cal O}_T\otimes^L_{\cal W}
R\Gam(X_{\bul \leq N}/{\cal W}(s))_{\mab Q})
={\cal H}^q({\cal K}_T\otimes^L_{K_0}
R\Gam(X_{\bul \leq N}/{\cal W}(s))_{\mab Q})\\
&={\cal K}_T\otimes_{K_0}
H^q_{\rm crys}(X_{\bul \leq N}/{\cal W}(s))_{\mab Q}. 
\end{align*} 
Hence we may assume that $T={\cal W}(s)$. 
\par 
Let ${\cal A}_0$ be a smooth ${\mab F}_p$-algebra contained in $\kap$ 
such that $X_{\bul \leq N}$ has a proper SNCL model ${\cal X}_{\bul \leq N}$ 
over 
${\cal S}_0:=({\rm Spec}({\cal A}_0),({\mab N}\oplus {\cal A}^*_0\lo {\cal A}_0))$ 
(cf.~the argument before \cite[(3.1)]{ndw}): 
$X_{\bul \leq N}={\cal X}_{\bul \leq N}\times_{{\cal S}_0}s$.   
Let ${\cal A}$ be a $p$-adically formally smooth algebra over ${\mab Z}_p$ 
which is a lift of ${\cal A}_0$. Endow ${\rm Spf}({\cal A})$ 
with a hollow log structure ${\mab N}\oplus {\cal A}^*\lo {\cal A}$ and 
let ${\cal S}$ be the resulting log formal scheme over 
${\rm Spf}({\mab Z}_p)=({\rm Spf}({\mab Z}_p), {\mab Z}_p^*)$. 
The log formal scheme ${\cal S}$ has a PD-ideal $p{\cal O}_{{\cal S}}$, 
which defines an exact closed immersion ${\cal S}_0 \os{\subset}{\lo} {\cal S}$.  
By  (\ref{coro:fenlt}),  
$H^q_{\rm crys}({\cal X}_{\bul \leq N}/{\cal S})_{\mab Q}$ is 
a locally free ${\cal K}_{\cal S}$-module 
(this locally freeness also follows from \cite[Lemma 36]{ollc}). 
For a nonempty log formal affine open subscheme ${\cal U}$ of ${\cal S}$,  set 
${\cal U}_1:={\cal U}{\otimes}_{{\mab Z}_p}{\mab F}_p$. 
We fix a lift $F_{\cal U} \col {\cal U} \lo {\cal U}$ of 
the Frobenius endomorphism(=$p$-th power endomorphism) of ${\cal U}_1$. 
(The lift $F_{\cal U}$ indeed exists.)
Then we have the Teichm\"{u}ller lift ${\cal O}_{\cal U}\lo {\cal W}$ of the inclusion 
${\cal O}_{{\cal U}_{1}}\os{\sus}{\lo} \kap$. 
This lift gives us a morphism ${\cal W}(s)\lo {\cal U}$ of log $p$-adic formal schemes. 
In particular,  we have a morphism ${\cal W}(s)\lo {\cal S}$ of 
log $p$-adic formal schemes. 
By the base change theorem of log crystalline cohomologies, 
we have an equality  
$$R\Gam(X_{\bul \leq N}/{\cal W}(s))=
{\cal W}\otimes^L_{{\cal O}_{\cal S}}R\Gam({\cal X}_{\bul \leq N}/{\cal S}).$$ 
Because $H^q_{\rm crys}({\cal X}_{\bul \leq N}/{\cal S})_{\mab Q}$ is a flat ${\cal K}_{\cal S}$-module 
and because 
$R\Gam({\cal X}_{\bul \leq N}/{\cal S})_{\mab Q}$ is a perfect complex of 
${\cal K}_{\cal S}$-modules (${\cal K}_S:={\cal O}_S\otimes_{\mab Z}{\mab Q}$), 
we have the following equalities as in (\ref{ali:obrgm}):    
\begin{align*} 
H^q_{\rm crys}(X_{\bul \leq N}/{\cal W}(s))_{\mab Q}& 
=H^q(R\Gam(X_{\bul \leq N}/{\cal W}(s))_{\mab Q})
=H^q({\cal W}(\kap)\otimes^L_{{\cal O}_{\cal S}}
R\Gam({\cal X}_{\bul \leq N}/{\cal S})_{\mab Q}) 
\tag{5.4.1.2}\label{ali:osrgm}\\
& =H^q(K_0\otimes^L_{{\cal K}_{\cal S}}
R\Gam({\cal X}_{\bul \leq N}/{\cal S})_{\mab Q})
=K_0\otimes_{{\cal K}_{\cal S}}H^q_{\rm crys}({\cal X}_{\bul \leq N}/{\cal S})_{\mab Q}. 
\end{align*} 
\par
The rest of the proof is the same as that of \cite[(3.6)]{ndw}. 
Indeed, take a closed point $\os{\circ}{u}$ of $\os{\circ}{\cal U}_1$. 
The point $\os{\circ}{u}$ is the spectrum of a finite field $\kap_u$. 
Then we have the Teichm\"{u}ller lift 
${\cal O}_{\cal U}\lo {\cal W}(\kap_u)$ of the morphism
${\cal O}_{{\cal U}_{1}}\lo \kap_u$. 
The ring ${\cal W}(\kap_u)$ becomes an ${\cal O}_{\cal U}$-algebra by this lift. 
Endow $\os{\circ}{u}$ and ${\rm Spf}({\cal W}(\kap_u))$
with the inverse images of the log structure of ${\cal U}$
and let $u$ and ${\cal W}(u)$ be the resulting log (formal) schemes.  
Let $K_{0,u}$ be the fraction field of ${\cal W}(\kap_u)$. 
Set ${\cal X}_{\bul \leq N,{\cal U}_1}:=
{\cal X}_{\bul \leq N}\times_{{\cal S}_1}{\cal U}_1$ and ${\cal X}_{\bul \leq N,u}
:={\cal X}_{\bul \leq N,T_1}{\times}_{{\cal U}_1}u$.  
Let $\{E^{\bul \bul}_r(X_{\bul \leq N}/K_0)\}_{r\geq 1}$ 
and $\{E^{\bul \bul}_r({\cal X}_{\bul \leq N,u}/K_{0,u})\}_{r\geq 1}$
be the $E_r$-terms of the weight spectral sequences (\ref{eqn:espssp}) 
of $X_{\bul \leq N}/K_0$ and ${\cal X}_{\bul \leq N,u}/K_{0,u}$, respectively. 
By the explicit description of the boundary morphisms between $E_1$-terms 
of (\ref{eqn:espssp}) ((\ref{prop:exbd})) and using the argument in the proof of \cite[(3.5)]{ndw}, 
we obtain the following: 
there exists a nonempty log formal affine open subscheme ${\cal U}$ of 
${\cal S}$ such that, for any closed point 
$\os{\circ}{u}$ of $\os{\circ}{\cal U}_1$ 
and for all $k$, $q\in {\mab Z}$, 
there exists a $($non-canonical$)$ isomorphism
$$E_2^{-k,q+k}(X_{\bul \leq N}/K_0){\otimes}_{{\mab Q}_p}K_{0,u} 
\os{\sim}{\lo}
E_2^{-k,q+k}({\cal X}_{\bul \leq N,u}/K_{0,u})
{\otimes}_{{\mab Q}_p}K_0
$$
of $K_0{\otimes}_{{\mab Q}_p}K_{0,u}\simeq 
K_{0,u}{\otimes}_{{\mab Q}_p}K_0$-modules. 
Hence it suffices to prove that 
the boundary morphisms 
$d_r^{-k,q+k}\col E^{-k,q+k}_r({\cal X}_{\bul \leq N,u}/K_{0,u}) 
\lo E_r^{-k+r,q+k-r+1}({\cal X}_{\bul \leq N,u}/K_{0,u})$ 
vanish for $r\geq 2$. 
These vanish because
$E^{-k,q+k}_1({\cal X}_{\bul \leq N,u}/K_{0,u})$ is the direct sum of 
the isocrystalline cohomologies of proper smooth schemes over a finite field 
with Tate twists and because $E^{-k,q+k}_1({\cal X}_{\bul \leq N,u}/K_{0,u})$ is 
of pure weight $q+k$ by  
\cite[Corollary 1. 2)]{kme}, \cite[(1.2)]{clpu} and \cite[(2.2) (4)]{ndw}.  
\end{proof}

\begin{rema}\label{rema:nllfcs} 
(1) For the proof of \cite[Corollary 1 2)]{kme}, the weak-Lefschetz conjecture 
for a hypersurface of a large degree has been used. 
This  weak-Lefschetz conjecture has been stated in \cite{bwl}. 
However, in \cite[(7.18)]{ny},  
I have pointed out the existence of a gap in the proof in \cite{bwl}
and I have filled the gap in \cite[(7.17)]{ny}. 
In the appendix of \cite{ny} I have given another short proof of 
the general weak Lefschetz conjecture by using theory of 
rigid cohomologies. 
\par 
(2) Let ${\cal U}$ be as in the proof of (\ref{theo:e2dam}). 
Consider the naive truncated simplicial case in \S\ref{sec:psc}. 
Let $d_r^{\bul \bul}({\cal X}_{\bul \leq N,{\cal U}_{1}}/{\cal K}_{\cal U})$ 
be the boundary morphism of the weight spectral sequence 
$(\ref{eqn:escssp})\otimes_{\mab Z}{\mab Q}$ 
of ${\cal X}_{\bul \leq N,{\cal U}_{1}}/{\cal U}$.  
Then we can reduce the vanishing of 
$d_r^{\bul \bul}(X_{\bul \leq N}/K_0)$ $(r \geq 2)$ 
to that of $d_r^{\bul \bul}({\cal X}_{\bul \leq N,{\cal U}_{1}}/{\cal K}_{\cal U})$ 
and then to that of $d_r^{\bul \bul}({\cal X}_u/K_{0,u})$ 
in a more standard way as in \cite[(2.15.4)]{nh2} because 
we have constructed the {\it integral} weight spectral sequence 
(\ref{eqn:escssp}) in the naive truncated simplicial case. 
This is an answer for the suggestive method in \cite[(3.7)]{ndw}.  
\end{rema}

\begin{theo}[{\bf $E_2$-degeneration II}]\label{theo:e2dgfam} 
Let ${\cal V}$ be as in {\rm \S\ref{sec:cfi}}. 
Let $T$ be an object of ${\rm Enl}_p(S/{\cal V})$ 
with structural morphism $T_1\lo S$. 
The spectral sequence $(\ref{eqn:espssp})$ for the case 
$E^{\bul \leq N}={\cal O}_{\os{\circ}{X}_{\bul \leq N,T_0}/\os{\circ}{T}}$ 
degenerates at $E_2$. 
Consequently the spectral sequence {\rm (\ref{eqn:getcpsp})} 
for the case $E^{\bul \leq N}_K={\cal O}_{\os{\circ}{X}_{\bul \leq N}/K}$ 
degenerates at $E_2$. 
\end{theo}
\begin{proof} 
(By virtue of (\ref{prop:e2tmc}) and (\ref{coro:finvliae}), 
the proof of this theorem is almost the same as that of \cite[(2.17.2)]{nh2}.) 
\par 
We may assume that $\os{\circ}{T}$ is 
a $p$-adic affine flat formal ${\cal V}$-scheme ${\rm Spf}(B)$.  
Let $\{E_{r}^{-k,q+k}(X_{\bul \leq N,\os{\circ}{T}_1}/S(T)^{\nat})_{K_0}\}$ be the $E_r$-terms 
of (\ref{eqn:espssp}). 
Consider the following boundary morphism: 
\begin{align*}
d_{r}^{-k, q+k} \col 
E_{r}^{-k,q+k}(X_{\bul \leq N,\os{\circ}{T}_1}/S(T)^{\nat})\lo 
E_{r}^{-k+r, q+k-r+1}(X_{\bul \leq N,\os{\circ}{T}_1}/S(T)^{\nat}) \quad (r \geq 2). 
\tag{5.4.3.1}\label{ali:bdafst}
\end{align*} 
We prove that $d_{r}^{-k, q+k}=0$  
$(r \geq 2)$. 
\par
Case A: First we consider a case where  $B$ is 
a topologically finitely generated ring  over ${\cal V}$ such that $B_K$ is an artinian local ring. 
Let ${\mathfrak m}$ be the maximal ideal of $B_K$. 
Set $K':=B_K/{\mathfrak m}$. 
Let ${\cal V}'$ be the integer ring of $K'$. 
Consider the following ideal of $B$: $I:={\rm Ker}(B\lo B_K/{\mathfrak m})$. 
Set $C=B/I$. Then $C_{K}=K'$ and  
${\cal V} \subset C\subset {\cal V}'$ (\cite[(2.17.1)]{nh2}). 
Let 
$\os{\circ}{\iota} \col {\rm Spf}(C) \os{\sus}{\lo} {\rm Spf}(B)$ 
be the nilpotent closed immersion. 
Since the characteristic of $K$ is 0, 
the morphism ${\rm Spec}(C_K) \lo {\rm Spec}(K)$ 
is smooth and hence there exists a retraction 
$\os{\circ}{\rho}_K \col {\rm Spec}(B_K) \lo {\rm Spec}(C_K)$ 
of the nilpotent closed immersion 
${\rm Spec}(C_K) \os{\sus}{\lo} {\rm Spec}(B_K)$. 
By \cite[(1.17)]{of} there exists a finite modification 
$\os{\circ}{\pi} \col {\rm Spf}(B') \lo {\rm Spf}(B)$, a nilpotent closed immersion 
$\os{\circ}{\iota}{}' \col {\rm Spf}(C) \os{\subset}{\lo} {\rm Spf}(B')$ with 
$\os{\circ}{\pi} \circ \os{\circ}{\iota}{}' = \os{\circ}{\iota}$ and a morphism 
$\os{\circ}{\rho} \col {\rm Spf}(B') \lo {\rm Spf}(C)$ such that $\os{\circ}{\rho}$ induces 
$\os{\circ}{\rho}_K$ and that $\os{\circ}{\rho} \circ \os{\circ}{\iota}{}' = {\rm id}_{{\rm Spf}(C)}$. 
Set $\os{\circ}{T}{}':={\rm Spf}(B')$, $\os{\circ}{T}{}'' := {\rm Spf}(C)$, 
$S(T')^{\nat}:=S(T)^{\nat}\times_{\os{\circ}{T}}\os{\circ}{T}{}'$ 
and $S(T'')^{\nat}:=S(T)^{\nat}\times_{\os{\circ}{T}}\os{\circ}{T}{}''$. 
Then we have the following morphisms
\begin{align*} 
\iota \col S(T'')^{\nat} \os{\sus}{\lo} S(T)^{\nat},~\iota' \col S(T'')^{\nat} \os{\sus}{\lo} S(T')^{\nat},~
\pi \col S(T')^{\nat} \lo S(T)^{\nat},~\rho \col S(T')^{\nat} \lo S(T'')^{\nat}  
\end{align*}
such that 
$\pi \circ \iota' =\iota$ and $\rho \circ \iota' = {\rm id}_{S(T'')^{\nat}}$. 
\par 
By (\ref{prop:e2tmc}) we have
$E_2^{-k, q+k}(X_{\bul \leq N,\os{\circ}{T}_1}/S(T)^{\nat})=
E_2^{-k, q+k}(X_{\bul \leq N,\os{\circ}{T}{}'_1}/S(T')^{\nat})$ since $B'_K=B_K$. 
Let $\{d_r'{}^{\bul \bul}\}$ $(r \geq 1)$ be the boundary morphism of 
(\ref{eqn:espssp}) for $X_{\bul \leq N,\os{\circ}{T}{}'_1}/S(T')^{\nat}$. 
Because $\{d_r^{\bul \bul}\}$ $(r \geq 2)$
are functorial with respect to a morphism of 
log $p$-adic enlargements,  
we have the following commutative diagram for $r\geq 2$:
\begin{equation*}
\begin{CD}
E_r^{-k, q+k}(X_{\bul \leq N,\os{\circ}{T}_1}/S(T)^{\nat}) @>>>  
E_r^{-k, q+k}(X_{\bul \leq N,\os{\circ}{T}{}'_1}/S(T')^{\nat})\\ 
@V{d_r^{{-k, q+k}}}VV  @VV{d_r'{}^{{-k, q+k}}}V \\
E_r^{-k+r, q+k-r+1}(X_{\bul \leq N,\os{\circ}{T}_1}/S(T)^{\nat}) @>>> 
E_r^{-k+r, q+k-r+1}(X_{\bul \leq N,\os{\circ}{T}{}'_1}/S(T')^{\nat}).
\end{CD}
\end{equation*}
Here, if $r=2$, then two horizontal morphisms above 
are isomorphisms by (\ref{prop:e2tmc}) since $B'_K=B_K$.
By induction on $r \geq 2$, we see that 
$d_r^{\bul \bul}$  vanishes 
if $d_r'{}^{\bul \bul}$ does. 
Hence it suffices to prove 
that the boundary morphism 
\begin{align*}
d_{r}'{}^{-k, q+k} \col & E_{r}^{-k,q+k}(X_{\bul \leq N,\os{\circ}{T}{}'_1}/S(T')^{\nat})
\lo E_{r}^{-k+r, q+k-r+1}(X_{\bul \leq N,\os{\circ}{T}{}'_1}/S(T')^{\nat})
\quad (r \geq 2) 
\tag{5.4.3.2}\label{eqn:bdsed} 
\end{align*} 
vanishes. 
Let $l(M)$ be the length of a finitely generated $B'_K=B_K$-module $M$.
Furthermore,  to prove the vanishing of $d_r'{}^{\bul \bul}$, 
it suffices to prove that 
\begin{equation*}
l(R^qf_{X_{\bul \leq N,\os{\circ}{T}{}'_1}/S(T')^{\nat}*}({\cal O}_{X_{\bul \leq N,\os{\circ}{T}{}'_1}/S(T')^{\nat}})_K)
= l({\bigoplus}_{k}E_2^{-k, q+k}(X_{\bul \leq N,\os{\circ}{T}{}'_1}/S(T')^{\nat})). 
\tag{5.4.3.3}
\label{eqn:legthe2ar}
\end{equation*}
For a log formal scheme $Z$ with $p$-adic topology, set $Z_1:=Z\mod p$. 
Then we have the morphisms 
$X_{\bul \leq N,\os{\circ}{T}{}''_1} \lo S(T'')^{\nat}$.  
Let us denote the pull-back of 
the morphism 
$X_{\bul \leq N,\os{\circ}{T}{}''_1} \lo S(T'')^{\nat}$ by 
$\rho \mod p: S(T')^{\nat}_1 \lo S(T'')^{\nat}_1$ by 
$X'_{\bul \leq N,\os{\circ}{T}{}'_1} \lo S(T')^{\nat}$. 
Then, since we have 
$\pi \circ \iota' = \iota$ and $\rho \circ \iota' = {\rm id}_{S(T'')^{\nat}}$, both 
$X_{\bul \leq N,\os{\circ}{T}{}'_1}$ and 
$X'_{\bul \leq N,\os{\circ}{T}{}'_1}$ are deformations of 
$X_{\bul \leq N,\os{\circ}{T}{}''_1}$ to $S(T')^{\nat}_1$.  
Hence, by (\ref{coro:finvcae}), the spectral sequence (\ref{eqn:espssp}) 
for $X_{\bul \leq N,\os{\circ}{T}{}'_1}/S(T')^{\nat}$ and 
that for $X'_{\bul \leq N,\os{\circ}{T}{}'_1}/S(T')^{\nat}$ are isomorphic. 
Therefore we have the following equalities  
\begin{align*}
E_2^{-k, q+k}(X_{\bul \leq N,\os{\circ}{T}{}'_1}/S(T')^{\nat})  &= 
E_2^{-k, q+k}(X'_{\bul \leq N,\os{\circ}{T}{}'_1}/S(T')^{\nat}) \\
&= B' \otimes_C E_2^{-k, q+k}(X'_{\bul \leq N,\os{\circ}{T}{}''_1}/S(T'')^{\nat}).  
\end{align*}  
(The last equality follows from (\ref{prop:e2tmc}).) 
Hence, to prove (\ref{eqn:legthe2ar}), it suffices 
to prove that 
\begin{align*}
{\rm dim}_{K'}(R^qf_{X'_{\bul \leq N,\os{\circ}{T}{}''_1}/S(T'')^{\nat}*} & 
({\cal O}_{X'_{\bul \leq N,\os{\circ}{T}{}''_1}/S(T'')^{\nat}})_K) 
= 
{\rm dim}_{K'}({\bigoplus}_{k}
E_2^{-k, q+k}(X'_{\bul \leq N,\os{\circ}{T}{}''_1}/S(T'')^{\nat})).   
\tag{5.4.3.4}\label{eqn:dimareq}  
\end{align*}
\par 
Set $S(T''')^{\nat}:=S(T'')^{\nat}\times_{{\rm Spf}(C)}{\rm Spf}({\cal V}')$. 
Because there exists a natural morphism $S(T''')^{\nat} \lo S(T'')^{\nat}$ 
of log $p$-adic enlargements of $S$, 
it suffices to prove that 
\begin{align*}
&{\rm dim}_{K'}(R^qf_{X'_{\bul \leq N,\os{\circ}{T}{}'''_1}/S(T''')^{\nat}*}
({\cal O}_{X'_{\bul \leq N,\os{\circ}{T}{}'''_1}/S(T''')^{\nat}*})_{K'})  \tag{5.4.3.5}\label{eqn:dimv}\\
&= 
{\rm dim}_{K'}({\bigoplus}_kE_2^{-k, q+k}(X'_{\bul \leq N,\os{\circ}{T}{}'''_1}/S(T''')^{\nat})). 
\end{align*}
We reduce (\ref{eqn:dimv}) to (\ref{theo:e2dam}) by using  (\ref{theo:definv}) as follows.  
\par 
Let $\kap'$ be the residue field of ${\cal V}'$. 
Since $\kap$ is perfect and since $\kap'$ is a finite extension of $\kap$, 
$\kap'$ is also perfect. Let $s'$ be the log point of $\kap'$. 
Let ${\cal W}'$ be the Witt ring of $\kap'$ and set $K'_0:={\rm Frac}({\cal W}')$. 
Set ${\cal V}'_1:={\cal V}'/p$. 
The ring ${\cal V}'_1$ is an artinian local $\kap'$-algebra with residue field $\kap'$ 
(\cite[II Proposition 8]{sec}). 
This gives us a retraction $S(T''')^{\nat}_1\lo s'$ of the immersion $s'\os{\sus}{\lo} S(T''')^{\nat}_1$. 
Set $X''_{\bul \leq N}:=X'_{\bul \leq N,\os{\circ}{T}{}'''_1}\times_{\os{\circ}{T}{}'''_1}\os{\circ}{s}{}'$. 
Then $X''_{\bul \leq N}\times_{\os{\circ}{s}{}'}\os{\circ}{T}{}'''_1$ and 
$X'_{\bul \leq N,\os{\circ}{T}{}''_1}$ are two log deformations of $X''_{\bul \leq N}$. 
We obtain the following isomorphisms 
by (\ref{coro:finvce}) and (\ref{coro:finvcae}): 
\begin{align*}
R^qf_{X'_{\bul \leq N,\os{\circ}{T}{}'''_1}/S(T''')^{\nat}*}
( {\cal O}_{X'_{\bul \leq N,\os{\circ}{T}{}'''_1}/S(T''')^{\nat}})\otimes_{{\cal V}'}K'  
\os{\sim}{\lo}
R^qf_{X''_{\bul \leq N}/{\cal W}(s')*}({\cal O}_{X''_{\bul \leq N}/{\cal W}(s')})
{\otimes}_{{\cal W}'}K' \tag{5.4.3.6}\label{eqn:dfvw}
\end{align*}
and 
\begin{equation*}
E_2^{-k, q+k}(X'_{\bul \leq N,\os{\circ}{T}{}'''_1}/S(T''')^{\nat})\otimes_{K'_0}K'
\os{\sim}{\lo}E_2^{-k, q+k}(X''_{\bul \leq N}/{\cal W}(s')){\otimes}_{K'_0}K'. 
\tag{5.4.3.7}\label{eqn:e2dfinv}
\end{equation*}
Hence it suffices to prove that 
$$E_2^{-k,q+k}(X''_{\bul \leq N}/{\cal W}(s')){\otimes}_{K'_0}K'=
E_{\infty}^{-k,q+k}(X''_{\bul \leq N}/{\cal W}(s')){\otimes}_{K'_0}K'.$$
We have already proved this in (\ref{theo:e2dam}).
\par
Case B: Next we consider the general case. 
Though the rest of the proof is almost the same as that of 
\cite[(2.17.2)]{nh2} by using (\ref{prop:e2tmc}), 
we repeat this for the completeness of this book. 
Let ${\mathfrak m}$ be a maximal ideal of $B_K$. 
Consider the following ideal $I^{(n)}$ and the following ring $B_{(n)}$ 
in \cite[p.~780]{of}:
$$I^{(n)}:={\rm Ker}(B \lo B_{K}/{\mathfrak m}^n), 
\quad B_{(n)} :=B/I^{(n)} \quad (n\in {\mab N}).$$
The ring $B_{(n)}$ defines log $p$-adic enlargements 
$S(T_{(n)})^{\nat}$ of $S(T)^{\nat}$. 
Let 
\begin{align*}
d_{r,(n)}^{{-k, q+k}} \col 
E_r^{-k, q+k}(X_{\bul \leq N,\os{\circ}{T}_{(n),1}}/S(T_{(n)})^{\nat}) 
\lo E_r^{-k+r, q+k-r+1}(X_{\bul \leq N,\os{\circ}{T}_{(n),1}}/S(T_{(n)})^{\nat})
\end{align*}
be the boundary morphism.
Because $\{d_r^{\bul \bul}\}$ is functorial, 
we have the following commutative diagram: 
\begin{equation*}
\begin{CD}
E_r^{-k, q+k}(X_{\bul \leq N,\os{\circ}{T}_1}/S(T)^{\nat})
{\otimes}_BB_{(n)} \!\!@>>>\!\! 
E_r^{-k, q+k}(X_{\bul \leq N,\os{\circ}{T}_{(n),1}}/S(T_{(n)})^{\nat})\\ 
@V{d_r^{-k, q+k}\otimes_{B}B_{(n)}}VV 
@VV{d_{r,(n)}^{{-k, q+k}}}V \\
E_r^{-k+r, q+k-r+1}(X_{\bul \leq N,\os{\circ}{T}_1}/S(T)^{\nat}){\otimes}_{B}B_{(n)} \!\!@>>>\!\!  
E_r^{-k+r, q+k-r+1}(X_{\bul \leq N,T_{(n),1}}/S(T_{(n)})^{\nat}).
\end{CD}
\end{equation*}
Because $E_2^{-k,q+k}(X_{\bul \leq N,\os{\circ}{T}_1}/S(T)^{\nat})$ is 
a convergent $F$-isocrystal ((\ref{prop:e2tmc})), 
two horizontal morphisms are 
isomorphisms if $r=2$.
By induction on $r$ and by the proof for the Case A, 
the boundary morphism 
$d_r^{\bul \bul}{\otimes}_{B}B_{(n)}$ 
$(r \geq 2)$ vanishes. Thus 
${\vpl}_n(d_r^{\bul \bul}{\otimes}_{B_{K}}
B_{K}/{\mathfrak m}^n)=0$. 
Because $B_K$ is a noetherian ring and 
$E_2^{-k, q+k}(X_{\bul \leq N,\os{\circ}{T}_1}/S(T)^{\nat})$ is a 
finitely generated $B_K$-module,
we have 
$$d_r^{\bul \bul}{\otimes}_{B_{K}}
({\vpl}_n B_{K}/{\mathfrak m}^n)=
\vpl_n(d_r^{\bul \bul}{\otimes}_{B_{K}}
B_{K}/{\mathfrak m}^n)=0.$$
Since $(B_{K})_{\mathfrak m}$ is a Zariski ring, 
$\vpl_n (B_{K})_{\mathfrak m}/
{\mathfrak m}^n(B_{K})_{\mathfrak m}$ 
is faithfully flat over $(B_{K})_{\mathfrak m}$  
(\cite[III \S3 Proposition 9]{bou2}).
Therefore 
$d_r^{\bul \bul}{\otimes}_{B_{K}}
(B_{K})_{\mathfrak m}=0$.
Since ${\mathfrak m}$ is 
an arbitrary maximal ideal of $B_{K}$,  
$d_r^{\bul \bul}=0$ $(r \geq 2)$.
\end{proof}

As in \cite[(4.7)]{ndw} we have the following corollary:

\begin{theo}\label{theo:edhw}
The spectral sequences $(\ref{eqn:espmuasp})$ modulo torsion and 
$(\ref{eqn:sesplasp})$ 
degenerate at $E_2$ for the case 
$E^{\bul \leq N}
={\cal O}_{\os{\circ}{X}_{\bul \leq N}/\os{\circ}{S}}$.
\end{theo}

\begin{coro}\label{coro:ne2}
The spectral sequence {\rm (\ref{eqn:getcpsp})} 
degenerates at $E_2$.
\end{coro}

\par 
Next we prove the strict compatibility of the pull-back of a morphism 
with respect to the weight filtration. 
We give the following strict compatibility 
whose proof is not the same as that of the strict compatibility 
for an open log scheme in \cite[(2.18.1)]{nh2}.  

\begin{theo}[{\bf Strict compatibility I}]\label{theo:stpfgb}
Let the notations be as in {\rm (\ref{theo:e2dam})} and the proof of it. 
Let $Y_{\bul \leq N}/s'$ be an analogous object to $X_{\bul \leq N}/s$. 
Let $f' \col Y_{\bul \leq N}\lo s'$ be the structural morphism. 
Let $h\col s\lo s'$ be a morphism of log schemes. 
Let $g_{\bul \leq N}\col 
X_{\bul \leq N,\os{\circ}{T}_0}\lo Y_{\bul \leq N,\os{\circ}{T}{}'_0}$ 
be the morphism in {\rm (\ref{eqn:xdbquss})} for the case $S=s$ and $S'=s'$ 
fitting into the commutative diagram {\rm (\ref{cd:xdjqbxy})}. 
Let ${\cal W}'$ be the Witt ring of $\Gam(s',{\cal O}_{s'})$ 
and set $K'_0:={\rm Frac}({\cal W}')$. 
Assume that $\os{\circ}{s}\lo \os{\circ}{s}{}'$ is finite. 
Let $q$ be an integer. 
Let us endow 
$H^q_{{\rm crys}}(Y_{\bul \leq N}/{\cal W}(s'))\otimes_{{\cal W}'}K_0$  
with the induced filtration by $P$ on 
$H^q_{{\rm crys}}(Y_{\bul \leq N}/{\cal W}(s'))_{K'_0}$.   
Then the induced morphism 
\begin{equation*}
g^* \col H^q_{{\rm crys}}(Y_{\bul \leq N}/{\cal W}(s'))\otimes_{{\cal W}'}K_0
\lo 
H^q_{{\rm crys}}(X_{\bul \leq N}/{\cal W}(s))\otimes_{\cal W}K_0
\tag{5.4.6.1}\label{eqn:gbwstn}
\end{equation*} 
is strictly compatible with the weight filtration. 
\end{theo}
\begin{proof}
Let the notations be as in the proof of (\ref{theo:e2dam}). 
Let ${\cal B}_0$, ${\cal B}$ and ${\cal Y}_{\bul \leq N}$ 
are analogous objects to ${\cal A}_0$, ${\cal A}$ and ${\cal X}_{\bul \leq N}$, 
respectively, such that $g$ has a model 
${\mathfrak g}_{\bul \leq N}\col {\cal X}_{\bul \leq N}\lo {\cal Y}_{\bul \leq N}$ 
over $({\rm Spec}({\cal A}_0),({\mab N}\oplus {\cal A}_0^*\lo {\cal A}_0))\lo 
({\rm Spec}({\cal B}_0),({\mab N}\oplus {\cal B}_0^*\lo {\cal B}_0))$. 
Here we may assume that the morphism ${\cal B}_0\lo {\cal A}_0$ is of finite type. 
There exists a morphism ${\cal B}\lo {\cal A}$ lifting the morphism above. 
Let $F_{\cal A}\col {\cal A}\lo {\cal A}$ and 
$F_{\cal B}\col {\cal B}\lo {\cal B}$ be lifts of the Frobenius endomorphisms  of 
$F_{{\cal A}_0}\col {\cal A}_0\lo {\cal A}_0$ 
and $F_{{\cal B}_0}\col {\cal B}_0\lo {\cal B}_0$, respectively. 
Then we have endomorphisms of 
$({\rm Spec}({\cal A}),({\mab N}\oplus {\cal A}^*\lo {\cal A}))$ 
and $({\rm Spec}({\cal B}),({\mab N}\oplus {\cal B}^*\lo {\cal B}))$ lifting  
the Frobenius endomorphisms of  
$({\rm Spec}({\cal A}_0),({\mab N}\oplus {\cal A}_0^*\lo {\cal A}_0))$ 
and $({\rm Spec}({\cal B}_0),({\mab N}\oplus {\cal B}_0^*\lo {\cal B}_0))$, 
respectively. 
\par 
Let $\{{\cal P}_k\}_{k=0}^{\infty}$ and $\{{\cal Q}_k\}_{k=0}^{\infty}$ 
be the weight filtrations on 
$R^qf_{{\cal X}_{\bul \leq N}/{\cal A}*}({\cal O}_{{\cal X}_{\bul \leq N}/{\cal A}})_{\mab Q}$ 
and $R^qf_{{\cal Y}_{\bul \leq N}/{\cal B}*}
({\cal O}_{{\cal Y}_{\bul \leq N}/{\cal B}})_{\mab Q}\otimes_{\cal B}{\cal A}$, respectively. 
We claim that 
\begin{align*} 
{\cal P}_k\cap {\rm Im}({\cal Q}_{\infty})={\rm Im}({\cal Q}_k). 
\tag{5.4.6.2}\label{ali:bh}
\end{align*}
Indeed, by (\ref{theo:fugennas}), 
the right hand side of (\ref{ali:bh}) is contained in 
the left hand side of (\ref{ali:bh}). 
By (\ref{exam:ofl}) and  (\ref{prop:nupf}),  
${\cal P}_k$, ${\cal P}_{\infty}$, ${\cal Q}_k$ and ${\cal Q}_{\infty}$ 
extend to objects ${\cal P}^{\rm conv}_k$ 
${\cal P}^{\rm conv}_{\infty}$, ${\cal Q}^{\rm conv}_k$ and 
${\cal Q}^{\rm conv}_{\infty}$ 
of $F{\textrm -}{\rm Isoc}({\rm Spec}({\cal A}_0)/{\mab Z}_p)$, respectively. 
Because the category $F{\textrm -}{\rm Isoc}({\rm Spec}({\cal A}_0)/{\mab Z}_p)$
is an abelian category (\cite[p.~795]{of}), 
both hand sides on (\ref{ali:bh}) extend to objects of 
$F{\textrm -}{\rm Isoc}({\rm Spec}({\cal A}_0)/{\mab Z}_p)$. 
By \cite[(4.1)]{of}, for a closed point $\os{\circ}{s}$ of ${\rm Spec}({\cal A}_0)$, 
the pull-back functor $\os{\circ}{s}{}^*\col {\rm Isoc}({\rm Spec}({\cal A}_0)/{\mab Z}_p)
\lo {\rm Isoc}(\os{\circ}{s}/{\mab Z}_p)$ is faithful. 
Hence it suffices to prove that 
\begin{align*} 
\os{\circ}{s}{}^*({\cal P}^{\rm conv}_k)\cap 
{\rm Im}(\os{\circ}{s}{}^*({\cal Q}'{}^{\rm conv}_{\infty}))
={\rm Im}(\os{\circ}{s}{}^*({\cal Q}'{}^{\rm conv}_k)). 
\tag{5.4.6.3}\label{ali:bsh}
\end{align*}
Let $\os{\circ}{t}$ be the image of $\os{\circ}{s}$ in ${\rm Spec}({\cal B}_0)$. 
Then $k(\os{\circ}{t})$ is a finite field. 
Consider the $p$-adic enlargements 
${\cal W}(\os{\circ}{s})$ and ${\cal W}(\os{\circ}{t})$ 
of ${\rm Spec}({\cal A}_0)/{\mab Z}_p$ and 
${\rm Spec}({\cal B}_0)/{\mab Z}_p$, respectively. 
Let $u$ be $s$ or $t$ and let ${\cal C}$ be ${\cal A}$ or ${\cal B}$, respectively. 
Because ${\rm Spf}({\cal C})$ is formally smooth over ${\rm Spf}({\mab Z}_p)$, 
there exists a morphism  ${\cal W}(\os{\circ}{u})\lo {\rm Spf}({\cal C})$ 
fitting into the following commutative diagram
\begin{equation*} 
\begin{CD} 
\os{\circ}{u}@>{\subset}>> {\cal W}(\os{\circ}{u}) \\
@VVV @VVV \\
{\rm Spec}({\cal C}_0) @>>> {\rm Spf}({\cal C}) \\
@VVV @VVV \\
{\rm Spec}({\mab F}_p) @>>> {\rm Spf}({\mab Z}_p).  
\end{CD} 
\end{equation*}  
Then the log structure of ${\cal W}(u)$ is 
the inverse image of the log structure of 
$({\rm Spf}({\cal C}),({\mab N}\oplus {\cal C}^*\lo {\cal C}))$. 
Let ${\cal Z}_{\bul \leq N}$ be ${\cal X}_{\bul \leq N}$ or 
${\cal Y}_{\bul \leq N}$. 
By the proof of (\ref{prop:nupf}), we see that 
\begin{align*} 
(\os{\circ}{u}{}^*({\cal P}{}^{\rm conv}_{\infty}))_{{\cal W}(\os{\circ}{u})}=
R^qf_{{\cal Z}_{\bul \leq N,u}/{\cal W}(u)*}
({\cal O}_{{\cal Z}_{\bul \leq N,u}/{\cal W}(u)})_{\mab Q}
\tag{5.4.6.4}\label{ali:bsznh}
\end{align*}
and 
\begin{align*} 
(\os{\circ}{u}{}^*({\cal P}{}^{\rm conv}_k))_{{\cal W}(\os{\circ}{u})}=
P_kR^qf_{{\cal Z}_{\bul \leq N,u}/{\cal W}(u)*}
({\cal O}_{{\cal Z}_{\bul \leq N,u}/{\cal W}(u)})_{\mab Q}.
\tag{5.4.6.5}\label{ali:bnush}
\end{align*}
Since the residue field $k(\os{\circ}{s})$ is a finite field, 
the equality of (\ref{ali:bsh}) at the value ${\cal W}(s)$ holds by 
the existence of the weight spectral sequence (\ref{eqn:espssp})
and by the purity of the weight as in the proof of \cite[(2.18.2)]{nh2}.   
This shows the equality (\ref{ali:bsh}) 
because ${\cal W}(\os{\circ}{u})$ is a formally smooth lift 
of $\os{\circ}{u}$ over ${\rm Spf}({\mab Z}_p)$. 
\end{proof}


\begin{theo}[{\bf Strict compatibility II}]\label{theo:stpbgb}
Let the notations and the assumption be as in {\rm (\ref{theo:e2dgfam})}. 
Let $Y_{\bul \leq N}/S'$ and $T'$ be analogous objects to $X_{\bul \leq N}/S$ 
and $T$, respectively. 
Let $g$ be a morphism {\rm (\ref{eqn:xdbquss})} 
fitting into the following commutative diagram {\rm (\ref{cd:xdjqbxy})}.
Let $q$ be an integer. 
Then the induced morphism 
\begin{equation*}
g^* \col u^*(R^qf'_{Y_{\bul \leq N,\os{\circ}{T}{}'_1}/S'(T')^{\nat}*}
({\cal O}_{Y_{\bul \leq N,\os{\circ}{T}{}'_1}/S'(T')^{\nat}}))_{\mab Q}
\lo R^qf_{X_{\bul \leq N,\os{\circ}{T}_1}/S(T)^{\nat}*}
({\cal O}_{X_{\bul \leq N,\os{\circ}{T}{}_1}/S(T)^{\nat}})_{\mab Q}
\tag{5.4.7.1}\label{eqn:gvbstn}
\end{equation*}  
is strictly compatible with the weight filtration.
Consequently the induced morphism 
\begin{equation*}
g^* \col u^*(R^qf'_*({\cal O}_{Y_{\bul \leq N}/K'})^{\nat})
\lo R^qf_*({\cal O}_{X_{\bul \leq N}/K})^{\nat} 
\tag{5.4.7.2}\label{eqn:gbstn}
\end{equation*} 
in $F{\textrm -}{\rm Isoc}^{\sq}(S/{\cal V})$ 
is strictly compatible with the weight filtration.
\end{theo}
\begin{proof}
Set 
$Y_{\bul \leq N,\os{\circ}{T}_1}:=Y_{\bul \leq N}\times_{\os{\circ}{S}{}'}\os{\circ}{T}_1$. 
Then we have the following equality by (\ref{coro:flft}) and 
the argument in (\ref{ali:osrgm}): 
\begin{align*} 
(R^qf_{Y_{\bul \leq N,\os{\circ}{T}_1}/S'(T)^{\nat}*}
({\cal O}_{Y_{\bul \leq N,\os{\circ}{T}_1}/S'(T)^{\nat}})_{\mab Q},P)
=(u^*(R^qf_{Y_{\bul \leq N,\os{\circ}{T}{}'_1}/S'(T')^{\nat}*}
({\cal O}_{Y_{\bul \leq N,\os{\circ}{T}{}'_1}/S'(T')^{\nat}})_{\mab Q}),P).
\end{align*} 
We may assume that $\os{\circ}{S}$ is the formal spectrum ${\rm Spf}(C)$ 
and the log structure of $S$ is split.  
It suffices to prove that 
\begin{align*} 
&P_kR^qf_{X_{\bul \leq N,\os{\circ}{T}_1}/S(T)^{\nat}*}
({\cal O}_{X_{\bul \leq N,\os{\circ}{T}_1}/S(T)^{\nat}})_K
\cap {\rm Im}(R^qf_{Y_{\bul \leq N,\os{\circ}{T}_1}/S'(T)^{\nat}*}
({\cal O}_{Y_{\bul \leq N,\os{\circ}{T}_1}/S'(T)^{\nat}})_K) 
\tag{5.4.7.3}\label{ali:bsth}\\
&={\rm Im}(P_kR^qf_{Y_{\bul \leq N,\os{\circ}{T}_1}/S'(T)^{\nat}*}
({\cal O}_{Y_{\bul \leq N,\os{\circ}{T}_1}/S'(T)^{\nat}})_K)). 
\end{align*}
By (\ref{exam:ofl}) and  (\ref{prop:nupf}),  
both hand sides on (\ref{ali:bsth}) extend to objects 
of $F{\textrm -}{\rm Isoc}(\os{\circ}{S}/{\cal V})$ as in (\ref{theo:stpfgb}). 
Let $\{t_i\}_{i=1}^j$ be a set of closed points of $\os{\circ}{S}$ 
such that, for any connected component 
$\os{\circ}{S}_i$ of $\os{\circ}{S}$, there exists a point $\os{\circ}{t}_i$  
such that $\os{\circ}{t}_i\in \os{\circ}{S}_i$.  
By \cite[(4.1)]{of} the pull-back functor $\prod_{i=1}^j\os{\circ}{t}{}^*_i\col 
{\rm Isoc}(\os{\circ}{S}/{\cal V})
\lo \prod_{i=1}^j{\rm Isoc}(\os{\circ}{t}_i/{\cal V})$ is faithful. 
By the proof of  \cite[(3.17)]{of}, we may assume that $\os{\circ}{T}$ is 
the formal spectrum of a finite extension 
${\cal V}'$ of ${\cal V}$. Let $\kap'$ be the residue field of ${\cal V}'$. 
As mentioned in the proof of (\ref{theo:e2dgfam}), 
${\cal V}'/p$ is an $\kap'$-algebra; the two pairs $X_{\bul \leq N,\os{\circ}{T}_1}$ and 
$X_{\bul \leq N,\os{\circ}{T}_1}\otimes_{{\cal V}'}\kap'\otimes_{\kap'}{\cal V}'/p$
are two deformations of $X_{\bul \leq N,\os{\circ}{T}_1}\otimes_{{\cal V}'}\kap'$; 
the obvious analogue for $Y_{\bul \leq N,\os{\circ}{T}_1}$ also holds. 
Hence, by the deformation invariance of log crystalline cohomologies with 
weight filtrations ((\ref{coro:finvliae})), we may assume 
that $\os{\circ}{T}={\rm Spf}(W(\kap'))$ and 
that $X_{\bul \leq N,\os{\circ}{T}_1}$ and $Y_{\bul \leq N,\os{\circ}{T}_1}$ 
are $N$-truncated simplicial base changes of SNCL schemes  
over the log point $({\rm Spec}(\kap'),{\mab N}\oplus \kap'{}^*\lo \kap')$. 
Hence the equality (\ref{ali:bsth}) follows from (\ref{theo:stpfgb}). 
\end{proof}

\section{Variational $p$-adic monodromy-weight conjecture II and 
variational log $p$-adic hard Lefschetz conjecture II}
\label{sec:vpmwc}
Let the notations be as in the previous section for the case $N=0$.  
In this section we give several results for 
the variational $p$-adic monodromy-weight conjecture 
((\ref{conj:rcpmc}), (\ref{eqn:grmmppd})) 
and the variational log $p$-adic hard Lefschetz conjecture 
((\ref{conj:lhlc}), (\ref{eqn:fcicpl})). 
\par 
We start with the following, which resembles Deligne's Principle B 
about the absolute Hodge cycle on an abelian variety 
over a field of characteristic $0$ (\cite{dhc})
(see also \cite[p.~103 (${}^{13}$)]{grdr}). 

\begin{lemm}\label{lemm:mwcprs} 
Let the notations be as in {\rm (\ref{conj:rcpmc})}. 
Then the following hold$:$ 
\par 
$(1)$ Let $s$ be the log point of a perfect field $\kap$ of characteristic $p>0$. 
Consider the case $S=s$. 
If the morphism {\rm (\ref{eqn:grmpd})} is injective or surjective 
for $X/S=X/s$ and $(T,{\cal J},\del)=({\cal W}(s),p{\cal W},[~])$, 
then {\rm (\ref{conj:rcpmc})} holds for $X/s$ and any log PD-enlargement 
$(T,{\cal J},\del)$ of $s$ over ${\rm Spf}({\cal W})$. 
\par 
$(2)$ Let $T$ be an object of ${\rm Enl}(S/{\cal V})$ in {\rm \S\ref{sec:cfi}}.
If the morphism {\rm (\ref{eqn:grmpd})} is injective or surjective for $X/S$ 
and a log PD-enlargement $(T,{\cal J},\del)$ of $S$, 
then {\rm (\ref{conj:rcpmc})} holds for $X_{\os{\circ}{T}_0,t}/{\cal W}(t)$ 
for any exact closed point $t\in S_{\os{\circ}{T}_0}$ 
with a lifting morphism ${\cal W}(t)\lo S(T)^{\nat}$ of $t\lo S_{\os{\circ}{T}_0}$.  
\par 
$(3)$ 
Let $T$ be an object of ${\rm Enl}(S/{\cal V})$ in {\rm \S\ref{sec:cfi}}.
Let $t$ be an exact closed point of $S_{\os{\circ}{T}_0}$.  
Assume that $\os{\circ}{T}_0$ is connected. 
If the morphism {\rm (\ref{eqn:grmpd})} 
is injective or surjective for $X_{\os{\circ}{T}_0,t}/{\cal W}(t)$, 
then the morphism {\rm (\ref{eqn:grmpd})} is an isomorphism for 
$X_{\os{\circ}{T}_0}/S_{\os{\circ}{T}_0}$ and $(T,{\cal J},\del)$. 
\end{lemm}
\begin{proof} 
(1): By the duality between $E_2^{-k,2d-q-k} $ and 
$E_2^{k,q+k}$ of the weight spectral sequence of 
$H^q_{{\rm crys}}(X/{\cal W}(s))_{\mab Q}$
((\ref{prop:dual})) and 
the $E_2$-degeneration of it (\cite[(3.6)]{ndw}, (\ref{theo:e2dgfam})), 
the injectivity (resp.~surjectivity) 
tells us the surjectivity (resp.~injectivity) of the morphism 
(\ref{eqn:grmpd}) for the case $T={\cal W}(s)$. 
By the filtered base change theorem in log crystalline cohomologies 
((\ref{theo:bccange})) and by the argument in the proof of (\ref{prop:cxst}), 
{\rm (\ref{conj:rcpmc})} holds for $X/s$ and any log PD-enlargement 
$(T,{\cal J},\del)$ of $s$ over ${\rm Spf}({\cal W})$. 
\par 
(2): 
Set $K_{0,t}:={\rm Frac}({\cal W}(\kap_t))$. 
The injectivity for $X_t/{\cal W}(t)$ follows from (\ref{coro:flft})  
and the filtered base change theorem in log crystalline cohomologies 
((\ref{theo:bccange})). 
Indeed, we obtain the following by (\ref{ali:offhfrgm}): 
\begin{align*} 
P_kH^q_{\rm crys}(X_t/{\cal W}(t))_{\mab Q}
=K_{0,t}\otimes_{{\cal K}_T}
P_kR^qf_{X_{\os{\circ}{T}_0}/S(T)^{\nat}*}
({\cal O}_{X_{\os{\circ}{T}_0}/S(T)^{\nat}})_{\mab Q}. 
\end{align*} 
Because $(P_kR^qf_{X_{\os{\circ}{T}_0}/S(T)^{\nat}*}
({\cal O}_{X_{\os{\circ}{T}_0}/S(T)^{\nat}})_{\mab Q},P)$ 
is a filteredly flat ${\cal K}_T$-module ((\ref{coro:fenlt})), 
$${\rm gr}^P_k(R^qf_{X_{\os{\circ}{T}_0}/S(T)^{\nat}*}
({\cal O}_{X_{\os{\circ}{T}_0}/S(T)^{\nat}})_{\mab Q})
\otimes_{{\cal K}_S}K_{0,t}=
{\rm gr}^P_kH^q_{\rm crys}(X_t/{\cal W}(t))_{\mab Q}.$$ 
Hence the morphism $(\ref{ali:incsx})\otimes_{{\cal K}_T}K_{0,t}$ 
is equal to the morphism 
\begin{align*} 
\nu^k_{\rm zar} \col 
{\rm gr}^P_{q+k}H^q_{\rm crys}(X_t/{\cal W}(t))_{\mab Q}  
\lo 
{\rm gr}^P_{q-k}H^q_{\rm crys}(X_t/{\cal W}(t))(-k,u)_{\mab Q}.  
\tag{5.5.1.2}\label{ali:incbcuqsx}
\end{align*}
Furthermore, if the morphism 
\begin{align*} 
\nu^k_{\rm zar} \col 
{\rm gr}^P_{q+k}R^qf_{X_{\os{\circ}{T}_0/S(T)^{\nat}}*}
({\cal O}_{X_{\os{\circ}{T}_0/S(T)^{\nat}}})_{\mab Q}  
\lo 
{\rm gr}^P_{q-k}R^qf_{X_{\os{\circ}{T}_0/S(T)^{\nat}}*}
({\cal O}_{X_{\os{\circ}{T}_0/S(T)^{\nat}}})(-k,u)_{\mab Q}  
\tag{5.5.1.1}\label{ali:incsx}
\end{align*}
is injective (resp.~surjective), 
then the morphism (\ref{ali:incbcuqsx}) is injective (resp.~surjective).  
\par 
(3): First we prove that the morphism 
(\ref{eqn:grmpd}) is injective for $X/S$ and $(T,{\cal J},\del)$. 
This is a local question on $S(T)^{\nat}$. 
Hence we may assume that $M_{S(T)^{\nat}}$ is split. 
By (\ref{rema:nlcfi}) and (\ref{prop:nupf}),  
${\rm gr}^P_kR^qf_*({\cal O}_{X_{\os{\circ}{T}_0}/S(T)^{\nat}})_{\mab Q}$ 
defines an object of ${\rm Isoc}(\os{\circ}{T}/{\cal V})$. 
Let ${\cal V}'$ be the integer ring of $K_{0,t}$.
Because 
\begin{equation*}
\os{\circ}{t}{}^* \col {\rm Isoc}(\os{\circ}{T}_0/{\cal V}) 
\lo {\rm Isoc}(\os{\circ}{t}/{\cal V}) 
\tag{5.5.1.3}\label{eqn:sppt}
\end{equation*} 
is faithful by \cite[(4.1)]{of}, to prove the injectivity of 
the morphism (\ref{eqn:grmpd}) for $X/S$ and $(T,{\cal J},\del)$, 
we may assume that $(T,{\cal J},\del)=({\cal V}',p{\cal V}',[~])$. 
In this case, the morphism is injective by the assumption and (1). 
Thus we have proved that 
the morphism (\ref{eqn:grmpd}) is injective for $X/S$ and $(T,{\cal J},\del)$.
\par 
Next we prove that 
the morphism (\ref{eqn:grmpd}) is surjective for $X/S$ and $(T,{\cal J},\del)$.
Set 
\begin{equation*} 
C:= {\rm Coker}
(\nu^e_{\rm zar}\col {\rm gr}^P_{q+e}R^qf_{*}({\cal O}_{X_{\os{\circ}{T}_0}/K})\lo 
{\rm gr}^P_{q-e}R^qf_{*}({\cal O}_{X_{\os{\circ}{T}_0}/K})(-e)).
\end{equation*}  
We would like to prove that $C_{S(T)^{\nat}}=0$. 
This is a local question on $S(T)^{\nat}$ 
and we may assume that $M_{S(T)^{\nat}}$ is split.  
By (\ref{rema:nlcfi}) and (\ref{prop:nupf}),  
$C$ gives an object of ${\rm Isoc}(\os{\circ}{T}/{\cal V})$. 
By the faithfulness of the functor (\ref{eqn:sppt}), 
it suffices to prove that $C_{{\cal V}'}=0$. 
This follows from (1) and the assumption on the injectivity. 
\end{proof}

\par 
The following is a $p$-adic analogue of Terasoma-Ito's result 
(\cite{ter}, \cite{itp}). 
As stated in the Introduction, Lazda and P\'{a}l recently 
have proved the coincidence of 
the monodromy filtration and the weight filtration 
on the log crystalline cohomology of 
a proper strict semistable family over   
a complete discrete valuation ring 
with finite residue field of equal characteristic $p>0$ 
by using their theory of rigid cohomologies over 
fields of Laurent series, their result 
``Hyodo-Kato isomorphism in equi-characteristic $p>0$'',  
the $p$-adic local monodromy theorem, 
Marmora's functor and Crew's result (\cite{lp}). 
We generalize this result in the case where 
the residue field is a general perfect field of characteristic $p>0$. 
Our proof is different from their proof and is natural 
because we have constructed the weight spectral sequence 
of the log crystalline cohomology sheaf of a proper SNCL family over 
a family of log points ((\ref{eqn:escssp})). 

\begin{theo}[{\bf Monodromy-weight conjecture in equal characteristic 
for log crystalline cohomologies}]\label{theo:aswc} 
Let $\wt{V}$ be a complete discrete valuation ring 
with perfect residue field of equal characteristic $p>0$. 
Set $\os{\circ}{D}:={\rm Spec}(\wt{V})$. 
Endow $\os{\circ}{D}$ with canonical log structure 
and let $D$ be the resulting log scheme. 
Let $s$ be the log special fiber of $D$.  
Let ${\cal X}$ be a proper strict semistable log scheme 
over $D$. Let $X$ be the log special fiber of ${\cal X}/D$.  
Then {\rm (\ref{conj:rcpmc})} holds for the case $S=s$ and 
$(T,{\cal J},\del)=({\cal W}(s),p{\cal W},[~])$. 
\end{theo} 
\begin{proof} 
Let $\tau$ be a uniformizer of $\wt{V}$. Set $\kap:=\wt{V}/\tau$. 
In this case, it is well-known (by using Neron's blow up) that 
there exists a smooth ${\mab F}_p[\tau]$-subalgebra 
$\wt{A}_1$ of $\wt{V}$ such that $\wt{V}$ is 
the henselization of $\wt{A}_1$ at the prime ideal $(\tau)$ in $\wt{A}_1$ 
and such that there exists a proper log smooth morphism 
$\wt{\cal X}\lo \wt{S}_1$ such that 
${\cal X}=\wt{\cal X}\times_{\wt{S}_1}D$ (see \cite[\S5,6]{itp}). 
Here $\wt{S}_1$ is the log scheme 
whose underlying scheme is ${\rm Spec}(\wt{A}_1)$ 
and whose log structure is associated to a morphism 
${\mab N}\owns 1\lom \tau \in \wt{A}_1$. 
Set $S_1:=\ul{\rm Spec}^{\log}_{\wt{S}_1}
({\cal O}_{\wt{S}_1}/\tau{\cal O}_{\wt{S}_1})$ and 
$A_1:=\Gam(S_1,{\cal O}_{S_1})=\wt{A}_1/\tau$. 
Then we have the following commutative diagram 
\begin{equation*} 
\begin{CD}
\wt{A}_1 @>{\sus}>> \wt{V} \\
@VVV @VVV \\
A_1 @>>> \wt{V}/\tau=\kap.  
\end{CD}
\tag{5.5.2.1}\label{cd:wtp} 
\end{equation*}
Set $X_{S_1}:=\wt{\cal X}{\times}_{\wt{S}_1}S_1$ by abuse of notation. 
We may assume that 
$X_{S_1}$ is a proper SNCL scheme 
over $S_1$ (\cite[\S4, \S6]{itp}). 
(Before \cite{itp} was published, K.~Fujiwara 
has essentially suggested to me this technique.)
Note that we have used the notion of the SNCL scheme. 
Take any exact closed point $t$ of $S_1$. 
Set $X_{S_1,t}:=X_{S_1}{\times}_{S_1}t$. 
Let $A$ be a $p$-adically formally smooth algebra over 
${\mab Z}_p$ which is a lift of $A_1$.
Endow ${\rm Spf}(A)$ with a log structure 
$({\mab N}\oplus A^*\owns (1,u)\lom pu\in A)^a$ 
and let $S$ be the 
resulting log formal scheme over 
${\rm Spf}({\mab Z}_p)=
({\rm Spf}({\mab Z}_p), {\mab Z}_p^*)$. 
The log formal scheme $S$ has a 
PD-ideal sheaf $p{\cal O}_S$, which defines 
an exact closed immersion 
$S_1 \os{\subset}{\lo} S$.
Assume that ${\rm dim}(\wt{A}_1)\geq 2$. 
In this case, there exists a smooth curve $\os{\circ}{C}_1$ in $\os{\circ}{\wt{S}}_1$ 
which passes through $\os{\circ}{t}$ and transversal to $(\tau=0)$ 
since $\wt{A}_1$ is smooth over ${\mab F}_p[\tau]$. 
In the case ${\rm dim}(\wt{A}_1)=1$, 
set 
$\os{\circ}{C}_1:={\rm Spec}(\wt{A}_1)$ 
(in \cite{itp} this case has been forgotten). 
Endow $\os{\circ}{C}_1$ with the log structure 
associated to the point $\os{\circ}{t}$. 
Let $C_1$ be the resulting log scheme. 
Then $t=C_1\times_{\wt{S}_1}S_1$ and the closed immersion 
$C_1\os{\sus}{\lo}\wt{S}_1$ is exact. 
As a summary, we have the following cartesian diagrams: 
\begin{equation*} 
\begin{CD} 
X @>{\subset}>> {\cal X} @>>> 
\wt{\cal X} @<{\supset}<< X_{S_1} @<{\supset}<< X_{S_1,t} \\ 
@VVV @VVV @VVV @VVV @VVV \\ 
s @>{\subset}>> D @>>> \wt{S}_1 
@<{\supset}<< S_1 @<{\supset}<< t\\
@. @. @A{\bigcup}AA @VV{\bigcap}V \\ 
@. @. C_1  @. S  @. .
\end{CD}
\tag{5.5.2.2}\label{cd:xssp} 
\end{equation*}  
Set $\wt{\cal X}_{C_1}:=\wt{\cal X}\times_{\wt{S}_1}C_1$.   
Because $\os{\circ}{C}_1$ is a smooth curve over ${\mab F}_p$ containing $t$ 
and because $\wt{\cal X}_{C_1}/C_1$ is smooth outside $t$, 
(\ref{conj:rcpmc}) holds for 
$H^q_{\rm crys}(X_{S_1,t}/{\cal W}(\kap_t))_{\mab Q}=
H^q_{\rm crys}(\wt{\cal X}_{C_1,t}/{\cal W}(\kap_t))_{\mab Q}$ 
as stated in \cite[p.~242 (e)]{fao} and \cite[(5.1)]{ctze} 
using Crew's result about the purity of 
the graded module of the monodromy filtration 
on $H^q_{\rm crys}(X_{S_1,t}/{\cal W}(\kap_t))_{\mab Q}$ 
(\cite[(10.8)]{cr}).  
\par 
For a connected affine log formal open subscheme 
$U$ of $S$ containing $t$,  
set $U_1:=U{\otimes}_{{\mab Z}_p}{\mab F}_p$. 
Set also $X_{U_1}:=X_{S_1}\times_{S_1}U_1$, 
$X_{U_1,t}:=X_{S_1}\times_{S_1}U_1\times_{U_1}t=X_{S_1,t}$
$\os{\circ}{X}{}^{(k)}_{U_1}:=
\os{\circ}{X}{}^{(k)}_{S_1}\times_{\os{\circ}{S_1}}
\os{\circ}{U}_1$ 
and 
$\os{\circ}{X}{}^{(k)}_{U_1,t}:=
\os{\circ}{X}{}^{(k)}_{U_1}\times_{\os{\circ}{U_1}}\os{\circ}{t}$. 
Using the lower horizontal morphism in (\ref{cd:wtp}), 
we have the natural inclusion morphism ${\cal O}_{U_1} \os{\sus}{\lo} \kap$. 
We fix an endomorphism of $\os{\circ}{U}$ lifting  
the Frobenius endomorphism(=$p$-th power endomorphism) 
of $\os{\circ}{U}_1$. 
Then we have the Teichm\"{u}ller lift 
${\cal O}_U\lo {\cal W}(\kap_t)$ 
(resp.~${\cal O}_U\lo {\cal W}$)  
of the morphism ${\cal O}_{U_1}\lo \kap_t$ 
(resp.~${\cal O}_{U_1}\lo \kap$). 
The rings ${\cal W}(\kap_t)$ and  
${\cal W}$ become ${\cal O}_U$-algebras by these lifts. 
Endow ${\rm Spf}({\cal W}(\kap_t))$ and $\os{\circ}{t}$
with the inverse images of the log structure of $U$.  
Then, by the convergence of the isocrystalline cohomology (\cite[(3.7)]{of}), 
the morphism 
\begin{equation*}
H^q_{\rm crys}(\os{\circ}{X}{}^{(k)}_{U_1}/U)_{\mab Q}
{\otimes}_{{\cal O}_U} {\cal W}(\kap_t) {\lo}
H^q_{\rm crys}
(\os{\circ}{X}{}^{(k)}_{U_1,t}/{\cal W}(\kap_t))_{\mab Q}
\tag{5.5.2.3}\label{eqn:noncosp}
\end{equation*} 
is an isomorphism for all $q,k\in {\mab N}$. Let the notations be as in 
the proof of (\ref{theo:e2dgfam}). 
By the description of 
$d^{-k,q+k}_1\col E_1^{-k,q+k}(X_{U_1}/U)_{\mab Q}\lo 
E_1^{-k+1,q+k}(X_{U_1}/U)_{\mab Q}$
((\ref{cd:gsmsd})), $d^{-k,q+k}_1$ extends to 
a morphism of convergent isocrystals. 
Hence we have the following isomorphism 
\begin{equation*} 
E_2^{-k,q+k}(X_{U_1}/U)_{\mab Q}{\otimes}_{{\cal O}_U}{\cal W}(\kap_t)
\os{\sim}{\lo} E_2^{-k,q+k}(X_{U_1,t}/{\cal W}(\kap_t))
_{\mab Q}.  
\tag{5.5.2.4}\label{eqn:usacp} 
\end{equation*}
The morphism $N$ is compatible with the base change ((\ref{prop:bcmop})). 
Hence we have the following commutative diagram: 
\begin{equation*} 
\begin{CD} 
E_2^{-k,q+k}(X_{U_1}/U)_{\mab Q}
{\otimes}_{{\cal O}_U}{\cal W}(\kap_t)
@>{N^k\otimes 1}>>
E_2^{k,q-k}(X_{U_1}/U)_{\mab Q}
{\otimes}_{{\cal O}_U}{\cal W}(\kap_t)\\
@V{\simeq}VV @VV{\simeq}V \\ 
E_2^{-k,q+k}(X_{U_1,t}/{\cal W}(t))_{\mab Q} 
@>{N^k}>> 
E_2^{k,q-k}(X_{U_1,t}/{\cal W}(t))_{\mab Q}.  
\end{CD}
\tag{5.5.2.5}\label{eqn:usuqcp}  
\end{equation*}  
By the $E_2$-degeneration 
(\cite[(3.6)]{ndw}, see also (\ref{theo:e2dgfam})) 
and the coincidence of the monodromy filtration 
and the weight filtration on $H^q_{\rm crys}(X_{U_1,t}/{\cal W}(t))_{\mab Q}$ 
as stated before, the lower morphism $N^k$ 
in (\ref{eqn:usuqcp}) is an isomorphism. 
Because the monodromy operator is a morphism of 
convergent isocrystals ((\ref{prop:convmon})) 
and because the faithfulness of the functor 
(\ref{eqn:sppt}), the morphism 
\begin{align*}
N^k \col    
E_2^{-k,q+k}(X_{U_1}/U)_{\mab Q} \lo 
E_2^{k,q-k}(X_{U_1}/U)_{\mab Q} 
\tag{5.5.2.6}\label{align:e2kqxs} 
\end{align*} 
is injective, hence bijective by (\ref{lemm:mwcprs}) (3). 
\par 
By Deligne's remark in \cite[(3.10)]{isp} (cf.~\cite[\S3]{ndw}), 
there exists a connected affine nonempty log formal open subscheme $U$ of $S$ 
such that the canonical morphism 
\begin{equation*}
H^q_{\rm crys}
(\os{\circ}{X}{}^{(k)}_{U_1}/U)
\otimes_{{\cal O}_U} {\cal W} \lo
H^q_{\rm crys}(\os{\circ}{X}{}^{(k)}/{\cal W})
\tag{5.5.2.7}\label{eqn:nonckuosp}
\end{equation*} 
is an isomorphism of ${\cal W}$-modules. 
In fact, there exists 
a connected affine nonempty log formal open subscheme $U$ of $S$ 
such that the canonical morphism 
\begin{align*} 
E_2^{-k,q+k}(X_{U_1}/U)
{\otimes}_{{\cal O}_U}{\cal W}\lo 
E_2^{-k,q+k}(X/{\cal W}(s))
\tag{5.5.2.8}\label{ali:uw}
\end{align*} 
is an isomorphism as in the proof of \cite[(3.5)]{ndw}. 
Hence we have the following diagram 
\begin{equation*} 
\begin{CD} 
E_2^{-k,q+k}(X_{U_1}/U)_{\mab Q}
{\otimes}_{{\cal O}_U}{\cal W}
@>{N^k\otimes 1,~\simeq}>>
E_2^{k,q-k}(X_{U_1}/U)_{\mab Q}
{\otimes}_{{\cal O}_U}{\cal W}\\
@V{\simeq}VV @VV{\simeq}V \\ 
E_2^{-k,q+k}(X/{\cal W}(s))
_{\mab Q} 
@>{N^k}>> E_2^{k,q-k}(X/{\cal W}(s))
_{\mab Q}.  
\end{CD}
\tag{5.5.2.9}\label{eqn:usuaqcp}  
\end{equation*}  
Now it is clear that the morphism 
\begin{equation*} 
N^k\col E_2^{-k,q+k}(X/{\cal W}(s))_{\mab Q} 
\lo 
E_2^{k,q-k}(X/{\cal W}(s))_{\mab Q} 
\tag{5.5.2.10}\label{eqn:usuzqcp}  
\end{equation*} 
is an isomorphism. We completes the proof.
\end{proof}

\par 
The following resembles Deligne's Principle B 
about the absolute Hodge cycle 
on an abelian variety over a field of characteristic $0$. 
This is a generalized $p$-adic version 
of Fujiwara's try to 
reduce the $l$-adic monodromy-weight conjecture in the case of 
a complete discrete valuation ring in mixed characteristics 
to the $l$-adic monodromy-weight conjecture in the case of 
a complete discrete valuation ring in equal characteristic. 
(See also Scholze's work \cite{sc}.)


\begin{coro}[{\bf Log ($p$-adically) convergent 
monodromy-weight conjecture in equal characteristic 
for log crystalline cohomologies}]\label{coro:anmwc} 
Let the notations be as in {\rm (\ref{lemm:mwcprs}) (3)}.  
Assume that $\os{\circ}{T}$ is connected.  
Assume also that there exists an exact closed point $t$ of 
$S_{\os{\circ}{T}_0}$ such that $X_{\os{\circ}{T}_0,t}$  
is the log special fiber of a proper strict semistable family over a 
complete discrete valuation ring of equal characteristic $p>0$. 
Then {\rm (\ref{conj:ccpwmc})} holds for 
$F{\textrm -}{\rm Isoc}^{\sq}({\rm Enl}(S/{\cal V}))$. 
\end{coro} 
\begin{proof}  
(\ref{coro:anmwc}) follows from (\ref{theo:aswc}) and (\ref{lemm:mwcprs}) (3). 
\end{proof}

\par 
In \cite{nd} Nakayama has reduced 
the log $l$-adic monodromy-weight conjecture 
for any base field to the conjecture for finite fields 
by using a specialization argument. 
The following is a variational $p$-adic analogue of his result. 

\begin{theo}\label{theo:nrlp}  
If {\rm (\ref{conj:rcpmc})} holds for the case where 
$(\os{\circ}{T},{\cal J},\del)$ is the formal spectrum of the Witt ring 
of any finite field with canonical PD-structure on $p{\cal O}_T$, 
then {\rm (\ref{conj:ccpwmc})} is true in $F{\textrm -}{\rm Isoc}^{\sq}(S/{\cal V})$. 
\end{theo}
\begin{proof}  
By using (\ref{lemm:mwcprs}) (3) and the argument in the proof of (\ref{theo:aswc}),  
the proof of this theorem is the same as that of \cite[(2.18.2)]{nh2}. 
\end{proof}

The following is Kajiwara's result (\cite{kaji}) using 
results of M.~Saito's result (\cite{samp}) and using hyperplane sections and 
Hodge index theorem for surfaces. 

\begin{theo}[{\rm {\bf \cite{kaji} (cf.~\cite[\S9]{nlpi}}}]\label{theo:kajia} 
Let $s$ be a log point of a field of characteristic $p>0$. 
Let $X$ be a projective SNCL scheme over $s$. 
Let $l\not= p$ be a prime. Then the monodromy filltration and the weight filtration 
on the Kummer log \'{e}tale cohomology 
$H^1_{\rm et}(X_{\ol{s}},{\mab Q}_l)$ {\rm (\cite{nd})} are the same. 
Moreover, if $\os{\circ}{X}$ is of pure dimension $d$, then the filtrations on 
$H^{2d-1}_{\rm et}(X_{\ol{s}},{\mab Q}_l)$ are the same. 
\end{theo}

\begin{coro}\label{coro:hhd}
Let the notations be as in {\rm (\ref{theo:kajia})}. 
Assume that the field in {\rm (\ref{theo:kajia})} is perfect. 
Then the monodromy filltration and the weight filtration 
on $H^1_{\rm crys}(X/{\cal W}(s))_{\mab Q}$ are the same. 
Moreover, if $\os{\circ}{X}$ is of pure dimension $d$, then the filtrations on 
$H^{2d-1}_{\rm crys}(X/{\cal W}(s))_{\mab Q}$ are the same. 
\end{coro}
\begin{proof}  
In \cite[\S9]{nlpi} we have proved that 
the $l$-adic quasi-monodromy operator 
$\nu_l \col E_2^{-1,1}\lo E_2^{1,0}$ and the 
the $p$-adic quasi-monodromy operator 
$\nu_p \col E_2^{-1,1}\lo E_2^{1,0}$ are defined over ${\mab Q}$ 
and they are the same. Hence (\ref{coro:hhd}) for 
$H^1_{\rm crys}(X/{\cal W}(s))_{\mab Q}$ follows from (\ref{theo:kajia}); 
(\ref{coro:hhd}) for 
$H^{2d-1}_{\rm crys}(X/{\cal W}(s))_{\mab Q}$ follows from the duality 
((\ref{prop:quasimon}), (\ref{prop:dual})).  
\end{proof}

\begin{coro}\label{coro:kjat}
The conjecture {\rm (\ref{conj:ccpwmc})} is true for the case $q=1$ and $q=2d-1$ 
if $\os{\circ}{X}$ is of pure dimension $d$.
\end{coro} 
\begin{proof} 
This follows from (\ref{lemm:mwcprs}) and (\ref{coro:hhd}). 
\end{proof} 

\begin{theo}\label{theo:22}
Let the notations be as in {\rm (\ref{conj:rcpmc})}. 
If the relative dimension of $X/S_0$ is less than or equal to $2$, 
then the conjecture {\rm (\ref{conj:rcpmc})} is true. 
\end{theo} 
\begin{proof} 
This follows from (\ref{coro:kjat}), 
\cite[Lemme 6.2.1, 6.2.2]{msemi} 
and the convergences of the weight filtration and 
the monodromy operator. 
\end{proof}

Next we discuss the log $p$-adic hard Lefschetz conjecture. 

\begin{lemm}\label{lemm:mlhs} 
Let the notations be as in {\rm (\ref{lemm:mwcprs})}.  
Assume that the morphism $\os{\circ}{X}\lo \os{\circ}{S}$ is projective and 
that the relative dimension is of pure dimension $d$. 
Then the following hold$:$ 
\par 
$(1)$ Consider the case $S=s$. 
If the morphism {\rm (\ref{eqn:fcpl})} is injective or surjective 
for $X/S=X/s$ and $(T,{\cal J},\del)=({\cal W}(s),p{\cal W},[~])$, 
then {\rm (\ref{conj:lhlc})} holds for $X/s$ and any log PD-enlargement 
$(T,{\cal J},\del)$ of $s$ over ${\rm Spf}({\cal W})$ of $t\lo S_{\os{\circ}{T}_0}$. 
\par 
$(2)$ Let $T$ be an object of ${\rm Enl}(S/{\cal V})$ in {\rm \S\ref{sec:cfi}}.
If the morphism {\rm (\ref{eqn:fcpl})} is injective or surjective for $X/S$ 
and a log PD-enlargement $(T,{\cal J},\del)$ of $S$, 
then {\rm (\ref{conj:lhlc})} holds for $X_{\os{\circ}{T}_0,t}/{\cal W}(t)$ 
for any exact closed point $t\in S_{\os{\circ}{T}_0}$ 
with a lifting morphism ${\cal W}(t)\lo S(T)^{\nat}$.  
\par 
$(3)$ 
Let $T$ be as in $(2)$. 
Let $t$ be an exact closed point of $S_{\os{\circ}{T}_0}$.  
Assume that $\os{\circ}{T}_0$ is connected. 
If the morphism {\rm (\ref{eqn:fcpl})} 
is injective or surjective for $X_{\os{\circ}{T}_0,t}/{\cal W}(t)$, 
then the morphism {\rm (\ref{eqn:fcpl})} is an isomorphism for 
$X_{\os{\circ}{T}_0}/S_{\os{\circ}{T}_0}$ and $(T,{\cal J},\del)$. 
\end{lemm}
\begin{proof} 
(1): (1) follows from the log Poincar\'{e} duality of Tsuji (\cite[Theorem (5.6)]{tsp}).  
\par 
(2): Assume that the morphism (\ref{eqn:fcpl}) is injective (resp.~surjective) 
for $X/S$ and $(T,{\cal J},\del)$.  
Then the injectivity (resp.~surjectivity) 
of the morphism {\rm (\ref{eqn:fcpl})} 
for $X_s/{\cal W}(s)$ follows from 
the flatness of $R^qf_{X_{\os{\circ}{T}_0}/S(T)^{\nat}*}
({\cal O}_{X_{\os{\circ}{T}_0}/S(T)^{\nat}})
\otimes_{\mab Z}{\mab Q}$ already 
used in the proof of (\ref{theo:e2dam}) and 
the base change theorem in log crystalline cohomologies 
(see (\ref{ali:osrgm})).  
\par 
(3): By using (1), the proof is the same as that of (\ref{lemm:mwcprs}) (3). 
\end{proof}

Let us recall the following Faltings' theorem: 

\begin{theo}[{\bf \cite[2 (f)]{fao}}]\label{theo:faltm}
Let $\kap$ be a perfect field of characteristic $p>0$. 
Let ${\cal W}$ be the Witt ring of $\kap$. 
Let $K_0$ be the fraction field of ${\cal W}$. 
Let $f\col (Y,D)\lo (Z,E)$ be a proper log smooth morphism of 
smooth schemes with SNCD's over $\kap$. 
Then there exists a log convergent $F$-isocrystal 
$R^qf_*({\cal O}_{(Y,D)/K_0})$ on $(Z,E)/{\cal W}$ 
such that, for any object $(T,z)\in {\rm Enl}_p((Z,E)/{\cal W})$, 
\begin{align*} 
R^qf_*({\cal O}_{(Y,D)/K_0})_T=R^qf_{(Y,D)_{T_1}/T*}({\cal O}_{(Y,D)_{T_1}/T})_{\mab Q}. 
\end{align*} 
For any object  $(T,z)\in {\rm Enl}((Z,E)/{\cal W})$, 
$R^qf_*({\cal O}_{(Y,D)/K_0})_T$ is a locally free ${\cal K}_T$-module of finite rank, 
where ${\cal K}_T:={\cal O}_T\otimes_{\mab Z}{\mab Q}$. 
\end{theo}

\begin{rema} 
(1) If $Z$ is a smooth curve over $\kap$, 
then (\ref{theo:faltm}) has also been proved in  
\cite[(3.12)]{ctcs}. 
\par
(2) The log structure of $(Z,E)$ is not constant. 
\end{rema} 

The following is an obvious log version of \cite[(1.11)]{of}: 

\begin{defi} 
Let $B$ be as in the beginning of \S\ref{sec:cfi}. 
Let $T$ be a fine log $p$-adic formal $B$-scheme. 
A coherent crystal of ${\cal O}_{T/B}\otimes_{\mab Z}{\mab Q}$ is the following datum:
\parno  
For any open log formal scheme $U$ of $T$ and any   
exact closed nilpotent immersion $U\os{\sus}{\lo} U'$ into 
a fine log $p$-adic formal ${\cal V}$-scheme, an object 
$E_{U'}\in {\rm Coh}({\cal K}_{U'})$ satisfying the following condition: 
\parno 
For any morphism $g \col  (U_1, U'_1)\lo (U_2, U'_2)$ of 
exact closed nilpotent immersions of  
fine log $p$-adic formal $B$-schemes, where $U_i$ $(i=1,2)$ 
is an open log formal scheme of $T$, 
we are given an isomorphism $\iota_g\col g^*(E_{U'_2})\lo E_{U'_1}$ 
satisfying the usual two relations 
$\iota_g\circ g^*(\iota_h)=\iota_{h\circ g}$ for a similar morphism 
$h \col  (U_2, U'_2)\lo (U_3, U'_3)$ to $g$ and 
$\iota_{{\rm id}_{(U,U')}}={\rm id}_{E_{U'}}$.  
\end{defi} 

The following is the obvious log version of \cite[(1.12)]{of}: 

\begin{prop}\label{prop:prodes} 
Let $T\os{\sus}{\lo} Y$ be an immersion 
into a log smooth $p$-adic formal $B$-scheme. 
Set $Y(1):=Y\times_BY$. 
Let $Y(1)^{\rm ex}$ be the exactification of the diagonal immersion 
$Y\os{\sus}{\lo} Y(1)$. 
Let $Y_n$ and $Y(1)_n$ be the $n$-th infinitesimal neighborhoods 
of the immersion $T\os{\sus}{\lo} Y$ and the composite immersion 
$T\os{\sus}{\lo} Y\os{\sus}{\lo} (Y\times_BY)^{\rm ex}$, respectively.
Let $p_i \col Y(1)_n \lo Y_n$ $(i=1,2)$ be the induced morphism 
from the $i$-th projection $Y(1)\lo Y$. 
Let $E$ be a coherent crystal of ${\cal O}_{T/B}\otimes_{\mab Z}{\mab Q}$-modules. 
Then we have the isomorphism 
$\eps_n \col p_2^*(E_{Y_n}) \os{\sim}{\lo} E_{Y(1)_n} \os{\sim}{\longleftarrow}  p_1^*(E_{Y_n})$ 
satisfying the standard cocycle condition. 
The correspondence $E\lom \{(E_{Y_n},\eps_n)\}_{n=1}^{\infty}$ is 
an equivalence of categories. 
\end{prop} 
\begin{proof} 
Because the proof is the obvious log version of  \cite[(1.12)]{of}, 
we leave the proof to the reader. 
\end{proof} 

\begin{lemm}\label{lemm:cpff}
Let ${\cal V}$ be a complete discrete valuation ring of mixed characteristics $(0,p)$. 
Let $T$ be a log smooth connected $p$-adic formal scheme over ${\rm Spf}({\cal V})$. 
Assume that $\os{\circ}{T}$ is smooth over ${\rm Spf}({\cal V})$. 
Let $E$ be a coherent crystal of ${\cal O}_{T/{\cal V}}\otimes_{\mab Z}{\mab Q}$-modules. 
Let ${\cal V}'$ be a finite extension of ${\cal V}$. 
Let $\tau \col S \os{\sus}{\lo} T$ be an exact closed immersion such that 
$\os{\circ}{S}={\rm Spf}({\cal V}')$. Let $e$ be a global section of $E$. 
If $\tau^*(e)=0$, then $e=0$. 
\end{lemm}
\begin{proof} 
By (\ref{prop:prodes}) it suffices to prove that $e_T=0$. 
Then the rest of the proof is the same as that of \cite[(1.18)]{of}. 
\end{proof}

\begin{theo}[{\bf Log hard Lefschetz conjecture of 
log isocrystalline cohomologies in equal characteristic}]\label{theo:llccec}
Let the notations be as in {\rm (\ref{theo:aswc})}.  
Assume that the structural morphism 
$\os{\circ}{f} \col \os{\circ}{\cal X} \lo \os{\circ}{D}$ 
is projective and the relative dimension of 
$\os{\circ}{f}$ is of pure dimension $d$. 
Let ${\cal L}$ be a relatively ample line bundle on 
$\os{\circ}{\cal X}/\os{\circ}{D}$. 
Let $L$ be the restriction of ${\cal L}$ to $\os{\circ}{X}/\os{\circ}{s}$.  
Then {\rm (\ref{conj:lhlc})} holds for the case $S=s$ and 
$(T,{\cal J},\del)=({\cal W}(s),p{\cal W},[~])$. 
\end{theo} 
\begin{proof} 
Let the notations be as in the proof of (\ref{theo:aswc}). 
We may assume that ${\cal L}$ comes from 
a relatively ample line bundle $\wt{\cal L}$ 
on $\os{\circ}{\wt{\cal X}}/\os{\circ}{\wt{S}}_1$ by using Neron's blow up 
(\cite[(4.6)]{a}, \cite[I (0.5)]{sga7}). 
Because $\os{\circ}{C}_1$ is a smooth curve over ${\mab F}_p$, 
$\os{\circ}{C}_1$ has a smooth lift $\os{\circ}{C}$ 
over ${\rm Spec}({\mab Z}_p)$. 
Let $({\os{\circ}{C}})^{\wh{}}$ be the $p$-adic completion of $\os{\circ}{C}$. 
Endow $({\os{\circ}{C}})^{\wh{}}$ with the PD-ideal sheaf $p{\cal O}_{({\os{\circ}{C}})^{\wh{}}}$ with 
canonical PD-structure $[~]$ 
(this PD-structure has been forgotten in \cite{ctcs}). 
Since $({\os{\circ}{C}})^{\wh{}}$ is formally smooth 
over ${\rm Spf}({\mab Z}_p)$, the section $\os{\circ}{t} \os{\sus}{\lo} \os{\circ}{C}_1$ 
extends to a section ${\cal W}(\os{\circ}{t}) \os{\sus}{\lo} ({\os{\circ}{C}})^{\wh{}}$. 
Endow $({\os{\circ}{C}})^{\wh{}}$ with the associated log structure to 
a formal local parameter of ${\cal W}(\os{\circ}{t})$ 
(cf.~\cite[(2.1.8.1), (2.1.9)]{nh2}; ${\cal W}(\os{\circ}{t})$ is a formal SNCD on $({\os{\circ}{C}})^{\wh{}}$)
and let $\wh{C}$ be the resulting log formal scheme. 
Then we have a natural exact closed immersion 
${\cal W}(t)\lo \wh{C}$. 
Let $f\col \wt{\cal X}_{C_1}\lo C_1\os{\sus}{\lo} \wh{C}$ 
be the structural morphism. 
Then we have the log crystalline cohomology sheaf 
$R^qf_{\wt{\cal X}_{C_1}/\wh{C}*}({\cal O}_{\wt{\cal X}_{C_1}/\wh{C}})$ 
$(q\in {\mab Z})$. 
Set ${\cal K}_{\wh{C}}:={\cal O}_{\wh{C}}\otimes_{{\mab Z}}{\mab Q}$. 
By (\ref{theo:faltm}),  
$R^qf_{\wt{\cal X}_{C_1}/\wh{C}*}({\cal O}_{\wt{\cal X}_{C_1}/\wh{C}})_{{\mab Q}}$ 
is a flat coherent sheaf of ${\cal K}_{\wh{C}}$-modules and this extends to 
a log convergent $F$-isocrystal 
$R^qf_{\wt{\cal X}_{C_1}*}({\cal O}_{\wt{\cal X}_{C_1}/{\mab Q}_p})$ 
on $C_1=(\os{\circ}{C}_1,\os{\circ}{t})/{\mab Z}_p$.  
This implies that the cohomology sheaf 
$R^qf_{\wt{\cal X}_{C_1}*}({\cal O}_{\wt{\cal X}_{C_1}/{\mab Q}_p})_{\wh{C}}$ for 
any object $T\in {\rm Enl}(C_1/{\mab Z}_p)$ 
defines a log crystal of ${\cal O}_{\wh{C}/{\mab Z}_p}\otimes_{\mab Z}{\mab Q}$-modules.  
Indeed, for a nilpotent exact closed immersion $U\os{\sus}{\lo} U'$ 
over ${\rm Spf}({\mab Z}_p)$ ($U$ is a log open formal subscheme of $\wh{C}$), 
$U_0=U'_0:=\ul{\rm Spec}^{\log}({\cal O}_{U'}/p{\cal O}_{U'})_{\rm red}$ 
as in \cite[(2.8.1)]{of}; $(U'_0,U')$ is an object of ${\rm Enl}(\wh{C}/{\mab Z}_p)$.  
\par 
Let $\wt{\cal L}_{C_1}$ be the restriction of $\wt{\cal L}$ to $\wt{\cal X}_{C_1}$. 
By the same argument as that of  
\cite[\S3]{boi} (see also \cite[p.~130]{nh2}), 
we obtain the log cohomology class 
$\wt{\eta}=c_{1,{\rm crys}}(\wt{\cal L}_{C_1})$ of $\wt{\cal L}_{C_1}$
in $R^2f_{\wt{\cal X}_{C_1}/\wh{C}*}({\cal O}_{\wt{\cal X}_{C_1}/\wh{C}})$.   
Because $\wt{\eta}$ commutes with base change, 
the morphism 
\begin{align*} 
\wt{\eta}^i 
\col 
R^{d-i}f_{\wt{\cal X}_{C_1}/\wh{C}*}({\cal O}_{\wt{\cal X}_{C_1}/\wh{C}})_{{\mab Q}} 
\lo R^{d+i}f_{\wt{\cal X}_{C_1}/\wh{C}*}({\cal O}_{\wt{\cal X}_{C_1}/\wh{C}})(i)_{{\mab Q}} 
\quad (i\in {\mab N})
\tag{5.5.15.1}\label{eqn:fcacpl} 
\end{align*}
prolongs to a global section  $\wt{\eta}^i_{{\mab Q}_p}:=\{\wt{\eta}^i_T\}_T$, where 
$\wt{\eta}^i_T$ is a section of  
\begin{align*} 
{\rm Hom}_{{\cal K}_T}(R^{d-i}f_{\wt{\cal X}_{T_1}/T*}({\cal O}_{\wt{\cal X}_{T_1}/T}), 
R^{d+i}f_{\wt{\cal X}_{T_1}/T*}({\cal O}_{\wt{\cal X}_{T_1}/T})) 
\quad (i\in {\mab N}).
\tag{5.5.15.2}\label{eqn:fcvcpl} 
\end{align*}
By (\ref{theo:faltm}) we can consider $\wt{\eta}^i_{{\mab Q}_p}$ 
as a global section of a crystal of ${\cal O}_{\wh{C}/{\mab Z}_p}\otimes_{\mab Z}{\mab Q}$-modules.  
(The Hom in (\ref{eqn:fcvcpl}) has an internal Hom.) 
Let $t'\not= \os{\circ}{t}$ be a closed point of $\os{\circ}{C}_1$. 
Set $\kap_{t'}:=\Gam(t',{\cal O}_{t'})$.  
Then we have a section ${\rm Spf}({\cal W}(\kap_{t'}))\os{\sus}{\lo} \wh{C}$ of log formal schemes 
over ${\rm Spf}({\mab Z}_p)$. 
The section gives a log $p$-adic enlargement of 
$\wh{C}$ over ${\rm Spf}({\mab Z}_p)$. 
Set $\wt{\cal X}_{t}:=\wt{\cal X}_{C_1}\times_{C_1}t=X_{S_1,t}$ 
and $\wt{\cal X}_{t'}:=\wt{\cal X}_{C_1}\times_{C_1}t'$. 
Consider the morphism 
\begin{align*} 
\wt{\eta}^i_{{\cal W}(t')} \col 
R^{d-i}f_{\wt{\cal X}_{t'}/{\cal W}(t')*}({\cal O}_{\wt{\cal X}_{t'}/{\cal W}(t')})_{{\mab Q}} 
\lo R^{d+i}f_{\wt{\cal X}_{t'}/{\cal W}(t')*}({\cal O}_{\wt{\cal X}_{t'}/{\cal W}(t')})(i)_{{\mab Q}}.  
\tag{5.5.15.3}\label{eqn:fctcpl} 
\end{align*}
This morphism is an isomorphism by a result of Katz-Messing 
(\cite[Corollary 1.~2)]{kme} (see also (\ref{rema:nllfcs}) (1))).   
Hence $({\rm Ker}(\wt{\eta}^i_{{\mab Q}_p}))_{{\cal W}(t')}=0
=({\rm Coker}(\wt{\eta}^i_{{\mab Q}_p}))_{{\cal W}(t')}$. 
By 
these equalities and (\ref{lemm:cpff}), 
${\rm Ker}(\wt{\eta}^i_{{\mab Q}_p})=0
={\rm Coker}(\wt{\eta}^i_{{\mab Q}_p})$. 
In particular, 
$({\rm Ker}(\wt{\eta}^i_{{\mab Q}_p}))_{{\cal W}(t)}=0
=({\rm Coker}(\wt{\eta}^i_{{\mab Q}_p}))_{{\cal W}(t)}$. 
This means that the morphism 
\begin{align*} 
\wt{\eta}^i_{{\cal W}(t)} \col 
&R^{d-i}f_{X_{S_1,t}/{\cal W}(t)*}({\cal O}_{X_{S_1,t}/{\cal W}(t)})_{{\mab Q}} = 
R^{d-i}f_{\wt{\cal X}_{t}/{\cal W}(t)*}({\cal O}_{\wt{\cal X}_{t}/{\cal W}(t)})_{{\mab Q}} 
\tag{5.5.15.4}\label{eqn:fctczpl} \\ 
&\lo R^{d+i}f_{\wt{\cal X}_{t}/{\cal W}(t)*}({\cal O}_{\wt{\cal X}_{t}/{\cal W}(t)})(i)_{{\mab Q}} 
=R^{d+i}f_{X_{S_1,t}/{\cal W}(t)*}({\cal O}_{X_{S_1,t}/{\cal W}(t)})(i)_{{\mab Q}}
\end{align*}
is an isomorphism. 
\par 
For an open log subscheme $U$ of $S$, set $U_1:=U\times_SS_1$ 
and $X_{U_1}:=X_{S_1}\times_{S_1}U_1$. 
By \cite[(3.4)]{ndw} there exists a connected affine nonempty open log formal 
subscheme $U$ of $S$ such that the canonical morphism 
\begin{equation*}
R^qf_{X_{U_1}/U*}({\cal O}_{X_{U_1}/U})
{\otimes}_{{\cal O}_U} {\cal W}(\kap_{t'}) \lo
R^qf_{X_{S_1,t'}/{\cal W}(t')*}({\cal O}_{X_{S_1,t'}/{\cal W}(t')})
\tag{5.5.15.5}\label{eqn:nonwcosp}
\end{equation*}
for all $q\in {\mab N}$ and any exact closed point $t'$ of $U_1$ 
is an isomorphism of ${\cal W}(\kap_{t'})$-modules. 
Hence, by (\ref{lemm:mlhs}) (3), 
the morphism 
\begin{align*} 
\wt{\eta}^i \col 
R^{d-i}f_{X_{U_1}/U*}({\cal O}_{X_{U_1}/U})_{{\mab Q}} 
\lo R^{d+i}f_{X_{U_1}/U*}({\cal O}_{X_{U_1}/U})(i)_{{\mab Q}}
\tag{5.5.15.6}\label{eqn:fctcoxzpl} 
\end{align*}
is an isomorphism. 
Since we have a natural morphism $s\lo S_1$, 
we have a natural morphism $s\lo U_1$. 
Because we can assume that $U$ has an endomorphism 
lifting the Frobenius endomorphism of $U_1$, 
we have the Teichm\"{u}ller lift ${\cal W}(s)\lo U$. 
By \cite[(3.4)]{ndw}, there exists a connected affine nonempty open log formal 
subscheme $U$ of $S$ such that the canonical morphism 
\begin{equation*}
R^qf_{X_{U_1}/U*}({\cal O}_{X_{U_{1}}/U})
{\otimes}_{{\cal O}_U} {\cal W} {\lo}
R^qf_{X/{\cal W}(s)*}({\cal O}_{X/{\cal W}(s)})
\tag{5.5.15.7}\label{eqn:nonlcosp}
\end{equation*}
for all $q\in {\mab N}$ is an isomorphism of ${\cal W}$-modules. 
Therefore  
the morphism 
\begin{align*} 
\wt{\eta}^i \col 
R^{d-i}f_{X/{\cal W}(s)*}({\cal O}_{X/{\cal W}(s)})_{{\mab Q}} 
\lo R^{d+i}R^qf_{X/{\cal W}(s)*}({\cal O}_{X/{\cal W}(s)})(i)_{{\mab Q}}
\tag{5.5.15.8}\label{eqn:fctuzpl} 
\end{align*}
is an isomorphism. 
\end{proof}

\begin{coro}[{\bf Log convergent log hard Lefschetz conjecture of 
log isocrystalline cohomologies in equal characteristic}]\label{coro:lclt}
Let the notations be as in {\rm (\ref{lemm:mlhs}) (3)}.  
Assume that $\os{\circ}{T}$ is connected.  
Assume also that there exists an exact closed point $t$ of 
$S_{\os{\circ}{T}_0}$ such that $X_{\os{\circ}{T}_0,t}$  
is the log special fiber of a projective strict semistable family of pure dimension 
over a complete discrete valuation ring of equal characteristic $p>0$. 
Then the morphism {\rm (\ref{eqn:fcicpl})} is an isomorphism.  
\end{coro} 
\begin{proof} 
(\ref{coro:lclt}) follows from (\ref{lemm:mlhs}) (3), (\ref{theo:llccec}) and (\ref{prop:oxlp}). 
\end{proof}

We recall the following result, which is stated in \cite[p.~167]{nlpi}: 

\begin{theo}[{\bf Log hard Lefschetz conjecture of 
log isocrystalline cohomologies in mixed characteristics}]\label{theo:llcmec}
Let the notations be as in {\rm (\ref{theo:llccec})}, where 
$\os{\circ}{D}$ is assumed to be a complete discrete valuation ring of 
mixed characteristics $(0,p)$.  
Then {\rm (\ref{conj:lhlc})} holds for the case 
$S=s$ and $T={\cal W}(s)$.   
\end{theo} 

As is remarked in \cite[p.~167]{nlpi}, 
this theorem follows from from \cite[(5.1)]{hk}, 
by proofs in \cite[\S3]{boi} and  
by the following, which gives  the compatibility 
of the $p$-adic Chern classes of invertible sheaves 
with Hyodo-Kato's isomorphism:

\begin{lemm}\label{lemm:hkcc}
Let $s$ be a log point such that 
$\os{\circ}{s}$ is the spectrum 
of a perfect field $\kap$ of characteristic $p>0$.
Let $Y$ be a proper log smooth scheme of Cartier type 
over $s$. Let $\widehat{W\langle t\rangle}$ be 
the $p$-adic completion of the 
divided power polynomial $W\langle t\rangle$ over $W$.  
Endow $\widehat{W\langle t\rangle}$ 
with log structure given by a morphism 
${\mab N} \owns 1 \lom t \in \widehat{W\langle t\rangle}$.
Let ${\cal L}$ be an invertible sheaf on $\os{\circ}{Y}$.
Then, under the isomorphism {\rm \cite[(4.13)]{hk}} 
\begin{equation*}
\widehat{W\langle t\rangle}\otimes_W
(H^2_{{\rm crys}}(Y/W(s))\otimes_W K_0)
\os{\sim}{\lo}
H^2_{{\rm crys}}(Y/\widehat{W\langle t\rangle})
\otimes_W K_0, 
\tag{5.5.18.1}\label{eqn:delis}
\end{equation*}
$1\otimes c_{1,p}({\cal L})$ corresponds to 
the log cohomology class 
$c_{1,p}({\cal L})_{\widehat{W\langle t\rangle}}$ 
with respect to $Y/\widehat{W\langle t\rangle}$.
\end{lemm}

\begin{coro}[{\bf Log convergent log hard Lefschetz conjecture of 
log crystalline cohomologies in mixed characteristics}]\label{coro:c1o} 
Let the notations be as in {\rm (\ref{lemm:mlhs}) (3)}.  
Assume that $\os{\circ}{T}$ is connected.  
Assume also that there exists an exact closed point $t$ of 
$S_{\os{\circ}{T}_0}$ such that $X_{\os{\circ}{T}_0,t}$  
is the log special fiber of a projective strict semistable family of pure dimension 
over a complete discrete valuation ring of mixed characteristics. 
Then the morphism {\rm (\ref{eqn:fcicpl})} is an isomorphism.  
\end{coro} 
\begin{proof} 
(\ref{coro:lclt}) follows from (\ref{lemm:mlhs}) (3), (\ref{theo:llcmec}) and (\ref{prop:oxlp}). 
\end{proof}

\begin{rema} 
In a future paper we would like to discuss 
filtered versions of (\ref{theo:llccec}), (\ref{coro:lclt}), (\ref{theo:llcmec}) and (\ref{coro:c1o}). 
\end{rema} 

We conclude this section by giving the following, which 
has a different flavor from (\ref{coro:anmwc}) and  (\ref{coro:c1o}):

\begin{theo}\label{theo:ctsncl}
Let ${\cal V}$ be a complete discrete valuation ring 
of mixed characteristics $(0,p)$ with perfect residue field $\kap$. 
Let $S$ be a family of log points 
such that $\os{\circ}{S}$ is a $p$-adic formal ${\cal V}$-scheme. 
Let $T$ be an object of ${\rm Enl}^{\sq}_p(S/{\cal V})$. 
Let $X/S$ be a projective SNCL scheme. 
For each point $\os{\circ}{s}\in \os{\circ}{S}_{\os{\circ}{T}_1}$, 
let $s$ be the exact closed point of $S_{\os{\circ}{T}_1}$. 
Assume that there exists a point $\os{\circ}{s}\in \os{\circ}{S}_{\os{\circ}{T}_1}$ 
such that the fiber $X_s$ of $X_{\os{\circ}{T}_1}$ at $s$ 
has a projective SNCL lift over 
a complete discrete valuation ring of mixed characteristics 
with canonical log structure. 
Assume that $\os{\circ}{S}_{\os{\circ}{T}_1}$ is connected. 
Then the following hold$:$
\par 
{\rm (1)} The conjecture {\rm (\ref{conj:remc})} holds 
for $X_{\os{\circ}{T}_1}/S_{\os{\circ}{T}_1}\os{\subset}{\lo} S(T)^{\nat}$. 
\par 
{\rm (2)} 
The analogous statement to {\rm (1)} 
in the $l$-adic case holds.  
\end{theo}
\begin{proof} 
(1): By (\ref{prop:convmon}) and (\ref{theo:e2dgfam}),  
we may assume that $\os{\circ}{T}$ consists of a point and hence that 
$S_{\os{\circ}{T}_1}$ is a log point $s$. 
Let ${\cal W}$ be the Witt ring of $\os{\circ}{s}$. 
Let $K_0$ be the fraction field of ${\cal W}$. 
We consider the following morphism
\begin{align*} 
N \col & E_1^{-k, q+k}=\us{\us{j\geq -k}{j\geq 0}}{\oplus}
H_{\rm crys}^{q-2j-k}((X^{(2j+k)}/{\cal W})_{\rm crys},
{\cal O}_{X^{(2j+k)}/{\cal W}}
\otimes_{\mab Z}\vp^{(2j+k)}_{\rm crys}(X/{\cal W}))(-j-k)_{K}  & &
\tag{5.5.21.1}\label{ali:xwjk} \\
& \lo E_1^{-k+2, q+k-2}(-1)& &
\\
& =\us{\us{j\geq -(k-2)}{j\geq 0}}{\oplus}
H_{\rm crys}^{q-2j-k+2}((X^{(2j+k-2)}/{\cal W})_{\rm crys},\\
&{\cal O}_{X^{(2j+k-2)}/{\cal W}}
\otimes_{\mab Z}\vp^{(2j+k-2)}_{\rm crys}(X/{\cal W}))(-j-k+2)(-1)_{K}. &\\
&  &
\end{align*} 
The target above is equal to 
$$
\us{\us{j\geq -(k-1)}{j\geq 0}}{\oplus}
H_{\rm crys}^{q-2j-k}((X^{(2j+k)}/{\cal W})_{\rm crys},
{\cal O}_{X^{(2j+k)}/{\cal W}}
\otimes_{\mab Z}\vp^{(2j+k)}_{\rm crys}(X/{\cal W}))(-j-k)_{K}$$
$$\oplus \left\{
\begin{array}{ll}
H_{\rm crys}^{q-k+2}((X^{(k-2)}/{\cal W})_{\rm crys},
{\cal O}_{X^{(k-2)}/{\cal W}}
\otimes_{\mab Z}\vp^{(k-2)}_{\rm crys}(X/{\cal W}))(-k+1)_{K}
&  \;\mbox{$(\,0\geq 
-(k-2)\,),$} \\ 
0 &  \;\mbox{$(\, 0< 
-(k-2)\,).$} 
\end{array}\right.$$
\parno
By (\ref{prop:qauasimon}), 
the morphism
(\ref{ali:xwjk}) 
is the sum of $\pm$(identities) and zero-morphisms. 
Because (\ref{eqn:escssp}) degenerates at $E_2$ modulo torsion 
(\cite[(3.6)]{ndw}, see also (\ref{theo:edhw})), 
it suffices to show that 
$N^e \col E_2^{-e,q+e} \lo E_2^{e,q-e}(-e)$ 
is an isomorphism. 
By (\ref{cd:gsmsd}) 
the boundary morphisms $d_1^{\bul \bul}$ of 
the spectral sequences (\ref{eqn:escssp}) are 
expressed by Gysin morphisms 
and the morphisms induced by inclusion morphisms 
(see also \cite[4.13]{msemi} and \cite[(10.3.1;$\star$)]{ndw}). 
By \cite[(2.4)]{boi} there exists a canonical isomorphism
$H^q_{\rm crys}(X^{(j)}/{\cal W})(k)\us{\cal W}
{\otimes}K_0\os{\simeq}{\lo}
H^q_{\rm dR}({\cal X}_{K_0}^{(j)}/K_0)
\otimes H^2_{\rm dR}({\mab P}^1_{K_0}/K_0)^{\otimes (-k)}$. 
The analogous boundary morphisms $d_1^{\bul \bul}$ of 
the analogous spectral sequence over ${\mab C}$ 
are also induced by Gysin morphisms and 
inclusion morphisms (\cite{nlpi}). 
If we take a subfield $F$ of the complex number field ${\mab C}$ 
such that 
there exists an embedding $F\os{\subset}{\lo} K$ and that 
all ${\cal X}^{(j)}_{K_0}$ are defined over $F$ and 
we make a scalar extension to 
${\mab C}$, then
$N_K \col E_1^{-k, q+k} \lo  E_1^{-k+2, q+k-2}(-1)$ 
induces a morphism which 
preserves the Hodge structures.
Hence (1) follows from 
the argument of \cite[(4.2)]{sam} ([loc.~cit.,~(4.2.2)]).
\par 
(2): The analogous $l$-adic spectral sequence to 
(\ref{eqn:escssp}) for the trivial coefficient 
degenerates at $E_2$ by \cite[(2.1)]{nd}. 
In the $l$-adic case, 
the monodromy operator $N$ and 
the quasi-monodromy operator $\nu$ are different 
in general. However, in order to prove (2), 
it suffices to show that 
$\nu^e \col {\rm gr}^P_{i+e}
H^i_{{\rm log}{\textrm -}{\rm \acute{e}t}}(\ol{X},{\mab Q}_l)
\os{\simeq}{\lo} {\rm gr}^P_{i-e}
H^i_{{\rm log}{\textrm -}{\rm \acute{e}t}}(\ol{X},{\mab Q}_l)(-e)$ because 
$N^e=\nu^e$. (This is a strategy of Rapoport-Zink (\cite[(2.11)]{rz}): 
let $T$ be a generator of ${\mab Z}_l(1)$.
$(T-1) \otimes 
\check{T}$ is realized by $\nu$ and the term $(T-1)^s \otimes 
\check{T}$ $(s\geq 2)$ 
in $N={\rm log}\,T\otimes \check{T}$  induces a 
zero-morphism in the morphism 
$N^e \col 
{\rm gr}^P_{i+e}H^i_{{\rm log}{\textrm -}{\rm \acute{e}t}}(\ol{X},{\mab Q}_l)
\lo 
{\rm gr}^P_{i-e}H^i_{{\rm log}{\textrm -}{\rm \acute{e}t}}
(\ol{X},{\mab Q}_l)(-e)$.) 
Therefore the same proof as that of (1) works. 
(We have only to use [SGA4$\dfrac{1}{2}]$ 
Arcata V (1.7) b) instead of 
\cite[(2.4)]{boi} and the 
comparison of the \'{e}tale cohomologies of 
a smooth variety over ${\mab C}$ 
with the Betti ones [SGA 4-3] XI (4.4).)
\end{proof}

\begin{theo}\label{theo:ctslncl}
Let the notations be and the assumptions be as in 
${\rm (\ref{theo:ctsncl})}$.  
Then the following hold$:$
\par 
{\rm (1)} The conjecture {\rm (\ref{conj:lhilc}) (1)} 
holds 
for $X_{\os{\circ}{T}_1}/S_{\os{\circ}{T}_1}\os{\subset}{\lo} S(T)^{\nat}$. 
\par 
{\rm (2)} 
The analogous statement to {\rm (1)} 
in the $l$-adic case holds.  
\end{theo}
\begin{proof} 
The same proof as that of (\ref{theo:ctsncl}) works. 
\end{proof}

\chapter{Limits of weight filtrations and limits of slope filtrations 
on infinitesimal cohomologies in mixed characteristics}

\section{Good and split proper hypercoverings}\label{sec:ph} 
Let ${\cal V}$ and $K$ be as in \S\ref{sec:cfi}. 
Fix an algebraic closure $\ol{K}$ of $K$. 
Let ${\cal X}$ be a proper flat scheme over ${\cal V}$. 
Set ${\cal X}_K:={\cal X}\otimes_{\cal V}K$ 
and $X:={\cal X}\otimes_{\cal V}\kap$.    
For an algebraic extension $L$ of $K$ in $\ol{K}$, 
let ${\cal O}_L$ be the integer ring of $L$. 
Set ${\cal X}_L:={\cal X}_K\otimes_KL$ and 
${\cal X}_{{\cal O}_L}:={\cal X}\otimes_{\cal V}{\cal O}_L$. 
Let $N$ be a nonnegative integer. 
We use the same notation for 
an $N$-truncated simplicial fine log scheme 
over ${\rm Spec}({\cal V})=({\rm Spec}({\cal V}),{\cal V}^*)$. 
For a proper strict semistable family  $\os{\circ}{\cal Y}$ 
over a complete discrete valuation ring, 
endow it with the canonical log structure and 
denote by ${\cal Y}$ the resulting log scheme.  
Set ${\cal V}_{-1}:={\cal V}$ and 
$\os{\circ}{S}_{-1}:={\rm Spec}({\cal V}_{-1})$. 

\begin{defi} 
Let $N$ be a nonnegative integer.  
Let ${\cal X}_{\bul \leq N}$ 
be an $N$-truncated simplicial fine log scheme over ${\rm Spec}({\cal V})$. 
If the generic fiber ${\cal X}_{\bul \leq N,K}$ of ${\cal X}_{\bul \leq N}$ 
is an $N$-truncated proper hypercovering 
of ${\cal X}_K$, then we call 
${\cal X}_{\bul \leq N}$ 
an $N$-truncated simplicial generically proper hypercovering 
of ${\cal X}$ over ${\cal V}$. 
\end{defi}

\begin{prop}\label{prop:vv0} 
There exists a sequence  
\begin{equation*} 
\cdots  \supset {\cal V}_N  \supset {\cal V}_{N-1}
\supset \cdots  \supset {\cal V}_0   \supset {\cal V}_{-1}
\tag{6.1.2.1}\label{eqn:nsv}
\end{equation*} 
of finite extensions of 
complete discrete valuation rings of ${\cal V}$ in $\ol{K}$
and a proper strict semistable family 
$\os{\circ}{\cal N}_m$ over ${\cal V}_m$ 
$(m\in {\mab N})$ 
and a log smooth 
$m$-truncated simplicial log scheme  
${\cal X}(m)_{\bul \leq m}$  
such that 
${\cal X}(m)_{m'}
=\coprod_{0\leq l \leq m'}
\coprod_{[m'] \twoheadrightarrow [l]}({\cal N}_l
\times_{S_l}S_m)$ for each $0\leq m' \leq m$, 
where $S_m$ is the spectrum 
${\rm Spec}({\cal V}_m)$ with canonical log structure, 
and such that 
${\cal X}(m)_{\bul \leq m}$ is an $m$-truncated 
generically proper hypercovering of 
${\cal X}_{{\cal V}_m}$ over ${\cal V}_m$. 
\end{prop}
\begin{proof}
By de Jong's theorem (\cite[(6.5)]{dj}) 
about the semistable reduction theorem 
by using the base change by an alteration 
(if one makes a finite extension 
of a complete discrete valuation ring)
and by a standard argument in 
\cite[${\rm V}^{\rm bis}$ \S5]{sga4-2} 
and \cite[(6.2.1.1)]{dh3}, 
we obtain (\ref{prop:vv0}). 
\end{proof}  
\parno 
The sequence (\ref{eqn:nsv}) 
gives us the following sequence of log schemes: 
\begin{equation*} 
\cdots \lo S_N \lo S_{N-1} \lo \cdots  \lo S_0\lo S_{-1}, 
\tag{6.1.2.2}\label{eqn:wtsnna}
\end{equation*} 
where $S_m$ is the log scheme whose underlying scheme is 
${\rm Spec}({\cal V}_m)$ and 
whose log structure is the canonical log structure 
on ${\rm Spec}({\cal V}_m)$. 
Let $\kap_m$ be the residue field of ${\cal V}_m$. 
Set $s_m:=S_m\times_{{\rm Spec}({\cal V}_m)}
{\rm Spec}(\kap_m)$.  
Then we have the following sequence 
\begin{equation*} 
\cdots \lo s_N \lo  s_{N-1}\lo \cdots \lo 
s_0\lo s_{-1}. 
\tag{6.1.2.3}\label{eqn:wtslpnn}
\end{equation*}
We also have the following sequence 
\begin{equation*} 
\cdots \lo {\cal W}_{\star}(s_N) \lo 
{\cal W}_{\star}(s_{N-1}) 
\lo \cdots \lo {\cal W}_{\star}(s_0) 
\lo 
{\cal W}_{\star}(s_{-1}) \quad 
({\star}\in {\mab Z}_{\geq 1}~{\rm   or }~\star
={\rm nothing}).   
\tag{6.1.2.4}\label{eqn:wtswpsnna}
\end{equation*}
Let $N$ be a nonnegative integer. 
Set ${\cal X}_{\bul \leq N}:={\cal X}(N)_{\bul \leq N}$ and $s:=s_N$. 
Set also $N_m:={\cal N}_m\otimes_{{\cal V}_m}\kap_m$ $(0\leq m\leq N)$. 
Set $X(m)_{m'}
=\coprod_{0\leq l \leq m'}
\coprod_{[m'] \twoheadrightarrow [l]}
(N_l\times_{s_l}s_m)$ for each $0\leq m' \leq m$. 
Then 
$X(m)_{\bul \leq m}=
{\cal X}(m)_{\bul \leq m}\times_{S_m}s_m$.   
We have the family 
$\{X(m)_{\bul \leq m}\}_{m=0}^N$ of successive 
split truncated simplicial log schemes 
with respect to the sequence (\ref{eqn:wtslpnn}). 
Consequently, for each $N\in {\mab N}$, we obtain 
the associated split $N$-truncated successive simplicial SNCL scheme 
$X_{\bul \leq N}:=X(N)_{\bul \leq N}$ to  
$\{X(m)_{\bul \leq m}\}_{m=0}^N$ 
with respect to the sequence (\ref{eqn:wtslpnn}). 

\begin{defi}\label{defi:gdefa}
We say that the $N$-truncated simplicial generically proper 
hypercovering ${\cal X}_{\bul \leq N}$ of 
${\cal X}_{{\cal V}_N}$ is 
{\it gs$($=good and split$)$}.  
\end{defi}

\begin{lemm}\label{lemm:cmd}
Let $U$ be a separated scheme of finite type over $K$. 
Let ${\cal U}$ and ${\cal U}'$ be two flat models 
over ${\cal V}$ of $U$. 
Then there exists a flat model ${\cal U}''$ of $U$ over ${\cal V}$ 
with morphisms  
${\cal U}''\lo {\cal U}$ and ${\cal U}''\lo {\cal U}'$ 
which induce ${\rm id} \col U \os{=}{\lo} U$.  
If $U$, ${\cal U}$ and ${\cal U}'$ is proper over $K$ and over ${\cal V}$, respectively, 
then ${\cal U}''$ is proper over ${\cal V}$.  
\end{lemm}
\begin{proof}
Since ${\cal X}'_K={\cal X}_K$, we have the following composite morphism 
\begin{equation*} 
{\cal X}_K\os{\rm diag.}{\lo} {\cal X}_K\times_K{\cal X}'_K
\os{\subset}{\lo}  {\cal X}\times_{\cal V}{\cal X}'. 
\end{equation*}
Let ${\cal X}''$ be the scheme theoretic closure of ${\cal X}_K$ 
in ${\cal X}\times_{\cal V}{\cal X}'$ (\cite[(2.1)]{raynd}). 
(Recall that, for a flat scheme ${\cal Y}$ over ${\cal V}$ and 
a closed subscheme $Z$ of ${\cal Y}_K$, 
the scheme theoretic closure of $Z$ in ${\cal Y}$ is, by definition,
$\ul{\rm Spec}_{\cal Y}({\cal O}_{\cal Y}
/{\rm Ker}({\cal O}_{\cal Y}\os{\subset}{\lo} {\cal O}_{{\cal Y}_K}\lo {\cal O}_Z))$). 
We easily see that ${\cal X}''$ is flat over ${\cal V}$, 
${\cal X}''_K={\cal X}_K$ and 
the following composite morphisms    
\begin{equation*} 
{\cal X}'' \os{\subset}{\lo} {\cal X}\times_{\cal V}{\cal X}'
\os{p_1}{\lo}  {\cal X} \quad 
{\rm and}   \quad 
{\cal X}'' \os{\subset}{\lo} {\cal X}\times_{\cal V}{\cal X}'
\os{p_2}{\lo}  {\cal X}' 
\end{equation*}
do the job. 
\end{proof}

\begin{prop}\label{prop:coif} 
Let $N$ be a nonnegative integer. 
Then the following hold$:$
\par
$(1)$ Two gs $N$-truncated simplicial generically proper hypercoverings of 
the base change of ${\cal X}$ 
over an extension of ${\cal V}$ 
are covered by a gs $N$-truncated simplicial
generically proper hypercovering of the base change of 
${\cal X}$ over an extension of ${\cal V}$. 
\par
$(2)$   
For a morphism ${\cal X}'\lo {\cal X}$ of proper schemes over ${\cal V}$ 
and for a gs $N$-truncated simplicial generically proper hypercovering 
${\cal X}_{\bul \leq N}$ of the base change of ${\cal X}$ over 
an extension of ${\cal V}$, 
there exist a finite extension ${\cal V}'$ of ${\cal V}$ 
and a gs $N$-truncated simplicial generically proper hypercovering 
${\cal X}'_{\bul \leq N}$ of ${\cal X}'_{{\cal V}'}$ 
and a morphism 
${\cal X}'_{\bul \leq N} \lo 
{\cal X}_{\bul \leq N,{\cal V}'}$ 
fitting into the following commutative diagram$:$
\begin{equation*}
\begin{CD}
{\cal X}'_{\bul \leq N} @>>>  
{\cal X}_{\bul \leq N,{\cal V}'}\\
@VVV @VVV  \\
{\cal X}'_{{\cal V}'}@>>> {\cal X}_{{\cal V}'}.
\end{CD}
\tag{6.1.5.1}\label{cd:smpbca}
\end{equation*}
\end{prop}
\begin{proof}
As in the proof of \cite[(9,4)]{nh2},  
by using a general formalism 
in \cite[${\rm V}^{\rm bis}$ \S5]{sga4-2} and 
de Jong's semistable reduction theorem, 
we obtain (\ref{prop:coif}).   
Here we give only the proof of (1). 
\par 
Let ${\cal X}'$ be another proper flat model of ${\cal X}_K$ over ${\cal V}$. 
Then we have a sequence  
\begin{equation*} 
\cdots  \supset {\cal V}'_N  \supset {\cal V}'_{N-1}
\supset \cdots  \supset {\cal V}'_0 
 \supset {\cal V}'_{-1}={\cal V} 
\tag{6.1.5.2}\label{eqn:npv}
\end{equation*} 
of finite extensions of complete discrete valuation rings 
of ${\cal V}$ which is similar to (\ref{eqn:nsv}), 
$\os{\circ}{\cal N}{}'_m$ over ${\cal V}'_m$ and 
${\cal X}(m)'_{\bul \leq m}$ 
which are similar to 
$\os{\circ}{\cal N}_m$ over ${\cal V}_m$ and 
${\cal X}(m)_{\bul \leq m}$ $(0\leq m\leq N)$, respectively. 
Let ${\cal X}''$ be another proper flat model of ${\cal X}_K$ 
in (\ref{lemm:cmd}). 
By de Jong's semistable reduction theorem, 
there exists a finite extension ${\cal V}''_0$ of 
${\cal V}_0{\cal V}'_0:=
\{a\in \ol{K}~\vert~a=\sum_ia_ia'_i~(a_i\in {\cal V}_i, a'_i\in {\cal V}_i')\}$ 
over which there exists a proper semistable family ${\cal N}_0''$ 
over ${\cal V}''_0$ with an alteration 
${\cal N}_0''\lo (\os{\circ}{\cal N}_{0,{\cal V}''_0}
\times_{{\cal X}_{{\cal V}''_0}}{\cal X}''_{{\cal V}''_0})
\times_{{\cal X}''_{{\cal V}''_0}}
(\os{\circ}{\cal N}{}'_{0,{\cal V}''_0}
\times_{{\cal X}'_{{\cal V}''_0}}{\cal X}''_{{\cal V}''_0})$. 
As a result, we have the following commutative diagram: 
\begin{equation*} 
\begin{CD} 
\os{\circ}{\cal N}_0@<<< \os{\circ}{\cal N}_{0,{\cal V}''_0}
@<<<\os{\circ}{\cal N}{}''_0@>>> \os{\circ}{\cal N}{}'_{0,{\cal V}''_0} @>>>
\os{\circ}{\cal N}{}'_0,\\
@VVV @VVV @VVV @VVV @VVV \\
{\cal X}@<<< {\cal X}_{{\cal V}''_0}@<<<{\cal X}''_{{\cal V}''_0}@>>> {\cal X}'_{{\cal V}''_0}
@>>> {\cal X}'.
\end{CD} 
\end{equation*} 
The middle vertical morphism is proper and surjective. 
Assume that we are given a sequence 
\begin{equation*} 
{\cal V}''_{m-1}\supset \cdots  \supset {\cal V}''_0 
\supset {\cal V}''_{-1}={\cal V}
\tag{6.1.5.3}\label{eqn:nppv}
\end{equation*} 
such that  
\smallskip 
\parno
(1): ${\cal V}''_l \supset {\cal V}_l{\cal V}'_l$ in $\ol{K}$ 
$(0\leq l\leq m-1)$ 
\parno 
and 
\parno 
(2): a desired $(m-1)$-truncated simplicial proper scheme 
${\cal X}''(m-1)_{\bul \leq m-1}$ with morphisms 
${\cal X}''(m-1)_{\bul \leq m-1}\lo {\cal X}'(m-1)_{\bul \leq m-1}$
and ${\cal X}''(m-1)_{\bul \leq m-1}\lo {\cal X}(m-1)_{\bul \leq m-1}$. 
\smallskip 
\parno 
Then, by using de Jong's theorem about an alteration 
(\cite[Theorem 6.5]{dj}),  
there exist a complete discrete valuation ring ${\cal V}''_m$ 
and a strict semistable family 
$\os{\circ}{\cal N}{}''_m$ over ${\cal V}''_m$ 
with morphisms 
$\os{\circ}{\cal N}{}''_m
\lo {\rm cosk}_m^{{\cal X}''_{{\cal V}''_m}}
({\cal X}''_{\bul \leq m-1,{\cal V}''_m})_m$, 
$\os{\circ}{\cal N}{}''_m\lo \os{\circ}{\cal N}_{m,{\cal V}''_m}$ 
and 
$\os{\circ}{\cal N}{}''_m\lo \os{\circ}{\cal N}{}'_{m,{\cal V}''_m}$   
of schemes over ${\cal V}''_m$ fitting into the following commutative diagram 
\begin{equation*} 
\begin{CD} 
\os{\circ}{\cal N}_{m,{\cal V}''_m}
@<<< 
\os{\circ}{\cal N}{}''_m
@>>> \os{\circ}{\cal N}{}'_{m,{\cal V}''_m} \\ 
@VVV @VVV @VVV \\
 {\rm cosk}_m^{{\cal X}_{{\cal V}''_m}}
({\cal X}_{\bul \leq m-1,{\cal V}''_m})_m
@<<< 
 {\rm cosk}_m^{{\cal X}''_{{\cal V}''_m}}
({\cal X}''_{\bul \leq m-1,{\cal V}''_m})_m@>>> 
 {\rm cosk}_m^{{\cal X}'_{{\cal V}''_m}}
({\cal X}'_{\bul \leq m-1,{\cal V}''_m})_m. 
\end{CD}
\end{equation*}
Let $S''_l$ $(0\leq l\leq m)$ be 
the log scheme ${\rm Spec}({\cal V}''_l)$ with 
canonical log structure. 
Set 
${\cal X}(m)''_{\bul \leq m}
:=\coprod_{0\leq l \leq m}
\coprod_{[m] \twoheadrightarrow [l]}
({\cal N}''_l\times_{S''_l}S''_m)$, which is a desired 
$m$-truncated generically proper hypercovering of ${\cal X}_{{\cal V}''_N}$. 
\end{proof}

\section{Log infinitesimal cohomologies}\label{lic}
In this section we prove fundamental properties of 
cosimplicial (log) infinitesimal cohomologies 
which are the cosimplicial version of (log) infinitesimal cohomologies
defined in \cite{grcr} and \cite{cm}. 
In this section we work only in characteristic $0$. 
Though we need only the infinitesimal cohomology 
of the trivial coefficient for the trivial log case for 
the application in the mixed characteristics case in the next section,  
it seems standard to consider the nontrivial coefficient 
and the nontrivial log structure. 
\par 
Almost all results in this section seems known. 
Indeed, in almost all parts, we give only easier 
analogues of results in \cite{oc} and \cite{s2} except   
the log Poincar\'{e} lemma.  
\par 
Let $K$ be a field of characteristic $0$. 
Let $S$ be a fine log locally noetherian scheme over $K$. 
Let $X$ be a fine log locally noetherian scheme over $S$.  
Let $f\col X\lo S$ be the structural morphism.
Then we have the log infinitesimal site ${\rm Inf}(X/S)$ 
of $X/S$ as in \cite{grcr} and \cite[(0.1)]{cm}.  
Let us recall this. 
An object of ${\rm Inf}(X/S)$ is 
an exact closed nilpotent immersion $U \os{\sus}{\lo} T$ 
into a fine log locally noetherian scheme over $S$ 
such that $U$ is a log open subscheme of $X$. 
A morphism in ${\rm Inf}(X/S)$ is defined in a usual way. 
In this book we make the assumption ``locally noetherianness'' of $\os{\circ}{T}$
which has not been made in \cite[Definition 3.2.12]{s1} and 
\cite{cm} because we use this assumption. 
We denote simply by $T$ 
the object $U \os{\sus}{\lo} T$ of ${\rm Inf}(X/S)$ as usual. 
The covering family of $T$ is defined by 
a Zariski open covering $\{T_i\}_i$ of $T$. 
Let $(X/S)_{\rm inf}$ be the associated topos to the site 
${\rm Inf}(X/S)$.  
Let ${\cal O}_{X/S}$ be the structure sheaf of $(X/S)_{\rm inf}$ 
($\Gam(T,{\cal O}_{X/S}):=\Gam(T,{\cal O}_T)$ for $T\in {\rm Inf}(X/S)$).   
Let $U$ be a log open subscheme of $X$ and 
let $T$ be a fine log locally noetherian scheme over $S$ 
admitting a closed immersion $U\os{\sus}{\lo} T$ over $S$.    
Following \cite{ollc} and \cite{s2}, 
we call $T=(U,T)$ a {\it prewidening} of $X/S$. 
If the closed immersion $U\os{\sus}{\lo} T$ is nilpotent, 
we call $T=(U,T)$ a {\it widening} of $X/S$. 
If the immersion $U\os{\sus}{\lo} T$ for a (pre)widening 
$(U,T)$ of $X/S$ is exact, 
we say that $(U,T)$ is an {\it exact {\rm (}pre{\rm )}widening} of $X/S$. 
More generally, let $T$ be a fine log locally noetherian scheme over $S$ 
admitting an immersion $U\os{\sus}{\lo} T$ over $S$.    
We call $T=(U,T)$ a {\it quasi-prewidening} of $X/S$. 
If the immersion $U\os{\sus}{\lo} T$ is nilpotent, 
we call $T=(U,T)$ a {\it quasi-widening} of $X/S$. 
\par 

\par 
Let $(U,T)$ be a (pre)widening of $X/S$. 
We consider $(U,T)$ as the following functor 
\begin{equation*} 
h_T:=h_{(U,T)} \col {\rm Inf}(X/S)\owns 
(U',T') \lom {\rm Hom}_S((U',T'),(U,T))\in 
{\rm (Sets)},  
\end{equation*} 
where an element of ${\rm Hom}_S((U',T'),(U,T))$ 
is the following commutative diagram 
\begin{equation*} 
\begin{CD}
U' @>>> U\\
@V{\bigcap}VV @VV{\bigcap}V \\
T'@>>> T 
\end{CD} 
\end{equation*}
over $S$.  
As in \cite[(5.8)]{klog1}, \cite[(3.2.2)]{s1} 
and \cite[(0.9)]{cm}, 
we have the inductive system of universal infinitesimal neighborhoods 
$\{{\mathfrak T}_{U,n}(T)\}_{n=1}^{\infty}$ 
with natural morphisms 
$\bet_n \col {\mathfrak T}_{U,n}(T) \lo T$ which is an analogue of 
the inductive system of universal log enlargements defined 
in \cite{of}, \cite{oc} and \cite{ollc}. 
The definition of  
${\mathfrak T}_{U,n}(T)$ is as follows.  
Take the exactification 
$U\os{\sus}{\lo} T^{\rm ex}$ 
of the closed immersion $U\os{\sus}{\lo} T$ defined before (\ref{rema:pds}).   
Note that $\os{\circ}{T}{}^{\rm ex}$ is a locally noetherian scheme.
Let ${\cal J}$ be the defining ideal sheaf of the closed immersion 
$U\os{\sus}{\lo} T^{\rm ex}$. 
Then ${\mathfrak T}_{U,n}(T)=
\ul{\rm Spec}^{\log}_{T^{\rm ex}}
({\cal O}_{T^{\rm ex}}/{\cal J}^n)$.  
For simplicity of notation,  
set $T_n:={\mathfrak T}_{U,n}(T)$.  
Then 
\begin{equation*}
h_{(U,T)}=\vil_{n}h_{T_n}.  
\tag{6.2.0.1}\label{eqn:hynhn}
\end{equation*} 
Since the morphisms 
${\cal O}_{T^{\rm ex}}/{\cal J}^{n+1}\lo {\cal O}_{T^{\rm ex}}/{\cal J}^n$ 
and $M_{T^{\rm ex}_{n+1}}\lo M_{T^{\rm ex}_n}$ 
are surjective, the morphism $h_{T_n}\lo h_{T_{n+1}}$ is injective. 
Since the morphisms 
${\cal O}_{T^{\rm ex}}\lo {\cal O}_{T^{\rm ex}}/{\cal J}^n$ 
and 
$M_{T^{\rm ex}}\lo M_{T^{\rm ex}_n}$ 
are surjective, we see that the morphism 
$h_{T_n}\lo h_T$ is also injective. 
More generally, let $T=(U,T)$ be a quasi-prewidening of $X/S$. 
Take an open log subscheme $T'$ of $T$ 
such that the immersion $U\os{\sus}{\lo} T$ factors through 
a closed immersion $U \os{\sus}{\lo} T'$ over $S$. 
Then, set 
$T_n:={\mathfrak T}_{U,n}(T):={\mathfrak T}_{U,n}(T')$.  
As in the case of the PD-envelope of an immersion in \cite{bb}, 
we can easily show that 
${\mathfrak T}_{U,n}(T)$ is independent of the choice of $T'$.    
Then we can define the inductive system of 
universal infinitesimal neighborhoods 
$\{{\mathfrak T}_{U,n}(T)\}_{n=1}^{\inf}$ of 
$U$ in $T$ as above.  
Set $T_n:={\mathfrak T}_{U,n}(T)$ as above. 
\par 
Let $T$ be a quasi-prewidening of $X/S$. 
Let $\os{\to}{T}$ be the site 
which is the obvious analogue of  
\cite[\S3]{oc} and \cite[Definition 2.1.28]{s2}: 
an object of $\os{\to}{T}$ is 
a log open subscheme $V_n$ of $T_n$ 
for a positive integer  $n$; let $V_m$ be 
a log open subscheme of $T_m$;  
for $n>m$, ${\rm Hom}_{\os{\to}{T}}(V_n,V_m):= \emptyset$, 
and, for $n\leq m$, 
${\rm Hom}_{\os{\to}{T}}(V_n,V_m)$ 
is the set of morphisms $V_n \lo V_m$'s of 
log open subschemes over the natural morphism $T_n \lo T_m$; 
a covering of $V_n$ is a Zariski open covering of $V_n$. 
Let ${\rm Top}(\os{\to}{T})$ be the topos associated to $\os{\to}{T}$
(${\rm Top}(\os{\to}{T})$ has been denoted by $\os{\to}{T}$ in \cite{oc}.). 
Let $\phi_n \col T_n \lo T_{n+1}$ be the natural morphism 
of infinitesimal neighborhoods of $X/S$. 
Then an object ${\cal F}_{\os{\to}{T}}$ in ${\rm Top}(\os{\to}{T})$
is a family $\{({\cal F}_n,\psi_n)\}_{n=1}^{\infty}$ of pairs, 
where ${\cal F}_n$ is a Zariski sheaf in $(T_n)_{\rm zar}$ 
and $\psi_n \col \phi^{-1}_n({\cal F}_{n+1}) \lo {\cal F}_n$ 
is a morphism of Zariski sheaves. 
Let ${\cal O}_{\os{\to}{T}}
=\{({\cal O}_{T_n},\phi_n^{-1})\}_{n=1}^{\infty}$ 
be the structure sheaf of $\os{\to}{T}$, where 
$\phi_n^{-1}$ is the natural pull-back morphism
$\phi_n^{-1}({\cal O}_{T_{n+1}}) \lo {\cal O}_{T_n}$.  
Imitating \cite[p.~141]{oc}, we say that 
an ${\cal O}_{\os{\to}{T}}$-module ${\cal F}_{\os{\to}{T}}$ is {\it crystalline} 
if the natural morphism 
${\cal O}_{T_n}\otimes_{\phi^{-1}_{n+1}({\cal O}_{T_{n+1}})}
\phi^{-1}_{n+1}({\cal F}_{n+1}) \lo {\cal F}_n$ is an isomorphism, 
and we say that 
${\cal F}_{\os{\to}{T}}$ is {\it coherent} 
if ${\cal F}_n \in {\rm Coh}({\cal O}_{T_n})$ 
for all $n\in {\mab Z}_{>0}$. 
Here ${\rm Coh}({\cal O}_{T_n})$ 
is the category of coherent ${\cal O}_{T_n}$-modules. 
For an object $E$ of 
$(X/S)_{\rm inf}\vert_{T}$ and 
for a morphism $T' \lo T$ of quasi-prewidenings, 
we have an associated object 
${\cal E}'=\{{\cal E}'_n\}_{n=1}^{\infty}$ 
$({\cal E}'_n=E_{T'_n})$ 
in ${\rm Top}(\os{\to}{T'})$. 
Imitaing \cite[p.~147]{oc}, we say that an 
${\cal O}_{X/S}\vert_T$-module $E$ is 
{\it coherent} and {\it crystalline} if ${\cal E}'$ is  coherent and crystalline  
for any morphism $T' \lo T$ of quasi-prewidenings.  
Let $\star$ be ex or nothing. 
Let 
${\rm Top}(T^{\star})$ be an analogue of ${\rm Top}(\os{\to}{T})$ 
defined by replacing the transition morphism $T_n\lo T_{n+1}$ 
$(n\in {\mab Z}_{\geq 1})$ 
by ${\rm id}_{T^{\star}}\col T^{\star}\lo T^{\star}$. 
Let 
$\bet_{T^{\star}}\col {\rm Top}(\os{\to}{T})\lo {\rm Top}(T^{\star}) 
\lo T^{\star}_{\rm zar}$ 
be the natural composite morphism of topoi, 
which is the log infinitesimal version of the morphism 
$\gam$ in \cite[p.~91]{s2} 
(=the log version of \cite[p.~141]{oc}):  
$\bet_{T^{\star}}$ is, by definition, the inverse limit of 
the system of direct images $\{T_n\lo T^{\star}\}_{n=1}^{\infty}$. 
Then, by the definition of $T_n$, we see that 
the functor $\bet_{T^{\rm ex}}$ induces 
the following equivalence of categories by \cite[II (9.6)]{hartsag}$:$
\begin{align*} 
\bet_{T^{\rm ex}*} \col & 
\{{\rm coherent \;crystalline\;}{\cal O}_{\os{\to}{T}}{\rm -modules}\}
\os{\sim}{\lo}  \{{\rm coherent\;}{\cal O}_{T^{\rm ex}}{\rm -modules}\}. 
\end{align*} 
Here we have used the assumption of the locally noetherianness of $\os{\circ}{T}$. 

\begin{prop}[{\rm {\bf cf.~\cite[(3.7)]{oc}}}]\label{prop:cv}
For a coherent crystalline ${\cal O}_{\os{\to}{T}}$-module $E_{\os{\to}{T}}$ and 
a positive integer $q$, $R^q\bet_{T^{\star}*}(E_{\os{\to}{T}})=0$. 
If $\os{\circ}{T}$ is an affine scheme, 
then $H^q(\os{\to}{T},E_{\os{\to}{T}})=0$. 
\end{prop} 
\begin{proof} 
(Though the proof is the same as that of \cite[(3.7)]{oc}, 
we give the proof for the completeness of this book.) 
We have only to prove the second statement. 
Consider the topos ${\mab N}$ as in [loc.~cit.] 
(${\mab N}$ is the topos of consisting of inverse systems of sets indexed by 
the set ${\mab N}$ of nonnegative integers)
and the morphism 
$\os{\to}{\bet}\col \os{\to}{T}\lo {\mab N}$ of topoi: 
$\os{\to}{\bet}_*({\cal F}_{\os{\to}{T}})_n:=\Gam(T_n,{\cal F}_n)$ 
(${\cal F}_{\os{\to}{T}}\in {\rm Top}(\os{\to}{T})$). 
Then we have the following spectral sequence: 
\begin{align*}  
E_2^{pq}=H^p({\mab N},R^q\os{\to}{\bet}_*(E_{\os{\to}{T}}))\Lo 
H^{p+q}(\os{\to}{T},E_{\os{\to}{T}}). 
\end{align*}
Since $R^q\os{\to}{\bet}_*(E_{\os{\to}{T}})_n=H^q(T_n,E_n)=0$ $(q>0)$, 
\begin{align*} 
H^p(\os{\to}{T},E_{\os{\to}{T}})=H^p({\mab N},\os{\to}{\bet}_*(E_{\os{\to}{T}}))
=R^p\vpl_n \Gam(T_n,E_n).
\end{align*}  
Since the morphism $\Gam(T_{n+1},E_{n+1})\lo \Gam(T_n,E_n)$ is surjective 
(the analogue of \cite[(3.8)]{oc} is obvious), 
$H^p(\os{\to}{T},E_{\os{\to}{T}})=R^p\vpl_n \Gam(T_n,E_n)=0$ for $p>0$. 
\end{proof}

\par 
Let $T=(U,T)$ be a widening of $X/S$. Then, as usual, 
we obtain the localized topos $(X/S)_{\rm inf}\vert_T$. 
Let $j_T\col (X/S)_{\rm inf}\vert_T\lo (X/S)_{\rm inf}$ 
be the localization of topoi. 
Let us define a functor $
\varphi_{\os{\to}{T}*}\col (X/S)_{\rm inf}\vert_T \lo \os{\to}{T}_{\rm zar}$
by the following 
\begin{align*} 
\varphi_{\os{\to}{T}*}(E)(V_n):= E((U\vert_{V_n},V_n)\lo (U,T)). 
\end{align*} 
Let 
\begin{align*} 
u_{X/S}^{\rm inf} \col (X/S)_{\rm inf}\lo X_{\rm zar} 
\end{align*} 
be the canonical projection of topoi: 
$u_{X/S}^{{\rm inf}*}(E)(U,T)=E(U)$ and 
$u_{X/S*}^{\rm inf}(E)(U)=\Gam((U/S)_{\rm inf},j^*_{\rm inf}(E))$, 
where $j_{\rm inf}\col (U/S)_{\rm inf}\lo (X/S)_{\rm inf}$ is the induced morphism of 
topoi by the open immersion $j\col U\os{\sus}{\lo} X$. 
It is straightforward to check that the following diagram is commutative: 
\begin{equation*} 
\begin{CD} 
(X/S)_{\rm inf}\vert_T@>{\varphi_{\os{\to}{T}*}}>> \os{\to}{T}_{\rm zar}
@>{\bet_{T*}}>> T_{\rm zar} \\
@V{j_{T*}}VV @.  @| \\
(X/S)_{\rm inf}@>{u_{X/S*}^{\rm inf}}>> X_{\rm zar}
@<<{j_*}< U_{\rm zar}. 
\end{CD}
\tag{6.2.1.1}\label{cd:xuz}
\end{equation*}   

\begin{prop}\label{prop:ij}
The functor 
$\varphi_{\os{\to}{T}*} \col (X/S)_{\rm inf}\vert_T\lo \os{\to}{T}_{\rm zar}$ 
transforms injective sheaves to flasque sheaves. 
\end{prop} 
\begin{proof} 
The proof is the same as that of \cite[(4.1)]{oc}. 
We give the proof of this proposition for the completeness of this book. 
\par 
Let $j_n\col T_n\lo T_{n+1}$ be the natural morphism in $\os{\to}{T}$. 
Then $j_n$ is a monomorphism. 
More generally, let $k\col V\lo V'$ be a morphism in $\os{\to}{T}$. 
Then $k$ is a monomorphism since $k$ is a composite morphism of 
the restrictions of $j_n$'s. 
Let $k_{{\rm inf}}\col (X/S)_{\rm inf}\vert_{V}\lo (X/S)_{\rm inf}\vert_{V'}$ 
be the induced morphism. 
Let $E$ be an abelian sheaf in $(X/S)_{\rm inf}\vert_{V}$. 
Let $F$ be a presheaf in $(X/S)_{\rm inf}\vert_{V'}$
defined by the following: 
for a morphism $h\col W\lo V'$ in $(X/S)_{\rm inf}$, set 
$F(h):=\us{g \col W'\lo V}{\oplus}\{E(g)~\vert~k\circ g=h\}$. 
Since $k$ is a monomorphism, $F(h)=0$ if $h$ does not factor through $k$ 
and $F(h)=E(g)$ if $h$ factors through $k$, 
where $g$ is a unique morphism such that $k\circ g=h$. 
Let $k_{{\rm inf}!}(E)$ be the sheafification of $F$. 
Then $k_{{\rm inf}!}$ is the left adjoint functor of $k_{{\rm inf}}^*$ and 
it is clear that the natural morphism $k_{{\rm inf}!}k_{{\rm inf}}^*(E)\lo E$ is injective 
as in the proof of \cite[(4.1)]{oc}. 
\par  
Let $I$ be an injective sheaf in $(X/S)_{\rm inf}\vert_T$. It suffices to prove that 
the morphism $I(V'\lo T)\lo I(V\lo T)$ is surjective. 
Consider the constant sheaf ${\mab Z}_{V^{\star}}$ in 
$(X/S)_{\rm inf}\vert_{V^{\star}}$ $(\star=\prime$ or nothing). 
Then $I(V^{\star}\lo T)=
{\rm Hom}_{(X/S)_{\rm inf}\vert_{V^{\star}}}
({\mab Z}_{V^{\star}},I\vert_{(X/S)_{\rm inf}\vert_{V^{\star}}})$. 
Since ${\mab Z}_{V}=j_{{\rm inf},!}j_{\rm inf}^*({\mab Z}_{V'})$, 
the morphism ${\mab Z}_{V}\lo {\mab Z}_{V'}$ is injective. 
Hence the morphism $I(V'\lo T)\lo I(V\lo T)$ is surjective. 
\end{proof} 

\begin{prop}\label{prop:ex}
The following hold$:$ 
\par 
$(1)$ The functor 
$\varphi_{\os{\to}{T}*} \col (X/S)_{\rm inf}\vert_T\lo \os{\to}{T}_{\rm zar}$ 
is exact. 
\par 
$(2)$ Let $E$ be an abelian sheaf in $(X/S)_{\rm inf}\vert_T$. 
Set $E_{\os{\to}{T}}:=\{E_{T_n\lo T}\}_{n=1}^{\infty}$.  
Then there exists a canonical isomorphism 
\begin{align*} 
H^q(\os{\to}{T}_{\rm zar},E_{\os{\to}{T}})
\os{\sim}{\lo} H^q((X/S)_{\rm inf}\vert_T,E). 
\end{align*} 
\par 
$(3)$ If $\os{\circ}{T}$ is affine and if 
$E$ is a crystal of ${\cal O}_{X/S}$-modules 
in $(X/S)_{\rm inf}\vert_T$, 
then $H^q((X/S)_{\rm inf}\vert_T,E)=0$ for $q>0$. 
\end{prop} 
\begin{proof} 
By using (\ref{prop:ij}), 
the proof is the same as that of \cite[(4.2)]{oc}. 
\end{proof}

\begin{prop}\label{prop:av} 
Set $u_T^{\rm inf}:=u_{X/S}^{\rm inf}\circ j_T$. 
If $\os{\circ}{T}$ is affine and if 
$E$ is a crystal of ${\cal O}_{X/S}$-modules in $(X/S)_{\rm inf}\vert_T$, 
then $R^qu_{T*}^{\rm inf}(E)$ and $R^qj_{T*}(E)$ vanish for $q>0$. 
\end{prop} 
\begin{proof} 
By using (\ref{prop:ex}) (3), the proof is the same as that of \cite[(4.3)]{oc}. 
\end{proof} 

\begin{coro}\label{coro:je} 
Let the notations be as in {\rm (\ref{prop:av})}.  
Then $j_{T*}(E)$ is acyclic for $u_{X/S*}^{\rm inf}$. 
\end{coro} 
\begin{proof}  
Indeed, as in the proof of \cite[Corollary 2.3.4]{s2}, 
we have the following by (\ref{prop:cv}), (\ref{prop:ex}) (1) and 
(\ref{prop:av}): 
\begin{align*} 
Ru_{X/S*}^{\rm inf}j_{T*}(E)& =Ru_{X/S*}^{\rm inf}Rj_{T*}(E)
=Rj_*R\bet_{T*}R\varphi_{\os{\to}{T}*}(E)
=Rj_*\bet_{T*}\varphi_{\os{\to}{T}*}(E) \\
&=j_*\bet_{T*}\varphi_{\os{\to}{T}*}(E)
=u_{X/S*}^{\rm inf}j_{T*}(E).
\end{align*} 
\end{proof}

\par 
Assume, for the moment, that there exists 
an immersion $\iota \col X \os{\sus}{\lo} P$ 
into a log smooth scheme over $S$.  
We can consider $\iota$ as a quasi-prewidening $P:=(X,P)$ of $X/S$. 

\par 
Let $j_{P} \col (X/S)_{\rm inf}\vert_{P} \lo (X/S)_{\rm inf}$ 
be the localization of topoi.

Let 
\begin{align*}
\varphi^* 
\col & \{\text{the category of coherent crystals of }
{\cal O}_{\os{\to}{P}}\text{-modules}\} 
\tag{6.2.5.1}\label{eqn:phyt} \\
{} & \lo \{\text{the category of coherent crystals of }
{\cal O}_{X/S}\vert_{P}\text{-modules}\}
\end{align*}
be a natural functor which is defined as follows: 
for a coherent ${\cal O}_{P_n}$-module ${\cal F}$ and 
for an object $T$ of $(X/S)_{\rm inf}\vert_P$ 
with a morphism $T\lo P_n$ for some $n\in {\mab N}$,  
$\varphi^*({\cal F}):=
{\cal O}_T\otimes_{{\cal O}_{P_n}}{\cal F}$ 
(cf.~\cite[p.~147]{oc}); 
$\varphi^*({\cal F})$ is independent of the choice of 
the morphism $T\lo P_n$.

\par 
For a coherent crystal of 
${\cal O}_{\os{\to}{P}}$-module 
${\cal E}_{\bul}=\{{\cal E}_n\}_{n=1}^{\infty}$, set 
\begin{equation*} 
L^{\rm inf}_{X/S}({\cal E})
:=j_{P*}\varphi^*({\cal E}_{\bul})
\tag{6.2.5.2}\label{eqn:luys}
\end{equation*} 
by abuse of notation (cf.~\cite[Remark 6.10.1]{bob}, 
\cite[p.~152]{oc}, \cite[p.~95]{s2}). 
Let $T:=(U,T)$ be a quasi-prewidening of $X/S$. 
The following diagram 
\begin{equation*} 
\begin{CD} 
U @>{\subset}>> X @>{\subset}>> P\\ 
@V{\iota}VV  @VVV @VVV\\\ 
T  @>>>  S@=S
\end{CD} 
\tag{6.2.5.3}\label{cd:uxp}
\end{equation*} 
gives a natural immersion 
$U \os{\subset}{\lo} T\times_SP$.  
Let 
$\{{\mathfrak T}_{U,n}(T\times_SP)\}_{n=1}^{\infty}$ 
be the system of the universal infinitesimal neighborhoods 
of this immersion. 
Then we have the following diagram 
of the inductive systems of 
infinitesimal neighborhoods of $X/S$:  
\begin{equation*}
\begin{CD}
\{{\mathfrak T}_{U,n}(T\times_SP)\}_{n=1}^{\infty} 
@>{\{p_n\}_{n=1}^{\infty}}>> 
\{P_n\}_{n=1}^{\infty}\\ 
@V{\{p'_n\}_{n=1}^{\infty}}VV \\
T.@. @. 
\end{CD}
\tag{6.2.5.4}\label{cd:ppp}
\end{equation*} 
Let ${\cal E}_{\bul}=\{{\cal E}_n\}_{n=1}^{\infty}$ 
be a coherent crystal of 
${\cal O}_{\os{\to}{P}}$-modules.  
Let $T:=(U,T)$ be an object of $(X/S)_{\rm inf}$. 
Then, as in the proof of \cite[(5.4)]{oc} and 
\cite[Theorem 2.3.5]{s2},   
we have the following equality   
\begin{equation*}
L^{\rm inf}_{X/S}({\cal E}_{\bul})_{T}
=\vpl_np'_{n*}p_n^*({\cal E}_n) 
\tag{6.2.5.5}\label{eqn:lueyset}
\end{equation*} 
by the definition of 
$L^{\rm inf}_{X/S}$ ((\ref{eqn:luys})). 
Set $P(1)=P\times_SP$ and let $P(1)^{\rm ex}$ be the exactification of  
the diagonal immersion $P\os{\sus}{\lo} P(1)$.

\par 
We would like to give the log infinitesimal version of 
the linearization functor in \cite[\S2]{bob} 
(see also \cite[(0.3)]{cm}). 
\par 
Let ${\cal E}$ be a coherent ${\cal O}_{P^{\rm ex}}$-module.  
Let $\{P(1)^{\rm ex}_n\}_{n=1}^{\infty}$ 
be the system of the universal log enlargements 
for the diagonal immersion $P\os{\subset}{\lo}P(1)$. 
Set 
\begin{align*} 
L_n({\cal E}):={\cal O}_{P(1)^{\rm ex}_n}
\otimes_{{\cal O}_{P^{\rm ex}}}{\cal E} 
\end{align*}  
and 
\begin{align*} 
L({\cal E}):=\vpl_n
({\cal O}_{P(1)^{\rm ex}_n})
\otimes_{{\cal O}_{P^{\rm ex}}}{\cal E}
=\vpl_n
({\cal O}_{P(1)^{\rm ex}_n}
\otimes_{{\cal O}_{P^{\rm ex}}}{\cal E}). 
\end{align*}  
Here we have used a fact that 
${\cal O}_{P(1)^{\rm ex}_n}$ is a locally free 
${\cal O}_{P^{\rm ex}}$-module 
((\cite[Lemma 3.27]{s1}); see also the proof of (\ref{prop:pl}) below)
to see that the natural morphism 
$\vpl_n
({\cal O}_{P(1)^{\rm ex}_n})
\otimes_{{\cal O}_{P^{\rm ex}}}{\cal E}
\lo \vpl_n
({\cal O}_{P(1)^{\rm ex}_n}
\otimes_{{\cal O}_{P^{\rm ex}}}{\cal E})$ is an isomorphism.

Let 
$\tau_n \col 
{\cal O}_{P(1)^{\rm ex}_n}\lo 
{\cal O}_{P(1)^{\rm ex}_n}$ 
be the induced morphism by the morphism 
$P(1)\owns (x,y)\lom (y,x)\in P(1)$. 
Then we have the following morphism 
\begin{align*} 
L_n({\cal E}) &=
{\cal O}_{P(1)^{\rm ex}_n}
\otimes_{{\cal O}_{P^{\rm ex}}}{\cal E} 
\owns 
a\otimes x\lom x\otimes \tau_n(a) \in {\cal E}\otimes_{{\cal O}_{P^{\rm ex}}}
{\cal O}_{P(1)^{\rm ex}_n}
\quad (a\in {\cal O}_{P(1)^{\rm ex}_n}, \; 
x\in {\cal E}), 
\end{align*} 
which we denote by $\tau_n$ again. 
Then we have the following morphism 
\begin{align*} 
L_n({\cal E})\os{\tau_n}{\lo} 
{\cal E}\otimes_{{\cal O}_{P^{\rm ex}}}{\cal O}_{P(1)^{\rm ex}_n} &
\owns x\lom 1\otimes x\in 
{\cal O}_{P(1)^{\rm ex}_n}
\otimes_{{\cal O}_{P^{\rm ex}}}
{\cal E}\otimes_{{\cal O}_{P^{\rm ex}}}{\cal O}_{P(1)^{\rm ex}_n},  
\tag{6.2.5.6}\label{eqn:strn}\\
\end{align*} 
which induces the following isomorphism:
\begin{align*} 
\eps_{n,n}\col 
{\cal O}_{P(1)^{\rm ex}_n}\otimes_{{\cal O}_{P^{\rm ex}}}L_{n}({\cal E}) 
\os{\sim}{\lo} 
L_n({\cal E})\otimes_{{\cal O}_{P^{\rm ex}}} {\cal O}_{P(1)^{\rm ex}_n}.   
\tag{6.2.5.7}\label{eqn:stirn}\\ 
\end{align*}  
Fix a positive integer $m$. Then, for $n\geq m$, 
we have the following isomorphism 
\begin{align*} 
\eps_{m,n}\col 
{\cal O}_{P(1)^{\rm ex}_m}\otimes_{{\cal O}_{P^{\rm ex}}}L_{n}({\cal E}) 
\os{\sim}{\lo} 
L_n({\cal E})\otimes_{{\cal O}_{P^{\rm ex}}} {\cal O}_{P(1)^{\rm ex}_m}.   
\tag{6.2.5.8}\label{eqn:stirmn}\\ 
\end{align*}  
Taking the projective limit $\vpl_n$ of the system above 
and noting that 
${\cal O}_{P(1)^{\rm ex}_n}$ is a locally free module, 
we have the following isomorphism 
\begin{align*} 
\eps_m\col 
{\cal O}_{P(1)^{\rm ex}_m}\otimes_{{\cal O}_{P^{\rm ex}}}L({\cal E}) 
\os{\sim}{\lo} 
L({\cal E})\otimes_{{\cal O}_{P^{\rm ex}}} {\cal O}_{P(1)^{\rm ex}_m}.   
\tag{6.2.5.9}\label{eqn:stirmin}\\ 
\end{align*}  
Hence the family $\{\eps_n\}_{n=1}^{\infty}$ defines a stratification 
on $L({\cal E})$. 

\begin{rema}\label{rema:fds}
(1) For each $m\in {\mab Z}_{\geq 1}$, consider 
the composite immersion 
$P\os{\sus}{\lo} P^{\rm ex}(1)_m\os{\sus}{\lo} P^{\rm ex}(1)$. 
We have two projections 
$p_1, p_2\col P(1)^{\rm ex}_m\lo P^{\rm ex}$ and 
$p_1, p_2\col P(1)^{\rm ex}\lo P^{\rm ex}$. 
It seems to me that we cannot have a morphism
$p_i\col {\mathfrak T}_{P(1)^{\rm ex}_m,n}(P(1)^{\rm ex})\lo 
{\mathfrak T}_{P^{\rm ex},n+m}(P(1)^{\rm ex})$ 
$(i=1,2)$ such 
that the composite morphism 
${\mathfrak T}_{P^{\rm ex},n}(P(1)^{\rm ex})
\lo {\mathfrak T}_{P(1)^{\rm ex},n}(P(1)^{\rm ex})\lo 
{\mathfrak T}_{P^{\rm ex},n+m}(P(1)^{\rm ex})$
is the transition morphism. 
It seems to me that we cannot imitate the analogous parts of the proofs of 
the (log) Poincar\'{e} lemma in \cite[(5.4)]{oc} and \cite[Theorem 2.3.5]{s2} 
for the proof of (\ref{prop:pl}) below. 
\par 
(2) We would not like to use the \v{C}ech-Alexander complexes used in 
\cite{grcr} and \cite{cm} because we would like to follow the linearization in  
\cite{bob} as possible. 
\end{rema}

\begin{defi} 
Let ${\cal E}$ and ${\cal F}$ be coherent ${\cal O}_{P^{\rm ex}}$-modules. 
A differential operator $u\col {\cal E}\lo {\cal F}$ of order $\leq k$ 
is an ${\cal O}_{P^{\rm ex}}$-linear morphism 
\begin{equation*} 
{\cal O}_{P(1)^{\rm ex}_k}
\otimes_{{\cal O}_{P^{\rm ex}}}{\cal E}
\lo {\cal F} 
\tag{6.2.7.1}\label{eqn:pk1} 
\end{equation*} 
for $k\in {\mab Z}_{\geq 1}$. 
\end{defi}

\parno 
Let 
\begin{equation*} 
\del_{m,n} \col {\cal O}_{P(1)^{\rm ex}_{m+n}}\lo 
{\cal O}_{P(1)^{\rm ex}_m}\otimes_{{\cal O}_{P^{\rm ex}}}{\cal O}_{P(1)^{\rm ex}_n}
\end{equation*} 
be the morphism obtained by the natural morphism 
$P(1)^{\rm ex}_m\times_{P^{\rm ex}}P(1)^{\rm ex}_n\lo P(1)^{\rm ex}_{m+n}$ 
(\cite[Lemma 3.2.3]{s1}). 
Let $h\col {\cal O}_{P(1)^{\rm ex}_k}\otimes_{{\cal O}_{P^{\rm ex}}}{\cal E}\lo {\cal F}$ 
be the differential operator of order $\leq k$.  
Then we have the following composite morphism 
\begin{align*}  
L_n(h)\col {\cal O}_{P(1)^{\rm ex}_n}\otimes_{{\cal O}_{P^{\rm ex}}}{\cal E}
\os{\del_{n-k,k}}{\lo}  
{\cal O}_{P(1)^{\rm ex}_{n-k}}\otimes_{{\cal O}_{P^{\rm ex}}}{\cal O}_{P(1)^{\rm ex}_k} 
\otimes_{{\cal O}_{P^{\rm ex}}}{\cal E} 
\os{1\otimes h}{\lo} {\cal O}_{P(1)^{\rm ex}_{n-k}}\otimes_{{\cal O}_{P^{\rm ex}}}{\cal F},  
\end{align*} 
which is called the linearization of $h$. 
Thus we have a morphism 
$$L(h):=\{\vpl_nL_n(h)\}_{n=1}^{\inf} \col L({\cal E})\lo L({\cal F}).$$ 
Because ${\cal O}_{P(1)^{\rm ex}_n}$ is a locally free ${\cal O}_{P^{\rm ex}}$-module  
and the morphism 
${\cal O}_{P(1)^{\rm ex}_{n+1}}\lo {\cal O}_{P(1)^{\rm ex}_n}$ is surjective, 
\begin{align*} 
L\col & \{{\rm category\;of}\;{\cal O}_{P^{\rm ex}}{\rm -modules\;and\; 
differential\;operators}\}\lo \\
&\{{\rm stratified}\;{\cal O}_{P^{\rm ex}}{\rm -modules\;and\;} 
{\cal O}_{P^{\rm ex}}{\rm -linear\;horizontal\;morphisms}\}  
\end{align*} 
is an exact functor. 
\par 
As is well-known, the integrable connection 
$\nabla \col L({\cal E})\lo L({\cal E})\otimes_{{\cal O}_P}\Om^1_{P/S}$ 
is obtained by the stratification (\ref{eqn:stirmin}) on $L({\cal E})$ as 
follows: 
\par 
Let $\theta \col L({\cal E})\lo 
L({\cal E})\otimes_{{\cal O}_P}{\cal O}_{P(1)^{\rm ex}_1}$ 
be the following composite morphism 
\begin{align*} 
L({\cal E})\os{\vpl_n(\tau_n),~\simeq}{\lo} 
\vpl_n({\cal E}\otimes_{{\cal O}_{P^{\rm ex}}}{\cal O}_{P(1)^{\rm ex}_n})
&\lo 
\vpl_n({\cal O}_{P(1)^{\rm ex}_n}\otimes_{{\cal O}_{P^{\rm ex}}}
\vpl_n({\cal E}\otimes_{{\cal O}_{P^{\rm ex}}}{\cal O}_{P(1)^{\rm ex}_n}))\\
&\lo 
\vpl_n({\cal O}_{P(1)^{\rm ex}_n}\otimes_{{\cal O}_{P^{\rm ex}}}
{\cal E}\otimes_{{\cal O}_{P^{\rm ex}}}{\cal O}_{P(1)^{\rm ex}_1})\\
&=L({\cal E})\otimes_{{\cal O}_{P^{\rm ex}}}{\cal O}_{P(1)^{\rm ex}_1}.   
\end{align*} 
(Here we have used the locally freeness of 
${\cal O}_{P(1)^{\rm ex}_1}$ over ${\cal O}_{P^{\rm ex}}$.)
Then $\nabla$ is defined by the formula 
$\nabla(e)=\theta(e)-e\otimes 1$ $(e\in L({\cal E}))$. 
\par
As usual (\cite[(2.14)]{bob} and \cite[Proposition 3.2.14]{s1}), 
we can consider $L({\cal E})$  
as a crystal of ${\cal O}_{P^{\rm ex}/S}$-modules 
by using the stratification (\ref{eqn:stirmin}) on $L({\cal E})$. 
Set $L_{X/S}'^{\rm inf}({\cal E}):=\iota^*_{\rm inf}(L({\cal E}))$.

\begin{prop} 
Let $\lam_n \col P_n\lo P$ be the natural morphism. 
Let ${\cal E}$ be a coherent ${\cal O}_{P^{\rm ex}}$-module. 
Then there exists an isomorphism 
\begin{align*} 
L_{X/S}'^{\rm inf}({\cal E})
\os{\sim}{\lo} 
L^{\rm inf}_{X/S}(\{\lam_n^*({\cal E})\}_{n=1}^{\inf}).
\tag{6.2.8.1}\label{ali:llxe} 
\end{align*} 
\end{prop} 
\begin{proof} 
(The proof is essentially the same as that of \cite[Proposition 6.10]{bob}.)
To give a morphism (\ref{ali:llxe}) 
is equivalent to giving a morphism  
\begin{align*} 
j^*_{P}L_{X/S}'^{\rm inf}({\cal E}) 
\lo 
\varphi^*(\{\lam_n^*({\cal E})\}_{n=1}^{\inf}).
\tag{6.2.8.2}\label{ali:len}
\end{align*} 
Let $u\col (U,T)\lo (X,P)=P$ be a morphism from 
an object of ${\rm Inf}(X/S)$.    
Then 
\begin{align*} 
(j^*_{P}L_{X/S}'^{\rm inf}({\cal E}))_{(U,T)\os{u}{\lo} (X,P)}
=L_{X/S}'^{\rm inf}({\cal E})_T=u^*L({\cal E}). 
\end{align*}  
We have a morphism $u_m\col T\lo P_m$ for some $m$. 
On the other hand, 
$$(\varphi^*(\{\lam_m^*({\cal E})\}_{m=1}^{\inf}))_{(U,T)\os{u}{\lo} (X,P)}
=u_m^*\lam^*_m({\cal E}).$$ 
Hence the composite 
${\cal O}_{P(1)^{\rm ex}_n}\lo {\cal O}_{P^{\rm ex}}
\lo {\cal O}_{P^{\rm ex}_m}$ 
$(n\geq m)$
of the multiplication morphism and the projection 
induces the morphism (\ref{ali:len}).  
We claim that the morphism (\ref{ali:llxe}) is an isomorphism. 
Indeed, let $(U,T)$ be an object of ${\rm Inf}(X/S)$.  
We may assume that there exists a morphism 
$T\lo P$. 
Then 
$$L^{\rm inf}_{X/S}(\{\lam_n^*({\cal E})\}_{n=1}^{\inf})_T
=\vpl_n({\cal O}_{{\mathfrak T}_{U,n}(T\times_SP)}
\otimes_{{\cal O}_{P^{\rm ex}_n}}{{\cal O}_{P^{\rm ex}_n}}
\otimes_{{\cal O}_{P^{\rm ex}}}{\cal E})
=\vpl_n({\cal O}_{{\mathfrak T}_{U,n}(T\times_SP)}
\otimes_{{\cal O}_{P^{\rm ex}}}{\cal E})$$ 
and 
$$L_{X/S}'^{\rm inf}({\cal E})_T
={\cal O}_T\otimes_{{\cal O}_{P^{\rm ex}}}\vpl_n
({\cal O}_{P(1)^{\rm ex}_n}\otimes_{{\cal O}_{P^{\rm ex}}}{\cal E})=
\vpl_n({\cal O}_T\otimes_{{\cal O}_{P^{\rm ex}}}{\cal O}_{P(1)^{\rm ex}_n}
\otimes_{{\cal O}_{P^{\rm ex}}}{\cal E}).$$  
(The last equality 
is obtained by the explicit local description of a base of 
${\cal O}_{P(1)^{\rm ex}_n}$ 
as a sheaf of ${\cal O}_{P^{\rm ex}}$-modules 
in the proof of (\ref{prop:pl}) below.) 
\par 
Now it suffices to prove that 
$T\times_{P^{\rm ex}}P(1)^{\rm ex}_n={\mathfrak T}_{U,n}(T\times_SP)$. 
By (\ref{lemm:lexl}) below, 
the morphism $U\lo T\times_{P^{\rm ex}}P(1)^{\rm ex}_n$ is an exact immersion. 
Let $(V,T')$ be an object of ${\rm Inf}(X/S)$ over 
$(U,T\times_SP)$. Then $(V,T')$ is over $(U,T)$. 
Assume that $n$-th power of the ideal of 
definition of the exact closed immersion $V\os{\sus}{\lo} T'$ is $0$. 
Then, by using the composite morphism 
$(V,T')\lo (V,T\times_SP)\lo (U,P\times_SP)\lo (P,P\times_SP)$ 
and the universality of $P(1)^{\rm ex}_n$, we have a morphism $T'\lo P(1)^{\rm ex}_n$. 
Hence we have a morphism $T'\lo T\times_{P^{\rm ex}}P(1)^{\rm ex}_n$. 
This shows that $(U,T\times_{P^{\rm ex}}P(1)^{\rm ex}_n)$ has the universality of that of 
$(U,{\mathfrak T}_{U,n}(T\times_SP))$. 
Consequently $T\times_{P^{\rm ex}}P(1)^{\rm ex}_n={\mathfrak T}_{U,n}(T\times_SP)$. 
\end{proof} 

\begin{lemm}\label{lemm:lexl}
Let $Z_1$ and $Z_2$ be fine log formal schemes over $Z_3$. 
Let $Y$ be a fine log formal scheme over $Z_1$ and $Z_2$. 
Assume that the morphism $Y\lo Z_1$ is exact and that the morphism $Z_2\lo Z_3$ is solid. 
Then the morphism $Y\lo Z_1\times_{Z_3}Z_2$ is exact. 
\end{lemm} 
\begin{proof}  
The base change morphism 
$Y\times_{Z_3}Z_2 \lo Y$ is exact because exact morphisms are stable 
under the base change. 
By taking a local chart of the solid morphism 
$Z_2\lo Z_3$, we see that the immersion $Y\os{\subset}{\lo} Y\times_{Z_3}Z_2$ 
is exact. 
Because the morphism $Y\lo Z_1$ is exact, the base change morphism 
$Y\times_{Z_3}Z_2\lo Z_1\times_{Z_3}Z_2$ is exact. 
Because the composite morphism of exact morphisms are exact, 
the morphism $Y\lo Z_1\times_{Z_3}Z_2$ is exact. 
\end{proof} 

\par 
Let $d\col {\cal O}_{P_1(1)}\otimes_{{\cal O}_P}\Om^i_{P/S}
\lo \Om^{i+1}_{P/S}$ be the usual differential operator of order $\leq 1$. 
Then we have the following morphism 
\begin{equation*} 
L(d)\col L^{\rm inf}_{X/S}({\cal O}_{\os{\to}{P}}\otimes_{{\cal O}_P}\Om^i_{P/S})
=L_{X/S}'^{\rm inf}(\Om^i_{P/S})
\lo L_{X/S}'^{\rm inf}(\Om^{i+1}_{P/S})
=L^{\rm inf}_{X/S}({\cal O}_{\os{\to}{P}}\otimes_{{\cal O}_P}\Om^{i+1}_{P/S}).  
\end{equation*} 
As usual, we can check $L(d)\circ L(d)=0$ 
(see the formula (\ref{ali:ngdf}) below).  
\par 
The following is essentially the same as that of 
\cite[(1.2)]{cm}: 

\begin{prop}[{\bf Poincar\'{e} lemma in the log infinitesimal topos}]\label{prop:pl}
The natural morphism 
\begin{equation*}
{\cal O}_{X/S}\lo L^{\rm inf}_{X/S}({\cal O}_{\os{\to}{P}}
\otimes_{{\cal O}_P}{\Om}^{\bul}_{P/S}).  
\tag{6.2.10.1}\label{eqn:pilyk}
\end{equation*} 
is a quasi-isomorphism.  
\end{prop} 
\begin{proof}
(The proof is easier than that of \cite[(2.2.7)]{nh2}.) 
The problem  is local;  
we may assume that ${\Om}_{P/S}^1$ has a basis 
$\{d\log \tau_j\}_{j=1}^d$, where $\tau_j$ 
is  a local section of $M_P$. 
We may assume that there exists a morphism 
$T \lo P$ fitting into the following commutative diagram 
\begin{equation*} 
\begin{CD} 
U@>{\subset}>> T \\ 
@VVV @VVV \\ 
X@>{\subset}>> P. 
\end{CD} 
\end{equation*} 
The immersion $U=U\times_XX\os{\sus}{\lo} T\times_SP$ is 
nothing but the composite morphism of 
the base change morphism 
of the diagonal immersion 
$P \os{\sus}{\lo} P\times_SP$ 
by the morphism 
$T\times_SP\lo P\times_SP$
and the immersion $U\os{\sus}{\lo} T$.  
Let $u_j$ be a local section of 
${\rm Ker}({\cal O}_{P_1(1)}^* 
\lo {\cal O}_P^*)$ such 
that $p_2^*(\tau_j)=p_1^*(\tau_j)u_j$, 
where $p_i \col P_1(1) \lo P$ 
is the induced morphism by the $i$-th projection 
$(i=1,2)$. 
Then, by the same proof as that of \cite[(6.6)]{klog1}, 
the following morphism 
$${\cal O}_T[[s_1,\ldots,s_d]]\owns s_j
\lom u_j-1 \in 
\vpl_n{\cal O}_{{\mathfrak T}_{U,n}(T\times_SP)}$$
is an isomorphism, where $s_j$'s are 
independent indeterminates.  
Hence 
\begin{equation*} 
L^{\rm inf}_{X/S}(\{{\cal O}_{P_n}\otimes_{{\cal O}_P}{\Om}^q_{P/S}\}_{n=1}^{\inf})_T
\simeq 
{\cal O}_T[[s_1,\ldots,s_d]]\otimes_{{\cal O}_P}
{\Om}^q_{P/S}. 
\tag{6.2.10.2}\label{eqn:lops}
\end{equation*} 
Let $[n]$ 
$(n\in {\mab Z}_{\geq 1})$ mean the divided power. 
As in \cite[(7)]{cm} and \cite[(2.2.7)]{nh2}, the boundary morphism 
$L^{\rm inf}_{X/S}({\Om}^q_{P/S})_T \lo 
L^{\rm inf}_{X/S}({\Om}^{q+1}_{P/S})_T$ is given by 
the following formula:  
\begin{equation*}
\nabla^q_T(as_1^{[i_1]}\cdots s_d^{[i_d]}\otimes 
\om)  = 
a(\sum_{j=1}^d s_1^{[i_1]}\cdots s_j^{[i_j-1]} 
\cdots s_d^{[i_d]}
(s_j+1) d\log 
\tau_j\wedge \om 
\tag{6.2.10.3}
\label{ali:ngdf} 
\end{equation*}
$$+s_1^{[i_1]}
\cdots s_d^{[i_d]}\otimes d\om) \quad (a\in 
{\cal O}_T,~i_1, \ldots, i_d 
\in {\mab N},~\om \in {\Om}^q_{P/S})
$$ 
(cf.~\cite[6.11 Lemma]{bob}). 
\par 
Now consider the case $d=1$ and set $s=s_1$ and 
$\tau=\tau_1$.  
Then the complex 
${\cal O}_T[[s_1,\ldots, s_d]]
\otimes_{{\cal O}_P}
{\Om}_{P/S}^{\bul}$ 
is equal to 
${\cal O}_T[[s]]
\os{\nabla^0_T}{\lo} {\cal O}_T[[s]]d\log t$. 
Because $\nabla^0_T(s^{[n]})=s^{[n-1]}(s+1)d\log \tau
=(ns^{[n]}+s^{[n-1]})d\log \tau$ 
for a positive integer $n$, we have the following formula    
\begin{align*}
\nabla^0_T(\sum_{n=0}^{\infty}a_ns^{[n]}) = & 
\sum_{n=1}^{\infty}(a_n+(n-1)a_{n-1})s^{[n-1]}d\log \tau
\quad (a_n \in {\cal O}_T).
\tag{6.2.10.4}
\label{ali:nbof}\\ 
\end{align*} 
Hence ${\rm Ker}(\nabla^0_T)={\cal O}_T$ 
and ${\rm Coker}(\nabla^0_T)=0$. 
Therefore we have 
checked that the morphism (\ref{eqn:pilyk}) is a quasi-isomorphism 
for the case $d=1$. 
\par 
The rest of the proof is the same as that 
of \cite[6.12 Theorem]{bob}.
\end{proof} 

\begin{prop}[{\bf Poincar\'{e} lemma 
for a coherent crystal}]\label{prop:lpe}
Let $E$ be a coherent crystal of 
${\cal O}_{X/S}$-modules. 
Set ${\cal E}_n:=E_{P_n}$ and 
${\cal E}:=\{{\cal E}_n\}_{n=1}^{\infty}$. 
Then the natural morphism 
\begin{equation*} 
E \lo 
L^{\rm inf}_{X/S}
({\cal E}\otimes_{{\cal O}_P}{\Om}^{\bul}_{P/S}) 
\tag{6.2.11.1}\label{exsop}
\end{equation*} 
is a quasi-isomorphism. 
\end{prop}
\begin{proof}
The proof is the same as that of \cite[Proposition 6.15]{bob}. 
\end{proof}


\begin{coro}\label{coro:ztx}
Let $X_{\rm zar}$ be the Zariski topos of $X$. 
Let $u^{\rm inf}_{X/S} \col (X/S)_{\rm inf} \lo X_{\rm zar}$ 
be the natural projection. 
Let $\bet_n \col P_{n,{\rm zar}} 
\os{\sim}{\lo} X_{\rm zar}$ be the natural equivalence of topoi. 
Using $\bet_n$, we identify 
$P_{n,{\rm zar}}$ with  $X_{\rm zar}$. 
Then 
\begin{equation*}
Ru^{\rm inf}_{X/S*}(E) =\vpl_n
({\cal E}_n\otimes_{{\cal O}_P}{\Om}^{\bul}_{P/S}).  
\tag{6.2.12.1}\label{eqn:olyzk}
\end{equation*}  
\end{coro}
\begin{proof}
By (\ref{exsop}) 
\begin{align*} 
Ru^{\rm inf}_{X/S*}(E)=Ru^{\rm inf}_{X/S*}(L^{\rm inf}_{X/S}
({\cal E}\otimes_{{\cal O}_P}{\Om}^{\bul}_{P/S})).
\tag{6.2.12.2}\label{ali:uhagh}
\end{align*}   
By (\ref{cd:xuz}) we have 
$(\vpl_n\bet_{n*}){\varphi}_{\os{\to}{P}*}
=u^{\rm inf}_{X/S*}j_{P*}$. 
Hence we have the following equalities 
for an ${\cal O}_{\os{\to}{P}}$-module 
${\cal F}=\{{\cal F}_n\}_{n=1}^{\infty}$ 
by the definition (\ref{eqn:luys}) and by (\ref{coro:je}): 
\begin{align*} 
Ru^{\rm inf}_{X/S*}(L^{\rm inf}_{X/S}({\cal F}))
&=
u^{\rm inf}_{X/S*}L^{\rm inf}_{X/S}({\cal F})
= u^{\rm inf}_{X/S*}j_{P*}{\varphi}^*({\cal F})
\tag{6.2.12.3}\label{ali:uhiph}\\
&= 
(\vpl_n\bet_{n*}){\varphi}_{\os{\to}{P}*}{\varphi}^*({\cal F}) 
= \vpl_n\bet_{n*}({\cal F}_n).  
\end{align*} 
By (\ref{ali:uhagh}) and (\ref{ali:uhiph}) 
we obtain (\ref{eqn:olyzk}).
\end{proof}

\par 
Now we consider the general case: the case where 
there does not 
necessarily exist an immersion 
$X\os{\sus}{\lo} P$ into a log smooth scheme over $S$. 
Let $X_{\bul}$ be the \v{C}ech diagram obtained by an affine open covering of $X$. 
Let $X_{\bul} \os{\sus}{\lo} P_{\bul}$ be 
a simplicial immersion into a log smooth simplicial scheme 
over $S$ (cf.~\cite[(2.18)]{hk}). 
Let 
\begin{equation*} 
\pi_{\rm inf} \col ((X_{\bul}/S)_{\rm inf},{\cal O}_{X_{\bul}/S})
\lo ((X/S)_{\rm inf},{\cal O}_{X/S})
\end{equation*} 
be the natural morphism of ringed topoi. 
Let 
\begin{equation*} 
\pi_{\rm zar} \col (X_{\bul})_{\rm zar}\lo X_{\rm zar}
\end{equation*} 
be also the natural morphism of ringed topoi. 
Set ${\cal E}^{\bul}_n
:=(\pi^{-1}_{\rm inf}(E))_{P_{\bul n}}$. 
We consider 
$\vpl_nR\pi_{{\rm zar}*}
({\cal E}^{\bul}_n
\otimes_{{\cal O}_{P_{\bul}}}{\Om}^{\bul}_{P_{\bul}/S})$ 
as in \cite[(0.14)]{cm},  
though Chiarellotto and Fornasiero 
have considered only the trivial coefficient in [loc.~cit.].  

The following is the relative and sheafied version of 
\cite[(39)]{cm} for a fine log scheme: 

\begin{coro}[{\bf {\rm {\bf cf.~\cite[(39)]{cm}}}}]\label{coro:co}
There exists the following canonical isomorphism$:$ 
\begin{equation*} 
Ru^{\rm inf}_{X/S*}(E) \os{\sim}{\lo} 
\vpl_nR\pi_{{\rm zar}*}({\cal E}^{\bul}_n 
\otimes_{{\cal O}_{P_{\bul}}}{\Om}^{\bul}_{P_{\bul}/S}). 
\tag{6.2.13.1}\label{eqn:ifxs} 
\end{equation*} 
\end{coro} 
\begin{proof}
(\ref{eqn:ifxs}) immediately follows from (\ref{eqn:olyzk}) 
and the cohomological descent.  
\end{proof} 

\begin{coro}\label{coro:vn} 
Set 
$f^{\rm inf}_{X/S}:=f\circ u^{\rm inf}_{X/S} 
\col ((X/S)_{\rm inf},{\cal O}_{X/S}) 
\lo (S_{\rm zar},{\cal O}_S)$. 
Assume that $\os{\circ}{S}$ is quasi-compact and that 
$\os{\circ}{f}$ is quasi-separated and quasi-compact. 
Then $Rf^{\rm inf}_{X/S*}(E)$ is bounded. 
\end{coro}
\begin{proof} 
By using a standard argument, 
(\ref{coro:vn}) follows  from (\ref{coro:co}) as in 
\cite[Theorem 7.6]{bob}. 
\end{proof}

The following is a log infinitesimal version of 
the base change theorem in the (log) crystalline cohomology 
(\cite[Theorem 7.8]{bob}, \cite[(6.10)]{klog1}):

\begin{prop}[{\bf Base change}]\label{prop:bcpif}
Let $u\col S'\lo S$ be a morphism of 
fine log schemes over $K$. 
Set $X':=X\times_SS'$ and let $p\col X'\lo X$ 
be the first projection. 
Assume that $\os{\circ}{S}$ is quasi-compact, that 
$\os{\circ}{f}$ is quasi-separated and quasi-compact 
and that $f$ is log smooth and integral. 
Assume that $E$ is a flat ${\cal O}_{X/S}$-module. 
Then the canonical  morphism 
\begin{equation*} 
Lu^*Rf^{\rm inf}_{X/S*}(E)
\lo Rf^{\rm inf}_{X'/S'*}(p^*_{\rm inf}(E)) 
\tag{6.2.15.1}\label{eqn:wtyif}
\end{equation*}
is an isomorphism. 
\end{prop}
\begin{proof} 
By using (\ref{eqn:ifxs}), we can give the same proof 
as that of \cite[Theorem 7.8]{bob} by using the simplical formal log scheme 
$P^{\rm ex}_{\bul}$. 
\end{proof}

\begin{rema}
(\ref{prop:bcpif}) is a generalization of \cite[III (5.2)]{hadr}. 
\end{rema}


\par 
Set 
\begin{equation*} 
Rf^{\rm inf}_{X/S*}(E):=Rf_*Ru^{\rm inf}_{X/S*}(E),   
\end{equation*} 
\begin{equation*} 
R^qf^{\rm inf}_{X/S*}(E):=
{\cal H}^q(Rf^{\rm inf}_{X/S*}(E)),   
\end{equation*} 
\begin{equation*} 
H^q_{\rm inf}(X/S,E):=
H^q(R\Gam(S,Rf^{\rm inf}_{X/S*}(E)))   
\end{equation*} 
and 
\begin{equation*} 
H^q_{\rm inf}(X/S):=H^q_{\rm inf}(X/S,{\cal O}_{X/S}). 
\end{equation*} 

\begin{prob}
The following problems seem interesting: 
\par 
(1) Is $R^qf^{\rm inf}_{X/S*}({\cal O}_{X/S})$ 
($q\in {\mab N}$)
a coherent ${\cal O}_S$-module if $X/S$ is of finite type?  
More generally, is $R^qf^{\rm inf}_{X/S*}(E)$ 
a coherent ${\cal O}_S$-module for a coherent ${\cal O}_{X/S}$-module 
$E$ if $X/S$ is of finite type?  
\par 
(2) 
Consider the case $K={\mab C}$. 
Does there exist a comparison theorem between 
$Rf_{X^{\log}_{\rm an}/S_{\rm an}*}({\mab C})$ and  
$Rf^{\rm inf}_{X/S*}({\cal O}_{X/S})_{\rm an}$ 
as in the absolute case in \cite[Theorem 5.1]{cm}? 
Or can one find a dense open subset of $S$ such that 
the problem above hold for $U$ and $f^{-1}(U)$ 
as in \cite[IV (4.3)]{hadr}?  
\end{prob}

\par 
Let $N$ be a nonnegative integer or $\infty$. 
Set $S=({\rm Spec}(K),K^*)$ and 
$R\Gam_{\rm inf}(X/K):=Rf^{\rm inf}_{X/S*}({\cal O}_{X/S})$. 
Let $X_{\bul \leq N}$ be 
an $N$-truncated proper hypercovering of $X$ over $K$. 

\begin{prop}\label{prop:hnn}
Assume that $M_X$ is trivial and that the structural morphism 
$f\col X\lo {\rm Spec}(K)$ is of finite type. 
Let $c$ be a nonnegative integer 
such that $H^q_{\rm inf}(X/K)=0$ for any $q>c$.  
$($The existence of $c$ will also be shown in the proof below.$)$
Let $N$ be a nonnegative integer satisfying 
the inequality in {\rm \cite[(2.2.1)]{nh3}} for $c$, that is, 
\begin{equation*}
N > 2^{-1}(c+1)(c+2). 
\tag{6.2.18.1}\label{eqn:shdf}
\end{equation*}
Set $R\Gam_{\rm inf}(X/K):=
R\Gam(X_{\rm zar},Ru^{\rm inf}_{X/K*}({\cal O}_{X/K}))$.  
Then the natural morphism 
\begin{equation*} 
R\Gam_{\rm inf}(X/K)
=\tau_cR\Gam_{\rm inf}(X/K) 
\lo \tau_cR\Gam_{\rm inf}(X_{\bul \leq N}/K)
\tag{6.2.18.2}\label{eqn:wcdyif}
\end{equation*} 
is an isomorphism. 
\end{prop}
\begin{proof} 
By (\ref{prop:bcpif}) we may assume that $K$ is algebraically closed. 
\par 
The first proof: 
Assume first that $K={\mab C}$. 
Then, by \cite{grcr} and \cite{hadr}, we have 
the following functorial isomorphisms: 
\begin{equation*} 
H^q_{\rm inf}(X/{\mab C}) \os{\sim}{\lo} 
H^q_{\rm dR}(X/{\mab C}) \os{\sim}{\lo} 
H^q(X_{\rm an},{\mab C}). 
\tag{6.2.18.3}\label{eqn:xcainc} 
\end{equation*} 
Since  $(X_{\bul \leq N})_{\rm an}$ is 
an $N$-truncated proper hypercovering of 
$X_{\rm an}$ over ${\mab C}$, 
we have the following commutative diagram 
by the proper cohomological descent: 
\begin{equation*} 
\begin{CD} 
H^q_{\rm inf}(X_{\bul \leq N}/{\mab C})
@>{\sim}>> H^q((X_{\bul \leq N})_{\rm an},{\mab C}) \\
@AAA @AA{\simeq}A \\
H^q_{\rm inf}(X/{\mab C}) @>{\sim}>> 
H^q(X_{\rm an},{\mab C}). 
\end{CD}
\tag{6.2.18.4}\label{eqn:xcanc} 
\end{equation*}  
(By (\ref{eqn:shdf}) and by \cite[(2.3)]{nh3}, 
we have the right vertical isomorphism.) 
Hence we have the isomorphism (\ref{eqn:wcdyif}) 
in the case $K={\mab C}$. 
Using (\ref{prop:bcpif}), 
we obtain the isomorphism (\ref{eqn:wcdyif}) 
by the Lefschetz principle. 
\par 
The second proof: 
(The following argument is much simpler 
than that of the proof of \cite[(4.4), (4.5)]{tzp}
and \cite[(10.9)]{nh3}.) 
\par  
Let $N$ be a positive integer. 
Later we assume that $N$ is 
a positive integer satisfying the inequality (\ref{eqn:shdf}).
\par 
Let $U$ be the disjoint union of 
an affine open covering of $X$. 
Let $U \os{\sus}{\lo} {\cal U}$ be a closed immersion 
into a smooth scheme over $K$. 
Set $U_{\bul}:={\rm cosk}^X_0(U)$ and 
${\cal U}_{\bul}:={\rm cosk}^K_0({\cal U})$. 
The fiber product $X_{\bul \leq N}\times_XU$  
is an $N$-truncated proper hypercovering of $U$.
Let $V_{\bul \leq N}$ be a refinement  of 
the proper hypercovering 
$X_{\bul \leq N}\times_XU$ (\cite[(4.2.1)]{tzp}) 
of $U$ such that there exists a closed immersion
$V_N \os{\sus}{\lo} P_N$
into a smooth scheme over $K$
(cf.~\cite[(4.2.3)]{tzp}).
Consider the $N$-truncated couple
$(V_{\bul \leq N},\Gam^K_N(P_N)_{\bul \leq N})$. 
The $N$-truncated simplicial $K$-scheme 
$\Gam^K_N(P_N)_{\bul \leq N}$ 
contains $V_{\bul \leq N}$ as 
an $N$-truncated simplicial closed subscheme over $K$. 
Set $Q_{\bul}:=\Gam^K_N(P_N)$.
\par
Following the idea in \cite[(4.4)]{tzp}, 
we consider the pair
$(V_{\bul \leq N, \bul},Q_{\bul \leq N, \bul})$
of $(N,\infty)$-truncated bisimplicial schemes
defined by 
\begin{align*}
(V_{mn},Q_{mn}) := & 
({\rm cosk}_0^{X_m}(V_m)_n,
{\rm cosk}_0^K(Q_m\times_K{\cal U})_n) 
\quad (0\leq m \leq N, n \in {\mab N}) 
\tag{6.2.18.5}\label{eqn:vbscq}
\end{align*}
with the natural morphisms which make
$(V_{\bul \leq N,\bul},Q_{\bul \leq N,\bul})$ 
a couple of $(N,\infty)$-truncated bisimplicial schemes.  
\par
We have a natural 
morphism 
$(V_{\bul \leq N, \bul},Q_{\bul \leq N, \bul}) 
\lo
(U_{\bul},{\cal U}_{\bul})$ 
of pairs and have 
the following commutative diagram of 
($N$-truncated) (bi)simplicial schemes:
\begin{equation*}
\begin{CD}
X_{\bul \leq N}
@<<< 
V_{\bul \leq N, \bul}\\
@V{}VV @VV{}V  \\
X @<<< U_{\bul}.
\end{CD}\label{cd:uvbnb}
\tag{6.2.18.6}
\end{equation*}
By the proof of \cite[(4.4.1) (1)]{tzp}
the morphism 
$V_{\bul \leq N,n}\lo U_n \quad (n\in {\mab N})$ 
is an $N$-truncated proper hypercovering of $U_n$. 
\par
Let $q$ be a fixed arbitrary nonnegative integer. 
Let $N$ be an integer satisfying 
the inequality (\ref{eqn:shdf}).
Then we claim that the morphism
\begin{equation*}
R\Gam_{\rm inf}(U_{\bul}/K)
\lo 
R\Gam_{\rm inf}(V_{\bul \leq N, \bul}/K)
\tag{6.2.18.7}
\end{equation*}
induces an isomorphism
\begin{equation*}
H^q_{\rm inf}(X/K)=
H^q(R\Gam_{\rm inf}(U_{\bul}/K))
\os{\sim}{\lo}
H^q(R\Gam_{\rm inf}(V_{\bul \leq N, \bul}/K)).
\tag{6.2.18.8}\label{eqn:chrchc}
\end{equation*}
Indeed, we have the following two spectral sequences:
\begin{equation*}
E^{ij}_1=H^j_{\rm inf}(U_i/K) \Lo H^{i+j}_{\rm inf}(U_{\bul}/K),
\tag{6.2.18.9}\label{eqn:eij1hdcz}
\end{equation*}
\begin{equation*}
E^{ij}_1=
H^j(R\Gam_{\rm inf}(V_{\bul \leq N,i}/K)) 
\Lo 
H^{i+j}(R\Gam_{\rm inf}(V_{\bul \leq N, \bul}/K)).
\tag{6.2.18.10}\label{eqn:vctrs}
\end{equation*}
By the analogue of \cite[(2.1.3)]{tzp}, 
the $E_r$-terms 
$(1 \leq r \leq \infty)$ 
of (\ref{eqn:eij1hdcz}) are 
isomorphic to those of (\ref{eqn:vctrs}) 
for $i+j \leq q$. 
Hence we have an isomorphism (\ref{eqn:chrchc}).
\par
So far we have proved that,
for any integer $q$, there exists a 
sufficiently large integer $N$ 
depending on $q$ such that
there exists an isomorphism 
\begin{equation*}
\tau_qR\Gam_{\rm inf}(X/K) \os{\sim}{\lo} 
\tau_qR\Gam_{\rm inf}(V_{\bul \leq N, \bul}/K).
\tag{6.2.18.11}\label{eqn:tahth}
\end{equation*}
\par
Let $N$ be any integer satisfying the inequality 
(\ref{eqn:shdf}).
Then, by (\ref{eqn:tahth}), we have
\begin{equation*}
R\Gam_{\rm inf}(X/K) =
\tau_cR\Gam_{\rm inf}(X/K)\os{\sim}{\lo}
\tau_cR\Gam_{\rm inf}(V_{\bul \leq N, \bul}/K). 
\tag{6.2.18.12}\label{eqn:rrgukc}
\end{equation*} 
By the  proof of 
\cite[(4.4.1) (2)]{tzp}, 
for a nonnegative integer $m \leq N$,  
\begin{equation*}
R\Gam_{\rm inf}(X_m/K)=
R\Gam_{\rm inf}(V_{m\bul}/K).
\tag{6.2.18.13}\label{eqn:rrgetb}
\end{equation*}
Hence the natural morphism 
$V_{\bul \leq N, \bul} \lo X_{\bul \leq N}$ 
induces the following isomorphism 
\begin{equation*} 
\tau_cR\Gam_{\rm inf}(X_{\bul \leq N}/K)
\os{\sim}{\lo} 
\tau_cR\Gam_{\rm inf}(V_{\bul \leq N, \bul}/K). 
\tag{6.2.18.14}\label{eqn:rrgietb}
\end{equation*}  
By (\ref{eqn:rrgukc}) and (\ref{eqn:rrgietb}), 
we have the isomorphism (\ref{eqn:wcdyif}). 
\end{proof}

\begin{coro}\label{coro:nxki}
Let $X_{\bul}\lo X$ be a proper hypercovering of $X$. 
Then the following hold$:$
\par 
$(1)$ 
The natural morphism 
\begin{equation*} 
R\Gam_{\rm inf}(X/K) \lo R\Gam_{\rm inf}(X_{\bul}/K)
\tag{6.2.19.1}\label{eqn:xkif}
\end{equation*} 
is an isomorphism. 
\par 
$(2)$ Assume that $X_{\bul}$ is smooth over $K$. 
Then there exists the following canonical isomorphism 
\begin{equation*} 
R\Gam_{\rm inf}(X/K) \os{\sim}{\lo} 
R\Gam(X_{\bul},\Om^{\bul}_{X_{\bul}/K}). 
\tag{6.2.19.2}\label{eqn:xbom}
\end{equation*} 
\end{coro}
\begin{proof} 
This is easy to prove. 
Indeed, 
because $q$ in (\ref{eqn:tahth}) is arbitrary, 
$H^q(R\Gam_{\rm inf}(X_{\bul}/K))=0$ for $q> c$. 
Hence $R\Gam_{\rm inf}(X_{\bul}/K)
=\tau_cR\Gam_{\rm inf}(X_{\bul}/K)$. 
Hence we have the isomorphism 
(\ref{eqn:xkif}) by (\ref{eqn:wcdyif}). 
By (1) and (\ref{coro:co})  
we have the following composite isomorphism 
\begin{equation*} 
R\Gam_{\rm inf}(X/K) \os{\sim}{\lo} 
R\Gam_{\rm inf}(X_{\bul}/K) 
\os{\sim}{\lo} R\Gam(X_{\bul},\Om^{\bul}_{X_{\bul}/K}). 
\tag{6.2.19.3}\label{eqn:ryfci} 
\end{equation*} 
\end{proof}

\begin{theo}\label{theo:fdyef} 
Let $Y_{\bul}$ be a proper smooth simplicial scheme  over $K$
and let $D_{\bul}$ be an SNCD  on $Y_{\bul}/K$. 
Set $X_{\bul}:=Y_{\bul}\setminus D_{\bul}$. 
Assume that $X_{\bul}$ is a proper smooth hypercovering of $X$ over $K$.  
Let $F$ be the filtration on  
$H^q_{\rm inf}(X/K)$ $(q\in {\mab N})$ defined by 
the filtered complex 
$(\Om^{\bul}_{Y_{\bul}/K}(\log D_{\bul}),
\{\Om^{\bul \geq i}_{Y_{\bul}/K}(\log D_{\bul})
\}_{i\in {\mab Z}})$ 
and the isomorphism {\rm (\ref{eqn:ryfci})}. 
Then $F$ is independent of the choice of 
$Y_{\bul}$, $D_{\bul}$ and $X_{\bul}$. 
$($It depends only on $X/K.)$  
\end{theo}
\begin{proof} 
Let $X'_{\bul}$ be another proper hypercovering of $X/K$ 
such that $X'_{\bul}$ is the complement of 
a simplicial SNCD $D'_{\bul}$ on 
a proper smooth simplicial scheme $Y'_{\bul}$ over $K$. 
Because two proper hypercoverings of $X/K$ is covered 
by another proper hypercovering of $X/K$ 
(\cite[${\rm V}^{\rm bis}$ (5.1.7), (5.1.3)]{sga4-2},  
\cite{dh3}), 
we may assume that there exists a morphism 
$(Y'_{\bul},D'_{\bul})\lo (Y_{\bul},D_{\bul})$ 
which gives the following commutative diagram
\begin{equation*} 
\begin{CD} 
X'_{\bul} @>>> X_{\bul} \\ 
@VVV @VVV \\ 
X@= X. 
\end{CD} 
\end{equation*} 
We may assume that $K={\mab C}$ by the Lefschetz principle. 
Let $(j_{\bul})_{\rm an} \col (X_{\bul})_{\rm an} 
\os{\sus}{\lo} (Y_{\bul})_{\rm an}$ 
be the natural open immersion. 
In \cite[(3.1.8)]{dh2} Deligne has proved 
that the natural inclusion morphism 
\begin{equation*} 
\Om^{\bul}_{Y_{\bul{\rm an}}/{\mab C}}
(\log D_{\bul{\rm an}}) 
\os{\sus}{\lo}  
j_{\bul {\rm an}*}\Om^{\bul}_{X_{\bul{\rm an}}/{\mab C}}
\end{equation*} 
is a quasi-isomorphism. 
Hence the natural morphism 
\begin{equation*} 
R\Gam(Y_{\bul{\rm an}},
\Om^{\bul}_{Y_{\bul {\rm an}}/{\mab C}}
(\log D_{\bul {\rm an}})) \lo 
R\Gam(X_{\bul {\rm an}},
\Om^{\bul}_{X_{\bul{\rm an}}/{\mab C}})
\end{equation*} 
is an isomorphism. 
Let $(R\Gam(X_{\bul{\star}},\Om^{\bul}_{X_{\bul{\star}}/K}),F)
\in {\rm D}^+{\rm F}(K)$  
be the filtered complex defined by 
the filtered complex 
$(\Om^{\bul}_{X_{\bul{\star}}/K},
\{\Om^{\bul \geq i}_{X_{\bul{\star}}/K}\}_{i\in {\mab Z}})$ (${\star}={\rm an}$ or nothing).  
By GAGA 
we have the following commutative diagram 
\begin{equation*} 
\begin{CD}
(R\Gam(Y_{\bul},
\Om^{\bul}_{Y_{\bul}/K}(\log D_{\bul})),F) @>>> 
(R\Gam(Y'_{\bul},
\Om^{\bul}_{Y'_{\bul}/K}(\log D'_{\bul})),F) \\ 
@V{\simeq}VV @VV{\simeq}V \\
(R\Gam(Y_{\bul{\rm an}},
\Om^{\bul}_{Y_{\bul{\rm an}}/K}(\log D_{\bul{\rm an}})),F) @>>> 
(R\Gam(Y'_{\bul{\rm an}},
\Om^{\bul}_{Y'_{\bul{\rm an}}/K}(\log D'_{\bul{\rm an}})),F). 
\end{CD}
\end{equation*} 
Consequently we have the following natural morphism 
$(H^q(Y_{\bul},
\Om^{\bul}_{Y_{\bul}/K}(\log D_{\bul})),F) 
\lo (H^q(Y'_{\bul},
\Om^{\bul}_{Y'_{\bul}/K}(\log D'_{\bul})),F)$.  
Hence we have an isomorphism 
$F^iH^q(Y_{\bul},
\Om^{\bul}_{Y_{\bul}/K}(\log D_{\bul}))
\os{\sim}{\lo} F^iH^q(Y'_{\bul},
\Om^{\bul}_{Y'_{\bul}/K}(\log D_{\bul}))$ 
by the theory of mixed Hodge structures by Deligne (\cite{dh3}).
\end{proof}

\begin{defi}\label{defi:hfi} 
We call the filtration $F$ on 
$H^q_{\rm inf}(X/K)$ $(q\in {\mab N})$ 
the {\it Hodge filtration} on $H^q_{\rm inf}(X/K)$.  
\end{defi} 

\begin{prop}\label{prop:lnadv} 
Let the notations and the assumptions be as in 
{\rm (\ref{theo:fdyef})}. 
Then  the boundary morphism 
\begin{equation*} 
H^j(Y_{\bul},\Om^i_{Y_{\bul}/K}(\log D_{\bul})) \lo 
H^j(Y_{\bul},\Om^{i+1}_{Y_{\bul}/K}(\log D_{\bul})) 
\quad (i,j\in {\mab N}) 
\end{equation*} 
vanishes. 
\end{prop} 
\begin{proof} 
\par
The first proof: By the Lefschetz principle, 
we may assume that $K={\mab C}$. 
This is a consequence of theory of mixed Hodge structures. 
\par 
The second proof: 
Fix $i$ and $j$. Let $N$ be a positive integer such that 
$N > 2^{-1}(i+j+2)(i+j+3)$. 
Then $H^j(Y_{\bul},\Om^k_{Y_{\bul}/K}(\log D_{\bul}))
=H^j(Y_{\bul \leq N},\Om^k_{Y_{\bul \leq N}/K}(\log D_{\bul \leq N}))$ 
for $k=i$ and $k=i+1$. 
By the standard technique explained in 
the proof of \cite[(2.7)]{di}, we may assume that $K$ is a field whose characteristic is a large prime number $p$ with respect to $i+j$ and $N$. 
In this case there exists the inverse log Cartier isomorphism 
$$C^{-1}\col \Om^k_{Y_{\bul \leq N}/K}(\log D_{\bul \leq N})
\lo {\cal H}^k(\Om^{\bul}_{Y_{\bul \leq N}/K}(\log D_{\bul \leq N}))$$
on $Y_{\bul \leq N}$.  
Hence the arguments in the proofs of \cite[(2.1)]{di} and \cite[(2.7)]{di} 
show the desired vanishing. 
\end{proof} 

\begin{theo}\label{theo:hysst}
Let $g\col X\lo Y$ be a morphism of schemes of finite type over $K$. 
Then the pull-back 
\begin{align*} 
g^*\col H^q_{\rm inf}(Y/K)\lo H^q_{\rm inf}(X/K) \quad (q\in {\mab N})
\end{align*} 
of $g$ is strictly compatible 
with the Hodge filtration. 
\end{theo}
\begin{proof}
We may assume that $K$ is algebraically closed. 
By the Lefschetz principle, we may assume that 
$K={\mab C}$. 
Then (\ref{theo:hysst}) is a consequence of theory of mixed Hodge structures. 
\end{proof}

\begin{rema}\label{rema:ll}
As in the $E_1$-degeneration of the Hodge de Rham spectral sequence 
of a separated scheme of finite type over a field of characteristic 0
((\ref{prop:lnadv})), it is a very interesting problem 
to prove (\ref{theo:fdyef}) and (\ref{theo:hysst}) algebraically. 
\end{rema}

\section{Monodromy operators on the infinitesimal 
cohomologies of proper schemes in mixed characteristics}\label{sec:pmm}
Let $K$ and $K_0$ be as in \S\ref{sec:cfi} and 
let ${\mathfrak X}$ be a proper scheme over $K$. 
Let $L/K$ be an extension of fields.  
Set ${\mathfrak X}_L:={\mathfrak X}\otimes_KL$.  
In this section we define a monodromy operator 
on the infinitesimal cohomology $H^q_{\rm inf}({\mathfrak X}_L/L)$ 
$(q\in {\mab N})$ for a certain finite extension $L$ of $K$. 
In \S\ref{sec:pff} below 
we shall define a  monodromy operator over $K_0$. 
\par 
We recall two Tsuji's results for the construction of 
the monodromy operator above. 
The first his result (\ref{theo:hkt}) below is 
an $N$-truncated simplicial version 
of the Hyodo-Kato isomorphism in \cite[(5.1)]{hk}. 
The second his result (\ref{prop:compb}) below is 
a result about  the behavior of 
the ($N$-truncated simplicial) Hyodo-Kato isomorphisms  
under a finite extension of complete discrete valuation rings 
of mixed characteristics with perfect residue fields. 
\par 
Let ${\cal V}$, ${\cal W}$ and $\kap$ be as in \S\ref{sec:cfi}. 
Let $N$ be a nonnegative integer.  
Endow ${\rm Spec}({\cal V})$ with the canonical log structure 
and let $S$ be the resulting log scheme. 
Let $e$ be the absolute ramification index of ${\cal V}$. 
Set $s:=S\otimes_{\cal V}\kap$. 
Let ${\cal Y}_{\bul \leq N}$ 
be a log smooth $N$-truncated simplicial log scheme over $S$.  
Let $Y_{\bul \leq N}/s$ be the log special fiber of ${\cal Y}_{\bul \leq N}/S$. 
Assume that 
$\os{\circ}{\cal Y}_{\bul \leq N}/\os{\circ}{S}$ 
is proper and that the structural morphism 
$Y_{\bul \leq N} \lo s$ is of Cartier type. 
Let ${\mathfrak Y}_{\bul \leq N}$ be 
the log generic fiber of ${\cal Y}_{\bul \leq N}$. 
\par

\begin{theo}[{\rm {\bf \cite[(6.3.1), (6.3.2)]{tsgep}}}]\label{theo:hkt} 
Let $\pi$ be a uniformizer of ${\cal V}$. 
Let $[~] \col \kap \lo {\cal W}$ 
be the Teim\"{u}ller representative. 
For an element $a\in {\cal V}^*$, set $\ol{a}:=a~{\rm mod}~\pi$ 
and $\log (a):=\sum_{n\geq 1}(-1)^{n-1}n^{-1}(a[\ol{a}]^{-1}-1)^n$. 
Then there exists the following $($natural$)$ isomorphism 
\begin{equation*} 
\Psi_{\pi} \col R\Gam_{{\rm crys}}(Y_{\bul \leq N}/{\cal W}(s))
\otimes_{{\cal W}}^LK 
\os{\sim}{\lo} R\Gam_{\rm dR}({\mathfrak Y}_{\bul \leq N}/K)
\tag{6.3.1.1}\label{eqn:rghk}
\end{equation*} 
depending on $\pi$ such that 
\begin{equation*} 
H^q(\Psi_{a\pi})=H^q(\Psi_{\pi})\circ H^q({\rm exp}(\log (a)N_{\rm zar})) \col 
H^q_{{\rm crys}}(Y_{\bul \leq N}/{\cal W}(s))
\otimes_{{\cal W}}K 
\os{\sim}{\lo} H^q_{\rm dR}({\mathfrak Y}_{\bul \leq N}/K)
\tag{6.3.1.2}\label{eqn:rghenk}
\end{equation*} 
for $a\in {\cal V}^*$. 
Here $N_{\rm zar}$ is the monodromy operator 
{\rm (\ref{ali:lqfnk})} of $Y_{\bul \leq N}/({\cal W}(s),p{\cal W},[~])$.  
\end{theo}

\begin{rema}\label{rema:nhk} 
First consider the case $N=0$. 
In \cite[(5.1)]{hk} only the existence of the Hyodo-Kato isomorphism 
\begin{equation*} 
H^q(\Psi_{\pi}) \col H^q(Y_0/{\cal W}(s))\otimes_{{\cal W}}K 
\os{\sim}{\lo} H^q_{\rm dR}({\mathfrak Y}_0/K) \quad (q\in {\mab N})
\tag{6.3.2.1}\label{eqn:rgyhk}
\end{equation*} 
has been claimed.  
\par 
Next consider the case where $N$ is general. 
Then, though the isomorphism \cite[(6.3.1)]{tsgep} is an isomorphism of 
complexes of a derived category $D^b(K)$ of bounded complexes of $K$-vector spaces, 
the isomorphism \cite[(6.3.2)]{tsgep} is an isomorphism of 
cohomologies. 
Hence the existence of the isomorphism of 
complexes of the derived category $D^b(K)$ 
has not been stated in \cite{hk} and \cite{tsgep}. 
However we can construct the morphism $\Psi_{\pi}$ and 
hence we see that $\Psi_{\pi}$ is an isomorphism by  \cite[(6.3.1), (6.3.2)]{tsgep} (see 
the proofs of \cite[(5.1)]{hk}, \cite[(6.3.1), (6.3.2)]{tsgep} and (\ref{theo:ctff}) below). 
\end{rema}

\begin{prob}\label{prob:nil} 
Let the notations and the assumptions be as in (\ref{theo:hkt}). 
Is the morphism 
\begin{equation*} 
N_{\rm zar}\col 
R\Gam_{{\rm crys}}(Y_{\bul \leq N}/{\cal W}(s))\otimes_{{\cal W}}^LK \lo
R\Gam_{{\rm crys}}(Y_{\bul \leq N}/{\cal W}(s))\otimes_{{\cal W}}^LK 
\tag{6.3.3.1}\label{eqn:rnlkhk}
\end{equation*} 
is nilpotent?  
As to this problem, we obtain the proposition (\ref{prop:dcis}) below. 
\end{prob} 

\begin{prop}\label{prop:dcis}
Let the notations be as in {\rm (\ref{theo:hkt})}.  
Furthermore assume that $Y_{\bul \leq N}/s$ is vertical and 
$Y_{\bul \leq N}$ is saturated. Assume also that $Y_{\bul \leq N}$ is split.   
Then the following hold$:$ 
\par 
$(1)$ The morphism 
\begin{equation*} 
N_{\rm zar}\col 
R\Gam_{{\rm crys}}(Y_{\bul \leq N}/{\cal W}(s))\otimes_{{\cal W}}^LK \lo
R\Gam_{{\rm crys}}(Y_{\bul \leq N}/{\cal W}(s))\otimes_{{\cal W}}^LK 
\tag{6.3.4.1}\label{eqn:rnnyhk}
\end{equation*} 
is nilpotent.  
Hence the isomorphism 
\begin{equation*} 
{\rm exp}(\log (a)N_{\rm zar})\col 
R\Gam_{{\rm crys}}(Y_{\bul \leq N}/{\cal W}(s))\otimes_{{\cal W}}^LK \os{\sim}{\lo} 
R\Gam_{{\rm crys}}(Y_{\bul \leq N}/{\cal W}(s))\otimes_{{\cal W}}^LK 
\tag{6.3.4.2}\label{eqn:rgmyhk}
\end{equation*} 
is well-defined.  
\par 
$(2)$ If $N=0$, then the morphism 
\begin{equation*} 
N_{\rm zar}\col 
R\Gam_{{\rm crys}}(Y_{\bul \leq N}/{\cal W}(s)) \lo
R\Gam_{{\rm crys}}(Y_{\bul \leq N}/{\cal W}(s))
\tag{6.3.4.3}\label{eqn:rnnws}
\end{equation*} 
is nilpotent.
\par 
$(3)$ $($for our memory$)$ 
If $Y_{\bul \leq N}/s$ is a split $N$-truncated SNCL scheme, 
then the morphism 
\begin{equation*} 
N_{\rm zar}\col 
R\Gam_{{\rm crys}}(Y_{\bul \leq N}/{\cal W}(s)) \lo
R\Gam_{{\rm crys}}(Y_{\bul \leq N}/{\cal W}(s))
\tag{6.3.4.4}\label{eqn:rnnwsn}
\end{equation*} 
is nilpotent.
\end{prop}
\begin{proof} 
(1): 
By \cite[(2.4.2.1)]{vi} (see also \cite{yoshi})
there exists a finite extension ${\cal V}(m)$ of ${\cal V}$ 
such that ${\cal Y}\times_SS(m)$  is a strictly semistable log scheme over $S(m)$, 
where $S(m)$ is the log scheme whose underlying scheme is ${\rm Spec}({\cal V}(m))$ 
and whose log structure is canonical.  Let $S'$ be a log scheme whose underlying scheme 
is the spectrum of the integral closure of a finite extension of ${\cal V}$ 
including ${\cal V}(0)\cdots {\cal V}(N)$ and whose log structure is canonical. 
Set ${\cal V}':=\Gam(S',{\cal O}_{S'})$ and 
let $\kap'$ be the residue field of ${\cal V}'$. 
Set $s':=S'\otimes_{{\cal V}'}\kap'$ and $Y'_{\bul \leq N}:=Y_{\bul \leq N}\times_ss'$. 
By the theorem of the existence of admissible immersion (\ref{theo:thenad}),  
we have the complex $A_{{\rm zar},{\mab Q}}(Y'_{\bul \leq N}/{\cal W}(s'))$. 
By (\ref{coro:nancilp}) we see that the morphism 
(\ref{eqn:rnnyhk}) is nilpotent. 
\par 
(2): (2) follows from the proof of (1) and (\ref{coro:nncilp}). 
\par 
\par 
(3): (3) is a special case of (\ref{coro:nncilp}). 
\end{proof}

\begin{rema}\label{rema:cncp} 
Let the notations and the assumptions be as in (\ref{prop:dcis}). 
Then, does the following equality 
\begin{equation*} 
\Psi_{a\pi}=\Psi_{\pi}\circ {\rm exp}(\log (a)N_{\rm zar})
\col R\Gam_{{\rm crys}}(Y_{\bul \leq N}/{\cal W}(s))
\otimes_{{\cal W}}^LK 
\os{\sim}{\lo} R\Gam_{\rm dR}({\mathfrak Y}_{\bul \leq N}/K)
\tag{6.3.5.1}\label{eqn:rghck}
\end{equation*} 
hold for $a\in {\cal V}^*$?
Because I do not know whether the equality (\ref{ali:sirbnq}) holds, 
I cannot answer this question in general. However we shall see that 
the equality (\ref{eqn:rghck}) holds when $a\in [\kap^*]$ only 
under the assumption of (\ref{theo:hkt}); we do not need to make the assumptions 
in (\ref{prop:dcis}). 
See (\ref{coro:can}) for the generalization of the last remark. 
\end{rema}

\begin{lemm}[{\bf cf.~\cite[p.~253]{tst}}]\label{lemm:ptt}
Consider the log scheme ${\rm Spec}^{\log}({\cal W}[t])$ 
whose underlying scheme is ${\rm Spec}({\cal W}[t])$ and 
whose log structure is associated to the morphism 
${\mab N}\owns 1\lom t\in {\cal W}[t]$. 
Consider the closed immersion 
$S\os{\sus}{\lo} {\rm Spec}^{\log}({\cal W}[t])$
defined by $t\lom \pi$.  
Let $D$ be the $p$-adic completion of the log PD-envelope of this closed immersion 
over $({\rm Spec}({\cal W}),p{\cal W},[~])$. 
Then $D$ has a lift $F_D$ of the Frobenius endomorphism of $D\mod p$.  
\end{lemm}
\begin{proof} 
Set  $B:=\Gam(D,{\cal O}_D)$. 
Let $f(t)$ be a minimal polynomial of $\pi$ over ${\cal W}$. 
Then $f(t)=t^e+pg(t)$ $(\deg g< e)$ is an Eisenstein polynomial over ${\cal W}$.
(This is a key point in the following argument, 
which has not been mentioned in references.) 
Then $B={\cal W}[t][f(t)^{[n]}~\vert~n\in {\mab N}]^{\wh{}}$ 
and $B$ has a PD-ideal 
$\langle f(t)B+pB\rangle =\langle t^eB+pB\rangle$. 
Let $\sig \col {\cal W}\lo {\cal W}$ be the Frobenius automorphism of ${\cal W}$. 
Let $\varphi \col {\cal W}[t]\lo {\cal W}[t]$ be an endomorphism 
defined by the following equalities $\varphi(a)=\sig(a)$ and $\varphi(t)=t^p$. 
Since 
$$\varphi(f(t){\cal W}[t]+p{\cal W}[t])=\varphi(t^e{\cal W}[t]+p{\cal W}[t])=
t^{ep}{\cal W}[t]+p{\cal W}[t]\subset t^e{\cal W}[t]+p{\cal W}[t],$$ 
$B$ has a lift $\varphi$ of the Frobenius endomorphism of $B/p$. 
We have proved (\ref{lemm:ptt}). 
\end{proof} 

\begin{rema}\label{rema:hkuso}
Though in \cite[p.~264]{hk}, Hyodo and Kato have used a fact
that $B/p^n$ has a lift of the Frobenius endomorphism, 
I cannot find the claim of this fact and the proof of it in \cite{hk}.  
\end{rema}

\begin{defi}\label{defi:frd} 
We call the morphism $F_D \col D\lo D$ 
constructed in the proof of (\ref{lemm:ptt}) 
the {\it Frobenius endomorphism} of $D$. 
\end{defi}

\begin{lemm}\label{lemm:np} 
Let $T$ be the $p$-adic formal completion of 
the log PD-envelope of a closed immersion 
${\cal W}(s)\os{\sus}{\lo} {\rm Spec}^{\log}({\cal W}[t])$ 
defined by $t\lom 0$ over $({\rm Spec}({\cal W}),p{\cal W},[~])$. 
Let $\wh{S}$ be the $p$-adic completion of $S$. 
Let $F_T\col T\lo T$ be the endomorphism of $T$ induced by the endomorphism 
${\mab N}\owns 1\lom p\in {\mab N}$ and the endomorphism of $\Gam(T,{\cal O}_T)$ 
defined by $t \lom t^p$ and $a\lom \sig(a)$. 
Then the following hold$:$
\par 
$(1)$ There exists a natural morphism 
$H\col T\lo D$ fitting into the following commutative diagram 
\begin{equation*} 
\begin{CD} 
s@>{\subset}>> \wh{S}\\
@V{\bigcap}VV @VV{\bigcap}V \\
T@>{H}>> D. 
\end{CD} 
\end{equation*} 
\par 
$(2)$ Let $m$ be a positive integer such that $p^m\geq e$.  
Then there exists a morphism $G^{\{m\}} \col D\lo T$ such that 
$G^{\{m\}}\circ H=F_T^m$.  
\par 
$(3)$ Set $S_1:=S~{\rm mod}~p$.  
Let the notations be as in $(2)$. 
Then the following diagram is commutative$:$ 
\begin{equation*} 
\begin{CD} 
S_1@>>> s@>>> S_1\\
@V{\bigcap}VV @V{\bigcap}VV @V{\bigcap}VV \\
D@>{G^{\{m\}}}>> T@>{H}>> D, 
\end{CD} 
\tag{6.3.9.1}\label{cd:sss1}
\end{equation*} 
where the two horizontal composite morphisms are the $m$-times iterations of 
the Frobenius endomorphisms. 
\end{lemm}
\begin{proof} 
(1): Note that $\Gam(T,{\cal O}_T)={\cal W}\langle t\rangle$. 
Hence we have the following inclusion 
\begin{align*} 
B=\Gam(D,{\cal O}_D)={\cal W}[t][f(t)^{[n]}~\vert~n\in {\mab N})^{\wh{}}
={\cal W}[t][(t^e)^{[n]}~\vert~n\in {\mab N})^{\wh{}}
\subset {\cal W}\langle t\rangle =\Gam(T,{\cal O}_T)
\end{align*} 
We also have the following commutative diagram 
\begin{equation*} 
\begin{CD} 
{\mab N}@={\mab N}\\
@VVV @VVV\\
\Gam(D,{\cal O}_D)@>{\subset}>> \Gam(T,{\cal O}_T)
\end{CD}
\end{equation*} 
Hence we have the desired morphism $H\col T\lo D$. 
\par 
(2): Since $p^m\geq e$, 
$F_T^*{}^m(t)=t^{p^m}\in t^e{\cal W}[t]+p{\cal W}[t]=f(t){\cal W}[t]+p{\cal W}[t]$. 
Hence the morphism $F_T^*{}^m\col \Gam(T,{\cal O}_T)\lo \Gam(T,{\cal O}_T)$ 
induces the morphism $\Gam(T,{\cal O}_T)\lo \Gam(D,{\cal O}_D)$. 
By noting that $t$ is mapped to $t^{p^m}$ by this morphism and 
that the following diagram 
\begin{equation*} 
\begin{CD} 
{\mab N}@>{p^m\times}>>{\mab N}\\
@VVV @VVV\\
\Gam(T,{\cal O}_T)@>>> \Gam(D,{\cal O}_D)
\end{CD}
\end{equation*} 
is commutative, we have the desired morphism 
$G^{\{m\}}\col D\lo T$ such that $G^{\{m\}}\circ H=F_T^m$. 
\par 
(3): Since $p^m\geq e$, the morphism $G^{\{m\}} \col D\lo T$ 
induces the morphism $S_1\lo s$. 
The commutativity of the right square is obvious.
\end{proof} 

\begin{rema}\label{rema:nm}
In \cite[pp.~263--264]{hk} the proof of a similar statement has not given; 
I think that the statement (\ref{lemm:np}) and the proof of it are necessary. 
\end{rema}

\begin{theo}[{\bf Contravariant functoriality of $\Psi_{\pi}$}]\label{theo:ctff}
Let ${\cal Y}'_{\bul \leq N}/{\cal V}$ be a similar proper log smooth scheme 
to ${\cal Y}_{\bul \leq N}/{\cal V}$. 
Let $Y'_{\bul \leq N}$ $($resp.~${\mathfrak Y}'_{\bul \leq N})$
be the log special fiber $($resp.~the log generic fiber$)$ of 
${\cal Y}'_{\bul \leq N}/S$. 
For a morphism 
$g\col {\cal Y}_{\bul \leq N}/{\cal V}\lo {\cal Y}'_{\bul \leq N}/{\cal V}$ over $S$,  
the following diagram is commutative$:$  
\begin{equation*} 
\begin{CD} 
R\Gam_{{\rm crys}}
(Y_{\bul \leq N}/{\cal W}(s))\otimes_{{\cal W}}^LK 
@>{\Psi_{\pi},\sim}>> R\Gam_{\rm dR} ({\mathfrak Y}_{\bul \leq N}/K)\\
@A{g^*}AA @AA{g^*}A\\
R\Gam_{{\rm crys}}
(Y'_{\bul \leq N}/{\cal W}(s))
\otimes_{{\cal W}}^LK@>{\Psi_{\pi},\sim}>> R\Gam_{\rm dR} 
({\mathfrak Y}'_{\bul \leq N}/K). 
\end{CD} 
\tag{6.3.11.1}\label{eqn:rgkhk}
\end{equation*} 
The following induced isomorphism 
\begin{align*} 
H^q(\Psi_{\pi})\col 
H^q_{\rm crys}(Y_{\bul \leq N}/{\cal W}(s))\otimes_{{\cal W}}K 
\os{\sim}{\lo} 
H^q_{\rm dR} ({\mathfrak Y}_{\bul \leq N}/K)\quad (q\in {\mab N})
\tag{6.3.11.2}\label{eqn:rgkyhk}
\end{align*} 
by $\Psi_{\pi}$ is compatible with the cup products of 
both hand sides on {\rm (\ref{eqn:rgkyhk})}. 
\end{theo} 
\begin{proof} 
We recall the definition of $\Psi_{\pi}$ following the definition of $\Psi_{\pi}$ 
in \cite{hk}. 
\par 
Let $T$ be the $p$-adic formal log scheme in (\ref{lemm:np}). 
Set $T_n:=T\otimes_{\cal W}{\cal W}_n$. 
Then $\{T_n\}_{n=1}^{\inf}$ satisfies the properties 
in the assumptions in {\rm \cite[(4.13)]{hk}}. 
The most nontrivial parts of the proof of (\ref{theo:ctff}) 
are the existence of the following isomorphism 
\begin{align*}
\iota  \col 
({\cal W}\langle t\rangle \otimes^L_{\cal W}Ru_{Y_{\bul \leq N}/{\cal W}(s)*}
({\cal O}_{Y_{\bul \leq N}/{\cal W}(s)}))_{\mab Q}
\os{\sim}{\lo} 
Ru_{Y_{\bul \leq N}/T*}({\cal O}_{Y_{\bul \leq N}/T})_{\mab Q}, 
\tag{6.3.11.3}\label{ali:rlru}
\end{align*} 
the contravariant functoriality of this isomorphism and 
the compatibility of the following induced isomorphism 
\begin{align*}
\iota  \col {\cal W}\langle t\rangle_{\mab Q}\otimes_{K_0}
H^q_{\rm crys}(Y_{\bul \leq N}/{\cal W}(s))_{\mab Q}
\os{\sim}{\lo} 
H^q_{\rm crys}(Y_{\bul \leq N}/T)_{\mab Q}.
\tag{6.3.11.4}\label{ali:rlhru}
\end{align*} 
with the cup products. 
Let 
$$\lam \col Ru_{Y_{\bul \leq N}/T*}({\cal O}_{Y_{\bul \leq N}/T})_{\mab Q}
\lo Ru_{Y_{\bul \leq N}/{\cal W}(s)*}({\cal O}_{Y_{\bul \leq N}/{\cal W}(s)})_{\mab Q}$$  
be the natural morphism. 
As in \cite[(4.13)]{hk}, the isomorphism $\iota$ is the unique morphism 
such that 
\parno 
(1) the composite 
$$\lam \circ \iota \col 
({\cal W}\langle t\rangle \otimes^L_{\cal W}
Ru_{Y_{\bul \leq N}/{\cal W}(s)*}({\cal O}_{Y_{\bul \leq N}/{\cal W}(s)}))_{\mab Q}\lo 
Ru_{Y_{\bul \leq N}/{\cal W}(s)*}({\cal O}_{Y_{\bul \leq N}/{\cal W}(s)})_{\mab Q}$$ 
is the induced morphism by the natural projection ${\cal W}\langle t\rangle  \lo {\cal W}$ 
\parno 
and 
\parno 
(2) $\iota$ commutes with $F\otimes^L F$ and $F$, where the three $F$'s 
on ${\cal W}\langle t\rangle$, $Ru_{Y_{\bul \leq N}/{\cal W}(s)*}({\cal O}_{Y_{\bul \leq N}/{\cal W}(s)})$ 
and $Ru_{Y_{\bul \leq N}/T*}({\cal O}_{Y_{\bul \leq N}/T})$ 
(with tensorization of $\otimes_{\mab Z}^L{\mab Q}$)  
are the induced endomorphisms by the lifts of Frobenius endomorphisms 
of $T_1$ (defined by $t\lom t^p$), 
$Y/s$ and $Y/T_1$, respectively. 
Here let us mean by $?_1$ the ``reduction mod $p$ of $?$. 
\par 
The construction of the morphism $\iota$ and 
the uniqueness of $\iota$ are 
the same as the isomorphism in \cite{hk} by using 
\cite[(4.14), (4.15)]{hk}. Thus we have only to prove that  
$\iota$ is an isomorphism. This follows from the lemma 
(\ref{lemm:fdrw}) below and the proof of \cite[(4.16)]{hk}. 
The contravariant functoriality of $\iota$ follows from the construction of $\iota$ 
and that of the morphism (\ref{ali:zotu}) below. 
Since 
$H^q_{\rm crys}(Y_{\bul \leq N}/T)_{\mab Q}$ has the pull-back of 
the Frobenius endomorphism of $Y_{\bul \leq N}/T_1$ and since 
it induces the Frobenius endomorphism of $Y_{\bul \leq N}/s$, 
we obtain the compatibility of $\iota$ with the products by \cite[(4.4.11)]{tsst} 
as in [loc.~cit., p.~356]: we have only to replace 
$H^m((X,M)/(E,M_E))$ and 
$H^m_{\rm crys}((X,M)$ in [loc.~cit.] 
by $H^q_{\rm crys}(Y_{\bul \leq N}/T)$ and 
$H^q_{{\rm crys}}(Y_{\bul \leq N}/{\cal W}(s))$, respectively. 
\par 
Now the contravariant functoriality of $\Psi_{\pi}$ and the compatibility of 
$\Psi_{\pi}$ with the cup products is easy to prove. 
\par 
Let the notations be as in the proofs of (\ref{lemm:ptt}) and (\ref{lemm:np}). 
Let $\wh{\cal Y}_{\bul \leq N}$ be the $p$-adic completion of ${\cal Y}_{\bul \leq N}$. 
Consider the following commutative diagram 
\begin{equation*} 
\begin{CD} 
{\cal Y}_{\bul \leq N,1}@>>>{\cal Y}_{\bul \leq N,1}\\
@VVV @VVV \\
S_1@>>> S_1\\
@VVV @VVV \\
D@>{F^m_D}>> D, 
\end{CD}
\end{equation*} 
where the two upper horizontal endomorphisms are the 
$p^m$-th power Frobenius endomorphisms of ${\cal Y}_{\bul \leq N,1}$ 
and $S_1$, respectively. 
By (\ref{lemm:np}) this morphism of the diagram 
$({\cal Y}_{\bul \leq N,1}\lo S_1\lo D)\lo ({\cal Y}_{\bul \leq N,1}\lo S_1\lo D)$
factors through 
a morphism $({\cal Y}_1\lo S_1\os{\sus}{\lo} D)\lo (Y\lo s\lo T)$ of diagrams.    
Then we have the following equalities as in \cite[(5.2.1)]{hk}: 
\begin{align*} 
Ru_{\wh{\cal Y}_{\bul \leq N}/D*}({\cal O}_{\wh{\cal Y}_{\bul \leq N}/D})_{\mab Q}
&=Ru_{{\cal Y}_{\bul \leq N,1}/D*}({\cal O}_{{\cal Y}_{\bul \leq N,1}/D})_{\mab Q}
\os{\sim}{\longleftarrow} 
B\otimes^L_{G^{\{m\}}{}^*,{\cal W}\langle t\rangle}
Ru_{Y_{\bul \leq N}/T*}({\cal O}_{Y_{\bul \leq N}/T})_{\mab Q}
\tag{6.3.11.5}\label{ali:rwto}\\
& =B\otimes^L_{\sig^m,{\cal W}}
Ru_{Y_{\bul \leq N}/{\cal W}(s)*}({\cal O}_{Y_{\bul \leq N}/{\cal W}(s)})_{\mab Q}\\
&\os{{\rm id}_R\otimes F^m,\sim}{\lo}
B\otimes^L_{\cal W}
Ru_{Y_{\bul \leq N}/{\cal W}(s)*}({\cal O}_{Y_{\bul \leq N}/{\cal W}(s)})_{\mab Q}. 
\end{align*}  
Here we have used the fact that the Frobenius operator 
\begin{align*} 
F \col 
Ru_{Y_{\bul \leq N}/{\cal W}(s)*}({\cal O}_{Y_{\bul \leq N}/{\cal W}(s)})_{\mab Q}
\lo 
Ru_{Y_{\bul \leq N}/{\cal W}(s)*}({\cal O}_{Y_{\bul \leq N}/{\cal W}(s)})_{\mab Q}
\end{align*} 
is bijective; this follows from \cite[(2.24)]{hk}; however see (\ref{rema:sot}) (2). 
The isomorphism (\ref{ali:rwto}) is contravariantly functorial. 
By taking the tensorization ${\cal V}\otimes^L_B$ of (\ref{ali:rwto}), 
we have the following isomorphism 
\begin{align*} 
{\cal V}\otimes^L_B
Ru_{\wh{\cal Y}_{\bul \leq N}/D*}({\cal O}_{\wh{\cal Y}_{\bul \leq N}/D})_{\mab Q}
\os{\sim}{\lo} 
{\cal V}\otimes^L_{\cal W}
Ru_{Y_{\bul \leq N}/{\cal W}(s)*}({\cal O}_{Y_{\bul \leq N}/{\cal W}(s)})_{\mab Q}. 
\tag{6.3.11.6}\label{ali:rwmto}\\
\end{align*}   
The source of this isomorphism is equal to 
$Ru_{\wh{\cal Y}_{\bul \leq N}/{\cal V}*}
({\cal O}_{\wh{\cal Y}_{\bul \leq N}/{\cal V}})_{\mab Q}$ 
by the base change theorem of log crystalline cohomologies. 
Here we have endowed ${\cal V}$ with the PD-ideal $p{\cal V}$. 
Hence we have the following isomorphism 
\begin{align*} 
Ru_{\wh{\cal Y}_{\bul \leq N}/{\cal V}*}
({\cal O}_{\wh{\cal Y}_{\bul \leq N}/{\cal V}})_{\mab Q}
\os{\sim}{\lo} 
{\cal V}\otimes^L_{\cal W}
Ru_{Y_{\bul \leq N}/{\cal W}(s)*}({\cal O}_{Y_{\bul \leq N}/{\cal W}(s)})_{\mab Q}. 
\tag{6.3.11.7}\label{ali:rwuto}
\end{align*} 
By taking $R\Gam(Y_{\bul \leq N},?)$ of this isomorphism, 
we have the following isomorphism 
\begin{align*} 
R\Gam_{\rm crys}({\wh{\cal Y}_{\bul \leq N}/{\cal V}})_{\mab Q}
\os{\sim}{\lo} 
{\cal V}\otimes^L_{\cal W}
R\Gam_{\rm crys}((Y_{\bul \leq N}/{\cal W}(s))_{\mab Q}. 
\tag{6.3.11.8}\label{ali:rwgtao}
\end{align*}
Because the source of this isomorphism is isomorphic to 
\begin{align*} 
R\Gam_{\rm crys}({\wh{\cal Y}_{\bul \leq N}/{\cal V}})_{\mab Q}
&=R\Gam(\wh{\cal Y}_{\bul \leq N},\Om^{\bul}_{\wh{\cal Y}_{\bul \leq N}}/{\cal V})_{\mab Q}
=R\Gam({\cal Y}_{\bul \leq N},\Om^{\bul}_{{\cal Y}_{\bul \leq N}}/{\cal V})_{\mab Q}\\
&=R\Gam_{\rm dR}({\cal Y}_{\bul \leq N}/{\cal V})_{\mab Q}
=R\Gam_{\rm dR}({\mathfrak Y}_{\bul \leq N}/K)
\tag{6.3.11.9}\label{ali:rwyto}
\end{align*}  
by \cite[(5.1.2)]{ega3}, 
we obtain the desired isomorphism 
\begin{align*} 
R\Gam_{\rm dR}({\mathfrak Y}_{\bul \leq N}/K)
\os{\sim}{\lo} 
K\otimes^L_{\cal W}
R\Gam_{\rm crys}(Y_{\bul \leq N}/{\cal W}(s)).
\tag{6.3.11.10}\label{ali:rwgto}
\end{align*}  
Because the isomorphisms (\ref{ali:rwto}) and 
(\ref{ali:rwyto}) are contravariantly functorial,  
the isomorphism (\ref{ali:rwgto}) 
is contravariantly functorial. 
\end{proof}

The following is a straightforward generalization of \cite[(4.8)]{hk} 
(however see \cite[(7.5), (7.6)]{ndw}): 

\begin{lemm}\label{lemm:fdrw}
Let $N$ be a nonnegative integer or $\infty$. 
Let $Z_{\bul \leq N}/s$ be 
an $N$-truncated simplicial log smooth scheme of Cartier type. 
Let $T$ be an object of the log crystalline site $(s/{\cal W}_n)_{\rm crys}$. 
Let $f_T \col Z_{\bul \leq N}\lo s\os{\subset}{\lo} T$ be the structural morphism. 
Then there exists a contravariantly functorial morphism 
\begin{align*} 
\bigoplus_{i\in {\mab Z}_{\geq 0}}{\cal O}_T
\otimes_{{\cal W}_n}{\cal W}_n\Om^i_{Z_{\bul \leq N}}\lo 
\bigoplus_{i\in {\mab Z}_{\geq 0}}
R^iu_{Z_{\bul \leq N}/T*}({\cal O}_{Z_{\bul \leq N}/T})
\tag{6.3.12.1}\label{ali:zotu}
\end{align*} 
of $f_T^{-1}({\cal O}_T)$-algebras. 
If $\os{\circ}{T}$ is flat over ${\cal W}_n$, then 
{\rm (\ref{ali:zotu})} is an isomorphism. 
\end{lemm}
\begin{proof} 
Let $0\leq m\leq N$ be a fixed integer. 
By the proof of \cite[(4.8)]{hk} 
(however see \cite[(7.5), (7.6)]{ndw}), 
we have a contravariantly functorial morphism 
\begin{align*} 
\bigoplus_{i\in {\mab Z}_{\geq 0}}{\cal O}_T
\otimes_{{\cal W}_n}{\cal W}_n\Om^i_{Z_m}\lo 
\bigoplus_{i\in {\mab Z}_{\geq 0}}R^iu_{Z_m/T*}({\cal O}_{Z_m/T}). 
\tag{6.3.12.2}\label{ali:zomtu}
\end{align*} 
If $\os{\circ}{T}$ is flat over ${\cal W}_n$, then 
{\rm (\ref{ali:zomtu})} is an isomorphism (\cite[(4.8)]{hk}). 
Because the morphism (\ref{ali:zomtu}) is contravariantly functorial, 
we have the morphism (\ref{ali:zotu}). 
This is an isomorphism if $\os{\circ}{T}$ is flat over ${\cal W}_n$. 
\end{proof} 


\begin{rema}\label{rema:istp}  
In (\ref{theo:ctff}) we have also constructed the isomorphism $\Psi_{\pi}$ in 
(\ref{theo:hkt}).  Our proof for the construction of $\Psi_{\pi}$ 
is different from Tsuji's proof in \cite{tsgep}.  
\par 
(2) The proof of \cite[(4.9)]{hk} 
(resp.~\cite[(4.19)]{hk}) is incomplete ([loc.~cit.] 
has been used in the proof of the Hyodo-Kato 
isomorphism in \cite[(5.1)]{hk}. 
In \cite[(7.5)]{ndw} (resp.~\cite[(7.19)]{ndw}) 
we have given a complete proof of \cite[(4.9)]{hk} 
(resp.~\cite[(4.19)]{hk}).  
See (\ref{theo:ccrw}) (2) and (\ref{rema:altn}) (4) for the functoriality of the isomorphism 
in \cite[(4.19)]{hk} and \cite[(7.19)]{ndw}. 
\end{rema}

By believing in the virtual existence of the Hyodo-Kato isomorphism by using 
the element $p$ instead of $\pi$ (cf.~(\ref{prop:compb}) below), 
we can get rid of 
the dependence of the Hyodo-Kato isomorphism on $\pi$ as follows; 
it seems that this is an unexpected result for anybody!:  

\begin{coro-defi}[{\bf Canonical Hyodo-Kato isomorphism}]\label{coro:indne} 
Let the notations be as in {\rm (\ref{theo:hkt})}. 
Let $e$ be the absolute ramification index of ${\cal V}$. 
Express $\pi^e=pa$, where $a\in {\cal V}^*$. 
Then the following isomorphism 
\begin{align*} 
H^q(\Psi):=H^q(\Psi_{p})&
:=H^q(\Psi_{\pi})\circ H^q({\rm exp}(-\log (a)e^{-1}N_{\rm zar})) 
\tag{6.3.14.1}\label{ali:rgrdhk} \\
& \col H^q_{{\rm crys}}(Y_{\bul \leq N}/{\cal W}(s))\otimes_{{\cal W}}K 
\os{\sim}{\lo} H^q_{\rm dR}({\mathfrak Y}_{\bul \leq N}/K)
\end{align*} 
is independent of the choice of the uniformizer $\pi$ of ${\cal V}$. 
This isomorphism is also contravariantly functorial as in 
{\rm (\ref{eqn:rgkhk})}. 
We call $H^q(\Psi)$ the {\it canonical Hyodo-Kato isomorphism}. 
\end{coro-defi} 
\begin{proof} 
Let $\pi'$ be another uniformizer of ${\cal V}$. 
Then there exists an element $a'\in {\cal V}^*$ such that $\pi'^e=pa'$. 
Express $\pi'=\pi b$ $(b\in {\cal V}^*)$. 
Then $b^e=a'/a$. 
By (\ref{theo:hkt}) we have the following equalities: 
\begin{align*} 
H^q(\Psi_{\pi'}) & =H^q(\Psi_{\pi}){\rm exp}(\log (b)N_{\rm zar})
=H^q(\Psi_{\pi}){\rm exp}(\log (b^e)e^{-1}N_{\rm zar})\\
&=H^q(\Psi_{\pi}){\rm exp}(\log (a'/a)e^{-1}N_{\rm zar}). 
\end{align*} 
It is easy to check that 
$\log (cc')=\log c+\log c'$ $(c,c'\in {\cal V}^*)$. 
Hence 
\begin{align*} 
H^q(\Psi_{\pi'}){\rm exp}(-\log (a')e^{-1}N_{\rm zar})
=H^q(\Psi_{\pi}){\rm exp}(-\log (a)e^{-1}N_{\rm zar}).  
\tag{6.3.14.2}\label{ali:ben}\\
\end{align*}  
\parno  
By (\ref{prop:otj}) we have the following commutative diagram: 
\begin{equation*} 
\begin{CD} 
R\Gam_{{\rm crys}}
(Y_{\bul \leq N}/{\cal W}(s)) 
@>{N_{\rm zar}}>> R\Gam_{{\rm crys}}
(Y_{\bul \leq N}/{\cal W}(s)) \\
@A{g^*}AA @AA{g^*}A\\
R\Gam_{{\rm crys}}
(Y'_{\bul \leq N}/{\cal W}(s))
\otimes_{{\cal W}}^LK@>{N_{\rm zar}}>> 
R\Gam_{{\rm crys}}(Y'_{\bul \leq N}/{\cal W}(s)). 
\end{CD} 
\end{equation*} 
Hence we have the following commutative diagram 
\begin{equation*} 
\begin{CD} 
H^q_{{\rm crys}}
(Y_{\bul \leq N}/{\cal W}(s)) 
@>{{\rm exp}(-\log (a)e^{-1}N_{\rm zar}),~\sim}>> 
H^q_{{\rm crys}}(Y_{\bul \leq N}/{\cal W}(s)) \\
@A{g^*}AA @AA{g^*}A\\
H^q_{{\rm crys}}
(Y'_{\bul \leq N}/{\cal W}(s))
\otimes_{{\cal W}}^LK@>{{\rm exp}(-\log (a)e^{-1}N_{\rm zar}),~\sim}>> 
H^q_{{\rm crys}}(Y'_{\bul \leq N}/{\cal W}(s)). 
\end{CD} 
\tag{6.3.14.3}\label{eqn:rnhk}
\end{equation*} 
The contravariant functoriality 
follows from this commutative diagram and (\ref{theo:hkt}). 
\end{proof} 

\begin{rema}\label{rema:un}
To construct 
a canonical Hyodo-Kato isomorphism in (\ref{coro:indne}) (1), (2), 
one can use $pu$ instead of $p$ for any $u\in {\cal W}^*$. 
That is, one can use 
\begin{align*} 
H^q(\Psi_{pu}):=H^q(\Psi_{\pi}){\rm exp}(-\log(b)e^{-1}N_{\rm zar})
=H^q(\Psi_p){\rm exp}(\log(u)e^{-1}N_{\rm zar})
\end{align*} 
instead of $H^q(\Psi_p)$, where $\pi^e=pub$. 
\end{rema}


\begin{prop}\label{prop:cwcup}
The following isomorphism 
\begin{equation*} 
H^q(\Psi) \col H^q_{\rm crys}(Y_{\bul \leq N}/{\cal W}(s))
\otimes_{\cal W}K 
\os{\sim}{\lo} H^q_{\rm dR}({\mathfrak Y}_{\bul \leq N}/K)
\tag{6.3.16.1}\label{eqn:rrhk} 
\end{equation*} 
is compatible with the cup prodcuts on both hand sides.  
\end{prop} 
\begin{proof} 
Because the Hyodo-Kato isomorphism $\Psi_{\pi}$ in (\ref{eqn:rgkyhk}) 
is compatible with the cup products, we have only to prove 
that the automorphism 
\begin{equation*} 
{\rm exp}(aN_{\rm zar}) \col 
\col H^q_{{\rm crys}}(Y_{\bul \leq N}/{\cal W}(s))\otimes_{{\cal W}}K 
\os{\sim}{\lo} H^q_{{\rm crys}}(Y_{\bul \leq N}/{\cal W}(s))\otimes_{{\cal W}}K 
\quad (a\in K)
\tag{6.3.16.2}\label{eqn:rgrmdhk} 
\end{equation*} 
is compatible with the cup product of $H^{\bul}_{{\rm crys}}(Y_{\bul \leq N}/{\cal W}(s))$; 
we have to prove that 
 \begin{equation*} 
{\rm exp}(aN_{\rm zar})(x\cup y)=
{\rm exp}(aN_{\rm zar})(x)\cup {\rm exp}(aN_{\rm zar})(y) 
\tag{6.3.16.3}\label{eqn:rgrahk} 
\end{equation*} 
for $x\in H^{q'}_{\rm crys}(Y_{\bul \leq N}/{\cal W}(s))$ and 
$y\in H^{q''}_{\rm crys}(Y_{\bul \leq N}/{\cal W}(s))$.
This follows from (\ref{eqn:bndd}), the bilinearity of the cup product 
and an elementary calculation.  
\end{proof} 

\begin{prop}[{\bf \cite[(4.4.17)]{tst}}]\label{prop:compb}
Let the notations be as in {\rm (\ref{theo:hkt})}. 
Let ${\cal V}'/{\cal V}$ be a finite extension of 
complete discrete valuation rings of mixed characteristics. 
Let $S'$ be the log scheme ${\rm Spec}({\cal V}')$ endowed with 
the canonical log structure. 
Let $s'$ be the log point of the residue field of ${\cal V}'$ 
and let ${\cal W}'$ be the Witt ring of $\Gam(s',{\cal O}_{s'})$. 
Let $e'$ be the ramification index of ${\cal V}'/{\cal V}$. 
Let $\pi'$ be a uniformizer of ${\cal V}'$. 
Set $b:=\pi'{}^{e'}\pi^{-1}\in {\cal V}'{}^*$.
Set ${\cal Y}'_{\bul \leq N}:={\cal Y}_{\bul \leq N}\times_SS'$ and 
${\mathfrak Y}'_{\bul \leq N}:={\mathfrak Y}_{\bul \leq N}\otimes_KK'$. 
Then the following hold$:$
Then the following diagram is commutative$:$
\begin{equation*} 
\begin{CD} 
H^q_{{\rm crys}}
(Y'_{\bul \leq N}/{\cal W}(s'))\otimes_{{\cal W}'}K' 
@>{H^q(\Psi_{\pi'}),~\sim}>> H^q_{\rm dR}({\mathfrak Y}'_{\bul \leq N}/K')\\
@| @|\\
(H^q_{{\rm crys}}(Y_{\bul \leq N}/{\cal W}(s))\otimes_{\cal W}K)\otimes_KK'
@>{(H^q(\Psi_{\pi})\otimes {\rm id}_{K'}){\rm exp}(\log (b) N_{\rm zar}),~\sim}>> 
H^q_{\rm dR}({\mathfrak Y}_{\bul \leq N}/K)\otimes_KK'. 
\end{CD} 
\tag{6.3.17.1}\label{eqn:rghkahk}
\end{equation*} 
\end{prop}
\begin{proof}  
The following proof is a generalization of that in \cite[(5.5)]{hk}. 
(See also the proof of \cite[(4.4.17)]{tst}.) 
As in the proof of [loc.~cit.], consider the case where 
$b\equiv 1~{\rm mod}~\pi'$ and the case $b\in [\kap'{}^*]$ separately. 
Let 
$$\varphi \col R\Gam_{{\rm crys}}
(Y_{\bul \leq N}/{\cal W}(s))\otimes^L_{\cal W}K \os{\sim}{\lo}  
R\Gam_{{\rm crys}}
(Y_{\bul \leq N}/{\cal W}(s))\otimes^L_{\cal W}K$$ 
and
$$\varphi' \col R\Gam_{{\rm crys}}
(Y'_{\bul \leq N}/{\cal W}(s'))\otimes^L_{{\cal W}'}K' \os{\sim}{\lo}  
R\Gam_{{\rm crys}}
(Y'_{\bul \leq N}/{\cal W}(s'))\otimes^L_{{\cal W}'}K' $$  
be the induced isomorphisms by the Frobenius endomorphisms  
of $Y_{\bul \leq N}$ and $Y'_{\bul \leq N}$, respectively. 
\par 
First consider the case where $b\equiv 1~{\rm mod}~\pi'$. 
Let $n$ be a positive integer and let $m$ be a positive integer 
such that $\pi'{}^m\in p{\cal V}'$. 
Let $t'$ and $t$ be indeterminates. 
Let $t'$ and $t$ be indeterminates. 
Let $T$ be as in (\ref{lemm:np}) and 
let $T'$ be a similar log formal scheme to $T$ by using the indeterminate $t'$. 
Let $r_T\col T'\lo T$ be a morphism of formal log scheme 
defined by the morphism ${\mab N}\owns 1\lom e\in {\mab N}$ and 
the morphism 
\begin{align*} 
{\cal W}\langle t\rangle \owns t \lom t'{}^{e'}\in {\cal W}'\langle t'\rangle. 
\end{align*} 
The morphism $r_T$ induces a morphism 
$r_{{\cal W}(s)}\col {\cal W}(s')\lo {\cal W}(s)$. 
They induce the following morphisms over $r_T$ and $r_{{\cal W}(s)}$, 
respectively: 
\begin{align*} 
r_{Y_{\bul \leq N}/T}\col Y'_{\bul \leq N}\lo Y_{\bul \leq N}
\end{align*} 
and 
\begin{align*} 
r_{Y_{\bul \leq N}/{\cal W}(s)}\col Y'_{\bul \leq N}\lo Y_{\bul \leq N}. 
\end{align*} 
The morphisms $r_{Y_{\bul \leq N}/T}$ and $r_{Y_{\bul \leq N}/{\cal W}(s)}$ induce the following morphisms
\begin{align*} 
r_{Y_{\bul \leq N}/T}^*\col R\Gam_{{\rm crys}}(Y_{\bul \leq N}/T)\otimes^L_{\cal W}K
\lo R\Gam_{{\rm crys}}(Y'_{\bul \leq N}/T')\otimes^L_{{\cal W}'}K'
\end{align*} 
and 
\begin{align*} 
r_{Y_{\bul \leq N}/{\cal W}(s)}^*\col R\Gam_{{\rm crys}}
(Y_{\bul \leq N}/{\cal W}(s))\otimes^L_{\cal W}K
\lo R\Gam_{{\rm crys}}
(Y'_{\bul \leq N}/{\cal W}(s'))\otimes^L_{{\cal W}'}K', 
\end{align*} 
respectively. 
Then we have the following commutative diagrams:
\begin{equation*} 
\begin{CD}
R\Gam_{{\rm crys}}
(Y_{\bul \leq N}/{\cal W}(s))\otimes^L_{\cal W}K
@>{r_{Y_{\bul \leq N}/{\cal W}(s)}^*}>> R\Gam_{{\rm crys}}
(Y'_{\bul \leq N}/{\cal W}(s'))\otimes^L_{{\cal W}'}K'\\
@V{\iota}VV @VV{\iota'}V \\
R\Gam_{{\rm crys}}(Y_{\bul \leq N}/T)\otimes^L_{\cal W}K
@>{r_{Y_{\bul \leq N}/T}^*}>> 
R\Gam_{{\rm crys}}(Y'_{\bul \leq N}/T')\otimes^L_{{\cal W}'}K'\\
@VVV @VVV \\
R\Gam_{{\rm crys}}
(Y_{\bul \leq N}/{\cal W}(s))\otimes^L_{\cal W}K
@>{r_{Y_{\bul \leq N}/{\cal W}(s)}^*}>> R\Gam_{{\rm crys}}
(Y'_{\bul \leq N}/{\cal W}(s'))\otimes^L_{{\cal W}'}K'. 
\end{CD}
\tag{6.3.17.2}\label{cd:rlrwyu}
\end{equation*} 
Consider the following (not necessarily commutative) diagram 
\begin{equation*} 
\begin{CD} 
{\cal W}\langle t\rangle @>{r_T^*}>> {\cal W}'\langle t'\rangle\\ 
@V{v^*_{1}}VV @VV{v^*_{2}}V\\
{\cal V}@>{{\rm inc}}>> {\cal V}', 
\end{CD} 
\end{equation*} 
where the morphism $v_1$ (resp.~$v_2$) is a ${\cal W}$-linear morphism 
(resp.~${\cal W}'$-linear morphism)  
defined by $t\lom \pi$ (resp.~$t'\lom \pi'$) and ${\rm inc}$ is the natural inclusion. 
Then $(v_{2}\circ r)(t)=\pi'{}^{e'}$ and $({\rm inc} \circ v_{1})(t)=\pi$. 
The diagram above induces the following diagram of formal log schemes: 
\begin{equation*} 
\begin{CD} 
T@<{r_T}<< T'\\ 
@A{v_1}AA @AA{v_2}A\\
S@<<< S'.  
\end{CD} 
\end{equation*} 
By this diagram we obtain the following diagram of formal log schemes: 
\begin{equation*} 
\begin{CD} 
Y_{\bul \leq N}@<{r_{Y_{\bul \leq N}/T}}<< Y'_{\bul \leq N}\\ 
@A{w_1}AA @AA{w_2}A\\
\wh{\cal Y}_{\bul \leq N}@<<< \wh{\cal Y}{}'_{\bul \leq N}.  
\end{CD} 
\tag{6.3.17.3}\label{cd:rlroyu}
\end{equation*} 
The morphism $w_i$ $(i=1,2)$ induces the following morphisms 
\begin{align*} 
w_1^*\col R\Gam_{{\rm crys}}(Y_{\bul \leq N}/T)\otimes^L_{\cal W}K\lo 
R\Gam_{\rm dR}({\mathfrak Y}_{\bul \leq N}/K)
\end{align*} 
and 
\begin{align*} 
w_2^*\col R\Gam_{{\rm crys}}(Y'_{\bul \leq N}/T')\otimes^L_{{\cal W}'}K'\lo 
R\Gam_{\rm dR}({\mathfrak Y}'_{\bul \leq N}/K'). 
\end{align*} 
We also have the following (not necessarily commutative) diagram: 
\begin{equation*} 
\begin{CD} 
R\Gam_{{\rm crys}}(Y_{\bul \leq N}/T)\otimes^L_{\cal W}K
@>{r_{Y_{\bul \leq N}/T}^*}>> R\Gam_{{\rm crys}}(Y'_{\bul \leq N}/T')\otimes^L_{{\cal W}'}K' \\
@V{w_1^*}VV @VV{w_2^*}V \\
R\Gam_{\rm dR}({\mathfrak Y}_{\bul \leq N}/K)
@>{``{\rm inc}{\textrm '}{\textrm '}}>>R\Gam_{\rm dR}({\mathfrak Y}'_{\bul \leq N}/K'),  
\end{CD} 
\end{equation*} 
where the lower horizontal morphism is the induced morphism by 
${\rm inc}\col {\cal V}\os{\sus}{\lo} {\cal V}'$. 
\par 
Let 
\begin{align*}
\iota'  \col 
({\cal W}'\langle t'\rangle \otimes^L_{\cal W}Ru_{Y'_{\bul \leq N}/{\cal W}(s')*}
({\cal O}_{Y'_{\bul \leq N}/{\cal W}(s')}))_{\mab Q}
\os{\sim}{\lo} 
Ru_{Y'_{\bul \leq N}/{\cal W}\langle t' \rangle*}({\cal O}_{Y_{\bul \leq N}/{\cal W}\langle t' \rangle})_{\mab Q}, 
\tag{6.3.17.4}\label{ali:rlroyu}
\end{align*} 
be the morphism (\ref{ali:rlru}) for $Y'_{\bul \leq N}/{\cal W}(s')$. 
Let $\star$ be nothing or $\prime$. 
By abuse of notation, we denote the following induced morphisms by $\iota^{\star}$ 
by the same letter: 
\begin{align*}
\iota^{\star}  \col 
({\cal W}^{\star}\langle t^{\star}\rangle \otimes^L_{\cal W}
R\Gam_{\rm crys}(Y^{\star}_{\bul \leq N}/{\cal W}(s^{\star}))_{\mab Q}
\os{\sim}{\lo} 
R\Gam_{\rm crys}(Y^{\star}_{\bul \leq N}/T^{\star})_{\mab Q}
\tag{6.3.17.5}\label{ali:rlrioyu}
\end{align*} 
Then $\Phi_{\pi'}=w_2^*\circ \iota' \circ \varphi'{}^{-m}$. 
On the other hand, $\Phi_{\pi}=w_1^*\circ \iota \circ \varphi^{-m}$. 
Let ${\cal N}_{T{\rm zar}}\col H^q_{\rm crys}(Y^{\star}_{\bul \leq N}/T^{\star})_{\mab Q}
\lo H^q_{\rm crys}(Y^{\star}_{\bul \leq N}/T^{\star})_{\mab Q}$ be 
the monodromy operator. 
Because 
$v_2\circ r_{{\cal W}\langle t\rangle}~{\rm mod}~p = {\rm inc} \circ v_1~{\rm mod}~p$, 
the reduction ${\rm mod}~p$ of the diagram (\ref{cd:rlroyu}) is commutative. 
Hence 
we obtain the following formula 
\begin{align*} 
H^q(w_2^*)\circ H^q(r_{Y_{\bul \leq N}/T}^*)&=\sum_{i=0}(b^{p^m}-1)^{[i]}
H^q(({\rm inc} \circ w^*_1))\circ 
\prod_{j=0}^{i-1}(N_{T{\rm zar}}-j)
\tag{6.3.17.6}\label{ali:rlbyu}\\
&=\sum_{i=0}^{\infty}\log (b^{p^m})^{[i]}H^q(({\rm inc} \circ w^*_1))\circ N^i_{T{\rm zar}}. 
\end{align*} 
as the equality of the morphisms  
\begin{align*} 
H^q_{{\rm crys}}(Y_{\bul \leq N}/T)\otimes_{\cal W}K\lo H^q({\mathfrak Y}'_{\bul \leq N}/K')
\end{align*} 
by using (\ref{ali:sinnq}) as in \cite[p.~266]{hk}. 
Hence 
\begin{align*} 
H^q(\Phi_{\pi'})\circ H^q(r^*_{Y_{\bul \leq N}/{\cal W}(s)})&=
H^q(w^*_2)\circ H^q(\iota') \circ H^q(\varphi'{}^{-m})\circ H^q(r^*_{Y_{\bul \leq N}/{\cal W}(s)})
\tag{6.3.17.7}\label{ali:rzroyu}\\
&=H^q(w^*_2)\circ H^q(\iota') \circ H^q(r^*_{Y_{\bul \leq N}/{\cal W}(s)})\circ H^q(\varphi^{-m})\\
&=
H^q(w^*_2)\circ  H^q(r^*_T)\circ H^q(\iota) \circ H^q(\varphi^{-m})\\
&=\sum_{i=0}^{\infty}\log (b^{p^m})^{[i]}H^q({\rm inc} \circ w_1)\circ N^i_{T{\rm zar}}
\circ  \iota \circ \varphi^{-m}\\
&=
\sum_{i=0}^{\infty}\log (b)^{[i]}H^q({\rm inc} \circ w_1)\circ  \iota \circ N^i_{\rm zar}
\circ \varphi^{-m}\\
&=\Psi_{\pi}{\rm exp}(\log (b) N_{\rm zar}) 
\end{align*} 
as in [loc.~cit.]. 
Here we have used the relations 
\begin{align*} 
H^q(\iota) \circ N_{\rm zar}=N_{T{\rm zar}}\circ H^q(\iota) 
\tag{6.3.17.8}\label{ali:rlrtyu}
\end{align*} 
and 
\begin{align*} 
N_{T{\rm zar}}H^q(\varphi) =pH^q(\varphi) N_{T{\rm zar}}. 
\end{align*} 
(We can prove that the first relation holds by \cite[(4.4.12)]{tst}.)
This implies that we have proved (2) in the case where 
$b\equiv 1~{\rm mod}~\pi'$. 
\par 
Next consider the case $b\in [\kap'{}^*]$. 
Let $g\col {\cal W}'\langle t'\rangle\lo {\cal W}'\langle t'\rangle$ 
be the endomorphism defined by 
$t'\lom b^{-1}t'$. 
Then we have the following commutative diagram 
\begin{equation*} 
\begin{CD} 
{\cal W}'\langle t'\rangle @>{g}>>{\cal W}'\langle t'\rangle \\
@V{\varphi'}VV @VV{\varphi'}V \\
{\cal W}'\langle t'\rangle@>{g}>>{\cal W}'\langle t'\rangle. 
\end{CD}
\tag{6.3.17.9}\label{cd:rlroayu}
\end{equation*} 
Hence we have the following commutative diagram 
\begin{equation*} 
\begin{CD} 
R\Gam_{{\rm crys}}
(Y'_{\bul \leq N}/{\cal W}(s'))\otimes_{{\cal W}'}^LK'
@= R\Gam_{{\rm crys}}(Y'_{\bul \leq N}/{\cal W}(s'))\otimes^L_{{\cal W}'}K' \\
@V{\iota'}VV @VV{\iota'}V \\
R\Gam_{{\rm crys}}
(Y'_{\bul \leq N}/{\cal W}'\langle t'\rangle)\otimes_{{\cal W}'}^LK'
@>{g}>> R\Gam_{{\rm crys}}
(Y'_{\bul \leq N}/{\cal W}'\langle t'\rangle)\otimes_{{\cal W}'}^LK'.  
\end{CD}
\tag{6.3.17.10}\label{cd:rlrtyu}
\end{equation*} 
We also have the following commutative diagram 
\begin{equation*} 
\begin{CD} 
R\Gam_{{\rm crys}}
(Y'_{\bul \leq N}/{\cal W}'\langle t'\rangle)\otimes_{{\cal W}'}^LK'
@>{g}>> R\Gam_{{\rm crys}}
(Y'_{\bul \leq N}/{\cal W}'\langle t'\rangle)\otimes_{{\cal W}'}^LK' \\
@AAA @VV{w_2}V \\
R\Gam_{{\rm crys}}
(Y_{\bul \leq N}/{\cal W}\langle t\rangle)\otimes_{\cal W}^LK
@>{w_1}>> R\Gam_{{\rm dR}}({\mathfrak Y}'_{\bul \leq N}/K')\otimes^L_KK'=
R\Gam_{{\rm dR}}({\mathfrak Y}'_{\bul \leq N}/K').  
\end{CD}
\tag{6.3.17.11}\label{cd:rlranu}
\end{equation*} 
By (\ref{cd:rlrtyu}) and (\ref{cd:rlranu}) we see that 
the diagram (\ref{eqn:rghkahk}) is commutative for $b\in[\kap'{}^*]$. 
\par 
We complete the proof of this proposition. 
\end{proof}

\begin{rema}\label{rema:nzdr}
(1) In \cite[(4.4.17)]{tst} Tsuji has stated (\ref{prop:compb}) 
only in the constant simplicial case. 
\par 
(2) I do not know whether the following diagram 
\begin{equation*} 
\begin{CD} 
R\Gam_{{\rm crys}}
(Y'_{\bul \leq N}/{\cal W}(s'))\otimes_{{\cal W}'}K' 
@>{\Psi_{\pi'},~\sim}>> R\Gam_{\rm dR}({\mathfrak Y}'_{\bul \leq N}/K')\\
@| @|\\
(R\Gam_{{\rm crys}}(Y_{\bul \leq N}/{\cal W}(s))\otimes_{\cal W}K)\otimes_KK'
@>{(\Psi_{\pi}\otimes {\rm id}_{K'}){\rm exp}(\log (b) N_{\rm zar}),~\sim}>> 
R\Gam_{\rm dR}({\mathfrak Y}_{\bul \leq N}/K)\otimes_KK'. 
\end{CD} 
\tag{6.3.18.1}\label{eqn:rgkahk}
\end{equation*}  
is commutative for $b\not\in [\kap'{}^*]$. 
If $b\in [\kap'{}^*]$, this diagram is commutative by 
the proof of (\ref{prop:compb}). 
\end{rema} 

\begin{coro}[for our memory]\label{coro:can}
Let the notations be as in {\rm (\ref{prop:compb})}. 
If $b\in [\kap'{}^*]$, then the following diagram is commutative$:$
\begin{equation*} 
\begin{CD} 
R\Gam_{{\rm crys}}(Y'_{\bul \leq N}/{\cal W}(s'))\otimes_{{\cal W}'}^LK' 
@>{\Psi_{\pi'},~\sim}>> R\Gam_{\rm dR}({\mathfrak Y}'_{\bul \leq N}/K')\\
@| @|\\
(R\Gam_{{\rm crys}}(Y_{\bul \leq N}/{\cal W}(s))\otimes_{\cal W}^LK)\otimes^L_KK'
@>{\Psi_{\pi}\otimes {\rm id}_{K'},~\sim}>> 
R\Gam_{\rm dR}({\mathfrak Y}_{\bul \leq N}/K)\otimes^L_KK'. 
\end{CD} 
\tag{6.3.19.1}\label{eqn:rgkaahk}
\end{equation*}  
\end{coro} 
\begin{proof} 
As remarked in (\ref{rema:nzdr}) (2), we have proved 
this corollary in the proof of (\ref{prop:compb}). 
\end{proof} 

\begin{coro}[{\bf Functoriality of the Hyodo-Kato isomorphism}]\label{coro:nchk}
Let the notations be as in {\rm (\ref{prop:compb})}. 
Let $g\col {\cal Z}_{\bul \leq N}\lo {\cal Y}_{\bul \leq N}$ be a morphism of 
proper log smooth $N$-truncated simplicial log schemes over $S'\lo S$. 
Let $Z_{\bul \leq N}$ be the log special fiber of ${\cal Z}_{\bul \leq N}$.  
Assume that $Z_{\bul \leq N}/s'$ is of Cartier type. 
Then the following diagram is commutative$:$  
\begin{equation*} 
\begin{CD} 
H^q_{\rm crys}(Z_{\bul \leq N}/{\cal W}(s'))\otimes_{{\cal W}'}K' 
@>{H^q(\Psi_{\pi'}),~\sim}>> H^q_{\rm dR} ({\mathfrak Z}_{\bul \leq N}/K')\\
@AAA @AAA\\
(H^q_{{\rm crys}}(Y_{\bul \leq N}/{\cal W}(s))\otimes_{{\cal W}}^LK)\otimes^L_KK'
@>{H^q(\Psi_{\pi}\otimes {\rm id}_{K'}){\rm exp}(\log (b) N_{\rm zar}),~\sim}>> 
H^q_{\rm dR}({\mathfrak Y}_{\bul \leq N}/K)\otimes_KK'. 
\end{CD} 
\tag{6.3.20.1}\label{eqn:ragkkhk}
\end{equation*} 
\end{coro} 
\begin{proof} 
This follows from (\ref{eqn:rgkhk}) and (\ref{eqn:rghkahk}). 
\end{proof}

\begin{lemm}\label{lemm:nne} 
Let the notations be as in {\rm (\ref{prop:compb})}. 
Then 
\begin{align*} 
N'_{\rm zar}=e'N_{\rm zar} \col & R\Gam_{{\rm crys}}
(Y'_{\bul \leq N}/{\cal W}(s'))=R\Gam_{{\rm crys}}
(Y_{\bul \leq N}/{\cal W}(s))\otimes^L_{\cal W}{\cal W}(s')
\tag{6.3.21.1}\label{cd:vswm}\\ 
& \lo R\Gam_{{\rm crys}}
(Y'_{\bul \leq N}/{\cal W}(s'))
=R\Gam_{{\rm crys}}
(Y_{\bul \leq N}/{\cal W}(s))\otimes^L_{\cal W}{\cal W}(s').
\end{align*} 
\end{lemm}
\begin{proof} 
Let $s'$ be the log point of  the residue field of ${\cal V}'$. 
Let $t$ be $s$ or $s'$.  Set $\kap_t:=\Gam(t,{\cal O}_t)$. 
Let $\varpi_t$ be the global section of $M_t$ 
whose image in $M_t/\kap^*_t$ is the generator.   
Let $v' \col {\cal V}'\setminus \{0\}\lo{\mab N}$ and 
$v\col {\cal V}\setminus \{0\}\lo{\mab N}$ be 
the normalized valuation of ${\cal V}'$ and ${\cal V}$, respectively. 
Because the following diagram 
\begin{equation*} 
\begin{CD} 
({\cal V}'\setminus \{0\})/{\cal V}'{}^*@>{v',\sim}>>{\mab N}\\
@A{\bigcup}AA @AA{e'\times}A\\
({\cal V}\setminus \{0\})/{\cal V}^*@>{v,\sim}>>{\mab N}
\end{CD}
\tag{6.3.21.2}\label{cd:vsim}
\end{equation*} 
is commutative, (\ref{lemm:nne}) follows from (\ref{eqn:ynmumss}). 
\end{proof} 

\parno 
Using (\ref{prop:compb}) and (\ref{lemm:nne}), 
we can generalize (\ref{coro:indne}) as follows: 

\begin{coro}[{\bf Compatibility of $H^q(\Psi)$ 
with a base extension of ${\cal V}$'s}]\label{coro:indene} 
Let the notations be as in {\rm (\ref{coro:indne})} and 
{\rm (\ref{prop:compb})}. 
Let 
\begin{equation*} 
H^q(\Psi')\col H^q_{{\rm crys}}(Y'_{\bul \leq N}/{\cal W}'(s'))\otimes_{{\cal W}'}K' 
\os{\sim}{\lo} H^q_{\rm dR}({\mathfrak Y}'_{\bul \leq N}/K')
\end{equation*} 
be the isomorphism {\rm (\ref{ali:rgrdhk})} for ${\cal Y}'_{\bul \leq N}/S'$.  
Then the following diagram 
\begin{equation*} 
\begin{CD} 
H^q_{{\rm crys}}
(Y'_{\bul \leq N}/{\cal W}(s'))\otimes_{{\cal W}'}K' 
@>{H^q(\Psi'),~\sim}>> H^q_{\rm dR}({\mathfrak Y}'_{\bul \leq N}/K')\\
@| @|\\
(H^q_{{\rm crys}}(Y_{\bul \leq N}/{\cal W}(s))\otimes_{{\cal W}}K)\otimes_KK'
@>{H^q(\Psi) \otimes {\rm id}_{K'},~\sim}>> 
H^q_{\rm dR}({\mathfrak Y}_{\bul \leq N}/K)\otimes_KK'. 
\end{CD} 
\tag{6.3.22.1}\label{eqn:rgkmhk}
\end{equation*} 
is commutative. 
\end{coro} 
\begin{proof} 
Obviously $ee'$ is the absolute ramification index of ${\cal V}'$ 
and we have the following equalities: 
$(\pi')^{ee'}=(\pi b)^e=pab^e$. 
Let $N'_{\rm zar} \col H^q_{\rm crys}(Y'_{\bul \leq N}/{\cal W}(s'))\lo 
H^q_{\rm crys}(Y'_{\bul \leq N}/{\cal W}(s'))$ be the monodromy operator for 
$Y'_{\bul \leq N}/{\cal W}(s')$. 
By (\ref{lemm:nne}) we have the following equalities: 
\begin{align*} 
H^q(\Psi')&=H^q(\Psi_{\pi'}){\rm exp}(-\log (ab^e)(ee')^{-1}N'_{\rm zar})
\tag{6.3.22.2}\label{ali:bebn}\\
&=H^q(\Psi_{\pi'}){\rm exp}(-\log (b)e'{}^{-1}N'_{\rm zar})
{\rm exp}(-\log (a)(ee')^{-1}N'_{\rm zar})\\
&=H^q(\Psi_{\pi}\otimes {\rm id}_{K'}){\rm exp}(-\log (a)e^{-1}N_{\rm zar}\otimes {\rm id}_{K'})\\
&=H^q(\Psi) \otimes {\rm id}_{K'}.  
\end{align*} 
\end{proof}

\begin{coro}[{\bf Functoriality of the canonical Hyodo-Kato isomorphism}]\label{coro:chk}
Let the notations be as in {\rm (\ref{coro:nchk})}. 
Then the following diagram is commutative$:$  
\begin{equation*} 
\begin{CD} 
H^q_{{\rm crys}}
(Z_{\bul \leq N}/{\cal W}(s'))\otimes_{{\cal W}'}K' 
@>{H^q(\Psi),~\sim}>> H^q_{\rm dR} ({\mathfrak Z}_{\bul \leq N}/K')\\
@AAA @AAA\\
(H^q_{{\rm crys}}(Y_{\bul \leq N}/{\cal W}(s))\otimes_{{\cal W}}K)\otimes_KK'
@>{H^q(\Psi) \otimes {\rm id}_{K'},~\sim}>> 
H^q_{\rm dR}({\mathfrak Y}_{\bul \leq N}/K)\otimes_KK'. 
\end{CD} 
\tag{6.3.23.1}\label{eqn:rgkkhk}
\end{equation*} 
\end{coro} 
\begin{proof} 
This follows from (\ref{coro:indne})  and (\ref{coro:indene}). 
\end{proof}

\par 
Henceforth we denote $S$, $s$ and $\kap$ in 
the beginning of this section by 
$S_{-1}$, $s_{-1}$  and $\kap_{-1}$, respectively. 
Let ${\mathfrak X}$ be a proper scheme over $K$. 
Let $q$ be a nonnegative integer. 
Let $N$ be an integer which is greater than or equal to 
$2^{-1}(q+1)(q+2)$. 
Let ${\cal V}_N$ be a finite extension of ${\cal V}$. 
Let $\kap$ be the residue field of ${\cal V}_N$. 
Set ${\cal W}:={\cal W}(\kap)$.   
Let $e:=e_N$ be the ramification index of ${\cal V}_N/{\cal V}$. 
Let $S$ be the log scheme ${\rm Spec}({\cal V}_N)$ 
endowed with the canonical log structure. 
Let $s$ be the log special fiber of $S$. 
Let ${\cal X}$ be a proper flat model of ${\mathfrak X}$ over ${\cal V}$. 
Let ${\cal X}_{\bul \leq N}$ be a proper log smooth scheme over $S$ 
such that ${\cal X}_{\bul \leq N}$ is 
a gs $N$-truncated simplicial generically proper hypercovering of 
${\cal X}_{{\cal V}_N}:={\cal X}\otimes_{\cal V}{\cal V}_N$. 
Let $X_{\bul \leq N}$ (resp.~${\mathfrak X}_{\bul \leq N}$) 
be the log special fiber (resp.~the log generic fiber) of ${\cal X}_{\bul \leq N}$. 
In the following the log structure of 
${\mathfrak X}_{\bul \leq N}$ is assumed to be trivial. 
Set ${\mathfrak X}_{K_N}:={\mathfrak X}\otimes_{K}K_N$. 
Then, by (\ref{eqn:olyzk})  and (\ref{prop:hnn}), 
we have the following composite isomorphism 
\begin{align*} 
H^q_{{\rm crys}}(X_{\bul \leq N}/{\cal W}(s))
\otimes_{\cal W}K_N
&\os{H^q(\Psi),~\sim}{\lo} 
H^q_{\rm dR}({\mathfrak X}_{\bul \leq N}/K_N)
=
H^q({\mathfrak X}_{\bul \leq N},\Om^{\bul}_{{\mathfrak X}_{\bul \leq N}/K_N}) 
\tag{6.3.23.2}\label{eqn:fkci}\\
&\os{\sim}{\lo} H^q_{\rm inf}({\mathfrak X}_{\bul \leq N}/K_N)
\os{\sim}{\longleftarrow}
H^q_{\rm inf}({\mathfrak X}_{K_N}/K_N). 
\end{align*} 
We denote this isomorphism by $\Psi$ by abuse of notation. 
Hence we have the following isomorphism 
\begin{align*} 
H^q(\Psi) \col 
H^q_{{\rm crys}}(X_{\bul \leq N}/{\cal W}(s))\otimes_{{\cal W}}K_N
\os{\sim}{\lo} 
H^q_{\rm inf}({\mathfrak X}_{K_N}/K_N). 
\tag{6.3.23.3}\label{eqn:fkysci}
\end{align*} 
Because the left hand side of (\ref{eqn:fkysci}) has the monodromy operator 
((\ref{defi:mdop})), 
we have the monodromy operator  
$N_{K_N/K}({\cal X}_{\bul \leq N}/{\cal X}/S/S_{-1})$ 
on $H^q_{\rm inf}({\mathfrak X}_{K_N}/K_N)$:  
\begin{align*}
N_{K_N/K}({\cal X}_{\bul \leq N}/{\cal X}/S/S_{-1})
\col H^q_{\rm inf}({\mathfrak X}_{K_N}/K_N)
\lo H^q_{\rm inf}({\mathfrak X}_{K_N}/K_N). 
\tag{6.3.23.4}\label{ali:xsaknif}
\end{align*} 
Here note that the monodromy operator in (\ref{eqn:fkysci}) is not equal 
$N_{\rm zar}$ in (\ref{theo:hkt}): it is equal to $e^{-1}N_{\rm zar}$. 
The morphism 
$N_{K_N/K}({\cal X}_{\bul \leq N}/{\cal X}/S/S_{-1})$ 
is independent of the choice of a uniformizer $\pi$ of $K_N$ because 
the isomorphism $\Psi$ in (\ref{eqn:fkysci}) is independent of $\pi$. 
(If we consider $\Psi_{\pi}$ instead of $\Psi$,  
we can define a similar morphism to 
$N_{K_N/K}({\cal X}_{\bul \leq N}/{\cal X}/S/S_{-1})$ 
and it is also independent of the choice of $\pi$ since 
$\Psi_{a\pi}N_{\rm zar}\Psi_{a\pi}^{-1}=
\Psi_{\pi}N_{\rm zar}\Psi_{\pi}^{-1}$ ((\ref{theo:hkt})).) 
In (\ref{theo:wdnmo}) below we prove that 
$N_{K_N/K}({\cal X}_{\bul \leq N}/{\cal X}/S/S_{-1})$ depends only on 
${\mathfrak X}_{K_N}/K_N$ and $K_N$. 
The morphism (\ref{ali:xsaknif}) induces the following morphism 
\begin{align*}
N_{\ol{K}/K}({\cal X}_{\bul \leq N}/{\cal X}/S/S_{-1})
\col H^q_{\rm inf}({\mathfrak X}_{\ol{K}}/\ol{K})
\lo H^q_{\rm inf}({\mathfrak X}_{\ol{K}}/\ol{K}). 
\tag{6.3.23.5}\label{ali:xsknoif}
\end{align*} 
In (\ref{theo:wdnmo}) below we prove that 
$N_{\ol{K}/K}({\cal X}_{\bul \leq N}/{\cal X}/S/S_{-1})$ depends only on 
${\mathfrak X}_{\ol{K}}/\ol{K}$ and $K$. 

\par
The following is a generalization 
of \cite[Definition A.~12 (i)]{miepl}, 
though our definition of a monodromy operator is 
different from that in [loc.~cit.]. 

\begin{prop}\label{prop:gemid}  
Let the notations be as above.  
Let ${\cal V}'/{\cal V}$ be a finite extension of 
complete discrete valuation rings of mixed characteristics. 
Set $K':={\rm Frac}({\cal V}')$. 
Let ${\mathfrak Y}$ be a proper scheme over $K'$. 
Let ${\cal Y}$ be a proper flat model of ${\mathfrak Y}$ over ${\cal V}'$. 
Let ${\mathfrak g}\col {\mathfrak Y}\lo {\mathfrak X}$ 
be a morphism over ${\rm Spec}(K')\lo {\rm Spec}(K)$. 
Let $g \col {\cal Y}\lo {\cal X}$ be a morphism 
over ${\rm Spec}({\cal V}')\lo {\rm Spec}({\cal V})$ which induces the morphism 
${\mathfrak g}$. 
Let $S'_{-1}$ be the log scheme ${\rm Spec}({\cal V}')$ with 
the canonical log structure. 
Let $s'_{-1}$ be the log point of the residue field of ${\cal V}'$. 
Let ${\cal V}'_N$ and ${\cal Y}_{\bul \leq N}/S'$ 
be a similar complete discrete valuation ring to ${\cal V}_N$
and a similar $N$-truncated simplicial log scheme to ${\cal X}_{\bul \leq N}/S$, 
respectively, fitting into the following commutative diagram$:$  
\begin{equation*} 
\begin{CD} 
{\cal Y}_{\bul \leq N} 
@>{g}>> {\cal X}_{\bul \leq N}\\
@VVV @VVV \\
{\cal Y}\times_{\os{\circ}{S}{}'_{-1}}\os{\circ}{S}{}' 
@>>> {\cal X}\times_{\os{\circ}{S}_{-1}}\os{\circ}{S}\\
@VVV @VVV \\
\os{\circ}{S}{}'@>>> \os{\circ}{S}
\end{CD}
\tag{6.3.24.1}\label{ali:xxpss}
\end{equation*} 
such that ${\cal Y}_{\bul \leq N}$ is a gs $N$-truncated simplicial 
generically proper hypercovering  of ${\mathfrak Y}_{{\cal V}'_N}$. 
Here we assume that the log structure of the generic fiber of 
${\cal Y}_{\bul \leq N}$ is trivial. 
Set $K'_N:={\rm Frac}({\cal V}'_N)$. 
Let $e_{{\cal V}'/{\cal V}}$ be the ramification index of the extension 
${\cal V}'/{\cal V}$. 
Then the following diagram is commutative$:$ 
\begin{equation*} 
\begin{CD} 
H^q_{\rm inf}({\mathfrak Y}_{K'_N}/K'_N) 
@>{e_{{\cal V}'/{\cal V}}^{-1}\cdot 
N_{K'_N/K'}({\cal Y}_{\bul \leq N}/{\cal Y}/S'/S'_{-1})}>>
H^q_{\rm inf}({\mathfrak Y}_{K'_N}/K'_N)\\
@A{{\mathfrak g}^*}AA @AA{{\mathfrak g}^*}A \\
H^q_{\rm inf}({\mathfrak X}_{K_N}/K_N)
@>{N_{K_N/K}({\cal X}_{\bul \leq N}/{\cal X}/S/S_{-1})}>>
H^q_{\rm inf}({\mathfrak X}_{K_N}/K_N). 
\end{CD}
\tag{6.3.24.2}\label{cd:lknl}
\end{equation*} 
\end{prop}
\begin{proof} 
It suffices to prove that the following diagram is commutative$:$ 
\begin{equation*} 
\begin{CD} 
H^q_{\rm inf}({\mathfrak Y}_{K'_N}/K'_N) 
@>{N_{K'_N/K'}({\cal Y}_{\bul \leq N}/{\cal Y}/S'/S'_{-1})}>>
H^q_{\rm inf}({\mathfrak Y}_{K'_N}/K'_N)\\
@A{{\mathfrak g}^*}AA @AA{{\mathfrak g}^*}A \\
H^q_{\rm inf}({\mathfrak X}_{K_N}/K_N)\otimes_{K_N}K'_N
@>{N_{K_N/K}({\cal X}_{\bul \leq N}/{\cal X}/S/S_{-1})\otimes {\rm id}_{K'_N}}>>
H^q_{\rm inf}({\mathfrak X}_{K_N}/K_N)\otimes_{K_N}K'_N. 
\end{CD}
\tag{6.3.24.3}\label{cd:lknxl}
\end{equation*} 
Let $Y_{\bul \leq N}$ be the log special fiber of ${\cal Y}_{\bul  \leq N}$. 
Let $X_{\bul \leq N,\bul} \os{\sus}{\lo} {\cal P}^{\rm ex}_{\bul \leq N,\bul}$ 
and 
$Y_{\bul \leq N, \bul} \os{\sus}{\lo} {\cal Q}^{\rm ex}_{\bul \leq N,\bul}$ 
be the admissible immersions constructed in (\ref{theo:thenad}) 
over ${\cal W}(s)$ and ${\cal W}(s')$, respectively, 
fitting into the following commutative diagram 
\begin{equation*} 
\begin{CD} 
Y_{\bul \leq N,\bul}@>>> {\cal Q}^{\rm ex}_{\bul \leq N,\bul}@>>>{\cal W}(s')\\
@VVV @VVV @VVV\\
X_{\bul \leq N,\bul}@>>> {\cal P}^{\rm ex}_{\bul \leq N,\bul}@>>>{\cal W}(s). 
\end{CD}
\tag{6.3.24.4}\label{eqn:gsmnflxd}
\end{equation*} 
Let ${\mathfrak D}_{\bul \leq N,\bul}$ and 
${\mathfrak E}_{\bul \leq N,\bul}$ be the log PD-envelopes of these immersions 
over $({\cal W}(s),p{\cal W}(\kap_s),[~])$ and 
$({\cal W}(s'),p{\cal W}(\kap_{s'}),[~])$, respectively. 
Let $h\col {\cal Q}^{\rm ex}_{\bul \leq N,\bul}\lo {\cal P}^{\rm ex}_{\bul \leq N,\bul}$ 
be the middle vertical morphism above. 
By abuse of notation, we denote the induced morphism 
${\mathfrak E}_{\bul \leq N,\bul}\lo 
{\mathfrak D}_{\bul \leq N,\bul}$ by $h$ again. 
Let $s'$ be the log point of  the residue field of ${\cal V}'_N$. 
Let $t$ be $s$ or $s'$.  Set $\kap_t:=\Gam(t,{\cal O}_t)$. 
Let $\varpi_t$ be the global section of $M_t$ 
whose image in $M_t/\kap^*_t$ is the generator.   
Let $v' \col {\cal V}'_N\setminus \{0\}\lo{\mab N}$ and 
$v\col {\cal V}_N\setminus \{0\}\lo{\mab N}$ be 
the normalized valuation of ${\cal V}'_N$ and ${\cal V}_N$, respectively. 
Let $e_{S'/S}$ be the ramification index of ${\cal V}'_N/{\cal V}_N$.  
Let $e_{S'}$ and $e_S$ be the ramification indexes of 
${\cal V}'_N$ and ${\cal V}_N$ over ${\cal V}'$ and ${\cal V}$, respectively. 
By (\ref{cd:vsim}) 
$d\log \varpi_{s}= e_{S'/S}d\log \varpi_{s'}$ in $\Om^1_{s'/\os{\circ}{s}{}'}$. 
Hence $e_Sd\log \varpi_{s}=e_{{\cal V}'/{\cal V}}(e_{S'}d\log \varpi_{s'})$ 
and 
$e_S\theta_{{\cal P}_{\bul \leq N,\bul/{\cal W}(s)}}=
e_{{\cal V}'/{\cal V}}e_{S'}\theta_{{\cal Q}_{\bul \leq N,\bul/{\cal W}(s')}}$ 
in ${\cal O}_{{\cal Q}_{\bul \leq N,\bul}}$ 
and the following diagram is commutative: 
\begin{equation*} 
\begin{CD}
0@>>> \Om^{\bul}_{{\mathfrak E}_{\bul \leq N,\bul}/{\cal W}(s'),[~]}[-1] 
@>{e_{{\cal V}'/{\cal V}}e_{S'}\theta_{{\cal Q}_{\bul \leq N,\bul/{\cal W}(s')}}}>>
\Om^{\bul}_{{{\mathfrak D}'{}}_{\bul \leq N,\bul}/{\cal W}(\os{\circ}{s}{}'),[~]}\\
@. @A{h^*}AA @A{h^*}AA\\
0@>>> h^*(\Om^{\bul}_{{{\mathfrak D}}_{\bul \leq N,\bul}/{\cal W}(s),[~]})[-1] 
@>{e_{S}\theta_{{\cal P}_{\bul \leq N,\bul/{\cal W}(s)}}}>>
h^*(\Om^{\bul}_{{\mathfrak D}_{\bul \leq N,\bul}/{\cal W}(\os{\circ}{s}),[~]})
\end{CD}
\tag{6.3.24.5}\label{eqn:gsmnfalxd}
\end{equation*} 
\begin{equation*} 
\begin{CD}
@>>>  
\Om^{\bul}_{{\mathfrak E}_{\bul \leq N,\bul}/{\cal W}(s'),[~]} \lo 0 \\
@. @A{h^*}AA \\
@>>>  
h^*(\Om^{\bul}_{{\mathfrak D}_{\bul \leq N,\bul}/{\cal W}(s),[~]}) \lo 0. 
\end{CD}
\end{equation*} 
Since $ee_{{\cal V}'/{\cal V}}$ is the absolute ramification index of ${\cal V}'$, 
(\ref{prop:gemid}) follows from (\ref{eqn:rgkahk}) and the definition of 
$N_{K_N/K}({\cal X}_{\bul \leq N}/{\cal X}/S/S_{-1})$ and 
$N_{K'_{S'}/K'}({\cal Y}_{\bul \leq N}/{\cal Y}/S'/S'_{-1})$ 
as in the proof of (\ref{lemm:nne}).  
\end{proof} 

\begin{rema}\label{rema:npp}
In \cite[Definition A.~12 (i)]{miepl} there is no proof for it. 
It is impossible to judge whether a precise proof for  
\cite[Definition A.~12 (i)]{miepl} has been given in [loc.~cit.]: 
I do not know whether (\ref{prop:compb}) for the constant simplicial case 
has been used in [loc.~cit.].  
\end{rema}

\begin{theo}\label{theo:wdnmo}  
The monodromy operator 
$N_{K_N/K}({\cal X}_{\bul \leq N}/{\cal X}/S/S_{-1})$ 
in {\rm (\ref{ali:xsaknif})} is independent of the choices of 
${\cal X}$, 
$N\geq 2^{-1}(q+1)(q+2)$ and 
the gs $N$-truncated simplicial generically proper hypercovering  
${\cal X}_{\bul \leq N}$ of ${\cal X}_{{\cal V}_N}$. 
It depends only on ${\mathfrak X}_{K_N}/K_N$ and $K$. 
The monodromy operator also 
$N_{\ol{K}/K}({\cal X}_{\bul \leq N}/{\cal X}/S/S_{-1})$ in 
{\rm (\ref{ali:xsknoif})} depends only on 
${\mathfrak X}_{\ol{K}}/\ol{K}$ and $K$.  
\end{theo} 
\begin{proof} 
Assume that we are given another proper flat model 
${\cal X}'$ of ${\mathfrak X}$, $N'\geq 2^{-1}(q+1)(q+2)$ and 
another gs $N'$-truncated generically proper hypercovering 
of ${\cal X}'_{\bul \leq N'}$ of 
${\cal X}_{{\cal V}_{N'}}$. 
We may assume that $N\leq N'$ and 
we may consider ${\cal X}'_{\bul \leq N}$ instead of 
${\cal X}'_{\bul \leq N'}$  by \cite[(2.2)]{nh3}. 
By (\ref{lemm:cmd}) 
we may assume that there exists another proper flat model ${\cal X}''$ over 
${\cal V}$ of ${\mathfrak X}$ with morphisms 
${\cal X}''\lo {\cal X}$ and ${\cal X}''\lo {\cal X}'$ over ${\cal V}$ 
which induce ${\rm id}_{\mathfrak X}$.  
By (\ref{prop:coif}) there exists a gs $N$-truncated simplicial
generically proper hypercovering 
of ${\cal X}''_{\bul \leq N}$ of ${\cal X}''_{{\cal V}'_N}$ over 
a finite extension ${\cal V}'_N$ of ${\cal V}_N$  
fitting into the following commutative diagram  
\begin{equation*} 
\begin{CD} 
{\cal X}_{\bul \leq N}@<<<
{\cal X}''_{\bul \leq N}\times_SS'@>>>{\cal X}'_{\bul \leq N}\\ 
@VVV @VVV @VVV\\
{\cal X}@<<< {\cal X}_{{\cal V}'_N} @>>> {\cal X}'. 
\end{CD} 
\tag{6.3.26.1}\label{cd:lynln}
\end{equation*}   
Here $S'$ is the log scheme ${\rm Spec}({\cal V}'_N)$ with the canonical log structure. 
Set $K'_N:={\rm Frac}({\cal V}'_N)$. 
It suffices to prove that  
$$N_{K_N/K}({\cal X}_{\bul \leq N}/{\cal X}_{{\cal V}_N}/S/S_{-1})
\otimes_{K_N}K'_N=
N_{K_N/K}({\cal X}'_{\bul \leq N}/{\cal X}_{{\cal V}_N}/S/S_{-1})
\otimes_{K_N}K'_N.$$ 
Let $S'$ be the log scheme whose underlying scheme is ${\rm Spec}({\cal V}'_N)$ 
and whose log structure is canonical. 
By (\ref{cd:lknl}) we have the following commutative diagram 
\begin{equation*} 
\begin{CD} 
H^q_{\rm inf}({\mathfrak X}_{K_N}/K_N) @>>>
H^q_{\rm inf}({\mathfrak X}''_{K'_N}/K'_N)@<<<
H^q_{\rm inf}({\mathfrak X}'_{K_N}/K_N)\\
@VV{N_{K_N/K}({\cal X}_{\bul \leq N}/{\cal X}/S/S_{-1})}V
@VV{N_{K_N/K}({\cal X}''_{\bul \leq N}/{\cal X}''/S'/S_{-1})}V
@VV{N_{K_N/K}({\cal X}'_{\bul \leq N}/{\cal X}'/S/S_{-1})}V \\
H^q_{\rm inf}({\mathfrak X}_{K_N}/K_N) @>>>
H^q_{\rm inf}({\mathfrak X}_{K'_N}/K'_N)@<<<
H^q_{\rm inf}({\mathfrak X}'_{K_N}/K_N). 
\end{CD} 
\tag{6.3.26.2}\label{cd:lkinln}
\end{equation*} 
This shows the desired independence of 
$N_{K_N/K}({\cal X}_{\bul \leq N}/{\cal X}/S/S_{-1})$.  
Now the proof of the independence of 
$N_{\ol{K}/K}({\cal X}_{\bul \leq N}/{\cal X}/S/S_{-1})$ is obvious by (\ref{cd:lkinln}). 
\end{proof}

\begin{defi}\label{defi:mop}
We call  the following well-defined operator
\begin{align*}
N_{{\mathfrak X}_{K_N}/K_N/K}
:= N_{K_N/K}({\cal X}_{\bul \leq N}/{\cal X}/S/S_{-1})
\col &H^q_{\rm inf}({\mathfrak X}_{K_N}/K_N)
\tag{6.3.27.1}\label{ali:xsxkoif}\\
&\lo H^q_{\rm inf}({\mathfrak X}_{K_N}/K_N). 
\end{align*}  
the {\it $p$-adic monodromy operator} 
of ${\mathfrak X}_{K_N}/K_N/K$. 
We also call the following well-defined operator 
\begin{align*}
N_{{\mathfrak X}_{\ol{K}}/\ol{K}/K}
:= N_{\ol{K}/K}({\cal X}_{\bul \leq N}/{\cal X}/S/S_{-1})
\col H^q_{\rm inf}({\mathfrak X}_{\ol{K}}/\ol{K})
\lo H^q_{\rm inf}({\mathfrak X}_{\ol{K}}/\ol{K}). 
\tag{6.3.27.2}\label{ali:xsknkkoif}
\end{align*}  
the {\it $p$-adic monodromy operator} 
of ${\mathfrak X}_{\ol{K}}/\ol{K}/K$. 
\end{defi}   

\begin{rema}\label{rema:opo} 
Our $p$-adic monodromy operator (\ref{ali:xsknkkoif}) 
of ${\mathfrak X}_{\ol{K}}/\ol{K}$ 
is a generalization of 
the $p$-adic monodromy operator in \cite[A.~11]{miepl}:  
in [loc.~cit.] Mieda has defined the $p$-adic monodromy 
\begin{align*}
N_{{\mathfrak X}_{\ol{K}}/\ol{K}}\col H^q_{\rm dR}(Y_{\ol{K}}/\ol{K})
\lo H^q_{\rm dR}(Y_{\ol{K}}/\ol{K})
\tag{6.3.28.1}\label{ali:xkykbf}
\end{align*} 
for a proper smooth scheme $Y/K$. 
Though he has used Fontaine's conjecture 
${\rm C}_{\rm st}$ which has been proved by Tsuji (\cite[(0.2)]{tst}), 
we have not used ${\rm C}_{\rm st}$. 
\end{rema}

\par 
The following is a generalizations 
of \cite[Definition A.~12 (ii)]{miepl} 
(see also \cite[(4.13.2) (2)]{tst}): 

\begin{prop}\label{prop:gemmid}  
Let $L$ be $K_N$ or $\ol{K}$. 
Let $x$ and $y$ be elements of 
$H^q_{\rm inf}({\mathfrak X}_{L}/L)$ and 
$H^{q'}_{\rm inf}({\mathfrak X}_{L}/L)$, respectively. 
Set $N_{{\mathfrak X}_{L}/K}^{[i]}:=(i!)^{-1}N^i_{{\mathfrak X}_{L}/K}$ $(i\in {\mab N})$. 
Then the Leibniz rule for $N_{{\mathfrak X}_{L}/K}^{[i]}$ holds$:$ 
\begin{align*} 
N^{[i]}_{{\mathfrak X}_{L}/K}(x\cup y)
=\sum_{j=0}^iN_{{\mathfrak X}_{L}/K}^{[j]}(x)\cup N^{[i-j]}_{{\mathfrak X}_{L}/K}(y).
\tag{6.3.29.1}\label{cd:lknxyl} 
\end{align*} 
\end{prop}
\begin{proof} 
This follows from (\ref{lemm:lbg}). 
\end{proof}

\section{Actions of crystalline Weil-Deligne groups on the infinitesimal 
cohomologies of proper schemes in mixed characteristics}\label{faiw}
Let $\ol{K}$ be as in the previous section and 
let ${\mathfrak X}$ be a proper scheme over $\ol{K}$. 
In this section we define 
a canonical action of the crystalline Weil-Deligne group 
on the infinitesimal cohomology $H^q_{\rm inf}({\mathfrak X}/\ol{K})$ $(q\in {\mab N})$.  
First let us recall the crystalline Weil-Deligne group 
defined by Ogus (\cite{ollc}). 
\par 
Let the notations be as in the previous section.  
Let $K^{\rm nr}_0$ be 
the maximal unramified extension of $K_0$ in $\ol{K}$. 
In \cite[(4.1)]{boi} Berthelot and Ogus have defined 
the crystalline Weil group of $\ol{K}$: this is, by definition, 
the group of automorphisms of $\ol{K}$ whose restrictions to 
$K^{\rm nr}_0$ are the integral powers of Frobenius automorphism of 
$K^{\rm nr}_0$. 
Set $I_{\rm crys}(\ol{K}):={\rm Gal}(\ol{K}/K^{\rm nr}_0)$.  
Then we have the following exact sequence 
\begin{align*} 
1\lo I_{\rm crys}(\ol{K})\lo 
W_{\rm crys}(\ol{K})\os{{\rm deg}}{\lo}{\mab Z}\lo 1.
\end{align*}  
In the introduction of \cite{ollc}, 
Ogus has defined the crystalline Weil-Deligne group ${\bf W}_{\rm crys}(\ol{K})$: 
we denote it by  ${\rm WD}_{\rm crys}(\ol{K})$: 
${\rm WD}_{\rm crys}(\ol{K})
:= \ol{K}(1)\rtimes W_{\rm crys}(\ol{K})$ of $\ol{K}/K_0$ 
with the action of $W_{\rm crys}(\ol{K})$ on $\ol{K}(1)=\ol{K}$ 
by the following formula: 
\begin{align*} 
g\cdot a:=p^{-{\rm deg}(g)}g(a) \quad (g\in W_{\rm crys}(\ol{K}), 
a\in \ol{K}). 
\end{align*} 
By the definition of ${\rm WD}_{\rm crys}(\ol{K})$,  
${\rm WD}_{\rm crys}(\ol{K})$ fits into the following exact sequence 
\begin{align*} 
0\lo \ol{K}(1) \lo {\rm WD}_{\rm crys}(\ol{K}){\lo} W_{\rm crys}(\ol{K})\lo 0. 
\end{align*}
Let $\pi \col {\rm WD}_{\rm crys}(\ol{K}){\lo} W_{\rm crys}(\ol{K})$ 
be the projection above. 
Let $V$ be a finite dimensional vector space over $\ol{K}$. 
Let $\rho \col {\rm WD}_{\rm crys}(\ol{K})
\lo {\rm GL}_{{\mab Q}_p}(V)$ be a morphism of groups. 
Following \cite[p.~188]{boi}, 
we say that $\rho$ is {\it semi-linear} if 
\begin{align*} 
&\rho((a,g))(\al v)=g(\al)\rho((a,g))(v) \\
&((a,g)\in \ol{K}(1)\times W_{\rm crys}(\ol{K})~(as~an~element~of~a~set), 
\al\in \ol{K},v\in V). 
\end{align*}
When $\rho$ is semi-linear, we denote $\rho \col {\rm WD}_{\rm crys}(\ol{K})
\lo {\rm GL}_{{\mab Q}_p}(V)$ by 
$\rho \col {\rm WD}_{\rm crys}(\ol{K})
\lo {\rm GL}^{\rm sl}(V)$ to indicate the semi-linearity. 

Let us recall the following Ogus' result: 

\begin{prop}[{\bf \cite[Proposition 39]{ollc}}]\label{prop:ce}
The category of log convergent $F$-isocrystals on 
$s_{-1}/{\rm Spf}({\cal W})$ is equivalent to the category of 
finite dimensional vector spaces over $\ol{K}$ with a semi-linear action of 
${\rm WD}_{\rm crys}(\ol{K})$. 
If one choose a uniformizer of ${\cal V}$, 
then this category is equivalent to the category of finite dimensional 
vector spaces over $K_0$ with a nilpotent endomorphism 
$N$ and the Frobenius-linear endomorphism $F$ such 
that $NF=pFN$.  
\end{prop} 
\par 
Let ${\mathfrak X}$ be a proper scheme over $K$. 
To define the canonical action of 
${\rm WD}_{\rm crys}(\ol{K})$ on 
$H^q_{\rm inf}({\mathfrak X}_{\ol{K}}/\ol{K})$ 
$(q\in {\mab N})$, it suffices to find a finite extension $L/K$ 
and a well-defined $L_0$-form 
$H^q$ of $H^q_{\rm inf}({\mathfrak X}_{\ol{K}}/\ol{K})$ 
(that is, $H^q\otimes_{L_0}\ol{K}\simeq 
H^q_{\rm inf}({\mathfrak X}_{\ol{K}}/\ol{K})$)
with a well-defined nilpotent endomorphism 
$N\col H^q\lo H^q$ and a well-defined $\sig$-linear endomorphism 
$F\col H^q\lo H^q$ if one fixes  a uniformizer $\pi_L$ of the integer ring of $L$.  
Here $L_0:=L\cap K_0^{\rm nr}$ and $\sig \col K_0^{\rm nr}\lo 
K_0^{\rm nr}$ is the Frobenius endomorphism of $K_0^{\rm nr}$. 
\par 
Now to find this $L_0$-form is not difficult. 
Indeed we have the following isomorphism 
\begin{align*} 
H^q(\Psi_{\pi})\col 
H^q_{\rm crys}(X_{\bul \leq N}/{\cal W}(s))\otimes_{\cal W}K_N
\os{\sim}{\lo} 
H^q_{\rm inf}({\mathfrak X}_{K_N}/K_N) 
\quad (q\in {\mab N}) 
\tag{6.4.1.1}\label{ali:xkns} 
\end{align*} 
((\ref{eqn:fkysci})) 
for a uniformizer $\pi$ of $K_N$.  
Because $H^q_{\rm crys}(X_{\bul \leq N}/{\cal W}(s))\otimes_{{\cal W}}K_{N,0}$ 
has the monodromy operator 
\begin{align*} 
N_{\rm zar}\col H^q_{\rm crys}(X_{\bul \leq N}/{\cal W}(s))\otimes_{{\cal W}}K_{N,0}
\lo H^q_{\rm crys}(X_{\bul \leq N}/{\cal W}(s))\otimes_{{\cal W}}K_{N,0}
\end{align*} 
and the Frobenius action 
\begin{align*}
F\col H^q_{\rm crys}(X_{\bul \leq N}/{\cal W}(s))\otimes_{{\cal W}}K_{N,0}
\lo H^q_{\rm crys}(X_{\bul \leq N}/{\cal W}(s))\otimes_{{\cal W}}K_{N,0}
\end{align*}
such that $N_{\rm zar}F=pFN_{\rm zar}$ ((\ref{eqn:minabv})). 
We would like to prove the following: 

\begin{theo}\label{theo:kn0} 
$(1)$ The $K_{N,0}$-form 
$${\rm Im}(H^q(\Psi_{\pi})\vert_{H^q_{\rm crys}(X_{\bul \leq N}/{\cal W}(s))
\otimes_{{\cal W}}K_{N,0}})$$ 
of $H^q_{\rm inf}({\mathfrak X}_{K_N}/K_N)$ is independent of the choices  
of ${\cal X}$, $N$ and ${\cal X}_{\bul \leq N}$.  
$($It depends only on ${\mathfrak X}/K$ and $K_{N,0}$.$)$
\par 
$(2)$ The $K_{N,0}$-form in $(1)$ has the Frobenius endomorphism 
which is independent of the choice of of ${\cal X}$, $N$ 
and ${\cal X}_{\bul \leq N}$. 
$($It depends only on ${\mathfrak X}/K$, $K_{N,0}$ and a uniformizer $\pi$.$)$
\end{theo}
\begin{proof}
(1): Let ${\cal X}'$ be another proper flat model of ${\mathfrak X}$ over ${\cal V}$ 
and let ${\cal X}'_{\bul \leq N'}$ be 
another gs $N'$-truncated simplicial generically proper hypercovering of 
${\cal X}'_{{\cal V}_N}$. 
We may assume that $N'\geq N$ and even $N'=N$ by considering 
${\cal X}'_{\bul \leq N}$. 
Let ${\cal X}''$ be another proper flat model of ${\mathfrak X}$ over ${\cal V}$ 
constructed by ${\cal X}$ and 
${\cal X}'$ in (\ref{lemm:cmd}). 
By (\ref{prop:coif})  
there exist a finite extension ${\cal V}'_N$ and 
a gs $N$-truncated simplicial generically proper hypercovering  
${\cal X}''_{\bul \leq N}$ of ${\cal X}''_{{\cal V}'_N}$ 
fitting into the following commutative diagram: 
\begin{equation*} 
\begin{CD} 
{\cal X}_{\bul \leq N}@<<< {\cal X}''_{\bul \leq N}@>>> {\cal X}'_{\bul \leq N}\\
@VVV @VVV @VVV \\
{\cal X}_{{\cal V}_N}@<<< {\cal X}''_{{\cal V}'_N}@>>> {\cal X}'_{{\cal V}_N}. 
\end{CD} 
\end{equation*} 
Let $\pi'$ be a uniformizer of ${\cal V}'_N$. 
Let $X'_{\bul \leq N}$ and $X''_{\bul \leq N}$ 
be the log special fibers of ${\cal X}'_{\bul \leq N}$ and ${\cal X}''_{\bul \leq N}$ 
over $S$ and $S'$, respectively, where $S'$ is the log scheme ${\rm Spec}({\cal V}'_N)$ 
endowed with the canonical log structure.  
Let ${\mathfrak X}'_{\bul \leq N}$ and ${\mathfrak X}''_{\bul \leq N}$ 
be the log special fibers of ${\cal X}'_{\bul \leq N}$ and ${\cal X}''_{\bul \leq N}$ 
over $S$ and $S'$, respectively. 
By (\ref{eqn:rgkahk}) we have the following commutative diagram: 
\begin{equation*} 
\begin{CD} 
H^q_{\rm crys}(X_{\bul \leq N}/{\cal W}(s))\otimes_{{\cal W}}K'_N 
@>{H^q(\Psi_{\pi})\otimes {\rm id}_{K'_N},~\sim}>> 
H^q_{\rm dR} ({\mathfrak X}_{\bul \leq N}/K_N)\otimes_{K_N}K'_N\\
@| @| \\
H^q_{\rm crys}(X''_{\bul \leq N}/{\cal W}(s'))\otimes_{{\cal W}'}K'_N 
@>{H^q(\Psi_{\pi'}){\rm exp}(-\log(a) N_{\rm zar}),~\sim}>> 
H^q_{\rm dR} ({\mathfrak X}''_{\bul \leq N}/K'_N)\\
@| @|\\
H^q_{\rm crys}(X'_{\bul \leq N}/{\cal W}(s))\otimes_{{\cal W}}K'_N 
@>{H^q(\Psi_{\pi})\otimes {\rm id}_{K'_N},~\sim}>> 
H^q_{\rm dR} ({\mathfrak X}_{\bul \leq N}/K_N)\otimes_{K_N}K'_N, 
\end{CD} 
\tag{6.4.2.1}\label{eqn:rgnwhk}
\end{equation*} 
where $a:=\pi'{}^e/\pi$. 
This commutative diagram implies (1).  
\par 
(2): (2) follows from the following commutative diagram: 
\begin{equation*} 
\begin{CD} 
H^q_{\rm crys}(X_{\bul \leq N}/{\cal W}(s)) 
@>{F}>> 
H^q_{\rm crys}(X_{\bul \leq N}/{\cal W}(s)) \\
@VVV @VVV \\
H^q_{\rm crys}(X''_{\bul \leq N}/{\cal W}(s'))
@>{F}>> H^q_{\rm crys}(X''_{\bul \leq N}/{\cal W}(s'))\\
@AAA @AAA\\
H^q_{\rm crys}(X'_{\bul \leq N}/{\cal W}(s)) 
@>{F}>> 
H^q_{\rm crys}(X'_{\bul \leq N}/{\cal W}(s)).  
\end{CD} 
\tag{6.4.2.2}\label{eqn:rgfwhk}
\end{equation*}  
\end{proof}

By (\ref{theo:kn0}) we obtain the following: 

\begin{coro}\label{coro:qns}
Let ${\mathfrak X}$ be a proper scheme over $K$. 
Let $q$ be an nonnegative integer. 
Then there exist a finite extension $L/K$ and 
an well-defined $L_0$ structure $H^q_{L_0}$ 
of $H^q_{\rm inf}({\mathfrak X}_L/L)$ 
such that $H^q_{L_0}$ has a well-defined nilpotent monodromy operator
$N\col H^q_{L_0}\lo H^q_{L_0}$ and a bijective well-defined 
Frobenius-linear endomorphism
$F\col H^q_{L_0}\lo H^q_{L_0}$ such that 
$NF=pFN$. 
\end{coro}

\begin{rema} 
(\ref{coro:qns}) is a generalizations of Tsuji's result 
stated in the introduction of \cite{tsst}, which says that 
for a proper semistable family ${\cal X}/{\cal V}$, 
the $K_0$-structure,   the Frobenius action and the monodromy operator 
on $H^q_{\rm dR}({\cal X}_K/K)$ is independent of 
the proper semistable model ${\cal X}/{\cal V}$ of ${\cal X}_K/K$. 
In the next section we give a stronger result than this result. 
\end{rema}

\begin{coro}\label{coro:cisc}
Let the notations be as in {\rm (\ref{coro:qns})}. 
Then $H^q_{\rm inf}({\cal X}_{\ol{K}}/\ol{K})$ has the following 
semi-linear action 
\begin{align*} 
\rho'{}^q_{{\mathfrak X}_{\ol{K}},{\rm crys}}\col {\rm WD}_{\rm crys}(\ol{K})\lo 
{\rm GL}^{\rm sl}(H^q_{\rm inf}({\mathfrak X}_{\ol{K}}/\ol{K}))
\tag{6.4.5.1}\label{ali:wckl} 
\end{align*}
of ${\rm WD}_{\rm crys}(\ol{K})$. 
\end{coro} 

\begin{rema}\label{rema:nzcdr}
Fix a uniformizer $\pi_{-1}$ of ${\cal V}$. 
Then we can see that the representation (\ref{ali:wckl}) depends only on 
${\mathfrak X}/K$ (and $\pi_{-1}$). 
Indeed, let ${\cal X}_{\bul \leq N}/{\cal V}_N$, ${\mathfrak X}_{\bul \leq N}$,  
$X_{\bul \leq N}$, $K_N$ and $K_{N,0}$ be as in the proof of (\ref{theo:kn0}). 
Let ${\cal X}'_{\bul \leq N}/{\cal V}'_N$, $X'_{\bul \leq N'}$, 
${\mathfrak X}'_{\bul \leq N}$, $K'_{N'}$ and $K'_{N',0}$ be 
analogous objects of the objects above. 
Then we may assume that $N'=N$ without of loss of generality. 
Let ${\cal X}''_{\bul \leq N}/{\cal V}''_N$, ${\mathfrak X}''_{\bul \leq N}$,  
$X''_{\bul \leq N'}$, $K''_{N}$ and $K''_{N,0}$ 
be analogous objects which cover the objects of two kinds above simultaneously. 
Let $\pi$, $\pi'$ and $\pi''$ be uniformizers of ${\cal V}_N$, 
${\cal V}'_N$ and ${\cal V}''_N$, respectively. 
Let $e$, $e'$ and $e''$ be the ramification indexes of 
${\cal V}_N/{\cal V}$, ${\cal V}'_N/{\cal V}$ and 
${\cal V}''_N/{\cal V}$, respectively. 
Set $a:=\pi^e/\pi_{-1}$, $a':=\pi'{}^{e'}/\pi_{-1}$ and $a'':=\pi''{}^{e''}/\pi_{-1}$. 
Let ${\cal W}'$ and ${\cal W}''$ be the Witt rings of the residue fields of 
${\cal V}'$ and ${\cal V}''$, respectively. 
By (\ref{eqn:rgkahk}) we have the following commutative diagram by using
(\ref{eqn:ragkkhk}): 
\begin{equation*} 
\begin{CD} 
H^q_{\rm crys}(X_{\bul \leq N}/{\cal W}(s))\otimes_{{\cal W}}K''_N 
@>{H^q(\Psi_{\pi}){\rm exp}(-\log(a) e^{-1}N_{\rm zar}),~\sim}>> 
H^q_{\rm dR} ({\mathfrak X}_{\bul \leq N}/K_N)\otimes_{K_N}K''_N\\
@| @| \\
H^q_{\rm crys}(X''_{\bul \leq N}/{\cal W}(s''))\otimes_{{\cal W}''}K''_N 
@>{H^q(\Psi_{\pi''}){\rm exp}(-\log(a'') e''{}^{-1}N_{\rm zar}),~\sim}>> 
H^q_{\rm dR} ({\mathfrak X}_{\bul \leq N}/K''_N)\\
@| @|\\
H^q_{\rm crys}(X'_{\bul \leq N}/{\cal W}(s'))\otimes_{{\cal W}'}K''_N 
@>{H^q(\Psi_{\pi'}){\rm exp}(-\log(a')e'{}^{-1}N_{\rm zar})~\sim}>> 
H^q_{\rm dR} ({\mathfrak X}'_{\bul \leq N}/K_N)\otimes_{K'_N}K''_N.  
\end{CD} 
\tag{6.4.6.1}\label{eqn:rcpwhk}
\end{equation*}  
This commutative diagram shows 
that the representation (\ref{ali:wckl}) depends only on 
${\mathfrak X}/K$ (and $\pi_{-1}$). 
Hence it is reasonable to denote the representation 
(\ref{ali:wckl}) by $\rho'{}^q_{{\mathfrak X}_{\ol{K}},{\rm crys}}$.
\end{rema} 

As a generalization of \cite[(4.3)]{boi}, we prove the following 
which we prefer to (\ref{coro:cisc}) because 
we need not to fix a uniformizer of $K_N$:

\begin{theo}[{\bf A generalization of \cite[(4.3)]{boi} 
and \cite[Theorem 5]{ollc}}]\label{theo:cwa}
Let ${\mathfrak X}$ be a proper scheme over $\ol{K}$ and 
let $q$ be a nonnegative integer. 
There exists a unique semi-linear action 
\begin{align*} 
\rho^q_{\rm crys}:=
\rho^q_{{\mathfrak X},{\rm crys}}\col 
{\rm WD}_{\rm crys}(\ol{K})\lo 
{\rm GL}^{\rm sl}(H^q_{\rm inf}({\mathfrak X}/\ol{K}))
\tag{6.4.7.1}\label{cd:wda}
\end{align*}  
satisfying the following properties$:$ 
\par 
$(1)$ This action is contravariantly functorial with respect to a morphism
${\mathfrak X}\lo {\mathfrak Y}$ of proper schemes over $\ol{K}$. 
That is, for a morphism 
${\mathfrak f}_{-1}\col {\mathfrak X}\lo {\mathfrak Y}$ of 
proper schemes over $\ol{K}$, 
the following diagram is commutative for an element 
$\gam\in {\rm WD}_{\rm crys}(\ol{K})\!:$
\begin{equation*} 
\begin{CD}
H^q_{\rm inf}({\mathfrak X}/\ol{K})
@>{\rho^q_{{\mathfrak X},{\rm crys}}(\gam)}>> 
H^q_{\rm inf}({\mathfrak X}/\ol{K})\\
@A{{\mathfrak f}_{-1}^*}AA @AA{{\mathfrak f}_{-1}^*}A \\
H^q_{\rm inf}({\mathfrak Y}/\ol{K})
@>{\rho^q_{{\mathfrak Y},{\rm crys}}(\gam)}>> 
H^q_{\rm inf}({\mathfrak Y}/\ol{K}). 
\end{CD}
\tag{6.4.7.2}\label{cd:xqr}
\end{equation*} 
\par 
$(2)$ Let $N$ be a nonnegative integer such that $N\geq \dfrac{(q+1)(q+2)}{2}$. 
For a finite extension $L$ of $K$, 
let ${\cal O}_L$ be the integer ring of $L$ and let 
$L^{\rm nr}$ be the maximal non-ramified extension of $L$ in $\ol{K}$. 
If ${\cal X}$ is a proper scheme over ${\cal O}_L$ such that 
${\cal X}\otimes_{{\cal O}_L}\ol{K}={\mathfrak X}$ and 
if ${\cal X}_{\bul \leq N}$ is a  
gs $N$-truncated simplicial generically proper hypercovering 
of ${\cal X}_{{\cal V}_N}$ for some finite extension ${\cal V}_N$ over ${\cal O}_L$ 
of complete discrete valuation rings of mixed characteristics 
and if $X_{\bul \leq N}$ is the log special fiber of ${\cal X}_{\bul \leq N}$, 
then the action of ${\rm I}_{\rm crys}(\ol{K})$ on the image of the following morphism 
\begin{align*} 
&H^q_{\rm crys}(X_{\bul \leq N}/{\cal W}(s))\otimes_{{\cal W}}
K^{\rm nr}_0 \subset 
H^q_{\rm crys}(X_{\bul \leq N}/{\cal W}(s))\otimes_{{\cal W}}(K_N)^{\rm nr}
\tag{6.4.7.3}\label{ali:xwsbn}\\
& \os{H^q(\Psi),~\sim}{\lo} H^q_{\rm inf}({\cal X}_{(K_N)^{\rm nr}}/(K_N)^{\rm nr})
\subset H^q_{\rm inf}({\mathfrak X}/\ol{K}) 
\end{align*}
is trivial;
the action of $W_{\rm crys}(\ol{K})/
I_{\rm crys}(\ol{K})\simeq {\mab Z}$ 
is induced by the actions of the powers of the 
absolute Frobenius endomorphism $F_{X_{\bul \leq N}}$ of 
$X_{\bul \leq N}$ and the action of $W_{\rm crys}(\ol{K})$ on $\ol{K};$ 
the action of $\ol{K}(1)\owns a$ $(a\in \ol{K})$ is defined by 
the exponential ${\rm exp}(ae^{-1}_NN_{\rm zar})$ of the monodromy operator 
$ae^{-1}_NN_{\rm zar}\col 
H^q_{\rm crys}(X_{\bul \leq N}/{\cal W}(s))\otimes_{\cal W}\ol{K}
\lo 
H^q_{\rm crys}(X_{\bul \leq N}/{\cal W}(s))(-1)\otimes_{\cal W}\ol{K}$, 
where $e_N$ is the absolute ramification index of $K_N$. 
More precisely speaking, in the situation above, the action of $\rho^q_{\rm crys}((a,g))$ 
$(a\in \ol{K}(1),g\in W_{\rm crys}(\ol{K}))$ is defined by the following formula$:$
\begin{align*} 
\rho^q_{\rm crys}((a,g))(H^q(\Psi)(v)):=
H^q(\Psi)({\rm exp}(ae^{-1}_NN_{\rm zar})(F^*_{X_{\bul \leq N}})^{{\rm deg}(g)}(v)) 
\end{align*} 
for $v\in H^q_{\rm crys}(X_{\bul \leq N}/{\cal W}(s))\otimes_{{\cal W}}K^{\rm nr}_0$. 
$($See the following proof for the precise action of $\rho^q_{\rm crys}((a,g))$ on 
$H^q_{\rm inf}({\mathfrak X}/\ol{K}).)$
\end{theo}
\begin{proof}
Let $g$ be an element of  $W_{\rm crys}(\ol{K})$. 
Consider $g$ as an element of ${\rm WD}_{\rm crys}(\ol{K})$: 
$(0,g)\in  {\rm WD}_{\rm crys}(\ol{K})=\ol{K}(1)\rtimes W_{\rm crys}(\ol{K})$. 
Set $d:={\rm deg}(g)$. 
Let ${}_L{\mathfrak X}$ be a model of ${\mathfrak X}$ over a finite extension 
$L$ of $K$: ${}_L{\mathfrak X}\otimes_L\ol{K}={\mathfrak X}$. 
Let ${\cal X}$ be a proper flat model of $_L{\mathfrak X}$ over the integer ring ${\cal O}_L$.  
Then there exists a gs $N$-truncated simplicial generically proper hypercovering 
of ${\cal X}_{{\cal V}_N}$ over some finite extension ${\cal V}_N$ of ${\cal O}_L$. 
Consider the following isomorphism induced by 
the isomorphism (\ref{ali:rgrdhk}): 
\begin{align*} 
H^q(\Psi) \col H^q_{\rm crys}(X_{\bul \leq N}/{\cal W}(s))\otimes_{{\cal W}}
K^{\rm nr}_{N,0}\otimes_{K^{\rm nr}_{N,0}}\ol{K}\os{\sim}{\lo} 
H^q_{\rm inf}({\mathfrak X}/\ol{K}). 
\tag{6.4.7.4}\label{ali:xwsn}
\end{align*} 
Let $\Phi_N \col H^q_{\rm crys}(X_{\bul \leq N}/W(s))\lo 
H^q_{\rm crys}(X_{\bul \leq N}/{\cal W}(s))$ be the induced morphism of 
the absolute Frobenius endomorphism of $X_{\bul \leq N}$. 
Let $\sig_N$ be the Frobenius endomorphism of ${\cal W}$. 
Because $\Phi_N$ is $\sig_N$-linear, 
\begin{align*} 
\Phi_N^d\otimes g &\col 
H^q_{\rm crys}(X_{\bul \leq N}/{\cal W}(s))\otimes_{{\cal W}}
K^{\rm nr}_0\otimes_{K^{\rm nr}_0}\ol{K} \tag{6.4.7.5}\label{ali:xphin}\\
&\lo H^q_{\rm crys}(X_{\bul \leq N}/{\cal W}(s))
\otimes_{\cal W}K^{\rm nr}_0\otimes_{K^{\rm nr}_0}\ol{K} 
\end{align*}
is a well-defined $g$-linear automorphism of 
$H^q_{\rm crys}(X_{\bul \leq N}/{\cal W}(s))\otimes_{{\cal W}}
K^{\rm nr}_0\otimes_{K^{\rm nr}_0}\ol{K}$. 
By using the isomorphism (\ref{ali:xwsn}), 
we have the action of $W_{\rm crys}(\ol{K})$ on 
$H^q_{\rm inf}({\mathfrak X}/\ol{K})$. 
Because $\ol{K}^{I_{\rm crys}(\ol{K})}=K_0^{\rm nr}$, 
it is clear that 
\begin{align*} 
H^q_{\rm crys}(X_{\bul \leq N}/W(s))\otimes_{{\cal W}}K^{\rm nr}_0
\os{\subset}{\lo} 
H^q_{\rm inf}({\mathfrak X}/\ol{K})^{{\rm I}_{\rm crys}(\ol{K})}.
\end{align*}  
Since we have an automorphism 
\begin{align*} 
{\rm exp}(ae^{-1}_NN_{\rm zar})\col 
H^q_{\rm crys}(X_{\bul \leq N}/W(s))\otimes_{{\cal W}}\ol{K}
\os{\sim}{\lo} 
H^q_{\rm crys}(X_{\bul \leq N}/W(s))\otimes_{{\cal W}}\ol{K}, 
\tag{6.4.7.6}\label{ali:xphnin}
\end{align*} 
the action of $\ol{K}(1)=\ol{K}\owns a$ is defined by the exponential of 
the monodromy operator 
\begin{align*} 
{\rm exp}(ae^{-1}_NN_{\rm zar})
\col H^q_{\rm inf}({\mathfrak X}/\ol{K})
\os{\sim}{\lo}  H^q_{\rm inf}({\mathfrak X}/\ol{K}).
\tag{6.4.7.7}\label{ali:xphiin}
\end{align*} 
Let $(a_i,g_i)$ $(i=1,2)$ be an element of 
${\rm WD}_{\rm crys}(\ol{K})=\ol{K}(1)\rtimes W_{\rm crys}(\ol{K})$. 
Set $d_i:={\rm deg}(g_i)$. 
To define the action (\ref{cd:wda}),  
we have to prove that the following equality holds:
\begin{align*} 
&{\rm exp}(a_1e_N^{-1}N_{\rm zar})(\Phi^{d_1}_N\otimes g_1)
{\rm exp}(a_2e_N^{-1}N_{\rm zar})(\Phi^{d_2}_N\otimes g_2)
\tag{6.4.7.8}\label{ali:xphggn}\\
&={\rm exp}(a_1+p^{-d_1}g_1(a_2))e_N^{-1}N_{\rm zar})
(\Phi^{d_1+d_2}_N\otimes g_1g_2). 
\end{align*} 
Because $\Phi_NN_{\rm zar}=p^{-1}N_{\rm zar}\Phi_N$, 
we can easily check that (\ref{ali:xphggn}) holds. 
Consequently we have an action of ${\rm WD}_{\rm crys}(\ol{K})$ on 
$H^q_{\rm inf}({\mathfrak X}/\ol{K})$. 
\par 
We have to prove that this action is independent of the choices of   
${\cal X}_{\bul \leq N}$, ${\cal V}_N$ and $L$. 
Let $L'$, ${\cal X}'$, ${\cal V}_N'$ and ${\cal X}'_{\bul \leq N}$ be 
similar objects to $L$, ${\cal X}$, ${\cal V}_N$ 
and ${\cal X}_{\bul \leq N}$, respectively. 
Then 
${\cal X}\otimes_{{\cal O}_L}\ol{K}={\mathfrak X}
={\cal X}'\otimes_{{\cal O}_{L'}}\ol{K}$. 
We may assume that $L\subset L'$. 
Then we obtain two models 
${\cal X}\otimes_{{\cal O}_L}{{\cal O}_{L'}}$ 
and ${\cal X}'$ over ${\cal O}_{L'}$ of ${\mathfrak X}$ 
and we have another model ${\cal X}''$ over 
${\cal O}_{L'}$ of ${\mathfrak X}$ with morphisms 
${\cal X}\otimes_{{\cal O}_L}{{\cal O}_{L'}}\longleftarrow 
{\cal X}'' \lo {\cal X}'$ such that they induce ${\rm id}_{\mathfrak X}$.  
We also have a gs $N$-truncated simplicial generically proper hypercovering 
of ${\cal X}'_{{\cal V}''_N}$ of ${\cal X}''\otimes_{{\cal O}_{L'}}{{\cal V}''_N}$ 
over some finite extension ${\cal V}''_N$ of ${\cal O}_{L'}$ covering 
${\cal X}_{\bul \leq N}$ and ${\cal X}'_{\bul \leq N}$.  
Let $X'_{\bul \leq N}$ and $X''_{\bul \leq N}$ be the log special fibers 
of ${\cal X}'_{\bul \leq N}$ and ${\cal X}''_{\bul \leq N}$, respectively. 
Then, as in (\ref{rema:nzcdr}), 
we have the following commutative diagram as in (\ref{eqn:rcpwhk}): 
\begin{equation*} 
\begin{CD} 
H^q_{\rm crys}(X_{\bul \leq N}/{\cal W}(s))\otimes_{{\cal W}}K''_N 
@>{H^q(\Psi),~\sim}>> 
H^q_{\rm dR} ({\mathfrak X}_{\bul \leq N}/K_N)\otimes_{K_N}K''_N\\
@| @| \\
H^q_{\rm crys}(X''_{\bul \leq N}/{\cal W}(s''))\otimes_{{\cal W}''}K''_N 
@>{H^q(\Psi),~\sim}>> 
H^q_{\rm dR} ({\mathfrak X}_{\bul \leq N}/K''_N)\\
@| @|\\
H^q_{\rm crys}(X'_{\bul \leq N}/{\cal W}(s'))\otimes_{{\cal W}'}K''_N 
@>{H^q(\Psi),~\sim}>> 
H^q_{\rm dR} ({\mathfrak X}'_{\bul \leq N}/K_N)\otimes_{K'_N}K''_N.  
\end{CD} 
\tag{6.4.7.9}\label{eqn:rcrwhk}
\end{equation*}
Moreover we have the following commutative diagram: 
\begin{equation*} 
\begin{CD} 
H^q_{\rm crys}(X_{\bul \leq N}/{\cal W}(s))\otimes_{{\cal W}}K''_{N,0}
@>{\subset}>> 
H^q_{\rm crys}(X_{\bul \leq N}/{\cal W}(s))\otimes_{{\cal W}}K'' 
\\
@VVV @| \\
H^q_{\rm crys}(X''_{\bul \leq N}/{\cal W}(s''))\otimes_{{\cal W}''}K''_{N,0}
@>{\subset}>> 
H^q_{\rm crys}(X''_{\bul \leq N}/{\cal W}(s''))\otimes_{{\cal W}''}K''_N \\
@AAA @|\\
H^q_{\rm crys}(X'_{\bul \leq N}/{\cal W}(s'))\otimes_{{\cal W}'}K''_{N,0}
@>{\subset}>> 
H^q_{\rm crys}(X'_{\bul \leq N}/{\cal W}(s'))\otimes_{{\cal W}'}K''_N.  
\end{CD} 
\tag{6.4.7.10}\label{eqn:rcrwxhk}
\end{equation*}
Since the vertical arrows in (\ref{eqn:rcrwxhk}) are injective and 
since the dimensions of the sources and the targets of these arrows 
are the same, these arrows are isomorphisms. 
This tells us that the action of ${\rm WD}_{\rm crys}(\ol{K})$ 
is independent of the choices of   
${\cal X}_{\bul \leq N}$, ${\cal V}_N$ and $L$. 
The uniqueness of (\ref{cd:wda}) is clear since it is 
$W_{\rm crys}(\ol{K})$-semi-linear. 
\par 
The rest we have to prove is the functoriality. 
If we take a bigger $L$ if necessary, then 
we may assume that the morphism 
${\mathfrak f}_{-1}\col \ol{\mathfrak X}\lo \ol{\mathfrak Y}$ 
over $\ol{K}$ is defined over a finite extension $L$ of $K$. 
Let ${}_L{\mathfrak f}_{-1}\col  {}_L{\mathfrak X}\lo {}_L{\mathfrak Y}$ be this morphism. 
By (\ref{lemm:cmd}) and (\ref{prop:coif}) (1), (2), 
we may assume that there exists a morphism 
${}_{{\cal O}_L}{\mathfrak f}_{-1} \col {\cal X}\lo {\cal Y}$ over ${\cal O}_L$ 
which is a morphism of proper flat models of ${}_L{\mathfrak f}$ 
and that  there exists a morphism 
${\cal X}_{\bul \leq N}\lo {\cal Y}_{\bul \leq N}$ 
over ${\cal X}\lo {\cal Y}$ 
fitting into the following commutative diagram  
\begin{equation*} 
\begin{CD} 
{\cal X}_{\bul \leq N}@>>>{\cal Y}_{\bul \leq N}\\ 
@VVV @VVV \\
{\cal X}_{{\cal V}_N}@>{{}_{{\cal O}_L}{\mathfrak f}_{-1}\otimes_{\cal V}{\cal V}_N}
>> {\cal Y}_{{\cal V}_N}
\end{CD} 
\tag{6.4.7.11}\label{cd:lynxln}
\end{equation*}   
such that the vertical morphisms are 
gs $N$-truncated simplicial generically proper hypercoverings. 
Let $X_{\bul \leq N}$ and $Y_{\bul \leq N}$ be the log special fibers of 
${\cal X}_{\bul \leq N}$ and ${\cal Y}_{\bul \leq N}$, respectively. 
Then we have the following morphism 
\begin{align*} 
H^q_{\rm crys}(Y_{\bul \leq N}/{\cal W}(s))\lo 
H^q_{\rm crys}(X_{\bul \leq N}/{\cal W}(s)). 
\end{align*}
This morphism is compatible with the monodromy operator and  
the Frobenius endomorphism. 
Let ${\mathfrak X}_{\bul \leq N}$ and ${\mathfrak Y}_{\bul \leq N}$ 
be the log generic fibers of 
${\cal X}_{\bul \leq N}$ and ${\cal Y}_{\bul \leq N}$, respectively. 
By (\ref{eqn:rgkkhk}) we obtain the following commutative diagram: 
\begin{equation*} 
\begin{CD} 
H^q_{\rm crys}(X_{\bul \leq N}/{\cal W}(s))\otimes_{\cal W}K_N
@>{H^q(\Psi),\sim}>>H^q_{\rm dR}({\mathfrak X}_{\bul \leq N}/K_N)
\\
@AAA @AAA \\
H^q_{\rm crys}(Y_{\bul \leq N}/{\cal W}(s))\otimes_{\cal W}K_N
@>{H^q(\Psi),\sim}>>H^q_{\rm dR}({\mathfrak Y}_{\bul \leq N}/K_N). 
\end{CD}
\end{equation*}  
This means the desired contravariant functoriality of the action of 
${\rm WD}_{\rm crys}(\ol{K})$.
\end{proof}

\begin{rema}\label{rema:ofccc}
(1) After \cite[Theorem 5]{ollc}, Ogus has conjectured 
the contravariant functoriality of  the action of 
${\rm WD}_{\rm crys}(\ol{K})$ on $H^q_{\rm dR}({\mathfrak X}/\ol{K})$ 
if ${\mathfrak X}$ is proper and smooth over $\ol{K}$ and 
${\mathfrak X}$ has potentially semistable reduction. 
In (\ref{theo:cwa}) (1) we have proved a generalization of his conjecture.  
\par 
(2) It is possible to generalize (\ref{theo:cwa}) for 
a separated scheme of finite type over $\ol{K}$ without difficulty. 
In the next chapter we prove this. 
\end{rema}

\begin{prop}\label{prop:caccw}
The representation is compatible with the cup product of 
$H^{\bul}_{\rm inf}({\mathfrak X}/\ol{K})$. 
That is, for $x\in H^q_{\rm inf}({\mathfrak X}/\ol{K})$, 
$y\in H^{q'}_{\rm inf}({\mathfrak X}/\ol{K})$ and 
$\gam \in {\rm WD}_{\rm crys}(\ol{K})$, 
\begin{align*} 
\rho^{q+q'}_{{\mathfrak X},{\rm crys}}(\gam)(x\cup y)
=
\rho^q_{{\mathfrak X},{\rm crys}}(\gam)(x)\cup 
\rho^{q'}_{{\mathfrak X},{\rm crys}}(\gam)(y).
\tag{6.4.9.1}\label{ali:xxxy}
\end{align*}  
\end{prop} 
\begin{proof} 
Let the notations be as in the proof of (\ref{theo:cwa}). 
Because $H^q(\Psi)$ in (\ref{ali:xwsn}) is compatible with the cup products of both hand sides 
((\ref{prop:cwcup})), we have only to note that the morphisms (\ref{ali:xphin}) 
and (\ref{ali:xphnin}) are compatible with the cup product of 
the log crystalline cohomologies. 
The former is obvious and the latter has already proved in 
the proof of (\ref{prop:cwcup}). 
\end{proof}

\par 
In the rest of this section, 
we consider a more general action than the Frobenius action. 
Let $K'/K$, ${\cal Y}\lo {\cal X}$, 
${\mathfrak g}\col {\mathfrak Y}\lo {\mathfrak X}$, 
$S'\lo S$ and $g\col {\cal Y}_{\bul \leq N}\lo {\cal X}_{\bul \leq N}$ 
be as in (\ref{prop:gemid}). 
Let $Y_{\bul \leq N}/s'$ be as in the proof of (\ref{prop:gemid}). 
Let $s'_{-1}$ be the log point obtained by 
${\rm Spec}({\cal V}')$. 
Let $s'_{-1}\lo s_{-1}$ be the induced morphism 
by ${\rm Spec}({\cal V}')\lo {\rm Spec}({\cal V})$. 
Set $\kap'_{-1}:=\Gam(s'_{-1},{\cal O}_{s'_{-1}})$. 
Let 
\begin{equation*} 
\begin{CD} 
Y_{\bul \leq N}@>{g_0}>>X_{\bul \leq N}\\
@VVV @VVV \\
s'@>{u}>>s\\
@V{\bigcap}VV @VV{\bigcap}V \\
{\cal W}(s')@>{{\cal W}(u)}>>{\cal W}(s) \\
@VVV @VVV \\
{\cal W}(s'_{-1})@>{{\cal W}(u_{-1})}>>{\cal W}(s_{-1}) \\
@A{\bigcup}AA @AA{\bigcup}A \\
s'_{-1}@>{u_{-1}}>>s_{-1}
\end{CD}
\tag{6.4.9.2}\label{eqn:gwpsp}
\end{equation*} 
be a commutative diagram of (formal) log schemes 
obtained by (\ref{ali:xxpss}). 
Here $g_0$ is the induced morphism by $g$. 

\par 
Let 
\begin{align*} 
\sig_{u} \col K_N \lo K'_N
\tag{6.4.9.3}\label{ali:ksss}
\end{align*} 
be a morphism of fields extending the morphisms
$K\lo K'$ and 
${\cal W}(u)^* \col {\cal W}\lo {\cal W}'$. 
Because we have a pull-back morphism 
\begin{align*} 
u^* \col 
H^q_{{\rm crys}}(X_{\bul \leq N}/{\cal W}(s))
\lo 
H^q_{\rm crys}(Y_{\bul \leq N}/{\cal W}(s'))
\tag{6.4.9.4}\label{ali:xknutif}
\end{align*} 
of $u$,  
\begin{align*} 
u^*\otimes \sig_{u}  \col 
H^q_{{\rm crys}}(X_{\bul \leq N}/{\cal W}(s))\otimes_{{\cal W}}K_N
\lo 
H^q_{\rm crys}(X_{\bul \leq N}/{\cal W}(s))\otimes_{{\cal W}'}K'_N
\tag{6.4.9.5}\label{ali:xkautif}
\end{align*} 
of $u$, 
we have an endomorphism by using the canonical Hyodo-Kato isomorphism
\begin{align*}
(u^*_{K_N/K},\sig_u)
:=(u^*_{K_N/K}({\cal X}_{\bul \leq N}/S,{\cal Y}_{\bul \leq N}/S'),\sig_u)
\col H^q_{\rm inf}({\mathfrak X}_{K_N}/K_N)
\lo 
H^q_{\rm inf}({\mathfrak Y}_{K'_N}/K'_N). 
\tag{6.4.9.6}\label{ali:xkutif}
\end{align*}   
\par 
Let 
\begin{align*} 
\ol{\sig}_{u_{-1}} \col \ol{K} \lo \ol{K}
\tag{6.4.9.7}\label{ali:kssols}
\end{align*} 
be a morphism of fields extending the morphisms
$K'\lo K$ and 
${\cal W}(u^*_{-1}) \col {\cal W}(\kap'_{-1})\lo {\cal W}(\kap_{-1})$. 
Then the morphism (\ref{ali:kssols}) induces the following morphism 
\begin{align*}
(u^*_{-1,\ol{K}/K},\ol{\sig}_{u_{-1}}):=
(u^*_{-1,\ol{K}/K}({\cal X}_{\bul \leq N}/S,{\cal Y}_{\bul \leq N}/S'),\ol{\sig}_{u_{-1}})
\col H^q_{\rm inf}({\mathfrak Y}_{\ol{K}}/\ol{K})
\lo H^q_{\rm inf}({\mathfrak X}_{\ol{K}}/\ol{K}).
\tag{6.4.9.8}\label{ali:kuoif}
\end{align*} 
\par 
If one considers the actions 
(\ref{ali:xkutif}) and (\ref{ali:kuoif}), 
then the morphisms (\ref{ali:xsxkoif}) and (\ref{ali:xsknkkoif}) 
are the following morphisms, respectively: 
\begin{align*}
N_{{\mathfrak X}_{K_N}/K_N}\col H^q_{\rm inf}({\mathfrak X}_{K_N}/K_N)
\lo H^q_{\rm inf}({\mathfrak X}_{K_N}/K_N)(-1,u_{K_N/K}), 
\tag{6.4.9.9}\label{ali:xkumif}
\end{align*} 
\begin{align*}
N_{{\mathfrak X}_{\ol{K}}/\ol{K}}
\col H^q_{\rm inf}({\mathfrak X}_{\ol{K}}/\ol{K})
\lo H^q_{\rm inf}({\mathfrak X}_{\ol{K}}/\ol{K})(-1,u_{K_N/K}).
\tag{6.4.9.10}\label{ali:xkmuoif}
\end{align*} 
\par 
The following result tells us that  
$u^*_{\ol{K}/K}({\cal X}_{\bul \leq N}/S,{\cal Y}_{\bul \leq N}/S')$ 
is independent of the choice of 
${\cal X}_{\bul \leq N}/S$ and ${\cal Y}_{\bul \leq N}/S'$ 
in a certain sense. 

\begin{theo}\label{theo:wdmo} 
Let $L$ be $K_N$ or $\ol{K}$.  
Then the morphism 
in {\rm (\ref{ali:xkutif})} $($resp.~{\rm (\ref{ali:kuoif}))} on 
$H^q_{\rm inf}({\mathfrak X}_L/L)$ 
is independent of the choices of ${\cal X}$
and the gs $N$-truncated simplicial generically proper hypercovering  
of ${\cal X}_{\bul \leq N}$ in 
{\rm \S\ref{sec:ph}}. 
\end{theo} 
\begin{proof} 
We have only to prove this theorem only in the case $L=K_N$. 
Assume that we are given another 
$K'_N \subset K_N$, ${\cal X}'$, 
$N'\geq 2^{-1}(q+1)(q+2)$ and 
another gs $N'$-truncated generically proper hypercovering  
of ${\cal X}'_{\bul \leq N'}$. 
We may assume that $N\leq N'$. 
We may consider ${\cal X}'_{\bul \leq N}$ instead of 
${\cal X}'_{\bul \leq N'}$. 
By (\ref{prop:coif}) (2) we have other ${\cal V}_N''$ containing 
${\cal V}_N$ and ${\cal X}''$
fitting into the following commutative diagram 
\begin{equation*} 
\begin{CD} 
{\cal X}_{\bul \leq N,{\cal V}_{S''}}@<<<
{\cal X}''_{\bul \leq N}@>>>{\cal X}'_{\bul \leq N,{\cal V}_{S''}}\\ 
@VVV @VVV @VVV\\
{\cal X}_{{\cal V}_N''}@<<<{\cal X}''_{{\cal V}_N''}@>>>{\cal X}'_{{\cal V}_N''}. 
\end{CD} 
\tag{6.4.10.1}\label{cd:lkuxnln}
\end{equation*} 
By (\ref{cd:lknl}) we have the following commutative diagram 
\begin{equation*} 
\begin{CD} 
H^q_{\rm inf}({\mathfrak X}_{K_N}/K_N)_{K_N''} @>{\sim}>>
H^q_{\rm inf}({\mathfrak X}_{K_N'}'/K_N'')@<{\sim}<< 
H^q_{\rm inf}({\mathfrak X}_{K_N'}/K_N')_{K_N''}\\
@A{u^*_{K_N/K}}AA @A{u^*_{K_N''/K}}AA 
@AA{u^*_{K_N'/K}}A \\
H^q_{\rm inf}({\mathfrak Y}_{K_N}/K_N)_{K_{S''}} @>{\sim}>>
H^q_{\rm inf}({\mathfrak Y}_{K_N''}/K_N'') @<{\sim}<< 
H^q_{\rm inf}({\mathfrak Y}_{K_N'}/K_N')_{K_N''}. 
\end{CD} 
\tag{6.4.10.2}\label{cd:lakuinln}
\end{equation*} 
This commutative diagram shows (\ref{theo:wdmo}). 
\end{proof}

\begin{theo}\label{theo:indnu}
$(1)$ There exists a unique morphism 
\begin{align*} 
u_K\col H^q_{\rm inf}({\mathfrak Y}/K)\lo H^q_{\rm inf}({\mathfrak X}/K)
\tag{6.4.11.1}\label{ali:ukd}
\end{align*} 
such that $u_K\otimes_K\ol{K}=u_{\ol{K}/K}$. 
\par 
$(2)$ There exists a unique morphism 
\begin{align*} 
N_{{\mathfrak X}/K}\col H^q_{\rm inf}({\mathfrak X}/K)\lo 
H^q_{\rm inf}({\mathfrak X}/K)(-1,u_K)
\tag{6.4.11.2}\label{ali:nkd}
\end{align*} 
such that $N_{{\mathfrak X}/K}
\otimes_K\ol{K}=N_{{\mathfrak X}_{\ol{K}}/\ol{K}}$. 
\end{theo}
\begin{proof} 
(1): Let $\sig$ be an element of ${\rm Gal}(\ol{K}/K)$. 
Then $\sig u_{\ol{K}/K}= u_{\ol{K}/K} \sig$. 
Since $H^q_{\rm inf}({\mathfrak X}/K)$ is fixed by $\sig$, 
$\sig u_{\ol{K}/K}= u_{\ol{K}/K}$ on 
$H^q_{\rm inf}({\mathfrak X}/K)$. 
Hence the image of $H^q_{\rm inf}({\mathfrak X}/K)$ by 
$u_{\ol{K}/K}$ is contained in $H^q_{\rm inf}({\mathfrak X}/K)$. 
\par 
(2): Let $\pi_{-1}$ be a uniformizer of ${\cal V}$. 
Since $\sig (\pi_{-1})=\pi_{-1}$, it is easy to check that 
$\sig N_{{\mathfrak X}_{\ol{K}}/\ol{K}}= N_{{\mathfrak X}_{\ol{K}}/\ol{K}} \sig$. 
The rest of the proof is the same as that of (1). 
\end{proof}

\section{Limits of the weight filtrations on the infinitesimal 
cohomologies of proper schemes in mixed characteristics}\label{sec:pff} 
Let the notations be as in the previous section. 
In this section we endow 
$H^q_{\rm inf}({\mathfrak X}/K)$ $(q\in {\mab N})$  
with the limit $P$ of the weight filtration. 
More strongly, we show that 
there exists a well-defined $K_0$-form 
$H^q_{\rm inf}({\mathfrak X}/K/K_0)$ of $H^q_{\rm inf}({\mathfrak X}/K)$ 
and a finite increasing filtration $P_{K_0}$ such that 
$P_{K_0}\otimes_{K_0}K=P$ on $H^q_{\rm inf}({\mathfrak X}/K)$. 
We also prove that $P_{K_0}$ is contravariantly functorial, in particukar, 
$H^q_{\rm inf}({\mathfrak X}/K/K_0)$ is  contravariantly functorial. 
\par 
Let ${\cal V}_N$, $K_N$, ${\cal W}$, ${\cal X}_{\bul \leq N}/S$,  
$X_{\bul \leq N}/s$ and  $\kap$ be as in \S\ref{sec:ph}.  
Set ${\cal W}:={\cal W}(\kap)$. 
By (\ref{theo:thenad}) and (\ref{eqn:espssp}) 
we have the following spectral sequence 
\begin{equation*} 
\begin{split} 
{} & E^{-k,q+k}_1(/K_0(\kap)):=
\bigoplus_{m=0}^N\bigoplus_{j\geq \max \{-(k+m),0\}} 
H^{q-2j-k-2m}((\os{\circ}{X}{}^{(2j+k+m)}_m
/{\cal W})_{\rm crys}, \\ 
{} & \phantom{R^{q-2j-k-k}f_{(\os{\circ}{X}^{(k)}, 
Z\vert_{\os{\circ}{X}^{(2j+k)}})/{\cal W}*} 
{\cal O}}
{\cal O}_{\os{\circ}{X}{}^{(2j+k+m)}_m/{\cal W}} 
\otimes_{\mab Z}\varpi^{(2j+k+m)}_{\rm crys}
(\os{\circ}{X}_m/{\cal W}))_{\mab Q} \\
&\Lo 
H^q_{{\rm crys}}(X_{\bul \leq N}/{\cal W}(s))_{\mab Q} 
\quad (q\in {\mab Z}). 
\end{split} 
\tag{6.5.0.1}\label{eqn:getepsp}
\end{equation*}   
Here we omit the action of $u^*$ in (\ref{eqn:espssp}) for the time being. 

\begin{theo}\label{theo:mrfwt}
$(1)$ Let $q$ be a nonnegative integer. 
Let $N$ be an integer which is greater than or equal to 
$2^{-1}(q+1)(q+2)$. 
The spectral sequence 
{\rm (\ref{eqn:getepsp})}$\otimes_{K_{N,0}}K_N$  
converges to $H^q_{\rm inf}({\mathfrak X}_{K_N}/K_N)$. 
Consequently there exists the following spectral sequence 
\begin{equation*} 
\begin{split} 
{} & E^{-k,q+k}_1(/K_N):=
\bigoplus_{m=0}^N\bigoplus_{j\geq \max \{-(k+m),0\}} 
H^{q-2j-k-2m}((\os{\circ}{X}{}^{(2j+k+m)}_m/{\cal W})_{\rm crys}, \\ 
{} & \phantom{
\vert_{\os{\circ}{X}^{(2j+k)}})/{\cal W}* 
{\cal O}}
{\cal O}_{\os{\circ}{X}{}^{(2j+k+m)}_m/{\cal W}} 
\otimes_{\mab Z}\varpi^{(2j+k+m)}_{\rm crys}
(\os{\circ}{X}_m/{\cal W}(\kap)))
_{\mab Q}\otimes_{K_{N,0}}K_N \\
&\Lo 
H^q_{{\rm inf}}({\mathfrak X}_{K_N}/K_N)
_{\mab Q} 
\quad (q\in {\mab Z}). 
\end{split} 
\tag{6.5.1.1}\label{eqn:getniepsp}
\end{equation*}    
and there exists a filtration $P$ on 
$H^q_{\rm inf}({\mathfrak X}_{K_N}/K_N)$ 
such that 
${\rm gr}_{q+k}^PH^q_{\rm inf}({\mathfrak X}_{K_N}/K_N)
\simeq E^{-k,q+k}_{\infty}(/K_N):=
E^{-k,q+k}_{\infty}(/K_{N,0})\otimes_{K_{N,0}}K_N$.  
The spectral sequence {\rm (\ref{eqn:getniepsp})} degenerates at $E_2$. 
\par 
$(2)$ 
For $N'\geq N$, there exists a canonical isomorphism 
$E^{-k,q+k}_{\infty}(/K_N)
\otimes_{K_N}K_{N'} \os{\sim}{\lo}
E^{-k,q+k}_{\infty}(/K_{N'})$. 
\end{theo}
\begin{proof} 
(1): By (\ref{eqn:fkysci}) 
we have the following isomorphism 
\begin{align*} 
H^q_{\rm inf}({\mathfrak X}_{K_N}/K_N) \os{\sim}{\lo} 
H^q_{{\rm crys}}(X_{\bul \leq N}/{\cal W}(s))\otimes_{{\cal W}}K_N. 
\tag{6.5.1.2}\label{ali:wkknn}
\end{align*} 
Hence we obtain the spectral sequence (\ref{eqn:getniepsp}) 
by (\ref{eqn:emtttsp}).  
The degeneration at $E_2$ of this spectral sequence 
follows from (\ref{theo:e2dam}). 
Hence we obtain (1). 
\par 
(2): The base change theorem of crystalline cohomologies 
and the existence of the spectral sequence 
(\ref{eqn:getepsp}) show (2). 
\end{proof} 

\begin{rema}\label{rema:wib}
The $E^{ij}_r$-terms 
$E_r^{ij}(/K_{N,0})$ 
$(1 \leq r \leq \infty, i+j=q)$ 
of the spectral sequence 
(\ref{eqn:getepsp}) are independent of the choice of $N$ 
(\cite[(2.2)]{nh3}) in the sense  
$E_r^{ij}(/K_{N,0})\otimes_{K_{N,0}}K_{N',0}
=E_r^{ij}(/K_{N',0})$ $(N\leq N')$ 
since $({\cal X}_{\bul \leq N'})_{\bul \leq N}=
{\cal X}_{\bul \leq N}\otimes_{{\cal V}_N}{\cal V}_{N'}$, 
and the graded objects of  
$H^q_{\rm inf}({\mathfrak X}_{K_N}/K_N)$ 
by the filtration induced by the spectral sequence 
(\ref{eqn:getepsp})$\otimes_{K_{N,0}}K_N$ 
are $E^{ij}_{\infty}$-terms $(i+j=q)$ of (\ref{eqn:getepsp}) 
with tensor $\otimes_{K_{N,0}}K_N$. 
\end{rema}

By (\ref{eqn:getniepsp}) 
we have the following spectral sequence 
\begin{align*} 
E_1^{-k,q+k}& = \bigoplus^N_{m\geq 0}
\bigoplus_{j\geq \max \{-(k+m),0\}}H^{q-2j-k-2m}
((\os{\circ}{X}{}^{(2j+k+m)}_m/{\cal W})_{\rm crys}, 
({\cal O}_{\os{\circ}{X}{}^{(2j+k+m)}_m/{\cal W}}
\otimes_{\mab Z} 
\tag{6.5.2.1}\label{eqn:espwiiffsp} \\
& \vp^{(2j+k+m)}_{\rm crys}(\os{\circ}{X}_m/S))
\otimes_{{\cal W}}\ol{K}
\Lo 
H^q_{\rm inf}({\mathfrak X}_{\ol{K}}/\ol{K}). 
\end{align*}

\begin{defi}\label{defi:psw}
We call (\ref{eqn:espwiiffsp}) the {\it limit weight spectral sequence} 
of ${\cal X}_{\bul \leq N}/{\cal X}_K/\ol{K}$. 
\end{defi} 

\par 
The following is one of main results in this book: 

\begin{theo}\label{theo:xkp} 
$(1)$ The induced filtration $P_{\ol{K}}$ on 
$H^q_{\rm inf}({\mathfrak X}_{\ol{K}}/\ol{K})$ 
by the spectral sequence {\rm (\ref{eqn:espwiiffsp})} 
is independent of the choice of 
${\cal X}$ and ${\cal X}_{\bul \leq N}$. 
\par 
$(2)$ There exists a well-defined filtration $P$ 
on $H^q_{\rm inf}({\mathfrak X}/K)$ such that 
$P\otimes_K\ol{K}=P_{\ol{K}}$ on 
$H^q_{\rm inf}({\mathfrak X}_{\ol{K}}/\ol{K})$.
\par 
$(3)$ $0=P_{-1} \subset \cdots \subset P_{2q}=
H^q_{\rm inf}({\mathfrak X}/K)$. 
\end{theo}
\begin{proof}
(1): 
Let the notations be as in the proofs of (\ref{rema:nzcdr}). 
Let $P^{{\cal X}_{\bul \leq N}}_{\ol{K}}$ and $P^{{\cal X}''_{\bul \leq N}}_{\ol{K}}$ 
be the filtrations on $H^q_{\rm inf}({\mathfrak X}_{\ol{K}}/\ol{K})$ 
obtained by ${\cal X}_{\bul \leq N}$ and  ${\cal X}''_{\bul \leq N}$, respectively. 
Then we have the natural injection  
$P^{{\cal X}_{\bul \leq N}}_{\ol{K}}\os{\sus}{\lo} P^{{\cal X}''_{\bul \leq N}}_{\ol{K}}$.  
By the strict compatibility  (\ref{theo:stpbgb}), 
we have the equality 
$P^{{\cal X}_{\bul \leq N}}_{\ol{K}}=P^{{\cal X}''_{\bul \leq N}}_{\ol{K}}$.  
This proves (1). 
\par 
(2): (T.~Tsuji kindly pointed out the following argument to the author.) 
We may assume that the finite extension $K_N/K$ is Galois. 
Let $\sig$ be an element of ${\rm Gal}(\ol{K}/K)$. 
Set 
\begin{align*} 
f_{\sig}:={\rm id}_{H^q_{\rm inf}({\mathfrak X}/K)}\otimes \sig 
\col &
H^q_{\rm inf}({\mathfrak X}/K)\otimes_K\ol{K}=
H^q_{\rm inf}({\mathfrak X}_{\ol{K}}/\ol{K})\\
&\lo 
H^q_{\rm inf}({\mathfrak X}/K)\otimes_K\ol{K}=
H^q_{\rm inf}({\mathfrak X}_{\ol{K}}/\ol{K}).
\end{align*}  
Then $f_{{\rm id}_{\ol{K}}}={\rm id}_{H^q_{\rm inf}({\mathfrak X}_{\ol{K}}/\ol{K})}$ 
and $f_{\tau \sig}=f_{\tau}\circ f_{\sig}$ for $\sig, \tau\in {\rm Gal}(\ol{K}/K)$. 
We have two $N$-truncated simplicial generically proper hypercoverings 
${\cal X}_{\bul \leq N}$ and $\sig({\cal X}_{\bul \leq N})$ 
of ${\mathfrak X}\otimes_KK_N$. 
The covering $\sig({\cal X}_{\bul \leq N})$ induces the filtration 
$\sig(P_{\ol{K}})$ on $H^q_{\rm inf}({\mathfrak X}_{\ol{K}}/\ol{K})$.  
By the well-definedness of $P_{\ol{K}}$ in (1), 
we see that $P_{\ol{K}}=\sig(P_{\ol{K}})$. 
Hence the $\sig$-linear map $f_{\sig} \col 
H^q_{\rm inf}({\mathfrak X}_{\ol{K}}/\ol{K})\lo H^q_{\rm inf}({\mathfrak X}_{\ol{K}}/\ol{K})$ 
induces a $\sig$-linear map $(P_{\ol{K}})_k H^q_{\rm inf}({\mathfrak X}_{\ol{K}}/\ol{K})
\lo (P_{\ol{K}})_k H^q_{\rm inf}({\mathfrak X}_{\ol{K}}/\ol{K})$ $(k\in {\mab Z})$. 
Now (2) follows from the Galois descent. 
\par 
(3): Assume that $k< -q$. 
Consider the spectral sequence 
(\ref{eqn:espwiiffsp}). 
Then $q-2j-k-2m<-2(k+m+j)\leq 0$. 
Hence $E^{-k,q+k}_1=0$.  
Assume that $k>q$. 
Then $q-2j-k-2m=q-k-2(j+m)<0$. Hence $E^{-k,q+k}_1=0$. 
As a result, we obtain (3). 
\end{proof}

The following is also one of main results in this book. 
The following tells us that $H^q_{\rm inf}({\mathfrak X}/K)$ 
has a well-defined deep ``hidden structure'' 
(cf.~Illusie's observation stated in the introduction of \cite{ollc}): 
 
\begin{theo}\label{theo:k0st}
There exists a well-defined $K_0$-form $H^q_{\rm inf}({\mathfrak X}/K/K_0)$ of 
$H^q_{\rm inf}({\mathfrak X}/K)$ and a well-defined filtration 
$P_{K_0}$ on $H^q_{\rm inf}({\mathfrak X}/K/K_0)$ such that 
$P_{K_0}\otimes_{K_0}K=P$ on $H^q_{\rm inf}({\mathfrak X}/K)$. 
\end{theo}
\begin{proof} 
Let the notations be as in the proof of (\ref{theo:xkp}). 
We may assume that the finite extension $K_N/K_0$ is Galois. 
Let $K_{N,0}$ be the maximal unramified extension of $K_0$ in $K_N$. 
Then $K_N$ is a totally ramified extension of $K_{N,0}K$ and 
$K_{N,0}\cap K=K_0$: 
\begin{equation*} 
\begin{CD} 
K_{N,0}@>{\subset}>> K_{N,0}K @>{\subset}>> K_N\\
@A{\bigcup}AA @AA{\bigcup}A \\
K_0@>{\subset}>> K. 
\end{CD} 
\end{equation*} 
Let $\sig$ be an element of ${\rm Gal}(K/K_0)$. 
Then we can consider $\sig$ as an element of ${\rm Gal}(K_{N,0}K/K_{N,0})$ 
since ${\rm Gal}(K_{N,0}K/K_{N,0})={\rm Gal}(K/K_0)$. 
Let $\sig({\mathfrak X})$ be the base change of ${\mathfrak X}$ by using 
$\sig$ (\cite[p.~191]{boi})
Let 
\begin{align*} 
f_{\sig} \col H^q_{\rm inf}({\mathfrak X}/K)\lo H^q_{\rm inf}(\sig({\mathfrak X})/K)
\end{align*} 
be the induced morphism by $\sig$. 
We would like to prove that 
\begin{align*} 
f_{\sig \tau}=\sig(f_{\tau})\circ f_{\sig}
\col  H^q_{\rm inf}({\mathfrak X}/K)\lo  H^q_{\rm inf}(\sig({\mathfrak X})/K)\lo 
H^q_{\rm inf}(\sig \tau({\mathfrak X})/K)
\end{align*} 
for
$\sig, \tau\in {\rm Gal}(K/K_0)$.
Set 
\begin{align*} 
f_{\sig, N}:=f_{\sig}\otimes_K\sig 
\col &H^q_{\rm inf}({\mathfrak X}_{KK_{N,0}}/KK_{N,0})
=H^q_{\rm inf}({\mathfrak X}/K)\otimes_K(KK_{N,0})\\
&\lo 
H^q_{\rm inf}({\mathfrak X}/K)\otimes_K(KK_{N,0})=
H^q_{\rm inf}({\mathfrak X}_{KK_{N,0}}/KK_{N,0}).
\end{align*}  
By the comparison isomorphism 
\begin{align*}
H^q(\Psi) \col H^q_{\rm crys}(X_{\bul \leq N}/{\cal W}(s))_{K_{N,0}}
\otimes_{K_{N,0}}KK_{N,0}\os{\sim}{\lo} 
H^q_{\rm inf}({\mathfrak X}_{KK_{N,0}}/KK_{N,0}), 
\tag{6.5.5.1}\label{ali:fscif} 
\end{align*} 
we see that $H^q(\Psi)(H^q_{\rm crys}(X_{\bul \leq N}/{\cal W}(s))_{K_{N,0}})$ 
is a $K_{N,0}$-form of $H^q_{\rm inf}({\mathfrak X}_{KK_{N,0}}/KK_{N,0})$. 
Set 
\begin{align*} 
g_{\sig,N}:={\rm id}_{H^q_{\rm crys}(X_{\bul \leq N}/{\cal W}(s))_{K_{N,0}}}\otimes \sig 
\col &H^q_{\rm crys}(X_{\bul \leq N}/{\cal W}(s))_{K_{N,0}}\otimes_{K_{N,0}}KK_{N,0}\\
&\lo 
H^q_{\rm crys}(\sig(X_{\bul \leq N})/{\cal W}(s))_{K_{N,0}}\otimes_{K_{N,0}}KK_{N,0}. 
\end{align*} 
Here $\sig$ induces the identity of $X_{\bul \leq N}$. 
Then the following diagram is commutative: 
\begin{equation*} 
\begin{CD} 
H^q_{\rm crys}(X_{\bul \leq N}/{\cal W}(s))_{K_{N,0}}\otimes_{K_{N,0}}KK_{N,0}
@>{H^q(\Psi),\sim} >>H^q_{\rm inf}({\mathfrak X}_{KK_{N,0}}/KK_{N,0})\\
@V{g_{\sig,N}}VV @VV{f_{\sig,N}}V \\
H^q_{\rm crys}(\sig(X_{\bul \leq N})/{\cal W}(s))_{K_{N,0}}\otimes_{K_{N,0}}KK_{N,0}
@>{H^q(\Psi),\sim} >>H^q_{\rm inf}(\sig({\mathfrak X}_{KK_{N,0}})/KK_{N,0}). 
\end{CD} 
\tag{6.5.5.2}\label{cd:fsxsf} 
\end{equation*} 
Hence $\sig(f_{\tau, N})\circ  f_{\sig, N}=f_{\sig \tau, N}$ 
by the Galois descent. 
This means that $f_{\sig \tau}\otimes_K(\sig \tau)=
(\sig(f_{\tau})\circ f_{\sig})\otimes_K(\sig \tau)$. 
Hence $f_{\sig \tau}=\sig(f_{\tau})\circ f_{\sig}$. 
Hence there exists a canonical $K_0$-form 
$H^q_{\rm inf}({\mathfrak X}/K/K_0)$ of $H^q_{\rm inf}({\mathfrak X}/K)$ 
by the Galois descent: 
\begin{align*} 
H^q_{\rm inf}({\mathfrak X}/K)=H^q_{\rm inf}({\mathfrak X}/K/K_0)\otimes_{K_0}K. 
\tag{6.5.5.3}\label{ali:fsxkf} 
\end{align*} 
\par 
Set $P(N):=P\otimes_K(KK_{N,0})$.  Let $P(N,0)$ (resp.~$\sig^*(P(N))$) be the filtration on 
$H^q_{\rm inf}({\mathfrak X}_{KK_{N,0}}/KK_{N,0})$ (resp.~
$H^q_{\rm inf}(\sig({\mathfrak X}_{KK_{N,0}})/KK_{N,0})$) defined by 
the weight filtration on $H^q_{\rm crys}(X_{\bul \leq N}/{\cal W}(s))_{K_{N,0}}$ 
and the isomorphism (\ref{ali:fscif}) (resp.~the lower isomorphism in (\ref{cd:fsxsf})).  
Then $P(N,0)$ is a $K_{N,0}$-form of $P(N)$. 
Hence $f_{\sig,N}(P(N))=\sig^*(P(N))$. This means that $f_{\sig}(P)=\sig^*(P)$ in 
$H^q_{\rm inf}(\sig({\mathfrak X})/K)$. Because the relation 
$f_{\sig \tau}=\sig(f_{\tau})\circ f_{\sig}$ holds as we have seen, 
the existence of $P_{K_0}$ on $H^q_{\rm inf}({\mathfrak X}/K/K_0)$ is clear 
by the Galois descent. 
\end{proof} 

\begin{defi}\label{defi:lmmmt}
We call the filtration $P_{K_0}$ on $H^q_{\rm inf}({\mathfrak X}/K/K_0)$
the {\it hidden limit weight filtration} of $H^q_{\rm inf}({\mathfrak X}/K/K_0)$. 
\end{defi} 

\begin{rema}
By the definition of 
$P_{K_0}$ on $H^q_{\rm inf}({\mathfrak X}/K/K_0)$, 
we have the following natural morphism 
\begin{align*} 
H^q(\Psi) \col P_kH^q_{{\rm crys}}(X_{\bul \leq N}/{\cal W}(s))\otimes_{{\cal W}}K_{N,0}
\os{\sim}{\lo} 
(P_{K_0})_kH^q_{\rm inf}({\mathfrak X}/K/K_0)\otimes_{K_0}K_{N,0} 
\quad (k\in {\mab Z}). 
\tag{6.5.7.1}\label{ali:ppkk}
\end{align*} 
\end{rema}

\begin{coro}\label{coro:xcs}
The log crystalline cohomology 
$H^q_{\rm crys}(X_{\bul \leq N}/{\cal W}(s))\otimes_{{\cal W}}K_{N,0}$ 
is independent of the choice of the $N$-truncated simplicial gs generically 
proper hypercovering of ${\mathfrak X}_{K_N}/K_N$.  
More strongly, $P_kH^q_{{\rm crys}}(X_{\bul \leq N}/{\cal W}(s))\otimes_{{\cal W}}K_{N,0}$ 
is independent of the choice of the $N$-truncated simplicial gs generically 
proper hypercovering of ${\mathfrak X}_{K_N}/K_N$.  
\end{coro}

\begin{rema}\label{rema:dj} 
Let $U$ be a separated scheme of finite type over $\kap_{-1}$.
Let $j \col U \os{\sus}{\lo} \ol{U}$ be 
an open immersion into a proper scheme over $\kap_{-1}$.
Let $(U_{\bul},X_{\bul})$ be a split proper hypercovering
of $(U,\ol{U})$ in the sense of 
\cite{tzp}, that is, $(U_{\bul},X_{\bul})$ is split, 
$U_{\bul}$ is a proper hypercovering of $U$, 
$X_{\bul}$ is a proper simplicial scheme over $\ol{U}$ 
and $U_{\bul}=X_{\bul}\times_{\ol{U}}U$. 
Moreover, using de Jong's alteration theorem (\cite{dj}),   
we can require that 
$X_{\bul}$ is a proper smooth simplicial scheme over $\kap$ 
and that $U_{\bul}$ is the complement of 
a simplicial SNCD $D_{\bul}$ 
on  $X_{\bul}$ over $\kap$. 
We have called such a split proper hypercovering 
$(U_{\bul}, X_{\bul})$ of $(U, \ol{U})$  
a gs$($=good and split$)$ proper  hypercovering of $(U, \ol{U})$. 
In \cite{nh3} we have proved
\begin{equation*}
R\Gam_{\rm rig}(U/K)=
R\Gam((X_{\bul},D_{\bul})/{\cal W})\otimes_{\cal W}K.  
\tag{6.5.9.1}\label{eqn:itrorlc}
\end{equation*} 
In particular, we see that $R\Gam((X_{\bul},D_{\bul})/{\cal W})_K$ is independent 
of the choice of a gs$($=good and split$)$ proper  hypercovering of $(U, \ol{U})$.
\end{rema}

By the remark above, we would like to conjecture the following: 

\begin{conj}\label{conj:indp}
Let the notations be as in {\rm (\ref{coro:xcs})}. 
Then we conjecture that 
\begin{align*} 
R\Gam_{\rm inf}({\mathfrak X}_{K_N}/K_N)=
R\Gam_{\rm crys}(X_{\bul \leq N}/{\cal W}(s))\otimes^L_{{\cal W}}K_N. 
\tag{6.5.10.1}\label{ali:mxkn} 
\end{align*} 
In particular, the complex 
$R\Gam_{\rm crys}(X_{\bul \leq N}/{\cal W}(s))\otimes^L_{{\cal W}}K_N$ 
is independent of 
the choice of the $N$-truncated simplicial gs generically 
proper hypercovering ${\cal X}_{\bul \leq N}$ of ${\mathfrak X}_{K_N}/K_N$.  
\end{conj}

\begin{theo}[{\bf Strict compatibility}]\label{theo:stinfgb}
Let $K\os{\subset}{\lo} K'$ be a morphism of 
the fraction fields of complete discrete valuation rings 
of mixed characteristics $(0,p)$. 
Let ${\mathfrak Y}$ be a proper scheme over $K'$. 
Let $f\col {\mathfrak Y}\lo {\mathfrak X}$ be a morphism over 
${\rm  Spec}(K')\lo  {\rm  Spec}(K)$. 
Let 
\begin{align*} 
f^*\col H^q_{\rm inf}({\mathfrak X}/K)\otimes_{K}K'\lo 
H^q_{\rm inf}({\mathfrak Y}/K')
\tag{6.5.11.1}\label{ali:yxxlk} 
\end{align*} 
be the induced morphism by $f$. 
Let $K'_0$ be the maximal unramified extension of $K_0$ in $K'$.  
Then the following hold$:$
\par 
$(1)$ There exists a morphism 
\begin{align*} 
f^*_0\col H^q_{\rm inf}({\mathfrak X}/K/K_0)\otimes_{K_0}K'_0
\lo 
H^q_{\rm inf}({\mathfrak Y}/K'/K'_0)
\tag{6.5.11.2}\label{ali:yxalk} 
\end{align*} 
such that $f^*_0\otimes_{K'_0}K'=f^*$. 
This morphism satisfies the transitive relation and 
${\rm id}^*_{\mathfrak X}=
{\rm id}_{H^q_{\rm inf}({\mathfrak X}/K/K_0)}$.  
\par 
$(2)$ The morphism $f^*_0$ is strictly compatible with 
$P_{K_0}\otimes_{K_0}K'_0$ and $P_{K'_0}$ 
on $H^q_{\rm inf}({\mathfrak X}/K/K_0)\otimes_{K_0}K'_0$ and 
$H^q_{\rm inf}({\mathfrak Y}/K'/K'_0)$, respectively. 
\end{theo}
\begin{proof}
(1): Let the notations be as in the proof of (\ref{theo:k0st})  and after 
(\ref{prop:caccw}). 
Let $K'_{N,0}$ be an analogous field of $K_{N,0}$ for $K_N'$. 
Consider the following morphism 
\begin{align*} 
f^* \col H^q_{\rm inf}({\mathfrak X}_{K_{N,0}K}/K_{N,0}K)
\otimes_{K_{N,0}K}K'_{N,0}K'\lo 
H^q_{\rm inf}({\mathfrak Y}_{K'_{N,0}K'}/K'_{N,0}K'), 
\tag{6.5.11.3}\label{ali:ykklk} 
\end{align*} 
which is a scalar extension of the morphism (\ref{ali:yxxlk}).  
Here we have denoted this morphism by $f^*$ by abuse of notation. 
Let the notations be as in the proof of (\ref{prop:gemid}). 
Then we have the following commutative diagram 
\begin{equation*} 
\begin{CD} 
H^q_{\rm inf}({\mathfrak X}_{K_{N,0}K}/K_{N,0}K)\otimes_{K_{N,0}K}K'_{N,0}K'
@>{f^*}>> H^q_{\rm inf}({\mathfrak Y}_{K'_{N,0}K'}/K'_{N,0}K')\\
@A{H^q(\Psi)}A{\simeq}A @A{\simeq}A{H^q(\Psi)}A \\
H^q_{{\rm crys}}(X_{\bul \leq N}/{\cal W}(s))
\otimes_{\cal W}K'_{N,0}K'@>{f^*}>>
H^q_{{\rm crys}}(Y_{\bul \leq N}/{\cal W}(s'))
\otimes_{{\cal W}'}K'_{N,0}K. 
\end{CD} 
\tag{6.5.11.4}\label{ali:qinfx}
\end{equation*} 
The lower morphism $f^*$ in (\ref{ali:qinfx}) is defined over $K'_{N,0}$. 
Let $\sig$ be an element of ${\rm Gal}(K'_{N,0}K'/K'_{N,0})$. 
Let $f^*$ be the upper morphism in (\ref{ali:qinfx}). 
Then $\sig \circ f^*=f^*\circ \sig$. 
Hence, for an element $\sig$ of ${\rm Gal}(K'/K'_{0})$ and 
for the morphism $f^*$ in (\ref{ali:yxxlk}), $\sig(f^*):=\sig \circ f\circ \sig^{-1}=f^*$. 
This tells us the existence of $f^*_0$ by the Galois descent. 
\par 
(2): 
It suffices to prove that 
the pull-back 
\begin{equation*} 
f^*\col H^q_{\rm inf}({\mathfrak X}_{\ol{K}}/\ol{K})
\lo H^q_{\rm inf}({\mathfrak Y}_{\ol{K}}/\ol{K})
\tag{6.5.11.5}\label{eqn:xxylk} 
\end{equation*} 
by $f$ is strictly compatible with the $($induced$)$ weight filtrations.  
This immediately follows from (\ref{prop:coif}), (\ref{coro:indne}) 
and (\ref{theo:stpbgb}). 
\end{proof}

\begin{theo}\label{theo:lmon}
Let the notations be as in the previous section. 
Then the following hold$:$ 
\par 
$(1)$ There exists a 
well-defined $p$-adic monodromy operator 
\begin{align*} 
N_{{\mathfrak X}/K/K_0}\col H^q_{\rm inf}({\mathfrak X}/K/K_0)
\lo H^q_{\rm inf}({\mathfrak X}/K/K_0)
\end{align*} 
such that 
$N_{{\mathfrak X}/K/K_0}\otimes_{K_0}K
=N_{{\mathfrak X}/K} \col H^q_{\rm inf}({\mathfrak X}/K)\lo 
H^q_{\rm inf}({\mathfrak X}/K)$. 
\par 
$(2)$ There exists an action $F_{/K_0}$ on 
$H^q_{\rm inf}({\mathfrak X}/K/K_0)$ 
such that 
$F_{/K_0}\otimes_{K_0}K_{N,0}=F$ on the 
$K_{N,0}$-form 
$${\rm Im}(H^q(\Psi) \vert_{H^q_{\rm crys}(X_{\bul \leq N}/{\cal W}(s))
\otimes_{{\cal W}}K_{N,0}}).$$ 
\end{theo}
\begin{proof} 
(1): The morphism $N_{{\mathfrak X}_{K_N}/K_N}$ is defined over $K_{N,0}$. 
Let the notations be as in the proof of (\ref{theo:k0st}). 
Then it is clear that $\sig(N_{{\mathfrak X}_{K_N}/K_N})=
N_{{\mathfrak X}_{K_N}/K_N}$ 
since $\sig(\pi)$ is also a uniformizer of $K_N$ 
(cf.~(\ref{cd:vsim})). 
Hence $N_{{\mathfrak X}_{K_N}/K_N}$ is defined over $K_0$. 
\par 
(2): Let the notations be as in the proof of (\ref{theo:k0st}). 
Then $\sig$ induces a morphism 
${\rm Spec}(\sig) \col \sig(X_{\bul \leq N})=X_{\bul \leq N}\lo X_{\bul \leq N}$. 
It is clear that we have the following commutative diagram 
\begin{equation*} 
\begin{CD} 
\sig(X_{\bul \leq N})@>\sig>> X_{\bul \leq N}\\
@V{F}VV @VV{F}V\\
\sig(X_{\bul \leq N})@>\sig>> X_{\bul \leq N}. 
\end{CD} 
\end{equation*} 
Since $\sig$ induces the identity morphism
of ${\mathfrak X}_{K}$ and since 
$\sig({\cal X}_{\bul \leq N})$ 
is a generically gs proper hypercovering 
of ${\mathfrak X}_{K_N}$, we obtain (2) 
by using the following three (commutative) diagrams: 
\begin{equation*} 
\begin{CD}
H^q_{\rm crys}(\sig(X_{\bul \leq N})/{\cal W}(s))
\otimes_{{\cal W}}K_{N,0}\otimes_{K_{N,0}}K_N
@>{H^q(\Psi),\sim}>> 
H^q_{\rm inf}({\mathfrak X}/K)\otimes_KK_N\\
@V{F}VV \\
H^q_{\rm crys}(\sig(X_{\bul \leq N})/{\cal W}(s))
\otimes_{{\cal W}}K_{N,0}
\otimes_{K_{N,0}}K_N@>{H^q(\Psi),\sim}>> 
H^q_{\rm inf}({\mathfrak X}/K)\otimes_KK_N
\end{CD}
\end{equation*} 
over 
\begin{equation*} 
\begin{CD}
H^q_{\rm crys}(X_{\bul \leq N}/{\cal W}(s))
\otimes_{{\cal W}}K_{N,0}
\otimes_{K_{N,0}}K_N@>{H^q(\Psi),\sim}>> 
H^q_{\rm inf}({\mathfrak X}/K)\otimes_KK_N\\
@V{F}VV \\
H^q_{\rm crys}(X_{\bul \leq N}/{\cal W}(s))
\otimes_{{\cal W}}K_{N,0} 
\otimes_{K_{N,0}}K_N@>{H^q(\Psi)s,\sim}>> 
H^q_{\rm inf}({\mathfrak X}/K)\otimes_KK_N 
\end{CD}
\end{equation*} 
and the following commutative diagram  
\begin{equation*} 
\begin{CD}
H^q_{\rm crys}(\sig(X_{\bul \leq N})/{\cal W}(s))
\otimes_{{\cal W}}K_{N,0}
\otimes_{K_{N,0}}K_N@>{\sig,\sim}>> 
H^q_{\rm crys}(X_{\bul \leq N}/{\cal W}(s))
\otimes_{{\cal W}}K_{N,0}\\
@V{F}VV @VV{F}V\\
H^q_{\rm crys}(\sig(X_{\bul \leq N})/{\cal W}(s))
\otimes_{{\cal W}}K_{N,0}
\otimes_{K_{N,0}}K_N@>{\sig,\sim}>> 
H^q_{\rm crys}(X_{\bul \leq N}/{\cal W}(s))
\otimes_{{\cal W}}K_{N,0}. 
\end{CD}
\end{equation*} 
\end{proof}

\par 
Let the notations be as in the beginning of this section. 
Then we have the following spectral sequence by (\ref{ali:ppkk}): 
\begin{align*} 
E_1^{-k,q+k}& = \bigoplus^N_{m\geq 0}
\bigoplus_{j\geq \max \{-(k+m),0\}}H^{q-2j-k-2m}
((\os{\circ}{X}{}^{(2j+k+m)}_m/{\cal W})_{\rm crys}, 
({\cal O}_{\os{\circ}{X}{}^{(2j+k+m)}_m/{\cal W}}
\otimes_{\mab Z} 
\tag{6.5.12.1}\label{eqn:eswniiffsp} \\
& \vp^{(2j+k+m)}_{\rm crys}(\os{\circ}{X}_m/S))
\otimes_{{\cal W}}K_{N,0}
\Lo 
H^q_{\rm inf}({\mathfrak X}/K/K_0)\otimes_{K_0}K_{N,0}. 
\end{align*}
More precisely, if we consider the action of the Frobenius endomorphism, 
the spectral sequence above is the following spectral sequence 
\begin{align*} 
E_1^{-k,q+k}& = \bigoplus^N_{m\geq 0}
\bigoplus_{j\geq \max \{-(k+m),0\}}H^{q-2j-k-2m}
((\os{\circ}{X}{}^{(2j+k+m)}_m/{\cal W})_{\rm crys}, 
({\cal O}_{\os{\circ}{X}{}^{(2j+k+m)}_m/{\cal W}}
\otimes_{\mab Z} 
\tag{6.5.12.2}\label{eqn:esfiiffsp} \\
& \vp^{(2j+k+m)}_{\rm crys}(\os{\circ}{X}_m/S))(-j-k-m)
\otimes_{{\cal W}}K_{N,0}
\Lo 
H^q_{\rm inf}({\mathfrak X}/K/K_0)\otimes_{K_0}K_{N,0}. 
\end{align*}

\begin{coro}\label{coro:frig} 
If $\kap$ is finite, then 
$P_{K_0}$ is the filtration defined by the weight of eigenvalues of 
the Frobenius endomorphism on $H^q_{\rm inf}(X_K/K/K_0)$. 
\end{coro}
\begin{proof} 
It suffices to prove that $P_{K_0}\otimes_{K_0}K_{N,0}$ on 
$H^q_{\rm inf}(X_K/K/K_0)\otimes_{K_0}K_{N,0}$
is the filtration defined by the weight of eigenvalues of 
the Frobenius endomorphism on
$H^q_{\rm inf}(X_K/K/K_0)\otimes_{K_0}K_{N,0}$. 
By the spectral sequence (\ref{eqn:esfiiffsp}) and 
the purity of the eigenvalues of the Frobenius endomorphism on 
the classical crystalline cohomology of a 
proper smooth scheme over a finite field 
(\cite[Corollary 1 2)]{kme}, \cite[(1.2)]{clpu}, \cite[(2.2) (4)]{ndw} 
(see also (\ref{rema:nllfcs}) (1))), 
we immediately obtain (\ref{coro:frig}).  
\end{proof} 

\begin{rema}\label{rema:nfcs}
See \S\ref{sec:padicar} below for the generalization of 
(\ref{coro:frig}) in the case where $\kap$ is not necessarily finite. 
\end{rema}

\begin{prop}\label{prop:cupp}
$(1)$ There exists the following product 
\begin{align*} 
(?\cup ?)_{K_0}\col 
H^q_{\rm inf}({\mathfrak X}/K/K_0)\otimes_{K_0}H^{q'}_{\rm inf}({\mathfrak X}/K/K_0)
\lo H^{q+q'}_{\rm inf}({\mathfrak X}/K/K_0)
\end{align*}
such that $(?\cup ?)_{K_0}\otimes_{K_0}K$ is the cup product 
of $H^q_{\rm inf}({\mathfrak X}/K)$.  
\par 
$(2)$ Set $N_{{\mathfrak X}/K/K_0}^{[i]}:=(i!)^{-1}N_{{\mathfrak X}/K/K_0}$ $(i\in {\mab N})$. 
Then the Leibniz rule for $N_{{\mathfrak X}/K/K_0}^{[i]}$ holds$:$ 
$N^{[i]}_{{\mathfrak X}/K/K_0}((x\cup y)_{K_0})=
\sum_{j=0}^i(N^{[j]}_{{\mathfrak X}/K/K_0}(x)\cup N^{[i-j]}_{{\mathfrak X}/K/K_0}(y))_{K_0}$   
$(x\in H^{q'}_{\rm inf}({\mathfrak X}/K/K_0),~y\in H^{q''}_{\rm inf}({\mathfrak X}/K/K_0))$. 
\par 
$(3)$ 
$N_{{\mathfrak X}/K/K_0}(P_kH^q_{\rm inf}({\mathfrak X}/K/K_0))
\subset P_{k-2}H^q_{\rm inf}({\mathfrak X}/K/K_0)\quad (k\in {\mab Z})$. 
\end{prop}
\begin{proof} 
(1): The proof is the same as that of (\ref{theo:stinfgb}). 
\par 
(2): (2) follows from (1) and (\ref{cd:lknxyl}). 
\par 
(3): (3) follows from (\ref{prop:nucsaa}). 
\end{proof}

\begin{prop}\label{prop:irs}
Let ${\mathfrak X}$ be a proper scheme over $\ol{K}$. 
The action {\rm (\ref{cd:wda})} preserves 
the limit weight filtration $P$ on $H^q_{\rm inf}({\mathfrak X}/\ol{K})$.  
\end{prop}
\begin{proof} 
Because (\ref{ali:xphin}) preserves the limit weight filtration, the action of 
the crystalline Weil group preserves the limit weight filtration. 
Because $N_K\col H^q_{\rm inf}({\mathfrak X}/\ol{K})\lo H^q_{\rm inf}({\mathfrak X}/\ol{K})$ 
induces the morphism 
$N_{{\mathfrak X}/\ol{K}}
\col P_kH^q_{\rm inf}({\mathfrak X}/\ol{K})\lo 
P_{k-2}H^q_{\rm inf}({\mathfrak X}/\ol{K})$, 
\begin{align*} 
{\rm exp}(aN_{{\mathfrak X}/\ol{K}}) \col H^q_{\rm inf}({\mathfrak X}/\ol{K})
\lo H^q_{\rm inf}({\mathfrak X}/\ol{K})
\quad (a\in \ol{K})
\end{align*}  
preserves the limit weight filtration. 
Hence the action of ${\rm WD}_{\rm crys}(\ol{K})$ on 
$H^q_{\rm inf}({\mathfrak X}/\ol{K})$ preserves the limit weight filtration. 
\end{proof} 

The following is a reformulation of Tsuji's result (\cite{tsgep}). 

\begin{theo}\label{theo:nicor}
The following equality 
\begin{align*} 
H^q_{\rm inf}({\mathfrak X}_{K_N}/K_N)=
D_{\rm st}(H^q_{{\rm {e}}{\rm t}}({\mathfrak X}_{\ol{K}},{\mab Q}_p))
\otimes_{L_0}L
\end{align*} 
holds. This equality is compatible with the monodromy operator and 
the Frobenius action. 
Here $L_0$ is the fraction field of the Witt ring of the residue field of $L$ 
and $D_{\rm st}(H^q_{{\rm {e}}{\rm t}}(X_{\ol{K}},{\mab Q}_p)):=
(B_{\rm st}\otimes_{L_0}
H^q_{{\rm {e}}{\rm t}}(X_{\ol{K}},{\mab Q}_p))^{{\rm Gal}(\ol{K}/L)}$. 
\end{theo}

\begin{rema}\label{rema:mj} 
(\ref{theo:nicor}) tells us that 
the cohomology $H^q_{\rm inf}({\cal X}_{K_N}/K_N)$ has a 
canonical $K_{N,0}$-structure and it is 
a filtered $(\phi,N)$-module:   
it has the Frobenius action 
(more generally the action $u^*_K\otimes_KK_N$) and 
the monodromy action with the relation $N\phi=p\phi N$ and 
$H^q_{\rm inf}({\mathfrak X}_{K_N}/K_N)$ has the Hodge filtration 
((\ref{defi:hfi})).  
More strongly we have proved that 
$H^q_{\rm inf}({\mathfrak X}/K)$ has 
a canonical $K_0$-structure ((\ref{theo:k0st})) 
and it is a filtered $(\phi,N)$-module ((\ref{theo:lmon}), (\ref{defi:hfi})):  
$H^q_{\rm inf}({\mathfrak X}/K/K_0)$  has 
the Frobenius action (more generally the action $u^*_K$) 
and the monodromy action and  
$H^q_{\rm inf}({\cal X}_{K}/K/K_0)\otimes_{K_0}K_N$ has the Hodge filtration. 
This fact is a generalization of Jannsen's conjecture 
stated in the introduction of \cite{hk}. 
Moreover we have proved that it has the limit weight filtration 
((\ref{defi:lmmmt})). 
The author does not check whether one can reconstruct 
the limit weight filtration from
$H^q_{{\rm {e}}{\rm t}}({\mathfrak X}_{\ol{K}},{\mab Q}_p)$.
\end{rema}

\begin{prop}\label{theo:rntr}  
Let the notations be as in {\rm \S\ref{sec:ph}} and {\rm (\ref{theo:wdmo})}.  
Let $H^q_{\rm inf}({\mathfrak X}/K/K_0)
=\bigoplus_{i\in {\mab N}}H^q_{\rm inf}({\mathfrak X}/K/K_0)_{[i,i+1)}$ 
be the slope decomposition. 
Then the monodromy operator 
$N_{{\mathfrak X}/K/K_0}\col H^q_{\rm inf}({\mathfrak X}/K/K_0)
\lo H^q_{\rm inf}({\mathfrak X}/K/K_0)$ 
induces the following morphism 
\begin{align*} 
N_{{\mathfrak X}/K/K_0}
\col H^q_{\rm inf}({\mathfrak X}/K/K_0)\lo 
H^q_{\rm inf}({\mathfrak X}/K/K_0)(-1). 
\tag{6.5.19.1}\label{ali:cylkn} 
\end{align*} 
More precisely, it induces the following morphism 
\begin{align*} 
N_{{\mathfrak X}/K/K_0}\col P_kH^q_{\rm inf}({\mathfrak X}/K/K_0)_{[i,i+1)}
\lo P_{k-2}H^q_{\rm inf}({\mathfrak X}/K/K_0)_{[i-1,i)}(-1). 
\tag{6.5.19.2}\label{ali:cylkun} 
\end{align*} 
\end{prop} 
\begin{proof} 
Because we have the following relation 
\begin{align*}
(N_{{\mathfrak X}/K/K_0}\otimes_{K_0}K_{N,0})F
=pF(N_{{\mathfrak X}/K/K_0}\otimes_{K_0}K_{N,0}), 
\tag{6.5.19.3}\label{ali:pfnk}
\end{align*} 
we have the relation 
$N_{{\mathfrak X}/K/K_0}F=pFN_{{\mathfrak X}/K/K_0}$. 
\end{proof}

\begin{conj}\label{conj:mwcmc}
The following is a $p$-adic analogue of 
the famous $l$-adic monodromy-weight conjecture: 
\par 
Assume that ${\cal X}_K$ is proper and smooth over $K$. 
Let $M$ be the monodromy filtration 
on $H^q_{\rm inf}({\cal X}_K/K)=
H^q_{\rm dR}({\cal X}_K/K)$. 
Then, is ${\rm gr}^M_kH^q_{\rm dR}({\cal X}_K/K)$ 
$(k\in {\mab N})$ pure of weight $?$
\end{conj}

\begin{rema} 
To prove (\ref{conj:mwcmc}), it suffices to prove 
(\ref{conj:mwcmc}) in the case where ${\cal X}_K$ 
has a projective strict semistable model over ${\cal V}$ by 
\cite[Definition A.~12 (i), (ii)]{miepl}, Chow's lemma 
and de Jong alteration theorem and \cite[(1.2.4)]{kl}. 
Consequently it suffices to prove (\ref{conj:rcpmc}). 
\end{rema}

Let ${\rm WD}(\ol{K}/K)$ be the inverse image of 
$W_{\rm crys}(\ol{K})\cap {\rm Gal}(\ol{K}/K)$ 
in ${\rm W}_{\rm crys}(\ol{K})$ by the morphism 
${\rm WD}_{\rm crys}(\ol{K})\lo {\rm W}_{\rm crys}(\ol{K})$.  
This is nothing but the usual Weil-Deligne group 
defined in \cite[8.4.1]{dhlc}. 
Let 
\begin{align*} 
\rho^q_{\rm inf}\col 
{\rm WD}(\ol{K}/K)\lo {\rm GL}^{\rm sl}(H^q_{\rm inf}({\mathfrak X}_{\ol{K}}/\ol{K}))
\end{align*} 
be a representation obtained by the natural action of 
${\rm Gal}(\ol{K}/K)$ on $H^q_{\rm inf}({\mathfrak X}_{\ol{K}}/\ol{K})
=H^q_{\rm inf}({\mathfrak X}/K)\otimes_K\ol{K}$. 
Set $\rho^q_p(\gam):=\rho^q_{\rm inf}(\gam)\rho^q_{\rm crys}(\gam^{-1})$ 
$(\gam\in {\rm WD}(\ol{K}/K))$. 
Then we have the following $\ol{K}$-linear action 
\begin{align*} 
\rho_p \col 
{\rm WD}(\ol{K}/K)\lo {\rm GL}_{\ol{K}}
(H^q_{\rm inf}({\mathfrak X}_{\ol{K}}/\ol{K}))
\tag{6.5.21.1}\label{ali:rhop}
\end{align*} 
of ${\rm WD}(\ol{K}/K)$.

In \cite[1]{sga7} Grothendieck has proved the following result 
(cf.~\cite[Appendix]{seta}). 

\begin{theo}\label{theo:qnil}
Let ${\mathfrak X}$ be a separated scheme of finite type over $K$. 
Let $l$ be a prime which is prime to $p$. 
Let $\rho_l \col {\rm Gal}(\ol{K}/K)\lo 
{\rm GL}_{{\mab Q}_l}(H^q_{\rm et}({\mathfrak X}_{\ol{K}},{\mab Q}_l))$ 
be the $l$-adic representation. Then $\rho_l$ is quasi-unipotent. 
\end{theo}

\parno 
By (\ref{theo:qnil}) it is well-known that there exists a nilpotent morphism 
$N_l \col H^q_{\rm et}({\mathfrak X}_{\ol{K}},{\mab Q}_l)\lo 
H^q_{\rm et}({\mathfrak X}_{\ol{K}},{\mab Q}_l)(-1)$ 
such that $\rho_l(g)={\rm exp}(t_l(g)N_l)$ for any element $g$ of some open subgroup of 
the inertia group of ${\rm Gal}(\ol{K}/K)$, where 
$t_l \col {\rm Gal}(\ol{K}/K) \lo {\mab Z}_l(1)$ is the $l$-adic tame character. 
\par 
Assume that $\kap$ is a finite field ${\mab F}_q$. 
Let $F$ be an element of the Weil group $W(\ol{K}/K)$ which is a lift of 
the geometric Frobenius morphism of ${\rm Gal}(\ol{\kap}/\kap)$. 
Because $NF=qFN$, we have the action of the Weil-Deligne group 
${\rm WD}(\ol{K}/K)$ on $H^q_{\rm et}({\mathfrak X}_{\ol{K}},{\mab Q}_l)$ 
with respect to $F$. 
\begin{rema}\label{rema:lwt}
By \cite{nd} and the results in \S\ref{sec:ph} and the cohomological descent 
of $l$-adic \'{e}tale cohomologies, 
we obtain the contravariantly functorial limit weight filtration $P$ on 
$H^q_{\rm et}({\mathfrak X}_{\ol{K}},{\mab Q}_l)$ 
such that $P_{-1}H^q_{\rm et}({\mathfrak X}_{\ol{K}},{\mab Q}_l)=0$ 
and $P_{2q}H^q_{\rm et}({\mathfrak X}_{\ol{K}},{\mab Q}_l)=
H^q_{\rm et}({\mathfrak X}_{\ol{K}},{\mab Q}_l)$.  
Moreover, the action of ${\rm WD}(\ol{K}/K)$ on 
$H^q_{\rm et}({\mathfrak X}_{\ol{K}},{\mab Q}_l)$ 
preserves the filtration $P$.  
\end{rema} 
\par 
The following is an immediate generalization of Fontaine's conjecture 
in \cite[2.4.3]{fos} for a proper smooth scheme over $K$ 
and a special case of his conjecture in [loc.~cit., (2.4.7)] for a motive over $K$. 

\begin{conj}\label{conj:gfc} 
Let ${\mathfrak X}$ be a proper scheme over $K$. 
For a nonnegative integer $k$, the representations $\rho_l$ of ${\rm WD}(\ol{K}/K)$ 
on $P_kH^q_{\rm et}({\mathfrak X}_{\ol{K}},{\mab Q}_l)$ 
and $\rho_p$ of ${\rm WD}(\ol{K}/K)$ on 
$P_kH^q_{\rm inf}({\mathfrak X}_{\ol{K}}/\ol{K})$ 
are compatible.
\end{conj}

\parno 
In the case where $P$ is trivial, 
see \cite[2.4.6]{fos} for the case where (\ref{conj:gfc}) is true. 

\section{Limits of the slope filtrations on the infinitesimal cohomologies of 
proper schemes in mixed characteristics}\label{sec:lcwgac}  
Let the notations be as in the previous section.  
In this section we express 
$H^q_{\rm inf}({\mathfrak X}/K/K_0)_{[i,i+1)}$ 
in the previous section 
in a geometric way after enlarging $K_0$. 
By (\ref{ali:ppkk}) we have the following isomorphism 
\begin{align*} 
H^q_{\rm inf}({\mathfrak X}/K/K_0)\otimes_{K_{0}}K_{N,0} 
\os{\sim}{\lo} 
H^q_{\rm crys}(X_{\bul \leq N}/{\cal W}(s))\otimes_{{\cal W}(\kap_s)}K_{N,0} 
\quad (q\in {\mab N}). 
\tag{6.6.0.1}\label{ali:xknms} 
\end{align*} 

\begin{theo}\label{theo:rtr}  
The following hold$:$
\par 
$(1)$ 
\begin{align*} 
(H^q_{\rm inf}({\mathfrak X}/K/K_0)\otimes_{K_0}K_{N,0})_{[i,i+1)}
=H^{q-i}(X_{\bul \leq N},{\cal W}\Om^i_{X_{\bul \leq N}})
\otimes_{\cal W}K_{N,0}.
\tag{6.6.1.1}\label{ali:cxkn} 
\end{align*} 
More generally, 
\begin{align*} 
(P_kH^q_{\rm inf}({\mathfrak X}/K/K_{N,0})\otimes_{K_0}K_{N,0})_{[i,i+1)}
=P_kH^{q-i}(X_{\bul \leq N},{\cal W}\Om^i_{X_{\bul \leq N}})
\otimes_{\cal W}K_{N,0}.
\tag{6.6.1.2}\label{ali:cykn} 
\end{align*} 
\end{theo}
\begin{proof} 
This follows from (\ref{ali:xknms}) and (\ref{theo:e1deg}).  
\end{proof}  

\begin{coro}\label{coro:msa}
$P_kH^{q-i}(X_{\bul \leq N},{\cal W}\Om^i_{X_{\bul \leq N}})
\otimes_{{\cal W}}K_{N,0}$ is independent of the choice of 
the $N$-truncated simplicial gs generically 
proper hypercovering of ${\mathfrak X}_{K_N}/K_N$. 
\end{coro}

\begin{rema}
Let the notations be as in (\ref{rema:dj}). 
In \cite{nh3} we have proved that 
$P_kH^{q-i}(X_{\bul},{\cal W}\Om^i_{X_{\bul}})\otimes_{\cal W}K$ 
is independent of the choice of 
a gs$($=good and split$)$ proper  hypercovering of $(U,\ol{U})$.
The corollary (\ref{coro:msa}) is a mixed characteristic of 
this theorem. 
\end{rema}



\par

\section{$p$-adic arithmetic weight filtrations}\label{sec:padicar}
In this section we define the notion of the 
$p$-adic arithmetic weight filtration associated to the 
geometric class of $H^h_{\rm inf}(X/K)$ $(h\in {\mab N})$
and show that the filtration $P$ on $H^h_{\rm inf}(X/K)$ 
in the previous section coincides with 
the $p$-adic arithmetic weight filtration 
on $H^h_{\rm inf}(X/K)$ to the geometric class.

\par 
Let $Z/B$ be as in the beginning of \S\ref{sec:cfi}.  
Assume that $B$ is the formal spectrum of 
the Witt ring of a perfect field endowed with trivial log structure. 
Let $Z'/B'$ be a similar log scheme. 
Let $B\lo B'$ be a morphism of $p$-adic formal schemes such that 
the induced morphism $\Gam(B',{\cal O}_{B'})\lo \Gam(B,{\cal O}_{B})$ is 
a morphism of Witt rings of perfect fields.  
In this section we do not assume that this morphism is finite. 

\begin{defi}
Let $((E,P),\Phi_Z)$ (resp.~$((E,P),\Phi_{Z'})$) 
be an object of $F{\textrm -}{\rm IsocF}^{\rm sld}(Z/B)$ 
(resp.~$F{\textrm -}{\rm IsocF}^{\rm sld}(Z'/B')$). 
Then we say that $((E,P),\Phi_Z)$ {\it comes from} 
$((E',P'),\Phi_{Z'})$ if, for any object $T\in {\rm Enl}(Z/B)$, 
there exists an object $T'\in {\rm Enl}(Z'/B')$ with a morphism
$T\lo T'$ of log formal schemes over $B\lo B'$ fitting into 
the following commutative diagram 
\begin{equation*} 
\begin{CD} 
T_0@>>> T_0'\\
@VVV @VVV \\
Z@>>> Z'
\end{CD} 
\end{equation*} 
and there exists a functorial isomorphism 
$\rho_T \col ((E,P),\Phi_{Z'})_{T'}\otimes_{{\cal K}_{T'}}{\cal K}_T
\os{\sim}{\lo} ((E,P),\Phi_Z)_T$. 
\end{defi} 

\par 
Assume that $M_Z$ is hollow and locally split.
Let $((E,P),\Phi_Z)$ be an object of $F{\textrm -}{\rm IsocF}^{\rm sld}(Z/B)$ 
such that there exits a formal scheme ${\cal B}$ of 
topologically finite type over the 
Witt ring ${\cal W}({\mab F}_q)$ of 
the finite field ${\mab F}_q$ with $q$-elements 
and a fine log formal scheme ${\cal Z}/{\cal B}$ topologically of 
finite type fitting into the following commutative diagram 
\begin{equation*}
\begin{CD} 
Z@>>> {\cal Z}\\
@VVV @VVV\\
B@>>> {\cal B}
\end{CD} 
\end{equation*} 
and an object $(({\cal E},{\cal P}),\Phi_{\cal E})$ of 
$F{\textrm -}{\rm IsocF}^{\rm sld}({\cal Z}/{\cal B})$ 
such that $((E,P),\Phi_Z)$ comes from $(({\cal E},{\cal P}),\Phi_{\cal Z})$. 
Here 
we assume that the upper horizontal morphism 
$Z\lo {\cal Z}$ is solid. 
In the following we consider a category ${\cal C}$ consisting of above triples  
$({\cal Z}/{\cal B},({\cal E},{\cal P}),\Phi_{\cal Z})$ such that,   
for any two objects 
$({\cal Z}_i/{\cal B}_i,({\cal E}_i,{\cal P}_i),\Phi_{{\cal Z}_i})$ 
$(i=1,2)$ of ${\cal C}$, 
there exists an object 
${\cal Z}_3/{\cal B}_3$ which covers 
${\cal Z}_i/{\cal B}_i$ $(i=1,2)$ fitting into 
the following commutative diagram 
\begin{equation*}
\begin{CD} 
Z@>>> {\cal Z}_3@>>> {\cal Z}_i\\
@VVV @VVV @VVV\\
B@>>> {\cal B}_3@>>> {\cal B}_i
\end{CD} 
\end{equation*} 
such that $(({\cal E}_3,{\cal P}_3),\Phi_{{\cal E}_3})$ comes from 
$(({\cal E}_i,{\cal P}_i),\Phi_i)$ $(i=1,2)$. 
Here the morphism in ${\cal C}$ is defined in an obvious way. 

\begin{defi}
We say that $((E,P),\Phi_Z)$ is a {\it log convergent $F$-isocrystal on} 
$Z/B$ {\it with arithmetic Frobenius weight filtration} with respect to ${\cal C}$ 
if, for any above object $({\cal Z}/{\cal B},(({\cal E},{\cal P}),\Phi_{\cal E}))$ of 
${\cal C}$ and 
for any exact point $x$ of ${\cal Z}$, the eigenvalue of $\Phi^e_{\cal E}(x)$ 
of the induced subvector space 
${\cal P}(x)_k{\cal E}(x)$of ${\cal E}(x)$ 
is a Weil number of weight $\leq k$.  Here $e:=\log_pq$. 
\end{defi}

\begin{rema}\label{rema:pxk} 
Since ${\cal P}(x)_k{\cal E}(x)$ is determined by $\Phi_{\cal E}(x)$, 
${\cal P}$ is determined by $\Phi_{\cal E}$ by (\ref{prop:csds}) and 
\cite[(4.1)]{of}. 
Consequently $P$ is determined by $\Phi_{\cal E}$.
\end{rema}

\begin{defi} 
Let ${\mathfrak X}/K$ be a proper scheme. 
Let $h$ be a nonnegative integer. 
Let ${\cal C}^h({\mathfrak X}/K)$ be the following class. 
\par 
The category ${\cal C}^h({\mathfrak X}/K)$ is, by definition, the following: 
\par 
Let ${\cal X}/{\cal V}$ be a proper flat model of ${\mathfrak X}/K$. 
Let ${\cal X}_{\bul \leq N}/{\cal V}_N$ be a generically 
good proper hypercovering of ${\cal X}/{\cal V}$ in (\ref{prop:vv0}), 
where $N$ is a positive integer such that 
$N\geq 2^{-1}(h+1)(h+2)$.
Set $X_{\bul \leq N}:={\cal X}_{\bul \leq N}
\otimes_{{\cal V}_N}\kap_N$, 
where $\kap_N$ is the residue field of ${\cal V}_N$. 
Let 
\begin{equation*} 
s_N \lo  s_{N-1}\lo \cdots \lo 
s_0\lo s_{-1}
\end{equation*}
be as in the part of (\ref{eqn:wtslpnn}). 
Then there exists a successive model 
${\cal X}_{\bul \leq N,S_1}/S_1$ of $X_{\bul \leq N}/s_N$
for a formally smooth algebra 
$A$ over the Witt ring of a finite field. 
Here $S$ is a family of log points such that $\os{\circ}{S}:={\rm Spf}(A)$ 
and 
$S_1:=S\otimes_{{\mab Z}_p}{\mab F}_p$ with an injective morphism 
$\Gam(S_1,{\cal O}_{S_1})\lo \kap$ of algebras. 
Then we have the log crystalline cohomology 
$H^q_{\rm crys}({\cal X}_{\bul \leq N,S_1}/S)$ with a
weight filtration $P$ and the Frobenius endomorphism $\Phi$. 
The objects of ${\cal C}^h({\mathfrak X}/{\cal V})$ are, by definition, 
$$(S/{\rm Spf}({\mab Z}_p),
((H^q_{\rm crys}({\cal X}_{\bul \leq N,S_1}/S),P),\Phi)
{\textrm '}{\rm s}$$ 
for various $A$, $N$, ${\cal V}_N$ and ${\cal X}_{\bul \leq N}/{\cal V}_N$. 
We call ${\cal C}^h({\mathfrak X}/{\cal V})$ 
the geometric class of $H^h_{\rm inf}({\mathfrak X}/K)$. 
\end{defi} 

\begin{theo}\label{theo:clgwt}
Let $q$ be a nonnegative integer. 
Then $(H^q_{\rm inf}({\mathfrak X}/K),P)$ 
is the arithmetic Frobenius weight filtration with respect to the 
the geometric class of $X/K$ with respect to ${\cal C}^h({\mathfrak X}/{\cal V})$.
\end{theo}
\begin{proof} 
In the proof of (\ref{theo:e2dam}), (\ref{theo:e2dgfam})  
and (\ref{theo:xkp}), we have already proved this.
\end{proof}

\chapter{Generalized Ogus' conjecture and Fontaine's conjecture}
Let ${\mathfrak U}$ be a separated scheme of finite type over $\ol{K}$ and 
let $H^q_{\rm inf}({\mathfrak U}/\ol{K})$ $(q\in {\mab N})$ 
be the infinitesimal cohomology of ${\mathfrak U}/\ol{K}$. 
\par 
In this chapter we define a canonical semi-linear action of 
${\rm WD}_{\rm crys}(\ol{K})$ on $H^q_{\rm inf}({\mathfrak U}/\ol{K})$ 
and we prove the contravariant functoriality of this action quickly. 
We also prove that this action is compatible with the cup product of 
$H^q_{\rm inf}({\mathfrak U}/\ol{K})$. 
\par 
Roughly speaking, all ideas to prove these have already appeared 
before this chapter and in this chapter we have only to 
generalize our result only by adding the base change of the log structure of 
horizontal truncated simplicial SNCD on a nice truncated simplicial proper hypercovering of 
the compactification $\ol{\mathfrak U}$ of ${\mathfrak U}$. 

\section{Generalized Ogus' conjecture and Fontaine's conjecture on the action of crystalline 
Weil-Deligne group for varieties}\label{sec:ofc}
Let ${\mathfrak U}$ be a separated scheme of finite type over $K$. 
In this section, by using the action of the crystalline Weil group 
${\rm WD}_{\rm crys}(\ol{K})$ on 
$H^q_{\rm inf}({\mathfrak U_{\ol{K}}}/\ol{K})$ 
and the action of ${\rm Gal}(\ol{K}/K)$ on 
$H^q_{\rm inf}({\mathfrak U_{\ol{K}}}/\ol{K})
=H^q_{\rm inf}({\mathfrak U}/K)\otimes_K\ol{K}$, 
we can define a linear action of the usual Weil-Deligne group 
${\rm WD}(\ol{K}/K)$ on 
$H^q_{\rm inf}({\mathfrak U_{\ol{K}}}/\ol{K})$. 
As is well-known, there exists a canonical action 
of ${\rm WD}(\ol{K}/K)$ on the $l$-adic \'{e}tale cohomology 
$H^q_{\rm et}({\mathfrak U_{\ol{K}}},{\mab Q}_l)$ $(l\not=p)$.
Using these two actions, we can conjecture that these two actions 
are compatible as in \cite{fos}.  
\par 
Let ${\cal U}$ be a flat model of ${\mathfrak U}$ over ${\cal V}$. 
Embed ${\cal U}$ into a proper scheme 
$\ol{\cal U}$ over ${\cal V}$ by using 
Nagata's embedding theorem (\cite{imno}, \cite{lu}). 
Set $\ol{\mathfrak U}:=\ol{\cal U}_K$ 
and $U:={\cal U}\otimes_{\cal V}\kap$.    
Embed ${\mathfrak U}$ into $\ol{\mathfrak U}$ over $K$. 

\par 
For an algebraic extension $L$ of $K$ in $\ol{K}$, 
let ${\cal O}_L$ be the integer ring of $L$. 
Set ${\cal U}_L:={\cal U}_K\otimes_KL$ and 
${\cal U}_{{\cal O}_L}:={\cal U}\otimes_{\cal V}{\cal O}_L$. 
Let $N$ be a nonnegative integer. 
We use the same notation for 
an $N$-truncated simplicial fine log scheme 
over ${\rm Spec}({\cal V})=({\rm Spec}({\cal V}),{\cal V}^*)$. 
\par 
Let us recall the following notion (cf.~\cite[(6.3)]{dj}): 

\begin{defi}\label{defi:sst}
Let ${\cal X}$ be a proper strict semistable family over ${\cal V}$ with 
the canonical log structure. 
(${\cal X}$ is a log scheme and not only a scheme.)  
Let ${\cal D}$ be a relative SNCD(=simple normal crossing divisor) on 
$\os{\circ}{\cal X}/{\cal V}$. 
Then we have a canonical fs(=fine and saturated) log structure $M({\cal D})$ 
associated to ${\cal D}$. 
$($See \cite[(2.1.7), (2.1.9)]{nh2} for the definitions of 
the relative SNCD and $M({\cal D})$.$)$ 
Let $M_{\cal X}$ be the log structure of ${\cal X}$. 
Set $M_{\cal X}({\cal D}):=M_{\cal X}\oplus_{{\cal O}^*_{\cal X}}M({\cal D})$ 
and consider a natural morphism 
$M_{\cal X}({\cal D})\lo {\cal O}_{\cal X}$ of sheaves of 
commutative monoids on ${\cal X}_{\rm zar}$. 
Then we call the log scheme $({\cal X},M_{\cal X}({\cal D}))$ 
a {\it  proper strict semistable family with horizontal SNCD} over ${\cal V}$. 
We call ${\cal U}:={\cal X}\setminus{\cal D}$ 
an {\it open semistable log scheme} over ${\cal V}$. 
\end{defi} 

\par  
Set ${\cal V}_{-1}:={\cal V}$ and 
$\os{\circ}{S}_{-1}:={\rm Spec}({\cal V}_{-1})$. 

\begin{defi} 
Let $N$ be a nonnegative integer.  
Let ${\cal U}_{\bul \leq N}$ 
be an $N$-truncated simplicial fine log scheme over ${\rm Spec}({\cal V})$. 
If the generic fiber ${\cal U}_{\bul \leq N,K}$ of ${\cal U}_{\bul \leq N}$ 
is an $N$-truncated proper hypercovering 
of ${\cal U}_K$, then we call ${\cal U}_{\bul \leq N}$ 
an {\it $N$-truncated simplicial generically proper hypercovering} 
of ${\cal U}$ over ${\cal V}$. 
\end{defi}

First we have to generalize (\ref{prop:vv0}) as follows: 

\begin{prop}\label{prop:vvg0} 
There exists a sequence  
\begin{equation*} 
\cdots  \supset {\cal V}_N  \supset {\cal V}_{N-1}
\supset \cdots  \supset {\cal V}_0   \supset {\cal V}_{-1}
\tag{7.1.3.1}\label{eqn:nasv}
\end{equation*} 
of finite extensions of complete discrete valuation rings of ${\cal V}$ in $\ol{K}$
and a proper strict semistable family 
${\cal N}_m$ over ${\cal V}_m$ with a relative horizontal SNCD 
${\cal D}_m$ 
$(m\in {\mab N})$ 
and a log smooth 
$m$-truncated simplicial log scheme  
${\cal X}(m)_{\bul \leq m}$  
such that 
\begin{align*} 
{\cal X}(m)_{m'}=
(\os{\circ}{\cal X}(m)_{m'},M_{{\cal X}(m)_{m'}})
=\coprod_{0\leq l \leq m'}
\coprod_{[m'] \twoheadrightarrow [l]}
(({\cal N}_l,M_{{\cal N}_l}({\cal D}_l))\times_{S_l}S_m)
\end{align*} 
for each $0\leq m' \leq m$, where $S_m$ is the spectrum 
${\rm Spec}({\cal V}_m)$ with canonical log structure, 
and such that 
${\cal U}(m)_{\bul \leq m}$ is an $m$-truncated 
generically proper hypercovering of 
${\cal U}_{{\cal V}_m}$ over ${\cal V}_m$.  
Here 
$${\cal U}(m)_{\bul \leq m}:=
\coprod_{0\leq l \leq m'}
\coprod_{[m'] \twoheadrightarrow [l]}({\cal N}_l\setminus {\cal D}_l)\times_{S_l}S_m.$$
\end{prop}
\begin{proof}
By de Jong's theorem (\cite[(6.5)]{dj}) 
about the semistable reduction theorem 
by using the base change by an alteration 
(if one makes a finite extension 
of a complete discrete valuation ring)
and by a standard argument in 
\cite[${\rm V}^{\rm bis}$ \S5]{sga4-2} 
and \cite[(6.2.1.1)]{dh3}, 
we obtain (\ref{prop:vv0}). 
\end{proof}  
\parno 
The sequence (\ref{eqn:nasv}) 
gives us the following sequence of log schemes: 
\begin{equation*} 
\cdots \lo S_N \lo S_{N-1} \lo \cdots  \lo S_0\lo S_{-1}, 
\tag{7.1.3.2}\label{eqn:wtsnn}
\end{equation*} 
where $S_m$ is the log scheme whose underlying scheme is 
${\rm Spec}({\cal V}_m)$ and 
whose log structure is the canonical log structure 
on ${\rm Spec}({\cal V}_m)$. 
Let $\kap_m$ be the residue field of ${\cal V}_m$. 
Set $s_m:=S_m\times_{{\rm Spec}({\cal V}_m)}
{\rm Spec}(\kap_m)$.  
Then we have the following sequence 
\begin{equation*} 
\cdots \lo s_N \lo  s_{N-1}\lo \cdots \lo 
s_0\lo s_{-1}. 
\tag{7.1.3.3}\label{eqn:wtslpnna}
\end{equation*}
We also have the following sequence 
\begin{equation*} 
\cdots \lo {\cal W}_{\star}(s_N) \lo 
{\cal W}_{\star}(s_{N-1}) 
\lo \cdots \lo {\cal W}_{\star}(s_0) 
\lo 
{\cal W}_{\star}(s_{-1}) \quad 
({\star}\in {\mab Z}_{\geq 1}~{\rm   or }~\star
={\rm nothing}).   
\tag{7.1.3.4}\label{eqn:wtswpsnn}
\end{equation*}
Let $N$ be a nonnegative integer. 
Set ${\cal X}_{\bul \leq N}:={\cal X}(N)_{\bul \leq N}$ and $s:=s_N$. 
Set also $(N_m,M_{N_m})
:=({\cal N}_m,M_{{\cal N}_m})\otimes_{{\cal V}_m}\kap_m$ $(0\leq m\leq N)$. 
Set 
$$X(m)_{m'}
=\coprod_{0\leq l \leq m'}
\coprod_{[m'] \twoheadrightarrow [l]}(N_l,M_{N_l})\times_{s_l}s_m$$ 
for each $0\leq m' \leq m$. 
Then $X(m)_{\bul \leq m}={\cal X}(m)_{\bul \leq m}\times_{S_m}s_m$; 
we have the family 
$\{X(m)_{\bul \leq m}\}_{m=0}^N$ of 
log smooth split truncated simplicial log schemes of Cartier type. 
Analogously we obtain the family 
$\{U(m)_{\bul \leq m}\}_{m=0}^N$. 
Consequently, for each $N\in {\mab N}$, 
we obtain a log smooth split $N$-truncated log scheme 
$X_{\bul \leq N}:=X(N)_{\bul \leq N}$ of Cartier type 
and 
$U_{\bul \leq N}:=U(N)_{\bul \leq N}$ over $s:=s_N$. 

\begin{defi}\label{defi:gdef}
We say that the $N$-truncated simplicial generically proper hypercovering 
${\cal U}_{\bul \leq N}$ of ${\cal U}_{{\cal V}_N}$ is 
{\it gs$($=good and split$)$}.   
We say also that ${\cal X}_{\bul \leq N}$ is 
a {\it gs $N$-truncated simplicial generically proper hypercovering} 
$\ol{\cal U}_{{\cal V}_N}$.  
\end{defi}


The following is a generalization of (\ref{prop:coif}): 

\begin{prop}\label{prop:cogif} 
Let $N$ be a nonnegative integer. 
Then the following hold$:$
\par
$(1)$ Two gs $N$-truncated simplicial generically proper hypercoverings of 
the base change of $\ol{\cal U}$ over an extension of ${\cal V}$ 
are covered by a gs $N$-truncated simplicial
generically proper hypercovering of the base change of 
$\ol{\cal U}$ over an extension of ${\cal V}$. 
\par
$(2)$   
For a morphism $\ol{\cal U}{}'\lo \ol{\cal U}$ of proper schemes over ${\cal V}$ 
and for a gs $N$-truncated simplicial generically proper hypercovering 
${\cal X}_{\bul \leq N}$ of the base change of  $\ol{\cal U}$ over 
an extension of ${\cal V}$, 
there exist a finite extension ${\cal V}'$ of ${\cal V}$ 
and a gs $N$-truncated simplicial generically proper hypercovering 
${\cal X}'_{\bul \leq N}$ of $\ol{\cal U}{}'_{{\cal V}'}$ 
and a morphism 
${\cal X}'_{\bul \leq N} \lo {\cal X}_{\bul \leq N,{\cal V}'}$ 
fitting into the following commutative diagram$:$
\begin{equation*}
\begin{CD}
{\cal X}'_{\bul \leq N} @>>>  
{\cal X}_{\bul \leq N,{\cal V}'}\\
@VVV @VVV  \\
\ol{\cal U}{}'_{{\cal V}'}@>>> \ol{\cal U}_{{\cal V}'}.
\end{CD}
\tag{7.1.5.1}\label{cd:smpbc}
\end{equation*}
\end{prop}
\begin{proof}
As in the proof \cite[(9,4)]{nh2}, 
by using Nagata's embedding theorem and  
by using a general formalism 
in \cite[${\rm V}^{\rm bis}$ \S5]{sga4-2} and 
de Jong's semistable reduction theorem (\cite[(6.5)]{dj}), 
we obtain (\ref{prop:cogif}).   
\end{proof}

\begin{prop}\label{prop:ci}
Let ${\cal U}' \lo {\cal U}$ be a morphism of 
separated flat schemes of finite type over a morphism 
${\rm Spec}({\cal V}')\lo {\rm Spec}({\cal V})$ 
of the spectrums of complete discrete valuation rings 
with perfect residue fields of mixed characteristics.  
Let ${\cal U}\os{\sus}{\lo} \ol{\cal U}$ be an immersion 
into a proper scheme over ${\rm Spec}({\cal V})$.  
Then there exists an immersion 
${\cal U}'\os{\sus}{\lo} \ol{\cal U}{}'$ fitting into 
the following commutative diagram 
\begin{equation*} 
\begin{CD} 
{\cal U}'@>{\sus}>> \ol{\cal U}{}' \\
@VVV @VVV \\
{\cal U}@>{\sus}>> \ol{\cal U}
\end{CD}
\end{equation*} 
over ${\rm Spec}({\cal V}')\lo {\rm Spec}({\cal V})$.  
\end{prop} 
\begin{proof} 
Let ${\cal U}'\os{\sus}{\lo} \ol{\cal U}{}'$ 
be an immersion into a proper scheme over ${\rm Spec}({\cal V}')$. 
By the same argument as that of (\ref{rema:uuc}), 
we obtain (\ref{prop:ci}). 
\end{proof}

\par 
Next we have to generalize (\ref{eqn:fkci}). 
\par 
Let the notations be as in (\ref{eqn:fkci}) and as in the beginning of this section. 
Let ${\cal X}_{\bul \leq N}$ be as in (\ref{defi:gdef}). 
Let $X_{\bul \leq N}$ (resp.~${\mathfrak X}_{\bul \leq N}$) 
be the log special fiber (resp.~the log generic fiber) of ${\cal X}_{\bul \leq N}$. 
Set ${\mathfrak U}_{K_N}:={\mathfrak U}\otimes_{K}K_N$. 
Then, by (\ref{ali:rgrdhk}), (\ref{eqn:olyzk})  and (\ref{prop:hnn}), 
we have the following composite isomorphism 
\begin{align*} 
&H^q_{{\rm crys}}(X_{\bul \leq N}/{\cal W}(s))
\otimes_{\cal W}K_N
\os{H^q(\Psi),~\sim}{\lo} 
H^q_{{\rm log}{\textrm -}{\rm dR}}({\mathfrak X}_{\bul \leq N}/K_N) 
\tag{7.1.6.1}\label{eqn:fkgci}\\
& =
H^q({\mathfrak X}_{\bul \leq N},
\Om^{\bul}_{{\mathfrak X}_{\bul \leq N}/K}(\log {\mathfrak D}_{\bul \leq N}))
\os{\sim}{\lo} H^q({\mathfrak U}_{\bul \leq N},\Om^{\bul}_{{\mathfrak U}_{\bul \leq N}/K})\\
&=H^q_{\rm inf}({\mathfrak U}_{\bul \leq N}/K_N)
\os{\sim}{\longleftarrow}
H^q_{\rm inf}({\mathfrak U}_{K_N}/K_N). 
\end{align*} 
We denote this isomorphism by $H^q(\Psi)$. 
The isomorphism $H^q(\Psi)$ induces the following isomorphism 
\begin{align*} 
H^q(\Psi) \col 
H^q_{{\rm crys}}(X_{\bul \leq N}/{\cal W}(s))\otimes_{{\cal W}}K_N
\os{\sim}{\lo} 
H^q_{\rm inf}({\mathfrak U}_{K_N}/K_N), 
\tag{7.1.6.2}\label{eqn:fkygsci}
\end{align*} 
which we denote by $H^q(\Psi)$ again. 
Because the left hand side of (\ref{eqn:fkygsci}) has the monodromy operator, 
we have the monodromy operator  
$N_{K_N/K}({\cal X}_{\bul \leq N}/\ol{\cal U}/S/S_{-1})$ 
on $H^q_{\rm inf}({\mathfrak U}_{K_N}/K_N)$:  
\begin{align*}
N_{K_N/K}({\cal X}_{\bul \leq N}/{\cal U}/S/S_{-1})
\col H^q_{\rm inf}({\mathfrak U}_{K_N}/K_N)
\lo H^q_{\rm inf}({\mathfrak U}_{K_N}/K_N). 
\tag{7.1.6.3}\label{ali:xsakgnif}
\end{align*} 
The morphism (\ref{ali:xsakgnif}) induces the following morphism 
\begin{align*}
N_{\ol{K}/K}({\cal X}_{\bul \leq N}/\ol{\cal U}/S/S_{-1})
\col H^q_{\rm inf}({\mathfrak U}_{\ol{K}}/\ol{K})
\lo H^q_{\rm inf}({\mathfrak U}_{\ol{K}}/\ol{K}). 
\tag{7.1.6.4}\label{ali:xskngoif}
\end{align*} 

\begin{rema}\label{rema:ukll}
In the future, we would like to discuss the limit weight spectral sequence 
of $H^q_{\rm inf}({\mathfrak U}_{K_N}/K_N)$ by using (\ref{eqn:fkygsci}) 
and discuss fundamental properties of the induced filtration 
by this spectral sequence. 
\end{rema}

\begin{prop}\label{prop:gegmid}  
Let the notations be as above and {\rm (\ref{prop:ci})}.  
Set ${\mathfrak U}':={\cal U}'\otimes_{{\cal V}'}K'$ and 
$\ol{\mathfrak U}{}':=\ol{\cal U}{}'\otimes_{{\cal V}'}K'$.  
Let ${\mathfrak g}\col {\mathfrak U}'\lo {\mathfrak U}$ 
be a morphism over ${\rm Spec}(K')\lo {\rm Spec}(K)$. 
Let $\ol{\mathfrak g}:=\ol{\mathfrak U}{}'\lo \ol{\mathfrak U}$ be a morphism 
over  ${\rm Spec}(K')\lo {\rm Spec}(K)$ which is an extension of ${\mathfrak g}$. 
Let $\ol{g} \col \ol{\cal U}{}'\lo \ol{\cal U}$ and $g \col {\cal U}'\lo {\cal U}$ be morphisms  
over ${\rm Spec}({\cal V}')\lo {\rm Spec}({\cal V})$ 
which are extensions of the morphisms $\ol{\mathfrak g}$ and ${\mathfrak g}$, respectively. 
Let $S'_{-1}$ be the log scheme ${\rm Spec}({\cal V}')$ with 
the canonical log structure. 
Let ${\cal V}'_N$ and ${\cal Y}_{\bul \leq N}/S'$ 
be a similar complete discrete valuation ring to ${\cal V}_N$
and a similar $N$-truncated simplicial log scheme to ${\cal X}_{\bul \leq N}/S$, 
respectively, fitting into the following commutative diagram$:$  
\begin{equation*} 
\begin{CD} 
{\cal Y}_{\bul \leq N} 
@>{g}>> {\cal X}_{\bul \leq N}\\
@VVV @VVV \\
\ol{\cal U}{}'\times_{\os{\circ}{S}{}'_{-1}}\os{\circ}{S}{}' 
@>>> \ol{\cal U}\times_{\os{\circ}{S}_{-1}}\os{\circ}{S}\\
@VVV @VVV \\
\os{\circ}{S}{}'@>>> \os{\circ}{S}, 
\end{CD}
\tag{7.1.8.1}\label{ali:xxpssa}
\end{equation*} 
where ${\cal Y}_{\bul \leq N}$ is a gs $N$-truncated simplicial 
generically proper hypercovering  of $\ol{\cal U}{}'_{{\cal V}'_N}$. 
Set $K'_N:={\rm Frac}({\cal V}'_N)$. 
Then the following diagram is commutative$:$ 
\begin{equation*} 
\begin{CD} 
H^q_{\rm inf}({\mathfrak U}'_{K'_N}/K'_N) 
@>{e_{{\cal V}'/{\cal V}}\cdot 
N_{K'_N/K'}({\cal Y}_{\bul \leq N}/{\cal Y}/S'/S'_{-1})}>>
H^q_{\rm inf}({\mathfrak U}'_{K'_N}/K'_N)\\
@A{{\mathfrak g}^*}AA @AA{{\mathfrak g}^*}A \\
H^q_{\rm inf}({\mathfrak U}_{K_N}/K_N)
@>{N_{K_N/K}({\cal X}_{\bul \leq N}/{\cal X}/S/S_{-1})}>>
H^q_{\rm inf}({\mathfrak U}_{K_N}/K_N). 
\end{CD}
\tag{7.1.8.2}\label{cd:lkgnl}
\end{equation*} 
\end{prop}
\begin{proof} 
By using (\ref{eqn:fkgci}), 
the proof is the same as that of (\ref{prop:gemid}). 
\end{proof}

\begin{theo}\label{theo:wdngmo}  
The monodromy operator 
$N_{K_N/K}({\cal X}_{\bul \leq N}/{\cal X}/S/S_{-1})$ 
in {\rm (\ref{ali:xsakgnif})} is independent of the choices of 
${\cal X}$, 
$N\geq 2^{-1}(q+1)(q+2)$ and 
the gs $N$-truncated simplicial generically proper hypercovering  
${\cal X}_{\bul \leq N}$ of ${\cal X}_{{\cal V}_N}$. 
It depends only on ${\mathfrak X}_{K_N}/K_N$ and $K$. 
The monodromy operator also 
$N_{\ol{K}/K}({\cal X}_{\bul \leq N}/{\cal X}/S/S_{-1})$ in 
{\rm (\ref{ali:xskngoif})} depends only on 
${\mathfrak U}_{\ol{K}}/\ol{K}$ and $K$.  
\end{theo} 
\begin{proof} 
By using (\ref{eqn:fkgci}), 
the proof is the same as that of (\ref{theo:wdnmo}).  
\end{proof}

\begin{defi}\label{defi:mogp}
We call  the following well-defined operator
\begin{align*}
N_{{\mathfrak U}_{K_N}/K_N}
:= N_{K_N/K}({\cal X}_{\bul \leq N}/{\cal X}/S/S_{-1})
\col &H^q_{\rm inf}({\mathfrak U}_{K_N}/K_N)
\tag{7.1.10.1}\label{ali:xsxgkoif}\\
&\lo H^q_{\rm inf}({\mathfrak U}_{K_N}/K_N). 
\end{align*}  
the {\it $p$-adic monodromy operator} 
of ${\mathfrak U}_{K_N}/K_N$. 
We also call the following well-defined operator 
\begin{align*}
N_{{\mathfrak U}_{\ol{K}}/\ol{K}}
:= N_{{\mathfrak X}_{\ol{K}}/K}({\cal X}_{\bul \leq N}/{\cal X}/S/S_{-1})
\col H^q_{\rm inf}({\mathfrak U}_{\ol{K}}/\ol{K})
\lo H^q_{\rm inf}({\mathfrak U}_{\ol{K}}/\ol{K}). 
\tag{7.1.10.2}\label{ali:xskgnkkoif}
\end{align*}  
the {\it $p$-adic monodromy operator} 
of ${\mathfrak U}_{\ol{K}}/\ol{K}$. 
\end{defi}


\begin{theo}\label{theo:cwga}
Let ${\mathfrak U}$ be a separated scheme of finite type over $\ol{K}$. 
There exists a unique semi-linear action 
\begin{align*} 
\rho^q_{\rm crys}:=
\rho^q_{{\mathfrak U},{\rm crys}}\col 
{\rm WD}_{\rm crys}(\ol{K})\lo 
{\rm GL}^{\rm sl}(H^q_{\rm inf}({\mathfrak U}/\ol{K}))
\end{align*}  
satisfying the following properties$:$ 
\par 
$(1)$ This action is contravariantly functorial with respect to a morphism
${\mathfrak U}\lo {\mathfrak V}$ of separated schemes of finite type over $\ol{K}$. 
That is, for a morphism 
${\mathfrak f}\col {\mathfrak U}\lo {\mathfrak V}$ of 
separated schemes of finite type over $\ol{K}$, 
the following diagram is commutative for an element 
$\gam\in {\rm WD}_{\rm crys}(\ol{K})\!:$
\begin{equation*} 
\begin{CD}
H^q_{\rm inf}({\mathfrak U}/\ol{K})
@>\rho^q_{{\mathfrak U},{\rm crys}}(\gam)>> 
H^q_{\rm inf}({\mathfrak U}/\ol{K})\\
@A{{\mathfrak f}^*}AA @AA{{\mathfrak f}^*}A \\
H^q_{\rm inf}({\mathfrak V}/\ol{K})
@>\rho^q_{{\mathfrak V},{\rm crys}}(\gam)>> 
H^q_{\rm inf}({\mathfrak V}/\ol{K}). 
\end{CD}
\tag{7.1.11.1}\label{cd:xqgr}
\end{equation*} 
\par 
$(2)$ For a finite extension of $K'$ of $K$, 
let $K'{}^{\rm nr}$ be the maximal non-ramified extension of $K'$ in $\ol{K}$. 
If ${\cal U}$ is a separated scheme of finite type over ${\cal V}$ such that 
${\cal U}\otimes_{\cal V}\ol{K}={\mathfrak U}$, 
if ${\cal U}\os{\sus}{\lo} \ol{\cal U}$ is an open immersion into a proper scheme over 
${\cal V}$ and if ${\cal X}_{\bul \leq N}$ is a  
gs $N$-truncated simplicial generically proper hypercovering 
of $\ol{\cal U}_{{\cal V}_N}$ for some finite extension ${\cal V}_N$ over ${\cal V}$ 
of complete discrete valuation rings of mixed characteristics 
and if $X_{\bul \leq N}$ is the log special fiber of ${\cal X}_{\bul \leq N}$, 
then the image of the following morphism 
\begin{align*} 
&H^q_{\rm crys}(X_{\bul \leq N}/{\cal W}(s))\otimes_{{\cal W}}
K^{\rm nr}_0 \subset 
H^q_{\rm crys}(X_{\bul \leq N}/{\cal W}(s))\otimes_{{\cal W}}(K_N)^{\rm nr}
\tag{7.1.11.2}\label{ali:xwsgbn}\\
& \os{H^q(\Psi),~\sim}{\lo} H^q_{\rm inf}({\cal U}_{(K_N)^{\rm nr}}/(K_N)^{\rm nr})
\subset H^q_{\rm inf}({\mathfrak U}/\ol{K}) 
\end{align*}
is equal to the invariant space 
$H^q_{\rm inf}({\mathfrak U}/\ol{K})^{{\rm I}_{\rm crys}(\ol{K})};$ 
the action of $W_{\rm crys}(\ol{K})/
I_{\rm crys}(\ol{K})\simeq {\mab Z}$ 
is induced by the action of the absolute Frobenius endomorphism of 
$X_{\bul \leq N};$ 
the action of $\ol{K}(1)\owns a$ $(a\in \ol{K})$ is defined by 
the exponential ${\rm exp}(ae^{-1}_NN_{\rm zar})$ of the monodromy operator 
$ae^{-1}_NN_{\rm zar}\col 
H^q_{\rm crys}(X_{\bul \leq N}/{\cal W}(s))\otimes_{\cal W}\ol{K}
\lo 
H^q_{\rm crys}(X_{\bul \leq N}/{\cal W}(s))(-1)\otimes_{\cal W}\ol{K}$, 
where $e_N$ is the absolute ramification index of $K_N$. 
More precisely speaking, in the situation above, the action of $\rho^q_{\rm crys}((a,g))$ 
$(a\in \ol{K}(1),g\in W_{\rm crys}(\ol{K}))$ is defined by the following formula$:$
\begin{align*} 
\rho^q_{\rm crys}((a,g))(H^q(\Psi)(v)):=
H^q(\Psi)({\rm exp}(ae^{-1}_NN_{\rm zar})(F^*_{X_{\bul \leq N}})^{{\rm deg}(g)}(v)) 
\end{align*} 
for $v\in H^q_{\rm crys}(X_{\bul \leq N}/{\cal W}(s))\otimes_{{\cal W}}K^{\rm nr}_0$. 
\end{theo}
\begin{proof}
By using (\ref{prop:ci}), (\ref{prop:gegmid}) and (\ref{eqn:fkygsci}), 
the proof is the same as that of (\ref{theo:cwa}). 
\end{proof}

\begin{prop}\label{prop:cauccw}
The representation is compatible with the cup product of 
$H^{\bul}_{\rm inf}({\mathfrak U}/\ol{K})$. 
That is, for $x\in H^q_{\rm inf}({\mathfrak U}/\ol{K})$, 
$y\in H^{q'}_{\rm inf}({\mathfrak U}/\ol{K})$ and 
$\gam \in {\rm WD}_{\rm crys}(\ol{K})$, 
\begin{align*} 
\rho^{q+q'}_{{\mathfrak U},{\rm crys}}(\gam)(x\cup y)
=
\rho^q_{{\mathfrak U},{\rm crys}}(\gam)(x)\cup 
\rho^{q'}_{{\mathfrak U},{\rm crys}}(\gam)(y).
\tag{7.1.12.1}\label{ali:xxuxy}
\end{align*}  
\end{prop} 
\begin{proof} 
The proof is the same as that of (\ref{prop:caccw}). 
\end{proof}

\begin{theo}\label{theo:k0gst}
$(1)$ There exists a well-defined $K_0$-form 
$H^q_{\rm inf}({\mathfrak U}/K/K_0)$ 
of $H^q_{\rm inf}({\mathfrak U}/K)$. 
\par 
$(2)$ 
Let $K\os{\subset}{\lo} K'$ be a morphism of 
the fraction fields of complete discrete valuation rings 
of mixed characteristics $(0,p)$. 
Let ${\mathfrak V}$ be a separated scheme of finite type over $K'$. 
Let $f\col {\mathfrak V}\lo {\mathfrak U}$ be a morphism over 
${\rm  Spec}(K')\lo  {\rm  Spec}(K)$. Then there exists a morphism 
\begin{align*} 
f^*_0\col H^q_{\rm inf}({\mathfrak U}/K/K_0)\otimes_{K_0}K'_0
\lo 
H^q_{\rm inf}({\mathfrak V}/K'/K'_0)
\tag{7.1.13.1}\label{ali:yxuvlk} 
\end{align*} 
such that $f^*_0\otimes_{K'_0}K'=f^*$. 
This morphism satisfies the transitive relation and 
${\rm id}^*_{\mathfrak X}=
{\rm id}_{H^q_{\rm inf}({\mathfrak X}/K/K_0)}$.  
\end{theo}
\begin{proof} 
(1): By using (\ref{eqn:fkygsci}), 
the proof is the same as that of (\ref{theo:k0st}). 
\par 
(2): By using (\ref{eqn:fkygsci}), the proof is the same  as that of (\ref{theo:stinfgb}). 
\end{proof} 

\begin{defi}\label{defi:lmmmta}
We call $H^q_{\rm inf}({\mathfrak U}/K/K_0)$ 
the {\it hidden} $K_0$-{\it structure} 
of $H^q_{\rm inf}({\mathfrak U}/K)$. 
\end{defi} 

\begin{rema}
By the definition of $H^q_{\rm inf}({\mathfrak U}/K/K_0)$, 
we have the following natural morphism 
\begin{align*} 
H^q(\Psi) \col H^q_{{\rm crys}}(X_{\bul \leq N}/{\cal W}(s))\otimes_{{\cal W}}K_{N,0}
\os{\sim}{\lo} 
H^q_{\rm inf}({\mathfrak U}/K/K_0)\otimes_{K_0}K_{N,0} 
\quad (k\in {\mab Z}). 
\tag{7.1.15.1}\label{ali:ppgkk}
\end{align*} 
\end{rema}

\begin{theo}\label{theo:lmgon}
Let the notations be as in the previous section. 
Then the following hold$:$ 
\par 
$(1)$ There exists a 
well-defined $p$-adic monodromy operator 
\begin{align*} 
N_{{\mathfrak U}/K/K_0}\col H^q_{\rm inf}({\mathfrak U}/K/K_0)
\lo H^q_{\rm inf}({\mathfrak U}/K/K_0)
\end{align*} 
such that 
$N_{{\mathfrak U}/K/K_0}\otimes_{K_0}K
=N_{{\mathfrak U}/K} \col H^q_{\rm inf}({\mathfrak U}/K)\lo 
H^q_{\rm inf}({\mathfrak U}/K)$. 
\par 
$(2)$ There exists an action $F_{/K_0}$ on 
$H^q_{\rm inf}({\mathfrak U}/K/K_0)$ 
such that 
$F_{/K_0}\otimes_{K_0}K_{N,0}=F$ on the 
$K_{N,0}$-form 
$${\rm Im}(H^q(\Psi) \vert_{H^q_{\rm crys}(X_{\bul \leq N}/{\cal W}(s))
\otimes_{{\cal W}}K_{N,0}}).$$ 
\end{theo}
\begin{proof} 
The proof is the same as that of (\ref{theo:lmon}). 
\end{proof}

\begin{prop}\label{prop:cugpp}
$(1)$ There exists the following product 
\begin{align*} 
(?\cup ?)_{K_0}\col 
H^q_{\rm inf}({\mathfrak U}/K/K_0)\otimes_{K_0}H^q_{\rm inf}({\mathfrak U}/K/K_0)
\lo H^q_{\rm inf}({\mathfrak U}/K/K_0)
\end{align*}
such that $(?\cup ?)_{K_0}\otimes_{K_0}K$ is the cup product 
of $H^q_{\rm inf}({\mathfrak U}/K)$.  
\par 
$(2)$  Set $N_{{\mathfrak U}/K/K_0}^{[i]}:=(i!)^{-1}N_{{\mathfrak U}/K/K_0}$ $(i\in {\mab N})$. 
Then the Leibniz rule for $N_{{\mathfrak U}/K/K_0}^{[i]}$ holds: 
\begin{align*} 
N^{[i]}_{{\mathfrak U}/K/K_0}((x\cup y)_{K_0})&=
\sum_{j=0}^iN^{[j]}_{{\mathfrak U}/K/K_0}(x)\cup N^{[i-j]}_{{\mathfrak U}/K/K_0}(y) \\
&(x\in H^{q'}_{\rm inf}({\mathfrak U}/K/K_0),~y\in H^{q''}_{\rm inf}({\mathfrak U}/K/K_0)).
\end{align*}  
\end{prop}
\begin{proof} 
The proof is the same as that of (\ref{prop:cugpp}). 
\end{proof} 


Let 
\begin{align*} 
\rho^q_{\rm inf}\col 
{\rm WD}(\ol{K}/K)\lo {\rm GL}^{\rm sl}(H^q_{\rm inf}({\mathfrak U}_{\ol{K}}/\ol{K}))
\end{align*} 
be a representation obtained by the natural action of 
${\rm Gal}(\ol{K}/K)$ on $H^q_{\rm inf}({\mathfrak U}_{\ol{K}}/\ol{K})
=H^q_{\rm inf}({\mathfrak U}/K)\otimes_K\ol{K}$. 
Set $\rho^q_p(\gam):=\rho^q_{\rm inf}(\gam)\rho^q_{\rm crys}(\gam^{-1})$ 
$(\gam\in {\rm WD}(\ol{K}/K))$. 
Then we have the following $\ol{K}$-linear action 
\begin{align*} 
\rho_p \col 
{\rm WD}(\ol{K}/K)\lo {\rm GL}_{\ol{K}}
(H^q_{\rm inf}({\mathfrak U}_{\ol{K}}/\ol{K}))
\tag{7.1.17.1}\label{ali:rhopa}
\end{align*} 
of ${\rm WD}(\ol{K}/K)$.

\begin{conj}[{\bf cf.~(\ref{conj:gfc})}]\label{conj:gfic} 
Let ${\mathfrak U}$ be a separated scheme of finite type over $K$. 
Then the representations $\rho_l$ on 
$H^q_{\rm et}({\mathfrak U}_{\ol{K}},{\mab Q}_l)$ ((\ref{theo:qnil}))
and $\rho_p$ of ${\rm WD}(\ol{K}/K)$ on 
$H^q_{\rm inf}({\mathfrak U}_{\ol{K}}/\ol{K})$ 
are compatible.
\end{conj}


\bigskip
\bigskip
\parno
Yukiyoshi Nakkajima 
\parno
Department of Mathematics,
Tokyo Denki University,
5 Asahi-cho Senju Adachi-ku,
Tokyo 120--8551, Japan. 
\parno
{\it E-mail address\/}: 
nakayuki@cck.dendai.ac.jp


\begin{thebibliography} {EGA IV-3}


\bibitem[A]{a}Artin, M.
\newblock{\em Algebraic approximation 
of structures over complete local rings}.
\newblock  Publ.~Math.~IH\'{E}S 36 (1969), 23--58.


\bibitem[BBD]{bbd}Beilinson, A.~A., 
Bernstein, J., Deligne, P.
\newblock{\em Faisceaux pervers}.
\newblock Ast\'{e}risque 100 (1982).

\bibitem[Be1]{bwl}Berthelot, P.
\newblock{\em Sur le $\ll$th\'{e}or\`{e}me de 
Lefschetz faible$\gg$ en cohomologie cristalline}. 
C. R. Acad. Sci. Paris S\'{e}r. A-B277 (1973), A955--A958. 


\bibitem[Be2]{bb}Berthelot, P.                                            
\newblock{\em Cohomologie cristalline 
des sch\'{e}mas de caract\'{e}ristique $p>0$}.
\newblock Lecture Notes in Math.~407, 
Springer-Verlag (1974).


\bibitem[Be3]{brc}Berthelot, P. 
\newblock{\em G\'{e}om\'{e}trie 
rigide et cohomologie des vari\'{e}t\'{e}s 
alg\'{e}briques de caract\'{e}ristique $p$}. 
\newblock In: Introductions aux 
cohomologies $p$-adiques (Luminy, 1984).
M\'{e}m.~Soc.~Math.~France 23 (1986), 7--32.  

\bibitem[Be4]{bdmi}Berthelot, P. 
\newblock{\em $D$-modules arithm\'{e}tiques {\rm I}$:$ 
Op\'{e}rateurs diff\'{e}rentiels de niveau fini}. 
\newblock Ann.~Scient.~\'{E}c.~Norm.~Sup.~29 (1996), 185--272.

\bibitem[Be5]{bfi}Berthelot, P.
\newblock{\em Finitude et puret\'{e} 
cohomologique en cohomologie rigide,
avec un appendice par Aise Johan de Jong}.
\newblock Invent.~Math.~128 (1997),  329--377.

\bibitem[Be6]{bdmii}Berthelot, P. 
\newblock{\em $D$-modules arithm\'{e}tiques {\rm II}$:$ Descente par Frobenius}. 
\newblock M\'{e}m. ~Soc. ~Math. ~France 81 
(2000), 1--136.

\bibitem[Be7]{bdmitr}Berthelot, P. 
\newblock{\em Introduction \`{a} 
la th\'{e}orie arithm\'{e}tique des D-modules}. 
In: Cohomologies $p$-adiques et applications arithm\'{e}tiques II, 
\newblock Ast\'{e}risque 279, Soc.~Math.~de France 
(2002), 1--80.


\bibitem[Be8]{bpre}Berthelot, P. 
\newblock{\em Cohomologie rigide et  
Cohomologie rigide \`{a} supports propres}. 
\newblock Premi\`{e}re partie, available from 
Berthelot's homepage:  
http://name.math.univ-rennes1.fr/pierre.berthelot/.  


\bibitem[Be9]{blec}Berthelot, P.
\newblock{\em }
\newblock Three lectures ``Cohomology 
operators, ${\cal D}$-modules 
and crystalline sheaves'' in Dwork Trimester 2001. 
(The title ``Crystals and D-modules'' in 
``http://www.math.unipd.it/~dwork01/schedweek1.html''  
was changed.)       


\bibitem[BI]{beill}Berthelot, P., Illusie, L.
\newblock{\em Classes de Chern en cohomologie cristalline}. 
\newblock C.~R.~Acad.~Sci.~Paris, 
S\'{e}rie I 270 (1970), 1695--1697.

\bibitem[Bl]{bl}Bloch, S.
\newblock{\em Algebraic $K$-theory and crystalline cohomology}. 
\newblock IHES Publ.~Math.~47 (1977), 187--268. 





\bibitem[BO1]{bob}Berthelot, P., Ogus, A.  
\newblock{\em Notes on crystalline cohomology}. 
\newblock Princeton University Press, 
University of Tokyo Press (1978).

\bibitem[BO2]{boi}Berthelot, P., Ogus, A.                                           
\newblock{\em $F$-isocrystals and de Rham cohomology. {\rm I}}. 
\newblock Invent.~Math.~72 (1983), 159-199. 

\bibitem[Bou]{bou2}Bourbaki, N.                                           
\newblock{\em Alg\`{e}bre commutative}. 
\newblock  Masson, Paris (1985).


\bibitem[Ca1]{lcara}Caraiani, A.  
\newblock{\em 
Local-global compatibility and 
the action of monodromy on nearby cycles}. 
\newblock Duke Math.~J.~161 (2012), 2311--2413. 


\bibitem[Ca2]{pcara}Caraiani, A.                                           
\newblock{\em 
Monodromy and local-global compatibility for $l=p$}. 
\newblock Algebra and Number Theory 8 (2014), 1597--1646. 

\bibitem[CCPI]{ccpi}Chiarellotto, B., Coleman, R., 
Proietto, V.-D., Iovita, A. 
\newblock{\em On a p-adic invariant cycles theorem}. 
J.~reine angew.~Math.~711 (2016), 55--74.


\bibitem[CF]{cm}Chiarellotto, B., Fornasiero, M.
\newblock{\em  Logarithmic de Rham, infinitesimal and Betti cohomologies}.
\newblock  J.~Math.~Sci.~Univ.~Tokyo 13 (2006), 205--257.   


\bibitem[Ch1]{cric}Chiarellotto, B.
\newblock{\em  Rigid cohomology and invariant cycles 
for a semistable log scheme}. 
\newblock Duke Math.~J.~97 (1999), 155--169. 


\bibitem[CLa]{cla}Chiarellotto, Lazda, C.
\newblock{\em Around $l$-independence}. 
\newblock  Compositio Math.~154 (2018), 223--248




\bibitem[CLe]{clpu}Chiarellotto, B., Le Stum, B.
\newblock{\em Sur la puret\'{e} de la cohomologie cristalline}. 
\newblock C.~R.~Acad.~Sci.~Paris, 
S\'{e}rie I 326 (1998), 961--963.






\bibitem[Cr]{cr}Crew, R. 
\newblock{\em Finiteness theorems for the cohomology 
of an overconvergent isocrystal on a curve}. 
\newblock Ann.~Scient.~\'{E}c.~Norm.~Sup.~
$4^e$ s\'{e}rie 31 (1998), 717--763. 

\bibitem[CT1]{ctze}Chiarellotto, B., Tsuzuki, N.  
\newblock{\em Cohomological descent of 
rigid cohomology for \'{e}tale coverings}.
\newblock  Rend.~Sem.~Mat.~Univ.~Padova 109 (2003), 63--215. 

\bibitem[CT2]{ctcs}Chiarellotto, B., Tsuzuki, N.  
\newblock{\em Clemens-Schmid exact sequence in characteristic $p$}.  
\newblock  Math.~Ann.~358 (2014), 971--1004. 



\bibitem[D1]{dh1}Deligne, P.                                               
\newblock{\em Th\'{e}orie  de Hodge ${\rm I}$}.
\newblock Actes du Congr\`es International 
des Math\'ematiciens (Nice, 1970)
Tome 1, Gauthier-Villars, Paris (1971), 425--430. 

\bibitem[D2]{dh2}Deligne, P.                                              
\newblock{\em Th\'{e}orie  de Hodge, ${\rm II}$}.
\newblock Publ.~Math.~IH\'{E}S 40 (1971), 5--57.

\bibitem[D3]{dhlc}Deligne, P.  
\newblock{\em Les constantes des \'{e}quations 
fonctionnelles des fonctions L}. 
\newblock In: Modular functions of one variable, II,  
Lecture Notes in Math.~349, Springer-Verlag (1973), 501--597. 


\bibitem[D4]{dh3}Deligne, P.                                                  
\newblock{\em Th\'{e}orie  de Hodge, ${\rm III}$}.
\newblock Publ.~Math.~IH\'{E}S 44 (1974), 5--77.

\bibitem[D5]{dp}Deligne, P.        
\newblock{\em Poids dans la 
cohomologie des vari\'{e}t\'{e}s 
alg\'{e}briques}. 
\newblock Proceedings of the International Congress
of Mathematicians (Vancouver, B. C., 1974) Vol.~1, 79--85. 
                
\bibitem[D6]{dw2}Deligne, P.                                       
\newblock{\em La conjecture de Weil, ${\rm II}$}.
\newblock Publ.~Math.~IH\'{E}S 52 (1980), 137--252.
                                 
\bibitem[D7]{dhc}Deligne, P.                                 
\newblock{\em Hodge cycles on abelian varieties $($notes by J. S. Milne$)$}.  
\newblock Lecture Notes in Math.~900, Springer-Verlag (1982), 9--100.                                 



\bibitem[dB]{db}du Bois, P.
\newblock{\em  Structure de Hodge mixte 
sur la cohomologie \'{e}vanscente}.
\newblock Ann.~Inst.~Fourier (Grenoble) 35 (1985), 
191--213.
                                 
\bibitem[DI]{di}Deligne, P., Illusie, L.
\newblock{\em Rel\`{e}vements 
modulo $p^2$ et d\'{e}compositon du complexe de de Rham}. 
\newblock Invent.~Math.~89 (1987), 247--270.


\bibitem[dJ]{dj}de Jong, A.~J.
\newblock{\em  Smoothness, semi-stability and alterations}.
\newblock Publ.~Math.~IH\'{E}S 83 (1996), 51--93.

\bibitem[EGA III-1]{ega3}Grothendieck, A., Dieudonn\'e, J. 
\newblock{\em \'El\'ements de g\'eom\'etrie alg\'ebrique 
${\rm III}$-$1$}.
\newblock Publ.~Math.~IH\'{E}S 11 (1961). 

\bibitem[EGA IV-3]{ega43}Grothendieck, A., Dieudonn\'e, J. 
\newblock{\em \'El\'ements de g\'eom\'etrie alg\'ebrique 
${\rm IV}$-$3$}.
\newblock Publ.~Math.~IH\'{E}S 28 (1966). 

\bibitem[EGA IV-4]{ega4}Grothendieck, A., Dieudonn\'e, J. 
\newblock{\em \'El\'ements de g\'eom\'etrie alg\'ebrique 
${\rm IV}$-$4$}.
\newblock  Publ.~Math.~ IH\'{E}S 32 (1967). 


\bibitem[Ek]{ek}Ekedahl, T.                                           
\newblock{\em Duality for the de Rham-Witt 
complex ${\rm I}$}.
\newblock Ark.~f\"{o}r Matematik 22 (1984), 185--239.

\bibitem[Et]{et}Etesse, J.-Y. 
\newblock{\em Complexe de de Rham-Witt \`{a} 
coefficients dans un cristal.} 
\newblock Compositio Math.~66 (1988), 57--120. 

\bibitem[Fa]{fao}Faltings, G. 
\newblock{\em $F$-isocrystals on open varieties: 
Results and conjectures}.
\newblock The Grothendieck Festschrift Vol.~II,    
Progr.~Math.~87, Birkh\"{a}user  (1990), 219--248.


\bibitem[Fu]{fup}Fujisawa, T. 
\newblock{\em Polarizatons on limiting mixed Hodge structures}.
\newblock  J.~of Singularities 8 (2014), 146--193.

\bibitem[FN]{fn}Fujisawa, T., Nakayama, C.                                            
\newblock{\em Mixed Hodge structures on log deformations}.
\newblock  Rend.~Sem.~Mat.~Univ.~Padova 110 
(2003), 221--268.

\bibitem[Fo]{fos}Fontaine, J.-M. 
\newblock{\em Repr\'{e}sentations $l$-adiques 
potentiellement semi-stables}. 
\newblock In:  P\'eriodes $p$-adiques, Seminaire de Bures, 1988. 
Ast\'{e}risque 223, Soc.~Math.~de France 
(1994),  321--347.

\bibitem[Fr]{fr1}Friedman, R. 
\newblock{\em Global smoothings 
of varieties with normal crossings}.
\newblock   Ann.~of Math.~118 (1983), 75--114. 



\bibitem[Fu]{fuhs}Fujiwara, K.
\newblock{\em Etale topology and philosophy of log}.
\newblock 
The report of algebraic geometry symposium at Kinosaki (1990),
116--123, (in Japanese).

\bibitem[GK1]{gkwf}Gro{\ss}e-Kl\"onne, E.
\newblock{\em Frobenius and monodromy operators 
in rigid analysis, and Drinfeld's symmetric space}.
\newblock  Algebraic Geom.~14 (2005), 391--437. 



\bibitem[GK2]{gkc}Gro{\ss}e-Kl\"onne, E. 
\newblock{\em 
The \v{C}ech filtration and monodromy in log crystalline cohomology}. 
\newblock Transactions of the AMS~359 (2007),  2945--2972. 


\bibitem[Gr1]{grdr}Grothendieck, A. 
\newblock{\em On the de Rham cohomology of algebraic varieties}. 
\newblock Publ.~Math.~IH\'{E}S 29 (1966), 95--103. 	

\bibitem[Gr2]{grcr}Grothendieck, A. 
\newblock{\em Crystals and the de Rham cohomology of schemes}. 
\newblock In Dix Expos\'{e}s 
sur la Cohomologie des Sch\'{e}mas, 
North-Holland, Amsterdam; Masson, Paris (1968), 306--358 	

\bibitem[Gr3]{grrs}Grothendieck, A. 
\newblock{\em R\'{e}coltes et semailles$:$ 
R\'{e}flexions et t\'{e}moignage 
sur un pass\'{e} de math\'{e}maticien,  
{\rm Parties I, II, IV}}.
\newblock  Gendai-Sugaku-sha, 
Japanese translation by Y.~Tsuji (1989), (1990).


\bibitem[GL]{greglan}Gregory, O., Langer, A. 
\newblock{\em 
Overconvergent de Rham-Witt cohomology for semistable varieties}. 
\newblock Preprint, available from https://arxiv.org/pdf/1711.09943.pdf.
 
 
\bibitem[GN]{gn}Guill\'{e}n, F., Navarro Aznar, V.
\newblock{\em Sur le th\'{e}or\`{e}me local des cycles invariants}. 
\newblock  Duke Math.~J. 61 (1990), 133--155.  
 

\bibitem[Ha1]{hartsag}Hartshorne, R.
\newblock{\em Algebraic geometry}.
\newblock Graduate Texts in Math.~52, Springer-Verlag (1977).

\bibitem[Ha2]{hadr}Hartshorne, R. 
\newblock{\em On the De Rham cohomology of algebraic varieties}. 
Publ.~Math.~IH\'{E}S 45 (1975), 5--99.


\bibitem[HK]{hk}Hyodo, O., Kato, K. 
\newblock{\em Semi-stable reduction and 
crystalline cohomology with logarithmic poles}. 
\newblock In:  P\'eriodes $p$-adiques, 
Seminaire de Bures, 1988. 
Ast\'{e}risque 223, Soc.~Math.~de France 
(1994),  221--268.


\bibitem[Hy]{hdw}Hyodo, O.
\newblock{\em On the de Rham-Witt 
complex attached to a semi-stable family}. 
\newblock Comp.~Math.~78 (1991), 241--260.


\bibitem[I1]{isp}Illusie, L.                                            
\newblock{\em Report on crystalline  cohomology}.
\newblock  Proc.~of Symposia in Pure Math.~Vol.~XXIX (1975), 459--478.


\bibitem[I2]{idw}Illusie, L.                                            
\newblock{\em Complexe de de Rham-Witt 
et cohomologie cristalline}.
\newblock  
Ann.~Scient.~\'{E}c.~Norm.~Sup.~$4^e$ s\'{e}rie 12 (1979), 
501--661.



\bibitem[IR]{ir}Illusie, L., Raynaud, M.
\newblock{\em Les suites spectrales 
 associ\'ees au complexe de de Rham-Witt}.  
\newblock  Publ.~Math.~IH\'{E}S  57 (1983), 73--212. 


\bibitem[It1]{itp}Ito, T. 
\newblock{\em Weight-monodromy-weight conjecture 
over equal characteristic local fields.}
\newblock Amer.~J.~Math.~127 (2005), 647--658. 


\bibitem[It2]{itup}Ito, T. 
\newblock{\em Weight-monodromy conjecture for 
p-adically uniformized varieties}
\newblock  Invent.~Math.~159 (2005), 607--656. 
 
\bibitem[Kaj]{kaji}Kajiwara, T.                                            
\newblock{\em The monodromy-weight conjecture and 
the hard Lefschetz theorem for $H^1$ of projective log semistable varieties}.
\newblock Preprint.

\bibitem[Kf]{fkato}Kato, F.                                            
\newblock{\em Log smooth deformation theory}.
\newblock T$\hat{{\rm o}}$hoku Math.~Jour.~48 (1996), 317-354.


\bibitem[Kk1]{kcl2}Kato, K.
\newblock{\em A generalization of 
local class field theory by using 
$K$-groups.~${\rm II}$}. 
\newblock J.~Fac.~Sci.~Univ.~Tokyo 
Sect.~IA Math.~27 (1980), 603--683. 

\bibitem[Kk2]{klog1}Kato, K.
\newblock{\em Logarithmic structures 
of Fontaine-Illusie}.
\newblock In: Algebraic analysis, 
geometry, and number theory, 
Johns Hopkins Univ. Press (1989), 191--224.

\bibitem[Kk3]{klog2}Kato, K.
\newblock{\em Toric singularities}.
\newblock Amer.~J.~Math.~116 (1994),  1073--1099.


\bibitem[Kk4]{kln}Kato, K.                                           
\newblock{\em }
\newblock A letter to Y.~Nakkajima (in Japanese) (1997).


\bibitem[KKMS]{kkms}Kempf, G., Knudsen, F., Mumford, D., Saint-Donat, B.                                        
\newblock{\em Toroidal Embeddings 1}
\newblock Lecture Notes in Math.~339, Springer Verlag  (1973).


\bibitem[KH]{kiha}Kim, M., Hain, R.~M.
\newblock{\em A de Rham-Witt approach 
to crystalline rational homotopy theory}.
\newblock  Comp.~Math.~140 (2004), 1245--1276. 

\bibitem[Kl]{kl}Kleiman, S.~L.
\newblock{\em Algebraic cycles and the Weil conjectures}.
\newblock In:~Dix Espos\'{e} sur 
la cohomologie des sch\'{e}mas, 
North-Holland, Amsterdam Math.~Vol.~55 Part 1 (1968), 
359--386.



\bibitem[KM]{kme}Katz,  N., Messing, W.                                           
\newblock{\em Some consequences  of 
the Riemann hypothesis for varieties 
over finite fields}.
\newblock Invent.~Math.~23 (1974), 73--77. 


\bibitem[KO]{ko}Katz, N.~M., Oda, T.
\newblock{\em On the differentiation of 
De Rham cohomology classes with respect to parameters}. 
\newblock J.~Math.~Kyoto Univ.~ 8 (1968), 199--213.

\bibitem[KtN]{kn}Kato, K., Nakayama, C.
\newblock{\em Log Betti cohomology, 
log \'{e}tale cohomology, and 
log de Rham cohomology of log schemes over ${\mab C}$}.
\newblock Kodai Math.~J.~22 (1999), 161--186.



\bibitem[KwN]{kwn}Kawamata, Y., Namikawa, Y.
\newblock{\em  Logarithmic deformations 
of normal crossing varieties 
and smoothing of degenerate Calabi-Yau varieties}.
\newblock Invent.~Math.~118 (1994), 395--409. 


\bibitem[Laz]{la}Lazard, M.
\newblock{\em Commutative formal groups}. 
\newblock  Lecture Notes in Math.~443, 
Springer-Verlag (1975).





\bibitem[Le1]{blr}Le Stum, B. 
\newblock{\em Rigid cohomology}. 
\newblock Cambridge Tracts in Math.~172, 
Cambridge University Press (2007). 


\bibitem[Le2]{lesc}Le Stum, B. 
\newblock{\em The overconvergent site}. 
\newblock M\'{e}m.~Soc.~Math.~France 127 (2011). 

\bibitem[Lo]{lodw}Lorenzon, P.
\newblock{\em  Logarithmic 
Hodge-Witt forms and Hyodo-Kato cohomology}. 
\newblock  J.~Algebra 249 (2002), 247--265.

\bibitem[LP]{lp}Lazda, C., P\'{a}l, A.
\newblock{\em Rigid cohomology over Laurent series fields}. 
\newblock  Algebra and Applications 21, Springer-Verlag (2016).





\bibitem[Lu]{lu}L\"{u}tkebohmert, W. 
\newblock{\em  On compactification of schemes}. 
\newblock Manuscripta Math.~80 (1993), 95--111.

\bibitem[LZ]{lazi}Langer, A and Zink, T. 
\newblock{\em De Rham-Witt cohomology for a proper and smooth morphism}. 
\newblock J.~of the Inst.~of Math.~Jussieu~3 (2004), 231--314.



\bibitem[Mar]{marmo}Marmora, A.
\newblock{\em Facteurs epsilon $p$-adiques}. 
\newblock Compositio Math.~144 (2008),  439--483



\bibitem[MatH]{maoro}Matsuue, H.
\newblock{\em On relative and overconvergent 
de Rham-Witt cohomology for log schemes}. 
\newblock Math.~Zeit.~286 (2017), 19--87.




\bibitem[MatY]{magk3}Matsumoto, Y.
\newblock{\em Good reduction criterion for K3 surfaces}. 
\newblock Math.~Zeitschrift (2015), 241--266


\bibitem[Mi]{miepl}Mieda, Y.
\newblock{\em The Picard--Lefschetz formula for $p$-adic cohomology}.  
\newblock Math.~Zeitschrift 257 (2007), 403--425.


\bibitem[Mo]{msemi}Mokrane, A.
\newblock{\em La suite spectrale 
des poids en cohomologie de Hyodo-Kato}.
\newblock Duke Math.~J.~72 (1993), 301--336. 


\bibitem[Nag]{imno}Nagata, M. 
\newblock{\em A generalization of the imbedding problem 
of an abstract variety in a complete variety}. 
\newblock J.~Math.~Kyoto Univ.~3 (1963) 89--102. 

\bibitem[Nak]{nd}Nakayama, C.
\newblock{\em Degeneration of 
$l$--adic weight spectral sequences}.
\newblock Amer.~J.~Math.~122 (2000), 721--733.



\bibitem[Nakk1]{nlk3}Nakkajima, Y.
\newblock{\em Liftings of log $K3$ surfaces 
and classical log Enriques surfaces in mixed characteristics}.
\newblock J.~of Algebraic Geometry 9 (2000), 355--393.


\bibitem[Nakk2]{ndw}Nakkajima, Y. 
\newblock{\em $p$-adic weight spectral sequences of log varieties}.
\newblock  J.~Math.~Sci.~Univ.~Tokyo 12 (2005), 513--661. 


\bibitem[Nakk3]{nlpi}Nakkajima, Y.                 
\newblock{\em Signs in weight spectral sequences, 
monodromy--weight conjectures, 
log Hodge symmetry and degenerations of surfaces}.
\newblock Rend.~Sem.~Mat.~Univ.~Padova 116 
(2006), 71--185. 




\bibitem[Nakk4]{nh3}Nakkajima, Y.
\newblock{\em Weight filtration 
and slope filtration on the rigid cohomology of 
a variety in characteristic $p>0$}.
\newblock M\'{e}m.~Soc.~Math.~France 130--131 (2012).


\bibitem[Nakk5]{nlfc}Nakkajima, Y.
\newblock{\em Artin-Mazur heights and Yobuko heights of 
proper log smooth schemes of Cartier type, and Hodge-Witt decompositions 
and Chow groups of quasi-$F$-split threefolds}. 
\newblock Preprint: available from http://arxiv.org/abs/1902.00185. 

\bibitem[Nakk6]{mwf}Nakkajima, Y.
\newblock {\em Monodromies and weight filtrations, 
and simple normal crossing log surfaces with torsion canonical sheaves}. 
\newblock Preprint.



\bibitem[Nakk7]{zflv}Nakkajima, Y.
\newblock{\em Zeta functions of log varieties}. 
\newblock Preprint. 






\bibitem[NS1]{nh2}Nakkajima, Y., Shiho, A.
\newblock{\em Weight filtrations on log crystalline cohomologies of 
families of open smooth varieties}.
\newblock Lecture Notes in Math.~1959, 
Springer-Verlag (2008).

\bibitem[NS2]{nhw}Nakkajima, Y., Shiho, A.
\newblock{\em Weight-filtered convergent complex}. 
\newblock  Preprint. 


\bibitem[NY]{ny}Nakkajima, Y.,  Yobuko, F. 
\newblock{\em Degenerations of log Hodge de Rham spectral sequences, 
log Kodaira vanishing theorem in characteristic $p>0$ 
and log weak Lefschetz conjecture for log crystalline cohomologies}. 
\newblock Preprint: available from https://arxiv.org/pdf/1902.09110.pdf. 

\bibitem[Oc]{ocltm}Ochiai, T.
\newblock{\em l-independence of traces of monodromy}. 
Math.~Ann.~315 (1999), 321--340.

\bibitem[Og1]{of}Ogus, A.
\newblock{\em $F$-isocrystals and de Rham cohomology.~${\rm II}$. 
Convergent isocrystals}.
\newblock Duke Math.~J.~51 (1984), 765--850.
   
\bibitem[Og2]{oc}Ogus, A.
\newblock{\em The convergent topos in characteristic $p$}.
\newblock The Grothendieck Festschrift Vol.~III, 
Progr.~Math.~88, Birkh\"{a}user (1990), 133--162.

\bibitem[Og3]{ofg}Ogus, A.
\newblock{\em F-crystals, Griffiths transversality, 
and the Hodge decomposition}. 
\newblock  Ast\'{e}risque 221, Soc.~Math.~de France (1994).
   
\bibitem[Og4]{ollc}Ogus, A.
\newblock{\em $F$-crystals on schemes with constant log structure}.
\newblock Compositio Math.~97 (1995), 
187--225.

\bibitem[Og5]{ob}Ogus, A.
\newblock{\em Lectures on logarithmic algebraic geometry}. 
\newblock Cambridge University Press (2018). 


\bibitem[Ol]{ols}Olsson, M.~C.  
\newblock{\em Crystalline Cohomology of Algebraic Stacks and Hyodo-Kato Cohomology}.
\newblock Ast\'{e}risque 316, Soc.~Math.~de France (2007).


\bibitem[P]{pvd}Proietto, V.-D.
\newblock{\em Kernel of the monodromy operator for semistable curves}.
\newblock  Algebraic Number Theory 
and Related Topics 2011, RIMS, Kokyuroku Bessatsu 44 (2013), 51--57. 


\bibitem[R]{raynd}Raynaud, M. 
\newblock{\em Sch\'{e}mas en groupes de type $(p,...,p)$.}
\newblock Bulletin de la Soc.~Math.~France~102 (1974), 241--280. 

\bibitem[RX]{rx}Raskind, W., Xarles, X.   
\newblock{\em On the \'{e}tale cohomology of algebraic varieties 
with totally degenerate reduction over $p$-adic fields}. 
J.~Math.~Sci.~Univ.~Tokyo 14 (2007), 261--291. 


\bibitem[RZ]{rz}Rapoport, M., Zink,T.                                            
\newblock{\em \"{U}ber die 
lokale Zetafunction von Shimuravariet\"{a}ten, 
Monodromiefiltration und verschwindende 
Zyklen in ungleicher Characteristik}.
\newblock Invent.~Math.~68 (1982), 21--101.


\bibitem[SaM1]{sam}Saito, M.
\newblock{\em Modules de Hodge polarisables}.
\newblock Publ.~Res.~Inst.~Math.~Sci.~24 (1988), 849-995.

\bibitem[SaM2]{samp}Saito, M.
\newblock{\em Monodromy filtration and positivity}.
\newblock Preprint, math.AG/0006162v10, 2003. 

\bibitem[SaT1]{stmfph}Saito, T.
\newblock{\em Modular forms and p-adic Hodge theory}.
\newblock Invent.~Math.~129 (1997), 607--620


\bibitem[SaT2]{stwsl}Saito, T.
\newblock{\em 
Weight spectral sequences and independence of $l$}.
\newblock J.~of the Inst.~of Math.~Jussieu~2 (2003), 1--52.

\bibitem[Sat]{stk}Sato, K.
\newblock{\em Abel-Jacobi mappings and finiteness of 
motivic cohomology groups}. 
\newblock Duke Math.~J.~104 (2000),  75--112. 

\bibitem[Sc]{sc}Scholze, P.
\newblock{\em Perfectoid spaces}. 
\newblock Publ.~Math.~IH\'{E}S  116 (2012), 245--313. 

\bibitem[Se]{sec}Serre, J.-P. 
\newblock{\em Corps locaux}.
\newblock  Hermann, Paris, (1968). 


\bibitem[ST]{seta}Serre, J.-P. 
\newblock{\em Good reduction of abelian varieties}.  
\newblock Ann.~of Math.~88 (1968), 492--517 


\bibitem[SGA 1]{sga1}Grothendieck, A., Raynaud, M.                                            
\newblock{\em Rev\^{e}tements 
\'{e}tales et groupe fondamental}.
\newblock 
Lecture Notes in Math.~224, Springer-Verlag (1971).



\bibitem[SGA 4-2]{sga4-2}Grothendieck, A.~et al.
\newblock{\em Th\'{e}orie des topos et cohomologie 
\'{e}tale des sch\'{e}mas}. 
\newblock Lecture Notes in Math.~270, 
Springer-Verlag (1972). 

\bibitem[SGA 5]{sga5-2}Grothendieck, A.~et al.
\newblock{\em Cohomologie $l$-adique et fonctions $L$}. 
\newblock Lecture Notes in Math.~589, 
Springer-Verlag (1977). 

\bibitem[SGA 7-I]{sga7}Grothendieck, A.
\newblock{\em Groupes de monodromie en 
g\'{e}om\'{e}trie alg\'{e}brique}.
\newblock Lecture Notes in Math.~288 (1972), 
Springer-Verlag.



\bibitem[Sh1]{s1}Shiho, A.
\newblock{\em Crystalline 
fundamental groups.~${\rm I}$-Isocrystals on log 
crystalline site and log convergent site}.
\newblock  J.~Math.~Sci.~Univ.~Tokyo 7 (2000),  
509--656. 


\bibitem[Sh2]{s2}Shiho, A.
\newblock{\em Crystalline 
fundamental groups ${\rm II}$-Log convergent 
cohomology and rigid cohomology}.
\newblock J.~Math.~Sci.~Univ.~Tokyo 9 (2002), 1--163.

\bibitem[Sh3]{s3}Shiho, A.
\newblock{\em Relative log convergent 
cohomology and relative rigid cohomology ${\rm I}$}.
\newblock Preprint, available from http://arxiv.org/abs/0707.1742 	
 
 
 
\bibitem[St1]{st1}Steenbrink, J.~H.~M.
\newblock{\em Limits of Hodge structures}. 
\newblock Invent.~Math.~31 (1976),  229--257. 


\bibitem[St2]{st2}Steenbrink, J.~H.~M.
\newblock{\em Logarithmic embeddings of varieties with 
normal crossings and mixed Hodge 
structures}.
\newblock Math.~Ann.~301 (1995), 105--118.


\bibitem[Te]{ter}Terasoma, T. 
\newblock{\em Monodormy weight filtration is independent of $l$}. 
\newblock Preprint: available from https://arxiv.org/abs/math/9802051.

\bibitem[Tsuj1]{tsgep}Tsuji, T.
\newblock{\em  $p$-adic Hodge theory in the semi-stable reduction case}.
\newblock In: Proceedings of the International Congress of Mathematicians, Vol.~II (Berlin, 1998),~Doc.~Math.~1998, Extra Vol.~II, 207--216  


\bibitem[Tsuj2]{tsp}Tsuji, T.
\newblock{\em Poincar\'{e} duality for logarithmic crystalline cohomology}.
\newblock Comp.~Math.~118 (1999), 11--41.


\bibitem[Tsuj3]{tst}Tsuji, T.
\newblock{\em $p$-adic \'{e}tale cohomology and crystalline cohomology 
in the semi-stable reduction case}.
\newblock Invent.~Math.~137 (1999), 233--411. 

\bibitem[Tsuj4]{tsst}Tsuji, T.
\newblock{\em Semi-stable conjecture of Fontaine-Jannsen: a survey}. 
\newblock In: Cohomologies $p$-adiques et applications arithm\'{e}tiques, II. 
Ast\'{e}risque 
279 (2002), 323--370.

\bibitem[Tsuj5]{tsd}Tsuji, T.
\newblock{\em On nearby cycles and  
D-modules of log schemes in charcteristic $p > 0$}. 
\newblock Comp.~Math.~146 (2010), 1552--1616. 

\bibitem[Tsuj6]{tsa}Tsuji, T.
\newblock{\em Saturated morphisms of logarithmic schemes}.  
\newblock Tunisian J.~Math.~1 (2019), 185--220.

\bibitem[Tsuz1]{tzp}Tsuzuki, N.
\newblock{\em Cohomological descent of rigid cohomology 
for proper coverings}.
\newblock  Invent.~Math.~151 (2003), 101--133.


\bibitem[Tsuz2]{tzcp}Tsuzuki, N.
\newblock{\em Cohomological descent in rigid cohomology}.
\newblock In: Geometric aspects of Dwork theory Vol.~II 
(Adolphson, A., Baldassari, F., Berthelot, P., Katz, N., 
Loeser, F. eds.), 
Walter de Gruyter (2004), 931--981. 


\bibitem[TY]{tay}Taylor, R and Yoshida, T.
\newblock{\em  Compatibility of 
local and global Langlands correspondences}.
\newblock In: J.~Amer.~Math.~Soc., 20-2 (2007), 467-493. 

\bibitem[V]{vi}Vidal, I.
\newblock{\em Monodromie locale et fonctions Z\^eta des log sch\'{e}mas}. 
\newblock  Geometric aspects of Dwork theory I Vol.~2, de Gruyter, Berlin (2004),  
983--1038. 


\bibitem[W]{weib}Weibel, C.~A.
\newblock{\em An introduction to homological algebra}.
\newblock  Cambridge Studies in Advanced Mathematics, 38.
Cambridge University Press, Cambridge (1994).

\bibitem[Yd]{ytgll}Yoshida, T.
\newblock{\em 
Global and local Langlands correspondences for $GL_n$}.
\newblock  (in Japanese), 
Proceedings of the 50th Symposium on Algebra, Tokushima Univ. 2005, 
available from Yoshida's homepage:  
http://www.dpmms.cam.ac.uk/~ty245/.

\bibitem[Yk]{yoshi}Yoshioka, H.
\newblock{\em 
Semistable reduction theorem for logarithmically smooth varieties}.
\newblock Masters thesis at Tokyo University (1998).

\bibitem[Z]{zvmh}Zucker, S.
\newblock{\em Variation of mixed Hodge structure.~{\rm II}}.
\newblock  Invent.~ Math.~80 (1985), 543--565.

\end{thebibliography}
\end{document}